\documentclass[a4paper,12pt,frenchb]{article}

\usepackage{amsmath,amsbsy,amsfonts,amssymb,amsthm}
\usepackage[french]{babel}
\newtheorem*{thm}{Th\'eor\`eme}
\newtheorem*{lem}{Lemme}
\newtheorem*{prop}{Proposition}

\newtheorem*{cor}{Corollaire}

\begin{document}

 \title{  Sur la stabilit\'e des int\'egrales orbitales nilpotentes }
\author{J.-L. Waldspurger}
\date{10 octobre 2024}
 \maketitle

   {\bf R\'esum\'e} Soit $G$ un groupe r\'eductif connexe d\'efini sur un corps local non archim\'edien $F$ de caract\'eristique nulle. On suppose que $G$ est quasi-d\'eploy\'e, adjoint et absolument simple. Notons $\mathfrak{g}$ l'alg\`ebre de Lie de $G$ et consid\'erons l'espace des distributions invariantes sur $\mathfrak{g}(F)$, qui sont stables et dont le support est form\'e d'\'el\'ements nilpotents. Magdy Assem a formul\'e des conjectures d\'ecrivant cet espace. Nous prouvons une partie de ces conjectures, sous l'hypoth\`ese que la caract\'eristique r\'esiduelle de $F$ est "tr\`es grande" relativement \`a $G$. 
   
   \bigskip
   
   {\bf English title: }Stable distributions and nilpotent orbital integrals
   
   \bigskip

  {\bf Abstract} Let $G$ be a connected reductive group defined over a non-archimedean local field $F$ of characteristic $0$. We assume $G$ is quasi-split, adjoint and absolutely simple. Let $\mathfrak{g}$ be the Lie algebra of $G$. We consider the space of the invariant distributions on $\mathfrak{g}(F)$, which are stable and supported by the set of nilpotent elements of $\mathfrak{g}(F)$. Magdy Assem has stated several conjectures which describe this space. We prove some of these conjectures, assuming that the residual characteristic of $F$ is "very large" relatively to $G$.

 \newpage

 {\bf Table des matières}
 
 \bigskip
 
 1. Introduction
 
 \noindent 1.1 Enoncé du théorème p.7
 
  \noindent 1.2 L'espace ${\cal D}^{st}(\mathfrak{g}(F))$ p.9
 
  \noindent 1.3 Paramétrages p.12
 
 \bigskip
 2. L'espace 
 ${\cal D}^{st}(\mathfrak{g}(F))$
  \bigskip
 
 2.1 Immeubles
 
  \noindent  2.1.1 Notations p.15
 
 \noindent   2.1.2 Les hypothèses sur $p$ p.17
 
  \noindent  2.1.3 Forme intérieure quasi-déployée p.17
 
  \noindent  2.1.4 Immeubles de Bruhat-Tits p.18
 
 \noindent   2.1.5 Immeubles et $F^{nr}$-Levi p.19
 
  \noindent  2.1.6 $F^{nr}$-Levi et groupes en réduction p.20
 
  \noindent  2.1.7 Un lemme sur les épinglages p.21
 
  \noindent  2.1.8 Description d'un appartement et d'une alc\^ove p.24
 
 \noindent   2.1.9 Racines affines  p.25
 
  \noindent  2.1.10 Un corollaire p.28
 
  \noindent  2.1.11 Discrétisation de l'immeuble; le groupe ${\bf G}$ p.29
 
  \noindent  2.1.12 Description des groupes $G_{{\cal F}}$ p.31
 
  \noindent  2.1.13 Adhérence d'une facette p.33
 
  \noindent  2.1.14 Action de $G(F^{nr})$ p.34
   \bigskip
 
 2.2 L'espace $FC(\mathfrak{g}(F))$
 
  \noindent  2.2.1 Fonctions sur les algèbres de Lie p.35
  
   \noindent  2.2.2 Faisceaux-caractères cuspidaux sur les algèbres de Lie p.25
   
    \noindent  2.2.3 Faisceaux-caractères  sur les algèbres de Lie p.36
    
     \noindent  2.2.4 Calcul de $\epsilon^{\flat}_{N}$ p.38
     
      \noindent  2.2.5 Fonctions caractéristiques p.40
      
       \noindent  2.2.6 Les espaces $I(\mathfrak{g}(F))$ et $FC(\mathfrak{g}(F))$ p.42
       
        \noindent  2.2.7 Les espaces $I^{st}_{cusp}(\mathfrak{g}(F))$ et $FC^{st}(\mathfrak{g}(F))$ p.43
        
         \noindent  2.2.8 Décomposition de l'espace $FC^{st}(\mathfrak{g}(F))$ p.43
 
  \noindent  2.2.9 Immeubles et restriction à la Weil p.44
  
   \noindent  2.2.10 Action d'automorphismes p.45
   
    \noindent  2.2.11 Démonstration du lemme 2.2.8 p.45
    
     \noindent  2.2.12 Extensions non ramifiées et espaces $FC^{st}(\mathfrak{g}(F))$ p.47
      \bigskip
     
     2.3 Etude des couples $(M,s)$
  
  \noindent    2.3.1 Réduction au cas absolument simple p.47
  
   \noindent 2.3.2 Construction de $F^{nr}$-Levi p.48
   
    \noindent 2.3.3 Couples formés d'un Levi et d'un sommet de son immeuble p.49
    
     \noindent 2.3.4 Normalisateur p.50
     
      \noindent 2.3.5 Excellents couples $(M,s)$ p.52
      
       \noindent   2.3.6 Excellents ensembles $\Lambda$ p.53
       
        \noindent  2.3.7 Un corollaire p.55
        
         \noindent 2.3.8 Cas d'un excellent couple $(M,s)$ p.55
         
          \noindent 2.3.9 Levi standard et excellents couples p.56
          
           \noindent  2.3.10 Description de $K_{s}^{M,nr,0}\backslash K_{s}^{M,nr,\dag}$ p.57
           
            \noindent  2.3.11 Couples $(M,s_{M})$ avec $s_{M}\in S^{nr,st}(M)$ p.60
             \bigskip
            
            2.4 $F^{nr}$-Levi et conjugaison stable
            
       \noindent      2.4.1 Bons éléments p.62
       
        \noindent  2.4.2 Eléments de réduction régulière p.64
        
         \noindent  2.4.3 Existence d'éléments entiers et de réduction régulière p.64
         
          \noindent  2.4.4 Description de l'ensemble $Imm^G(H_{ad})$ p.66
          
           \noindent  2.4.5 Classes de conjugaison stable dans ${\cal L}^{nr}_{F}$ p.67
           
            \noindent  2.4.6 Passage au groupe $G^*$ pour l'ensemble ${\cal L}^{nr}_{F}$ p.68
            
             \noindent  2.4.7 Passage à $G^*$, le cas stable p.69
             
              \noindent  2.4.8 L'ensemble ${\cal L}^{*,st}_{F,sd}$ p.71
              
               \noindent  2.4.9 Les ensembles ${\cal L}_{F}^{nr,st}$ et ${\cal L}_{ell}^{nr,st}$ p.72
               
                \noindent  2.4.10 Sur l'action de certains automorphismes p.74
                
                 \noindent  2.4.11 Le groupe $H_{s_{H}}$ p.75 
                  \bigskip
                 
                 2.5 Stabilité et espace ${\cal D}(\mathfrak{g}(F))$ 
    
     \noindent 2.5.1 L'espace  ${\cal D}(\mathfrak{g}(F))$  p.76
     
     \noindent  2.5.2 Réduction aux éléments elliptiques p.77
     
      \noindent  2.5.3 Décomposition selon les Levi p.78
      
       \noindent 2.5.4 Une première description de l'espace ${\cal D}_{cusp}(\mathfrak{g}(F))$ p.78
       
        \noindent  2.5.5 Un calcul de mesures p.81
        
         \noindent 2.5.6 Description des distributions associées aux éléments de ${\cal K}^*_{cusp}(\mathfrak{g}(F))$ p.83
         
          \noindent 2.5.7 Expression intégrale de $D_{X_{H},\varphi}^G$ p.87
          
           \noindent  2.5.8 Une deuxième description de l'espace $D^G({\cal D}_{cusp}(\mathfrak{g}(F)))$ p.89
           
            \noindent  2.5.9 Un premier résultat d'instabilité p.89
            
             \noindent  2.5.10 Les espaces ${\cal K}_{cusp}^{st}(\mathfrak{g}(F))$ et ${\cal K}^{\flat}_{cusp}(\mathfrak{g}(F))$ p.91
             
              \noindent  2.5.11 Description de l'espace ${\cal D}^{st}_{cusp}(\mathfrak{g}(F))$ p.94
              
               \noindent  2.5.12 Description des distributions associées aux éléments de ${\cal D}(\mathfrak{g}(F))$ p.95
               
                \noindent  2.5.13 L'espace ${\cal K}^{st}(\mathfrak{g}(F))$ p.96
                
                 \noindent  2.5.14 L'espace ${\cal K}^{st}(\mathfrak{g}(F))$ et les Levi de $G$ p.97
                 
                  \noindent  2.5.15 Description de ${\cal D}^{st}(\mathfrak{g}(F))$ p.99
                  
                   \noindent  2.5.16 Dimension de l'espace ${\cal D}^{st}(\mathfrak{g}(F))$ p.99
                   
                  \bigskip
                   
                   3. Intégrales orbitales nilpotentes
                    \bigskip
                   
    3.1   Orbites nilpotentes
    
     \noindent 3.1.1 Rappels sur les représentations de $\mathfrak{sl}(2)$ p.100
     
      \noindent  3.1.2 $\mathfrak{sl}(2)$-triplets et orbites nilpotentes p.100
      
       \noindent  3.1.3 Orbites engendrées p.101
       
        \noindent  3.1.4 Les groupes $A(N)$ et $\bar{A}(N)$ p.102
        
         \noindent  3.1.5 Partitions p.102
         
          \noindent  3.1.6 Paramétrages des orbites nilpotentes p.103
          
           \noindent  3.1.7 Action galoisienne p.104
            \bigskip
           
           3.2 Orbites nilpotentes et groupe ${\bf G}$
           
            \noindent  3.2.1 Orbites nilpotentes dans $\boldsymbol{\mathfrak{g}}$ p.105
            
             \noindent  3.2.2 Les ensembles $A^+(N)$,  $\tilde{A}(N)$, $\bar{A}(N)$ et 
             $\tilde{\bar{A}}(N)$ p.106

\noindent 3.2.3 Les ensembles ${\cal C}^{\sharp}_{F}$ et $\bar{{\cal C}}^{\sharp}_{F}$ p.107

\noindent  3.2.4 Les applications $\iota_{{\cal F},nil}$ et $c_{{\cal F}}$ p.109

\noindent  3.2.5 Les définitions sur $F^{nr}$ p.110

\noindent  3.2.6 Description de $\iota_{{\cal F},nil}$ dans le cas ramifié p.110

\noindent  3.2.7 Preuve du lemme 3.2.3 p.112

\noindent  3.2.8 Sur le relèvement des commutants p.113

\noindent  3.2.9 Une construction auxiliaire p.116

\noindent  3.2.10 Description de $c_{s}$ dans le cas ramifié p.117

\noindent  3.2.11 Preuve de la proposition 3.2.4  p.120
 \bigskip 

3.3 Orbites nilpotentes de $\mathfrak{g}(F)$

\noindent  3.3.1 La construction de DeBacker  p.122

\noindent  3.3.2 Un diagramme d'applications p.126

\noindent  3.3.3 L'ensemble $\bar{{\cal C}}^{\natural}_{F}$ p.127

\noindent  3.3.4 Comparaison des commutants p.129

\noindent  3.3.5 Commutants et action galoisienne p.132

\noindent  3.3.6 L'isomorphisme ${\cal C}^{\natural}_{F}\to {\cal C}^{\sharp}_{F}$ p.133

\noindent  3.3.7 Interprétation cohomologique de l'application $\bar{c}_{{\cal F},F}$ p.133

\bigskip

4. Normalisation des actions de Frobenius
 \bigskip

4.1 Isomorphismes d'espaces $FC^{st}(\mathfrak{g}(F))$

\noindent 4.1.1 Position du problème p.135

\noindent 4.1.2 Familles admissibles p.135

\noindent 4.1.3 Résolution du problème p.137

\noindent 4.1.4 Un sous-ensemble de ${\cal O}_{{\cal E}}$ p.140

\noindent 4.1.5 Construction de la famille admissible $(\epsilon_{N})_{N\in {\cal O}}$ où ${\cal O}={\cal O}_{{\cal E}}$ p.142

\noindent  4.1.6 Equivariance par automorphisme p.144

\noindent 4.1.7 Prolongement du caractère $\epsilon_{N}$ p.144
             
 \noindent 4.1.8 Composition d'isomorphismes $\iota_{\phi}$ p.145
  \bigskip
 
 4.2 Isomorphismes d'espaces $FC^{st}(\mathfrak{g}(F))$ et transfert endoscopique
 
 \noindent 4.2.1 Enoncé du théorème p.146
 
 \noindent   4.2.2 Préparatifs à la preuve du théorème 4.2.1 p.146
 
 \noindent 4.2.3 Début de la preuve p.148
 
 \noindent 4.2.4 Preuve dans le cas  $B_{n}$ p.150
 
 \noindent 4.2.5 Preuve dans le cas     $C_{n}$ p.153
 
 \noindent 4.2.6 Preuve dans le cas $(D_{n},nr)$ p.155
 
 \noindent 4.2.7 Preuve dans le cas $(D_{n},ram)$ p.157
 
 \noindent 4.2.8 Preuve dans le cas où $G$ est de type $E_{6}$ déployé ou $E_{7}$ p.158
 
 \noindent 4.2.9 Preuve dans le cas où $G$ est de type $(A_{n-1},ram)$ avec $n$ pair p.160
  \bigskip
 
 4.3 Utilisation de la construction de Lusztig
 
 \noindent 4.3.1 Définition d'un caractère $\epsilon_{N}^{\flat}$ p.162
 
 \noindent 4.3.2 Cas stable p.166
 
 \noindent 4.3.3 Les cas faciles p.168
 
 \noindent 4.3.4 Preuve du théorème 4.3.2 dans le cas $(A_{n-1},ram)$ avec $n$ impair    p.168
                                      
   \noindent 4.3.5 Preuve du théorème 4.3.2 dans le cas $(A_{n-1},ram)$ avec $n$ pair    p.170
    
     \noindent 4.3.6 Preuve du théorème 4.3.2 dans le cas $B_{n}$ p.172
    
   \noindent 4.3.7 Début de la preuve du théorème 4.3.2 pour $G$ de type $E_{8}$ p.174
   
     \noindent  4.3.8 Sommet $s_{\alpha_{6}}$ pour $G$ de type $E_{8}$   p. 175

       \noindent  4.3.9 Sommet $s_{\alpha_{3}}$ pour $G$ de type $E_{8}$   p.177
    
     \noindent 4.3.10 Début de la preuve du théorème 4.3.2  dans le cas $(E_{6},ram)$ p.179
     
       \noindent  4.3.11 Sommet $s_{\beta_{1,6}}$ pour $G$ de type $(E_{6},ram)$ p.179
       
          \noindent  4.3.12 Sommet $s_{\beta_{4}}$ pour $G$ de type $(E_{6},ram)$ p.181
           \bigskip
          
          4.4 Isomorphismes d'espaces $FC^{st}(\mathfrak{h}_{SC}(F))$
          
             \noindent   4.4.1 Définition des isomorphismes p.183
             
                \noindent   4.4.2 Composition d'isomorphismes $\iota_{H',H}$ p.184

  \noindent  4.4.3 Comparaison avec le transfert endoscopique p.185
  
    \noindent  4.4.4 Comparaison de deux éléments de $FC^{st}(\mathfrak{h}_{SC}(F))$ p.186
    
    \bigskip
    
    5. Relation avec les faisceaux-caractères unipotents sur un groupe fini
     \bigskip
    
    5.1 Quelques constructions
    
   \noindent   5.1.1 Le groupe $G$ p.186
   
    \noindent   5.1.2 Conjugaison des éléments $t_{\alpha}$ p.187

     \noindent  5.1.3 Construction de Levi p.189
      \bigskip
     
     5.2 Faisceaux-caractères unipotents
     
      \noindent   5.2.1 Rappels sur la représentation de Springer généralisée p.189
      
       \noindent   5.2.2 Faisceaux-caractères unipotents cuspidaux sur $G$ p.191
       
        \noindent   5.2.3 Faisceaux-caractères unipotents sur $G$   p.192
        
         \noindent   5.2.4 Le caractère $\mu_{t}$ p.193
         
          \noindent   5.2.5 L'ensemble $\Lambda$ associé à un couple $(M,{\cal E})$ p.194
          
           \noindent   5.2.6 Restriction au commutant d'un élément semi-simple p.195
           
            \noindent   5.2.7 Un lemme de restriction p.195
            
             \noindent   5.2.8 Preuve de la proposition 5.2.6 p.198
              \bigskip
             
             5.3 Paramétrages
             
              \noindent   5.3.1 L'orbite spéciale associée à un triplet $(M,{\cal E},\rho)$ p.199
              
               \noindent   5.3.2 Calcul de l'orbite spéciale associée à un triplet $(M,{\cal E},\rho)$ p.200
                              
         \noindent  5.3.3 Une propriété des familles admissibles canoniques p.201
         
         \noindent  5.3.4 Majoration d'orbites p.203
         
         \noindent  5.3.5 Paramétrage de ${\cal A}(G,{\cal O})$ p.203
         
         \noindent  5.3.6 Relèvement d'éléments de $\bar{A}(N)$ p.205
         
         \noindent  5.3.7 Séparation des éléments de ${\cal A}(G)$ p.206
         
         \bigskip
         
         6. L'espace des combinaisons linéaires stables d'intégrales orbitales nilpotentes     
          \bigskip
    
    6.1 Intégrales orbitales nilpotentes et espace ${\cal D}(\mathfrak{g}(F))$
    
    \noindent 6.1.1 Normalisation des intégrales orbitales nilpotentes p.210

     \noindent 6.1.2 La fonction-test associée à un couple $(s,{\cal O})$ p.210
     
      \noindent  6.1.3 Calcul de la fonction $\underline{h}_{s,{\cal O}}$ p.211
      
       \noindent 6.1.4 Description de $SI(\mathfrak{g}(F))_{nil}^*$ p.213
       
        \noindent  6.1.5 L'ensemble ${\cal B}(G)$ p.213
        
         \noindent  6.1.6 Deux bases de ${\cal D}^{st}(\mathfrak{g}(F))$ p.214
         
          \noindent  6.1.7 Evaluation de la distribution $\hat{k}(\Lambda,\varphi,w)$ p.215
          
           \noindent 6.1.8 Evaluation de la distribution $\hat{k}(\Lambda,\varphi,\rho)$ p.221
           
            \noindent  6.1.9 Calcul de $\hat{k}(\Lambda,\varphi,\rho)(h_{s,{\cal O}})$ p.222
             \bigskip
            
            6.2 Deux théorèmes de maximalité
            
             \noindent  6.2.1 Paramétrages de l'ensemble ${\cal B}_{Irr}(G)$ p.224
             
              \noindent  6.2.2 Le cas de type $(A_{n-1},ram)$ p.226
              
               \noindent  6.2.3 Le cas de type $(D_{n},ram)$ p.226
               
                \noindent  6.2.4 Le cas $(D_{4},3-ram)$ p.228
                
                 \noindent  6.2.5 Le cas $(E_{6},ram)$ p.228
                 
                  \noindent  6.2.6 Un premier théorème p.230
                  
                   \noindent  6.2.7 Un deuxième théorème p.231

       \noindent  6.2.8 Preuve des théorèmes 6.2.6 et 6.2.7 dans le cas où $G$ est déployé sur $F^{nr}$ p.231
       
        \noindent  6.2.9 Rappels sur la correspondance de Springer généralisée pour les groupes classiques p.232
        
         \noindent  6.2.10        Preuve des théorèmes 6.2.6 et 6.2.7 dans le cas $(D_{n},ram)$ p.235
         
         \noindent 6.2.11 Preuve des théorèmes 6.2.6 et 6.2.7 dans le cas $(A_{n-1},ram)$   p.245  
         
         \noindent 6.2.12 Preuve des théorèmes 6.2.6 et 6.2.7 dans le cas $(D_{4},3-ram)$ p.253
         
         \noindent 6.2.13 Preuve des théorèmes 6.2.6 et 6.2.7 dans le cas $(E_{6},ram)$ p.255
          \bigskip
         
         6.3 Une nouvelle base de $SI(\mathfrak{g}(F))_{nil}^*$ 
         
         \noindent 6.3.1 Un calcul plus précis de $\hat{k}(\Lambda,\varphi,\rho)(h_{s,{\cal O}_{s}})$ p.260
         
         \noindent 6.3.2 Construction d'une famille de distributions p.261
         
         \noindent 6.3.3 Indépendance linéaire des éléments de la famille p.261
         
         \noindent 6.3.4 Nombres d'éléments p.262
         
         \noindent 6.3.5 L'énoncé final p.264
         
         \bigskip
         
         Références p.266

\newpage

\section{Introduction}

\subsection{Enonc\'e du th\'eor\`eme}

Soient $F$ un corps local non archim\'edien de caract\'eristique nulle et $G$ un groupe r\'eductif connexe d\'efini sur $F$. On note ${\mathbb F}_{q}$ le corps r\'esiduel de $F$, $q$ \'etant son nombre d'\'el\'ements,  $p$ la  caract\'eristique de ${\mathbb F}_{q}$ et $val_{F}$ la valuation usuelle de $F$. 
Notons $\mathfrak{g}$ l'alg\`ebre de Lie de $G$. On identifie $G$ et $\mathfrak{g}$ aux ensembles de leurs points dans une cl\^oture alg\'ebrique $\bar{F}$ de $F$.

 On note $C_{c}^{\infty}(\mathfrak{g}(F))$ l'espace des fonctions sur $\mathfrak{g}(F)$ \`a valeurs complexes, localement constantes et \`a support compact. On note
$I(\mathfrak{g}(F))$ le quotient de $C_{c}^{\infty}(\mathfrak{g}(F))$ par le sous-espace des fonctions dont toutes les int\'egrales orbitales sont nulles. On note $I(\mathfrak{g}(F))^*$ l'espace dual de $I(\mathfrak{g}(F))$. C'est l'espace des distributions invariantes  sur $\mathfrak{g}(F)$ (invariantes par conjugaison par $G(F)$, o\`u on appelle conjugaison l'action adjointe). On d\'efinit le sous-espace $SI(\mathfrak{g}(F))^*$ des distributions stables, cf. \ref{lespacecalD}. 

Notons $\mathfrak{g}_{nil}$ l'ensemble des \'el\'ements nilpotents de $\mathfrak{g}$. On note $I(\mathfrak{g}(F))^*_{nil}\subset I(\mathfrak{g}(F))^*$ le sous-espace des distributions invariantes \`a support dans $\mathfrak{g}_{nil}(F)$. Cet espace se d\'ecrit ais\'ement. Notons $\mathfrak{g}_{nil}(F)/conj$ l'ensemble des orbites nilpotentes de $\mathfrak{g}(F)$, c'est-\`a-dire des classes de conjugaison par $G(F)$ dans $\mathfrak{g}_{nil}(F)$. A toute orbite ${\mathbb  O}\in \mathfrak{g}_{nil}(F)/conj$, on associe l'int\'egrale orbitale $I_{{\mathbb O}}$ sur ${\mathbb  O}$, cf. \ref{integralesnilpotentes}. C'est un \'el\'ement de $I(\mathfrak{g}(F))^*_{nil}$ et la famille $(I_{{\mathbb  O}})_{{\mathbb O}\in \mathfrak{g}_{nil}(F)/conj}$ est une base de cet espace. Pour tout sous-espace, disons $I^*_{\square}$ de $I(\mathfrak{g}(F))^*$, on pose $SI^*_{\square}=SI(\mathfrak{g}(F))^*\cap I^*_{\square}$. En particulier, on
 pose $SI(\mathfrak{g}(F))^*_{nil}=SI(\mathfrak{g}(F))^*\cap I(\mathfrak{g}(F))^*_{nil}$. Cet espace est assez myst\'erieux et on n'en conna\^{\i}t pas actuellement de description g\'en\'erale. Le but de l'article est d'en donner une description partielle sous des hypoth\`eses restrictives.

Magdy Assem a formul\'e des conjectures \`a ce propos (\cite{MA}, conjectures A et C de l'introduction). En fait, il consid\'erait des distributions invariantes sur $G(F)$ \`a support unipotent mais la traduction est imm\'ediate. Rappelons ces conjectures, du moins partiellement. Notons $\mathfrak{g}_{nil}/conj$ l'ensemble des orbites nilpotentes alg\'ebriques de $\mathfrak{g}$, c'est-\`a-dire l'ensemble des classes de conjugaison par $G=G(\bar{F})$ dans $\mathfrak{g}_{nil}=\mathfrak{g}_{nil}(\bar{F})$. Notons $\Gamma_{F}$ le groupe de Galois de $\bar{F}/F$. Ce groupe agit sur $\mathfrak{g}_{nil}/conj$ et on  note $(\mathfrak{g}_{nil}/conj)^{\Gamma_{F}}$ le sous-ensemble des points fixes. 

{\bf Remarque.} Nous noterons ${\mathbb O}$  une orbite dans $\mathfrak{g}_{nil}(F)$ et $\bar{{\mathbb O}}$ une orbite dans $\mathfrak{g}_{nil}$. Dans cette  dernière notation, le  $\,\,\bar{}\,\,$  ne doit pas \^etre confondu avec une quelconque adhérence.
\bigskip

Tout orbite ${\mathbb O}\in \mathfrak{g}_{nil}(F)/conj$ est contenue dans une unique orbite alg\'ebrique $\bar{{\mathbb  O}}\in (\mathfrak{g}_{nil}/conj)^{\Gamma_{F}}$. Pour $\bar{{\mathbb  O}}\in (\mathfrak{g}_{nil}/conj)^{\Gamma_{F}}$, notons $I(\mathfrak{g}(F))^*_{\bar{{\mathbb  O}}}$ l'espace de distributions engendr\'e par les int\'egrales orbitales $I_{{\mathbb  O}}$ pour ${\mathbb O}\in \mathfrak{g}_{nil}(F)/conj$, ${\mathbb  O}\subset \bar{{\mathbb O}}$.   La premi\`ere conjecture est

(1) $SI(\mathfrak{g}(F))^*_{nil}$ est la somme directe des $SI(\mathfrak{g}(F))^*_{\bar{{\mathbb  O}}}$ quand $\bar{{\mathbb  O}}$ d\'ecrit $(\mathfrak{g}_{nil}/conj)^{\Gamma_{F}}$.

\noindent Bien que tr\`es naturelle, cette assertion n'est absolument pas \'evidente. 

Lusztig a d\'efini le sous-ensemble des orbites nilpotentes (alg\'ebriques) sp\'eciales. Bien que leur d\'efinition soit tr\`es peu parlante pour moi, ces orbites jouent un r\^ole particulier et essentiel dans de nombreuses situations. Elles interviennent ici par la deuxi\`eme conjecture

(2) soit $\bar{{\mathbb  O}}\in (\mathfrak{g}_{nil}/conj)^{\Gamma_{F}}$; si $\bar{{\mathbb  O}}$ n'est pas sp\'eciale, $SI(\mathfrak{g}(F))^*_{\bar{{\mathbb  O}}}=\{0\}$. 

La troisi\`eme conjecture est plus subtile. Magdy Assem la pose en supposant que $G$ est quasi-d\'eploy\'e et classique. Supposons plut\^ot, puisque ce seront nos hypoth\`eses, que 
  $G$ est adjoint, quasi-d\'eploy\'e  et absolument simple. Soit $\bar{{\mathbb  O}}\in (\mathfrak{g}_{nil}/conj)^{\Gamma_{F}}$, supposons que  cette orbite soit sp\'eciale. Parce que $G$ est quasi-d\'eploy\'e, l'ensemble $\bar{{\mathbb  O}}\cap \mathfrak{g}(F)$ est non vide, cf. \cite{K1} th\'eor\`eme 4.2. Fixons un \'el\'ement $e$ de cet ensemble. Notons $Z_{G}(e)$ le commutant de $e$ dans $G$, $Z_{G}(e)^0$ sa composante neutre et posons $A(e)=Z_{G}(e)/Z_{G}(e)^0$. Lusztig a d\'efini un quotient $\bar{A}(e)$ de ce groupe, que l'on appelle comme il se doit le quotient de Lusztig. L'action galoisienne sur $A(e)$ se quotiente en une action sur $\bar{A}(e)$ et il y a une application naturelle $d_{e}:H^1(\Gamma_{F},Z_{G}(e))\to H^1(\Gamma_{F},\bar{A}(e))$. 
 L'ensemble des orbites ${\mathbb  O}\in \mathfrak{g}_{nil}(F)/conj$ telles que ${\mathbb  O}\subset \bar{{\mathbb O}}$ s'identifie au noyau $H^1(\Gamma_{F},Z_{G}(e))_{1}$ de l'application naturelle entre ensembles point\'es $H^1(\Gamma_{F},Z_{G}(e))\to H^1(\Gamma_{F},G)$.  Pour un \'el\'ement $h\in H^1(\Gamma_{F},\bar{A}(e))$, notons $E(\bar{{\mathbb  O}},h)$ l'ensemble des  orbites ${\mathbb O}\in \mathfrak{g}_{nil}(F)/conj$ telles que ${\mathbb  O}\subset \bar{{\mathbb  O}}$ et dont le param\`etre dans $H^1(\Gamma_{F},Z_{G}(e))_{1}$ appartient \`a $d_{e}^{-1}(h)$.
 Remarquons que cet ensemble  peut a priori \^etre vide puisqu'il n'est pas clair que l'application  $d_{e}$ soit surjective, ni a fortiori sa restriction \`a $H^1(\Gamma_{F},Z_{G}(e))_{1}$. Notons  
$I(\mathfrak{g}(F))^*_{\bar{{\mathbb  O}},h}$ l'espace de distributions engendr\'e par les int\'egrales orbitales $I_{{\mathbb O}}$ pour   ${\mathbb O}\in E(\bar{{\mathbb O}},h)$. L'assertion que Magdy Assem pose comme conjecture dans le cas o\`u $G$ est classique est

(3) pour $h\in H^1(\Gamma_{F},\bar{A}(e))$ tel que $E(\bar{{\mathbb O}},h)\not=\emptyset$, l'espace $SI(\mathfrak{g}(F))^*_{\bar{{\mathbb  O}},h}$ est une droite; l'espace $SI(\mathfrak{g}(F))^*_{\bar{{\mathbb  O}}}$ est la somme directe de ces droites. 

Enon\c{c}ons maintenant notre th\'eor\`eme. On suppose que $G$ est adjoint, quasi-d\'eploy\'e  et absolument simple.  Supposons que $p$ est "tr\`es grand" relativement \`a $G$, pr\'ecis\'ement que l'hypoth\`ese $(Hyp)_{2}(p)$ de \ref{leshypothesessurp} est v\'erifi\'ee. Cela signifie que $p\geq c_{1}val_{F}(p)+c_{2}$, o\`u $c_{1}$ et $c_{2}$ sont des entiers d\'ependant de $G$. On doit d\'efinir les "orbites exceptionnelles". Il n'y en a pas si $G$ n'est pas de type $E_{7}$ ou $E_{8}$. Si $G$ est de type $E_{7}$ ou $E_{8}$, on utilise la classification des \'el\'ements de $\mathfrak{g}_{nil}/conj$ de \cite{C} pages 403 \`a 407. On dit qu'une orbite $\bar{{\mathbb  O}}\in \mathfrak{g}_{nil}/conj$ est exceptionnelle dans les cas suivants:

$G$ est de type $E_{7}$ et $\bar{{\mathbb O}}$ est de type $A_{4}+A_{1}$;

$G$ est de type $E_{8}$ et $\bar{{\mathbb O}}$ est de type $A_{4}+A_{1}$ ou $E_{6}(a_{1})+A_{1}$.

\begin{thm}{(i) Les assertions (1) et (2) ci-dessus sont v\'erifi\'ees. 

(ii) Soit $\bar{{\mathbb  O}}\in (\mathfrak{g}_{nil}/conj)^{\Gamma_{F}}$, supposons que $\bar{{\mathbb  O}}$ soit sp\'eciale et fixons $e\in \mathfrak{g}(F)\cap \bar{{\mathbb O}}$. L'application $d_{e}$ se restreint en une surjection $H^1(\Gamma_{F},Z_{G}(e))_{1}\to H^1(\Gamma_{F},\bar{A}(e))$. Supposons que $\bar{{\mathbb  O}}$ ne soit pas exceptionnelle ou que $4$ divise $q-1$. Alors la dimension de $SI(\mathfrak{g}(F))^*_{\bar{{\mathbb  O}}}$ est \'egale au nombre d'\'el\'ements de $H^1(\Gamma_{F},\bar{A}(e))$. Supposons que $\bar{{\mathbb O}}$ soit exceptionnelle et que $4$ ne divise pas $q-1$. Alors la dimension de $SI(\mathfrak{g}(F))^*_{\bar{{\mathbb  O}}}$ est  $2$ tandis que le nombre d'\'el\'ements de $H^1(\Gamma_{F},\bar{A}(e))$ est $4$.} \end{thm}

{\bf Remarques.} (4) Ces assertions sont d\'ej\`a connues dans le cas des groupes classiques non ramifi\'es, c'est-\`a-dire d\'eploy\'es sur une extension non ramifi\'ee, cf. \cite{W1} IX.16, ainsi que dans le cas o\`u $G$ est de type $G_{2}$, cf. \cite{DK}.

(5) Dans le cas o\`u $\bar{{\mathbb O}}$ est exceptionnelle et o\`u $4$ ne divise pas $q-1$, l'assertion (ii) interdit \`a (3) d'\^etre v\'erifi\'ee. Dans les autres cas, l'assertion (ii) ne d\'emontre nullement (3) mais les deux assertions ont  \'evidemment un air de famille.

(6) L'hypoth\`ese que $p$ est "tr\`es grand" n'est utilis\'ee de fa\c{c}on essentielle qu'en un seul point, cf. \ref{reductionelliptique}. Elle assure une propri\'et\'e que l'on a d\'emontr\'ee dans \cite{W5} dont la preuve utilise l'exponentielle. Il est probable qu'en utilisant les diff\'erents substituts de cette application, par exemple les quasi-logarithmes de Kazhdan et Varshavsky, on peut r\'eduire la condition sur $p$ \`a "$p$ grand", c'est-\`a-dire $p\geq c$, o\`u $c$ est une constante d\'ependant de $G$ (c'est par exemple l'hypoth\`ese de \cite{W1}). Il est par contre certain que la d\'emonstration ne marche pas pour $p$ trop petit.

(7) L'hypoth\`ese que $G$ est quasi-d\'eploy\'e est restrictive. Par contre, celle qu'il  est adjoint et absolument simple  ne l'est gu\`ere. En g\'en\'eral, l'espace $SI(\mathfrak{g}(F))^*_{nil}$ est isomorphe au m\^eme espace relatif au groupe adjoint $G_{AD}$ de $G$. Un groupe adjoint $G$ se d\'ecompose en produit de composantes simples et l'espace $SI(\mathfrak{g}(F))^*_{nil}$  est produit tensoriel des espaces analogues pour chaque composante. Enfin, si $G$ est adjoint et simple, c'est la restriction \`a la Weil d'un  groupe adjoint absolument simple $G'$ d\'efini sur une extension $F'$ de $F$ de degr\'e fini et on se ram\`ene ais\'ement au groupe $G'$ en changeant de corps de base. 

Dans les deux paragraphes suivants, on explique la structure de la preuve.

\subsection{L'espace ${\cal D}^{st}(\mathfrak{g}(F))$}

Posons quelques d\'efinitions dans une situation  plus g\'en\'erale: $G$ est un groupe r\'eductif connexe d\'efini sur $F$ et on suppose que $p$ est tr\`es grand relativement \`a $G$. Notons $\mathfrak{g}_{ell}(F)$ le sous-ensemble des \'el\'ements semi-simples r\'eguliers et elliptiques de $\mathfrak{g}(F)$. On note $I_{cusp}(\mathfrak{g}(F))$ le sous-espace des \'el\'ements cuspidaux de $I(\mathfrak{g}(F))$, c'est-\`a-dire des $f\in I(\mathfrak{g}(F))$ telles que les int\'egrales orbitales $I^G(X,f)$ sont nulles pour tout $X\in \mathfrak{g}(F)$ qui est semi-simple r\'egulier mais pas elliptique. On note $I^{st}_{cusp}(\mathfrak{g}(F))$ le sous-espace des $f\in I_{cusp}(\mathfrak{g}(F))$ telles que les fonctions $X\mapsto I^G(X,f)$ soient constantes sur les classes de conjugaison stable  dans $\mathfrak{g}_{ell}(F)$. 

On note $G_{AD}$ le groupe adjoint de $G$ et $G_{SC}$ son rev\^etement simplement connexe. Introduisons    l'immeuble de Bruhat-Tits $Imm(G_{AD})$ de $G$. Il se d\'ecompose en facettes, on  note $S(G)$ l'ensemble des sommets. Pour $s\in S(G)$, on dispose du sous-groupe parahorique $K_{s}\subset G(F)$, de son plus grand sous-groupe distingu\'e pro-$p$-unipotent $K_{s}^+$ et du groupe r\'eductif connexe $G_{s}$ d\'efini sur ${\mathbb F}_{q}$ tel que $G_{s}({\mathbb F}_{q})=K_{s}/K_{s}^+$. 
Il y a des objets correspondants pour les alg\`ebres de Lie: $\mathfrak{k}_{s}\subset \mathfrak{g}(F)$, $\mathfrak{k}_{s}^+\subset \mathfrak{k}_{s}$ et  l'alg\`ebre de Lie $\mathfrak{g}_{s}$ de $G_{s}$. On introduit une transformation de Fourier $f\mapsto \hat{f}$ dans $C_{c}^{\infty}(\mathfrak{g}(F))$ poss\'edant les propri\'et\'es habituelles et telle que, pour tout $s\in S(G)$, on ait $\hat{{\bf 1}}_{s}=q^{dim(\mathfrak{g}_{s})/2}{\bf 1}_{s}^+$, o\`u ${\bf  1}_{s}$ et ${\bf  1}_{s}^+$ sont les fonctions caract\'eristiques de $\mathfrak{k}_{s}$ et $\mathfrak{k}_{s}^+$. Cette transformation de Fourier se descend en une transformation de l'espace $I(\mathfrak{g}(F))$. Notons $\mathfrak{g}_{tn}$ le sous-ensemble des \'el\'ements topologiquement nilpotents de $\mathfrak{g}$. 
On introduit le sous-espace $FC(\mathfrak{g}(F))\subset I(\mathfrak{g}(F))$ form\'e des \'el\'ements $f\in I(\mathfrak{g}(F))$ tels que les int\'egrales orbitales $I^G(X,f)$ et $I^G(X,\hat{f})$ sont nulles pour tout $X\in \mathfrak{g}(F)-\mathfrak{g}_{tn}(F)$. L'espace $FC(\mathfrak{g}(F))$ est contenu dans $I_{cusp}(\mathfrak{g}(F))$.  On pose $FC^{st}(\mathfrak{g}(F))=FC(\mathfrak{g}(F))\cap I^{st}_{cusp}(\mathfrak{g}(F))$.  On consid\`ere aussi l'espace $FC(\mathfrak{g}_{SC}(F))$. Le groupe $G(F)$ agit par conjugaison sur cet espace et on note $FC^{G}(\mathfrak{g}_{SC}(F))$ le sous-espace des invariants. On a l'inclusion $FC^{st}(\mathfrak{g}_{SC}(F))\subset FC^G(\mathfrak{g}_{SC}(F))$. 

L'espace $FC(\mathfrak{g}(F))$ se d\'ecrit de la fa\c{c}on suivante. Pour un groupe
r\'eductif connexe ${\bf G}$ d\'efini sur ${\mathbb F}_{q}$, Lusztig a d\'efini la notion de faisceau-caract\`ere cuspidal \`a support nilpotent dans l'alg\`ebre de Lie $\boldsymbol{\mathfrak{g}}$ de ${\bf G}$. Le groupe de Galois $\Gamma_{{\mathbb F}_{q}}$ de $\bar{{\mathbb F}}_{q}/{\mathbb F}_{q}$ agit sur ces faisceaux. Notons ${\bf FC}_{{\mathbb F}_{q}}(\boldsymbol{\mathfrak{g}})$ le sous-ensemble des faisceaux invariants par cette action. Pour ${\cal E}\in {\bf FC}_{{\mathbb F}_{q}}(\boldsymbol{\mathfrak{g}})$, on peut munir ${\cal E}$ d'une action de Frobenius (c'est-\`a-dire que l'on fixe un isomorphisme entre $Fr^*({\cal E})$ et ${\cal E}$, o\`u $Fr\in \Gamma_{{\mathbb F}_{q}}$ est l'\'el\'ement de Frobenius) et il  s'en d\'eduit une fonction caract\'eristique $f_{{\cal E}}$ sur $\boldsymbol{\mathfrak{g}}({\mathbb F}_{q})$, qui est \`a support nilpotent. Soit $s\in S(G)$, supposons $Z(G_{s})^0=\{1\}$. Pour ${\cal E}\in {\bf FC}_{{\mathbb F}_{q}}(\mathfrak{g}_{s})$, on d\'efinit la fonction $f_{{\cal E}}$ sur $\mathfrak{g}_{s}({\mathbb F}_{q})=\mathfrak{k}_{s}/\mathfrak{k}_{s}^+$. On l'identifie \`a une fonction sur $\mathfrak{k}_{s}$ invariante par translations par $\mathfrak{k}_{s}^+$ et on l'\'etend \`a $\mathfrak{g}(F)$ tout entier par $0$ hors de $\mathfrak{k}_{s}$. Notons $f_{{\cal E}}^{\mathfrak{g}}$ la fonction obtenue. Alors $FC(\mathfrak{g}(F))$ est engendr\'e par les images dans $I(\mathfrak{g}(F))$ de ces fonctions 
$f_{{\cal E}}^{\mathfrak{g}}$ quand $s$ d\'ecrit les sommets de $Imm(G_{AD})$ tels que $Z(G_{s})^0=\{1\}$ et ${\cal E}$ d\'ecrit ${\bf FC}_{{\mathbb F}_{q}}(\mathfrak{g}_{s})$. Dans la suite de cette introduction, on identifie un \'el\'ement de $FC(\mathfrak{g}(F))$ \`a un rel\`evement dans $C_{c}^{\infty}(\mathfrak{g}(F))$ qui est  combinaison lin\'eaire de fonctions $f_{{\cal E}}^{\mathfrak{g}}$. 
On a prouv\'e dans \cite{W7} que l'on pouvait d\'efinir un sous-ensemble $\underline{S}^{st}(G)$ de l'ensemble des sommets pr\'ec\'edents et, pour tout $s\in \underline{S}^{st}(G)$, un sous-ensemble ${\bf FC}^{st}_{{\mathbb F}_{q}}(\mathfrak{g}_{s})\subset {\bf FC}_{{\mathbb F}_{q}}(\mathfrak{g}_{s})$ de sorte que, quand $s$ d\'ecrit $\underline{S}^{st}(G)$ et ${\cal E}$ d\'ecrit ${\bf FC}^{st}_{{\mathbb F}_{q}}(\mathfrak{g}_{s})$, les images dans $I(\mathfrak{g}(F))$ des fonctions $f_{{\cal E}}^{\mathfrak{g}}$ forment une base de $FC^{st}(\mathfrak{g}(F))$.

On introduit un certain espace ${\cal D}(\mathfrak{g}(F))$ et une application antilin\'eaire injective $D^G:{\cal D}(\mathfrak{g}(F))\to I(\mathfrak{g}(F))^*$. Pour cette introduction, d\'ecrivons seulement les \'el\'ements de l'image de $D^G$. Notons $F^{nr}$ l'extension non ramifi\'ee maximale de $F$. Notons ${\cal L}^{nr}_{F}$ l'ensemble des sous-groupes $H\subset G$ qui sont d\'efinis sur $F$ et qui sont des $F^{nr}$-Levi, c'est-\`a-dire qu'il existe un sous-groupe parabolique $P$ de $G$ d\'efini sur $F^{nr}$ tel que $H$ soit une composante de Levi de $P$. Soit $H\in {\cal L}^{nr}_{F}$. Notons $A_{H}^{nr}$ le plus grand tore central dans $H$ qui soit déployé sur $F^{nr}$. Fixons un \'el\'ement $X_{H}\in \mathfrak{a}_{H}^{nr}(F)$ qui est entier et de r\'eduction r\'eguli\`ere. Cela signifie que pour tout caract\`ere alg\'ebrique $x^*$ de $A_{H}^{nr}$ (identifi\'e \`a une forme lin\'eaire sur $\mathfrak{a}_{H}^{nr}$), $x^*(X_{H})$ est un entier de $\bar{F}$ et que, pour toute racine $\alpha$ de $A_{H}^{nr}$ dans $\mathfrak{g}$, $\alpha(X_{H})$ est une unit\'e. Soit $\varphi\in FC^H(\mathfrak{h}_{SC}(F))$. Introduisons la distribution $D^G_{X_{H},\varphi}$ d\'efinie par
$$D^G_{X_{H},\varphi}(f)=\int_{  G(F)}\int_{\mathfrak{h}_{SC}(F)}f(g^{-1}(X_{H}+Z)g)\bar{\varphi}(Z)\,dZ\,dg$$
pour tout $f\in C_{c}^{\infty}(\mathfrak{g}(F))$ (on note l'action adjointe de $G$ sur $\mathfrak{g}$ comme une conjugaison). Les mesures doivent \^etre convenablement normalis\'ees. L'int\'egrale ci-dessus converge dans l'ordre indiqu\'e. Il s'av\`ere que la restriction \`a $\mathfrak{g}_{tn}(F)$ de la transformation de Fourier $\hat{D}^G_{X_{H},\varphi}$ ne d\'epend pas du choix de $X_{H}$. Notons $k(H,\varphi)$ la distribution \`a support dans $\mathfrak{g}_{tn}(F)$ qui co\"{\i}ncide avec cette restriction. 
Alors $D^G({\cal D}(\mathfrak{g}(F)))$ est l'espace de distributions  engendr\'e par les  $k(H,\varphi)$ quand $H$ d\'ecrit ${\cal L}^{nr}_{F}$ et $\varphi$ d\'ecrit $FC^H(\mathfrak{h}_{SC}(F))$. 

On introduit le sous-espace ${\cal H}\subset C_{c}^{\infty}(\mathfrak{g}(F))$ engendr\'e par les fonctions $f$ telles qu'il existe une sous-alg\`ebre d'Iwahori $\mathfrak{j}\subset \mathfrak{g}(F)$ de sorte que $f$ soit invariante par translations par $\mathfrak{j}$. Notons ${\cal H}^*$ son dual. On dispose de l'application de restriction $res_{{\cal H}}:I(\mathfrak{g}(F))^*\to {\cal H}^*$. Notons $\hat{D}^G({\cal D}(\mathfrak{g}(F)))$ l'espace  des transform\'ees de Fourier de $D^G({\cal D}(\mathfrak{g}(F)))$. En utilisant un r\'esultat de DeBacker, cf.  \cite{D}, on montre que les restrictions de $res_{{\cal H}}$ \`a $I(\mathfrak{g}(F))^*_{nil}$ et \`a $\hat{D}^G({\cal D}(\mathfrak{g}(F)))$ sont injectives et ont m\^eme image. Il existe donc un unique isomorphisme (antilin\'eaire) $\delta:I(\mathfrak{g}(F))_{nil}^*\to {\cal D}(\mathfrak{g}(F))$ de sorte que $res_{{\cal H}}(J)=res_{{\cal H}}(\hat{D}^G(\delta(J)))$ pour tout $J\in I(\mathfrak{g}(F))^*_{nil}$. Notons ${\cal D}^{st}(\mathfrak{g}(F))$ le sous-espace des \'el\'ements $d\in {\cal D}(\mathfrak{g}(F))$ tels que $D^G(d)$ est stable. On prouve facilement, cf. lemme \ref{descriptionnil}, que 

(1) $\delta(SI(\mathfrak{g}(F))^*_{nil})={\cal D}^{st}(\mathfrak{g}(F))$. 

Pour exploiter cette \'egalit\'e, on doit d\'ecrire l'espace ${\cal D}^{st}(\mathfrak{g}(F))$.   
  Soient $H,H'\in {\cal L}_{F}^{nr}$. Disons que ces deux groupes sont stablement conjugu\'es s'il existe $g\in G$ tel que $H'=g^{-1}Hg$ et, pour tout $\sigma\in \Gamma_{F}$, $g\sigma(g)^{-1}$ appartient \`a $H$. Supposons que $H$ et $H'$ soient stablement conjugu\'es et fixons $g$ v\'erifiant les propri\'et\'es pr\'ec\'edentes. Alors l'application $Ad(g^{-1})$ se restreint en un torseur int\'erieur  $\psi:H\to H'$. Par une propri\'et\'e  basique de l'endoscopie, un tel torseur d\'efinit un isomorphisme $transfert_{\psi}:I_{cusp}^{st}(\mathfrak{h}_{SC}(F))\to I_{cusp}^{st}(\mathfrak{h}'_{SC}(F))$. Une cons\'equence des r\'esultats de \cite{W6} est que cet isomorphisme se restreint en un isomorphisme de $FC^{st}(\mathfrak{h}_{SC}(F))$ sur $FC^{st}(\mathfrak{h}_{SC}'(F))$. On montre que celui-ci ne d\'epend pas du choix de l'\'el\'ement $g$, on le note $\iota_{H',H}$. Notons ${\cal L}^{nr,st}_{F}$ l'ensemble des $H\in {\cal L}_{F}^{nr}$ tels que $FC^{st}(\mathfrak{h}_{SC}(F))\not=\{0\}$. Il est conserv\'e par conjugaison stable. Notons
  ${\cal L}^{nr,st}_{F}/st-conj$ l'ensemble des   classes de conjugaison stable dans ${\cal L}^{nr,st}_{F}$. Soit ${\cal C}\in {\cal L}^{nr,st}_{F}/st-conj$. Fixons des repr\'esentants $(H_{i})_{i=1,...,n}$ des classes de conjugaison par $G(F)$ dans  ${\cal C}$.  Pour $\varphi\in FC^{st}(\mathfrak{h}_{1,SC}(F))$, posons
 $$k^{st}({\cal C},\varphi)=\sum_{i=1,...,n}\vert W^{nr}(H_{i})^{\Gamma_{F}^{nr}}\vert \vert W_{F}(H_{i})\vert ^{-1} k(H_{i},\iota_{H_{i},H_{1}}(\varphi)).$$
 Les constantes sont des nombres d'\'el\'ements de groupes de Weyl qui seront d\'efinis plus tard. Dans cette formule, le groupe $H_{1}$ joue un r\^ole  particulier mais c'est insignifiant, les espaces $FC^{st}(\mathfrak{h}_{SC}(F))$ \'etant canoniquement isomorphes pour $H\in {\cal C}$. Le r\'esultat principal de la premi\`ere section est que l'espace ${\cal D}^{st}(\mathfrak{g}(F))$ est le sous-espace de ${\cal D}(\mathfrak{g}(F))$ tel que $D^G({\cal D}^{st}(\mathfrak{g}(F))$ soit engendr\'e par les distributions $k^{st}({\cal C},\varphi)$ quand ${\cal C}$ d\'ecrit ${\cal L}^{nr,st}_{F}/st-conj$ et $\varphi$ d\'ecrit $FC^{st}(\mathfrak{h}_{1,SC}(F))$ ($H_{1}$ \'etant un \'el\'ement de ${\cal C}$). Plus pr\'ecis\'ement ces distributions forment une base de $D^G({\cal D}^{st}(\mathfrak{g}(F))$ si $\varphi$ d\'ecrit seulement une base de $FC^{st}(\mathfrak{h}_{1,SC}(F))$.

Imposons d\'esormais les hypoth\`eses du th\'eor\`eme: $G$ est quasi-d\'eploy\'e sur $F$, adjoint et absolument simple. Fixons une paire de Borel \'epingl\'ee $\mathfrak{E}=(B,T,(E_{\alpha})_{\alpha\in \Delta})$ de $G$ conserv\'ee par l'action de $\Gamma_{F}$ ($B$ est un sous-groupe de Borel, $T\subset B$ est un sous-tore maximal et $(E_{\alpha})_{\alpha\in \Delta}$ est un \'epinglage). On note $W$ le groupe de Weyl de $G$ relatif \`a $T$. Le groupe $\Gamma_{F}$ agit sur $W$, on note $W^{I_{F}}$ le sous-groupe des invariants par le groupe d'inertie $I_{F}\subset \Gamma_{F}$.  Notons ${\cal L}_{F,sd}^{st}$ l'ensemble des $F$-groupes de Levi standard $M$ tels que $FC^{st}(\mathfrak{m}_{SC}(F))\not=\{0\}$. Pour un tel groupe, notons $Norm_{W^{I_{F}}}(M)$ le sous-groupe des \'el\'ements de $W^{I_{F}}$ qui conservent $M$, $W^{M,I_{F}}$ l'analogue de $W^{I_{F}}$ quand on remplace $G$ par $M$ et $W^{I_{F}}(M)=Norm_{W^{I_{F}}}(M)/W^{M,I_{F}}$. Posons $\Gamma_{F}^{nr}=\Gamma_{F}/ I_{F}$. Ce groupe est isomorphe au groupe de Galois $\Gamma_{{\mathbb F}_{q}}$  et est topologiquement engendr\'e par  l'\'el\'ement de Frobenius $Fr$.  Il  agit sur $W^{I_{F}}(M)$. On dit  que deux \'el\'ements $w,w'\in W^{I_{F}}(M)$ sont $Fr$-conjugu\'es s'il existe $u\in W^{I_{F}}(M)$ tel que $w'=Fr(u)wu^{-1}$ (la d\'efinition usuelle est plut\^ot $uwFr(u)^{-1}$, des raisons techniques nous ont conduit \`a la modifier). On note $W^{I_{F}}(M)/Fr-conj$ l'ensemble des classes de $Fr$-conjugaison. On montre en \ref{parametrage} qu'il existe une bijection naturelle de ${\cal L}^{nr,st}_{F}/st-conj$ sur l'ensemble des paires $(M,Cl_{Fr}(w))$ o\`u $M\in {\cal L}_{F,sd}^{st}$ et $Cl_{Fr}(w)\in W^{I_{F}}(M)/Fr-conj$. Plus pr\'ecis\'ement, soient ${\cal C}\in {\cal L}^{nr,st}_{F}/st-conj$, notons $(M,Cl_{Fr}(w))$ son image par cette bijection et soit $H\in {\cal C}$. Alors $H$ est conjugu\'e \`a $M$ par un \'el\'ement de $G(F^{nr})$, il y a une bijection naturelle de $\underline{S}^{st}(H_{AD})$ sur $\underline{S}^{st}(M_{AD})$ et, pour tout $s_{H}\in \underline{S}^{st}(H_{AD})$, d'image $s_{M}$ dans $\underline{S}^{st}(M_{AD})$, il y a une bijection naturelle de ${\bf FC}_{{\mathbb F}_{q}}^{st}(\mathfrak{h}_{SC,s_{H}})$ sur ${\bf FC}_{{\mathbb F}_{q}}^{st}(\mathfrak{m}_{SC,s_{M}})$. Cela  entra\^{\i}ne que  les espaces $FC^{st}(\mathfrak{h}_{SC}(F))$ et $FC^{st}(\mathfrak{m}_{SC}(F))$ sont de m\^eme dimension. Mais cela ne suffit pas pour identifier ces espaces. En effet, les fonctions caract\'eristiques $f_{{\cal E}}$ introduites ci-dessus d\'ependent du choix d'une action de Frobenius sur le faisceau ${\cal E}$ et, \`a notre connaissance, il n'y a pas de choix canonique d'une telle action. Les groupes $H$ et $M$ sont isomorphes sur $F^{nr}$ puisqu'ils sont conjugu\'es par un \'el\'ement de $G(F^{nr}) $ mais l'isomorphisme n'est pas en g\'en\'eral un torseur int\'erieur. Il n'y a pas de transfert endoscopique qui puisse nous aider. 
 Le probl\`eme est r\'esolu dans la section 4, qui est la plus technique de l'article. Soient $H$, $s_{H}$, $M$, $s_{M}$ comme pr\'ec\'edemment. Soit ${\cal E}_{H}\in {\bf FC}_{{\mathbb F}_{q}}^{st}(\mathfrak{h}_{SC,s_{H}})$, notons ${\cal E}_{M}$ son image dans $ {\bf FC}_{{\mathbb F}_{q}}^{st}(\mathfrak{m}_{SC,s_{M}})$. Munissons ${\cal E}_{M}$ d'une action de Frobenius. On d\'efinit en \ref{isomorphismesFC}  un proc\'ed\'e pour en d\'eduire une action de Frobenius sur ${\cal E}_{H}$.  On en d\'eduit un isomorphisme canonique $\iota_{H,M}:FC^{st}(\mathfrak{m}_{SC}(F))\to FC^{st}(\mathfrak{h}_{SC}(F))$. On montre en \ref{comparaison} que cet isomorphisme est compatible au transfert endoscopique, c'est-\`a-dire que si $H'\in {\cal L}^{nr,st}_{F}$ est stablement conjugu\'e \`a $H$, le diagramme suivant est commutatif:
 $$\begin{array}{ccc}&&FC^{st}(\mathfrak{h}_{SC}(F))\\ &\iota_{H,M}\nearrow\,\,&\\ FC^{st}(\mathfrak{m}_{SC}(F))&&\,\,\downarrow \iota_{H',H}\\&\iota_{H',M}\searrow\,\,\\ && FC^{st}(\mathfrak{h}'_{SC}(F))\\ \end{array}$$

  Soient $M\in {\cal L}_{F,sd}^{st}$ et $s_{M}\in \underline{S}^{st}(M_{AD})$. Pour tout ${\cal E}\in {\bf FC}^{st}_{{\mathbb F}_{q}}(\mathfrak{m}_{SC,s_{M}})$, on munit ${\cal E}$ d'une action de Frobenius. On note  ${\cal B}(M,s_{M})$
l'ensemble des fonctions  caract\'eristiques $f_{{\cal E}}$   quand ${\cal E}$ d\'ecrit $ {\bf FC}^{st}_{{\mathbb F}_{q}}(\mathfrak{m}_{SC,s_{M}})$. Notons ${\cal B}_{W}(G)$ l'ensemble des quadruplets $(M,s_{M},\varphi,w)$ o\`u $M$ et $s_{M}$ sont comme ci-dessus, $\varphi\in {\cal B}(M,s_{M})$ et $w$ parcourt $W^{I_{F}}(M)$ \`a $Fr$-conjugaison pr\`es. Consid\'erons un tel quadruplet. Au couple $(M,w)$ est associ\'e une classe de conjugaison stable ${\cal C}\subset {\cal L}^{nr,st}_{F}$. Fixons $H_{1}\in {\cal C}$. Posons simplement $\varphi_{H_{1}}=\iota_{H_{1},M}(\varphi^{\mathfrak{m}_{SC}})$. On a d\'efini ci-dessus la distribution $k^{st}({\cal C},\varphi_{H_{1}})$. On pose $k(M,s_{M},\varphi,w)=c(w)k^{st}({\cal C},\varphi_{H_{1}})$, o\`u $c(w)$ est une constante non nulle qui sera d\'efinie dans l'article. Alors la famille $(k(M,s_{M},\varphi,w))_{(M,s_{M},\varphi,w)\in {\cal B}_{W}(G)}$ est une base de l'espace $D^G({\cal D}^{st}(\mathfrak{g}(F)))$. Soit $M\in {\cal L}_{F,sd}^{st}$. Notons $Irr(W^{I_{F}}(M))$ l'ensemble des classes de repr\'esentations irr\'eductibles de $W^{I_{F}}(M)$. Le groupe $\Gamma_{F}^{nr}$ agit sur $W^{I_{F}}(M)$ donc sur $Irr(W^{I_{F}}(M))$. Soit $Irr_{F}(W^{I_{F}}(M))$ le sous-ensemble des points fixes. Pour $\rho\in Irr_{F}(W^{I_{F}}(M))$, on peut prolonger $\rho$ en une repr\'esentation irr\'eductible du produit semi-direct $W^{I_{F}}(M)\rtimes \Gamma_{F}^{nr}$ et on fixe un tel prolongement $\rho^{\flat}$. On note ${\cal B}_{Irr}(G)$ l'ensemble des quadruplets $(M,s_{M},\varphi,\rho)$ o\`u $M$ et $s_{M}$ sont comme ci-dessus, $\varphi\in {\cal B}(M,s_{M})$ et $\rho\in Irr_{F}(W^{I_{F}}(M))$. Pour un tel quadruplet, on pose
$$k(M,s_{M},\varphi,\rho)=\vert W^{I_{F}}(M)\vert ^{-1}\sum_{w\in W^{I_{F}}(M)}trace(\rho^{\flat}(wFr^{-1})k(M,s_{M},\varphi,w).$$
 Alors 
 
 (2) la famille $(k(M,s_{M},\varphi,\rho))_{(M,s_{M},\varphi,\rho)\in {\cal B}_{Irr}(G)}$ est  encore une base de $D^G({\cal D}^{st}(\mathfrak{g}(F)))$. 
 
 Pour $(M,s_{M},\varphi,\rho)\in {\cal B}_{Irr}(G)$, notons $I(M,s_{M},\varphi,\rho)$ l'\'el\'ement de $I(\mathfrak{g}(F))^*_{nil}$ tel que $D^G\circ \delta(I(M,s_{M},\varphi,\rho))=k(M,s_{M},\varphi,\rho)$. La famille $(I(M,s_{M},\varphi,\rho))_{(M,s_{M},\varphi,\rho)\in {\cal B}_{Irr}(G)}$ est une base de $SI(\mathfrak{g}(F))^*_{nil}$ d'après (1) et (2). Pour en d\'eduire le th\'eor\`eme, on a besoin de quelques  renseignements sur l'\'ecriture d'une distribution $I(M,s_{M},\varphi,\rho)$ dans la base $(I_{{\mathbb O}})_{{\mathbb O}\in \mathfrak{g}_{nil}(F)/conj}$ de $I(\mathfrak{g}(F))_{nil}^*$. 
 
 \subsection{Param\'etrages}
 
 Pour cela, on doit commencer par param\'etrer l'ensemble ${\cal B}_{Irr}(G)$. Introduisons l'extension mod\'er\'ement ramifi\'ee maximale $F^{mod}$ de $F$ contenue dans $\bar{F}$. Notons $\Gamma_{F}^{mod}$, resp. $I_{F}^{mod}$,  le groupe de Galois de $F^{mod}/F$, resp. $F^{mod}/F^{nr}$. Modulo certains choix, on peut identifier $\Gamma_{F}^{mod}$ au produit semi-direct $I_{F}^{mod}\rtimes \Gamma_{F}^{nr}$. On a fix\'e un \'epinglage $\mathfrak{E}$ de $G$ conserv\'e par $\Gamma_{F}$. On en d\'eduit de fa\c{c}on usuelle une "donn\'ee de racines", c'est-\`a-dire un sextuplet $(X^*(T),X_{*}(T),\Sigma,\check{\Sigma},\Delta,\check{\Delta})$ ($\Sigma$ est l'ensemble des racines de $T$ dans $\mathfrak{g}$, $\check{\Sigma}$ l'ensemble des coracines etc...). Le groupe $\Gamma_{F}$ agit sur cette donn\'ee et l'hypoth\`ese sur $p$ assure que cette action se factorise en une action de $\Gamma_{F}^{mod}$. En particulier, puisqu'on a identifi\'e $\Gamma_{F}^{nr}\simeq \Gamma_{{\mathbb F}_{q}}$ \`a un sous-groupe de $\Gamma_{F}^{mod}$, il y a une action de $\Gamma_{{\mathbb F}_{q}}$ sur la donn\'ee de racines. On peut alors introduire un groupe r\'eductif connexe ${\bf G}$ d\'efini sur ${\mathbb F}_{q}$ dont la donn\'ee de racines est isomorphe \`a celle de $G$, cet isomorphisme \'etant \'equivariant pour les actions de $\Gamma_{{\mathbb F}_{q}}$ (si $G$ est d\'eploy\'e sur $F$, ${\bf G}$ est simplement le groupe adjoint simple et d\'eploy\'e sur ${\mathbb F}_{q}$ de m\^eme type que $G$). L'action de $I_{F}^{mod}$ sur la donn\'ee de racines d\'efinit une action alg\'ebrique de $I_{F}^{mod}$ sur ${\bf G}$. Remarquons que cette action n'est pas forc\'ement d\'efinie sur ${\mathbb F}_{q}$. D'autre part, on conna\^{\i}t la structure de  $I_{F}^{mod}$. On peut fixer un pro-g\'en\'erateur $\gamma$ de ce groupe et l'action de $I_{F}^{mod}$ est d\'etermin\'ee par l'unique automorphisme associ\'e \`a $\gamma$.

  Comme plus haut, on note $\boldsymbol{\mathfrak{g}}_{nil}/conj$ l'ensemble des orbites nilpotentes dans l'alg\`ebre de Lie $\boldsymbol{\mathfrak{g}}=\boldsymbol{\mathfrak{g}}(\bar{{\mathbb F}}_{q})$ de ${\bf G}={\bf G}(\bar{{\mathbb F}}_{q})$.  Le groupe $I_{F}^{mod}$ agit sur ${\bf G}$ via un quotient fini que l'on note $\Gamma_{F^G/F^{nr}}$. En g\'en\'eral, ce quotient est d'ordre $1$ ou $2$ mais il peut \^etre d'ordre $3$ dans le cas o\`u $G$ est un groupe trialitaire de type $D_{4}$.  Introduisons le produit semi-direct ${\bf G}\rtimes \Gamma_{F^G/F^{nr}}$ et son sous-ensemble $\tilde{{\bf G}}={\bf G}\times\{\gamma\}$ (o\`u on note encore $\gamma$ l'image dans $\Gamma_{F^G/F^{nr}}$ du pro-g\'en\'erateur fix\'e de $I_{F}^{mod}$). Soit $N\in \boldsymbol{\mathfrak{g}}_{nil}$. On suppose que son orbite est conservée par l'action de $I_{F}^{mod}$. 
   On d\'efinit comme plus haut les groupes $Z_{{\bf G}}(N)$, $Z_{{\bf G}}(N)^0$, $A(N)=Z_{{\bf G}}(N)/Z_{{\bf G}}(N)^0$ et le quotient de Lusztig $\bar{A}(N)$.  On note $Z_{\tilde{{\bf G}}}(N)$ le sous-ensemble des \'el\'ements $(g,\gamma)\in \tilde{{\bf G}}$ telles que $g\gamma(N)g^{-1}=N$. Le groupe $Z_{{\bf G}}(N)$ agit sur cet ensemble par multiplication \`a gauche ou \`a droite et $Z_{\tilde{{\bf G}}}(N)$ est un bitorseur pour ces actions, c'est-\`a-dire que, pour chacune d'elle, $Z_{\tilde{{\bf G}}}(N)$ est un espace principal homog\`ene sous $Z_{{\bf G}}(N)$.  On pose $\tilde{A}(N)=Z_{{\bf G}}(N)^0\backslash Z_{\tilde{{\bf G}}}(N)$. Les actions de $Z_{{\bf G}}(N)$ se quotientent en des actions de $A(N)$ sur $\tilde{A}(N)$ et $\tilde{A}(N)$ est un bitorseur pour ces actions.  On montrera qu'on peut  quotienter  $\tilde{A}(N)$ en un ensemble   $\tilde{\bar{A}}(N)$ qui est un bitorseur pour des actions \`a droite et \`a gauche du groupe $\bar{A}(N)$. Notons $\bar{{\cal C}}^{Irr}$ l'ensemble des classes de conjugaison par ${\bf G}$ dans l'ensemble des triplets $(N,d,\mu)$, où $N$ est un élément de $\boldsymbol{\mathfrak{g}}_{nil}$ dont l'orbite est conservée par $I_{F}^{mod}$, $d
   \in \tilde{\bar{A}}(N)$ et $\mu$ est une repr\'esentation irr\'eductible du commutant $Z_{\bar{A}(N)}(d)$ de $d$ dans $\bar{A}(N)$. Notons $\bar{{\cal C}}^{Irr}_{sp,F}$ le sous-ensemble où on se limite aux $N$ dont l'orbite est spéciale et est conservée par le groupe  $\Gamma_{F}^{mod}$ tout entier. 
     En \ref{parametragecalB}, nous d\'efinirons une application injective
 $$\boldsymbol{\nabla}_{F}:{\cal B}_{Irr}(G)\to \bar{{\cal C}}^{Irr}_{sp,F}.$$
 
 Indiquons l'origine de cette d\'efinition. Supposons que $G$ soit d\'eploy\'e. L'action de $\gamma$ sur ${\bf G}$ est triviale, on peut l'oublier et supprimer les $\tilde{}$ dans les d\'efinitions ci-dessus. Toute orbite nilpotente dans $\boldsymbol{\mathfrak{g}}_{nil}$ est conservée par $\Gamma_{F}^{mod}$.    Dans \cite{L10} th\'eor\`eme 2.4, Lusztig \'etablit une bijection $\nabla:{\cal A}({\bf G})\to \bar{{\cal C}}^{Irr}_{sp,F}$, o\`u ${\cal A}({\bf G})$ est l'ensemble des faisceaux-caract\`eres unipotents sur ${\bf G}$. On sait d\'ecrire l'ensemble ${\cal A}({\bf G})$ en d\'ecomposant les complexes induits \`a partir de faisceaux-caract\`eres unipotents et cuspidaux sur  des groupes de Levi standard de ${\bf G}$. Un tel groupe ${\bf M}$ correspond \`a un groupe de Levi standard $M$ de $G$. On sait d\'ecrire les faisceaux-caract\`eres unipotents et cuspidaux sur ${\bf M}$.   Les descriptions de \cite{W7} paragraphe 9 montrent qu'il y a une injection naturelle qui, \`a un couple $(s_{M},\varphi)$ associe un tel faisceau-caract\`ere unipotent et cuspidal sur ${\bf M}$, $(s_{M},\varphi)$ parcourant les couples ci-dessus, c'est-\`a-dire $s_{M}\in \underline{S}^{st}(M_{AD})$ et $\varphi\in {\cal B}(M,s_{M})$.  Cette injection est bijective si $G$ est classique ou si $q-1$ est divisible par $3\times 4\times 5$. On en d\'eduit facilement une injection naturelle de ${\cal B}_{Irr}(G)$ dans ${\cal A}({\bf G})$. Alors $\boldsymbol{\nabla}_{F}$ est la compos\'ee de cette injection et de la bijection $\nabla$ de Lusztig. Cette construction se g\'en\'eralise au cas o\`u $G$ est d\'eploy\'e sur $F^{nr}$. 
 Dans le cas o\`u $G$ n'est pas d\'eploy\'e  sur $F^{nr}$, on donnera une autre d\'efinition de $\boldsymbol{\nabla}_{F}$. Il est possible qu'elle se rattache \`a une analogue de la bijection $\nabla$ pour le groupe non connexe ${\bf G}\rtimes \Gamma_{F^G/F^{nr}}$.  Mais nous avons pr\'ef\'er\'e utiliser des moyens plus terre \`a terre. 
 
 Pour calculer l'\'ecriture d'une distribution appartenant \`a $ I(\mathfrak{g}(F))_{nil}^*$ dans la base $(I_{{\mathbb  O}})_{{\mathbb  O}\in \mathfrak{g}_{nil}(F)/conj}$, on doit introduire une famille de fonctions dans $I(\mathfrak{g}(F))$ qui s\'eparent les \'el\'ements de la base. La construction d'une telle famille est connue, cf. \cite{W1} IX pour les groupes classiques non ramifi\'es et \cite{DK} pour un groupe de type $G_{2}$. Notons ${\cal J}_{F}^*$ l'ensemble des couples $(s,{\cal O}_{s})$ o\`u  $s\in S(G)$ et ${\cal O}_{s} \in \mathfrak{g}_{s,nil}({\mathbb F}_{q}) /conj$ (avec une notation compr\'ehensible). Dans \cite{D2}, DeBacker associe \`a $(s,{\cal O}_{s})\in {\cal J}^*_{F}$ une orbite ${\mathbb  O}\in \mathfrak{g}_{nil}(F)/conj$, qui est caract\'eris\'ee de la fa\c{c}on suivante. Rappelons que $\mathfrak{g}_{s}({\mathbb F}_{q})=\mathfrak{k}_{s}/\mathfrak{k}_{s}^+$, il y a donc une application de r\'eduction $\mathfrak{k}_{s}\to \mathfrak{g}_{s}({\mathbb F}_{q})$. Consid\'erons l'ensemble $E(s,{\cal O}_{s})$ des $e\in \mathfrak{k}_{s}\cap \mathfrak{g}_{nil}(F)$ dont la r\'eduction appartient \`a ${\cal O}_{s}$. Alors ${\mathbb  O}$ est la plus petite orbite nilpotente dans $\mathfrak{g}_{nil}(F)$ qui coupe $E(s,{\cal O}_{s})$. C'est-\`a-dire que  ${\mathbb  O}\cap E(s,{\cal O}_{s})\not=\emptyset$ et, si ${\mathbb  O}'\in \mathfrak{g}_{nil}(F)/conj$  v\'erifie ${\mathbb  O}'\cap E(s,{\cal O}_{s})\not=\emptyset$, alors ${\mathbb  O}$ est contenue dans l'adh\'erence de ${\mathbb O}'$ pour la topologie $p$-adique.  A l'aide de la th\'eorie des $\mathfrak{sl}(2)$-triplets, on associe \`a $(s,{\cal O}_{s})$ une fonction $h_{s,{\cal O}_{s}}\in C_{c}^{\infty}(\mathfrak{g}(F))$, \`a support dans $\mathfrak{k}_{s}$ et invariante par translations par $\mathfrak{k}_{s}^+$, qui v\'erifie les propri\'et\'es suivantes:
 
 $h_{s,{\cal O}_{s}}\in {\cal H}$;
 
 soit ${\mathbb O}'\in \mathfrak{g}_{nil}(F)/conj$, supposons que $I_{{\mathbb  O}'}(h_{s,{\cal O}_{s}})\not=0$; alors ${\mathbb  O}$ est contenue dans l'adh\'erence de ${\mathbb O}'$;
 
 $I_{{\mathbb  O}}(h_{s,{\cal O}_{s}})=1$.

 Pour tout ${\mathbb  O}\in \mathfrak{g}_{nil}(F)/conj$, il existe au moins un couple $(s,{\cal O}_{s})\in {\cal J}^*_{F}$  qui s'envoie sur ${\mathbb  O}$ par l'application de DeBacker. La famille $(h_{s,{\cal O}_{s}})_{(s,{\cal O}_{s})\in {\cal J}^*_{F}}$ s\'epare donc la famille de distributions $(I_{{\mathbb  O}})_{{\mathbb O}\in \mathfrak{g}_{nil}(F)/conj}$. 

Pour $\boldsymbol{{\cal O}}\in (\boldsymbol{\mathfrak{g}}_{nil}/conj)^{\Gamma_{F}^{mod}}$, notons $\bar{{\cal C}}^{\sharp}_{F}(\boldsymbol{{\cal O}})$ l'ensemble des classes de conjugaison par ${\bf G}$    dans l'ensemble des triplets  $(N,d,v)$, où $N\in \boldsymbol{{\cal O}}$ et $(d,v)\in\tilde{\bar{A}}(N)\times \bar{A}(N)$ vérifient la relation $vdv^{-1}=d^q$ (cette d\'efinition doit \^etre modifi\'ee dans le cas d'un groupe trialitaire de type $D_{4}$).   Un \'el\'ement de $\bar{{\cal C}}_{F}^{\sharp}(\boldsymbol{{\cal O}})$ sera not\'e comme un triplet le repr\'esentant.   Notons $\bar{{\cal C}}^{\sharp}_{F}$ la réunion des $\bar{{\cal C}}_{F}^{\sharp}(\boldsymbol{{\cal O}})$ sur les orbites $\boldsymbol{{\cal O}}\in (\boldsymbol{\mathfrak{g}}_{nil}/conj)^{\Gamma_{F}^{mod}}$.  Dans \ref{iotanil}, nous d\'efinirons une application surjective $\bar{c}_{F}:{\cal J}^*_{F}\to \bar{{\cal C}}^{\sharp}_{F}$. La premi\`ere composante de cette application se d\'efinit facilement. Soit $(s,{\cal O}_{s})\in {\cal J}^*_{F}$. On lui a associ\'e une orbite ${\mathbb  O}\in \mathfrak{g}_{nil}(F)/conj$. Cette orbite est contenue dans une orbite $\bar{{\mathbb O}}$
sur la cl\^oture alg\'ebrique $\bar{F}$ de $F$.   Il y a une bijection naturelle entre $\mathfrak{g}_{nil}/conj$ et $\boldsymbol{\mathfrak{g}}_{nil}/conj$: les deux groupes $G$ et ${\bf G}$ sont de m\^eme type, l'un sur $\bar{F}$, l'autre sur $\bar{{\mathbb F}}_{q}$,  et on sait que la classification des orbites nilpotentes ne d\'epend pas du corps de base, pourvu que celui-ci soit alg\'ebriquement clos de caract\'eristique nulle ou positive et assez grande. Ainsi, \`a $\bar{{\mathbb O}}$ correspond une orbite $\boldsymbol{{\cal O}}$  dans $\boldsymbol{\mathfrak{g}}_{nil}/conj$.   Cette orbite $\boldsymbol{{\cal O}}$ est conservée par $\Gamma_{F}^{mod}$. On a $\bar{c}_{F}(s,{\cal O}_{s})=(N,d,v)$ avec $N\in \boldsymbol{{\cal O}}$.    Les  composantes  $d$ et $v$ sont reli\'ees \`a l'ensemble de cohomologie $H^1(\Gamma_{F},\bar{A}(e))$, o\`u $e$ est un \'el\'ement de ${\mathbb O}$. On montrera d'ailleurs que cet ensemble est en bijection naturelle avec $\bar{{\cal C}}^{\sharp}_{F}(\boldsymbol{{\cal O}})$. 

 On peut maintenant \'enoncer le calcul cl\'e de l'article. Soient $(M,s_{M},\varphi,\rho)\in {\cal B}_{Irr}(G)$ et $(s,{\cal O}_{s})\in {\cal J}^*_{F}$. Posons $\boldsymbol{\nabla}_{F}(M,s_{M},\varphi,\rho)=(N,d,\mu)$ et $\bar{c}_{F}(s,{\cal O}_{s})=(N',d',v')$. On note $\boldsymbol{{\cal O}}$, resp. $\boldsymbol{{\cal O}}'$, l'orbite de $N$, resp. $N'$. On sait que $\boldsymbol{{\cal O}}$ est spéciale.  On note $\hat{k}(M,s_{M},\varphi,\rho)$ la transform\'ee de Fourier de $k(M,s_{M},\varphi,\rho)$. On a alors (cf. \ref{uncalculplusprecis}):

(1)(i) si $\hat{k}(M,s_{M},\varphi,\rho)(h_{s,{\cal O}_{s}})\not=0$, alors ou bien $dim(\boldsymbol{{\cal O}}')<dim(\boldsymbol{{\cal O}})$,  ou  bien $\boldsymbol{{\cal O}}'=\boldsymbol{{\cal O}}$ et  les couples $(N,d)$ et $(N',d')$ sont conjugués par un élément de ${\bf G}$; 

(ii) si $(N,d)=(N',d')$, $\hat{k}(M,s_{M},\varphi,\rho)(h_{s,{\cal O}_{s}})$ est non nul et s'exprime par une formule explicite.

Pour d\'emontrer ce r\'esultat, on se ram\`ene  \`a un calcul dans l'espace $\mathfrak{g}_{s}({\mathbb F}_{q})$ et on utilise de fa\c{c}on essentielle les constructions de Lusztig (correspondance de Springer g\'en\'eralis\'ee, fonctions de Green etc...), en particulier les r\'esultats de \cite{L8}. 

Soit $(M,s_{M},\varphi,\rho)\in {\cal B}_{Irr}(G)$. Posons $\boldsymbol{\nabla}_{F}(M,s_{M},\varphi,\rho)=(N,d,\mu)$ et notons $\boldsymbol{{\cal O}}$ l'orbite de $N$. Ecrivons $I(M,s_{M},\varphi,\rho)$ dans la base $(I_{{\mathbb O}})_{{\mathbb  O}\in \mathfrak{g}_{nil}(F)/conj}$:
$$I(M,s_{M},\varphi,\rho)=\sum_{{\mathbb O}\in \mathfrak{g}_{nil}(F)/conj}c(M,s_{M},\varphi,\rho,{\mathbb O})I_{{\mathbb  O}}.$$
 On a expliqu\'e ci-dessus qu'\`a une orbite ${\mathbb  O}\in \mathfrak{g}_{nil}(F)/conj$, on pouvait associer une orbite dans $\boldsymbol{\mathfrak{g}}_{nil}$, notons-la $\boldsymbol{{\cal O}}_{{\mathbb O}}$. 
  A l'aide de (1)(i), on montre que
 
 (2) soit ${\mathbb O}\in \mathfrak{g}_{nil}(F)/conj$, supposons $c(M,s_{M},\varphi,\rho,{\mathbb  O})\not=0$; alors $dim(\boldsymbol{{\cal O}}_{{\mathbb O}})<dim(\boldsymbol{{\cal O}})$ ou $ \boldsymbol{{\cal O}}_{{\mathbb O}}=\boldsymbol{{\cal O}}$. 
 
 Posons alors 
 $$I^{max}(M,s_{M},\varphi,\rho)=\sum_{{\mathbb O}\in \mathfrak{g}_{nil}(F)/conj, \boldsymbol{{\cal O}}_{{\mathbb  O}}=\boldsymbol{{\cal O}}}c(M,s_{M},\varphi,\rho,{\mathbb  O})I_{{\mathbb  O}}.$$
 C'est la somme des termes "maximaux" de la d\'ecomposition de $I(M,s_{M},\varphi,\rho)$. En utilisant l'homog\'en\'eit\'e des int\'egrales orbitales nilpotentes, on v\'erifie que la stabilit\'e de $I(M,s_{M},\varphi,\rho)$ entra\^{\i}ne celle de $I^{max}(M,s_{M},\varphi,\rho)$. En utilisant de nouveau (1)(i) et maintenant (1)(ii), on prouve que $(I^{max}(M,s_{M},\varphi,\rho))_{(M,s_{M},\varphi,\rho)\in {\cal B}_{Irr}(G)}$ est une base de $SI(\mathfrak{g}(F))^*_{nil}$. La premi\`ere assertion du th\'eor\`eme s'en d\'eduit. En effet, une distribution $I^{max}(M,s_{M},\varphi,\rho)$ est port\'ee par une unique orbite $\bar{{\mathbb  O}}\in (\mathfrak{g}_{nil}/conj)^{\Gamma_{F}}$ et celle-ci est sp\'eciale: si l'on \'ecrit $\boldsymbol{\nabla}_{F}(M,s_{M},\varphi,\rho)=(N,d,\mu)$, $\bar{{\mathbb O}}$ est l'orbite dans $\mathfrak{g}_{nil}$ qui correspond \`a l'orbite de $N$ dans $\boldsymbol{\mathfrak{g}}_{nil}$ par la bijection \'evoqu\'ee plus haut. 
 
  Soit $\bar{{\mathbb  O}}$ une orbite spéciale dans $\mathfrak{g}_{nil}$, conservée par $\Gamma_{F}$. Notons $\boldsymbol{{\cal O}}$ l'orbite dans $\boldsymbol{\mathfrak{g}}_{nil}$ qui lui correspond.     Notons $\bar{{\cal C}}^{Irr}(\boldsymbol{{\cal O}})$ le sous-ensemble des \'el\'ements de $\bar{{\cal C}}^{Irr}$  de la forme $(N,d,\mu)$ avec $N\in \boldsymbol{{\cal O}}$. D'apr\`es le r\'esultat ci-dessus, la dimension de $SI(\mathfrak{g}(F))^*_{\bar{{\mathbb O}}}$ est \'egale au nombre d'\'el\'ements de  $Im(\boldsymbol{\nabla}_{F})\cap \bar{{\cal C}}^{Irr}(\boldsymbol{{\cal O}})$, o\`u $Im(\boldsymbol{\nabla}_{F})$ est l'image de l'application $\boldsymbol{\nabla}_{F}$.    Pour $e\in \mathfrak{g}(F)\cap \bar{{\mathbb O}}$, on a dit ci-dessus que $H^1(\Gamma_{F},\bar{A}(e))$ a m\^eme nombre d'\'el\'ements que $\bar{{\cal C}}^{\sharp}_{F}(\boldsymbol{{\cal O}})$. L'assertion (ii) du th\'eor\`eme dit qu'hormis les cas particuliers indiqu\'es, les ensembles $Im(\boldsymbol{\nabla}_{F})\cap \bar{{\cal C}}^{Irr}(\boldsymbol{{\cal O}})$ et $\bar{{\cal C}}^{\sharp}_{F}(\boldsymbol{{\cal O}})$ ont m\^eme nombre d'\'el\'ements. Si $G$ est classique ou si $q-1$ est divisible par $3\times 4\times 5$, l'application $\boldsymbol{\nabla}_{F}$ est surjective et les deux ensembles sont d'une nature combinatoire simple. Le calcul est alors facile.  Si $G$ est exceptionnel et $q-1$ n'est pas divisible par $3\times 4\times 5$, $\boldsymbol{\nabla}_{F}$ n'est plus surjective. Cela vient du fait qu'il y a des couples $(M,s_{M})$ avec $M\in {\cal L}^{st}_{F,sd}$ et $s_{M}\in \underline{S}^{st}(M_{AD})$ pour lesquels certains \'el\'ements de ${\bf FC}^{st}_{{\mathbb F}_{q}}(\mathfrak{m}_{SC,s_{M}})$ "disparaissent": les faisceaux-caract\`eres existent bien mais ne sont plus invariants par l'action de Frobenius. Du c\^ot\'e de $\bar{{\cal C}}^{\sharp}_{F}(\boldsymbol{{\cal O}})$, la situation est aussi diff\'erente du cas o\`u $q-1$ est assez divisible \`a cause de l'\'egalit\'e $vdv^{-1}=d^q$ qui intervient dans la d\'efinition de cet ensemble: cette \'egalit\'e \'equivaut \`a $vdv^{-1}=d$ si $q-1$ est assez divisible mais ne l'est pas en g\'en\'eral. On montre qu'hormis les cas particuliers indiqu\'es dans l'assertion (ii) du th\'eor\`eme, ces modifications se compensent et que les ensembles $Im(\boldsymbol{\nabla}_{F})\cap  \bar{{\cal C}}^{Irr}(\boldsymbol{{\cal O}})$ et $\bar{{\cal C}}^{\sharp}_{F}(\boldsymbol{{\cal O}})$ ont encore m\^eme nombre d'\'el\'ements.   Ce n'est toutefois pas le cas si $\boldsymbol{{\cal O}}$ est exceptionnelle et que $q-1$ n'est pas divisible par $4$.  Ce point est un peu myst\'erieux mais on peut dire pour se rassurer qu'il est habituel dans la th\'eorie des groupes finis que ces orbites exceptionnelles cr\'eent des perturbations. 
  \bigskip
  
  Comme on le voit, cet article doit beaucoup aux travaux de Lusztig. Je remercie celui-ci pour les explications qu'il m'a fournies.

 \section{L'espace ${\cal D}^{st}(\mathfrak{g}(F))$}

\subsection{Immeubles}

 \subsubsection{Notations}\label{notations}
 
 Nous adoptons le principe suivant. Quand nous aurons défini un objet $X$ dépendant d'un objet $G$, par exemple ci-dessous $X=\mathfrak{g}_{nil}$ ou $X=W$, nous utiliserons la notation analogue ou enrichie de la lettre $H$ pour désigner l'objet analogue relatif à un objet $H$, par exemple $\mathfrak{h}_{nil}$ ou $W^H$. 
 
 Soit $F$ un corps local non-archim\'edien de caract\'eristique nulle. On note  ${\mathbb F}_{q}$ le corps r\'esiduel de $F$, $q$ \'etant son nombre d'\'el\'ements, $p$ la caract\'eristique  de ${\mathbb F}_{q}$, $\mathfrak{o}_{F}$ l'anneau d'entiers de $F$ et $\mathfrak{p}_{F}$ son id\'eal maximal.    On fixe des cl\^otures alg\'ebriques $\bar{F}$ de $F$ et $\bar{{\mathbb F}}_{q}$ de ${\mathbb F}_{q}$.  Toutes les extensions alg\'ebriques de $F$, resp. ${\mathbb F}_{q}$,  que l'on consid\'erera seront suppos\'ees incluses dans $\bar{F}$, resp. $\bar{{\mathbb F}}_{q}$.   On note $val_{F}$ la valuation usuelle de $F$ et on la prolonge \`a $\bar{F}$ en une valuation \`a valeurs dans ${\mathbb Q}\cup\{\infty\}$. On note $\vert .\vert _{F}$ la valeur absolue usuelle de $F$, que l'on prolonge aussi \`a $\bar{F}$. 
 
 On note ${\mathbb N}_{>0}$ l'ensemble des entiers strictement positifs. 
 
   Consid\'erons un groupe $G$ agissant sur un ensemble $U$. Pour $u\in U$, on note $Z_{G}(u)$ le fixateur de $u$ dans $G$. Pour un sous-ensemble $V\subset U$, on note $Norm_{G}(V)$ le stabilisateur de $V$ dans $G$. 
   
    Pour tout groupe ab\'elien $R$, on note $R^{\vee}=Hom(R,{\mathbb C}^{\times})$ son groupe des caract\`eres.        Soit $k$ un corps.    Pour tout espace vectoriel $V$ sur  $k$, on note $V^*$ son dual. 
      Pour tout ${\mathbb Z}$-module $R$, on note $R_{k}=R\otimes_{{\mathbb Z}}k$.

   Soit $k$ une extension alg\'ebrique de $F$ ou ${\mathbb F}_{q}$.  On note $\Gamma_{k}$ le groupe de Galois de $\bar{k}/k$. Pour toute extension galoisienne $k'$ de $k$, on note $\Gamma_{k'/k}$ le groupe de Galois de $k'/k$.  On note $F^{nr}$ l'extension non ramifi\'ee maximale de $F$. Son corps r\'esiduel s'identifie \`a $\bar{{\mathbb F}}_{q}$. On note $I_{F}$ le sous-groupe d'inertie de $\Gamma_{F}$, c'est-\`a-dire le groupe $\Gamma_{F^{nr}}$. On pose simplement $\Gamma_{F}^{nr}=\Gamma_{F^{nr}/F}\simeq \Gamma_{{\mathbb F}_{q}}$. 
 On note $Fr$ l'\'el\'ement de Frobenius de  l'un ou l'autre de ces deux  groupes.

  Pour tout groupe alg\'ebrique $H$ d\'efini sur $k$, on note $H^0$ sa composante neutre et $Z(H)$ son centre. Pour tout tore $T$ d\'efini sur $k$, on note $X_{*}(T)$, resp. $X^*(T)$, son groupe de cocaract\`eres alg\'ebriques, resp. caract\`eres alg\'ebriques, d\'efinis sur $\bar{k}$.

 Soit $G$ un groupe r\'eductif connexe d\'efini sur $k$. On note  $\mathfrak{g}$ l'alg\`ebre de Lie de $G$,  $G_{AD}$ le groupe adjoint et $G_{SC}$ le rev\^etement simplement connexe de $G_{AD}$. On note $\underline{\pi}:G\to G_{AD}$ la projection naturelle. Pour un sous-groupe alg\'ebrique $H$ de $G$, on note $H_{ad}=\underline{\pi}(H)$ son image dans $G_{AD}$ et $H_{sc}$ l'image r\'eciproque de $H_{ad}$ dans $G_{SC}$.
 On identifie $G$ \`a $G(\bar{k})$. Ce groupe et beaucoup d'objets qui lui sont reli\'es sont alors munis d'une action de $\Gamma_{k}$ que l'on note $\sigma\mapsto \sigma_{G}$. S'il n'y a pas d'ambigu\"{\i}t\'e, on note simplement $\sigma=\sigma_{G}$. Pour $g\in G$, on note $Ad(g)$ la conjugaison par $g$ dans $G$ ou l'action adjointe de $g$ dans $\mathfrak{g}$. On notera aussi ces actions $x\mapsto gxg^{-1}$ ou $X\mapsto gXg^{-1}$.   On note $\mathfrak{g}_{reg}$ l'ensemble des \'el\'ements semi-simples r\'eguliers de $\mathfrak{g}$ et $\mathfrak{g}_{nil}$ l'ensemble des \'el\'ements nilpotents.   Pour  $x\in G$, resp. $X\in \mathfrak{g}$, on note  $G_{x}$, resp. $G_{X}$,  la composante neutre de $Z_{G}(x)$, resp. $Z_{G}(X)$. Pour tout sous-groupe parabolique $P$ de $G$, on note $U_{P}$ son radical unipotent. 
   Un Levi de $G$ est un sous-groupe de $G$ qui est la composante de Levi d'un sous-groupe parabolique de $G$. La notion de Levi "sur $k$" est ambigu\"e, pr\'ecisons-la. Nous dirons qu'un sous-groupe $M\subset G$ est un Levi d\'efini sur $k$ si $M$ est un Levi et qu'il est d\'efini sur $k$. Nous dirons que c'est un $k$-Levi  si $M$ est d\'efini sur $k$ et qu'il existe un sous-groupe parabolique $P$ de $G$ d\'efini sur $k$ tel que $M$ soit une composante de Levi de $P$. On appelle paire parabolique un couple $(P,M)$ tel que $P$ soit un sous-groupe parabolique de $G$ et $M$ une composante de Levi de $P$. Si $P$ est un sous-groupe de Borel, on parlera d'une paire de Borel. On dit que la paire est d\'efinie sur $k$ si $P$ et $M$ le sont. Fixons une paire parabolique $(P,M)$ d\'efinie sur $k$ et minimale, c'est-\`a-dire minimale parmi les paires paraboliques d\'efinies sur $k$. On appelle alors
 sous-groupe parabolique standard un sous-groupe parabolique $Q$ tel que $P\subset Q$ et Levi standard une composante de Levi $L$ contenant $M$ d'un sous-groupe parabolique standard.      Si $L$ est un Levi standard, alors $L$ est d\'efini sur $k$ si et seulement si c'est un $k$-Levi (le groupe parabolique standard $Q$ associ\'e \`a $L$ est unique et est donc conserv\'e par $\Gamma_{k}$ si $L$ l'est).  Pour un Levi $M$ de $G$, on pose $W(M)=Norm_{G}(M)/M$.

   Fixons une paire de Borel $(B,T)$ de $G$. On note $\Sigma$ l'ensemble des racines de $T$ dans $\mathfrak{g}$ et $\Delta$ le sous-ensemble de racines simples. A toute racine $\alpha$ est associ\'ee une sous-alg\`ebre radicielle $\mathfrak{u}_{\alpha}$. Pour toute $\alpha\in \Delta$, fixons un \'el\'ement non nul $E_{\alpha}\in \mathfrak{u}_{\alpha}$. Le triplet $(B,T,(E_{\alpha})_{\alpha\in \Delta})$ est appel\'e un \'epinglage de $G$. On dit qu'il est d\'efini sur $k$ si $(B,T)$ l'est et que l'action de $\Gamma_{k}$ permute les \'el\'ements $E_{\alpha}$ pour $\alpha\in \Delta$. A tout Levi $M$ de $G$ contenant $T$ est associ\'e un sous-ensemble $\Sigma^M\subset \Sigma$. Si $M$ est standard, il lui est associ\'e un sous-ensemble $\Delta^M\subset \Delta$ qui est une base de $\Sigma^M$.  On utilise pour le syst\`eme de racines $\Sigma$ les notations usuelles: \`a toute racine $\alpha\in \Sigma$ est associ\'ee une coracine $\check{\alpha}\in X_{*}(T_{sc})$, que l'on peut identifier \`a un \'el\'ement de $\mathfrak{t}_{sc}=X_{*}(T_{sc})\otimes_{{\mathbb Z}}\bar{k}$; on note $\check{\Delta}=\{\check{\alpha},\alpha\in \Delta\}$;   \`a toute racine simple $\alpha\in \Delta$ est associ\'e un poids $\varpi_{\alpha}\in X^*(T_{sc})$ et un copoids $\check{\varpi}_{\alpha}\in X_{*}(T_{ad})$.  On note $W(T) $ ou simplement $W$ le groupe de Weyl relatif \`a $T$, c'est-\`a-dire $W=Norm_{G}(T)/T$. Le groupe $W$ agit naturellement sur l'ensemble des Levi de $G$ contenant $T$. Pour un tel Levi $M$,  on note 
   $Norm_{W}(M)$ le groupe des $w\in W$ tel que $w(M)=M$. On a alors 
   $W(M)=Norm_{W}(M)/W^M$.       
    Supposons que $G_{SC}$ soit absolument quasi-simple. On utilise alors pour le syst\`eme de racines de $G$ les notations de Bourbaki \cite{B} p. 250-275, en particulier la num\'erotation des \'el\'ements de $\Delta$. 
       
   On note $A_{G}$ le plus grand sous-tore de $Z(G)$ qui soit d\'eploy\'e sur $k$ et on pose  ${\cal A}_{G}= X_{*,{\mathbb R}}(A_{G})=X_{*}(A_{G})\otimes_{{\mathbb Z}}{\mathbb R}$. On note $\mathfrak{g}_{ell}(k)$ le sous-ensemble des \'el\'ements elliptiques de $\mathfrak{g}_{reg}(k)$, c'est-\`a-dire des $X\in \mathfrak{g}_{reg}(k)$ tels que $A_{G_{X}}=A_{G}$.  
   Si $k\subset F^{nr}$, on note $A_{G}^{nr}$ le plus grand sous-tore de $Z(G)$ qui soit d\'eploy\'e sur $F^{nr}$ et on pose  ${\cal A}_{G}^{nr}=X_{*}(A^{nr}_{G})\otimes_{{\mathbb Z}}{\mathbb R}$. 
   
   On note $Aut(G)$ le groupe d'automorphismes alg\'ebriques de $G$. Il contient le groupe d'automorphismes int\'erieurs, lequel s'identifie \`a $G_{AD}$. On pose $Out(G)=Aut(G)/G_{AD}$.

\subsubsection{Les hypoth\`eses sur $p$\label{leshypothesessurp}}
Dans chaque paragraphe de l'article, il est donn\'e un groupe r\'eductif connexe $G$ sur une extension alg\'ebrique $k$ de $F$ ou de ${\mathbb F}_{q}$. Notons $rg(G)$ le rang de $G$ sur $\bar{k}$ et $h(G)$ le plus grand nombre de Coxeter des composantes irr\'eductibles de $G$ sur $\bar{k}$. On  impose toujours l'hypoth\`ese

$(Hyp)_{1}(p):$ $p\geq sup(3h(G),rg(G)+2,5)$. 

Dans les paragraphes \ref{lespaceFC} \`a \ref{descriptionfinale} et dans la section 6, le groupe $G$ est d\'efini sur une extension de $F$ et on impose  l'hypoth\`ese plus forte

$(Hyp)_{2}(p):$ $p\geq5$ et il existe un entier $N\geq1$ tel qu'il existe un plongement d\'efini sur $F$ de $G$ dans $GL(N)$ et que $p\geq (2+val_{F}(p))N$. 

Signalons quelques cons\'equences de ces hypoth\`eses. Supposons que $G$ soit d\'efini sur $\bar{{\mathbb F}}_{q}$.  L'hypoth\`ese $(Hyp)_{1}(p)$  implique que  l'exponentielle est d\'efinie sur les \'el\'ements nilpotents de $\mathfrak{g}$, cf. \cite{C} proposition 5.5.5. Pr\'ecis\'ement, soit $X\in \mathfrak{g}$ un \'el\'ement nilpotent. Notons $ad(X)$ l'application $Y\mapsto [X,Y]$ de $\mathfrak{g}$ dans lui-m\^eme. Alors on a $ad(X)^p=0$ et il existe un unique homomorphisme alg\'ebrique 
$$\begin{array}{ccc}\bar{{\mathbb F}}_{q}&\to&G_{SC}\\t&\mapsto& exp(tX)\\ \end{array}$$
tel que, pour tout $t\in \bar{{\mathbb F}}_{q}$, on ait l'\'egalit\'e
$$Ad(exp(tX))=\sum_{i=0,...,p-1}\frac{t^{i}ad(X)^{i}}{i!}.$$

Supposons que $G$ soit d\'efini sur $F$. L'hypoth\`ese $(Hyp)_{1}(p)$ implique que $G$ est d\'eploy\'e sur une extension galoisienne finie mod\'er\'ement ramifi\'ee de $F$. L'hypoth\`ese $(Hyp)_{2}(p)$ implique que l'exponentielle  est bien définie sur l'ensemble des \'el\'ements topologiquement nilpotents de $\mathfrak{g}(F)$ et est un isomorphisme de cet ensemble  sur  celui des \'el\'ements topologiquement unipotents de $G(F)$, cf. \cite{DR} lemme B.0.3.

 \subsubsection{Forme int\'erieure quasi-d\'eploy\'ee\label{formeinterieure}}
 Pour la suite de la section, sauf mention explicite du contraire,  $G$ est un groupe r\'eductif connexe d\'efini sur $F$. 
 
Introduisons la forme int\'erieure quasi-d\'eploy\'ee $G^*$ de $G$, fixons un épinglage $\mathfrak{E}^*=(B^*,T^*,(E_{\alpha})_{\alpha\in \Delta})$ de $G^*$ défini sur $F$  et fixons un torseur int\'erieur $\psi:G\to G^*$ sur $\bar{F}$. On introduit le cocycle $u_{G}:\Gamma_{F}\to G_{AD}$ de sorte que $\psi(\sigma_{G}(g))=u_{G}(\sigma)\sigma_{G^*}(\psi(g))u_{G}(\sigma)^{-1}$  pour tous $\sigma\in \Gamma_{F}$ et $g\in G$. Fixons une paire parabolique $(P_{min},M_{min})$ de $G$ d\'efinie sur $F$ et minimale. Quitte \`a composer $\psi$ avec un automorphisme int\'erieur de $G^*$, on suppose que $(\psi(P_{min}),\psi(M_{min}))$ est standard. On note cette paire $(P^*_{min},M^*_{min})$. On sait que cette derni\`ere paire est uniquement d\'etermin\'ee, qu'elle est d\'efinie sur $F$ et que $u_{G}$ (ou plut\^ot sa classe) provient d'un \'el\'ement de $H^1(\Gamma_{F},M^*_{min,ad})$. Comme on le sait, l'inflation $H^1(\Gamma_{F}^{nr},M^*_{min,ad})\to H^1(\Gamma_{F},M^*_{min,ad})$ est bijective. On peut donc supposer   que $u_{G}$ se factorise en une application de $\Gamma_{F}^{nr}$ dans $M^*_{min,ad}$. 
  Alors $\psi$ est un isomorphisme d\'efini sur $F^{nr}$. Le groupe $M^*_{min}$ poss\`ede la propri\'et\'e suivante:
  
  (1) soit $w\in W^{\Gamma_{F}}$, supposons que $w(M^*_{min})$ soit un Levi standard; alors $w(M^*_{min})=M^*_{min}$.
  
  C'est une cons\'equence du fait que toutes les paires paraboliques de $G$ d\'efinies sur $F$ et minimales sont conjugu\'ees.

On dit qu'un $F$-Levi $M^*$ de $G^*$ se transf\`ere \`a $G$ si et seulement si $\psi^{-1}(M^*)$ est conjugu\'e \`a un $F$-Levi de $G$. Montrons que

(2) un $F$-Levi standard $M^*$ de $G^*$ se transf\`ere \`a $G$ si et seulement si $M^*$ contient $M^*_{min}$. 

Notons $P^*$ le sous-groupe parabolique standard de composante de Levi $M^*$. Supposons $M^*_{min}\subset M^*$. Parce que $u_{G}$ prend ses valeurs dans $M^*_{min}$, on v\'erifie que la paire $(\psi^{-1}(P^*),\psi^{-1}(M^*))$ de $G$ est d\'efinie sur $F$, donc $\psi^{-1}(M^*)$ est un $F$-Levi de $G$. Inversement, supposons que $\psi^{-1}(M^*)$ soit conjugu\'e \`a un $F$-Levi de $G$. Alors il est conjugu\'e \`a un $F$-Levi standard $L$. Posons $L^*=\psi(L)$. C'est un Levi standard de $G^*$ qui contient $M^*_{min}$. Pour la m\^eme raison que ci-dessus, puisque $L$ est un $F$-Levi, $L^*$ est aussi  un $F$-Levi. Puisque  $\psi^{-1}(M^*)$ est conjugu\'e \`a $L$, $M^*$ est conjugu\'e \`a $L^*$. Les deux groupes $M^*$ et $L^*$ sont des $F$-Levi conjugu\'es. D'apr\`es \cite{S} th\'eor\`eme A, ils sont conjugu\'es par un \'el\'ement de $G(F)$. Soit $g\in G(F)$ tel que $gL^*g^{-1}=M^*$. Le groupe $gM^*_{min}g^{-1}$ est un  $F$-Levi  de $M^*$.
  Quitte \`a multiplier $g$ \`a gauche par un  \'el\'ement de $M^*(F)$, on peut  supposer que $gM^*_{min}g^{-1}$ est un $F$-Levi standard. Le groupe $gT^*g^{-1}$ est un $F$-Levi minimal de $gM_{min}^*g^{-1}$. Quitte \`a multiplier $g$ \`a gauche par un \'el\'ement de $gM^*_{min}(F)g^{-1}$, on  peut  supposer que $ gT^*g^{-1}=T^*$.  Donc $g$ d\'efinit un \'el\'ement $w\in W^{\Gamma_{F}}$ tel que $M^*=w(L^*)$. D'apr\`es (1), on a $w(M^*_{min})=M^*_{min}$. Puisque $M^*_{min}\subset L^*$,  cela entra\^{\i}ne que $M^*_{min}\subset M^*$. Cela prouve (2).

    \subsubsection{Immeubles de Bruhat-Tits \label{immeubles}} 
   
 On note $Imm(G_{AD})$ l'immeuble de Bruhat-Tits du groupe adjoint $G_{AD}$ sur $F$. Le groupe $G(F)$, et plus g\'en\'eralement le groupe $G_{AD}(F)$, agissent sur cet immeuble. D'autre part, celui-ci est d\'ecompos\'e en facettes. On note $Fac(G)$ l'ensemble des facettes et  $S(G)\subset Fac(G)$ l'ensemble des sommets.   Pour tout ${\cal F}\in Fac(G)$, on note $\bar{{\cal F}}$ son adh\'erence et $S(\bar{{\cal F}})$ l'ensemble des sommets appartenant \`a $\bar{{\cal F}}$. 
 On note $K_{{\cal F}}^{\dag}$ le stabilisateur de ${\cal F}$ dans $G(F)$,  $K_{{\cal F}}^0\subset K_{{\cal F}}^{\dag}$ le sous-groupe parahorique, $K_{{\cal F}}^+$ son plus grand sous-groupe distingu\'e pro-$p$-unipotent et $G_{{\cal F}}$ le groupe r\'eductif connexe sur ${\mathbb F}_{q}$ d\'efini par Bruhat et Tits tel que $K_{{\cal F}}^0/K_{{\cal F}}^+\simeq G_{{\cal F}}({\mathbb F}_{q})$. Il y a des objets correspondants dans l'alg\`ebre de Lie: $\mathfrak{k}_{{\cal F}}^+\subset \mathfrak{k}_{{\cal F}}\subset \mathfrak{g}(F)$, avec $\mathfrak{k}_{{\cal F}}/\mathfrak{k}_{{\cal F}}^+\simeq \mathfrak{g}_{{\cal F}}({\mathbb F}_{q})$. Pour $x\in Imm(G_{AD})$, on pose $K_{x}^0=K_{{\cal F}}^0$, $\mathfrak{k}_{x}=\mathfrak{k}_{{\cal F}}$ etc... o\`u ${\cal F}$ est la facette \`a laquelle appartient $x$. 
 On ajoute des indices $AD$, resp. $SC$, pour les m\^emes objets relatifs au groupe $G_{AD}$, resp. $G_{SC}$: $K_{{\cal F},AD}^{\dag}$ etc...  
 
 On consid\'erera aussi l'immeuble de $G_{AD}$ sur des extensions alg\'ebriques $K$ de $F$. On ajoutera alors la lettre $K$ dans la notation: $Imm_{K}(G_{AD})$, $S_{K}(G)$ etc...  
 Si $K/F$ est galoisienne, le groupe $\Gamma_{K/F}$ agit sur $Imm_{K}(G_{AD})$  et, si $K/F$ est mod\'er\'ement ramifi\'ee, $Imm(G_{AD})$ s'identifie \`a l'ensemble des points fixes par $\Gamma_{K/F}$ dans $Imm_{K}(G_{AD})$. Supposons $K/F$ non ramifiée. Une facette ${\cal F}$ de $Imm(G_{AD})$ ne reste pas une facette de $Imm_{K}(G_{AD})$ mais est incluse dans une unique facette ${\cal F}_{K}$ de cet immeuble et on pose simplement $K_{{\cal F},K}^0=K_{{\cal F}_{K},K}^0$, etc... Si $K=F^{nr}$, on simplifiera certaines notations en rempla\c{c}ant $F^{nr}$ par un exposant $^{nr}$, par exemple on notera  ${\cal F}^{nr}={\cal F}_{F^{nr}}$ et $K_{{\cal F}}^{0,nr}=K_{{\cal F}^{nr},F^{nr}}^0$. 
 
 Notons ${\cal T}_{max}$ l'ensemble des sous-tores de $G$  d\'efinis et d\'eploy\'es sur $F$ maximaux (c'est-\`a-dire maximaux parmi les sous-tores  d\'efinis et d\'eploy\'es sur $F$). Notons de m\^eme ${\cal T}_{max}^{nr}$ l'ensemble des sous-tores de $G$  d\'efinis et d\'eploy\'es sur $F^{nr}$ maximaux. A $S\in {\cal T}_{max}$ est associ\'e un appartement $App(S)\subset Imm(G_{AD})$. C'est un espace affine sous le groupe ${\cal A}_{S}/{\cal A}_{G}$.  De m\^eme, soit $F'$ une extension galoisienne de $F$ mod\'er\'ement ramifi\'ee et soit $T\in {\cal T}_{F',max}$. A $T$ 
  est associ\'e un appartement $App_{F'}(T)\subset Imm_{F'}(G_{AD})$. Supposons que $S\subset T$ et que $T$ soit d\'efini sur $F$. Alors $App(S)$ est l'ensemble des points fixes par $\Gamma_{F'/F}$ dans $App_{F'}(T)$, cf. \cite{T} 2.6.1. Rappelons que, pour tout $S\in {\cal T}_{max}$, il existe un  tore $T\in {\cal T}^{nr}_{max}$ tel que  $S\subset T$ et que $T$ soit d\'efini sur $F$, cf. \cite{T} 1.10.
 
 Pour tout $x\in Imm(G_{AD})$ et tout r\'eel $r\in {\mathbb R}$, resp. tel que $r\geq0$, Moy et Prasad ont d\'efini un $\mathfrak{o}_{F}$-r\'eseau $\mathfrak{k}_{x,r}\subset \mathfrak{g}(F)$, resp. un sous-groupe $K_{x,r}\subset G(F)$.    L'application $r\mapsto \mathfrak{k}_{x,r}$, resp. $r\mapsto K_{x,r}$,  est d\'ecroissante et on pose $\mathfrak{k}_{x,r+}=\cup_{r'>r}\mathfrak{k}_{x,r'}$, resp. $K_{x,r+}=\cup_{r'>r}K_{x,r'}$. En particulier $\mathfrak{k}_{x,0}=\mathfrak{k}_{x}$ et $\mathfrak{k}_{x,0+}=\mathfrak{k}_{x}^+$, resp. $K_{x,0}=K_{x}^0$, $K_{x,0+}=K_{x}^+$. Pour $X\in \mathfrak{g}(F)$, on dit que $X$ est entier, resp. topologiquement nilpotent, si et seulement s'il existe $x\in Imm(G_{AD})$ tel que $X\in \mathfrak{k}_{x,0}$, resp. $X\in \mathfrak{k}_{x,0+}$. On note $\mathfrak{g}_{ent}(F)$, resp. $\mathfrak{g}_{tn}(F)$, le sous-ensemble des \'el\'ements entiers, resp. topologiquement nilpotents. 
 
 Soit $x\in Imm_{F^{nr}}(G_{AD})$. En g\'en\'eral, on a $G_{SC,x}\not=G_{x,SC}$ et m\^eme $\mathfrak{g}_{SC,x}\not=\mathfrak{g}_{x,SC}$: $G_{x,SC}$ s'identifie au rev\^etement simplement connexe du groupe d\'eriv\'e de $G_{SC,x}$. On a toutefois
 
 (1) si $x$ est un sommet de $Imm_{F^{nr}}(G_{AD})$, alors $\mathfrak{g}_{SC,x}=\mathfrak{g}_{x,SC}$.
 
 En effet, si $x$ est un sommet de $Imm_{F^{nr}}(G_{AD})$, $Z(G_{SC,x})^0$ est trivial et $G_{x,SC}$ est un rev\^etement de $G_{SC,x}$. 
 
  \subsubsection{Immeubles et $F^{nr}$-Levi \label{immeublesetLevi}}
 
 Soit  $H$ un $F^{nr}$-Levi de $G$. Remarquons que l'ensemble des tores $T\in {\cal T}_{max}^{nr}$ tels que $T\subset H$ est \'egal \`a l'analogue ${\cal T}_{max}^{H,nr}$ pour $H$ de l'ensemble ${\cal T}_{max}^{nr}$.  On note $Imm^G_{F^{nr}}(H_{ad}) $
    la r\'eunion des appartements $App_{F^{nr}}(T)$ dans $Imm_{F^{nr}}(G_{AD})$ sur  cet ensemble de sous-tores.    Il y a une projection naturelle $p_{H}:Imm_{F^{nr}}^G(H_{ad})\to Imm_{F^{nr}}(H_{AD})$. Pour un sous-tore $T\subset H$ appartenant \`a ${\cal T}_{max}^{nr}$, on a
    
    (1) $p_{H}^{-1}(App_{F^{nr}}^H(T))=App_{F^{nr}}(T)$,
    
    \noindent o\`u $App_{F^{nr}}^H(T)$ est l'appartement associ\'e \`a $T$ dans $Imm_{F^{nr}}(H_{AD})$. 
    
  Supposons que $H$ soit d\'efini sur $F$.  L'ensemble  $Imm^G_{F^{nr}}(H_{ad}) $ est stable par l'action de $\Gamma_{F}^{nr}$. On note $Imm^G(H_{ad})$ l'ensemble des points fixes.  Il s'identifie comme ensemble \`a l'immeuble \'etendu  de $H_{ad}$ mais il faut prendre garde que les d\'ecompositions en facettes ne sont pas les m\^emes.  La projection $p_{H}$ est \'equivariante pour les actions de $\Gamma_{F}^{nr}$, on note encore   $p_{H}:Imm^G(H_{ad})\to Imm(H_{AD})$ sa restriction, qui est encore surjective. L'ensemble $Imm^G(H_{ad})$ est une r\'eunion de facettes de $Imm(G_{AD})$ et, pour toute telle facette ${\cal F}\subset Imm^G(H_{ad})$, $p_{H}({\cal F})$ est incluse dans une facette de $Imm(H_{AD})$ que l'on note ${\cal F}_{H}$.
La d\'ecomposition de $Imm^G(H_{ad})$ est pour nous induite par celle de $Imm(G_{AD})$.  
Pour $x\in Imm^G(H_{ad})$, on a

(2) $\mathfrak{k}^H_{p_{H}(x)}=\mathfrak{k}_{x}\cap \mathfrak{h}(F)$; plus g\'en\'eralement,  $\mathfrak{k}^H_{p_{H}(x),r}=\mathfrak{k}_{x,r}\cap \mathfrak{h}(F)$ pour tout $r\in {\mathbb R}$;

(3) $K^{H,0}_{p_{H}(x)}= K^0_{x}\cap H(F)$.

Soit $F'$ une extension galoisienne  finie de $F$, mod\'er\'ement ramifi\'ee. Montrons que

(4) $Imm^G(H_{ad})=(Imm_{F'}^G(H_{ad}))^{\Gamma_{F'/F}}$.

Il suffit de prouver que $Imm_{F^{nr}}^G(H_{ad})=(Imm_{F^{'nr}}^G(H_{ad}))^{\Gamma_{F^{'nr}/F^{nr}}}$. Le membre de gauche est clairement inclus dans celui de droite. Soit $x\in (Imm_{F^{'nr}}^G(H_{ad}))^{\Gamma_{F^{'nr}/F^{nr}}}$. La projection $p_{H}: Imm_{F^{'nr}}^G(H_{ad})\to Imm_{F^{'nr}}(H_{AD})$ est \'equivariante pour les actions galoisiennes. Donc
$p_{H}(x)$ appartient \`a $Imm_{F^{nr}}(H_{AD})$. Fixons un tore $T\in {\cal T}_{max}^{H,nr}$ tel que $p_{H}(x)\in App_{F^{nr}}^H(T)$. Le commutant $Z_{H}(T)$ est un tore. Notons $T'$ le plus grand sous-tore d\'eploy\'e sur $F^{'nr}$ de $Z_{H}(T)$. On a aussi $p_{H}(x)\in App_{F^{'nr}}^H(T')$. En appliquant (1), on a $x\in App_{F^{'nr}}(T')$. Mais, par hypoth\`ese,  $x$ est fix\'e par $\Gamma_{F^{'nr}/F^{nr}}$. Donc $x\in App_{F^{'nr}}(T')^{\Gamma_{F^{'nr}/F^{nr}}}=App_{F^{nr}}(T)$. Cela prouve que $x$ appartient \`a $Imm_{F^{nr}}^G(H_{ad})$, d'o\`u (4). 

  Supposons  que $H$ soit un   $F$-Levi de $G$. Dans ce cas
  
 (5) $Imm^G(H_{ad})$ est la r\'eunion des appartements $App(S)$ de $Imm(G_{AD})$ pour les sous-tores $S$ de $H$ qui appartiennent \`a ${\cal T}_{max}$.
 
 La preuve est similaire \`a celle de (4).  Soit $x\in Imm^G(H_{ad})$. Alors $p_{H}(x)\in Imm(H_{AD})$. Fixons un sous-tore $S$ de $H$ d\'efini et d\'eploy\'e sur $F$, maximal, de sorte que $p_{H}(x)\in App^H(S)$. Fixons un tore $T\in {\cal T}_{max}^{H,nr}$ qui contient $S$ et est d\'efini sur $F$. Alors $p_{H}(x)\in App_{F^{nr}}^H(T)$. En appliquant (1), $x$ appartient \`a $App_{F^{nr}}(T)$. Parce que $H$ est un $F$-Levi, $S$ est encore d\'eploy\'e maximal dans $G$. Alors $App_{F^{nr}}(T)^{\Gamma_{F}^{nr}}
=App(S)$. Le point $x$ appartient \`a cet ensemble, ce qui prouve (5).

Remarquons que l'espace ${\cal A}_{H}/{\cal A}_{G}$ agit naturellement sur $Imm^G(H_{ad})$ et que la projection $p_{H}$ n'est autre que la projection sur l'espace des orbites pour cette action. Dans le cas o\`u $H$ est elliptique, c'est-\`a-dire $A_{H}=A_{G}$, $p_{H}$ est un isomorphisme.  En cons\'equence:

(6) supposons $H$ elliptique, soit $s$ un sommet de $Imm(H_{AD})$; alors $p_{H}^{-1}(s)$ est un sommet de $Imm^G(H_{ad})$.

Cela r\'esulte du fait que la d\'ecomposition en facettes de $Imm^G(H_{ad})$ est plus fine que celle de $Imm(H_{AD})$.

\subsubsection{$F^{nr}$-Levi et groupes en r\'eduction}\label{groupesenreduction}
Soient $T\in {\cal T}_{max}^{nr}$ et $x\in App_{F^{nr}}(T)\subset Imm_{F^{nr}}(G_{AD})$. Du tore $T$ se d\'eduit un sous-tore maximal $T_{x}$ de $G_{x}$. On a $X_{*}(T)\simeq X_{*}(T_{x})$ et il y a une bijection $T'\mapsto T'_{x}$ entre les sous-tores de $T$ et ceux de $T_{x}$:  $X_{*}(T'_{x})$ correspond \`a $X_{*}(T')$ par l'isomorphisme pr\'ec\'edent.  

 {\bf Remarque.} Les m\^emes d\'efinitions et propri\'et\'es valent si l'on consid\`ere une facette ${\cal F}\subset App_{F^{nr}}(T)$. Par exemple, du tore $T$ se d\'eduit un tore maximal $T_{{\cal F}}$ de $G_{{\cal F}}$. 
 
 \bigskip
 
Soit $H$ un $F^{nr}$-Levi de $G$. Consid\'erons un tore $T\in {\cal T}^{H,nr}_{max}$, un point $x\in App_{F^{nr}}(T)\subset Imm_{F^{nr}}^G(H_{ad})$ et sa projection $x_{H}=p_{H}(x)\in App^H_{F^{nr}}(T)\subset Imm_{F^{nr}}(H_{AD})$.  On a $A_{H}^{nr}\subset T$ et de $A_{H}^{nr}$ se d\'eduit un sous-tore $A_{H,x}^{nr}$ de $T_{x}$. Notons $H_{x}$ l'image dans $G_{x}$ de $K_{x}^{0,nr}\cap H(F^{nr})$. Alors $H_{x}$ s'identifie au groupe $H_{x_{H}}$ associ\'e au point $x_{H}$ de $Imm_{F^{nr}}(H_{AD})$. Le tore $A_{H,x}^{nr}$ est contenu dans $Z(H_{x})^0$. Montrons que 

(1) $H_{x}$ est le commutant de $A_{H,x}^{nr}$ dans $G_{x}$; on a l'\'egalit\'e $A_{H,x}^{nr}=Z(H_{x})^0$ si et seulement si $x_{H}$ est un sommet de $Imm_{F^{nr}}(H_{AD})$. 

Notons $\Sigma^{nr}$ l'ensemble des racines de $T$ dans $\mathfrak{g}$ et $\Sigma_{x}$ celui des racines de $T_{x}$ dans $\mathfrak{g}_{x}$. Par l'isomorphisme $X^*(T)\simeq X^*(T_{x})$, $\Sigma_{x}$ s'identifie \`a un sous-ensemble $\Sigma^{nr}(G_{x})\subset \Sigma^{nr}$, notons plus pr\'ecis\'ement $\beta\mapsto \beta^G$ cette bijection. Notons $\Sigma_{x}^{H_{x}}$, resp. $\Sigma^{nr,H}$, l'ensemble des racines de $T_{x}$ dans $H_{x}$,  resp. de $T$ dans $H$. Pour $\beta\in \Sigma_{x}$, on a $\beta\in \Sigma_{x}^{H_{x}}$ si et seulement si $\beta^G\in \Sigma^{nr,H}$, c'est-\`a-dire si et seulement si $\beta^G=1$ sur $A_{H}^{nr}$. Cette condition \'equivaut \`a $\beta=1$ sur $A_{H,x}^{nr}$, c'est-\`a-dire $\beta$ intervient dans le commutant de $A_{H,x}^{nr}$. Cela d\'emontre la premi\`ere assertion. Le point $x_{H}$ est un sommet de $Imm_{F^{nr}}(H_{AD})$ si et seulement si $Z(H_{x_{H}})^0$ et $A_{H}^{nr}$ ont m\^eme dimension. Cela entra\^{\i}ne la seconde assertion. 

Soit $H$ un $F$-Levi de $G$. Consid\'erons un tore $S\in {\cal T}^{H}_{max}$, un point $x\in App(S)\subset Imm^G(H_{ad})$ et sa projection $x_{H}=p_{H}(x)\in App^H(S)\subset Imm_{F^{nr}}(H_{AD})$.   De $A_{H}$ se d\'eduit un sous-tore $A_{H,x}$ de $S_{x}$. On a

(2) $H_{x}$ est le commutant de $A_{H,x}$ dans $G_{x}$; $A_{H,x}$ est le plus grand sous-tore de $Z(H_{x})^0$ d\'eploy\'e  sur ${\mathbb F}_{q}$ si et seulement si $x_{H}$ est un sommet de $Imm(H_{AD})$. 

On fixe un tore $T\in {\cal T}^{H,nr}_{max}$ contenant $S$ et d\'efini sur $F$ et on reprend les notations pr\'ec\'edentes. Parce que $H$ est un $F$-Levi, c'est le commutant de $A_{H}$ dans $G$. Une racine $\beta\in \Sigma^{nr}$ vaut donc $1$ sur $A_{H}$ si et seulement si elle vaut $1$ sur $A_{H}^{nr}$. Par la correspondance $\beta\mapsto \beta^G$, on en d\'eduit qu'une racine $\beta\in \Sigma_{x}$ vaut $1$ sur $A_{H,x}$ si et seulement si elle vaut $1$ sur $A_{H,x}^{nr}$. Les commutants dans $G_{x}$ de $A_{H,x}$ et $A_{H,x}^{nr}$ sont donc \'egaux. En appliquant (1), on en d\'eduit la premi\`ere assertion  de (2). 
Le tore $A_{H,x}$ est en tout cas un sous-tore de  $Z(H_{x})^0$ d\'eploy\'e  sur ${\mathbb F}_{q}$. C'est le plus grand tel sous-tore si et seulement si ce  plus grand sous-tore a la m\^eme dimension que $A_{H,x}$, ou encore que $A_{H}$. Mais cette condition \'equivaut \`a ce que $x_{H}$ soit un sommet de $Imm(H_{AD})$. Cela d\'emontre (2).

Soient $T\in {\cal T}_{max}^{nr}$ et $x\in App_{F^{nr}}(T)\subset Imm_{F^{nr}}(G_{AD})$. Consid\'erons un $\bar{{\mathbb F}}_{q}$-Levi $\bar{H}$ de $G_{x}$ et supposons que $Z(\bar{H})^0\subset T_{x}$. Introduisons le sous-tore $T_{\bar{H}}\subset T$ tel que $T_{\bar{H},x}=Z(\bar{H})^0$. Notons $H$ la composante neutre du commutant de $T_{\bar{H}}$ dans $G$. C'est un $F^{nr}$-Levi contenant $T$ donc $x\in Imm^G_{F^{nr}}(H_{ad})$. Posons $x_{H}=p_{H}(x)$. Alors

(3) $x_{H}$ est un sommet de $Imm_{F^{nr}}(H_{AD})$; on a $A_{H}^{nr}=T_{\bar{H}}$ et $H_{x}=\bar{H}$. 

Le m\^eme calcul que dans la preuve de (1) montre que $H_{x}$ et $\bar{H}$ ont m\^emes ensembles de racines donc sont \'egaux.  Par construction, $T_{\bar{H}}$ est un sous-tore d\'eploy\'e sur $F^{nr}$ de $Z(H)^0$. Donc $T_{\bar{H}}\subset A_{H}^{nr}$ et, par r\'eduction, $Z(\bar{H})^0\subset A_{H,x}^{nr}$. Mais $A_{H,x}^{nr}$ est central dans $H_{x}=\bar{H}$, d'o\`u $A_{H,x}^{nr}\subset Z(\bar{H})^0$. Ces deux tores sont donc \'egaux, d'o\`u $A_{H}^{nr}=T_{\bar{H}}$. Puisque $A_{H,x}^{nr}=Z(\bar{H})^0=Z(H_{x})^0$, 
$x_{H}$ est un sommet de $Imm_{F^{nr}}(H_{AD})$ d'apr\`es (1). 

 {\bf Remarque.} Soient $S\in {\cal T}_{max}$, $T\in {\cal T}_{max}^{nr}$, supposons que $S\subset T$ et que $T$ soit d\'efini sur $F$. Soit $x\in App(S)\subset  Imm(G_{AD})$ et soit $\bar{H}$ soit un ${\mathbb F}_{q}$-Levi de $G_{x}$ tel que  $Z(\bar{H})^0\subset T_{x}$. Notons $\bar{S}$ le plus grand sous-tore de $Z(\bar{H})^0$ d\'eploy\'e sur ${\mathbb F}_{q}$. On a $\bar{S}\subset S_{x}$ et on peut introduire le sous-tore $S_{\bar{H}}\subset S$ tel que $S_{\bar{H},x}=\bar{S}$. On prendra garde qu'en g\'en\'eral, les commutants de $T_{\bar{H}}$ et de $S_{\bar{H}}$ ne sont pas \'egaux.

 \subsubsection{Un lemme sur les \'epinglages\label{epinglages}}
 Pour ce paragraphe, on modifie nos hypoth\`eses en posant $k=\bar{F}$ ou $k=\bar{{\mathbb F}}_{q}$ et en supposant que $G$ est un groupe r\'eductif connexe d\'efini  sur $k$.

 On fixe une paire de Borel \'epingl\'ee $\mathfrak{E}=(B,T,(E_{\alpha})_{\alpha\in \Delta})$ de $G$.   L'ensemble $\Delta$ plus les relations de produit scalaire entre les racines forment le diagramme de Dynkin ${\cal D}$ de $G$.   Si $G$ est absolument quasi-simple,   
  on note $\alpha_{0}$ l'oppos\'ee de la plus grande racine dans $\Sigma$ et on pose $\Delta_{a}=\{\alpha_{0}\}\cup \Delta$.  L'ensemble $\Delta_{a}$ plus les relations de produit scalaire entre les racines forment le diagramme de Dynkin  \'etendu ${\cal D}_{a}$ de $G$. En g\'en\'eral, l'ensemble $\Delta$,  resp. le diagramme ${\cal D}$, se d\'ecompose en  union d'ensembles $\Delta'$, resp. de diagrammes ${\cal D}'$, analogues associ\'es aux composantes quasi-simples de $G_{SC}$.  On pose $\Delta_{a}=\cup \Delta'_{a}$, ${\cal D}_{a}=\cup {\cal D}'_{a}$, la r\'eunion portant sur les composantes quasi-simples de $G_{SC}$. 
 On note $Aut({\cal D})$ et $Aut({\cal D}_{a})$ les groupes d'automorphismes des diagrammes ${\cal D}$ et ${\cal D}_{a}$.

 Rappelons  que  de l'\'epinglage  $\mathfrak{E}$ se   d\'eduit une section  $S:W\to G$. Pour tout $\alpha\in \Sigma$, notons $w_{\alpha}\in W$ la sym\'etrie \'el\'ementaire associ\'ee \`a $\alpha$. Pour $\alpha\in \Delta$, notons $E_{-\alpha}$ l'unique élément de $\mathfrak{u}_{-\alpha}$ tel que $[E_{\alpha},E_{-\alpha}]=\check{\alpha}$. Alors, 
 par d\'efinition,  $S(w_{\alpha})=exp(E_{\alpha})exp(-E_{-\alpha})exp(E_{\alpha})$ pour $\alpha\in \Delta$.   Pour deux \'el\'ements $w_{1},w_{2}\in W$, on a l'\'egalit\'e 
 
 (1) $S(w_{1}w_{2})=t(w_{1},w_{2})S(w_{1})S(w_{2})$, 
 
 \noindent o\`u 
 $$t(w_{1},w_{2})=\prod_{\alpha\in \Sigma, \alpha>0, w_{1}^{-1}(\alpha)<0,w_{2}^{-1}w_{1}^{-1}(\alpha)>0}\check{\alpha}(-1),$$
 cf. \cite{LS} lemme 2.1.A. Bien s\^ur, la notion de positivit\'e est relative au Borel $B$.  Remarquons que l'on a aussi une section de Springer $S^{G_{SC}}:W\to G_{SC}$ et que $S$ est la composée de $S^{G_{SC}}$ et de la projection de $G_{SC}$ dans $G$. 
  On a aussi
 
 (2) soient $\alpha,\beta\in \Delta$ et $u\in W$; supposons $u(\alpha)=\beta$; alors $Ad(S(u))(E_{\alpha})=E_{\beta}$.
 
 C'est bien connu, rappelons bri\`evement l'argument.  L'assertion résulte de la m\^eme assertion pour $S^{G_{SC}}$. On peut donc supposer $G$ simplement connexe.  On a en tout cas $Ad(S(u))(E_{\alpha})=\lambda E_{\beta}$ pour un $\lambda\in \bar{F}^{\times}$. Les d\'efinitions de $S(w_{\alpha} )$ et $ S(w_{\beta})$ entra\^{\i}nent que $S(u)S(w_{\alpha})S(u)^{-1}=\check{\beta}(\lambda)S(w_{\beta})$. Pour d\'emontrer que $\lambda=1$, il suffit de prouver que $S(u)S(w_{\alpha})=S(w_{\beta})S(u)$. A l'aide de (1), on voit que les deux membres sont \'egaux \`a $S(uw_{\alpha})$. Cela prouve (2).
 
 Introduisons le sous-groupe fini $T_{S}$ de $T$ engendré par les éléments $\check{\alpha}(-1)$ pour $\alpha\in \Sigma$. La relation (1) entraîne que l'ensemble $T_{S}S(W)$ est un groupe dont $T_{S}$ est un sous-groupe distingué. Notons $\pi_{S}:T_{S}S(W)\to W$ la projection naturelle. On déduit de (1) et (2) la propriété suivante
 
 (3) pour tout $\alpha\in \Sigma$, il existe un unique sous-ensemble à deux éléments $\{\pm \dot{E}_{\alpha}\}$ de $\mathfrak{u}_{\alpha}$ de sorte que, pour tout $\beta\in \Delta$ et tout $g\in T_{S}S(W)$ tels que $\pi_{S}(g)(\beta)=\alpha$, on a $Ad(g)(E_{\beta})\in \{\pm \dot{E}_{\alpha}\}$.
 
 On vérifie que, un élément $\dot{E}_{\alpha}$ étant fixé, on peut choisir $\dot{E}_{-\alpha}$ de sorte que $[\dot{E}_{\alpha},\dot{E}_{-\alpha}]=\check{\alpha}$. De tels éléments étant fixés, posons $\dot{S}(w_{\alpha})=exp(\dot{E}_{\alpha})exp(-\dot{E}_{-\alpha})exp(\dot{E}_{\alpha})$. Montrons que
 
 (4) $\dot{S}(w_{\alpha})\in T_{S}S(w_{\alpha})$. 
 
  Remplacer $\dot{E}_{\alpha}$ par $-\dot{E}_{\alpha}$ multiplie $\dot{S}(w_{\alpha})$ par $\check{\alpha}(-1)$. On peut donc supposer qu'il existe $\beta\in \Delta$ et $g\in T_{S}S(W)$ de sorte que $\dot{E}_{\alpha}=Ad(g)(E_{\beta})$. On a alors aussi $\dot{E}_{-\alpha}=Ad(g)(E_{-\beta})$, puis $\dot{S}(w_{\alpha})=gS(w_{\beta})g^{-1}$. Donc $\dot{S}(w_{\alpha})\in T_{S}S(W)$. Alors (4) résulte du fait que  $\dot{S}(w_{\alpha})$ et $S(w_{\alpha})$ ont m\^eme image dans $W$. 
  
 Pour les propriétés ci-dessus, l'hypothèse $(Hyp)_{1}(p)$ n'est pas nécessaire mais nous l'imposons désormais. Pour $\alpha\in \Sigma$, notons   $\mathfrak{u}_{ \alpha, {\mathbb Z}_{p}}$ le  sous-${\mathbb Z}_{p}$-module de $\mathfrak{u}_{ \alpha}$ engendré par $\dot{E}_{ \alpha}$ (le choix de cet élément ou de son opposé ne change rien).  La théorie des bases de Chevalley  implique la propriété suivante:
 
 (5) pour $\alpha,\alpha'\in \Sigma$ tels que $\alpha+\alpha'\in \Sigma$, $[\mathfrak{u}_{\alpha,{\mathbb Z}_{p}},\mathfrak{u}_{\alpha',{\mathbb Z}_{p}}]=\mathfrak{u}_{\alpha+\alpha',{\mathbb Z}_{p}}$. 
 
 Pour tout sous-${\mathbb Z}_{p}$-module $R\subset \bar{F}$, on note $\mathfrak{u}_{\alpha,R}$ le sous-$R$-module de $\mathfrak{u}_{\alpha}$ engendré par $\mathfrak{u}_{\alpha,{\mathbb Z}_{p}}$.

 {\bf On suppose pour la fin du paragraphe que} $G$ {\bf est quasi-simple.}

 On a introduit ci-dessus  la racine  $\alpha_{0}$ oppos\'ee de la plus grande racine dans $\Sigma$.  Il existe une unique relation lin\'eaire
 
 (6) $\sum_{\alpha\in \Delta_{a}}d(\alpha)\alpha=0$
 
 \noindent avec $d(\alpha_{0})=1$. On a $d(\alpha)\in {\mathbb N}_{>0}$ pour tout $\alpha\in \Delta$. Fixons un \'el\'ement non nul $E_{\alpha_{0}}$ de $\mathfrak{u}_{\alpha_{0}}$. 
 Appelons \'epinglage affine la famille $\mathfrak{E}_{a}=(T,(E_{\alpha})_{\alpha\in \Delta_{a}})$.      On note $\Omega$ le sous-groupe des \'el\'ements de $W$ dont l'action sur $\Sigma$ conserve $\Delta_{a}$. Il s'identifie \`a un sous-groupe de $Aut({\cal D}_{a})$. D'apr\`es \cite{B} VI.2.3, on a
 
 (7) $\Omega$ est ab\'elien; l'application $\omega\mapsto \omega(\alpha_{0})$ est une bijection de $\Omega$ sur le sous-ensemble des $\alpha\in \Delta_{a}$ telles que $d(\alpha)=1$; on a l'\'egalit\'e  
   $Aut({\cal D}_{a})=\Omega\rtimes Aut({\cal D})$. 
   
   Montrons que
   
   (8) pour $\omega\in \Omega$ et $\alpha\in \Delta_{a}$, on a $d(\omega(\alpha))=d(\alpha)$.
   
   En appliquant $\omega^{-1}$ \`a la relation (6), on obtient
   $$\sum_{\alpha\in \Delta_{a}}d(\omega(\alpha))\alpha=0.$$
   Cette relation est proportionnelle \`a (6). Puisque, dans  les deux relations, le plus petit coefficient vaut $1$,   le coefficient de proportionnalit\'e ne peut qu'\^etre \'egal \`a $1$. D'o\`u (8).

   Notons  $\boldsymbol{\Omega}$ le sous-groupe des \'el\'ements de $ G_{AD}$ qui conservent $\mathfrak{E}_{a}$ (en particulier qui normalisent $T$).

 \begin{prop}{L'application naturelle $\boldsymbol{\Omega}\to \Omega$ est bijective.} \end{prop}
 
 Preuve.   
L'application de l'\'enonc\'e  est \'evidemment injective (un \'el\'ement de $G_{AD}$ qui  conserve $T$ et fixe $E_{\alpha}$ pour tout $\alpha\in \Delta_{a}$, et m\^eme seulement $\alpha\in \Delta$, est l'\'el\'ement neutre puisque $G$ est adjoint).     D'apr\`es (7),  il nous suffit de fixer $\beta\in \Delta$ tel que $d(\beta)=1$ et de construire un \'el\'ement $\boldsymbol{\omega}\in \boldsymbol{\Omega}$ tel que $\boldsymbol{\omega}(\beta)=\alpha_{0}$. En tout cas, l'\'el\'ement $\omega\in \Omega$ qui envoie $\beta$ sur $\alpha_{0}$ se construit de la fa\c{c}on suivante, cf. \cite{B} VI.2.3 proposition 6. Notons $\Sigma^{\beta}$ le sous-syst\`eme de racines de $\Sigma$ engendr\'e par $\Delta-\{\beta\}$ et $W^{\beta}$ son groupe de Weyl. Notons $w_{0}$ l'\'el\'ement de plus grande longueur de $W$ et $w^{\beta}$ celui de $W^{\beta}$. Alors $\omega=w_{0}w^{\beta}$. 
  L'image de $E_{\beta}$ par $Ad(S(\omega))$ est de la forme $\lambda^{-1}E_{0}$ pour un $\lambda\in k^{\times}$. Posons $\boldsymbol{\omega}=S(\omega)\check{\varpi}_{\beta}(\lambda)$. Alors $\boldsymbol{\omega}$ est un \'el\'ement de $Norm_{G}(T)$ qui envoie $E_{\beta}$ sur $E_{\alpha_{0}}$. Soit $\alpha\in \Delta-\{\beta\}$.   On a $Ad(\check{\varpi}_{\beta}(\lambda))(E_{\alpha})=E_{\alpha}$ donc $Ad(\boldsymbol{\omega})(E_{\alpha})=Ad(S(\omega))(E_{\alpha})$. Mais $\omega(\alpha)\in \Delta$ et, d'apr\`es (2), on a  $Ad(S(\omega))(E_{\alpha})=E_{\omega(\alpha)}$.
 Posons $\gamma=\omega(\alpha_{0})$. Il reste seulement \`a prouver que $Ad(\boldsymbol{\omega})(E_{\alpha_{0}})=E_{\gamma}$. On a $\alpha_{0}(\check{\varpi}_{\beta}(\lambda))=\beta(\check{\varpi}_{\beta}(\lambda))^{-d(\beta)}=\lambda^{-1}$ puisque $d(\beta)=1$. Donc $Ad(\boldsymbol{\omega})(E_{\alpha_{0}})=\lambda^{-1}Ad(S(\omega))(E_{\alpha_{0}})$. Puisque $E_{\alpha_{0}}=\lambda Ad(S(\omega))(E_{\beta})$, on obtient $Ad(\boldsymbol{\omega})(E_{\alpha_{0}})=Ad(S(\omega)S(\omega))(E_{\beta})$. En utilisant (1), on obtient $Ad(\boldsymbol{\omega})(E_{\alpha_{0}})=Ad(t \,S(\omega^2))(E_{\beta})$, o\`u $t=t(\omega,\omega)$. On a $\omega^2(\beta)=\gamma$, d'o\`u, comme ci-dessus, $Ad(S(\omega^2))(E_{\beta})=E_{\gamma}$.  Il reste \`a prouver que $Ad(t)(E_{\gamma})=E_{\gamma}$, autrement dit $\gamma(t)=1$. Puisque $\gamma=\omega(\alpha_{0})$, cela \'equivaut \`a $\alpha_{0}(t_{0})=1$, o\`u $t_{0}=\omega^{-1}(t)$. On calcule
$$t_{0}=\prod_{\alpha\in \Sigma; \alpha>0, \omega(\alpha)<0,\omega^{-1}(\alpha)<0}\check{\alpha}(-1).$$
L'ensemble des $\alpha\in \Sigma$ tels que $\alpha>0$ et $\omega(\alpha)<0$ est $\Sigma^+-\Sigma^{\beta,+}$, o\`u les exposants $+$ d\'esignent les sous-ensembles de racines positives. De m\^eme, l'ensemble des $\alpha\in \Sigma$ tels que $\alpha>0$ et $\omega^{-1}(\alpha)<0$ est $\Sigma^+-\Sigma^{\gamma,+}$. Notons $\Sigma^{\beta,\gamma}$ le sous-ensemble de racines de $\Sigma$ engendr\'e par $\Sigma-\{\beta,\gamma\}$. En notant ${\bf 1}_{\Sigma^+}$, ${\bf 1}_{\Sigma^{\beta,+}}$ etc... les fonctions caract\'eristiques des ensembles $\Sigma^+$, $\Sigma^{\beta,+}$ etc..., on obtient qu'un \'el\'ement $\alpha\in \Sigma$ v\'erifie les conditions $\alpha>0$, $\omega(\alpha)<0$, $\omega^{-1}(\alpha)<0$ si et seulement si ${\bf 1}_{\Sigma^+}(\alpha)-{\bf 1}_{\Sigma^{\beta,+}}(\alpha)-{\bf 1}_{\Sigma^{\gamma,+}}(\alpha)+{\bf 1}_{\Sigma^{\beta,\gamma,+}}(\alpha)=1$. D\'efinissons de la fa\c{c}on usuelle $2\check{\rho}=\sum_{\alpha>0}\check{\alpha}$ et de m\^eme $2\check{\rho}^{\beta}$ etc.... Posons $X=2\check{\rho}-2\check{\rho}^{\beta}-2\check{\rho}^{\gamma}+2\check{\rho}^{\rho,\gamma}$. Alors $t_{0}=X(-1)$ et $\alpha_{0}(t_{0})=(-1)^{<\alpha_{0},X>}$. Il nous suffit de prouver que $<\alpha_{0},X>\in 2{\mathbb Z}$. 
 Notons que les \'el\'ements $2\check{\rho}$ etc... appartiennent \`a ${\mathbb Z}[\check{\Delta}]$ et que seules comptent leurs images dans ${\mathbb Z}[\check{\Delta}]/2{\mathbb Z}[\check{\Delta}]$. 

On va calculer $<\alpha_{0},X>$ cas par cas. 

Supposons $G$ de type $A_{n-1}$, avec $n\geq2$. On identifie de fa\c{c}on usuelle $\Delta_{a}$ \`a ${\mathbb Z}/n{\mathbb Z}$, $\alpha_{0}$ s'identifiant \`a $0$. Alors $\alpha_{0}=-\varpi_{1}-\varpi_{n-1}$. Soit $i\in {\mathbb Z}/n{\mathbb Z}$ l'\'el\'ement s'identifiant \`a $\beta$. Alors $\omega$ est la translation par $-i$. Il en r\'esulte que $\gamma$ s'identifie \`a $n-i$. Autrement dit, $\gamma$ est l'image de $\beta$ par l'automorphisme ext\'erieur habituel $\theta$ de $G$ qui pr\'eserve $\mathfrak{E}$ et envoie $\alpha_{j}$ sur $\alpha_{n-j}$. Il r\'esulte des d\'efinitions qu'alors $ X$ est invariant par cet automorphisme.  On a $X\in {\mathbb Z}[\check{\Delta}]$ donc $<\varpi_{1},X>\in {\mathbb Z}$. Puisque $\varpi_{n-1}=\theta(\varpi_{1})$, on obtient $<\alpha_{0},X>=-2<\varpi_{1},X>\in 2{\mathbb Z}$.

Supposons $G$ de type $B_{n}$, avec $n\geq2$. On a $\Omega\simeq {\mathbb Z}/2{\mathbb Z}$. L'unique \'el\'ement $\beta\in \Delta$ tel que $d(\beta)=1$ est   $\alpha_{1}$. Si $n=2$, $\alpha_{0}=-2\varpi_{\alpha_{2}}$ et forc\'ement $<\alpha_{0},X>\in 2{\mathbb Z}$. Supposons $n\geq3$. Alors $\alpha_{0}=-\varpi_{\alpha_{2}}$. Le syst\`eme $\check{\Sigma}$ est de type $C_{n}$ avec pour "racines simples" $\check{\alpha}_{1},...,\check{\alpha}_{n}$ et le syst\`eme $\check{\Sigma}^{\alpha_{1}}$ est de type $C_{n-1}$ avec pour "racines simples" $\check{\alpha}_{2},...,\check{\alpha}_{n}$. On voit que le coefficient de $\check{\alpha}_{2}$ dans $2\check{\rho}$ comme dans $2\check{\rho}^{\alpha_{1}}$ est pair. Donc $<\alpha_{0},X>\in 2{\mathbb Z}$.

Supposons $G$ de type $C_{n}$ avec $n\geq2$. On a $\Omega\simeq {\mathbb Z}/2{\mathbb Z}$.  L'unique \'el\'ement $\beta\in \Delta$ tel que $d(\beta)=1$ est   $\alpha_{n}$.   On a $\alpha_{0}\equiv -\alpha_{n}\,\,mod\,\, 2{\mathbb Z}[\Delta]$. Mais $<\alpha_{n},\check{\alpha}>\in 2{\mathbb Z}$ pour tout $\alpha\in \Delta$, donc  $<\alpha_{0},X>\in 2{\mathbb Z}$.

 Supposons $G$ de type $D_{n}$ avec $n\geq4$. On a $\Omega\simeq ({\mathbb Z}/2{\mathbb Z})^2$ si $n$ est pair et $\Omega\simeq {\mathbb Z}/4{\mathbb Z}$ si $n$ est impair. Les \'el\'ements $\beta\in \Delta$ tels que $d(\beta)=1$ sont  $\alpha_{1},\alpha_{n-1},\alpha_{n}$ et on a $\alpha_{0}=-\varpi_{\alpha_{2}}$. Le coefficient de $\check{\alpha}_{2}$ dans $2\check{\rho}$ est pair. Le syst\`eme $\Sigma^{\alpha_{1}}$ est de type $D_{n-1}$ et le coefficient de $\check{\alpha}_{2}$ dans $2\check{\rho}^{\alpha_{1}}$ est pair. Les syst\`emes $\Sigma^{\alpha_{n-1}}$ et $\Sigma^{\alpha_{n}}$ sont de type $A_{n-1}$ et les coefficients de $\check{\alpha}_{2}$ dans $2\check{\rho}^{\alpha_{n-1}}$ et $2\check{\rho}^{\alpha_{n}}$  sont pairs. On ne peut avoir $\beta\not=\gamma$ que si $n$ est impair et $\{\beta,\gamma\}=\{\alpha_{n-1},\alpha_{n}\}$. Dans ce cas, $\Sigma^{\alpha_{n-1},\alpha_{n}}$ est de type $A_{n-2}$. Les hypoth\`eses $n\geq4$ et $n$ impair entra\^{\i}nent $n-2\geq 3$ et on voit encore que le coefficient de $\check{\alpha}_{2}$ dans $2\check{\rho}^{\alpha_{n-1},\alpha_{n}}$ est pair.  Donc $<\alpha_{0},X>\in 2{\mathbb Z}$.
 
 Supposons $G$ de type $E_{6}$. On a $\Omega\simeq {\mathbb Z}/3{\mathbb Z}$.  Les \'el\'ements $\beta\in \Delta$ tels que $d(\beta)=1$ sont  $\alpha_{1},\alpha_{6}$ et on a $\alpha_{0}=-\varpi_{\alpha_{2}}$. On a $2\check{\rho}\in 2{\mathbb Z}[\check{\Delta}]$. Les syst\`emes $\Sigma^{\alpha_{1}}$ et $\Sigma^{\alpha_{6}}$ sont de type $D_{5}$, le syst\`eme $\Sigma^{\alpha_{1},\alpha_{6}}$ est de type $D_{4}$ et on voit que $2\check{\rho}^{\alpha_{1}}$, $2\check{\rho}^{\alpha_{6}}$ et $2\check{\rho}^{\alpha_{1},\alpha_{6}}$ appartiennent tous \`a $2{\mathbb Z}[\check{\Delta}]$. Donc $<\alpha_{0},X>\in 2{\mathbb Z}$.
 
 Supposons $G$ de type $E_{7}$. On a $\Omega\simeq {\mathbb Z}/2{\mathbb Z}$.  L'unique \'el\'ement $\beta\in \Delta$ tel que $d(\beta)=1$ est  $\alpha_{7}$ et on a $\alpha_{0}=-\varpi_{\alpha_{1}}$. Le coefficient de $\check{\alpha}_{1}$ dans $2\check{\rho}$ est pair. Le syst\`eme $\Sigma^{\alpha_{7}}$ est de type $E_{6}$ et on a $2\check{\rho}^{\alpha_{7}}\in 2{\mathbb Z}[\check{\Delta}]$. Donc $<\alpha_{0},X>\in 2{\mathbb Z}$.
 
 Pour les groupes de type $E_{8}$, $F_{4}$ ou $G_{2}$, $\Omega=\{1\}$ et le lemme est trivial. $\square$
 
   On identifie le  groupe $Out(G)$ au sous-groupe $\boldsymbol{Out}(G)$ des \'el\'ements de $Aut(G)$ qui conservent l'\'epinglage $\mathfrak{E}$.  On a:
     
  (9) si $G$ n'est pas de type $A_{n-1}$ avec $n$ impair,  tout \'el\'ement de $\boldsymbol{ Out}(G)$ fixe $E_{0}$. 
  
  Cf. \cite{KS} 1.3.5.  
  
  Supposons que $G$ n'est pas de type $A_{n-1}$ avec $n$ impair. Alors $\boldsymbol{Out}(G)$ normalise $\boldsymbol{\Omega}$. 
  Evidemment, $Out(G)$ s'identifie \`a $Aut({\cal D})$. Le lemme ci-dessus permet d'identifier   $Aut({\cal D}_{a})=\Omega\rtimes Aut({\cal D})$  au sous-groupe $\boldsymbol{Aut}({\cal D}_{a})=Ad(\boldsymbol{\Omega})\rtimes \boldsymbol{Out}(G)$ de $Aut(G)$.

 \subsubsection{Description d'un appartement et d'une alc\^ove}\label{alcove}

On revient à nos hypothèses habituelles, le groupe $G$ étant défini sur $F$. Supposons d'abord $G$ quasi-d\'eploy\'e. On fixe une paire de Borel \'epingl\'ee  $\mathfrak{E}=(B,T,(E_{\alpha})_{\alpha\in \Delta})$ conserv\'ee par $\Gamma_{F}$. 
  Il y a un unique point $s_0\in Imm_{F^{nr}}(G_{AD})$ de sorte que $\mathfrak{k}_{s_0}^{nr}$ soit l'ensemble des points fixes par $I_{F}$ dans le $\mathfrak{o}_{\bar{F}}$-r\'eseau  somme des $\mathfrak{u}_{\alpha,\mathfrak{o}_{\bar{F}}}$ et   du réseau engendré par $X_{*}(T)\subset \mathfrak{t}$.  On note $T^{nr}$ le plus grand sous-tore de $T$ d\'eploy\'e sur $F^{nr}$. On a $X_{*}(T^{nr})=X_{*}(T)^{I_{F}}$. On note $C^{nr}$ l'alc\^ove  de $Imm_{F^{nr}}(G_{AD})$ qui v\'erifie les conditions suivantes: elle est contenue dans  l'appartement  $App_{F^{nr}}(T^{nr})$ associ\'e \`a $T^{nr}$, $s_0$ appartient \`a la cl\^oture de $C^{nr}$ et $\mathfrak{k}_{C^{nr}}^{nr}\cap \mathfrak{b}(F^{nr})=\mathfrak{k}_{s_0}^{nr}\cap \mathfrak{b}(F^{nr})$. Le point $s_0$ est fix\'e par l'action galoisienne de $\Gamma_{F}^{nr}$ et l'alc\^ove $C^{nr}$ est conserv\'ee par cette action.  A l'alc\^ove  $C^{nr}$ est associ\'e un ensemble de racines affines: les annulateurs de ces racines affines sont les murs de l'alc\^ove et elles sont positives sur celle-ci. Notons $\Delta_{a}^{nr}$ l'ensemble de racines sous-jacent \`a cet ensemble de racines affines. Ce sont des formes lin\'eaires sur $X_{*}(T)^{I_{F}}=X_{*}(T^{nr})$.  Cet ensemble de racines $\Delta_{a}^{nr}$ plus les relations de produit scalaire entre elles forment un diagramme ${\cal D}_{a}^{nr}$. Si $G$ est irréductible sur $F^{nr}$,  les \'el\'ements de $S(\bar{C}^{nr})$ correspondent aux \'el\'ements de $\Delta_{a}^{nr}$. On note $\beta\mapsto s_{\beta}$ cette bijection.  Pour toute facette ${\cal F}$ de $Imm_{F^{nr}}(G_{AD})$ contenue dans  $\bar{C}^{nr}$, le groupe $G_{{\cal F}}$ est muni d'une paire de Borel $(B_{{\cal F}},T_{{\cal F}})$: $B_{{\cal F}}$ est l'image par r\'eduction de $K^{0,nr}_{C^{nr}}\subset K_{{\cal F}}^{0,nr}$ et 
  $T_{{\cal F}}$ est l'image de $T(F^{nr})\cap K_{{\cal F}}^{0,nr}$ ou encore de  $ T^{nr}(\mathfrak{o}_{F^{nr}})\subset K_{{\cal F}}^{0,nr}$.

Posons $N^{nr}=K_{C^{nr},AD}^{\dag,nr}\cap Norm_{G_{AD}(F^{nr})}(T)$, $N_{0}^{nr}=K_{C^{nr},AD}^{0,nr}\cap Norm_{G_{AD}(F^{nr})}(T)$ et $\bar{N}^{nr}=N^{nr}/N^{nr}_{0}\simeq K_{C^{nr},AD}^{\dag,nr}/K_{C^{nr},AD}^{0,nr}$. 
 D'apr\`es \cite{HR} lemme 14, on a $\bar{N}^{nr}\simeq (Z(\hat{G}_{SC})^{I_{F}})^{\vee}$, en particulier, c'est un groupe ab\'elien. On sait que $H^1(\Gamma_{F},G_{AD})\simeq H^1(\Gamma_{F}^{nr},\bar{N}^{nr})$ (cf. \cite{W2} p. 441) et ce dernier groupe n'est autre que le quotient $\bar{N}^{nr}_{\Gamma_{F}^{nr}}$ de $\bar{N}^{nr}$ par le sous-groupe des $Fr(n)n^{-1}$ pour $n\in \bar{N}^{nr}$.  Remarquons que $\bar{N}^{nr}_{\Gamma^{nr}_{F}}$ est aussi \'egal \`a $(Z(\hat{G}_{SC})^{\Gamma_{F}})^{\vee}$ et on retrouve l'isomorphisme de Kottwitz $H^1(\Gamma_{F},G_{AD})\simeq (Z(\hat{G}_{SC})^{\Gamma_{F}})^{\vee}$, cf. \cite{K2} th\'eor\`eme 1.2.

Supprimons l'hypoth\`ese que $G$ est quasi-d\'eploy\'e. Comme en \ref{formeinterieure}, on  fixe une forme quasi-d\'eploy\'ee $G^*$ de $G$, un  torseur int\'erieur $\psi:G\to G^*$ et le cocycle associ\'e $u_{G}$. On effectue  les constructions pr\'ec\'edentes pour ce groupe $G^*$ (on ajoute si besoin est  des exposants $*$). Au cocycle $u_{G}$ est associ\'e un \'el\'ement de $\bar{N}^{nr,*}_{\Gamma_{F}^{nr}}$. Relevons ce dernier en un \'el\'ement de $\bar{N}^{nr,*}$. On peut ensuite relever celui-ci en un \'el\'ement $n_{G}\in N^{nr,*}$ de sorte que $n_{G}$  se prolonge en un cocycle de $\Gamma_{F}^{nr}$ dans $G_{AD}^*(F^{nr})$, c'est-\`a-dire un cocycle $\underline{n}_{G}$ tel que $\underline{n}_{G}(Fr)=n_{G}$.  On note encore $\underline{n}_{G}$ le cocycle de $\Gamma_{F}$ dans $G^*(F^{nr})$ obtenu par inflation. On peut alors supposer que $G(F^{nr})=G^*(F^{nr})$, que $\psi$ est l'identit\'e et que l'on a les \'egalit\'es d'actions galoisiennes $\sigma_{G}(g)=\underline{n}_{G}(\sigma)\sigma_{G^*}(g)\underline{n}_{G}(\sigma)^{-1}$ pour tous $g\in G(F^{nr})$ et $\sigma\in \Gamma_{F}$. Le tore $T$ est conserv\'e par cette action. Remarquons que $T^{nr}$ est encore le plus grand sous-tore de $T$ d\'eploy\'e sur $F^{nr}$ puisque les deux actions sur $G$ et $G^*$ co\"{\i}ncident sur $I_{F}$.  On note $S$ le plus grand sous-tore de $T$ d\'eploy\'e sur $F$ pour l'action $\sigma\mapsto \sigma_{G}$. C'est un élément de ${\cal T}_{max}$. Le groupe $\Gamma_{F}^{nr}$ agit sur $Imm_{F^{nr}}(G_{AD})$ par $(\sigma,x)\mapsto \underline{n}_{G}(\sigma)\circ\sigma(x)$ et $Imm(G_{AD})$ est le sous-ensemble des points fixes. L'alc\^ove $C^{nr}$ est conserv\'ee par l'action galoisienne. On note $C$ son sous-ensemble de points fixes. C'est une alc\^ove de $Imm(G_{AD})$, qui appartient \`a l'appartement $App(S)$ associ\'e \`a $S$. L'ensemble $S(\bar{C})$ des sommets contenus dans l'adh\'erence de $C$ s'identifie \`a celui des orbites de l'action $\sigma\mapsto \sigma_{G}$ de $\Gamma_{F}^{nr}$ dans $S(\bar{C}^{nr})$. 

{\bf Remarque.} On a d\'ecrit deux façons  de repr\'esenter la classe de cocycle $u_{G}$: ci-dessus, on la repr\'esente par un cocycle qui conserve $C^{nr}$; en \ref{formeinterieure}, on l'a repr\'esent\'ee par un cocycle qui conserve $M_{min}^*$. Ces deux repr\'esentations ne sont en g\'en\'eral pas compatibles.

\subsubsection{Racines affines}\label{racinesaffines}

  Dans ce paragraphe,  {\bf on suppose que} $G$ {\bf est absolument quasi-simple}. On introduit les objets d\'efinis en \ref{alcove}. La paire de Borel \'epingl\'ee $\mathfrak{E}$ est conserv\'ee par $I_{F}$.

  D\'ecrivons plus pr\'ecis\'ement l'appartement $App_{F^{nr}}(T^{nr})$ et la chambre $C^{nr}$. On pose ${\cal A}^{nr}=X_{*,{\mathbb R}}(T^{nr}) $. L'appartement $App_{F^{nr}}(T^{nr})$ est isomorphe \`a ${\cal A}^{nr}$, le point $s_0$ s'identifiant \`a $0\in {\cal A}^{nr}$. Toute racine $\alpha\in \Sigma$ se restreint en une forme lin\'eaire $\alpha^{res}$ sur ${\cal A}^{nr}$. On note $\Sigma^{nr}$ l'ensemble des $\alpha^{res}$ pour $\alpha\in \Sigma$. C'est un syst\`eme de racines irr\'eductible, pas toujours r\'eduit. On note $\Delta^{nr}$ le sous-ensemble des $\alpha^{res}$ pour $\alpha\in \Delta$. C'est une base du sous-syst\`eme des racines indivisibles de $\Sigma^{nr}$.

  Pour $\beta\in \Sigma^{nr}$ et $b\in {\mathbb R}$, notons $\beta[b]$ la forme affine $H\mapsto \beta(H)+b$ sur ${\cal A}^{nr}$. Pour tout $\beta\in \Sigma^{nr}$, il existe un sous-ensemble $\Gamma_{\beta}\subset {\mathbb R}$ tel que l'ensemble des racines affines $\Sigma^{aff}$ sur ${\cal A}^{nr}$ soit celui des $\beta[b]$, o\`u $\beta\in \Sigma^{nr}$ et $b\in \Gamma_{\beta}$. Pour $\beta\in \Delta^{nr}$, posons $\beta^{aff}=\beta[0]=\beta$. Il existe une racine suppl\'ementaire $\beta_{0}\in \Sigma^{nr}$ et un \'el\'ement $b_{0}\in \Gamma_{\beta_{0}}$ de sorte qu'en posant $\beta_{0}^{aff}=\beta_{0}[b_{0}]$ et $\Delta_{a}^{nr}=\{\beta_{0}\}\cup \Delta^{nr}$, on ait:
  
  l'ensemble $C^{nr}$ est celui des $H\in {\cal A}^{nr}$ tels que $\beta^{aff}(H)>0$ pour tout $\beta\in \Delta_{a}^{nr}$.

  Notons $E$ la plus petite extension de $F^{nr}$ telle que $G$ soit d\'eploy\'e sur $E$ et posons $e=[E:F^{nr}]$. On a en fait $b_{0}=\frac{1}{e}$, cf. \cite{W7} 1.6. Il existe une unique relation
  $$(1) \qquad \sum_{\beta\in \Delta_{a}^{nr}}d(\beta)\beta^{aff}=1,$$
  o\`u $d(\beta_{0})=e$. On a alors $d(\beta)\in {\mathbb N}_{>0}$ pour tout $\beta\in \Delta_{a}^{nr}$. 
  
  {\bf Notations.} L'ensemble $\Delta$ est num\'erot\'e comme dans \cite{B}.  Soit $\beta\in \Delta^{nr}$, notons $\alpha_{i_{1}}$, $\alpha_{i_{2}}$,... les \'el\'ements de $\Delta$ dont $\beta$ est la restriction. On note alors $\beta= \beta_{i_{1},i_{2},...} $. L'ensemble $S(\bar{C}^{nr})$ des sommets adh\'erents \`a $C^{nr}$ est naturellement en bijection avec $\Delta_{a}^{nr}$. Pour $\beta=\beta_{i_{1},i_{2},...}\in \Delta_{a}^{nr}$, on note $s_{\beta}$ ou $s_{i_{1},i_{2},...}$ le sommet associ\'e \`a $\beta$. Remarquons que le sommet $s_{0}$ d\'ej\`a introduit est associ\'e \`a $\beta_{0}$.
  \bigskip
  
  Le groupe $\Gamma_{E/F^{nr}}$ est cyclique    
    et agit sur ${\cal D}$ par un homomorphisme injectif $\Gamma_{E/F^{nr}}\to Aut({\cal D})$.  On  fixe un g\'en\'erateur $\gamma$ 
 de $\Gamma_{E/F^{nr}}$ et on note  $\theta\in Aut({\cal D})$  son image. On note encore $\theta$ l'unique automorphisme de $G$ qui conserve l'\'epinglage et dont l'action d\'eduite sur ${\cal D}$ est $\theta$. Les actions de $\gamma$ et $\theta$ coïncident sur $E_{\alpha}$ pour tout $\alpha\in \Delta$. Plus généralement, elles coïncident sur les $F$-droites $\mathfrak{u}_{\alpha,F}$ pour tout $\alpha\in \Sigma$.  Soit $\beta\in \Sigma^{res}$.  L'action de $\theta$ permute cycliquement l'ensemble des $\alpha\in \Sigma$ tels que $\alpha^{res}=\beta$. Notons $e(\beta)$ le nombre d'éléments de cet ensemble.   Notons $\mathfrak{u}_{\beta}$ la somme des $\mathfrak{u}_{\alpha}$ sur cet ensemble de racines $\alpha$.  Plus généralement, pour tout sous-${\mathbb Z}_{p}$-module $R\subset \bar{F}$, on note $\mathfrak{u}_{\beta,R}$ la somme des $\mathfrak{u}_{\alpha,R}$.   Si $\beta$ est indivisible, $\theta^{e(\beta)}$ se restreint en l'identité de $\mathfrak{u}_{\alpha}$. Pour $\beta\in \Delta^{nr}$, $\beta$ est toujours indivisible et on pose $E_{\beta}=\sum_{n=0,...,e(\beta)-1}E_{\theta^n(\alpha)}$. C'est un générateur du $\mathfrak{o}_{F^{nr}}$-module $\mathfrak{u}_{\beta,\mathfrak{o}_{F^{nr}}}^{I_{F}}$. 
 Fixons   une uniformisante $\varpi_{E}$ de $E$ telle que $\varpi_{E}^{e}\in F^{nr,\times}$. On fixe un générateur $E_{\beta_{0}}$ du $\mathfrak{o}_{F^{nr}}$-module $\mathfrak{u}_{\beta_{0},\varpi_{E}\mathfrak{o}_{F^{nr}}}^{I_{F}}$. Concrètement, deux cas se produisent. En général, $\beta_{0}$ est indivisible et $e(\beta_{0})=e$. On fixe alors  $\alpha\in \Sigma$ tel que $\alpha^{res}=\beta_{0}$, puis un générateur $E_{\alpha}$ du $\mathfrak{o}_{F^{nr}}$-module $\mathfrak{u}_{\alpha,\mathfrak{o}_{F^{nr}}}$ et on peut choisir $E_{\beta_{0}}=\sum_{n=0,...,e-1}\gamma^{n}(\varpi_{E})\theta^n(E_{\alpha})$. 
 Il y a un cas particulier: c'est celui o\`u $G$ est de type $A_{n-1}$ avec $n\geq 3$, $n$  impair et o\`u $E$ est l'extension quadratique de $F^{nr}$. Dans ce cas, $e=2$ et $e(\beta_{0})=1$. Il n'y a qu'une racine $\alpha\in \Sigma$ telle que $\alpha^{res}=\beta_{0}$.  L'opérateur $\theta$  agit sur  $\mathfrak{ u}_{\alpha}$ par moins l'identité. On fixe un générateur $E_{\alpha}$ du $\mathfrak{o}_{F^{nr}}$-module $\mathfrak{u}_{\alpha,\mathfrak{o}_{F^{nr}}}$  et on peut choisir $E_{\beta_{0}}=\varpi_{E}E_{\alpha}$. 
 
   On pose $\mathfrak{E}_{a}=(T^{nr},(E_{\beta})_{\beta\in \Delta_{a}^{nr}})$ et on appelle cette famille un \'epinglage affine. Remarquons que l'on a supprim\'e le groupe $B$ et que,  par construction,  $\mathfrak{E}_{a}$ est conserv\'e par $I_{F}$.
 
 {\bf Remarque} (2). Supposons que $G$ soit quasi-d\'eploy\'e sur $F$ et que $\mathfrak{E}$ soit conserv\'e par l'action de $\Gamma_{F}$. Il est clair que la famille $(E_{\beta})_{\beta\in \Delta^{nr}}$ est conserv\'ee par cette action. On peut supposer que $E_{\beta_{0}}$ est fix\'ee par cette action. En effet l'\'el\'ement construit ci-dessus est fix\'e par $I_{F}$. Pour un \'el\'ement $\sigma\in \Gamma_{F}$, on v\'erifie qu'il existe $x(\sigma)\in \mathfrak{o}_{F^{nr}}^{\times}$ tel que $\sigma(E_{\beta_{0}})=x(\sigma)E_{\beta_{0}}$ et $x(\sigma)$ ne d\'epend que de la projection de $\sigma$ dans $\Gamma_{F}^{nr}$. L'application $\sigma\mapsto x(\sigma)$ est un cocycle de $\Gamma_{F}^{nr}$ dans $\mathfrak{o}^{\times}_{F^{nr}}$. Un tel cocycle est un cobord donc il existe $x\in \mathfrak{o}^{\times}_{F^{nr}}$ tel que $x(\sigma)=x\sigma(x)^{-1}$. L'\'el\'ement $xE_{\beta_{0}}$ est alors fix\'e par $\Gamma_{F}$. Mais $xE_{\beta_{0}}$ est issu de la m\^eme construction que $E_{\beta_{0}}$ o\`u l'on a remplac\'e l'élément $E_{\alpha}$ par $xE_{\alpha}$. 
 
 \bigskip
 
 On note $\boldsymbol{\Omega}^{nr}$ le sous-groupe des \'el\'ements de $G_{AD}(F^{nr})$ qui conservent $\mathfrak{E}_{a}$. Il est inclus dans le groupe $N^{nr}$ d\'efini en \ref{alcove} et  il y a un homomorphisme naturel $N^{nr}\to Aut({\cal D}_{a}^{nr})$. 
 
 \begin{lem}{ L'homomorphisme $\boldsymbol{\Omega}^{nr}\to Aut({\cal D}_{a}^{nr})$ est injectif et les groupes $\boldsymbol{\Omega}^{nr}$ et $N^{nr}$ ont m\^eme image dans $Aut({\cal D}_{a}^{nr})$. Si $G$ est déployé sur $F^{nr}$, cette image est le groupe $\Omega$ d\'efini en \ref{epinglages}. Si $G$ n'est pas déployé sur $F^{nr}$, cette image est  $Aut({\cal D}_{a}^{nr})$ tout entier.}\end{lem}
 
 Preuve. 
 Supposons que $G$ soit   d\'eploy\'e sur $F^{nr}$. On a alors $e=1$, $\Sigma^{nr}=\Sigma$, $\beta_{0}=\alpha_{0}$, $\Delta_{a}^{nr}=\Delta_{a}$, ${\cal D}_{a}^{nr} ={\cal D}_{a}$,  $\boldsymbol{\Omega}^{nr}=\boldsymbol{\Omega}$.  Il est clair que l'image de $N^{nr}$ dans $Aut({\cal D}_{a})$ est incluse dans  $\Omega$.  Le lemme \ref{epinglages} implique alors les assertions de l'énoncé.
 
  Supposons  que $G$  ne soit pas d\'eploy\'e sur $F^{nr}$.  On ne perd rien à supposer qu'il est adjoint. Soit $g\in \boldsymbol{\Omega}^{nr}$ un \'el\'ement du noyau de l'homomorphisme $\boldsymbol{\Omega}^{nr}\to Aut({\cal D}_{a}^{nr})$. Puisque $g\in \boldsymbol{\Omega}^{nr}$, $Ad(g)$ conserve $T^{nr}$. Mais le commutant de $T^{nr}$ dans $G$ est $T$ donc $g$ normalise $T$. Il permute donc les espaces radiciels $\mathfrak{u}_{\alpha}$ pour $\alpha\in \Sigma$. Soit $\beta\in \Delta^{res}$. Alors $Ad(g)$ fixe $E_{\beta}$. Il permute donc les projections de $E_{\beta}$ sur les espaces radiciels, donc permute les $E_{\alpha}$ pour $\alpha\in \Delta$ tel que $\alpha^{res}=\beta$. Cette propriété étant vraie pour tous les $\beta\in \Delta^{res}$, $Ad(g)$ permute les $E_{\alpha}$ pour $\alpha\in \Delta$, autrement dit  $Ad(g)$ conserve l'\'epinglage $\mathfrak{E}$. Puisqu'on a suppos\'e $G$ adjoint, cela entra\^{\i}ne $g=1$. D'o\`u la première assertion d'injectivité  de l'\'enonc\'e.

En examinant tous les cas possibles, on voit que $Aut({\cal D}_{a}^{nr})$ est r\'eduit \`a $\{1\}$ sauf dans les deux cas suivants:

(3) $G$ est de type $A_{n-1}$ avec $n\geq4$, $n$ est pair et le corps $E$ est l'extension quadratique de $F^{nr}$;

(4) $G$ est de type $D_{n}$ avec $n\geq4$ et le corps $E$ est l'extension quadratique de $F^{nr}$.

Dans ces cas, $Aut({\cal D}_{a}^{nr})\simeq {\mathbb Z}/2{\mathbb Z}$. Pour  achever la preuve du lemme, on peut supposer que l'on est dans l'un de ces deux cas et relever dans $\boldsymbol{\Omega}^{nr}$ l'unique \'el\'ement non trivial $\omega\in Aut({\cal D}_{a}^{nr})$. 

Consid\'erons le cas (3). L'\'el\'ement non trivial de $\Gamma_{E/F^{nr}}$ agit sur $\Delta$ par l'automorphisme $\theta$ habituel: $\theta(\alpha_{i})=\alpha_{n-i}$.  On a $\Delta^{nr}=\{\beta_{i,n-i}; i=1,...,n/2-1\}\cup \{\beta_{n/2}\}$. La racine $\beta_{0}$ est la restriction commune de $\alpha'=-(\alpha_{1}+...+\alpha_{n-2})$ et de $\alpha''=-(\alpha_{2}+...+\alpha_{n-1})$. On a $E_{\beta_{i,n-i}}=E_{\alpha_{i}}+E_{\alpha_{n-i}}$, $E_{\beta_{n/2}}=E_{\alpha_{n/2}}$. On peut fixer des générateurs $E_{\alpha'}$, resp. $E_{\alpha''}$, de $\mathfrak{u}_{\alpha',{\mathbb Z}_{p}}$, resp. $\mathfrak{u}_{\alpha'',{\mathbb Z}_{p}}$, de sorte que $\theta(E_{\alpha'})=E_{\alpha''}$. Alors 
  $E_{\beta_{0}}=\lambda\varpi_{E}(E_{\alpha'}-E_{\alpha''})$ avec $\lambda\in \mathfrak{o}_{F^{nr}}^{\times}$. L'\'el\'ement $\omega$ permute $\beta_{0}$ et $\beta_{1,n-1}$ et fixe les autres \'el\'ements de $\Delta_{a}^{nr}$. Introduisons le Levi standard $M$ de $G$ associ\'e au sous-ensemble $\Delta^M\subset \Delta$ d\'efini par $\Delta^M=\{\alpha_{i}; i=2,...,n-2\}$. Notons $w_{0}$, resp. $w_{0}^M$, l'\'el\'ement de plus grande longueur de $W$, resp.  $W^M$, posons $w=w_{0}w_{0}^M$. On v\'erifie que $w$ est fix\'e par $\theta$,  que $w$ fixe $\alpha_{i}$ pour $i=2,...,n-2$ et que $w$ permute $\alpha_{1}$ et $\alpha''$ ainsi que $\alpha_{n-1}$ et $\alpha'$. Donc
 l'action de $w$ sur $\Sigma$ se descend en une action sur $\Sigma^{res}$ qui conserve $\Delta_{a}^{nr}$ et agit sur cet ensemble par $\omega$. Il existe un \'el\'ement $\mu\in \mathfrak{o}_{F^{nr}}^{\times}$ tel que $Ad(S(w))(E_{\alpha_{1}})=\mu E_{\alpha''}$. L'élément $S(w)$ est fix\'e par $\theta$ et on en d\'eduit $Ad(S(w))(E_{\alpha_{n-1}})=\mu E_{\alpha'}$. 
 Posons $\boldsymbol{\omega}=S(w)\check{\varpi}_{\alpha_{1}}(-\lambda\mu^{-1}\varpi_{E})\check{\varpi}_{\alpha_{n-1}}(\lambda\mu^{-1}\varpi_{E})$. C'est un \'el\'ement de $G(F^{nr})$. La d\'efinition est faite pour  que $Ad(\boldsymbol{\omega})(E_{\beta_{1,n-1}})=E_{\beta_{0}}$. Les m\^emes arguments qu'en \ref{epinglages} montrent que $Ad(\boldsymbol{\omega})$ fixe $E_{\beta_{n/2}}$ et $E_{\beta_{i,n-i}}$ pour $i=2,...,n/2-1$. Il reste \`a prouver que $Ad(\boldsymbol{\omega})(E_{\beta_{0}})=E_{\beta_{1,n-1}}$. On a
 $$(5) \qquad Ad(\boldsymbol{\omega})(E_{\beta_{0}})=\lambda\varpi_{E}Ad(S(w)\check{\varpi}_{\alpha_{1}}(-\lambda\mu^{-1}\varpi_{E})\check{\varpi}_{\alpha_{n-1}}(\lambda\mu^{-1}\varpi_{E}))(E_{\alpha'}-E_{\alpha''})$$
 $$=\lambda\varpi_{E}Ad(S(w))(-\lambda^{-1}\mu\varpi_{E}^{-1}E_{\alpha'}-\lambda^{-1}\mu\varpi_{E}^{-1}E_{\alpha''})=- \mu Ad(S(w))(E_{\alpha'}+E_{\alpha''})$$
 $$=-Ad(S(w)^2)(E_{\alpha_{1}}+E_{\alpha_{n-1}}).$$
 Notons $\Sigma^+$, resp. $\Sigma^{M+}$,  l'ensemble des racines positives dans $\Sigma$, resp. $\Sigma^M$. L'ensemble des $\alpha>0$ telles que $w(\alpha)<0$ est $\Sigma^+-\Sigma^{M+}$. 
 Gr\^ace \`a \ref{epinglages}(4), on a $S(w)^2=t$, o\`u $t=\prod_{\alpha\in \Sigma^+-\Sigma^{M,+}}\check{\alpha}(-1)$.  Alors (5) devient
 $Ad(\boldsymbol{\omega})(E_{\beta_{0}})=-Ad(t)(E_{\beta_{1,n-1}})$. Pour obtenir l'\'egalit\'e cherch\'ee $Ad(\boldsymbol{\omega})(E_{\beta_{0}})=E_{\beta_{1,n-1}}$, il nous reste \`a prouver que $\alpha_{1}(t)=\alpha_{n-1}(t)=-1$. On calcule facilement la somme des \'el\'ements de $\Sigma-\Sigma^{M+}$: c'est $(n-1)(\alpha_{1}+...+\alpha_{n-1})$. D'o\`u $t=(\check{\alpha}_{1}+...+\check{\alpha}_{n-1})((-1)^{n-1})$. D'o\`u $\alpha_{1}(t)=\alpha_{n-1}(t)=(-1)^{n-1}$. Puisque $n$ est pair, on a bien $\alpha_{1}(t)=\alpha_{n-1}(t)=-1$, ce qui ach\`eve la preuve dans le cas (3).
 
  Consid\'erons le cas (4). L'\'el\'ement non trivial de $\Gamma_{E/F^{nr}}$ agit sur $\Delta$ par l'automorphisme $\theta$ habituel: $\theta$ permute $\alpha_{n-1}$ et $\alpha_{n}$ et fixe les autres \'el\'ements de $\Delta$. On a $\Delta^{nr}=\{\beta_{i}; i=1,...,n-2\}\cup \{\beta_{n-1,n}\}$. Posons
  $$(\alpha',\alpha'')=\left\lbrace\begin{array}{cc}(-(\alpha_{1}+...+\alpha_{n-2}+\alpha_{n-1}),-(\alpha_{1}+...+\alpha_{n-2}+\alpha_{n}))&\text{ si }n \text{ est pair }\\ (-(\alpha_{1}+...+\alpha_{n-2}+\alpha_{n}),-(\alpha_{1}+...+\alpha_{n-2}+\alpha_{n-1}))&\text{ si }n \text{ est impair }\\ \end{array}\right.$$
 La racine $\beta_{0}$ est la restriction commune de $\alpha' $ et de $\alpha'' $. On a $E_{\beta_{i}}=E_{\alpha_{i}}$ pour $i=1,...,n-2$, $E_{\beta_{n-1,n}}=E_{\alpha_{n-1}}+E_{\alpha_{n}}$. On peut fixer des éléments $E_{\alpha'}$ et $E_{\alpha''}$ comme dans le cas (3), de sorte que   $E_{\beta_{0}}=\lambda\varpi_{E}(E_{\alpha'}-E_{\alpha''})$ avec $\lambda\in \mathfrak{o}_{F^{nr}}^{\times}$. L'\'el\'ement $\omega$ permute $\beta_{n,n-1}$ et $\beta_{0}$ et permute  $\beta_{i}$ et $\beta_{n-1-i}$ pour $i=1,...,n-2$. Introduisons le Levi standard $M$ de $G$ associ\'e au sous-ensemble $\Delta^M\subset \Delta$ d\'efini par $\Delta^M=\{\alpha_{i}; i=1,...,n-2\}$. Notons $w_{0}$, resp. $w_{0}^M$, l'\'el\'ement de plus grande longueur de $W$, resp.  $W^M$, posons $w=w_{0}w_{0}^M$. On v\'erifie que $w$ est fix\'e par $\theta$, que $w$ permute $\alpha_{i}$ et $\alpha_{n-1-i}$ pour $i=1,...,n-2$ et permute $\alpha_{n-1}$ et $\alpha'$ ainsi que $\alpha_{n}$ et $\alpha''$. Donc
 l'action de $w$ sur $\Sigma$ se descend en une action sur $\Sigma^{res}$ qui conserve $\Delta_{a}^{nr}$ et agit sur cet ensemble par $\omega$. Il existe un \'el\'ement $\mu\in \mathfrak{o}_{F^{nr}}^{\times}$ tel que $Ad(S(w))(E_{\alpha_{n-1}})=\mu E_{\alpha'}$. Comme ci-dessus, on a aussi $Ad(S(w))(E_{\alpha_{n}})=\mu E_{\alpha''}$. On pose $\boldsymbol{\omega}=S(w)\check{\varpi}_{n-1}(\lambda\mu^{-1}\varpi_{E})\check{\varpi}_{n}(-\lambda\mu^{-1}\varpi_{E})$. C'est un \'el\'ement de $G(F^{nr})$. Comme dans le cas (3), il faut prouver que  $Ad(\boldsymbol{\omega})(E_{\beta_{0}})=E_{\beta_{n-1,n}}$. Un calcul analogue \`a (5) conduit \`a l'\'egalit\'e
 $$Ad(\boldsymbol{\omega})(E_{\beta_{0}})=(-1)^nAd(S(w)^2)(E_{\beta_{n-1,n}}).$$
 On a  $S(w)^2=t$, o\`u $t$ est d\'efini comme dans le cas (3). Il reste \`a prouver que $\alpha_{n-1}(t)=\alpha_{n}(t)=(-1)^n$. On note $2\check{\rho}$, resp. $2\check{\rho}^M$, la somme des coracines $\check{\alpha}$ pour $\alpha\in \Sigma$, resp. $\alpha\in \Sigma^M$. On a $t=(2\check{\rho})(-1)(2\check{\rho}^M)(-1)$. Puisque $\Sigma$ est de type $D_{n}$, $2\check{\rho}$ appartient \`a $2X_{*}(T)$, donc $(2\check{\rho})(-1)=1$. L'\'el\'ement $2\check{\rho}^M$ est combinaison lin\'eaire des $\check{\alpha}_{i}$ pour $i=1,...,n-2
$ et seule $\check{\alpha}_{n-2}$ n'est pas annul\'ee par $\alpha_{n-1}$ et $\alpha_{n}$. Le coefficient de $\check{\alpha}_{n-2}$ est $n-2$. On obtient $\alpha_{n-1}(t)=\alpha_{n}(t)=(-1)^{n-2}=(-1)^n$, ce qui ach\`eve la d\'emonstration. $\square$

\subsubsection{Un corollaire}\label{uncorollaire}
On garde les hypothèses du paragraphe précédent. Le groupe $T$ est naturellement défini sur $\mathfrak{o}_{F^{nr}}$: on a $T(\mathfrak{o}_{\bar{F}})=X_{*}(T)\otimes_{{\mathbb Z}}\mathfrak{o}_{\bar{F}}$ et $T(\mathfrak{o}_{F^{nr}})=T(\mathfrak{o}_{\bar{F}})^{I_{F}}$. De m\^eme pour les tores $T_{ad}$ etc... Notons $\pi:G_{SC}\to G_{AD}$ l'homomorphisme naturel. 

\begin{cor}{(i) On a l'égalité $G_{AD}(F^{nr})=\boldsymbol{\Omega}^{nr}T_{ad}(\mathfrak{o}_{F^{nr}})\pi(G_{SC}(F^{nr}))$.

(ii)  Soient $x,y\in \bar{C}^{nr}$. Supposons qu'il existe $g\in G_{SC}(F^{nr})$, resp. $g\in G_{AD}(F^{nr})$, tel que $gx=y$. Alors $x=y$, resp.  il existe $\omega\in \boldsymbol{\Omega}^{nr}$ tel que $\omega x=y$.

(iii) De l'inclusion $\boldsymbol{\Omega}^{nr}\to G_{AD}(F^{nr})$ se déduit un isomorphisme $H^1(\Gamma_{F}^{nr},\boldsymbol{\Omega}^{nr})\to H^1(\Gamma_{F}^{nr},G_{AD}(F^{nr}))=H^1(F,G_{AD})$.
} \end{cor}

{\bf Remarque.} Dans le membre de droite de l'égalité du (i), les trois ensembles sont permutables car $\pi(G_{SC}(F^{nr})$ et $T_{ad}(\mathfrak{o}_{F^{nr}})\pi(G_{SC}(F^{nr})$ sont des sous-groupes distingués de $G_{AD}(F^{nr})$.

\bigskip
Preuve du corollaire. Soit $g\in G_{AD}(F^{nr})$. Alors $g(C^{nr})$ est une alc\^ove de $Imm_{F^{nr}}(G_{AD})$. Ainsi qu'il est bien connu, il existe $h\in G_{SC}(F^{nr})$ tel que $g(C^{nr})= h(C^{nr})$. Posons $g'=\pi(h)^{-1}g$. Cet élément conserve $C^{nr}$. Il envoie $T^{nr}$ sur un élément de ${\cal T}_{max}^{nr}$ dont l'appartement associé contient $C^{nr}$. Un tel tore est conjugué de $T^{nr}$ par un élément de $K_{C^{nr},SC}^{0,nr}$. Quitte à modifier $h$ on peut donc supposer que $g'$ normalise $T^{nr}$. Alors $g'\in N^{nr}$ donc de $g$ se déduit une action sur ${\cal D}_{a}^{nr}$. Le lemme \ref{racinesaffines} dit que l'on peut trouver $\omega\in \boldsymbol{\Omega}^{nr}$ qui ait la  m\^eme action. En posant $t=\omega^{-1}g'$, $t$ est alors un élément de $N^{nr}$ dont l'action sur ${\cal D}_{a}^{nr}$ est triviale, donc dont l'action sur $App_{F^{nr}}(T^{nr})$ est triviale. Un tel élément appartient à $T_{ad}(\mathfrak{o}_{F^{nr}})$.  L'égalité $g=\pi(h)\omega t$ démontre le (i) de l'énoncé.

L'assertion de (ii) concernant le groupe $G_{SC}$ est bien connue. Supposons $g\in G_{AD}(F^{nr})$ et $gx=y$. On écrit comme ci-dessus $g=\pi(h)\omega t$. Alors $h\omega x=y$. On a $\omega(x)\in \bar{C}^{nr}$ par définition de $\boldsymbol{\Omega}^{nr}$. Appliquant l'assertion concernant le groupe $G_{SC}$, l'égalité $h(\omega x)=y$ implique $y=\omega x$. Cela prouve (ii). 

Reprenons les notations de \ref{alcove}. Considérons  le diagramme
$$ \begin{array}{ccc} H^1(\Gamma_{F}^{nr}, N^{nr})&\to H^1(\Gamma_{F}^{nr},G_{AD}(F^{nr}))\\ \downarrow&&\\ H^1(\Gamma_{F}^{nr},\bar{N}^{nr})&&\\ \end{array}$$
D'après \cite{W2} p. 441, les trois  flèches sont bijectives. D'après le lemme \ref{racinesaffines}, l'injection $\boldsymbol{\Omega}^{nr}\to N^{nr}$ se pousse en une bijection $\boldsymbol{\Omega}^{nr}\to \bar{N}^{nr}$. Donc de  l'injection $\boldsymbol{\Omega}^{nr}\to G_{AD}(F^{nr})$ se déduit un isomorphisme  $H^1(\Gamma_{F}^{nr},\boldsymbol{\Omega}^{nr})\to H^1(\Gamma_{F}^{nr},G_{AD}(F^{nr}))$. On a déjà dit que $H^1(\Gamma_{F}^{nr},G_{AD}(F^{nr}))=H^1(F,G_{AD})$. Cela prouve le (iii) de l'énoncé. $\square$

   \subsubsection{Discr\'etisation de l'immeuble; le groupe ${\bf G}$\label{discretisation}}

 Notons $F^{mod}$ la plus grande extension   mod\'er\'ement ramifi\'ee de $F$ contenue dans $\bar{F}$. On note $\Gamma_{F}^{mod}$, resp.   $I_{F}^{mod}$, le groupe de Galois de l'extension $F^{mod}/F$, resp.   $F^{mod}/F^{nr}$. 

Notons ${\mathbb P}$ l'ensemble des entiers $n\geq1$ premiers \`a $p$. Soit $\varpi_{F}$ une uniformisante de $F$. On sait que l'on peut fixer, et on fixe, une famille $(\varpi_{1/n})_{n\in {\mathbb P}}$ d'\'el\'ements de $\bar{F}^{\times}$ de sorte que les conditions suivantes soient v\'erifi\'ees:

$\varpi_{1}=\varpi_{F}$;

pour tous $n,m\in {\mathbb P}$, $\varpi_{1/nm}^m=\varpi_{1/n}$.

Alors $F^{mod}$ est l'extension de $F$ engendr\'ee par $F^{nr}$ et les $\varpi_{1/n}$ pour $n\in {\mathbb P}$. Le groupe $\Gamma_{F}^{nr}$ s'identifie au sous-groupe des \'elements $\sigma\in \Gamma_{F}^{mod}$ tels que $\sigma(\varpi_{1/n})=\varpi_{1/n}$ pour tout $n\in {\mathbb P}$. On obtient ainsi une d\'ecomposition en produit semi-direct $\Gamma_{F}^{mod}=I_{F}^{mod}\rtimes \Gamma_{F}^{nr}$. En particulier, le Frobenius se rel\`eve en  un \'el\'ement de $\Gamma_{F}^{mod}$.

Pour $n\in {\mathbb P}$, notons $\boldsymbol{\zeta}_{1/n}$, resp. $\bar{\boldsymbol{\zeta}}_{1/n}$, le groupe des racines de l'unit\'e dans $\bar{F}^{\times}$, resp. $\bar{{\mathbb F}}_{q}^{\times}$, d'ordre divisant $n$. On a $\cup_{n\in {\mathbb P}}\bar{\boldsymbol{\zeta}}_{1/n}=\bar{{\mathbb F}}_{q}^{\times}$. On pose $\boldsymbol{\zeta}=\cup_{n\in {\mathbb P}}\boldsymbol{\zeta}_{1/n}$.  Notons ${\mathbb Z}_{(p)}$ le localis\'e de ${\mathbb Z}$ en $p$, c'est-\`a-dire l'anneau des $r\in {\mathbb Q}$ tels qu'il existe $n\in {\mathbb P}$ de sorte que $rn\in {\mathbb Z}$. On sait que ${\mathbb Z}_{(p)}/{\mathbb Z}$ est isomorphe \`a $ \boldsymbol{\zeta}$. On fixe un homomorphisme $r\mapsto \zeta_{r}$ de ${\mathbb Z}_{(p)}$ sur $\boldsymbol{\zeta}$ qui se quotiente en un isomorphisme ${\mathbb Z}_{(p)}/{\mathbb Z}\simeq \boldsymbol{\zeta}$. On note $r\mapsto \bar{\zeta}_{r}$ son compos\'e avec l'application de r\'eduction de $\boldsymbol{\zeta}$ dans $\bar{{\mathbb F}}_{q}^{\times}$. Ce compos\'e est un isomorphisme ${\mathbb Z}_{(p)}/{\mathbb Z}\to \bar{{\mathbb F}}_{q}^{\times}$. Pour tout $\sigma\in I_{F}^{mod}$ et tout $n\in {\mathbb P}$, il existe un unique \'el\'ement de $\boldsymbol{\zeta}_{1/n}$ not\'e $\zeta(\sigma)_{1/n}$ tel que $\sigma(\varpi_{1/n})=\zeta(\sigma)^{-1}_{1/n}\varpi_{1/n}$  (on a choisi $\zeta(\sigma)^{-1}_{1/n}$ plut\^ot que $\zeta(\sigma)_{1/n}$ car cela simplifiera certaines formules ci-dessous). Pour $n,m\in {\mathbb P}$, on a $\zeta(\sigma)_{1/nm}^m=\zeta(\sigma)_{1/n}$. On note $\gamma$ l'unique \'el\'ement de $I_{F}^{mod}$ tel que $\zeta(\gamma)_{1/n}=\zeta_{1/n}$ pour tout $n\in {\mathbb P}$. L'\'el\'ement $\gamma$ est un prog\'en\'erateur de $I_{F}^{mod}$.

 Supposons $G$  quasi-d\'eploy\'e sur $F$.  Soit $\mathfrak{ E}=(B,T,(E_{\alpha})_{\alpha\in \Delta})$   une paire de Borel \'epingl\'ee de $G$ conserv\'ee par l'action galoisienne.  Le groupe $\Gamma_{F}$ agit sur $X_{*}(T)$ et  notre  hypothèse $(Hyp)_{1}(p)$, cf. \ref{leshypothesessurp},  implique que cette action est mod\'er\'ement ramifi\'ee, c'est-\`a-dire se quotiente en une action de $\Gamma_{F}^{mod}=\Gamma_{F^{mod}/F}$.
  Le sextuplet $(X^*(T),X_{*}(T),\Sigma,\check{\Sigma},\Delta,\check{\Delta})$ est une donn\'ee de racines munie d'une action de $\Gamma_{F}$. Comme on le sait, on peut introduire un groupe r\'eductif connexe ${\bf G}$ d\'efini sur $\bar{{\mathbb F}}_{q}$ qui a les m\^emes donn\'ees de racines. C'est-\`a-dire que l'on peut fixer une paire de Borel \'epingl\'ee $\boldsymbol{\mathfrak{E}}=({\bf B},{\bf T},({\bf E}_{\alpha})_{\alpha\in \Delta})$  de ${\bf G}$ et des isomorphismes compatibles $X_{*}({\bf T})\simeq X_{*}(T)$ et $X^*({\bf T})\simeq X^*(T)$ identifiant les racines et coracines des groupes ${\bf G}$ et $G$ ainsi que leurs racines et coracines simples. De l'homomorphisme de r\'eduction $\mathfrak{o}_{\bar{F}}^{\times}\to \bar{{\mathbb F}}_{q}^{\times}$ se d\'eduit un homomorphisme de $T(\mathfrak{o}_{\bar{F}})$ dans ${\bf T}$ que l'on note $t\mapsto \bar{t}$. 
 Les groupes de Weyl de $G$ relativement \`a $T$ et de ${\bf G}$ relativement \`a ${\bf T}$ sont isomorphes, on les note tous deux $W$.  L'action de $I_{F}$ sur la donn\'ee de racines devient une action par automorphismes alg\'ebriques de $I_{F}$ sur ${\bf G}$ pr\'eservant la paire de Borel \'epingl\'ee $\boldsymbol{\mathfrak{ E}}$. De l'action de $\Gamma_{F}^{nr}\subset \Gamma_{F}^{mod}$ sur $\mathfrak{E}$ se d\'eduit une action sur $\boldsymbol{\mathfrak{E}}$ qui munit ${\bf G}$ d'une structure de groupe sur ${\mathbb F}_{q}$. En g\'en\'eral, les automorphismes d\'efinissant l'action de $I_{F}$ ne sont pas d\'efinis sur ${\mathbb F}_{q}$. 
 
 {\bf Remarque.} La structure de ${\bf G}$ sur ${\mathbb F}_{q}$ d\'epend en g\'en\'eral du choix  de notre famille $(\varpi_{1/n})_{n\in {\mathbb P}}$ qui nous a permis de relever $\Gamma_{F}^{nr}$ en un sous-groupe de $\Gamma_{F}^{mod}$. Toutefois, si $G$ est d\'eploy\'e sur $F^{nr}$, cette structure est canonique. 
 \bigskip

 Posons ${\mathbb T}=X_{*}(T)\otimes_{{\mathbb Z}}{\mathbb Z}_{(p)}$. Le groupe $\Gamma_{F}$ agit sur ${\mathbb T}$ via son action sur $X_{*}(T)$.  On d\'efinit deux homomorphismes
 $$v:T(F^{mod})\to  {\mathbb T},\,\, \underline{v}: {\mathbb T}\to T(F^{mod})$$
 de la fa\c{c}on suivante. On a $T(F^{mod})=X_{*}(T)\otimes F^{mod,\times}$ et $v$ est le produit tensoriel de l'identit\'e de $X_{*}(T)$ et de la valuation $val_{F}:F^{mod,\times}\to {\mathbb Z}_{(p)}$. Pour $y\in {\mathbb T}$, on choisit $n\in {\mathbb P}$ tel que $ny\in X_{*}(T)$ et on d\'efinit $\underline{v}(y)=(ny)(\varpi_{1/n})$. Cela ne d\'epend pas du choix de $n$. L'application $v\circ \underline{v}$ est l'identit\'e de $ {\mathbb T}$. L'application $v$ est \'equivariante pour les actions de $I_{F}$ mais $\underline{v}$ ne l'est pas. 
 
 On d\'efinit un homomorphisme $j_{T}:{\mathbb T}\to {\bf T}$ qui est le produit tensoriel de l'isomorphisme $X_{*}(T)\simeq X_{*}({\bf T})$ et de l'homomorphisme $r\mapsto \zeta_{r}$ de ${\mathbb Z}_{p}$ sur $\bar{{\mathbb F}}_{q}^{\times}$. Cet homomorphisme $j_{T}$ se quotiente en un isomorphisme ${\mathbb T}/X_{*}(T)\simeq {\bf T}$.

A tout sous-tore $T'$ de $T$ sont associ\'es un sous-tore ${\bf T}'$ de ${\bf T}$ et un sous-groupe ${\mathbb T}'$ de ${\mathbb T}$: $X_{*}({\bf T}')$ correspond \`a $X_{*}(T')$ et ${\mathbb T}'=X_{*}(T')\otimes_{{\mathbb Z}}{\mathbb Z}_{(p)}$. Pour toute extension $F'$ de $F$ contenue dans $F^{mod}$, on note $T_{F'}$ le plus grand sous-tore de $T$ d\'eploy\'e sur $F'$ et $A_{G,F'}$ le plus grand tore central dans $G$ d\'eploy\'e sur $F'$. Si $F'=F^{nr}$, on note plut\^ot ces tores $T^{nr}$ et $A_{G}^{nr}$. Soit $F'$ une extension finie de $F^{nr}$ contenue dans $F^{mod}$. Au tore $T$ est associ\'e un appartement $App_{F'}(T_{F'})$ dans $Imm_{F'}(G_{AD})$, qui s'identifie  \`a $X_{*,{\mathbb R}}(T_{F'}) /X_{*,{\mathbb R}}(A_{G,F'})$. Pour \'eviter les confusions, on note $s_{0}$ le point $0$ de ce dernier groupe. Le groupe $\Gamma_{F'/F}$ agit sur $App_{F'}(T_{F'})$ et sur $X_{*,{\mathbb R}}(T_{F'})/X_{*,{\mathbb R}}(A_{G,F'})$ et l'identification est \'equivariante pour ces actions. Le groupe ${\mathbb T}_{F'}$ est contenu dans $X_{*}(T_{F'})\otimes_{{\mathbb Z}}{\mathbb R}$ et il s'en d\'eduit un homomorphisme $a_{T}:{\mathbb T}_{F'}\to App_{F'}(T_{F'})$ qui se quotiente en une injection ${\mathbb T}_{F'}/{\mathbb A}_{G,F'}\to App_{F'}(T_{F'})$ d'image dense. On v\'erifie que, si $x,x'\in App_{F'}(T_{F'})$ et s'il existe $g\in G(F')$ tel que $gx=x'$, alors $x$ appartient \`a l'image de $a_{T}$ si et seulement s'il en est de m\^eme de $x'$. Pour $t\in T(F')$, on a l'\'egalit\'e $t^{-1}s_{0}=a_{T}\circ v(t)$ (remarquons que l'\'equivariance de l'homomorphisme $v$ implique que $v(t)$ appartient \`a ${\mathbb T}_{F'}$). 

Soit $y\in {\mathbb T}^{nr}$. Fixons $n\in {\mathbb P}$ tel que $ny\in X_{*}(T^{nr})$. Pour $\sigma\in I_{F}$, on pose $t_{y}(\sigma)=(ny)(\zeta(\sigma)_{1/n})$. C'est un \'el\'ement de $T^{nr}(\mathfrak{o}_{F^{nr}})$ dont on note $\bar{t}_{y}(\sigma)$ la r\'eduction dans ${\bf T}^{nr}$. On a $\underline{v}(y)\in T^{nr}(F^{mod})$ et on v\'erifie l'\'egalit\'e $\underline{v}(y)\sigma(\underline{v}(y))^{-1}=t_{y}(\sigma)$. Pour $\sigma=\gamma$, on a $\bar{t}_{y}(\gamma)=j_{T}(y)$.

Soit $F'$ une extension finie de $F^{nr}$ contenue dans $F^{mod}$ et supposons $G$ d\'eploy\'e sur $F'$. Alors le groupe r\'esiduel $G_{s_{0},F'}$ s'identifie \`a ${\bf G}$. Pr\'ecis\'ement, de la paire de Borel \'epingl\'ee $\mathfrak{E}$ se d\'eduit une telle paire $\mathfrak{E}_{s_{0},F'}=(B_{s_{0},F'},T_{s_{0},F'},(E_{\alpha,s_{0},F'})_{\alpha\in \Delta})$ de $G_{s_{0},F'}$. L'élément $E_{\alpha,s_{0},F'}$ est la r\'eduction naturelle de $E_{\alpha}\in \mathfrak{k}_{s_{0},F'}$. On peut identifier  et on identifie $G_{s_{0},F'}$ \`a ${\bf G}$ de sorte que $\mathfrak{ E}_{s_{0},F'}$ s'identifie \`a $\boldsymbol{\mathfrak{E}}$, l'identification \'etant compatible aux isomorphismes $X_{*}(T_{s_{0},F'})\simeq X_{*}(T)\simeq X_{*}({\bf T})$. Dans le cas o\`u $G$ est non ramifi\'e, on peut choisir $F'=F^{nr}$ et $G_{s_{0}}$ s'identifie \`a ${\bf G}$. 

On a suppos\'e $G$ quasi-d\'eploy\'e sur $F$. On peut lever cette restriction en consid\'erant que $G$ est d\'efini sur $F^{nr}$ (ou sur une extension finie non ramifi\'ee de $F$ sur laquelle $G$ est quasi-d\'eploy\'e). Les m\^emes constructions s'appliquent, en oubliant les actions de $\Gamma_{F}^{nr}$. C'est-\`a-dire que ${\bf G}$ est seulement d\'efini sur $\bar{{\mathbb F}}_{q}$, pas sur ${\mathbb F}_{q}$.

\subsubsection{Description des groupes $G_{{\cal F}}$\label{lesgroupesGF}}
 Soit $x\in a_{T}({\mathbb T}^{nr})\subset App_{F^{nr}}(T^{nr})$. Fixons $\tilde{x}\in {\mathbb T}^{nr}$ tel que $a_{T}(\tilde{x})=x$ et  une extension finie $F'$ de $F^{nr}$ contenue dans $F^{mod}$ telle que $G$ soit d\'eploy\'e sur $F'$ et que $\underline{v}(\tilde{x})$ appartienne \`a $T^{nr}(F')$. On a $x=a_{T}(\tilde{x})=a_{T}\circ v\circ\underline{v}(\tilde{x})=\underline{v}(\tilde{x})^{-1}s_{0}$. Alors $Ad(\underline{v}(\tilde{x}))$ se restreint en un isomorphisme de $\mathfrak{k}_{x,F'}$ sur $\mathfrak{k}_{s_{0},F'}$ qui entrelace l'action naturelle de $I_{F}$ sur $\mathfrak{k}_{x,F'}$ et l'action $\sigma\mapsto Ad(t_{\tilde{x}}(\sigma))\circ\sigma$ sur $\mathfrak{k}_{s_{0},F'}$. Donc cet isomorphisme se restreint en un isomorphisme de $\mathfrak{k}_{x}^{nr}$ sur le sous-ensemble des points fixes dans $\mathfrak{k}_{s_{0},F'}$ par l'action pr\'ec\'edente. On a identifi\'e l'espace r\'eduit $\mathfrak{g}_{s_{0},F'}$ \`a $\boldsymbol{\mathfrak{g}}$. Par r\'eduction, $\mathfrak{g}_{x}$ s'identifie au sous-espace des points fixes dans $\boldsymbol{\mathfrak{g}}$ par l'action $ \sigma\mapsto Ad(\bar{t}_{\tilde{x}}(\sigma))\circ\sigma$ de $I_{F}$. Le groupe $G_{x}$ s'identifie au sous-groupe de ${\bf G}$ d\'etermin\'e par cette alg\`ebre de Lie, c'est-\`a-dire la composante neutre du sous-espace des points fixes dans ${\bf G}$ par l'action pr\'ec\'edente. Cette description ne d\'epend pas du choix du rel\`evement $\tilde{x}$ de $x$: on ne peut changer $\tilde{x}$ que par un \'el\'ement central, cela ne modifie pas  $Ad(\underline{v}(\tilde{x}))$.

Remarquons que l'on peut aussi bien remplacer ci-dessus l'action de $I_{F}$ par celle du g\'en\'erateur $\gamma$. On a $\bar{t}_{\tilde{x}}(\gamma)=j_{T}(\tilde{x})$ et l'action de $\gamma$ est  $Ad(j_{T}(\tilde{x}))\circ\gamma$.  Dans le cas o\`u $G$ est d\'eploy\'e sur $F^{nr}$, l'action de $\gamma$ sur ${\bf G}$ est triviale et on obtient que l'image de $G_{x}$ dans ${\bf G}$ est ${\bf G}_{j_{T}(\tilde{x})}$ (dont on rappelle qu'il s'agit de $Z_{{\bf G}}(j_{T}(\tilde{x}))^0$). 

 Soit ${\cal F}$ une facette de $App_{F^{nr}}(T^{nr})$.  L'intersection ${\cal F}\cap a_{T}({\mathbb T}^{nr})$ est non vide. On choisit $x$ dans cette intersection. On a $\mathfrak{k}_{{\cal F}}^{nr}=\mathfrak{k}_{x}^{nr}$ et la construction pr\'ec\'edente fournit des homomorphismes   $ \mathfrak{g}_{{\cal F}}\to \boldsymbol{ \mathfrak{g}}$ et $G_{{\cal F}}\to {\bf G}$ que l'on  note tous deux $\iota_{x,{\cal F}}$.

\begin{lem}{Soit ${\cal F}$ une facette de $App_{F^{nr}}(T^{nr})$. Les homomorphismes $\iota_{x,{\cal F}}$ ne d\'ependent pas du choix de    $x\in {\cal F}\cap a_{T}({\mathbb T}^{nr})$.}\end{lem}

 Preuve. Soient $x,y\in {\cal F}\cap a_{T}({\mathbb T}^{nr})$. On rel\`eve $x$ et $y$ en des \'el\'ements $\tilde{x}$ et $\tilde{y}$ de ${\mathbb T}$. 
  On peut fixer un \'el\'ement $n\in {\mathbb P}$ commun tel que $n\tilde{x},n\tilde{y}\in X_{*}(T^{nr})$. Fixons une extension finie $F'$ de $F^{nr}$ contenue dans $F^{mod}$, contenant   $\varpi_{1/n}$ et telle que $G$ soit d\'eploy\'e sur $F'$. 
  Notons $e'$ l'indice de ramification de $F'/F^{nr}$. Remarquons que $n$ divise $e'$.  
  
  Parce que $x\in {\cal F}$, on a les \'egalit\'es $\mathfrak{k}_{{\cal F}}^{nr}=\mathfrak{k}_{x}^{nr}$, $\mathfrak{k}_{{\cal F}}^{+,nr}=\mathfrak{k}_{x}^{+,nr}$, $K_{{\cal F}}^{0,nr}=K_{x}^{0,nr}$ et $K_{{\cal F}}^{+,nr}=K_{x}^{+,nr}$. Les homomorphismes $\iota_{x,{\cal F}}$ sont par d\'efinition les compos\'es
  $$\mathfrak{g}_{{\cal F}}=\mathfrak{k}_{{\cal F}}^{nr}/ \mathfrak{k}_{{\cal F}}^{+,nr}=\mathfrak{k}_{x}^{nr}/ \mathfrak{k}_{x}^{+,nr}\subset \mathfrak{k}_{x,F'}/\mathfrak{p}_{F'}\mathfrak{k}_{x,F'}=\mathfrak{g}_{x,F'}\stackrel{\tau(\tilde{x})}{\to }\boldsymbol{ \mathfrak{g}}$$
  et
  $$G_{{\cal F}}=K_{{\cal F}}^{0,nr}/K_{{\cal F}}^{+,nr}=K_{x}^{0,nr}/K_{x}^{+,nr}\subset K_{x,F'}^0/K_{x,F'}^+=G_{x,F'}\stackrel{\tau(\tilde{x})}{\to }{\bf G},$$
  o\`u on a not\'e $\tau(\tilde{x})$ l'isomorphisme d\'eduit par r\'eduction de $Ad(\underline{v}(\tilde{x}))$. 
  On a une description analogue des homomorphismes $\iota_{y,{\cal F}}$. 
  
   Posons $t= \underline{v}(\tilde{y})^{-1} \underline{v}(\tilde{x})$. L'action de $t$ sur l'immeuble $Imm_{F'}(G_{AD})$ envoie $x$ sur $y$. Donc $Ad(t)$ envoie $\mathfrak{k}_{x,F'}$ sur $\mathfrak{k}_{y,F'}$ et  $K_{x,F'}^0$ sur $K_{y,F'}^0$ et on en d\'eduit par r\'eduction des isomorphismes $\mathfrak{g}_{x,F'}\to \mathfrak{g}_{y,F'}$ et $G_{x,F'}\to G_{y,F'}$ que l'on note tous deux $\tau$. Les diagrammes suivants sont commutatifs:
   
   $$\begin{array}{ccc}\mathfrak{g}_{x,F'}&&\\ &\searrow \tau(\tilde{x})&\\ \downarrow \tau&&\boldsymbol{\mathfrak{g}}\\ &\nearrow  \tau(\tilde{y})&\\ \mathfrak{g}_{y,F'}&&\\ \end{array}$$
    $$\begin{array}{ccc}G_{x,F'}&&\\ &\searrow \tau(\tilde{x})&\\ \downarrow \tau&&{\bf G}\\ &\nearrow  \tau(\tilde{y})&\\ G_{y,F''}&&\\ \end{array}$$
    
    Il nous suffit donc de prouver que les diagrammes suivants sont eux-aussi commutatifs:
    $$(1) \qquad \begin{array}{ccc}\mathfrak{g}_{x}&\to &\mathfrak{g}_{x,F'}\\ \parallel&&\\ \mathfrak{g}_{{\cal F}}&&\downarrow \tau\\ \parallel&&\\ \mathfrak{g}_{y}&\to&\mathfrak{g}_{y,F'}\\ \end{array}$$
    $$(2) \qquad \begin{array}{ccc}G_{x}&\to&G_{x,F'}\\ \parallel&&\\ G_{{\cal F}}&&\downarrow \tau\\ \parallel&&\\G_{y} &\to &G_{y,F'}\\ \end{array}$$

    Pour $\alpha\in \Sigma$ et $r\in \frac{1}{e'}{\mathbb Z}$, notons $\mathfrak{u}_{\alpha,r}$ le $\mathfrak{o}_{F'}$-module $\mathfrak{u}_{\alpha,\mathfrak{p}_{F'}^{e'r}}$. 
   Le $\mathfrak{o}_{F'}$-module $\mathfrak{k}_{x,F'}$ est la somme de $\mathfrak{t}(\mathfrak{o}_{F'})$ et des  $\mathfrak{u}_{\alpha,-\alpha(x)}$ pour $\alpha\in \Sigma$ (remarquons que par construction, $\alpha(x)\in \frac{1}{n}{\mathbb Z}\subset \frac{1}{e'}{\mathbb Z}$).  
  Pour $\beta\in \Sigma^{nr}$, notons  $\Sigma[\beta]$ l'ensemble des $\alpha\in \Sigma$ telles que $\alpha^{res}=\beta$. C'est une orbite pour l'action de $I_{F}$.  Notons $\mathfrak{u}_{\beta}$ la somme des $\mathfrak{u}_{\alpha}$ sur les $\alpha\in \Sigma[\beta]$ et $\mathfrak{u}_{\beta,r}$ la somme des $\mathfrak{u}_{\alpha,r}$.   Ce module est conservé par l'action de $I_{F}$.   Puisque $x\in App_{F^{nr}}(T^{nr})$, $\beta(x)$ est bien défini et on a $\alpha(x)=\beta(x)$  pour tout  $\alpha\in \Sigma[\beta]$. Le module $\mathfrak{k}_{x}^{nr}$ est la somme de $\mathfrak{t}(\mathfrak{o}_{F'})^{I_{F}}$ et des $\mathfrak{u}_{\beta,-\beta(x)}^{I_{F}}$ pour $\beta\in \Sigma^{nr}$. De m\^eme, $\mathfrak{k}_{x}^{+,nr}$  est la somme de $\mathfrak{t}(\mathfrak{p}_{F'})^{I_{F}}$ et des $\mathfrak{u}_{\beta,-\beta(x)+1/e'}^{I_{F}}$ pour $\beta\in \Sigma^{nr}$. Notons $\Sigma^{nr}(x)$ l'ensemble des $\beta\in \Sigma$ tels que $\mathfrak{u}_{\beta,-\beta(x)}^{I_{F}}\not=\mathfrak{u}_{\beta,-\beta(x)+1/e'}^{I_{F}}$. Alors $\mathfrak{g}_{x}^{nr}$ est la somme des images par r\'eduction de $\mathfrak{t}(\mathfrak{o}_{F'})^{I_{F}}$ et des $\mathfrak{u}_{\beta,-\beta(x)}^{I_{F}}$ pour $\beta\in \Sigma^{nr}(x)$.

   On a des descriptions analogues de $\mathfrak{k}_{y,F'}$, $\mathfrak{k}_{y}^{nr}$ et $\mathfrak{k}_{y}^{+,nr}$. Puisque
  $\mathfrak{k}_{x}^{nr}=\mathfrak{k}_{y}^{nr}$ et $\mathfrak{k}_{x}^{+,nr}=\mathfrak{k}_{y}^{+,nr}$, on a les \'egalit\'es  $\mathfrak{u}_{\beta,-\beta(x)}^{I_{F}}=\mathfrak{u}_{\beta,-\beta(y)}^{I_{F}}$ et $\mathfrak{u}_{\beta,-\beta(x)+1/e'}^{I_{F}}=\mathfrak{u}_{\beta,-\beta(y)+1/e'}^{I_{F}}$ pour tout $\beta\in \Sigma^{nr}$. Les \'egalit\'es pr\'ec\'edentes entra\^{\i}nent $\Sigma^{nr}(x)=\Sigma^{nr}(y)$. 
 Pour $r\in \frac{1}{e'}{\mathbb Z}$, l'ensemble $\mathfrak{u}_{\beta,r}^{I_{F}}$ ne d\'etermine pas $r$. Mais, si $\mathfrak{u}_{\beta,r}^{I_{F}}\not=\mathfrak{u}_{\beta,r+1/e'}^{I_{F}}$, l'ensemble $\mathfrak{u}_{\beta,r}^{I_{F}}$ d\'etermine $r$: c'est le plus grand \'el\'ement $r'\in\frac{1}{e'}{\mathbb Z}$ tel qu'un \'el\'ement de $\mathfrak{u}_{\beta,r}^{I_{F}}-\mathfrak{u}_{\beta,r+1/e'}^{I_{F}}$ appartienne \`a $\mathfrak{u}_{\beta,r'}$. Il en r\'esulte que, pour $\beta\in \Sigma^{nr}(x)=\Sigma^{nr}(y)$, on a $\beta(x)=\beta(y)$. 

Il r\'esulte de ces descriptions que, pour prouver la commutativit\'e du diagramme (1), il suffit de prouver que $Ad(t)$ se restreint en l'identit\'e de $\mathfrak{t}(\mathfrak{o}_{F'})^{I_{F}}$ et des $\mathfrak{u}_{\beta,-\beta(x)}^{I_{F}}$ pour $\beta\in \Sigma^{nr}(x)=\Sigma^{nr}(y)$. Puisque $t\in T^{nr}(F')$, $Ad(t)$ est bien l'identit\'e sur $\mathfrak{t}(\mathfrak{o}_{F'})^{I_{F}}$. Soit $\beta\in \Sigma^{nr}(x)$. Puisque $t\in T^{nr}(F')$ et que $\beta\in \Sigma^{nr}\subset X^*(T^{nr})$, $\beta(t)$ est bien défini et on a $\alpha(t)=\beta(t)$  pour tout  $\alpha\in \Sigma[\beta]$. L'automorphisme $Ad(t)$ agit sur $\mathfrak{u}_{\beta}$ par multiplication par $\beta(t)$. Par d\'efinition de $t$, on a $\beta(t)=\varpi_{1/n}^{\beta(n\tilde{x}-n\tilde{y})}=\varpi_{1/n}^{\beta(nx-ny)}$. Mais $\beta\in \Sigma^{nr}(x)=\Sigma^{nr}(y)$ donc $\beta(x)=\beta(y)$. Donc $\beta(t)=1$ et $Ad(t)$ est bien l'identit\'e sur $\mathfrak{u}_{\beta,-\beta(x)}^{I_{F}}$. Cela prouve la commutativit\'e du diagramme (1).

 Le groupe $G_{{\cal F}}$ est engendr\'e par son sous-tore maximal image par r\'eduction de $T^{nr}(\mathfrak{o}_{F^{nr}})$ et par les exponentielles  des \'el\'ements nilpotents de $\mathfrak{g}_{{\cal F}}$. Quand on se restreint au tore maximal, le diagramme (2) est commutatif pour la m\^eme raison que ci-dessus: $t$ appartient \`a $T(F')$. Pour les exponentielles des \'el\'ements nilpotents de $\mathfrak{g}_{{\cal F}}$, l'identit\'e requise r\'esulte de la commutativit\'e du diagramme (1). Cela d\'emontre la commutativit\'e de (2) et prouve le lemme. $\square$

Ce lemme nous autorise \`a noter $\iota_{{\cal F}}$ l'homomorphisme $\iota_{x,{\cal F}}$ pour un \'el\'ement quelconque $x\in {\cal F}\cap a_{T}({\mathbb T}^{nr})$. 

Supposons que $G$ soit quasi-d\'eploy\'e. Soit ${\cal F}$ une facette de $App(T_{F})$. Elle est contenue dans une unique facette ${\cal F}^{nr}$ de $App_{F^{nr}}(T^{nr})$. On a $\mathfrak{g}_{{\cal F}}=\mathfrak{g}_{{\cal F}^{nr}}$ et $G_{{\cal F}}=G_{{\cal F}^{nr}}$. 
On vient de définir les homomorphismes $\iota_{{\cal F}^{nr}}$. Par les égalités précédentes, ils deviennent des homomorphismes notés $\iota_{{\cal F}}:\mathfrak{g}_{{\cal F}}\to \boldsymbol{\mathfrak{g}}$ et $\iota_{{\cal F}}:G_{{\cal F}}\to {\bf G}$. Montrons que 

(3) ces homomorphismes sont équivariants pour les actions de $\Gamma_{F}^{nr}$.

Dans la construction de $\iota_{{\cal F}^{nr}}$, on peut choisir pour élément $x$ un élément de $ a_{T}({\mathbb T}_{F})\cap {\cal F}$. L'\'el\'ement $\underline{v}(\tilde{x})$ est alors  fix\'e par $\Gamma_{F}^{nr}$. L'assertion (3) s'en déduit.

\subsubsection{Adh\'erence d'une facette\label{adherence}}
Soient ${\cal F}$ et ${\cal F}'$ deux facettes de $App_{F^{nr}}(T^{nr})$.  Supposons ${\cal F}'\subset \bar{{\cal F}}$. On a alors l'inclusion $\mathfrak{k}_{{\cal F}}^{nr}\subset \mathfrak{k}_{{\cal F}'}^{nr}$. La r\'eduction de $\mathfrak{k}_{{\cal F}}^{nr}$ dans $\mathfrak{g}_{{\cal F}'}$ est une sous-alg\`ebre parabolique $\mathfrak{p}_{{\cal F}}$ et $\mathfrak{g}_{{\cal F}}$ est canoniquement son quotient de Levi, c'est-\`a-dire le quotient de $\mathfrak{p}_{{\cal F}}$ par son radical nilpotent. On peut l'identifier  \`a la sous-alg\`ebre de Levi de $\mathfrak{p}_{{\cal F}}$ contenant   l'alg\`ebre de Lie $\mathfrak{t}_{{\cal F}'}$ du tore maximal $T_{{\cal F}'} $ de $G_{{\cal F}'}$ d\'eduit par r\'eduction de $T^{nr}$.  On a alors un diagramme
 $$\begin{array}{ccc}\mathfrak{g}_{{\cal F}}&&\\ &\searrow \iota_{{\cal F}}&
\\ \downarrow&&\bar{\mathfrak{g}}\\&\nearrow \iota_{{\cal F}'}&\\\mathfrak{g}_{{\cal F}'}&&\\ \end{array}$$

\begin{lem}{Le diagramme ci-dessus est commutatif.}\end{lem} 

Preuve. On fixe $x\in {\cal F}\cap a_{T}({\mathbb T}^{nr})$ et $y\in {\cal F}'\cap a_{T}({\mathbb T}^{nr})$. On reprend les constructions de la preuve pr\'ec\'edente. On est encore ramen\'e \`a prouver que $Ad(t)$ se restreint en l'identit\'e de $\mathfrak{u}_{\beta}$ pour $\beta\in \Sigma^{nr}(x)$. On n'a plus $\Sigma^{nr}(x)=\Sigma^{nr}(y)$ mais tout-de-m\^eme  l'inclusion $\Sigma^{nr}(x)\subset \Sigma^{nr}(y)$. Pour $\beta\in \Sigma^{nr}(x)$, on a encore $\beta(x)=\beta(y)$. Cela suffit pour conclure. $\square$

\subsubsection{Action de $G(F^{nr})$\label{conjugaison}}
Soient ${\cal F}$ une facette de $App_{F^{nr}}(T^{nr})$ et $g\in G(F^{nr})$. Notons ${\cal F}'=g{\cal F}$ l'image de ${\cal F}$ par l'action de $g$ sur $Imm_{F^{nr}}(G_{AD})$ et supposons ${\cal F}'\subset App_{F^{nr}}(T^{nr})$. De $Ad(g)$ se d\'eduisent  des isomorphismes de $\mathfrak{g}_{{\cal F}}$ sur $\mathfrak{g}_{{\cal F}'}$ et de $G_{{\cal F}}$ sur $G_{{\cal F}'}$ que l'on note tous deux $\overline{Ad(g)}$. Fixons $x\in {\cal F}\cap a_{T}({\mathbb T}^{nr})$, posons $x'=gx$. On a $x'\in {\cal F}'\cap a_{T}({\mathbb T}^{nr})$. Fixons des \'el\'ements $\tilde{x}$ et $\tilde{x}'$  dans ${\mathbb T}^{nr}$ tels que $a_{T}(\tilde{x})=x$ et $a_{T}(\tilde{x}')=x'$.

\begin{lem}{Il existe $z\in Z({\bf G})$ et  $u\in {\bf G}$ de sorte que
  
  (i) les diagrammes suivants soient commutatifs:
    
   $$\begin{array}{ccccc}&G_{{\cal F}}&\stackrel{\iota_{{\cal F}}}{\to}&{\bf G}&\\   \overline{Ad(g)} &\downarrow&&\downarrow &Ad(u)\\     &G_{{\cal F}'}&\stackrel{\iota_{{\cal F}'}}{\to}&{\bf G}&\\ \end{array}$$
     $$\begin{array}{ccccc}&\mathfrak{g}_{{\cal F}}&\stackrel{\iota_{{\cal F}}}{\to}&\boldsymbol{\mathfrak{g}}&\\  \overline{Ad(g)} &\downarrow &&\downarrow &Ad(u)\\     &\mathfrak{g}_{{\cal F}'}&\stackrel{\iota_{{\cal F}'}}{\to}&\boldsymbol{\mathfrak{g}}&\\ \end{array};$$

    (ii) on ait l'\'egalit\'e $zuj_{T}(\tilde{x})\gamma(u)^{-1}=j_{T}(\tilde{x}')$. }\end{lem}
    
    Preuve.   Puisque que $x$ et $x'$ appartiennent \`a $App_{F^{nr}}(T^{nr})$ et qu'il existe $g\in G(F^{nr})$ de sorte que $g x=x'$, on sait que
    $g\in K_{x'}^{0,nr}Norm_{G(F^{nr})}(T)$. Choisissons $g'\in Norm_{G(F^{nr})}(T)$ et $k_{1}\in K_{x'}^{0,nr}$ tels que $g=k_{1}g'$. On a encore $g'x=x'$.    L'image de $g'$ dans $W$ est fixe par $I_{F}$.   Le groupe de Weyl de $G_{s_{0}}$ est $W^{I_{F}}$. Ce groupe se rel\`eve en un sous-groupe de $K_{s_{0}}^{0,nr}\cap Norm_{G(F^{nr})}(T)$. On peut donc fixer un \'el\'ement $w\in K_{s_{0}}^{0,nr}\cap Norm_{G(F^{nr})}(T)$ qui a m\^eme image que $g'$ dans $W$. En posant $\tau=g'w^{-1}$, on a alors $\tau\in T(F^{nr})$. Posons $t=\underline{v}(\tilde{x})$ et $t'=\underline{v}(\tilde{x}')$.   Fixons une extension finie $F'$ de $F^{nr}$ contenue dans $F^{mod}$ de sorte que $t,t'\in T^{nr}(F')$ et que $G$ soit d\'eploy\'e sur $F'$.  Ainsi qu'on l'a d\'ej\`a dit, on a $x=t^{-1}s_{0}$ et $x'={t'}^{-1}s_{0}$. L'\'egalit\'e $g'x=x'$ devient $\tau wt^{-1} s_{0}={t'}^{-1}s_{0}$, c'est-\`a-dire $t'\tau wt^{-1}\in K_{s_{0},F'}^{\dag}$. 
    Puisque $w$ appartient \`a $K_{s_{0}}^{0,nr}$, il fixe $s_{0}$ et on a aussi bien $t'\tau wt^{-1}w^{-1} \in K_{s_{0},F'}^{\dag}$. Remarquons que $t'\tau wt^{-1}w^{-1}\in T(F')$.
 Puisque $s_{0}$ est hypersp\'ecial dans $Imm_{F'}(G_{AD})$, on a l'\'egalit\'e $K_{s_{0},F'}^{\dag}=Z(G)(F')K_{s_{0},F'}^0$. L'intersection de ce groupe avec $T(F')$ est $Z(G)(F')T(\mathfrak{o}_{F'})$. On peut donc fixer $z_{1}\in Z(G)(F')$ et $u_{1}\in T(\mathfrak{o}_{F'})$ de sorte que $t'\tau wt^{-1}w^{-1}=z_{1}u_{1}$, ou encore $t'=z_{1}u_{1}wtw^{-1}\tau^{-1}$. Remarquons que $z_{1}$, $u_{1}$, $wtw^{-1}$ et $\tau$ appartiennent tous \`a $T(F')$ donc commutent. Par d\'efinition, on a $t_{\tilde{x}'}(\gamma)=t'\gamma(t')^{-1}$, $t_{\tilde{x}}(\gamma)=t\gamma(t)^{-1}$. On utilise l'\'egalit\'e pr\'ec\'edente en se rappelant que $w$ et $\tau$ sont fixes par $I_{F}$. On obtient une \'egalit\'e que l'on peut \'ecrire au choix sous les formes
 $$t_{\tilde{x}'}(\gamma)=z_{1}\gamma(z_{1})^{-1}u_{1} wt_{\tilde{x}}(\gamma)\gamma(w)^{-1}\gamma(u_{1})^{-1},$$
  $$t_{\tilde{x}'}(\gamma)=z_{1}\gamma(z_{1})^{-1}u_{1} \gamma(u_{1})^{-1}wt_{\tilde{x}}(\gamma)w^{-1}.$$
     Les \'el\'ements   $t_{\tilde{x}}(\gamma)$ et $t_{\tilde{x}'}(\gamma)$ appartiennent \`a $T^{nr}(\mathfrak{o}_{F'})$. On a aussi $u_{1}\gamma(u_{1})^{-1}\in T(\mathfrak{o}_{F'})$. Il en r\'esulte que $z_{1}\gamma(z_{1})^{-1}$ appartient aussi \`a $T(\mathfrak{o}_{F'})$. Posons $k_{2}=t'k_{1}t^{'-1}$. Puisque $k_{1}\in K_{x'}^{0,nr}$, $k_{2}$ appartient \`a $K_{s_{0},F'}^0$ et est fixe par $Ad(t_{\tilde{x}'}(\gamma))\circ\gamma$. On a donc $t_{\tilde{x}'}(\gamma)=k_{2}t_{\tilde{x}'}(\gamma)\gamma(k_{2})^{-1}$ et on obtient
    $$t_{\tilde{x}'}(\gamma)=z_{1}\gamma(z_{1})^{-1}u_{2} t_{\tilde{x}}(\gamma) \gamma(u_{2})^{-1},$$
    o\`u $u_{2}=k_{2}u_{1}w$.   
     On r\'eduit cette \'egalit\'e  dans ${\bf G}$. Les r\'eductions de $t_{\tilde{x}'}(\gamma)$ et $t_{\tilde{x}}(\gamma)$ sont $j_{T}(\tilde{x}')$ et $j_{T}(\tilde{x})$. Il est clair que $z_{1}\gamma(z_{1})^{-1}$ se r\'eduit en un \'el\'ement $z\in Z({\bf G})$. On note $u$ la r\'eduction de $u_{2}$. On obtient alors $j_{T}(\tilde{x}')=zuj_{T}(\tilde{x})\gamma(u)^{-1}$, ce qui est l'assertion (ii) de l'\'enonc\'e. 
   
   L'application $\iota_{{\cal F}'}\circ \overline{Ad(g)} $ se d\'eduit par r\'eduction de l'application $Ad(t'k_{1}g')$. On a $t'k_{1}g'=k_{2}t'\tau w=z_{1}k_{2}u_{1}wt=z_{1}u_{2}t$. Puisque $z_{1}$ est central, on a $Ad(t'k_{1}g')=Ad(u_{2})Ad(t)$. Mais l'application qui se d\'eduit en r\'eduction de cet op\'erateur est pr\'ecis\'ement $Ad(u)\circ \iota_{{\cal F}}$. D'o\`u l'assertion (i). $\square$
   
   \subsection{L'espace $FC(\mathfrak{g}(F))$}

    \subsubsection{Fonctions sur les alg\`ebres de Lie\label{fonctionssurlesalgebresdeLie}}
 Pour les cinq premiers paragraphes de cette sous-section, on consid\`ere un groupe $G$ r\'eductif connexe d\'efini sur ${\mathbb F}_{q}$.  On  impose comme toujours l'hypoth\`ese $(Hyp)_{1}(p)$, cf. \ref{leshypothesessurp}.

  On note  ${\bf C}(\mathfrak{g}({\mathbb F}_{q}))$  l'espace des fonctions sur $\mathfrak{g}({\mathbb F}_{q})$ \`a valeurs complexes et $C(\mathfrak{g}({\mathbb F}_{q}))$  le sous-espace des fonctions invariantes par conjugaison par $G({\mathbb F}_{q})$ (on notera aussi $C^G(\mathfrak{g}({\mathbb F}_{q}))$ ce dernier espace s'il y a ambigu\"{\i}t\'e sur le groupe $G$). On ajoute un indice "nil", resp. "cusp", resp. "nil,cusp", pour d\'esigner le sous-espace des fonctions \`a support nilpotent, resp. des fonctions cuspidales, resp. des fonctions cuspidales et \`a support nilpotent.  On munit l'espace ${\bf C}(\mathfrak{g}({\mathbb F}_{q}))$  du produit hermitien d\'efini positif
  $$(f,f')=\vert G({\mathbb F}_{q})\vert ^{-1}\sum_{X\in \mathfrak{g}({\mathbb F}_{q})}\bar{f}(X)f'(X)$$
  pour tous $f,f'\in {\bf C}(\mathfrak{g}({\mathbb F}_{q}))$. Fixons un caract\`ere non trivial $\psi:{\mathbb F}_{q}\to {\mathbb C}^{\times}$ et une forme bilin\'eaire sym\'etrique $<.,.>$ sur $\mathfrak{g}({\mathbb F}_{q})$ non d\'eg\'en\'er\'ee et invariante par conjugaison par $G({\mathbb F}_{q})$. On d\'efinit la transformation de Fourier $f\mapsto \hat{f}$ de $C(\mathfrak{g}({\mathbb F}_{q}))$ par
  $$\hat{f}(X)=q^{-dim(\mathfrak{g})/2}\sum_{Y\in \mathfrak{g}({\mathbb F}_{q})}\psi(<X,Y>)f(Y).$$
  Le carr\'e de cette transformation est l'application $f\mapsto f^-$ d\'efinie par $f^-(X)=f(-X)$. 
  
   On consid\'erera aussi l'espace $C^G(\mathfrak{g}_{SC}({\mathbb F}_{q}))$   des fonctions sur $\mathfrak{g}_{SC}({\mathbb F}_{q}))$ qui sont invariantes par conjugaison par $G({\mathbb F}_{q})$

  \subsubsection{Faisceaux-caract\`eres cuspidaux  sur les alg\`ebres de Lie\label{faisceauxcaracterescuspidaux}}
 La notion de faisceau-caract\`ere n\'ec\'essite de fixer un nombre premier $l\not=p$, ces faisceaux sont alors des $\bar{{\mathbb Q}}_{l}$-faisceaux.  On fixe  aussi un isomorphisme entre $\bar{{\mathbb Q}}_{l}$ et ${\mathbb C}$.

Lusztig  a d\'efini la notion de faisceau-caract\`ere cuspidal \`a support nilpotent sur $\mathfrak{g}$. Il n'y en a que si $Z(G)$ est fini.    Remarquons que cette notion d\'epend non seulement de $\mathfrak{g}$ mais du groupe $G$: en supposant $G$ semi-simple, ces faisceaux-caract\`eres   ne sont en g\'en\'eral pas les m\^emes pour les groupes $G_{SC}$ et $G_{AD}$, alors que leurs alg\`ebres de Lie sont les m\^emes. On note ${\bf FC}^G(\mathfrak{g})$ l'ensemble de ces faisceaux-caract\`eres cuspidaux et \`a support nilpotent, ou simplement ${\bf FC}(\mathfrak{g})$ s'il n'y a pas d'ambigu\"{\i}t\'e sur le groupe $G$. 
Un tel faisceau ${\cal E}$ est port\'e par une orbite nilpotente ${\cal O}_{{\cal E}}$. Fixons $N\in {\cal O}_{{\cal E}}$. Alors le faisceau est d\'etermin\'e par une repr\'esentation irr\'eductible $\xi_{N}$ de $Z_{G}(N)/Z_{G}(N)^0$ \`a valeurs dans un  $\bar{{\mathbb Q}}_{l}$-espace vectoriel  $V$ de dimension finie. Le groupe $Z_{G}(N)/Z_{G}(N)^0$ agit par multiplication \`a droite sur  $G/Z_{G}(N)^0$. On note $G/Z_{G}(N)^0\times_{\xi_{N}} V$ le quotient de $G/Z_{G}(N)^0\times V$ par la relation d'\'equivalence $(gZ_{G}(N)^0z,v)\equiv (gZ_{G}(N)^0,\xi_{N}(z)v)$ pour tout $z\in Z_{G}(N)/Z_{G}(N)^0$. L'application 
$$G/Z_{G}(N)^0\times_{\xi_{N}} V\quad \stackrel{(gZ_{G}(N)^0,v)\mapsto gNg^{-1}}{\to}\quad {\cal O}_{{\cal E}}$$
fait appara\^{\i}tre $G/Z_{G}(N)^0\times_{\xi_{N}} V$ comme un syst\`eme local sur ${\cal O}_{{\cal E}}$. Alors ${\cal E}$ est, \`a un d\'ecalage pr\`es, ce syst\`eme local prolong\'e au-dessus de $\mathfrak{g}$ par $0$ hors de ${\cal O}_{{\cal E}}$. On a fix\'e un \'el\'ement $N\in {\cal O}_{{\cal E}}$. Soit $x\in G$, posons $N'=x^{-1}Nx$. En rempla\c{c}ant ci-dessus $N$ par $N'$ et $\xi_{N}$ par la repr\'esentation $\xi_{N'}=\xi_{N}\circ Ad(x)$ de  $Z_{G}(N')/Z_{G}(N')^0$ dans $V$, on obtient un faisceau isomorphe au pr\'ec\'edent. 
 
 Le groupe $\Gamma_{{\mathbb F}_{q}}$ agit sur l'ensemble ${\bf FC}(\mathfrak{g})$. On note ${\bf FC}_{{\mathbb F}_{q}}(\mathfrak{g})$ le sous-ensemble des points fixes. Pour ${\cal E}\in {\bf FC}_{{\mathbb F}_{q}}(\mathfrak{g})$, on fixe une action de Frobenius sur ${\cal E}$. Il s'en d\'eduit des actions sur les fibres de ${\cal E}$ au-dessus des \'el\'ements de ${\cal O}_{{\cal E}}$ fix\'es par $\Gamma_{{\mathbb F}_{q}}$ et on suppose que ces actions sont d'ordre fini. 
   De cette action de Frobenius sur ${\cal E}$ se d\'eduit une fonction caract\'eristique , à valeurs dans $\bar{{\mathbb Q}}_{l}$. Par l'isomorphisme que l'on a fixé entre ce corps et ${\mathbb C}$, cette fonction devient une fonction $f_{{\cal E}}$ à valeurs complexes qui appartient à $ C_{nil,cusp}^G(\mathfrak{g}({\mathbb F}_{q}))$. Elle ne d\'epend du choix de l'action de Frobenius que par une homoth\'etie. Dans \cite{W6}, on  a d\'efini l'espace $FC(\mathfrak{g}({\mathbb F}_{q}))$ comme  le sous-espace de $C(\mathfrak{g}({\mathbb F}_{q}))$ engendr\'e par les fonctions caract\'eristiques $f_{{\cal E}}$ pour ${\cal E}\in {\bf FC}_{{\mathbb F}_{q}}(\mathfrak{g})$. D'apr\`es Lusztig, cf. \cite{L6}, l'espace $FC(\mathfrak{g}({\mathbb F}_{q}))$ est le sous-espace des fonctions $f\in C(\mathfrak{g}({\mathbb F}_{q}))$ telles que $f$ et $\hat{f}$ sont \`a support dans $\mathfrak{g}_{nil}({\mathbb F}_{q})$.

 On consid\'erera aussi  l'ensemble ${\bf FC}^G(\mathfrak{g}_{SC})$ des faisceaux-caract\`eres  cuspidaux et \`a support nilpotent sur $\mathfrak{g}_{SC}$ qui sont \'equivariants par l'action de $G$ sur $\mathfrak{g}_{SC}$.   On appelera ces faisceaux des $G$-faisceaux.  Cela revient \`a remplacer dans les constructions pr\'ec\'edentes le groupe $G$ par $G/Z(G)^0$. 
Un $G$-faisceau caract\`ere cuspidal \`a support nilpotent sur $\mathfrak{g}_{SC}$ est un tel faisceau pour $G_{SC}$ pour lequel la repr\'esentation $\xi_{N}$ ci-dessus de $Z_{G_{SC}}(N)/Z_{G_{SC}}(N)^0$  se factorise par la surjection $Z_{G_{SC}}(N)/Z_{G_{SC}}(N)^0\to Z_{G}(N)/Z_{G}(N)^0$. On d\'efinit les variantes \'evidentes ${\bf FC}_{{\mathbb F}_{q}}^G(\mathfrak{g}_{SC})$ et  $FC^G(\mathfrak{g}_{SC}({\mathbb F}_{q}))$ des objets introduits ci-dessus. 

\subsubsection{Faisceaux-caract\`eres sur les alg\`ebres de Lie\label{faisceauxcaracteres}}
Rappelons deux propri\'et\'es importantes. Soient $M$ un Levi de $G$ et ${\cal E}\in {\bf FC}^M(\mathfrak{m}_{SC})$. D'apr\`es \cite{L2} th\'eor\`eme 9.2, on a

(1) les sous-groupes paraboliques de $G$  de composante de Levi $M$ sont tous  contenus dans une m\^eme classe de conjugaison par $G $;

(2) pour tout $n\in Norm_{G}(M)$, l'action $Ad(n)$ sur $\mathfrak{m}_{SC}$ conserve ${\cal E}$. 

Posons ${\cal O}={\cal O}_{{\cal E}}$ et fixons $N\in {\cal O}$. Le faisceau ${\cal E}$ est d\'etermin\'e par une repr\'esentation irr\'eductible $\xi_{N}$ de $Z_{M}(N)/Z_{M}(N)^0$ dans un $\bar{{\mathbb Q}}_{l}$-espace $V$.   Notons $\mathfrak{z}_{reg}(M)$ le sous-ensemble des \'el\'ements de $\mathfrak{z}(M)$ dont le commutant dans $G$ est \'egal \`a $M$. Posons  
$${\cal V}_{reg}=\{(X,gZ_{M}(N))\in \mathfrak{g}\times G/Z_{M}(N); g^{-1}Xg\in N+\mathfrak{z}_{reg}(M)\},$$ 
 $$\tilde{{\cal V}}_{reg}=\{(X,gZ_{M}(N)^0)\in \mathfrak{g}\times G/Z_{M}(N)^0; g^{-1}Xg\in N+\mathfrak{z}_{reg}(M)\}.$$
   Le groupe $Z_{M}(N)/Z_{M}(N)^0$ agit par multiplication \`a droite sur $\tilde{{\cal V}}_{reg}$ et le quotient de $\tilde{{\cal V}}_{reg}$ par cette action est ${\cal V}_{reg}$. On d\'efinit comme en \ref{faisceauxcaracterescuspidaux} le syst\`eme local ${\cal C}_{{\cal V}_{reg}} =\tilde{{\cal V}}_{reg}\times_{\xi_{N}}V$ sur ${\cal V}_{reg}$. Notons $p_{reg}:{\cal V}_{reg}\to \mathfrak{g}$ l'application $(X,gZ_{M}(N))\mapsto X$. Notons $\nabla_{reg}$ son image et $\nabla$ l'adh\'erence de $\nabla_{reg}$. On introduit le syst\`eme local $p_{reg,*}({\cal C}_{{\cal V}_{reg}})$ sur $\nabla_{reg}$. On note $K$ le complexe prolongement d'intersection de  $p_{reg,*}({\cal C}_{{\cal V}_{reg}})[dim(\nabla)]$. C'est un faisceau pervers de support  $\nabla$.

  Le groupe $Norm_{G}(M)$ agit sur $\mathfrak{m}_{SC}$ en conservant le faisceau ${\cal E}$. Lusztig a prouv\'e que l'alg\`ebre d'endomorphismes du complexe $K$ \'etait isomorphe \`a l'alg\`ebre du groupe $W(M)$, cf. \cite{L2} th\'eor\`eme 9.2. Concr\`etement, cela signifie ce qui suit. Notons $Norm_{G}(M,N)$ le sous-groupe des $n\in Norm_{G}(M)$ tels que $nNn^{-1}=N$. Un \'el\'ement de $n\in Norm_{G}(M)$ conserve ${\cal E}$ comme on vient de le dire, donc conserve ${\cal O}$. Il existe donc $m\in M$ tel que $nm\in Norm_{G}(M,N)$. Donc l'inclusion $Norm_{G}(M,N)\to Norm_{G}(M)$ se quotiente en une bijection $Norm_{G}(M,N)/Z_{M}(N)\simeq W(M)$. Alors la repr\'esentation $\xi_{N}$ de $Z_{M}(N)$ dans $V$ se prolonge en une repr\'esentation irr\'eductible $\epsilon_{N}$ de $Norm_{G}(M,N)$. Fixons un tel prolongement.   On munit ${\cal V}_{reg}$ de l'action $n\mapsto \tau_{n}$ de $Norm_{G}(M,N)$ d\'efinie par $\tau_{n}(X,gZ_{M}(N))\mapsto (X,gZ_{M}(N)n^{-1})$ et on munit ${\cal C}_{{\cal V}_{reg}}$ de l'action compatible $n\mapsto \tilde{\tau}_{n}$ de $Norm(M,N)$ d\'efinie par $\tilde{\tau}_{n}(X,gZ_{M}(N)^0,v)\mapsto (X,gZ_{M}(N)^0n^{-1},\epsilon_{N}(n)v)$. Ces actions se quotientent en des actions de $W(M)$. On a $p_{reg}\circ \tau_{n}=p_{reg}$. Donc $\tilde{\tau}_{n}$ se descend en un endomorphisme de $p_{reg,*}({\cal C}_{{\cal V}_{reg}})[dim(\nabla)]$, qui se prolonge en un endomorphisme de $K$. On obtient ainsi une action de $W(M)$ sur $K$.

  A l'orbite nilpotente ${\cal O}\subset \mathfrak{m}_{SC}\subset \mathfrak{m}$  est associ\'ee une orbite induite $\underline{{\cal O}}$ dans $\mathfrak{g}$: si $P$ est un sous-groupe parabolique de $G$ de composante de Levi $M$, $\underline{{\cal O}}$ est l'unique orbite nilpotente de $\mathfrak{g}$ telle que $\underline{{\cal O}}\cap ({\cal O}+\mathfrak{u}_{P})$ soit un ouvert dense de  ${\cal O}+\mathfrak{u}_{P}$. La restriction de $K$ \`a $\underline{{\cal O}}$ est un syst\`eme local irr\'eductible, cf. \cite{L2} th\'eor\`eme 9.2. Donc, pour 
  $\underline{N}\in \underline{{\cal O}}$, $\tau_{n,K}$ agit sur la fibre de $K$ au-dessus de $\underline{N}$ par un scalaire ind\'ependant de $\underline{N}$. L'application qui \`a $n$ associe ce scalaire est un caract\`ere de $W(M)$. On a choisi ci-dessus une repr\'esentation $\epsilon_{N}$. Il y a un unique choix pour lequel le caract\`ere pr\'ec\'edent est trivial. Nous noterons $\epsilon_{N}^{\flat}$ cette repr\'esentation et c'est celle que nous utilisons dans la suite. 
  
   Par notre isomorphisme de $\bar{{\mathbb Q}}_{l}$ sur 
   ${\mathbb C}$, les représentations de $W(M)$ dans des espaces complexes s'identifient à des représentations dans des $\bar{{\mathbb Q}}_{l}$-espaces. Pour tout $\rho\in Irr(W(M))$, fixons une r\'ealisation de $\rho$ dans un $\bar{{\mathbb Q}}_{l}$-espace $V_{\rho}$. Il existe un unique  faisceau pervers  irr\'eductible  $K_{M,{\cal E}\rho}$ de sorte que 
  $$(3) \qquad K\simeq \oplus_{\rho\in Irr(W(M))}V_{\rho}\otimes K_{M,{\cal E},\rho},$$
  l'isomorphisme entrela\c{c}ant l'action de $W(M)$ sur $K$ en l'action \'evidente sur le membre de droite.  
    
 On a choisi   un \'el\'ement $N\in {\cal O}$.  Soit 
   $N'$ un autre \'el\'ement de ${\cal O}$. Fixons $m$ tel que $N'=m^{-1}Nm$. Reprenons la construction en rempla\c{c}ant $N$ par $N'$  et la repr\'esentation $\xi_{N}$ de $Z_{M}(N)$ par la repr\'esentation $\xi_{N'}=\xi_{N}\circ Ad(m)$ de $Z_{M}(N')$.  On obtient le m\^eme complexe $K$. L'action de $W(M)$ est inchang\'ee, en utilisant bien s\^ur la repr\'esentation canonique $\epsilon_{N'}^{\flat}$ de  $Norm_{G}(M,N')$.     On a  $Norm_{G}(M,N')=
  m^{-1}Norm_{G}(M,N)m$ et on v\'erifie que
  
  (4)  $\epsilon_{N'}^{\flat} =\epsilon^{\flat}_{N}\circ Ad(m)$.  
  
  Indiquons quelques propri\'et\'es de la repr\'esentation $\epsilon_{N}^{\flat}$ qui r\'esultent de sa construction.
  
  (5) L'application $\epsilon_{N}^{\flat}$ co\"{\i}ncide avec $\xi_{N}$ sur $Z_{M}(N)$. 
  
   (6) Soit $L$ un Levi de $G$ contenant $M$. En rempla\c{c}ant $G$ par $L$, on d\'efinit  l'application $\epsilon^{L,\flat}_{N}$ de $Norm_{L}(M,N)$.  Alors la restriction de $\epsilon_{N}^{\flat}$ \`a $Norm_{L}(M,N)$ est \'egale \`a $\epsilon_{N}^{L,\flat}$.

 (7) Supposons que $G$ soit un produit $G_{1}\times G_{2}$ de deux groupes r\'eductifs connexes d\'efinis sur $\bar{{\mathbb F}}_{q}$. Tous nos objets se d\'ecomposent: $M=M_{1}\times M_{2}$, $N=N_{1}\oplus N_{2}$, $Norm_{G}(M,N)=Norm_{G_{1}}(M_{1},N_{1})\times Norm_{G_{2}}(M_{2},N_{2})$. Alors $\epsilon_{N}^{\flat}=\epsilon_{N_{1}}^{G_{1},\flat}\otimes \epsilon_{N_{2}}^{G_{2},\flat}$.
 
 (8) Soit $A$ un automorphisme de $V$. Si l'on remplace la représentation $\xi_{N}$ par la représentation équivalente $z\mapsto A\xi_{N}(z)A^{-1}$, $\epsilon_{N}^{\flat}$ est remplacé par la représentation $n\mapsto A\epsilon_{N}^{\flat}(n)A^{-1}$. 
 
  (9) Consid\'erons un  groupe r\'eductif connexe $G'$ d\'efini sur $\bar{{\mathbb F}}_{q}$ et un homomorphisme $\delta:G'\to G$. On  
  suppose que le noyau de $\delta$ est central dans $G'$ et que l'homomorphisme compos\'e $G'\stackrel{\delta}{\to }G\to G_{AD}$ est surjectif. 
  Notons $M'=\delta^{-1}(M)$, $N'=\delta^{-1}(N)$. Par $\delta^{-1}$, on transporte ${\cal E}$ en un faisceau-caractère cuspidal ${\cal E}'$ sur $\mathfrak{m}'_{SC}$.  La construction appliqu\'ee \`a $G'$ et $M'$ fournit une représentation  $\epsilon_{N'}^{G',\flat}$ de $Norm_{G'}(M',N')$. On a alors $\epsilon_{N'}^{G',\flat}=\epsilon_{N}^{\flat}\circ \delta$. 
  
  Supposons que $M$ soit un tore, donc $N=0$. Le syst\`eme local ${\cal E}$ est trivial et on v\'erifie que
  
  (10) si $M$ est un tore, $\epsilon_{0}^{\flat}=1$.

  \subsubsection{Calcul de $\epsilon_{N}^{\flat}$  \label{calculepsilonflat}}
  Soient $M$ un Levi de $G$ et ${\cal E}\in {\bf FC}^M(\mathfrak{m}_{SC})$.   
 Posons ${\cal O}={\cal O}_{{\cal E}}$. On a d\'efini l'orbite induite $\underline{{\cal O}}$.  Notons $\bar{{\cal O}}$ l'adh\'erence de ${\cal O}$ dans $\mathfrak{m}_{SC}$. Soit $P$ un sous-groupe parabolique de $G$ de composante de Levi $M$. Fixons deux  \'el\'ements $N\in {\cal O}$ et  $\underline{N}\in (N+\mathfrak{u}_{P})\cap \underline{{\cal O}}$.   Montrons que
   
   (1)  pour $g\in G$, les trois conditions suivantes sont équivalentes;
   
    (a) $g^{-1}\underline{N}g\in\bar{{\cal O}} +\mathfrak{u}_{P}$;
    
    (b) $g^{-1}\underline{N}g\in {\cal O}+\mathfrak{u}_{P}$;
    
    (c) $g\in P$.

    Soit ${\cal O}'\in \mathfrak{m}_{SC,nil}/conj$ telle que ${\cal O}'\subset \bar{{\cal O}}$ et $\underline{{\cal O}}\cap ({\cal O}'+\mathfrak{u}_{P})\not=\emptyset$. Introduisons l'orbite induite $\underline{{\cal O}}'$. Par d\'efinition de cette induite, on a $({\cal O}'+\mathfrak{u}_{P})\subset \bar{\underline{{\cal O}}}'$, donc $\underline{{\cal O}}\subset \bar{\underline{{\cal O}}}'$ puis $dim(\underline{{\cal O}})\leq dim(\underline{{\cal O}}')$. D'apr\`es \cite{LSp} th\'eor\`eme  1.3(a), on a $dim(\underline{{\cal O}})=dim({\cal O})+2dim(U_{P})$, $dim(\underline{{\cal O}}')=dim({\cal O}')+2dim(U_{P})$. Donc $dim({\cal O})\leq dim({\cal O}')$. Avec l'hypoth\`ese ${\cal O}'\subset \bar{{\cal O}}$, cela entra\^{\i}ne ${\cal O}'={\cal O}$. Cela démontre que (a) est équivalent à (b). D'après \cite{L??} proposition 2.14, (b) implique (c). Il est clair que (c) implique (b). Cela prouve (1).

 Le complexe $K$ admet la description suivante. Posons ${\cal W}=\{(X,gZ_{M}(N)U_{P})\in  \mathfrak{g}\times G/Z_{M}(N)U_{P}; g^{-1}Xg\in N+\mathfrak{z}(M)+\mathfrak{u}_{P}\}$.    Posons $\tilde{{\cal W}}= \{(X,gZ_{M}(N)^0U_{P})\in  \mathfrak{g}\times G/Z_{M}(N)^0U_{P}; g^{-1}Xg\in N+\mathfrak{z}(M)+\mathfrak{u}_{P}\}$. Le groupe $Z_{M}(N)/Z_{M}(N)^0$ agit par multiplication \`a droite sur $\tilde{{\cal W}}$ et le quotient de $\tilde{{\cal W}}$ par cette action est ${\cal W}$. On d\'efinit le syst\`eme local ${\cal C}=\tilde{{\cal W}}\otimes_{\xi_{N}}V$ sur ${\cal W}$.  Introduisons le sous-ensemble ${\cal W}_{reg}$ de ${\cal W}$ form\'e des $(X,gZ_{M}(N)U_{P})$ tels que $g^{-1}Xg\in N+\mathfrak{z}_{reg}(M)+\mathfrak{u}_{P}$ et le sous-ensemble similaire $\tilde{{\cal W}}_{reg}$ de $\tilde{{\cal W}}$. Les applications $\iota:(X,gZ_{M}(N))\mapsto (X,gZ_{M}(N)U_{P})$ et $\tilde{\iota}:(X,gZ_{M}(N)^0)\mapsto (X,gZ_{M}(N)^0U_{P})$  sont des isomorphismes de ${\cal V}_{reg}$ sur ${\cal W}_{reg}$ et de $\tilde{{\cal V}}_{reg}$ sur $\tilde{{\cal W}}_{reg}$. Le syst\`eme local $\iota^*({\cal C})$ sur ${\cal V}_{reg}$ est le syst\`eme  ${\cal C}_{{\cal V}_{reg}}$. 
 On note $p_{1}:{\cal W}\to \mathfrak{g}$ la projection $(X,gZ_{M}(N)U_{P})\mapsto X$.  Son image est la cl\^oture de $p_{1}({\cal W}_{reg})$, qui est \'egal \`a $p_{reg}({\cal V}_{reg})=\nabla_{reg}$. Donc l'image de $p_{1}$ est $\nabla$. Lusztig prouve que $K$ est \'egal \`a 
  $p_{1!}({\cal C})[dim(\nabla)]$, cf. \cite{L2} proposition 4.5.   
  
   {\bf On suppose maintenant que } $G$ {\bf  est semi-simple et que} $M$ {\bf est un Levi propre maximal de} $G$. 
 Le groupe $Norm_{G}(M)/M$ a au plus $2$ \'el\'ements. La propri\'et\'e \ref{faisceauxcaracteres}(1)  implique que ce groupe a bien $2$ \'el\'ements.   La restriction de $\epsilon_{N}^{\flat}$ \`a $Z_{M}(N)$ est $\xi_{N}$. Pour d\'eterminer compl\`etement $\epsilon_{N}^{\flat}$, il suffit de fixer $n\in Norm_{G}(M,N)-Z_{M}(N)$ et de calculer $\epsilon_{N}^{\flat}(n)$.  
 
           Les hypothèses sur $M$ et $G$ impliquent que $\mathfrak{z}(M)$ est une droite. Fixons-en un g\'en\'erateur $Z$. Il appartient \`a  $ \mathfrak{z}_{reg}(M)$. Soit $\nu\in \bar{{\mathbb F}}_{q}^{\times}={\mathbb G}_{m}(\bar{{\mathbb F}}_{q})$. Il existe un unique $u_{\nu}\in U_{P}$ tel que $u_{\nu}^{-1}(\nu Z+\underline{N})u_{\nu}=\nu Z+N$.  L'application $\nu\mapsto u_{\nu}$ est alg\'ebrique, c'est-\`a-dire est un morphisme alg\'ebrique de ${\mathbb G}_{m}$ dans $U_{P}$.  Puisque la vari\'et\'e $G/P$ est projective, le morphisme $\nu\mapsto u_{\nu}n^{-1}P$ de ${\mathbb G}_{m}$ dans $G/P$   se prolonge en un morphisme de ${\mathbb A}^1$ dans $G/P$, autrement dit se prolonge en $0$. Notons $g_{0}P$ sa valeur en $0$. Pour $\nu\not=0$, on
 a $n u_{\nu}^{-1}(\nu Z+\underline{N})u_{\nu}n^{-1}=-\nu Z+N$.  On en d\'eduit que   $g_{0}^{-1}\underline{N}g_{0}\in \bar{{\cal O}}+\mathfrak{u}_{P}$. D'apr\`es (1), on a $g_{0}P=P$.  
  Notons $P^{op}$ le sous-groupe parabolique de $G$ de composante de Levi $M$ et oppos\'e \`a $P$. Puisque l'application $\underline{v}\mapsto \underline{v}P$ de $U_{P^{op}}$ dans $G/P$ est une carte de la vari\'et\'e $G/P$ au voisinage de $P\in G/P$, il existe un voisinage de Zariski ${\cal U}_1$ de $0$ dans ${\mathbb A}^1$ tel que notre morphisme de ${\mathbb A}^1$ dans $G/P$ se rel\`eve en un morphisme de ${\cal U}_1$ dans $U_{P^{op}}$, que l'on note $\nu\mapsto \underline{v}'_{\nu}$. On a $\underline{v}'_{0}=1$. Pour $\nu\not=0$, on a 
  $u_{\nu}n^{-1}P=\underline{v}'_{\nu}P$ et  on peut \'ecrire  de fa\c{c}on unique $u_{\nu}n^{-1}=\underline{v}'_{\nu}m_{\nu}v'_{\nu}$ avec $m_{\nu}\in M$ et $v'_{\nu}\in U_{P}$. Les applications $\nu\mapsto m_{\nu}$ et $\nu\mapsto v'_{\nu}$ sont alg\'ebriques. 
  L'\'egalit\'e $n u_{\nu}^{-1}(\nu Z+\underline{N})u_{\nu}n^{-1}=-\nu Z+N$ \'equivaut \`a
 $$\underline{v}_{\nu}^{'-1}(\nu Z+\underline{N})\underline{v}'_{\nu}=m_{\nu}v'_{\nu}(-\nu Z+N)v^{'-1}_{\nu}m_{\nu}^{-1}.$$ 
 On d\'ecompose chaque terme selon la d\'ecomposition $\mathfrak{g}=\mathfrak{u}_{P^{op}}\oplus\mathfrak{m}\oplus\mathfrak{u}_{P}$.  Le membre de gauche  se prolonge en $\nu=0$, sa valeur en $0$ est $\underline{N}$, dont la composante dans $\mathfrak{m}$ est $N$. La composante dans $\mathfrak{m}$ du membre de droite est $-\nu Z+m_{\nu}Nm_{\nu}^{-1}$. Donc  l'application $\nu\mapsto m_{\nu}Nm_{\nu}^{-1}$ se prolonge en un morphisme de ${\cal U}_1$ dans ${\cal O}$ dont la valeur en $\nu=0$   est $N$. L'application $m\mapsto mNm^{-1}$ de $M$ dans ${\cal O}$ est submersive en $m=1$. Il existe donc un voisinage ${\cal U}_2$ de $0$ dans ${\cal U}_1$ tel que notre application $\nu\mapsto m_{\nu}Nm_{\nu}^{-1}$ se rel\`eve en un morphisme de ${\cal U}_2$ dans $M$ dont la valeur en $0$ est $1$. On note $\nu\mapsto m'_{\nu}$ ce morphisme. Pour $\nu\not=0$, on a alors 
   $m_{\nu}=m'_{\nu}z_{\nu}$ o\`u   $z_{\nu}\in Z_{M}(N)$. L'application $\nu\mapsto z_{\nu}$ est alg\'ebrique et il en est de m\^eme de $\nu\mapsto z_{\nu}Z_{M}(N)^0\in Z_{M}(N)/ Z_{M}(N)^0$. Puisque ce dernier groupe est fini, l'application $\nu\mapsto z_{\nu}Z_{M}(N)^0$ est constante. On peut donc fixer $z\in Z_{M}(N)^0$ de sorte que $z_{\nu}=zz^0_{\nu}$, avec $z^0_{\nu}\in Z_{M}(N)^0$. 
   Posons $B=\epsilon_{N}^{\flat}(n)\in Aut(V)$.  Soit $v\in V$. 
 Pour $\nu\not=0$, le triplet $(\nu Z+\underline{N}, u_{\nu}Z_{M}(N)^0,v)$ appartient \`a $\tilde{{\cal V}}_{reg}\times_{\xi_{N}}V$. Son image par $\tilde{\tau}_{n}$   est $(\nu Z+\underline{N}, u_{\nu}n^{-1}Z_{M}(N)^0,B(v))=(\nu Z+\underline{N}, \underline{v}'_{\nu}m'_{\nu}zz^0_{\nu}v'_{\nu}Z_{M}(N)^0,B(v))$. En poussant ces \'el\'ements dans $\tilde{{\cal W}}_{reg}\otimes_{\xi_{N}}V$, on a $\tilde{\iota}(\nu Z+\underline{N}, u_{\nu}Z_{M}(N)^0,v)=(\nu Z+\underline{N}, Z_{M}(N)^0U_{P},v)$ et $\tilde{\iota}(\nu Z+\underline{N}, \underline{v}'_{\nu}m'_{\nu}zz^0_{\nu}v'_{\nu}Z_{M}(N)^0,B(v))=(\nu Z+\underline{N}, \underline{v}'_{\nu}m'_{\nu}Z_{M}(N)^0U_{P},\xi_{N}(z)B(v))$. Les composantes   dans $\tilde{{\cal W}}_{reg}$ de ces deux termes se prolongent en $\nu=0$ et leur valeur commune  en ce point est $(\underline{N},Z_{M}(N)^0U_{P})$. Donc $\tau_{n,K}$ agit sur la fibre de $K$ au-dessus de $\underline{N}$ par l'automorphisme  $\xi_{N}(z)B$. 
Par d\'efinition de $\epsilon_{N}^{\flat}$,  cet automorphisme est l'identité.   D'o\`u l'\'egalit\'e $\xi_{N}(z)B=1$, c'est-à-dire

 (2) $\epsilon_{N}^{\flat}(n)= \xi_{N}(z)^{-1}$. 
 
 Consid\'erons le {\bf  cas particulier o\`u}  ${\cal O}$ {\bf  est l'orbite r\'eguli\`ere de} $\mathfrak{m}_{SC}$. Dans ce cas, 
 on a $Z_{M}(N)/Z_{M}(N)^0=Z(M)/Z(M)^0$, $V$ est de dimension $1$ et $\xi_{N}$ est un caractère.
Fixons un sous-tore maximal $T$ de $M$ et un sous-groupe de Borel $B^M$ de $M$  contenant $T$  de sorte  que $N\in \mathfrak{u}_{B^M}$. Il existe alors un sous-groupe $R\subset U_{B^M}$ tel que $Z_{M}(N)=Z(M)R$. Pour $\nu\in {\cal U}_2$, on \'ecrit $z^0_{\nu}=z^{''0}_{\nu}r''_{\nu}$ avec $z^{''0}_{\nu}\in Z(M)^0$ et $r''_{\nu}\in R$. 
  Posons $B=B^MU_{P}$. Le groupe $B$ est un sous-groupe de Borel de $G$. On note $B^{M,opp}$ et $B^{opp}$ les sous-groupes de Borel oppos\'es de $B^M$ et $B$. 
  Puisque $m'_{0}=1$, il existe un voisinage ${\cal U}_3$ de $0$ dans  ${\cal U}_2$ de sorte que, pour $\nu\in {\cal U}_3$, on puisse \'ecrire $m'_{\nu}=\underline{v}''_{\nu}t''_{\nu}v''_{\nu}$ avec $\underline{v}''_{\nu}\in U_{B^{M,opp}}$, $t''_{\nu}\in T$, $v''_{\nu}\in U_{B^M}$. Les applications $\nu\mapsto \underline{v}''_{\nu}$, $\nu\mapsto t''_{\nu}$ et $\nu\mapsto v''_{\nu}$ sont alg\'ebriques dans ${\cal U}_3$ et valent $1$ en $\nu=0$. Pour $\nu\not=0$, posons $\underline{v}_{\nu}=\underline{v}'_{\nu}\underline{v}''_{\nu}$, $t_{\nu}=t''_{\nu}zz^{''0}_{\nu}$, $v_{\nu}=v''_{\nu}r''_{\nu}v'_{\nu}$. On a $\underline{v}_{\nu}\in U_{B^{opp}}$, $t_{\nu}\in T$, $v_{\nu}\in U_{B}$ et on a 
   l'\'egalit\'e
    
 (3)  $u_{\nu}n^{-1}=\underline{v}_{\nu}t_{\nu}v_{\nu}$.
 
 L'application $\nu\mapsto \underline{v}_{\nu}$ se prolonge alg\'ebriquement \`a ${\cal U}_3$ et sa valeur en $0$ est $1$. L'application $\nu\mapsto t_{\nu}$ ne se prolonge pas forc\'ement \`a ${\cal U}_3$. Mais notons ${\bf t}_{\nu}$, resp. ${\bf z}$,  l'image de $t_{\nu}$, resp. $z$,  dans $T/Z(M)^0$. Alors l'application $\nu\mapsto {\bf t}_{\nu}$ se prolonge \`a ${\cal U}_3$ et sa valeur en $0$ est ${\bf z}$. Cela d\'etermine ${\bf z}$, donc aussi $\xi_{N}(z)$.

    \subsubsection{Fonctions caract\'eristiques\label{fonctionscaracteristiques}}

 On reprend les constructions du paragraphe \ref{faisceauxcaracteres} en supposant que   $M$  est un ${\mathbb F}_{q}$-Levi de $G$ et  que ${\cal E}\in {\bf FC}_{{\mathbb F}_{q}}^M(\mathfrak{m}_{SC})$. On munit ${\cal E}$ d'une action de Frobenius, c'est-\`a-dire que l'on fixe un isomorphisme $Fr^*({\cal E})\simeq {\cal E}$. Concr\`etement, l'orbite ${\cal O}:={\cal O}_{{\cal E}}$ est conserv\'ee par l'action galoisienne et on choisit l'\'el\'ement $N$ dans ${\cal O}^{\Gamma_{{\mathbb F}_{q}}}$. La repr\'esentation $\xi_{N}$ de $Z_{M}(N)$ dans $V$ est alors conserv\'ee par l'action de Frobenius. On fixe un automorphisme $r_{N}$ de $V$ tel que $\xi_{N}(z)r_{N}=r_{N}\xi_{N}(Fr(z))$ pour tout $z\in Z_{M}(N)$. On suppose $r_{N}$ d'ordre fini.  L'isomorphisme $Fr^*({\cal E})\simeq {\cal E}$, c'est-\`a-dire 
   $$Fr^*(M/Z_{M}(N)^0\times_{\xi_{N}}V)\to M/Z_{M}(N)^0\times_{\xi_{N}}V$$
   est alors $(mZ_{M}(N)^0,v)\mapsto (Fr^{-1}(m)Z_{M}(N)^0,r_{N}(v))$. Soit $N'$ un autre \'el\'ement de ${\cal O}^{\Gamma_{{\mathbb F}_{q}}}$, fixons $m\in M$ tel que $N'=m^{-1}Nm$ et rempla\c{c}ons dans la construction ci-dessus l'\'el\'ement $N$ par $N'$ et $\xi_{N}$ par $\xi_{N'}=\xi_{N}\circ Ad(m)$. Pour obtenir la m\^eme action de Frobenius sur ${\cal E}$, on doit choisir l'isomorphisme 
   
   (1) $r_{N'}=r_{N}\xi_{N}(Fr(m)m^{-1})$.

   Remarquons que le prolongement canonique $\epsilon_{N}^{\flat}$ de $\xi_{N}$ \`a $Norm_{G}(M,N)$ v\'erifie la relation
   
   (2) $\epsilon_{N}^{\flat}(n)r_{N}=r_{N}\epsilon_{N}^{\flat}(Fr(n))$ pour tout $n\in Norm_{G}(M,N)$.

 Le groupe $\Gamma_{{\mathbb F}_{q}}$ agit sur $W(M)$. Il agit aussi sur $Irr(W(M))$. Notons $Irr_{{\mathbb F}_{q}}(W(M))$ le sous-ensemble des repr\'esentations conserv\'ees par cette action. Pour tout $\rho\in Irr_{{\mathbb F}_{q}}(W(M))$, fixons un prolongement continu $\rho^{\flat}$ de $\rho$ au produit semi-direct $W(M)\rtimes \Gamma_{{\mathbb F}_{q}}$, "continu" signifiant qu'il est trivial sur $\Gamma_{{\mathbb F}_{q^n}}$ pour un certain  entier $n\geq1$. Evidemment, si $\Gamma_{{\mathbb F}_{q}}$ agit trivialement sur $W(M)$, on suppose que $\rho^{\flat}$ est triviale sur $\Gamma_{{\mathbb F}_{q}}$. Pour une repr\'esentation $\rho\in Irr(W(M))-Irr_{{\mathbb F}_{q}}(W(M))$, posons $\rho'=\rho\circ Fr$ et fixons  un isomorphisme $\rho^{\flat}(Fr^{-1}):V_{\rho}\to V_{\rho'}$ tel que $\rho^{\flat}(Fr^{-1})\rho(Fr(w))=\rho'(w)\rho^{\flat}(Fr^{-1})$ pour tout $w\in W(M)$.  De l'isomorphisme $Fr^*({\cal E})\simeq {\cal E}$ se d\'eduit naturellement un isomorphisme $Fr^*(K)\simeq K$.   Pour tout $\rho\in Irr(W(M))$, il existe un unique isomorphisme  $Fr^*(K_{M,{\cal E},\rho})\simeq K_{M,{\cal E},\rho\circ Fr}$ de sorte que, par l'isomorphisme \ref{faisceauxcaracteres}(3), notre isomorphisme $Fr^*(K)\simeq K$ s'identifie \`a la somme sur $\rho\in Irr(W(M))$ des produits tensoriels des isomorphismes  $Fr^*(K_{M,{\cal E},\rho})\simeq K_{M,{\cal E},\rho\circ Fr}$ et $\rho^{\flat}(Fr^{-1}):V_{\rho}\to V_{\rho\circ Fr}$. En particulier, pour $\rho\in Irr_{{\mathbb F}_{q}}(W(M))$, $K_{M,{\cal E},\rho}$ est ainsi muni d'une action de Frobenius.     On note $\boldsymbol{\chi}_{M,{\cal E},\rho}$ sa fonction caract\'eristique. Cette fonction, comme celles qui interviendront ci-dessous, est identifiée à une fonction à valeurs complexes. 
  
   On dit que deux \'el\'ements $w,w'\in W(M)$   sont $Fr$-conjugu\'es si et seulement s'il existe $u\in W(M)$ tel que $w=Fr(u)w'u^{-1}$.  
  
  {\bf Remarque.} La place du Frobenius dans cette relation est  inhabituelle.  On l'a choisie par souci de coh\'erence avec certaines d\'efinitions de Lusztig.
  \bigskip
   
   On note $W(M)/Fr-conj$ l'ensemble des classes de $Fr$-conjugaison.       Soit $w\in W(M)$.  Fixons un repr\'esentant $n\in Norm_{G}(M,N) $ de $w$ et un \'el\'ement  $g\in G $ tel que $gFr(g)^{-1}= n^{-1}$. Posons $M_{w}=g^{-1}Mg$. Alors $M_{w}$ est d\'efini sur ${\mathbb F}_{q}$. De ${\cal E}$ se d\'eduit naturellement un $M_{w}$-faiceau-caract\`ere cuspidal ${\cal E}_{w}$ sur $\mathfrak{m}_{w,SC}$. On pose ${\cal O}_{w}=g^{-1}{\cal O}g$, $N_{w}=g^{-1}Ng$ et on d\'efinit la repr\'esentation $\xi_{N_{w}}$ de $Z_{M_{w}}(N_{w})$ dans $V$ par $\xi_{N_{w}}(z)=\xi_{N}(gzg^{-1})$ pour tout $z\in Z_{M_{w}}(N_{w})$. Alors ${\cal E}_{w}$ est le faisceau-caract\`ere port\'e par l'orbite ${\cal O}_{w}$ qui est associ\'e \`a la repr\'esentation $\xi_{N_{w}}$ de $Z_{M_{w}}(N_{w})$. On a $N_{w}\in {\cal O}_{w}^{\Gamma_{{\mathbb F}_{q}}}$. On d\'efinit l'automorphisme $r_{N_{w}}=r_{N}\epsilon_{N}^{\flat}(n)$ de $V$ et on munit ${\cal E}_{w}$ de l'action de Frobenius d\'efinie comme plus haut par l'automorphisme $r_{N_{w}}$. On reprend alors la construction de \ref{faisceauxcaracteres} o\`u l'on remplace $M$ et ${\cal E}$ par $M_{w}$ et ${\cal E}_{w}$. On obtient un complexe $K_{w}$ muni d'une action de Frobenius. On note $\boldsymbol{\chi}_{M,{\cal E},w}$ sa fonction caract\'eristique. Elle ne d\'epend pas des choix effectu\'es (except\'e bien s\^ur celui de l'action de Frobenius de d\'epart sur ${\cal E}$) et ne d\'epend que de la classe de $Fr$-conjugaison de $w$.  Pour $\rho \in Irr_{{\mathbb F}_{q}}(W(M))$, on a alors l'\'egalit\'e
   
  (3)  $\boldsymbol{\chi}_{M,{\cal E},\rho}=\vert W(M)\vert ^{-1}\sum_{w\in W(M)}trace(\rho^{\flat}(w^{-1}Fr))\boldsymbol{\chi}_{M,{\cal E},w}$,
  
  \noindent cf. \cite{L3} 10.4.5.   

    Pour $w\in W(M)$, notons $f_{{\cal E}_{w}}$ la fonction caract\'eristique de ${\cal E}_{w}$ muni de l'action de Frobenius d\'efinie ci-dessus. 
 C'est une fonction \`a support nilpotent sur $\mathfrak{m}_{w,SC}({\mathbb F}_{q})$, que l'on peut consid\'erer comme une fonction \`a support nilpotent sur $\mathfrak{m}({\mathbb F}_{q})$. Elle est invariante par conjugaison par $M_{w}({\mathbb F}_{q})$. On dispose de l'induction de Deligne-Lusztig $R_{M_{w}}^G:C(\mathfrak{m}_{w}({\mathbb F}_{q}))\to C(\mathfrak{g}({\mathbb F}_{q}))$. On note  $Q_{M,{\cal E},w}$  l'image de $f_{{\cal E}_{w}}$ par $R_{M_{w}}^G$.  D'apr\`es \cite{L7} th\'eor\`eme 1.14, on a l'\'egalit\'e
 
 (4) $Q_{M,{\cal E},w}(X)=(-1)^{dim(Z(M)^0)}\boldsymbol{\chi}_{M,{\cal E},w}(X)$  pour tout $X\in \mathfrak{g}({\mathbb F}_{q})$.

 Soit $g\in G({\mathbb F}_{q})$, posons $M'=g^{-1}Mg$. L'application $Ad(g)^{-1}$ transporte ${\cal E}$ en un faisceau ${\cal E}'$ sur $\mathfrak{m}'_{SC}$, qui est muni d'une action de Frobenius. Elle \'etablit aussi une bijection $w\mapsto w'$ de $W(M)$ sur $W(M')$. Il est clair que $Q_{M',{\cal E}',w'}=Q_{M,{\cal E},w}$.

Pour tout Levi $H$ de $G$ d\'efini sur ${\mathbb F}_{q}$, posons $W_{{\mathbb F}_{q}}(H)=Norm_{G({\mathbb F}_{q})}(H)/H({\mathbb F}_{q})$. Les deux relations suivantes r\'esultent du corollaire 9.11 de \cite{L3}.   Pour tout $w,w'\in W(M)$, on a
 
   $$(5)  (Q_{M,{\cal E},w},Q_{M,{\cal E},w'})=\left\lbrace\begin{array}{cc}\vert W_{{\mathbb F}_{q}}(M_{w})\vert  \vert Z(M_{w})^0({\mathbb F}_{q})\vert^{-1} q^{dim({\cal O}_{{\cal E}})-dim(M_{SC})}&\text{ si }w\text{ et }w'\text{ sont }\\ &Fr-\text{ conjugu\'es}\\ 0&\text{ sinon}\end{array}\right.$$ 
  
  D'autre part, soient  $(M',{\cal E}',w') $ v\'erifiant les m\^emes propri\'et\'es que $(M,{\cal E},w)$. On a l'\'egalit\'e
  
 (6) $(Q_{M,{\cal E},w},Q_{M',{\cal E}',w'})=0$  si $(M',{\cal E}')$ n'est pas conjugu\'e \`a $(M,{\cal E})$ par un \'el\'ement de $G({\mathbb F}_{q})$.

D'autres fonctions nous seront utiles. Soit $M$ un ${\mathbb F}_{q}$-Levi de $G$ tel que ${\bf FC}_{{\mathbb F}_{q}}^M(\mathfrak{m}_{SC})$ soit non vide. Soit $w\in W(M)$.   De l'application ${\cal E}\mapsto {\cal E}_{w}$ de ${\bf FC}_{{\mathbb F}_{q}}^M(\mathfrak{m}_{SC})$ dans ${\bf FC}_{{\mathbb F}_{q}}^{M_{w}}(\mathfrak{m}_{w,SC})$ se d\'eduit une application $f_{{\cal E}}\mapsto f_{{\cal E}_{w}}$ qui se prolonge par lin\'earit\'e en un isomorphisme $\iota^{\flat}_{M_{w},M}:FC^{M}(\mathfrak{m}_{SC}({\mathbb F}_{q}))\to  FC^{M_{w}}(\mathfrak{m}_{w,SC}({\mathbb F}_{q}))$.  Soit $\varphi\in FC^{M_{w}}(\mathfrak{m}_{w,SC}({\mathbb F}_{q}))$.  
 Fixons un \'el\'ement $X_{w}\in \mathfrak{z}_{reg}(M_{w},{\mathbb F}_{q})$.   Notons $\varphi_{X_{w}}^{\natural}$ la fonction sur $\mathfrak{g}({\mathbb F}_{q})$, \`a support dans $X_{w}+\mathfrak{m}_{w,SC}({\mathbb F}_{q})$, telle que, pour $Y\in \mathfrak{m}_{w,SC}({\mathbb F}_{q})$, on ait $\varphi_{X_{w}}^{\natural}(X_{w}+Y)=\varphi(Y)$. Notons $\varphi_{X_{w}}$ la fonction sur $\mathfrak{g}({\mathbb F}_{q})$ d\'efinie par
$$\varphi_{X_{w}}(Z)=\vert M_{w}({\mathbb F}_{q})\vert ^{-1}\sum_{g\in G({\mathbb F}_{q})}\varphi_{X_{w}}^{\natural}(g^{-1}Zg).$$
 Notons enfin
  ${\cal Q}_{\varphi}$ la restriction \`a $\mathfrak{g}_{nil}({\mathbb F}_{q})$ de $\hat{\varphi}_{X_{w}}$. Cette fonction ne d\'epend pas du choix de $X_{w}$.  Soit ${\cal E}\in {\bf FC}_{{\mathbb F}_{q}}^M(\mathfrak{m}_{SC})$, appliquons la construction \`a $\varphi=f_{{\cal E}_{w}}$. On a
  
  (7) ${\cal Q}_{f_{{\cal E}_{w}}}$ est la restriction \`a $\mathfrak{g}_{nil}({\mathbb F}_{q})$ de $\hat{Q}_{M,{\cal E},w}$. 
  
  \noindent Cette \'egalit\'e a \'et\'e prouv\'ee par Kazhdan dans le cas o\`u $M_{w}$ est un tore, cf. \cite{Kazh}. Dans le cas g\'en\'eral, la proportionnalit\'e des deux fonctions r\'esulte de \cite{L6} 7(a). On a donn\'e une preuve de l'\'egalit\'e en \cite{W1} proposition II.8. Dans  cette r\'ef\'erence, le groupe est suppos\'e classique mais cette hypoth\`ese ne sert pas pour la preuve en question.

  L'espace  $C_{nil}(\mathfrak{g}({\mathbb F}_{q}))$ a pour base la famille des fonctions ${\cal Q}_{\varphi}$, quand $M$ d\'ecrit les  ${\mathbb F}_{q}$-Levi de $G$, \`a conjugaison par $G({\mathbb F}_{q})$ pr\`es, tels que ${\bf FC}_{{\mathbb F}_{q}}^M(\mathfrak{m}_{SC})\not=\emptyset$, $w$ d\'ecrit $W(M)/Fr-conj$   et $\varphi$ d\'ecrit une base de $FC^{M_{w}}(\mathfrak{m}_{w,SC}({\mathbb F}_{q}))$. 

 Notons $W(M)_{Fr-reg}$ le sous-ensemble des $w\in W(M)$ tels que l'action de $ wFr^{-1}$ dans $X_{*,{\mathbb Q}}(Z(M)^0)/X_{*,{\mathbb Q}}(Z(G)^0)$ n'ait pas de point fixe non nul. Il est conserv\'e par $Fr$-conjugaison. On note $W(M)_{Fr-reg}/Fr-conj$ l'ensemble des classes de $Fr$-conjugaison dans $W(M)_{Fr-reg}$. Alors $C_{nil,cusp}(\mathfrak{g}({\mathbb F}_{q}))$ a pour base la famille des fonctions ${\cal Q}_{\varphi}$, quand $M$ d\'ecrit les  ${\mathbb F}_{q}$-Levi de $G$, \`a conjugaison par $G({\mathbb F}_{q})$ pr\`es, tels que ${\bf FC}_{{\mathbb F}_{q}}^M(\mathfrak{m}_{SC})\not=\emptyset$, $w$ d\'ecrit $W(M)_{Fr-reg}/Fr-conj$   et $\varphi$ d\'ecrit une base de $FC^{M_{w}}(\mathfrak{m}_{w,SC}({\mathbb F}_{q}))$, cf. \cite{L3} th\'eor\`eme 1.14.

\subsubsection{Les espaces $I(\mathfrak{g}(F))$ et $FC(\mathfrak{g}(F))$\label{lespaceFC}}
 Pour la fin de la section 2,  $G$ est un groupe r\'eductif connexe d\'efini sur $F$ et  {\bf on suppose v\'erifi\'ee l'hypoth\`ese $(Hyp)_{2}(p)$}, cf. \ref{leshypothesessurp}.
 
 Fixons un caract\`ere non trivial $\psi$ de $F$  et une forme bilin\'eaire sym\'etrique et non d\'eg\'en\'er\'ee $<.,.>$ sur $\mathfrak{g}(F)$, invariante par l'action par conjugaison de $G(F)$. On d\'efinit la transformation de Fourier $f\mapsto \hat{f}$ dans $C_{c}^{\infty}(\mathfrak{g}(F))$ par la formule usuelle
 $$\hat{f}(X)=\int_{\mathfrak{g}(F)}f(Y)\psi(<X,Y>)\,dY,$$
 o\`u $dY$ est la mesure autoduale. Le carr\'e de cette transformation est l'application $f\mapsto f^-$, o\`u $f^-(X)=f(-X)$. 
 On peut choisir et on choisit la forme $<.,.>$ de sorte que, pour tout ${\cal F}\in Fac(G)$, on ait l'\'egalit\'e $\hat{{\bf 1}}_{\mathfrak{k}_{{\cal F}}}=\vert \mathfrak{g}_{{\cal F}}({\mathbb F}_{q})\vert ^{1/2}{\bf 1}_{\mathfrak{k}_{{\cal F}}^+}$, o\`u l'on a not\'e ${\bf 1}_{\mathfrak{k}_{{\cal F}}}$ et ${\bf 1}_{\mathfrak{k}_{{\cal F}}^+}$ les fonctions caract\'eristiques de $\mathfrak{k}_{{\cal F}}$ et $\mathfrak{k}_{{\cal F}}^+$.  La mesure sur $\mathfrak{g}(F)$ v\'erifie alors l'\'egalit\'e $mes(\mathfrak{k}_{{\cal F}}^+)=\vert \mathfrak{g}_{{\cal F}}({\mathbb F}_{q})\vert ^{-1/2}$ pour tout ${\cal F}\in Fac(G)$. On peut relever cette mesure en une mesure de Haar sur $G(F)$ telle que  $mes(K_{{\cal F}}^+)=mes(\mathfrak{k}_{{\cal F}}^+)$ pour tout ${\cal F}$.  Des mesures analogues seront choisies sur tout autre groupe r\'eductif connexe. Pour  ${\cal F}\in Fac(G)$, la transformation de Fourier se descend naturellement en une transformation de l'espace $C(\mathfrak{g}_{{\cal F}}({\mathbb F}_{q}))$ et c'est celle-ci que nous utilisons. 

Pour $f\in C_{c}^{\infty}(\mathfrak{g}(F))$ et $X\in \mathfrak{g}_{reg}(F)$, on d\'efinit l'int\'egrale orbitale
$$I^G(X,f)=d^G(X)^{1/2}\int_{G_{X}(F)\backslash G(F)}f(g^{-1}Xg)\,dg,$$
o\`u $d^G$ est le discriminant de Weyl. On note $C_{cusp}^{\infty}(\mathfrak{g}(F))$ le sous-espace des $f\in C_{c}^{\infty}(\mathfrak{g}(F))$ telles que   $I^G(X,f)=0$ pour tout $X\in \mathfrak{g}_{reg}(F)$ tel que $X\not\in\mathfrak{g}_{ell}(F)$. 
On note $I(\mathfrak{g}(F))$ le quotient de $C_{c}^{\infty}(\mathfrak{g}(F))$ par le sous-espace des fonctions $f$ telles que $I^G(X,f)=0$ pour tout $X\in \mathfrak{g}_{reg}(F)$. Les int\'egrales orbitales se descendent en des formes lin\'eaires sur $I(\mathfrak{g}(F))$. On note $I_{cusp}(\mathfrak{g}(F))$  l'image de $C_{cusp}^{\infty}(\mathfrak{g}(F))$ dans $I(\mathfrak{g}(F))$.  

 Pour $X\in \mathfrak{g}_{reg}(F)$ et $f\in C_{c}^{\infty}(\mathfrak{g}(F))$, on d\'efinit l'int\'egrale orbitale stable $S^G(X,f)$ comme la somme des $I^G(X',f)$ o\`u $X'$ parcourt un ensemble de repr\'esentants des classes de conjugaison par $G(F)$ dans la classe de conjugaison stable de $X$. On note $SI(\mathfrak{g}(F))$ le quotient de $C_{c}^{\infty}(\mathfrak{g}(F))$ par le sous-espace des fonctions $f$ telles que $S^G(X,f)=0$ pour tout $X\in \mathfrak{g}_{reg}(F)$. C'est aussi un quotient de $I(\mathfrak{g}(F))$. On note $SI_{cusp}(\mathfrak{g}(F))$ l'image de $I_{cusp}(\mathfrak{g}(F))$ dans $SI(\mathfrak{g}(F))$.

  Pour une facette ${\cal F}\in Fac(G)$, l'espace des fonctions sur $ \mathfrak{g}_{{\cal F}}({\mathbb F}_{q})$  s'identifie \`a un sous-espace de $C_{c}^{\infty}(\mathfrak{g}(F))$: une fonction sur $\mathfrak{g}_{{\cal F}}({\mathbb F}_{q})\simeq \mathfrak{k}_{{\cal F}}/\mathfrak{k}_{{\cal F}}^+$ se rel\`eve en une fonction sur $\mathfrak{k}_{{\cal F}}$ puis s'\'etend \`a $\mathfrak{g}(F)$ par $0$ hors de $\mathfrak{k}_{{\cal F}}$.  On note $FC(\mathfrak{g}(F))$ l'image  dans $I(\mathfrak{g}(F))$ de l'espace
  $$\sum_{s\in S(G)}FC(\mathfrak{g}_{s}({\mathbb F}_{q})).$$
   On sait qu'en fait, cette image est contenue dans le sous-espace $I_{cusp}(\mathfrak{g}(F))$, cf. \cite{W6} proposition 10.

    La transformation de Fourier se descend en une transformation de l'espace $I(\mathfrak{g}(F))$. L'espace $FC(\mathfrak{g}(F))$ est le sous-espace des $f\in I(\mathfrak{g}(F))$ telles que $I^G(X,f)=I^G(X,\hat{f})=0$ pour tout $X\in \mathfrak{g}(F)-\mathfrak{g}_{tn}(F)$., cf. \cite{W6} proposition 10.
    
    Le groupe $G(F)$ agit par conjugaison sur l'espace $FC(\mathfrak{g}_{SC}(F))$. On note $FC^G(\mathfrak{g}_{SC}(F))$ le sous-espace des invariants. 
   
Soit $s\in Imm(G_{AD})$. En g\'en\'eral $s$ n'est pas un sommet de $Imm_{F^{nr}}(G_{AD})$. C'est seulement l'unique point fixe par l'action de Frobenius dans une facette de cet immeuble qui est conserv\'ee par cette action. Toutefois, on a la propri\'et\'e suivante.
\begin{lem}{Soit $s\in Imm(G_{AD})$. Supposons que   ${\bf FC}(\mathfrak{g}_{s})\not=\emptyset$. Alors $s$ est un sommet de $Imm_{F^{nr}}(G_{AD})$.}\end{lem}
Preuve.   Si $s $ appartient \`a une facette de dimension strictement positive, le groupe $Z(G_{s})^0$ est non trivial et ${\bf FC}(\mathfrak{g}_{s})=\emptyset$ comme on l'a dit en \ref{faisceauxcaracterescuspidaux}.  $\square$

   \subsubsection{Les espaces $I_{cusp}^{st}(\mathfrak{g}(F))$ et $FC^{st}(\mathfrak{g}(F))$\label{FCstable}}
    
   Fixons  un ensemble de repr\'esentants   $Endo_{ell}(G)$ des classes d'\'equivalence de donn\'ees endoscopiques elliptiques de $G$, cf. \cite{MW} I.1.5. A  chaque telle donn\'ee   ${\bf G}'$ est associ\'e un groupe endoscopique $G'$ et  on fixe un facteur de transfert $\Delta^{{\bf G}'}$ sur $\mathfrak{g}'(F)\times \mathfrak{g}(F)$. Cela permet de d\'efinir le transfert endoscopique
   $$transfert^{{\bf G}'}:I(\mathfrak{g}(F))\to SI(\mathfrak{g}'(F)).$$
   On note $I_{cusp}(\mathfrak{g}(F),{\bf G}')$ le sous-espace des $f\in I_{cusp}(\mathfrak{g}(F))$ tels que $transfert^{{\bf G}''}(f)=0$ pour tout ${\bf G}''\in Endo_{ell}(G)$, ${\bf G}''\not={\bf G}'$. D'apr\`es \cite{MW} proposition I.4.11, on a l'\'egalit\'e
   $$I_{cusp}(\mathfrak{g}(F))=\oplus_{{\bf G}'\in Endo_{ell}(G)}I_{cusp}(\mathfrak{g}(F),{\bf G}').$$
   Cette d\'ecomposition est stable par transformation de Fourier. 
  
  Pour ${\bf G}'\in Endo_{ell}(G)$, on pose $FC(\mathfrak{g}(F),{\bf G}')=FC(\mathfrak{g}(F))\cap I_{cusp}(\mathfrak{g}(F),{\bf G}')$. D'apr\`es \cite{W6} proposition 12, on a l'\'egalit\'e
$$FC(\mathfrak{g}(F))=\oplus_{{\bf G}'\in Endo_{ell}(G)}FC(\mathfrak{g}(F),{\bf G}').$$
   
   Il y a une unique donn\'ee endoscopique elliptique de $G$, \`a \'equivalence pr\`es, dont le groupe endoscopique associ\'e est une forme int\'erieure quasi-d\'eploy\'ee $G^*$ de $G$.   On l'appelle la donn\'ee principale et on la note ${\bf G}$. On pose simplement $I^{st}_{cusp}(\mathfrak{g}(F))=I_{cusp}(\mathfrak{g}(F),{\bf G})$ et $FC^{st}(\mathfrak{g}(F))=FC(\mathfrak{g}(F),{\bf G})$. Si $G$ est semi-simple,  l'espace $FC^{st}(\mathfrak{g}(F))$ est insensible aux isog\'enies, c'est-\`a-dire est le m\^eme pour $G$, $G_{SC}$ et $G_{AD}$. L'application naturelle $I_{cusp}^{st}(\mathfrak{g}(F))\to SI_{cusp}(\mathfrak{g}(F))$ est bijective. 
   
       \subsubsection{D\'ecomposition de l'espace $FC^{st}(\mathfrak{g}(F))$\label{decompositionFCstable}}
       Fixons un ensemble $\underline{S}(G)$ de repr\'esentants des orbites de l'action de $G(F)$ dans $S(G)$. D'apr\`es \cite{W6} 10(2), l'application naturelle 
   $$(1) \qquad \oplus_{s\in \underline{S}(G)}FC(\mathfrak{g}_{s}({\mathbb F}_{q}))^{K_{s}^{\dag}}\to FC(\mathfrak{g}(F))$$
   est bijective. Soit $s\in S(G)$. Puisqu'on peut choisir $\underline{S}(G)$ contenant $s$, l'application naturelle
  $FC(\mathfrak{g}_{s}({\mathbb F}_{q}))^{K_{s}^{\dag}} \to FC(\mathfrak{g}(F))$ est injective.  On note $FC^{st}(\mathfrak{g}_{s}({\mathbb F}_{q}))$ le sous-espace des \'el\'ements de $FC(\mathfrak{g}_{s}({\mathbb F}_{q}))^{K_{s}^{\dag}}$ dont l'image dans $FC(\mathfrak{g}(F))$ appartient \`a $FC^{st}(\mathfrak{g}(F))$. 
  
   Soit $s\in S(G)$. Consid\'erons l'orbite de $s$ pour l'action de $G(F)$ sur $Imm(G_{AD})$. Elle se d\'ecompose en un nombre fini d'orbites pour l'action de $G_{SC}(F)$. Fixons un ensemble de repr\'esentants $\underline{S}(s)$  de ces derni\`eres orbites. Pour $s'\in \underline{S}(s)$, fixons $g_{s'}\in G(F)$ tel que $g_{s'}s=s'$. De $Ad(g_{s'})$ se d\'eduit un isomorphisme $C^{G_{SC,s}}_{nil,cusp}(\mathfrak{g}_{SC,s}({\mathbb F}_{q}))\to C^{G_{SC,s'}}_{nil,cusp}(\mathfrak{g}_{SC,s'}({\mathbb F}_{q}))$ encore not\'e $Ad(g_{s'})$. Posons $\underline{S}(G_{SC})=\cup_{s\in \underline{S}(G)}\underline{S}(s)$. Alors $\underline{S}(G_{SC})$ est un ensemble de repr\'esentants des orbites de l'action de $G_{SC}(F)$ dans $S(G)$. On a l'isomorphisme analogue \`a (1):
  $$(2) \qquad  \oplus_{s'\in \underline{S}(G_{SC})}FC(\mathfrak{g}_{SC,s'}({\mathbb F}_{q})) \to FC(\mathfrak{g}_{SC}(F)).$$
  On a supprim\'e les exposants $K_{SC,s'}^{\dag}$ car, $G_{SC}$ \'etant simplement connexe, on a $K_{SC,s'}^{\dag}=K_{SC,s'}^0$ et l'invariance par ce groupe est automatique.   
   
   Supposons $A_{G}^{nr}=\{1\}$.  Pour tout $s\in S(G)$, on a l'\'egalit\'e $\mathfrak{g}_{s}=\mathfrak{g}_{SC,s}$. L'espace $C^{G_{s}}_{nil,cusp}(\mathfrak{g}_{s}({\mathbb F}_{q}))$ s'identifie \`a l'espace d'invariants $(C^{G_{SC,s}}_{nil, cusp}(\mathfrak{g}_{SC,s}({\mathbb F}_{q})))^{K^0_{s}}$. De m\^eme, d'apr\`es les d\'efinitions de ces espaces, $FC(\mathfrak{g}_{s}({\mathbb F}_{q}))$ s'identifie \`a l'espace d'invariants $FC(\mathfrak{g}_{SC,s}({\mathbb F}_{q}))^{K^0_{s}}$, donc $FC(\mathfrak{g}_{s}({\mathbb F}_{q}))^{K_{s}^{\dag}}$ s'identifie \`a $FC(\mathfrak{g}_{SC,s}({\mathbb F}_{q}))^{K^{\dag}_{s}}$. Pour $s\in \underline{S}(G)$, d\'efinissons
   $$(3) \qquad j_{s}:FC(\mathfrak{g}_{s}({\mathbb F}_{q}))^{K_{s}^{\dag}}\to \oplus_{s'\in \underline{S}(s)}FC(\mathfrak{g}_{SC,s'}({\mathbb F}_{q}))$$
   comme la compos\'ee de l'application pr\'ec\'edente et de l'application $\oplus_{s'\in \underline{S}(s)}Ad(g_{s'})$, le tout multipli\'e par $\vert \underline{S}(s)\vert ^{-1}$. On d\'efinit
   $$j:\oplus_{s\in \underline{S}(G)}FC(\mathfrak{g}_{s}({\mathbb F}_{q}))^{K_{s}^{\dag}}\to  \oplus_{s'\in \underline{S}(G_{SC})}FC(\mathfrak{g}_{SC,s'}({\mathbb F}_{q})) $$
   comme la somme des $j_{s}$ pour $s\in \underline{S}(G)$. Par (1) et (2), $j$ s'identifie \`a un homomorphisme $FC(\mathfrak{g}(F))\to FC(\mathfrak{g}_{SC}(F))$. L'espace $G(F)$ agit naturellement sur $FC(\mathfrak{g}_{SC}(F))$ et, par construction, l'image de $j$ est \'egale au sous-espace des invariants $FC(\mathfrak{g}_{SC}(F))^{G(F)}$. 
    
 \begin{lem}{(i) Soit $s\in S(G)$ tel que $FC^{st}(\mathfrak{g}_{s}({\mathbb F}_{q}))\not=\{0\}$. Alors les orbites de $s$  pour  les actions de $G(F)$ et de  $G_{AD}(F)$ sur $Imm(G_{AD})$ sont \'egales. 
 
 (ii) On a l'\'egalit\'e $FC^{st}(\mathfrak{g}(F))=\oplus_{s\in \underline{S}(G)}FC^{st}(\mathfrak{g}_{s}({\mathbb F}_{q}))$.
 
 (iii) Supposons $A_{G}^{nr}=\{1\}$. Alors l'homomorphisme $j$ se restreint en un isomorphisme $FC^{st}(\mathfrak{g}(F))\to FC^{st}(\mathfrak{g}_{SC}(F))$. }\end{lem}
 
 On d\'emontrera ce lemme en \ref{demonstration}. 
 
   \subsubsection{Immeubles et  restriction \`a la Weil}\label{Weil}
  Consid\'erons une extension finie mod\'er\'ement ramifi\'ee $F'$ de $F$ et un  groupe r\'eductif connexe $G'$ d\'efini sur $F$. Supposons que $G=Res_{F'/F}(G')$, c'est-\`a-dire que $G$ est la restriction \`a la Weil de $G'$.  
  
  Supposons d'abord que $F'/F$ soit totalement ramifi\'ee. On a l'égalité $Imm_{F^{nr}}(G_{AD})=Imm_{F^{'nr}}(G'_{AD})$, cette  identification \'etant \'equivariante pour les actions de $\Gamma_{F}^{nr}\simeq \Gamma_{F'}^{nr}$. Donc aussi $Imm(G_{AD})=Imm_{F'}(G'_{AD})$. Pour tout sommet $s$ de cet immeuble commun, on a $\mathfrak{g}_{s}=\mathfrak{g}'_{s}$, cette  identification \'etant \'equivariante pour les actions de $\Gamma_{{\mathbb F}_{q}}$. L'\'egalit\'e $FC(\mathfrak{g}_{s}({\mathbb F}_{q}))=FC(\mathfrak{g}'_{s}({\mathbb F}_{q}))$ r\'esulte de la construction de ces espaces. On a aussi $\mathfrak{g}(F)=\mathfrak{g}'(F')$ et  $SI(\mathfrak{g}(F))=SI(\mathfrak{g}'(F'))$. On en d\'eduit que $FC(\mathfrak{g}(F))=FC(\mathfrak{g}'(F'))$ et $FC^{st}(\mathfrak{g}(F))=FC^{st}(\mathfrak{g}'(F'))$. Pour tout $s\in S(G)$,  on a $FC^{st}(\mathfrak{g}_{s}({\mathbb F}_{q}))=FC^{st}(\mathfrak{g}'_{s}({\mathbb F}_{q}))$. Plus g\'en\'eralement, soit $E$ une extension finie non ramifi\'ee de $F$ de degr\'e $d$. On a $FC^{st}(\mathfrak{g}(E))=FC^{st}(\mathfrak{g}'(EF'))$ et, pour tout $s\in S_{E}(G)$, on a $FC^{st}(\mathfrak{g}_{s}({\mathbb F}_{q^d}))=FC^{st}(\mathfrak{g}'_{s}({\mathbb F}_{q^d}))$. 
  
  Passons au cas g\'en\'eral.  Notons $f$ le degr\'e r\'esiduel de l'extension $F'/F$. L'immeuble $Imm_{F^{'nr}}(G'_{AD})$ est muni d'une action, non pas du Frobenius $Fr$, mais de $Fr^f$. L'immeuble $Imm_{F^{nr}}(G_{AD})$ s'identifie \`a $Imm_{F^{'nr}}(G'_{AD})^f$, l'action du Frobenius \'etant
  $$(x_{1},...,x_{f})\mapsto (Fr^f(x_{f}),x_{1},...,x_{f-1}).$$
  Un sommet $s$ de $Imm(G_{AD})$ s'identifie \`a $(s',...,s')$, o\`u $s'$ est un sommet de $ Imm_{F'}(G'_{AD})$. On a $\mathfrak{g}_{s}=(\mathfrak{g}'_{s'})^f$, l'action du Frobenius \'etant donn\'ee par la m\^eme formule que ci-dessus. D'o\`u $\mathfrak{g}_{s}({\mathbb F}_{q})=\mathfrak{g}'_{s'}({\mathbb F}_{q^f})$. De nouveau, l'\'egalit\'e $FC(\mathfrak{g}_{s}({\mathbb F}_{q}))=FC(\mathfrak{g}'_{s'}({\mathbb F}_{q^f}))$ r\'esulte de la construction de ces espaces. On a encore $SI(\mathfrak{g}(F))=SI(\mathfrak{g}'(F'))$ et on en d\'eduit l'\'egalit\'e  $FC^{st}(\mathfrak{g}_{s}({\mathbb F}_{q}))=FC^{st}(\mathfrak{g}'_{s'}({\mathbb F}_{q^f}))$. Plus g\'en\'eralement, soit $E$ une extension finie non ramifi\'ee de $F$ de degr\'e $d$. Notons $b$ le pgcd de $f$ et $d$. On voit que l'application $x=(x_{1},...,x_{f})\mapsto (x_{1},...,x_{b})$ identifie $Imm_{E}(G_{AD})$ \`a $Imm_{F'E}(G'_{AD})^b$. Pour un sommet $s=(s_{1},...,s_{b})$ de $Imm_{E}(G_{AD})$, on a $\mathfrak{g}_{s}({\mathbb F}_{q^d})=\oplus_{i=1,...,b}\mathfrak{g}'_{s_{i}}({\mathbb F}_{q^{fd/b}})$. On en d\'eduit  $FC^{st}(\mathfrak{g}_{s}({\mathbb F}_{q^d}))=\otimes_{i=1,...,b}FC^{st}(\mathfrak{g}'_{s_{i}}({\mathbb F}_{q^{fd/b}}))$. L'espace  $FC^{st}(\mathfrak{g}(E))$ est le produit tensoriel de $b$ copies de $FC^{st}(\mathfrak{g}'(EF'))$.

\subsubsection{Action d'automorphismes\label{lemmeautomorphismes}}
Dans ce paragraphe, on suppose que $G$ est absolument quasi-simple. Soit $\delta$ un automorphisme de $G$ d\'efini sur $F^{nr}$. Il s'en d\'eduit un automorphisme de $Imm_{F^{nr}}(G)$ encore not\'e $\delta$.

\begin{lem}{Soit $s\in S(G)$ tel que $FC^{st}(\mathfrak{g}_{s}({\mathbb F}_{q}))\not=\{0\}$. Supposons que $\delta(s)$ appartienne \`a $Imm(G_{AD})$. Alors il existe $g\in G_{SC}(F)$ tel que $\delta(s)=gs$.}\end{lem}

Preuve. On fixe des objets comme en \ref{alcove}. On ne perd rien \`a supposer que $s\in S(\bar{C})$. Fixons une alc\^ove $C'$ de $Imm(G_{AD})$ telle que $\delta(s)\in \bar{C}'$. On peut fixer $g\in G_{SC}(F)$ tel que $C'=g(C)$. L'alc\^ove $C'$ est l'ensemble des points fixes par $\Gamma_{F}^{nr}$ d'une alc\^ove $C^{'nr}$ de $Imm_{F^{nr}}(G_{AD})$ et on a aussi $C^{'nr}=g(C^{nr})$. Posons $C^{''nr}=\delta(C^{nr})$. C'est encore une alc\^ove de $Imm_{F^{nr}}(G_{AD})$ et, puisque $s\in \bar{C}^{nr}$, on a $\delta(s)\in \bar{C}^{''nr}$. L'ensemble des alc\^oves de $Imm_{F^{nr}}(G_{AD})$ dont l'adh\'erence contient $\delta(s)$ est en bijection avec celui des sous-groupes de Borel de $G_{\delta(s)}$. Ce dernier forme une unique orbite pour l'action de $G_{\delta(s)}$ et on en d\'eduit que l'ensemble d'alc\^oves pr\'ec\'edent forme une unique orbite pour l'action de $K_{\delta(s)}^{0,nr}$. On peut donc fixer $k\in K_{\delta(s)}^{0,nr}$ tel que $k(C^{''nr})=C^{'nr}$. Posons $\delta'=Ad(g^{-1}k)\circ \delta$. Alors $\delta'$ conserve $C^{nr}$. L'action de $\delta'$ sur $C^{nr}$ se fait par un \'el\'ement de $Aut({\cal D}_{a}^{nr})$. D'apr\`es l'hypoth\`ese que $FC^{st}(\mathfrak{g}_{s}({\mathbb F}_{q}))\not=\{0\}$ et d'apr\`es \cite{W7} 9(4), $s$ est fixe par cette action. Il en r\'esulte que $Ad(k)\circ \delta(s)=gs$. Puisque $k\in K_{\delta(s)}^{0,nr}$, on a aussi $Ad(k)\circ \delta(s)=\delta(s)$, d'o\`u $\delta(s)=gs$. $\square$

  \subsubsection{D\'emonstration du lemme \ref{decompositionFCstable}}\label{demonstration}

  Supposons d'abord $G$ simplement connexe.  L'assertion (iii) est triviale. On peut identifier $G$ \`a un produit de groupes de la forme $Res_{F'/F}(G')$, o\`u $F'/F$ est une extension  finie mod\'er\'ement ramifi\'ee de $F$ (gr\^ace \`a l'hypoth\`ese sur $p$) et $G'$ est un groupe r\'eductif d\'efini sur $F'$ qui est simplement connexe et absolument quasi-simple. Il r\'esulte de \ref{Weil} que l'immeuble $Imm(G_{AD})$ s'identifie au produit  des $Imm_{F'}(G'_{AD})$ et $FC(\mathfrak{g}(F))$, resp. $FC^{st}(\mathfrak{g}(F))$, s'identifie au produit tensoriel des $FC(\mathfrak{g}'(F'))$, resp. $FC^{st}(\mathfrak{g}'(F'))$.  On voit qu'il suffit de d\'emontrer le lemme pour chaque composante $G'$ o\`u l'on remplace $F$ par $F'$. Cela nous ram\`ene au cas o\`u $G$ est absolument quasi-simple. L'assertion (ii) est alors \cite{W7} 9.1. L'assertion (i) r\'esulte du lemme \ref{lemmeautomorphismes} appliqu\'e aux automorphismes adjoints.

 Passons au cas g\'en\'eral. 
  Si $A_{G}^{nr}\not=\{1\}$, on a $Z(G_{s})^0\not=\{1\}$ pour tout sommet $s\in S(G)$ donc $FC(\mathfrak{g}_{s}({\mathbb F}_{q}))=\{0\}$ d'apr\`es \ref{faisceauxcaracterescuspidaux} et tout est nul. On suppose d\'esormais  $A_{G}^{nr}=\{1\}$. De la d\'ecomposition $  \mathfrak{g}_{SC}\oplus  \mathfrak{z}(G)=\mathfrak{g}$ se d\'eduit un isomorphisme $ C_{c}^{\infty}(\mathfrak{g}_{SC}(F))\otimes C_{c}^{\infty}(\mathfrak{z}(G)(F))\simeq C_{c}^{\infty}(\mathfrak{g}(F))$. Le groupe $G(F)$ agit sur $ C_{c}^{\infty}(\mathfrak{g}_{SC}(F))$ et les espaces qui s'en d\'eduisent et de l'isomorphisme ci-dessus se d\'eduit un  isomorphisme
  $$(1) \qquad I_{cusp}(\mathfrak{g}_{SC}(F))^{G(F)}\otimes C_{c}^{\infty}(\mathfrak{z}(G)(F))\simeq I_{cusp}(\mathfrak{g}(F)).$$
  L'action de $G(F)$ sur $\mathfrak{g}_{SC}(F)$ pr\'eservant les classes de conjugaison stable, l'action de $G(F)$ sur $SI_{cusp}(\mathfrak{g}_{SC}(F))$ est triviale et   l'isomorphisme (1) se restreint en un isomorphisme
  $$(2) \qquad I^{st}_{cusp}(\mathfrak{g}_{SC}(F))\otimes C_{c}^{\infty}(\mathfrak{z}(G)(F))\simeq I^{st}_{cusp}(\mathfrak{g}(F)).$$
  Soit $s\in S(G)$. Comme on l'a dit en \ref{decompositionFCstable}, parce que $A_{G}^{nr}=\{1\}$,  l'espace $C^{G_{s}}_{nil,cusp}(\mathfrak{g}_{s}({\mathbb F}_{q}))$ s'identifie \`a l'espace d'invariants $C^{G_{SC,s}}_{nil, cusp}(\mathfrak{g}_{SC,s}({\mathbb F}_{q}))^{K^0_{s}}$.  Notons $f_{Z}$ la fonction caract\'eristique de l'ensemble $\mathfrak{z}(G)_{tn}(F)$, notons
  $proj: I_{cusp}(\mathfrak{g}_{SC}(F))\to I_{cusp}(\mathfrak{g}_{SC}(F))^{G(F)}$ la projection naturelle et $\iota:I_{cusp}(\mathfrak{g}_{SC}(F))^{G(F)}\to I_{cusp}(\mathfrak{g}(F))$ la compos\'ee   de l'application $f\mapsto f\otimes f_{Z}$ et de l'isomorphisme (1). L'application $\iota$ est injective. 
   Le diagramme suivant est commutatif:
  $$\begin{array}{ccc}C^{G_{s}}_{nil,cusp}(\mathfrak{g}_{s}({\mathbb F}_{q}))& \to& I_{cusp}(\mathfrak{g}(F))\\ \downarrow &&\,\uparrow \iota\circ proj\\ C^{G_{SC,s}}_{nil,cusp}(\mathfrak{g}_{SC,s}({\mathbb F}_{q}))&\to&I_{cusp}(\mathfrak{g}_{SC}(F))\\ \end{array}$$
  Il r\'esulte alors des constructions que le diagramme suivant est aussi commutatif
  $$\begin{array}{ccccc}&&FC(\mathfrak{g}(F))&\to&I_{cusp}(\mathfrak{g}(F))\\ &&\,\downarrow j&&\,\uparrow \iota\\ FC(\mathfrak{g}_{SC}(F))&\supset&FC(\mathfrak{g}_{SC}(F))^{G(F)}&\to&I_{cusp}(\mathfrak{g}_{SC}(F))^{G(F))}\\ \end{array}$$
  Par d\'efinition de $FC^{st}(\mathfrak{g}(F))$ et $FC^{st}(\mathfrak{g}_{SC}(F))$ et puisque l'isomorphisme (1) se restreint en l'isomorphisme (2), l'isomorphisme $j$ se restreint en un isomorphisme de $FC^{st}(\mathfrak{g}(F))$ sur $FC(\mathfrak{g}_{SC}(F))^{G(F)}\cap FC^{st}(\mathfrak{g}_{SC}(F))$. Mais tout \'el\'ement de $FC^{st}(\mathfrak{g}_{SC}(F))$ est forc\'ement invariant par $G(F)$ donc l'espace pr\'ec\'edent n'est autre que $FC^{st}(\mathfrak{g}_{SC}(F))$. Cela d\'emontre l'assertion (iii) du lemme \ref{decompositionFCstable}. 
  
  Soit $f\in FC(\mathfrak{g}(F))$. Ecrivons $f=\oplus_{s\in \underline{S}(G)}f_{s}$ conform\'ement \`a la d\'ecomposition \ref{decompositionFCstable}(1). Pour $s\in \underline{S}(G)$, \'ecrivons $j_{s}(f_{s})=\oplus_{s'\in \underline{S}(s)}f'_{s'}$ selon \ref{decompositionFCstable}(3). Alors $j(f)=\oplus_{s'\in \underline{S}(G_{SC})}f'_{s'}$. Appliquons le (ii) de l'\'enonc\'e d\'ej\`a d\'emontr\'e  pour $G_{SC}$. Alors 
  $$f\in FC^{st}(\mathfrak{g}(F)) \iff j(f)\in FC^{st}(\mathfrak{g}_{SC}(F))\iff\,  \forall s'\in \underline{S}(G_{SC}), \,
   f_{s'}\in FC^{st}(\mathfrak{g}_{SC}(F)) $$
  $$ \iff  \,\forall s\in \underline{S}(G),\, j(f_{s})\in FC^{st}(\mathfrak{g}_{SC}(F))  \iff \,\forall s\in \underline{S}(G),\,f_{s}\in FC^{st}(\mathfrak{g}(F)) $$
  $$\iff \,\forall s\in \underline{S}(G),\, f_{s}\in FC^{st}(\mathfrak{g}_{s}({\mathbb F}_{q})).$$
  Cela prouve (ii). 
  
  Enfin, soit $s\in \underline{S}(G)$, supposons que les orbites de $s$ pour les actions de $G(F)$ et de $G_{AD}(F)$ sur $Imm(G_{AD})$ ne sont pas \'egales. Alors, pour tout $s'\in \underline{S}(s)$, les orbites de $s'$   pour les actions de $G_{SC}(F)$ et de $G_{AD}(F)$ sur $Imm(G_{AD})$ ne sont pas \'egales. D'apr\`es le (i) de l'\'enonc\'e d\'ej\`a d\'emontr\'e pour $G_{SC}$, on a $FC^{st}(\mathfrak{g}_{SC,s'}({\mathbb F}_{q}))=\{0\}$ pour tout $s'\in \underline{S}(s)$. Les m\^emes \'equivalences que ci-dessus montrent que $FC^{st}(\mathfrak{g}_{s}({\mathbb F}_{q}))=\{0\}$. Cela prouve (i). $\square$

  \subsubsection{Extensions non ramifi\'ees et espaces $FC^{st}(\mathfrak{g}(F))$}\label{extensionsnonram}

\begin{lem}{Soit   $s$ un sommet de $Imm(G_{AD})$ et soit $E$ une extension   non ramifi\'ee de $F$ de degr\'e fini $d$. Si $FC^{st}(\mathfrak{g}_{s}({\mathbb F}_{q}))\not=\{0\}$, alors $s$ reste un sommet dans $Imm_{E}(G_{AD})$ et on a $FC^{st}(\mathfrak{g}_{s}({\mathbb F}_{q^d}))\not=\{0\}$. A fortiori, si $FC^{st}(\mathfrak{g}(F))\not=\{0\}$, alors $FC^{st}(\mathfrak{g}(E))\not=\{0\}$. }\end{lem}

Preuve. La derni\`ere assertion r\'esulte de la premi\`ere et du lemme \ref{decompositionFCstable}.

D\'emontrons la premi\`ere assertion.  L'hypoth\`ese   $FC^{st}(\mathfrak{g}_{s}({\mathbb F}_{q}))\not=\{0\}$ implique  que $s$ reste un sommet dans $Imm_{F^{nr}}(G_{AD})$, cf. lemme \ref{lespaceFC}, et que $A_{G}^{nr}=\{1\}$. D'apr\`es le (iii) du lemme \ref{decompositionFCstable}, on peut remplacer $G$ par $G_{SC}$. On peut d\'ecomposer $G$ en produit $G_{1}\times...\times G_{n}$ o\`u, pour tout $i=1,...,n$, $G_{i}=Res_{F_{i}/F}(G'_{i})$, $F_{i}$ \'etant une extension finie mod\'er\'ement ramifi\'ee de $F$ et $G'_{i}$ \'etant un groupe d\'efini sur $F'_{i}$, simplement connexe et absolument quasi-simple. L'immeuble $Imm_{F^{nr}}(G_{AD})$ s'identifie 	au produit des $Imm_{F^{nr}}(G_{i,AD})$ et $s$ s'identifie \`a une famille de sommets $s_{i}\in Imm_{F^{nr}}(G_{i,AD})$, qui sont fix\'es par $\Gamma_{F}^{nr}$. On a $FC^{st}(\mathfrak{g}_{s}({\mathbb F}_{q}))=\otimes_{i=1,...,n}FC^{st}(\mathfrak{g}_{s_{i}}({\mathbb F}_{q}))$ et une \'egalit\'e analogue vaut en rempla\c{c}ant $q$ par $q^d$. Alors il suffit   de traiter chaque groupe $G_{i}$. Autrement dit, on peut supposer $i=1$. On pose simplement $F'=F_{1}$ et $G'=G'_{1}$. Supposons d'abord $F'=F$. Alors l'assertion du lemme r\'esulte des descriptions de  \cite{W7} paragraphe 9. Traitons le cas g\'en\'eral.  Notons $f$ le degr\'e r\'esiduel de $F'/F$ et $b$ le pgcd de $d$ et $f$. On a vu en \ref{Weil} que $s$ s'identifiait \`a un sommet $s'\in Imm_{F'}(G'_{AD})$ et que l'on avait des isomorphismes
$FC^{st}(\mathfrak{g}_{s}({\mathbb F}_{q}))\simeq FC^{st}(\mathfrak{g}'_{s'}({\mathbb F}_{q^f}))$ et $FC^{st}(\mathfrak{g}_{s}({\mathbb F}_{q^d}))\simeq FC^{st}(\mathfrak{g}'_{s'}({\mathbb F}_{q^{df/b}}))^{\otimes b}$. L'hypoth\`ese entra\^{\i}ne que $FC^{st}(\mathfrak{g}'_{s'}({\mathbb F}_{q^f}))\not=\{0\}$. Puisque $G'$ est absolument quasi-simple, on vient de voir que cela entra\^{\i}ne $FC^{st}(\mathfrak{g}'_{s'}({\mathbb F}_{q^{df/b}}))\not=\{0\}$, donc aussi $FC^{st}(\mathfrak{g}_{s}({\mathbb F}_{q^d}))\not=\{0\}$. $\square$

 On note $S^{st}(G)$ l'ensemble des sommets $s\in S(G)$ tels que $FC^{st}(\mathfrak{g}_{s}({\mathbb F}_{q}))\not=\{0\}$.
On a $S^{st}(G)\subset S_{F^{nr}}(G)$. Le lemme ci-dessus montre que, si $E$ est une extension finie non ramifi\'ee de $F$, on a $S^{st}(G)\subset S_{E}^{st}(G)$. On note $S^{nr,st}(G)$ la r\'eunion des $S^{st}_{E}(G)$ quand $E$ d\'ecrit les 
extensions finies non ramifi\'ees de $F$.  

{\bf Remarque.} On peut définir $S^{nr,st}(G)$ si $G$ est défini seulement sur $F^{nr}$ au lieu de $F$. En effet, la structure de $G$ sur $F^{nr}$ se prolonge en une structure de $G$ sur une extension finie non ramifiée  $E$ de $F$. En remplaçant $F$ par $E$, l'ensemble $S^{nr,st}(G)$  est alors défini. L'extension $E$ et la structure de $G$ sur $E$ ne sont pas uniques mais les ensembles $S^{nr,st}(G)$  déduits de choix différents sont égaux. 
 
\bigskip

Supposons $G$ absolument quasi-simple. On a

(1) soit $s\in S^{nr,st}(G)$; alors il existe $g\in G_{SC}(F^{nr})$ tel que $gs\in Imm(G_{AD})$. 

Comme en \ref{alcove}, on fixe un torseur int\'erieur $\psi:G\to G^*$ dont le cocycle associ\'e $\underline{n}_{G}$ conserve l'alc\^ove $C^{nr}$. Alors $C^{nr}$ est une alc\^ove de $Imm_{F^{nr}}(G_{AD})$ et l'action galoisienne sur $S(\bar{C}^{nr})$ se fait par des \'el\'ements de $Aut({\cal D}_{a}^{nr})$. A conjugaison pr\`es par un \'el\'ement de $G_{SC}(F^{nr})$, on peut supposer que $s\in S(\bar{C}^{nr})$. Alors, d'apr\`es \cite{W7} 9(4), l'hypoth\`ese $s\in S^{nr,st}(G)$
entra\^{\i}ne que $s$ est fix\'e par tout \'el\'ement de $Aut({\cal D}_{a}^{nr})$, donc par l'action galoisienne. Cela prouve (1).

\subsection{Etude des couples $(M,s)$}

\subsubsection{Réduction au cas absolument simple}\label{reduction}

Pour cette sous-section, on suppose que $G$ est défini sur  $F^{nr}$. On fixe un \'epinglage $\mathfrak{E}=(B,T,(E_{\alpha})_{\alpha\in \Delta})$ de $G$ défini sur $F^{nr}$.  On utilise les constructions de \ref{alcove} et \ref{racinesaffines}.     Dans certains paragraphes, on supposera (et cela sera dit explicitement) que $G$   est défini et quasi-déployé sur $F$. On suppose alors que $\mathfrak{E}$ est défini sur $F$ et on pose $C=C^{nr}\cap Imm(G_{AD})$. C'est une alc\^ove de cet immeuble. 

 Dans ce paragraphe, $G$ est défini sur $F^{nr}$.   
Soient $M$ un $F^{nr}$-Levi et $s$ un sommet de $Imm_{F^{nr}}(M_{AD})$.     
Le groupe $G_{SC}$ se décompose en produit $\prod_{i=1,...,n}G_{i}$, où, pour tout $i$, $G_{i}$ est un groupe simplement connexe défini et irréductible  sur $F^{nr}$. Pour tout $i$, il existe une extension finie et modérément ramifiée $E_{i}$ de $F^{nr}$  et un groupe $G'_{i}$ simplement connexe défini sur $E_{i}$, absolument irréductible, de sorte que $G_{i}=Res_{E_{i}/F^{nr}}(G'_{i})$. Le Levi $M_{sc}$ se décompose conformément en $M_{sc}=\prod_{i=1,...,n}M_{i}$ et, pour tout $i$, il existe un $E_{i}$-Levi $M'_{i}$ de sorte que $M_{i}=Res_{E_{i}/F^{nr}}(M'_{i})$. On a $Imm_{F^{nr}}(M_{AD})=\prod_{i=1,...,n}Imm_{F^{nr}}(M_{i,AD})$ et, pour tout $i$, $Imm_{F^{nr}}(M_{i,AD})=Imm_{E_{i}}(M'_{i,AD})$. Ainsi, $s$ s'identifie à un produit $\prod_{i=1,...,n}s_{i}$, où $s_{i}\in Imm_{E_{i}}(M'_{i,AD})$. Remarquons que, pour tout $i$, il existe une extension finie modérément ramifiée $F_{i}$ de $F$ de sorte que $E_{i}=F_{i}^{nr}$. Les triplets $(G'_{i},M'_{i},s_{i})$ vérifient donc les m\^emes hypothèses que $(G,M,s)$, au changement près du corps de base $F$ en $F_{i}$. 

Rappelons la projection $p_{M}:Imm^G_{F^{nr}}(M_{ad})\to Imm_{F^{nr}}(M_{AD})$.  
 Posons $P(s)=p_{M}^{-1}(s)$. Cet ensemble est un sous-espace affine de $Imm_{F^{nr}}(G_{AD})$ d'espace vectoriel associé ${\cal A}_{M}^{nr}/{\cal A}_{G}^{nr}$.    L'espace $P(s)$ se décompose en $\prod_{i=1,...,n}P(s_{i})$. Pour $x\in P(s)$, on dispose du groupe $G_{x}$, dont $M_{s}$ est un $\bar{{\mathbb F}}_{q}$-Levi. Ecrivons $x=\prod_{i=1,...,n}x_{i}$, où $x_{i}\in P(s_{i})$ pour tout $i$. On a alors l'égalité $G_{x,SC}=\prod_{i=1,...,n}G'_{i,x_{i},SC}$.  

           \subsubsection{Construction de $F^{nr}$-Levi\label{constructiondeLevi}}

 {\bf On suppose   que} $G$ {\bf est absolument quasi-simple}.     
 
Soit $\Lambda$ un sous-ensemble  propre de $\Delta^{nr}_{a}$. Notons $A_{\Lambda}^{nr}$  la composante connexe du groupe des $t\in T^{nr}$ tels que $\beta(t)=1$ pour tout $\beta\in \Lambda$. Notons $M_{\Lambda}$ le commutant de $A_{\Lambda}^{nr}$ dans $G$. C'est un  $F^{nr}$-Levi de $G$. Notons ${\mathbb Q}[\Lambda]$ le sous-${\mathbb Q}$-espace vectoriel de $X^*_{{\mathbb Q}}(T^{nr})$ engendré par $\Lambda$. 
 Montrons que

(1) l'ensemble des racines de $T^{nr}$ dans $\mathfrak{m}_{\Lambda}$ est $\Sigma^{nr}\cap {\mathbb Q}[\Lambda]$; le rang  de $M_{\Lambda,SC}$ sur $F^{nr}$ est \'egal \`a $\vert \Lambda\vert $; on a $A_{\Lambda}^{nr}=A_{M_{\Lambda}}^{nr}$. 

Puisque l'espace des relations lin\'eaires entre \'el\'ements de $\Delta_{a}^{nr}$ est une droite, $A_{\Lambda}^{nr}$ est de dimension $\vert \Delta_{a}^{nr}\vert -\vert \Lambda\vert -1=dim(T^{nr}) -\vert \Lambda\vert$. Son annulateur dans $X^*_{{\mathbb Q}}(T^{nr})$ est un espace de dimension $\Lambda$. Puisque cet espace contient $\Lambda$ par construction, cet annulateur est ${\mathbb Q}[\Lambda]$. La premi\`ere assertion s'en d\'eduit et la deuxi\`eme est cons\'equence de la premi\`ere.  On a  $A_{\Lambda}^{nr}\subset A_{M_{\Lambda}}^{nr}$ par d\'efinition de $M_{\Lambda}$. Ces deux tores \'etant de m\^eme dimension, ils sont \'egaux. Cela prouve (1). 

 L'appartement $App_{F^{nr}}(T^{nr})\subset Imm_{F^{nr}}(G_{AD})$ est contenu dans  $Imm_{F^{nr}}^{G}(M_{\Lambda,ad})$. L'application $p_{M_{\Lambda}}:App_{F^{nr}}(T^{nr})\to Imm_{F^{nr}}(M_{AD})$ est donc bien définie. 

 Pour $x\in \bar{C}^{nr}$, notons $\Lambda(x)$ l'ensemble des $\beta\in \Delta_{a}^{nr}$ telles que $\beta^{aff}(x)=0$. On a dit que de  $T$ et de $C^{nr}$ se d\'eduisait une paire de Borel $(B_{x},T_{x})$ de $G_{x}$. L'ensemble des racines simples  de $G_{x}$ pour cette paire s'identifie \`a $\Lambda(x)$. Posons $\Lambda^{aff}=\{\beta^{aff}; \beta\in \Lambda\}$ et notons $P(\Lambda^{aff})$ l'ensemble des $x\in App_{F^{nr}}(T^{nr})$ tels que $\beta^{aff}(x)=0$ pour tout $\beta^{aff}\in \Lambda^{aff}$. L'intersection $P(\Lambda^{aff})\cap \bar{C}^{nr}$ est l'ensemble des $x\in \bar{C}^{nr}$ tels que $\Lambda\subset \Lambda(x)$. Cette intersection contient les sommets $s_{\beta}$ associ\'es aux \'el\'ements $\beta\in \Delta_{a}^{nr}-\Lambda$: on a $\Lambda(s_{\beta})=\Delta_{a}^{nr}-\{\beta\}$. 
   Montrons que

(2) la projection $p_{M}(P(\Lambda^{aff}))$ est   un unique point $s_{\Lambda}\in Imm_{F^{nr}}(M_{\Lambda,AD})$;
 ce point $s_{\Lambda}$ est un sommet de $Imm_{F^{nr}}(M_{\Lambda,AD})$;  on a l'égalité $P(\Lambda^{aff})=P(s_{\Lambda})$;  pour $x\in P(\Lambda^{aff})\cap \bar{C}^{nr}$, le $\bar{{\mathbb F}}_{q}$-Levi standard de $G_{ x}$  associ\'e au sous-ensemble $\Lambda$ de $\Lambda(x)$ est le groupe $M_{\Lambda,s_{\Lambda}}$.
      
Posons    ${\cal A}^{nr}= X_{*,{\mathbb R}}(T^{nr})$ et ${\cal A}_{\Lambda}^{nr}=X_{*,{\mathbb R}}(A_{M_{\Lambda}}^{nr})$. 
L'appartement $App_{F^{nr}}(T^{nr})$ est un espace affine sous le groupe ${\cal A}^{nr}$.  L'appartement associ\'e \`a $T^{nr}$ dans $Imm_{F^{nr}}(M_{\Lambda,AD} )$ est le quotient de $App_{F^{nr}}(T^{nr})$ par l'action du sous-groupe ${\cal A}_{\Lambda}^{nr}$. D'apr\`es la troisi\`eme assertion de (1), ce groupe est l'annulateur de $\Lambda$ dans ${\cal A}^{nr}$. Consid\'erons deux \'el\'ements $ x,x'\in P(\Lambda^{aff})$. Soit $X\in {\cal A}^{nr}$ tel que $ x'=x+ X$. Pour $\beta\in \Lambda$, on $\beta^{aff}( x)=\beta^{aff}(x')=0$. Donc $\beta(X)=0$. Cela entra\^{\i}ne que $X\in {\cal A}_{\Lambda}^{nr}$.  Donc $x$ et $x'$ ont m\^eme projection dans $Imm_{F^{nr}}(M_{\Lambda,AD})$. Cela d\'emontre la premi\`ere assertion de (2). Soit $x\in P(\Lambda^{aff})\cap \bar{C}^{nr}$, posons $s_{\Lambda}=p_{M_{\Lambda}}(x)$.   D'apr\`es \ref{groupesenreduction}(2), le Levi $M_{\Lambda,s_{\Lambda}}$ de $G_{x}$ est le commutant de la r\'eduction de $A_{\Lambda}^{nr}$ dans $G_{x}$. Il est clair que les \'el\'ements de $\Lambda$ annulent cette r\'eduction. Donc $M_{\Lambda,s_{\Lambda}}$ contient le Levi standard de $G_{x}$ associ\'e au sous-ensemble $\Lambda$ de $ \Lambda(x)$. Donc le rang de $M_{\Lambda,s_{\Lambda},SC}$ est au moins \'egal \`a $\vert \Lambda\vert $. Ce rang est au plus \'egal \`a $\vert \Lambda\vert $ d'apr\`es (1), donc $M_{\Lambda,s_{\Lambda}}$ est \'egal au   Levi standard de $G_{x}$ associ\'e \`a l'ensemble $\Lambda$. Puisque le rang de $M_{\Lambda,s_{\Lambda},SC}$ est \'egal au rang de 
  $M_{\Lambda,SC}$ sur $F^{nr}$,   $s_{\Lambda}$ est un sommet de $Imm_{F^{nr}}(M_{\Lambda,AD})$. Enfin,   par définition de $s_{\Lambda}$, on a $P(\Lambda^{aff})\subset P(s_{\Lambda})$. Ces deux ensembles sont des espaces affines de m\^eme dimension, égale à celle de $A_{\Lambda}^{nr}$. Ils sont donc égaux. Cela prouve (2).

  {\bf Remarque.} Il y a un certain conflit de notations: dans le cas o\`u $\Lambda$ est r\'eduit \`a une racine $\beta$,  le $s_{\Lambda}$ de l'assertion (2) n'est pas le sommet $s_{\beta}$. 
  
    \subsubsection{Couples formés d'un Levi et d'un sommet de son immeuble}\label{couples}

 {\bf On suppose  que }$G$ {\bf est défini   sur} $F$ {\bf et qu'il est absolument quasi-simple.}   Notons $T_{F}$ le plus grand sous-tore de $T$ déployé sur $F$. 
 
\begin{prop}{Soient $M$ un $F$-Levi de $G$, $s$ un sommet de $Imm(M_{AD})$ et $x\in Imm^G(M_{ad})$. Supposons que $s$ reste un sommet dans $Imm_{F^{nr}}(M_{AD})$ et que $p_{M}(x)=s$.  Soit $S\in {\cal T}_{max}^M$ tel que $s$ appartienne \`a l'appartement associ\'e \`a $S$ dans $Imm(M_{AD})$. Alors il existe  un sous-ensemble  propre $\Lambda\subsetneq \Delta^{nr}_{a}$  et un \'el\'ement de $G_{SC}(F)$ qui conjugue $(S,M,s)$ en $(T_{F},M_{\Lambda},s_{\Lambda})$ et qui conjugue $x$ en un \'el\'ement $x'\in P(\Lambda^{aff})\cap \bar{C}$. L'ensemble $\Lambda$ est conserv\'e par l'action galoisienne et  $M_{\Lambda}$ est un $F$-Levi. }\end{prop}

Preuve.  Puisque $M$ est un $F$-Levi, ${\cal T}_{max}^M$ est le sous-ensemble des $T'\in {\cal T}_{max}$ tels que $T'\subset M$.  Le point $x$ appartient \`a l'appartement $App(S)$ associ\'e \`a $S$ dans $Imm(G_{AD})$.    Les \'el\'ements de ${\cal T}_{max}$ forment une unique orbite sous l'action de $G_{SC}(F)$. Quite \`a effectuer une conjugaison par un \'el\'ement de ce groupe, on peut supposer $S=T_{F}$. Tout \'el\'ement de  l'appartement $App(T_{F})$ est conjugu\'e \`a un \'el\'ement de $\bar{C}$ par un \'el\'ement de $Norm_{G_{SC}(F)}(T_{F})$. Quitte \`a effectuer une conjugaison par un \'el\'ement de ce groupe, on peut supposer non seulement que $S=T_{F}$ mais que, de plus, $x\in \bar{C}$.  Puisque $M$ est le commutant de $A_{M}\subset T_{F}$ et que $T$ commute à $T_{F}$,  on a $T\subset M$. En appliquant \ref{groupesenreduction}(2), on introduit le sous-tore $A_{M,x}$ de $T_{x}$ et son commutant   dans $G_{x}$. Celui-ci est \'egal \`a $M_{s}$.      Tout ${\mathbb F}_{q}$-Levi de $G_{x}$ contenant $T_{x}$ est conjugu\'e \`a un ${\mathbb F}_{q}$-Levi standard par un \'el\'ement de $Norm_{G_{x,SC}({\mathbb F}_{q})}(T_{x})$. Un \'el\'ement de ce groupe se rel\`eve en un \'el\'ement de $K_{SC,x}^{0}$ qui normalise $T_{F}$. Quitte \`a effectuer une conjugaison par un \'el\'ement de ce groupe, on peut donc supposer de plus que $M_{s}$ est un  ${\mathbb F}_{q}$-Levi standard de $G_{x}$. Celui-ci correspond \`a un sous-ensemble $\Lambda$ de l'ensemble $ \Lambda(x)$ des racines simples de $G_{x}$.  Puisque $M_{s}$ est un ${\mathbb F}_{q}$-Levi,   $\Lambda$ est conserv\'e par l'action galoisienne. On introduit le sous-tore $A^{nr}_{\Lambda}$ de $T^{nr}$, cf. \ref{constructiondeLevi}. Le groupe $Z(M_{s})^0$ est la composante neutre de l'annulateur de $\Lambda$ dans $T_{x}$ et est donc \'egal \`a $A_{\Lambda,x}^{nr}$. Puisque $A_{M,x}$ est contenu dans $ Z(M_{s})^0$, on en d\'eduit $A_{M,x}\subset A_{\Lambda,x}^{nr}$, ce qui se rel\`eve en $A_{M}\subset A_{\Lambda}^{nr}$. Alors le commutant $M$ de $A_{M}$ contient le commutant $M_{\Lambda}$ de $A_{\Lambda}^{nr}$. Puisque $s$ est un sommet de $Imm_{F^{nr}}(M_{AD})$, le rang de $M_{SC}$ sur $F^{nr}$ est \'egal au rang de $M_{s,SC}$, c'est-\`a-dire $\vert \Lambda\vert $, qui est \'egal au rang de $M_{\Lambda,SC}$ sur $F^{nr}$. L'inclusion $M_{\Lambda}\subset M$ est donc une \'egalit\'e. On a alors $s=p_{M}(x)=p_{M_{\Lambda}}(x)=s_{\Lambda}$. 
 Cela achève la preuve. $\square$ 
 
 Remarquons que la première assertion de l'énoncé vaut pour un $F^{nr}$-Levi $M$, un sommet $s$ de $Imm_{F^{nr}}(M_{AD})$ et un point $x\in Imm^G_{F^{nr}}(M_{ad})$ (en remplaçant  $F$ par $F^{nr}$ dans cette assertion). Il suffit d'appliquer la proposition sur une extension finie $F'$ non ramifiée de $F$ telle que les données $G$, $M$, $s$ et $x$ vérifient les hypothèses de l'énoncé quand on remplace $F$ par $F'$. La m\^eme remarque vaudra pour plusieurs énoncés ci-dessous. 
 
 \subsubsection{Normalisateur}\label{normalisateur}
 
 Posons $T(\mathfrak{o}_{\bar{F}})^{I_{F}}_{0}=T(\mathfrak{o}_{\bar{F}})^{I_{F}}\cap K_{C^{nr}}^{nr,0}$.  Prouvons que:
 
 (1) si $\hat{T}^{I_{F}}$ est connexe, $T(\mathfrak{o}_{\bar{F}})^{I_{F}}_{0}=T(\mathfrak{o}_{\bar{F}})^{I_{F}}$.
 
 Le groupe $T(\mathfrak{o}_{\bar{F}})^{I_{F}}$ agit trivialement sur l'appartement $App_{F^{nr}}(T^{nr})$, donc est inclus dans $K_{C^{nr}}^{nr,\dag}$. L'homomorphisme de Kottwitz $w_{G}: G(F^{nr})\to X^*(Z(\hat{G})^{I_{F}})$ se restreint à $K_{C^{nr}}^{nr,\dag}$ en un homomorphisme 
 dont le noyau est $K_{C^{nr}}^{nr,0}$ (cf. \cite{W4} paragraphe 3). Donc $T(\mathfrak{o}_{\bar{F}})^{I_{F}}_{0}=T(\mathfrak{o}_{\bar{F}})^{I_{F}}\cap Ker(w_{G})$. Sur $T(\bar{F})^{I_{F}}$, $w_{G}$ se factorise par l'homomorphisme $w_{T}:T(\bar{F})^{I_{F}}\to X^*(\hat{T}^{I_{F}})$. Si $\hat{T}^{I_{F}}$ est connexe, ce dernier groupe est sans torsion, donc le groupe compact $T(\mathfrak{o}_{\bar{F}})^{I_{F}}$ est inclus dans le noyau de $w_{T}$, a fortiori dans celui de $w_{G}$. Cela conclut. $\square$
 
 Remarquons que l'hypothèse de (1) est vérifiée si $G$ est simplement connexe ou si $G$ est adjoint, le groupe $\hat{T}$ étant alors induit.

   Introduisons le groupe de Weyl affine $W_{aff}=N_{G(F^{nr})}(T^{nr})/T(\mathfrak{o}_{\bar{F}})^{I_{F}}_{0}$. Il agit proprement sur $App_{F^{nr}}(T^{nr})$ en conservant l'ensemble des racines affines. Le noyau de cette action est $T(\mathfrak{o}_{\bar{F}})^{I_{F}}/T(\mathfrak{o}_{\bar{F}})^{I_{F}}_{0}$.

Soit $M$ un $F^{nr}$-Levi de $G$ contenant $T$ et $s$ un sommet de $Imm_{F^{nr}}(M_{AD})$ appartenant à l'appartement $App^M_{F^{nr}}(T^{nr})$ de cet immeuble associé  au tore $T^{nr}$. Le tore $T^{nr}$ détermine un sous-tore maximal $T_{s}$ de $M_{s}$. Fixons un sous-groupe de Borel  $B_{s}$ de ce groupe contenant $T_{s}$.  Le groupe $Norm_{G(F^{nr})}(M)$ agit naturellement sur $Imm_{F^{nr}}(M_{AD})$. On note $Norm_{G(F^{nr})}(M,s)$ le sous-groupe des éléments qui conservent $s$. Tout élément $g$ de  ce groupe détermine un automorphisme $Ad(g)_{M_{s}}$ de $M_{s}$.  

Posons $Norm_{G(F^{nr})}(T^{nr},M,s)=Norm_{G(F^{nr})}(T^{nr})\cap Norm_{G(F^{nr})}(M,s)$. Pour $n\in Norm_{G(F^{nr})}(T^{nr},M,s)$, $Ad(n)_{ M_{s}}$ conserve $T_{s}$. Notons  $Norm_{G(F^{nr})}(T^{nr},M,s,B_{s})$ le sous-groupe des $n\in  Norm_{G(F^{nr})}(T^{nr},M,s)$ tels que $Ad(n)_{ M_{s}}$ conserve $B_{s}$.  Ce groupe contient $T(\mathfrak{o}_{\bar{F}})^{I_{F}}$.  Notons $W_{aff}(M,s)$ et $W_{aff}(M,s,B_{s})$   les images de $Norm_{G(F^{nr})}(T^{nr},M,s)$ et $Norm_{G(F^{nr})}(T^{nr},M,s,B_{s})$ dans $W_{aff}$.

Soit $x\in  P(s)$. Le tore $T^{nr}$ détermine un sous-tore maximal $T_{x}$ de $G_{x}$. Le groupe $M_{s}$ est un Levi de $G_{x}$ et on a $T_{s}=T_{x}$. On pose $Norm_{G_{x}}(T_{x},M_{s})=Norm_{G_{x}}(T_{x})\cap Norm_{G_{x}}(M_{s})$. Il est  clair que l'image naturelle de $Norm_{G(F^{nr})}(T^{nr},M,s)\cap K_{x}^{nr,0}$ dans $G_{x}$ appartient à $Norm_{G_{x}}(T_{x},M_{s})$.

\begin{lem} {(i) Le groupe $K_{s}^{M,nr,0}$ est distingué dans $Norm_{G(F^{nr})}(M,s)$. On a les égalités 
$$Norm_{G(F^{nr})}(M,s)=K_{s}^{M,nr,0}Norm_{G(F^{nr})}(T^{nr},M,s,B_{s})$$ et $$K_{s}^{M,nr,0}\cap Norm_{G(F^{nr})}(T^{nr},M,s,B_{s})=T(\mathfrak{o}_{\bar{F}})^{I_{F}}_{0}.$$

(ii) L'action sur $Imm_{F^{nr}}(G_{AD})$ de $Norm_{G(F^{nr})}(M,s)$ conserve $P(s)$. Si $G$ est simplement connexe, $W_{aff}(M,s,B_{s})$  agit fidèlement sur $P(s)$. 

(iii) Soit $x\in P(s)$. Le groupe $Norm_{G(F^{nr})}(T^{nr},M,s)\cap K_{x}^{nr,0}$ s'envoie surjectivement sur $Norm_{G_{x}}(T_{x},M_{s})$. 

 }\end{lem}

Preuve. La première assertion de (i) est évidente. Soit $g\in Norm_{G(F^{nr})}(M,s)$. Alors $gT^{nr}g^{-1}$ est un tore déployé sur $F^{nr}$ de $M$ et $s$ appartient à l'appartement associé à ce tore dans $Imm_{F^{nr}}(M_{AD})$ Un tel tore est conjugué à $T^{nr}$ par un élément de $K_{s}^{M,nr,0}$. Quiite à multiplier $g$ à gauche par un tel élément, on peut supposer que $Ad(g)$ conserve $T^{nr}$, autrement dit que $g\in 
Norm_{G(F^{nr})}(T^{nr})$.  Plus précisément, on a $g \in Norm_{G(F^{nr})}(T^{nr},M,s)$.  Il existe un élément $n_{s}\in Norm_{M_{s}}(T_{s})$ tel que $Ad(n_{s})Ad(g)_{M_{s}}$  conserve $B_{s}$. Puisque $n_{s}$ se relève en un élément de $K_{s}^{M,nr,0}$ normalisant $T^{nr}$, on peut aussi bien multiplier $g$ à gauche par un tel élément et supposer que $Ad(g)_{M_{s}}$ conserve $B_{s}$. Alors $g$ appartient à $Norm_{G(F^{nr})}(T^{nr},M,s,B_{s})$ et cela démontre la première égalité du (i) de l'énoncé. On sait bien que l'intersection $K_{s}^{M,nr, 0}\cap Norm_{G(F^{nr})}(T^{nr})$ s'inscrit dans une suite exacte
  $$1\to T(\mathfrak{o}_{\bar{F}})^{I_{F}}_{0}\to K_{s}^{M,nr, 0}\cap Norm_{G(F^{nr})}(T^{nr})\to W^{M_{s}}\to 1$$
  où $W^{M_{s}}$ est le groupe de Weyl de $M_{s}$ relatif au tore $T_{s}$. Un élément du groupe  $K_{s}^{M,0}\cap N_{G(F^{nr})}(T^{nr}, M,s,B_{s})$ s'envoie sur $1\in W^{M_{s}}$ puisque l'action d'un tel élément conserve $B_{s}$. La deuxième égalité de (i) s'en déduit.
  
    Pour un élément $g\in Norm_{G(F^{nr})}(M)$, on a une égalité $g_{M_{AD}}\circ p_{M}=p_{M}\circ g$ sur $Imm_{F^{nr}}^G(M_{ad})$, où on note $g_{M_{AD}}$ l'action sur $Imm_{F^{nr}}(M_{AD})$ déduite de l'action de $g$. On en déduit que, si $g\in Norm_{G(F^{nr})}(M,s)$ l'action de $g$ sur $Imm_{F^{nr}}(G_{AD})$ conserve $p_{M}^{-1}(s)=P(s)$. D'où  la première assertion de (ii). Supposons $G$ simplement connexe. Soit $g\in Norm_{G(F^{nr})}(T^{nr},M,s,B_{s})$. Supposons que son action sur $P(s)$ soit triviale. Puisque $P(s)$ est un espace affine sous l'espace vectoriel ${\cal A}_{M}^{nr}$,  l'action linéaire de $g$ sur ${\cal A}_{M}^{nr}$ est triviale. Autrement dit, l'action de $Ad(g)$ sur $A_{M}^{nr}$ est triviale. Un élément de $Norm_{G(F^{nr})}(M)$ dont l'action sur ce tore est triviale appartient à $M(F^{nr})$. Puisque $g$ fixe $s$, on a $g\in K_{s}^{M,nr,\dag}$.  Comme on l'a rappelé ci-dessus, de  l'homomorphisme de Kottwitz $w_{M}$ se  déduit un homomorphisme injectif
   $$ w_{M}:K_{s}^{M,nr,0}\backslash K_{s}^{M,nr,\dag}\to X^*(Z(\hat{M})^{I_{F}}).$$
   En utilisant (1) et (i), de $w_{M}$ se déduit encore un homomorphisme injectif
   $$w_{M}:  T(\mathfrak{o}_{\bar{F}})^{I_{F}}\backslash (K_{s}^{M,nr,\dag}\cap Norm_{G(F^{nr})}(T^{nr},M,s,B_{s}))\to X^*(Z(\hat{M})^{I_{F}}).$$
   Parce que $G$ est simplement connexe, le groupe $Z(\hat{M})^{I_{F}}$ est un tore de m\^eme dimension que $A_{M}^{nr}$.  Donc $X^*(Z(\hat{M})^{I_{F}})$ est un ${\mathbb Z}$-module libre de cette dimension. Il en résulte que $ T(\mathfrak{o}_{\bar{F}})^{I_{F}}\backslash (K_{s}^{M,nr,\dag}\cap Norm_{G(F^{nr})}(T^{nr},M,s,B_{s}))$ est sans torsion. 
   Le groupe $A_{M}^{nr}$ est inclus dans $K_{s}^{M,nr,\dag}\cap Norm_{G(F^{nr})}(T^{nr},M,s,B_{s})$.  L'image par   $w_{M}$ de  $(T(\mathfrak{o}_{\bar{F}})^{I_{F}}\cap A_{M}^{nr})\backslash A_{M}^{nr}$ est encore un ${\mathbb Z}$-module libre de m\^eme  dimension que $X^*(Z(\hat{M})^{I_{F}})$ . Il existe donc un entier $N$ tel que $N X^*(Z(\hat{M})^{I_{F}})$ soit inclus dans cette image. En conséquence $g^N$ appartient à $ T(\mathfrak{o}_{\bar{F}})^{I_{F}}A_{M}^{nr}$.  C'est-à-dire que $g^N=ta$, avec $t\in T(\mathfrak{o}_{\bar{F}})^{I_{F}}$ et $a\in A_{M}^{nr}$. Mais $g$ agit trivialement sur $P(s)$ donc $a$ aussi. L'action de $a\in A_{M}^{nr}$ sur $App_{F^{nr}}(T^{nr})$ est une translation par l'image naturelle de $a$ dans $X_{*}(A_{M}^{nr})$. Cette translation étant triviale, $a$ appartient au sous-groupe   $A_{M}^{nr}(\mathfrak{o}_{F^{nr}})$, lequel est inclus dans $T(\mathfrak{o}_{\bar{F}})^{I_{F}}$. On conclut que $g^N$ appartient aussi à ce groupe, donc $g$ aussi puisque $ T(\mathfrak{o}_{\bar{F}})^{I_{F}}\backslash (K_{s}^{M,nr,\dag}\cap Norm_{G(F^{nr})}(T^{nr},M,s,B_{s}))$ est sans torsion. 
   Donc l'image de $g$ dans $W_{aff}$ est triviale. Cela prouve la deuxième assertion de (ii).  
    
   Prouvons (iii). Le groupe $T(\mathfrak{o}_{\bar{F}})^{I_{F}}_{0}$ est contenu dans $Norm_{G(F^{nr})}(T^{nr},M,s)\cap K_{x}^{nr,0}$  et s'envoie surjectivement sur $T_{x}$. Il suffit donc de prouver que l'homomorphisme composé
   $$(2) \qquad Norm_{G(F^{nr})}(T^{nr},M,s)\cap K_{x}^{nr,0} \to Norm_{G_{x}}(T_{x},M_{s})/T_{x}$$
   est surjectif. L'espace d'arrivée est égal à $Norm_{W^{G_{x}}}(M_{s})$. Ce groupe ne change pas quand on remplace $G$ par $G_{SC}$. Cela nous ramène à démontrer l'assertion quand $G$ est simplement connexe. Le procédé de \ref{reduction} nous ramène au cas où $G$ est  absolument quasi-simple. En appliquant la proposition \ref{couples}, on peut supposer que $(M,s)=(M_{\Lambda},s_{\Lambda})$ pour un sous-ensemble propre $\Lambda$ de $\Delta_{a}^{nr}$ et que $x\in \bar{C}^{nr}\cap P(s)$. On note $B_{s}$ le sous-groupe de Borel de $M_{s}$ déterminé par l'alc\^ove $C^{nr}$. Soit $w\in Norm_{W^{G_{x}}}(M_{s})$.  Il existe un unique $w'\in W^{M_{s}}$ tel que $w'w$ conserve $B_{s}$. L'élément $w'$ se relevant en un élément de $K_{s}^{M,nr,0}\cap Norm_{M(F^{nr})}(T^{nr})\subset Norm_{G(F^{nr})}(T^{nr},M,s)\cap K_{x}^{nr,0}$, il suffit de prouver que $w'w$ appartient à l'image de l'homomorphisme (2). Cela nous ramène au cas où $w$ conserve $B_{s}$.  Relevons $w$ en un élément 
$n\in Norm_{G(F^{nr})}(T^{nr})\cap K_{x}^{nr,0}$. Parce que $w$ conserve $B_{s}$ qui est associé au sous-ensemble $\Lambda$ de l'ensemble de racines simples $\Lambda(x)$ de $G_{x}$, l'action linéaire de $n$ sur ${\cal A}^{nr}$ conserve $\Lambda$ et aussi  ${\mathbb Q} [\Lambda]$. Alors elle conserve l'ensemble de racines ${\mathbb Q}[\Lambda]\cap \Sigma$ de $M$. Donc $Ad(n)$ conserve $M$, autrement dit $n\in Norm_{G(F^{nr})}(T^{nr})\cap Norm_{G(F^{nr})}(M)$. Puisque $n$ conserve $x$ et $M$, $n$ conserve aussi $p_{M}(x)=s$. Donc $n\in  Norm_{G(F^{nr})}(T^{nr},M,s)\cap K_{x}^{nr,0}$, ce qui prouve (iii). $\square$

   Il résulte du (i) du lemme une suite exacte
    $$(3) \qquad 1 \to K_{s}^{M,nr, 0}\to Norm_{G(F^{nr})}(M,s)\to W_{aff}( M,s,B_{s})\to 1. $$

Considérons le cas où $G$ est absolument quasi-simple et où $(M,s)=(M_{\Lambda},s_{\Lambda})$ pour un sous-ensemble propre $\Lambda$ de $\Delta_{a}^{nr}$. On suppose alors que $B_{s_{\Lambda}}$ est le sous-groupe de $M_{\Lambda,s_{\Lambda}}$ déterminé par l'alc\^ove $C^{nr}$.  
   On note $W_{aff}(\Lambda^{aff})$ le sous-ensemble des éléments de $W_{aff}$ qui conservent l'ensemble de racines affines $\Lambda^{aff} $.   On note $Norm_{G(F^{nr})}(T^{nr},\Lambda^{aff})$ l'image réciproque de $W_{aff}(\Lambda^{aff})$ dans $N_{G(F^{nr})}(T^{nr})$. Montrons que

    (4) on a l'égalité $W_{aff}(\Lambda^{aff})=W_{aff}(M_{\Lambda},s_{\Lambda},B_{s_{\lambda}})$. 
    
    Soit $g\in N_{G(F^{nr})}(T^{nr},\Lambda^{aff})$. Puisque l'action de $g$ sur $App_{F^{nr}}(T^{nr})$ conserve $\Lambda^{aff}$, l'action linéaire de $g$ sur $ {\cal A}^{nr}$ conserve $\Lambda$ donc aussi l'ensemble $\Sigma^{nr}\cap {\mathbb Q}[\Lambda]$ des racines de $M_{\Lambda}$. Donc $Ad(g)$ normalise $M_{\Lambda}$.  Puisque l'action de $g$ sur $App_{F^{nr}}(T^{nr})$ conserve  $P(\Lambda^{aff})=P(s_{\Lambda})$, l'action de $g$ sur $Imm_{F^{nr}}(M_{AD})$  conserve  aussi $p_{M}(P( s_{\Lambda}))$, c'est-à-dire $s_{\Lambda}$. Donc $g\in  Norm_{G(F^{nr})}(M_{\Lambda},s_{\Lambda})$. Puisque $\Lambda$ est un ensemble de racines simples associé au sous-groupe de Borel $B_{s_{\Lambda}}$ (cf. \ref{constructiondeLevi}(2)) et que   l'action linéaire de $g$ sur $ {\cal A}^{nr}$ conserve $\Lambda$, l'action de $g$ sur $M_{s_{\Lambda}}$ conserve $B_{s_{\Lambda}}$. Inversement, soit $g\in N_{G(F^{nr})}(T^{nr},M_{\Lambda},s_{\Lambda},B_{s_{\Lambda}})$. Pour la m\^eme raison que ci-dessus,  l'action linéaire de $g$ sur $ {\cal A}^{nr}$ conserve $\Lambda$. Soit $\beta\in \Lambda$, notons $\beta'$ son image par cette action. Puisque  $Ad(g)$ conserve $s_{\Lambda}$, l'action de $g$ sur l'ensemble des racines affines envoie $\beta^{aff}$ sur une racine affine  de partie vectorielle $\beta'$ et qui s'annule en $s_{\Lambda}$. L'unique telle racine affine est $\beta^{'aff}$. Donc l'action de $g$ conserve $\Lambda^{aff}$, c'est-à-dire $g\in N_{G(F^{nr})}(T^{nr},\Lambda^{aff})$. Cela prouve (4).

\subsubsection{Excellents couples $(M,s)$}\label{excellentscouples}

Soient $M$ un $F^{nr}$-Levi et $s$ un sommet de $Imm_{F^{nr}}(M_{AD})$.   Supposons que $M\not=G$, soit ${\cal F}$ une facette de $P(s)$ de dimension $dim(A_M^{nr})-dim(A_G^{nr})-1$. Alors le groupe $M_{s}$ est un Levi maximal de $G_{{\cal F}}$. Le quotient $Norm_{G_{{\cal F}}}(M_{s})/M_{s}$ a un ou deux éléments. 

{\bf Définition.} On dit que $(M,s)$ est excellent si l'une des conditions suivantes est vérifiée:

(a) $M=G$;

(b) $M\not=G$ et, pour toute facette ${\cal F}$ de $P(s)$ de dimension $dim(A_M^{nr})-dim(A_G^{nr})-1$, $Norm_{G_{{\cal F}}}(M_{s})/M_{s}$ a deux éléments. 

\bigskip
Cette définition dépend du groupe ambiant. Le cas échéant, on dira plut\^ot que $(M,s)$ est $G$-excellent.
 Utilisons la décomposition de $G_{SC}$ introduite en \ref{reduction} et les objets $G'_{i},M'_{i},s_{i}$ afférents. 
 
\begin{lem}{(i) Le couple $(M,s)$ est $G$-excellent si et seulement si $(M_{sc},s)$ est $G_{SC}$-excellent. 

(ii) Le couple $(M,s)$ est $G$-excellent si et seulement si $(M'_{i},s_{i})$ est $G'_{i}$-excellent pour tout $i=1,...,n$.

(iii) Soit $L$ un $F^{nr}$-Levi de $G$ tel que $M\subset L$. Si  le couple $(M,s)$ est $G$-excellent, alors il est $L$-excellent. }\end{lem}

Preuve. L'assertion (i) est évidente si $M=G$. Supposons $M\not=G$. Passer de $G$ à $G_{SC}$ ne change pas $P(s)$. Soit ${\cal F}$ une facette de $P(s)$ de dimension $dim(A_M^{nr})-dim(A_G^{nr})-1$. Les groupes dérivés de $G_{{\cal F}}$ et $G_{SC,{\cal F}}$ ont  m\^eme rev\^etement simplement connexe. Les images réciproques de $M_{s}$ et $M_{sc,s}$  dans ce groupe $G_{{\cal F},SC}$ sont égales, notons-les $M_{s,sc}$. Alors $Norm_{G_{{\cal F}}}(M_{s})/M_{s}\simeq Norm_{G_{{\cal F},SC}}(M_{s}/M_{s,sc})\simeq Norm_{G_{SC,{\cal F}}}(M_{sc,s})/M_{sc,s}$. L'assertion (i) s'ensuit.

 En vertu de (i), pour démontrer (ii), on peut supposer que $G$ est simplement connexe. L'assertion (ii) est encore évidente si $M=G$.  Supposons $M\not=G$. On peut supposer qu'il existe $m\in \{1,...,n\}$ tel que $M'_{i}=G'_{i}$ si et seulement si $i> m$. On a $P(s)=\prod_{i=1,...,n}P(s_{i}) $ et $P(s_{i})$  est réduit à un point pour $i>m$. Une facette ${\cal F}$ de $P(s)$   est un produit $\prod_{i=1,...,n}{\cal F}_{i}$ où, pour tout $i$, ${\cal F}_{i}$ est une facette de $P(s_{i})$. On a ${\cal F}_{i}=P(s_{i})$ pour $i> m$. La facette ${\cal F}$ est
de dimension $dim(A_M^{nr})-dim(A_G^{nr})-1$ si et seulement s'il existe $j\in \{1,...,m\}$ de sorte que $dim({\cal F}_{j})=dim(A_{M'_{j}}^{nr})-dim(A_{G'_{j}}^{nr})-1$ tandis que $dim({\cal F}_{i})=dim(A_{M'_{i}}^{nr})-dim(A_{G'_{i}}^{nr})$ pour tout $i\not=j$ (les "$nr$"  se réfèrent ici au fait qu'un groupe $G'_{j}$ est défini sur    $F_{j}^{nr}$ où $F_{j}$ est une extension finie modérément ramifiée de $F$). On a alors $G_{{\cal F}}=\prod_{i=1,...,n}G'_{{\cal F}_{i}}$;  $M'_{j,s_{j}}$ est un $E_{j}$-Levi maximal de $G'_{j,{\cal F}_{j}}$; pour $i\not=j$, $G'_{i,{\cal F}_{i}}=M'_{i,s_{i}}$. Donc $Norm_{G_{{\cal F}}}(M_{s})/M_{s}$ s'identifie à $Norm_{G'_{j,{\cal F}_{j}}}(M'_{j,s_{j}})/M'_{j,s_{j}}$. Si $(M'_{j},s_{j})$ est $G'_{j}$-excellent, ce groupe a deux éléments et cela démontre l'assertion "si" du (ii) de l'énoncé. La réciproque résulte du m\^eme calcul car, pour $j\in \{1,...,m\}$ et pour une facette ${\cal F}_{j}$ de $P(s_{j})$ telle que $dim({\cal F}_{j})=dim(A_{M'_{j}}^{nr})-dim(A_{G'_{j}}^{nr})-1$, on peut trouver une facette ${\cal F}$ comme ci-dessus dont la composante  en $j$ soit ${\cal F}_{j}$. Cela prouve le (ii) de l'énoncé. 

Prouvons (iii). On peut supposer $L\not=M$.  On a deux ensembles $P(s)\subset Imm_{F^{nr}}(G_{AD})$ et $P^L(s)\subset Imm_{F^{nr}}(L_{AD})$. L'ensemble $P(s) $ est contenu dans $Imm^G_{F^{nr}}(L_{ad})$ et on a l'égalité $P^L(s)=p_{L}(P(s))$. Soit ${\cal F}^L$ une facette de $P^L(s)$ de dimension $dim(A_M^{nr})-dim(A_L^{nr})-1$. On peut fixer une facette ${\cal F}$ de $P(s)$, de dimension $dim(A_M^{nr})-dim(A_G^{nr})-1$, telle que $p_{L}({\cal F})\subset {\cal F}^L$. On a alors $G_{{\cal F}}=L_{{\cal F}^L}$ et $Norm_{G_{{\cal F}}}(M_{s})/M_{s}=Norm_{L_{{\cal F}^L}}(M_{s})/M_{s}$. Si $(M,s)$ est $G$-excellent, ce groupe a deux éléments et cela prouve que $(M,s)$ est $L$-excellent. $\square$

Montrons que  

(1) supposons que l'ensemble ${\bf FC}^{M_{s}}(\mathfrak{m}_{SC,s})$ soit non vide; alors $(M,s)$ est excellent.

C'est clair si $M=G$. Sinon, la condition (b) ci-dessus est vérifiée d'après \ref{faisceauxcaracteres}(1).

    \subsubsection{Excellents ensembles $\Lambda$}\label{excellentsensemblesLambda}
    
     {\bf On suppose   que} $G$ {\bf est absolument quasi-simple et simplement connexe}.   
     
   Soit $\Lambda$ un sous-ensemble de $\Delta_{a}^{nr}$.  Supposons que $\Delta_{a}^{nr}-\Lambda$ ait au moins deux éléments. Soit $\gamma\in \Delta_{a}^{nr}-\Lambda$. Notons ${\cal F}_{\gamma}$ la facette de $\bar{C}^{nr}$ formée des $x\in \bar{C}^{nr}$ vérifiant les égalités $\beta^{aff}(x)=0$ pour $\beta\in \{\gamma\}\cup \Lambda$. D'après \ref{constructiondeLevi}(2), $G_{{\cal F}_{\gamma}}$ a pour ensemble de racines simples $ \{\gamma\}\cup \Lambda$ et $M_{\Lambda,s_{\Lambda}}$ en est le Levi maximal standard associé au sous-ensemble $\Lambda$.  Le groupe $Norm_{G_{{\cal F}_{\gamma}}}(M_{\Lambda,s_{\Lambda}})/M_{\Lambda,s_{\Lambda}}$ est alors d'ordre $1$ ou $2$. Les arguments de la preuve du lemme \ref{couples}  montrent que ce groupe  est isomorphe à  l'image de $K_{{\cal F}_{\gamma}}^{0,nr}\cap  N_{G(F^{nr})}(T^{nr},\Lambda^{aff})$ dans $W_{aff}(\Lambda^{aff})$. S'il est d'ordre $2$, on note $w_{\gamma}$ l'élément non trivial de cette image. 
   
   Supposons seulement que $\Lambda$ est un sous-ensemble propre de $\Delta_{a}^{nr}$. La définition suivante imite Lusztig, cf. \cite{L?} 2.4.
   
   {\bf Définition.} on dit que $\Lambda$ est excellent si l'une des conditions suivantes est vérifiée:
   
 (a)  $\Delta_{a}^{nr}-\Lambda$ a un élément;
   
  (b)  $\Delta_{a}^{nr}-\Lambda$ a au moins deux éléments et le groupe $Norm_{G_{{\cal F}_{\gamma}}}(M_{\Lambda,s_{\Lambda}})/M_{\Lambda,s_{\Lambda}}$  est d'ordre $2$ pour tout $\gamma\in \Delta_{a}^{nr}-\Lambda$.   
  
  \bigskip
   
   Supposons (b) vérifiée. On pose $S=\{w_{\gamma}; \gamma\in \Delta_{a}^{nr}-\Lambda\}$. Notons $Aut(P(\Lambda^{aff}))$ le groupe des tranformations affines de $P(\Lambda^{aff})$. D'après \ref{normalisateur}(4) et le (ii) du lemme \ref{normalisateur}, l'action de $W_{aff}(\Lambda^{aff})$ sur $App_{F^{nr}}(T^{nr})$ conserve $P(\Lambda^{aff})$ et l'application $res:W_{aff}(\Lambda^{aff})\to Aut(P(\Lambda^{aff}))$ qui s'en déduit est injective. Notons ${\bf W}_{aff}(\Lambda^{aff})$ son image.  L'espace affine $App_{F^{nr}}(T^{nr})$ est un espace euclidien et $P(\Lambda^{aff})$ en  est un sous-espace   euclidien. On munit $P(\Lambda^{aff})$ de la décomposition en facettes induite par celle de $App_{F^{nr}}(T^{nr})$: une facette de $P(\Lambda^{aff})$ est une facette de cet appartement qui est incluse dans $P(\Lambda^{aff})$. Appelons $P(\Lambda^{aff})$-facette une telle facette.  On appelle $P(\Lambda^{aff})$-chambre une $P(\Lambda^{aff})$-facette  de dimension maximale et $P(\Lambda^{aff})$-mur un hyperplan affine de $P(\Lambda^{aff})$ qui est un mur d'une $P(\Lambda^{aff})$-chambre. Notons ${\cal W}_{aff}(\Lambda^{aff})$ le sous-groupe de  $Aut(P(\Lambda^{aff}))$ engendré par les symétries orthogonales relatives aux $P(\Lambda^{aff})$-murs. 
  
  \begin{prop}{On suppose que $\Delta_{a}^{nr}-\Lambda$ a au moins deux éléments et que $\Lambda$ est excellent. Alors ${\bf W}_{aff}(\Lambda^{aff})={\cal W}_{aff}(\Lambda^{aff})$. Le couple $(W_{aff}(\Lambda^{aff}),S)$ est un système  de Coxeter, cf. \cite{B} IV.1.3. Pour tout $x\in P(\Lambda^{aff})$, il existe $w\in W_{aff}(\Lambda^{aff})$ tel que $w(x)\in P(\Lambda^{aff})\cap \bar{C}^{nr}$. }\end{prop}
  
  Preuve  Pour $\gamma\in \Delta_{a}^{nr}-\Lambda$, posons ${\bf w}_{\gamma}=res(w_{\gamma})$. L'ensemble  $P(\Lambda^{aff})\cap \bar{C}^{nr}$ est une $P(\Lambda^{aff})$-chambre et l'ensemble  ${\bf S}=\{{\bf w}_{\gamma}; \gamma\in \Delta_{a}^{nr}-\Lambda\}$ est l'ensemble des symétries orthogonales relatives aux $P(\Lambda^{aff})$-murs de cette $P(\Lambda^{aff})$-chambre. Notons ${\cal W}_{{\bf S}}$ le sous-groupe de ${\cal W}_{aff}(\Lambda^{aff})$ engendré par ${\bf S}$. On applique les résultats de \cite{B} paragraphe V.3.1. Les hypothèses de Boubaki sont que ${\cal W}_{aff}(\Lambda^{aff})$ conserve la décomposition en facettes de $P(\Lambda^{aff})$ et que ${\cal W}_{aff}(\Lambda^{aff})$ agit proprement sur $P(\Lambda^{aff})$. Mais en examinant la preuve de  lemme 2 de loc.cit., on s'aperçoit qu'il suffit que  ${\cal W}_{{\bf S}}$ conserve la décomposition en facettes de $P(\Lambda^{aff})$ et que ${\cal W}_{{\bf S}}$ agisse proprement sur $P(\Lambda^{aff})$. Or ces hypothèses sont vérifiées. En effet, notons $W_{S}$ le sous-groupe de $W_{aff}(\Lambda^{aff})$ engendré par  l'ensemble $S $.   On a l'égalité ${\cal W}_{{\bf S}}=res(W_{S})$. Les propriétés requises de l'action de ${\cal W}_{{\bf S}}$ se déduisent alors des m\^emes propriétés de l'action de $W_{aff}$ sur $App_{F^{nr}}(T^{nr})$. Appliquons donc le lemme 2 de loc. cit. Il affirme d'abord que, pour tout $x\in P(\Lambda^{aff})$, il existe ${\bf w}\in  {\cal W}_{{\bf S}}$ tel que ${\bf w}(x)\in P(\Lambda^{aff})\cap \bar{C}^{nr}$. La dernière assertion de la proposition en résulte. Le lemme affirme aussi que ${\cal W}_{aff}(\Lambda^{aff})={\cal W}_{{\bf S}}$.
   Il en résulte que les hypothèses de Bourbaki évoquées ci-dessus sont vérifiées. On peut donc appliquer le théorème 1 de \cite{B} V.3.2. Il affirme que ${\cal W}_{{\bf S}}={\cal W}_{aff}(\Lambda^{aff})$ est un groupe de Coxeter, plus exactement que $({\cal W}_{{\bf S}},{\bf S})$ est un système de Coxeter. 
   Il en résulte que $(W_{S},S)$ est un système de Coxeter. Pour achever la preuve de la proposition, il reste à prouver que $W_{aff}(\Lambda^{aff})=W_{S}$. Soit $w\in W_{aff}(\Lambda^{aff})$. Comme on l'a dit, son action sur $P(\Lambda^{aff})$ conserve la décomposition en facette. Notons ${\bf C}=w(C^{nr}\cap P(\Lambda^{aff}))$. Le théorème de Bourbaki affirme aussi qu'il existe  ${\bf w}'\in {\cal W}_{{\bf S}}$ tel que ${\bf w}'({\bf C})=C^{nr}\cap P(\Lambda^{aff})$. Autrement dit, il existe $w'\in W_{S}$ tel que $w'({\bf C})=C^{nr}\cap P(\Lambda^{aff})$. Alors $w'w$ conserve $C^{nr}\cap P(\Lambda^{aff})$.  On sait que deux éléments de $C^{nr}$ ne sont conjugués par un élément de $W_{aff}$ que s'ils sont égaux (parce que $G$ est simplement connexe). Donc $w'w$ fixe tout point de $C^{nr}\cap P(\Lambda^{aff})$. Mais cet ensemble engendre l'espace affine $P(\Lambda^{aff})$. Donc $w'w$ agit trivialement sur $P(\Lambda^{aff})$, c'est-à-dire $res(w'w)=1$, d'où $w'w=1$.  Donc $w$ appartient à $W_{S}$. Cela achève la preuve. $\square$
   
   \subsubsection{Un corollaire}\label{uncorollaire2}
   On. conserve les m\^emes hypothèses. 
   
   \begin{cor}{Soit $\Lambda$ un sous-ensemble propre de $\Delta_{a}^{nr}$. Alors le couple $(M_{\Lambda},s_{\Lambda})$ est excellent si et seulement si $\Lambda$ est excellent.}\end{cor}
   
   Preuve. Si $\Delta_{a}^{nr}-\Lambda$ n'a qu'un élément, $\Lambda$ est excellent et $M_{\Lambda}=G$, donc $(M_{\Lambda},s_{\Lambda})$ est excellent. Supposons que $\Delta_{a}^{nr}-\Lambda$  a au moins deux éléments.   La condition (b) de \ref{excellentsensemblesLambda} est  un cas particulier de la condition (b) de \ref{excellentscouples}. Cela entraîne l'implication "seulement si" de l'énoncé. Inversement, supposons que  $\Lambda$ est excellent. Soit ${\cal F}$ une facette de $P(\Lambda^{aff})=P(s_{\Lambda})$ de dimension $dim(A_{M_{\Lambda}}^{nr})-dim(A_{G}^{nr})-1=dim(P(\Lambda^{aff}))-1$. D'après la proposition \ref{excellentsensemblesLambda}, il existe $w\in W_{aff}(\Lambda^{aff})$ tel que $w({\cal F})$ soit incluse dans $P(\Lambda^{aff})\cap \bar{C}^{nr}$. Une facette de cet ensemble  qui est de dimension  $dim(P(\Lambda^{aff}))-1$ est de la forme ${\cal F}_{\gamma}$ pour un $\gamma\in \Delta_{a}^{nr}-\Lambda$. Puisque $\Lambda$ est excellent, le groupe $Norm_{G_{{\cal F}_{\gamma}}}(M_{\Lambda,s_{\Lambda}})/ M_{\Lambda,s_{\Lambda}}$ a deux éléments. L'action de $w^{-1}$ se quotiente en une bijection de ce groupe sur $Norm_{G_{{\cal F}}(}M_{\Lambda,s_{\Lambda}})/M_{\Lambda,s_{\Lambda}}$. Ce dernier groupe a donc deux éléments, ce qui prouve que  le couple $(M_{\Lambda},s_{\Lambda})$ est excellent. $\square$

 \subsubsection{Cas  d'un excellent couple $(M,s)$}\label{casdunexcellentcouple}
 
  {\bf On suppose  que }$G$ {\bf est défini et quasi-déployé sur} $F$ {\bf et qu'il est absolument quasi-simple.}  
 
 \begin{prop}{Soient $M$ un $F$-Levi de $G$ et $s$ un sommet de $Imm(M_{AD})$. Supposons que $s$ reste un sommet dans $Imm_{F^{nr}}(M_{AD})$ et que le couple $(M,s)$ soit excellent.

 (i)  Il existe un unique sous-ensemble $\Lambda\subsetneq\Delta_{a}^{nr}$ de sorte que les couples $(M,s_{M})$ et $(M_{\Lambda},s_{\Lambda})$ soient conjugués par un élément de $G_{SC}(F^{nr})$. Ces couples sont en fait conjugués par un élément de $G_{SC}(F)$. 
 
 (ii) Soit $\Lambda'\subsetneq \Delta_{a}^{nr}$. Les couples $(M,s)$ et $(M_{\Lambda'},s_{\Lambda'})$ sont conjugués par un élément de $G_{AD}(F^{nr})$ si et seulement si $\Lambda'$ est conjugué à $\Lambda$ par un élément de $\Omega^{nr}$.}\end{prop}
 
 Preuve.   L'existence de l'ensemble $\Lambda$ vérifiant (i) résulte donc de la proposition \ref{couples}. Fixons un tel $\Lambda$.  Supposons d'abord que $\Delta_{a}^{nr}-\Lambda$ n'a qu'un élément. Alors $A_{\Lambda}^{nr}=\{1\}$, $M_{\Lambda}=G$ et $s_{\Lambda}$ est le sommet de $\bar{C}^{nr}$ associé à l'unique élément de $\Delta_{a}^{nr}-\Lambda$. L'assertion  d'unicité de (i) et l'assertion (ii) résultent des faits que deux éléments de $\bar{C}^{nr}$  sont conjugués par un élément de $G_{SC}(F^{nr})$, resp. $G_{AD}(F^{nr})$, si et seulement s'ils sont égaux, resp. conjugués par un élément de ${\bf \Omega}^{nr}$, cf. \ref{racinesaffines}(6).  On suppose maintenant que $\Delta_{a}^{nr}-\Lambda$  a au moins deux éléments. L'ensemble $\Lambda$ est excellent d'après le corollaire \ref{uncorollaire2}. 
  Soit $\Lambda'\subsetneq \Delta_{a}^{nr}$, supposons qu'il existe $g\in G_{SC}(F^{nr})$ qui conjugue $(M_{\Lambda'},s_{\Lambda'})$ en $(M_{\Lambda},s_{\Lambda})$. L'action de $g$ sur $Imm_{F^{nr}}(G)$ envoie $ P(s_{\Lambda'})$ sur $P(s_{\Lambda})$, c'est-à-dire $P(\Lambda^{'aff})$ sur $P(\Lambda^{aff})$. Soit $\gamma\in \Delta_{a}^{nr}-\Lambda'$. Le sommet $s_{\gamma}$ de $\bar{C}^{nr}$ associé à $\gamma$ appartient à $P(\Lambda^{'aff})$. Donc $g(s_{\gamma})\in P(\Lambda^{aff})$. D'après   la proposition \ref{excellentsensemblesLambda}, on peut trouver $n\in Norm_{G_{SC}(F^{nr})}(T^{nr},\Lambda^{aff})$ de sorte que $ng(s_{\gamma})$ appartienne à $\bar{C}^{nr}\cap P(\Lambda^{aff})$. Mais deux éléments de $\bar{C}^{nr}$ ne sont conjugués par l'action d'un élément de $G_{SC}(F^{nr})$ que s'ils sont égaux. Donc $s_{\gamma}=ng(s_{\gamma})\in 
 P(\Lambda^{aff})$. Les sommets de $\bar{C}^{nr}$ qui appartiennent à $P(\Lambda^{aff})$ sont ceux qui sont associés à des éléments de $\Delta_{a}^{nr}-\Lambda$. Donc $\gamma$ appartient à $\Delta_{a}^{nr}-\Lambda$. Cela prouve l'inclusion $\Delta_{a}^{nr}-\Lambda'\subset \Delta_{a}^{nr}-\Lambda$. Mais $\Lambda$ et $\Lambda'$ ont m\^eme nombre d'éléments pour des raisons de dimension. L'inclusion précédente équivaut alors à $\Lambda=\Lambda'$. Cela prouve le (i) de l'énoncé. 
 
   Soit $\Lambda'\subsetneq \Delta_{a}^{nr}$, supposons qu'il existe $g\in G_{AD}(F^{nr})$ qui conjugue $(M_{\Lambda'},s_{\Lambda'})$ en $(M_{\Lambda},s_{\Lambda})$.   Notons $\pi:G_{SC}\to G_{AD}$ l'homomorphisme naturel. D'après \ref{racinesaffines}(6), on peut trouver $\omega\in \boldsymbol{\Omega}^{nr}$, $t\in T_{ad}(\mathfrak{o}_{\bar{F}})^{I_{F}}$ et  $g'\in G_{SC}(F^{nr})$ de sorte que $g=\omega t\pi(g')$.  Il est clair que  l'image de $(M_{\Lambda},s_{\Lambda})$ par $t^{-1}\omega^{-1}$ est $(M_{\Lambda''},s_{\Lambda''})$, où $\Lambda''=\omega^{-1}(\Lambda)$. Puisque $g$ conjugue $(M_{\Lambda'},s_{\Lambda'})$ en $(M_{\Lambda},s_{\Lambda})$, $g'$ conjugue  $(M_{\Lambda'},s_{\Lambda'})$ en $t^{-1}\omega^{-1}(M_{\Lambda},s_{\Lambda})$, c'est-à-dire en $(M_{\Lambda''},s_{\Lambda''})$. En appliquant (i), cela entraîne $\Lambda'=\Lambda''$, d'où le (ii) de l'énoncé. $\square$

  \subsubsection{Levi standard et excellents couples}\label{Levistandard}
  \begin{prop}{Soit $M$ un $F^{nr}$-Levi standard de $G$.  
 Supposons qu'il existe un sommet $s\in Imm_{F^{nr}}(M_{AD})$ de sorte que $ (M,s)$ soit excellent.  Les propri\'et\'es suivantes sont v\'erifi\'ees:

(i) soit $w\in W^{I_{F}}$, supposons que $w(M)$ soit encore un Levi standard;  alors $w(M)=M$;

(ii) soient $P'$, $P''$ deux sous-groupes paraboliques de $G$ d\'efinis sur $F^{nr}$ de composantes de Levi $M$; alors $P'$ et $P''$ sont conjugu\'es par un \'el\'ement de $W^{I_{F}}$.
 }\end{prop} 

Preuve.  Commen\c{c}ons par d\'eduire l'assertion (ii)   de (i) (la preuve est standard). Le groupe $W^{I_{F}}$ est le groupe de Weyl de $G(F^{nr})$ relativement au tore $T(F^{nr})$. Donc tout sous-groupe parabolique de $G$ défini sur $F^{nr}$ et contenant $T$ est conjugué par un élément de $W^{I_{F}}$ à un unique sous-groupe parabolique standard, lequel est défini sur $F^{nr}$. Notons ${\cal P}^{nr}(M)$ l'ensemble des sous-groupes paraboliques de $G$ d\'efinis sur $F^{nr}$ dont $M$ est une composante de Levi et notons ${\cal M}$ l'ensemble des Levi standard qui sont conjugués à $M$ par un élément de $W^{I_{F}}$. A $P\in {\cal P}^{nr}(M)$, associons le Levi standard de l'unique sous-groupe parabolique standard conjugué à $P$ par un élément de $W^{I_{F}}$. On obtient une application ${\cal P}^{nr}(M)\to {\cal M}$. Elle se quotiente en une bijection    
$${\cal P}^{nr}(M)/W^{I_{F}}\simeq {\cal M}.$$
Chacune des assertions (i) et (ii) est équivalente au fait que l'un de ces ensembles n'a qu'un élément. Elles sont donc équivalentes.

D\'emontrons maintenant (i). Evidemment, cette assertion est triviale si $M=G$. On suppose $M\not=G$. Soit $w\in W^{I_{F}}$ tel que $w(M)$ soit un Levi standard. D'apr\`es \cite{Sha} lemme 2.1.2, on peut trouver des suites $L_{1},...,L_{n}$ et $M_{1},...,M_{n+1}$ de $F^{nr}$-Levi standard de $G$ et, pour tout $j\in \{1,...,n\}$,  un \'el\'ement $w_{j}\in W^{L_{j},I_{F}}$ de sorte que les conditions suivantes soient v\'erifi\'ees:

- pour tout $j\in \{1,...,n\}$, $M_{j}$ et $M_{j+1}$ sont des $F^{nr}$-Levi propres et maximaux de $L_{j}$;

- pour tout $j\in \{1,...,n\}$, $w_{j}(M_{j})=M_{j+1}$;

- $M=M_{1}$ et $w=w_{n}...w_{1}$.

    Les Levi $M_{j}$ \'etant tous conjugu\'es \`a $M$ v\'erifient la m\^eme hypoth\`ese que celui-ci. Pour prouver que $w(M)=M$, il suffit de prouver que $M_{j+1}=M_{j}$ pour tout $j$. Pour cela, on peut remplacer le groupe $G$ par $L_{j}$.  L'hypothèse sur $M_{j}$ se conserve d'après le (iii) du lemme \ref{excellentscouples}. 
Cela nous ram\`ene au cas o\`u $M$ est un $F^{nr}$-Levi standard propre et maximal de $G$. D'après l'équivalence de (i) et (ii), on doit alors prouver que le groupe $W^{nr}(M)=Norm_{G(F^{nr})}(M)/M(F^{nr})$ a deux éléments. 

On utilise les groupes $G'_{i}$, $M'_{i}$ et les sommets $s_{i}$ introduits en \ref{reduction} relatifs à un sommet $s$ vérifiant l'hypothèse de l'énoncé. Le (ii) du lemme  \ref{excellentscouples} nous ramène au cas où $(G,M,s)$ est l'un de ces triplets $(G'_{i},M'_{i},s_{i})$. On peut donc supposer que $G$ est simplement connexe et absolument simple.

    Appliquons la proposition \ref{casdunexcellentcouple}, introduisons l'ensemble $\Lambda$ tel que $(M,s)$ soit conjugué à $(M_{\Lambda},s_{\Lambda})$ par un élément de $G_{SC}(F^{nr})$. Il suffit de prouver que $W^{nr}(M_{\Lambda})$ a deux éléments. Remarquons que $\Delta_{a}^{nr}-\Lambda$ a au moins deux éléments parce que $M\not=G$.  L'ensemble  $\Lambda$ est excellent. Pour $\gamma\in \Delta_{a}^{nr}-\Lambda$, on dispose donc d'une symétrie $w_{\gamma}$ que l'on relève en un élément de $ Norm_{G_{SC}(F^{nr})}(T^{nr},\Lambda^{aff})$. Cet élément normalise $M_{\Lambda}$. Il n'appartient pas à ce groupe car, par construction, il n'agit pas trivialement sur $A_{\Lambda}^{nr}$. Donc son image dans $W^{nr}(M_{\Lambda})$ est non triviale, ce qui achève la preuve. $\square$ 
    
    \subsubsection{Description de $K_{s}^{M,nr,0}\backslash K_{s}^{M,nr,\dag}$}\label{descriptionKdag}
    
    {\bf On suppose que} $G$  {\bf  est simplement connexe. } Soit $M$ un $F^{nr}$-Levi de $G$ et $s$ un sommet de $Imm_{F^{nr}}(M_{AD})$. On suppose que le couple $(M,s)$ est excellent. On se propose de décrire le quotient $K_{s}^{M,nr,0}\backslash K_{s}^{M,nr,\dag}$.
    Si $M=G$, on a  $K_{s}^{M,nr,0}= K_{s}^{M,nr,\dag}$ puisque $G$ est simplement connexe. On suppose désormais que $M\not=G$.
    
    Soit $L$ un $F^{nr}$-Levi de $G$. Supposons que $M\subset L$ et que $M$ est un $F^{nr}$-Levi propre  maximal de $L$. Introduisons le groupe $L_{SC}$. Notons $M_{L,sc}$ l'image réciproque de $M$ dans ce groupe. 
    Parce que $G$ est simplement connexe, $L_{SC}$ est le groupe dérivé de $L$ et on a $M_{L,sc}=M\cap L_{SC}$. On a l'égalité    $K_{s}^{M_{L,sc},nr,\dag}=M_{L,sc}(F^{nr})\cap K_{s}^{M,nr,\dag}$. Montrons que l'on a aussi l'égalité
    
    (1) $K_{s}^{M_{L,sc},nr,0}=M_{L,sc}(F^{nr})\cap K_{s}^{M,nr,0}$.
    
    Comme on l'a dit dans \ref{normalisateur}, on dispose des homomorphismes de Kottwitz qui s'inscrivent dans un diagramme commutatif
    $$(2) \qquad \begin{array}{ccc} K_{s}^{M_{L,sc},nr,\dag}&\stackrel{w_{M_{L,sc}}}{\to}&X^*(Z(\hat{M}_{L,sc})^{I_{F}})\\ \downarrow&&\downarrow\\ K_{s}^{M,nr,\dag}&\stackrel{w_{M}}{\to}&X^*(Z(\hat{M})^{I_{F}})\\ \end{array}$$
    Admettons un instant que la   flèche verticale de droite soit  injective. Alors (1) résulte du fait que $K_{s}^{M_{L,sc},nr,0}$, resp. $K_{s}^{M,nr,0}$, est le noyau de $w_{M_{L,sc}}$, resp. $w_{M}$. Prouvons l'injectivité requise. On a $Z(\hat{M}_{L,sc})=Z(\hat{M})/Z(\hat{L})$, dont on déduit la suite exacte
    $$(3) \qquad 1\to Z(\hat{L})^{I_{F}}\to Z(\hat{M})^{I_{F}}\to Z(\hat{M}_{L,sc})^{I_{F}}.$$
    Parce que $G$ et $L_{SC}$ sont simplement connexes, ces trois groupes sont des tores de dimensions respectives $dim(A_L^{nr})-dim(A_G^{nr})$, $dim(A_M^{nr})-dim(A_G^{nr})$, $dim(A_{M_{L,sc}}^{nr})$. La dimension du terme central est la somme des dimensions des termes extr\^emes. Donc la   dernière flèche de (3) est surjective. Par dualité, la flèche verticale de droite de (2) est injective. Cela prouve (1). 
    
    On déduit de (1) un homomorphisme  injectif
    $$j_{L}:K_{s}^{M_{L,sc},nr,0}\backslash K_{s}^{M_{L,sc},nr,\dag}\to K_{s}^{M,nr,0}\backslash K_{s}^{M,nr,\dag}.$$
    
    \begin{lem}{Le groupe $K_{s}^{M,nr,0}\backslash K_{s}^{M,nr,\dag}$ est engendré par les images des homomorphismes $j_{L}$ quand $L$ décrit les $F^{nr}$-Levi de $G$ contenant $M$ comme $F^{nr}$-Levi propre maximal.}\end{lem}
    
    Preuve.
          Avec les notations de \ref{reduction}, on a l'égalité
    $$(4) \qquad K_{s}^{M,nr,0}\backslash K_{s}^{M,nr,\dag}=\prod_{i=1,...,n}K_{s_{i}}^{M_{i}',nr,0}\backslash K_{s_{i}}^{M_{i}',nr,\dag}. $$
     Les "$nr$" du membre de droite se réfèrent au fait qu'un groupe $G'_{i}$ est défini sur    $F_{i}^{nr}$ où $F_{i}$ est une extension finie modérément ramifiée de $F$. Il suffit de démontrer l'assertion du lemme pour chacun des facteurs. Autrement dit, on peut supposer que  $G$ est absolument quasi-simple. 
    
    On peut alors identifier le couple $(M,s)$ à $(M_{\Lambda},s_{\Lambda})$ pour un excellent ensemble $\Lambda\subset \Delta_{a}^{nr}$.  D'après \ref{normalisateur}(3) et (4),  on a la suite exacte
  $$1\to K_{s}^{M,nr,0}\to Norm_{G(F^{nr})}(M,s)\to W_{aff}(\Lambda^{aff})\to 1.$$
  On a aussi  l'homomorphisme naturel $W_{aff}(\Lambda^{aff})\to W^{nr}(M)=Norm_{G(F^{nr})}(M)/M(F^{nr})$. Notons $t_{aff}$ son noyau. De la suite exacte ci-dessus se déduit une suite exacte
  $$(5) \qquad 1\to K_{s}^{M,nr,0}\to K_{s}^{M,nr,\dag}\to t_{aff}\to 1.$$
 Reprenons les notations introduites dans \ref{excellentsensemblesLambda}. Posons 
 ${\bf  t}_{aff}=res(t_{aff})$.  D'après la proposition \ref{excellentsensemblesLambda} et la définition de $t_{aff}$, ${\bf t}_{aff}$ est le sous-groupe des éléments de ${\cal W}_{aff}(\Lambda^{aff})$ dont la partie linéaire est triviale, autrement dit c'est le groupe des translations contenues dans ${\cal W}_{aff}(\Lambda^{aff})$. Considérons deux $P(\Lambda^{aff})$-murs $H$ et $H'$ de $P(\Lambda^{aff})$, supposons qu'ils soient parallèles. On dispose dans ${\cal W}_{aff}(\Lambda^{aff})$ des deux symétries orthogonales relatives à ces murs, notons-les ${\bf w}_{H}$ et ${\bf w}_{H'}$. Le produit ${\bf w}_{H}{\bf w}_{H'}$ est une translation. Notons $\boldsymbol{\tau}_{aff}$ le groupe engendré par ces éléments quand $(H,H')$ décrit l'ensemble des couples de $P(\Lambda^{aff})$-murs parallèles. Montrons que
 
  (6) on a l'égalité ${\bf t}_{aff}=\boldsymbol{\tau}_{aff}$. 
  
  On a vu dans la preuve de  \ref{excellentsensemblesLambda}  que les hypothèses de \cite{B} V.3 étaient vérifiées. On applique la proposition 10 de \cite{B} V.3.10: $P(\Lambda^{aff})$ possède un point spécial. Fixons-en un, que l'on note $x$. Notons ${\cal W}_{x}$ son fixateur dans ${\cal W}_{aff}(\Lambda^{aff})$. Soit $H$ un $P(\Lambda^{aff})$-mur. La proposition 9 de loc.cit. nous dit qu'il existe un $P(\Lambda^{aff})$-mur $H_{x}$ qui est parallèle à $H$ et contient $x$. On a ${\bf w}_{H_{x}}\in {\cal W}_{x}$ et ${\bf w}_{H_{x}}{\bf w}_{H}\in \boldsymbol{\tau}_{aff}$. Donc ${\bf w}_{H}\in {\cal W}_{x}\boldsymbol{\tau}_{aff}$. Il résulte de la définition de $\boldsymbol{\tau}_{aff}$ que ce groupe est un sous-groupe distingué de ${\cal W}_{aff}(\Lambda^{aff})$. Puisque ${\cal W}_{aff}(\Lambda^{aff})$ est engendré par les symétries ${\bf w}_{H}$, la relation précédente entraîne que ${\cal W}_{aff}(\Lambda^{aff})={\cal W}_{x}\boldsymbol{\tau}_{aff}$. Soit $u\in {\bf t}_{aff}$. On écrit $u=wu'$ avec $w\in {\cal W}_{x}$ et $u'\in \boldsymbol{\tau}_{aff}$. Alors $w=u{u'}^{-1}$ est une translation qui fixe $x$ et est donc triviale. Donc $u=u'$ appartient à $\boldsymbol{\tau}_{aff}$. Cela prouve (6).
  
  En vertu de (6), il nous reste à prouver l'assertion suivante:
  
  (7) soient $H$, $H'$ deux $P(\Lambda^{aff})$ murs parallèles; alors il existe un $F^{nr}$-Levi $L$ de $G$ contenant $M$ comme $F^{nr}$-Levi maximal propre tel que ${\bf w}_{H}{\bf w}_{H'}$ appartienne à l'image de l'homomorphisme composé
  $$K_{s}^{M_{L,sc},nr,\dag}\to K_{s}^{M,nr,\dag}\to t_{aff}.$$
  
  L'ensemble $P(\Lambda^{aff})$ est un espace affine sous l'espace vectoriel ${\cal A}_{M}^{nr}$.   Par définition des $P(\Lambda^{aff})$-murs, il existe une racine affine $\beta^{aff}$ dont la racine sous-jacente $\beta$ n'appartient pas à $\Sigma^M$ et telle que $H$ soit l'annulateur de $\beta^{aff}$ dans $P(\Lambda^{aff})$. Notons $A(H)$ le sous-tore de $A_{M}^{nr}$ tel que $X_{*}(A(H))$ soit l'annulateur de $\beta$ dans $X_{*}(A_{M}^{nr})$. Posons ${\cal A}(H)=X_{*,{\mathbb R}}(A(H))$. Alors ${\cal A}(H)$ est l'espace vectoriel sous-jacent à l'espace affine $H$ et aussi à l'espace affine $H'$ puisque ces espaces affines sont parallèles. Notons $L$ le commutant de $A(H)$ dans $G$. C'est un $F^{nr}$-Levi de $G$ qui contient $M$. Il est distinct de $M$ car la racine $\beta$ intervient dans $L$ mais pas dans $M$. Le groupe $M$ est maximal propre dans $L$ car $dim(A(H))=dim(A_{M}^{nr})-1$ par définition. L'image $s_{H}=p_{L}(H)$ dans $Imm_{F^{nr}}(L_{AD})$ appartient à l'ensemble $P^L(s)$ (analogue de $P(s)$ quand on remplace $G$ par $L$).  Cet ensemble $P^L(s)$ est un espace affine de dimension $1$. La racine affine $\beta^{aff}$ n'est pas identiquement nulle sur $P^L(s)$ mais s'annule en $s_{H}$. Donc $s_{H}$ est un sommet de $Imm_{F^{nr}}(L_{AD})$. 
    Le groupe $Norm_{L_{SC,s_{H}}}(M_{L,sc,s})/M_{L,sc,s}$ a deux éléments puisque $(M_{L,sc},s)$ est $L_{SC}$-excellent d'après le lemme \ref{excellentscouples}. D'après le (iii) du lemme \ref{normalisateur}, on peut relever l'élément non trivial de ce groupe en un élément   $n_{H}\in  K_{s_{H}}^{L_{SC},nr,0}\cap Norm_{L_{SC}(F^{nr})}(M_{L,sc},s)$. Alors ${\bf w}_{H}$ est l'image de $n_{H}$ par l'homomorphisme composé
  $$Norm_{L_{SC}(F^{nr})}(M_{L,sc},s)\to Norm_{G(F^{nr})}(M,s)\to W_{aff}(\Lambda^{aff})\to {\cal W}_{aff}(\Lambda^{aff}).$$
  On construit de m\^eme un élément $n_{H'}$. Posons $n=n_{H}n_{H'}$. Son image par l'homomorphisme précédent est ${\bf w}_{H}{\bf w}_{H'}$. Comme ce dernier élément est une translation, l'action de $n$ sur ${\cal A}_{M_{L,sc}}^{nr}$ est triviale. Cela entraîne que $n\in K_{s}^{M_{L,sc},nr,\dag}$. Cela prouve (7) et le lemme. $\square$
  
  {\bf Remarque.} Supposons que $T\subset M$ et que $s$ appartienne à l'appartement associé à $T^{nr}$ dans $Imm_{F^{nr}}(M_{AD})$. La preuve ci-dessus montre que l'on peut restreindre l'ensemble des Levi $L$ de l'énoncé à ceux qui contiennent $T$.

  \bigskip
  
  {\bf Supposons de plus que} $G$ {\bf soit absolument quasi-simple}. Supposons que $(M,s)=(M_{\Lambda},s_{\Lambda})$ pour un excellent ensemble $\Lambda\subset \Delta_{a}^{nr}$ tel que $\Delta_{a}^{nr}-\Lambda$ ait deux éléments. Le Levi $M$ est un $F^{nr}$-Levi maximal propre de $G$.  Avec les notations introduites dans la preuve ci-dessus, on a l'égalité $K_{s}^{M,nr,0}\backslash K_{s}^{M,nr,\dag}\simeq t_{aff}$, cf. (5).  On a vu dans la preuve du lemme \ref{normalisateur} que ce groupe était un ${\mathbb Z}$-module libre de rang la dimension de $Z(\hat{M})^{I_{F}}$ qui est ici égale à $1$. 
 Notons   $\gamma_{1}$ et $\gamma_{2}$ les éléments de $\Delta_{a}^{nr}-\Lambda$.     A ces éléments sont associées deux symétries $w_{\gamma_{i}}\in W_{aff}(\Lambda^{aff})$ pour $i=1,2$.  
  
  (8) $t_{aff}$ est engendré par $w_{\gamma_{1}}w_{\gamma_{2}}$. 
  
  Le groupe $W_{aff}(\Lambda^{aff})$ est un groupe de Coxeter engendré par les symétries $w_{\gamma_{i}}$ pour $i=1,2$. Il est donc libre car sinon, il serait fini et ne pourrait pas contenir $t_{aff}$ qui est un ${\mathbb Z}$-module libre de rang $1$.  Tout élément $w\in W_{aff}(\Lambda^{aff})$ s'écrit donc de façon unique comme produit 
  $$ w=w_{\gamma_{i_{1}}}w_{\gamma_{i_{2}}}...w_{\gamma_{i_{m}}},$$
  où $i_{1},... i_{m}$ est une suite d'éléments de $\{1,2\}$, deux éléments successifs étant distincts. 
 Comme on l'a déjà vu,  les images  de $w_{\gamma_{i}}$ dans $W^{nr}(M)$ pour $i=1,2$ sont non triviales. Puisque ce groupe n'a que deux éléments,  un élément $w$ écrit comme ci-dessus appartient  à $t_{aff}$ si et seulement si $m$ est pair. 
 En particulier $w_{\gamma_{1}}w_{\gamma_{2}}\in t_{aff}$. Le groupe engendré par cet élément contient tous les $w$ tels que $m$ est pair et $i_{1}=1$. Pour un élément $w$ tel que $m$ est pair et $i_{1}=2$, on introduit l'élément $w'=w_{\gamma_{i_{2}}}...w_{\gamma_{i_{m}}}w_{\gamma_{i_{1}}}$. Il appartient au groupe engendré par $w_{\gamma_{1}}w_{\gamma_{2}}$.  On a $w=w_{\gamma_{2}}w'w_{\gamma_{2}}^{-1}$. La conjugaison par $w_{\gamma_{2}}$ conserve $t_{aff}$ par définition de cet ensemble. Puisque celui-ci est un ${\mathbb Z}$-module libre de rang $1$, ses seuls automorphismes sont $1$ et l'inversion. On a donc $w=w'$ ou $w={w'}^{-1}$ (en fait, il est clair que cette dernière égalité est la bonne), en tout cas $w$ appartient au groupe engendré par $w_{\gamma_{1}}w_{\gamma_{2}}$. Cela démontre (1).

 Relevons les symétries $w_{\gamma_{i}}$ en deux éléments $n_{i}\in Norm_{G(F^{nr})}(M,s)$ pour $i=1,2$. Alors $n_{1}n_{2}$ appartient à $M(F^{nr})$, donc à $K_{s}^{M,nr,\dag}$.  L'assertion (8) équivaut à
 
 (9) le groupe     $K_{s}^{M,nr,0}\backslash K_{s}^{M,nr,\dag}$ est le ${\mathbb Z}$-module libre engendré par l'image de  $n_{1}n_{2}$ dans ce groupe.

\subsubsection{Couples $(M,s_{M})$ avec $s_{M}\in S^{nr,st}(M)$}\label{couplesstables}
{\bf On suppose que } $\boldsymbol{G}$ {\bf est absolument quasi-simple.}    Soit $M$ un $F^{nr}$-Levi standard. De la paire de Borel \'epingl\'ee  $\mathfrak{E}$ se d\'eduit une telle paire $\mathfrak{E}^M=(B\cap M,T,(E_{\alpha})_{\alpha\in \Delta^M}))$, dont on d\'eduit comme en \ref{alcove} une alc\^ove $C^{nr,M}\subset Imm_{F^{nr}}(M_{AD})$. On pose $S^{st}(\bar{C}^{nr,M})=S^{nr,st}(M)\cap S(\bar{C}^{nr,M})$ (il s'agit ici de sommets de $Imm_{F^{nr}}(M_{AD})$). Dans le cas $M=G$, on a donné dans \cite{W7} paragraphe 9  une description cas pas cas de $S^{st}(\bar{C}^{nr})$. Il est aisé d'en déduire une description des couples $(M,s_{M})$ où $M$ est standard et $s_{M}\in S^{st}(\bar{C}^{nr,M})$. Remarquons que le couple $(T,s_{T})$, où $s_{T}$ est l'unique élément de $Imm(T_{AD})$ vérifie ces conditions. On suppose désormais $M\not=T$.

  On d\'ecrit ci-dessous les conditions  portant sur le diagramme ${\cal D}^M$ qui sont n\'ecessaires et suffisantes pour que  $S^{st}(\bar{C}^{nr,M})$  ne soit pas vide. On indique le nombre d'\'el\'ements de  cet ensemble et les diagrammes des groupes $M_{s_{M},SC}$ quand $s_{M}$ le d\'ecrit. On note  ${\cal D}_{s_{M}}$  le diagramme associé à $M_{s_{M},SC}$ pour simplifier la notation. 
   Pour les groupes de type $A_{n-1}$, $D_{n}$ avec $n\geq5$, $E_{6}$, on note $\theta$ l'unique automorphisme d'ordre $2$ de ${\cal D}$. Pour un groupe de type $D_{4}$, on note $\theta$ l'automorphisme de ${\cal D}$ qui permute $\alpha_{3}$ et $\alpha_{4}$ et fixe $\alpha_{1}$ et $\alpha_{2}$ et on note $\theta_{3}$ celui qui fixe $\alpha_{2}$ et envoie $(\alpha_{1},\alpha_{3},\alpha_{4})$ sur $(\alpha_{3},\alpha_{4},\alpha_{1})$. Par abus de notation, on pose $B_{1}=C_{1}=A_{1}$ et $D_{2}=A_{1}\times A_{1}$. 
 
 $(A_{n-1},nr)$ Supposons que $G$ est de type $A_{n-1}$ avec $n\geq2$ et que $G$ est d\'eploy\'e sur $F^{nr}$. Aucun $M$ ne convient.
 
    $(A_{n-1},ram)$ Supposons que $G$ est de type $A_{n-1}$ avec $n\geq3$ et que $I_{F}$ agit sur ${\cal D}$ via un homomorphisme surjectif $I_{F}\to \{1,\theta\}$.  Alors ${\cal D}^M$ est de type $A_{m-1}$ avec $3\leq m\leq n$, $m\equiv n\,mod\,2{\mathbb Z}$, et il existe $h,k\in {\mathbb N}$ avec $h=k$ ou $h=k+1$ tels que $m=h^2+k(k+1)$. On a  $\vert S^{st}(\bar{C}^{nr,M})\vert =1$ et, pour l'unique \'el\'ement $s_{M}\in  S^{st}(\bar{C}^{nr,M})$, ${\cal D}_{s_{M}}$ est de type $B_{(h^2-1)/2}\times C_{k(k+1)/2}$ si $n$ est impair, $D_{h^2/2}\times C_{k(k+1)/2}$ si $n$ est pair.

     $(B_{n})$ Supposons que $G$ est de type $B_{n}$ avec $n\geq2$. Alors ${\cal D}^M$ est de type $B_{m}$ avec $2\leq m\leq n$ et il existe $h,k\in {\mathbb N}$ tels que $k$ est pair, $h$ est impair, $\vert k-h\vert =1$ et  $2m+1=k^2+h^2$. On a  $\vert S^{st}(\bar{C}^{nr,M})\vert =1$ et, pour l'unique \'el\'ement $s_{M}\in  S^{st}(\bar{C}^{nr,M})$, ${\cal D}_{s_{M}}$ est de type $D_{k^2/2}\times B_{(h^2-1)/2}$. 
     
           $(C_{n})$ Supposons que $G$ est de type $C_{n}$ avec $n\geq2$. Alors ${\cal D}^M$ est de type $C_{m}$ avec $2\leq m\leq n$ et il existe $k\in {\mathbb N}$ tel  que   $m=k(k+1)$. On a  $\vert S^{st}(\bar{C}^{nr,M})\vert =1$ et, pour l'unique \'el\'ement $s_{M}\in  S^{st}(\bar{C}^{nr,M})$, ${\cal D}_{s_{M}}$ est de type $C_{k(k+1)/2}\times C_{k(k+1)/2}$.

$(D_{n},nr)$ Supposons que $G$ est de type $D_{n}$ avec $n\geq4$ et que $G$ est d\'eploy\'e sur $F^{nr}$. Alors ${\cal D}^M$ est de type $D_{m}$ avec $4\leq m\leq n$ et il existe $k\in {\mathbb N}$ tel  $k$ est pair et   $m=k^2$. On a  $\vert S^{st}(\bar{C}^{nr,M})\vert =1$ et, pour l'unique \'el\'ement $s_{M}\in  S^{st}(\bar{C}^{nr,M})$, ${\cal D}_{s_{M}}$ est de type $D_{k^2/2}\times D_{k^2/2}$.

$(D_{n},ram)$ Supposons que $G$ est de type $D_{n}$ avec $n\geq4$ et que $I_{F}$ agit sur ${\cal D}$ via un homomorphisme surjectif $I_{F}\to \{1,\theta\}$. Alors ${\cal D}^M$ est de type $(D_{m},ram)$ avec $9\leq m\leq n$ et il existe $k\in {\mathbb N}$ tel  $k$ est impair et   $m=k^2$. On a  $\vert S^{st}(\bar{C}^{nr,M})\vert =1$ et, pour l'unique \'el\'ement $s_{M}\in  S^{st}(\bar{C}^{nr,M})$, ${\cal D}_{s_{M}}$ est de type $B_{(k^2-1)/2}\times B_{(k^2-1)/2}$.

$(D_{4},3-ram)$ Supposons que $G$ est de type $D_{4}$ et que $I_{F}$ agit sur ${\cal D}$ via un homomorphisme surjectif $I_{F}\to \{1,\theta_{3},\theta_{3}^2\}$. Alors $M=G$. On a  $\vert S^{st}(\bar{C}^{nr,M})\vert =2$ et les deux diagrammes ${\cal D}_{s_{M}}$ sont  de type $G_{2}$ et $A_{1}\times A_{1}$.

$(E_{6},nr)$ Supposons que $G$ est de type $E_{6}$ et que $G$ est d\'eploy\'e sur $F^{nr}$. Alors ${\cal D}^M$ est de type $E_{6}$ ou $D_{4}$. Si ${\cal D}^M$ est de type $E_{6}$, c'est-\`a-dire $M=G$, 
on a  $\vert S^{st}(\bar{C}^{nr,M})\vert =1$ et, pour l'unique \'el\'ement $s_{M}\in  S^{st}(\bar{C}^{nr,M})$, ${\cal D}_{s_{M}}$ est de type  $A_{2}\times A_{2}\times A_{2}$. Si ${\cal D}^M$ est de type $D_{4}$, on a  $\vert S^{st}(\bar{C}^{nr,M})\vert =1$ et, pour l'unique \'el\'ement $s_{M}\in  S^{st}(\bar{C}^{nr,M})$, ${\cal D}_{s_{M}}$ est de type $A_{1}\times A_{1}\times A_{1}\times A_{1}$.

$(E_{6},ram)$ Supposons que $G$ est de type $E_{6}$ et que $I_{F}$ agit sur ${\cal D}$ via un homomorphisme surjectif $I_{F}\to \{1,\theta\}$. Alors ${\cal D}^M$ est de type $(E_{6},ram)$ ou $(A_{5},ram)$. Si ${\cal D}^M$ est de type $(E_{6},ram)$, c'est-\`a-dire $M=G$, 
on a  $\vert S^{st}(\bar{C}^{nr,M})\vert =2$ et les deux diagrammes ${\cal D}_{s_{M}}$ sont de type $F_{4}$ et $A_{2}\times A_{2}$.  Si ${\cal D}^M$ est de type $(A_{5},ram)$, on a  $\vert S^{st}(\bar{C}^{nr,M})\vert =1$ et, pour l'unique \'el\'ement $s_{M}\in  S^{st}(\bar{C}^{nr,M})$, ${\cal D}_{s_{M}}$ est de type $A_{1}\times A_{1}\times A_{1}$.

$(E_{7})$ Supposons que $G$ est de type $E_{7}$. Alors ${\cal D}^M$ est de type $E_{7}$, $E_{6}$ ou $D_{4}$. Si ${\cal D}^M$ est de type $E_{7}$, c'est-\`a-dire $M=G$, 
 on a  $\vert S^{st}(\bar{C}^{nr,M})\vert =1$ et, pour l'unique \'el\'ement $s_{M}\in  S^{st}(\bar{C}^{nr,M})$, ${\cal D}_{s_{M}}$ est de type $A_{3}\times A_{3}\times A_{1}$. Les autres $M$ sont comme dans le type $(E_{6},nr)$.

   $(E_{8})$ Supposons que $G$ est de type $E_{8}$. Alors ${\cal D}^M$ est de type $E_{8}$, $E_{7}$, $E_{6}$ ou $D_{4}$. Si ${\cal D}^M$ est de type $E_{8}$, c'est-\`a-dire $M=G$, on a $\vert S^{st}(\bar{C}^{nr,M})\vert =7$ et les sept diagrammes ${\cal D}_{s_{M}}$ sont de type $E_{8}$, $D_{8}$, $A_{1}\times A_{2}\times A_{5}$, $A_{4}\times A_{4}$, $D_{5}\times A_{3}$, $E_{6}\times A_{2}$, $E_{7}\times A_{1}$. Les autres $M$ sont comme dans le type $E_{7}$. 
     
   $(F_{4})$ Supposons que $G$ est de type $F_{4}$. Alors ${\cal D}^M$ est de type $F_{4}$ ou $B_{2}=C_{2}$. Si ${\cal D}^M$ est de type $F_{4}$, c'est-\`a-dire $M=G$, on a $\vert S^{st}(\bar{C}^{nr,M})\vert =5$   et les cinq diagrammes ${\cal D}_{s_{M}}$ sont de type $F_{4}$, $A_{1}\times C_{3}$, $A_{2}\times A_{2}$, $A_{3}\times A_{1}$, $B_{4}$. Si ${\cal D}^M$ est de type $B_{2}$, on a $\vert S^{st}(\bar{C}^{nr,M})\vert =1$ et, pour l'unique \'el\'ement $s_{M}\in  S^{st}(\bar{C}^{nr,M})$, ${\cal D}_{s_{M}}$ est de type $A_{1}\times A_{1}$.

   $(G_{2})$ Supposons que $G$ est de type $G_{2}$. Alors  $M=G$.  On a  $\vert S^{st}(\bar{C}^{nr,M})\vert =3$  et les trois diagrammes ${\cal D}_{s_{M}}$ sont de type $G_{2}$, $A_{1}\times A_{1}$, $A_{2}$. 
   
   Remarquons que 
   
   (1) $M_{SC}$ est toujours absolument quasi-simple.
   
   D'autre part, les diagrammes ${\cal D}^M$ sont tous conserv\'es par le groupe $Aut({\cal D})$. Il en r\'esulte que
   
   (2) si $G$ est quasi-d\'eploy\'e sur $F$, tout $F^{nr}$-Levi standard $M$ tel que $S^{nr,st}(M)$  n'est pas vide   est un $F$-Levi. 
   
   Soient $M$ un $F^{nr}$-Levi standard et $s_{M}\in S^{st}(\bar{C}^{nr,M})$. On va prouver
   
   (3) il existe un unique $\Lambda\subsetneq\Delta_{a}^{nr}$ de sorte que les couples $(M,s_{M})$ et $(M_{\Lambda},s_{\Lambda})$ soient conjugués par un élément de $G_{AD}(F^{nr})$.

   Un tel ensemble existe d'après le (i) de la proposition \ref{casdunexcellentcouple} et \ref{excellentscouples}(1).  Fixons un tel $\Lambda$. D'après le (ii) de la proposition \ref{excellentscouples}, il suffit de prouver que $\Lambda$ est conservé par l'action de $\Omega^{nr}$. On va en fait prouver qu'il est conservé par tout $Aut({\cal D}_{a}^{nr})$. Pour cela, on va déterminer $\Lambda$ dans chaque cas. Supposons d'abord que $M=T$. Alors $\Lambda=\emptyset$ et l'assertion est claire. Supposons $M=G$. On peut supposer que $\Lambda$ est le complémentaire de l'unique élément de $\Delta_{a}^{nr}$ dont est issu le sommet $s_{G}=s_{M}$. D'après \cite{W7} 9(4), $s_{G}$ est fixé par $Aut({\cal D}_{a}^{nr})$ donc son complémentaire est conservé par l'action de ce groupe.  On suppose maintenant $M\not=T$, $M\not=G$. On note ${\cal D}_{\Lambda}$ le diagramme obtenu en supprimant de ${\cal D}_{a}^{nr}$ les sommets associés aux racines hors de $\Lambda$. C'est le diagramme de Dynkin de $M_{\Lambda,s_{\Lambda}}$ donc aussi celui de $M_{s_{M}}$, que l'on a noté ${\cal D}_{s_{M}}$ ci-dessus. En supposant fixée une structure euclidienne sur $X^*(T)\otimes_{{\mathbb Z}}{\mathbb R}$, vérifiant les propriétés habituelles, les racines dans ${\cal D}_{\Lambda}$ comme dans ${\cal D}_{s_{M}}$ se retrouvent affectées d'une longueur et l'isomorphisme entre ces deux diagrammes préserve les longueurs. D'autre part, l'ensemble $\Lambda$ doit \^etre excellent. Cela a la conséquence suivante. Soit $\Gamma$ un ensemble non vide de $\Delta_{a}^{nr}-\Lambda$. Considérons une  composante connexe  $\Theta$ de $\Delta_{a}^{nr}-\Gamma$  et supposons qu'elle  soit de type $A_{n-1}$. On note $\alpha_{1},...,\alpha_{n-1}$ les racines intervenant dans cette composante. Alors  il existe deux entiers $m,d\geq1$ tels que $md=n$ et  $\Theta-\Lambda\cap \Theta=\{\alpha_{m},\alpha_{2m},...,\alpha_{(d-1)m}\}$ (cela résulte du fait qu'un  Levi standard d'un groupe $GL(n)$ qui n'est conjugué à aucun autre Levi standard que lui-m\^eme est de la forme $GL(m)^d$ avec $m,d$ comme ci-dessus). On voit alors cas par cas qu'il y a un unique ensemble $\Lambda$ possible qui soit excellent et pour lequel ${\cal D}_{\Lambda}$ soit isomorphe à ${\cal D}_{s_{M}}$, l'isomorphisme préservant les longueurs. On les décrit ci-dessous, en supposant encore que $M$ n'est ni $T$, ni $G$.  Rappelons que l'on a num\'erot\'e les \'el\'ements de $\Delta_{a}^{nr}$ soit par des entiers, soit par des familles d'entiers dans le cas o\`u $G$ n'est pas d\'eploy\'e sur $F^{nr}$.  Notons symboliquement $\underline{i}$ de tels indices.  On note $I(\Lambda)$ l'ensemble des indices $\underline{i}$ tels que $\beta_{\underline{i}}\in \Lambda$.  Les entiers $h$ et $k$ figurant ci-apr\`es sont ceux de la description donnée plus haut.
    
  Supposons $G$ de type $(A_{n-1},ram)$ et $M$ est de type $A_{m-1}$, avec $3\leq m<n$. Si $m$ est impair, $I(\Lambda)=\{0,(1,n-1),...,(k(k+1)/2-1,n+1-k(k+1)/2)\}\cup\{( (n-h^2)/2+1,(n+h^2)/2-1),...,((n-1)/2,(n+1)/2)\}$. Si $m$ est pair, $I(\Lambda)=\{0,(1,n-1),...(h^2/2-1,n+1-h^2/2)\}\cup\{((n-k(k+1))/2+1,(n+k(k+1))/2-1),...,(n/2-1,n/2+1),n/2\}$.
    
  Supposons $G$ de type $B_{n}$  et $M$ de type $B_{m}$ avec $2\leq m<n$.  Alors $I(\Lambda)=\{ 0,...,k^2/2-1\}\cup\{n+1-(h^2-1)/2,...,n\}$.  
  
Supposons $G$ de type $C_{n}$  et $M$ de type $C_{m}$ avec $2\leq m<n$.  Alors $I(\Lambda)=\{0,...,k(k+1)/2-1\}\cup\{n+1-k(k+1)/2,...,n\}$.  

Supposons $G$ de type $(D_{n},nr)$  et  $M$  de type $D_{m}$ avec $4\leq m<n$.  Alors $I(\Lambda)=\{0,...,k^2/2-1\}\cup\{n+1-k^2/2,...,n\}$.  

 Supposons $G$ de type $(D_{n},ram)$ et $M$ est de type $D_{m}$, avec $9\leq m<n$.  Alors $I(\Lambda)=\{0,...,(k^2-3)/2\}\cup \{n-(k^2-1)/2,...,n-2,(n-1,n)\}$. 

Supposons $G$ de type $(E_{6},nr)$ et $M$  de type $D_{4}$. Alors $I(\Lambda)=\{0,1,4,6\}$.  

Supposons $G$ de type $(E_{6},ram)$ et $M$ de type $A_{5}$. Alors $I(\Lambda)=\{0,2,(3,5)\}$. 

Supposons $G$ de type $E_{7}$. Si  $M$ est  de type $E_{6}$,  $I(\Lambda)=\{0,1,2,4,6,7\}$. Si $M$ est de type $D_{4}$, $I(\Lambda)=\{0,3,5,7\}$. 
 
Supposons $G$ de type $E_{8}$. Si  $M$ est  de type $E_{7}$, $I(\Lambda)=\{0,1,2,4,5,7,8\}$. Si $M$ est de type $E_{6}$, $I(\Lambda)=\{0,1,3,5,6,8\}$. Si $M$ est de type $D_{4}$, $I(\Lambda)=\{0,2,5,7\}$.

Supposons $G$ de type $F_{4}$ et $M$ est de type $C_{2}$. Alors $I(\Lambda)=\{0,2\}$.

   \subsection{  $F^{nr}$-Levi et conjugaison stable}
 
 \subsubsection{Bons \'el\'ements\label{bonselements}}
 Nous reprenons la d\'efinition des bons \'el\'ements de \cite{KM}, en nous limitant au cas de profondeur $0$. Pour tout sous-tore $T$ de $G$ (pas forc\'ement maximal), on note $\Sigma(T)$ l'ensemble des racines de $T$ dans $\mathfrak{g}$. 
 
 Soit $X\in \mathfrak{g}(F)$. Supposons $X$ semi-simple et fixons un sous-tore maximal $T$ de $G$ d\'efini sur $F$ tel que $X\in \mathfrak{t}(F)$.   Consid\'erons la condition suivante: pour toute racine $\beta\in \Sigma(T)$, on a $val_{F}(\beta(X))=0$ ou $\beta(X)=0$. Que cette condition soit v\'erifi\'ee ne d\'epend pas du choix de $T$. On dit que $X$ est un bon \'el\'ement si et seulement si elle l'est.

 Soit $X$ un bon \'el\'ement, posons $H=G_{X}$. Les propri\'et\'es suivantes sont v\'erifi\'ees:
 
 (1) soient $Y,Y'\in \mathfrak{h}_{tn}(F)$; alors l'ensemble des $g\in G(F)$ tels que $g(X+Y)=X+Y'$ co\"{\i}ncide avec celui des $h\in H(F)$ tels que $hYh^{-1}=Y'$; 
cf. \cite{KM} 2.1.5(5).

 Dans les deux assertions suivantes, les propri\'et\'es des \'el\'ements $Y$ sont relatives au groupe ambiant $H$, celles de $X+Y$ sont relatives au groupe ambiant $G$: 

(2) pour $Y\in \mathfrak{h}_{tn}(F)$, on a  $Z_{H}(Y)=Z_{G}(X+Y)$; $Y$ est semi-simple, resp. r\'egulier,  si et seulement si il en est de m\^eme de $X+Y$;  si $A_{H}=A_{G}$, $Y$ est semi-simple r\'egulier elliptique si et seulement s'il en est de m\^eme de $X+Y$;

(3) soit $Y\in \mathfrak{h}_{tn}(F)$ un \'el\'ement semi-simple r\'egulier; alors $d^G(X+Y)=d^H(Y)$;  

(4) soient $Y,Y'\in \mathfrak{h}_{tn}(F)$; on suppose que $X+Y$ et $X+Y'$ sont semi-simples; alors $Y$ et $Y'$ sont stablement conjugu\'es si et seulement $X+Y$ et $X+Y'$ sont stablement conjugu\'es. 

La propri\'et\'e (2) r\'esulte imm\'ediatement de (1) et la propri\'et\'e (4) aussi, en \'etendant le corps de base. La propri\'et\'e (3) r\'esulte d'un calcul facile. La propri\'et\'e suivante est bien connue, on en donne une preuve faute d'avoir trouv\'e une r\'ef\'erence:

(5) soit $\omega\subset \mathfrak{g}(F)$ un sous-ensemble compact; alors l'ensemble des $g\in G(F)$ tels que $g\omega g^{-1}\cap (X+\mathfrak{h}_{tn}(F))\not=\emptyset$ a une projection compacte dans $H(F)\backslash G(F)$.

Preuve. Pour tout sous-ensemble $E$ de $\mathfrak{g}(F)$, notons $\Omega(E)$ l'ensemble des $g\in G(F)$ tels que $g\omega g^{-1}\cap E\not=\emptyset$.  On doit prouver que la projection de $\Omega(X+\mathfrak{h}_{tn}(F))$ dans $H(F)\backslash G(F)$ est compacte. 

Soit $Y\in \mathfrak{g}(F)$ un \'el\'ement semi-simple. D'apr\`es \cite{HC} lemme 19, on peut fixer un voisinage $V_{Y}$ de $Y$ dans $\mathfrak{g}_{Y}(F)$, qui est  ouvert et compact et qui v\'erifie la condition:  l'image de $\Omega(V_{Y})$   dans $G_{Y}(F)\backslash G(F)$ est relativement compacte. Notons $\mathfrak{h}_{tn,ss}(F)$ l 'ensemble des \'el\'ements semi-simples de $\mathfrak{h}_{tn}(F)$. 
Pour $Y\in X+\mathfrak{h}_{tn,ss}(F)$, on a $G_{Y}=H_{Y}\subset H$ d'apr\`es (2) donc l'image de $\Omega(V_{Y})$   dans $H(F)\backslash G(F)$  est relativement compacte. Soit $W$ un voisinage de $1$ dans $H(F)$. L'ensemble $\{wZw^{-1}; w\in W, Z\in V_{Y}\}$ contient un voisinage de $Y$ dans $\mathfrak{h}(F)$. Fixons un voisinage ouvert $U_{Y}$ de $Y$ dans $\mathfrak{h}(F)$ qui est contenu dans l'ensemble pr\'ec\'edent. L'image de  $\Omega(U_{Y})$ dans $H(F)\backslash G(F)$ est contenue dans celle de $\Omega(V_{Y})$, donc est relativement compacte.  

 L'ensemble  $\mathfrak{h}_{tn,ss}(F)$ est compact modulo conjugaison par $H(F)$: il n'y a qu'un nombre fini de classes de conjugaison par $H(F)$ de sous-tores maximaux de $H$ d\'efinis sur $F$ et, pour tout tel sous-tore $T$, $\mathfrak{t}(F)\cap \mathfrak{h}_{tn}(F)$ est compact. Il en est de m\^eme de $X+\mathfrak{h}_{tn,ss}(F)$.  Fixons donc un ensemble compact $C\subset X+ \mathfrak{h}_{tn,ss}(F)$    tel que tout \'el\'ement de $X+ \mathfrak{h}_{tn,ss}(F)$ soit conjugu\'e \`a un \'el\'ement de $C$ par un \'el\'ement de $H(F)$.   Puisque $C$ est compact, on peut fixer une famille finie ${\cal Y}\subset C$ telle qu'en posant $U=  \cup_{Y\in {\cal Y}}U_{Y}$, on ait $C\subset U$. Montrons que

(6) tout \'el\'ement de $X+\mathfrak{h}_{tn}(F)$ est conjugu\'e \`a un \'el\'ement de $U$ par un \'el\'ement de $H(F)$. 

En effet, soit $Z\in X+\mathfrak{h}_{tn}(F)$,  notons $Z_{s}$   sa partie semi-simple. On a  $Z_{s}\in X+\mathfrak{h}_{tn,ss}(F)$. A conjugaison pr\`es par $H(F)$, on peut supposer $Z_{s}\in C\subset U$. La classe de conjugaison de $Z$ par $G_{Z_{s}}(F)=H_{Z_{s}-X}(F)$ contient un point arbitrairement proche de $Z_{s}$. 
Puisque $U$ est  ouvert dans $\mathfrak{h}(F)$, cette classe contient un \'el\'ement de $U$. D'o\`u (6).

 D'apr\`es (6), on a $\Omega(X+\mathfrak{h}_{tn}(F))\subset H(F)\Omega(U)=\cup_{Y\in {\cal Y}}H(F)\Omega(U_{Y})$. L'image de $\Omega(X+\mathfrak{h}_{tn}(F))$ dans $H(F)\backslash G(F)$ est donc contenue dans la r\'eunion sur $Y\in {\cal Y}$ des images de $\Omega(U_{Y})$. Ce dernier ensemble est r\'eunion finie d'ensembles relativement compacts, donc l'image de $\Omega(X+\mathfrak{h}_{tn}(F))$ dans $H(F)\backslash G(F)$ est relativement compacte. 
 
  Il reste \`a prouver qu'elle est ferm\'ee. Soit $(g_{i})_{i\in {\mathbb N}}$ une suite d'\'el\'ements de $\Omega(X+\mathfrak{h}_{tn}(F))$ telle que les projections des $g_{i}$ dans $H(F)\backslash G(F)$ aient une limite dans ce quotient. Puisque l'application $G(F)\to H(F)\backslash G(F)$ est partout submersive, on peut fixer des d\'ecompositions $g_{i}=h_{i}z_{i}$ de sorte que $h_{i}\in H(F)$ et la suite $(z_{i})_{i\in {\mathbb N}}$ soit convergente. Notons $z$ sa limite. Il suffit de prouver que $z$ appartient \`a $\Omega(X+\mathfrak{h}_{tn}(F))$. Pour tout $i\in {\mathbb N}$, fixons $Z_{i}\in \omega$ tel que $g_{i}Z_{i}g_{i}^{-1}\in X+\mathfrak{h}_{tn}(F)$. Quitte \`a extraire une sous-suite, on peut supposer que la suite $(Z_{i})_{i\in {\mathbb N}}$ converge vers un point $Z\in \omega$. La suite $(z_{i}Z_{i}z_{i}^{-1})_{i\in {\mathbb N}}$  converge alors vers $zZz^{-1}$. Puisque les \'el\'ements de cette suite appartiennent \`a l'ensemble ferm\'e $X+\mathfrak{h}_{tn}(F)$,  sa limite $zZz^{-1}$ appartient aussi \`a cet ensemble, donc $z\in \Omega(X+\mathfrak{h}_{tn}(F))$. Cela ach\`eve de prouver (5).

\subsubsection{El\'ements de r\'eduction r\'eguli\`ere \label{reductionreguliere}}
On note ${\cal L}_{F}^{nr}$ l'ensemble des $F^{nr}$-Levi de $G$ qui sont définis sur $F$. 
   Soient $H\in {\cal L}_{F}^{nr}$ et $X\in \mathfrak{a}_{H}^{nr}(F)$.  Disons que $X$ est de r\'eduction r\'eguli\`ere si et seulement si   $val_{F}(\beta(X))=0$ pour tout $\beta\in \Sigma(A_{H}^{nr})$.   

\begin{lem}{Soit  $H\in {\cal L}_{F}^{nr}$.

(i) Soit  $X\in \mathfrak{a}_{H}^{nr}(F)$ un \'el\'ement de r\'eduction r\'eguli\`ere. Alors $X$ est un bon \'el\'ement et $G_{X}=H$.

(ii) Soient  $X\in \mathfrak{a}_{H}^{nr}(F)$ un \'el\'ement de r\'eduction r\'eguli\`ere et $F'$ une extension finie de $F$.  Alors $X$ est encore un \'el\'ement de r\'eduction r\'eguli\`ere quand on remplace le corps de base $F$ par $F'$.

(iii) Supposons que $H$ soit un $F$-Levi. Soit $X\in \mathfrak{a}_{H}(F)$, supposons que, pour tout $\beta\in \Sigma(A_{H})$, on ait $val_{F}(\beta(X))=0$. Alors $X$ est un \'el\'ement de r\'eduction r\'eguli\`ere. }\end{lem}

Preuve. Soit $X $  comme en (i). Fixons un sous-tore maximal $T$ de $H$ d\'efini sur $F$. On a les inclusions $A_{H}^{nr}\subset Z(H)^0\subset T$. Pour $\beta\in \Sigma(T)$, $\beta$ peut intervenir dans $\mathfrak{h}$, alors $\beta$ est nulle sur $Z(H)^0$ et $\beta(X)=0$. Sinon, puisque $H$ est le commutant de $A_{H}^{nr}$, la restriction de $\beta$ \`a $A_{H}^{nr}$ n'est pas triviale et est donc un \'el\'ement de $\Sigma(A_{H}^{nr})$. Alors $val_{F}(\beta(X))=0$.   Cela d\'emontre que $X$ est un bon \'el\'ement et que la composante neutre de  son commutant est $H$. 

Le m\^eme argument vaut pour (ii) et (iii): pour (ii), toute racine $\beta\in \Sigma(A_{H,F'}^{nr})$ se restreint en une racine de $A_{H}^{nr}$; pour (iii), toute racine $\beta\in \Sigma(A_{H}^{nr})$ se restreint en une racine de $A_{H}$. $\square$

\subsubsection{Existence d'éléments entiers et de réduction régulière}\label{existence}

\begin{lem}{Soit $H\in {\cal L}_{F}^{nr}$. Il existe des éléments de $\mathfrak{a}_{H}^{nr}(F)$ qui sont entiers et de réduction régulière.}\end{lem}

Preuve. Soit $X\in \mathfrak{a}_{H}^{nr}(F)$, écrivons $X=Z+X_{sc}$, avec $Z\in \mathfrak{a}_{G}^{nr}(F)$ et $X_{sc}\in \mathfrak{a}_{H_{sc}}^{nr}(F)$. Alors $X$ est entier et de réduction régulière si et seulement si $X_{sc}$ a les m\^emes propriétés et que $Z$ est entier. De plus, que $X_{sc}$ soit de réduction régulière implique qu'il est entier. On peut donc supposer que $G$ est simplement connexe et prouver qu'il existe des éléments de réduction régulière de $\mathfrak{a}_{H}^{nr}(F)$.

Considérons les ${\mathbb Z}$-modules $X_{*}(A_{H}^{nr})$ et  $X^*(A_{H}^{nr})$. Ils sont munis d'une action de Frobenius. Introduisons le tore ${\bf A}$ défini sur ${\mathbb F}_{q}$  tel que $X_{*}({\bf A})$ et $X^*({\bf A})$, munis de leurs actions de Frobenius, s'identifient aux modules précédents. Disons qu'un élément de $\boldsymbol{\mathfrak{a}}({\mathbb F}_{q})$ est régulier s'il n'est annulé par aucun élément de $\Sigma(A_{H}^{nr})\subset X^{*}(A_{H}^{nr})=X^*({\bf A})$. 
Le groupe $\mathfrak{a}_{H,ent}^{nr}(F)$ des éléments entiers de $\mathfrak{a}_{H}^{nr}(F)$ est celui des points fixes de Frobenius dans $X_{*}(A_{H}^{nr})\otimes_{{\mathbb Z}}\mathfrak{o}_{F^{nr}}$. L'homomorphisme de réduction de $\mathfrak{o}_{F^{nr}}$ dans $\bar{{\mathbb F}}_{q}$ envoie surjectivement $\mathfrak{a}_{H,ent}^{nr}(F)$ sur $\boldsymbol{\mathfrak{a}}({\mathbb F}_{q})$. Un élément de $\mathfrak{a}_{H,ent}^{nr}(F)$ est de réduction régulière si et seulement si sa réduction dans $\boldsymbol{\mathfrak{a}}({\mathbb F}_{q})$ est régulière. Il nous suffit  donc de prouver qu'il existe des éléments réguliers dans $\boldsymbol{\mathfrak{a}}({\mathbb F}_{q})$. 

Nous utiliserons la propriété générale suivante:

(1) soit $d\in {\mathbb N}$ et soit $P(X_{1},...,X_{d})$ un polyn\^ome en $d$ variables, à coefficients dans ${\mathbb F}_{q}$; supposons que $P$ soit non nul et que le degré total de $P$ soit strictement inférieur à $q$; alors il existe $(x_{1},...,x_{d})\in {\mathbb F}_{q}^{d}$ tel que $P(x_{1},...,x_{d})\not=0$. 

Preuve de (1). On raisonne par récurrence sur $d$. Si $d=0$, $P$ est une constante non nulle. Supposons $d\geq1$. On considère $P$ comme un polyn\^ome en $X_{d}$, à coefficients dans ${\mathbb F}_{q}[X_{1},...,X_{d-1}]$. Il y a au moins un coefficient de ce polyn\^ome qui est non nul. Appliquant l'hypothèse de récurrence, on peut fixer $(x_{1},...,x_{d-1})\in {\mathbb F}_{q}^{d-1}$ de sorte que ce coefficient soit non nul en ce point, a fortiori  que le polyn\^ome $P(x_{1},...,x_{d-1},X_{d})$ ne soit pas identiquement nul. Le nombre de racines de celui-ci est inférieur ou égal à son degré, donc est strictement inférieur à $q$ d'après l'hypothèse. Il y a donc un élément $x_{d}\in {\mathbb F}_{q}$ en lequel ce polyn\^ome ne s'annule pas. La famille $(x_{1},...,x_{d})$ vérifie les conditions requises. Cela prouve (1). 

Démontrons maintenant:

(2) il existe une famille linéairement indépendante $(\beta_{1},...,\beta_{t})$ d'éléments de $\Sigma(A_{H}^{nr})$ telle que tout $\beta\in \Sigma(A_{H}^{nr})$ s'écrive $\beta=\sum_{i=1,...,t}n_{i}\beta_{i}$, où les $n_{i}$ sont des entiers tous positifs ou nuls ou tous négatifs ou nuls et où
$1\leq \sum_{i=1,...,t}\vert n_{i}\vert \leq h(G)-1$.

Fixons un épinglage $\mathfrak{E}=(B,T,(E_{\alpha})_{\alpha\in \Delta})$ de $G$ défini sur $F^{nr}$. La 
propriété (2) ne dépend que de la classe de conjugaison de $H$ par $F^{nr}$. On peut donc supposer que $H$ est un $F^{nr}$-Levi standard associé à un sous-ensemble $\Delta^H\subset \Delta$. Pour tout $\alpha\in \Sigma$, notons $\tau_{H}(\alpha)$ la restriction de $\alpha$ à $ X_{*}(A_{H}^{nr})\subset X_{*}(T)$.  Les éléments de $\Sigma(A_{H}^{nr})$ sont les éléments non nuls  dans l'ensemble $\{\tau_{H}(\alpha); \alpha\in \Sigma\}$. On a $\tau_{H}(\alpha)=0$ si $\alpha\in \Delta^H$.   Par contre, si l'on note $(\beta_{1},...,\beta_{t})$ les éléments  de l'ensemble $\{\tau_{H}(\alpha); \alpha\in \Delta-\Delta^H\}$, la famille $(\beta_{1},...,\beta_{t})$ est linéairement indépendante. Soit $\beta\in \Sigma(A_{H}^{nr})$, fixons $\gamma\in \Sigma$ tel que $\tau_{H}(\gamma)=\beta$. Ecrivons $\gamma=\sum_{\alpha\in \Delta}n_{\alpha}\alpha$. On déduit de ce qui précède que $\beta=\sum_{i=1,...,t}n_{i}\beta_{i}$ où, pour tout $i$, $n_{i}=\sum_{\alpha\in \Delta; \tau_{H}(\alpha)=\beta_{i}}n_{\alpha}$. Les propriétés requises des $n_{i}$ résultent des m\^emes propriétés  bien connues des $n_{\alpha}$, à savoir que ce sont des entiers tous positifs ou nuls ou tous négatifs ou nuls et que $\sum_{\alpha\in \Delta}\vert n_{\alpha}\vert \leq h(G)-1$. 
Cela prouve (2). 

 On considère maintenant les éléments de $\Sigma(A_{H}^{nr})$ comme des formes linéaires sur $\boldsymbol{\mathfrak{a}}(\bar{{\mathbb F}}_{q})$. 
On fixe une famille $(\beta_{1},...,\beta_{t})$ vérifiant (2). L'hypothèse que $G$ est  semi-simple implique que $t$ est égal  à la dimension du $\bar{{\mathbb F}}_{q}$-espace vectoriel $\boldsymbol{\mathfrak{a}}(\bar{{\mathbb F}}_{q})$. Fixons une base $\{Y_{i}; i=1,...,t\}$ de $\boldsymbol{\mathfrak{a}}({\mathbb F}_{q})$ sur ${\mathbb F}_{q}$. C'est aussi une base de $\boldsymbol{\mathfrak{a}}(\bar{{\mathbb F}}_{q})$ sur $\bar{{\mathbb F}}_{q}$.  Notons $M$ la matrice carrée $(\beta_{i}(Y_{j}))_{i,j=1,...,t}$. Puisque la famille $(\beta_{1},...,\beta_{t})$  est linéairement indépendante, on a $det(M)\not=0$. 
Pour $i=1,...,t$, notons $\beta_{i,F}$ la restriction de $\beta_{i}$ à $\boldsymbol{\mathfrak{a}}({\mathbb F}_{q})$.  Ce sont des applications ${\mathbb F}_{q}$-linéaires de $\boldsymbol{\mathfrak{a}}({\mathbb F}_{q})$ dans $\bar{{\mathbb F}}_{q}$. On peut évidemment fixer un sous-corps fini ${\mathbb F}_{q^d}$ de sorte qu'elles prennent leurs valeurs dans ce sous-corps.   Fixons une base $\{\xi_{1},...,\xi_{d}\}$ de ${\mathbb F}_{q^d}$ sur ${\mathbb F}_{q}$.
    Pour tout $i,j=1,...,t$, écrivons $\beta_{i,F}(Y_{j})=\sum_{k=1,...,d} b_{i,j,k}\xi_{k}$  avec des $b_{i,j,k}\in {\mathbb F}_{q}$ et introduisons le polyn\^ome à $d$ variables $m_{i,j} (X_{1},...,X_{d})=\sum_{k=1,...,d}b_{i,j,k}X_{k}$. On introduit la matrice carrée $M(X_{1},...,X_{d})=(m_{i,j}(X_{1},...,X_{d}))_{i,j=1,...,t}$. Pour tous $i,j$, on a $\beta_{i,F}(Y_{j})=\beta_{i}(Y_{j})=m_{i,j}(\xi_{1},...,\xi_{d})$, donc $M=M(\xi_{1},...,\xi_{d})$. Le déterminant $det(M(X_{1},...,X_{d}))$ est un polyn\^ome à coefficients dans ${\mathbb F}_{q}$. Son degré total est $t$, qui est strictement inférieur à $q$ d'après l'hypothèse $(Hyp)_{1}(p)$.  Ce polyn\^ome n'est pas nul car $det(M(\xi_{1},...,\xi_{d}))=det(M)\not=0$. Appliquant (1), on peut  fixer des éléments $x_{1},...,x_{d}\in {\mathbb F}_{q}$ de sorte que $det(M(x_{1},...,x_{d}))\not=0$. On définit l'application  ${\mathbb F}_{q}$-linéaire $\ell:{\mathbb F}_{q^d}\to {\mathbb F}_{q}$ par les égalités $\ell(\xi_{k})=x_{k}$ pour tout $k=1,...,d$. Alors le déterminant de la matrice carrée $(\ell\circ\beta_{i,F}(Y_{j}))_{i,j=1,...,t}$ est non nul, c'est-à-dire que les formes linéaires $\ell\circ \beta_{i,F}$ pour $i=1,...,t$ sont linéairement indépendantes. Il s'agit alors d'une base du dual de $\boldsymbol{\mathfrak{a}}({\mathbb F}_{q})$.  Il existe donc un élément $Y\in \boldsymbol{\mathfrak{a}}({\mathbb F}_{q})$ tel que $\ell\circ\beta_{i,F}(Y)=1$ pour tout $i$. Soit $\beta\in \Sigma(A_{H}^{nr})$. Ecrivons $\beta=\sum_{i=1,...,t}n_{i}\beta_{i}$ comme dans (2). Alors $\ell\circ\beta_{F}(Y)=\sum_{i=1,...,t}n_{i}$. La propriété (2) et l'hypothèse $(Hyp)_{1}(p)$ entraînent que $\ell\circ\beta_{F}(Y)$ n'est pas nul (on est ici dans ${\mathbb F}_{q}$). A fortiori $\beta(Y)\not=0$. Ceci étant vrai pour tout $\beta$, $Y$ est un élément régulier de $\boldsymbol{\mathfrak{a}}({\mathbb F}_{q})$. Cela achève la preuve. $\square$

\subsubsection{Description de l'ensemble $Imm^G(H_{ad})$ \label{ImmGHad}}

  \begin{lem}{ Soit $H\in {\cal L}_{F}^{nr}$.  Soit  $X\in  \mathfrak{a}_{H}^{nr}(F)$ un \'el\'ement  entier et de r\'eduction r\'eguli\`ere et soit $Z\in \mathfrak{h}_{tn}(F)$. Alors l'ensemble des points $x\in Imm(G_{AD})$ tels que $X+Z\in \mathfrak{k}_{x}$ est \'egal \`a celui des $x\in Imm^G(H_{ad})$ tels que $Z\in \mathfrak{k}^H_{p_{H}(x)}$. En particulier, 
$Imm^G(H_{ad})$ est l'ensemble des $x\in Imm(G_{AD})$ tels que $X\in \mathfrak{k}_{x}$. }\end{lem}

Preuve. Notons ${\cal X}$ l'ensemble des points $x\in Imm(G_{AD})$ tels que $X+Z\in \mathfrak{k}_{x}$ et ${\cal Y}$ celui des $x\in Imm^G(H_{ad})$ tels que $Z\in \mathfrak{k}^H_{p_{H}(x)}$.  Remarquons que ${\cal Y}$ n'est pas vide. En effet, 
puisque $Z\in \mathfrak{h}_{tn}(F)$, on peut fixer un point $y\in Imm(H_{AD})$ tel que $Z\in \mathfrak{k}^{H,+}_{y}$.  Soit $x\in Imm^G(H_{ad})$ tel que $p_{H}(x)=y$. Alors $x\in {\cal Y}$.

Supposons d'abord que $G$ soit d\'eploy\'e et que  $H$ soit un $F$-Levi.   Soit $x\in  {\cal Y}$.   L'\'el\'ement $X$ \'etant central  dans $\mathfrak{h}$ et entier est contenu dans $\mathfrak{k}^{H}_{z}$ pour tout $z\in Imm(H_{AD})$.  Donc $X\in \mathfrak{k}^H_{p_{H}(x)}$ puis $X+Z\in \mathfrak{k}^H_{p_{H}(x)}$. D'apr\`es \ref{immeublesetLevi}(1), on a $X+Z\in \mathfrak{k}_{x}$, c'est-\`a-dire que $x$ appartient \`a ${\cal X}$.     Inversement, soit $x\in {\cal X}$. On va montrer que $x\in Imm^G(H_{ad})$. Fixons un point $y\in {\cal Y}$. Comme on vient de le voir, on a $y\in {\cal X}$.  Tra\c{c}ons la g\'eod\'esique $[x,y]$ reliant $x$ \`a $y$. Il est bien connu que $\mathfrak{k}_{x}\cap \mathfrak{k}_{y}\subset \mathfrak{k}_{z}$ pour tout $z\in [x,y]$. En particulier $X+Z\in \mathfrak{k}_{z}$. Notons $u$ le point le plus proche de $x$ dans l'ensemble ferm\'e $[x,y]\cap Imm^{G}(H_{ad})$. Si $u=x$, $x\in Imm^G(H_{ad})$ comme on veut le prouver. Supposons $u\not=x$. Notons ${\cal F}$ la facette de $Imm(G_{AD})$ contenant $u$. Il existe une unique facette ${\cal F}'$ telle que tout point de $[x,u[$ assez proche de $u$ appartienne \`a ${\cal F}'$. Puisque $u\in Imm^G(H_{ad})$, on a $X\in \mathfrak{k}_{u}=\mathfrak{k}_{{\cal F}}$ (car, comme ci-dessus, $X$ appartient \`a $\mathfrak{k}^H_{p_{H}(u)}$). On a aussi $X+Z\in \mathfrak{k}_{{\cal F}}$ donc $Z\in \mathfrak{k}_{{\cal F}}$. Notons $X_{{\cal F}}$ et $Z_{{\cal F}}$ les r\'eductions de $X$ et $Z$ dans $\mathfrak{g}_{{\cal F}}({\mathbb F}_{q})$. Utilisons les d\'efinitions de \ref{groupesenreduction} o\`u l'on remplace les indices $u$ par ${\cal F}$. Ainsi, on  note  $H_{{\cal F}}$ le sous-groupe de $G_{{\cal F}}$  qui est l'image naturelle de $K_{{\cal F}}^{0,nr}\cap H(F^{nr})$. C'est un Levi de $G_{{\cal F}}$.  
Montrons que

(1) $X_{{\cal F}}$ est semi-simple et son commutant dans $G_{{\cal F}}$ est $H_{{\cal F}}$. 

On a $Z(H)^0=A_{H}$ puisque $G$ est d\'eploy\'e et que $H$ est un $F$-Levi. Donc $X$ appartient \`a $\mathfrak{a}_{H}(F)\cap \mathfrak{k}_{{\cal F}}$ et la r\'eduction $X_{{\cal F}}$ appartient \`a $\mathfrak{a}_{H,{\cal F}}({\mathbb F}_{q})$.   A fortiori, $X_{{\cal F}}$ est semi-simple.  Compte tenu de \ref{groupesenreduction}(2), il reste \`a prouver que le commutant de $A_{H,{\cal F}}$ dans $G_{{\cal F}}$ est \'egal au commutant de $X_{{\cal F}}$. Puisque nos groupes sont d\'eploy\'es, il suffit de montrer qu'aucune racine  de $A_{H,{\cal F}}$ dans $\mathfrak{g}_{{\cal F}}$  n'est nulle en $X_{{\cal F}}$. Une telle racine $\beta$ provient d'une racine $\beta^G$ de $A_{H}$ dans $\mathfrak{g}$. Puisque $X$ est de r\'eduction r\'eguli\`ere, on a $\beta^G(X)\in \mathfrak{o}_{F}^{\times}$. Le terme $\beta(X_{{\cal F}})$ en est la r\'eduction dans ${\mathbb F}_{q}^{\times}$ et n'est donc pas nulle. Cela prouve (1).

 La facette ${\cal F}$ est contenue dans l'adh\'erence de ${\cal F}'$. Cela entra\^{\i}ne qu'il existe un sous-groupe parabolique $P_{{\cal F}'}$ de $G_{{\cal F}}$ de sorte que la r\'eduction de $\mathfrak{k}_{{\cal F}'}\cap \mathfrak{k}_{{\cal F}}$ dans $\mathfrak{g}_{{\cal F}}({\mathbb F}_{q})$ soit \'egale \`a $\mathfrak{p}_{{\cal F}'}({\mathbb F}_{q})$. Puisque $X+Z\in \mathfrak{k}_{{\cal F}'}\cap \mathfrak{k}_{{\cal F}}$, $X_{{\cal F}}+Z_{{\cal F}}$ appartient \`a $\mathfrak{p}_{{\cal F}'}({\mathbb F}_{q})$. Les \'el\'ements $X_{{\cal F}}$ et $Z_{{\cal F}}$ commutent. Le premier est semi-simple d'apr\`es (1), le second est nilpotent puisque $Z\in \mathfrak{h}_{tn}(F)\subset \mathfrak{g}_{tn}(F)$. Donc $X_{{\cal F}}$ est la composante semi-simple de $X_{{\cal F}}+Z_{{\cal F}}$. Puisque $X_{{\cal F}}+Z_{{\cal F}}$ appartient \`a $\mathfrak{p}_{{\cal F}'}({\mathbb F}_{q})$, sa composante semi-simple $X_{{\cal F}}$ appartient aussi \`a $ \mathfrak{p}_{{\cal F}'}({\mathbb F}_{q})$.  Il existe donc un sous-tore maximal $T_{{\cal F}}$ de 
   $P_{{\cal F}'}$ d\'efini sur ${\mathbb F}_{q}$ tel que $X_{{\cal F}}\in \mathfrak{t}_{{\cal F}}({\mathbb F}_{q})$. On a $T_{{\cal F}}\subset H_{{\cal F}}$ d'apr\`es (1). Donc $H_{{\cal F}}\cap P_{{\cal F}'}$ contient un sous-tore maximal de $G_{{\cal F}}$. Cela entra\^{\i}ne que $H_{{\cal F}}\cap P_{{\cal F}'}$
 est un sous-groupe parabolique de $H_{{\cal F}}$. Etant d\'efini sur ${\mathbb F}_{q}$, il contient un sous-tore d\'eploy\'e maximal $S_{{\cal F}}$.  On rel\`eve $S_{{\cal F}}$ en un sous-tore d\'eploy\'e maximal $S$ de $H$.    Puisque $P_{{\cal F}'}$ contient $S_{{\cal F}}$,  on sait que ${\cal F}'$ est  contenue dans $App(S)$. Mais $App(S)\subset Imm^G(H_{ad})$, donc ${\cal F}'\subset Imm^G(H_{ad})$. Cela contredit la d\'efinition de $u$. On a ainsi prouv\'e que $x\in Imm^G(H_{ad})$. On a $X+Z\in \mathfrak{k}_{x}$ donc $X+Z\in \mathfrak{k}^H_{p_{H}(x)}$ d'apr\`es \ref{immeublesetLevi}(1). On a d\'ej\`a dit que $X\in \mathfrak{k}^H_{p_{H}(x)}$, donc aussi $Z\in \mathfrak{k}^H_{p_{H}(x)}$. Cela prouve que $x\in {\cal Y}$.  
 
 Revenons au cas g\'en\'eral o\`u  $G$ et $H$ ne sont pas suppos\'es d\'eploy\'es. On fixe une extension finie $F'$ de $F$ galoisienne et mod\'er\'ement ramifi\'ee, de sorte que $G$ et $H$ soient d\'eploy\'es sur $F'$.  D'apr\`es le lemme \ref{reductionreguliere}, les hypoth\`eses restent v\'erifi\'ees si l'on remplace le corps de base $F$ par $F'$. On applique ce que l'on vient de d\'emontrer: l'ensemble des points $x\in Imm_{F'}(G_{AD})$ tels que $X+Z\in \mathfrak{k}_{x,F'}$ est \'egal \`a celui des $x\in Imm_{F'}^G(H_{ad})$ tels que $Z\in \mathfrak{k}^H_{p_{H}(x),F'}$.  On prend les invariants par $\Gamma_{F'/F}$. Alors, gr\^ace \`a \ref{immeublesetLevi}(3),  l'ensemble des points $x\in Imm(G_{AD})$ tels que $X+Z\in \mathfrak{k}_{x,F'}$ est \'egal \`a celui des $x\in Imm^G(H_{ad})$ tels que $Z\in \mathfrak{k}^H_{p_{H}(x),F'}$.  Mais, pour $x\in Imm(G_{AD})$, $\mathfrak{k}_{x}=\mathfrak{g}(F)\cap \mathfrak{k}_{x,F'}$ et $\mathfrak{k}^H_{p_{H}(x)}=\mathfrak{h}(F)\cap \mathfrak{k}^H_{p_{H}(x),F'}$. Les conditions $X+Z\in \mathfrak{k}_{x,F'}$ et   $Z\in \mathfrak{k}^H_{p_{H}(x),F'}$ \'equivalent donc \`a $X+Z\in \mathfrak{k}_{x}$ et $Z\in \mathfrak{k}^H_{p_{H}(x)}$. Cela d\'emontre l'\'egalit\'e ${\cal X}={\cal Y}$ qui est  la premi\`ere assertion de l'\'enonc\'e. La seconde est le cas particulier $Z=0$.  $\square$
 
 \subsubsection{ Classes de conjugaison stable dans ${\cal L}^{nr}_{F}$ \label{classesdeconjstable}}
   Soient $H,H'\in {\cal L}^{nr}_{F}$. Disons que $H$ et $H'$ sont stablement conjugu\'es si et seulement s'il existe $g\in G(\bar{F})$ tel que $H'=g^{-1}Hg$ et $g\sigma_{G}(g)^{-1}\in H(\bar{F})$ pour tout $\sigma\in \Gamma_{F}$. Fixons $H$ et un \'el\'ement $X_{H}\in \mathfrak{a}_{H}^{nr}(F)$ entier et de r\'eduction r\'eguli\`ere.  Alors $H=Z_{G}(X_{H})^0$ d'apr\`es le (i) du lemme \ref{reductionreguliere}. On d\'efinit la classe de conjugaison stable de $X_{H}$ de la fa\c{c}on habituelle: $X\in \mathfrak{g}(F)$ est stablement conjugu\'e \`a $X_{H}$ si et seulement s'il existe $g\in G$ tel que $X=g^{-1}X_{H}g$ et $g\sigma_{G}(g)^{-1}\in H$ pour tout $\sigma\in \Gamma_{F}$. Notons $ker^1(H,G)$ le noyau de l'application naturelle $H^1(F,H)\to H^1(F,G)$. 

{\bf Remarque.} Les ensembles $H^1(F,H)$ et $H^1(F,G)$ ne sont pas naturellement des groupes. Le noyau est pris au sens des ensembles point\'es. Toutefois, on peut bel et bien munir ces ensembles de structures de groupes, par exemple en les identifiant aux duaux de Pontryagin de $\pi_{0}(Z(\hat{H})^{\Gamma_{F}})$ et $\pi_{0}(Z(\hat{G})^{\Gamma_{F}})$. L'application $H^1(F,H)\to H^1(F,G)$ est alors un homomorphisme de groupes et $ker^1(H,G)$ est un sous-groupe de $H^1(F,H)$. 
\bigskip

Notons $Cl_{st}(X_{H})$ la classe de conjugaison stable de $X_{H}$ et $Cl_{st}(X_{H})/conj$ l'ensemble des  classes de conjugaison par $G(F)$ dans $Cl_{st}(X_{H})$. Soit $X\in Cl_{st}(X_{H})$. On choisit $g$ comme ci-dessus. L'application $\sigma\mapsto g\sigma_{G}(g)^{-1}$ est un cocycle \`a valeurs dans $H$ dont la classe ne d\'epend pas du choix de $g$. L'\'el\'ement de $H^1(F,H)$ qu'il d\'efinit appartient \'evidemment \`a $ker^1(H,G)$.   On obtient une application de $Cl_{st}(X_{H})$ dans $ ker^1(H,G)$. Elle se quotiente en une bijection de l'ensemble $Cl_{st}(X_{H})/conj$  sur $ker^1(H,G)$. Notons $Cl_{st}(H)$ la classe de conjugaison stable de $H$ et $Cl_{st}(H)/conj$ l'ensemble des classes de conjugaison par $G(F)$ dans $Cl_{st}(H)$. 

  \begin{lem}{(i) L'application $X\mapsto G_{X}$ est une application surjective de $Cl_{st}(X_{H})$ sur $Cl_{st}(H)$ qui se quotiente en une surjection de $Cl_{st}(X_{H})/conj$ sur $Cl_{st}(H)/conj$.
  
  (ii) Pour $X\in Cl_{st}(X_{H})$, $X$ est un \'el\'ement entier de r\'eduction r\'eguli\`ere de $\mathfrak{a}^{nr}_{G_{X}}(F)$. 
  
  (iii) Pour $X\in Cl_{st}(X_{H})$, resp. $H'\in Cl_{st}(H)$, il existe un \'el\'ement $g\in G(F^{nr})$ tel que $X=g^{-1}X_{H}g$, resp. $H'=g^{-1}Hg$, et $gFr(g)^{-1}\in H(F^{nr})$. Pour un tel $g$, l'application $Ad(g)^{-1}$ se restreint en un isomorphisme d\'efini sur $F$ de $Z(H)$ sur $Z(G_{X})$, resp. sur $Z(H')$. 
  
  (iv) Le groupe $ker^1(H,G)$ est l'image naturelle de $H^1(F,H_{sc})$ dans $H^1(F,H)$. }\end{lem}
  
  Preuve. La d\'emonstration de (iv) est standard. Le diagramme
  $$\begin{array}{ccc}H^1(F,H_{sc})&\to&H^1(F,G_{SC})\\ \downarrow&&\downarrow\\ H^1(F,H)&\to& H^1(F,G)\\ \end{array}$$
  est commutatif. Puisque $G_{SC}$ est simplement connexe, on a $H^1(F,G_{SC})=\{0\}$. Donc $H^1(F,H_{sc})$ s'envoie dans $ker^1(H,G)$. Inversement, soit $u$ un cocycle repr\'esentant un \'el\'ement de $ ker^1(H,G)$. Par hypoth\`ese, il existe $g\in G$ tel que $u(\sigma)=g\sigma_{G}(g)^{-1}$ pour tout $\sigma$. Notons $\pi:G_{SC}\to G$ l'homomorphisme naturel. On peut \'ecrire $g=z\pi(g_{sc})$ avec $z\in Z(G)$ et $g_{sc}\in G_{SC}$. 
  Puisque $Z(G)\subset H$, on peut remplacer $u$ par le cocycle $\sigma\mapsto u'(\sigma)=z^{-1}u(\sigma)\sigma_{G}(z)$. D\'efinissons le cocycle $u'_{sc}$ \`a valeurs dans $G_{SC}$ par $u'_{sc}(\sigma)=g_{sc}\sigma_{G}(g_{sc})^{-1}$. Alors $u'(\sigma)=\pi(u'_{sc}(\sigma))$. Cela entra\^{\i}ne que $u'_{sc}$ prend ses valeurs dans $H_{sc}$ et $u'$ d\'efinit un \'el\'ement de $H^1(F,H_{sc})$. L'image de $u$ dans $H^1(F,H)$ est l'image de cet \'el\'ement.

  D\'emontrons  (iii). Soit $X\in Cl_{st}(X_{H})$, resp. $H'\in Cl_{st}(H)$. Fixons un \'el\'ement $g\in G$ tel que $X=g^{-1}X_{H}g$, resp. $H'=g^{-1}Hg$, et $g\sigma_{G}(g)^{-1}\in H$ pour tout $\sigma\in \Gamma_{F}$. Alors le cocycle $\sigma\mapsto g\sigma_{G}(g)^{-1}$ d\'efinit un \'el\'ement de $ker^1(H,G)$. On sait qu'il existe une extension finie $F'$ de $F$, non ramifi\'ee, telle que l'application de restriction $H^1(F,H)\to H^1(F',H)$ soit nulle. Fixons une telle extension. Le cocycle pr\'ec\'edent devient trivial dans $H^1(F',H)$ donc on peut fixer $h\in H$ tel que  $hg\sigma_{G}(hg)^{-1}=1$ pour tout $\sigma\in \Gamma_{F'}$ ou encore $hg\in G(F')$. Alors l'\'el\'ement $hg$ v\'erifie la condition de (iii). Pour unifier la notation, on pose $H'=G_{X}$ quand la donn\'ee de d\'epart est un \'el\'ement  $X\in Cl_{st}(X_{H})$.  L'application $Ad(g)^{-1}$ se restreint en un isomorphisme $\varphi:Z(H)\to Z(H')$ qui est \'evidemment d\'efini sur $F^{nr}$.   Pour $z\in Z(H)$, on a $Fr_{G}(\varphi(z))=\varphi\circ Ad(gFr_{G}(g)^{-1})(Fr_{G}(z))$. Mais $gFr_{G}(g)^{-1}\in H(F^{nr})$ donc $Ad(gFr_{G}(g)^{-1})$ est l'identit\'e sur $Z(H)$ et on obtient $Fr_{G}(\varphi(z))=\varphi(Fr_{G}(z))$. Autrement dit $\varphi$ est d\'efini sur $F$. Cela prouve (iii). 
  
    L'assertion (ii) est \'evidente, puisque les couples $(H,X_{H})$ et $(G_{X},X)$ sont conjugu\'es par un \'el\'ement de $G(F^{nr})$ d'apr\`es (iii).

  Soit $X\in Cl_{st}(X_{H})$. Appliquons (iii) et fixons $g\in G(F^{nr})$ tel que $X=g^{-1}X_{H}g$  et $gFr_{G}(g)^{-1}\in H(F^{nr})$.
Posons $H'=G_{X}=g^{-1}Hg$. Puisque $H$ est un $F^{nr}$-Levi de $G$ et que $g\in G(F^{nr})$, $H'$ est encore un $F^{nr}$-Levi de $G$. Il est stablement conjugu\'e \`a $H$.  L'application $X\mapsto G_{X}$ est  donc une application  de $Cl_{st}(X_{H})$ dans $Cl_{st}(H)$. Montrons qu'elle est surjective. Soit $H'\in Cl_{st}(H)$. On fixe $g\in G(F^{nr})$ tel que $H'=g^{-1}Hg$  et $gFr_{G}(g)^{-1}\in H(F^{nr})$. On pose $X=g^{-1}X_{H}g$.   Il est clair que $X$ est stablement conjugu\'e \`a $X_{H}$ et que $G_{X}=H'$. Cela d\'emontre (i). $\square$

 L'application de $Cl_{st}(X_{H})/conj$ sur $Cl_{st}(H)/conj$ n'est pas injective en g\'en\'eral car un \'el\'ement de $(H\backslash Norm_{G}(H))^{\Gamma_{F}}$ n'est pas forc\'ement repr\'esent\'e par un \'el\'ement de $Norm_{G(F)}(H)$. Pour un tel \'el\'ement que l'on rel\`eve en un \'el\'ement $n\in  Norm_{G}(H)$, l'\'el\'ement $X=n^{-1}X_{H}n$ est stablement conjugu\'e \`a $X_{H}$ mais pas conjugu\'e par un \'el\'ement de $G(F)$. Pourtant $G_{X}=H$. 
 
Notons ${\cal L}^{nr}_{ell}$ le sous-ensemble des $H\in {\cal L}_{F}^{nr}$ tels que $A_{H}=A_{G}$. On a

(1) pour $H\in {\cal L}^{nr}_{ell}$, $Cl_{st}(H)$ est contenu dans ${\cal L}^{nr}_{ell}$.

En effet,  d'apr\`es le (iii) du lemme ci-dessus, les groupes $Z(H')^0$ et $Z(H)^0$ sont isomorphes sur $F$. L'\'egalit\'e $A_{G}=A_{H}$ entra\^{\i}ne donc $A_{G}=A_{H'}$.

 \subsubsection{Passage au groupe $G^*$ pour l'ensemble ${\cal L}_{F}^{nr}$\label{passage}}
 On reprend  pour  la fin de la sous-section 2.4  les notations de \ref{formeinterieure}. En particulier, on  fixe  un épinglage $\mathfrak{E}^*=(B^*,T^*,(E_{\alpha})_{\alpha\in \Delta})$ de $G^*$ défini sur $F$. 
   Soit $H\in {\cal L}_{F}^{nr}$. Le groupe $H$ est un $F^{nr}$-Levi de $G$   donc $\psi(H)$ est un $F^{nr}$-Levi de $G^*$. Il est conjugu\'e par un \'el\'ement de $G^*(F^{nr})$ \`a un $F^{nr}$-Levi standard $M^*$ de $G^*$. On peut fixer $x\in G^*(F^{nr})$ de sorte que $Ad(x)\circ \psi(H)=M^*$. On pose $\psi_{x}=Ad(x)\circ \psi$, dont le cocycle associ\'e $u_{x}$ est d\'efini par $u_{x}(\sigma)=xu_{G}(\sigma)\sigma_{G^*}(x)^{-1}$ pour tout $\sigma\in \Gamma_{F}$.  On  a
$$\sigma_{G^*}(M^*)= \sigma_{G^*}(\psi_{x}(H)) = u_{x}(\sigma)^{-1}\psi_{x}(\sigma_{G}(H))u_{x}(\sigma) =u_{x}(\sigma)^{-1}\psi_{x}(H)u_{x}(\sigma) ,$$
puisque $H$ est d\'efini sur $F$. Donc
$$(1) \qquad \sigma_{G^*}(M^*)=u_{x}(\sigma)^{-1}M^*u_{x}(\sigma).$$

Soit $s_{H}$ un sommet de $Imm(H_{AD})$, supposons  que ${\bf FC}^{H_{s_{H}}}(\mathfrak{h}_{SC,s_{H}})\not=\emptyset$. 
 Puisque $\psi_{x}:H\to M^*$  est un  isomorphisme d\'efini sur $F^{nr}$,  il s'en d\'eduit un isomorphisme not\'e $\psi^{Imm}_{x}:Imm_{F^{nr}}(H_{AD})\to Imm_{F^{nr}}(M^*_{AD})$. Notons $s_{M^*}$ l'image de $s_{H}$. On a $M^*_{SC,s_{M^*}}\simeq H_{SC,s_{H}}$ et l'hypoth\`ese implique que ${\bf FC}^{M^*_{s_{M^*}}}(\mathfrak{m}^*_{SC,s_{M^*}})\not=\emptyset$. 
 Donc  $M^*$ v\'erifie les conclusions de la proposition \ref{Levistandard}. Montrons que

(2) $M^*$ est d\'efini sur $F$ et est uniquement d\'etermin\'e.

Soit $\sigma\in \Gamma_{F}$. Posons $g=u_{x}(\sigma)^{-1}$. C'est un \'el\'ement de $ G_{AD}^*(F^{nr})$ tel que $g^{-1}\sigma_{G^*}(M^*)g=M^*$ d'apr\`es (1). Notons $P^*$ le sous-groupe parabolique standard de $G^*$ de composante de Levi $M^*$. Alors $g^{-1}\sigma_{G^*}(P^*)g$ est un sous-groupe parabolique d\'efini sur $F^{nr}$ et de composante de Levi $M^*$ d'apr\`es (1). D'apr\`es le (ii) de la proposition \ref{Levistandard}, il est conjugu\'e \`a $P^*$ par un \'el\'ement de $W^{I_{F}}$. Tout \'el\'ement de $W^{I_{F}}$ est repr\'esent\'e par un \'el\'ement de $G^*(F^{nr})$. Il en r\'esulte que $\sigma_{G^*}(P^*)$ est conjugu\'e \`a $P^*$ par un \'el\'ement de $G^*_{AD}(F^{nr})$. Mais $P^*$ et $\sigma_{G^*}(P^*)$ sont standard. Leur conjugaison entra\^{\i}ne qu'ils sont \'egaux. Donc $\sigma_{G^*}(P^*)=P^*$ puis $\sigma_{G^*}(M^*)=M^*$. Donc $M^*$ est d\'efini sur $F$. La d\'efinition de $M^*$ ne d\'epend que du choix de l'\'el\'ement $x$. Choisissons un autre \'el\'ement $\underline{x}\in G^*(F^{nr})$ de sorte que $Ad(\underline{x})\circ \psi(H)$ soit encore un Levi standard, que l'on note $\underline{M}^*$. Posons $y=\underline{x}x^{-1}$. Alors $y\in G^*(F^{nr})$ et $yM^*y^{-1}=\underline{M}^*$. Les paires $(B^*\cap M^*, T^*)$ et $(y^{-1}(B^*\cap \underline{M}^*)y, y^{-1}T^*y)$ sont deux paires de Borel de $M^*$ d\'efinies sur $F^{nr}$. Quitte \`a multiplier $y$ \`a droite  par un \'el\'ement de $M^*(F^{nr})$, on peut supposer qu'elles sont \'egales. Alors $y$ normalise $T^*$ et d\'efinit un \'el\'ement  $w\in W^{I_{F}}$ tel que $w(M^*)=\underline{M}^*$. Le (i) de la proposition \ref{Levistandard}  implique que $M^*=\underline{M}^*$. Cela prouve que $M^*$  est uniquement d\'efini, d'o\`u (2).

\subsubsection{Passage \`a $G^*$, le cas stable\label{passagestable}}
Soit $H\in {\cal L}_{F}^{nr}$. On lui associe comme en \ref{passage} un Levi standard $M^*$ de $G^*$ d\'efini sur $F^{nr}$.   Supposons que $M^*$ soit d\'efini sur $F$.   Fixons des alc\^oves $C^{nr,H}$, resp. $C^{nr,M^*}$, de $Imm_{F^{nr}}(H_{AD})$, resp.  de $Imm_{F^{nr}}(M_{AD})$, conserv\'ees par $\Gamma_{F}^{nr}$.   L'isomorphisme $\psi_{x}:H\to M^*$ est d\'efini 
  gr\^ace au choix d'un \'el\'ement $x\in G^*(F^{nr})$ que l'on peut, si l'on veut, multiplier \`a gauche par un \'el\'ement de $M^*(F^{nr})$. On peut donc supposer que l'isomorphisme d'immeubles  $\psi^{Imm}_{x}$ envoie $C^{nr,H}$ sur $C^{nr,M^*}$.  Fixons un sommet $s_{H}\in  S(\bar{C}^{nr,H})$ et posons $s_{M^*}=\psi_{x}^{Imm}(s_{H})$.  On a $s_{M^*}\in S( \bar{C}^{nr,M^*})$.

 \begin{lem}{(i) Les propri\'et\'es suivantes sont \'equivalentes:
 
 (a) $s_{H}\in Imm(H_{AD})$ et  $FC^{st}(\mathfrak{h}_{SC,s_{H}}({\mathbb F}_{q}))\not=\{0\}$;
 
 (b) $M^*$ est d\'efini sur $F$, $s_{M^*}\in Imm(M^*_{AD})$ et $FC^{st}(\mathfrak{m}^*_{SC,s_{M^*}}({\mathbb F}_{q}))\not=\{0\}$.
 
 Si elles sont v\'erifi\'ees, les espaces  $FC^{st}(\mathfrak{h}_{SC,s_{H}}({\mathbb F}_{q}))$ et $FC^{st}(\mathfrak{m}^*_{SC,s_{M^*}}({\mathbb F}_{q}))$ sont isomorphes. 
 
 (ii) Les propri\'et\'es suivantes sont \'equivalentes:
 
 (a) $FC^{st}(\mathfrak{h}_{SC}(F))\not=\{0\}$;
 
 (b) $M^*$ est d\'efini sur $F$ et $FC^{st}(\mathfrak{m}^*_{SC}(F))\not=\{0\}$.
 
 Si elles sont v\'erifi\'ees, les espaces $FC^{st}(\mathfrak{h}_{SC}(F))$ et $FC^{st}(\mathfrak{m}^*_{SC}(F))$ sont isomorphes.
 
  (iii) Supposons $G$  absolument quasi-simple. Si l'une des propri\'et\'es ci-dessus est v\'erifi\'ee,  $M^*_{SC}$ est lui aussi absolument quasi-simple. }\end{lem}
 
 {\bf Remarque.} On peut supposer que $M^*$ est d\'efini sur $F$: c'est une hypoth\`ese si (i)(b) ou (ii)(b) est v\'erifi\'ee; cela r\'esulte de  \ref{passage} (2) si (i)(a) ou (ii)(a) l'est.  
 \bigskip
 
  Preuve. En vertu du lemme \ref{decompositionFCstable}, l'assertion (ii) r\'esulte de (i) appliqu\'ee \`a tous les sommets $s_{H}\in S(\bar{C}^{nr,H})$. L'assertion (iii) r\'esulte de \ref{couplesstables}(1). 
  
  D\'emontrons (i). Si $G^*$ est simplement connexe et absolument quasi-simple, l'assertion resulte  de \cite{W7} 9(3). Remarquons que l'on a de plus:  $s_{H}$, resp. $s_{M^*}$ est fix\'e par $Aut({\cal D}_{a}^{nr,H})$, resp. $Aut({\cal D}_{a}^{nr,M^*})$, cf. \cite{W7} 9(4). On va se ramener \`a ce cas.

  On ne perd rien \`a supposer que $G$ et $G^*$ sont simplement connexes. On d\'ecompose $G^*$ en produit de groupes $Res_{F'/F}(G^{'*})$, o\`u $F'$ est une extension finie  mod\'er\'ement ramifi\'ee  de $F$ et $G^{'*}$ est un groupe quasi-d\'eploy\'e sur $F'$. Alors $G$ se d\'ecompose conform\'ement en produit de groupes $Res_{F'/F}(G')$, o\`u  $G'$ est une forme int\'erieure de $G^{'*}$. Le groupe $H$, resp. $M^*$, se d\'ecompose lui-aussi en produit  de groupes $Res_{F'/F}(H')$, resp. $Res_{F'/F}(M^{'*})$. On voit qu'il suffit de d\'emontrer l'\'enonc\'e pour chacun des facteurs. On peut  donc supposer   $G^*=Res_{F'/F}(G^{'*})$. 
  
  Si $F'=F$,  $G^*$ est absolument quasi-simple et on a vu ci-dessus que le lemme \'etait v\'erifi\'e. 
  
  Supposons $F'/F$ totalement ramifi\'e. On a alors  $Imm_{F^{nr}}(H_{AD})=Imm_{F^{'nr}}(H'_{AD})$ et $Imm_{F^{nr}}(M^*_{AD})=Imm_{F^{'nr}}(M^{'*}_{AD})$ et ces identifications sont \'equivariantes pour les actions de $\Gamma_{F}^{nr}=\Gamma_{F'}^{nr}$. L'\'enonc\'e se d\'eduit alors du m\^eme \'enonc\'e (d\'ej\`a d\'emontr\'e) pour les groupes $H'$ et $M^{'*}$. 
  
  En g\'en\'eral, soit $E$ la plus grande extension non ramifi\'ee contenue dans $F'$. Posons $d=[E:F]$.  Posons $G^{''*}=Res_{F'/E}(G^{'*})$, $G''=Res_{F'/E}(G')$, $H''=Res_{F'/E}(H')$, $M^{''*}=Res_{F'/E}(M^{'*})$. Comme en \ref{Weil}, on  identifie $Imm_{F^{nr}}(H_{AD})$ \`a $Imm_{F^{nr}}(H''_{AD})^d$, l'action de $Fr_{H}$ \'etant
  $$(x_{1},...,x_{d})\mapsto (Fr^d_{H''}(x_{d}),x_{1},...,x_{d-1}).$$
  La chambre $C^{nr,H}$ s'identifie \`a $(C^{nr,H''})^d$ et le sommet $s_{H}$ \`a un $d$-uple de sommets $(s_{H,1},...,s_{H,d})$. La condition $s_{H}\in Imm(H_{AD})$ \'equivaut \`a $s_{H,1}=...=s_{H,d}$ et $s_{H,1}\in Imm_{E}(H''_{AD})$. Si ces conditions sont v\'erifi\'ees, on a $\mathfrak{h}_{SC,s_{H}}=Res_{{\mathbb F}_{q^d}/{\mathbb F}_{q}}(\mathfrak{h}''_{SC,s_{H,1}})$. En particulier $\mathfrak{h}_{SC,s_{H}}({\mathbb F}_{q})=\mathfrak{h}''_{SC,s_{H,1}}({\mathbb F}_{q^d})$. De m\^emes propri\'et\'es valent pour le groupe $M^*$. L'isomorphisme $\psi_{x}^{Imm}$ se d\'ecompose en produit $\prod_{i=1,...,d}\psi_{i}^{Imm}$, o\`u $\psi_{i}^{Imm}$ est un isomorphisme de $ Imm_{E}(H''_{AD})$ sur $ Imm_{E}(M^{''*}_{AD})$ qui envoie $ C^{nr,H''}$ sur $C^{nr,M^{''*}}$. On a $\psi_{i}^{Imm}(s_{H,i})=s_{M^*,i}$ pour tout $i$. L'action des Frobenius n'\'etablit pas de relation entre les $\psi_{i}^{Imm}$ car $\psi_{x}^{Imm}$ n'est pas d\'efini sur $F$. Toutefois, notons $\psi_{i,  {\cal D}_{a}^{nr,H''}}$ la bijection de ${\cal D}_{a}^{nr,H''}$ sur ${\cal D}_{a}^{nr,M^{''*}}$ d\'eduite de $\psi_{i}^{Imm}$. Il existe alors $\delta_{i}\in Aut({\cal D}_{a}^{nr,M^{''*}})$ tel que    $\psi_{i, {\cal D}_{a}^{nr,H''}}=\delta_{i}\circ \psi_{1, {\cal D}_{a}^{nr,H''}}$. Supposons (i)(a) v\'erifi\'ee. Alors, pour tout $i$, on a $s_{H,i}=s_{H,1}\in Imm_{E}(H''_{AD})$ et $FC^{st}(\mathfrak{h}''_{SC,s_{H,i}}({\mathbb F}_{q^d}))\not=\{0\}$. Appliquons l'\'enonc\'e d\'ej\`a d\'emontr\'e pour $H''$ et $M^{''*}$ sur $E$. On en d\'eduit que, pour tout $i$,  $s_{M^*,i}\in Imm_{E}(M^{''*}_{AD})$ et $FC^{st}(\mathfrak{m}^{''*}_{SC,s_{M^*,i}}({\mathbb F}_{q^d}))\not=\{0\}$. En identifiant $s_{M^*,i}$, resp $s_{H,i}$, \`a un sommet du diagramme ${\cal D}_{a}^{nr,M^{''*}}$, resp. ${\cal D}_{a}^{nr,H''}$, on a aussi $s_{M^*,i}=\psi_{i, {\cal D}_{a}^{nr,H''}}(s_{H,i})=\delta_{i}\circ \psi_{1, {\cal D}_{a}^{nr,H''}}(s_{H,1})=\delta_{i}(s_{M^*,1})$. Mais, comme on l'a rappel\'e ci-dessus, les sommets qui peuvent intervenir sont fixes par l'action de $Aut({\cal D}_{a}^{nr,M^{''*}})$. Donc $s_{M^*,i}=s_{M^*,1}$. Jointe \`a la relation $s_{M^*,1}\in Imm_{E}(M^{''*}_{AD})$, cette \'egalit\'e entra\^{\i}ne que $s_{M^*}\in Imm(M^*_{AD})$. On a alors $\mathfrak{m}^*_{SC,s_{M^*}}({\mathbb F}_{q})=\mathfrak{m}^{''*} _{SC,s_{M^*,1}}({\mathbb F}_{q^d})$ et la relation $FC^{st}(\mathfrak{m}^{''*}_{SC,s_{M^*,1}}({\mathbb F}_{q^d}))\not=\{0\}$ entra\^{\i}ne $FC^{st}(\mathfrak{m}^*_{SC,s_{M^*}}({\mathbb F}_{q}))\not=\{0\}$. Cela d\'emontre (i)(b). De plus, $FC^{st}(\mathfrak{m}^*_{SC,s_{M^*}}({\mathbb F}_{q}))=FC^{st}(\mathfrak{m}^{''*}_{SC,s_{M^*,1}}({\mathbb F}_{q^d}))\simeq FC^{st}(\mathfrak{h}''_{SC,s_{H,1}}({\mathbb F}_{q^d}))=FC^{st}(\mathfrak{h}_{SC,s_{H}}({\mathbb F}_{q}))$. 
  On d\'emontre de la m\^eme fa\c{c}on que (i)(b) implique (i)(a). $\square$
  
  \subsubsection{L'ensemble ${\cal L}^{*,st}_{F,sd}$\label{L*stFsd}}
  
  On note ${\cal L}_{F,sd}^{*,st}$ l'ensemble des $F$-Levi standard $M^*$ de $G^*$  tels  que $FC^{st}(\mathfrak{m}^*_{SC}(F))\not=\{0\}$ (l'indice $sd$ de la notation indique "standard"). Remarquons que tout \'el\'ement de ${\cal L}_{F,sd}^{*,st}$ v\'erifie la proposition \ref{Levistandard}. Pour $M^*\in {\cal L}_{F,sd}^{*,st}$, notons $W^{I_{F}}(M^*)=W^{M^*,I_{F}}\backslash Norm_{W^{I_{F}}}(M^*)$. Le groupe $W^{I_{F}}(M^*)$ agit sur lui-m\^eme par $Fr$-conjugaison: l'action d'un \'el\'ement $w$ est $w'\mapsto Fr(w)w'w^{-1}$. Pour $w\in W^{I_{F}}(M^*)$, on note $Cl_{Fr}(w)$ la classe de $Fr$-conjugaison de $w$. 

Soit $w\in W^{I_{F}}(M^*)$. Cet  \'el\'ement $w$ agit sur $X_{*}(Z(M^*)^0)^{I_{F}}$.  Notons $X_{*}(Z(M^*)^0)^{I_{F},w^{-1}\circ Fr}$ le sous-espace des points fixes de $w^{-1}\circ Fr$ dans $X_{*}(Z(M^*)^0)^{I_{F}}$. Notons $A_{w}$ le sous-tore de $Z(M^*)^0$ dont le groupe de cocaract\`eres est  $ X_{*}(Z(M^*)^0)^{I_{F},w^{-1}\circ Fr}$. Notons $M^*_{w}$ le commutant de $A_{w}$ dans $G^*$. Le groupe $A_{w}$ est un tore d\'eploy\'e sur $F^{nr}$ et  $M^*_{w}$ est un $F^{nr}$-Levi de $G^*$  contenant $M^*$. Disons que $w$ est $G$-admissible, resp. $Fr$-r\'egulier, si et seulement si $\psi^{-1}(M^*_{w})$ est conjugu\'e \`a un $F$-Levi de $G$, resp. $M^*_{w}=G$. Un \'el\'ement $Fr$-r\'egulier est \'evidemment $G$-admissible. On note $W^{I_{F}}(M^*)_{G-adm}$, resp. $W^{I_{F}}(M^*)_{Fr-reg}$, l'ensemble des \'el\'ements $G$-admissibles, resp. $Fr$-r\'eguliers. 

\begin{lem}{Soit $w\in W^{I_{F}}(M^*)$. 

(i) Soit $w'\in Cl_{Fr}(w)$. Alors $w$ est $G$-admissible, resp. $Fr$-r\'egulier, si et seulement s'il en est de m\^eme de $w'$.

(ii) L'\'el\'ement $w$ est $G$-admissible si et seulement s'il existe un $F$-Levi standard $L^*$ de $G^*$ tel que $M^*\subset L^*$, $M^*_{min}\subset L^*$ et   $ W^{L^*,I_{F}}(M^*)_{Fr-reg}\cap Cl_{Fr}(w)\not=\emptyset$. }\end{lem}

Preuve.   Soit $u\in W^{I_{F}}(M^*)$, posons $w'=Fr(u)wu^{-1}$. On v\'erifie que le groupe $X_{*}(Z(M^*)^0)^{I_{F},w^{'-1}\circ Fr}$
est l'image naturelle  de $X_{*}(Z(M^*)^0)^{I_{F},w^{-1}\circ Fr}$  par l'action de $u$. Donc $M^*_{w'}=uM^*_{w}u^{-1}$ (ici et dans la suite, on s'autorise \`a  identifier $u$ \`a un rel\`evement dans $G^*(F^{nr})$ qui conserve $T^*$ et $M^*$). L'assertion (i) en r\'esulte imm\'ediatement. 

Gr\^ace \`a (i), les deux propri\'et\'es de $w$ dont (ii) affirme l'\'equivalence sont invariantes par $Fr$-conjugaison. On s'autorise donc \`a remplacer $w$ par un \'el\'ement $Fr$-conjugu\'e. 
Fixons un \'el\'ement $X\in X_{*}(Z(M^*)^0)^{I_{F},w^{-1}\circ Fr}$  qui est r\'egulier en ce sens que  $\alpha(X)\not=0$ pour toute racine $\alpha\in \Sigma(A_{w})$. Notons $P_{w}^*$ le sous-groupe parabolique de $G^*$ engendr\'e par $M_{w}^*$ et les sous-espaces radiciels associ\'es aux racines $\alpha\in \Sigma(A_{w})$ telles que $\alpha(X)>0$. La paire parabolique $(P^*_{w},M^*_{w})$ est d\'efinie sur $F^{nr}$ et $M^*_{w}$ contient  $T^*$. Il existe donc une paire parabolique standard $(Q^*,L^*)$ de $G^*$ et un \'el\'ement 
$u\in W^{I_{F}}$ tels que  $(uP^*_{w}u^{-1},uM_{w}^*u^{-1})=(Q^*,L^*)$. Le groupe $uM^*u^{-1}$ est un  $F^{nr}$-Levi de $ L^*$ qui contient $T^*$. Quitte \`a multiplier $u$ \`a gauche par un \'el\'ement de $W^{L^*,I_{F}}$, on peut donc supposer que $uM^*u^{-1}$ est standard. D'apr\`es la proposition \ref{Levistandard}, on a alors $uM^*u^{-1}=M^*$ et $u$ d\'efinit un \'el\'ement de $W^{I_{F}}(M^*)$. 
Posons $w'=Fr(u)wu^{-1}$ et recommen\c{c}ons la construction en rempla\c{c}ant $w$ par $w'$ et $X$ par $X'=u(X)\in X_{*}(Z(M^*)^0)^{I_{F},w^{'-1}\circ Fr}$. On voit que la paire $(P_{w}^*,M_{w}^*)$ est remplac\'ee par $(uP^*_{w}u^{-1},uM_{w}^*u^{-1})$, c'est-\`a-dire $(Q^*,L^*)$. Quitte \`a effectuer ces remplacements, on peut donc aussi bien supposer que $(P_{w}^*,M_{w}^*)=(Q^*,L^*)$. Relevons $Fr$ en un \'el\'ement de $\Gamma_{F}$. Par construction, la paire $(P_{w}^*,M_{w}^*)=(Q^*,L^*) $ est conserv\'ee par $w^{-1}\circ Fr$. Puisqu'elle est standard, un raisonnement habituel entra\^{\i}ne qu'elle est conserv\'ee par $Fr$ et par $w$. Cela entra\^{\i}ne que $(Q^*,L^*)$ est d\'efinie sur $F$, autrement dit que $L^*$ est un $F$-Levi, et que $w\in W^{L^*,I_{F}}(M^*)$. L'\'egalit\'e $M_{w}^*=L^*$ et la d\'efinition de $M^*_{w}$ entra\^{\i}nent alors que $w\in W^{L^*,I_{F}}(M^*)_{Fr-reg}$. Supposons que $w$ soit $G$-admissible. Alors $M_{w}^*=L^*$ se transf\`ere en un $F$-Levi de $G$. D'apr\`es \ref{formeinterieure}(2), cela implique $M^*_{min}\subset L^*$. La conclusion de (ii) est donc v\'erifi\'ee.

Inversement, supposons que $w$ v\'erifie la conclusion de (ii).  En rempla\c{c}ant $w$ par un \'el\'ement $Fr$-conjugu\'e. on peut  supposer que $w\in  W^{L^*,I_{F}}(M^*)_{Fr-reg}$. Alors $M_{w}^*=L^*$. Puisque $M^*_{min}\subset L^*$, $\psi^{-1}(L^*)$ est un Levi standard de $G$ (pour la paire $(P_{min},M_{min})$, cf. \ref{formeinterieure}) d\'efini sur  $F$, donc un $F$-Levi de $G$. Cela prouve que $w$ est $G$-admissible. $\square$

 On note $W^{I_{F}}(M^*)_{G-adm}/Fr-conj$, resp.  $W^{I_{F}}(M^*)_{Fr-reg}/ Fr-conj$, l'ensemble des classes de $Fr$-conjugaison dans $W^{I_{F}}(M^*)_{G-adm}$, resp. $W^{I_{F}}(M^*)_{Fr-reg}$.

 \subsubsection{Les ensembles  ${\cal L}_{F}^{nr,st}$ et ${\cal L}^{nr,st}_{ell}$\label{parametrage}}
 
On a d\'efini les ensembles ${\cal L}_{F}^{nr}$ et ${\cal L}_{ell}^{nr}$ en \ref{classesdeconjstable}.  On note   ${\cal L}_{F}^{nr,st}$ l'ensemble des $H\in {\cal L}_{F}^{nr}$ tels que $FC^{st}(\mathfrak{h}_{SC}(F))\not=\{0\}$. On pose 
${\cal L}^{nr,st}_{ell}={\cal L}_{ell}^{nr}\cap {\cal L}_{F}^{nr,st}$.

{\bf Remarque.} Les $F^{nr}$-Levi qui sont des tores sont les éléments de ${\cal T}_{max}^{nr}$. Tout $H\in {\cal T}_{max}^{nr}$ qui est défini sur $F$ appartient à ${\cal L}_{F}^{nr,st}$: on a $\mathfrak{h}_{SC}=\{0\}$ et $FC^{st}(\mathfrak{h}_{SC}(F))={\mathbb C}$. 

\bigskip

D'apr\`es le lemme \ref{extensionsnonram}, on a ${\cal L}_{F}^{nr,st}\subset {\cal L}_{F'}^{nr,st}$ pour toute extension finie non ramifi\'ee $F'$ de $F$. On note ${\cal L}^{nr}$, resp. ${\cal L}^{nr,st}$,   la r\'eunion des ${\cal L}_{F'}^{nr}$, resp. ${\cal L}_{F'}^{nr,st}$,  quand $F'$ parcourt les extensions finies non ramifi\'ees de $F$. Ces deux ensembles sont conserv\'es par conjugaison par $F^{nr}$.   L'ensemble ${\cal L}^{nr}$ est celui des $F^{nr}$-Levi de $G$. Remarquons que ${\cal L}^{nr,st}$ ne change pas par une extension finie non ramifi\'ee du corps de base et que l'application $H\mapsto H_{sc}$ est une bijection de ${\cal L}^{nr,st}$ sur ${\cal L}^{nr,st,G_{SC}}$.

 Rappelons la propri\'et\'e g\'en\'erale suivante. Soient $G_{1}$ et  $G_{2}$ deux groupes r\'eductifs connexes d\'efinis sur $F$. Deux torseurs int\'erieurs $\Psi,\Psi':G_{1}\to G_{2}$ sont dits dans la m\^eme classe modulo automorphismes int\'erieurs si et seulement s'il existe $g_{2}\in G_{2}$ tel que $\Psi'=Ad(g_{2})\circ \Psi$. Soit  $\Psi:G_{1}\to G_{2}$ un torseur int\'erieur. Il se d\'eduit de $\Psi$ un isomorphisme $I^{st}_{cusp}(\mathfrak{g}_{1}(F))\to I^{st}_{cusp}(\mathfrak{g}_{2}(F))$ qui ne  d\'epend que de la classe de $\Psi$ modulo automorphismes int\'erieurs. Cet isomorphisme se restreint en un isomorphisme $FC^{st}(\mathfrak{g}_{1}(F))\to FC^{st}(\mathfrak{g}_{2}(F))$ (c'est une cons\'equence de la proposition 12 de \cite{W6}).

Soit $H\in {\cal L}_{F}^{nr,st}$ et soit $H'\in {\cal L}_{F}^{nr}$ un \'el\'ement stablement conjugu\'e \`a $H$. Fixons $g\in G(F^{nr})$ tel que $H'=g^{-1}Hg$ et $gFr(g)^{-1}\in H(F^{nr})$. L'automorphisme $Ad(g)^{-1}$ se restreint en un torseur int\'erieur de $H$ sur $H'$. Il s'en d\'eduit un isomorphisme $\iota_{g}:FC^{st}(\mathfrak{h}_{SC}(F))\to FC^{st}(\mathfrak{h}'_{SC}(F))$. A fortiori, ce dernier espace est non nul et $H'$ appartient \`a ${\cal L}_{F}^{nr,st}$.  Donc
l'ensemble ${\cal L}_{F}^{nr,st}$   est conserv\'e par la conjugaison stable. D'apr\`es \ref{classesdeconjstable}(1), il en est de m\^eme de ${\cal L}_{ell}^{nr,st}$.

 Soit $H\in {\cal L}_{F}^{nr,st}$. On a associ\'e \`a $H$ un $F$-Levi $M^*$ de $G^*$  et un isomorphisme $\psi_{x}:H\to M^*$ d\'efini sur $F^{nr}$. Le lemme \ref{passagestable} montre que $M^*$ appartient \`a ${\cal L}_{F,sd}^{*,st}$.    Posons    $n=Fr_{G^*}(x)u(Fr)^{-1}x^{-1}\in G^*_{AD}(F^{nr})$. On a l'\'egalit\'e 

(1) $Fr_{G^*}(\psi_{x}(h))=n\psi_{x}(Fr_{G}(h))n^{-1}$ pour tout $h\in H(F^{nr})$.

\noindent Il en r\'esulte que $n\in Norm_{G_{AD}^*(F^{nr})}(M^*)$. Les paires de Borel $(n(B^*\cap M^*)n^{-1},nT^*n^{-1})$ et $(B^*\cap M^*,T^*)$ sont toutes deux des paires de Borel de $M^*$ d\'efinies sur $F^{nr}$. Elles sont conjugu\'ees par un \'el\'ement de $M^*(F^{nr})$.  On peut donc fixer $m\in M^*(F^{nr})$ de sorte que $mn\in Norm_{G_{AD}^*(F^{nr})}(T^*)$. Alors $mn$ d\'efinit un \'el\'ement de $W^{I_{F}}$, plus pr\'ecis\'ement de $Norm_{W^{I_{F}}}(M^*)$. Notons $w$ son image dans $W^{I_{F}}(M^*)$. Montrons que

(2) l'\'el\'ement $w$ est  $G$-admissible  et sa classe de $Fr$-conjugaison est uniquement d\'etermin\'ee; si $H\in {\cal L}^{nr,st}_{ell}$, $w$ est $Fr$-r\'egulier.

L'homomorphisme $\psi_{x}$ se restreint en un isomorphisme de $Z(H)^0$ sur $Z(M^*)^0$ qui est d\'efini  sur $F^{nr}$. Il s'en d\'eduit un isomorphisme de $X_{*}(Z(H)^0)^{I_{F}}$ sur $X_{*}(Z(M^*)^0)^{I_{F}}$. D'apr\`es (1) et la d\'efinition de $w$, cet isomorphisme entrelace l'action de $Fr_{G}$ sur l'espace de d\'epart avec celle de $w^{-1}\circ Fr_{G^*}$ sur celui d'arriv\'ee. A fortiori, il se restreint en un isomorphisme de  $X_{*}(Z(H)^0)^{\Gamma_{F}}$ sur $X_{*}(Z(M^*)^0)^{I_{F},w^{-1}\circ Fr_{G^*}}$. On  a d\'efini le groupe $M^*_{w}$ en  \ref{L*stFsd}.
Notons $L_{H}$ le commutant de $A_{H}$ dans $G$. C'est un $F$-Levi de $G$ et la propri\'et\'e pr\'ec\'edente entra\^{\i}ne que $\psi_{x}(L_{H})=M^*_{w}$. Donc  $M^*_{w}$ se tranf\`ere en un $F$-Levi de $G$, c'est-\`a-dire que $w$ est $G$-admissible. Si $H\in {\cal L}^{nr,st}_{ell}$, on a $L_{H}=G$ et $M^*_{w}= G^*$ donc $w$ est $Fr$-r\'egulier. 
 L'\'el\'ement $w$ d\'epend uniquement du choix de $x$. On ne peut changer $x$ qu'en le multipliant \`a gauche par un \'el\'ement de $Norm_{G^*(F^{nr})}(M^*)$. On voit qu'une telle multiplication remplace $w$ par un \'el\'ement qui lui est $Fr$-conjugu\'e. Cela prouve (2).

  A un \'el\'ement $H\in {\cal L}_{F}^{nr,st}$, on vient d'associer    un couple $(M^*,Cl_{Fr}(w))$, o\`u $M^*\in {\cal L}_{F,sd}^{*,st}$ et $Cl_{Fr}(w)\in W^{I_{F}}(M^*)_{G-adm}/ Fr-conj$.

\begin{prop}{Cette application se quotiente en une bijection entre l'ensemble des classes de conjugaison stable dans l'ensemble  ${\cal L}_{F}^{nr,st}$, resp.  ${\cal L}^{nr,st}_{ell}$, et l'ensemble des couples 
 de la forme $(M^*,Cl_{Fr}(w))$, où $M^*$ est un élément de $ {\cal L}_{F,sd}^{*,st}$ et     $Cl_{Fr}(w)$ est un élément de $ W^{I_{F}}(M^*)_{G-adm}/Fr-conj$, resp.  de $W^{I_{F}}(M^*)_{Fr-reg}/ Fr-conj$.} \end{prop}

Preuve. Conservons notre \'el\'ement $H$ et les notations aff\'erentes. Soit $H'\in {\cal L}_{F}^{nr,st}$ un \'el\'ement stablement conjugu\'e \`a $H$.  On peut  fixer $g\in G(F^{nr})$ tel que $H'=g^{-1}Hg$ et $gFr_{G}(g)^{-1}\in H(F^{nr})$. Posons $x'=x\psi(g)$. Alors $Ad(x')\circ \psi(H')=M^*$.  Notre application associe donc \`a $H'$ le m\^eme groupe $M^*$. L'analogue   de $n$ pour $H'$ est $n'=Fr_{G^*}(x')u(Fr)^{-1}{x'}^{-1}=Fr_{G^*}(x)Fr_{G^*}\circ\psi(g)u(Fr)^{-1} \psi(g)^{-1} x^{-1}=Fr_{G^*}(x)u(Fr)^{-1}\psi\circ Fr_{G}(g)\psi(g)^{-1}x^{-1}=Fr_{G^*}(x)u(Fr)^{-1}x^{-1}y=ny$, o\`u $y=x\psi(Fr_{G}(g)g^{-1})x^{-1}$. Puisque $Fr_{G}(g)g^{-1}\in H(F^{nr})$, on a $y\in M^*(F^{nr})$. Donc l'\'el\'ement de $W^{I_{F}}(M^*)$ d\'eduit de $n'$ est le m\^eme que celui d\'eduit de $n$. Notre application envoie donc $H'$ et $H$ sur le m\^eme couple $(M^*,Cl_{Fr}(w))$. 

Inversement, soit $H'\in {\cal L}_{F}^{nr,st}$. Supposons que notre application envoie $H'$ sur le couple $(M^*,Cl_{Fr}(w))$ associ\'e \`a $H$. Fixons $x'\in G^*(F^{nr})$ tel que $Ad(x')\circ \psi(H')=M^*$,  posons $n'=Fr_{G^*}(x')u(Fr)^{-1}{x'}^{-1}$ et notons $w'$ son image dans $W^{I_{F}}(M^*)$. L'hypoth\`ese signifie que  l'on peut fixer $v\in W^{I_{F}}(M^*)$ tel que $w'=Fr_{G^*}(v)wv^{-1}$.     L'homomorphisme  $Norm_{G^*(F^{nr})}(M^*)\to W^{I_{F}}(M^*)$ est surjectif. Relevons $v$ en un \'el\'ement $y\in Norm_{G^*(F^{nr})}(M^*)$.   L'\'egalit\'e $w'=Fr_{G^*}(v)wv^{-1}$  \'equivaut \`a $n'\in  Fr_{G^*}(y)n y^{-1}M^*_{ad}(F^{nr})$, ou encore
$Fr_{G^*}(x')u(Fr)^{-1}{x'}^{-1}\in  Fr_{G^*}(yx)u(Fr)^{-1}(yx)^{-1}M^*_{ad}(F^{nr})$.  Posons $x''=y^{-1}x'$. On a encore $Ad(x'')\circ\psi(H')=M^*$ et maintenant $Fr_{G^*}(x'')u(Fr)^{-1}{x''}^{-1}$ appartient à $ Fr_{G^*}(x)u(Fr)^{-1}x^{-1}M^*_{ad}(F^{nr})$. Posons $g=\psi^{-1}(x^{-1}x'')$. On a $g\in G(F^{nr})$,  $H'=g^{-1}Hg$ et on voit que l'\'egalit\'e pr\'ec\'edente entra\^{\i}ne $gFr_{G}(g)^{-1}\in H$. Puisque c'est un \'el\'ement de $G(F^{nr})$, cela  \'equivaut \`a $gFr_{G}(g)^{-1}\in  H(F^{nr})$. Donc $H'$ est stablement conjugu\'e \`a $H$. Cela d\'emontre que notre application se quotiente en une application injective  de l'ensemble des classes de conjugaison stable dans l'ensemble   ${\cal L}_{F}^{nr,st}$ dans l'ensemble des couples $(M^*,Cl_{Fr}(w))$ indiqu\'es.

 Il reste \`a prouver que cette application est surjective. Fixons $M^*\in {\cal L}_{F,sd}^{*,st}$ et un \'el\'ement $w\in W^{I_{F}}(M^*)_{Fr-reg}$.  
 On rel\`eve $w$ en un \'el\'ement $\dot{w}\in W^{I_{F}}$. Notons   $\sigma\mapsto \sigma_{T^*}$ l'action galoisienne sur $T^*$. On peut munir le tore $T^*$ d'une nouvelle structure galoisienne not\'ee $\sigma\mapsto \sigma_{\dot{w}}$ d\'efinie de la fa\c{c}on suivante. Il existe un unique cocycle ${\bf w}:\Gamma_{F}^{nr}\to W^{I_{F}}$ tel que ${\bf w}(Fr)=\dot{w}^{-1}$. On l'identifie par inflation \`a un cocycle d\'efini sur $\Gamma_{F}$. 
   Alors, pour $\sigma\in \Gamma_{F}$, on  pose $\sigma_{\dot{w}}= {\bf w}(\sigma)\circ\sigma_{T^*}$. 
  Notons $T^*_{\dot{w}}$ le tore ainsi d\'efini. L'identit\'e $j:T^*_{\dot{w}}\to T^*$ est un plongement de $T^*$ dans $G^*$ tel que $\sigma(j):=\sigma_{T^*}\circ j\circ \sigma_{\dot{w}}^{-1}$ est conjugu\'e \`a $j$ pour tout $\sigma\in \Gamma_{F}$. D'apr\`es \cite{K1} corollaire 2.2, il existe un \'el\'ement $y\in G^*$ tel que $Ad(y)^{-1}\circ j$ soit un plongement de $T^*_{\dot{w}}$ dans $G^*$ d\'efini sur $F$.   On peut supposer $y=\pi(y_{sc})$ avec $y_{sc}\in G_{SC}^*(\bar{F})$, o\`u $\pi:G^*_{SC}\to G^*$ est l'homomorphisme naturel. Fixons une extension finie non ramifi\'ee de $F$ telle que ${\bf w}$ soit trivial sur $\Gamma_{F'}$.  Par construction de $T^*_{\dot{w}}$,  $j$ est \'equivariant pour les actions de $\Gamma_{F'}$. Cela entra\^{\i}ne que $y_{sc}\sigma_{G^*}(y_{sc})^{-1}\in T^*_{sc}$ pour tout $\sigma\in \Gamma_{F'}$. Mais $T^*_{sc}$ est un tore induit donc tout cocycle de $\Gamma_{F'}$ dans $T^*_{sc}$ est cohomologue au cocycle trivial. Quitte \`a multiplier $y_{sc}$ \`a gauche par un \'el\'ement de $T^*_{sc}$, on peut donc supposer $y_{sc}\sigma_{G^*}(y_{sc})^{-1}=1$ pour tout $\sigma\in \Gamma_{F'}$. Alors $y_{sc}\in G_{SC}^*(F')$ et aussi $y\in G^*(F')$ a fortiori $y\in G^*(F^{nr})$. Dire que 
$Ad(y)^{-1}\circ j$ est un plongement de $T^*_{\dot{w}}$ dans $G^*$ d\'efini sur $F$ \'equivaut \`a dire que  $yFr_{G^*}(y)^{-1}$ est un \'el\'ement de $Norm_{G^*(F^{nr})}(T^*)$ dont l'image dans $W^{I_{F}}$ est $\dot{w}^{-1}$. 
Posons $H^*=y^{-1}M^*y$.  Ce groupe est  un $F^{nr}$-Levi de $G^*$ d\'efini sur $F$.  De l'application $Ad(y)$ se d\'eduit un isomorphisme de $X_{*}(Z(H^*)^0)^{I_{F}}$ sur $X_{*}(Z(M^*)^0)^{I_{F}}$ qui entrelace les actions $Fr_{G^*}$ et $w^{-1}\circ Fr_{G^*}$. Parce qu'on a supposé que  $w$ était $Fr$-r\'egulier, on a 
 $X_{*}(Z(H^*)^0)^{\Gamma_{F}}=X_{*}(Z(G^*)^0)^{\Gamma_{F}}$. Fixons un sous-tore elliptique maximal $T_{1}^*$ de $H^*$ d\'efini sur $F$. L'\'egalit\'e pr\'ec\'edente implique que $T_{1}^*$ est aussi un sous-tore elliptique maximal de $G^*$. Un tel tore se tranf\`ere \`a toute forme int\'erieure.  On peut donc fixer un sous-tore elliptique maximal $T_{1}$ de $G$ d\'efini sur $F$ et un \'el\'ement $z\in L^*$ de sorte qu'en posant  $\psi_{z}=Ad(z)\circ\psi$, l'isomorphisme $\psi_{z}$ se restreigne en un isomorphisme d\'efini sur $F$ de $T_{1}$ sur $T^*_{1}$.  Cela entra\^{\i}ne que 
 
 (3) $zu(\sigma)\sigma_{G^*}(z)^{-1}\in T_{1,ad}^*$ pour tout $\sigma\in \Gamma_{F}$.  
 
 Posons $H=\psi_{z}^{-1}(H^*)$.  D'apr\`es (3),  $H$ est d\'efini sur $F$ et $\psi_{z}$ se restreint en un isomorphisme d\'efini sur $F$ de $Z(H)^0$ sur $Z(H^*)^0$. Donc aussi de $A^{nr}_{H}$ sur $A^{nr}_{H^*}$. Puisque $H^*$ est le commutant de $A^{nr}_{H^*}$, $H$ est le commutant de $A^{nr}_{H}$. Le commutant d'un tore d\'eploy\'e sur $F^{nr}$ est un $F^{nr}$-Levi  donc $H$ est un $F^{nr}$-Levi de $G$ d\'efini sur $F$.  
L'\'egalit\'e $X_{*}(Z(H^*)^0)^{\Gamma_{F}}=X_{*}(Z(G^*)^0)^{\Gamma_{F}}$ se transf\`ere par $\psi_{z}^{-1} $ en l'\'egalit\'e $X_{*}(Z(H)^0)^{\Gamma_{F}}=X_{*}(Z(G)^0)^{\Gamma_{F}}$. Donc  le groupe $H$ appartient \`a ${\cal L}^{nr}_{ell}$. Posons $x=yz$. On a $Ad(x)\circ \psi(H)=M^*$. Posons $n=Fr_{G^*}(x)u(Fr)^{-1}x^{-1}$. En utilisant (3) et la définition de $y$, on voit que 
$n$ appartient à $Norm_{G^*_{AD}(F^{nr})}(M^*)$ et que son image dans $W^{I_{F}}(M^*)$  est $w$. Le lemme \ref{passagestable}  et  l'hypoth\`ese $FC^{st}(\mathfrak{m}^*_{SC}(F))\not=\{0\}$ entraînent que $FC^{st}(\mathfrak{h}_{SC}(F))\not=\{0\}$. Alors $H$ est un antécédent du couple $(M^*,Cl_{Fr}(w))$ dans  ${\cal L}_{ell}^{nr,st}$.
  
  On a supposé que $w$ était $Fr$-régulier. Il faut considérer le cas général où $w\in W^{I_{F}}(M^*)_{G-adm}$. D'après le lemme \ref{L*stFsd}, quitte à remplacer $w$ par un élément $Fr$-conjugué, on peut fixer un $F$-Levi standard $L^*$ de $G^*$ contenant $M^*_{min}$ et $M^*$ de sorte que $w\in W^{L^*,I_{F}}(M^*)_{Fr-reg}$. Posons $L=\psi^{-1}(L^*)$. C'est un $F$-Levi de $G$. On peut définir une application similaire à celle de l'énoncé en remplaçant le triplet $(G,G^*,\psi)$ par $(L,L^*,\psi_{\vert L})$. D'après ce que l'on vient de prouver, le couple $(M^*,Cl_{Fr}^{L^*}(w))$ a un antécédent $H\in {\cal L}_{ell}^{L,nr,st}$ par l'application relative au triplet $(L,L^*,\psi_{\vert L})$. Il est facile de vérifier que les applications relatives aux deux triplets sont compatibles en un sens facile à préciser et que $H$ est aussi un antécédent de $(M^*,Cl_{Fr}(w))$ par l'application relative à  $(G,G^*,\psi)$. Cela achève la preuve. $\square$

      \subsubsection{Sur l'action de certains automorphismes\label{automorphismes}}

Soit $H\in {\cal L}_{F}^{nr,st}$. Soit $g\in Norm_{G(F^{nr})}(H)$. Supposons $gFr(g)^{-1}\in H$. L'application $Ad(g)$ se restreint alors en  un torseur int\'erieur de $H$ dans lui-m\^eme. Il s'en d\'eduit un automorphisme $\iota$ de $FC^{st}(\mathfrak{h}_{SC}(F))$.

\begin{lem}{Cet automorphisme $\iota$ est l'identit\'e.}\end{lem}

Preuve.   Comme dans plusieurs preuves pr\'ec\'edentes, on se ram\`ene au cas o\`u  $G^*$ est simplement connexe et absolument quasi-simple.   Introduisons comme en \ref{passage} le torseur int\'erieur $\psi_{x}:G\to G^*$ dont on a suppos\'e que $\psi_{x}(H)=M^*$.  D'apr\`es le lemme \ref{passagestable}, $M^*_{SC}$ et $H_{SC}$ sont eux-aussi absolument quasi-simples et on a $FC^{st}(\mathfrak{m}^*_{SC}(F))\not=\{0\}$.  Introduisons la forme quasi-d\'eploy\'ee $H^*$ de $H$ et un torseur int\'erieur $\psi_{H}:H\to H^*$, dont on peut supposer qu'il est un  isomorphisme sur $F^{nr}$. Il s'en d\'eduit un isomorphisme $\iota_{\psi_{H}}:FC^{st}(\mathfrak{h}_{SC}(F))\to FC^{st}(\mathfrak{h}_{SC}^*(F))$. Posons $\iota^*=\iota_{\psi_{H}}\circ \iota\circ \iota_{\psi_{H}}^{-1}$. Il s'agit de voir que $\iota^*$ est l'identit\'e. Posons $\delta^*=\psi_{H}\circ Ad(g)\circ \psi_{H}^{-1}$. C'est un   torseur int\'erieur   de $H^*$ dans lui-m\^eme et $\iota^*$ est d\'eduit de $\delta^*$. Fixons une paire de Borel \'epingl\'ee de $H^*$ conserv\'ee par l'action galoisienne.  On peut \'ecrire $\delta^*=Ad(h^*)\circ \theta^*$ o\`u $h^*\in H^*$ et $\theta^*$ conserve la paire de Borel \'epingl\'ee.    Parce que $\delta^*$ est un torseur intérieur,  $\theta^*$ est un automorphisme de $H^*$ d\'efini sur $F$. L'isomorphisme $\iota^*$ est aussi celui d\'efini par $\theta^*$.  D'apr\`es \cite{W7} 9(5), ce dernier  agit trivialement sur $FC^{st}(\mathfrak{h}^*_{SC}(F))$ sauf dans le cas o\`u $H^*_{SC}$ est de type $(A_{n-1},ram)$, cf. \ref{couplesstables} pour cette notation. Il suffit donc de v\'erifier que, dans le cas o\`u $H^*_{SC}$  est de ce type, on a $\theta^*=1$, autrement dit $\delta^*$ est un automorphisme int\'erieur.  Le compos\'e $\psi_{x}\circ \psi_{H}^{-1}$ est un isomorphisme de $H^*$ sur $M^*$ d\'efini sur $F^{nr}$. Donc $M^*_{SC}$ est aussi du type $(A_{n-1},ram)$. Posons $\delta^{M^*}=\psi_{x}\circ \psi_{H}^{-1}\circ \delta^*\circ\psi_{H}\circ\psi_{x}^{-1}$.  Il suffit encore de montrer que $\delta^{M^*}
$ est un automorphisme int\'erieur. On a $\delta^{M^*}=\psi_{x}\circ Ad(g)\circ \psi_{x}^{-1}=Ad(g^*)$, 
 o\`u $g^*=\psi_{x}(g)$. L'\'el\'ement $g^*$ appartient \`a $G^*(F^{nr})$ et normalise $M^*$. Quitte \`a le multiplier par un \'el\'ement de $M^*(F^{nr})$, on peut supposer qu'il normalise $T^*$. Il lui est associ\'e un \'el\'ement de $ W^{I_{F}}(M^*)$.    Il s'agit de prouver que l'action sur $M^*$ d'un tel \'el\'ement est int\'erieure.   En utilisant les  descriptions de \ref{couplesstables}, on voit que les deux conditions $FC^{st}(\mathfrak{m}^*_{SC}(F))\not=\{0\}$ et $M^*$ est de type $(A_{n-1},ram)$ ne peuvent se produire que  dans les deux cas suivants:

(1) $G^*$ est de type $(A_{m-1},ram)$ avec $m\geq n$ et $m\equiv n\,\,mod\,\,2{\mathbb Z}$;  

(2) $G^*$ est de type $(E_{6},ram)$; dans ce cas $M^*$ est de type $A_{5}$.

  Dans le cas (1), on voit que le groupe $W^{I_{F}}(M^*)$ est engendr\'e par les sym\'etries \'el\'ementaires associ\'ees aux sous-groupes de Levi standard $L^*$ d\'efinis sur $F^{nr}$, contenant strictement $M^*$ et minimaux parmi ceux v\'erifiant ces conditions. Pour un tel groupe, on a soit $L^*_{SC}=Res_{E/F}SL(2) \times M^*_{SC}$, o\`u $E$ est l'extension quadratique ramifi\'ee de $F$ telle que $\Gamma_{E}$ agisse trivialement sur ${\cal D}$, soit $L^*_{SC}$ est de type $A_{n+1}$. Dans le premier cas, la sym\'etrie \'el\'ementaire est celle du premier facteur, elle agit trivialement sur $M^*$. Dans le deuxi\`eme cas, la sym\'etrie \'el\'ementaire est le produit des \'el\'ements de plus grande longueur des groupes de Weyl de  $L^*$ et $M^*$. On voit que cet \'el\'ement fixe les points du diagramme de Dynkin de $M^*$, donc agit sur $M^*$ par automorphisme int\'erieur. Dans le cas (2), $M^*$ est un $F^{nr}$-Levi propre maximal et $W^{I_{F}}(M^*)$  est r\'eduit \`a deux \'el\'ements. L'\'el\'ement non trivial est l'image dans $W^{I_{F}}(M^*)$ du produit des \'el\'ements de plus grande longueur des groupes de Weyl de  $G^*$ et $M^*$. En utilisant la description de \cite{B} de ces \'el\'ements, on voit que cette sym\'etrie fixe encore les points du diagramme de Dynkin de $M^*$, donc agit sur $M^*$ par automorphisme int\'erieur. Cela ach\`eve la d\'emonstration. $\square$

Soient $H,H'\in {\cal L}_{F}^{nr,st}$, supposons que ces groupes sont stablement conjugu\'es. Fixons $g\in G(F^{nr})$ tel que $H'=g^{-1}Hg$ et $gFr(g)^{-1}\in H(F^{nr})$. Comme on l'a dit en \ref{parametrage}, de  l'automorphisme $Ad(g)^{-1}$ se d\'eduit un isomorphisme $\iota_{g}:FC^{st}(\mathfrak{h}_{SC}(F))\to FC^{st}(\mathfrak{h}'_{SC}(F))$. Il est ind\'ependant   du choix de $g$. En effet, on ne peut modifier $g$ qu'en le multipliant \`a gauche par un \'el\'ement $n\in Norm_{G(F^{nr})}(H)$ tel que $nFr(n)^{-1}\in H(F^{nr})$. On a alors $\iota_{ng}=\iota_{g}\iota_{n}$. Or le lemme pr\'ec\'edent dit que $\iota_{n}$ est l'identit\'e. On note simplement $\iota_{H',H}$ cet isomorphisme $\iota_{g}$ pour un choix quelconque de $g$. 

\subsubsection{Le groupe $H_{s_{H}}$\label{HsH}}

{\bf On suppose dans ce paragraphe que} $\boldsymbol{G}$ {\bf est adjoint, simple et quasi-d\'eploy\'e sur} $\boldsymbol{F}$. Cette derni\`ere hypoth\`ese permet d'identifier $G$ \`a $G^*$. On utilise les constructions de \ref{racinesaffines}. 
Soient $H\in {\cal L}_{F}^{nr,st}$ et $x\in Imm^G(H_{ad})$. Posons $s_{H}=p_{H}(x)$ et supposons que $s_{H}$ appartient \`a $S^{st}(H_{AD})$, c'est-\`a-dire que $s_{H}$ est un sommet de $Imm(H_{AD})$ et $FC^{st}(\mathfrak{h}_{SC,s_{H}}({\mathbb F}_{q}))\not=\{0\}$. Puisque $G=G^*$, la construction de \ref{passage} associe au couple $(H,s_{H})$ un couple $(M,s_{M})$, o\`u $M$ est un $F$-Levi standard de $G$ et $s_{M}\in S^{st}(M_{AD})$. Ces deux couples sont conjugu\'es par un \'el\'ement de $G(F^{nr})$. On introduit le sous-ensemble $\Lambda\subset \Delta_{a}^{nr}$ associ\'e en \ref{couples} au couple $(M,s_{M})$. 
Ainsi que tout \'el\'ement de $Imm(G)$, $x$ est conjugu\'e par un \'el\'ement de $G(F)$ \`a un \'el\'ement de $\bar{C}$. 

\begin{lem}{(i) Soit $x'\in \bar{C}$ un \'el\'ement conjugu\'e \`a $x$ par un \'el\'ement de $G(F)$. Alors $\Lambda\subset \Lambda(x')$. 

(ii) Supposons $x\in \bar{C}$. Alors il existe $k\in K_{x}^{0,nr}$ tel que $Ad(k)$ transporte $(H,s_{H})$ sur $(M_{\Lambda},s_{\Lambda})$. En notant $\bar{k}$ la r\'eduction de $k$ dans $G_{x}$, on a $Ad(\bar{k})(H_{s_{H}})=M_{\Lambda,s_{\Lambda}}$.}\end{lem}

Preuve. On a dit que les deux couples $(H,s_{H})$ et $(M,s_{M})$ \'etaient conjugu\'es par un \'el\'ement de $G(F^{nr})$. Ils sont en fait conjugu\'es par un \'el\'ement de $G_{SC}(F^{nr})$: il suffit d'appliquer la construction de \ref{passage} au groupe $G_{SC}$ et aux $F^{nr}$-Levi $H_{sc}$ et $M_{sc}$. Fixons $g\in G_{SC}(F^{nr})$ tel que $Ad(g)$ transporte $(H,s_{H})$ en $(M,s_{M})$. Fixons une extension finie non ramifi\'ee $F'$ de $F$ telle que $g\in G_{SC}(F')$. On a alors $gx\in Imm_{F'}^G(M_{ad})$ et $p_{M}(gx)=s_{M}$. D'apr\`es les propositions \ref{couples} et \ref{casdunexcellentcouple}, on peut fixer $g'\in G_{SC}(F')$ tel que $Ad(g')$ transporte $(M,s_{M})$ en $(M_{\Lambda},s_{\Lambda})$, que $g'gx\in \bar{C}^{nr}$ et que $\Lambda\subset \Lambda(g'gx)$. Les deux points $g'gx$ et $x'$ sont conjugu\'es par un \'el\'ement de $G(F^{nr})$ et appartiennent tous deux \`a $\bar{C}^{nr}$. D'apr\`es le corollaire  \ref{uncorollaire}, ils sont conjugu\'es par un \'el\'ement de $\boldsymbol{\Omega}^{nr}$. L'ensemble $\Lambda(x')$ est l'image  de $\Lambda(g'gx)$ par l'action  d'un tel \'el\'ement.  D'apr\`es \ref{couplesstables}(3), $\Lambda$ est conserv\'e par cette action. Puisque $\Lambda\subset \Lambda(g'gx)$, on a donc $\Lambda\subset \Lambda(x')$, ce qui prouve (i).

Supposons $x\in \bar{C}$. Reprenons la construction pr\'ec\'edente et posons $k_{sc}=g'g$. On a $k_{sc}\in G_{SC}(F')$ et $k_{sc}x\in \bar{C}^{nr}$. Deux \'el\'ements de $\bar{C}^{nr}$ ne sont conjugu\'es par un \'el\'ement de $G_{SC}(F^{nr})$ que s'ils sont \'egaux. Donc $k_{sc}x=x$, d'o\`u $k_{sc}\in K_{SC,x}^{\dag,nr}=K_{SC,x}^{0,nr}$. Notons $k$ l'image de $k_{sc}$ dans $G(F^{nr})$. Alors $k\in K_{x}^{0,nr}$. Par construction de $g'$ et $g$, $Ad(k)$ transporte $(H,s_{H})$ sur $(M_{\Lambda},s_{\Lambda})$. Puisque $k$ fixe $x$, l'\'egalit\'e $p_{H}(x)=s_{H}$ entra\^{\i}ne $p_{M_{\Lambda}}(x)=s_{\Lambda}$ et $Ad(\bar{k})$ transporte $H_{x}=H_{s_{H}}$ sur $M_{\Lambda,x}=M_{\Lambda,s_{\Lambda}}$. Cela prouve (ii). $\square$

 \subsection{Stabilit\'e et espace ${\cal D}(\mathfrak{g}(F))$}

\subsubsection{L'espace  ${\cal D}(\mathfrak{g}(F))$\label{lespacecalD}}

  Posons
$$\boldsymbol{{\cal D}}_{cusp}(\mathfrak{g}(F))=\oplus_{s\in S(G)}C_{nil,cusp}(\mathfrak{g}_{s}({\mathbb F}_{q})).$$
Comme on l'a dit en \ref{lespaceFC}, l'espace $C(\mathfrak{g}_{s}({\mathbb F}_{q}))$ s'identifie \`a un sous-espace de $C_{c}^{\infty}(\mathfrak{g}(F))$. L'espace $\boldsymbol{{\cal D}}_{cusp}(\mathfrak{g}(F))$ est lui-aussi un sous-espace de $C_{c}^{\infty}(\mathfrak{g}(F))$ (la somme reste directe dans cet espace, cf. \cite{W6} paragraphe 8).  Notons ${\cal L}_{F}$ l'ensemble des $F$-Levi de $G$. Posons
$$\boldsymbol{{\cal D}}(\mathfrak{g}(F))=\oplus_{M\in {\cal L}_{F}}\boldsymbol{{\cal D}}_{cusp}(\mathfrak{m}(F)).$$
Le groupe $G(F)$ agit naturellement sur ces espaces par conjugaison, on note ${\cal D}_{cusp}(\mathfrak{g}(F))$ et ${\cal D}(\mathfrak{g}(F))$ leurs quotients de coinvariants. 

On d\'efinit un homomorphisme antilin\'eaire $D^G:\boldsymbol{{\cal D}}_{cusp}(\mathfrak{g}(F))\to I(\mathfrak{g}(F))^*$ par la formule suivante, pour $\varphi\in \boldsymbol{{\cal D}}_{cusp}(\mathfrak{g}(F))$ et $f\in C_{c}^{\infty}(\mathfrak{g}(F))$:
$$D^G_{\varphi}(f)=\int_{A_{G}(F)\backslash G(F)}\int_{\mathfrak{g}(F)}f(g^{-1}Xg)\bar{\varphi}(X)\,dX\, dg.$$
On \'etend cette d\'efinition en un homomorphisme antilin\'eaire $D^G:\boldsymbol{{\cal D}}(\mathfrak{g}(F))\to I(\mathfrak{g}(F))^*$ de la fa\c{c}on suivante. Pour $M\in {\cal L}_{F}$, on dispose de l'homomorphisme d'induction $Ind_{M}^G:I(\mathfrak{m}(F))^*\to I(\mathfrak{g}(F))^*$.  Pour $\varphi\in \boldsymbol{{\cal D}}_{cusp}(\mathfrak{m}(F))$, on pose $D^G_{\varphi}=Ind_{M}^G(D^M_{\varphi})$. Cet homomorphisme $D^G$ se quotiente en un homomorphisme antilin\'eaire $D^G:{\cal D}(\mathfrak{g}(F))\to I(\mathfrak{g}(F))^*$. On a

(1) cet homomorphisme $D^G:{\cal D}(\mathfrak{g}(F))\to I(\mathfrak{g}(F))^*$ est injectif.

Cf. \cite{W5} proposition 5.5. 
On note $\hat{D}^G$ le compos\'e de  $D^G$ et de la transformation de Fourier dans $I(\mathfrak{g}(F))^*$.  

Notons $I^{inst}(\mathfrak{g}(F))$ le sous-espace de $I(\mathfrak{g}(F))$ form\'e des \'el\'ements dont toutes les int\'egrales orbitales stables r\'eguli\`eres sont nulles. On a $SI(\mathfrak{g}(F))=I(\mathfrak{g}(F))/I^{inst}(\mathfrak{g}(F))$. Le dual 
  $SI(\mathfrak{g}(F))^*$ est l'espace des distributions stables. Il s'identifie au sous-espace des \'el\'ements de $I(\mathfrak{g}(F))^*$ qui annulent $I^{inst}(\mathfrak{g}(F))$. Rappelons que $I^{inst}(\mathfrak{g}(F))$ est invariant par transformation de Fourier et qu'il en est donc de m\^eme de $SI(\mathfrak{g}(F))^*$. Notons ${\cal D}^{st}(\mathfrak{g}(F))$ le sous-ensemble des $d\in {\cal D}(\mathfrak{g}(F))$ tels que $D^G(d)\in SI(\mathfrak{g}(F))^*$.

  On  note ${\cal H}$ le sous-espace des $f\in C_{c}^{\infty}(\mathfrak{g}(F))$ telles que le support de $\hat{f}$ soit contenu dans $\mathfrak{g}_{tn}(F)$. On note $res_{{\cal H}}$ l'homomorphisme de restriction $I(\mathfrak{g}(F))^*\to {\cal H}^*$. Notons $S{\cal H}^*$ l'image de $SI(\mathfrak{g}(F))^*$ par $res_{{\cal H}}$. Pour tout sommet $s\in Imm(G_{AD})$, fixons une sous-$\mathfrak{o}_{F}$-alg\`ebre d'Iwahori $\mathfrak{i}_{s}\subset \mathfrak{k}_{s}$. Posons
$${\cal H}^0=\sum_{s\in S(G)}C(\mathfrak{k}_{s}/\mathfrak{i}_{s}).$$
On a ${\cal H}^0\subset {\cal H}$ et on note $res_{{\cal H}^0,{\cal H}}:{\cal H}^*\to {\cal H}^{0,*}$ l'homomorphisme de restriction.  On note  $I(\mathfrak{g}(F))^*_{ent}$ l'espace des distributions \`a support dans $\mathfrak{g}_{ent}(F)$.

\begin{lem}{ (i) La restriction \`a l'espace $res_{{\cal H}}(I(\mathfrak{g}(F))^*_{ent})$ de l'homomorphisme $res_{{\cal H}^0,{\cal H}}
$ est injective.

(ii) Les homomorphismes    $res_{{\cal H}}\circ\hat{D}^G$ et $res_{{\cal H}^0,{\cal H}}\circ res_{{\cal H}}\circ\hat{D}^G$ sont injectifs.
 
(iii) Soit  $d\in {\cal D}(\mathfrak{g}(F))$. Alors $res_{{\cal H}}\circ\hat{D}^G(d)\in S{\cal H}^*$ si et seulement si $d\in {\cal D}^{st}(\mathfrak{g}(F))$.
  }\end{lem}
  
  Preuve. Le (i) est l'une des assertions du th\'eor\`eme 2.1.5 de \cite{D} telle que nous l'avons reformul\'ee en \cite{W6} proposition 7. 
  
  Notons $I(\mathfrak{g}(F))^*_{tn}$ l'espace des distributions \`a support dans $\mathfrak{g}_{tn}(F)$ et $res_{tn}:I(\mathfrak{g}(F))^*\to I(\mathfrak{g}(F))^*_{tn}$ l'endomorphisme de restriction \`a $\mathfrak{g}_{tn}(F)$. Par transformation de Fourier, l'injectivit\'e de $res_{{\cal H}}\circ \hat{D}^G$ \'equivaut \`a celle de $res_{tn}\circ D^G$. Or par d\'efinition, l'image de $D^G$ est form\'ee de distributions \`a support dans $\mathfrak{g}_{tn}(F)$. Donc $res_{tn}$ est l'identit\'e sur l'image de $D^G$. Puisque $D^G$ est injectif, il en est de m\^eme de $res_{tn}\circ D^G$. 
 
  Par construction, l'image de $\hat{D}^G$ est contenue dans $I(\mathfrak{g}(F))^*_{ent}$. D'apr\`es (i),  $res_{{\cal H}^0,{\cal H}}$ est injective sur l'image de $res_{{\cal H}}\circ\hat{D}^G$ et la deuxi\`eme assertion de (ii) r\'esulte de la premi\`ere.

Pour (iii), le sens "si" est \'evident. Soit  $d\in {\cal D}(\mathfrak{g}(F))$. Supposons  $res_{{\cal H}}\circ\hat{D}^G(d)\in S{\cal H}^*$. Fixons $S\in SI(\mathfrak{g}(F))^*$ tel que $res_{{\cal H}}\circ\hat{D}^G(d)=res_{{\cal H}}(S^-)$ (rappelons que, pour $f\in C_{c}^{\infty}(\mathfrak{g}(F))$, on note $f^-$ la fonction $X\mapsto f(-X)$; alors $S^-$ est la distribution telle que $S^-(f)=S(f^-)$ pour tout $f$). Alors $D^G(d)$ co\"{\i}ncide avec $\hat{S}$ sur $\mathfrak{g}_{tn}(F)$. Soit $f\in I^{inst}(\mathfrak{g}(F))$. Puisque $\mathfrak{g}_{tn}(F)$ est un sous-ensemble ouvert et ferm\'e de $\mathfrak{g}(F)$ invariant par conjugaison stable, on peut \'ecrire $f=f_{tn}+f'$, o\`u $f_{tn}$ est \`a support dans $\mathfrak{g}_{tn}(F)$, le support de $f'$ ne coupe pas $\mathfrak{g}_{tn}(F)$ et $f_{tn}, f'\in I^{inst}(\mathfrak{g}(F))$. Par construction, le support de $D^G(d)$ est contenu dans $\mathfrak{g}_{tn}(F)$. Donc $D^G(d)(f)=D^G(d)(f_{tn})=\hat{S}(f_{tn})$ et $\hat{S}(f_{tn})=0$ puisque $\hat{S}$ est stable. Donc $D^G(d)(f)=0$, ce qui entraine que $D^G(d)$ est stable. $\square$

\subsubsection{R\'eduction aux \'el\'ements elliptiques \label{reductionelliptique}}

Posons $I_{cusp}^{inst}(\mathfrak{g}(F))=I_{cusp}(\mathfrak{g}(F))\cap I^{inst}(\mathfrak{g}(F))$. Disons qu'une distribution $J\in I(\mathfrak{g}(F))^*$ est stable sur les elliptiques si elle annule $I_{cusp}^{inst}(\mathfrak{g}(F))$.  

\begin{prop}{Soit $d\in {\cal D}_{cusp}(\mathfrak{g}(F))$. Alors $ D^G(d)$ est stable si et seulement si $D^G(d)$ est stable sur les elliptiques.}\end{prop}

Sous notre hypoth\`ese $(Hyp)_{2}(p)$, l'exponentielle est une bijection de $\mathfrak{g}_{tn}(F)$ sur le sous-ensemble des \'el\'ements topologiquement unipotents de $G(F)$. Alors la proposition r\'esulte de la proposition 8.3 de \cite{W5} descendue \`a l'alg\`ebre de Lie par l'exponentielle. $\square$

\subsubsection{D\'ecomposition selon les Levi\label{decompositionLevi}}
 Pour tout $M\in {\cal L}_{F}$, le normalisateur $Norm_{G(F)}(M)$ agit naturellement dans ${\cal D}_{cusp}(\mathfrak{m}(F))$. On note ${\cal D}_{cusp}(\mathfrak{m}(F))^{Norm_{G(F)}(M)}$ le sous-espace des invariants. Fixons un ensemble de repr\'esentants $\underline{{\cal L}}_{F}$ des classes de conjugaison de   $F$-Levi.  Posons
 $$(1) \qquad \underline{{\cal D}}(\mathfrak{g}(F))=\oplus_{M\in \underline{{\cal L}}_{F}}{\cal D}_{cusp}(\mathfrak{m}(F))^{Norm_{G(F)}(M)}.$$
 L'application naturelle $\underline{{\cal D}}(\mathfrak{g}(F))\to {\cal D}(\mathfrak{g}(F))$ est un isomorphisme.
 
 \begin{lem}{Soit $d\in \underline{{\cal D}}(\mathfrak{g}(F))$, que l'on \'ecrit $d=\oplus_{M\in \underline{{\cal L}}}d_{M}$ conform\'ement \`a la d\'ecomposition (1). Alors $d\in {\cal D}^{st}(\mathfrak{g}(F))$ si et seulement si $d_{M}\in {\cal D}^{st}(\mathfrak{m}(F))$ pour tout $M\in \underline{{\cal L}}$.}\end{lem}
 
 Preuve. Le sens "si" est \'evident, en tenant compte que, pour tout $M\in {\cal L}$,  l'induction $Ind_{M}^G:I(\mathfrak{m}(F))^*\to I(\mathfrak{g}(F))^*$ pr\'eserve la stabilit\'e. 
 
 Dans l'autre sens, supposons $d\in {\cal D}^{st}(\mathfrak{g}(F))$. Notons $\underline{{\cal L}}_{F}(d)$ l'ensemble des
  $M\in \underline{{\cal L}}_{F}$ tels que $d_{M}\not=0$ et notons $l(d)$ son nombre d'\'el\'ements. On raisonne par r\'ecurrence sur $l(d)$.   Si $l(d)=0$, $d_{M}=0$ pour tout $M$ et la conclusion est claire. Supposons $l(d)>0$ et, parmi les $M\in \underline{{\cal L}}_{F}(d)$, fixons un \'el\'ement maximal $L$. Pour tout $M\in {\cal L}_{F}$, notons $Res_{M}^G:I(\mathfrak{g}(F))\to I(\mathfrak{m}(F))$ l'application "terme constant". Soit $f_{L}\in I_{cusp}(\mathfrak{l}(F))^{Norm_{G(F)}(L)}$. Parce que $L$ est maximal parmi les  \'el\'ements de $\underline{{\cal L}}_{F}(d)$, la filtration habituelle de $I(\mathfrak{g}(F))$, cf. \cite{MW} I.4.2, montre qu'il existe $f\in I(\mathfrak{g}(F))$ tel que $Res^G_{M}(f)=0$ pour tout $M\in \underline{{\cal L}}_{F}(d)$, $M\not=L$, et $Res_{L}^G(f)=f_{L}$. Supposons de plus $f_{L}\in I^{inst}_{cusp}(\mathfrak{l}(F))^{Norm_{G(F)}(L)}$. Alors on peut supposer $f\in I^{inst}(\mathfrak{g}(F))$, cf. \cite{MW} I.4.6. Puisque  $d\in {\cal D}^{st}(\mathfrak{g}(F))$, on a $D^G(d)(f)=0$. C'est-\`a-dire $\sum_{M\in \underline{{\cal L}}_{F}(d)}D^M(d_{M})\circ Res_{M}^G(f)=0$. Cette expression se r\'eduit \`a $D^L(d_{L})(f_{L})=0$. On a suppos\'e 
 $f_{L}\in I^{inst}_{cusp}(\mathfrak{l}(F))^{Norm_{G(F)}(L)}$ mais, puisque $d^L$ est suppos\'ee invariante par $Norm_{G(F)}(L)$, le r\'esultat s'\'etend imm\'ediatement \`a toute $f_{L}\in I^{inst}_{cusp}(\mathfrak{l}(F))$. Donc $D^L(d_{L})$ est stable sur les elliptiques. Par la proposition \ref{reductionelliptique}, $D^L(d_{L})$ est stable. Posons $d'=\oplus_{M\in \underline{{\cal L}}_{F}(d), M\not=L} d_{M}$. Puisque $D^G(d')=D^G(d)-D^G(d_{L})$, $D^G(d')$ est encore stable. On a $l(d')=l(d)-1$ et on applique l'hypoth\`ese de r\'ecurrence \`a $d'$: $D^M(d_{M})$ est stable pour tout $M\in \underline{{\cal L}}_{F}(d)$ tel que $M\not=L$. Joint au r\'esultat d\'ej\`a d\'emontr\'e pour $M=L$, c'est la conclusion du lemme. $\square$

\subsubsection{Une premi\`ere description de l'espace ${\cal D}_{cusp}(\mathfrak{g}(F))$\label{premieredescription}}

Notons ${\cal L}^{nr,\star}_{ell}$ l'ensemble des couples $(H,s_{H})$ tels que $H\in {\cal L}^{nr}_{ell}$, $s_{H}$ est un sommet de $Imm(H_{AD})$ et $FC^{H_{s_{H}}}(\mathfrak{h}_{SC,s_{H}}({\mathbb F}_{q}))\not=\{0\}$.
    
 Soit $(H,s_{H})\in {\cal L}^{nr,\star}_{ell}$. Posons $s=p_{H}^{-1}(s_{H})$. C'est un sommet de $Imm^G(H_{ad})$, cf. \ref{immeublesetLevi}(6).  D'apr\`es \ref{groupesenreduction}, il se d\'eduit de $H$ un $\bar{{\mathbb F}}_{q}$-Levi $H_{s}$ de $G_{s}$ qui est isomorphe \`a $H_{s_{H}}$. Puisque $s_{H}$ reste un sommet dans $Imm_{F^{nr}}(H_{AD})$, le plus grand tore central de $H_{s}$ est $A_{H,s}^{nr}$. Puisque $H$ est d\'efini sur $F$, $H_{s}$ est d\'efini sur ${\mathbb F}_{q}$. Puisque $FC^{H_{s_{H}}}(\mathfrak{h}_{SC,s_{H}}({\mathbb F}_{q}))\not=\{0\}$, on a $\mathfrak{h}_{SC,s_{H}}=\mathfrak{h}_{s_{H},SC}=\mathfrak{h}_{s,SC}$ et ${\bf FC}_{{\mathbb F}_{q}}^{H_{s}}(\mathfrak{h}_{s,SC})\not=\emptyset$. 
  Il en r\'esulte que les sous-groupes paraboliques de $G_{s}$ de composante de Levi $H_{s}$ sont  tous conjugu\'es dans $G_{s}$, cf. \ref{faisceauxcaracteres}(1). Alors, le fait que $H_{s}$ soit d\'efini sur ${\mathbb F}_{q}$ implique qu'il existe un ${\mathbb F}_{q}$-Levi $\bar{M}$ de $G_{s}$   et un \'el\'ement $g\in G_{s}$ tel que $H_{s}=g^{-1}\bar{M}g$.
    Posons $n=Fr(g)g^{-1}$. Alors $n\in Norm_{G_{s}}(\bar{M})$. Notons $w$ son image dans $W(\bar{M})$.  Parce que $H$ appartient \`a ${\cal L}^{nr}_{ell}$,  le plus grand sous-tore de $A_{H,s}^{nr}$ d\'eploy\'e sur ${\mathbb F}_{q}$ est central dans $G_{s}$. Cela entra\^{\i}ne que $w\in W(\bar{M})_{Fr-reg}$. On est alors dans la situation de \ref{fonctionscaracteristiques}.  Fixons un \'el\'ement $X_{H}\in A_{H}^{nr}(F)$  entier et de r\'eduction r\'eguli\`ere. Notons $X_{H,s}$ sa r\'eduction dans $A_{H,s}^{nr}({\mathbb F}_{q})$. C'est un élément régulier au sens    que son commutant dans $G_{s}$ est $H_{s}$. Alors, pour $\varphi\in FC^{H_{s_{H}}}(\mathfrak{h}_{SC,s_{H}}({\mathbb F}_{q}))$, on construit comme en \ref{fonctionscaracteristiques} les fonctions $\varphi_{X_{H,s}}$ (que l'on notera simplement $\varphi_{X_{H}}$) et ${\cal Q}_{\varphi}$ sur $\mathfrak{g}_{s}({\mathbb F}_{q})$. Ce sont des fonctions cuspidales. L'application $\varphi\mapsto {\cal Q}_{\varphi}$ est injective. 
 
 Ainsi, au couple $(H,s_{H})\in {\cal L}^{nr,\star}_{ell}$, on a associ\'e $s\in S(G)$ et une application lin\'eaire injective  $\varphi\mapsto {\cal Q}_{\varphi}$ de $FC^{H_{s_{H}}}(\mathfrak{h}_{SC,s_{H}}({\mathbb F}_{q}))$ dans $ C_{nil,cusp}(\mathfrak{g}_{s}({\mathbb F}_{q}))\subset \boldsymbol{{\cal D}}_{cusp}(\mathfrak{g}(F))$.   Posons
 $$\boldsymbol{{\cal K}}^{\star}_{cusp}(\mathfrak{g}(F))=\oplus_{(H,s_{H})\in {\cal L}^{nr,\star}_{ell}}FC^{H_{s_{H}}}(\mathfrak{h}_{SC,s_{H}}({\mathbb F}_{q})).$$
 Des applications lin\'eaires d\'efinies ci-dessus se d\'eduit 
     une application lin\'eaire
 $${\bf j}:\boldsymbol{{\cal K}}_{cusp}^{\star}(\mathfrak{g}(F))\to \boldsymbol{{\cal D}}_{cusp}(\mathfrak{g}(F)).$$
 
 Cette application n'est pas injective mais v\'erifie la propri\'et\'e suivante:
 
 (1) soient $(H,s_{H}),(H',s_{H'})\in {\cal L}_{ell}^{nr,\star}$, $\varphi\in FC^{H_{s_{H}}}(\mathfrak{h}_{SC,s_{H}}({\mathbb F}_{q}))$ et $\varphi'\in FC^{H'_{s_{H'}}}(\mathfrak{h}'_{SC,s_{H'}}({\mathbb F}_{q}))$; supposons ${\bf j}(\varphi)={\bf j}(\varphi')$; alors il existe $g\in G(F)$ tel que $Ad(g)(H)=H'$ et que l'isomorphisme $Ad(g)$ de $H$ sur $H'$ transporte $s_{H}$ et $\varphi$ en $s_{H'}$ et $\varphi'$.
 
 Preuve. L'\'egalit\'e ${\bf j}(\varphi)={\bf j}(\varphi')$ implique que le sommet $s$ de la construction pr\'ec\'edente est le m\^eme pour les deux couples $(H,s_{H})$ et $(H',s_{H'})$. En utilisant \ref{fonctionscaracteristiques}, elle implique aussi qu'il existe $x\in G_{s}({\mathbb F}_{q})$ tel que $Ad(x)(H_{s})=H'_{s}$ et que l'isomorphisme $Ad(x)$ de $H_{s}$ sur $H'_{s}$ transporte $\varphi$ sur $\varphi'$. Fixons un sous-tore $S$ de $H$, déployé maximal sur $F$, tel que $s_{H}$ appartienne \`a l'appartement de $Imm(H_{AD})$ associ\'e \`a $S$ et fixons un sous-tore $T$ de $H$, d\'efini sur $F$, d\'eploy\'e maximal sur $F^{nr}$ et contenant $S$. Il s'en d\'eduit des sous-tores $S_{s}\subset T_{s}$ de $H_{s}$ d\'efinis sur ${\mathbb F}_{q}$. Le tore $S_{s}$ est un sous-tore d\'eploy\'e sur ${\mathbb F}_{q}$ maximal de $H_{s}$ et $T_{s}$ est le commutant de $S_{s}$ dans $H_{s}$. On fixe des sous-tores $S'\subset T'$ de $H'$ v\'erifiant des propri\'et\'es analogues. Puisque  $Ad(x)(H_{s})=H'_{s}$, on peut multiplier $x$ \`a gauche par un \'el\'ement de $H'_{s}({\mathbb F}_{q})$ de sorte que $Ad(x)(S_{s})=S'_{s}$. Alors $Ad(x)$ envoie le commutant $T_{s}$ de $S_{s}$ dans $H_{s}$ sur le commutant $T'_{s}$ de $S'_{s}$ dans $H'_{s}$. En appliquant le lemme 2.2.2 de \cite{D3}, on peut relever $x$ en un \'el\'ement $g\in K_{s}^0$ tel que $Ad(g)(T)=T'$. On a $A_{H}^{nr}\subset T$ et $A_{H'}^{nr}\subset T'$. Puisque $Ad(x)(H_{s})=H'_{s}$, on a aussi $Ad(x)(A_{H,s}^{nr})=A_{H',s}^{nr}$. Cette \'egalit\'e se rel\`eve en $Ad(g)(A_{H}^{nr})=A_{H'}^{nr}$. 
   Puisque $H$, resp. $H'$, est le commutant de $A_{H}^{nr}$, resp. $A_{H'}^{nr}$, on a aussi $Ad(g)(H)=H'$. L'isomorphisme $Ad(g)$ transporte l'image $s_{H}$ de $s$ dans $Imm(H_{AD})$ en l'image de $gs$ dans $Imm(H'_{AD})$. Puisque $g\in K_{s}^0$, on a $gs=s$ et cette image est $s_{H'}$. Enfin, parce que $Ad(x)$ transporte $\varphi$ sur $\varphi'$, $Ad(g)$ a la m\^eme propri\'et\'e. Cela prouve (1).

 Inversement, soient $s\in S(G)$, $\bar{M}$   un ${\mathbb F}_{q}$-Levi  de $G_{s}$ tel que $FC^{\bar{M}}(\bar{\mathfrak{m}}_{SC}({\mathbb F}_{q}))\not=\{0\}$ et  $w$ un \'el\'ement  de $ W(\bar{M})_{Fr-reg}$.     On peut fixer un tore $T\in {\cal T}_{max}^{nr}$ d\'efini sur $F$  tel que $s\in App_{F^{nr}}(T)$ et que la r\'eduction $T_{s}$ de $T$ dans $G_{s}$ soit un sous-tore maximal et maximalement d\'eploy\'e de ce groupe contenu dans $\bar{M}$. On peut relever $w$ en un \'el\'ement de $Norm_{G_{s}}(\bar{M})\cap Norm_{G_{s}}(T_{s})$.    Fixons $g_{w}\in G_{s}$ tel que $g_{w}Fr(g_{w})^{-1}$ appartienne à $ Norm_{G_{s}}(T_{s})$ et se projette sur $w^{-1}$.  Posons   $T'_{s}=g_{w}^{-1}T_{s}g_{w}$. D'apr\`es \cite{D3} lemmes 2.2.2 et 2.3.1, on peut relever $T'_{s}$ en un sous-tore $T'\in {\cal T}_{max}^{nr}$  d\'efini sur $F$ et $g_{w}$ en un \'el\'ement $k\in K^{0,nr}_{s}$ de sorte que $k^{-1}Tk=T'$.  On a encore $s\in App_{F^{nr}}(T')$.    Posons 
 $\bar{M}_{w}=g_{w}^{-1}\bar{M}g_{w}$.  C'est un $\bar{{\mathbb F}}_{q}$-Levi de $G_{s}$ d\'efini sur ${\mathbb F}_{q}$. Le tore $Z(\bar{M}_{w})^0$ se rel\`eve naturellement en un sous-tore $T'_{w}$ de  $T'$ qui est d\'efini sur $F$. Notons $H$ le commutant de $T'_{w}$. 
   C'est un $F^{nr}$- Levi de $G$ qui est d\'efini sur $F$. 
    On a $T'\subset H$. Puisque $s\in App_{F^{nr}}(T')$, $s$ appartient \`a $Imm^G_{F^{nr}}(H_{ad})$ et m\^eme \`a $Imm^G(H_{ad})$ puisque $s\in S(G)$ est fix\'e par l'action galoisienne. Notons $s_{H}$ son image dans $Imm(H_{AD})$. 
   D'apr\`es \ref{groupesenreduction}(1), $s_{H}$ est un sommet de $Imm_{F^{nr}}(H_{AD})$ et on a $H_{s_{H}}\simeq H_{s}=\bar{M}_{w}$. Puisque $s$ est un  sommet de $Imm_{F^{nr}}(H_{AD})$, on a $\mathfrak{h}_{SC,s_{H}}=\mathfrak{h}_{s_{H},SC}=\bar{\mathfrak{m}}_{w,SC}$. D'apr\`es la construction de Lusztig, cf. \ref{fonctionscaracteristiques}, l'hypoth\`ese $FC^{\bar{M}}(\bar{\mathfrak{m}}_{SC}({\mathbb F}_{q}))\not=\{0\}$ implique $FC^{H_{s_{H}}}(\mathfrak{h}_{SC,s_{H}}({\mathbb F}_{q}))\not=\{0\}$. On a 
    
      (2) le plus grand sous-tore de $A^{nr}_{H}$ d\'eploy\'e sur $F$ est central dans $G$.
 
 Ce sous-tore se r\'eduit en le plus grand sous-tore de $Z(\bar{M}_{w})^0$ d\'eploy\'e sur ${\mathbb F}_{q}$. A cause de l'hypoth\`ese que $w$ est $Fr$-r\'egulier, ce sous-tore est central dans $G_{s}$. Puisque $s$ est un sommet de $Imm(G_{AD})$, le plus grand sous-tore central dans $G_{s}$ qui est d\'eploy\'e sur ${\mathbb F}_{q}$ est la r\'eduction de $A_{G}$.   L'assertion (2) en r\'esulte. 
 
 Alors le couple $(H,s_{H})$  appartient \`a ${\cal L}^{nr,\star}_{ell}$. L'homomorphisme $\varphi\mapsto {\cal Q}_{\varphi}$ est une injection de $FC^{H_{s_{H}}}(\mathfrak{h}_{SC,s_{H}}({\mathbb F}_{q}))=FC^{\bar{M}_{w}}(\bar{\mathfrak{m}}_{w,SC}({\mathbb F}_{q}))$ dans $C_{nil,cusp}(\mathfrak{g}_{s}({\mathbb F}_{q}))$. Notons $C(s,\bar{M},w)$ son image.  Ainsi, \`a $\varphi\in C(s,\bar{M},w)$, on peut associer son image inverse  $\varphi^H$ qui appartient à 
  $FC^{H_{s_{H}}}(\mathfrak{h}_{SC,s_{H}}({\mathbb F}_{q})\subset \boldsymbol{{\cal K}}^{\star}_{cusp}(\mathfrak{g}(F))$. Puisque $C_{nil,cusp}(\mathfrak{g}_{s}({\mathbb F}_{q}))$ est somme directe des espaces $C(s,\bar{M},w)$ quand $(\bar{M},w)$ d\'ecrit un certain ensemble, cf. \ref{fonctionscaracteristiques}, on 
  d\'eduit de cette construction une application lin\'eaire
 $${\bf j}':\boldsymbol{{\cal D}}_{cusp}(\mathfrak{g}(F))\to \boldsymbol{{\cal K}}^{\star}_{cusp}(\mathfrak{g}(F)).$$
 Elle d\'epend de plusieurs choix et  n'est pas   un inverse de $j$ (et de fait les espaces $\boldsymbol{{\cal D}}_{cusp}(\mathfrak{g}(F))$ et $\boldsymbol{{\cal K}}_{cusp}^{\star}(\mathfrak{g}(F))$ ne sont pas isomorphes). Toutefois, il r\'esulte des constructions que
 
 (3) ${\bf j}\circ {\bf j}'$ est l'identit\'e de $\boldsymbol{{\cal D}}_{cusp}(\mathfrak{g}(F))$. 
 
  Le groupe $G(F)$ agit naturellement par conjugaison sur $\boldsymbol{{\cal K}}_{cusp}^{\star}(\mathfrak{g}(F))$ et $\boldsymbol{{\cal D}}_{cusp}(\mathfrak{g}(F))$. On note ${\cal K}_{cusp}^{\star}(\mathfrak{g}(F))$ le quotient des coinvariants de $\boldsymbol{{\cal K}}_{cusp}^{\star}(\mathfrak{g}(F))$ et $\pi:\boldsymbol{{\cal K}}_{cusp}^{\star}(\mathfrak{g}(F))\to  {\cal K}_{cusp}^{\star}(\mathfrak{g}(F))$ la projection. Par construction, ${\bf j}$ est \'equivariante pour ces actions. On en d\'eduit une application lin\'eaire quotient 
 $$j:{\cal K}^{\star}_{cusp}(\mathfrak{g}(F))\to {\cal D}_{cusp}(\mathfrak{g}(F)).$$
 
 En notant  $id$ l'identit\'e de $\boldsymbol{{\cal K}}_{cusp}^{\star}(\mathfrak{g}(F))$ , montrons  que
 
 (4) l'image de ${\bf j}'\circ{\bf j}-id$  est contenue dans $Ker(\pi)$.

Soient $(H,s_{H})\in {\cal L}_{ell}^{nr,\star}$ et $\varphi\in FC^{H_{s_{H}}}(\mathfrak{h}_{SC,s_{H}}({\mathbb F}_{q}))$. En  introduisant les termes $\bar{M}$ et $w$ de la construction de ${\bf j}(\varphi)$,  on a ${\bf j}(\varphi)\in C(s,\bar{M},w)$. La construction de son image par ${\bf j}'$ utilise le m\^eme sommet $s$ mais peut utiliser un autre Levi $\bar{M}'$ et un autre \'el\'ement $w'$ tels que $C(s,\bar{M},w)=C(s,\bar{M}',s')$. En tout cas, il existe un couple  $(H',s_{H'})\in {\cal L}_{ell}^{nr,\star}$ et un \'el\'ement  $\varphi'\in FC^{H'_{s_{H'}}}(\mathfrak{h}'_{SC,s_{H'}}({\mathbb F}_{q}))$ de sorte que ${\bf j}'\circ{\bf j}(\varphi)=\varphi'$. D'apr\`es (3), on a ${\bf j}(\varphi)={\bf j}(\varphi')$. Alors (1) nous dit que $\varphi'-\varphi\in Ker(\pi)$, c'est-\`a-dire $({\bf j}'\circ{\bf j}-id)(\varphi)\in Ker(\pi)$. Cela prouve (4). 

Soient $\phi\in \boldsymbol{{\cal D}}_{cusp}(\mathfrak{g}(F))$ et $g\in G(F)$. D'apr\`es (3), on a ${\bf j}'(g(\phi))={\bf j}'(g({\bf j}\circ{\bf j}'(\phi)))$. Parce que ${\bf j}$ est \'equivariante, on a aussi ${\bf j}'(g(\phi))={\bf j}'\circ{\bf j}(g({\bf j}'(\phi)))$. En appliquant (4), on obtient ${\bf j}'(g(\phi))\in Ker(\pi)+g({\bf j}'(\phi))$. Mais $g({\bf j}'(\phi))-{\bf j}'(\phi)$ appartient \`a $Ker(\pi)$ par d\'efinition de cet espace, d'o\`u ${\bf j}'(g(\phi))-{\bf j}'(\phi)\in Ker(\pi)$. Cela entra\^{\i}ne que l'application $\pi\circ {\bf j}'$ 
 se quotiente en une application lin\'eaire $j':{\cal D}_{cusp}(\mathfrak{g}(F))\to {\cal K}^{\star}_{cusp}(\mathfrak{g}(F))$.
 Les assertions (3) et (4) entra\^{\i}nent que $j\circ j'$, resp.  $j'\circ j$ sont les identit\'es de ${\cal D}_{cusp}(\mathfrak{g}(F))$, resp. $ {\cal K}^{\star}_{cusp}(\mathfrak{g}(F))$. C'est-\`a-dire que

(5) les applications $j$ et $j'$ sont des  isomorphismes inverses l'un de l'autre.

\subsubsection{Un calcul de mesures\label{uncalculdemesures}}
Le groupe $G(F)$ agit naturellement par conjugaison dans $\mathfrak{g}_{SC}(F)$. Le sous-groupe $A_{G}(F)$ agit trivialement. Notons $\pi:G_{SC}\to G$ l'homomorphisme naturel. Le quotient $A_{G}(F) \pi(G_{SC}(F))\backslash G(F)$ est compact. On le munit de la mesure telle que, pour $f\in C_{c}^{\infty}(A_{G}(F)\backslash G(F))$, on ait l'\'egalit\'e
$$(1) \qquad \int_{A_{G}(F)\pi(G_{SC}(F))\backslash G(F)}\int_{G_{SC}(F)}f(\pi(g_{sc})x)\,dg_{sc}\,dx=\int_{A_{G}(F)\backslash G(F)}f(g)\,dg.$$
On pose $m(G)=mes(A_{G}(F)\pi(G_{SC}(F))\backslash G(F))$.

Un groupe complexe diagonalisable est par d\'efinition un sous-groupe alg\'ebrique d'un tore complexe. Pour tout tel groupe complexe diagonalisable $Y$, notons $X^*(Y)$ le groupe des homomorphismes alg\'ebriques de $Y$ dans ${\mathbb C}^{\times}$. De l'injection $A_{G}\to G$ se d\'eduit un homomorphisme $Z(\hat{G})\to \hat{A}_{G}$, puis un homomorphisme $b:X^*(\hat{A}_{G})\to X^*(Z(\hat{G})^{I_{F}})^{\Gamma_{F}^{nr}}$. Son conoyau est fini, notons $c_{0}(G)$ le nombre d'\'el\'ements de ce conoyau.  Notons $A_{G}(F)_{c}$, resp. $A_{G}^{nr}(F)_{c}$, le plus grand sous-groupe compact de $A_{G}(F)$, resp. $A_{G}^{nr}(F)$.

\begin{lem}{On a l'\'egalit\'e $m(G)=c_{0}(G)mes(A_{G}^{nr}(F)_{c}) mes(A_{G}(F)_{c})^{-1}$.} \end{lem}

Preuve. Fixons un sommet $s\in Imm(G_{AD})$. Notons $f_{0}$ la fonction caract\'eristique de $K_{s}^0$ dans $G(F)$. D\'efinissons $f\in C_{c}^{\infty}(A_{G}(F)\backslash G(F))$ par 
$$f(g)=\int_{A_{G}(F)}f_{0}(ag)\,da.$$
     L'intersection de $A_{G}(F)$ et de $K_{s}^0$ est $A_{G}(F)_{c}$. 
Alors $f$ est la fonction caract\'eristique de $A_{G}(F)K_{s}^0$, multipli\'ee par $mes(A_{G}(F)_{c})$. Appliquons \`a $f$ la relation (1). Le membre de droite vaut $mes(K_{s}^0)$.  Du c\^ot\'e gauche, l'int\'egrale en $x$ est \`a support dans l'image naturelle de $K_{s}^0$ dans $A_{G}(F)\pi(G_{SC}(F))\backslash G(F)$. Pour $x\in K_{s}^0$, le support de la  fonction $g_{sc}\mapsto f(\pi(g_{sc})x)$ est $K_{SC,s}^0$: ce groupe est clairement dans le support et un \'el\'ement dans le support agit sur $Imm(G_{AD})$ en fixant le sommet $s$, donc appartient \`a $K_{SC,s}^0$. Cette fonction est donc 
la fonction caract\'eristique de $K_{SC,s}^0$ multipli\'ee par $mes(A_{G}( F)_{c})$ et son int\'egrale vaut $mes(K_{SC,s}^0)mes(A_{G}(F)_{c})$. On en d\'eduit l'\'egalit\'e
$$(2) \qquad mes(A_{G}(F)\pi(G_{SC}(F))\backslash A_{G}(F)\pi(G_{SC}(F))K_{s}^0)mes(K_{SC,s}^0)mes(A_{G}(F)_{c})=mes(K_{s}^0),$$
la premi\`ere mesure \'etant prise dans le groupe $A_{G}(F) \pi(G_{SC}(F))\backslash G(F)$. 

 Introduisons les homomorphismes de Kottwitz $w_{G}:G(F)\to X^*(Z(\hat{G})^{I_{F}})^{\Gamma_{F}^{nr}}$ et $w_{A_{G}}:A_{G}(F)\to X^*(\hat{A}_{G})$. Ils sont surjectifs et s'inscrivent dans un diagramme commutatif
 $$(3) \qquad\begin{array}{ccc}A_{G}(F)&\stackrel{w_{A_{G}}}{\to}&X^*(\hat{A}_{G})\\ \downarrow&&\,\downarrow b\\ G(F)&\stackrel{w_{G}}{\to}& X^*(Z(\hat{G})^{I_{F}})^{\Gamma_{F}^{nr}}\\ \end{array}$$
 Montrons que
 
 (4)   $\pi(G_{SC}(F))K_{s}^0=Ker(w_{G})$. 
 
 Notons $\hat{F}^{nr}$ le compl\'et\'e de $F^{nr}$. L'homomorphisme $w_{G}$ est la restriction \`a $G(F)$ d'un homomorphisme $\hat{w}_{G}^{nr}:G(\hat{F}^{nr})\to X^*(Z(\hat{G})^{I_{F}})$. Le groupe $K_{s}^0$ a un analogue sur le corps de base $\hat{F}^{nr}$ que l'on note ${\hat{K}}_{s}^{0,nr}$. D'apr\`es \cite{HR} lemme 17, on a l'\'egalit\'e $\pi(G_{SC}(\hat{F}^{nr})){\hat{K}}_{s}^{0,nr}=Ker(\hat{w}_{G}^{nr})$. Il suffit donc de d\'emontrer que $(\pi(G_{SC}(\hat{F}^{nr})){\hat{K}}_{s}^{0,nr})\cap G(F)=\pi(G_{SC}(F))K_{s}^0$. Soient $g_{sc}\in G_{SC}(\hat{F}^{nr})$ et $k\in {\hat{K}}_{s}^{0,nr}$, supposons $\pi(g_{sc})k\in G(F)$. Pour $\sigma\in \Gamma_{F}^{nr}$, on a $\pi(\sigma(g_{sc}))\sigma(k)=\pi(g_{sc})k$, d'o\`u $\pi(g_{sc}^{-1}\sigma(g_{sc}))=k\sigma(k)^{-1}$. L'image r\'eciproque dans $G_{SC}(\hat{F}^{nr})$ de ${\hat{K}}_{s}^{0,nr}$ est le groupe analogue ${\hat{K}}_{SC,s}^{0,nr}$ et $g_{sc}^{-1}\sigma(g_{sc})$ appartient \`a ce groupe. Alors l'application $\sigma\mapsto g_{sc}^{-1}\sigma(g_{sc})$ est un cocycle \`a valeurs dans ${\hat{K}}_{SC,s}^{0,nr}$. Un tel cocycle est un cobord, essentiellement d'apr\`es le th\'eor\`eme de Lang. On peut donc fixer $k_{sc}\in {\hat{K}}_{SC,s}^{0,nr}$ de sorte que $g_{sc}^{-1}\sigma(g_{sc})=k_{sc}\sigma(k_{sc})^{-1}$ pour tout $\sigma\in \Gamma_{F}^{nr}$. Cela \'equivaut \`a $g_{sc}k_{sc}=\sigma(g_{sc}k_{sc})$, ce qui implique $g_{sc}k_{sc}\in G_{SC}(F)$. On a $\pi(g_{sc})k=\pi(g'_{sc})k'$, o\`u $g'_{sc}=g_{sc}k_{sc}$ et $k'=\pi(k_{sc})^{-1}k$. On a $k'\in {\hat{K}}_{s}^{0,nr}$. Par hypoth\`ese $\pi(g_{sc})k$ appartient \`a $G(F)$ et on vient de voir que $g'_{sc}$ appartient \`a $G_{SC}(F)$. Donc $k'$ appartient \`a $K_{s}^{0}$. L'\'egalit\'e $\pi(g_{sc})k=\pi(g'_{sc})k'$ montre que cet \'el\'ement appartient \`a $\pi(G_{SC}(F))K_{s}^0$, ce qui ach\`eve la preuve de (4).

 On d\'eduit ais\'ement de (4)  que le conoyau de l'injection
 $$A_{G}(F)\pi(G_{SC}(F))\backslash A_{G}(F)\pi(G_{SC}(F))K_{s}^0\to A_{G}(F)\pi(G_{SC}(F))\backslash G(F)$$
 a m\^eme nombre d'\'el\'ements que $X^*(Z(\hat{G})^{I_{F}})^{\Gamma_{F}^{nr}}/w_{G}(A_{G}(F))$. D'apr\`es le diagramme (3), c'est aussi le nombre d'\'el\'ements du conoyau de $b$, c'est-\`a-dire $c_{0}(G)$. On en d\'eduit l'\'egalit\'e
 $$(5) \qquad mes(A_{G}(F)\pi(G_{SC}(F))\backslash A_{G}(F)\pi(G_{SC}(F))K_{s}^0)=m(G)c_{0}(G)^{-1}.$$

  L'immeuble du groupe adjoint de $A_{G}^{nr}$ est r\'eduit \`a un unique point, notons-le $\xi$. On en d\'eduit des sous-groupes parahoriques dans $A_{G}^{nr}(F^{nr})$ et dans $A_{G}^{nr}(F)$ que l'on note simplement $K_{\xi}^{0,nr}$ et $K^0_{\xi}$,   un groupe $A_{G,\xi}^{nr}$ sur ${\mathbb F}_{q}$ et son alg\`ebre de Lie $\mathfrak{a}_{G,\xi}^{nr}$. On a l'\'egalit\'e $K^0_{\xi}=A_{G}^{nr}(F)_{c}$. On a aussi
  $K_{\xi}^{0,nr}= A_{G}^{nr}(F^{nr})\cap K_{s}^{0,nr}$ et ce groupe a m\^eme projection dans $G_{s}$ que $ Z(G)^0(F^{nr})\cap K_{s}^{0,nr}$. On en d\'eduit l'\'egalit\'e
  $$(6) \qquad \mathfrak{g}_{s}=\mathfrak{g}_{SC,s}\oplus \mathfrak{a}^{nr}_{G,\xi}.$$
  Par d\'efinition de nos mesures, on  a les \'egalit\'es
 $$mes(K_{s}^0)=\vert G_{s}({\mathbb F}_{q})\vert \vert \mathfrak{g}_{s}({\mathbb F}_{q}) \vert ^{-1/2},$$
 $$mes(K_{SC,s}^0)=\vert G_{SC,s}({\mathbb F}_{q})\vert \vert \mathfrak{g}_{SC,s} ({\mathbb F}_{q})\vert ^{-1/2}.$$
  Fixons une paire de Borel $(B_{s},T_{s})$ de $G_{s}$ conserv\'ee par l'action galoisienne. On a $\vert G_{s}({\mathbb F}_{q})\vert=\vert T_{s}({\mathbb F}_{q})\vert P_{G_{s}}(q)$, o\`u $P_{G_{s}}$ est le polyn\^ome de Poincar\'e de $G_{s}$. Il y a un homomorphisme naturel $\pi_{s}:G_{SC,s}\to G_{s}$, notons $T_{SC,s}$ l'image r\'eciproque de $T_{s}$.  On a de m\^eme $\vert G_{SC,s}({\mathbb F}_{q})\vert=\vert T_{SC,s}({\mathbb F}_{q})\vert P_{G_{SC, s}}(q)$. Les polyn\^omes $P_{G_{s}}$ et $P_{G_{SC,s}}$ sont \'egaux. D'apr\`es (6), on a
$$X_{*,{\mathbb Q}}(T_{s})=  X_{*,{\mathbb Q}}(T_{SC,s})\oplus X_{*,{\mathbb Q}}(A_{G,\xi}^{nr}).$$
Pour tout tore $T$ d\'efini sur ${\mathbb F}_{q}$, $\vert T({\mathbb F}_{q})\vert $ est \'egal \`a la valeur absolue du d\'eterminant de l'action de $1-q Fr$ agissant sur $X_{*,{\mathbb Q}}(T)$, cf. \cite{C} proposition 3.2.2. On en d\'eduit l'\'egalit\'e $\vert T_{s}({\mathbb F}_{q})\vert = \vert T_{SC,s}({\mathbb F}_{q})\vert \vert A_{G,\xi}^{nr}({\mathbb F}_{q})\vert $. Toujours d'apr\`es (6), on a
$\vert \mathfrak{g}_{s}({\mathbb F}_{q})\vert =\vert \mathfrak{g}_{SC,s}({\mathbb F}_{q})\vert \vert \mathfrak{a}_{G,\xi}^{nr}({\mathbb F}_{q})\vert $.
 En mettant ces calculs bout-\`a-bout, on obtient
 $$(7) \qquad mes(K_{s}^0)=mes(K_{SC,s}^0)\vert A_{G,\xi}^{nr}({\mathbb F}_{q})\vert  \vert \mathfrak{a}_{G,\xi}^{nr}({\mathbb F}_{q})\vert ^{-1/2}$$
 $$=mes(K_{SC,s}^0)mes(K_{\xi}^0)=mes(K_{SC,s}^0)mes(A_{G}^{nr}(F)_{c}).$$ 
 Le  lemme r\'esulte de (2), (5) et (7). $\square$
 
 Pour $J\in I(\mathfrak{g}_{SC}(F))^*$ et $f\in C_{c}^{\infty}(\mathfrak{g}_{SC}(F))$, posons
$$inv^{G}(J)(f)=\int_{A_{G}(F)\pi(G_{SC}(F))\backslash G(F)}J({^gf})\,dg,$$
o\`u $^gf$ est la fonction $X\mapsto f(g^{-1}Xg)$. 
Cela d\'efinit une application lin\'eaire $inv^G$ de 
 $I(\mathfrak{g}_{SC}(F))^* $ dans lui-m\^eme dont l'image est form\'ee de distributions invariantes par l'action de $G(F)$. Remarquons que, pour une distribution $J$ invariante par cette action, on a $inv^G(J)=m(G)J$. 
 
 En particulier, pour $ \varphi\in {\cal D}_{cusp}(\mathfrak{g}_{SC}(F))$, on a d\'efini une distribution $D^{G_{SC}}_{\varphi}\in I(\mathfrak{g}_{SC}(F))^*$. On note $D^G_{\varphi}$ son image par $inv^G$.

    \subsubsection{Description des distributions associ\'ees aux \'el\'ements de $\boldsymbol{{\cal K}}^{\star}_{cusp}(\mathfrak{g}(F))$\label{descriptiondistributions}}

Soient $(H,s_{H})\in {\cal L}^{nr,\star}_{ell}$ et $\varphi\in FC^{H_{s_{H}}}(\mathfrak{h}_{SC,s_{H}}({\mathbb F}_{q}))$. En \ref{premieredescription}, on a associ\'e \`a ces donn\'ees un sommet $s\in Imm(G_{AD})$ et une fonction ${\cal Q}_{\varphi}\in C_{nil,cusp}(\mathfrak{g}_{s}({\mathbb F}_{q}))$. On en d\'eduit une distribution $D^G_{Q_{\varphi}}\in I(\mathfrak{g}(F))^*$. Nous allons d\'ecrire cette distribution de sorte que le point $s$ n'apparaisse pas. Fixons un \'el\'ement $X_{H}\in A_{H}^{nr}(F)$ entier et de r\'eduction r\'eguli\`ere. La fonction $\varphi$ peut  \^etre consid\'er\'ee comme une fonction sur $\mathfrak{k}^{H_{SC}}_{s_{H}}$ invariante par translations par $\mathfrak{k}^{H_{SC},+}_{s_{H}}$. On l'\'etend en une fonction sur $\mathfrak{h}_{SC}(F)$ par $0$ hors de $\mathfrak{k}^{H_{SC}}_{s_{H}}$. On obtient un \'el\'ement de $FC(\mathfrak{h}_{SC}(F))$. D'o\`u une distribution $D_{\varphi}^{H_{SC}}\in I(\mathfrak{h}_{SC}(F))^*$ puis, par le proc\'ed\'e  d\'ecrit en \ref{uncalculdemesures}, une distribution $D^H_{\varphi}\in I(\mathfrak{h}_{SC}(F))^*$. Pour $f\in C_{c}^{\infty}(\mathfrak{g}(F))$, notons $f_{\vert X_{H}+\mathfrak{h}_{SC}}$ la fonction sur $\mathfrak{h}_{SC}(F)$ d\'efinie par $Y\mapsto f(X_{H}+Y)$.  On pose
$$D^G_{X_{H},\varphi}(f)=\int_{H(F)\backslash G(F)}D_{\varphi}^{H}((^gf)_{\vert X_{H}+  \mathfrak{h}_{SC}})\,dg.$$
L'int\'egrale est \`a support compact d'apr\`es \ref{bonselements}(5). Posons
$$c= q^{-dim(A_{H}^{nr})/2}mes(K_{s}^0) mes(K^{H,0}_{s_{H}})^{-1}.$$

\begin{prop}{ La distribution $D^G_{{\cal Q}_{\varphi}^-}$ est \'egale \`a la restriction \`a $\mathfrak{g}_{tn}(F)$ de    $c\hat{D}^G_{X_{H},\varphi} $.}\end{prop}

Preuve.  
On reprend la construction de \ref{fonctionscaracteristiques} pour le groupe $G_{s}$ et $M_{w}=H_{s}$. On peut supposer que l'\'el\'ement $X_{w}$ est la r\'eduction $\bar{X}_{H}$ de $X_{H}$ dans $\mathfrak{g}_{s}({\mathbb F}_{q})$. On a construit les fonctions $\varphi_{X_{w}}^{\natural}$ et $\varphi_{X_{w}}$ que l'on note $\varphi_{\bar{X}_{H}}^{\natural}$ et $\varphi_{\bar{X}_{H}}$. La fonction ${\cal Q}_{\varphi}$ est la restriction \`a $\mathfrak{g}_{s,nil}({\mathbb F}_{q})$ de $\hat{\varphi}_{\bar{X}_{H}}$, cf. \ref{fonctionscaracteristiques}(7). A cause de la conjugaison complexe que l'on a introduite dans la d\'efinition de l'application $D^G$,    $D^G_{{\cal Q}^-_{\varphi}}$ est \'egale \`a la restriction \`a $\mathfrak{g}_{tn}(F)$ de $\hat{D}^G_{\varphi_{\bar{X}_{H}}}$. Calculons la distribution $D^G_{\varphi_{\bar{X}_{H}}}$. 

Pour plus de pr\'ecision, notons $\varphi^{\mathfrak{g}}_{\bar{X}_{H}}$ la fonction sur $\mathfrak{g}(F)$ issue de $\varphi_{\bar{X}_{H}}$ qui, elle,  vit sur $\mathfrak{g}_{s}({\mathbb F}_{q})$. Pour toute $f\in C_{c}^{\infty}(\mathfrak{g}(F))$, posons
$$\Lambda(f)=\int_{\mathfrak{g}(F)}f(Y)\bar{\varphi}^{\mathfrak{g}}_{\bar{X}_{H}}(Y)\,dY.$$
Par d\'efinition,  
 $$(1) \qquad D^G_{\varphi_{\bar{X}_{H}}}( f)=\int_{A_{G}(F)\backslash G(F)}\int_{\mathfrak{g}(F)} f(g^{-1}Y g)\bar{\varphi}^{\mathfrak{g}}_{\bar{X}_{H}}(Y)\, dY\, dg=\int_{A_{G}(F)\backslash G(F)}\Lambda(^gf)\,dg.$$
 
 Notons $f_{s}$ la fonction sur $\mathfrak{g}_{s}({\mathbb F}_{q})$ telle que, pour $Z\in \mathfrak{k}_{s}$, on ait $f_{s}(\bar{Z})=\int_{\mathfrak{k}_{s}^{+}}f(Z+Y)\,dY$, o\`u $\bar{Z}$ est la r\'eduction de $Z$. Par d\'efinition de $\varphi^{\mathfrak{g}}_{\bar{X}_{H}}$, 
on a
$$ \Lambda(f)= \sum_{Z\in \mathfrak{g}_{s}({\mathbb F}_{q})}f_{s}(Z)\bar{\varphi}_{\bar{X}_{H}}(Z),$$
 $$ \vert H_{s_{H}}({\mathbb F}_{q})\vert ^{-1} \sum_{Z\in \mathfrak{g}_{s}({\mathbb F}_{q})}\sum_{g\in G_{s}({\mathbb F}_{q})}f_{s}(Z)\bar{\varphi}^{\natural}_{\bar{X}_{H}}(g^{-1}Zg)$$
$$=\vert H_{s_{H}}({\mathbb F}_{q})\vert ^{-1} \sum_{Z\in \mathfrak{g}_{s}({\mathbb F}_{q})}\sum_{g\in G_{s}({\mathbb F}_{q})}f_{s}(gZg^{-1})\bar{\varphi}^{\natural}_{\bar{X}_{H}}(Z).$$
Par d\'efinition de $\varphi^{\natural}_{\bar{X}_{H}}$, on obtient
$$ \Lambda(f)=\vert H_{s_{H}}({\mathbb F}_{q})\vert ^{-1} \sum_{Z\in \mathfrak{h}_{SC,s_{H}}({\mathbb F}_{q})}\sum_{g\in G_{s}({\mathbb F}_{q})}f_{s}(g(\bar{X}_{H}+Z)g^{-1})\bar{\varphi}(Z).$$
 En utilisant l'\'egalit\'e $mes(K_{s}^0)=\vert G_{s}({\mathbb F}_{q})\vert mes(\mathfrak{k}_{s}^+)$, on voit que, pour tout $Y\in \mathfrak{g}_{s}({\mathbb F}_{q})$, on a l'\'egalit\'e 
$$\sum_{g\in G_{s}({\mathbb F}_{q})}f_{s}(gYg^{-1})=mes(\mathfrak{k}_{s}^+)^{-1}\int_{K_{s}^0}({^kf})_{s}(Y)\,dk.$$
L'\'egalit\'e pr\'ec\'edente devient 
$$ \Lambda(f)=\vert H_{s_{H}}({\mathbb F}_{q})\vert ^{-1}mes(\mathfrak{k}_{s}^+)^{-1} \sum_{Z\in \mathfrak{h}_{SC,s_{H}}({\mathbb F}_{q})}\int_{K_{s}^0} (^kf)_{s}(\bar{X}_{H}+Z)\bar{\varphi}(Z)\,dk,$$
ou encore
$$(2) \qquad   \Lambda(f)=\vert H_{s_{H}}({\mathbb F}_{q})\vert ^{-1}mes(\mathfrak{k}_{s}^+)^{-1} \sum_{Z\in \mathfrak{k}^{H_{SC}}_{s_{H}}/  \mathfrak{k}^{H_{SC},+}_{s_{H}}}\int_{K_{s}^0} (^kf)_{s}(\bar{X}_{H}+\bar{Z})\bar{\varphi}(\bar{Z})\,dk.$$

 Soit $Z\in \mathfrak{k}^{H}_{s_{H}}\cap \mathfrak{h}_{tn}(F) $.  Posons
 $$J_{Z}(f)=\int_{K_{s}^+}\int_{\mathfrak{k}^{H,+}_{s_{H}}} {^kf}(X_{H}+Z+Y))\,dY\,dk.$$
  Pour tout r\'eel $r$, introduisons la fonction $f[r]$ d\'efinie par 
$$f[r](X)=mes(\mathfrak{k}_{s,r})^{-1}\int_{\mathfrak{k}_{s,r}}f(X+U)\,dU$$
pour tout $X\in \mathfrak{g}(F)$. Montrons que

(3) on a l'\'egalit\'e $J_{Z}(f)=J_{Z}(f[r])$ pour tout r\'eel $r>0$. 

Pour $r$ assez grand, $f$ est invariante par translations par $\mathfrak{k}_{s,r}$, donc $f[r]=f$ d'o\`u l'\'egalit\'e ci-dessus. On sait que les sauts de l'application $r\mapsto \mathfrak{k}_{s,r}$ forment un sous-ensemble discret de ${\mathbb R}$. On peut raisonner par r\'ecurrence descendante sur ces sauts et il suffit de prouver que $J_{Z}(f[r])=J_{Z}(f[r+])$ pour tout $r>0$, o\`u on note $f[r+]=f[r+\epsilon]$ pour $\epsilon>0$ assez petit. Fixons $r>0$. On peut aussi bien remplacer $f$ par $f[r+]$. On est ramen\'e \`a prouver que

(4)  si $f$ est invariante par translations par $\mathfrak{k}_{s,r+}$, on a $J_{Z}(f[r])=J_{Z}(f)$.

Pour $k\in K_{s}^0$, on a 
$$^k(f[r])(X)=f[r](k^{-1}Xk)=\int_{\mathfrak{k}_{s,r}}f(k^{-1}Xk+U)\,dU.$$
Par le changement de variables $U\mapsto k^{-1}Uk$, c'est \'egal \`a
$$\int_{\mathfrak{k}_{s,r}}f(k^{-1}Xk+k^{-1}Uk)\,dU=\int_{\mathfrak{k}_{s,r}}{^kf}(X+U)\,dU=(^kf)[r](X).$$
On en d\'eduit
$$J_{Z}(f[r])=mes(\mathfrak{k}_{s,r})^{-1}\int_{K_{s}^+}\int_{\mathfrak{k}^{H,+}_{s_{H}}} \int_{\mathfrak{k}_{s,r}}{^kf}(X_{H}+Z+Y+U)\,dU\,dY\,dk.$$
La propri\'et\'e suivante se d\'eduit de \cite{KM} lemme 2.2.4:
  
  soit  $U\in \mathfrak{k}_{s,r}$; alors il existe $k_{U}\in K_{s,r}$ et $V_{U}\in \mathfrak{k}^{H}_{s_{H},r}$ de sorte que $X_{H}+Z+U=k_{U}^{-1}(X_{H}+Z+V_{U})k_{U}$. 
  
  \noindent On applique cette propri\'et\'e \`a l'\'el\'ement $U$ intervenant dans l'expression ci-dessus de $J_{Z}(f[r])$.   On obtient
  $${^kf}(X_{H}+Z+Y+U)={^kf}(k_{U}^{-1}(X_{H}+Z+V_{U})k_{U}+Y)={^{k_{U}k}f}(X_{H}+Z+V_{U}+k_{U}Yk_{U}^{-1}).$$
  Parce que $Y\in \mathfrak{k}^{H,+}_{s_{H}}\subset \mathfrak{k}_{s,0+}$ et $k_{U}\in K_{s,r}$, on a $k_{U}Yk_{U}^{-1}-Y\in \mathfrak{k}_{s,r+}$. Puisque $f$ et donc aussi ${^{k_{U}k}f}$, est invariante par translation par $\mathfrak{k}_{s,r+}$, on peut remplacer $k_{U}Yk_{U}^{-1}$ par $Y$ dans l'expression ci-dessus:
  $${^kf}( X_{H}+Z+Y+U)={^{k_{U}k}f}(X_{H}+Z+V_{U}+Y).$$
  D'o\`u 
$$J_{Z}(f[r])=mes(\mathfrak{k}_{s,r})^{-1}\int_{K_{s}^+}\int_{\mathfrak{k}^{H,+}_{s_{H}}} \int_{\mathfrak{k}_{s,r}}{^{k_{U}k}f}(X_{H}+Z+V_{U}+Y) \,dU\,dY\,dk$$
$$=mes(\mathfrak{k}_{s,r})^{-1}\int_{\mathfrak{k}_{s,r}}\int_{K_{s}^+}\int_{\mathfrak{k}^{H,+}_{s_{H}}}  {^{k_{U}k}f}(X_{H}+Z+V_{U}+Y) \,dY\,dk\,dU.$$
Par les changements de variables $Y\mapsto Y-V_{U}$ et $k\mapsto k_{U}^{-1}k$, on en d\'eduit
$$J_{Z}(f[r])=mes(\mathfrak{k}_{s,r})^{-1}\int_{\mathfrak{k}_{s,r}}\int_{K_{s}^+}\int_{\mathfrak{k}^{H,+}_{s_{H}}}  {^kf}(X_{H}+Z+Y) \,dY\,dk\,dU$$
$$=mes(\mathfrak{k}_{s,r})^{-1}\int_{\mathfrak{k}_{s,r}} J_{Z}(f)\,dU=J_{Z}(f).$$
Cela prouve (4), d'o\`u (3). 

 Montrons que
$$(5) \qquad  f_{s}(\bar{X}_{H}+\bar{Z})=mes(\mathfrak{k}^{H,+}_{s_{H}})^{-1} J_{Z}(f).$$

D'apr\`es (3), on a 
$$J_{Z}(f)= J_{Z}(f[0+])=\int_{K_{s}^+}\int_{\mathfrak{k}^{H,+}_{s_{H}}} {^k(f[0+])}(X_{H}+Z+Y))\,dY\,dk.$$
Parce que $X_{H}+Z\in \mathfrak{k}_{s}$ et $Y\in \mathfrak{k}^{H,+}_{s_{H}}\subset \mathfrak{k}_{s}^+$, et parce que $f[0+]$ est invariante par translations par $\mathfrak{k}_{s}^+$, on voit que $ {^k(f[0+])}(X_{H}+Z+Y))=f[0+](X_{H}+Z)$ pour tous $k\in K_{s}^+$ et $Y\in \mathfrak{k}^{H,+}_{s_{H}}$. D'o\`u
$$J_{Z}(f)=mes(K_{s}^+)mes(\mathfrak{k}^{H,+}_{s_{H}}) f[0+](X_{H}+Z).$$
Il r\'esulte des d\'efinitions que 
$$ f[0+](X_{H}+Z)=mes(\mathfrak{k}_{s}^+)^{-1}f_{s}(\bar{X}_{H}+\bar{Z})=mes(K_{s}^+)^{-1}f_{s}(\bar{X}_{H}+\bar{Z}).$$
L'assertion (5) en r\'esulte.

Dans le membre de droite de (2), on peut limiter la somme en $Z$ aux \'el\'ements de $\mathfrak{k}^{H_{SC}}_{s_{H}}\cap \mathfrak{h}_{tn}(F)$ puisque $\varphi$ est \`a support nilpotent. On applique alors (5) \`a la fonction ${^kf}$. D'o\`u
 $$ \Lambda(f)=c_{1}\sum_{Z\in \mathfrak{k}^{H_{SC}}_{s_{H}}/\mathfrak{k}^{H_{SC},+}_{s_{H}}}\int_{K_{s}^0}J_{Z}(^kf)\bar{\varphi}(\bar{Z})\,dk$$
 $$=c_{1}\sum_{Z\in \mathfrak{k}^{H_{SC}}_{s_{H}}/\mathfrak{k}^{H_{SC},+}_{s_{H}}}\int_{K_{s}^0}\int_{K_{s}^+}\int_{\mathfrak{k}_{s_{H}}^{H,+}}{^{k'k}f}(X_{H}+Z+Y)\bar{\varphi}(\bar{Z})\,dY\,dk'\,dk.$$
 o\`u
 $$c_{1}=\vert H_{s_{H}}({\mathbb F}_{q})\vert ^{-1}mes(\mathfrak{k}_{s}^+)^{-1}mes(\mathfrak{k}^{H,+}_{s_{H}})^{-1}=mes(\mathfrak{k}_{s}^+)^{-1}mes(K_{s_{H}}^{H,0})^{-1}.$$
 L'int\'egrale sur $K_{s}^+$ est absorb\'ee par celle sur $K_{s}^0$, elle est remplac\'ee par la multiplication par $mes(K_{s}^+)=mes(\mathfrak{k}_{s}^+)$. On peut remplacer $\bar{\varphi}(\bar{Z})$ par $\bar{\varphi}^{\mathfrak{h}_{SC}}(Z)$, o\`u $\varphi^{\mathfrak{h}_{SC}}$ est la fonction sur $\mathfrak{h}_{SC}(F)$ d\'eduite de $\varphi$. On a l'\'egalit\'e $\mathfrak{k}^{H,+}_{s_{H}}=\mathfrak{z}(H)_{tn}(F)\oplus \mathfrak{k}^{H_{SC},+}_{s_{H}}$. On peut remplacer l'int\'egrale en $Y$ par une double int\'egrale en $X'\in \mathfrak{z}(H)_{tn}(F)$ et $Z'\in \mathfrak{k}^{H_{SC},+}_{s_{H}}$. La somme en $Z$ et l'int\'egrale en $Z'$ se reconstituent en une int\'egrale sur $\mathfrak{k}^{H_{SC}}_{s_{H}}$. On peut aussi bien int\'egrer sur tout $\mathfrak{h}^{SC}(F)$ puisque $\varphi^{\mathfrak{h}_{SC}}(Z)$ est nul si $Z\not\in \mathfrak{k}^{H_{SC}}_{s_{H}}$.On obtient
  $$ \Lambda(f)=mes(K_{s_{H}}^{H,0})^{-1}\int_{K_{s}^0}\int_{\mathfrak{z}(H)_{tn}(F)}\int_{\mathfrak{h}^{SC}(F)}{^kf}(X_{H}+X'+Z)\bar{\varphi}^{\mathfrak{h}_{SC}}(Z)\,dZ\,dX'\,dk.$$
  
 Dans le membre de droite de (1), on applique l'\'egalit\'e pr\'ec\'edente \`a la fonction ${^gf}$. L'int\'egrale sur $K_{s}^0$ est absorb\'ee par celle sur $A_{G}(F)\backslash G(F)$ et on obtient
  $$D_{\varphi_{\bar{X}_{H}}}^G(f)=c_{2}\int_{A_{G}(F)\backslash G(F)}\int_{\mathfrak{z}(H)_{tn}(F)}\int_{\mathfrak{h}^{SC}(F)}{^gf}(X_{H}+X'+Z)\bar{\varphi}^{\mathfrak{h}_{SC}}(Z)\,dZ\,dX'\,dg,$$
  o\`u
  $$c_{2}=mes(K_{s}^0) mes(K^{H,0}_{s_{H}})^{-1}.$$
  
  On peut d\'ecomposer l'int\'egrale sur $A_{G}(F)\backslash G(F)$ en une double int\'egrale sur les ensembles $H(F)\backslash G(F)$ et $A_{H}(F)\backslash H(F)$ (rappelons que $A_{H}=A_{G}$ puisque $H$ est elliptique). Parce que $\varphi^{\mathfrak{h}_{SC}}$ est tr\`es cuspidale, l'int\'egrale sur $A_{H}(F)\backslash H(F)$ est \`a support compact et il en est manifestement de m\^eme de la double int\'egrale sur  $A_{H}(F)\backslash H(F)$ et $\mathfrak{z}(H)_{tn}(F)$. On peut donc permuter ces deux int\'egrales. Alors, 
  puisque $\varphi^{\mathfrak{h}_{SC}}$ est \`a support dans $\mathfrak{h}_{SC,tn}(F)$, la propri\'et\'e (5) de \ref{bonselements} montre que la double int\'egrale sur $H(F)\backslash G(F)$ et $\mathfrak{z}(H)_{tn}(F)$ est \`a support compact. On les permute, puis on reconstitue les deux int\'egrales int\'erieures en une int\'egrale sur $A_{G}(F)\backslash G(F)$.  On obtient
 $$(6) \qquad D_{\varphi_{\bar{X}_{H}}}^G(f)=c_{2}\int_{\mathfrak{z}(H)_{tn}(F)}D_{X_{H}+X',\varphi}^G(f),$$
 o\`u
 $$D_{X_{H}+X',\varphi}^G(f)=\int_{A_{G}(F)\backslash G(F)} \int_{\mathfrak{h}^{SC}(F)}{^gf}(X_{H}+X'+Z)\bar{\varphi}^{\mathfrak{h}_{SC}}(Z)\,dZ\,dg.$$
 Dans le cas $X'=0$,  on peut r\'ecrire
 $$D_{X_{H},\varphi}^G(f)=\int_{H(F)\backslash G(F)}\int_{A_{H}(F)\backslash H(F)}\int_{\mathfrak{h}^{SC}(F)}{^{hg}f}(X_{H}+Z)\bar{\varphi}^{\mathfrak{h}_{SC}}(Z)\,dZ\,dh\,dg,$$
 $$=\int_{H(F)\backslash G(F)}D_{\varphi}^H((^gf)_{X_{H}+\mathfrak{h}_{SC}})\,dg.$$
 C'est bien la d\'efinition que l'on a donn\'ee au d\'ebut du paragraphe. 
 
 Montrons que:
 
 (7) pour $f\in {\cal H}$, $D_{X_{H}+X',\varphi}^G(f)$ ne d\'epend pas de $X'\in \mathfrak{z}(H)_{tn}(F)$.
 
 Puisque $\varphi^{\mathfrak{h}_{SC}}$ est \`a support topologiquement nilpotent, $D_{X_{H}+X',\varphi}^G$ est une distribution \`a support entier. D'apr\`es le (i) de la proposition \ref{lespacecalD}, il suffit de prouver que  $D_{X_{H}+X',\varphi}^G(f)$ ne d\'epend pas de $X'$ quand $f\in {\cal H}^0$. On peut donc fixer un sommet $x\in Imm(G_{AD})$ et supposer que $f$ est \`a support dans $\mathfrak{k}_{x}$ et est invariante par une sous-$\mathfrak{o}_{F}$-alg\`ebre d'Iwahori de $\mathfrak{k}_{x}$, a fortiori par $\mathfrak{k}_{x}^+$. Soit $g\in G(F)$ tel que l'int\'egrale int\'erieure dans la d\'efinition de $D_{X_{H}+X',\varphi}^G(f)$ soit non nulle. Alors le support de ${^gf}$ doit contenir un \'el\'ement de la forme $X_{H}+Z'$ avec $Z'\in \mathfrak{h}_{tn}(F)$. Le support de ${^gf}$ est contenu dans $\mathfrak{k}_{gx}$. Le lemme \ref{ImmGHad} entra\^{\i}ne que $gx$ appartient \`a $Imm^G(H_{ad})$. Pour un tel point, on a $\mathfrak{z}(H)_{tn}(F)\subset \mathfrak{k}_{gx}^+$. Mais ${^gf}$ est invariante par ce r\'eseau, donc invariante par translation par $X'$. Cet \'el\'ement $X'$ dispara\^{\i}t donc de l'int\'egrale. Cela d\'emontre (7).

 En cons\'equence, la formule (6) entra\^{\i}ne  

(8) on a l'\'egalit\'e  $ D_{\varphi_{\bar{X}_{H}}}^G(f)=cD_{X_{H},\varphi}^G(f)$ pour tout $f\in {\cal H}$, 

\noindent o\`u
 $$c=mes(\mathfrak{z}(H)_{tn}(F))mes(K_{s}^0) mes(K_{s_{H}}^{H,0})^{-1}.$$
 Pour tout tore $T$ d\'efini sur $F$, la d\'efinition de nos mesures implique que $mes(\mathfrak{t}_{tn}(F))=q^{-dim(T^{nr})/2}$, o\`u $T^{nr}$ est le plus grand sous-tore de $T$ d\'eploy\'e sur $F^{nr}$. Donc $c$ co\"{\i}ncide avec la d\'efinition que l'on a donn\'ee avant l'\'enonc\'e. 
 Dire que deux distributions co\"{\i}ncident sur ${\cal H}$ revient \`a dire que leurs transform\'ees de Fourier co\"{\i}ncident sur $\mathfrak{g}_{tn}(F)$. Puisque $D_{Q_{\varphi}^-}^G$ est pr\'ecis\'ement la restriction \`a $\mathfrak{g}_{tn}(F)$ de $\hat{D}^G_{\varphi_{\bar{X}_{H}}}$,  (8) entra\^{\i}ne l'\'egalit\'e  de l'\'enonc\'e. $\square$

 \subsubsection{Expression int\'egrale de $D^G_{X_{H},\varphi}$\label{expressionintegrale}}
On note $\mathfrak{g}_{reg}(F)/conj$ l'ensemble des classes de conjugaison par $G(F)$ dans l'ensemble $\mathfrak{g}_{reg}(F)$. Il est naturellement muni d'une structure de vari\'et\'e analytique sur $F$ et d'une mesure de sorte que, pour tout sous-tore maximal $T$ de $G$ d\'efini sur $F$, l'application naturelle $\mathfrak{t}(F)\cap \mathfrak{g}_{reg}(F)\to \mathfrak{g}_{reg}(F)/conj$ soit localement un isomorphisme analytique  et pr\'eserve localement les mesures. Pour $f\in C_{c}^{\infty}(\mathfrak{g}(F))$, l'application $X\mapsto I^G(X,f)$ d\'efinie sur $\mathfrak{g}_{reg}(F)$ se quotiente en une application d\'efinie sur $\mathfrak{g}_{reg}(F)/conj$. La formule suivante est une fa\c{c}on d'\'enoncer la classique formule de Weyl:
 $$(1) \qquad \int_{\mathfrak{g}(F)}f(X)\,dX=\int_{\mathfrak{g}_{reg}(F)/conj}I^G(X,f)d^G(X)^{1/2}\,dX.$$
  
On note $\mathfrak{g}_{SC,reg}(F)/G-conj$  l'ensemble des classes de conjugaison par $G(F)$ dans l'ensemble $\mathfrak{g}_{SC,reg}(F)$.  Il est lui aussi naturellement muni d'une structure de vari\'et\'e analytique sur $F$. La projection naturelle $\tau:\mathfrak{g}_{SC,reg}(F)/conj\to \mathfrak{g}_{SC,reg}(F)/G-conj$ est un isomorphisme local \`a fibres finies. On munit l'ensemble d'arriv\'ee de la mesure  telle que $\tau$ pr\'eserve localement les mesures. Pour $f\in C_{c}^{\infty}(\mathfrak{g}_{SC}(F))$ et $X\in  \mathfrak{g}_{SC,reg}(F)/G-conj$, on pose $I^G(X,f)=\sum_{X'\in \tau^{-1}(X)}I^{G_{SC}}(X',f)$. Soit $X\in \mathfrak{g}_{SC,reg}(F)$, notons $T$ son commutant dans $G$. Il r\'esulte de la d\'efinition de nos mesures que l'injection $T_{sc}(F)\backslash G_{SC}(F)\to T(F)\backslash G(F)$ pr\'eserve les mesures. On en d\'eduit la formule \`a laquelle on s'attend:
 $$I^G(X,f)=d^G(X)^{1/2}\int_{T(F)\backslash G(F)}f(g^{-1}Xg)\, dg.$$

Pour $f\in C_{c}^{\infty}(\mathfrak{g}_{SC}(F))$,  la formule (1) appliqu\'ee au groupe $G_{SC}$ et les d\'efinitions entra\^{\i}nent
$$\int_{\mathfrak{g}_{SC}(F)}f(X)\,dX=\int_{\mathfrak{g}_{SC,reg}/G-conj}I^G(X,f)d^G(X)^{1/2}\,dX.$$

 Soit $H\in {\cal L}^{nr}_{ell}$ et soit $\varphi\in FC^H(\mathfrak{h}_{SC}(F))$. On a d\'efini la distribution $D^{H_{SC}}_{\varphi}$. Elle est localement int\'egrable d'apr\`es Harish-Chandra. On note $e^{H_{SC}}_{\varphi}$ la fonction associ\'ee. Puisque, par hypoth\`ese,  $\varphi$ est invariante par conjugaison par $H(F)$, cette fonction $e^{H_{SC}}_{\varphi}$ est invariante par conjugaison par $H(F)$.  Pour toute $f\in C_{c}^{\infty}(\mathfrak{h}_{SC}(F))$, on a l'\'egalit\'e
$$D^{H_{SC}}_{\varphi}(f)=\int_{\mathfrak{h}_{SC,reg}(F)/H-conj}I^H(Y,f)e^{H_{SC}}_{\varphi}(Y)d^{H}(Y)^{1/2}\, dY.$$
 
En se rappelant la d\'efinition de $D^H_{\varphi}(f)$, on obtient
$$D^{H}_{\varphi}(f)=\int_{\mathfrak{h}_{SC,reg}(F)/H-conj}I^H(Y,f)e^{H}_{\varphi}(Y)d^{H}(Y)^{1/2}\, dY,$$
o\`u $e^H_{\varphi}(Y)=m(H)e^{H_{SC}}_{\varphi}(Y)$. 
Puisque $\varphi$ est \`a support dans $\mathfrak{h}_{SC,tn}(F)$, la fonction $e^{H}_{\varphi}$ est \`a support dans l'image de  $\mathfrak{h}_{SC,tn}(F)\cap \mathfrak{h}_{SC,reg}(F)$ dans $\mathfrak{h}_{SC,reg}(F)/H-conj$. Notons cette image $\mathfrak{h}_{SC,reg,tn}(F)/H-conj$.

\begin{lem}{Pour toute $f\in C_{c}^{\infty}(\mathfrak{g}(F))$, on a l'\'egalit\'e
$$D^G_{X_{H},\varphi}(f)=\int_{\mathfrak{h}_{SC,reg}(F)/H-conj}I^G(X_{H}+Y,f)e^{H}_{\varphi}(Y)d^{G}(X_{H}+Y)^{1/2}\, dY.$$}\end{lem}
Preuve. On a 
$$D^G_{X_{H},\varphi}(f)=\int_{H(F)\backslash G(F)}D_{\varphi}^{H}((^gf)_{\vert X_{H}+ \mathfrak{h}_{SC}})\,dg$$
d'o\`u
$$(1) \qquad D^G_{X_{H},\varphi}(f)=\int_{H(F)\backslash G(F)}\int_{\mathfrak{h}_{SC,reg}(F)/H-conj}I^H(Y,(^gf)_{\vert  X_{H}+\mathfrak{h}_{SC}})e^{H}_{\varphi}(Y)d^{H}(Y)^{1/2}\, dY\,dg.$$
D'apr\`es la propri\'et\'e de support de $e^{H}_{\varphi}$ et d'apr\`es \ref{bonselements}(5), il existe un sous-ensemble compact $\omega\subset H(F)\backslash G(F)$ tel que $I^H(Y,(^gf)_{\vert X_{H}+ \mathfrak{h}_{SC}})=0$ pour tout $Y\in \mathfrak{h}_{SC,reg,tn}(F)/H-conj$  et pour tout $g\not\in \omega$. La double int\'egrale ci-dessus est donc absolument convergente et on peut intervertir les int\'egrations. Pour $Y\in \mathfrak{h}_{SC,reg,tn}(F)/H-conj$, on a
$$\int_{H(F)\backslash G(F)}I^H(Y,(^gf)_{\vert X_{H}+ \mathfrak{h}_{SC}}) \,dg=d^H(Y)^{1/2}\int_{H(F)\backslash G(F)}\int_{H_{Y}(F)\backslash H(F)}(^gf)(X_{H}+h^{-1}Yh)\,dh\,dg$$
$$=d^H(Y)^{1/2}\int_{H(F)\backslash G(F)}\int_{H_{Y}(F)\backslash H(F)}f(g^{-1}h^{-1}(X_{H}+Y)hg)\,dh\,dg.$$
On a $H_{Y}=G_{X_{H}+Y}$. La double int\'egrale se reconstitue en une int\'egrale sur $G_{X_{H}+Y}(F)\backslash G(F)$. On a aussi $d^H(Y)^{1/2}=d^G(X_{H}+Y)^{1/2}$ d'apr\`es \ref{bonselements}(3). On obtient
$$\int_{H(F)\backslash G(F)}I^H(Y,(^gf)_{\vert X_{H}+ \mathfrak{h}_{SC}}) \,dg=I^G(X_{H}+Y,f).$$
En utilisant encore une fois l'\'egalit\'e $d^H(Y)^{1/2}=d^G(X_{H}+Y)^{1/2}$, l'\'egalit\'e (1) devient celle de l'\'enonc\'e. $\square$

 \subsubsection{Une deuxi\`eme  description de l'espace $D^G({\cal D}_{cusp}(\mathfrak{g}(F)))$\label{deuxiemedescription}}
On a d\'efini l'espace
$$\boldsymbol{{\cal K}}^{\star}_{cusp}(\mathfrak{g}(F))=\oplus_{(H,s_{H})\in {\cal L}^{nr,\star}_{ell}}FC^{H_{s_{H}}}(\mathfrak{h}_{SC,s_{H}}({\mathbb F}_{q})),$$
et son quotient des coinvariants ${\cal K}^{\star}_{cusp}(\mathfrak{g}(F))$ par l'action de $G(F)$. 
 Posons
$$\boldsymbol{{\cal K}}_{cusp}(\mathfrak{g}(F))=\oplus_{H\in {\cal L}^{nr}_{ell}}FC^{H}(\mathfrak{h}_{SC}(F)).$$
De nouveau, le groupe $G(F)$ agit naturellement sur cet espace et on note ${\cal K}_{cusp}(\mathfrak{g}(F))$ le quotient des coinvariants. Il y a un homomorphisme naturel  surjectif $\boldsymbol{{\cal K}}^{\star}_{cusp}(\mathfrak{g}(F))\to \boldsymbol{{\cal K}}_{cusp}(\mathfrak{g}(F))$: on envoie chaque sous-espace $FC^{H_{s_{H}}}(\mathfrak{h}_{SC,s_{H}}({\mathbb F}_{q}))$ dans $FC(\mathfrak{h}_{SC}(F))$ puis on projette sur le sous-espace $FC^{H}(\mathfrak{h}_{SC}(F))$ des invariants par $H(F)$. Il est clair que, par passage aux coinvariants par $G(F)$, cet homomorphisme devient un  isomorphisme de ${\cal K}^{\star}_{cusp}(\mathfrak{g}(F))$ sur ${\cal K}_{cusp}(\mathfrak{g}(F))$.

 Pour tout $H\in {\cal L}^{nr}_{ell}$, fixons un \'el\'ement $X_{H}\in A^{nr}_{H}(F)$ entier et de r\'eduction r\'eguli\`ere.   
 On d\'efinit  l'application antilin\'eaire
$$\boldsymbol{\Delta}^G:\boldsymbol{{\cal K}}_{cusp}(\mathfrak{g}(F))\to I(\mathfrak{g}(F))^*$$
telle que, pour $H\in {\cal L}^{nr}_{ell}$ et $\varphi\in FC^H(\mathfrak{h}_{SC}(F))$, $\boldsymbol{\Delta}^G(\varphi)=D_{X_{H},\varphi}^G$. 
Cette application d\'epend du choix des $X_{H}$.  Pour $\varphi\in \boldsymbol{{\cal K}}(\mathfrak{g}(F))$, notons $k^G(\varphi)$ la restriction \`a $\mathfrak{g}_{tn}(F)$ de la transform\'ee de Fourier $\hat{\boldsymbol{\Delta}}^G(\varphi)$ de $\boldsymbol{\Delta}^G(\varphi)$. On obtient une application antilin\'eaire
 $$k^G:\boldsymbol{{\cal K}}_{cusp}(\mathfrak{g}(F))\to I(\mathfrak{g}(F))^*.$$
   La proposition \ref{descriptiondistributions} et l'isomorphisme \ref{premieredescription}(1) entra\^{\i}nent que  $k^G$ ne d\'epend plus des choix et qu'elle se quotiente en un isomorphisme antilin\'eaire encore not\'e
   
 (1) $k^G:  
      {\cal K}_{cusp}(\mathfrak{g}(F))\simeq D^G({\cal D}_{cusp}(\mathfrak{g}(F)))$.

\subsubsection{Un premier r\'esultat d'instabilit\'e \label{instabilite1}}
Pour ${\bf G}'\in Endo_{ell}(G)$,  notons $I(\mathfrak{g}(F))^*[{\bf G}']$ l'espace des $d\in I(\mathfrak{g}(F))^*$ tels que $d$ annule $I_{cusp}(\mathfrak{g}(F),{\bf G}'')$ pour tout ${\bf G}''\not={\bf G}'$. Notons $I(\mathfrak{g}(F))^*_{ind}$ la somme sur tous les Levi propres des espaces de distributions induites $Ind_{M}^G(I(\mathfrak{m}(F))^*)$. C'est aussi le sous-espace des \'el\'ements de $I(\mathfrak{g}(F))^*$ qui annulent $I_{cusp}(\mathfrak{g}(F))$. On a une suite exacte
$$0\to I(\mathfrak{g}(F))^*_{ind}\to I(\mathfrak{g}(F))^*[{\bf G}']\to  I_{cusp}(\mathfrak{g}(F),{\bf G}')^*\to 0.$$
L'espace $I(\mathfrak{g}(F))^*[{\bf G}']$ est stable par transformation de Fourier. Posons $I(\mathfrak{g}(F))^*[\not={\bf G}]=\sum_{{\bf G}'\not={\bf G}}I(\mathfrak{g}(F))^*[{\bf G}']$. Soulignons que les espaces $I(\mathfrak{g}(F))^*[{\bf G}']$ ne sont pas en somme directe: leur intersection est l'espace $I(\mathfrak{g}(F))^*_{ind}$. 

  Pour ${\bf G}'\in Endo_{ell}(G)$, on pose ${\cal D}_{cusp}(\mathfrak{g}(F),{\bf G}')=\{d\in{\cal D}_{cusp}(\mathfrak{g}(F)); D^G(d)\in I(\mathfrak{g}(F))^*[{\bf G}']\}$.   En \cite{W5}, on a d\'efini des espaces de distributions analogues sur le groupe $G(F)$ au lieu de son alg\`ebre de Lie $\mathfrak{g}(F)$. Dans le cas des groupes, les assertions qui suivent r\'esultent du lemme 6.4 et de la proposition 7.2 de \cite{W5}. Notre hypoth\`ese sur $p$ permet de les descendre \`a nos distributions \`a support dans $\mathfrak{g}_{tn}(F)$. On obtient ainsi les \'egalit\'es 
 $${\cal D}_{cusp}(\mathfrak{g}(F))=\oplus_{{\bf G}'\in Endo_{ell}(G)}{\cal D}_{cusp}(\mathfrak{g}(F),{\bf G}')$$
et     ${\cal D}^{st}_{cusp}(\mathfrak{g}(F))={\cal D}_{cusp}(\mathfrak{g}(F),{\bf G})$. 
On note ${\cal D}^{inst}_{cusp}(\mathfrak{g}(F))$ la somme des ${\cal D}_{cusp}(\mathfrak{g}(F),{\bf G}')$ sur les ${\bf G}'\in Endo_{ell}(G)$ tels que ${\bf G}'\not={\bf G}$.

   On a une inclusion naturelle $Endo_{ell}(G)\subset Endo_{ell}(G_{SC})$. En notant $I_{cusp}^G(\mathfrak{g}_{SC}(F))$ l'espace des invariants par $G(F)$ dans $I_{cusp}(\mathfrak{g}_{SC}(F))$, on a 
$$I_{cusp}^G(\mathfrak{g}_{SC}(F))=\oplus_{{\bf G}'\in Endo_{ell}(G)}   I(\mathfrak{g}_{SC}(F),{\bf G}').$$
   On a d\'efini l'espace $FC^G(\mathfrak{g}_{SC}(F))$. D'apr\`es l'\'egalit\'e pr\'ec\'edente, il 
    est \'egal \`a la somme des $FC^{G_{SC}}(\mathfrak{g}_{SC}(F),{\bf G}')$ o\`u ${\bf G}'$ parcourt  $Endo_{ell}(G)$. On a l'\'egalit\'e $FC^{G,st}(\mathfrak{g}_{SC}(F))=FC^{G_{SC},st}(\mathfrak{g}_{SC}(F))$, c'est pourquoi on a not\'e simplement $FC^{st}(\mathfrak{g}_{SC}(F))$ cet espace. On pose 
$$FC^{G,inst}(\mathfrak{g}_{SC}(F))=\oplus_{{\bf G}'\in Endo_{ell}(G),{\bf G}'\not={\bf G}}FC^{G_{SC}}(\mathfrak{g}_{SC}(F),{\bf G}').$$

Posons 
$$\boldsymbol{{\cal K}}_{cusp}^{\natural}(\mathfrak{g}(F))=\oplus_{H\in {\cal L}^{nr}_{ell}}FC^{H,inst}(\mathfrak{h}_{SC}(F)),$$
$$\boldsymbol{{\cal K}}_{cusp}^{\sharp}(\mathfrak{g}(F))=\oplus_{H\in {\cal L}^{nr}_{ell}}FC^{st}(\mathfrak{h}_{SC}(F)).$$
  
 Le groupe $G(F)$ agit naturellement sur ces espaces.  On note ${\cal K}_{cusp}^{\natural}(\mathfrak{g}(F))$ et ${\cal K}_{cusp}^{\sharp}(\mathfrak{g}(F))$ les quotients des coinvariants. On a l'\'egalit\'e
 $${\cal K}_{cusp}(\mathfrak{g}(F))={\cal K}_{cusp}^{\natural}(\mathfrak{g}(F))\oplus {\cal K}_{cusp}^{\sharp}(\mathfrak{g}(F)).$$
 
 \begin{lem}{ Soit   $\varphi\in {\cal K}_{cusp}^{\natural}(\mathfrak{g}(F))$. Alors  $k^G(\varphi)$ appartient \`a $D^G({\cal D}_{cusp}^{inst}(\mathfrak{g}(F)))$.
}\end{lem}

Preuve. On peut fixer $H\in {\cal L}^{nr}_{ell}$ et supposer   que $\varphi$ est l'image d'un \'el\'ement de $ FC^{H,inst}(\mathfrak{h}_{SC}(F))$, que l'on note encore $\varphi$. 
Notons $\mathfrak{h}_{SC,reg}(F)/st-conj$ l'ensemble des classes de conjugaison stable dans $\mathfrak{h}_{SC,reg}(F)$. Comme en \ref{expressionintegrale}, c'est une vari\'et\'e analytique sur $F$ que l'on peut munir d'une mesure.  L'application $proj:\mathfrak{h}_{SC,reg}(F)/H-conj\to \mathfrak{h}_{SC,reg}(F)/st-conj$ est un isomorphisme local qui pr\'eserve  localement  les mesures. 
Pour toute $f\in C_{c}^{\infty}(\mathfrak{g}(F))$,  la formule du lemme \ref{expressionintegrale} se r\'ecrit
$$(1) \qquad  \boldsymbol{\Delta}^G(\varphi)(f)=D^G_{X_{H},\varphi}(f)=\int_{\mathfrak{h}_{SC,reg}(F)/st-conj} \Phi(Y)d^G(X_{H}+Y)^{1/2}\, dY,$$
o\`u
$$\Phi(Y)=\sum_{Y'\in proj^{-1}(Y)}I^G(X_{H}+Y',f)e^{H}_{\varphi}(Y').$$
Dans le calcul qui suit, on identifie les \'el\'ements des espaces quotients de $\mathfrak{h}_{SC,reg}(F)$ \`a des repr\'esentants dans ces espaces. Supposons $f\in I_{cusp}^{st}(\mathfrak{g}(F))$. Soit  $Y\in \mathfrak{h}_{SC,reg}(F)/st-conj$. Si $Y\not\in \mathfrak{h}_{SC,tn}(F)$, $e^{H}_{\varphi}(Y')=0$ pour tout $Y'\in proj^{-1}(Y)$. Supposons $Y\in \mathfrak{h}_{SC,tn}(F)$. Puisque $f$ est cuspidale, 
on a $I^G(X_{H}+Y',f)=0$ si $X_{H}+Y'$ n'est pas elliptique, ou encore si $Y'$ n'est pas elliptique, ou encore si $Y$ n'est pas elliptique. Supposons $Y$ elliptique. Les \'el\'ements  $X_{H}+Y'$ pour $Y'\in proj^{-1}(Y)$ sont stablement conjugu\'es, cf. \ref{bonselements}(4). Puisque  $f\in I_{cusp}^{st}(\mathfrak{g}(F))$, les valeurs $I^G(X_{H}+Y',f)$ sont constantes. On note  $I^G(X_{H}+Y,f)$ cette valeur constante. On obtient
$$\Phi(Y)=I^G(X_{H}+Y,f)\sum_{Y'\in proj^{-1}(Y)}e^{H}_{\varphi}(Y').$$
Mais cette somme est nulle parce que $\varphi\in FC^{H,inst}(\mathfrak{h}_{SC}(F))$. Donc $\Phi(Y)=0$. Alors (1) entra\^{\i}ne que $\boldsymbol{ \Delta}^G(\varphi)(f)=0$. Puisque cela est vrai pour toute $f\in I_{cusp}^{st}(\mathfrak{g}(F))$, on a $\boldsymbol{ \Delta}^G(\varphi)\in I(\mathfrak{g}(F))^*[\not={\bf G}]$. 
 La transformation de Fourier pr\'eserve $ I(\mathfrak{g}(F))^*[\not={\bf G}]$. Il en est de m\^eme de la restriction \`a $\mathfrak{g}_{tn}(F)$. Donc $k^G(\varphi)$ appartient \`a $I(\mathfrak{g}(F))^*[\not={\bf G}]$. Cet \'el\'ement appartenant \`a  $D^G({\cal D}_{cusp}(\mathfrak{g}(F)))$, il r\'esulte des d\'efinitions qu'il appartient \`a $D^G({\cal D}_{cusp}^{inst}(\mathfrak{g}(F)))$.  $\square$

\subsubsection{Les espaces ${\cal K}_{cusp}^{st}(\mathfrak{g}(F))$ et ${\cal K}_{cusp}^{\flat}(\mathfrak{g}(F))$\label{lesespacesstflat}}

Notons ${\cal L}^{nr,st}_{ell}/conj$, resp. ${\cal L}^{nr,st}_{ell}/st-conj$, les ensembles de classes d'\'equivalence dans ${\cal L}^{nr,st}_{ell}$ pour la conjugaison par $G(F)$, resp. la conjugaison stable. La premi\`ere \'equivalence \'etant plus fine que la seconde, il y a une application naturelle surjective $\Pi:{\cal L}^{nr,st}_{ell}/conj \to {\cal L}^{nr,st}_{ell}/st-conj$ et la notion   de conjugaison stable a un sens dans l'ensemble ${\cal L}^{nr,st}_{ell}/conj$.  On a dit en \ref{automorphismes} que, pour deux \'el\'ements stablement conjugu\'es $H,H'\in {\cal L}^{nr,st}_{ell}$, il y avait un isomorphisme canonique $\iota_{H',H}:FC^{st}(\mathfrak{h}_{SC}(F))\to FC^{st}(\mathfrak{h}'_{SC}(F))
$. Cela entra\^{\i}ne que l'on peut d\'efinir un espace $FC^{st}(\mathfrak{h}_{SC}(F))$ pour un \'el\'ement $H\in {\cal L}^{nr,st}_{ell}/st-conj$. Pour \^etre bourbakiste,  l'espace en question est la limite inductive des  $FC^{st}(\mathfrak{h}'_{SC}(F))$ sur les $H'\in {\cal L}^{nr,st}_{ell}$ appartenant \`a la classe $H$, les applications de transition \'etant les isomorphismes canoniques ci-dessus.  A fortiori, $FC^{st}(\mathfrak{h}_{SC}(F))$ est d\'efini pour un \'el\'ement $H\in {\cal L}^{nr,st}_{ell}/conj$.

La d\'efinition de l'espace  ${\cal K}_{cusp}^{\sharp}(\mathfrak{g}(F))$ peut se r\'ecrire
$$(1) \qquad {\cal K}_{cusp}^{\sharp}(\mathfrak{g}(F))=\oplus_{H\in {\cal L}^{nr,st}_{ell}/conj }FC^{st}(\mathfrak{h}_{SC}(F)).$$
En effet, il est clair que l'homomorphisme naturel $\boldsymbol{{\cal K}}_{cusp}^{\sharp}(\mathfrak{g}(F))\to {\cal K}_{cusp}^{\sharp}(\mathfrak{g}(F))$ se factorise par l'espace de droite de l'\'egalit\'e ci-dessus. Que cet espace s'envoie injectivement dans l'espace de coinvariants ${\cal K}_{cusp}^{\sharp}(\mathfrak{g}(F))$ r\'esulte du lemme \ref{automorphismes}: l'action de $Norm_{G(F)}(H)$ sur $FC^{st}(\mathfrak{h}_{SC}(F))$ est triviale.

Posons
$${\cal K}_{cusp}^{st}(\mathfrak{g}(F))=\oplus_{H'\in {\cal L}^{nr,st}_{ell}/st-conj }FC^{st}(\mathfrak{h}'_{SC}(F)).$$
 Posons $W^{nr}(H)=H(F^{nr})\backslash Norm_{G(F^{nr})}(H)$ et $W_{F}(H)=H(F)\backslash Norm_{G(F)}(H)$ pour tout $H\in {\cal L}_{F}^{nr}$. On d\'efinit une application lin\'eaire
$$\begin{array}{cccc}\tau:&{\cal K}_{cusp}^{st}(\mathfrak{g}(F))&\to& {\cal K}^{\sharp}_{cusp}(\mathfrak{g}(F))\\ &(\phi_{H'})_{H'\in {\cal L}_{ell}^{nr,st}/st-conj}&\mapsto& (\varphi_{H})_{H\in {\cal L}_{ell}^{nr,st}/conj}\\ \end{array}$$
par la formule $\varphi_{H}=\vert W^{nr}(H)^{\Gamma_{F}^{nr}}\vert \vert W_{F}(H)\vert ^{-1}\phi_{\Pi(H)}$ ($\Pi$ a \'et\'e d\'efinie au d\'ebut du paragraphe). L'application $\tau$ est injective.

On  note ${\cal K}^{\flat}_{cusp}(\mathfrak{g}(F))$ le sous-espace des \'el\'ements  $(\varphi_{H})_{H\in {\cal L}_{ell}^{nr,st}/conj}\in  {\cal K}_{cusp}^{\sharp}(\mathfrak{g}(F))$ tels que, pour tout $H'\in {\cal L}^{nr,st}_{ell}/st-conj$, on ait l'\'egalit\'e $\sum_{H\in \Pi^{-1}(H')}\varphi_{H}=0$. Montrons que

(2) ${\cal K}_{cusp}^{\sharp}(\mathfrak{g}(F))={\cal K}_{cusp}^{\flat}(\mathfrak{g}(F))\oplus \tau({\cal K}_{cusp}^{st}(\mathfrak{g}(F)))$. 

Pour $H'\in {\cal L}^{nr,st}_{ell}/st-conj$, notons ${\cal K}_{cusp}^{\sharp}(\mathfrak{g}(F);H')$ la sous-somme de la formule (1) o\`u on se limite aux $H\in \Pi^{-1}(H')$. Par construction, nos espaces sont sommes directes de leurs intersections avec ces sous-espaces ${\cal K}_{cusp}^{\sharp}(\mathfrak{g}(F);H')$ quand $H'$ d\'ecrit ${\cal L}^{nr,st}_{ell}/st-conj$. Il suffit de fixer $H'$ et de voir ce qui se passe dans ce sous-espace ${\cal K}_{cusp}^{\sharp}(\mathfrak{g}(F);H')$. On a  $ \tau({\cal K}_{cusp}^{st}(\mathfrak{g}(F)))\cap {\cal K}_{cusp}^{\sharp}(\mathfrak{g}(F);H')=\tau(FC^{st}(\mathfrak{h}'_{SC}(F)))$. On pose ${\cal K}^{\flat}_{cusp}(\mathfrak{g}(F);H')={\cal K}_{cusp}^{\flat}(\mathfrak{g}(F))\cap {\cal K}_{cusp}^{\sharp}(\mathfrak{g}(F);H')$.  Pour tout ensemble $Y$, notons ${\mathbb C}[Y]$ l'espace vectoriel complexe de base $Y$. Pour \'eviter les ambigu\"{\i}t\'es, pour $y\in Y$, on note $[y]$ l'\'el\'ement de base correspondant. On note ${\mathbb C}[Y]_{0}$ le sous-espace des $\sum_{y\in Y}a_{y}[y]$ tels que $\sum_{y\in Y}a_{y}=0$. L'espace ${\cal K}_{cusp}^{\sharp}(\mathfrak{g}(F);H')$ s'identifie \`a $FC^{st}(\mathfrak{h}'_{SC}(F))\otimes {\mathbb C}[\Pi^{-1}(H')]$. Introduisons l'\'el\'ement  $\delta$ de ${\mathbb C}[\Pi^{-1}(H')]$ d\'efini par
$$\delta=\sum_{H\in \Pi^{-1}(H')}\vert W^{nr}(H)^{\Gamma_{F}^{nr}}\vert \vert W_{F}(H)\vert ^{-1} [H].$$
L'espace ${\cal K}_{cusp}^{\flat}(\mathfrak{g}(F);H')$  s'identifie \`a $FC^{st}(\mathfrak{h}'_{SC}(F))\otimes {\mathbb C}[\Pi^{-1}(H')]_{0}$ tandis que l'espace $ \tau(FC^{st}(\mathfrak{h}'_{SC}(F))) $ s'identifie \`a $FC^{st}(\mathfrak{h}'_{SC}(F))\otimes {\mathbb C}\delta$. Or les sous-espaces ${\mathbb C}[\Pi^{-1}(H')]_{0}$ et ${\mathbb C}\delta$ de ${\mathbb C}[\Pi^{-1}(H')]$ sont en somme directe. Cela prouve (2). 

\begin{lem}{(i) Soit $\varphi\in {\cal K}^{\flat}_{cusp}(\mathfrak{g}(F))$. Alors $k^G(\varphi)$ appartient \`a $D^G({\cal D}_{cusp}^{inst}(\mathfrak{g}(F)))$.

(ii) Soit $\phi\in {\cal K}^{st}_{cusp}(\mathfrak{g}(F))$. Alors $k^G\circ \tau (\phi)$ appartient \`a $D^G({\cal D}_{cusp}^{st}(\mathfrak{g}(F)))$.}\end{lem}

Preuve. On peut encore fixer $H'\in {\cal L}^{nr,st}_{ell}/st-conj$ et supposer $\varphi\in {\cal K}^{\flat}_{cusp}(\mathfrak{g}(F);H')$ pour l'assertion (i), $\phi\in FC^{st}(\mathfrak{h}'_{SC}(F))$ pour l'assertion (ii). Identifions les classes de conjugaison ou de conjugaison stable dans ${\cal L}^{nr,st}_{ell}$ \`a des repr\'esentants dans cet ensemble. En particulier, on identifie 
  $H'$ \`a un \'el\'ement de ${\cal L}_{ell}^{nr,st}$.  Notons $\sim$ la relation d'\'equivalence "conjugaison par $G(F)$" dans ${\cal L}^{nr,st}_{ell}$. Fixons un \'el\'ement entier et de r\'eduction r\'eguli\`ere $X_{H'}\in \mathfrak{a}_{H'}^{nr}(F)$. La classe de conjugaison stable de $X_{H'}$ est param\'etr\'ee par le groupe $ker^1(H',G)$, cf. \ref{classesdeconjstable}. Pour $\nu\in ker^1(H',G)$, fixons un \'el\'ement $X_{\nu}$ param\'etr\'e par $\nu$ dans cette classe et posons $H_{\nu}=G_{X_{\nu}}$.  L'application qui, \`a $\nu$, associe l'unique \'el\'ement  $H\in {\cal L}^{nr,st}_{ell}/conj$ tel que $H\sim H_{\nu}$ est une surjection de $ker^1(H',G)$ sur $\Pi^{-1}(H')$. On en d\'eduit une surjection
  $$\beta:  {\mathbb C}[ker^1(H',G)]\to    {\mathbb C}[\Pi^{-1}(H')] .$$
 La transformation de Fourier sur $ker^1(H',G)$ fournit un isomorphisme
  $$\alpha:  {\mathbb C}[ker^1(H',G)^{\vee}]\to   {\mathbb C}[ker^1(H',G)].$$
  Pr\'ecis\'ement, pour $\kappa\in ker^1(H',G)^{\vee}$,   $\alpha( [\kappa])=\sum_{\nu\in ker^1(H',G)}\kappa(\nu)  [\nu]$.
  On calcule
  $$(3) \qquad \beta\circ \alpha([\kappa])=\sum_{H\in \Pi^{-1}(H')}(\sum_{\nu\in ker^1(H',G); H_{\nu}\sim H}\kappa(\nu)) [H].$$
  Notons $1^{\vee}$ l'\'el\'ement neutre de $ker^1(H',G)^{\vee}$. 
   On a introduit ci-dessus un \'el\'ement $\delta\in {\mathbb C}[\Pi^{-1}(H')]$.  Montrons que
$$(4) \qquad \beta\circ \alpha([1^{\vee}])=\delta.$$
En vertu de (3) et de la d\'efinition de $\delta$, il s'agit de prouver que, pour tout $H\in \Pi^{-1}(H')$, on a l'\'egalit\'e
$$(5) \qquad \vert \{\nu\in ker^1(H',G); H_{\nu}\sim H\}\vert =\vert W^{nr}(H)^{\Gamma_{F}^{nr}}\vert \vert W_{F}(H)\vert ^{-1} .$$
Pour tout $H''\in {\cal L}^{nr,st}_{ell}$, posons  ${\cal Y}_{H''}=\{g\in G(F^{nr}); gFr(g)^{-1}\in H''(F^{nr})\}$. L'application qui, \`a $g\in {\cal Y}_{H'}$ associe la classe du  cocycle de $\Gamma_{F}^{nr}$ dans $H'(F^{nr})$ qui  envoie $Fr$ sur $gFr(g)^{-1}$ se quotiente en une bijection de $H'(F^{nr})\backslash {\cal Y}_{H'}/G(F)$ sur $ker^1(H',G)$. Notons ${\cal Y}_{H'}[H]$ le sous-ensemble des $g\in {\cal Y}_{H'}$  
  tels que $g^{-1}H'g$ soit conjugu\'e \`a $H$ par un \'el\'ement de $G(F)$. Le membre de gauche de (5) est \'egal au nombre d'\'el\'ements de l'image de ${\cal Y}_{H'}[H]$ dans $H'(F^{nr})\backslash {\cal Y}_{H'}/G(F)$. 
  Puisque $H$ est par d\'efinition stablement conjugu\'e \`a $H'$,  on peut fixer un \'el\'ement $g_{0}\in {\cal Y}_{H'}$ tel que $g_{0}^{-1}H'g_{0}=H$. L'application $g\mapsto g_{0}g$ est une bijection de ${\cal Y}_{H}$ sur ${\cal Y}_{H'}$ qui se quotiente en une bijection de $H(F^{nr})\backslash {\cal Y}_{H}/G(F)$ sur $H'(F^{nr})\backslash {\cal Y}_{H'}/G(F)$. 
   On voit que  
$${\cal Y}_{H'}[H]={\cal Y}_{H'}\cap (Norm_{G(F^{nr})}(H')g_{0}G(F))=g_{0}({\cal Y}_{H}\cap Norm_{G(F^{nr})}(H))G(F).$$
L'image de ${\cal Y}_{H'}[H]$ dans $H'(F^{nr})\backslash {\cal Y}_{H'}/ G(F)$ a donc m\^eme nombre d'\'el\'ements que celle de ${\cal Y}_{H}\cap Norm_{G(F^{nr})}(H)$ dans $H(F^{nr})\backslash {\cal Y}_{H}/ G(F)$. Il y a une surjection naturelle de  ${\cal Y}_{H}\cap Norm_{G(F^{nr})}(H)$  sur $W^{nr}(H)^{\Gamma_{F}^{nr}}$ et on voit que deux \'el\'ements de ${\cal Y}_{H}\cap Norm_{G(F^{nr})}(H)$ ont m\^eme image dans $H(F^{nr})\backslash {\cal Y}_{H}/ G(F)$ si et seulement si leurs images dans $W^{nr}(H)^{\Gamma_{F}^{nr}}$  ont m\^eme image dans $W^{nr}(H)^{\Gamma_{F}^{nr}}/W_{F}(H)$. D'o\`u les \'egalit\'es (5) puis (4). 

 De $\beta$, resp. $\alpha$, se d\'eduisent une surjection, resp.  un isomorphisme
  $$b:FC^{st}(\mathfrak{h}'_{SC}(F))\otimes {\mathbb C}[ker^1(H',G)]\to  FC^{st}(\mathfrak{h}'_{SC}(F))\otimes {\mathbb C}[\Pi^{-1}(H')]={\cal K}_{cusp}^{\sharp}(\mathfrak{g}(F);H').$$
  $$a:FC^{st}(\mathfrak{h}'_{SC}(F))\otimes {\mathbb C}[ker^1(H',G)^{\vee}]\to FC^{st}(\mathfrak{h}'_{SC}(F))\otimes {\mathbb C}[ker^1(H',G)].$$
  D'apr\`es (4) et le calcul de la preuve de (2), on a $b\circ a(FC^{st}(\mathfrak{h}'_{SC}(F))\otimes {\mathbb C}[1^{\vee}])=\tau(FC^{st}(\mathfrak{h}'_{SC}(F)))$. 
 Notons ${\mathbb C}[ker^1(H',G)^{\vee}]^{\flat}$ le sous-espace engendr\'e par les $[\kappa]$ pour $\kappa\not=1^{\vee}$. On voit gr\^ace \`a (3) que $b\circ a( FC^{st}(\mathfrak{h}'_{SC}(F))\otimes {\mathbb C}[ker^1(H',G)^{\vee}]^{\flat})$ est contenu dans ${\cal K}_{cusp}^{\flat}(\mathfrak{g}(F);H')$. Puisque $b\circ a$ est surjectif, (2) implique que ces deux espaces sont \'egaux. Alors, le lemme est \'equivalent \`a l'assertion suivante. Soient $\phi\in FC^{st}(\mathfrak{h}'_{SC}(F))$ et $\kappa\in ker^1(H',G)$. Alors
 
 (6) $k^G\circ b\circ a(\phi\otimes [1^{\vee}])$ appartient \`a  $D^G({\cal D}_{cusp}^{st}(\mathfrak{g}(F)))$; si $\kappa\not=1^{\vee}$, $k^G\circ b\circ a(\phi\otimes [\kappa])$ appartient \`a  $D^G({\cal D}_{cusp}^{inst}(\mathfrak{g}(F)))$.
  
 Pour $\nu\in ker^1(H',G)$, notons $\varphi_{\nu}\in FC^{st}(\mathfrak{h}_{\nu,SC}(F))$ l'image de $\phi$ par l'application canonique $\iota_{H_{\nu},H'}$. Introduisons la distribution
 $$\Theta_{\kappa}=\sum_{\nu\in ker^1(H',G)}\kappa(\nu) D_{X_{\nu},\varphi_{\nu}}^G.$$
 Il r\'esulte des d\'efinitions que $k^G\circ b\circ a(\phi\otimes [\kappa])$ est \'egal \`a la restriction \`a $\mathfrak{g}_{tn}(F)$ de la transform\'ee de Fourier de $\Theta_{\kappa}$. Soit $f\in I_{cusp}(\mathfrak{g}(F))$.  Le m\^eme calcul qu'en \ref{instabilite1} conduit \`a l'\'egalit\'e
$$ \Theta_{\kappa}(f)=\sum_{\nu\in ker^1(H',G)}\kappa(\nu)\int_{\mathfrak{h}_{\nu,SC,reg}(F)/st-conj}\Phi_{\nu}(Y)d^G(X_{\nu}+Y)^{1/2}\, dY,$$
o\`u
$$\Phi_{\nu}(Y)=\sum_{Y'\in proj_{\nu}^{-1}(Y)}I^G(X+Y',f)e_{ \varphi_{\nu}}^{H_{\nu}}(Y').$$
On a ajout\'e un indice $\nu$ \`a $proj$ pour plus de pr\'ecision. Puisque les $ \varphi_{\nu}$ sont stables, les fonctions $e_{ \varphi_{\nu}}^{H_{\nu,SC}}$ sont invariantes par conjugaison stable. Il en est de m\^eme des fonctions $e_{ \varphi_{\nu}}^{H_{\nu}}=m(H_{\nu})e_{ \varphi_{\nu}}^{H_{\nu,SC}}$. 
 On note $e_{ \varphi_{\nu}}^{H_{\nu}}(Y)$ la valeur de $e_{\varphi_{\nu}}^{H_{\nu}}(Y')$ pour $Y'\in proj_{\nu}^{-1}(Y)$. On obtient
$$\Phi_{\nu}(Y)=e_{\varphi_{\nu}}^{H_{\nu}}(Y)\sum_{Y'\in proj_{\nu}^{-1}(Y)}I^G(X+Y',f).$$
Puisque $f$ est cuspidale et que $e_{\varphi_{\nu}}^{H_{\nu}}$ est \`a support topologiquement nilpotent, on peut limiter l'int\'egration \`a $\mathfrak{h}_{\nu,SC,ell,tn}(F)/st-conj$ avec une d\'efinition \'evidente de cet ensemble.  La bijection canonique $\iota_{ H_{\nu},H'}$ a \'et\'e d\'efinie \`a l'aide d'un torseur int\'erieur de $H'$ vers $H_{\nu}$ dont on peut supposer qu'il envoie $X_{H'}$ sur $X_{\nu}$. On d\'eduit de ce torseur une bijection $\iota'_{H_{\nu},H'}:\mathfrak{h}'_{SC,ell}(F)/st-conj\to \mathfrak{h}_{\nu,SC,ell}(F)/st-conj$ . Par d\'efinition du transfert, on a $e_{ \varphi_{\nu}}^{H_{\nu,SC}}(\iota'_{H_{\nu},H'}(Y))=e_{\phi}^{H'_{SC}}(Y)$ pour tout $Y\in \mathfrak{h}'_{SC,ell}(F)/st-conj$. Cela entra\^{\i}ne $e_{\varphi_{\nu}}^{H_{\nu}}(\iota'_{H_{\nu},H'}(Y))=e_{\phi}^{H'}(Y)$  car $m(H_{\nu})$ est ind\'ependant de $\nu$ d'apr\`es le lemme \ref{uncalculdemesures}.  
On obtient donc
$$(7) \qquad  \Theta_{\kappa}(f)=\int_{\mathfrak{h}'_{SC,ell,tn}(F)/st-conj}e_{\phi}^{H'}(Y)d^G(X_{H'}+Y)^{1/2}I^{G,\kappa}(X_{H'}+Y,f)\,dY,$$
o\`u
$$I^{G,\kappa}(X_{H'}+Y,f)=\sum_{\nu\in  ker^1(H',G)}\kappa(\nu)\sum_{Y'\in proj_{\nu}^{-1}(\iota'_{H_{\nu},H'}(Y))}I^G(X_{\nu}+Y',f).$$
Fixons $Y\in \mathfrak{h}'_{SC,ell,tn}(F)/st-conj$ et un rel\`evement $Y_{0}$ dans $\mathfrak{h}'_{SC,ell}(F)\cap \mathfrak{h}'_{SC,tn}(F)$. Notons $T_{0}$ le centralisateur de $X_{H'}+Y_{0}$ dans $G$. On sait que l'ensemble $Cl_{st}(X_{H'}+Y_{0})/conj$ des classes de conjugaison par $G(F)$ dans la classe de conjugaison stable de $X_{H'}+Y_{0}$ s'identifie au noyau $ker^1(T_{0},G)$ de l'application $H^1(F,T_{0})\to H^1(F,G)$. En utilisant \ref{bonselements}(4), on voit que cet ensemble de classes de conjugaison s'identifie aussi \`a
$$\cup_{\nu\in ker^1(H',G)}proj_{\nu}^{-1}(\iota'_{H_{\nu},H'}(Y))$$
et que l'application naturelle de cet ensemble dans $ ker^1(H',G)$ s'identifie \`a l'application naturelle de $ker^1(T_{0},G) $ dans ce groupe. Ces applications sont surjectives. Le caract\`ere $\kappa$ se rel\`eve en un caract\`ere $\kappa_{0}$ de $ker^1(T_{0},G)$ et alors $I^{G,\kappa}(X_{H'}+Y,f)$ est \'egale \`a la $\kappa_{0}$-int\'egrale orbitale habituelle $I^{G,\kappa_{0}}(X_{H}+Y_{0},f)$. Si $\kappa\not=1^{\vee}$ et $f\in I^{st}_{cusp}(\mathfrak{g}(F))$, cette   $\kappa_{0}$-int\'egrale orbitale est nulle donc $ \Theta_{\kappa}(f)=0$ d'apr\`es (7). De m\^eme, si $\kappa=1^{\vee}$  (donc $\kappa_{0}=1$) et $f\in I_{cusp}^{inst}(\mathfrak{g}(F))$, cette   $\kappa_{0}$-int\'egrale orbitale est nulle donc  $\Theta_{1^{\vee}}(f)=0$.  Il en r\'esulte que

$ \Theta_{\kappa}$ appartient \`a $I(\mathfrak{g}(F))^*[\not={\bf G}]$ si $\kappa\not=1^{\vee}$;

$ \Theta_{1^{\vee}}$ appartient \`a $I(\mathfrak{g}(F))^*[{\bf G}]$.  

On en d\'eduit (6) comme dans la fin de la preuve de \ref{instabilite1}. Cela ach\`eve la d\'emonstration. $\square$

\subsubsection{Description de l'espace ${\cal D}^{st}_{cusp}(\mathfrak{g}(F))$\label{descriptiondeDstable}}

Notons $k^{G,st}$ l'application   $k^G\circ \tau :{\cal K}^{st}_{cusp}(\mathfrak{g}(F))\to I(\mathfrak{g}(F))^*$.

\begin{prop}{L'application $k^{G,st}$ est un isomorphisme antilin\'eaire de ${\cal K}^{st}_{cusp}(\mathfrak{g}(F))$ sur $D^G({\cal D}_{cusp}^{st}(\mathfrak{g}(F)))$.}\end{prop}

Preuve. D'apr\`es \ref{deuxiemedescription}(1), $k^G$ est un isomorphisme de ${\cal K}_{cusp}(\mathfrak{g}(F))$ sur $D^G({\cal D}_{cusp}(\mathfrak{g}(F)))$. D'apr\`es les d\'efinitions et \ref{lesespacesstflat}(2), on a l'\'egalit\'e
 $$(1) \qquad {\cal K}_{cusp}(\mathfrak{g}(F))={\cal K}^{\natural}_{cusp}(\mathfrak{g}(F))\oplus  {\cal K}^{\flat}_{cusp}(\mathfrak{g}(F))\oplus \tau({\cal K}^{st}_{cusp}(\mathfrak{g}(F))).$$
 On a aussi
$$(2) \qquad D^G({\cal D}_{cusp}(\mathfrak{g}(F)))=D^G({\cal D}_{cusp}^{inst}(\mathfrak{g}(F)))\oplus D^G({\cal D}^{st}_{cusp}(\mathfrak{g}(F))).$$
Les lemmes \ref{instabilite1} et \ref{lesespacesstflat} disent que les deux premiers facteurs de (1) s'envoient dans le premier facteur de (2) tandis que le troisi\`eme facteur de (1) s'envoie dans le second facteur de (2). Donc $k^G$ se restreint en un isomorphisme du troisi\`eme facteur de (1) sur le second facteur de (2).  Puisque $\tau$ est injective, cela \'equivaut \`a l'\'enonc\'e. $\square$.

 \subsubsection{Description  des distributions associ\'ees aux \'el\'ements de ${\cal D}(\mathfrak{g}(F))$\label{descriptioncasgeneral}}
 
 Soit $L $ un $F$-Levi de $G$, consid\'erons des donn\'ees $(H,s_{H})$ et $\varphi$ comme en \ref{descriptiondistributions}, le groupe ambiant $G$ \'etant remplac\'e par $L$. C'est-\`a-dire $(H,s_{H})\in {\cal L}^{L,nr,\star}_{ell}$ et $\varphi\in FC^{H_{s_{H}}}(\mathfrak{h}_{SC,s_{H}}({\mathbb F}_{q}))$.  On a associ\'e \`a ces donn\'ees un sommet $s\in Imm(L_{AD})$ et une fonction ${\cal Q}_{\varphi}\in C_{nil,cusp}(\mathfrak{l}_{s}({\mathbb F}_{q}))$. On en d\'eduit une distribution $D^L_{{\cal Q}_{\varphi}}\in I(\mathfrak{l}(F))^*$ puis, par induction, une distribution $D^G_{{\cal Q}_{\varphi}}\in I(\mathfrak{g}(F))^*$. Fixons un \'el\'ement $X_{H}\in A_{H}^{nr}(F)$ entier et de r\'eduction r\'eguli\`ere. Comme  en \ref{uncalculdemesures}, on associe \`a $\varphi$ une distribution $D^H_{\varphi}\in I(\mathfrak{h}_{SC}(F))^*$.  Pour $f\in C_{c}^{\infty}(\mathfrak{g}(F))$, notons $f_{\vert X_{H}+\mathfrak{h}_{SC}}$ la fonction sur $\mathfrak{h}_{SC}(F)$ d\'efinie par $Y\mapsto f(X_{H}+Y)$.  On pose
$$D^G_{X_{H},\varphi}(f)=\int_{H(F)\backslash G(F)}D_{\varphi}^{H}((^gf)_{\vert X_{H}+  \mathfrak{h}_{SC}})\,dg.$$
L'int\'egrale est \`a support compact d'apr\`es \ref{bonselements}(5). Posons
$$c= q^{-dim(A_{H}^{nr})/2}mes(K_{s}^L)mes(K_{s_{H}}^H)^{-1}.$$

 \begin{prop}{(i) On a l'\'egalit\'e $D^G_{X_{H},\varphi}=Ind_{L}^G(D^L_{X_{H},\varphi})$.
 
 (ii) La distribution $D^G_{{\cal Q}_{\varphi}^-}$ est \'egale \`a la restriction \`a $\mathfrak{g}_{tn}(F)$ de    $c\hat{D}^G_{X_{H},\varphi} $.}\end{prop}

Preuve.  Posons $\Theta=Ind_{L}^G(D^L_{X_{H},\varphi})$. 
Fixons un sous-groupe parabolique $P$ de $G$ de composante de Levi $L$ et d\'efini sur $F$. 
Pour $f\in C_{c}^{\infty}(\mathfrak{g}(F))$, on note $f_{U_{P}}$ la fonction sur $\mathfrak{l}(F)$ d\'efinie par
$$f_{U_{P}}(Y)=\int_{\mathfrak{u}_{P}(F)}f(Y+N)\,dN.$$
 Par d\'efinition,
 $$  \Theta(f)=\int_{P(F)\backslash G(F)}D^L_{ X_{H},\varphi}((^gf)_{U_{P}})\,dg$$
 $$=\int_{P(F)\backslash G(F)}\int_{H(F)\backslash L(F)}D^H_{\varphi}(f[g,l])\,dl\,dg,$$
 o\`u
 $f[g,l]=(^l((^gf)_{U_{P}}))_{\vert X_{H}+\mathfrak{h}_{SC}}$. Calculons cette fonction $f[g,l]$. Puisque la distribution $D^H_{\varphi}$ est \`a support dans $\mathfrak{h}_{SC,tn}(F)$, on peut se limiter \`a la calculer sur cet ensemble. 
 Fixons $g$, $l$ et $Y\in \mathfrak{h}_{SC,tn}(F)$.  Alors
 $$f[g,l](Y)=\int_{\mathfrak{u}_{P}(F)}f(g^{-1}(l^{-1}(X_{H}+Y)l+N)g)\,dN.$$
 Par le changement de variables $N\mapsto l^{-1}Nl$, on obtient
 $$f[g,l](Y)=\delta_{P}(l)^{-1}\int_{\mathfrak{u}_{P}(F)}f(g^{-1}l^{-1}(X_{H}+Y+N)lg)\,dN.$$
 Parce $X_{H}+Y$ est un \'el\'ement de $L(F)$ dont le commutant est  contenu dans $H$ d'apr\`es \ref{bonselements}(1), donc dans $ L$, l'application $u\mapsto u^{-1}(X_{H}+Y)u-(X_{H}+Y)$ est un isomorphisme de $U_{P}(F)$ sur $\mathfrak{u}_{P}(F)$.  Son jacobien est \'egal \`a $d^G(X_{H}+Y)^{1/2}d^L(X_{H}+Y)^{-1/2}$ qui vaut $1$ d'apr\`es \ref{bonselements}(3). D'o\`u
 $$f[g,l](Y)=\delta_{P}(l)^{-1}\int_{U_{P}(F)}f(g^{-1}l^{-1}u^{-1}(X_{H}+Y)ulg)\,du=\delta_{P}(l)^{-1}\int_{U_{P}(F)}(^{ulg}f)_{\vert X_{H}+\mathfrak{h}_{SC}}(Y)\,du.$$
 La propri\'et\'e \ref{bonselements}(5) entra\^{\i}ne qu'il existe un sous-ensemble compact $\omega\subset U_{P}(F)$ tel que $ (^{ulg}f)_{\vert X_{H}+\mathfrak{h}_{SC}}(Y)=0$ pour tout $Y\in \mathfrak{h}_{SC,tn}(F)$ si $u\not\in \omega$. Il en r\'esulte que
 $$D^H_{\varphi}(f[g,l])=\delta_{P}(l)^{-1}\int_{U_{P}(F)}D^H_{\varphi}((^{ulg}f)_{\vert X_{H}+\mathfrak{h}_{SC}})\,du.$$
 On a alors
 $$\Theta(f)=\int_{P(F)\backslash G(F)}\int_{H(F)\backslash L(F)}\int_{U_{P}(F)}\delta_{P}(l)^{-1}D^H_{\varphi}((^{ulg}f)_{\vert X_{H}+\mathfrak{h}_{SC}})\,du\,dl\,dg.$$
 Cette int\'egrale est absolument convergente toujours d'apr\`es \ref{bonselements}(5). Elle se
  reconstitue alors  en
 $$\Theta(f)=\int_{H(F)\backslash G(F)}D_{\varphi}^{H}((^gf)_{\vert X_{H}+  \mathfrak{h}_{SC}})\,dg,$$
  c'est-\`a-dire $\Theta(f)=D^G_{X_{H},\varphi}(f)$. Cela d\'emontre le (i) de l'\'enonc\'e.
  
   On applique \`a $D^L_{{\cal Q}_{\varphi}^-}$ la proposition \ref{descriptiondistributions}: cette distribution est la restriction \`a $\mathfrak{l}_{tn}(F)$ de $c\hat{D}^L_{X_{H},\varphi}$. Par induction, $D^G_{{\cal Q}_{\varphi}^-}$ est la restriction \`a $\mathfrak{g}_{tn}(F)$ de $c\hat{\Theta}$. Le (ii) de l'\'enonc\'e r\'esulte alors de (i). $\square$

\subsubsection{L'espace ${\cal K}^{st}(\mathfrak{g}(F))$\label{lespaceKstablegeneral}}
Notons ${\cal L}_{F}^{nr,st}/conj$, resp. ${\cal L}_{F}^{nr,st}/st-conj$, les ensembles de classes de conjugaison par $G(F)$, resp. de conjugaison stable, dans ${\cal L}_{F}^{nr,st}$.  Comme en \ref{lesespacesstflat}, il y a une surjection naturelle $\Pi:{\cal L}_{F}^{nr,st}/conj \to {\cal L}_{F}^{nr,st}/st-conj$. 

Posons
$${\cal K}^{\sharp}(\mathfrak{g}(F))=\oplus_{H\in {\cal L}_{F}^{nr,st}/conj}FC^{st}(\mathfrak{h}_{SC}(F)).$$
 Soit $H\in {\cal L}_{F}^{nr,st}$. Fixons un \'el\'ement $X_{H}\in \mathfrak{a}^{nr}_{H}(F)$ entier et de r\'eduction r\'eguli\`ere. Pour $\varphi\in FC^{st}(\mathfrak{h}_{SC}(F))$, notons $k^G(\varphi)$ la restriction \`a $\mathfrak{g}_{tn}(F)$ de la transform\'ee de Fourier de $D^G_{X_{H},\varphi}$. D'apr\`es la proposition \ref{descriptioncasgeneral}, cette distribution ne d\'epend pas du choix de $X_{H}$. Si $H'\in {\cal L}_{F}^{nr,st}$ est conjugu\'e \`a $H$ par un \'el\'ement $g\in G(F)$, c'est-\`a-dire $H'=g^{-1}Hg$, l'application $Ad(g)^{-1}$ induit un isomorphisme de $FC^{st}(\mathfrak{h}_{SC}(F))$ sur $FC^{st}(\mathfrak{h}'_{SC}(F))$ qui n'est autre que l'application canonique $\iota_{H',H}$ et il est clair par simple transport de structure que $k^G_{\iota_{H',H}(\varphi)}=k^G_{\varphi}$. Autrement dit, l'application $k^G$ d\'efinie sur $FC^{st}(\mathfrak{h}_{SC}(F))$ ne d\'epend que de la classe de conjugaison de $H$ par $G(F)$. En faisant varier cette classe de conjugaison, on obtient une application antilin\'eaire
 $$k^G:{\cal K}^{\sharp}(\mathfrak{g}(F))\to I(\mathfrak{g}(F))^*.$$
 
  Posons 
 $${\cal K}^{st}(\mathfrak{g}(F))=\oplus_{H'\in {\cal L}_{F}^{nr,st}/st-conj}FC^{st}(\mathfrak{h}'_{SC}(F)).$$
 En imitant la d\'efinition de \ref{lesespacesstflat}, on d\'efinit une application lin\'eaire
 $$\begin{array}{cccc}\tau:&{\cal K}^{st}(\mathfrak{g}(F))&\to& {\cal K}^{\sharp}(\mathfrak{g}(F))\\ &(\phi_{H'})_{H'\in {\cal L}^{nr,st}/st-conj}&\mapsto& (\varphi_{H})_{H\in {\cal L}^{nr,st}/conj}\\ \end{array}$$
par la formule $\varphi_{H}=\vert W^{nr}(H)^{\Gamma_{F}^{nr}}\vert \vert W_{F}(H)\vert ^{-1}\phi_{\Pi(H)}$. L'application $\tau$ est injective.

 On d\'efinit $k^{G,st}:{\cal K}^{st}(\mathfrak{g}(F))\to I(\mathfrak{g}(F))^*$ par $ k^{G,st}=k^G\circ \tau$.
 
 Soit $H\in {\cal L}_{F}^{nr,st}/st-conj$ que l'on identifie \`a un repr\'esentant dans ${\cal L}_{F}^{nr,st}$ et soit $\varphi\in FC^{st}(\mathfrak{h}_{SC}(F))$. Les classes de conjugaison par $G(F)$ dans la classe de conjugaison stable de $X_{H}$ sont param\'etr\'ees par 
  $ker^1(H,G)$. Pour $\nu\in ker^1(H,G)$, fixons un \'el\'ement $X_{\nu}$ stablement conjugu\'e \`a $X$, dont la classe de conjugaison par $G(F)$ est param\'etr\'ee par $\nu$. Posons $H_{\nu}=G_{X_{\nu}}$ et $\varphi_{\nu}=\iota_{H_{\nu},H}(\varphi)$. Le m\^eme calcul qu'en \ref{lesespacesstflat}(5) montre que  
   
 (1) $k^{G,st}(\varphi)$ est la restriction \`a $\mathfrak{g}_{tn}(F)$ de la transform\'ee de Fourier de la distribution 
 $$\sum_{\nu\in ker^1(H,G)}D^G_{X_{\nu},\varphi_{\nu}}.$$
 
 \subsubsection{L'espace ${\cal K}^{st}(\mathfrak{g}(F))$ et les Levi de $G$\label{lespaceKstableetlesLevi}}
 Soit $L$ un $F$-Levi de $G$. Il y a une application naturelle $p_{L}:{\cal L}^{L,nr,st}_{ell}/st-conj\to {\cal L}_{F}^{nr,st}/st-conj$. On prendra garde que la conjugaison stable dans le premier ensemble est relative au groupe ambiant $L$ tandis que celle dans le second ensemble est relative au groupe ambiant $G$. Il s'en d\'eduit une application lin\'eaire
 $$j_{L}:{\cal K}^{st}_{cusp}(\mathfrak{l}(F))\to {\cal K}^{st}(\mathfrak{g}(F)).$$ 
  Le groupe $Norm_{G(F)}(L)$ agit naturellement dans ${\cal L}^{L,nr,st}_{ell}$ et on v\'erifie que cette action pr\'eserve  la  conjugaison stable pour le groupe ambiant $L$ (pour $n\in Norm_{G(F)}(L)$ et $X,X'\in \mathfrak{l}_{reg}(F)$, $X$ et $X'$ sont stablement conjugu\'es dans $L$ si et seulement si $n(X)$ et $n(X')$ le sont).  Il s'en d\'eduit une action sur ${\cal L}^{L,nr,st}_{ell}/st-conj$. De plus,  pour  $n\in Norm_{G(F)}(L)$ et $H\in {\cal L}^{L,nr,st}_{ell}$,  les groupes $H$ et $H'=nHn^{-1}$ sont stablement conjugu\'es pour le groupe ambiant $G$ et les espaces $FC^{st}(\mathfrak{h}_{SC}(F))$ et $FC^{st}(\mathfrak{h}'_{SC}(F))$ sont canoniquement isomorphes. De l'action de $Norm_{G(F)}(L)$ dans ${\cal L}^{L,nr,st}_{ell}$ se d\'eduit donc une action lin\'eaire de ce groupe dans l'espace ${\cal K}^{st}_{cusp}(\mathfrak{l}(F))$. Notons ${\cal K}^{st}_{cusp}(\mathfrak{l}(F))/Norm_{G(F)}(L)$ le quotient des coinvariants.

 \begin{lem}{(i) Pour $\varphi\in {\cal K}^{st}_{cusp}(\mathfrak{l}(F))$, on a l'\'egalit\'e $k^{G,st}(j_{L}(\varphi))=Ind_{L}^G(k^{L,st}(\varphi))$.

(ii)  L'application $j_{L}$ se quotiente en une application injective 
$$\bar{j}_{L}:{\cal K}^{st}_{cusp}(\mathfrak{l}(F))/Norm_{G(F)}(L) \to {\cal K}^{st}(\mathfrak{g}(F)).$$}\end{lem}
 
 Preuve. Pour (i), on 
  peut fixer $H\in {\cal L}_{ell}^{L,nr,st}/st-conj $, que l'on identifie \`a un repr\'esentant dans ${\cal L}^{L,nr,st}_{ell}$, et supposer $\varphi\in FC^{st}(\mathfrak{h}_{SC}(F))$.  Utilisons la formule \ref{lespaceKstablegeneral}(1) et son analogue dans $L$. On obtient
  que $k^{G,st}(j_{L}(\varphi))$, resp. $Ind_{L}^G(k^{L,st}(\varphi))$, est la  restriction \`a $\mathfrak{g}_{tn}(F)$ de la transform\'ee de Fourier de la distribution 
 $\sum_{\nu\in ker^1(H,G)}D^G_{X_{\nu},\varphi_{\nu}}$, resp. $\sum_{\nu'\in ker^1(H,L)}Ind_{L}^G(D^L_{X_{\nu'},\varphi_{\nu'}})$. D'apr\`es la proposition \ref{descriptioncasgeneral}, cette derni\`ere distribution est \'egale \`a $\sum_{\nu'\in ker^1(H,L)}D^G_{X_{\nu'},\varphi_{\nu'}}$. Il y a une application naturelle $ker^1(H,L)\to ker^1(H,G)$. Celle-ci   est bijective car, $L$ \'etant un $F$-Levi de $G$, l'application $H^1(F,L)\to H^1(F,G)$ est injective. Pour $\nu'\in ker^1(H,L)$ s'identifiant \`a un \'el\'ement $\nu\in ker^1(H,G)$, on peut supposer que les \'el\'ements $X_{\nu'}$ et $X_{\nu}$ intervenant ci-dessus sont \'egaux. Alors les deux distributions ci-dessus sont \'egales, ce qui prouve (i).

 Notons $^L{\cal L}_{F}^{nr,st}/st-conj$ l'image de l'application $p_{L}$ d\'efinie si-dessus. 
 Pour $H\in{^L{\cal L}}_{F}^{nr,st}/st-conj$, notons $ {\cal K}^{st}_{cusp}(\mathfrak{l}(F))[H]$ la somme des $FC^{st}(\mathfrak{h}'_{SC}(F))$ sur les $H'\in p_{L}^{-1}(H)$. On a
 $${\cal K}^{st}_{cusp}(\mathfrak{l}(F))=\oplus_{H\in {^L{\cal L}}_{F}^{nr,st}/st-conj}{\cal K}^{st}_{cusp}(\mathfrak{l}(F))[H].$$
 L'application $j_{L}$ est la somme directe sur $H\in  {^L{\cal L}}_{F}^{nr,st}/st-conj$ de ses restrictions
 $$j_{L}[H]:{\cal K}^{st}_{cusp}(\mathfrak{l}(F))[H]\to FC^{st}(\mathfrak{h}_{SC}(F)).$$
 L'action de $Norm_{G(F)}(L)$ conserve les sous-espaces ${\cal K}^{st}_{cusp}(\mathfrak{l}(F))[H]$. Il suffit donc de prouver les analogues  de (ii) pour les applications $j_{L}[H]$, pour $H\in {^L{\cal L}}_{H}^{nr,st}/st-conj$. Fixons un tel $H$. Pour $H'\in p_{L}^{-1}(H)$, l'espace $FC^{st}(\mathfrak{h}'_{SC}(F))$ s'identifie canoniquement \`a $FC^{st}(\mathfrak{h}_{SC}(F))$. Alors $ {\cal K}^{st}_{cusp}(\mathfrak{l}(F))[H]$ s'identifie \`a $  FC^{st}(\mathfrak{h}_{SC}(F))\otimes {\mathbb C}[p_{L}^{-1}(H)]$.D\'efinissons une forme lin\'eaire 
  $$\begin{array}{cccc}\bar{j}_{L}[H]:&{\mathbb C}[p_{L}^{-1}(H)]&\to &{\mathbb C}\\ &\sum_{H'\in p_{L}^{-1}(H)}c_{H'}[H']&\mapsto&\sum_{H'\in p_{L}^{-1}(H)}c_{H'}.\\ \end{array}$$ 
  Alors 
   $j_{L}[H]$ s'identifie  au produit tensoriel de l'identit\'e de $FC^{st}(\mathfrak{h}_{SC}(F))$ et de cette forme lin\'eaire $\bar{j}_{L}[H]$. 
  L'action de $Norm_{G(F)}(L)$ sur ${\cal K}^{st}_{cusp}(\mathfrak{l}(F))[H]   $ est le produit tensoriel de l'identit\'e de $FC^{st}(\mathfrak{h}_{SC}(F))$ et de l'action sur  $ {\mathbb C}[p_{L}^{-1}(H)]$ d\'eduite de l'action par permutations de $Norm_{G(F)}(L)$ sur $p_{L}^{-1}(H)$. Cela nous ram\`ene \`a prouver l'analogue de (ii)   pour la forme lin\'eaire  $\bar{j}_{L}[H]$. Que celle-ci se quotiente par l'action de $Norm_{G(F)}(L)$ est \'evident: pour $\sum_{H'\in p_{L}^{-1}(H)}c_{H'}[H']\in {\mathbb C}[p_{L}^{-1}(H)]$, $\sum_{H'\in p_{L}^{-1}(H)}c_{H'}$ est insensible \`a toute permutation de $p_{L}^{-1}(H)$. Puisque $Norm_{G(F)}(L)$ permute la base $p_{L}^{-1}(H)$ de ${\mathbb C}[p_{L}^{-1}(H)]$, l'espace des coinvariants ${\mathbb C}[p_{L}^{-1}(H)]/ Norm_{G(F)}(L)$ s'identifie \`a la somme de droites index\'ees par les orbites de l'action de  $Norm_{G(F)}(L)$ dans $p_{L}^{-1}(H)$. La deuxi\`eme assertion de (ii) r\'esulte alors de 
 
 (1) $Norm_{G(F)}(L)$ agit transitivement dans $p_{L}^{-1}(H)$.
 
 Soit $H',H''\in p_{L}^{-1}(H)$, que l'on identifie \`a des repr\'esentants dans ${\cal L}^{L,nr,st}_{ell}$.  Fixons $g\in G(F^{nr})$ tel que $H''=g^{-1}H'g$ et $gFr(g)^{-1}\in H'(F^{nr})$. L'application $Ad(g)^{-1}$ se restreint en un isomorphisme $\Gamma_{F}$-\'equivariant de $Z(H')^0$ sur $Z(H'')^0$, a fortiori de $A_{H'}$ sur $A_{H''}$. Puisque $H'$ et $H''$ sont elliptiques dans $L$, ces tores $A_{H'}$ et $A_{H''}$ sont tous deux \'egaux \`a $A_{L}$. Donc $Ad(g)^{-1}$ conserve $A_{L}$ et conserve aussi son commutant $L$. Autrement dit $g\in Norm_{G(F^{nr})}(L)$. Puisque $gFr(g)^{-1}\in H'(F^{nr})\subset L(F^{nr})$, l'image de $g$ dans $L(F^{nr})\backslash Norm_{G(F^{nr})}(L)$ est fix\'ee par $\Gamma_{F}^{nr}$. Puisque $L$ est un $F$-Levi de $G$, l'application naturelle
 $$L(F)\backslash Norm_{G(F)}(L)\to (L(F^{nr})\backslash Norm_{G(F^{nr})}(L))^{\Gamma_{F^{nr}}}$$
 est surjective. On peut donc fixer $n\in Norm_{G(F)}(L)$ et $l\in L(F^{nr})$ tels que $g=ln$. La relation $gFr(g)^{-1}\in H'(F^{nr})$ entra\^{\i}ne $lFr(l)^{-1}\in H'(F^{nr})$. Posons $H'_{0}=l^{-1}H'l$. Alors $H'_{0}$ est un \'el\'ement de ${\cal L}^{L,nr,st}_{ell}$
 qui est stablement conjugu\'e dans $L$ \`a $H'$, autrement dit les images de $H'$ et $H'_{0}$ dans ${\cal L}^{L,nr,st}_{ell}/st-conj$ sont \'egales. L'\'egalit\'e $H''=n^{-1}H'_{0}n$ montre que l'image de $H''$ dans ${\cal L}^{L,nr,st}_{ell}/st-conj$ se d\'eduit de celle de $H$ par l'action de $n^{-1}\in Norm_{G(F)}(L)$. Cela prouve (1) et le lemme. $\square$
 
 \subsubsection{Description de ${\cal D}^{st}(\mathfrak{g}(F))$\label{descriptionfinale}}
 
\begin{prop}{L'application $k^{G,st}$ est un isomorphisme antilin\'eaire de ${\cal K}^{st}(\mathfrak{g}(F))$ sur $D^G({\cal D}^{st}(\mathfrak{g}(F))$.}\end{prop}

Preuve. Fixons un ensemble $\underline{{\cal L}}_{F}$ de repr\'esentants des classes de conjugaison par $G(F)$ dans l'ensemble des $F$-Levi de $G$. Des applications $\bar{j}_{L}$ d\'efinies en \ref{lespaceKstableetlesLevi} se d\'eduit une application lin\'eaire
$$\bar{j}:\oplus_{L\in \underline{{\cal L}}_{F}}{\cal K}^{st}_{cusp}(\mathfrak{l}(F))/Norm_{G(F)}(L)\to {\cal K}^{st}(\mathfrak{g}(F)).$$
On a

(1) $\bar{j}$ est un isomorphisme.

Puisque les applications $\bar{j}_{L}$ sont injectives pour tout $L$, il suffit de montrer que ${\cal K}^{st}(\mathfrak{g}(F))$ est la somme directe de leurs images. D'apr\`es la d\'efinition de ces applications, il suffit de prouver que, pour tout $H\in {\cal L}_{F}^{nr,st}$, il existe un et un seul $L\in \underline{{\cal L}}_{F}$ tel que $H$ soit stablement conjugu\'e \`a un \'el\'ement de ${\cal L}^{L,nr,st}_{ell}$. On note $L_{H}$ le commutant de $A_{H}$ dans $G$. C'est un Levi de $G$, qui est donc conjugu\'e par un \'el\'ement de $G(F)$ \`a un unique $L\in \underline{{\cal L}}_{F}$. Par des calculs d\'ej\`a faits plusieurs fois, on voit que ce Levi $L$ r\'epond \`a la question et que c'est le seul. Cela prouve (1). 

Fixons $L\in \underline{{\cal L}}_{F}$. On a prouv\'e  en \ref{descriptiondeDstable} que l'application $k^{L,st}_{cusp}:{\cal K}^{st}_{cusp}(\mathfrak{l}(F))\to I(\mathfrak{l}(F))^*$ \'etait un isomorphisme (antilin\'eaire) de l'espace de d\'epart sur $D^L({\cal D}^{st}_{cusp}(\mathfrak{l}(F)))$. D'autre part, $D^L$ est injective sur ${\cal D}^{st}_{cusp}(\mathfrak{l}(F))$. Notons $ h^{L,st}_{cusp}: {\cal K}^{st}_{cusp}(\mathfrak{l}(F)) \to {\cal D}^{st}_{cusp}(\mathfrak{l}(F))$ l'isomorphisme (lin\'eaire)  tel que $k^{L,st}_{cusp}=D^L\circ h^{L,st}_{cusp}$. 
 Le groupe $Norm_{G(F)}(L)$ agit sur tous  les espaces intervenant et l'application $h^{L,st}_{cusp}$ est \'equivariante (par simple transport de structure). Cette application se quotiente donc en un isomorphisme  
 $$\bar{h}_{cusp}^{L,st}:{\cal K}^{st}_{cusp}(\mathfrak{l}(F))/Norm_{G(F)}(L) \to {\cal D}^{st}_{cusp}(\mathfrak{l}(F))/Norm_{G(F)}(L).$$
L'assertion (i) du lemme \ref{lespaceKstableetlesLevi} signifie que, pour $\varphi\in {\cal K}^{st}_{cusp}(\mathfrak{l}(F))/Norm_{G(F)}(L)$, on a l'\'egalit\'e $k^{G,st}(\bar{j}(\varphi))=D^G(\bar{h}^{L,st}_{cusp}(\varphi))$.

 On obtient un diagramme commutatif
 $$\begin{array}{ccc}\oplus_{L\in \underline{{\cal L}}}{\cal K}^{st}_{cusp}(\mathfrak{l}(F))/Norm_{G(F)}(L)&\stackrel{\bar{j}}{\to}& {\cal K}^{st}(\mathfrak{g}(F))\\ \downarrow \oplus_{L}\bar{h}^{L,st}_{cusp}&&\downarrow k^{G,st}\\ \oplus_{L\in \underline{{\cal L}}}{\cal D}^{st}_{cusp}(\mathfrak{l}(F))/Norm_{G(F)}(L)&\stackrel{D^G}{\to}& I(\mathfrak{g}(F))^*\\ \end{array}$$
 Les applications horizontales du haut et verticales de gauche sont des isomorphismes et l'application horizontale du bas est injective. D'apr\`es le lemme \ref{decompositionLevi}, l'espace sud-ouest est \'egal \`a ${\cal D}^{st}(\mathfrak{g}(F))$. La commutativit\'e du diagramme entra\^{\i}ne l'\'enonc\'e. $\square$

\subsubsection{ Dimension de l'espace ${\cal D}^{st}(\mathfrak{g}(F))$}\label{dimension}
Pour tout $M^*\in {\cal L}_{F,sd}^{*,st}$ et tout $Cl_{Fr}(w)\in W^{I_{F}}(M^*)_{G-adm}/Fr-conj$, fixons un \'el\'ement $H(M^*,w)\in {\cal L}_{F}^{nr,st}$ dont la classe de conjugaison stable est param\'etr\'ee par $(M^*,Cl_{Fr}(w))$, cf. proposition \ref{parametrage}. Alors
$${\cal K}^{st}(\mathfrak{g}(F))=\oplus_{M^*\in {\cal L}^{*,st}_{F,sd}}\oplus_{Cl_{Fr}(w)\in W^{I_{F}}(M^*)_{G-adm}/Fr-conj}FC^{st}(\mathfrak{h}(M^*,w)_{SC}(F)). $$
D'apr\`es le lemme \ref{passagestable}, les espaces $FC^{st}(\mathfrak{h}(M^*,w)_{SC}(F))$ et $FC^{st}(\mathfrak{m}^*_{SC}(F))$ sont isomorphes. On d\'eduit de la proposition \ref{descriptionfinale} l'\'egalit\'e
$$dim({\cal D}^{st}(\mathfrak{g}(F)))=\sum_{M^*\in {\cal L}_{F,sd}^{*,st}}\vert W^{I_{F}}(M^*)_{G-adm}/Fr-conj\vert  dim(FC^{st}(\mathfrak{m}^*_{SC}(F))).$$

 \section{  Int\'egrales orbitales nilpotentes  }

\subsection{Orbites nilpotentes}

\subsubsection{Rappels sur les repr\'esentations de $\mathfrak{sl}(2)$\label{sl2}}

 Dans cette sous-section 3.1, le corps de base $k$ est une extension alg\'ebrique de  ${\mathbb F}_{q}$ ou de $F$. On note $\{\pmb{ f},\pmb{ h},\pmb{ e}\}$ la base habituelle de $\mathfrak{sl}(2,k)$:
 $$\pmb{f}=\left(\begin{array}{cc}0&0\\1&0\\ \end{array}\right),\, \pmb{ h}=\left(\begin{array}{cc}1&0\\0&-1\\ \end{array}\right),\, \pmb{ e}=\left(\begin{array}{cc}0&1\\0&0\\ \end{array}\right).$$
 
  Soit $n$ un entier, supposons 

(1) $1\leq n< p$, $2\not=p$.

 Alors, \`a isomorphisme pr\`es,  il existe une unique repr\'esentation irr\'eductible $\rho_{n}$ de $\mathfrak{sl}(2,k)$ dans un $k$-espace vectoriel  ${\cal W}_{n}$   de dimension $n$. On peut choisir une base $\{w_{1},...,w_{n}\}$ de ${\cal W}_{n}$ de sorte que
$\rho_{n}(\pmb{ h})(w_{i})= (2i-n-1)w_{i}$  et $\rho_{n}(\pmb{ e})(w_{i})=w_{i+1}$ pour tout $i=1,...,n$, avec la convention $w_{n+1}=0$. On a alors $\rho_{n}(\pmb{ f})(w_{i})=c_{n,i}w_{i-1}$ avec la convention $w_{0}=0$, o\`u $c_{n,i}$ est un entier premier \`a $p$. Les valeurs propres de $\rho_{n}(\pmb{h})$ sont des entiers premiers \`a $p$ et  la diff\'erence entre deux valeurs propres distinctes est encore un entier premier \`a $p$. Deux telles bases sont homoth\'etiques, \'etant d\'etermin\'es par le vecteur non nul $w_{1}$ appartenant \`a la droite des \'el\'ements $w\in {\cal W}_{n}$ tels que $\rho_{n}(\pmb{ h})(w)=(1-n)w$. On fixe une telle base. 

Supposons pour la fin du paragraphe que $k$ soit  \'egal \`a $F$ ou \`a $F^{nr}$. On note  ${\cal L}_{n}$ le sous-$\mathfrak{o}_{k}$-module de ${\cal W}_{n}$ engendr\'e par les \'el\'ements de base. Alors ${\cal L}_{n}$ est un sous-$\mathfrak{sl}(2,\mathfrak{o}_{k})$-module de $W_{n}$. Tout $\mathfrak{o}_{k}$-r\'eseau de $W_{n}$ qui est conserv\'e par $\mathfrak{sl}(2,\mathfrak{o}_{k})$ est homoth\'etique \`a ${\cal L}_{n}$. 

Soit $\rho$ une repr\'esentation de $\mathfrak{sl}(2,k)$ dans un $k$-espace $V$.   Consid\'erons une cha\^{\i}ne de $\mathfrak{o}_{k}$-r\'eseaux $L_{.}=(L_{0}\supset L_{1}\supset... \supset L_{m}=\mathfrak{p}_{k}L_{0})$ de $V$. 

\begin{prop}{Supposons  que toute sous-repr\'esentation irr\'eductible de $\rho$ soit de dimension strictement inf\'erieure  \`a $p$. Supposons que $L_{i}$ soit conserv\'e par $\mathfrak{sl}(2,\mathfrak{o}_{k})$ pour tout $i=0,...,m-1$. Alors il existe un ensemble $J$, muni d'une d\'ecomposition en r\'eunion disjointe $J=\sqcup_{i=0,...m-1}J_{i}$ et d'une application $d:J\to \{1,...,p-1\}$ et 
  il existe un isomorphisme de $\mathfrak{sl}(2,k)$-modules
$\phi:V\to  \oplus_{j\in J}{\cal W}_{d(j)}$ de sorte que, pour tout $i=0,...,m-1 $, on ait l'\'egalit\'e
$$\phi(L_{i})=\left(\oplus_{i'=0,...,i-1}\oplus_{j\in J_{i'}}\mathfrak{p}_{k}{\cal L}_{d(j)}\right)\oplus\left(\oplus_{i'=i,...,m-1}\oplus_{j\in J_{i'}}{\cal L}_{d(j)}\right).$$}\end{prop}

La d\'emonstration est \'el\'ementaire, on la laisse au lecteur.

\subsubsection{$\mathfrak{sl}(2)$-triplets et orbites nilpotentes\label{sl2triplets}}

 Pour la suite de la sous-section 3.1, on suppose que $G$ est un groupe r\'eductif connexe d\'efini sur $k$.    On impose l'hypoth\`ese $(Hyp)_{1}(p)$.

On appelle orbite nilpotente dans $\mathfrak{g}(k)$ une classe de conjugaison par $G(k)$ dans $\mathfrak{g}_{nil}(k)$. On  note $\mathfrak{g}_{nil}(k)/conj$ l'ensemble des orbites nilpotentes dans $\mathfrak{g}(k)$.  On appelle $\mathfrak{sl}(2)$-triplet de $\mathfrak{g}(k)$ un triplet $(f,h,e)$ d'\'el\'ements de $\mathfrak{g}(k)$ vérifiant l'une des conditions suivantes:

$(f,h,e)$ engendrent une sous-alg\`ebre isomorphe \`a $\mathfrak{sl}(2,k)$, les trois \'el\'ements s'identifiant aux \'el\'ements $\pmb{ f}$, $\pmb{h}$ et $\pmb{e}$ de $\mathfrak{sl}(2,k)$;

$(f,h,e)=(0,0,0)$. 

 Les \'el\'ements $f$, $h$ et $e$ appartiennent forc\'ement \`a la sous-alg\`ebre $\mathfrak{g}_{SC}(k)$ de $\mathfrak{g}(k)$.

Pour un tel $\mathfrak{sl}(2)$-triplet, il existe un unique homomorphisme $\phi:SL(2)\to G_{SC}$ d\'efini sur $k$ tel que $\phi(exp(\lambda{\bf e}))=exp(\lambda e)$ et $\phi(exp(\lambda{\bf f}))=exp(\lambda f)$ pour tout $\lambda\in k$, cf. \cite{C} paragraphe 5.5. L'application $Ad\circ \phi$ est une repr\'esentation de $SL(2)$ dans $\mathfrak{g}_{SC}$. L'hypoth\`ese $(Hyp)_{1}(p)$ entra\^{\i}ne que  cette  repr\'esentation est semi-simple et que les dimensions de toutes ses composantes irr\'eductibles   sont strictement inférieures à $p$.   

Notons $x_{*,h}$ l'homomorphisme 
 $$\lambda\mapsto \phi(\left(\begin{array}{cc}\lambda&0\\ 0&\lambda^{-1}\\ \end{array}\right))$$
 de $GL(1)$ dans $G_{SC}$. C'est un sous-groupe \`a un param\`etre de $G_{SC}$.  On a $h=\partial x_{*,h}(1)$, o\`u  $\partial x_{*,h}:\mathfrak{gl}(1)\to \mathfrak{g}_{SC}$ est  l'homomorphisme d\'eriv\'e de $x_{*,h}$.       Si $k$ est de caract\'eristique nulle, les  sous-espaces propres de $\mathfrak{g}_{SC}$ pour l'action de $x_{*,h}$ correspondent bijectivement aux sous-espaces propres pour l'action de $h$. Si $k$ est de caract\'eristique $p$, cela serait vrai sous une hypoth\`ese plus forte que $(Hyp)_{1}(p)$ mais ne l'est pas toujours sous cette hypoth\`ese. Quel que soit $k$, on a toutefois, cf. \cite{C}, proposition 5.5.8:

(1) le commutant de $h$ dans $G$ est \'egal au commutant du groupe $x_{*,h}(\bar{k}^{\times})$; c'est un $k$-Levi de $G$.

Dans le cas o\`u $k\subset F^{nr}$, on  a m\^eme

(2) le commutant de $h$ dans $G$ est \'egal au commutant du groupe $x_{*,h}(\mathfrak{o}_{F^{nr}}^{\times})$.

 Notons $Trip(\mathfrak{g}(k))$ l'ensemble des $\mathfrak{sl}(2)$-triplets dans $\mathfrak{g}(k)$.  On sait que l'application
$(f,h,e)\mapsto e$ se quotiente en une bijection entre l'ensemble des classes de conjugaison par $G(k)$ dans $Trip(\mathfrak{g}(k))$ et l'ensemble des orbites nilpotentes   dans $\mathfrak{g}(k)$, cf. \cite{C} th\'eor\`eme 5.5.11. Conform\'ement \`a notre convention, si $k=\bar{k}$ est alg\'ebriquement clos, on supprime le $k$ de la notation: on  appelle orbite nilpotente dans $\mathfrak{g}(=\mathfrak{g}(\bar{k}))$ une classe de conjugaison par $G(=G(\bar{k}))$ dans $\mathfrak{g}_{nil}(=\mathfrak{g}_{nil}(\bar{k}))$; on note $\mathfrak{g}_{nil}/conj$ l'ensemble  des orbites nilpotentes dans $\mathfrak{g}$ et $Trip(\mathfrak{g})$ l'ensemble des $\mathfrak{sl}(2)$-triplets dans $\mathfrak{g}$.

  Supposons $k$ alg\'ebriquement clos.  
  Soit ${\cal O}\in \mathfrak{g}_{nil}/conj$. Fixons un $\mathfrak{sl}(2)$-triplet $(f,h,e)$ de $\mathfrak{g}$ tel que $e\in {\cal O}$. La classe de conjugaison de $x_{*,h}$ ne d\'epend pas du choix du $\mathfrak{sl}(2)$-triplet. On note $x_{*,{\cal O}}$ cette classe de conjugaison. Fixons une paire de Borel $(B,T)$ de $G$.   Il existe un unique \'el\'ement dans cette classe, notons-le encore $x_{*,{\cal O}}$, qui appartient \`a  la chambre positive fermée  de $X_{*}(T_{sc})$ relative à $B$. D'après Dynkin,  $x_{*,{\cal O}}$ v\'erifie les relations $<\alpha,x_{*,{\cal O}}>\in \{0,1,2\}$ pour toute racine simple $\alpha\in \Delta$.   L'application ${\cal O}\mapsto x_{*,{\cal O}}$ est injective.   Munissons $X_{*}(T_{sc})$ d'un produit scalaire invariant par l'action du groupe de Weyl, que l'on peut supposer \`a valeurs dans ${\mathbb Q}$. La norme $\vert \vert x_{*,{\cal O}}\vert \vert $ est alors d\'efinie et on a ainsi associ\'e \`a ${\cal O}$ un r\'eel positif (strictement positif si ${\cal O}\not=\{0\}$). 
  
    Levons l'hypoth\`ese que $k$ est alg\'ebriquement clos. A une orbite nilpotente ${\cal O}\in \mathfrak{g}_{nil}(k)/conj$ est associ\'ee l'unique  orbite ${\cal O}_{\bar{k}}\in \mathfrak{g}_{nil}/conj$ telle que ${\cal O}\subset {\cal O}_{\bar{k}}$. On pose $x_{*,{\cal O}}=x_{*,{\cal O}_{\bar{k}}}$. Evidemment, l'application ${\cal O}\mapsto x_{*,{\cal O}}$ n'est plus injective en g\'en\'eral.

    \subsubsection{Orbites engendr\'ees\label{orbitesengendrees}}
      Soit $H$ un sous-groupe r\'eductif connexe de $G$. Fixons des sous-tores maximaux $T$ de $G$ et $T_{H}$ de $H$ et munissons $X_{*}(T)$ d'un produit scalaire invariant par l'action du groupe de Weyl $W$ de $G$. On en d\'eduit un tel produit scalaire sur $X_{*}(T_{H})$. En effet,  \`a conjugaison pr\`es, on peut supposer $T_{H}\subset T$ et le produit cherch\'e est la restriction de celui sur $X_{*}(T)$. Cela ne d\'epend pas du choix de la conjugaison. 
      
      Soit ${\cal O}\in \mathfrak{h}_{nil}/conj$. Il s'en d\'eduit une orbite ${\cal O}^G\in \mathfrak{g}_{nil}/conj$: c'est l'unique orbite pour l'action de $G$ qui contient ${\cal O}$.  Fixons un $\mathfrak{sl}(2)$-triplet $(f,h,e)$ de $\mathfrak{h}$ tel que $e\in {\cal O}$. C'est aussi un $\mathfrak{sl}(2) $-triplet dans $\mathfrak{g}$ tel que $e\in {\cal O}^G$. On en d\'eduit que $x_{*,{\cal O}^G}$ est la $G$-orbite de sous-groupes \`a un param\`etre engendr\'ee par la $H$-orbite $x_{*,{\cal O}}$.  
  A fortiori, 
  
  (1) $\vert \vert x_{*,{\cal O}^G}\vert \vert =\vert \vert x_{*,{\cal O}}\vert \vert $.
  
  L'ensemble $\mathfrak{h}_{nil}/conj$ est muni de l'ordre usuel: ${\cal O} \leq {\cal O}'$ si ${\cal O}$ est incluse dans l'adh\'erence de Zariski de ${\cal O}'$. Il en est de m\^eme de $\mathfrak{g}_{nil}/conj$. L'application ${\cal O}\mapsto {\cal O}^G$ n'est pas injective en g\'en\'eral. On a toutefois le lemme suivant.
  
  \begin{lem}{L'application ${\cal O}\mapsto {\cal O}^G$ est strictement croissante.}\end{lem}
  
  Preuve. Si ${\cal O}$ est contenue dans l'adh\'erence de ${\cal O}'$, il est clair que ${\cal O}^G$ est contenue dans l'adh\'erence de ${\cal O}^{'G}$. On doit prouver que si, de plus ${\cal O}\not={\cal O}'$ alors ${\cal O}^G\not={\cal O}^{'G}$. Comme ci-dessus, on fixe un $\mathfrak{sl}(2)$-triplet $(f,h,e)$ de $\mathfrak{h}$ tel que $e\in {\cal O}$.  On pose $\mathfrak{z}_{\mathfrak{h}}(f)=\{X\in \mathfrak{h}; [f,X]=0\}$ et on d\'efinit de m\^eme $\mathfrak{z}_{\mathfrak{g}}(f)$. On consid\`ere l'application
  $$\begin{array}{ccc}H\times \mathfrak{z}_{\mathfrak{h}}(f)&\to&\mathfrak{h}\\ (h,X)&\mapsto& h(e+X)h^{-1}.\\ \end{array}$$
 D'apr\`es \cite{Slo} paragraphe 7.4, son image contient un voisinage de $e$. Il existe donc un \'el\'ement $X\in \mathfrak{z}_{\mathfrak{h}}(f)$ tel que $e+X\in {\cal O}'$. On a $\mathfrak{z}_{\mathfrak{h}}(f)\subseteq\mathfrak{z}_{\mathfrak{g}}(f)$. Supposons ${\cal O}^{'G}={\cal O}^G$. Alors $e+X\in {\cal O}^G$. Or l'ensemble des $Y\in \mathfrak{z}_{\mathfrak{g}}(f)$ tels que $e+Y\in {\cal O}^G$ est r\'eduit \`a $\{0\}$, cf. \cite{W1} V.7(9). Donc $X=0$ et $e\in {\cal O}'$. Cela contredit l'hypoth\`ese ${\cal O}'\not={\cal O}$. $\square$
 
 {\bf Remarque.} Slodowy \'enonce ses r\'esultats en imposant \`a $p$ une hypoth\`ese plus forte que la n\^otre, parce qu'il \'enonce le th\'eor\`eme de Jacobson-Morozov sous cette hypoth\`ese plus forte. En fait, d'apr\`es \cite{C}, th\'eor\`eme 5.5.11, notre hypoth\`ese $(Hyp)_{1}(p)$ suffit.

   \subsubsection{Les groupes $A(N)$ et $\bar{A}(N)$\label{lesgroupesA(N)}}

  Pour $N\in \mathfrak{g}_{nil}$, on pose $A(N)=Z_{G}(N)/Z_{G}(N)^0$. Fixons un triplet $(f,h,N)\in Trip(\mathfrak{g})$ dont le troisi\`eme terme soit notre nilpotent $N$.  Notons $Z_{G}(f,h,N)$ le commutant commun de $f$, $h$ et $N$. Il est connu que $Z_{G}(f,h,N)$ est une composante de Levi de $Z_{G}(N)$, c'est-\`a-dire que $Z_{G}(f,h,N)$ est r\'eductif (en g\'en\'eral non connexe) et que $Z_{G}(N)$ est produit semi-direct de $Z_{G}(f,h,N)$ et du radical unipotent de $Z_{G}(N)$. De m\^eme, $Z_{G}(f,h,N)^0$ est une composante de Levi de $Z_{G}(N)^0$. On a l'\'egalit\'e $A(N)=Z(f,h,N)/Z(f,h,N)^0$.

 Lusztig a d\'efini un quotient $\bar{A}(N)$ de $A(N)$ (Lusztig donne la d\'efinition en \cite{L1} 13.1 en supposant que l'orbite de $N$ est sp\'eciale; la d\'efinition se g\'en\'eralise \`a toute orbite nilpotente, cf. \cite{Sommers} paragraphe 5).  On note $A^1(N)$ le noyau de l'homomorphisme naturel $A(N)\to \bar{A}(N)$. Il contient l'image de $Z(G)$ dans 
 $A(N)$. 

\subsubsection{Partitions\label{partitions}}

 Pour un ensemble fini $Y$ et pour $\eta\in \{\pm 1\}$, on note $\{\pm 1\}^{Y}_{\eta}$ le sous-ensemble des \'el\'ements $(\epsilon_{y})_{y\in Y}\in \{\pm 1\}^{Y}$ tels que $\prod_{y\in Y}\epsilon_{y}=\eta$. 
 
 On appelle partition une classe d'\'equivalence de suites finies d\'ecroissantes d'entiers positifs ou nuls $\lambda=(\lambda_{1}, \lambda_{2},...)$, deux telles suites \'etant \'equivalentes si et seulement si elles ne diff\`erent que par des termes nuls. Pour une telle partition $\lambda$, on utilise les notations suivantes:
 
 $S(\lambda)=\sum_{j\geq1}\lambda_{j}$; 
 
 pour $n\in {\mathbb N}$, $S_{n}(\lambda)=\sum_{j;1\leq j\leq n}\lambda_{j}$;
 
  pour $i\in {\mathbb N}_{>0}$,  $mult_{\lambda}(i)=\vert \{j\geq1; \lambda_{j}=i\}$ et $mult_{\lambda}(\geq i)=\sum_{i'\in {\mathbb N},i'\geq i}mult_{\lambda}(i')$; 
  
  pour un sous-ensemble $E$ de ${\mathbb N}_{>0}$, $mult_{\lambda}(E)=\vert \{j\geq 1; \lambda_{j}\in E\}$; 
  
  $l(\lambda)=mult_{\lambda}(\geq1)$, autrement dit $l(\lambda)$ est le nombre d'indices $j\geq1$ tels que $\lambda_{j}\not=0$;
  
  $Jord(\lambda)=\{i\in {\mathbb N}_{>0}; mult_{\lambda}(i)\geq1\}$.
  
  \noindent En notant $i_{1}>...>i_{r}$ les \'el\'ements de $Jord(\lambda)$, on notera aussi la partition $\lambda=(i_{1}^{mult_{\lambda}(i_{1})},...,i_{r}^{mult_{\lambda}(i_{r})})$. Pour $N\in {\mathbb N}$, on note  
 ${\cal P}(N)$ l'ensemble des partitions $\lambda$ telles que $S(\lambda)=N$.  Pour $\lambda'\in {\cal P}(N')$ et $\lambda''\in {\cal P}(N'')$, on d\'efinit leur r\'eunion $\lambda=\lambda'\cup \lambda''\in {\cal P}(N'+N'')$: c'est l'unique partition telle que $mult_{\lambda}(i)=mult_{\lambda'}(i)+mult_{\lambda''}(i)$ pour tout $i\geq1$. On d\'efinit la somme $\lambda=\lambda'+\lambda''\in {\cal P}(N'+N'')$ par $\lambda_{j}=\lambda'_{j}+\lambda''_{j}$ pour tout $j\geq1$. On introduit sur ${\cal P}(N)$ l'ordre usuel: $\lambda'\leq \lambda''$ si et seulement si $S_{n}(\lambda')\leq S_{n}(\lambda'')$ pour tout entier $n\geq0$. 
 
 On  note ${\cal P}_{2}(N)$ l'ensemble des couples de partitions $(\lambda',\lambda'')$ tels que $S(\lambda')+S(\lambda'')=N$.

On note ${\cal P}^{symp}(N)$, resp. ${\cal P}^{orth}(N)$, l'ensemble des $\lambda\in {\cal P}(N)$ telles que $mult_{\lambda}(i)$ est pair pour tout $i\geq1$ impair, resp. pair. Evidemment ${\cal P}^{symp}(N)=\emptyset$ si $N$ est impair. Pour $\lambda\in {\cal P}^{symp}(N)$, resp. $\lambda\in{\cal P}^{orth}(N)$, on note $Jord_{bp}(\lambda)=\{i\in Jord(\lambda); i\text{ est pair}\}$, resp. $Jord_{bp}(\lambda)=\{i\in Jord(\lambda); i\text{ est impair}\}$.    Notons $\xi(\lambda)$ l'\'el\'ement $(\xi(\lambda)_{i})_{i\in Jord_{bp}(\lambda)}\in \{\pm 1\}^{Jord_{bp}(\lambda)}$ tel que $\xi(\lambda)_{i}=(-1)^{mult_{\lambda}(i)}$. Pour $\lambda\in {\cal P}^{orth}(N)$, on a $\xi(\lambda)\in \{\pm 1\}^{Jord_{bp}(\lambda)}_{\eta}$, o\`u $\eta=(-1)^N$.

\subsubsection{Param\'etrage des orbites nilpotentes\label{parametrageorbitesnilpotentes}}
   
 Pour ce paragraphe, {\bf on suppose que} $\boldsymbol{G}$ {\bf est adjoint et simple.}
 
 Si $G$ est de type exceptionnel, le param\'etrage des  orbites nilpotentes de $\mathfrak{g}$ se trouve dans diff\'erentes tables.  Nous prenons comme r\'ef\'erence \cite{C}, pages 401-407, dont nous utiliserons les notations. 
Pour un \'el\'ement $N\in \mathfrak{g}_{nil}$, ces tables d\'ecrivent le groupe $A(N)$. Il est toujours isomorphe \`a un groupe de permutations $\mathfrak{S}_{m}$ pour un entier $m\in \{1,...,5\}$.   Le groupe $\bar{A}(N)$ est d\'ecrit en \cite{Sommers} paragraphe 9.   On a toujours $\bar{A}(N)=A(N)$ ou $\bar{A}(N)=\{1\}$, sauf dans le cas o\`u $G$ est de type $E_{8}$ et l'orbite de $N$ est de type $E_{8}(b_{6})$. Dans ce dernier cas, on a $A(N)=\mathfrak{S}_{3}$ et $\bar{A}(N)=\mathfrak{S}_{2}$.

 Si $G$ est classique, $\mathfrak{g}_{nil}/conj$ est param\'etr\'e par un ensemble de partitions (avec une petite entorse dans le cas d'un groupe de type $D_{n}$, cf. ci-dessous). Pour une partition $\lambda$ dans cet ensemble, on note $A(\lambda)$ le groupe $A(N)$ pour un \'el\'ement $N\in \mathfrak{g}_{nil}$ param\'etr\'e par $\lambda$.
 
 Si $G$ est de type $A_{n-1}$, $\mathfrak{g}_{nil}/conj$ est param\'etr\'e par ${\cal P}(n)$. On a $A(\lambda)=\bar{A}(\lambda)=\{1\}$ pour tout $ \lambda\in {\cal P}(n)$. 
 
 Si $G$ est de type $B_{n}$, $\mathfrak{g}_{nil}/conj$ est param\'etr\'e par ${\cal P}^{orth}(2n+1)$. Soit $\lambda\in {\cal P}^{orth}(2n+1)$. On a $A(\lambda)=\{\pm 1\}^{Jord_{bp}(\lambda)}_{1}$. Notons $i_{1}>i_{2}>...>i_{t}$ les \'el\'ements de $Jord_{bp}(\lambda)$. On d\'efinit par r\'ecurrence une suite strictement  croissante d'entiers $(r_{j})_{j=0,...,k}$ telle que $0\leq r_{j}\leq t$  de la fa\c{c}on suivante. On pose $r_{0}=0$. Supposons d\'efini $r_{j}$. Si $r_{j} =t$, on s'arr\^ete et on pose $k=j$. Si $r_{j}<t$,  il existe forc\'ement un entier $m\in \{r_{j}+1,...,t\}$ tel que $mult_{\lambda}(\geq i_{m})$ soit impair: l'entier $m=t$ v\'erifie cette condition.  On note $r_{j+1}$ le plus petit tel entier. Cela d\'efinit notre suite $(r_{j})_{j=0,...,k}$. A tout $j=1,...,k$, on associe l'intervalle $\Delta_{j}\subset Jord_{bp}(\lambda)$ form\'e des $i_{m}$ pour $r_{j-1}<m\leq r_{j}$. Remarquons que $\sum_{i\in \Delta_{j}}mult_{\lambda}(i)$ est impair pour $j=1$ et est pair pour $j\geq2$. On note $Int(\lambda)=\{\Delta_{j};j=1,...,k\}$ l'ensemble de ces intervalles. Cet ensemble est naturellement ordonn\'e: $\Delta_{1}>\Delta_{2}>... >\Delta_{k}$ (les entiers appartenant \`a $\Delta_{1}$ sont strictement sup\'erieurs aux entiers appartenant \`a $\Delta_{2}$ etc...). On a 
  $\bar{A}(\lambda)=\{\pm 1\}^{Int(\lambda)}_{1}$. L'homomorphisme $A(\lambda)\to \bar{A}(\lambda)$ envoie $(\epsilon_{i})_{i\in Jord_{bp}(\lambda)}$ sur $(\eta_{\Delta})_{\Delta\in Int(\lambda)}$ o\`u $\eta_{\Delta}=\prod_{i\in \Delta}\epsilon_{i}$.

  Si $G$ est de type $C_{n}$, $\mathfrak{g}_{nil}/conj$ est param\'etr\'e par ${\cal P}^{symp}(2n)$. Soit $\lambda\in {\cal P}^{symp}(2n)$. On a $A(\lambda)=\{\pm 1\}^{Jord_{bp}(\lambda)}/\{1,\xi(\lambda)\}$. Remarquons que l'\'el\'ement $\xi(\lambda)$ peut \^etre \'egal \`a $1$.   Notons $i_{1}>i_{2}>...>i_{t}$ les \'el\'ements de $Jord_{bp}(\lambda)$. On d\'efinit par r\'ecurrence une suite strictement  croissante d'entiers $(r_{j})_{j=0,...,k}$ telle que $0\leq r_{j}\leq t$  de la fa\c{c}on suivante. On pose $r_{0}=0$. Supposons d\'efini $r_{j}$. Si $r_{j} =t$, on s'arr\^ete et on pose $k=j$. Si $r_{j}<t$, supposons qu'il existe $m\in \{r_{j}+1,...,t\}$ tel que $mult_{\lambda}(\geq i_{m})$ soit pair. Dans ce cas on note $r_{j+1}$ le plus petit tel entier. Supposons qu'il n'existe pas de tel entier.  Alors la suite s'arr\^ete, on pose $k=j$. La suite $(r_{j})_{j=0,...,k}$ \'etant d\'efinie, on d\'efinit l'ensemble d'intervalles $Int(\lambda)$ comme dans le cas $B_{n}$. Remarquons que, cette fois, $\sum_{i\in \Delta_{j}}mult_{\lambda}(i)$ est pair pour tout $j$.   
  On a
  $\bar{A}(\lambda)=\{\pm 1\}^{Int(\lambda)}$. L'homomorphisme $A(\lambda)\to \bar{A}(\lambda)$ se construit de la fa\c{c}on suivante.  On d\'efinit comme dans le cas $B_{n}$ un homomorphisme de $\{\pm 1\}^{Int(\lambda)}$ dans $\bar{A}(\lambda)$. On v\'erifie que son noyau contient $\xi(\lambda)$. L'homomorphisme se quotiente donc en un homomorphisme $A(\lambda)\to \bar{A}(\lambda)$ qui est l'homomorphisme cherch\'e.

  Si $G$ est de type $D_{n}$, il y a une application surjective de $\mathfrak{g}_{nil}/conj$ sur ${\cal P}^{orth}(2n)$. Les fibres ont un \'el\'ement sauf au-dessus d'une partition dont tous les termes sont pairs. Dans ce cas, la fibre a deux \'el\'ements. Pour un \'el\'ement $N$ d'une telle fibre exceptionnelle, on a $A(N)=\bar{A}(N)=\{1\}$. Soit $\lambda\in {\cal P}^{orth}(2n)$ poss\'edant au moins un terme impair. Alors $A(\lambda)=\{\pm 1\}^{Jord_{bp}(\lambda)}_{1}/\{1,\xi(\lambda)\}$. Remarquons que l'\'el\'ement $\xi(\lambda)$ peut \^etre \'egal \`a $1$. On d\'efinit l'ensemble d'intervalles $Int(\lambda)$ de la m\^eme fa\c{c}on que dans le cas $C_{n}$. On a $\bar{A}(\lambda)=\{\pm 1\}^{Int(\lambda)}_{1}$.  L'homomorphisme $A(\lambda)\to \bar{A}(\lambda)$ se construit comme dans le cas $C_{n}$.

  \subsubsection{ Action galoisienne\label{actiongaloisienne}}
  On suppose  que $k=F$ ou $k={\mathbb F}_{q}$, que $G_{AD}$ {\bf est simple et que} $G$ {\bf est quasi-d\'eploy\'e} (cette deuxième condition est automatique si $k={\mathbb F}_{q}$).   Le groupe $\Gamma_{k}$ agit sur $\mathfrak{g}_{nil}/conj$. Toute orbite nilpotente conserv\'ee par l'action galoisienne contient un point de $\mathfrak{g}(k)$ (cf. \cite{K1} th\'eor\`eme 4.2 si $k=F$;  si $k={\mathbb F}_{q}$, l'assertion  est v\'erifi\'ee pour toute orbite nilpotente ou pas, cela r\'esulte du th\'eor\`eme de Lang). Autrement dit, l'application naturelle $\mathfrak{g}_{nil}(k)/conj\to (\mathfrak{g}_{nil}/conj)^{\Gamma_{k}}$ est surjective. 
Montrons que

(1) si $G$ est d\'eploy\'e sur $k$ ou si $G$ n'est pas de type $D_{n}$, toute orbite nilpotente dans $\mathfrak{g}$ est conserv\'ee par l'action galoisienne.

Fixons une paire de Borel $(B,T)$ d\'efinie sur $k$. Pour ${\cal O}\in \mathfrak{g}_{nil}/conj$, on identifie $x_{*,{\cal O}}$ \`a l'unique \'el\'ement de cette classe appartenant \`a $X_{*}(T_{sc})$ et v\'erifiant les conditions de Dynkin. Pour $\sigma\in \Gamma_{k}$, on a \'evidemment  $x_{*,\sigma({\cal O})}=\sigma(x_{*,{\cal O}})$. Puisque $x_{*,{\cal O}}$ d\'etermine ${\cal O}$, il suffit de prouver que, sous nos hypoth\`eses, $x_{*,{\cal O}}$ est fix\'e par l'action galoisienne. Si $G$ est d\'eploy\'e, cette action est triviale sur $X_{*}(T_{sc})$ et le r\'esultat est clair. Si $G$ n'est pas d\'eploy\'e et n'est pas de type $D_{n}$, il ne peut \^etre que de type $A_{n-1}$ ou $E_{6}$. Si $G$ est de type $E_{6}$, $\Gamma_{k}$ agit sur $X_{*}(T_{sc})$ par un homomorphisme surjectif $\Gamma_{k}\to Aut({\cal D})={\mathbb Z}/2{\mathbb Z}$.   L'inspection de la table \cite{C} p. 402 montre que, pour toute orbite  ${\cal O}$, $x_{*,{\cal O}}$ est  fix\'e par  $Aut({\cal D})$.  Dans le cas o\`u $G$ est de type $A_{n-1}$, il est plus simple de faire appel \`a la description \'el\'ementaire des orbites nilpotentes par l'alg\`ebre lin\'eaire, qui conduit \`a la conclusion. D'o\`u (1).

 Supposons $G$ de type $D_{n}$ avec $n\geq4$ et fixons comme ci-dessus une paire de Borel $(B,T)$ d\'efinie sur $k$. Le groupe d'automorphismes  $Aut({\cal D})$  contient toujours l'\'el\'ement $\theta$ qui \'echange les racines $\alpha_{n-1}$ et $\alpha_{n}$ et fixe $\alpha_{i}$ pour $i=1,...,n-2$. Si $n\geq5$, ce groupe est r\'eduit \`a $\{1,\theta\}$. Supposons $n=4$. Alors ce groupe d'automorphismes contient l'\'el\'ement $\theta_{3}$ qui fixe $\alpha_{2}$ et envoie $\alpha_{1},\alpha_{3},\alpha_{4}$ respectivement sur $\alpha_{3},\alpha_{4},\alpha_{1}$. Le groupe est engendr\'e par $\theta$ et $\theta_{3}$ et est isomorphe \`a $\mathfrak{S}_{3}$. 
 Supposons que $\Gamma_{k}$ agisse par un homomorphisme surjectif $\Gamma_{k}\to \{1,\theta\}$. 
   On v\'erifie qu'une orbite nilpotente param\'etr\'ee par une partition $\lambda\in {\cal P}^{orth}(2n)$ poss\'edant au moins un terme impair est conserv\'ee par cette action. Par contre, les deux orbites param\'etr\'ees par une partition $\lambda$ dont les termes sont tous pairs sont \'echang\'es par $\theta$. Supposons maintenant que $n=4$ et que $\Gamma_{k}$ agisse par un homomorphisme $\Gamma_{k}\to Aut({\cal D})$ dont l'image contient un \'el\'ement d'ordre $3$.  En calculant les termes $x_{*,{\cal O}}$, l'argument de la preuve de (1) montre que les orbites conserv\'ees par l'action galoisienne sont celles param\'etr\'ees par les partitions
$$(2) \qquad (71),(53),(3^21^2),(32^21), (2^21^4), (1^8).$$

Revenons à un groupe quelconque vérifiant les hypothèses posées au début du paragraphe. Soit $N\in \mathfrak{g}_{nil}(k)$. Le groupe $\Gamma_{k}$ agit naturellement sur $A(N)$. Montrons que

(2) cette action se quotiente en une action sur $\bar{A}(N)$.

Preuve. Supposons $G$ est classique (on exclut par là le cas où $G$ est de type $D_{4}$ et où l'image de $\Gamma_{k}$ dans $Aut({\cal D})$ contient un élément d'ordre $3$). L'algèbre élémentaire permet de décrire $A(N)$ et on voit que l'action galoisienne sur ce groupe est triviale. L'assertion est alors triviale. Supposons que $G$ est exceptionnel. Comme on l'a dit en \ref{parametrageorbitesnilpotentes}, on a presque toujours $\bar{A}(N)=A(N)$ ou $\bar{A}(N)=\{1\}$, auquels cas l'assertion est triviale. Si $G$ est de type $E_{8}$ et que l'orbite de $N$ est de type $E_{8}(b_{6})$, on a $A(N)=\mathfrak{S}_{3}$ et $\bar{A}(N)=\mathfrak{S}_{2}$. Le noyau de la projection $A(N)\to \bar{A}(N)$ est l'unique sous-groupe d'ordre $3$ de $A(N)$. Il est forcément conservé par l'action galoisienne. Enfin, supposons que $G$ est de type $D_{4}$ et que  l'image de $\Gamma_{k}$ dans $Aut({\cal D})$ contient un élément d'ordre $3$.  D'après la description de \ref{parametrageorbitesnilpotentes}, on voit que l'on a toujours $\bar{A}(N)=A(N)$ ou $\bar{A}(N)=\{1\}$. On conclut comme ci-dessus. $\square$

\subsection{Orbites nilpotentes et groupe ${\bf G}$}

\subsubsection{Orbites nilpotentes dans $\boldsymbol{\mathfrak{g}}$\label{dansboldsymbolg}}
Dans la sous-section 3.2, on suppose que $G$ {\bf est d\'efini  et quasi-d\'eploy\'e sur} $F$, {\bf qu'il est adjoint et absolument simple}. On impose comme toujours la condition $(Hyp)_{1}(p)$.

On fixe un épinglage $\mathfrak{E}=(B,T,(E_{\alpha})_{\alpha\in \Delta})$ défini sur $F$ et on  effectue les constructions de \ref{discretisation}, on introduit en particulier le groupe ${\bf G}$.  
Il est d\'efini sur ${\mathbb F}_{q}$ et est muni d'une action alg\'ebrique de $I_{F}$ qui se quotiente en une action de  $I_{F}^{mod}$. Les groupes $\Gamma_{F}$ et $I_{F}$ agissent sur $\boldsymbol{\mathfrak{g}}_{nil}/conj$ et on d\'efinit les sous-ensembles d'invariants $(\boldsymbol{\mathfrak{g}}_{nil}/conj)^{\Gamma_{F}}$ et $(\boldsymbol{\mathfrak{g}}_{nil}/conj)^{I_{F}}$. On se gardera de croire qu'une orbite $\boldsymbol{{\cal O}}\in (\boldsymbol{\mathfrak{g}}_{nil}/conj)^{I_{F}}$ contient toujours un \'el\'ement fix\'e par $I_{F}$. Par contre, on a

(1) pour $\boldsymbol{{\cal O}}\in \boldsymbol{\mathfrak{g}}_{nil}/conj$, $\boldsymbol{{\cal O}}$ est conserv\'ee par $I_{F}$ si et seulement si la classe de sous-groupe \`a un param\`etre $x_{*,\boldsymbol{{\cal O}}}$ contient un \'el\'ement fix\'e par $I_{F}$.

Si $x_{*,\boldsymbol{{\cal O}}}$ contient un \'el\'ement fix\'e par $I_{F}$, cette classe est conserv\'ee par $I_{F}$. Puisque $\boldsymbol{{\cal O}}$ est d\'etermin\'ee par $x_{*,\boldsymbol{{\cal O}}}$, $\boldsymbol{{\cal O}}$ est elle-aussi conserv\'ee par $I_{F}$. Inversement, supposons $\boldsymbol{{\cal O}}$ conserv\'ee par $I_{F}$. Identifions $x_{*,\boldsymbol{{\cal O}}}$ \`a son unique repr\'esentant dans  la chambre positive fermée de $X_{*}({\bf T}_{sc})$. La classe de conjugaison de $x_{*,\boldsymbol{{\cal O}}}$ est conserv\'ee par $I_{F}$. Pour $\sigma\in I_{F}$, $\sigma(x_{*,\boldsymbol{{\cal O}}})$ est donc conjugu\'e \`a $x_{*,\boldsymbol{{\cal O}}}$. Puisque $\Delta$ est conserv\'e par $I_{F}$, $\sigma(x_{*,\boldsymbol{{\cal O}}})$  appartient aussi à la chambre positive fermée. Les deux \'el\'ements $x_{*,\boldsymbol{{\cal O}}}$ et $\sigma(x_{*,\boldsymbol{{\cal O}}})$ sont  deux \'el\'ements de cette chambre  qui sont conjugu\'es. Ils sont donc \'egaux. Donc $x_{*,\boldsymbol{{\cal O}}}$ est fix\'e par $I_{F}$. Cela prouve (1).

L'interpr\'etation en termes d'\'el\'ements $x_{*,\boldsymbol{{\cal O}}}$ permet d'appliquer les m\^emes arguments qu'en \ref{actiongaloisienne}. On obtient

(2) si $G$ est d\'eploy\'e sur $F$ ou si $G$ n'est pas de type $D_{n}$, toute orbite nilpotente dans $\boldsymbol{\mathfrak{g}}$ est conserv\'ee par l'action de $\Gamma_{F}$;

(3) supposons que $G$ est de type $D_{n}$ et que $\Gamma_{F}$, resp. $I_{F}$, agit sur ${\cal D}$ par un homomorphisme surjectif $\Gamma_{F}\to \{1,\theta\}$, resp. $I_{F}\to \{1,\theta\}$; alors $(\boldsymbol{\mathfrak{g}}_{nil}/conj)^{\Gamma_{F}}$, resp. $(\boldsymbol{\mathfrak{g}}_{nil}/conj)^{I_{F}}$, est form\'e des orbites param\'etr\'ees par une partition $\lambda\in {\cal P}^{orth}(2n)$ poss\'edant au moins un terme impair;

(4) supposons que $G$ est de type $D_{4}$ et  que $\Gamma_{F}$, resp. $I_{F}$, agit sur ${\cal D}$ par un homomorphisme $\Gamma_{F}\to Aut({\cal D})$, resp. $I_{F}\to Aut({\cal D})$ dont l'image contient un \'el\'ement d'ordre $3$; alors $(\boldsymbol{\mathfrak{g}}_{nil}/conj)^{\Gamma_{F}}$, resp. $(\boldsymbol{\mathfrak{g}}_{nil}/conj)^{I_{F}}$, est l'ensemble d\'ecrit en \ref{actiongaloisienne}(2).

  \subsubsection{Les ensembles $A^+(N)$, $\tilde{A}(N)$, $\bar{A}^+(N)$ et $\tilde{\bar{A}}(N)$\label{lesensemblesA+(N)} }
  
     On note $F^{G}$ l'extension de $F^{nr}$ telle que $\Gamma_{F^G}$ soit le noyau de l'action de $I_{F}$ sur ${\cal D}$. On pose $e^G=[F^G:F^{nr}]$. On a introduit un g\'en\'erateur  $\gamma$ de $ I_{F}^{mod}$, cf. \ref{discretisation}. Son image dans $\Gamma_{F^G/F^{nr}}$ est un g\'en\'erateur de ce groupe fini. Pour ne pas compliquer les notations, on note encore $\gamma$ cette image.

   Le groupe $\Gamma_{F^G/F^{nr}}$ agit par automorphismes alg\'ebriques sur ${\bf G}$, on peut d\'efinir le groupe ${\bf G}^+={\bf G}\rtimes \Gamma_{F^G/F^{nr}}$, non connexe si $F^G\not=F^{nr}$. On note $\tilde{{\bf G}}$ le sous-ensemble ${\bf G}\gamma$ de ${\bf G}^+$. 
   
   Soit $N\in \boldsymbol{\mathfrak{g}}_{nil}$. Le groupe $Z_{{\bf G}^+}(N)$ coupe toute composante connexe de ${\bf G}^+$ si et seulement si l'orbite de $N$ est conserv\'ee par $I_{F}$. On suppose qu'il en est ainsi. Remarquons que la structure de $Z_{{\bf G}^+}(N)$ ne change pas si l'on remplace $N$ par un \'el\'ement conjugu\'e \`a $N$ par un \'el\'ement de ${\bf G}$: cet \'el\'ement conjugue aussi les commutants. Les groupes  $Z_{{\bf G}}(N)$ et $Z_{{\bf G}}(N)^0$ sont des sous-groupes distingu\'es de $Z_{{\bf G}^+}(N)$. 
   On pose   $A^+(N)=Z_{{\bf G}^+}(N)/Z_{{\bf G}}(N)^0$.  On pose $Z_{\tilde{{\bf G}}}(N) =\tilde{{\bf G}}\cap Z_{{\bf G}^+}(N)$ et $\tilde{A}(N)=Z_{\tilde{{\bf G}}}(N)/Z_{{\bf G}}(N)^0$. 
 On a une suite exacte
 $$ 1\to A(N)\to A^+(N)\to \Gamma_{F^G/F^{nr}}\to 1$$

Le groupe $\Gamma_{F}^{nr}\simeq \Gamma_{{\mathbb F}_{q}}$ agit sur ${\bf G}$ et sur $\Gamma_{F^G/F^{nr}}$. Il en r\'esulte une action sur ${\bf G}^+$. Supposons comme ci-dessus que l'orbite de $N$ soit conserv\'ee par l'action de $I_{F}$ et que $N$ soit fix\'e par l'action de $\Gamma_{F}^{nr}$. Alors l'action de $\Gamma_{F}^{nr}$ sur ${\bf G}^+$ conserve le sous-groupe   $Z_{{\bf G}^+}(N)$ et cette action se descend en une action sur $A^+(N)$.
 
 {\bf Remarque.} L'action sur $I_{F}^{mod}$ de l'\'el\'ement de Frobenius $Fr\in \Gamma_{F}^{nr}$ est l'\'el\'evation \`a la puissance $q$. Si $G$ n'est pas du type $(D_{4},3-ram)$, on a $e_{G}=1$ ou $2$ et cette action se quotiente en l'action triviale sur $\Gamma_{F^G/F^{nr}}$. Alors l'action de $\Gamma_{F}^{nr}$ sur $A^+(N)$ conserve le sous-ensemble $\tilde{A}(N)$. Il en est de m\^eme si $G$ est de type $(D_{4},3-ram)$ et que $q\equiv 1\,mod\,3{\mathbb Z}$. Par contre, si $G$ est de type $(D_{4},3-ram)$ et que $q\equiv 2\, mod\,3{\mathbb Z}$, le sous-ensemble $\tilde{A}(N)$ n'est pas conserv\'e par l'action de $\Gamma_{F}^{nr}$. 
 
 \bigskip
 
   Soit $N\in \boldsymbol{\mathfrak{g}}_{nil}$ un \'el\'ement dont l'orbite est conserv\'ee par $I_{F}$. On a introduit en \ref{lesgroupesA(N)} le sous-groupe distingu\'e $A^1(N)\subset A(N)$. Montrons que
   
   (1) $A^1(N)$ est aussi distingué dans $A^+(N)$.
   
   Preuve. Si $G$ est déployé sur $F^{nr}$, on a $\boldsymbol{G}^+=\boldsymbol{G}$ et l'assertion est évidente. Si $G$ est de type $A_{n-1}$, resp. $E_{6}$ ou  $D_{4}$, on a $A(N)=\{1\}$, resp. $\bar{A}(N)=A(N)$ ou $\bar{A}(N)=\{1\}$, et l'assertion est tout aussi évidente. Reste le cas noté $(D_{n,ram})$ en  \ref{couplesstables}, c'est-à-dire $G$ est de type $D_{n}$ et $\gamma$ agit sur ${\cal D}$ par l'automorphisme $\theta$ d'ordre $2$. Puisqu'on a déjà traité le cas de $D_{4}$, on peut supposer que $n\geq5$.  On peut choisir un espace sur ${\mathbb F}_{q}$ muni d'une forme quadratique de sorte qu'en notant ${\bf SO}$ son groupe sp\'ecial orthogonal, on ait ${\bf SO}_{AD}={\bf G}$ et que l'action de $\gamma\in I_{F}$ sur ${\bf G}$ se d\'eduise par passage au quotient de l'action d'une sym\'etrie simple $S(\gamma)$ qui est d\'efinie sur ${\mathbb F}_{q}$.
   Introduisons le groupe orthogonal tout entier ${\bf O}$ de notre forme quadratique. Alors ${\bf G}^+$ s'identifie au quotient de ${\bf O}$ par le centre de ${\bf SO}$, l'\'el\'ement $S(\gamma)$ s'envoyant sur $\gamma\in \Gamma_{F^G/F^{nr}}$. Cette 
   identification est compatible aux actions de $\Gamma_{F}^{nr}$. Pour tout $N\in \boldsymbol{\mathfrak{g}}_{nil}$, le groupe $A^+(N)$ s'identifie \`a $Z_{{\bf O}}(N)/Z({\bf SO})Z_{{\bf SO}}(N)^0$. Des constructions d'alg\`ebre \'el\'ementaire montrent que ce groupe est commutatif. Plus précisément, soit $\lambda\in {\cal P}^{orth}(2n)$ la partition paramétrant l'orbite de $N$. .  Puisque cette  orbite est  conserv\'ee par $I_{F}$, $\lambda$ poss\`ede au moins un terme impair. Alors le groupe $A^+(\lambda)$  est \'egal \`a $\{\pm 1\}^{Jord_{bp}(\lambda)}/\{1,\xi(\lambda)\}$ et  on a $\tilde{A}(\lambda)=\{\pm 1\}^{Jord_{bp}(\lambda)}_{-1}/\{1,\xi(\lambda)\}$. L'assertion (1)  résulte de la commutativité de$A^+(\lambda)$. $\square$
   
    D'après (1), on peut  d\'efinir sans ambigu\"{\i}t\'e le quotient $\bar{A}^+(N)=A^+(N)/A^1(N)$ et le sous-ensemble $\tilde{\bar{A}}(N)=\tilde{A}(N)/A^1(N)$. 
    
     \subsubsection{Les ensembles ${\cal C}^{\sharp}_{F}$ et $\bar{{\cal C}}^{\sharp}_{F}$}\label{calCF}
   Le groupe $\Gamma_{F}^{mod}$ agit sur $\boldsymbol{G}$ via son quotient $\Gamma_{F^G/F}$. Introduisons le produit semi-direct
   $\boldsymbol{G}^{++}=\boldsymbol{G}\rtimes \Gamma_{F}^{mod}$. Il se projette sur $\Gamma_{F}^{mod}$. Pour $\sigma\in \Gamma_{F}^{mod}$, notons $\boldsymbol{G}^{++}_{\sigma}$ la fibre au-dessus de $\sigma$.  Le groupe $\boldsymbol{G}^{++}$ agit par conjugaison sur $\boldsymbol{G}$, donc aussi sur $\boldsymbol{\mathfrak{g}}$.    Pour $N\in \boldsymbol{\mathfrak{g}}_{nil}$, notons $A^{++}(N)=Z_{\boldsymbol{G}}(N)^0\backslash Z_{\boldsymbol{G}^{++}}(N)$. Pour $\sigma\in \Gamma_{F}^{mod}$, notons $A^{++}_{\sigma}(N)_=Z_{\boldsymbol{G}}(N)^0\backslash Z_{\boldsymbol{G}^{++}_{\sigma}}(N)$. Dire que l'orbite de $N$ est conservée par $\Gamma_{F}$ équivaut dire que $A^{++}_{\sigma}(N)$ est  non  vide pour tout  $\sigma\in \Gamma_{F}^{mod}$. Puisque $\Gamma_{F}^{mod}$ est engendré par $\gamma$ et $Fr$, il revient au m\^eme que $A^{++}_{\gamma}(N)$ et $A^{++}_{Fr}(N)$ soient non vides. 
   Remarquons que $A^{++}_{\gamma}(N)=\tilde{A}(N)$. 
   
   {\bf Remarque.} (1) Soit $(f,h,N)$ un $\mathfrak{sl}(2)$-triplet de $\boldsymbol{\mathfrak{g}}$. On a les isomorphismes $A^{++}(N)=Z_{{\bf G}}(f,h,N)^0\backslash Z_{{\bf G}^{++}}(f,h,N)$ et $A^{++}_{\sigma}(N)=Z_{{\bf G}}(f,h,N)^0\backslash Z_{{\bf G}_{\sigma}^{++}}(f,h,N)$ pour tout $\sigma\in \Gamma_{F}^{mod}$. 
   \bigskip
   
   Considérons  l'ensemble des triplets $(N,b,u)$ où:
   
   $N$ est un élément de $\boldsymbol{\mathfrak{g}}_{nil}$ dont l'orbite est conservée par $\Gamma_{F}$;
   
   $b\in A^{++}_{\gamma}(N)$ et $u\in A^{++}_{Fr}(N)$;
   
   $ub=b^qu$.

   Le groupe $\boldsymbol{G}$ agit par conjugaison sur cet ensemble de triplets. On note ${\cal C}^{\sharp}_{F}$ l'ensemble des classes de conjugaison. 
  Soit  ${\boldsymbol{\cal O}}\in (\boldsymbol{\mathfrak{g}}_{nil}/conj)^{\Gamma_{F}}$. On note ${\cal C}^{\sharp}_{F}(\boldsymbol{{\cal O}})$ le sous-ensemble des classe de conjugaison de triplets $(N,b,u)$ pour lesquels $N$ appartient à $\boldsymbol{{\cal O}}$.

     {\bf Remarque.}  (2) La  condition $ub=b^qu$  équivaut à l'existence d'un homomorphisme $\rho:\Gamma_{F}^{mod}\to A^{++}(N)$ tel que son composé avec la projection de $A^{++}(N)$ sur $\Gamma_{F}^{mod}$ soit l'identite et tel que $\rho(Fr)=u$, $\rho(\gamma)=b$. 
   
   \bigskip
   
   {\bf Notation.} On notera un élément de ${\cal C}^{\sharp}_{F}$  comme un triplet $(N,b,u)$ la représentant. Nous utiliserons le m\^eme abus de notation pour différents ensembles de classes de conjugaison  que nous introduirons dans la suite.
   
   \bigskip

Pour décrire plus concrètement l'ensemble ${\cal C}^{\sharp}_{F}$, on utilisera le lemme suivant, que nous démontrerons en \ref{preuvedecalCF}. On y utilise la terminologie introduite en \ref{couplesstables}.

 \begin{lem}{Soit ${\boldsymbol{\cal O}}\in (\boldsymbol{\mathfrak{g}}_{nil}/conj)^{\Gamma_{F}}$.  
 
  (i) Supposons que $G$ est d\'eploy\'e sur $F^{nr}$.  Il existe $N\in {\boldsymbol{\cal O}}$ qui est fix\'e par $\Gamma_{F}$ et tel  que l'action de $\Gamma_{F}^{nr}$ sur $A^+(N)=A(N)$ soit triviale.  
 
 (ii) Supposons que $G$ soit des types $(A_{n-1},ram)$, $(D_{4},3-ram)$ ou $(E_{6},ram)$. Si $A(N)=\{1\}$, il existe $N\in \boldsymbol{{\cal O}}$ qui est fix\'e par $\Gamma_{F}^{nr}$ et on a $A^+(N)\simeq  \Gamma_{F^G/F^{nr}}$. Si $A(N)\not=\{1\}$, il existe $N\in {\boldsymbol{\cal O}}$ qui est fix\'e par $\Gamma_{F}$ et  pour lequel   l'action de $\Gamma_{F}$ sur $A(N)$ est triviale. 
 
 (iii) Supposons que $G$ soit de type $(D_{n},ram)$. Alors il existe $N\in {\boldsymbol{\cal O}}$ qui est fix\'e par $\Gamma_{F}^{nr}$. 
  Pour un tel $N$, l'action de $\Gamma_{F}^{nr}$ sur $A^+(N)$ est triviale. Les groupes 
  $A(N)$ et $A^+(N)$ sont des produits de ${\mathbb Z}/2{\mathbb Z}$ et on a non canoniquement $A^+(N)\simeq A(N)\times {\mathbb Z}/2{\mathbb Z}$. }\end{lem}
  
  On fixe un ensemble de représentants ${\cal N}(\Gamma_{F})$ de l'ensemble des orbites nilpotentes conservées par l'action de $\Gamma_{F}$ formé d'éléments $N$ vérifiant les propriétés du lemme. Pour $N\in{\cal N}(\Gamma_{F})$, $N$ est fixé par $\Gamma_{F}^{nr}$ et 
 on peut identifier $A^{++}_{Fr}(N)$ à $A(N)$ par $(u,Fr)\mapsto u$.  Considérons l'ensemble des couples $(b,u)\in \tilde{A}(N)\times A(N)$ tels que $u Fr(b)=b^qu$. Le groupe $A(N)$ agit par conjugaison sur cet ensemble (parce que l'action de $Fr$ sur $A(N)$ est triviale). On note ${\cal C}^{\sharp}_{F}(N)$ l'ensemble quotient. Notons $\boldsymbol{{\cal O}}$ l'orbite contenant $N$. L'application qui à   $(b,u)\in {\cal C}^{\sharp}_{F}(N)$ associe le triplet $(N,b,(u,Fr))$ se quotiente en une bijection
 $$(3) \qquad  {\cal C}^{\sharp}_{F}(N)\to {\cal C}^{\sharp}_{F}(\boldsymbol{{\cal O}}).$$
 
 {\bf Remarques.}  (4) Supposons que $G$ soit  classique ou du type $(D_{4},3-ram)$. Alors la condition $uFr(b)=b^qu$ est toujours vérifiée. En effet, si $G$ est classique, $A^+(N)$ est produit de groupes ${\mathbb Z}/2{\mathbb Z}$. C'est un groupe commutatif et l'élévation à la puissance $q$ est l'identité. De plus l'action de $\Gamma_{F}^{nr}$ sur $A^+(N)$ est triviale (toutes ces propriétés résultent du lemme ci-dessus). L'assertion s'ensuit. Supposons $G$ de type $(D_{4},3-ram)$. Si $A(N)=\{1\}$, on a $u=1$ et $b$ est l'unique élément de $\tilde{A}(N)$. Alors $Fr(b)$ et $b^q$ sont égaux car ce sont les uniques éléments de $A^+(N)$ qui se projettent sur $\gamma^q$. Si $A(N)\not=\{1\}$, on a $A(N)={\mathbb Z}/2{\mathbb Z}$ et $A^+(N)={\mathbb Z}/2{\mathbb Z}\times {\mathbb Z}/3{\mathbb Z}$, l'action de Frobenius étant la multiplication par $q$ dans le deuxième facteur. Un calcul immédiat prouve l'assertion. 
 
  (5) Supposons que $G$ soit  exceptionnel.   Alors $\Gamma_{F}^{nr}$ agit trivialement sur $\tilde{A}(N)$. La condition $uFr(b)=b^qu$ équivaut donc à $ub=b^qu$. Dans le cas où $A(N)=\{1\}$, elle est automatique pour la m\^eme raison que dans le cas ci-dessus d'un groupe de type $(D_{4},3-ram)$. Si $A(N)\not=\{1\}$, on a $A^+(N)=A(N)\times \Gamma_{F^G/F}$. En écrivant $b=a\gamma$ avec $a\in A(N)$, la condition $ub=b^qu$ équivaut à $ua=a^qu$. 
  
  (6) Pour $N\in {\cal N}(\Gamma_{F})$, l'ensemble ${\cal C}^{\sharp}_{F}(N)$ n'est pas vide et l'application naturelle ${\cal C}^{\sharp}_{F}(N)\to \tilde{A}(N)$ est surjective. Cela résulte de la remarque (4) si $G$ est  classique ou du type $(D_{4},3-ram)$.  Supposons $G$ exceptionnel. D'après la remarque (5), il suffit de prouver que, pour tout $a\in A(N)$, il existe $u\in A(N)$ tel que $ua=a^qu$, autrement dit $a$ est conjugué à $a^q$ dans $A(N)$. Cela résulte d'un calcul facile compte tenu des faits que  le groupe $A(N)$ est isomorphe à $\mathfrak{S}_{i}$ pour un $i\in \{1,...,5\}$ et que $q$ est premier à $2,\,3,\,5$. 
  
  \bigskip
  
  Il y a une application naturelle ${\cal C}^{\sharp}_{F}\to (\boldsymbol{\mathfrak{g}}_{nil}/conj)^{\Gamma_{F}}$. La remarque (6) montre qu'elle est surjective.

 Soit $N$ un élément de $\boldsymbol{\mathfrak{g}}_{nil}$ dont l'orbite est conservée par $\Gamma_{F}$. Le sous-groupe $A^1(N)$ de 
 $A(N)$ est distingué dans $A^{++}(N)$. En effet, cette assertion est insensible à la conjugaison de $N$ par un élément de $\boldsymbol{G}$. On peut donc supposer $N\in {\cal N}(\Gamma_{F})$, auquel cas l'assertion équivaut à \ref{lesensemblesA+(N)}(1). On peut donc définir les quotients $\bar{A}^{++}(N)=A^{++}(N)/A^1(N)$, $\bar{A}^{++}_{\sigma}(N)=A^{++}_{\sigma}(N)/A^1(N)$ etc...En  remplaçant dans les constructions précédentes tous les ensembles $A^{++}(N)$ etc... par les quotientés $\bar{A}^{++}(N)$  etc..., on définit les ensembles $\bar{{\cal C}}^{\sharp}_{F}$, $\bar{{\cal C}}^{\sharp}_{F}(\boldsymbol{{\cal O}})$ et, pour $N\in {\cal N}(\Gamma_{F})$, l'ensemble $\bar{{\cal C}}^{\sharp}_{F}(N)$. On a pour ces ensembles l'analogue de la bijection (1). On a aussi des applications naturelles ${\cal C}^{\sharp}_{F}\to \bar{{\cal C}}^{\sharp}_{F}$, ${\cal C}^{\sharp}_{F}(N)\to \bar{{\cal C}}^{\sharp}_{F}(N)$. Montrons que
 
 (7) l'application  ${\cal C}^{\sharp}_{F}\to \bar{{\cal C}}^{\sharp}_{F}$ est surjective; pour $N\in {\cal N}(\Gamma_{F})$, l'application ${\cal C}^{\sharp}_{F}(N)\to \bar{{\cal C}}^{\sharp}_{F}(N)$ est surjective.
 
 Preuve. La première assertion résulte de la seconde. Soit $N\in {\cal N}(\Gamma_{F})$. Si $\bar{A}(N)=A(N)$, on a $\bar{{\cal C}}^{\sharp}_{F}(N)={\cal C}^{\sharp}_{F}(N)$ et l'assertion est triviale. Supposons $\bar{A}(N)\not=A(N)$. 
 Soit $(d,v)\in {\cal C}^{\sharp}_{F}(N)$. Si $G$ est classique ou de type $(D_{4},3-ram)$, on relève $(d,v)$ en un couple $(b,u)\in \tilde{A}(N)\times A(N)$. Ce couple appartient à ${\cal C}^{\sharp}_{F}(N)$ d'après la remarque (4). Supposons $G$ exceptionnel.  Si $\bar{A}(N)=\{1\}$, on relève $d$ en un élément $b\in \tilde{A}(N)$.  On a $A(N)\not=\{1\}$ puisqu'on a supposé $\bar{A}(N)\not=A(N)$. On écrit $b=a\gamma$, cf. remarque (5).   D'après cette remarque (5), on peut  choisir $u\in A(N)$ tel que $ua=a^qu$. Alors le couple $(b,u)$ appartient à $ {\cal C}^{\sharp}_{F}(N)$ et relève $(d,v)$. Supposons enfin que $G$ soit du type $E_{8}$, que $A(N)=\mathfrak{S}_{3}$ et $\bar{A}(N)=\mathfrak{S}_{2}$. On a $F^G=F$ et $\tilde{A}(N)=A(N)$, $\tilde{\bar{A}}(N)=\bar{A}(N)$.   En fixant une section $\bar{A}(N)\to A(N)$, on peut identifier $(d,v)$ à un couple $(b,u)\in \tilde{A}(N)\times A(N)$ dont on voit qu'il appartient à ${\cal C}^{\sharp}_{F}(N)$ et qu'il relève $(d,v)$. Cela prouve (7).

\subsubsection{Les applications $\iota_{{\cal F},nil}$ et $c_{{\cal F}}$\label{iotanil}}
L'hypoth\`ese que $G$ est adjoint implique que l'application $a_{T}:{\mathbb T}^{nr}\to App_{F^{nr}}(T^{nr})$ est injective. On a d\'efini une application $j_{T}:{\mathbb T}^{nr}\to {\bf T}$. On note ${\bf j}_{T}:a_{T}({\mathbb T}^{nr})\to {\bf T}$ l'application $j_{T}\circ a_{T}^{-1}$.

 Notons $T_{F}$ le plus grand sous-tore de $T$ d\'eploy\'e sur $F$.  Soit ${\cal F}$ une facette de l'appartement $App(T_{F})\subset Imm(G)$. 
 On a défini en \ref{discretisation}(3) l'homomorphisme 
 $\iota_{{\cal F}}:\mathfrak{g}_{{\cal F}}\to \boldsymbol{\mathfrak{g}}$. On en. déduit une application $\iota_{{\cal F},nil}:\mathfrak{g}_{{\cal F},nil}({\mathbb F}_{q})/conj\to \boldsymbol{\mathfrak{g}}_{nil}/conj$. Fixons un \'el\'ement  $x\in a_{T}({\mathbb T}_{F})\cap {\cal F}$.  L'image de $\iota_{{\cal F}}$ est l'ensemble des points fixes de $Ad({\bf j}_{T}(x))\circ\gamma$ dans $\boldsymbol{\mathfrak{g}}$. On en d\'eduit que $\iota_{{\cal F},nil}$ prend ses valeurs dans $(\boldsymbol{\mathfrak{g}}_{nil}/conj)^{I_{F}}$. D'autre part, $\iota_{{\cal F},nil}$ est équivariant 
  pour les actions de $\Gamma_{F}^{nr}$. Donc $\iota_{{\cal F},nil}$ prend ses valeurs dans $(\boldsymbol{\mathfrak{g}}_{nil}/conj)^{\Gamma_{F}}$.  
   
  Soit ${\cal O}\in \mathfrak{g}_{{\cal F},nil}({\mathbb F}_{q})/conj$. Fixons un élément $N'\in {\cal O}$ et un point $x\in a_{T}({\mathbb T}_{F})\cap {\cal F}$. Posons $N=\iota_{{\cal F}}(N')$.  Pour la m\^eme raison d'équivariance que ci-dessus,  $N$ est fixé par $\Gamma_{F}^{nr}$.  Il est aussi fixé par $Ad({\bf j}_{T}(x))\circ \gamma$.   On définit les éléments  $\dot{b}$ et $\dot{u}$ de $\boldsymbol{G}^{++}$ par $\dot{b}={\bf j}_{T}(x)\gamma$ et $\dot{u}=Fr$. Les propriétés précédentes entraînent que $\dot{b}\in Z_{\boldsymbol{G}^{++}_{\gamma}}(N)$ et $\dot{u}\in Z_{{\bf G}^{++}_{Fr}}(N)$. Notons $b$ et $u$ les images de ces éléments dans $A^{++}_{\gamma}(N)$ et $A^{++}_{Fr}(N)$. Montrons que
  
  (1) le triplet $(N,b,u)$ représente un élément de ${\cal C}^{\sharp}_{F}$; cet élément est indépendant des choix effectués. 
  
  Preuve. Parce que $x\in  a_{T}({\mathbb T}_{F})$, on a ${\bf j}_{T}(x)\in {\bf T}_{F}=X_{*}(T_{F})\otimes \bar{{\mathbb F}}_{q}^{\times}$. Donc ${\bf j}_{T}(x)$ est fixé par $I_{F}$ et on a $Fr({\bf j}_{T}(x))={\bf j}_{T}(x)^q$. Alors $\dot{u}\dot{b}=Fr({\bf j}_{T}(x))Fr\, \gamma={\bf j}_{T}(x)^q\gamma^q Fr=\dot{b}^q\dot{u}$. D'où $ub=b^qu$, ce qui prouve la première assertion. Les seuls choix sont ceux de $N'$ et $x$. Remplaçons $N'$ par un autre élément $\underline{N}'\in {\cal O}$. Posons $\underline{N}=\iota_{{\cal F}}(\underline{N}')$. Notre nouveau triplet est $(\underline{N},b,u)$. L'élément $\underline{N}'$ est  conjugué à $N'$ par un élément de $G_{{\cal F}}({\mathbb F}_{q})$. D'après la définition de $\iota_{{\cal F}}$ et son équivariance pour les actions de $\Gamma_{F}^{nr}$,  il existe  $g\in Z_{{\bf G}}({\bf j}_{T}(x)\gamma)^0$ tel que  $g$ soit  fixe par $Fr$ et $\underline{N}=gNg^{-1}$.  Alors $g$ conjugue le triplet $(N,b,u)$ en $(\underline{N},b,u)$ et les images dans ${\cal C}^{\sharp}_{F}$ de ces triplets sont égales. Remplaçons maintenant $x$   par un autre \'el\'ement $y\in a_{T}({\mathbb T}_{F})\cap {\cal F}$. On peut conserver le m\^eme \'el\'ement $N'$. Les éléments $N$ et $u$ ne changent pas. L'élément $\dot{b}$ est remplacé par $\underline{\dot{b}}={\bf j}_{T}(y)\gamma$. Posons $\tilde{x}=a_{T}^{-1}(x)$, $\tilde{y}=a_{T}^{-1}(y)$.  Choisissons un \'el\'ement $n\in {\mathbb P}$ tel que $n\tilde{x},n\tilde{y}\in X_{*}(T_{F})$ et $e^G$ divise $n$.  Soit $\lambda\in \mathfrak{o}_{F^{nr}}^{\times}$. Le m\^eme calcul que dans la preuve du lemme \ref{lesgroupesGF} montre que la conjugaison par $n(\tilde{x}-\tilde{y})(\lambda)$ se r\'eduit en  un automorphisme de $\boldsymbol{\mathfrak{g}}$ qui fixe l'image de $\iota_{{\cal F}}$.  Il en r\'esulte que, pour tout $\zeta\in \bar{{\mathbb F}}_{q}^{\times}$, l'\'el\'ement $n(\tilde{x}-\tilde{y})(\zeta)\in {\bf T}$ commute \`a $N$.  L'application $\zeta\mapsto n(\tilde{x}-\tilde{y})(\zeta)$ est  un sous-groupe \`a un param\`etre de  $Z_{{\bf G}}(N)$. Il est donc inclus dans $Z_{{\bf G}}(N)^0$.  On a les relations 
  $$\underline{\dot{b}}={\bf j}_{T}(y){\bf j}_{T}(x)^{-1}\dot{b}=(n\tilde{x}-n\tilde{y})(\zeta_{1/n})\dot{b}\in Z_{{\bf G}}(N)^0\dot{b}.$$
  Donc les images de $\underline{\dot{b}}$ et $\dot{b}$ dans $A^{++}_{\gamma}(N)$ sont égales et le triplet $(N,b,u)$ n'a pas changé. Cela achève la preuve de (1). 
  
  On définit l'application $c_{{\cal F},F}:\mathfrak{g}_{{\cal F},nil}({\mathbb F}_{q})/conj\to {\cal C}^{\sharp}_{F}$ qui, avec les notations ci-dessus, envoie ${\cal O}$ sur $(N,b,u)$ (ou plus exactement l'image de ce triplet dans ${\cal C}^{\sharp}_{F}$).  On définit l'application $\bar{c}_{{\cal F},F}:\mathfrak{g}_{{\cal F},nil}({\mathbb F}_{q})/conj\to \bar{{\cal C}}^{\sharp}_{F}$ qui envoie ${\cal O}$ sur l'image $(N,d,v)$ de $(N,b,u)$ dans $\bar{{\cal C}}^{\sharp}_{F}$.

 Rappelons que l'on a fixé un épinglage $\mathfrak{E}$ et effectué les constructions afférentes. On dispose donc d'une chambre $C^{nr}$ de $App_{F^{nr}}(T^{nr})$ et de la chambre $C\in App(T_{F})$ formée des éléments de $C^{nr}$ fixés par $\Gamma_{F}^{nr}$.   La   proposition suivante sera prouvée en \ref{surjectivite}.
 
 \begin{prop}{Quand $s$ décrit $S(\bar{C})$, la réunion des images des applications $\bar{c}_{s,F}$ est $\bar{{\cal C}}^{\sharp}_{F}$ tout entier.}\end{prop}
 
 \subsubsection{Les définitions sur $F^{nr}$}\label{surFnr}
  
  Les définitions des paragraphes précédents s'adaptent en se simplifiant si l'on remplace le corps de base $F$ par $F^{nr}$.  Il s'agit simplement d'y supprimer ce qui concerne les actions de Frobenius. Ainsi, on note ${\cal C}$ l'ensemble des classes de conjugaison par ${\bf G}$ dans l'ensemble des couples $(N,b)$ où $N$ est un élément de $\boldsymbol{\mathfrak{g}}_{nil}$ dont l'orbite est conservée par l'action de $I_{F}$ et $b$ appartient à $\tilde{A}(N)$. On fixe un ensemble de représentants ${\cal N}(I_{F})$ de l'ensemble des orbites  dans $\boldsymbol{\mathfrak{g}}_{nil}$ qui sont conservées par $I_{F}$, ces représentants vérifiant les conditions du lemme  \ref{calCF} où l'on oublie 
 celles concernant l'action de $\Gamma_{F}^{nr}$.  
  On suppose que cet ensemble ${\cal N}(I_{F})$  contient l'ensemble ${\cal N}(\Gamma_{F})$ déjà fixé. 
  
  {\bf Remarque.} Rappelons qu'un élément de ${\cal N}(I_{F})$ n'est pas forcément fixé par $I_{F}$.
  \bigskip
  
  Pour $N\in {\cal N}(I_{F})$, on note ${\cal C}(N)$ l'ensemble des classes de conjugaison par $A(N)$ dans $\tilde{A}(N)$. On a une bijection
  $$\bigsqcup_{N\in {\cal N}(I_{F})}{\cal C}(N)\to {\cal C}.$$
  Soit ${\cal F}$ une facette de $App_{F^{nr}}(T^{nr})$. On définit l'application $\iota_{{\cal F},nil}:\mathfrak{g}_{{\cal F},nil}/conj\to (\boldsymbol{\mathfrak{g}}_{nil}/conj)^{I_{F}}$, puis 
 une application $c_{{\cal F}}: \mathfrak{g}_{{\cal F},nil}/conj\to {\cal C}$ de m\^eme  qu'en \ref{iotanil},  en oubliant le terme $u$. 
 
 Soit ${\cal O}\in \mathfrak{g}_{{\cal F},nil}/conj$, posons $\boldsymbol{{\cal O}}=\iota_{{\cal F},nil}({\cal O})$. Notons $W_{{\cal F}}$ le groupe de Weyl de $G_{{\cal F}}$, qui s'identifie par $\iota_{{\cal F}}$ \`a un sous-groupe de $W$. 
Puisque $X_{*}(T_{{\cal F}})\simeq X_{*}(T^{nr})\simeq X_{*}({\bf T}^{nr})\subseteq X_{*}({\bf T})$,  l'\'el\'ement $x_{*,{\cal O}}$  peut \^etre identifi\'e \`a  une classe de conjugaison par $W_{{\cal F}}$ dans $X_{*}({\bf T}^{nr})$. L'\'el\'ement $x_{*,\boldsymbol{{\cal O}}}$  peut \^etre identifi\'e \`a une classe de conjugaison par $W$ dans $X_{*}({\bf T})$. 
 En appliquant  ce que l'on a dit en \ref{orbitesengendrees}, $x_{*,\boldsymbol{{\cal O}}}$ est la classe  de conjugaison par $W$ engendr\'ee par $x_{*,{\cal O}}$.

 En remplaçant les groupes $A(N)$ par $\bar{A}(N)$, on définit les variantes $\bar{{\cal C}}$, $\bar{{\cal C}}(N)$ et $\bar{c}_{{\cal F}}$ des objets précédents.

  \subsubsection{Description de $\iota_{{\cal F},nil}$ dans le cas ramifi\'e\label{iotanilram}}
 Dans ce paragraphe, on suppose que $G$ {\bf n'est pas d\'eploy\'e sur} $F^{nr}$. On consid\`ere un sommet $s\in S(\bar{C}^{nr})$ et on se propose de d\'ecrire l'application $\iota_{s,nil}$.  
 
   Supposons d'abord que $G$ est de type $(A_{n-1},ram)$.   L'ensemble $\boldsymbol{\mathfrak{g}}_{nil}/conj$  est en bijection avec   ${\cal P}(n)$. Tout \'el\'ement est conserv\'e par $I_{F}$. Si $n$ est impair, le  groupe $G_{s,SC}$ est produit de deux groupes de type $C_{n'}$ et $B_{n''}$ avec $2(n'+n'')+1=n$. L'ensemble $\mathfrak{g}_{s,nil}/conj$ est en bijection avec ${\cal P}^{symp}(2n')\times {\cal P}^{orth}(2n''+1)$. L'application $\iota_{s,nil}$ est la r\'eunion des deux  partitions. Si $n$ est pair, le  groupe $G_{s,SC}$ est produit de deux groupes de type $C_{n'}$ et $D_{n''}$ avec $2(n'+n'')=n$. L'ensemble $\mathfrak{g}_{s,nil}/conj$ s'envoie surjectivement sur ${\cal P}^{symp}(2n')\times {\cal P}^{orth}(2n'')$.    L'application $\iota_{s,nil}$ est compos\'ee de cette surjection et de la r\'eunion des partitions.
  
  Supposons maintenant que $G$ est de type $(D_{n},ram)$. L'ensemble $\boldsymbol{\mathfrak{g}}_{nil}/conj$  s'envoie surjectivement sur   ${\cal P}^{orth}(2n)$ et $(\boldsymbol{\mathfrak{g}}_{nil}/conj)^{I_{F}}$  est l'image r\'eciproque du sous-ensemble des partitions $\lambda\in {\cal P}^{orth}(2n)$ ayant au moins un terme impair.  En fait,  $(\boldsymbol{\mathfrak{g}}_{nil}/conj)^{I_{F}}$  est  en bijection avec cet ensemble de partitions.   Le groupe $G_{s,SC}$ est produit de deux groupes de type $B_{n'}$ et $B_{n''}$ avec $n'+n''=n-1$. L'ensemble $\mathfrak{g}_{s,nil}/conj$ est en bijection avec ${\cal P}^{orth}(2n'+1)\times {\cal P}^{orth}(2n''+1)$. L'application $\iota_{s,nil}$  est la r\'eunion des deux  partitions. 
  
  Les deux derniers cas sont $(D_{4},3-ram)$ et $(E_{6},ram)$. La description est moins \'el\'ementaire. On munit $X_{*}(T)$ du produit scalaire habituel pour lequel le carr\'e de la  norme d'une coracine vaut $2$ (ces coracines sont de m\^eme norme dans les deux cas). Il s'en d\'eduit des produits scalaires sur $X_{*}(T_{s})$ et $X_{*}({\bf T})$. Pour une orbite nilpotente ${\cal O}$ de $\mathfrak{g}_{s}$, resp. $\boldsymbol{\mathfrak{g}}$, la norme $\vert \vert x_{*,{\cal O}}\vert \vert $ est donc d\'efinie. D'apr\`es  \ref{orbitesengendrees}(1), l'application $\iota_{s,nil}$ conserve cette norme.   Il  s'av\`ere que, dans les deux cas, pour $\boldsymbol{{\cal O}}\in \boldsymbol{\mathfrak{g}}_{nil}/conj$, la norme de $x_{*,\boldsymbol{{\cal O}}}$ suffit \`a d\'eterminer $\boldsymbol{{\cal O}}$. Le calcul des \'el\'ements $x_{*,{\cal O}}$ et de leurs normes suffit donc \`a d\'eterminer la correspondance. 
  
  Supposons que $G$ est de type $(D_{4},3-ram)$. Pour $\boldsymbol{{\cal O}}\in \boldsymbol{\mathfrak{g}}_{nil}/conj$, on calcule $x_{*,\boldsymbol{{\cal O}}}$ sous la forme de Dynkin. L'orbite $\boldsymbol{{\cal O}}$ est conserv\'ee par $I_{F}$ si et seulement si $x_{*,\boldsymbol{{\cal O}}}$ est fix\'e par $\theta_{3}$, cf. \ref{actiongaloisienne}.   Dans le tableau suivant, la premi\`ere colonne d\'ecrit les \'el\'ements de $(\boldsymbol{\mathfrak{g}}_{nil}/conj)^{I_{F}}$ en termes de partitions. La deuxi\`eme colonne d\'ecrit les valeurs du carr\'e de la norme des $x_{*,\boldsymbol{{\cal O}}}$. Les trois colonnes suivantes d\'ecrivent les \'el\'ements de $\mathfrak{g}_{s,nil}/conj$ pour les trois sommets possibles. Ces sommets sont indiqu\'es \`a la premi\`ere ligne par l'\'el\'ement de $\Delta_{a}^{nr}$ auquel ils correspondent et, \`a la deuxi\`eme ligne, par le type du  groupe $G_{s,SC}$. Selon l'usage, la diff\'erence entre $\tilde{A}_{1}$ et $A_{1}$ signifie que la racine de la composante $\tilde{A}_{1}$ est courte tandis que celle de la composante $A_{1}$ est longue (les racines d'une composante $A_{2}$  sont courtes mais il n'y a pas d'ambigu\"{\i}t\'e possible). Les orbites nilpotentes de $\mathfrak{g}_{s,nil}/conj $ sont plac\'ees sur la ligne de   leur image dans  $\boldsymbol{\mathfrak{g}}_{nil}/conj$  par $\iota_{s,nil}$. Dans les deux derni\`eres colonnes, les composantes de $G_{s,SC}$ sont classiques et les orbites nilpotentes sont param\'etr\'ees par des partitions, plus exactement par des paires de partitions en ce qui concerne la quatri\`eme colonne. Pour la troisi\`eme colonne, le groupe $G_{s,SC}$ est de type $G_{2}$ et on utilise pour ses orbites les notations de \cite{C} p. 427. 
$$\begin{array}{ccccc}\boldsymbol{{\cal O}}&\vert \vert x_{*,\boldsymbol{{\cal O}}}\vert \vert ^2&\beta_{0}&\beta_{134}&\beta_{2}\\ &
&G_{2}&\tilde{A}_{1}\times A_{1}&A_{2}\\ &&&&\\ 
71&56&G_{2}&&\\53&24&&&3\\3^21^2&8&G_{2}(a_{1})&(2,2)&\\32^21&6&\tilde{A}_{1}&(2,1^2)&21\\2^21^4&2&A_{1}&(1^2,2)&\\ 1^8&0&1&(1^2,1^2)&1^3\\
\end{array}$$

Supposons que $G$ est de type $(E_{6},ram)$. On obtient alors le tableau suivant, qui s'explique comme le pr\'ec\'edent. 

$$\begin{array}{ccccccc}\boldsymbol{{\cal O}}&\vert \vert x_{*,\boldsymbol{{\cal O}}}\vert \vert ^2&\beta_{0}&\beta_{16}&\beta_{35}&\beta_{4}&\beta_{2}\\ &&

F_{4}&A_{1}\times B_{3}&\tilde{A}_{2}\times A_{2}&A_{3}\times A_{1}&C_{4}\\&&&&&&\\

E_{6}&312&F_{4}&&&&

\\E_{6}(a_{1})&168&&&&&8

\\D_{5}&120&F_{4}(a_{1})&&&&

\\E_{6}(a_{3})&72&F_{4}(a_{2})&&&&62

\\A_{5}&70&C_{3}&&&&61^2

\\D_{5}(a_{1})&60&&(2,7)&&&

\\ D_{4}&56&B_{3}&(1^2,7)&&&

\\ A_{4}+A_{1}&42&&&&(4,2)&

\\A_{4}&40&&&&(4,1^2)&4^2

\\D_{4}(a_{1})&24&F_{4}(a_{3})&(2,51^2)&(3,3)&&42^2

\\A_{3}+A_{1}&22&C_{3}(a_{1})&&&&421^2

\\ A_{3}&20&B_{2}&(1^2,51^2)&&&41^4

\\2A_{2}+A_{1}&18&\tilde{A}_{2}+A_{1}&&(3,21)&(31,2)&3^22

\\2A_{2}&16&\tilde{A}_{2}&&(3,1^3)&(31,1^2)&3^21^2

\\A_{2}+2A_{1}&12&A_{2}+\tilde{A}_{1}&(2,3^21)&(21,3)&&

\\A_{2}+A_{1}&10&&(2,32^2)&&(2^2,2)&

\\A_{2}&8&A_{2}&(2,31^4),(1^2,3^21)&(1^3,3)&(2^2,1^2)&2^4

\\3A_{1}&6&A_{1}+\tilde{A}_{1}&(2,2^21^3),(1^2,32^2)&(21,21)&(21^2,2)&2^31^2

\\2A_{1}&4&\tilde{A}_{1}&(2,1^7),(1^2,31^4)&(21,1^3)&(21^2,1^2)&2^21^4

\\A_{1}&2&A_{1}&(1^2,2^21^3)&(1^3,21)&(1^4,2)&21^6

\\1&0&1&(1^2,1^7)&(1^3,1^3)&(1^4,1^2)&1^8 \\ \end{array}$$

 \subsubsection{Preuve du lemme \ref{calCF} }\label{preuvedecalCF}

  Supposons que $G$ soit d\'eploy\'e sur $F^{nr}$. On a $F^G=F^{nr}$ et ${\bf G}^+={\bf G}$, d'o\`u  $A^+(N)=A(N)$ pour tout $N\in {\boldsymbol{\cal O}}$. Il est connu que, ${\boldsymbol{\cal O}}$ \'etant conserv\'e par $\Gamma_{F}^{nr}=\Gamma_{{\mathbb F}_{q}}$, il existe un \'el\'ement $N\in {\boldsymbol{\cal O}}$ qui est fix\'e par l'action de ce groupe et tel que l'action sur $A(N)$ est triviale, cf. \cite{Taylor} proposition 2.4. Cela prouve l'assertion (i) du lemme.
 
Supposons que $G$ soit des types $(A_{n-1},ram)$, $(D_{4},3-ram)$ ou $(E_{6},ram)$.   Si $A(N)=\{1\}$, on a \'evidemment $A^+(N)\simeq  \Gamma_{F^G/F^{nr}}$ pour tout $N\in \boldsymbol{{\cal O}}$. Puisque $\boldsymbol{{\cal O}}$ est fix\'ee par $\Gamma_{F}^{nr}$, cette orbite contient un \'el\'ement $N$ fix\'e par ce groupe. On traite maintenant les   ${\boldsymbol{\cal O}}$ tels que $A(N)\not=\{1
\}$ pour $N\in {\boldsymbol{\cal O}}$. Cela exclut le cas o\`u $G$ est de type $(A_{n-1},ram)$.

 Supposons que $G$ soit de type $(D_{4},3-ram)$.    Il existe une unique   orbite ${\boldsymbol{\cal O}}$ pour laquelle $A(N)\not=\{1\}$ pour $N\in {\boldsymbol{\cal O}}$:  elle est param\'etr\'ee par la partition $(3^21^2)$.  On a vu  en \ref{iotanilram} que l'orbite provenait par $\iota_{s_{0},nil}$ d'une unique  orbite de $\mathfrak{g}_{s_{0},nil}$.   On a identifi\'e ${\bf G}$ \`a $G_{s_{0},F^G}$ et  $\iota_{s_{0}}$ est le plongement naturel $G_{s_{0}}\to G_{s_{0},F^G}$. Le groupe  $G_{s_{0}}$ est la composante neutre du groupe des points fixes par $I_{F}$ dans ${\bf G}$. L'intersection ${\boldsymbol{\cal O}}\cap \mathfrak{g}_{s_{0}}$ est une orbite pour le groupe $G_{s_{0}}$ et est conserv\'ee par $\Gamma_{F}^{nr}$. Comme ci-dessus, elle contient un \'el\'ement $N$ fix\'e par ce groupe. L'\'el\'ement $N$ est aussi fix\'e par $I_{F}$, donc par tout le groupe $\Gamma_{F}$. On a $A(N)\simeq {\mathbb Z}/2{\mathbb Z}$ et l'action de $\Gamma_{F}$ sur ce groupe est forc\'ement triviale. 
 D'o\`u l'assertion  (ii) dans ce cas. 
  
  Supposons que $G$ soit de type $(E_{6},ram)$.  Il existe seulement trois   orbites ${\boldsymbol{\cal O}}$ pour lesquelles $A(N)\not=\{1\}$ pour $N\in {\boldsymbol{\cal O}}$:  ce sont les orbites  $A_{2}$ et $E_{6}(a_{3})$ pour lesquelles $A(N)={\mathbb Z}/2{\mathbb Z}$ pour $N\in {\boldsymbol{\cal O}}$ et  $D_{4}(a_{1})$ pour laquelle $A(N)=\mathfrak{S}_{3}$. On a vu  en \ref{iotanilram} que chacune de ces orbites provenait  par $\iota_{s_{0},nil}$ d'une unique orbite de $\mathfrak{g}_{s_{0},nil}$.  Dans le cas des orbites  $A_{2}$ et $E_{6}(a_{3})$, le m\^eme raisonnement que dans le cas $(D_{4},3-ram)$ conclut. Consid\'erons le cas de l'orbite $D_{4}(a_{1})$. Elle provient par $\iota_{s_{0},nil}$ de l'orbite ${\boldsymbol{\cal O}}_{0}$ de type $F_{4}(a_{3})$ de $\mathfrak{g}_{s_{0}}$.    Comme toujours, on affecte d'un exposant $G_{s_{0}}$ les objets relatifs \`a ce groupe. Par le m\^eme argument que dans le cas o\`u $G$ est d\'eploy\'e sur $F^{nr}$, on peut fixer $N\in {\boldsymbol{\cal O}}_{0}$ qui est fix\'e par $\Gamma_{F}^{nr}$ et tel que ce groupe agisse trivialement sur $A^{G_{s_{0}}}(N)$. L'action de $I_{F}$ sur ces objets est triviale puisqu'ils vivent dans $G_{s_{0}}$.
   Fixons un $\mathfrak{sl}(2)$-triplet $(f,h,N)$ de $\mathfrak{g}_{s_{0}}$. Le groupe $G_{s_{0}}$ est de type $F_{4}$ donc adjoint et on a $A^{G_{s_{0}}}(N)=Z_{G_{s_{0}}}(N)/Z_{G_{s_{0}}}(N)^0=Z_{G_{s_{0}}}(f,h,N)/Z_{G_{s_{0}}}(f,h,N)^0$. On a aussi $A(N)=Z_{{\bf G}}(f,h,N)/Z_{{\bf G}}(f,h,N)^0$ (en identifiant $f,h,N$ à des éléments de $\boldsymbol{\mathfrak{g}}$ par $\iota_{s_{0}}$).  
  De l'injection $Z_{G_{s_{0}}}(f,h,N)\to Z_{{\bf G}}(f,h,N)$ se d\'eduit un homomorphisme $A^{G_{s_{0}}}(N)\to A(N)$, dont le noyau est $(Z_{{\bf G}}(f,h,N)^0)\cap G_{s_{0}})/Z_{G_{s_{0}}}(f,h,N)^0$.  D'apr\`es \cite{C} pages 401-402, $A^{G_{s_{0}}}(N)$ est isomorphe \`a $\mathfrak{S}_{4}$ et $Z_{{\bf G}}(f,h,N)^0$ est commutatif. Il en r\'esulte que le noyau de l'homomorphisme $A^{G_{s_{0}}}(N)\to A(N)$ est un sous-groupe commutatif de $A^{G_{s_{0}}}(N)$. Or un sous-groupe commutatif de $\mathfrak{S}_{4}$ a au plus $4$ \'el\'ements. Il en r\'esulte que l'homomorphisme $A^{G_{s_{0}}}(N)\to A(N)$ est surjectif. Puisque l'action de $\Gamma_{F}$ sur $A^{G_{s_{0}}}(N)$ est triviale, il en est de m\^eme de l'action de ce groupe sur $A(N)$. Alors $N$ v\'erifie les conclusions du (ii) de l'\'enonc\'e.

   Supposons  que $G$ est de type $(D_{n},ram)$. On  a déjà décrit en \ref{lesensemblesA+(N)}(1) les ensembles $A(N)$ et $A^+(N)$. Cette description  montre que l'action de $\Gamma_{F}^{nr}$ sur $A^+(N)$ est triviale pour   tout  élément $N\in \boldsymbol{{\cal O}}$ qui  est fixé par $\Gamma_{F}^{nr}$ . Cela achève la preuve.

    \subsubsection{Sur le rel\`evement des commutants\label{relevementdescommutants}}
  
  Fixons $N \in \boldsymbol{\mathfrak{g}}_{nil}$ dont l'orbite est conserv\'ee par $I_{F}$. Fixons  un $\mathfrak{sl}(2)$-triplet $(f,h,N)$ de $\boldsymbol{\mathfrak{g}}$  dont le troisi\`eme terme est $N$. Posons ${\bf Z}=Z_{{\bf G}}(f,h,N)$, ${\bf Z}^+=Z_{{\bf G}^+}(f,h,N)$ et $\tilde{{\bf Z}}=Z_{\tilde{{\bf G}}}(f,h,N)$. On  a $A(N)={\bf Z}/{\bf Z}^0$, $A^+(N)={\bf Z}^+/{\bf Z}^0$ et $\tilde{A}(N)=\tilde{{\bf Z}}/{\bf Z}^0$. 

  \begin{lem}{(i) Supposons que $G$ ne soit pas de type $B_{n}$, $(D_{n}, nr)$ ou $(D_{n},ram)$. Soit $b\in \tilde{A}(N)$. Alors il existe un \'el\'ement semi-simple $\tau\in \tilde{{\bf Z}}$ dont l'image dans $\tilde{A}(N)$ est \'egale \`a $b$ et tel que l'image dans $A(N)$ de ${\bf Z}\cap  {\bf G}_{\tau}$ soit \'egale \`a $Z_{A(N)}(b)$ (on a pos\'e ${\bf G}_{\tau}=Z_{{\bf G}}(\tau)^0$).
  
  (ii)   Soit $d\in \tilde{\bar{A}}(N)$. Alors il existe un \'el\'ement semi-simple $\tau\in \tilde{{\bf Z}}$ dont l'image dans $\tilde{\bar{A}}(N)$ est \'egale \`a $d$ et tel que l'image dans $\bar{A}(N)$ de ${\bf Z}\cap {\bf G}_{\tau}$ soit \'egale \`a $Z_{\bar{A}(N)}(d)$.}\end{lem}
  
  Preuve.  Les assertions sont insensibles \`a la conjugaison de notre $\mathfrak{sl}(2)$-triplet (et cons\'equemment de $N$) par un \'el\'ement de ${\bf G}$. On se r\'eserve donc le droit de choisir convenablement ces objets. 

Traitons (i).  Il existe un \'el\'ement semi-simple $\tau\in \tilde{{\bf Z}}$ dont l'image dans $\tilde{A}(N)$ est \'egale \`a $b$. En effet, 
  fixons une paire de Borel $({\bf B}_{N},{\bf T}_{N})$ de $ {\bf Z}^0$. Puisque toutes les paires de Borel de ce groupe sont conjugu\'ees, toute composante connexe de ${\bf Z}^+$ contient un \'el\'ement dont l'action par conjugaison sur ${\bf Z}^0$  conserve $({\bf B}_{N},{\bf T}_{N})$. Un tel \'el\'ement est semi-simple. Si $A(N)=\{1\}$, cela suffit \`a prouver (i). Cela r\'esout le cas o\`u $G$ est de type $A_{n-1}$. 
  
  Supposons que $G$ soit de type $C_{n}$.   L'\'el\'ement $N$ est param\'etr\'e par une partition $\lambda\in {\cal P}^{symp}(2n)$ et $b$ appartient \`a $\{\pm 1\}^{Jord_{bp}(\lambda)}/\{1,\xi(\lambda)\}$. On rel\`eve $b$ en un \'el\'ement 
  $\dot{b}=(\dot{b}_{i})_{i\in Jord_{bp}(\lambda)}\in \{\pm 1\}^{Jord_{bp}(\lambda)}$. Notons $\lambda'$ la partition telle que, pour tout entier $i\geq1$, $mult_{\lambda'}(i)=1$ si $i\in Jord_{bp}(\lambda)$ et $\dot{b}_{i}=-1$ et $mult_{\lambda'}(i)=0$ sinon. Notons $\lambda''$ l'unique partition telle que $\lambda'\cup \lambda''=\lambda$. Les entiers $S(\lambda')$ et $S(\lambda'')$ sont  pairs, notons-les $2n'$ et $2n''$. On a $\lambda'\in {\cal P}^{symp}(2n')$ et $\lambda''\in {\cal P}^{symp}(2n'')$. On peut fixer un \'el\'ement semi-simple $\tau\in {\bf G}$ tel que ${\bf G}_{\tau,SC}$  soit isomorphe \`a $Sp(2n')\times Sp(2n'')$. Notons cette d\'ecomposition ${\bf G}'\times {\bf G}''$. Fixons des $\mathfrak{sl}(2)$-triplets $(f',h',N')$ de $\boldsymbol{\mathfrak{g}}'$ et $(f'',h'',N'')$ de $\boldsymbol{\mathfrak{g}}''$ tels que $N'$, resp. $N''$, soit param\'etr\'e par $\lambda'$, resp. $\lambda''$. On peut supposer que $(f,h,N)$ est la somme de ces deux triplets. On  a alors $\tau\in {\bf Z}$ et on v\'erifie que son image dans $A(N)$ est \'egale \`a $b$. Le groupe $Z_{{\bf  G}'}(f',h',N')\times Z_{{\bf  G}''}(f'',h'',N'')$ s'envoie dans ${\bf Z}$, son image est \'egale \`a ${\bf Z}\cap {\bf G}_{\tau}$ et on v\'erifie qu'il s'envoie
   surjectivement sur $A(N)$. Cela v\'erifie (i). 
  
    Supposons que $G$ soit exceptionnel et d\'eploy\'e sur $F^{nr}$. Alors $A(N)=\tilde{A}(N)$ est isomorphe \`a $\mathfrak{S}_{i}$ pour un entier $i=1,...,5$. On a d\'ej\`a traité le cas $i=1$. Si $A(N)=\mathfrak{S}_{4}$ ou $\mathfrak{S}_{5}$, on sait que  $G$ est de type $E_{8}$ ou $F_{4}$ donc est \`a la fois adjoint et simplement connexe. On sait aussi que ${\bf Z}^0=\{1\}$, cf. \cite{C} pages 401 et 406,  donc $\tilde{{\bf Z}}={\bf Z}=A(N)$. On pose $\tau=b$.  Il est \'evident que cet \'el\'ement est semi-simple et que ${\bf Z}\cap Z_{{\bf G}}(\tau)=Z_{A(N)}(b)$. Puisque $\tau$ est semi-simple et que $G$ est simplement connexe, le commutant $Z_{{\bf G}}(\tau)$ est connexe. L'\'egalit\'e pr\'ec\'edente  \'equivaut donc \`a ${\bf Z}\cap {\bf G}_{\tau}=Z_{A(N)}(b)$.  Supposons $A(N)=\mathfrak{S}_{i}$, avec $i=2,3$. Si $b=1$, l'\'el\'ement $\tau=1$ v\'erifie les conditions requises. Si $b\not=1$, $Z_{A(N)}(b)$ est le sous-groupe de $A(N)$ engendr\'e par $b$. On choisit un \'el\'ement semi-simple $\tau\in \tilde{{\bf Z}}={\bf Z}$ dont l'image dans $A(N)$ est \'egale \`a $b$. Parce que $\tau$ est semi-simple, il appartient \`a ${\bf G}_{\tau}$. Le groupe engendr\'e par $\tau$ est donc contenu dans ${\bf Z}\cap {\bf G}_{\tau}$ et s'envoie surjectivement sur 
   $Z_{A(N)}(b)$. 
   
    Supposons  que $G$ soit de type $(E_{6},ram)$.   On a d\'ej\`a trait\'e le cas $A(N)=\{1\}$. Il nous reste \`a consid\'erer les cas o\`u l'orbite de $N$ est de type $A_{2}$, $E_{6}(a_{3})$ ou $D_{4}(a_{1})$.  
      On a vu  en \ref{iotanilram} que l'on pouvait supposer $N\in \mathfrak{g}_{s_{0}}$, en identifiant $\mathfrak{g}_{s_{0}}$ à son image par $\iota_{s_{0}}$, autrement dit à l'ensemble des points fixes par $\gamma$ dans $\boldsymbol{\mathfrak{g}}$. On a alors $A^+(N)=A(N)\times \Gamma_{F^G/F^{nr}}$. On peut  m\^eme supposer que le $\mathfrak{sl}(2)$-triplet $(f,h,e)$ est contenu dans $\mathfrak{g}_{s_{0}}$.  On note par des exposants $G_{s_{0}}$ les objets relatifs au groupe $G_{s_{0}}$.  Démontrons  que
  
    (1) l'homomorphisme naturel ${\bf Z}^{G_{s_{0}}}\to A(N)$ est surjectif.

   Si l'orbite de $N$ dans $\boldsymbol{\mathfrak{g}}$ est de type $D_{4}(a_{1})$, on a prouvé cette assertion en \ref{preuvedecalCF}. 
  Supposons que cette orbite soit de type $E_{6}(a_{3})$. L'orbite de $N$ dans $\mathfrak{g}_{s_{0}}$ est  de type $F_{4}(a_{2})$. Il r\'esulte   de \cite{C} p. 401-402 que ${\bf Z}^0={\bf Z}^{G_{s_{0}},0}=\{1\}$ et que $A(N)=A^{G_{s_{0}}}(N)={\mathbb Z}/2{\mathbb Z}$.  Donc  $A(N)={\bf Z}$, $A^{G_{s_{0}}}(N)={\bf Z}^{G_{s_{0}}}$  et  ces groupes  ont deux \'el\'ements. Puisque l'homomorphisme ${\bf Z}^{G_{s_{0}}}\to {\bf Z}$ est injectif, il est surjectif, ce qui prouve (1). Supposons enfin que l'orbite de $N$ soit de type $A_{2}$. On a $A(N)={\mathbb Z}/2{\mathbb Z}$ et il suffit de prouver que ${\bf Z}^{G_{s_{0}}}$ n'est pas contenu dans ${\bf Z}^0$. Simplifions les notations en posant $L={\bf Z}^{G_{s_{0}}}$, $M={\bf Z}^0$. L'action de $\gamma$ se restreint en un automorphisme de $M$, noté encore $\gamma$, d'ordre au plus $2$ et on note $M^{\gamma}$ le sous-groupe des points fixes. Supposons par l'absurde que $L\subset M$.   Puisque $G_{s_{0}}$ est la composante neutre de ${\bf G}^{\gamma}$, on a l'inclusion $L\subset M^{\gamma}$ et l'égalité $L^0=M^{\gamma,0}$.  On a $L/L^0=A^{G_{s_{0}}}(N)={\mathbb Z}/2{\mathbb Z}$, donc $M^{\gamma}/M^{\gamma,0}$ contient un sous-groupe ${\mathbb Z}/2{\mathbb Z}$.   Les tables de \cite{C} p.401-402 décrivent la dimension de   $M$ et le type de  $M_{AD}$. On voit que    $M$ est semi-simple de type  $A_{2}\times A_{2}$. Notons $\boldsymbol{\zeta}_{3}$ le centre du groupe $SL(3)$. On peut identifier $M_{SC}$ à $SL(3)\times SL(3)$ et $M$ à $M_{SC}/\Xi$,  où $\Xi$ est un sous-groupe de $ \boldsymbol{\zeta}_{3}\times \boldsymbol{\zeta}_{3}$. L'automorphisme $\gamma$ se relève en un automorphisme de $ M_{SC}$, dont le groupe des points fixes $M_{SC}^{\gamma}$ est connexe car $SL(3)$  est simplement connexe. Le groupe $M^{\gamma,0}$ est l'image dans $M$ de $M_{SC}^{\gamma}$. Le groupe $M^{\gamma}$ est l'image dans $M$ du groupe des $h\in M_{SC}$  tels que $\gamma(h)\in h\Xi$. Parce que  $\xi^3=1$ pour tout $\xi\in \Xi$, tout élément $h$ de ce groupe  vérifie $h^3\in M_{SC}^{\gamma}$. Il en résulte que $m^3=1$ pour tout $m\in M^{\gamma}/M^{\gamma,0}$.
  Cela interdit à $M^{\gamma}/M^{\gamma,0}$ de contenir un sous-groupe isomorphe à ${\mathbb Z}/2{\mathbb Z}$. Cette contradiction  démontre (1). 
  
     Ecrivons $b=a\gamma$ avec $a\in A(N)$. Si $a=1$, on pose $\tau=\gamma$. Par d\'efinition, on a ${\bf G}_{\gamma}=Z_{{\bf G}}(\gamma)^0=G_{s_{0}}$, donc ${\bf Z}^{G_{s_{0}}}$ est  \'egal \`a ${\bf Z}\cap  {\bf G}_{\tau}$. D'apr\`es (1), $\tau$ v\'erifie la propri\'et\'e requise. Si $a\not=1$, on rel\`eve $a $ gr\^ace \`a (1) en un \'el\'ement  $x\in {\bf Z}^{G_{s_{0}}}$, que l'on peut supposer semi-simple par le m\^eme argument qu'au d\'ebut de la preuve.  On pose $\tau=x\gamma$. Cet \'el\'ement est semi-simple car son carr\'e vaut $x^2$ qui l'est. Notons $X$ le groupe engendr\'e par $x$. Parce que $x$ est semi-simple, $X$ est contenu dans  $(G_{s_{0}})_{x}$. Ce dernier groupe  commute \`a la fois \`a $x$ et \`a $\gamma$ (car il est contenu dans $G_{s_{0}}$) donc aussi \`a $\tau$.  Puisqu'il est connexe, on a $(G_{s_{0}})_{x}\subset {\bf G}_{\tau}$ puis
   $X\subset {\bf Z}\cap (G_{s_{0}})_{x}\subset  {\bf Z}\cap  {\bf G}_{\tau}$. Mais l'image de $X$ dans $A(N)$ contient le groupe engendr\'e par $a$. Ce groupe est le commutant de $b$ dans $A(N)$ puisque $a\not=1$ et que $A(N)=\mathfrak{S}_{2}$ ou $\mathfrak{S}_{3}$. Cela prouve l'assertion (i) de l'énoncé.

 Supposons  que $G$ soit de type $(D_{4},3-ram)$.  On a d\'ej\`a trait\'e le cas $A(N)=\{1\}$. Il nous reste \`a consid\'erer le cas o\`u l'orbite de $N$ est param\'etr\'ee par la partition $3^21^2$.   On a vu en \ref{iotanilram} que l'on pouvait supposer $N\in \mathfrak{g}_{s_{0}}$ et que l'orbite de $N$ pour l'action de $G_{s_{0}}$ est de type $G_{2}(a_{1})$.  On peut aussi supposer que le $\mathfrak{sl}(2)$-triplet $(f,h,N)$ est contenu dans $\mathfrak{g}_{s_{0}}$.    On sait que ${\bf Z}^{G_{s_{0}},0}=\{1\}$  et $A^{G_{s_{0}},0}(N)=\mathfrak{S}_{3}$, cf. \cite{C} p. 401, donc $ {\bf Z}^{G_{s_{0}}}=\mathfrak{ S}_{3}$. Le groupe ${\bf Z}^0$ est commutatif. Donc l'image de ${\bf Z}^{G_{s_{0}}}$ dans ${\bf Z}$ n'est pas contenue dans ${\bf Z}^0$, donc l'application naturelle ${\bf Z}^{G_{s_{0}}}\to A(N)={\mathbb Z}/2{\mathbb Z}$ est surjective. On peut alors construire $\tau$ comme dans le cas  o\`u $G$ est de type $(E_{6},ram)$.  Cela ach\`eve la preuve de (i).  
 
 Avant de traiter l'assertion (ii), prouvons
 
    (2) soit $d\in \tilde{\bar{A}}(N)$; alors il existe $b\in \tilde{A}(N)$ dont l'image dans $\tilde{\bar{A}}(N)$ soit $d$ et tel que l'image dans $\bar{A}(N)$ de $Z_{A(N)}(b)$ soit \'egale \`a $Z_{\bar{A}(N)}(d)$. 
   
   Preuve. Si $G$ est classique, le lemme \ref{calCF} implique que $A^+(N)$ est commutatif. Alors tout rel\`evement $b$ de $d$ dans $\tilde{A}(N)$ vérifie la condition requise. Supposons $G$ exceptionnel. Si $\bar{A}(N)=A(N)$, l'assertion est tautologique. Si $\bar{A}(N)=\{1\}$, pour tout rel\`evement $b$ de $d$ dans $\tilde{A}(N)$, $Z_{A(N)}(b)$ se projette surjectivement sur $\bar{A}(N)=\{1\}$. Comme on l'a dit en \ref{parametrageorbitesnilpotentes}, l'une ou l'autre des hypoth\`eses ci-dessus est toujours v\'erifi\'ee sauf dans le cas o\`u $G$ est de type $E_{8}$. Dans ce cas, il existe un $N$ tel que $A(N)=\mathfrak{S}_{3}$ et $\bar{A}(N)=\mathfrak{S}_{2}$. Le groupe $G$ est forc\'ement d\'eploy\'e donc $\tilde{A}(N)=A(N)$. On peut relever $d$ en un \'el\'ement $b$ d'ordre au plus $2$. Alors $Z_{A(N)}(b)$ se projette surjectivement sur $\bar{A}(N)$. Cela prouve (2).

 Traitons (ii). Supposons que $G$ ne soit pas de type $B_{n}$, $(D_{n},nr)$ ou $(D_{n},ram)$. On rel\`eve $d$ en un \'el\'ement $b\in \tilde{A}(N)$ v\'erifiant  (2), puis on applique (i) \`a cet \'el\'ement. L'\'el\'ement $\tau$ que l'on obtient v\'erifie la condition requise.
 
 Supposons que $G$ soit de type $B_{n}$. Identifions $G$ \`a un groupe  sp\'ecial orthogonal $SO(2n+1)$. L'\'el\'ement $N$ est param\'etr\'e par une partition $\lambda\in {\cal P}^{orth}(2n+1)$ et on a $d=(d_{i})_{i\in Int(\lambda)}\in \{\pm 1\}^{Int(\lambda)}_{1}$. Si $d=1$, l'\'el\'ement $\tau=1$ convient.  Supposons $d\not=1$. Notons $\delta:Jord_{bp}(\lambda)\to Int(\lambda)$ l'application qui, \`a $i\in Jord_{bp}(\lambda)$, associe l'intervalle qui contient $i$. Notons $K'$  l'ensemble des $\Delta\in Int(\lambda)$ tels que $d_{i}=-1$ et fixons un sous-ensemble $J'\subset Jord_{bp}(\lambda)$ tel que $\delta$ se restreigne en une bijection de $J'$ sur $K'$.    Notons $\lambda'$ la partition telle que, pour tout entier $i\geq1$, $mult_{\lambda'}(i)=1$ si $i\in J'$ et $mult_{\lambda'}(i)=0$ sinon. Notons $\lambda''$ l'unique partition telle que $\lambda'\cup \lambda''=\lambda$.  Parce que $d$ appartient \`a  $\{\pm 1\}^{Int(\lambda)}_{1}$, on voit que $S(\lambda')$ est pair, donc $S(\lambda'')$ est impair. On note ces entiers $2n'$ et $2n''+1$.   On a $\lambda'\in {\cal P}^{orth}(2n')$ et $\lambda''\in {\cal P}^{orth}(2n''+1)$. On peut fixer $\tau\in {\bf G}$ tel que ${\bf G}_{\tau}$  soit isomorphe \`a $SO(2n')\times SO(2n''+1)$. Notons cette d\'ecomposition ${\bf G}'\times {\bf G}''$. Fixons des $\mathfrak{sl}(2)$-triplets $(f',h',N')$ de $\boldsymbol{\mathfrak{g}}'$ et $(f'',h'',N'')$ de $\boldsymbol{\mathfrak{g}}''$ tels que $N'$, resp. $N''$, soit param\'etr\'e par $\lambda'$, resp. $\lambda''$. On peut supposer que $(f,h,N)$ est la somme de ces deux triplets. On  a alors $\tau\in {\bf Z}$ et on v\'erifie que son image dans $\tilde{\bar{A}}(N)=\bar{A}(N)$ est \'egale \`a $d$. On introduit les groupes ${\bf Z}'$ et ${\bf Z}''$ analogues \`a ${\bf Z}$, relatifs aux triplets $(f',h',N')$ et $(f'',h'',N'')$. On a ${\bf Z}'\times {\bf Z}''\subset {\bf Z}\cap  {\bf G}_{\tau}$. Il reste \`a prouver que
 
 (3) ${\bf Z}'\times {\bf Z}''$ s'envoie surjectivement sur $\bar{A}(N)$.
 
 On a $Jord_{bp}(\lambda')=J'$ et $\delta(J')=K'$. Posons $K''=\delta(Jord_{bp}(\lambda''))$. L'application de ${\bf Z}'$ dans $\bar{A}(N)$ se d\'ecompose en une suite d'applications
 $${\bf Z}'\to {\bf Z}'/{\bf Z}^{'0}=\{\pm 1\}^{J'}_{1}\to \{\pm 1\}^{K'}_{1}\subset \{\pm 1\}^{Int(\lambda)}_{1}.$$
 La derni\`ere inclusion consiste \`a prolonger tout \'el\'ement de $ \{\pm 1\}^{K'}_{1}$ par des $1$ pour les indices appartenant \`a $Int(\lambda)-K'$. C'est ce que l'on fera implicitement dans la suite. 
 L'image de ${\bf Z}'$ dans $\bar{A}(N)$ est $\{\pm 1\}^{K'}_{1}$. De m\^eme, l'image de ${\bf Z}''$ dans $\bar{A}(N)$ est $\{\pm 1\}^{K''}_{1}$.  On doit prouver que $\bar{A}(N)$ est produit de ces deux images. Par construction, on a $K'\cup K''=Int(\lambda)$. Montrons que l'on a
 
 (4) $K'\cap K''\not=\emptyset$.
 
 Rappelons que l'on pose $mult_{\lambda}(J)=\sum_{i\in J}mult_{\lambda}(i)$  pour tout sous-ensemble $J$ de $Jord_{bp}(\lambda)$.  Soit $\Delta\in Int(\lambda)$. Par construction, $mult_{\lambda}(\Delta)$ est pair sauf si $\Delta$ est le plus grand intervalle de $\lambda$. 
 Le nombre d'\'el\'ements de $K'$ est pair et non nul. Donc $K'$ contient un intervalle $\Delta$  tel que $mult_{\lambda}(\Delta)$ est pair et non nul. Par construction de $\lambda'$, on a $mult_{\lambda'}(\Delta)=1$.  On a donc $mult_{\lambda''}(\Delta)=mult_{\lambda}(\Delta)-mult_{\lambda'}(\Delta)>0$, ce qui \'equivaut \`a $\Delta\in K''$. Cela prouve (4). 
 
 Fixons $\Delta_{0}\in K'\cap K''$. Soit $\xi=(\xi_{\Delta})_{\Delta\in Int(\lambda)}\in \{\pm 1\}^{Int(\lambda)}_{1}$. Définissons  l'élément $\xi^{0}\in  \{\pm 1\}^{Int(\lambda)-\{\Delta_{0}\}}$ par  $\xi^0=(\xi_{\Delta})_{\Delta\in Int(\lambda)-\{\Delta_{0}\}} $. On peut \'evidemment fixer $\xi^{'0}=(\xi'_{\Delta})_{\Delta\in K'-\{\Delta_{0}\}}\in \{\pm 1\}^{K'-\{\Delta_{0}\}}$ et $\xi^{''0}=(\xi''_{\Delta})_{\Delta\in K''-\{\Delta_{0}\}}\in \{\pm 1\}^{K''-\{\Delta_{0}\}}$  de sorte que $\xi^{'0}\xi^{''0}=\xi^0$. En adjoignant un terme $\xi'_{\Delta_{0}}$ convenable, on peut prolonger la famille $\xi^{'0}$ en une famille $\xi'\in \{\pm 1\}^{K'}_{1}$. De m\^eme, on prolonge la famille $\xi^{''0}$ en une famille $\xi''\in \{\pm 1\}^{K''}_{1}$. Pour prouver (3), il suffit de  montrer que l'on a encore $\xi'\xi''=\xi$. Cela \'equivaut \`a l'\'egalit\'e $\xi'_{\Delta_{0}}\xi''_{\Delta_{0}}=\xi_{\Delta_{0}}$. Par hypoth\`ese, on a $\xi_{\Delta_{0}}=\prod_{\Delta\in Int(\lambda)-\{\Delta_{0}\}}\xi_{\Delta}$. Par construction, on a aussi $\xi'_{\Delta_{0}}\xi''_{\Delta_{0}}=\prod_{\Delta\in Int(\lambda)-\{\Delta_{0}\}}\xi'_{\Delta}\xi''_{\Delta}$. L'\'egalit\'e $\xi'_{\Delta_{0}}\xi''_{\Delta_{0}}=\xi_{\Delta_{0}}$ r\'esulte alors de l'\'egalit\'e $\xi^{'0}\xi^{''0}=\xi^0$. Cela ach\`eve la preuve de (3) et celle du (ii) de l'\'enonc\'e dans le cas o\`u $G$ est de type $B_{n}$. Les cas o\`u $G$ est de type $(D_{n},nr)$ ou $(D_{n},ram)$ se traitent de la m\^eme fa\c{c}on. $\square$

    \subsubsection{ Une construction auxiliaire\label{auxiliaire}} 
    
    Considérons l'ensemble des triplets $(N,b,R)$, où
    
    $N$ est un élément de $\boldsymbol{\mathfrak{g}}_{nil}$ dont l'orbite est conservée par $I_{F}$;
    
    $b\in \tilde{A}(N)$ et $R$ est un sous-groupe de $Z_{A(N)}(b)$.
    
    Le groupe ${\bf G}$ agit par conjugaison sur cet ensemble, notons ${\cal R}$ l'ensemble des classes de conjugaison. On note ${\cal R}_{max}$ le sous-ensemble des $(N,b,R)\in {\cal R}$ tels que $R=Z_{A(N)}(b)$.

 Soit ${\cal F}$ une facette de $App_{F^{nr}}(T^{nr})$ et soit ${\cal O}\in \mathfrak{g}_{{\cal F},nil}/conj$. On choisit un point $x\in a_{T}({\mathbb T}^{nr})\cap {\cal F}$.    Fixons 
$N'\in {\cal O}$. Posons $N=\iota_{{\cal F}}(N')$, $\dot{R}=\iota_{{\cal F}}(Z_{G_{{\cal F}}}(N'))$ et $\dot{b}={\bf j}_{T}(x)\gamma$. On a $b\in Z_{\tilde{{\bf G}}}(N)$, on note $b$ son image dans $\tilde{A}(N)$ et $R$ l'image de $\dot{R}$ dans $A(N)$. Le triplet $(N,b,R)$ représente un élément de ${\cal R}$. Comme en \ref{iotanil}, on voit que cet élément ne dépend pas des choix effectués.   On d\'efinit une application $r_{{\cal F}}:\mathfrak{g}_{{\cal F},nil}/conj\to {\cal R}$ par $r_{{\cal F}}({\cal O})=(N,b,R)$. 

On a les propri\'et\'es suivantes. 

(1) Il y a une application "oubli de $R$" de ${\cal R}$ dans ${\cal C}$. L'application $c_{{\cal F}}$ est la compos\'ee de $r_{{\cal F}}$ et de cette application d'oubli. La restriction de l'application d'oubli de ${\cal R}_{max}$ dans ${\cal C}$ est bijective: un élément $(N,b)\in {\cal C}$ se relève en l'unique triplet $(N,b,R)$ tel que $R=Z_{A(N)}(b)$. 

 (2) Soit $g\in G(F^{nr})$. Posons ${\cal F}'=g{\cal F}$ et supposons que ${\cal F}'$ soit elle-aussi contenue dans $App_{F^{nr}}(T^{nr})$. Comme en \ref{conjugaison}, on a un isomorphisme $\overline{Ad(g)}:\mathfrak{g}_{{\cal F}}\to \mathfrak{g}_{{\cal F}'}$.  Il r\'esulte du lemme \ref{conjugaison}  et de la construction que 
 $ r_{{\cal F}'}\circ \overline{ad(g)}=r_{{\cal F}}$. 

(3) Soit ${\cal F}'$ une  facette  de $App_{F^{nr}}(T^{nr})$   telle que ${\cal F}'\subset \bar{{\cal F}}$.  Comme en  \ref{adherence},  on peut identifier $\mathfrak{g}_{{\cal F}}$ \`a une sous-alg\`ebre de $\mathfrak{g}_{{\cal F}'}$. Il en r\'esulte une application $\mathfrak{g}_{{\cal F},nil}/conj\to \mathfrak{g}_{{\cal F}',nil}/conj$ qui est canonique. On voit que $c_{{\cal F}}$ est la compos\'ee de cette application et de $c_{{\cal F}'}$.  Il n'en est pas de m\^eme pour les applications $r_{{\cal F}}$ et $r_{{\cal F}'}$. En passant de ${\cal F}$ \`a ${\cal F}'$, les commutants peuvent augmenter. Pr\'ecis\'ement, soit ${\cal O}\in \mathfrak{g}_{{\cal F},nil}/conj$, notons ${\cal O}'$ son image dans $\mathfrak{g}_{{\cal F}',nil}/conj$. Posons $r_{{\cal F}}({\cal O})=(N,b,R)$. Alors $r_{{\cal F}'}({\cal O}')$ est de la forme $(N,b, R')$, o\`u $R\subset R'\subset Z_{A(N)}(b)$. 

\bigskip

En rempla\c{c}ant les ensembles $\tilde{A}(N)$ et $A(N)$ par $\tilde{\bar{A}}(N)$ et $\bar{A}(N)$, on d\'efinit de m\^eme les ensembles $\bar{{\cal R}}$ et $\bar{{\cal R}}_{max}$, ainsi qu'une application $\bar{r}_{{\cal F}}:\mathfrak{g}_{{\cal F},nil}/conj\to \bar{{\cal R}}$ pour toute facette ${\cal F}$ de $App_{F^{nr}}(T^{nr})$.   Il y a une application naturelle ${\cal R}\to \bar{{\cal R}}$: \`a $(N,b,R)\in {\cal R}$, elle associe $(N,d,V)$, o\`u $d$ est l'image de $b$ dans $\tilde{\bar{A}}(N)$ et $V$ est l'image de $R$ dans $\bar{A}(N)$. L'application d'oubli $\bar{{\cal R}}_{max}\to \bar{{\cal C}}$ est encore bijective.

 \begin{lem}{  
 (i) Supposons que $G$ n'est pas du type $B_{n}$, $(D_{n},nr)$ ou $(D_{n},ram)$. Quand $s$ d\'ecrit $S(\bar{C}^{nr})$, la r\'eunion des images des applications $r_{s}$ contient ${\cal R}_{max}$.
 
 (ii)  Quand $s$ d\'ecrit $S(\bar{C}^{nr})$, la r\'eunion des images des applications $\bar{r}_{s}$ contient $\bar{{\cal R}}_{max}$.

 (iii) Quand $s$ d\'ecrit $S(\bar{C}^{nr})$, la r\'eunion des images des applications $ c_{s}$ est ${\cal C}$ tout entier.

  }\end{lem}
 
 Preuve.  Commen\c{c}ons par prouver (i).    Fixons $(N,b,R)\in {\cal R}_{max}$. 
 Fixons  un $\mathfrak{sl}(2)$-triplet $(f,h,e)$ de $\boldsymbol{\mathfrak{g}}$ tel que $e=N$. Avec les notations du lemme \ref{relevementdescommutants}, celui-ci nous dit que l'on peut fixer un \'el\'ement 
  semi-simple $\tau\in \tilde{{\bf Z}}$ dont l'image dans $\tilde{A}(N)$ est \'egale \`a $b$ et tel que l'image dans $A(N)$ de ${\bf Z}\cap  {\bf G}_{\tau}$ soit \'egale \`a $Z_{A(N)}(b)=R$. 
  Dans la composante $\tilde{{\bf G}}$,   tout \'el\'ement semi-simple est conjugu\'e \`a un \'el\'ement $t\gamma$, o\`u $t\in {\bf T}^{I_{F}} $. Quitte à conjuguer $(N,b,R)$, on peut supposer que $\tau$ est un tel élément  $t\gamma$. 
   Puisque ${\bf G}$ est adjoint, ${\bf T}^{I_{F}}$ est connexe et donc \'egal \`a ${\bf T}^{nr}$.   Fixons $ x\in a_{T}({\mathbb T}^{nr})$ tel que ${\bf j}_{T}(x)=t$.  Notons ${\cal F}$ la facette de  l'appartement $App_{F^{nr}}(T^{nr})$ contenant $ x$. L'image de $\mathfrak{g}_{{\cal F}}$ par l'application $ \iota_{{\cal F}}$ est le sous-espace des \'el\'ements de $\boldsymbol{\mathfrak{g}}$ fix\'es par $Ad(\tau)$. L'image de $G_{{\cal F}}$ par cette application est ${\bf G}_{\tau}$. 
   En particulier, notre triplet $(f,h,e)$ appartient \`a l'image de $\mathfrak{g}_{{\cal F}}$ par  $ \iota_{{\cal F}}$. Notons $(f_{{\cal F}},h_{{\cal F}},e_{{\cal F}})$ son image r\'eciproque et $Z_{{\cal F}}$ le commutant de ce triplet dans $G_{{\cal F}}$. L'image de $Z_{{\cal F}}$ par $\iota_{{\cal F}}$ est \'egale \`a ${\bf Z}\cap  {\bf G}_{\tau}$.     Notons ${\cal O}$  l'image de $e_{{\cal F}}$ dans $\mathfrak{g}_{{\cal F},nil}/conj$.  D'après la construction de l'application $r_{{\cal F}}$, on a $ r_{{\cal F}}({\cal O})=(N,b, R)$. La facette ${\cal F}$ n'a pas de raison d'\^etre contenue dans $\bar{C}^{nr}$. Mais  son image par l'action d'un \'el\'ement convenable de $G(F^{nr})$ est contenue dans cette cl\^oture. La propri\'et\'e (2) nous permet de remplacer ${\cal F}$ par son image.   La propri\'et\'e (3) nous permet ensuite de remplacer ${\cal F}$ par un sommet adh\'erent \`a cette facette: cela ne peut qu'augmenter la troisi\`eme composante de $ r_{{\cal F}}({\cal O})$ mais, puisque celle-ci est d\'ej\`a maximale, elle ne change pas. 
  Un tel sommet appartient \`a $S(\bar{C}^{nr})$. Cela prouve le (i) de l'\'enonc\'e.  Le (ii) se prouve de fa\c{c}on analogue.
  Le (iii)  aussi: pour $(N,b)\in {\cal C}$, on choisit simplement un élément semi-simple $\tau\in \tilde{{\bf Z}}$ dont l'image dans $\tilde{A}(N)$ est égale à $b$, ce qui est toujours possible par l'argument donné au début de la preuve de 3.2.5.  
  $\square$

   \subsubsection{Description de $ c_{s}$ dans le cas ramifi\'e\label{csram}}
  {\bf On  suppose que }$G$ {\bf n'est pas d\'eploy\'e sur} $F^{nr}$. Soit $s\in \bar{C}^{nr}$. On note $\beta_{s}$ l'\'el\'ement de $\Delta_{a}^{nr}$ qui lui est associ\'e.  On  a d\'ecrit en \ref{iotanilram} l'application $\iota_{s,nil}$. 
    
  Dans le cas o\`u $G$ est de type $(A_{n-1},ram)$, $I_{F}$ agit trivialement sur $\boldsymbol{\mathfrak{g}}_{nil}/conj$. On a $A(N)=\{1\}$ pour  tout $N\in \boldsymbol{\mathfrak{g}}_{nil}$. Il n'y a qu'une classe de conjugaison dans $\tilde{A}(N)$, qui est r\'eduite au seul \'el\'ement de cet ensemble. Donc ${\cal C}\simeq  \boldsymbol{\mathfrak{g}}_{nil}/conj$ et $ c_{s}$ s'identifie \`a $\iota_{s,nil}$.
  
  Supposons que $G$ soit de type $(D_{n},ram)$.   Le groupe $G_{s,SC}$ est produit de deux groupes de types $B_{n'}$ et $B_{n''}$ avec $n'+n''=n-1$. On doit distinguer les deux composantes. Celles-ci correspondent aux composantes connexes du diagramme ${\cal D}_{a}^{nr}-\{\beta_{s}\}$. Si $\beta_{s}\not=\beta_{0}$,   on suppose que $B_{n''}$ correspond \`a la composante qui contient $\beta_{0}$.  Si  $\beta_{s}=\beta_{0}$, cette composante dispara\^{\i}t, $G_{s,SC}$ est quasi-simple de type $B_{n'}$ avec $n'=n-1$. On identifie  les orbites nilpotentes intervenant dans chacun de nos groupes  \`a des partitions.  
   Soient $\lambda'\in {\cal P}^{orth}(2n'+1)$ et $\lambda''\in {\cal P}^{orth}(2n''+1)$. Posons $\lambda=\iota_{s,nil}(\lambda',\lambda'')$. On a d\'ej\`a dit que $ \lambda=\lambda'\cup \lambda''$. On a $\lambda\in {\cal P}^{orth}(n)$. 
   L'ensemble $Jord_{bp}(\lambda)$ est non vide: chacune des deux partitions $\lambda'$ et $\lambda''$ est non vide et poss\`ede  un terme impair.    L'ensemble $\tilde{A}(\lambda)$ s'identifie \`a   $\{\pm 1\}^{Jord_{bp}(\lambda)}_{-1}/\{1,\xi(\lambda)\}$, cf. \ref{lesensemblesA+(N)}(2). Notons $\dot{b}(\lambda',\lambda'')=(\dot{b}(\lambda',\lambda'')_{i})_{i\in Jord_{bp}(\lambda)}$ l'\'el\'ement  de $ \{\pm 1\}^{Jord_{bp}(\lambda)}$ 
   d\'efini par $\dot{b}(\lambda',\lambda'')_{i}=(-1)^{mult_{\lambda''}(i)}$. On v\'erifie que $\dot{b}(\lambda',\lambda'')\in \{\pm 1\}^{Jord_{bp}(\lambda)}_{-1}$ et on note $b(\lambda',\lambda'')$ son image dans $\tilde{A}(\lambda)$. Il est facile de voir  que 
$ c_{s}(\lambda',\lambda')=(\lambda,b(\lambda',\lambda''))$. 

Consid\'erons maintenant les deux cas restants: $G$ de type $(D_{4},3-ram)$  ou $(E_{6},ram)$. On identifie $(\boldsymbol{\mathfrak{g}}_{nil}/conj)^{I_{F}}$ à l'ensemble de représentants ${\cal N}(I_{F})$ fixé en \ref{surFnr}. 
Soit   ${\cal O}_{s}\in \mathfrak{g}_{s,nil}/conj$, posons $N=\iota_{s,nil}({\cal O}_{s})$.   Comme dans le cas o\`u $G$ est de type $(A_{n-1},ram)$, si $A(N)=\{1\}$, on n'a rien \`a ajouter pour d\'ecrire $ c_{s}({\cal O}_{s})$ puisqu'il n'y a qu'une classe de conjugaison dans $\tilde{A}(N)$, qui est r\'eduite  \`a l'unique \'el\'ement de cet ensemble. On n'a donc \`a se pr\'eoccuper que des $N\in {\cal N}(I_{F})$ tels que  $A(N)\not=\{1\}$. D'apr\`es  \ref{iotanilram} et \cite{C} p. 402,  $N$ est fix\'e par $I_{F}$ et appartient \`a l'image de $\iota_{s_{0},nil}$. Plus pr\'ecis\'ement, il existe une  unique orbite ${\cal O}_{0}\in  \mathfrak{g}_{s_{0},nil}/conj$ telle que $\iota_{s_{0},nil}({\cal O}_{0})=N$.   On note   $b_{1}$ la classe de conjugaison de $\gamma$ dans $\tilde{A}(N)=A(N)\gamma$.  On a $s_{0}=a_{T}(0)$, ${\bf j}_{T}(s_{0})=1$ et il r\'esulte des constructions que  $ c_{s_{0}}({\cal O}_{0})=(N,b_{1})$.

 Supposons $G$ de type $(D_{4},3-ram)$.  Il n'y a  qu'un \'el\'ement $N\in {\cal N}(I_{F})$ tel que $A(N)\not=\{1\}$. Il est param\'etr\'e par  la partition $3^21^2$ et on a $A(N)={\mathbb Z}/2{\mathbb Z}$.  D'apr\`es \ref{iotanilram}, il y a  deux sommets $s$ tels que $N$ appartienne \`a l'image de $ \iota_{s,nil}$, \`a savoir $s_{0}$ et $s_{134}$.  De plus, pour $i=0,134$, il y a une seule orbite ${\cal O}_{i}\in \mathfrak{g}_{s_{i},nil}/conj$ telle que $\iota_{s_{i},nil}({\cal O}_{i})=N$.   Il n'y a que deux classes de conjugaison dans $\tilde{A}(N)$: la classe de  $(0,\gamma)\in \tilde{A}(N)={\mathbb Z}/2{\mathbb Z}\times \{\gamma\} $ que l'on a not\'ee $b_{1}$ et celle de $(1,\gamma) $ que l'on note $b_{2}$. On a calcul\'e ci-dessus $c_{s_{0}}({\cal O}_{0})=(N,b_{1})$. L'assertion (iii) du lemme \ref{auxiliaire} implique que $c_{s_{134}}({\cal O}_{134})=(N,b_{2})$. 
   
 Supposons maintenant $G$ de type $(E_{6},ram)$.   Rappelons que $\Delta_{a}^{nr}=\{\beta_{0},\beta_{16},\beta_{35},\beta_{4},\beta_{2}\}$, avec $\beta_{0}=-(2\beta_{16}+3\beta_{35}+2\beta_{4}+\beta_{2})$.   On introduit la base  $\{\check{\varpi}_{i}; i=1,...,6\}$ de $X_{*}(T)$ duale de $\Delta$ et les \'el\'ements $\check{\varpi}_{16}=\check{\varpi}_{1}+\check{\varpi}_{6}$, $\check{\varpi}_{35}=\check{\varpi}_{3}+\check{\varpi}_{5}$. Alors $\check{\Pi}^{nr}=\{\check{\varpi}_{16},\check{\varpi}_{35},\check{\varpi}_{4},\check{\varpi}_{2}\}$ est une base de $X_{*}(T^{nr})$. La racine affine $\beta_{0}^{aff}$ est $\beta_{0}+1/2$.

  Il y a trois \'el\'ements $N\in {\cal N}(I_{F})$ tels que $A(N)\not=\{1\}$, dont les orbites sont  not\'ees $E_{6}(a_{3})$, $D_{4}(a_{1})$ et $A_{2}$, cf. \cite{C} page 402. Pour $N$ de type $E_{6}(A_{3})$, on a $A(N)={\mathbb Z}/2{\mathbb Z}$ et la situation est la m\^eme que dans le cas $(D_{4},3-ram)$: $N$ provient des deux sommets $s_{0}$ et $s_{2}$  et d'une seule classe pour chaque sommet.  Avec les m\^emes notations que ci-dessus, on obtient $ c_{s_{2}}({\cal O}_{2})=(N,b_{2})$.
 
  Supposons maintenant $N$ de type $D_{4}(a_{1})$. Alors $A(N)=S_{3}$ et $A^+(N)= A(N)\times \{1,\gamma\}$. Il y a donc trois classes de conjugaison dans $\tilde{A}(N)$: les classes $b_{1}$ de $(a_{1},\gamma)$, $b_{2}$ de $(a_{2},\gamma)$ et $b_{3}$ de $(a_{3},\gamma)$, o\`u $a_{1}$, resp. $a_{2}$, resp. $a_{3}$, est un \'el\'ement d'ordre $1$ resp. $2$, resp. $3$, de $A(N)$.  
 Il y a quatre sommets  $s$ tels que $N$ appartienne \`a l'image de $\iota_{s,nil}$: $s_{0}$, $s_{16}$, $s_{35}$, $s_{2}$. On introduit les \'el\'ements suivants de ${\mathbb T}^{nr}$:  $\tilde{x}_{0}=0$, $\tilde{x}_{16}=\frac{1}{4}\check{\varpi}_{16}$, $\tilde{x}_{35}=\frac{1}{6}\check{\varpi}_{35}$, $\tilde{x}_{2}=\frac{1}{2}\check{\varpi}_{2}$. On voit que $a_{T}(\tilde{x}_{0})=s_{0}$, $a_{T}(\tilde{x}_{16})=s_{16}$, $a_{T}(\tilde{x}_{35})=s_{35}$, $a_{T}(\tilde{x}_{2})=s_{2}$. Pour chacun des indices $i=0,16,35,2$, on pose $t_{i}={\bf j}_{T}(s_{i})$. On calcule $t_{0}=1$,    $t_{16}=\check{\varpi}_{16}(\zeta_{1/4})^{-1}$, $t_{35}=\check{\varpi}_{35}(\zeta_{1/6})^{-1}$,   $t_{2}=\check{\varpi}_{2}(\zeta_{1/2})^{-1}$. Il 
 y a une unique orbite ${\cal O}_{i}\in \mathfrak{g}_{s_{i},nil}/conj$ telle que $ \iota_{s_{i},nil}({\cal O}_{i})=N$.  Posons $c_{s_{i}}({\cal O}_{i})=(N,b(s_{i}))$. Par construction, l'image r\'eciproque de $b(s_{i})$ dans $Z_{\tilde{{\bf G}}}(N)$ contient un \'el\'ement conjugu\'e \`a $t_{i}\gamma$.  L'\'el\'ement $t_{i}\in {\bf T}^{nr}$ commute \`a $\gamma$. D'o\`u $(t_{i}\gamma)^2=t_{i}^2$. On en d\'eduit que l'\'el\'ement $t_{i}\gamma$ est d'ordre $2$ dans ${\bf G}^+$ si $i=0$ ou $i=2$, d'ordre $4$ si $i=16$ et d'ordre $6$ si $i=35$.  L'assertion (iii) du lemme \ref{auxiliaire} nous dit qu'il y a au moins un indice $i$  tel que $b(s_{i})=b_{3}$. Pour cet indice,   l'ordre de $t_{i}\gamma$ est divisible par $3$.  Ce ne peut \^etre que  l'indice $35$, donc $b(s_{35})=b_{3}$. On a d\'ej\`a calcul\'e $b(s_{0})=b_{1}$. Pour les deux autres indices $i=16, 2$,   on a forc\'ement $b(s_{i})=b_{1}$ ou $b_{2}$.  Prouvons que

 (1) $b(s_{16})=b(s_{2})=b_{2}$.
 
 Pour chacun des indices $i=0,16,2$, fixons un $\mathfrak{sl}(2)$-triplet $(f'_{i},h'_{i},N'_{i})$ de $\mathfrak{g}_{s_{i}}$ tel que $N'_{i}\in {\cal O}_{i}$. Posons $(f_{i},h_{i},N_{i})=\iota_{s_{i}}(f'_{i},h'_{i},N'_{i})$. Par construction, $t_{i}\gamma$ appartient \`a $Z_{\tilde{{\bf G}}}(f_{i},h_{i},N_{i})$ et $\iota_{s_{i}}$ identifie $Z_{G_{s_{i}}}(f'_{i},h'_{i},N'_{i})^0$ \`a la composante neutre du commutant de $t_{i}\gamma$ dans $Z_{{\bf G}}(f_{i},h_{i},N_{i})^0$.  On ne perd rien \`a supposer que $N=N_{0}$. On a $t_{0}=1$, donc $Z_{G_{s_{0}}}(f'_{0},h'_{0},N'_{0})^0$ s'identifie \`a la composante neutre de l'action de $\gamma$ dans $Z_{{\bf G}}(f_{0},h_{0},N_{0})^0$.  Consid\'erons un indice $i=16,2$ et supposons $b(s_{i})=b_{1}$. Les ensembles $Z_{\tilde{{\bf G}}}(f_{i},h_{i},N_{i})$ et $Z_{\tilde{{\bf G}}}(f_{0},h_{0},N_{0})$ sont conjugu\'es par un \'el\'ement de ${\bf G}$ et l'hypoth\`ese pr\'ec\'edente signifie que l'on peut imposer \`a  cette conjugaison d'envoyer $t_{i}\gamma$ sur un \'el\'ement $z\gamma$, o\`u  $z\in Z_{{\bf G}}(f_{0},h_{0},N_{0})^0$. Alors $Z_{G_{s_{i}}}(f'_{i},h'_{i},N'_{i})^0$ s'identifie \`a la composante neutre du commutant de $z\gamma$ dans  $Z_{{\bf G}}(f_{0},h_{0},N_{0})^0$. Or, d'apr\`es \cite{C} page 402,  $Z_{{\bf G}}(f_{0},h_{0},N_{0})^0$ est 
  un tore de dimension $2$ et est donc commutatif. Le commutant de $z\gamma$ est donc le m\^eme que celui de $\gamma$ et on obtient un isomorphisme  $Z_{G_{s_{i}}}(f'_{i},h'_{i},N'_{i})^0\simeq Z_{G_{s_{0}}}(f'_{0},h'_{0},N'_{0})^0$. L'orbite ${\cal O}_{0}$ est de type $F_{4}(a_{3})$. D'apr\`es \cite{C} page 401, on a $Z_{G_{s_{0}}}(f'_{0},h'_{0},N'_{0})^0=\{1\}$. Si $i=16$, l'orbite ${\cal O}_{16}$ est param\'etr\'ee par $(2,51^2)$; si $i=2$, ${\cal O}_{2}$ est param\'etr\'ee par $(42^2)$. Dans les deux cas, $Z_{G_{s_{i}}}(f'_{i},h'_{i},N'_{i})^0$ est un tore de dimension $1$. Cela contredit l'\'egalit\'e pr\'ec\'edente. Notre hypoth\`ese $b(s_{i})=b_{1}$ est donc contradictoire, ce qui prouve (1).

  Supposons enfin $N$ de type $A_{2}$. On a $A(N)={\mathbb Z}/2{\mathbb Z}$ et $A^+(N)=A(N)\times \{1,\gamma\}$. Il y a deux classes de conjugaison dans $\tilde{A}(N)$:  la classe $b_{1}$ de $(0,\gamma)$ et la classe $b_{2}$ de $(1,\gamma)$. L'\'el\'ement  $N$ appartient \`a l'image de $\iota_{s,nil}$ pour tous les sommets $s\in S(\bar{C}^{nr})$. Sauf pour  l'indice $16$, il y a  une unique orbite ${\cal O}_{i}\in \mathfrak{g}_{s_{i},nil}/conj$ telle que $\iota_{s_{i},nil}({\cal O}_{i})=N$.  On pose alors $c_{s_{i}}({\cal O}_{i})=(N,b(s_{i}))$.  Pour l'indice $16$, il y en a deux: notons ${\cal O}'_{16}$ l'orbite param\'etr\'ee par  $(2,31^4)$ et ${\cal O}''_{16}$ l'orbite param\'etr\'ee par  $(1^2,3^21)$. On pose $c_{s_{16}}({\cal O}'_{16})=(N,b'(s_{16}))$ et $c_{s_{16}}({\cal O}''_{16})=(N,b''(s_{16}))$. Montrons que
 
 (2) $b(s_{0})=b''(s_{16})=b(s_{35})=b_{1}$, $b'(s_{16})=b(s_{4})=b(s_{2})=b_{2}$. 
 
 Il suffit de prouver que $b(s_{0})=b''(s_{16})=b(s_{35})$ et $b'(s_{16})=b(s_{4})=b(s_{2})$. En effet, on sait d\'ej\`a que $b(s_{0})=b_{1}$, d'o\`u les premi\`eres \'egalit\'es de (2). L'assertion (iii) du lemme \ref{auxiliaire} entra\^{\i}ne alors que les trois autres classes sont forc\'ement \'egales \`a $b_{2}$. 
 
 Introduisons la facette ${\cal F}\subset \bar{C}^{nr}$ telle que les \'el\'ements de $S(\bar{C}^{nr})$ adh\'erents \`a ${\cal F}$ soient $s_{0},s_{16},s_{35}$. Le groupe $G_{{\cal F},SC}$ est de type $A_{2}$. Consid\'erons l'orbite nilpotente r\'eguli\`ere ${\cal O}_{{\cal F}}\subset \mathfrak{g}_{{\cal F},nil}$ param\'etr\'ee par la partition $3$. On a des plongements 
 
 $$\begin{array}{ccc}&&\mathfrak{g}_{s_{0}}\\ &\nearrow&\\\mathfrak{g}_{{\cal F}}&\to&\mathfrak{g}_{s_{16}}\\ &\searrow&\\ &&\mathfrak{g}_{s_{35}}\\ \end{array}$$
 On v\'erifie que l'image de ${\cal O}_{{\cal F}}$ par le deuxi\`eme plongement est ${\cal O}''_{16}$ et, par le troisi\`eme plongement, est ${\cal O}_{35}$. En vertu de la compatibilit\'e des applications $\iota_{s,nil}$, l'image de ${\cal O}_{{\cal F}}$ par le premier plongement est une orbite ${\cal O}'_{0}\in \mathfrak{g}_{s_{0},nil}/conj$ telle que $$\iota_{s_{0},nil}({\cal O}'_{0})=\iota_{{\cal F},nil}({\cal O}_{{\cal F}})=\iota_{s_{35},nil}({\cal O}_{35})=N.$$
 Il n'y a qu'une seule telle orbite \`a savoir ${\cal O}_{0}$, donc ${\cal O}'_{0}={\cal O}_{0}$. D'apr\`es \ref{auxiliaire} (3), on obtient alors
 $$\begin{array}{ccccc}\\ &&&&c_{s_{16}}({\cal O}''_{16})\\ &&&&\parallel\\ (N,b_{1})&=&c_{s_{0}}({\cal O}_{0})&=&c_{{\cal F}}({\cal O}_{{\cal F}})\\ &&&&\parallel\\&&&&c_{s_{35}}({\cal O}_{35})\\ \end{array}$$
 Cela prouve les premi\`eres \'egalit\'es voulues.
 
 Oublions la facette ${\cal F}$ pr\'ec\'edente et notons maintenant ${\cal F}$ la facette  telle que les \'el\'ements de $S(\bar{C}^{nr})$ adh\'erents \`a ${\cal F}$ soient  $s_{16}$, $s_{4}$, $s_{2}$. Le groupe $G_{{\cal F},SC}$ est de type $A_{1}\times A_{1}$. Consid\'erons l'orbite ${\cal O}_{{\cal F}}$ dans $\mathfrak{g}_{{\cal F},nil}$ form\'ee des couples d'\'el\'ements r\'eguliers dans les deux composantes. On a comme ci-dessus des plongements
  $$\begin{array}{ccc}&&\mathfrak{g}_{s_{16}}\\ &\nearrow&\\\mathfrak{g}_{{\cal F}}&\to&\mathfrak{g}_{s_{4}}\\ &\searrow&\\ &&\mathfrak{g}_{s_{2}}\\ \end{array}$$
  On v\'erifie que les images de ${\cal O}_{{\cal F}}$ par ces plongements sont respectivement ${\cal O}'_{16}$, ${\cal O}_{4}$, ${\cal O}_{2}$. Le m\^eme raisonnement que ci-dessus conduit aux \'egalit\'es  
  $$\begin{array}{ccc}\\ &&c_{s_{16}}({\cal O}'_{16})\\ &&\parallel\\ c_{s_{4}}({\cal O}_{4})&=&c_{{\cal F}}({\cal O}_{{\cal F}})\\ &&\parallel\\&&c_{s_{2}}({\cal O}_{2})\\ \end{array}$$
   Cela prouve les deuxi\`emes \'egalit\'es voulues, d'o\`u (2). 
   
     \subsubsection{ Preuve de la proposition \ref{iotanil}} \label{surjectivite}
  
    Soit $(N,d,v)\in \bar{{\cal C}}^{\sharp}_{F}$.  Posons $V=Z_{\bar{A}(N)}(d)$. D'apr\`es l'assertion (ii) du lemme \ref{auxiliaire}, on peut fixer $s^{nr}\in S(\bar{C}^{nr})$ et ${\cal O}^{nr}\in \mathfrak{g}_{s^{nr},nil}/conj$ de sorte que $\bar{r}_{s^{nr}}({\cal O}^{nr})=(N,d,V)$. Supposons les conditions suivantes v\'erifi\'ees: 
 
 (1)  $s^{nr}$ est fix\'e par l'action de $\Gamma_{F}^{nr}$ et   ${\cal O}^{nr}$ est fix\'ee par l'action de $\Gamma_{{\mathbb F}_{q}}$ sur $\mathfrak{g}_{s^{nr},nil}/conj$. 

  Alors $s^{nr}$ appartient \`a $S(\bar{C})$, on pose $s=s^{nr}$.   Puisque ${\cal O}^{nr}$ est fix\'ee par l'action de $\Gamma_{{\mathbb F}_{q}}$, l'ensemble ${\cal O}^{nr}\cap  \mathfrak{g}_{s,nil}({\mathbb F}_{q})$ n'est pas vide et on peut fixer une orbite ${\cal O}_{1}\in \mathfrak{g}_{s,nil}({\mathbb F}_{q})/conj$ contenue dans cette intersection.  Fixons  $N'_{1}\in {\cal O}_{1}$. Par définition de $\bar{r}_{s^{nr}}$, on peut supposer, quitte à conjuguer $(N,d,v)$, que $N=\iota_{s}(N'_{1})$ et que $d$ est l'image dans $\tilde{\bar{A}}(N)$ de l'élément $\dot{d}={\bf j}_{T}(s)\gamma$. Par construction, $\bar{c}_{s,F}({\cal O}'_{1})=(N,d,v_{1})$, où $v_{1}$ est l'image de $Fr$ dans $\bar{A}^{++}_{Fr}(N)$. L'élément $v$ appartient au m\^eme ensemble $\bar{A}^{++}_{Fr}(N)$. Donc $v_{1}^{-1}v$ appartient à $\bar{A}(N)$. On a aussi les égalités $vdv^{-1}=d^q=v_{1}dv_{1}^{-1}.$. Donc $v_{1}^{-1}v$ appartient à $Z_{\bar{A}(N)}(d)=V$. L'\'egalit\'e $\bar{r}_{s^{nr}}({\cal O}^{nr})=(N,d,V)$ signifie que $V$ est l'image dans $\bar{A}(N)$ du groupe $\iota_{s}(Z_{G_{s}}(N'_{1}))$. Soit $z\in Z_{G_{s}}(N'_{1})$ tel que $v_{1}^{-1}v$ soit l'image dans $\bar{A}(N)$ de $ \iota_{s}(z) $. En appliquant le th\'eor\`eme de Lang dans $G_{s}$, on \'ecrit $Fr(z)=x Fr(x)^{-1}$ pour un $x\in G_{s}$. Posons $N'=x^{-1}N'_{1}x$. Parce que $Fr(z)$ appartient \`a $Z_{G_{s}}(N'_{1})$, on a encore $N'\in \mathfrak{g}_{s,nil}({\mathbb F}_{q})$. On note ${\cal O}'$ l'orbite de $N'$ et  $y=\iota_{s}(x)$. On a alors $\bar{c}_{s,F}({\cal O}')=(y^{-1}Ny,d',v')$, où  $d'$ est l'image dans $\tilde{\bar{A}}(y^{-1}Ny)$ de $\dot{d}$ et $v'$ est celle de $Fr$ dans $\bar{A}^{++}_{Fr}(y^{-1}Ny)$. Montrons que
  
  (2) ce triplet est l'image par $Ad(y^{-1})$ de $(N,d,v)$. 
  
  En effet, $y$ commute à $\dot{d}$ car $y$ appartient à l'image de $\iota_{s}$. L'élément $Ad(y)(v')$ est l'image dans $\bar{A}^{++}_{Fr}(N)$ de $yFr\, y^{-1}=Fr \,Fr^{-1}(y)y^{-1}$. Puisque $\iota_{s}$ est $Fr$-équivariant, $Fr^{-1}(y)y^{-1}=\iota_{s}(Fr^{-1}(x)x^{-1})=\iota_{s}(z)$. C'est-à-dire que $Ad(y)(v')$ est l'image dans $\bar{A}^{++}_{Fr}(N)$ du produit de $v_{1}$ et de l'image dans $\bar{A}(N)$ de $\iota_{s}(z)$, autrement dit $Ad(y)(v')=v$. Cela prouve (2). 
  
  D'après (2), on a  l'égalité $\bar{c}_{s,F}({\cal O}')=(N,d,v)$ dans $\bar{{\cal C}}^{\sharp}_{F}$, c'est-à-dire modulo conjugaison des triplets. Donc $(N,d,v)$ appartient à l'image de $\bar{c}_{s,F}$, ce qui démontre l'assertion de  la proposition \ref{iotanil} sous l'hypothèse (1). 
  
     On va montrer que, quitte \`a changer de couple $(s^{nr}, {\cal O}^{nr})$, on peut supposer que la condition (1) est v\'erifi\'ee. 
   
   Supposons que $G$ est d\'eploy\'e sur $F$. Alors toutes les actions de $\Gamma_{F}^{nr}$  sont triviales et (1) est \'evident.

  Supposons  que $G$ n'est pas d\'eploy\'e sur $F^{nr}$. Par construction, $\Gamma_{F}^{nr}$ agit sur ${\cal D}_{a}^{nr}$ en conservant la racine $\beta_{0}$. Or l'identit\'e est le seul automorphisme de ${\cal D}_{a}^{nr}$ fixant $\beta_{0}$ donc cette action est triviale.  En particulier $s^{nr}$ est fix\'e par l'action de $\Gamma_{F}^{nr}$. La trivialit\'e de l'action de $\Gamma_{F}^{nr}$ implique aussi que $G_{s^{nr}}$ est d\'eploy\'e. Toute orbite nilpotente est alors conserv\'ee par l'action galoisienne. Donc (1) est v\'erifi\'ee.

  Supposons que $G$ est d\'eploy\'e sur $F^{nr}$ mais pas sur $F$. Alors $G$ est de type $A_{n-1}$, $D_{n}$ avec une action d'ordre $2$ de $\Gamma_{F}^{nr}$ sur ${\cal D}$,  $D_{4}$ trialitaire ou $E_{6}$.  Le sous-groupe $\boldsymbol{\Omega}$ de $G(F^{nr})$ introduit en \ref{epinglages} agit sur l'ensemble des couples $(s^{nr},{\cal O}^{nr})$ sans changer l'image $\bar{r}_{s^{nr}}({\cal O}^{nr})$  d'apr\`es \ref{auxiliaire}(2). 
    On v\'erifie cas par cas que, pour tout sommet $s'\in  S(\bar{C}^{nr})$, il existe $\omega\in \boldsymbol{\Omega}$ tel que $\omega(s')$ soit fix\'e par $\Gamma_{F}^{nr}$. On peut donc supposer que 
    $s^{nr}$ est fixe par l'action de $\Gamma_{F}^{nr}$. On pose simplement $s=s^{nr}$. Supposons v\'erifi\'ees les conditions
  
  (3)  l'action de $\Gamma_{F}^{nr}$ sur le diagramme de Dynkin de $G_{s}$ conserve chaque composante irr\'eductible et ces composantes ne sont pas de type $D_{n}$. 
  
  D'apr\`es \ref{actiongaloisienne}(1), ${\cal O}^{nr}$ est alors conserv\'ee par l'action galoisienne  donc (1) est v\'erifi\'ee.
  
  Il reste \`a consid\'erer les cas o\`u (3) n'est pas v\'erifi\'ee. Ce sont les cas (a) \`a (d) suivants.
  
  (a) On suppose que $G$ est de type $D_{n}$, le Frobenius agissant sur $\Delta_{a}$ par l'automorphisme $\theta$ qui permute $\alpha_{n-1}$ et $\alpha_{n}$ et fixe les autres racines. On param\`etre l'orbite de $N$ par une partition $\lambda\in {\cal P}^{orth}(2n)$. 
  Puisque $s$ est un \'el\'ement de $S(\bar{C})$ qui reste un sommet  dans $Imm_{F^{nr}}(G_{AD})$,  il est associ\'e \`a une racine $\alpha_{n''}$, avec $n''\in \{0,..., n-2\}$. Puisqu'il existe un \'el\'ement de $\boldsymbol{\Omega}$ qui permute $\alpha_{0}$ et $\alpha_{1}$, on peut supposer $n''\not=1$. Alors 
   $G_{s,SC}$ est de type $D_{n'}\times D_{n''}$, o\`u $n'=n-n''\geq2$.   Le groupe $\Gamma_{F}^{nr}$ agit trivialement sur le diagramme de $D_{n''}$ et non trivialement sur celui de $D_{n'}$. L'orbite ${\cal O}^{nr}$ est param\'etr\'ee par un couple $(\lambda',\lambda'')\in {\cal P}^{orth}(2n')\times {\cal P}^{orth}(2n'')$ tel que $\lambda'\cup\lambda''=\lambda$.  Si $\lambda'$ poss\`ede au moins un terme impair, ${\cal O}^{nr}$ est
  conserv\'ee par l'action galoisienne  et (1) est v\'erifi\'ee. Supposons que $\lambda'$ soit form\'ee de termes pairs. Par hypoth\`ese,  l'orbite de $N$  est conservée par $\Gamma_{F}$ donc $\lambda$ poss\`ede au moins un terme impair et $\lambda''$  aussi. A fortiori, $n''\not=0$ donc $n''\geq2$. Mais il y a un \'el\'ement $\delta\in \boldsymbol{\Omega}$ qui agit par $\alpha_{i}\mapsto \alpha_{n-i}$ sur le sous-ensemble $\{\alpha_{i};i=2,...,n-2\}$ de $\Delta_{a}$.   On  remplace $s$ par son image $\delta(s)$, qui est associ\'ee \`a $\alpha_{n'}$. Cela permute $(\lambda',\lambda'')$. Apr\`es cette permutation, $\lambda'$ poss\`ede au moins un terme impair et (1) est v\'erifi\'ee.

(b) On suppose que $G$ est de type $D_{4}$ trialitaire et que $s=s_{0}$. Dans ce cas, $G_{s_{0}}$ est de type $D_{4}$. Le groupe $G_{s_{0}}$ s'identifie \`a ${\bf G}$ et $\iota_{s_{0},nil}$ est l'identit\'e.  Puisque $\iota_{s_{0},nil}({\cal O}^{nr})$ est l'orbite de $N$ qui est   par hypoth\`ese conserv\'ee par $\Gamma_{F}^{nr}$, il en est de m\^eme de ${\cal O}^{nr}$  et (1) est v\'erifi\'ee.

(c) On suppose que $G$ est de type $D_{4}$ trialitaire et que $s=s_{2}$ (c'est-à-dire le sommet associ\'e \`a $\alpha_{2}$). Dans ce cas, $G_{s_{2}}$ est de type $A_{1}^4$. On note plus pr\'ecis\'ement ses composantes $A_{1,j}$, o\`u $j=0,1,3,4$, chacune d'elles correspondant \`a la racine $\alpha_{j}\in \Delta_{a}$. Le groupe $\Gamma_{{\mathbb F}_{q}}$  fixe la composante $A_{1,0}$ et permute  cycliquement les trois autres. Dans chaque composante $A_{1}$, il n'y a que deux orbites possibles sur $\bar{{\mathbb F}}_{q}$, l'orbite r\'eguli\`ere et l'orbite nulle.  L'orbite ${\cal O}^{nr}$ s'identifie \`a un quadruplet de telles orbites. Supposons que le quadruplet contienne trois termes \'egaux. L'action du groupe $\boldsymbol{\Omega}$ fixe $s_{2}$ et permute transitivement les composantes $A_{1}$.   On peut remplacer ${\cal O}^{nr}$ par son image par un \'el\'ement convenable et supposer que les composantes de ${\cal O}^{nr}$ dans $A_{1,1}$, $A_{1,3}$ et $A_{1,4}$ sont \'egales. Alors ${\cal O}^{nr}$ est conserv\'ee par l'action galoisienne et  (1) est v\'erifi\'ee. Il reste le cas o\`u le quadruplet d\'etermin\'e par ${\cal O}^{nr}$ a deux \'el\'ements r\'eguliers et deux \'el\'ements nuls. On raisonne alors comme dans la preuve de \ref{iotanilram} en fixant le produit scalaire sur $X_{*}(T)$ pour lequel $\vert \vert \check{\alpha}\vert \vert ^2=2$ pour toute racine $\alpha$. On calcule $\vert \vert x_{*,{\cal O}^{nr}}\vert \vert^2 =4$. Cette norme se conserve par l'application $\iota_{s_{2},nil}$. En notant $x_{*,N}$ le terme associ\'e \`a l'orbite de $N$, on a donc $\vert \vert x_{*,N}\vert \vert ^2=4$. Mais, d'apr\`es notre tableau de \ref{iotanilram}, c'est impossible puisque cette orbite est conservée par $\Gamma_{F}$. 

 (d) On suppose que $G$ est de type $E_{6}$ et que $s=s_{4}$. Dans ce cas, $G_{s_{4}}$ est de type $A_{2}^3$. Comme ci-dessus, on note plus pr\'ecis\'ement ses composantes $A_{2,0,2}$, $A_{2,1,3}$, $A_{2,6,5}$. Le groupe $\Gamma_{{\mathbb F}_{q}}$ fixe la composante  $A_{2,0,2}$ et permute les deux autres. Dans chaque composante $A_{2}$, il n'y a que trois orbites possibles sur $\bar{{\mathbb F}}_{q}$,  param\'etr\'ees par les partitions $(3)$, $(21)$ et $(111)$.  L'orbite ${\cal O}^{nr}$ s'identifie \`a un triplet de telles orbites. Supposons que le triplet contienne deux termes \'egaux. L'action du groupe $\boldsymbol{\Omega}$ fixe $s_{4}$ et permute transitivement les composantes $A_{2}$.   On peut remplacer ${\cal O}^{nr}$ par son image par un \'el\'ement convenable et supposer que les composantes de ${\cal O}^{nr}$ dans $A_{2,1,3}$ et $A_{2,6,5}$   sont \'egales. Alors ${\cal O}^{nr}$ est conserv\'ee par l'action galoisienne et (1) est v\'erifi\'ee. Reste le cas o\`u les trois termes du triplet sont distincts. Comme dans le cas (c), on voit que cela entra\^{\i}ne $\vert \vert x_{*,N}\vert \vert ^2=10$. D'apr\`es notre tableau de \ref{iotanilram}, l'orbite de $N$ est de type $A_{2}+A_{1}$. Pour cette orbite, $A(N)=\{1\}$ d'apr\`es \cite{C} p. 402 donc $d=1$ et $V=\{1\}$. On peut alors remplacer notre sommet $s_{4}$ par $s_{0}$.   En effet, on a
   $G_{s_{0}}={\bf G}$ et $\iota_{s_{0},nil}$ est l'identit\'e. Notons ${\cal O}_{0}^{nr}$ l'\'el\'ement de $\mathfrak{g}_{s_{0},nil}/conj$  qui s'identifie \`a l'orbite de $N$.    Puisque $A(N)=\{1\}$, on a forc\'ement $r_{s_{0}}({\cal O}_{0}^{nr})=(N,d,V)$.  Pour le couple $(s_{0},{\cal O}^{nr}_{0})$, l'assertion (3) est v\'erifi\'ee puisque 
   $G_{s_{0}}$ est de type $E_{6}$.
   Cela ach\`eve la preuve de la proposition. $\square$

   \subsection{Orbites nilpotentes de  $\mathfrak{g}(F)$}

\subsubsection{La construction de DeBacker\label{debacker}}
Dans ce paragraphe, on  suppose
  que $G$ est un groupe r\'eductif connexe d\'efini sur $F$.  On impose comme toujours l'hypothèse $(Hyp)_{1}(p)$. On a la propri\'et\'e suivante:
  
(1) soient $x\in Imm(G_{AD})$, $r\in {\mathbb R}$ avec $r\geq0$ et $e\in \mathfrak{g}_{nil}(F)\cap \mathfrak{k}_{x,r}$; alors $exp(e)\in K_{x,r}$. 

Preuve. 
Puisque $e$ est nilpotent, les \'el\'ements $e$ et $exp(e)$ sont les images naturelles d'\'el\'ements similaires appartenant respectivement \`a $\mathfrak{g}_{SC}(F)$ et $G_{SC}(F)$. L'assertion (1) pour ces \'el\'ements entra\^{\i}ne (1) pour nos \'el\'ements initiaux. Autrement dit,  on peut remplacer $G$ par $G_{SC}$ et supposer $G$ simplement connexe.    On peut aussi remplacer le corps de base $F$ par $F^{nr}$. Le groupe $G$ est alors quasi-d\'eploy\'e et on peut utiliser les constructions de \ref{discretisation}. On peut supposer $x\in App(T^{nr})$.  
La transformation $ad(e)$ de $\mathfrak{g}(F^{nr})$ conserve $\mathfrak{k}^{nr}_{x}$ donc $Ad(exp(e))$ aussi d'apr\`es la formule de \ref{leshypothesessurp} d\'ecrivant $exp(e)$. Un \'el\'ement de $G(F^{nr})$ qui conserve $\mathfrak{k}^{nr}_{x}$ conserve aussi la facette contenant $x$ (celle-ci est l'unique facette ${\cal F}$ telle que $\mathfrak{k}_{{\cal F}}^{nr}=\mathfrak{k}_{x}^{nr}$). Puisque $G$ est simplement connexe, le sous-groupe des \'el\'ements de $G(F^{nr})$ qui conservent cette facette est $K_{x}^{0,nr}$. Donc $exp(e)\in K_{x}^{0,nr}$. Si $r=0$, on a $K_{x}^{0,nr}=K^{nr}_{x,0}$ et on a termin\'e. Supposons $r>0$. On peut fixer $\tilde{y}\in {\mathbb T}^{nr}$ tel qu'en posant $y=a_{T}(\tilde{y})$, le point $y$ appartienne \`a la m\^eme facette que $x$.  On a alors  $\mathfrak{k}_{x,0+}^{nr}=\mathfrak{k}_{y,0+}^{nr}$, donc $e\in \mathfrak{k}_{y,0+}^{nr}$. On fixe une extension finie $F'$ de $F^{nr}$ mod\'er\'ement ramifi\'ee telle que $G$ soit d\'eploy\'e sur $F'$ et que, en notant $n=[F':F^{nr}]$, on ait $n\tilde{y}\in X_{*}(T^{nr})$.  Le raisonnement que l'on vient de faire montre que $exp(\lambda e)\in K_{y,F'}^0$ pour tout $\lambda\in \mathfrak{o}_{F'}$.  Donc $exp(\lambda e)$ se r\'eduit en un \'el\'ement de $G_{y,F'}$ que l'on note $u(\lambda)$. Parce que $e\in \mathfrak{k}_{y,0+}$, $Ad(u(\lambda))$ agit trivialement sur $\mathfrak{g}_{y,F'}$. Donc $u(\lambda)$ appartient au centre de $G_{y,F'}$. Parce que $n\tilde{y}\in X_{*}(T^{nr})$, $y$ est un point sp\'ecial dans $Imm_{F'}(G_{AD})$ et, puisque $G$ est semi-simple, $G_{y,F'}$ l'est aussi. Le centre de $G_{y,F'}$ est fini. L'application $\lambda\mapsto u(\lambda)$ se r\'eduit en une application alg\'ebrique d\'efinie sur $\bar{{\mathbb F}}_{q}$ qui ne prend qu'un nombre fini de valeurs. C'est donc forc\'ement l'application constante de valeur $u(0)=1$. En particulier $exp(e)$ a pour r\'eduction $1$ donc appartient \`a $K_{y,0+,F'}$.  Puisque c'est un \'el\'ement de $G(F^{nr})$, on a en fait $exp(e)\in K_{y,0+}^{nr}=K_{x,0+}^{nr}$. Puisque l'application $s\mapsto K_{x,s}^{nr}$ est continue \`a gauche, on peut d\'efinir le plus grand r\'eel $s\in [0,r]$ tel que $exp(e)\in K_{x,s}^{nr}$.  La relation $exp(e)\in K_{x,0+}^{nr}$ entra\^{\i}ne que $s>0$.  Supposons $s<r$. D'apr\`es \cite{A}, proposition 1.6.3, il y a un isomorphisme $\phi:K_{x,s}^{nr}/K_{x,s+}^{nr}\simeq \mathfrak{k}_{x,s}^{nr}/\mathfrak{k}_{x,s+}^{nr}$ v\'erifiant les conditions suivantes. Soient $g\in K_{x,s}^{nr}$ et $X\in \mathfrak{k}_{x,s}^{nr}$. Notons $\bar{g}$, resp. $\bar{X}$, l'image de   $g$ dans $K_{x,s}^{nr}/K_{x,s+}^{nr}$, resp. $X$ dans $\mathfrak{k}_{x,s}^{nr}/\mathfrak{k}_{x,s+}^{nr}$. Supposons $\phi(\bar{g})=\bar{X}$. Soient $z\in {\mathbb R}$ et $Z\in \mathfrak{k}^{nr}_{x,z}$. Alors $Ad(g)(Z)-Z-[X,Z]$ appartient \`a $\mathfrak{k}^{nr}_{x,(z+s)+}$. On applique cela \`a $g=exp(e)$. Par d\'efinition de $s$ et d'apr\`es l'hypoth\`ese $s<r$, on a $\bar{g}\not=1$ d'o\`u $\bar{X}\not=0$. D'apr\`es \cite{DR} lemme B.6.1, on peut trouver $z\in {\mathbb R}$ et $Z\in \mathfrak{k}^{nr}_{x,z}$ de sorte que $[X,Z]\not\in \mathfrak{k}^{nr}_{x,(z+s)+}$ (les hypoth\`eses de l'article de DeBacker et Reeder sont plus fortes que les n\^otres mais ces derni\`eres suffisent pour le lemme cit\'e). On a alors $Ad(exp(e))(Z)-Z\not\in \mathfrak{k}^{nr}_{x,(z+s)+}$. Mais la d\'efinition de $exp(e)$ et l'hypoth\`ese $e\in \mathfrak{k}_{x,r}$ entra\^{\i}nent que $Ad(exp(e))(Z)-Z\in \mathfrak{k}^{nr}_{x,z+r}\subset \mathfrak{k}^{nr}_{x,(z+s)+}$.  Cette contradiction \`a laquelle  conduit notre hypoth\`ese $s<r$ prouve que $s=r$, c'est-\`a-dire que $exp(e)\in K_{x,r}^{nr}$. Cela prouve (1).

Soit ${\cal F}$ une facette de l'immeuble $Imm(G_{AD})$. On note $x\mapsto \bar{x}$ les applications de r\'eduction de $\mathfrak{k}^{nr}_{{\cal F}}$ dans $\mathfrak{g}_{{\cal F}} $ ou de $K^{0,nr}_{{\cal F}}$ dans $G_{{\cal F}}$. On a d\'efini les ensembles $\mathfrak{g}_{{\cal F},nil}({\mathbb F}_{q})/conj$ et $ \mathfrak{g}_{nil}(F)/conj$. Pour une orbite ${\mathbb  O}\in \mathfrak{g}_{nil}(F)/conj$, on note $Cl_{F}({\mathbb O})$ son adh\'erence pour la topologie $p$-adique. 
 Dans \cite{D2}, DeBacker  d\'efinit une application 
$$rel_{{\cal F},F}:\mathfrak{g}_{{\cal F},nil}({\mathbb F}_{q})/conj\to \mathfrak{g}_{nil}(F)/conj.$$
Rappelons la construction. Soit ${\cal O}\in \mathfrak{g}_{{\cal F},nil}({\mathbb F}_{q})/conj$. Choisissons un $\mathfrak{sl}(2)$-triplet $(\underline{f},\underline{ h},\underline{ e})$ de $\mathfrak{g}_{{\cal F}}({\mathbb F}_{q})$ tel que $\underline{e}\in {\cal O}$. DeBacker d\'emontre qu'il existe un $\mathfrak{sl}(2)$-triplet $(f,h,e)$ de $\mathfrak{g}(F)$ v\'erifiant les conditions suivantes:

$f,h,e\in \mathfrak{k}_{{\cal F}}$;

$\bar{f}=\underline{ f}$, $\bar{h}=\underline{ h}$, $\bar{e}=\underline{ e}$.

On dira qu'un tel triplet rel\`eve $(\underline{f},\underline{ h},\underline{e})$ dans $\mathfrak{k}_{{\cal F}}$. 
DeBacker d\'emontre que l'orbite nilpotente ${\mathbb O}\in \mathfrak{g}_{nil}(F)/conj$ contenant $e$ ne d\'epend pas du choix des $\mathfrak{sl}(2)$-triplets. Alors $rel_{{\cal F},F}({\cal O})={\mathbb O}$. 

{\bf Remarque.} Inversement, si $(f,h,e)$ est un $\mathfrak{sl}(2)$-triplet de $\mathfrak{g}(F)$ tel que $f,h,e\in \mathfrak{k}_{{\cal F}}$, alors $(\bar{f},\bar{h},\bar{e})$ est un $\mathfrak{sl}(2)$-triplet de $\mathfrak{g}_{{\cal F}}({\mathbb F}_{q})$.
\bigskip

Pr\'ecisons quelques points de la construction qui sont implicites dans l'article de DeBacker mais  que l'on pr\'ef\`ere \'enoncer. 

\begin{lem}{(i) Soient $(f,h,e)$ et $(f',h',e')$ deux $\mathfrak{sl}(2)$-triplets de $\mathfrak{g}(F)$. Supposons que les \'el\'ements $f,h,e,f',h',e'$ appartiennent \`a $ \mathfrak{k}_{{\cal F}}$ et supposons que les r\'eductions $(\bar{f},\bar{h},\bar{e})$ et $(\bar{f}',\bar{h}',\bar{e}')$ soient deux $\mathfrak{sl}(2)$-triplets de $\mathfrak{g}_{{\cal F}}({\mathbb F}_{q})$ qui sont conjugu\'es par un \'el\'ement de $G_{{\cal F}}({\mathbb F}_{q})$, resp. qui sont \'egaux. Alors $(f,h,e)$ et $(f',h',e')$ sont conjugu\'es par un \'el\'ement de $K_{{\cal F}}^0$, resp. de $K_{{\cal F}}^+$.

(ii) Soit $(\underline{f},\underline{h},\underline{e})$ un $\mathfrak{sl}(2)$-triplet de $\mathfrak{g}_{{\cal F}}({\mathbb F}_{q})$, soit $S\in {\cal T}_{max}$ un tore tel que ${\cal F}$ appartienne \`a l'appartement de $Imm(G_{AD})$ associ\'e \`a $S$ et que $\underline{h}$ appartienne \`a $\mathfrak{s}_{{\cal F}}({\mathbb F}_{q})$. Alors il existe un $\mathfrak{sl}(2)$-triplet $(f,h,e)$ de $\mathfrak{g}(F)$ qui  rel\`eve $(\underline{f},\underline{ h},\underline{e})$ dans $\mathfrak{k}_{{\cal F}}$ et tel que $h\in \mathfrak{s}(\mathfrak{o}_{F})$. On a alors $x_{*,h}\in X_{*}(S)$, $x_{*,\underline{h}}\in X_{*}(S_{{\cal F}})$ et ces \'el\'ements s'identifient par l'isomorphisme $X_{*}(S)\simeq X_{*}(S_{{\cal F}})$.
}\end{lem}

Preuve de (i). L'assertion concernant le cas o\`u les r\'eductions des deux $\mathfrak{sl}(2)$-triplets sont \'egales implique celle o\`u elles sont seulement conjugu\'ees. En effet si ces r\'eductions sont conjugu\'ees, on choisit un \'el\'ement $k\in K_{{\cal F}}^0$ se r\'eduisant en un \'el\'ement $\bar{k}$ conjuguant $(\bar{f}',\bar{h}',\bar{e}')$ en $(\bar{f},\bar{h},\bar{e})$. On remplace $(f',h',e')$ par son image par la conjugaison par $k$. On obtient un $\mathfrak{sl}(2)$-triplet dans $\mathfrak{k}_{{\cal F}}$ dont la r\'eduction co\"{\i}ncide avec $(\bar{f},\bar{h},\bar{e})$ et il reste \`a appliquer \`a ces triplets l'assertion concernant le cas o\`u les r\'eductions des deux $\mathfrak{sl}(2)$-triplets sont \'egales.

On suppose donc que $(\bar{f},\bar{h},\bar{e})=(\bar{f}',\bar{h}',\bar{e}')$. D'apr\`es le r\'esultat de DeBacker mentionn\'e ci-dessus, $e$ et $e'$ appartiennent \`a la m\^eme orbite dans $\mathfrak{g}_{nil}(F)$. 
  Puisque $\bar{e}'=\bar{e}$, on a $e'\in e+\mathfrak{k}_{{\cal F}}^+$. Le lemme 5.2.1 de \cite{D2} dit qu'il existe $k\in K_{{\cal F}}^+$ et $Y\in \mathfrak{z}_{G}(f)\cap \mathfrak{k}_{{\cal F}}^+$ de sorte que $ad(k)(e')=Y+e$.  Parce que $e$ et $e'$ appartiennent \`a la m\^eme orbite,   on a $Y=0$ d'apr\`es \cite{W1} V.7(9). On conjugue le triplet $(f',h',e')$ par $k$ et on est ramen\'e au cas o\`u $e=e'$, ce que l'on suppose d\'esormais. Fixons un point $x\in {\cal F}$. 
  On sait que l'ensemble des r\'eels $r$ tels que $\mathfrak{k}_{x,r}\not=\mathfrak{k}_{x,r+}$ est un sous-ensemble de ${\mathbb Q}$ invariant par translations par ${\mathbb Z}$ et qu'il est d'image finie dans ${\mathbb Q}/{\mathbb Z}$. Notons ses \'el\'ements positifs ou nuls $r_{0}=0<r_{1}<r_{2}<...$. On note $m$ l'entier tel que $\mathfrak{k}_{x,r_{i+m}}=\mathfrak{p}_{F}\mathfrak{k}_{x,r_{i}}$ pour tout $i$.  On a $f'\in f+\mathfrak{k}_{{\cal F}}^+=f+\mathfrak{k}_{x,r_{1}}$, $h'\in h+\mathfrak{k}_{{\cal F}}^+=h+\mathfrak{k}_{x,r_{1}}$ et $K_{x,r}\subset K_{{\cal F}}^+$ pour tout $r>0$.  On voit qu'il suffit de d\'emontrer l'assertion suivante:

(2) soit $i\geq1$ un entier et soit $(f_{i},h_{i},e)$ un $\mathfrak{sl}(2)$-triplet de $\mathfrak{g}(F)$ tel que $f_{i}\in f+\mathfrak{k}_{x,r_{i}}$, $h_{i}\in h+\mathfrak{k}_{x,r_{i}}$; alors il existe $k\in K_{x,r_{i}}\cap Z_{G}(e)(F)$ tel qu'en posant $f_{i+1}=ad(k)(f_{i})$ et $h_{i+1}=ad(k)(h_{i})$, on ait $f_{i+1}\in f+\mathfrak{k}_{x,r_{i+1}}$, $h_{i+1}\in h+\mathfrak{k}_{x,r_{i+1}}$. 

Pour tout $j\in {\mathbb Z}$, posons $\mathfrak{g}_{x,r_{j}}=\mathfrak{k}_{x,r_{j}}/\mathfrak{k}_{x,r_{j+1}}$. L'action adjointe de l'alg\`ebre $\mathfrak{k}_{x}=\mathfrak{k}_{{\cal F}}$ conserve chaque $\mathfrak{k}_{x,r_{j}}$ et cette action se quotiente en une action de $\mathfrak{g}_{{\cal F}}({\mathbb F}_{q})$ dans $\mathfrak{g}_{x,r_{j}}$. En particulier, l'action du $\mathfrak{sl}(2)$-triplet $(f,e,h)$ dans $\mathfrak{k}_{x,r_{i}}$ se quotiente en une action de $(\bar{f},\bar{h},\bar{e})$ dans $\mathfrak{g}_{x,r_{i}}$. En identifiant $(f,h,e)$ au triplet standard dans $\mathfrak{sl}(2,F)$, on applique la proposition \ref{sl2} \`a l'action du triplet $(f,h,e)$ dans $\mathfrak{g}(F)$, cet espace \'etant muni de la cha\^{\i}ne de r\'eseaux $\mathfrak{k}_{x,r_{i}}\supset \mathfrak{k}_{x,r_{i+1}}\supset...\supset \mathfrak{k}_{x,r_{i+m}}=\mathfrak{p}_{F}\mathfrak{k}_{x,r_{i}}$ (l'hypothèse de cette proposition est vérifiée d'après $(Hyp)_{1}(p)$).  On peut fixer un ensemble fini $J$ muni d'une application $d:J\to \{1,...,p-1\}$ et un morphisme injectif de $\mathfrak{sl}(2,F)$-modules $\phi:\oplus_{j\in J}{\cal W}_{d(j)}\to \mathfrak{g}(F)$ de sorte que $\phi(\oplus_{j\in J}{\cal L}_{d(j)})\subset \mathfrak{k}_{x,r_{i}}$ et que cette derni\`ere application se r\'eduise en un isomorphisme 

  $$(3) \qquad\oplus_{j\in J}\bar{{\cal L}}_{d(j)} \simeq \mathfrak{g}_{x,r_{i}},$$
 o\`u on a not\'e $\bar{{\cal L}}_{d(j)}={\cal L}_{d(j)}/\mathfrak{p}_{F}{\cal L}_{d(j)}$. On fait dispara\^{\i}tre $\phi$ de la notation en identifiant son ensemble de d\'epart \`a son image. Posons $y=f_{i}-f$ et $z=h_{i}-h$. Ce sont des \'el\'ements de $\mathfrak{k}_{x,r_{i}}$ par hypoth\`ese et on note $\bar{y}$ et $\bar{z}$ leurs r\'eductions dans $\mathfrak{g}_{x,r_{i}}$. Les relations $[e,f]=h$, $[e,f_{i}]=h_{i}$, $[h,e]=2e$ et $[h_{i},e]=2e$ impliquent $[e,y]=z$ et $[z,e]=0$, d'o\`u $[\bar{e},\bar{y}]=\bar{z}$ et $[\bar{z},\bar{e}]=0$. Ecrivons $\bar{y}=\sum_{j\in J}\bar{y}_{j}$ et $\bar{z}=\sum_{j\in J}\bar{z}_{j}$ selon l'isomorphisme (3). On a encore $[\bar{e},\bar{y}_{j}]=\bar{z}_{j}$ et $[\bar{z}_{j},\bar{e}]=0$ pour tout $j$. Ces relations impliquent $[\bar{h},\bar{z}_{j}]=(d(j)-1)\bar{z}_{j}$ et   $\bar{z}_{j}=0$ si $d(j)=1$.  Si $d(j)=1$, posons $z'_{j}=0$. Si $d(j)\geq2$, il est loisible de fixer un \'el\'ement $z'_{j}\in{\cal L}_{d(j)}$ qui v\'erifie $[h,z'_{j}]=(d(j)-1)z'_{j}$, $[e,z'_{j}]=0$ et $\bar{z}'_{j}=\frac{1}{d(j)-1}\bar{z}_{j}$. Posons $z'=\sum_{j\in J}z'_{j}$. Les relations impos\'ees impliquent que $z'$ est un \'el\'ement du radical nilpotent de $\mathfrak{z}_{G}(e)(F)$.  Posons $k=exp(z')$.  On a $k\in Z_{G}(e)(F)$. Puisque $z'\in \mathfrak{k}_{x,r_{i}}$, on a $k\in K_{x,r_{i}}$ d'apr\`es (1). Par construction de $k$, on a $Ad(k)(h_{i})\equiv h_{i}+[z',h_{i}]\, mod\, \mathfrak{k}_{x, r_{i+1}}$. On a aussi $h_{i}+[z',h_{i}]\equiv h_{i}+[z',h]=h+z+[z',h]\,mod\, \mathfrak{k}_{x, r_{i+1}}$. Mais, d'apr\`es les d\'efinitions, on a $\bar{z}+[\bar{z}',\bar{h}]=0$. Donc $Ad(k)(h_{i})\in h+\mathfrak{k}_{x,r_{i+1}}$. Notons $(f_{i+1},h_{i+1},e)$ l'image de $(f_{i},h_{i},e)$ par la conjugaison par $k$. Alors $h_{i+1}\in h+\mathfrak{k}_{x,r_{i+1}}$. Posons $y''=f_{i+1}-f$, $z''=h_{i+1}-h$, notons encore $\bar{y}''$ et $\bar{z}''$ leurs r\'eductions dans $\mathfrak{g}_{x,r_{i}}$ et introduisons leurs composantes $\bar{y}''_{j}$ et $\bar{z}''_{j}$ comme ci-dessus. On a $\bar{z}''=0$. 
  Les \'egalit\'es $[h_{i+1},f_{i+1}]=-2f_{i+1}$ et $[h,f]=-2f$ entra\^{\i}nent $[\bar{h},\bar{y}'']=-2\bar{y}''$ d'o\`u $[\bar{h},\bar{y}''_{j}]=-2\bar{y}''_{j}$ pour tout $j\in J$. On a aussi comme ci-dessus  $[\bar{e},\bar{y}''_{j}]=\bar{z}''_{j}=0$. Ces deux \'egalit\'es entra\^{\i}nent $\bar{y}''_{j}=0$ pour tout $j\in J$, d'o\`u $\bar{y}''=0$. Alors $f_{i+1}\in f+\mathfrak{k}_{x,r_{i+1}}$. Cela d\'emontre (2) et l'assertion (i) de l'\'enonc\'e.

  Preuve de (ii). Fixons un $\mathfrak{sl}(2)$-triplet $(f',h',e')$ relevant $(\underline{f},\underline{h},\underline{e})$ dans $\mathfrak{k}_{{\cal F}}$. Notons $M$ le commutant de $h'$ dans $G$. D'apr\`es \ref{sl2triplets} (1), c'est un $F$-Levi de $G$ qui est \'egal au commutant dans $G$ de l'image du groupe \`a un param\`etre $X_{*,h'}$. 
   D'apr\`es \cite{D2} corollaire 4.5.9,  la facette ${\cal F}$ est contenue dans $Imm^G(M_{ad})$. Notons ${\cal F}_{M}$ la facette de $Imm(M_{AD})$ qui contient $p_{M}({\cal F})$. En appliquant les constructions de  \ref{groupesenreduction} \`a un point $x\in {\cal F}$ quelconque, on d\'eduit de $M$ 
   un ${\mathbb F}_{q}$-Levi $M_{{\cal F}}$ de $G_{{\cal F}}$,  qui s'identifie \`a $M_{{\cal F}_{M}}$. Montrons que
  
  (4) $M_{{\cal F}}$ est le commutant de $\underline{h}$ dans $G_{{\cal F}}$.
  
  Pour tout $i\in {\mathbb Z}$, notons $\mathfrak{g}_{i}$  le sous-espace des $X\in \mathfrak{g}$ tels que $[h',X]=iX$. Pour $i\in {\mathbb Z}/p{\mathbb Z}$, notons  $\mathfrak{g}_{{\cal F},i}$ le sous-espace des $X\in \mathfrak{g}_{{\cal F}}$ tels que $[\underline{h},X]=iX$. Notons $\pi:\mathfrak{k}_{{\cal F}}\to \mathfrak{g}_{{\cal F}}$ la projection naturelle. En appliquant la proposition \ref{sl2} \`a la suite de r\'eseaux $\mathfrak{k}_{{\cal F}}\supset \mathfrak{k}_{{\cal F}}^+\supset \mathfrak{p}_{F}\mathfrak{k}_{{\cal F}}$ de $\mathfrak{g}(F)$, on voit que $\mathfrak{g}_{{\cal F}}$ est somme directe des sous-espaces $\pi(\mathfrak{g}_{i}\cap \mathfrak{k}_{{\cal F}})$ quand $i$ d\'ecrit ${\mathbb Z}$. Il est clair que $\pi(\mathfrak{g}_{i}\cap \mathfrak{k}_{{\cal F}})\subset \mathfrak{g}_{{\cal F},\bar{i}}$ o\`u $\bar{i}$ est l'image de $i$ dans ${\mathbb Z}/p{\mathbb Z}$. D'apr\`es notre hypoth\`ese  $(Hyp)_{1}(p)$, l'ensemble des $i$ tels que $\mathfrak{g}_{i}\not=\{0\}$ est contenu dans un intervalle $[-a,a]$ avec $a<p$.  Il y a un seul \'el\'ement $i$ de cet intervalle tel que $\bar{i}=0$, \`a savoir $i=0$.   Il en r\'esulte que $\mathfrak{g}_{{\cal F},0}=\pi(\mathfrak{g}_{0}\cap \mathfrak{k}_{{\cal F}})$. L'espace $\mathfrak{g}_{0}$ n'est autre que $\mathfrak{m}$ donc $\pi(\mathfrak{g}_{0}\cap \mathfrak{k}_{{\cal F}})$ est \'egal \`a $\mathfrak{m}_{{\cal F}}$. L'espace $\mathfrak{g}_{{\cal F},0}$ est \'egal \`a l'alg\`ebre de Lie  de $Z_{G_{{\cal F}}}(\underline{h})$. On obtient donc que cette derni\`ere alg\`ebre de Lie est \'egale \`a $\mathfrak{m}_{{\cal F}}$. Puisque $M_{{\cal F}}$  et $Z_{G_{{\cal F}}}(\underline{h})$ sont d\'etermin\'es par leurs alg\`ebres de Lie, cela implique (4).
  
  Puisque $\underline{h}\in \mathfrak{s}_{{\cal F}}({\mathbb F}_{q})$,  (4) implique que $S_{{\cal F}}\subset M_{{\cal F}}$. On peut alors relever $S_{{\cal F}}$ en un tore $S'\in {\cal T}_{max}$ contenu dans $M$ et tel que ${\cal F}_{M}$ appartienne \`a l'appartement de $Imm(M_{AD})$ associ\'e \`a $S'$ (dire que $S'$ rel\`eve $S_{{\cal F}}$ signifie que $S_{{\cal F}}=S'_{{\cal F}}$). Ces propri\'et\'es restent valables en rempla\c{c}ant $M$ par $G$: $S'$ rel\`eve $S_{{\cal F}}$ et ${\cal F}$ appartient \`a l'appartement de $Imm(G_{AD})$ associ\'e \`a $S'$. Le tore $S$ v\'erifie les m\^emes propri\'et\'es. Il existe alors $k\in K_{{\cal F}}^+$ tel que $kS'k^{-1}=S$. Fixons un tel $k$ et notons $(f,h,e)$ l'image de $(f',h',e')$ par $Ad(k)$. Puisque $k$ appartient \`a $K_{{\cal F}}^+$, le triplet $(f,h,e)$ rel\`eve encore $(\underline{f},\underline{h},\underline{e})$ dans $\mathfrak{k}_{{\cal F}}$.   Reprenons les constructions ci-dessus en rempla\c{c}ant $(f',h',e')$ par $(f,h,e)$, c'est-\`a-dire en d\'efinissant $M$ comme le commutant de $h$. On a maintenant $S\subset M$. L'image du groupe \`a un param\`etre $X_{*,h}$ est un sous-tore d\'eploy\'e de $G$ qui est contenu dans le centre de $M$, donc commute \`a $S$. Puisque $S$ est un sous-tore d\'eploy\'e maximal, cette image est contenue dans $S$. Cela implique $x_{*,h}\in X_{*}(S)$ et $h\in \mathfrak{s}(\mathfrak{o}_{F})$. Introduisons l'homomorphisme $\phi:SL(2)\to G$ associ\'e \`a $(f,h,e)$, cf. \ref{sl2triplets} (en \ref{sl2triplets}, $\phi$ \'etait \`a valeurs dans $G_{SC}$, on prend ici son compos\'e avec l'homomorphisme naturel $G_{SC}\to G$). Introduisons l'homomorphisme similaire $\underline{\phi}:SL(2)\to G_{{\cal F}}$ associ\'e \`a $(\underline{f},\underline{h},\underline{e})$. Le groupe $SL(2,\mathfrak{o}_{F^{nr}})$ est engendr\'e par les \'el\'ements $exp(\lambda \pmb{e})$ et $exp(\lambda \pmb{f})$ pour $\lambda\in \mathfrak{o}_{F^{nr}}$. Les images par $\phi$ de ces \'el\'ements sont $exp(\lambda e)$ et $exp(\lambda f)$ qui appartiennent \`a $K_{{\cal F}}^{0,nr}$ d'apr\`es (1)  et se r\'eduisent en $exp(\bar{\lambda}\underline{e})=\underline{\phi}(\bar{\lambda}{\bf e})$ et  $exp(\bar{\lambda}\underline{f})=\underline{\phi}(\bar{\lambda}{\bf f})$ . On en d\'eduit que  $\phi(SL(2,\mathfrak{o}_{F^{nr}}))\subset K_{{\cal F}}^{0,nr}$ et que le diagramme suivant est commutatif
  $$\begin{array}{ccc}SL(2,\mathfrak{o}_{F^{nr}})&\stackrel{\phi}{\to}&K_{{\cal F}}^{0,nr}\\ \downarrow&&\downarrow\\ SL(2,\bar{{\mathbb F}}_{q})&\stackrel{\underline{\phi}}{\to}&G_{{\cal F}}\\ \end{array}$$
  Alors $x_{*,\underline{h}}$ est la r\'eduction de $x_{*,h}$, c'est-\`a-dire que $x_{*,\underline{h}}(\bar{\lambda})$ est la r\'eduction de $x_{*,h}(\lambda)$ pour tout $\lambda\in \mathfrak{o}_{F^{nr}}^{\times}$. Puisque $x_{*,\underline{h}}$ prend ses valeurs dans $S$, $x_{*,\underline{h}}$ prend ses valeurs dans $S_{{\cal F}}$. Autrement dit,  $x_{*,\underline{h}}\in X_{*}(S_{{\cal F}})$ et  cet \'el\'ement s'identifie \`a  $x_{*,h}$ par l'isomorphisme $X_{*}(S_{{\cal F}})\simeq X_{*}(S)$.   
 $\square$

   Il r\'esulte du lemme ci-dessus que
   
   (5) $rel_{{\cal F},F}$ conserve les termes $x_{*,{\cal O}}$.
   
    Pr\'ecisons ce que l'on entend par l\`a. Fixons un tore $T^{nr}\in {\cal T}^{nr}_{max}$ tel que ${\cal F}$ soit contenue dans l'appartement de $Imm_{F^{nr}}(G_{AD})$  associ\'e \`a $T^{nr}$ et  notons $T$ le commutant de $T^{nr}$ dans $G$.  Soit ${\cal O}\in \mathfrak{g}_{{\cal F},nil}/conj$, posons ${\mathbb O}=rel_{{\cal F},F}({\cal O})$. L'\'el\'ement $x_{*,{\cal O}}$ s'identifie \`a la classe de conjugaison par le groupe de Weyl  $W_{{\cal F}}$ de $G_{{\cal F}}$ d'un \'el\'ement de $X_{*}(T^{nr})$ et $x_{*,{\mathbb O}}$ s'identifie \`a la classe de conjugaison par $W$ d'un \'el\'ement de $X_{*}(T)$. Alors $x_{*,{\mathbb O}}$ est l'unique  classe de conjugaison par $W$ contenant la classe $x_{*,{\cal O}}$. 

 On a aussi, cf.  \cite{D2} corollaire 5.2.5: 

(6) soit ${\mathbb O}'\in \mathfrak{g}_{nil}/conj$; supposons qu'il existe $e'\in {\mathbb O}'\cap \mathfrak{k}_{{\cal F}}$ tel que $\bar{e}'\in {\cal O}$; alors $rel_{{\cal F},F}({\cal O})\subset Cl_{F}({\mathbb O}')$; 
  
   Soit ${\cal F}'$ une facette contenue dans  $\bar{{\cal F}}$. On a alors $\mathfrak{k}_{{\cal F}}\subset \mathfrak{k}_{{\cal F}'}$.  D'apr\`es \ref{adherence}, on peut identifier $\mathfrak{g}_{{\cal F}}$ \`a une sous-alg\`ebre de $\mathfrak{g}_{{\cal F}'}$. Il en r\'esulte    une application    $\mathfrak{g}_{{\cal F},nil}({\mathbb F}_{q})/conj\to \mathfrak{g}_{{\cal F}',nil}({\mathbb F}_{q})/conj$ qui est canonique. Alors
  
  (7) $rel_{{\cal F},F}$ est la compos\'ee  de $rel_{{\cal F}',F}$ et de cette application.
  
  Soit $g\in G(F)$, posons ${\cal F}'=g{\cal F}$. Comme en \ref{conjugaison}, on a un isomorphisme $\overline{Ad(g)}:\mathfrak{g}_{{\cal F},nil}({\mathbb F}_{q})/conj\to \mathfrak{g}_{{\cal F}',nil}({\mathbb F}_{q})/conj$. Il est clair que
  
  (8)  $rel_{{\cal F},F}$ est la compos\'ee de   $rel_{{\cal F}',F}$ et de $\overline{Ad(g)}$.

  Enfin, et surtout,  on a, cf. \cite{D2} th\'eor\`eme p. 297:    
  
  (9) $\mathfrak{g}_{nil}(F)/conj$ est r\'eunion des images des applications $rel_{{\cal F},F}$ quand ${\cal F}$ d\'ecrit les facettes de l'immeuble $Imm(G_{AD})$.
  
  D'apr\`es (7), on peut se limiter dans (9) \`a des facettes qui sont des sommets et, d'apr\`es (8),  on peut m\^eme se limiter \`a des sommets dans un ensemble de repr\'esentants des orbites pour l'action de $G(F)$ dans l'ensemble des sommets.

\subsubsection{Un diagramme d'applications\label{undiagramme}}
Pour la fin de la sous-section 3.3, {\bf on reprend les hypoth\`eses de \ref{dansboldsymbolg}:} $G$ {\bf est quasi-d\'eploy\'e, adjoint et absolument simple.} 
   
Soit ${\cal F}$ une facette de l'immeuble $Imm(G_{AD})$.  Soit $F'$ une extension finie non  ramifi\'ee de $F$, notons ${\mathbb F}_{q^f}$ son corps r\'esiduel. Il existe une unique facette ${\cal F}_{F'}$ de $Imm_{F'}(G_{AD})$ telle que ${\cal F}$ soit le sous-ensemble des points fixes dans ${\cal F}_{F'}$ par l'action de $\Gamma_{F'/F}$. On a $\mathfrak{g}_{{\cal F}}=\mathfrak{g}_{{\cal F}_{F'}}$. Il r\'esulte des d\'efinitions que le diagramme suivant est commutatif:
  $$\begin{array}{ccc}\mathfrak{g}_{{\cal F},nil}({\mathbb F}_{q})/conj&\stackrel{rel_{{\cal F},F}}{\to}& \mathfrak{g}_{nil}(F)/conj\\ \downarrow&&\downarrow\\ \mathfrak{g}_{{\cal F},nil}({\mathbb F}_{q^f})/conj&\stackrel{rel_{{\cal F}_{F'},F'}}{\to}& \mathfrak{g}_{nil}(F')/conj\ \end{array}$$
  Avec des notations \'evidentes, on en d\'eduit une suite d'applications
  $$\mathfrak{g}_{{\cal F},nil}/conj=\mathfrak{g}_{{\cal F},nil}(\bar{{\mathbb F}}_{q})/conj\to \mathfrak{g}_{nil}(F^{nr})/conj\to (\mathfrak{g}_{nil}(\bar{F})/conj)^{I_{F}}=(\mathfrak{g}_{nil}/conj)^{I_{F}},$$
  dont on note $Rel_{{\cal F}}$ la compos\'ee. On a suppos\'e que ${\cal F}$ \'etait une facette de $Imm(G_{AD})$ mais on peut remplacer le corps de base $F$ par une extension finie non ramifi\'ee. L'application pr\'ec\'edente est donc d\'efinie pour toute facette ${\cal F}$ de $ Imm_{F^{nr}}(G_{AD})$. 
  
   Le r\'esultat de DeBacker \ref{debacker}(9) nous dit que $\mathfrak{g}_{nil}(F^{nr})/conj$ est r\'eunion des images des applications $rel_{{\cal F},F^{nr}}$ quand ${\cal F}$ d\'ecrit l'ensemble des facettes de $ Imm_{F^{nr}}(G)$. Le groupe $G$ est suppos\'e quasi-d\'eploy\'e.  Pour toute extension finie $F'$ de $F$, toute orbite dans $\mathfrak{g}_{nil}$ qui est conserv\'ee par $\Gamma_{F'}$ contient un point de $\mathfrak{g}(F')$. Il en r\'esulte que 
 l'application $\mathfrak{g}_{nil}(F^{nr})/conj\to (\mathfrak{g}_{nil}/conj)^{I_{F}}$ est surjective. Donc $(\mathfrak{g}_{nil}/conj)^{I_{F}}$ est r\'eunion des images des applications $Rel_{{\cal F}}$ quand ${\cal F}$ d\'ecrit l'ensemble des facettes de  $Imm_{F^{nr}}(G)$. 
  
Reprenons  les constructions de \ref{discretisation}. Rappelons que l'on  note $F^G$ la plus petite extension de $F^{nr}$ telle que $G$ soit d\'eploy\'e sur $F^G$. Fixons une extension finie mod\'er\'ement ramifi\'ee $F'$ de $F^G$. On a introduit le point $s_{0}\in Imm(G)$. Le groupe ${\bf G}$  s'identifie \`a $G_{s_{0},F'}$. On peut fixer une extension finie $E$ de $F$ de sorte que $F'=E^{nr}$. En se pla\c{c}ant sur le corps de base $E$, on a d\'efini l'application $Rel_{s_{0}}:\mathfrak{g}_{s_{0},F',nil}/conj\to \mathfrak{g}_{nil}/conj$. Elle s'identifie \`a une application que nous notons
 $${\bf Rel}:\boldsymbol{\mathfrak{g}}_{nil}/conj\to \mathfrak{g}_{nil}/conj.$$
 Les algèbres $\boldsymbol{\mathfrak{g}}$ et $\mathfrak{g}$ sont des algèbres de Lie simples de m\^eme type définies, la première sur $\bar{{\mathbb F}}_{q}$, la seconde sur $\bar{F}$. On sait que la classification des  orbites nilpotentes d'une alg\`ebre de Lie ne d\'epend pas du corps de base pourvu que celui-ci soit alg\'ebriquement clos et de caract\'eristique nulle ou  de caract\'eristique $p$  assez grande, ce qui est assuré par notre hypothèse $(Hyp)_{1}(p)$. Autrement dit, le paramétrage décrit en \ref{parametrageorbitesnilpotentes} vaut pour nos deux algèbres. Modulo ce paramétrage, l'application ${\bf Rel}$ est l'identité. Cela résulte simplement du fait que cette 
 application "conserve les termes $x_{*,{\cal O}}$", cf. \ref{debacker}(5), et que ces termes appartiennent au m\^eme ensemble de classes de conjugaison car les groupes de Weyl de ${\bf G}$ et de $G$ sont les m\^emes.
  
   Le groupe $\Gamma_{F}$ agit sur $\boldsymbol{\mathfrak{g}}_{nil}/conj$ et $ \mathfrak{g}_{nil}/conj$. Il r\'esulte de sa construction que l'application ${\bf Rel}$ est \'equivariante pour ces actions. En particulier ${\bf Rel}((\boldsymbol{\mathfrak{g}}_{nil}/conj)^{I_{F}})=( \mathfrak{g}_{nil}/conj)^{I_{F}}$. 
 
  Soit ${\cal F}$ une facette de l'appartement $App_{F^{nr}}(T^{nr})$. On a d\'efini un diagramme d'applications
 $$\begin{array}{ccc}&&(\boldsymbol{\mathfrak{g}}_{nil}/conj)^{I_{F}}\\ &\qquad \nearrow \iota_{{\cal F},nil}&\\ \mathfrak{g}_{{\cal F},nil}/conj&&\,\, \downarrow {\bf Rel}\\ &\qquad \searrow Rel_{{\cal F}}&\\&&(\mathfrak{g}_{nil}/conj)^{I_{F}}\\ \end{array}$$
 Ce diagramme est commutatif.
 En effet, chaque application conserve les termes $x_{*,{\cal O}}$, en un sens convenable, et ces termes d\'eterminent les orbites nilpotentes.

   \subsubsection{L'ensemble $\bar{{\cal C}}_{F}^{\natural}$}\label{lensemblecalFnatural}
   
   Soit $e\in \mathfrak{g}_{nil}(F^{mod})$. Montrons que 
   
   (1) l'homomorphisme $Z_{G}(e)(F^{mod})/Z_{G}(e)^0(F^{mod})\to A(e)=Z_{G}(e)/Z_{G}(e)^0$ est bijectif.
   
   Fixons un $\mathfrak{sl}(2)$ triplet $(f,h,e)$ formé d'éléments de $\mathfrak{g}(F^{mod})$. Posons $Z=Z_{G}(f,h,e)$. On ne perd rien à remplacer $Z_{G}(e)$ par $Z$. 
   On peut remplacer $F$ par une extension finie modérément ramifiée et supposer que $f,h,e\in \mathfrak{g}(F)$. En vertu de notre hypothèse $(Hyp)_{1}(p)$, le groupe $G$  se déploie sur une extension finie et modérément ramifiée de $F$. Quitte à remplacer $F$ par une telle extension, on peut supposer  que $G$ est  déployé. Si $G$ est classique, des constructions d'algèbre élémentaire montrent que l'homomorphisme $ Z(F)/Z^0(F)\to A(e)$ est bijectif.  Supposons que $G$ soit exceptionnel. Le groupe $\Gamma_{F}$ agit sur $A(e)$ par un homomorphisme de $\Gamma_{F}$ dans le groupe $Aut(A(e))$ des automorphismes de $A(e)$.  Le groupe $A(e)$ est isomorphe à $\mathfrak{S}_{i}$ pour $i=1,...,5$. Tout automorphisme d'un tel groupe est intérieur. Notre hypothèse $(Hyp)_{1}(p)$ implique donc que l'ordre de $Aut(A(e))$ est premier à $p$. Quitte à remplacer $F$ par une extension finie modérément ramifiée, on peut donc supposer que $\Gamma_{F}$ agit trivialement sur $A(e)$. Soit $g\in Z$. Alors le cocycle galoisien $\sigma\mapsto g\sigma(g)^{-1}$ prend ses valeurs dans $Z^0$. On sait que l'application $H^1(\Gamma_{F}^{nr},Z^0(F^{nr}))\to H^1(\Gamma_{F},Z^0)$ est bijective. Le cocycle précédent est donc cohomologue à un cocycle trivial sur $I_{F}$. Autrement dit, on peut fixer $h\in Z^0$ de sorte que $hg\sigma(hg)^{-1}=1$ pour tout $\sigma\in I_{F}$. Alors $hg\in Z(F^{nr})$ et $hg$ a m\^eme image que $g$ dans $A(e)$. Cela prouve (1).
   
  Supposons $e\in \mathfrak{g}_{nil}(F)$.  L'assertion (1) entraîne que le groupe d'inertie sauvage $\Gamma_{F^{mod}}$ agit trivialement sur $A(e)$. Tout cocycle de $\Gamma_{F^{mod}}$ dans $A(e)$ est trivial puisque $A(e)$ est d'ordre premler à $p$. En conséquence,
   
   (2)  l'inflation $H^1(\Gamma_{F}^{mod},A(e))\to H^1(\Gamma_{F},A(e))$ est bijective.
   
   Montrons que
   
   (3) l'application $\mathfrak{g}_{nil}(F^{mod})/conj \to \mathfrak{g}_{nil}/conj$ est bijective. 
   
   Pour la m\^eme raison que ci-dessus, on peut supposer que $G$ est déployé sur $F$. Alors $\Gamma_{F}$ agit trivialement sur $\mathfrak{g}_{nil}/conj$, cf. \ref{actiongaloisienne}(1) et toute orbite dans cet ensemble contient un élément de $\mathfrak{g}_{nil}(F)$. Il reste à prouver que, si deux éléments $e,e'\in \mathfrak{g}_{nil}(F^{mod})$ sont conjugués par un élément de $G$, ils le sont par un élément de $G(F^{mod})$. Introduisons des $\mathfrak{sl}(2)$-triplets $(f,h,e)$ et $(f',h',e')$ formés d'éléments de $\mathfrak{g}(F^{mod})$. Quitte à remplacer $F$ par une extension finie et modérément ramifiée, on peut supposer que les éléments de ces triplets appartiennent à $\mathfrak{g}(F)$. 
  Ces triplets  sont conjugués par un élément de $G$ et on fixe $g\in G$ tel que $Ad(g)$ envoie $(f',h',e')$ sur $(f,h,e)$.  On pose $Z=Z_{G}(f,h,e)$.      L'application $\sigma\mapsto g\sigma(g)^{-1}$ est un cocycle de $\Gamma_{F}$ dans $Z$. Elle se quotiente en un cocycle de $\Gamma_{F}$ dans $A(e)=Z/Z^0$.  En appliquant (2), quitte \`a multiplier $g$ \`a gauche  par un \'el\'ement de $Z$, on peut supposer que ce cocycle est trivial sur $\Gamma_{F^{mod}}$. Autrement dit $g\sigma(g)^{-1}\in Z^0$ pour $\sigma\in \Gamma_{F^{mod}}$. Fixons une extension galoisienne finie $E$ de $F$ telle que $g\in G(E)$, posons $E'=E\cap F^{mod}$.  L'application $\sigma\mapsto g\sigma(g)^{-1}$ prend ses valeurs dans $Z^0$ pour $\sigma\in \Gamma_{F^{mod}}$ et pour $\sigma\in \Gamma_{E}$. Puisque c'est un cocycle et que ces deux groupes engendrent $\Gamma_{E'}$, cette application prend ses valeurs dans $Z^0$ pour tout $\sigma\in \Gamma_{E'}$. Puisque $Z^0$ est connexe, tout cocycle de $\Gamma_{E'}$ dans $Z^0$ devient cohomologue \`a $0$ dans $\Gamma_{E^{'nr}}$. Cela signifie que, quitte \`a multiplier encore $g$ \`a gauche par un \'el\'ement de $Z^0$, on a $g\sigma(g)^{-1}=1$ pour $\sigma\in \Gamma_{E^{'nr}}$. Alors $g\in G(E^{'nr})$. Puisque $E'\subset F^{mod}$, on a $E^{'nr}\subset F^{mod}$.  Donc $g\in G(F^{mod})$, ce qui d\'emontre (3). 
  
  L'application (3) est compatible aux actions galoisiennes. On en déduit la bijection
  $$(\mathfrak{g}_{nil}(F^{mod})/conj )^{\Gamma_{F}^{mod}}\to( \mathfrak{g}_{nil}/conj)^{\Gamma_{F}}.$$

    Notons $G^{\natural}$ le produit semi-direct $G(F^{mod})\rtimes \Gamma_{F}^{mod}$. Il se projette sur $\Gamma_{F}^{mod}$. Pour $\sigma\in \Gamma_{F}^{mod}$, on note $G_{\sigma}^{\natural}$ la fibre au-dessus de $\sigma$.    
Le groupe $G^{\natural}$ agit sur $G(F^{mod})$ et aussi sur $\mathfrak{g}(F^{mod})$. Soit $e\in \mathfrak{g}_{nil}(F^{mod})$.  L'orbite de $e$ est conservée par $\Gamma_{F}^{mod}$ si et seulement si $Z_{G^{\natural}_{\sigma}}(e)$ est non vide pour tout $\sigma\in \Gamma_{F}^{mod}$. Comme en \ref{calCF}, cela équivaut à ce que si $Z_{G^{\natural}_{\sigma}}(e)$  soit non vide pour $\sigma=Fr$ et $\sigma=\gamma$.  

{\bf Remarque.} (4) Si $e\in \mathfrak{g}(F)$, $Z_{G^{\natural}}(e)$ est le produit semi-direct $Z_{G}(e)(F^{mod})\rtimes \Gamma_{F}^{mod}$. 

\bigskip

L'ensemble $G^{\natural}_{1}$ est égal à $G(F^{mod})$.  Son intersection avec $Z_{G^{\natural}}(e)$ est $Z_{G}(e)(F^{mod})$, qui est donc un sous-groupe distingué de $Z_{G^{\natural}}(e)$. Il est clair que $Z_{G}(e)^0(F^{mod})$ est lui-aussi un sous-groupe distingué de $Z_{G^{\natural}}(e)$. On pose $A^{\natural}(e)=Z_{G}(e)^0(F^{mod})\backslash Z_{G^{\natural}}(e)$. Ce groupe  contient $A(e)$ comme sous-groupe distingué. Pour $\sigma\in \Gamma_{F}^{mod}$, on pose $A^{\natural}_{\sigma}(e)=Z_{G}(e)^0(F^{mod})\backslash Z_{G_{\sigma}^{\natural}}(e)$. 

Fixons un $\mathfrak{sl}(2)$-triplet $(f,h,e)$ d'éléments de $\mathfrak{g}(F^{mod})$. On peut reprendre les définitions ci-dessus en remplaçant $Z_{G}(e)$, $Z_{G^{\natural}}(e)$ etc... par $Z_{G}(f,h,e)$, $Z_{G^{\natural}}(f,h,e)$ etc... On obtient le m\^eme quotient $A^{\natural}(e)$.

Considérons les triplets $(e,b,u)$ où

$e$ est un élément de $\mathfrak{g}_{nil}(F^{mod})$ dont l'orbite est conservée par $\Gamma_{F}^{mod}$;

$b\in A^{\natural}_{\gamma}(e)$, $u\in A^{\natural}_{Fr}(e)$;

$ub=b^qu$. 

Le groupe $G(F^{mod})$ agit par conjugaison sur cet ensemble de triplets. On note ${\cal C}_{F}^{\natural}$ l'ensemble des classes de conjugaison. Pour une orbite $\bar{{\mathbb  O}}\in ( \mathfrak{g}_{nil}/conj)^{\Gamma_{F}}$, on note ${\cal C}_{F}^{\natural}(\bar{{\mathbb  O}})$ l'ensemble des classes de conjugaison de triplets $(e,b,u)$ pour lesquels $e\in \bar{{\mathbb  O}}$.

 On a défini en \ref{lesgroupesA(N)}    le quotient $\bar{A}(e)$ de $A(e)$ et le noyau $A(e)^1$ de l'homomorphisme $A(e)\to \bar{A}(e)$. On a

(5) supposons que l'orbite de $e$ soit conservée par $\Gamma_{F}^{mod}$; alors  le groupe $A(e)^1$ est un sous-groupe distingué dans $A^{\natural}(e)$. 

La preuve est la m\^eme qu'en \ref{actiongaloisienne}(2).

Sous l'hypothèse de (5), on pose $\bar{A}^{\natural}(e)=A^1(e)\backslash A^{\natural}(e)$. En remplaçant ci-dessus $A^{\natural}(e)$ par $\bar{A}^{\natural}(e)$, on définit l'ensemble $\bar{{\cal C}}_{F}^{\natural}$ et $\bar{{\cal C}}_{F}^{\natural}(\bar{{\mathbb  O}})$ pour $\bar{{\mathbb  O}}\in ( \mathfrak{g}_{nil}/conj)^{\Gamma_{F}}$.

\subsubsection{Comparaison des commutants\label{comparaison}}

Fixons une extension finie mod\'er\'ement ramifi\'ee $F'$ de $F^{nr}$ contenant $F^G$, autrement dit sur laquelle $G$ est d\'eploy\'e. 
 Pour simplifier la notation, on pose ${\bf K}=K_{s_{0},F'}^0$, $\boldsymbol{\mathfrak{k}}=\mathfrak{k}_{s_{0},F'}$ etc...  On note $k\mapsto \bar{k}$ l'application de r\'eduction naturelle de ${\bf K}$ dans ${\bf G}\simeq G_{s_{0,F'}}$. On note de m\^eme $X\mapsto \bar{X}$ l'application de r\'eduction de $\boldsymbol{\mathfrak{k}}$ dans $\boldsymbol{\mathfrak{g}}$. 
 
 Soit $(\bf{f},\bf{h},{\bf e})$ un $\mathfrak{sl}(2)$ triplet de $\boldsymbol{\mathfrak{g}}$. Comme on l'a dit en \ref{undiagramme}, on peut appliquer la construction de DeBacker en rempla\c{c}ant le corps de base $F$ par une extension finie convenable contenue dans $F'$. Ainsi, on fixe     un $\mathfrak{sl}(2)$-triplet $(f,e,h)$ de $\mathfrak{g}(F')$ relevant $(\bf{f},\bf{h},{\bf e})$ dans $\boldsymbol{\mathfrak{k}}$. On pose $Z=Z_{G}(f,h,e)$, ${\bf Z}=Z_{{\bf G}}({\bf f},{\bf h},{\bf e})$. 
 On a $Z/Z^0=A(e)$ d'o\`u un homomorphisme surjectif $Z\to \bar{A}(e)$ dont on note $Z^1$ le noyau. On d\'efinit de m\^eme ${\bf Z}^1$.

  \begin{prop}{ 
 (i) L'homomorphisme de r\'eduction $Z(F')\cap {\bf K}\to {\bf Z}$ est surjectif et se quotiente en un isomorphisme  
 $$(Z(F')\cap {\bf K})/(Z^0(F')\cap {\bf K})\simeq {\bf Z}/{\bf Z}^0=A({\bf e}).$$
 
 (ii) L'injection naturelle $Z(F')\cap {\bf K}\to Z$ se quotiente en un isomorphisme 
 $$(Z(F')\cap {\bf K})/(Z^0(F')\cap {\bf K})\simeq Z/Z^0=A(e).$$
 
 (iii) L'isomorphisme $A({\bf e})\to A(e)$ d\'eduit des isomorphismes de (i) et (ii) se quotiente en un isomorphisme $\bar{A}({\bf e})\to \bar{A}(e)$.  
 
   }\end{prop}

 Preuve de (i). Quitte \`a \'etendre le corps de base $F$, on peut supposer $F'=F^{nr}$, cela simplifie les notations. Soit $k\in {\bf K}$. Rempla\c{c}ons les triplets $(f,h,e)$,  resp. $({\bf f},{\bf h},{\bf e})$, par leurs images par $Ad(k)$, resp. $Ad(\bar{k})$.  L'assertion (i) pour ces triplets est \'equivalente \`a la m\^eme assertion pour les triplets d'origine. 
 Cela nous autorise \`a effectuer de telles conjugaison, ce qui nous permettra  d'imposer des conditions suppl\'ementaires \`a nos triplets.

 Soit ${\bf x}\in {\bf Z}$. On rel\`eve ${\bf x}$ en un \'el\'ement $k\in {\bf K}$. Notons $(f',h',e')$ l'image de $(f,h,e)$ par $Ad(k)$. Alors les triplets $(f,h,e)$ et $(f',h',e')$ rel\`event tous deux $(\bf{f},\bf{h},{\bf e})$ dans $\boldsymbol{\mathfrak{k}}$. On applique le lemme \ref{debacker}: il existe $k'\in {\bf K}^+$ tel que l'image de $(f',h',e')$ par $Ad(k')$ soit \'egale \`a $(f,h,e)$. Alors $Ad(k'k)$ conserve $(f,h,e)$, c'est-\`a-dire $k'k\in Z$. On a plus pr\'ecisement $k'k\in Z(F^{nr})\cap {\bf K}$ et l'image dans ${\bf Z}$ de $k'k$ est ${\bf x}$. Donc l'homomorphisme  de r\'eduction $Z (F^{nr})\cap {\bf K}\to  {\bf Z}$ est surjectif. 
 
 Montrons que
 
 (1) $\boldsymbol{\mathfrak{z}}$ est l'image par r\'eduction de $\mathfrak{z}(F^{nr})\cap \boldsymbol{\mathfrak{k}}$; a fortiori, $dim(Z^0)=dim({\bf Z}^0)$.
 
 On applique le lemme \ref{sl2} \`a l'action du triplet $(f,h,e)$ dans $\mathfrak{g}(F^{nr})$, cet espace \'etant muni de la cha\^{\i}ne de r\'eseaux r\'eduite \`a $\boldsymbol{\mathfrak{k}}\supset \mathfrak{p}_{F^{nr}}\boldsymbol{\mathfrak{k}}$. Ce lemme fournit une d\'ecomposition de $\mathfrak{g}(F^{nr})$ en somme directe de sous-espaces irr\'eductibles pour l'action du triplet $(f,h,e)$ et une d\'ecomposition correspondante de $\boldsymbol{\mathfrak{g}}$ en  somme directe de sous-espaces irr\'eductibles pour l'action du triplet $(\bf{f},\bf{h},{\bf e})$. 
  L'alg\`ebre $\mathfrak{z}(F^{nr})$ est la somme des sous-repr\'esentations de dimension $1$. Le lemme entra\^{\i}ne que $\mathfrak{z}(F^{nr})\cap \boldsymbol{\mathfrak{k}}$ se r\'eduit en la somme des sous-repr\'esentations de dimension $1$ de $\boldsymbol{\mathfrak{g}}$ et cette somme est justement $\boldsymbol{\mathfrak{z}}$. D'o\`u (1). 
  
  Fixons un sous-tore maximal ${\bf S}$ de ${\bf Z}^0$. Quitte \`a conjuguer le triplet $(\bf{f},\bf{h},{\bf e})$,on peut supposer que ${\bf S}\subset {\bf T}$. Notons ${\bf M}$ le commutant de ${\bf S}$ dans ${\bf G}$. C'est un Levi de ${\bf G}$. Puisque ${\bf S}\subset {\bf Z}$, on a $\bf{f},\bf{h},{\bf e}\in \boldsymbol{\mathfrak{m}}$. Donc $Z({\bf M})^0$ est contenu dans ${\bf Z}^0$. Par construction, on a ${\bf S}\subset Z({\bf M})^0$. La maximalit\'e de ${\bf S}$ entra\^{\i}ne que ${\bf S}=Z({\bf M})^0$. Par la bijection $X_{*}({\bf T})\simeq X_{*}(T)$, ${\bf S}$ se rel\`eve en un sous-tore $S$ de $T$. On note $M$ son commutant dans $G$. D'apr\`es \ref{groupesenreduction}(3), on a encore $S=Z(M)^0$, le point $s_{0}$ appartient \`a $Imm_{F^{nr}}^G(M_{ad})$ et, en notant $s_{0,M}$ l'image de $s_{0}$ dans $Imm_{F^{nr}}(M_{ad})$, on a $M_{s_{0,M}}={\bf M}$. Affectons d'un exposant $M$ les objets analogues pour $M$ \`a ceux introduits pour $G$, en particulier ${\bf K}^M=K_{s_{0,M}}^{M,0,nr}$, $\boldsymbol{\mathfrak{k}}^M=\mathfrak{k}_{s_{0,M}}^{M,nr}$. On a $\boldsymbol{\mathfrak{k}}^M=\boldsymbol{\mathfrak{k}}\cap \mathfrak{m}(F^{nr})$. Puisque le triplet  $(\bf{f},\bf{h},{\bf e})$ est contenu dans $\boldsymbol{\mathfrak{m}}=\mathfrak{m}_{x_{0,M}}^{nr}$, on peut le relever en un triplet contenu dans $\boldsymbol{\mathfrak{k}}^M$ et on peut supposer que notre triplet $(f,h,e)$ est un tel triplet, c'est-\`a-dire $f,h,e\in \mathfrak{m}(F^{nr})$. On a alors $S\subset Z^0$. 
 Montrons que
 
 (2) $S$ est un sous-tore maximal (sur $\bar{F}$) de $Z^0$, a fortiori $Z^0$ est d\'eploy\'e sur $F^{nr}$.
 
Il existe en tout cas un sous-tore maximal $S'$ de $Z^0$ qui est d\'efini sur $F^{nr}$ et contient $S$. Alors $\mathfrak{s}'(F^{nr})\cap \boldsymbol{\mathfrak{k}}$ est un sous-r\'eseau de $\mathfrak{z}(F^{nr})\cap\boldsymbol{ \mathfrak{k}}$ qui se r\'eduit en une sous-alg\`ebre  $\boldsymbol{\mathfrak{s}}'$ de  $\boldsymbol{\mathfrak{z}}$ d'apr\`es (1). Puisque $S\subset S'$ et que $S'$ est commutatif, $\boldsymbol{\mathfrak{s}}'$ est contenu dans le commutant de $\boldsymbol{\mathfrak{s}}$ dans $\boldsymbol{\mathfrak{z}}$. Puisque ${\bf S}$ est un sous-tore maximal de ${\bf Z}^0$, cela implique $\boldsymbol{\mathfrak{s}}'=\boldsymbol{\mathfrak{s}}$. Mais les dimensions de $S$ et $S'$ sont \'egales aux dimensions de $\boldsymbol{\mathfrak{s}}$ et $\boldsymbol{\mathfrak{s}}'$, donc $S$ et $S'$ sont de m\^eme dimension et l'inclusion $S\subset S'$ est une \'egalit\'e. Cela prouve (2). 

Notons $\Sigma^Z$ l'ensemble des racines de $S$ dans $\mathfrak{z}$. Pour $\alpha\in \Sigma^Z$, notons $\mathfrak{z}_{\alpha}$ l'espace radiciel associ\'e \`a $\alpha$ et $\boldsymbol{\mathfrak{z}}_{\alpha}$ l'image dans $\boldsymbol{\mathfrak{g}}$ de $\mathfrak{z}_{\alpha}(F^{nr})\cap \boldsymbol{\mathfrak{k}}$. L'espace $\boldsymbol{\mathfrak{z}}_{\alpha}$ est non nul et inclus dans $\boldsymbol{\mathfrak{z}}$ d'apr\`es (1). Il est clair que, modulo l'isomorphisme $X^*(S)\simeq X^*({\bf S})$, ${\bf S}$ agit dans $\boldsymbol{\mathfrak{z}}_{\alpha}$ par la racine $\alpha$. Cela implique que $\alpha$ est une racine de ${\bf S}$ dans $\boldsymbol{\mathfrak{z}}$ et que $\boldsymbol{\mathfrak{z}}_{\alpha}$ en est le sous-espace radiciel associ\'e. Par comparaison des dimensions, cf. (1), on voit que $\Sigma^Z$ s'identifie \`a l'ensemble de racines de ${\bf S}$ dans $\boldsymbol{\mathfrak{z}}$. 

Fixons une paire de Borel \'epingl\'ee $(B^Z,S,(H_{\alpha})_{\alpha\in \Delta^Z})$  de $Z$  conserv\'ee par l'action de $I_{F}$ et telle que, pour tout $\alpha\in \Delta^Z$,  $H_{\alpha}$ soit un g\'en\'erateur du $\mathfrak{o}_{F^{nr}}$-module $\mathfrak{z}_{\alpha}(F^{nr})\cap \boldsymbol{\mathfrak{k}}$. Comme en \ref{alcove}, on déduit de cette paire de Borel épinglée un sous-${\mathbb Z}_{p}$-module $\mathfrak{z}_{\alpha,{\mathbb Z}_{p}}$ de $\mathfrak{z}_{\alpha}$ pour tout $\alpha\in \Sigma^Z$ puis le $\mathfrak{o}_{F^{nr}}$-module $\mathfrak{z}_{\alpha,\mathfrak{o}_{F^{nr}}}$ engendré par $\mathfrak{z}_{\alpha,{\mathbb Z}_{p}}$. 
 Montrons que

(3) pour tout $\alpha\in \Sigma^Z$, on a l'égalité $\mathfrak{z}_{\alpha,\mathfrak{o}_{F^{nr}}}= \mathfrak{z}_{\alpha}(F^{nr})\cap \boldsymbol{\mathfrak{k}}$. 

Commençons par prouver (3) quand $-\alpha\in \Delta^Z$. On note  $H_{\alpha}$ l'unique élément de $\mathfrak{z}_{\alpha}$ tel que $[H_{\alpha},H_{-\alpha}]=\check{\alpha}^Z$, où $\check{\alpha}^Z$ est la coracine dans $X_{*}(S)\subset \mathfrak{z}$ associée à $\alpha$. Par définition, $\mathfrak{z}_{\alpha,\mathfrak{o}_{F^{nr}}}$ est le $\mathfrak{o}_{F^{nr}}$-module engendré par $H_{\alpha}$. L'élément $H_{\alpha}$  appartient à $\mathfrak{z}(F^{nr})$ car $Z^0$ est déployé sur $F^{nr}$. Donc il existe $\lambda\in F^{nr,\times}$ tel qu'en posant 
$H'_{\alpha}=\lambda H_{\alpha}$, on ait $\mathfrak{z}(F^{nr})\cap \boldsymbol{\mathfrak{k}}=\mathfrak{o}_{F^{nr}}H'_{\alpha}$. On a $[H'_{\alpha},H_{-\alpha}]=\lambda\check{\alpha}^Z$.  Puisque $H_{\alpha}$ et $ H'_{-\alpha}$ appartiennent \`a $\boldsymbol{\mathfrak{k}}$,  on a $\lambda\check{\alpha}^Z\in \boldsymbol{\mathfrak{k}}$. Par d\'efinition de $S$, on a $\mathfrak{s}(F^{nr})\cap \boldsymbol{\mathfrak{k}}=\mathfrak{s}(\mathfrak{o}_{F^{nr}})=X_{*}(S)\otimes_{{\mathbb Z}}\mathfrak{o}_{F^{nr}}$. La coracine $\check{\alpha}^Z$ n'est pas forc\'ement un \'el\'ement primitif de $X_{*}(S)$ mais elle l'est  dans $X_{*}(S_{sc})$, o\`u $S_{sc}$ est l'image r\'eciproque de $S$ dans $Z^0_{SC}$. Notre hypoth\`ese $(Hyp)_{1}(p)$ implique que $X_{*}(S_{sc})$ est un sous-groupe de $X_{*,{\mathbb Q}}(S_{sc})\cap X_{*}(S)$ d'indice premier \`a $p$.  La relation $\lambda\check{\alpha}^Z\in \boldsymbol{\mathfrak{k}}$ implique donc $\lambda\in \mathfrak{o}_{F^{nr}}$. Par r\'eduction de la relation $[H'_{\alpha},H_{-\alpha}]=\lambda\check{\alpha}^Z$, on obtient $[\bar{H}'_{\alpha},\bar{H}_{-\alpha}]=\bar{\lambda}\check{\alpha}^Z$, o\`u  $\check{\alpha}^Z$ est identifi\'e \`a un \'el\'ement de $X_{*}({\bf S})$. Puisque $\bar{H}'_{\alpha}$ et $\bar{H}_{-\alpha}$ sont des \'el\'ements non nuls de $\boldsymbol{\mathfrak{z}}_{\alpha}$ et $\boldsymbol{\mathfrak{z}}_{-\alpha}$, $[\bar{H}'_{\alpha},\bar{H}_{-\alpha}]$ est non nul. Cela implique $\bar{\lambda}\not=0$, ce qui \'equivaut \`a $\lambda\in \mathfrak{o}_{F^{nr}}^{\times}$. Mais alors, $H'{\alpha}$  et $H_{\alpha}$ engendrent le m\^eme $\mathfrak{o}_{F^{nr}}$-module,  ce qui prouve (3) pour $\alpha$. 

 Prouvons maintenant (3) pour toute racine. Rappelons qu'une racine $\alpha$ a une longueur $l(\alpha)$:  on écrit $\alpha=\epsilon\sum_{\beta\in \Delta^Z}m(\beta)\beta$, avec des entiers $m(\beta)\geq0$ et un signe $\epsilon\in \{\pm 1\}$; alors  $l(\alpha)=\sum_{\beta\in \Delta^Z}m(\beta)$.  On raisonne par r\'ecurrence sur $l(\alpha)$. L'assertion (3) est v\'erifi\'ee si $l(\alpha)=1$: c'est une 
 hypoth\`ese si $\alpha\in \Delta^Z$ et on vient de la vérifier si $-\alpha\in \Delta^Z$. Supposons $l(\alpha)>1$. 
   On  peut d\'ecomposer $\alpha$ en $\alpha_{1}+\alpha_{2}$ o\`u $\alpha_{1},\alpha_{2}\in \Sigma^Z$ v\'erifient $l(\alpha_{1})=l(\alpha)-1$, $l(\alpha_{2})=1$. Par l'hypothèse de récurrence, on a $\mathfrak{z}_{\alpha_{i},\mathfrak{o}_{F^{nr}}}=\mathfrak{z}_{\alpha_{i}}(F^{nr})\cap \boldsymbol{\mathfrak{k}}$ pour $i=1,2$. On a par définition $[\mathfrak{z}_{\alpha_{1},\mathfrak{o}_{F^{nr}}},\mathfrak{z}_{\alpha_{2},\mathfrak{o}_{F^{nr}}}]=\mathfrak{z}_{\alpha,\mathfrak{o}_{F^{nr}}}$. Donc $\mathfrak{z}_{\alpha,\mathfrak{o}_{F^{nr}}}\subset \mathfrak{z}(F^{nr})\cap \boldsymbol{\mathfrak{k}}$. Supposons que ces modules ne soient pas égaux, autrement dit $\mathfrak{z}_{\alpha,\mathfrak{o}_{F^{nr}}}\subset\mathfrak{z}(F^{nr})\cap \mathfrak{p}_{F^{nr}}\boldsymbol{\mathfrak{k}}$. Par réduction de la relation ci-dessus, on obtient $[\boldsymbol{\mathfrak{z}}_{\alpha_{1}},\boldsymbol{\mathfrak{z}}_{\alpha_{2}}]=0$. C'est contradictoire car $\alpha_{1}+\alpha_{2}=\alpha$ est une racine de ${\bf S}$. Cela prouve (3).

 Puisque $Z^0$ est d\'eploy\'e sur $F^{nr}$, la paire de Borel \'epingl\'ee que l'on a fix\'ee d\'etermine un point hypersp\'ecial $x^Z\in Imm_{F^{nr}}(Z^0_{AD})$.  Par d\'efinition, la sous-alg\`ebre parahorique $\boldsymbol{\mathfrak{k}}^Z$ associ\'ee \`a $x^Z$ est la somme de $\mathfrak{s}(\mathfrak{o}_{F^{nr}})$ et des $\mathfrak{z}_{\alpha,\mathfrak{o}_{F^{nr}}}$ pour $\alpha\in \Sigma^Z$. Il résulte de (3) et du fait que la réduction de $\mathfrak{z}_{\alpha}(F^{nr})\cap \boldsymbol{\mathfrak{k}}$ est $\boldsymbol{\mathfrak{z}}_{\alpha}$ que $\boldsymbol{\mathfrak{k}}^Z\subset  \mathfrak{z}(F^{nr})\cap\boldsymbol{\mathfrak{k}}$ et que $\boldsymbol{\mathfrak{k}}^Z$ s'envoie surjectivement sur $\boldsymbol{\mathfrak{z}}$.   Le lemme de Nakayama entra\^{\i}ne alors $\boldsymbol{\mathfrak{k}}^Z=\mathfrak{z}(F^{nr})\cap \boldsymbol{\mathfrak{k}}$.

Montrons que

(4) l'application $Z(F^{nr})\cap {\bf K}\to {\bf Z}$ se restreint en une surjection $Z^0(F^{nr})\cap {\bf K}\to {\bf Z}^0$. 

Introduisons la complétion $\hat{F}^{nr}$ de $F^{nr}$. Notons $\hat{{\bf K}}^Z$, resp. ${\bf K}^Z$, le sous-groupe parahorique de $Z^0(\hat{F}^{nr})$, resp. $Z^0(F^{nr})$, associ\'e \`a $x^Z$.  Notons aussi $\hat{{\bf K}}$ l'analogue de ${\bf K}$ sur $\hat{F}^{nr}$.  Etant associé à un point hypersp\'ecial, le groupe $\hat{{\bf K}}^Z$ est un sous-groupe compact maximal de $Z^0(\hat{F}^{nr})$. Il est engendr\'e par $S(\mathfrak{o}_{\hat{F}^{nr}})$ et les exponentielles $exp(H)$ pour $H\in \mathfrak{z}_{\alpha,\mathfrak{o}_{\hat{F}^{nr}}}$, $\alpha$ décrivant $\Sigma^Z$.   Tous ces \'el\'ements appartiennent \`a $\hat{{\bf K}}$ d'apr\`es (3) et \ref{debacker}(1) (qui s'étend sur le corps $\hat{F}^{nr}$ par continuité). 
Il en r\'esulte que $\hat{{\bf K}}^Z\subset Z^0(\hat{F}^{nr})\cap \hat{{\bf K}}$. Cette intersection est un sous-groupe compact de $Z^0(\hat{F}^{nr})$. D'après  la maximalit\'e de $\hat{{\bf K}}^Z$, cela implique l'\'egalit\'e  $\hat{{\bf K}}^Z= Z^0(\hat{F}^{nr})\cap \hat{{\bf K}}$. Par intersection avec $G(F^{nr})$, on obtient ${\bf K}^Z=Z^0(F^{nr})\cap {\bf K}$.  L'image de $Z^0(F^{nr})\cap {\bf K}$ dans $ {\bf Z}^0$ est donc l'image de ${\bf K}^Z$ par r\'eduction  dans ${\bf G}$. Cette image est engendr\'ee par ${\bf S}$ et les exponentielles $exp(\bar{\lambda} \bar{H})$ pour $\lambda\in \mathfrak{o}_{F^{nr}}$ et $\bar{H}\in \boldsymbol{\mathfrak{z}}_{\alpha}$, $\alpha$ décrivant $\Sigma^Z$.  Ces \'el\'ements engendrent ${\bf Z}^0$. Cela prouve (4). 

Il r\'esulte de (4) que l'image r\'eciproque de ${\bf Z}^0$ dans $Z(F^{nr})\cap {\bf K}$ est \'egale \`a  $(Z^0(F^{nr})\cap {\bf K})(Z(F^{nr})\cap {\bf K}^+)$. Pour achever la preuve de l'assertion (i) de l'\'enonc\'e, il suffit de prouver que $Z(F^{nr})\cap {\bf K}^+\subset Z^{0}(F^{nr})$. Or $Z(F^{nr})\cap {\bf K}^+$ est un pro-$p$-groupe. Il suffit donc de prouver que  le nombre d'\'el\'ements de $Z(F^{nr})/Z^0(F^{nr})$ est premier \`a $p$. Evidemment, ce quotient s'injecte dans $Z/Z^0$. Or $Z/Z^0=A(e)$ et ce groupe  est soit un produit de groupes ${\mathbb Z}/2{\mathbb Z}$, soit $\mathfrak{S}_{3}$, $\mathfrak{S}_{4}$ ou $\mathfrak{S}_{5}$. En tout cas, il est d'ordre premier \`a $p$. Cela ach\`eve la preuve de l'assertion (i) de l'\'enonc\'e. 

Preuve de (ii). L'injection $Z(F')\cap {\bf K}\to Z$ se quotiente en un 
homomorphisme injectif $(Z(F')\cap {\bf K})/(Z^0(F')\cap {\bf K})\to Z/Z^0$. Par composition avec l'isomorphisme de l'assertion (i), on obtient un homorphisme injectif ${\bf Z}/{\bf Z}^0\to Z/Z^0$, autrement dit $A({\bf e})\to A(e)$.  Les orbites de $e$ et ${\bf e}$ se correspondent selon les param\'etrages habituels de $\mathfrak{g}_{nil}/conj$ et $\boldsymbol{\mathfrak{g}}_{nil}/conj$ comme on l'a dit en \ref{undiagramme}.  Deux orbites se correspondant ont m\^eme groupe $A$ associ\'e, celui-ci ne d\'ependant pas du corps de base, pourvu que la caract\'eristique de celui-ci soit assez grande si elle n'est pas nulle. Donc $A({\bf e})$ et $A(e)$ sont isomorphes. Puisqu'ils sont finis, l'injection pr\'ec\'edente est bijective. Cela entra\^{\i}ne que l'injection $(Z(F')\cap {\bf K})/(Z^0(F')\cap {\bf K})\to Z/Z^0$ est elle-aussi bijective. D'o\`u (ii).

  Preuve de (iii).   Dans le cas des groupes classiques, on peut expliciter l'isomorphisme  $A({\bf e})\simeq A(e)$ en utilisant des constructions d'alg\`ebre linaire. On constate qu'elle devient l'identit\'e modulo les descriptions combinatoires que l'on a donn\'ees de ces deux groupes. L'assertion (iii) en r\'esulte. On laisse la mise au point de la construction au lecteur. 
  Dans le cas des groupes exceptionnels, on montre que tout isomorphisme $A({\bf e})\simeq A(e)$ se quotiente en un isomorphisme $\bar{A}({\bf e})\simeq \bar{A}(e)$. Cela resulte du fait que les paires $(A({\bf e}),\bar{A}({\bf e}))$ et $(A(e),\bar{A}(e))$ sont isomorphes et forc\'ement isomorphes \`a l'une des paires $(A(e),A(e))$, $(A(e),\{1\})$ ou $(\mathfrak{S}_{3},\mathfrak{S}_{2})$. $\square$

  \subsubsection{Commutants et action galoisienne}\label{comparaisongalois}
  On conserve les hypothèses et notations du paragraphe précédent. 
  Notons  $\boldsymbol{K}^{\natural}$  le sous-groupe $\boldsymbol{K}\rtimes \Gamma_{F}^{mod}$ de $G^{\natural}$. L'application de réduction de $\boldsymbol{K}$ sur ${\bf G}$ est équivariante pour les actions de $\Gamma_{F}^{mod}$. Elle se prolonge naturellement en un homomorphisme surjectif $\boldsymbol{K}^{\natural}\to {\bf G}^{++}$.   
  Pour $\sigma\in \Gamma_{F}^{mod}$, posons $Z_{\sigma}=Z_{G_{\sigma}^{\natural}}(f,h,e)$ et ${\bf Z}_{\sigma}= Z_{{\bf G}^{++}_{\sigma}}({\bf f},{\bf h},{\bf e})$.

  \begin{prop}{Supposons que l'orbite de ${\bf e}$ est conservée par l'action de $\Gamma_{F}^{mod}$. Soit $\sigma\in \Gamma_{F}^{mod}$. 
  
  (i)  L'application de réduction se quotiente en des bijections
  $$(Z^0\cap \boldsymbol{K})\backslash (Z_{\sigma}\cap\boldsymbol{K}^{\natural})\to {\bf Z}^0\backslash {\bf Z}_{\sigma}.$$
  
  $$(Z^1\cap \boldsymbol{K})\backslash (Z_{\sigma}\cap\boldsymbol{K}^{\natural})\to {\bf Z}^1\backslash {\bf Z}_{\sigma}.$$
  
  (ii) L'injection naturelle $Z_{\sigma}\cap \boldsymbol{K}^{\natural}\to Z_{\sigma}$ se quotiente en des bijections
  $$(Z^0\cap \boldsymbol{K})\backslash (Z_{\sigma}\cap\boldsymbol{K}^{\natural})\to Z^0\backslash Z_{\sigma}.$$
  $$(Z^1\cap \boldsymbol{K})\backslash (Z_{\sigma}\cap\boldsymbol{K}^{\natural})\to Z^1\backslash Z_{\sigma}.$$}\end{prop}

  Preuve. On utilise la proposition précédente.  L'isomorphisme du (ii) se quotiente en un isomorphisme  
  $$(Z^1\cap {\bf K})\backslash(Z\cap {\bf K})\to Z^1\backslash Z=\bar{A}(e).$$
  Il r\'esulte de (iii) que l'isomorphisme du (i) se quotiente lui-aussi en un isomorphisme
  $$(Z^1\cap {\bf K})\backslash(Z\cap {\bf K})\to {\bf Z}^1\backslash{\bf Z}=\bar{A}({\bf e}).$$
  Si $Z_{\sigma}\cap {\bf K}^{\natural}$ n'est pas vide, c'est 
 une unique classe pour la multiplication \`a gauche par  $Z\cap {\bf K}$. Alors  l'énoncé   se d\'eduit imm\'ediatement de ces   isomorphismes et de ceux de la proposition \ref{comparaison}.  Tout revient \`a prouver que $Z_{\sigma}\cap {\bf K}^{\natural}$ n'est pas vide. Le triplet $(\sigma(f),\sigma(h),\sigma(e))$ est contenu dans $\boldsymbol{\mathfrak{k}}$ et rel\`eve le triplet $( \sigma(\bf{f}),\sigma(\bf{h}),\sigma({\bf e}))$ de $\boldsymbol{\mathfrak{g}}$. Par hypoth\`ese, l'orbite de ${\bf e}$ est conserv\'ee par l'action galoisienne. Donc $( \sigma(\bf{f}),\sigma(\bf{h}),\sigma({\bf e}))$ est conjugu\'e \`a $({\bf f},{\bf h},{\bf e})$ par un \'el\'ement de ${\bf G}$. Le lemme \ref{debacker} implique qu'il existe $k\in {\bf K}$ tel que  $Ad(k)$ envoie $(\sigma(f),\sigma(h),\sigma(e))$ sur $(f,e,h)$.  Pour un tel $k$, on a $(k,\sigma)\in Z_{\sigma}\cap {\bf K}^{\natural}$. Cela ach\`eve la preuve. $\square$
 
 \subsubsection{L'isomorphisme ${\cal C}^{\natural}_{F}\to {\cal C}_{F}^{\sharp}$}\label{lisomorphisme}
 
 Soit $\boldsymbol{{\cal O}}\in (\boldsymbol{\mathfrak{g}}_{nil}/conj)^{\Gamma_{F}}$, posons $\bar{{\mathbb O}}={\bf Rel}(\boldsymbol{{\cal O}})$. Soit $N\in \boldsymbol{{\cal O}}$. Complétons $N$ en un $\mathfrak{sl}(2)$-triplet $({\bf f},{\bf h},N)$ et appliquons les constructions de \ref{comparaison}. En particulier, on fixe un $\mathfrak{sl}(2)$-triplet $(f,h,e)$ de $\mathfrak{g}(F')$ relevant $({\bf f},{\bf h},N)$. On a $e\in \bar{{\mathbb O}}$. 
 Posons ${\bf Z}^{++}=Z_{{\bf G}^{++}}({\bf f},{\bf h},N)$, $Z^{\natural}=Z_{G^{\natural}}(f,h,e)$. De  la proposition \ref{comparaisongalois} se déduisent deux isomorphismes
 $$(1) \qquad A^{\natural}(e)\simeq (Z^0\cap \boldsymbol{K}^{\natural})\backslash (Z^{\natural}\cap \boldsymbol{K}^{\natural})\simeq A^{++}(N),$$
 où l'on a utilisé les notations de \ref{comparaisongalois}. L'ensemble ${\cal C}^{\sharp}_{F}(\boldsymbol{{\cal O}})$ est celui des classes de conjugaison par $A(N)$ dans l'ensemble des couples $(b,u)\in A^{++}_{\gamma}(N)\times A^{++}_{Fr}(N)$ tels que $ub=b^qu$. De m\^eme, ${\cal C}_{F}^{\natural}(\bar{{\mathbb O}})$ est celui  des classes de conjugaison par $A(e)$ dans l'ensemble des couples $(b,u)\in A^{\natural}_{\gamma}(e)\times A^{\natural}_{Fr}(e)$ tels que $ub=b^qu$. La bijection ci-dessus identifie ces deux ensembles. On en déduit une bijection
 $$(2) \qquad  C:  {\cal C}_{F}^{\natural}(\bar{{\mathbb O}})\simeq {\cal C}^{\sharp}_{F}(\boldsymbol{{\cal O}}).$$
 Cette bijection est indépendante des choix effectués. En effet, remplaçons nos triplets $({\bf f},{\bf h},N)$ et $(f,h,e)$  par d'autres, notés $({\bf f}',{\bf h}',N')$ et $(f',h',e')$. Il existe alors $g\in \boldsymbol{G}$ tel que $Ad(k)$ envoie $({\bf f},{\bf h},N)$ sur $({\bf f}',{\bf h}',N')$. D'après le lemme \ref{debacker}, il. existe un relèvement $k\in \boldsymbol{K}$ de $g$ tel que $Ad(k)$ envoie $(f,h,e)$ sur $(f',h',e')$. Les isomorphismes (1) pour nos deux couples de triplets s'insèrent dans un diagramme commutatif
 $$\begin{array}{ccc} A^{\natural}(e)&\simeq& A^{++}(N)\\ &&\\\,\,\downarrow Ad(k)&&\,\,\downarrow Ad(g)\\ &&\\ A^{\natural}(e')&\simeq& A^{++}(N')\\ \end{array}$$
 On en déduit que l'application (2) n'a pas changé.
 
 En faisant varier $\boldsymbol{{\cal O}}$, on obtient une bijection que l'on note encore
 $$C:{\cal C}_{F}^{\natural} \simeq {\cal C}^{\sharp}_{F} .$$
 
 La m\^eme construction fournit des isomorphismes
 $$\bar{C}:\bar{{\cal C}}^{\natural}_{F}(\bar{{\mathbb O}})\simeq \bar{{\cal C}}_{F}^{\sharp}(\boldsymbol{{\cal O}}),\,\, \bar{C}:\bar{{\cal C}}^{\natural}_{F}\simeq \bar{{\cal C}}_{F}^{\sharp}.$$

  \subsubsection{Interpr\'etation cohomologique de l'application $ \bar{c}_{{\cal F},F}$\label{cohomologie}}
 Soit $\bar{{\mathbb O}}\in (\mathfrak{g}_{nil}/conj)^{\Gamma_{F}}$. L'intersection  $\bar{{\mathbb O}}\cap \mathfrak{g}(F)$ n'est pas vide puisque $G$ est quasi-d\'eploy\'e. Notons $(\bar{{\mathbb O}}\cap \mathfrak{g}(F))/conj$ l'ensemble des classes de conjugaison par $G(F)$ dans $\bar{{\mathbb O}}\cap \mathfrak{g}(F)$. Soit $e\in \bar{{\mathbb O}}\cap \mathfrak{g}(F)$. On d\'efinit l'ensemble de cohomologie $H^1(\Gamma_{F},Z_{G}(e))$. 
  Il s'envoie dans $H^1(\Gamma_{F},G)$, on note $H^1(\Gamma_{F},Z_{G}(e))_{1}$ le noyau de cette application (au sens des applications entre ensembles point\'es).  
  Soit $e'\in \bar{{\mathbb O}}\cap \mathfrak{g}(F)$. Choisissons un \'el\'ement $r\in G$ tel que $e'=r^{-1}er$. Soit $u'$ un cocycle de $\Gamma_{F}$ dans $Z_{G}(e')$. Pour $\sigma\in \Gamma_{F}$, posons $u(\sigma)=ru'(\sigma)\sigma(r)^{-1}$. C'est 
  un cocycle de $\Gamma_{F}$ dans $Z_{G}(e)$. L'application $u'\mapsto u$ se quotiente en une bijection $r_{e,e'}:H^1(\Gamma_{F},Z_{G}(e'))\to H^1(\Gamma_{F},Z_{G}(e))$ qui ne d\'epend pas du choix de $r$. Celle-ci    se restreint en une bijection de  $H^1(\Gamma_{F},Z_{G}(e'))_{1}$ sur $ H^1(\Gamma_{F},Z_{G}(e))_{1}$. Notons $ 1_{e'}$ la classe triviale dans $H^1(\Gamma_{F},Z_{G}(e'))$ et $1_{e,e'}$ son image par $r_{e,e'}$.   L'application $e'\mapsto 1_{e,e'}$ se quotiente en une bijection de $(\bar{{\mathbb O}}\cap \mathfrak{g}(F))/conj$ sur $H^1(\Gamma_{F},Z_{G}(e))_{1}$.

  D'après \ref{actiongaloisienne}(1), de l'homomorphisme naturel $Z_{G}(e)\to \bar{A}(e)$ se d\'eduit une application 
  
   $$d_{e}:H^1(\Gamma_{F},Z_{G}(e))\to H^1(\Gamma_{F},\bar{A}(e)).$$
   Pour $e'\in \bar{{\mathbb O}}\cap \mathfrak{g}(F)$, l'application $r_{e,e'}$ se quotiente en une application
   $$\bar{r}_{e,e'}:H^1(\Gamma_{F},\bar{A}(e')) \to H^1(\Gamma_{F},\bar{A}(e)).$$
   
   L'ensemble $\bar{{\cal C}}_{F}^{\natural}(\bar{{\mathbb O}})$ s'identifie à celui des classes de conjugaison par $\bar{A}(e)$ dans l'ensemble des couples $(d,v)\in \bar{A}^{\natural}_{\gamma}(e)\times \bar{A}^{\natural}_{Fr}(e)$ tels que $vd=d^qv$. L'inflation $H^1(\Gamma_{F}^{mod},\bar{A}(e))\to H^1(\Gamma_{F},\bar{A}(e))$ est bijective (c'est l'analogue de  \ref{lensemblecalFnatural}(2)). A un cocycle  $u:\Gamma_{F}^{mod}\to \bar{A}(e)$, associons le couple $(d,v)=(u(\gamma)\gamma,u(Fr)Fr)$. Il appartient à $\bar{{\cal C}}_{F}^{\natural} (\bar{{\mathbb O}})$ (ou plut\^ot il représente une classe dans cet ensemble) et l'application $u\mapsto (d,v)$ se quotiente en une bijection $\bar{c}_{e}: H^1(\Gamma_{F}^{mod},\bar{A}(e))\to \bar{{\cal C}}_{F}^{\natural}(\bar{{\mathbb O}})$. Notons $\boldsymbol{{\cal O}}$ l'élément de $(\boldsymbol{\mathfrak{g}}_{nil}/conj)^{\Gamma_{F}}$ tel que ${\bf Rel}(\boldsymbol{{\cal O}})=\bar{{\mathbb O}}$. Pour $e'\in \bar{{\mathbb O}}\cap \mathfrak{g}(F)$, on obtient un diagramme commutatif

 $$(1) \qquad \begin{array}{ccccccc}H^1(\Gamma_{F},Z_{G}(e'))&\stackrel{d_{e'}}{\to}&H^1(\Gamma_{F},\bar{A}(e'))&&&&\\ &&&\searrow \bar{c}_{e'}&&&\\
 
 \downarrow r_{e,e'}&&\downarrow  {\bf r}_{e,e'}&&\bar{{\cal C}}_{F}^{\natural}(\bar{{\mathbb O}})&\stackrel{\bar{C}}{\to}&\bar{{\cal C}}^{\sharp}_{F}(\boldsymbol{{\cal O}})\\&&&\nearrow \bar{c}_{e}&&&\\ H^1(\Gamma_{F},Z_{G}(e))&\stackrel{d_{e}}{\to}&H^1(\Gamma_{F},\bar{A}(e))&&&&\\ \end{array}$$

  Dans la proposition suivante, on utilise les applications définies en \ref{iotanil} et \ref{undiagramme}.

  \begin{prop}{  (i) 
  Soient ${\cal F}$ une facette de $App(T_{F})$ et ${\cal O}\in \mathfrak{g}_{{\cal F},nil}({\mathbb F}_{q})/conj$. Supposons que $Rel_{{\cal F}}({\cal O})=\bar{{\mathbb O}}$. 
  Soit $e$ un \'el\'ement de l'orbite $rel_{{\cal F},F}({\cal O})$.  Alors $\bar{C}\circ \bar{c}_{e}\circ d_{e}( 1_{e})=\bar{c}_{{\cal F},F}({\cal O})$.

  (ii) Pour tout $e\in \bar{{\mathbb O}}\cap \mathfrak{g}(F)$, l'application $d_{e}$ se restreint en une surjection $$H^1(\Gamma_{F},Z_{G}(e))_{1}\to H^1(\Gamma_{F},\bar{A}(e)).$$}\end{prop} 
  
  Preuve.  Prouvons (i). L'assertion ne dépend pas du choix de $e\in rel_{{\cal F},F}({\cal O})$: cela résulte du diagramme   (1) et du fait que, si $e'$ et $e$ sont conjugués par un élément de $G(F)$, on a l'égalité $1_{e,e'}=1_{e}$. On peut donc fixer un $\mathfrak{sl}(2)$-triplet $(\underline{f},\underline{h},\underline{e})$ de $\mathfrak{g}_{{\cal F}}({\mathbb F}_{q})$ tel que $\underline{e}\in {\cal O}$, le relever en un $\mathfrak{sl}(2)$-triplet $(f,h,e)$ dans 
  $\mathfrak{k}_{{\cal F}}$ et supposer que $e$ est l'élément ainsi noté de ce triplet. Fixons $x\in a_{T}({\mathbb T}_{F})\cap {\cal F}$. 
  Posons $\iota_{{\cal F}}(\underline{f},\underline{h},\underline{e})=({\bf f},{\bf h},{\bf e})$. C'est un $\mathfrak{sl}(2)$-triplet de $\boldsymbol{\mathfrak{g}}$. On a ${\bf e}\in \boldsymbol{{\cal O}}$. Par définition de $\bar{c}_{{\cal F},F}$, $\bar{c}_{{\cal F},F}({\cal O})$ est l'image naturelle dans $\bar{{\cal C}}^{\sharp}_{F}(\boldsymbol{{\cal O}})$ du triplet $({\bf e},{\bf j}_{T}(x)\gamma,Fr)$. Notons $\tilde{x}$ l'élément de ${\mathbb T}_{F}$ tel que $a_{T}(\tilde{x})=x$. Fixons une extension finie modérément ramifiée de $F^{nr}$, contenant $F^G$ et telle que $\underline{v}(\tilde{x})$ appartiennne à $T_{F}(F')$. Posons $Ad(\underline{v}(\tilde{x}))(f,h,e)=(f_{0},h_{0},e_{0})$. C'est un $\mathfrak{sl}(2)$-triplet dans $\boldsymbol{\mathfrak{k}}=\mathfrak{k}_{s_{0},F'}$. Par définition de $\iota_{{\cal F}}$, il relève le triplet $({\bf f},{\bf h},{\bf e})$. L'élément $\bar{c}_{e}\circ d_{e}( 1_{e})$ est l'image naturelle dans $\bar{{\cal C}}_{F}^{\natural}(\bar{{\mathbb O}})$ du triplet $(e,\gamma,Fr)$. En conjuguant par $\underline{v}(\tilde{x})$, c'est aussi l'image de $(e_{0},\underline{v}(\tilde{x})\gamma \underline{v}(\tilde{x})^{-1}, \underline{v}(\tilde{x})Fr\,\underline{v}(\tilde{x})^{-1})$. On a $\underline{v}(\tilde{x})\gamma \underline{v}(\tilde{x})^{-1}=t_{\tilde{x}}(\gamma)\gamma$. On a aussi  $\underline{v}(\tilde{x})Fr\,\underline{v}(\tilde{x})^{-1})=Fr$ car, parce que $\tilde{x}\in {\mathbb T}_{F}$, $\underline{v}(\tilde{x})$ est fixé par $Fr$. Donc $\bar{c}_{e}\circ d_{e}( 1_{e})$ est l'image naturelle dans $\bar{{\cal C}}_{F}^{\natural}({\cal O})$ du triplet $(e_{0},t_{\tilde{x}}(\gamma)\gamma, Fr)$. L'élément $t_{\tilde{x}}(\gamma)\gamma$ appartient à l'ensemble $\boldsymbol{K}^{\natural}$ de \ref{comparaisongalois} et la réduction de $t_{\tilde{x}}(\gamma)$ dans ${\bf G}$ est ${\bf j}_{T}(x)$. Il résulte alors de la définition de l'application $\bar{C}$ que l'image de  $\bar{c}_{e}\circ d_{e}( 1_{e})$  par $\bar{C}$ est l'image naturelle dans $\bar{{\cal C}}^{\sharp}_{F}(\boldsymbol{{\cal O}})$ du triplet $({\bf e},{\bf j}_{T}(x)\gamma,Fr)$. Cela prouve l'assertion (i) de l'énoncé.  
  
  Prouvons (ii). Soit $u\in H^1(\Gamma_{F},\bar{A}(e))$, posons ${\bf u}= \bar{C}\circ \bar{c}_{e}(u)$. C'est un élément de $\bar{{\cal C}}_{F}^{\sharp}(\boldsymbol{{\cal O}})$. D'après la proposition \ref{iotanil}, on peut fixer un sommet $s\in S(\bar{C})$ et une orbite ${\cal O}\in \mathfrak{g}_{s,nil}({\mathbb F}_{q})/conj$ de sorte que $\bar{c}_{{\cal F},F}({\cal O})={\bf u}$. Soit $e'\in rel_{s,F}({\cal O})$. On a $e'\in \bar{{\mathbb O}}$ et, d'après l'assertion (i), on a l'égalité $\bar{C}\circ \bar{c}_{e'}\circ d_{e'}( 1_{e'})=\bar{c}_{{\cal F},F}({\cal O})={\bf u}$. En appliquant le diagramme (1), on a aussi $\bar{C}\circ \bar{c}_{e}\circ d_{e}( 1_{e,e'})= {\bf u}$. En se rappelant la définition de ${\bf u}$, et parce que les applications $\bar{C}$ et $\bar{c}_{e}$ sont bijectives, on obtient l'égalité $u=d_{e}(1_{e,e'})$. Mais $1_{e,e'}$ appartient à $H^1(\Gamma_{F},Z_{G}(e))_{1}$ donc $u$ appartient à l'image de cet ensemble dans $H^1(\Gamma_{F},\bar{A}(e))$. Cela achève la preuve de (ii). $\square$

\section{Normalisation des actions de Frobenius}

\subsection{Isomorphismes d'espaces $FC^{st}(\mathfrak{g}(F))$}

\subsubsection{Position du problème}\label{position}

Soit $G$ un groupe réductif connexe défini sur $F$. On impose la condition $(Hyp)_{1}(p)$. {\bf On suppose que} $G$ {\bf est adjoint et absolument simple}.   En \ref{couplesstables}, on a raffiné selon la structure de $G$ sur $F^{nr}$ la classification habituelle du système de racines de $G$ et nous utiliserons cette classification raffinée. 

Soit $H$ un groupe réductif connexe défini sur $F$  qui est une forme intérieure de $G$ et fixons un torseur intérieur $\phi:H\to G$. On a dit en \ref{parametrage} que de ce torseur se déduisait un isomorphisme de transfert endoscopique de $FC^{st}(\mathfrak{h}(F))$ sur $ FC^{st}(\mathfrak{g}(F))$. Sa définition n'est pas reliée  aux bases de ces espaces déduites de faisceaux-caractères cuspidaux à support nilpotent sur des algèbres réduites $\mathfrak{g}_{s}$.  
On se propose dans cette sous-section de définir un isomorphisme 
$$ \iota_{\phi}:FC^{st}(\mathfrak{h}(F))\to FC^{st}(\mathfrak{g}(F))$$
qui, lui, se lit sur ces bases. On prouvera dans la sous-section suivante qu'il co\"{\i}ncide avec le transfert endoscopique. 
 
 Dans les cas où $G$ est de l'un des types   $(D_{n},nr)$ avec $n$ impair,  $(D_{n},ram)$ avec $n$ pair  ou $(A_{n-1},nr)$, la question est triviale puisque les deux espaces intervenant ci-dessus sont nuls. On exclut ces cas dans la suite. 

\subsubsection{Familles admissibles}\label{familles}

  Pour un \'el\'ement   $s\in S(G)$, c'est-\`a-dire que $s$ est un sommet de $Imm(G)$, on a not\'e ${\bf FC}(\mathfrak{g}_{s})$ l'ensemble des  faisceaux-caract\`eres cuspidaux \`a support nilpotent sur $\mathfrak{g}_{s}$ et ${\bf FC}_{{\mathbb F}_{q}}(\mathfrak{g}_{s})$ le sous-ensemble des faisceaux  conserv\'es  par l'action de Frobenius. Pour simplifier les notations, nous notons ces ensembles ${\bf FC}_{s}$ et ${\bf FC}_{s,{\mathbb F}_{q}}$. 
 L'espace $FC^{st}(\mathfrak{g}_{s}({\mathbb F}_{q}))$  est engendr\'e par les fonctions caract\'eristiques de certains \'el\'ements de ${\bf FC}_{s,{\mathbb F}_{q}}$. On  note ${\bf FC}_{s,{\mathbb F}_{q}}^{st}$ le sous-ensemble de ces faisceaux.  On a noté $S^{st}(G)$ l'ensemble des $s\in S(G)$ tels que $FC^{st}(\mathfrak{g}_{s}({\mathbb F}_{q}))\not=\{0\}$.  Soit $E$ une extension finie non ramifi\'ee de degr\'e $n$ de $F$. D'apr\`es \ref{extensionsnonram}, si $s\in S^{st}(G)$, on a aussi $s\in S_{E}^{st}(G)$ et ${\bf FC}_{s,{\mathbb F}_{q}}^{st}={\bf FC}_{s,{\mathbb F}_{q^n}}^{st}\cap {\bf FC}_{s,{\mathbb F}_{q}}$.  On a déjà défini l'ensemble  
   $S^{nr,st}(G)$  des sommets $s$ de $Imm_{F^{nr}}(G)$ pour lesquels il existe une extension $E$ comme ci-dessus de sorte que $s\in S_{E}^{st}(G)$. Pour 
 $s\in S^{nr,st}(G)$, on note  ${\bf FC}^{st}_{s}$  l'ensemble des ${\cal E}\in {\bf FC}_{s}$ pour lesquels il existe une extension $E$ comme ci-dessus de sorte que ${\cal E}\in {\bf FC}^{st}_{s,{\mathbb F}_{q^n}}$.

  Soit $s\in S^{nr,st}(G)$. Un  faisceau-caract\`ere cuspidal \`a support nilpotent  ${\cal E}$ sur $\mathfrak{g}_{s}$ est d\'etermin\'e par une orbite nilpotente  ${\cal O}$ (ou plus pr\'ecis\'ement ${\cal O}_{{\cal E}}$) et, pour tout $N\in {\cal O}$, par une repr\'esentation irr\'eductible $\xi_{N}$ (ou plus pr\'ecis\'ement $\xi_{N,{\cal E}}$) de $Z_{G_{s}}(N)/Z_{G_{s}}(N)^0$, cf. \ref{faisceauxcaracterescuspidaux}. D'apr\`es les r\'esultats de \cite{W7}, si ${\cal E}\in {\bf FC}_{s}^{st}$, $\xi_{N}$ est de dimension $1$, c'est-\`a-dire que c'est un caract\`ere. 
  
  {\bf Remarque.} Pour \^etre corrects, $\xi_{N}$ doit prendre ses valeurs dans $\bar{{\mathbb Q}}_{l}^{\times}$. Mais on a fixé un isomorphisme de $\bar{{\mathbb Q}}_{l}$ sur ${\mathbb C}$ et, pour toute la section, on considère que $\xi_{N}$ est à valeurs complexes.
  
  \bigskip
  
  Quand $N$ parcourt ${\cal O}$, les caractères $\xi_{N}$ v\'erifient  la condition de coh\'erence suivante:
  
  (1) $\xi_{g^{-1}Ng}(z)= \xi_{N}(gzg^{-1})$  pour tous $g\in G_{s}$ et $z\in Z_{G_{s}}(g^{-1}Ng)/Z_{G_{s}}(g^{-1}Ng)^0$. 
  
  Supposons que $s\in S^{st}(G)$ et que ${\cal E}$ soit conserv\'e par l'action de Frobenius.  Alors ${\cal O}$ est stable par l'action de $\Gamma_{{\mathbb F}_{q}}$. Munir ${\cal E}$ d'une action de Frobenius revient \`a se donner pour tout $N\in {\cal O}^{\Gamma_{{\mathbb F}_{q}}}$ un scalaire $r_{N}\in {\mathbb C}^{\times}$ (ou plus pr\'ecis\'ement $r_{N,{\cal E}}$), cf. \ref{fonctionscaracteristiques}. Ces scalaires v\'erifient la condition de coh\'erence
  
  (2) $r_{g^{-1}Ng}=r_{N}\xi_{N}(Fr(g)g^{-1})$ pour tout $g\in G_{s}$  tel que $g^{-1}Ng\in {\cal O}^{{\mathbb F}_{q}}$.

   Soient $s\in S^{nr,st}(G)$ et ${\cal E}\in {\bf FC}_{s}^{st}$.   Notons $Aut_{s}(G,F^{nr})$ le sous-groupe des éléments    $\delta\in Aut(G,F^{nr})$ dont  l'action  sur $Imm_{F^{nr}}(G)$ conserve $s$. Un élément $\delta\in Aut_{s}(G,F^{nr})$ se descend en un automorphisme de $G_{s}$ que l'on note $\delta_{s}$. 
  D'apr\`es \cite{W7} 9(5),  le faisceau ${\cal E}$  est conserv\'e par $\delta_{s}$.    
   Cela entra\^{\i}ne que, pour $N\in {\cal O}$, $\delta_{s}(N)$ appartient \`a ${\cal O}$ et  que l'on a l'\'egalit\'e
  
  (3) $ \xi_{\delta_{s}(N)}(\delta_{s}(z))=\xi_{N}(z)$ pour tout $z\in Z_{G_{s}}(N)/Z_{G_{s}}(N)^0$. 
  
  Pour $N\in \mathfrak{g}_{s}$, on note $K_{s,N}^{\dag,nr}$ le groupe des $k\in K_{s}^{\dag,nr}$ tels que $Ad(k)_{s}(N)=N$. On pose $K_{s,N}^{0,nr}=K_{s,N}^{\dag,nr}\cap K_{s}^{0,nr}$. 
On note $g\mapsto \bar{g}$ l'homomorphisme de r\'eduction de $K_{s}^{0,nr}$ dans $G_{s}$.

  Consid\'erons une famille $(\epsilon_{N})_{N\in {\cal O}}$ d'applications $\epsilon_{N}: K_{s,N}^{\dag,nr}\to {\mathbb C}^{\times}$.    
    Pour tous $N,N'\in {\cal O}$, introduisons les conditions  suivantes:
    
      $Hyp(1,N)$ $\epsilon_{N}$ est un caract\`ere de $ K_{s,N}^{\dag,nr}$;
  
  $Hyp(2,N)$    pour $k\in K_{s,N}^{0,nr}$, on a l'égalité $\epsilon_{N}(k)=\xi_{N}(\bar{k})$;

  $Hyp(3,N,N')$ $\epsilon_{N'}( h^{-1}kh)=\epsilon_{N}(k)$ pour tout $k\in K_{s,N}^{\dag,nr}$ et tout $h\in K_{s}^{\dag,nr}$ tel que $Ad(h)_{s}(N)=N'$.  
  
  Nous dirons que $(\epsilon_{N})_{N\in {\cal O}}$  est une famille admissible  si ces trois conditions sont v\'erifi\'ees pour tous $N,N'\in {\cal O}$. 
  
  Supposons que $s\in S^{st}(G)$ et que ${\cal E}$ soit conserv\'e par l'action galoisienne.
  Pour $N\in {\cal O}^{\Gamma_{{\mathbb F}_{q}}}$, introduisons l'hypoth\`ese
  
  $Hyp(4,N)$  on a l'\'egalit\'e $\epsilon_{N}=\epsilon_{N}\circ Fr$.
  
  Nous dirons que $(\epsilon_{N})_{N\in {\cal O}}$  est une famille admissible   d\'efinie sur $F$ si c'est une famille admissible  et que $Hyp(4,N)$ est v\'erifi\'ee pour tout $N\in {\cal O}^{\Gamma_{{\mathbb F}_{q}}}$. 
  
  Les deux propriétés ci-dessous sont immédiates:
  
  (4) soit $N\in {\cal O}$ et soit $\epsilon_{N}$ une application de $K_{s,N}^{\dag,nr} $ dans ${\mathbb C}^{\times}$; supposons v\'erifi\'ees les conditions $Hyp(1,N)$ et $Hyp(2,N)$; alors $\epsilon_{N}$ est la composante en $N$ d'une unique famille admissible $(\epsilon_{N'})_{N'\in {\cal O}}$; 
  
  (5) supposons que $s\in S^{st}(G)$ et que ${\cal E}$ soit conserv\'e par l'action galoisienne;  soit $N\in {\cal O}^{\Gamma_{{\mathbb F}_{q}}}$ et soit $\epsilon_{N}$ une application de $ K_{s,N}^{\dag,nr}$ dans ${\mathbb C}^{\times}$; supposons v\'erifi\'ees les conditions $Hyp(1,N)$, $Hyp(2,N)$ et $Hyp(4,N)$;  alors $\epsilon_{N}$ est la composante en $N$ d'une unique famille admissible $(\epsilon_{N'})_{N'\in {\cal O}}$  d\'efinie sur $F$.

       Faisons maintenant varier 
  le sommet $s$ et le faisceau ${\cal E}\in {\bf FC}_{s}^{st}$.   On ajoute si besoin est des indices $s$ et ${\cal E}$ aux notations. Pour tout sommet $s\in  S^{nr,st}(G)$  et pour tout faisceau-caract\`ere ${\cal E}\in {\bf FC}^{st}_{s}$,  supposons donn\'ee une famille  admissible $(\epsilon_{s,{\cal E},N})_{N\in {\cal O}_{{\cal E}}}$.

 Soit $g\in G(F^{nr})$. L'action de $g$ sur $Imm_{F^{nr}}(G)$ conserve  l'ensemble $S^{nr,st}(G)$.  Fixons $s\in S^{nr,st}(G)$, notons $s'$ son image par cette action .  L'application $Ad(g)$ se descend en un isomorphisme $\overline{Ad(g)}$   de $G_{s}$ sur $G_{s'}$. Cet isomorphisme transporte un faisceau-caract\`ere ${\cal E}\in {\bf FC}_{s}$   en un faisceau-caract\`ere ${\cal E}'=\overline{Ad(g)}({\cal E})\in {\bf FC}_{s'}$.  On a $\overline{Ad(g)}({\cal O}_{{\cal E}})={\cal O}_{{\cal E}'}$. Le premier faisceau appartient \`a ${\bf FC}^{st}_{s}$    si et seulement si le second  appartient \`a ${\bf FC}^{st}_{s'}$. Supposons que ces deux faisceaux v\'erifient ces conditions.   Soit $N\in {\cal O}_{{\cal E}}$, posons $N'=\overline{Ad(g)}(N)$. On a l'\'egalit\'e $ K_{s',N'}^{\dag,nr}=\{gkg^{-1}; k\in K_{s,N}^{\dag,nr}\}$.  Consid\'erons la condition
 
$Hyp(5)$ sous les hypoth\`eses ci-dessus, on a l'\'egalit\'e $\epsilon_{s',{\cal E}',N'}( gkg^{-1})=\epsilon_{s,{\cal E},N}(k)$ pour tout $ k\in K_{s,N}^{\dag,nr}$. 

Nous dirons que la famille  $(\epsilon_{s,{\cal E},N})_{s\in S^{nr,st}(G),{\cal E}\in {\bf FC}^{st}_{s},N\in {\cal O}_{{\cal E}}}$ est une $G$-famille admissible si cette condition est v\'erifi\'ee pour tous $s$, ${\cal E}$, $N$ et $g$ comme ci-dessus. 

On dit que $(\epsilon_{s,{\cal E},N})_{s\in S^{nr,st}(G),{\cal E}\in {\bf FC}^{st}_{s},N\in {\cal O}_{{\cal E}}}$ est une  $G$-famille admissible d\'efinie sur $F$ si c'est une $G$-famille admissible et que, pour tous 
 $s\in S^{st}(G)$ et ${\cal E}\in {\bf FC}_{s,{\mathbb F}_{q}}^{st}$, la famille $(\epsilon_{s,{\cal E},N})_{N\in {\cal O}_{{\cal E}}}$ est d\'efinie sur $F$. 
 
{\bf Remarque.} (6) Supposons que les familles $(\epsilon_{s,{\cal E},N})_{N\in {\cal O}_{{\cal E}}}$ soient canoniques en ce sens que leurs d\'efinitions ne d\'ependent d'aucun objet auxiliaire (ce sera le cas pour celles que nous d\'efinirons). Alors $Hyp(5)$ est v\'erifi\'ee par simple transport de structure. 

\subsubsection{R\'esolution du probl\`eme\label{resolution}}
 On suppose fix\'ee une $G$-famille admissible $(\epsilon_{s,{\cal E},N})_{s\in S^{nr,st}(G),{\cal E}\in {\bf FC}^{st}_{s},N\in {\cal O}_{{\cal E}}}$ d\'efinie sur $F$.

Consid\'erons un groupe $H$  et un torseur intérieur $\phi:H\to G$ comme en \ref{position}.      Fixons un ensemble de représentants $\underline{S}^{st}(H)$ des orbites de l'action de $H(F)$ dans $S^{st}(H)$. D'après le lemme \ref{decompositionFCstable}, on a l'égalité
$$(1) \qquad FC^{st}(\mathfrak{h}(F))= \oplus_{s_{H}\in \underline{S}^{st}(H)}FC^{st}(\mathfrak{h}_{s_{H}}({\mathbb F}_{q})).$$ 

{\bf Remarque.} Il s'agit en v\'erit\'e d'un isomorphisme. L'espace de droite est un sous-espace de $C_{c}^{\infty}(\mathfrak{g}(F))$ tandis que celui de gauche est un sous-espace de $I(\mathfrak{g}(F))$. La projection naturelle $C_{c}^{\infty}(\mathfrak{g}(F))\to I(\mathfrak{g}(F))$ envoie bijectivement le membre de droite sur celui de gauche.
\bigskip

On pose des définitions analogues pour le groupe $G$.  Puisque l'homomorphisme d'inflation $H^1(\Gamma_{F}^{nr},G(F^{nr}))\to H^1(\Gamma_{F},G)$ est un isomorphisme, on peut supposer que $\phi$ est un isomorphisme sur $F^{nr}$, quitte à le remplacer par $Ad(y) \phi$ pour un $y\in G$ convenable.  
  Notons 
$\phi^{Imm}: Imm_{F^{nr}}(H)\to  Imm_{F^{nr}}(G)$ l'isomorphisme d'immeubles qui se d\'eduit de $\phi$.  D'apr\`es les r\'esultats de \cite{W7} paragraphe 9, $\phi^{Imm}$ envoie bijectivement $S^{nr,st}(H)$ sur $S^{nr,st}(G)$ et, pour $s_{H}\in S^{st}(H)$, $\phi^{Imm}(s_{H})$ est conjugué à un élément de $S^{st}(G)$ par un élément de $G(F^{nr})$. D'autre part deux éléments de $S^{st}(G)$ ne sont conjugués par un élément de $G(F^{nr})$ que s'ils sont conjugués par un élément de $G(F)$, cf. lemme \ref{lemmeautomorphismes}. On en déduit une bijection $\underline{\phi}:\underline{S}^{st}(H)\to \underline{S}^{st}(G)$: pour $s_{H}\in \underline{S}^{st}(H)$, $\underline{\phi}(s_{H})$ est l'unique élément de $\underline{S}^{st}(G)$ qui est conjugué par un élément de $G(F^{nr})$ à $\phi^{Imm}(s_{H})$. Fixons $s_{H}\in \underline{S}^{st}(H)$, posons $s=\underline{\phi}(s_{H})$. Nous allons définir un isomorphisme
 $$\iota_{s}:FC^{st}(\mathfrak{h}_{s_{H}}({\mathbb F}_{q}))\to FC^{st}(\mathfrak{g}_{s}({\mathbb F}_{q})).$$
En tenant compte de l'égalité (1) et de son analogue pour le groupe $G$, la somme des $\iota_{s}$ définira l'isomorphisme $\iota_{\phi}$ de \ref{position}. 

  Les sommets $s_{H}$ et $s$ étant maintenant fixés, quitte  \`a remplacer encore $\phi$ par $Ad(y) \phi$ pour un $y\in G(F^{nr})$ convenable, on peut supposer que  $\phi^{Imm}(s_{H})=s$.    L'isomorphisme $\phi$ se descend en un isomorphisme $\phi_{s}: H_{s_{H}}\to G_{s}$.   L'isomorphisme $\phi_{s}^{-1}$ transporte un faisceau ${\cal E}\in {\bf FC}^{st}_{s,{\mathbb F}_{q}}$ en un faisceau ${\cal E}_{H}$ sur $\mathfrak{h}_{s_{H}}$.  On pose ${\cal O}={\cal O}_{{\cal E}}$ etc... et ${\cal O}_{H}={\cal O}_{{\cal E}_{H}}$ etc...  On a  ${\cal O}_{H}=\phi_{s}^{-1}({\cal O})$. Pour $N_{H}\in {\cal O}_{H}$, on a l'\'egalit\'e $\xi_{N_{H}}(z)=\xi_{\phi_{s}(N_{H})}(\phi_{s}(z))$ pour tout $z\in Z_{H_{s_{H}}}(N_{H})/Z_{H_{s_{H}}}(N_{H})^0$.   Puisque $\phi$ est un torseur intérieur  et un isomorphisme sur $F^{nr}$ et que $G$ est adjoint, il existe un unique $k_{\phi}\in G(F^{nr})$ tel que $Fr(\phi)\phi^{-1}=Ad(k_{\phi})$. Puisque  $\phi^{Imm}(s_{H})=s$  et que les deux sommets $s_{H}$ et $s$ sont fixes par l'action galoisienne, on a $k_{\phi}\in K_{s}^{\dag,nr}$. Le faisceau ${\cal E}$ est conservé par $Ad(k_{\phi})_{s}=Fr(\phi_{s})\phi_{s}^{-1}$. Cela \'equivaut \`a:  le faisceau ${\cal E}_{H}$ est conserv\'e par l'action de Frobenius.
Il r\'esulte des descriptions de \cite{W7} que ${\cal E}_{H}$ appartient \`a ${\bf FC}^{st}_{s_{H},{\mathbb F}_{q}}$ et que l'application ainsi d\'efinie de ${\bf FC}^{st}_{s,{\mathbb F}_{q}}$ sur ${\bf FC}^{st}_{s_{H},{\mathbb F}_{q}}$ est bijective. 

Fixons   un \'el\'ement $N\in {\cal O}^{\Gamma_{{\mathbb F}_{q}}}$. 
Montrons que

(2) quitte à remplacer $\phi$ par $Ad(y)\circ\phi$ pour un $y\in K_{s}^{0,nr}$ convenable, on peut supposer $k_{\phi}\in K_{s,N}^{\dag,nr}$. 

Le faisceau ${\cal E}$ est conservé par $Ad(k_{\phi})_{s}$. Il  existe donc $g\in G_{s}$ tel que $ Ad(k_{\phi})_{s}(N)=gNg^{-1}$. Munissons un instant le groupe $G_{s}$ de l'action de Frobenius $x\mapsto Ad(k_{\phi})_{s}^{-1}\circ Fr(x)$. On applique le th\'eor\`eme de Lang: il existe $x\in G_{s}$ tel que $x Ad(k_{\phi})_{s}^{-1}\circ Fr(x)^{-1}=Ad(k_{\phi})_{s}^{-1}(g)$. Relevons $x$ en  un \'el\'ement $k\in K_{s}^{0,nr}$ et rempla\c{c}ons $\phi$ par $Ad(k) \phi$.    On voit que, maintenant, on a $Ad(k_{\phi})_{s}(N)=N$, donc  $k_{\phi}\in K_{s,N}^{\dag,nr}$. Cela prouve (2). 

Supposons donc $k_{\phi}\in K_{s,N}^{\dag,nr}$.  Munissons 
   ${\cal E}$ d'une action de Frobenius,  qui est d\'etermin\'ee par un scalaire  $r_{N}\in {\mathbb C}^{\times}$. Posons comme ci-dessus $N_{H}=\phi_{s}^{-1}(N)$. Parce que $Ad(k_{\phi})_{s}=Fr(\phi_{s})\phi_{s}^{-1}$ fixe $N$, on voit que $N_{H}\in {\cal O}_{H}^{\Gamma_{{\mathbb F}_{q}}}$. On d\'efinit une action de Frobenius sur ${\cal E}_{H}$ \`a l'aide du scalaire $r_{N_{H}}=r_{N}\epsilon_{s,{\cal E},N}(k_{\phi})$.
 
 Notons $f_{{\cal E}}$ la fonction caract\'eristique de ${\cal E}$ et $f_{{\cal E}_{H}}$ celle de ${\cal E}_{H}$. Faisons maintenant varier  le faisceau ${\cal E}$. La famille $(f_{{\cal E}})_{{\cal E}\in {\bf FC}^{st}_{s,{\mathbb F}_{q}}}$ est une base de $FC^{st}(\mathfrak{g}_{s}({\mathbb F}_{q}))$ et la famille $(f_{{\cal E}_{H}})_{{\cal E}\in {\bf FC}^{st}_{s,{\mathbb F}_{q}}}$ est une base de $FC^{st}(\mathfrak{h}_{s_{H}}({\mathbb F}_{q}))$. On d\'efinit l'isomorphisme  $\iota_{s}$ par l'\'egalit\'e
$$(3) \qquad \iota_{s}(f_{{\cal E}_{H}})=\vert G_{s}({\mathbb F}_{q}) \vert ^{-1}\vert H_{s_{H}}({\mathbb F}_{q})\vert  f_{{\cal E}}$$
pour tout ${\cal E}\in {\bf FC}^{st}_{s,{\mathbb F}_{q}}$.

\begin{prop}{L'isomorphisme $\iota_{\phi}$ ainsi défini ne d\'epend pas des choix effectu\'es.
}\end{prop}

Preuve. On a choisi les  ensembles $\underline{S}^{st}(H)$ et $\underline{S}^{st}(G)$, puis, les sommets $s$ et $s_{H}$ et un faisceau ${\cal E}$ \'etant fix\'e, on a choisi un \'el\'ement $N$,  un scalaire $r_{N}$ et on a remplacé $\phi$ par $Ad(y)\circ\phi$ pour un $y\in G$ convenable. Le  choix de $r_{N}$ n'influe pas: en vertu de la d\'efinition  $r_{N_{H}}=r_{N}\epsilon_{s,{\cal E},N}(\delta)$, changer  $r_{N}$ multiplie les deux fonctions $f_{{\cal E}}$ et $f_{{\cal E}_{H}}$ par le m\^eme scalaire et ne modifie pas $\iota_{\phi}$. 

Supposons fix\'es  $\underline{S}^{st}(H)$, $\underline{S}^{st}(G)$, $ s$, $s_{H}$ et  ${\cal E}$.   Montrons que la d\'efinition de l'action du Frobenius sur ${\cal E}_{H}$ ne d\'epend pas des choix de $N$ et $y$. Choisissons un autre élément $N'\in {\cal O}^{\Gamma_{{\mathbb F}_{q}}}$ et un autre élément $y'\in G$ tel que $Ad(y')\phi$ vérifie les conditions requises relativement à $N'$. Pour simplifier les notations, on  peut supposer que $y=1$ et on pose   $\phi'=Ad(y')\phi$.  On déduit des nouvelles données un faisceau ${\cal E}'_{H}={\phi'_{s}}^{-1}({\cal E})$ sur $\mathfrak{h}_{s_{H}}$, muni d'une action de Frobenius. On affecte d'un $'$ les objets relatifs à ces nouvelles données. Puisque $\phi$ et $\phi'$ sont des isomorphismes sur $F^{nr}$, on a $y'\in G(F^{nr})$. 
Puisque $\phi^{Imm}(s_{H})=s={\phi'}^{Imm}(s_{H})$ et que ${\phi'}^{Imm}$ est le composé de l'action de $y'$ sur $Imm_{F^{nr}}(G)$ et de $\phi^{Imm}$, on a $y'\in K_{s}^{\dag,nr}$. Puisque $Ad(y')_{s}$ conserve ${\cal E}$, on a l'égalité ${\cal E}'_{H}={\cal E}_{H}$. Posons $N'_{H}={\phi'}_{s}^{-1}(N')$. On doit prouver l'égalité

(4) $r'_{N'_{H}}=r_{N'_{H}}$. 

On simplifie les notations en supprimant les termes $s$ et ${\cal E}$ de notre famille admissible $(\epsilon_{s,{\cal E},N''})_{N''\in {\cal O}}$. 
Par définition, on a $r'_{N'_{H}}=\epsilon_{N'}(k_{\phi'})r_{N'}$. Le terme $k_{\phi'}$ est défini par l'égalité $Ad(k_{\phi'})=Fr(\phi'){\phi'}^{-1}=Ad(Fr(y'))Fr(\phi)\phi^{-1}Ad(y')^{-1}=Ad(Fr(y')k_{\phi}{y'}^{-1})$, d'où $k_{\phi'}=Fr(y')k_{\phi}{y'}^{-1}$. Puisque $N'\in {\cal O}$, on peut fixer $g\in K_{s}^{0,nr}$ tel que $N'=\bar{g}^{-1}N\bar{g}$. D'après $Hyp(3,N,N')$, on. a $\epsilon_{N'}(k_{\phi'})=\epsilon_{N}(gk_{\phi'}g^{-1})$. D'après \ref{familles}(2)  et $Hyp(2,N)$, on a aussi $r_{N'}=r_{N}\xi_{N}(Fr(\bar{g})\bar{g}^{-1})=r_{N}\epsilon_{N}(Fr(g)g^{-1})$. En utilisant $Hyp(1,N)$, on obtient  l'égalité

(5)  $r'_{N'_{H}}=r_{N}\epsilon_{N}(Fr(gy')k_{\phi}(gy')^{-1})$.

Puisque $Ad(y')_{s}$ conserve ${\cal E}$, on peut fixer $u\in K_{s}^{nr,0}$ tel que $Ad(y')_{s}^{-1}(N')=\bar{u}^{-1}N\bar{u}$. Posons $h=\phi^{-1}(u)$. On a $h\in K_{s_{H}}^{H,0,nr}$ et les égalités $N'_{H}={\phi'}_{s}^{-1}(N')=\phi_{s}^{-1}Ad(y')_{s}^{-1}(N')=\phi_{s}^{-1}Ad(u)_{s}^{-1}(N)=Ad(h)_{s_{H}}^{-1}\phi_{s}^{-1}(N)=\bar{h}^{-1}N_{H}\bar{h}$, où $\bar{h}$ est la réduction de $h$ dans $H_{s_{H}}$.  D'après  \ref{familles}(2), on a $r_{N'_{H}}=r_{N_{H}}\xi_{N_{H}}(Fr(\bar{h})\bar{h}^{-1})$, d'où, en appliquant les définitions de $r_{N_{H}}$ et $\xi_{N_{H}}$, $r_{N'_{H}}=r_{N}\epsilon_{N}(k_{\phi})\xi_{N}\circ\phi_{s}(Fr(\bar{h})\bar{h}^{-1})$. Le terme $\phi_{s}(Fr(\bar{h})\bar{h}^{-1})$ est la réduction dans $G_{s}$ de $\phi(Fr(h)h^{-1})$, d'où $\xi_{N}\circ\phi_{s}(Fr(\bar{h})\bar{h}^{-1})=\epsilon_{N}(\phi(Fr(h)h^{-1}))$. 
On a $Fr(h)=Fr(\phi^{-1}(u))=Fr(\phi)^{-1}(Fr(u))$. Par définition de $k_{\phi}$, $Fr(\phi)^{-1}=\phi^{-1}Ad(k_{\phi})^{-1}$. D'où l'égalité $\phi(Fr(h))=k_{\phi}^{-1}Fr(u)k_{\phi}$, puis $\phi(Fr(h)h^{-1})=k_{\phi}^{-1}Fr(u)k_{\phi}u^{-1}$. D'où $r_{N'_{H}}=r_{N}\epsilon_{N}(k_{\phi})\epsilon_{N}(k_{\phi}^{-1}Fr(u)k_{\phi}u^{-1})$, puis, en
 utilisant $Hyp(1,N)$, 

(6) $ r_{N'_{H}}=r_{N}\epsilon_{N}(Fr(u)k_{\phi}u^{-1})$.

On a les égalités $Ad(y')_{s}^{-1}(N')=\bar{u}^{-1}N\bar{u}$ et $N'=\bar{g}^{-1}N\bar{g}$. Posons $v=gy'u^{-1}$. On a $v\in K_{s}^{\dag,nr}$ et les égalités précédentes impliquent que $Ad(v)_{s}$ fixe $N$, autrement dit $v\in K_{s,N}^{\dag,nr}$. L'égalité (6) se récrit $r_{N'_{H}}=r_{N}\epsilon_{N}(Fr(v)^{-1}Fr(gy')k_{\phi}(gy')^{-1}v^{-1})$. Parce que $v\in K_{s,N}^{\dag,nr}$, en utilisant $Hyp(1,N)$ et $Hyp(4,N)$, cela se simplifie en  $r_{N'_{H}}=r_{N}\epsilon_{N}(Fr(gy')k_{\phi}(gy')^{-1})$. En comparant avec (5), on obtient (4).

Il nous reste maintenant \`a faire varier les  ensembles  $\underline{S}^{st}(H)$ et $\underline{S}^{st}(G)$.  Un changement de ces ensembles remplace le couple $(s_{H},s)$ de la construction ci-dessus par un couple $(\underline{s}_{H},\underline{s})=(h^{-1}s_{H},gs)$, où  $h\in H(F)$ et $g\in G(F)$. 
    On a un diagramme
$$ \qquad \begin{array}{ccccccc}&&FC^{st}(\mathfrak{h}_{\underline{s}_{H}}({\mathbb F}_{q}))&\stackrel{\iota_{\underline{s}}}{\to}&FC^{st}(\mathfrak{g}_{\underline{s}}({\mathbb F}_{q}))&&\\ &\swarrow&&&&\searrow&\\ I(\mathfrak{h}(F))&& \downarrow Ad(h)
&&\uparrow Ad(g)&&I(\mathfrak{g}(F))\\&\nwarrow&&&&\nearrow&\\&&  FC^{st}(\mathfrak{h}_{s_{H}}({\mathbb F}_{q}))&\stackrel{\iota_{s}}{\to}&FC^{st}(\mathfrak{g}_{s}({\mathbb F}_{q}))&&\\ \end{array}$$
    Les deux triangles sont commutatifs.  Pour d\'emontrer que le changement de sommets ne modifie pas $\iota$, il suffit donc de d\'emontrer que le  carr\'e int\'erieur est commutatif. Soit ${\cal E}\in {\bf FC}^{st}_{s,{\mathbb F}_{q}}$, avec les notations aff\'erentes. Fixons $N\in {\cal O}^{\Gamma_{{\mathbb F}_{q}}}$.  De $Ad(g)$ se d\'eduit l'isomorphisme $\overline{Ad(g)}:G_{s}\to G_{\underline{s}}$. Posons $\underline{{\cal E}}=\overline{Ad(g)}({\cal E})$, $\underline{N}=\overline{Ad(g)}(N)$. On a $\underline{{\cal E}}\in {\bf FC}^{st}_{\underline{s},{\mathbb F}_{q}}$. Puisque $g\in G(F)$, on a encore $\underline{N}\in \underline{{\cal O}}^{\Gamma_{{\mathbb F}_{q}}}$ (o\`u \'evidemment $\underline{{\cal O}}$ est l'orbite supportant $\underline{{\cal E}}$).  L'action $\overline{Ad(g)}$ transporte l'action fix\'ee de Frobenius sur ${\cal E}$ en l'action sur $\underline{{\cal E}}$ d\'efinie par le scalaire $\underline{r}_{\underline{N}}=r_{N}$. Quitte à  remplacer $\phi$ par $Ad(y)\phi$ pour un $y\in G$ convenable, on suppose que $\phi$ vérifie les conditions imposées dans la construction de $\iota_{s}$ relativement à $s$, $s_{H}$ et $N$. Posons 
   $\underline{\phi}=Ad(g)\phi Ad(h)^{-1}=Ad(g\phi(h)^{-1})\phi$. On voit que $\underline{\phi}$ vérifie les conditions imposées relativement à $\underline{s}$, $\underline{s}_{H}$ et $\underline{N}$.  Appliquons les constructions de $\iota_{s}$  et $\iota_{\underline{s}}$.  On en d\'eduit  des faisceaux ${\cal E}_{H}\in {\bf FC}_{s_{H},{\mathbb F}_{q}}^{st}$ et $\underline{{\cal E}}_{H}\in {\bf FC}_{\underline{s}_{H},{\mathbb F}_{q}}^{st}$ et des \'el\'ements $N_{H}\in {\cal O}_{H}^{\Gamma_{{\mathbb F}_{q}}}$ et $\underline{N}_{H}\in \underline{{\cal O}}_{H}^{\Gamma_{{\mathbb F}_{q}}}$. Il est clair que l'isomorphisme $\overline{Ad(h)}:H_{s_{H}}\to H_{\underline{s}_{H}}$ transporte ${\cal E}_{H}$ et $N_{H}$ sur $\underline{{\cal E}}_{H}$ et $\underline{N}_{H}$.  Il suffit de prouver que cet isomorphisme entrelace les actions de Frobenius. En reprenant les d\'efinitions, cela  \'equivaut \`a l'\'egalit\'e $\epsilon_{\underline{s},\underline{{\cal E}},\underline{N}}(k_{\underline{\phi}})=\epsilon_{s,{\cal E},N}(k_{\phi})$. Par définition, $Ad(k_{\underline{\phi}})=Fr(\underline{\phi})\underline{\phi}^{-1}=Ad(g)Fr(\phi)Ad(h)Ad(h)^{-1}\phi^{-1}Ad(g)^{-1}$ puisque $g\in G(F)$ et $h\in H(F)$, d'où $Ad(k_{\underline{\phi}})=Ad(g)Ad(k_{\phi})Ad(g)^{-1}$, puis $k_{\underline{\phi}}=gk_{\phi}g^{-1}$. L'égalité à démontrer est donc $\epsilon_{\underline{s},\underline{{\cal E}},\underline{N}}(gk_{\phi}g^{-1})=\epsilon_{s,{\cal E},N}(k_{\phi})$. Elle résulte de $Hyp(5)$.  Cela ach\`eve la  démonstration.  $\square$

{\bf Remarque}(7).  Supposons que $\phi$ soit équivalent à un isomorphisme défini sur $F$, c'est-à-dire qu'il existe $y\in G$ tel que $Ad(y)\phi$ soit équivariant pour les actions de $\Gamma_{F}$.   En remplaçant $\phi$ par un tel $Ad(y)\phi$,  on obtient un isomorphisme $FC^{st}(\mathfrak{h}(F))\to FC^{st}(\mathfrak{g}(F))$ par transport de structure.  Celui-ci ne dépend pas du choix de $y$: on ne peut changer $y$ qu'en le multipliant à gauche par un élément $g\in G(F)$ et l'action de $Ad(g) $ sur  $FC^{st}(\mathfrak{g}(F))$ est triviale.   Il coïncide avec l'isomorphisme $\iota_{\phi}$ construit ci-dessus. En effet, partant de données $s$, ${\cal E}$, $N$, on peut supposer que $s_{H}={\phi^{Imm}}^{-1}(s)$ et on voit que $\phi$ vérifie les conditions requises relativement à $s$, $s_{H}$ et $N$ avec $k_{\phi}=1$. 

\bigskip

Nous devons donc   construire une $G$-famille admissible
$(\epsilon_{s,{\cal E},N})_{s\in S^{nr,st}(G), {\cal E}\in {\bf FC}_{s}^{st}, N\in {\cal O}_{{\cal F}}}$ d\'efinie sur $F$. Nous le ferons dans les deux paragraphes suivants. Dans les paragraphes ultérieurs, on supposera que l'isomorphisme $\iota_{\phi}$ est défini à l'aide de cette famille. 

\subsubsection{Un sous-ensemble de ${\cal O}_{{\cal E}}$}\label{unsousensemble}
  
  On veut construire une  $G$-famille admissible $(\epsilon_{s,{\cal E},N})_{s\in S^{nr,st}(G), {\cal E}\in {\bf FC}_{s}^{st}, N\in {\cal O}_{{\cal F}}}$ d\'efinie sur $F$. Il nous suffit de fixer $s\in S^{nr,st}(G)$ et ${\cal E}\in {\bf FC}_{s}^{st}$ et de définir une famille admissible  $(\epsilon_{s,{\cal E},N})_ { N\in {\cal O}_{{\cal F}}}$ "canonique" au sens qu'elle ne dépende pas de choix auxiliaires.   Comme on l'a dit dans la remarque (6) de \ref{familles}, l'hypoth\`ese $Hyp(5)$ sera alors automatiquement v\'erifi\'ee.  On fixe donc $s$ et ${\cal E}$ et on simplifie la notation: on va construire une famille admissible   $(\epsilon_{N})_{N\in {\cal O}}$ où ${\cal O}={\cal O}_{{\cal E}}$. Seule la condition  $Hyp(4,N)$ d\'epend de la structure de $G$ sur $F$. 
    Nous supposerons donc que le corps de base est $F^{nr}$, en ne revenant \`a la situation sur $F$ que pour traiter cette condition $Hyp(4,N)$. Bien s\^ur,  pour traiter cette derni\`ere condition, il faut imposer les hypoth\`eses suppl\'ementaires que $s\in S^{st}(G)$ et que ${\cal E}\in {\bf FC}^{st}_{s,{\mathbb F}_{q}}$.

   On reprend les constructions de \ref{alcove} et \ref{racinesaffines}. On fixe donc une paire de Borel \'epingl\'ee $\mathfrak{E}=(B,T,(E_{\alpha})_{\alpha\in \Delta})$  et  un \'epinglage affine $\mathfrak{E}_{a}=(T,(E_{\alpha})_{\alpha\in \Delta_{a}})$  tous deux  fix\'es par $I_{F}$. On a d\'eduit de $\mathfrak{E}$ une alc\^ove $C^{nr}\subset Imm_{F^{nr}}(G)$ et on suppose que $s\in S(\bar{C}^{nr})$ (le sommet $s$ \'etant fix\'e, on peut toujours construire des donn\'ees telles que cette condition soit v\'erifi\'ee). Le sommet $s$ est associ\'e \`a une racine $\beta_{s}\in \Delta_{a}^{nr}$. 
   On introduit le groupe $\boldsymbol{\Omega}^{nr}$. Montrons que
   
   (1) $K_{s}^{\dag,nr}$ est égal au produit semi-direct $K_{s}^{0,nr}\rtimes \boldsymbol{\Omega}^{nr}$.
   
  Parce que $s\in S^{nr,st}(G)$, $\beta_{s}$ est fixé par l'action de $\boldsymbol{\Omega}^{nr}$ sur ${\cal D}_{a}^{nr}$. Donc $\boldsymbol{\Omega}^{nr}\subset K_{s}^{\dag,nr}$. Soit $g\in K_{s}^{\dag,nr}$. D'après le corollaire \ref{uncorollaire}, on peut écrire $g=\omega t\pi(g_{sc})$, avec $\omega\in \boldsymbol{\Omega}^{nr}$, $t\in T(\mathfrak{o}_{F^{nr}})$ et $g_{sc}\in G_{SC}(F^{nr})$. On a $\omega\in K_{s}^{\dag,nr}$ et $ T(\mathfrak{o}_{F^{nr}})\subset K_{s}^{0,nr}$ (cf.  \ref{normalisateur}(1)). Donc aussi $\pi(g_{sc})\in K_{s}^{\dag,nr}$, autrement dit l'action de $g_{sc}$ sur $Imm_{F^{nr}}(G)$ fixe $s$. Puisque  $G_{SC}$ est simplement connexe, cela entraîne que $g_{sc} \in K_{SC, s}^{0,nr}$, d'où $\pi(g_{sc})\in K_{s}^{0,nr}$. Cela prouve que $K_{s}^{\dag,nr}$ est bien le produit de $K_{s}^{0,nr}$ et de $\boldsymbol{\Omega}^{nr}$. Prouvons que l'intersection de ces groupes est triviale. Soit $k\in K_{s}^{0,nr}\cap \boldsymbol{\Omega}^{nr}$. Alors $k$ conserve $C^{nr}$. Un élément de $K_{s}^{0,nr}$ qui conserve $C^{nr}$ appartient à $K_{C^{nr}}^{0,nr}$. Un tel élément agit trivialement sur ${\cal D}_{a}^{nr}$. L'action de  $\boldsymbol{\Omega}^{nr}$ sur ${\cal D}_{a}^{nr}$ est fidèle, donc $k=1$. Cela démontre (1). 
  
  Tout élément $ k\in K_{s}^{\dag,nr}$ définit un automorphisme $Ad(k)_{s}$ de $G_{s}$. Notons $R$ la composante neutre du sous-groupe des éléments de $G_{s}$ fixés par $Ad(\omega)_{s}$ pour tout $\omega\in \boldsymbol{\Omega}^{nr}$. 
  
  \begin{lem}{L'ensemble ${\cal O}\cap \mathfrak{r}$ est une unique orbite sous l'action de $R$.}\end{lem}
  
  Preuve.  On rappelle que l'on a exclu plusieurs cas en \ref{position}. On peut aussi exclure ici les cas où $\boldsymbol{\Omega}^{nr}=\{1\}$, le lemme étant trivial dans ce cas. 
On traite tous les cas restants, en se r\'ef\'erant \`a \cite{W7} o\`u on a d\'ecrit le sommet $s$ et l'orbite ${\cal O}$ associ\'ee au  faisceau ${\cal F}$ (l'ensemble ${\bf FC}^{st}_{s}$ peut avoir plusieurs \'el\'ements mais l'orbite ${\cal O}$ est commune \`a tous ceux-ci).   

Supposons $G$ de type $B_{n}$ avec  $2n+1=k^2+h^2$, o\`u $k,h\in {\mathbb N}$, $k$ est pair et $ \vert k-h\vert =1$.  Alors $G_{s,SC}=Spin(k^2)\times Spin(h^2)$. 
 Le groupe $ \boldsymbol{\Omega}^{nr}$ s'identifie à $Aut({\cal D}_{a}) $. Il  a deux \'el\'ements, l'\'el\'ement non trivial agissant  sur $G_{s,SC}$ par un automorphisme ext\'erieur non trivial de la premi\`ere composante $Spin(k^2)$. 
 On a ${\cal O}={\cal O}_{1}\times {\cal O}_{2}$, o\`u ${\cal O}_{1}$, resp. ${\cal O}_{2}$, sont param\'etr\'ees par les partitions $(2k-1,...,3,1)$, resp. $(2h-1,...,3,1)$. On a $R_{SC}=Spin(k^2-1)\times Spin(h^2)$, l'homomorphisme de $R_{SC}$ dans $G_{s,SC}$ \'etant \'evident.  Consid\'erons un \'el\'ement nilpotent  de $\mathfrak{spin}(k^2-1)$ param\'etr\'e par une partition $(\lambda_{1},...,\lambda_{l})$, avec $\lambda_{l}>0$. Alors son image dans $\mathfrak{spin}(k^2)$ est param\'etr\'ee par la partition $(\lambda_{1},..,\lambda_{l},1)$. Donc  ${\cal O}\cap \mathfrak{r}$ est \'egal \`a ${\cal O}'_{1}\times {\cal O}_{2}$, o\`u ${\cal O}'_{1}$ est param\'etr\'ee par la partition $(2k-1,...,5,3)$. 
 
 Supposons $G$ de type $C_{n}$ o\`u   $n=k(k+1)$ avec $k\in {\mathbb N}$. Alors $G_{s,SC}=Sp(k(k+1))\times Sp(k(k+1))$. Le groupe $ \boldsymbol{\Omega}^{nr}$ s'identifie à $Aut({\cal D}_{a}) $. Il  a deux \'el\'ements,  l'\'el\'ement non trivial agissant sur $G_{s,SC}$ par permutation des facteurs. On a ${\cal O}={\cal O}_{1}\times {\cal O}_{2}$, o\`u ${\cal O}_{1}$, resp. ${\cal O}_{2}$, sont param\'etr\'ees par la m\^eme partition $(2k,...,4,2)$. On a $R_{SC}=Spin(k(k+1))$,  l'homomorphisme de $R_{SC}$ dans $G_{s,SC}$ \'etant  diagonal. Alors ${\cal O}\cap \mathfrak{r}$ est \'egal \`a  l'orbite param\'etr\'ee par cette m\^eme partition $(2k,...,4,2)$.
 
 Supposons $G$ de type $(D_{n},nr)$ avec $n=h^2$ o\`u $h\in {\mathbb N}$ est pair.  Alors $G_{s,SC}=Spin(h^2)\times Spin(h^2)$.Le groupe $ \boldsymbol{\Omega}^{nr}$ a $4$ \'el\'ements. Son action sur $G_{s,SC}$ est engendr\'ee par le produit des automorphismes ext\'erieurs habituels de chacun des facteurs $Spin(h^2)$ et par une permutation des deux facteurs.  On a ${\cal O}={\cal O}_{1}\times {\cal O}_{2}$, o\`u ${\cal O}_{1}$, resp. ${\cal O}_{2}$, sont param\'etr\'ees par la m\^eme partition $(2h-1,...,3,1)$. On a $R_{SC}=Spin(h^2-1)$, l'homomorphisme de $R_{SC}$ dans $G_{s,SC}$ \'etant compos\'e du plongement \'evident $Spin(h^2-1)\to Spin(h^2)$ et du plongement diagonal. Comme dans le cas de type $B_{n}$, on voit que ${\cal O}\cap \mathfrak{r}$ est \'egal \`a  l'orbite param\'etr\'ee par la partition $(2h-1,...,5,3)$.

 Supposons $G$  est  de type $(E_{6},nr)$. Alors $G_{s,SC}=SL(3)^3$. Le groupe $ \boldsymbol{\Omega}^{nr}$ a $3$ éléments. Son action sur $G_{s,SC}$ est engendrée par une permutation cyclique des trois facteurs. 
  L'orbite ${\cal O}$ est le produit des orbites r\'eguli\`eres. On a $R_{SC}=SL(3)$,  l'homomorphisme de $R_{SC}$ dans $G_{s,SC}$ \'etant  diagonal.  Alors ${\cal O}\cap \mathfrak{r}$ est l'orbite r\'eguli\`ere de $\mathfrak{r}$. 

 Supoosons $G$ de type $E_{7}$. Alors $G_{s,SC}=SL(4)\times SL(4)\times SL(2)$.   Le groupe $ \boldsymbol{\Omega}^{nr}$ s'identifie à $Aut({\cal D}_{a}) $. Il  a deux \'el\'ements. L'action sur $G_{s,SC}$ de l'élément non trivial   permute les deux facteurs $SL(4)$. L'orbite  ${\cal O}$ est le produit des orbites r\'eguli\`eres. On a $R_{SC}=SL(4)\times SL(2)$, l'homomorphisme de $R_{SC}$ dans $G_{s,SC}$ \'etant le produit de l'homomorphisme diagonal $SL(4)\to SL(4)\times SL(4)$ et de l'identit\'e de $SL(2)$.  Alors ${\cal O}\cap\mathfrak{r}$ est  
  le produit des orbites r\'eguli\`eres. 

  Supposons $G$ de type $(A_{n-1},ram)$ avec $n$ pair. On a $n=h^2+k(k+1)$ avec $h,k\in {\mathbb N}$ et $h=k$ ou $h=k+1$. Puisque $n$ est pair, $h$ l'est aussi.    Alors $G_{s,SC}=Spin(h^2)\times Sp(k(k+1))$. On a ${\cal O}={\cal O}_{1}\times {\cal O}_{2}$, o\`u ${\cal O}_{1}$, resp. ${\cal O}_{2}$, sont param\'etr\'ees par les partitions $(2h-1,...,3,1)$, resp. $(2k,...,4,2)$.  Le groupe  $ \boldsymbol{\Omega}^{nr}$ a deux éléments et l'élément non trivial  agit sur $G_{s,SC}$ par un automorphisme ext\'erieur non trivial du facteur $Spin(h^2)$. On a $R_{SC}=Spin(h^2-1)\times Sp(k(k+1))$, l'homomorphisme de $R_{SC}$ dans $G_{s,SC}$ \'etant \'evident. Comme dans le cas d'un groupe de type $B_{n}$, ${\cal O}\cap \mathfrak{r}$ est \'egal \`a ${\cal O}'_{1}\times {\cal O}_{2}$, o\`u ${\cal O}'_{1}$ est param\'etr\'e par la partition $(2h-1,...,5,3)$.

   Supposons $G$ de type $(D_{n},ram)$ avec $n=h^2$ où $h\geq 3$ est un entier impair. On a $G_{s}=SO(n)\times SO(n)$. Le groupe  $ \boldsymbol{\Omega}^{nr}$ a deux éléments et l'élément non trivial  agit sur $G_{s}$ par permutation des facteurs. On a ${\cal O}={\cal O}_{1}\times {\cal O}_{2}$, o\`u ${\cal O}_{1}$  et ${\cal O}_{2}$ sont toutes deux param\'etr\'ees par la partition $(2h-1,...,3,1)$. Alors $R=SO(n)$, l'homomorphisme de $R$ dans $G_{s}$ \'etant diagonal. L'ensemble ${\cal O}\cap \mathfrak{r}$ est l'orbite param\'etr\'ee par la m\^eme partition.  $\square$
   
   \subsubsection{Construction de la famille admissible $(\epsilon_{N})_{N\in {\cal O}}$ où ${\cal O}={\cal O}_{{\cal E}}$}\label{constructionfamilleadmissible}
   On conserve les données du paragraphe précédent. On définit un caractère $\boldsymbol{\epsilon}$  de $\boldsymbol{\Omega}^{nr}$ de la façon suivante, où on utilise les notations introduites dans la preuve de ce paragraphe:
   
   si $G$ est de type $(B_{n})$, le groupe $\boldsymbol{\Omega}^{nr}$ a $2$ éléments; on note $sgn$ son caractère non trivial; alors $\boldsymbol{\epsilon}=sgn^{(k+h+1)/2}$;
    
  si $G$ est de type $(A_{n-1},ram)$ avec $n$ pair, le groupe $\boldsymbol{\Omega}^{nr}$ a $2$ éléments; on note $sgn$ son caractère non trivial; alors $\boldsymbol{\epsilon}=sgn^{k+1}$;
  
  dans les autres cas, $\boldsymbol{\epsilon}=1$. 
  
  Soit $N\in {\cal O}\cap \mathfrak{r}$. Montrons que
  
  (1) $K_{s,N}^{\dag,nr}$ est le produit semi-direct $K_{s,N}^{0,nr}\rtimes \boldsymbol{\Omega}^{nr}$.
  
  La condition $N\in \mathfrak{r}$ assure que $\boldsymbol{\Omega}^{nr}$ est contenu dans $K_{s,N}^{\dag,nr}$. L'assertion résulte alors de \ref{unsousensemble}(1).
  
  Soit $k\in K_{s,N}^{\dag,nr}$. Ecrivons $k=k_{0}\omega$ avec $k_{0}\in K_{s,N}^{0,nr}$ et $\omega\in \boldsymbol{\Omega}^{nr}$. Posons $\epsilon_{N}(k)=\xi_{N}(\bar{k}_{0})\boldsymbol{\epsilon}(\omega)$. Cela définit l'application $\epsilon_{N}$. L'hypothèse $Hyp(1,N)$ est vérifiée en vertu de \ref{familles}(3). L'hypothèse $Hyp(2,N)$ est vérifiée par construction. D'après \ref{familles}(3), $\epsilon_{N}$ est la composante en $N$ d'une unique famille admissible $(\epsilon_{ N'})_{N'\in {\cal O}}$. 
  
  \begin{lem}{Cette famille admissible ne dépend pas des choix effectués. Dans le cas où $s\in S^{st}(G)$ et où ${\cal E}$ est conservé par l'action galoisienne, cette famille admissible est définie sur $F$.}\end{lem}
  
  Preuve. Les choix sont ceux de l'épinglage $\mathfrak{E}$, de l'épinglage affine $\mathfrak{E}_{a}$, autrement dit de l'élément supplémentaire $E_{\beta_{0}}$,  et de l'élément $N\in {\cal O}\cap \mathfrak{r}$. 
  
  Conservons l'épinglage $\mathfrak{E}$ et l'élément $E_{\beta_{0}}$ et remplaçons $N$ par un autre élément  $\tilde{N}\in {\cal O}\cap \mathfrak{r}$. On en déduit une autre famille admissible $(\tilde{\epsilon}_{N'})_{N'\in {\cal O}}$. Il nous suffit de prouver que $\tilde{\epsilon}_{\tilde{N}}=\epsilon_{\tilde{N}}$. D'après le lemme \ref{unsousensemble}, on peut fixer $\bar{r}\in R$ de sorte que $\tilde{N}=\bar{r}^{-1}N\bar{r}$. Relevons $\bar{r}$ en un élément $r\in K_{s}^{0,nr}$. Soit $k\in K_{s,N}^{\dag,nr}$. D'après $Hyp(3,N,\tilde{N})$, on a l'égalité 
  
  (2) $\epsilon_{\tilde{N}}(r^{-1}kr)=\epsilon_{N}(k)$.
  
   Ecrivons $k=k_{0}\omega$ avec $k_{0}\in K_{s,N}^{0,nr}$ et $\omega\in \boldsymbol{\Omega}^{nr}$. On a $r^{-1}kr=\tilde{k}'_{0}\tilde{h}\omega$, où $\tilde{k}'_{0}=r^{-1}k_{0}r$  et $\tilde{h}=r^{-1}\omega r\omega^{-1}$. On a $\tilde{k}'_{0}\in K_{s,\tilde{N}}^{0,nr}$,  $\tilde{h}\in K_{s}^{0,nr}$ et la réduction $\bar{\tilde{h}}$ est égale à $\bar{r}^{-1}\omega_{s}(\bar{r})$. Or $r\in R$, donc $\omega_{s}(\bar{r})=\bar{r}$. Donc $\bar{\tilde{h}}=1$, a fortiori $\tilde{h}$ appartient à $K_{s,\tilde{N}}^{0,nr}$. En posant $\tilde{k}_{0}=\tilde{k}'_{0}\tilde{h}$, l'égalité $r^{-1}kr=\tilde{k}_{0}\omega$ est la décomposition de $r^{-1}kr$ en produit d'un élément de $K_{s,\tilde{N}}^{0,nr}$ et d'un élément de $\boldsymbol{\Omega}^{nr}$. Par définition, on a $\tilde{\epsilon}_{\tilde{N}}(r^{-1}kr)=\xi_{\tilde{N}}(\bar{\tilde{k}}_{0})\boldsymbol{\epsilon}(\omega)$. On a $\bar{\tilde{k}}_{0}=\bar{r}^{-1}\bar{k}_{0}\bar{r}$ d'où, d'après \ref{familles}(1), $\xi_{\tilde{N}}(\bar{\tilde{k}}_{0})=\xi_{N}(\bar{k}_{0})$. Alors $\tilde{\epsilon}_{\tilde{N}}(r^{-1}kr)=\xi_{N}(\bar{k}_{0})\boldsymbol{\epsilon}(\omega)=\epsilon_{N}(k)$. En comparant avec (2), on obtient l'égalité requise $\tilde{\epsilon}_{\tilde{N}}(r^{-1}kr)=\epsilon_{\tilde{N}}(r^{-1}kr)$. 

Conservons l'épinglage $\mathfrak{E}$ et changeons $E_{\beta_{0}}$ en un autre élément $\tilde{E}_{\beta_{0}}$. On veut montrer que ce changement ne modifie pas notre famille admissible. L'élément $E_{\beta_{0}}$   n'est intervenu
dans nos constructions que via la définition du groupe  $\boldsymbol{\Omega}^{nr}$. Si ce groupe est égal à $1$, $E_{\beta_{0}}$ ne joue donc pas de rôle et la conclusion s'ensuit.  Supposons $\boldsymbol{\Omega}^{nr}\not=\{1\}$. On a introduit la racine $\beta_{s}$ associée à $s$. Un examen cas par cas montre que $\beta_{s}\not=\beta_{0}$. On affecte de\, $\tilde{}$ \,les objets construits à l'aide du nouveau choix $\tilde{E}_{\beta_{0}}$. Il est commode de noter $\tilde{E}_{\beta}=E_{\beta}$ pour $\beta\in \Delta_{a}^{nr}-\{\beta_{0}\}$. On a $\tilde{E}_{\beta_{0}}=\lambda E_{\beta_{0}}$ pour un  $\lambda\in \mathfrak{o}_{F^{nr}}^{\times}$.   L'ensemble $\Delta_{a}^{nr}-\{\beta_{s}\}$ est un ensemble linéairement indépendant de formes linéaires sur $X_{*}(T^{nr})$. Il existe $t\in T^{nr}(\mathfrak{o}_{F^{nr}})$ tel que $\beta_{0}(t)=\lambda$ et $\beta(t)=1$ pour $\beta\in \Delta_{a}^{nr}-\{\beta_{0},\beta_{s}\}$. On fixe un tel $t$. On a $Ad(t)(E_{\beta})=\tilde{E}_{\beta}$ pour tout $\beta\in \Delta_{a}^{nr}-\{\beta_{s}\}$, tandis qu'il existe $\mu\in \mathfrak{o}_{F^{nr}}^{\times}$ tel que $\mu Ad(t)(E_{\beta_{s}})=\tilde{E}_{\beta_{s}}$. 
 Montrons que

(3) $\tilde{\boldsymbol{\Omega}}^{nr}=\{t\omega t^{-1}; \omega\in \boldsymbol{\Omega}^{nr}\}$. 

Soit $\omega\in \boldsymbol{\Omega}^{nr}$ et $\beta\in \Delta_{a}^{nr}$. Posons $\gamma=\omega(\beta)$. On doit prouver que $Ad(t\omega t^{-1})(\tilde{E}_{\beta})=\tilde{E}_{\gamma}$. On sait que $\omega$ conserve $s$, donc fixe $\beta_{s}$. Si $\beta\not=\beta_{s}$, on a aussi $\gamma\not=\beta_{s}$. Alors $Ad(t\omega t^{-1})(\tilde{E}_{\beta})=Ad(t\omega)(E_{\beta})=Ad(t)(E_{\gamma})$, car $\omega\in \boldsymbol{\Omega}^{nr}$, d'où $Ad(t\omega t^{-1})(\tilde{E}_{\beta})=\tilde{E}_{\gamma}$. Si $\beta=\beta_{s}$, on a $\gamma=\beta_{s}$. Alors  $Ad(t\omega t^{-1})(\tilde{E}_{\beta_{s}})=\mu Ad(t\omega)(E_{\beta_{s}})=\mu Ad(t)(E_{\beta_{s}})=\tilde{E}_{\beta_{s}}$. Cela prouve (3).

Evidemment, l'application $\omega\mapsto t\omega t^{-1}$ identifie les caractères $\boldsymbol{\epsilon}$ de $\boldsymbol{\Omega}^{nr}$ et son analogue $\tilde{\boldsymbol{\epsilon}}$ pour $\tilde{\boldsymbol{\Omega}}^{nr}$.

En vertu de (3), on a l'égalité $\tilde{R}=Ad(\bar{t})(R)$. Soit $N\in {\cal O}\cap \mathfrak{r}$. Posons $\tilde{N}=Ad(\bar{t})(N)$. Alors $\tilde{N}\in {\cal O}\cap \tilde{\mathfrak{r}}$. Soit $k\in K_{s,N}^{\dag,nr}$. Posons $\tilde{k}=tkt^{-1}$. Alors $\tilde{k}\in K_{s,\tilde{N}}^{\dag,nr}$. D'après $Hyp(3,N,\tilde{N})$, on a l'égalité 

(4) $\epsilon_{\tilde{N}}(\tilde{k})=\epsilon_{N}(k)$. 

Ecrivons $k=k_{0}\omega$ avec $k_{0}\in K_{s,N}^{0,nr}$ et $\omega\in \boldsymbol{\Omega}^{nr}$. 
Alors $\tilde{k}=\tilde{k}_{0}\tilde{\omega}$, où $\tilde{k}_{0}=tk_{0}t^{-1}\in K_{s,\tilde{N}}^{0,nr}$ et $\tilde{\omega}=t\omega t^{-1}\in \tilde{\boldsymbol{\Omega}}^{nr}$. 
Par définition $\tilde{\epsilon}_{\tilde{N}}(\tilde{k})=\xi_{\tilde{N}}(\bar{\tilde{k}}_{0})\tilde{\boldsymbol{\epsilon}}(\tilde{\omega})$. On a $\xi_{\tilde{N}}(\bar{\tilde{k}}_{0})=\xi_{\tilde{N}}(\bar{t}\bar{k}_{0}\bar{t}^{-1})=\xi_{N}(\bar{k}_{0})$ d'après \ref{familles}(1). On a aussi $\tilde{\boldsymbol{\epsilon}}(\tilde{\omega})=\boldsymbol{\epsilon}(\omega)$. D'où $\tilde{\epsilon}_{\tilde{N}}(\tilde{k})=\xi_{N}(\bar{k}_{0})\boldsymbol{\epsilon}(\omega)=\epsilon_{N}(k)$. Avec (4), cela démontre l'égalité requise $\tilde{\epsilon}_{\tilde{N}}=\epsilon_{\tilde{N}}$.  

  Consid\'erons maintenant un autre \'epinglage $\tilde{\mathfrak{E}}=(\tilde{B},\tilde{T},(\tilde{E}_{\alpha})_{\alpha\in \Delta})$ de $G$ v\'erifiant les m\^emes hypoth\`eses que $\mathfrak{E}$: il est fix\'e  par l'action de $I_{F}$; en notant $\tilde{C}^{nr}$ l'alc\^ove de $Imm_{F^{nr}}(G)$ qui s'en d\'eduit, $s$ appartient \`a $S(\bar{\tilde{C}}^{nr})$. Puisque $\mathfrak{E}$ et $\tilde{\mathfrak{ E}}$ sont tous deux  d\'efinis sur $F^{nr}$ et que $G$ est adjoint, il existe un unique $g\in G(F^{nr})$ tel que $Ad(g)(\mathfrak{ E})=\tilde{\mathfrak{ E}}$. Montrons que
  
  (5) on a $g\in K_{s}^{\dag,nr}$. 
  
  L'action de $g$ sur $Imm_{F^{nr}}(G)$ envoie $C^{nr}$ sur $\tilde{C}^{nr}$. Puisque $s\in S(\bar{\tilde{C}}^{nr})$, on a $g^{-1}s\in S(\bar{C}^{nr})$. Mais on a aussi $s\in S(\bar{C}^{nr})$. En appliquant le corollaire \ref{uncorollaire}, il existe $\omega\in \boldsymbol{\Omega}^{nr}$ tel que $g^{-1}s=\omega s$. Parce que $s\in S^{nr,st}(G)$, $s$ est fixé par $\boldsymbol{\Omega}^{nr}$, d'où $g^{-1}s=s$, ce qui démontre (5).
  
  On peut compléter l'épinglage $\tilde{\mathfrak{E}}$ en un épinglage affine $\tilde{\mathfrak{E}}_{a}$ en posant $\tilde{E}_{\tilde{\beta}_{0}}=Ad(g)(E_{\beta_{0}})$. A l'aide de ces nouvelles données, on construit de nouvelles applications $\tilde{\epsilon}_{N}$ dont on doit prouver qu'elles sont égales aux applications $\epsilon_{N}$ initiales. De nouveau, on affecte d'un  $\,\,\tilde{}\,\,$  les nouveaux objets. Il est clair que $\tilde{\boldsymbol{\Omega}}^{nr}=\{g\omega g^{-1}; \omega\in \boldsymbol{\Omega}^{nr}\}$ et que l'application $\omega\mapsto g\omega g^{-1}$ identifie les caractères $\boldsymbol{\epsilon}$ de $\boldsymbol{\Omega}^{nr}$ et $\tilde{\boldsymbol{\epsilon}}$ de $\tilde{\boldsymbol{\Omega}}^{nr}$. On a aussi $\tilde{R}=Ad(g)_{s}(R)$. Soit $N\in {\cal O}\cap \mathfrak{r}$. Posons $\tilde{N}=Ad(g)_{s}(N)$. Alors $\tilde{N}\in {\cal O}\cap \tilde{\mathfrak{r}}$. Soit $k\in K_{s,N}^{\dag,nr}$. Posons $\tilde{k}=gkg^{-1}$. On a $\tilde{k}\in K_{s,\tilde{N}}^{\dag,nr}$. D'après $Hyp(3,N,\tilde{N})$, on a 
  
  (6) $\epsilon_{\tilde{N}}(\tilde{k})=\epsilon_{N}(k)$.  
  
  Ecrivons $k=k_{0}\omega$ avec $k_{0}\in K_{s,N}^{0,nr}$ et $\omega\in \boldsymbol{\Omega}^{nr}$. En posant $\tilde{k}_{0}=gk_{0}g^{-1}$ et $\tilde{\omega}=g\omega g^{-1}$, on a $\tilde{k}=\tilde{k}_{0}\tilde{\omega}$ et $\tilde{k}_{0}\in K_{s,\tilde{N}}^{0,nr}$, $\tilde{\omega}\in \tilde{\boldsymbol{\Omega}}^{nr}$. Par définition, $\tilde{\epsilon}_{\tilde{N}}(\tilde{k})=\xi_{\tilde{N}}(\bar{\tilde{k}}_{0})\tilde{\boldsymbol{\epsilon}}(\tilde{\omega})$. On a $\tilde{\boldsymbol{\epsilon}}(\tilde{\omega})=\boldsymbol{\epsilon}(\omega)$. On a $\bar{\tilde{k}}_{0}=Ad(g)_{s}(\bar{k}_{0})$, d'où, en utilisant \ref{familles}(3), $\xi_{\tilde{N}}(\bar{\tilde{k}}_{0})=\xi_{N}(\bar{k}_{0})$. Alors $\tilde{\epsilon}_{\tilde{N}}(\tilde{k})=\xi_{N}(\bar{k}_{0})\boldsymbol{\epsilon}(\omega)=\epsilon_{N}(k)$.  Avec (6), cela démontre l'égalité requise $\tilde{\epsilon}_{\tilde{N}}=\epsilon_{\tilde{N}}$. Cela achève la preuve de la première assertion de l'énoncé.

    Supposons que $s\in S^{st}(G)$ et que ${\cal E}\in {\bf FC}^{st}_{s,{\mathbb F}_{q}}$. Pour d\'emontrer que notre famille admissible est d\'efinie sur $F$, il suffit d'apr\`es  \ref{familles} (4) que, pour un seul \'el\'ement  $N\in {\cal O}^{\Gamma_{{\mathbb F}_{q}}}$, la condition $Hyp(4,N)$ soit v\'erifi\'ee.  
    On vient de  construire  la famille admissible $(\epsilon_{N'})_{N'\in {\cal O}}$ à l'aide  d'un épinglage $\mathfrak{E}$ et d'un épinglage affine $\mathfrak{E}_{a}$. On vient de démontrer qu'elle est en fait indépendante de ces choix. Parce que $s\in S^{st}(G)$, l'épinglage $Fr(\mathfrak{E})$ vérifie les m\^emes conditions requises de $\mathfrak{E}$. On peut donc utiliser aussi l'épinglage $Fr(\mathfrak{E})$ et l'épinglage affine $Fr(\mathfrak{E}_{a})$.
     Pour ces nouveaux choix, le groupe $R$ est remplacé par $Fr(R)$. Soit $N'\in {\cal O}\cap \mathfrak{r}$. On a $Fr(N')\in {\cal O}\cap Fr(\mathfrak{r})$. En utilisant les choix initiaux pour construire $\epsilon_{N'}$ et les nouveaux choix pour construire $\epsilon_{Fr(N')}$, il est clair que l'on a l'égalité
     
     (7)  $\epsilon_{Fr(N')}(Fr(k))=\epsilon_{N'}(k)$ pour tout $k\in K_{s,N'}^{\dag,nr}$. 
     
     Fixons $g\in K_{s}^{0,nr}$ tel que $N=\bar{g}N'\bar{g}^{-1}$. Parce que $Fr(N)=N$,  on a aussi $N=Fr(\bar{g})Fr(N')Fr(\bar{g})^{-1}$.  Soit $k\in K_{s,N}^{\dag,nr}$. En appliquant $Hyp_{3}(N,N')$ et $Hyp_{3}(N,Fr(N'))$, on a
 $\epsilon_{N}(k)=\epsilon_{N'} (g^{-1}kg)$ et $\epsilon_{N}(Fr(k))=\epsilon_{Fr(N')}(Fr(g)^{-1}Fr(k)Fr(g))$. Il reste à appliquer (7) pour obtenir l'égalité $\epsilon_{N}(k)=\epsilon_{N}(Fr(k))$, c'est-à-dire   $Hyp(4,N)$.  Cela achève la preuve du lemme. $\square$

{\bf Remarques} (7). Dans le cas où $\boldsymbol{\Omega}^{nr}=\{1\}$, il est clair d'après (1) qu'il existe une unique famille admissible. 

(8) On a supposé que $G$ était défini sur $F$ mais la définition de la famille admissible n'utilise que la structure de $G$ sur $F^{nr}$.  Elle est donc bien définie si $G$ est seulement défini sur $F^{nr}$. 

\subsubsection{Equivariance par automorphisme}\label{equivariance}
Soit $G'$ un groupe réductif défini sur $F^{nr}$ et soit $\delta:G\to G'$ un isomorphisme défini sur $F^{nr}$. Soient $s$ et ${\cal E}$ comme dans le paragraphe précédent. On en déduit une famille admissible $(\epsilon_{N})_{N\in {\cal O}}$. De $\delta$ se déduit un isomorphisme $\delta^{Imm}:Imm_{F^{nr}}(G)\to Imm_{F^{nr}}(G')$. Posons $s'=\delta^{Imm}(s)$. Alors de $\delta$ se déduit un isomorphisme $\delta_{s}:G_{s}\to G'_{s'}$. Cet isomorphisme transporte le faisceau ${\cal E}$ en un faisceau ${\cal E}'$ sur $\mathfrak{g}'_{s'}$, qui appartient à ${\bf FC}^{st}_{s'}$. Posons  ${\cal O}'={\cal O}_{{\cal E}'}=\delta_{s}({\cal O})$. On construit la famille admissible $(\epsilon_{N'})_{N'\in {\cal O}'}$. Soit $N\in {\cal O}$, posons $N'=\delta_{s}(N)$. On a  $K_{s',N'}^{G',\dag,nr}=\delta(K_{s,N}^{\dag,nr})$. Parce que nos familles sont canoniques, il résulte d'un simple transport de structure que l'on a l'égalité

(1) $\epsilon_{N'}(\delta(k))=\epsilon_{N}(k)$ pour tout $k\in K_{s,N}^{\dag,nr}$.

\subsubsection{Prolongement du caractère $\epsilon_{N}$}\label{prolongement}
On reprend les données  $s$ et ${\cal E}$ de \ref{unsousensemble}.    On a défini le groupe  $Aut_{s}(G,F^{nr})$ en \ref{familles}.   On fixe $N\in {\cal O}$ et on  note $Aut_{s,N}(G,F^{nr})$ le sous-groupe des $\delta\in  Aut_{s}(G,F^{nr})$ tels que $\delta_{s}(N)=N$. Pour $\delta\in Aut_{s,N}(G,F^{nr})$, on a l'égalité $\xi_{N}=\xi_{N}\circ \delta_{s}$. 
     Le   groupe  $Aut_{s,N}(G,F^{nr})$ contient $K_{s,N}^{\dag,nr}=G(F^{nr})\cap Aut_{s,N}(G,F^{nr})$ comme sous-groupe distingué et le quotient $K_{s,N}^{\dag,nr}\backslash Aut_{s,N}(G,F^{nr})$ s'envoie injectivement dans le groupe $Out(G)^{I_{F}}$ des points fixes par $I_{F}$ dans le groupe des automorphismes extérieurs de $G$ (en fait bijectivement). En particulier, ce quotient est cyclique d'ordre au plus $3$ sauf dans le cas ou $G$ est de type $(D_{4},nr)$.

   \begin{lem}{ Le caractère $\epsilon_{N}$ de $K_{s,N}^{\dag,nr}$ se prolonge en un caractère de $Aut_{s,N}(G,F^{nr})$.}\end{lem} 
    
  Preuve.  Supposons d'abord que $G$ n'est pas du type $(D_{4},nr)$. Puisque $Aut_{s,N}(G,F^{nr})$ est une extension d'un groupe cyclique d'ordre premier par $K_{s,N}^{\dag,nr}$,  il suffit de prouver que l'action par conjugaison de $Aut_{s,N}(G,F^{nr})$ sur $K_{s,N}^{\dag,nr}$ préserve $\epsilon_{N}$. Mais c'est un cas particulier de l'assertion (1) de \ref{equivariance}. 
  
  Supposons maintenant que $G$ est du type $(D_{4},nr)$. On identifie $Out(G)$ au groupe des automorphismes de $G$ qui conservent l'épinglage $\mathfrak{E}$. Ils conservent aussi  l'épinglage affine $\mathfrak{E}_{a}$ et $Aut({\cal D}_{a})$ s'identifie au produit semi-direct $\boldsymbol{\Omega}^{nr}\rtimes Out(G)$. Le groupe $G_{s}$ a pour épinglage l'ensemble des réductions  $\bar{E}_{\alpha_{i}}$ dans $\mathfrak{g}_{s}$ des éléments $E_{\alpha_{i}}$ pour $i\in \{0,1,3,4\}$. L'élément $N_{0}=\sum_{i=0,1,3,4}\bar{E}_{i}$ appartient à ${\cal O}$. On vérifie aisément que l'assertion de l'énoncé ne dépend pas de l'élément $N\in {\cal O}$ (ceci étant bien s\^ur vrai pour tout $G$). On peut donc supposer que $N=N_{0}$. Alors $Out(G)$ est inclus dans $Aut_{s,N}(G,F^{nr})$ et ce  dernier groupe s'identifie au produit semi-direct $K_{s,N}^{\dag,nr}\rtimes Out(G)$, ou encore  $K_{s,N}^{0,nr}\rtimes (\boldsymbol{\Omega}^{nr}\rtimes Out(G))$.   
   On définit une application $\tilde{\epsilon}_{N}:Aut_{s,N} (G,F^{nr})\to {\mathbb C}^{\times}$  par $\tilde{\epsilon}_{N}(k\times \delta)
=\epsilon_{N}(k)$ pour $k\in K_{s,N}^{\dag,nr}$ et $\delta\in Out(G)$. On voit que c'est un caractère qui prolonge $\epsilon_{N}$.   $\square$
  
  \subsubsection{Composition d'isomorphismes $\iota_{\phi}$}\label{composition}
  On considère maintenant trois groupes $G,H,H'$ définis sur $F$, adjoints et absolument simples, et trois torseurs intérieurs s'insérant dans un diagramme
   $$\begin{array}{ccccc} H&&\stackrel{\psi}{\to}&&H'\\ &\phi \searrow\,\,&&\,\,\swarrow \phi'&\\ && G&&\\ \end{array}$$
   \begin{lem}{On a l'égalité $\iota_{\phi}=\iota_{\phi'}\circ\iota_{\psi}$.}\end{lem}
   Preuve. Fixons des sommets $s_{H}\in \underline{S}^{st}(H)$, $s_{H'}\in \underline{S}^{st}(H')$, $s\in \underline{S}^{st}(G)$, supposons que $\underline{\psi}(s_{H})=s_{H'}$ et $\underline{\phi}'(s_{H'})=s$. Cela entraîne $\underline{\phi}(s_{H})=s$. On peut modifier $\phi$, $\phi'$ et $\psi$ de sorte que $\psi^{Imm}(s_{H})=s_{H'}$, ${\phi'}^{Imm}(s_{H'})=s$, d'où $\phi^{Imm}(s_{H})=s$. On a alors des isomorphismes
   $$\begin{array}{ccccc} H_{s_{H}}&&\stackrel{\psi_{s_{H'}}}{\to}&&H'_{s_{H'}}\\&\phi _{s}\searrow\,\,&&\,\,\swarrow \phi'_{s}&\\ && G_{s}&&\\ \end{array}$$
   Fixons ${\cal E}\in {\bf FC}^{st}_{{\mathbb F}_{q}}(\mathfrak{g}_{s})$, posons ${\cal E}_{H'}={\phi_{s}'}^{-1}({\cal E})$, ${\cal E}_{H}=\psi_{s_{H'}}^{-1}({\cal E}_{H'})=\phi_{s}^{-1}({\cal E})$. On fixe aussi des éléments $N\in {\cal O}_{{\cal E}}$, $N_{H}\in {\cal O}_{{\cal E}_{H}}$ et $N_{H'}\in {\cal O}_{{\cal E}_{H'}}$, tous trois fixés par $\Gamma_{{\mathbb F}_{q}}$. Quitte à modifier encore nos torseurs, on peut supposer que $\psi_{s_{H'}}(N_{H})=N_{H'}$, $\phi'_{s}(N_{H'})=N$, d'où $\phi_{s}(N_{H})=N$. Fixons une action de Frobenius sur ${\cal E}$ définie par un scalaire $r_{N}$. En appliquant la définition de \ref{resolution} au couple $(H',G)$, on en déduit une action de Frobenius sur ${\cal E}_{H'}$ définie par un scalaire $r_{N_{H'}}$. En appliquant la m\^eme définition au couple $(H,H')$, on en déduit une action de Frobenius sur ${\cal E}_{H}$ définie par un scalaire que l'on note $\tilde{r}_{N_{H}}$. En appliquant la m\^eme définition au couple $(H,G)$, on munit ${\cal E}_{H}$ d'une action de Frobenius, a priori différente de la précédente, définie par un scalaire que l'on note $r_{N_{H}}$. Le lemme équivaut à dire que ces deux actions de Frobenius sont les m\^emes, autrement dit que l'on a l'égalité $\tilde{r}_{N_{H}}=r_{N_{H}}$. On a
   $$\tilde{r}_{N_{H}}=r_{N_{H'}}\epsilon_{s_{H'},{\cal E}_{H'}, N_{H'}}(k_{\psi})=r_{N}\epsilon_{s,{\cal E},N}(k_{\phi'})\epsilon_{s_{H'},{\cal E}_{H'}, N_{H'}}(k_{\psi}).$$
   En utilisant \ref{equivariance}(1) pour l'isomorphisme $\delta=\phi'$, on a 
   $$\epsilon_{s_{H'},{\cal E}_{H'}, N_{H'}}(k_{\psi})=\epsilon_{s,{\cal E},N}(\phi'(k_{\psi})).$$
   D'où
   $$(1) \qquad \tilde{r}_{N_{H}}=r_{N}\epsilon_{s,{\cal E},N}(k_{\phi'}\phi'(k_{\psi}) ).$$
   Par définition, on a $Ad(k_{\phi'})=Fr(\phi'){\phi'}^{-1}$ et $Ad(k_{\psi})=Fr(\psi)\psi^{-1}$. D'où $Ad(\phi'(k_{\psi}))=\phi' Ad(k_{\psi}) {\phi'}^{-1}=\phi'Fr(\psi)\psi^{-1}{\phi'}^{-1}$, puis 
   $$Ad(k_{\phi'}\phi'(k_{\psi}))=Fr(\phi'){\phi'}^{-1}\phi'Fr(\psi)\psi^{-1}{\phi'}^{-1}=Fr(\phi')Fr(\psi)\psi^{-1}{\phi'}^{-1}=Fr(\phi)\phi^{-1}=Ad(k_{\phi}).$$
   Cela démontre l'égalité $k_{\phi'}\phi'(k_{\psi})=k_{\phi}$. Alors (1) devient
   $$ \tilde{r}_{N_{H}}=r_{N}\epsilon_{s,{\cal E},N}(k_{\phi})=r_{N_{H}}.$$
   C'est ce qu'il fallait démontrer. $\square$

\subsection{Isomorphismes d'espaces $FC^{st}(\mathfrak{g}(F))$ et transfert endoscopique}

\subsubsection{Enoncé du théorème}\label{enonce}
On considère deux groupes $G$ et $H$ comme en \ref{position}. On définit l'isomorphisme $\iota_{\phi}:FC^{st}(\mathfrak{h}(F))\to FC^{st}(\mathfrak{g}(F))$ comme en   \ref{resolution} où l'on utilise la $G$-famille admissible definie dans \ref{constructionfamilleadmissible}. Puisque $\phi:H\to G$ est un torseur intérieur, il s'en  déduit  un isomorphisme de  transfert endoscopique $transfert_{\phi}:
I^{st}_{cusp}(\mathfrak{h}(F))\to I^{st}_{cusp}(\mathfrak{g}(F))$ qui 
 se restreint en un isomorphisme $transfert_{\phi}:FC^{st}(\mathfrak{h}(F))\to FC^{st}(\mathfrak{g}(F))$.
 
 \begin{thm}{On a l'égalité $\iota_{\phi}=transfert_{\phi}$.}\end{thm}
 
 La sous-section est consacrée à la preuve de ce théorème. Introduisons une forme intérieure $G^*$ de $G$ quasi-déployée sur $F$ et un torseur intérieur $\psi:G\to G^*$. Posons $\psi'=\psi\circ\phi$ qui est encore un torseur intérieur. D'après le lemme \ref{composition}, on a l'égalité $\iota_{\phi}=\iota_{\psi}^{-1}\circ\iota_{\psi'}$. On a l'égalité analogue pour les transferts endoscopiques: $transfert_{\phi}=transfert_{\psi}^{-1}\circ transfert_{\psi'}$. Il suffit donc de démontrer le théorème pour les torseurs intérieurs $\psi$ et $\psi'$. Autrement dit, en oubliant cette construction, on peut supposer que $G$ est quasi-déployé sur $F$ et on pose cette hypothèse pour la suite de la sous-section.

 \subsubsection{Pr\'eparatifs \`a la preuve du théorème  \ref{enonce}\label{preparatifs}}
  On se place dans la situation de \ref{enonce}. Supposons d'abord que $\phi$ soit équivalent à un isomorphisme d\'efini sur $F$, c'est-\`a-dire qu'il existe $x\in G(F^{nr})$ tel que $Ad(x) \phi$ soit d\'efini sur $F$. L'isomorphisme $transfert_{\phi}$ est insensible au remplacement de $\phi$ par $Ad(x) \phi$, on peut donc aussi bien supposer que $\phi$ est un isomorphisme d\'efini sur $F$. Donc  $transfert_{\phi}$ n'est autre que l'isomorphisme d\'eduit de $\phi$ par transport de structure. Il en est de m\^eme de $\iota_{\phi}$  d'après  \ref{resolution}(6), donc $\iota_{\phi}=transfert_{\phi}$. On suppose d\'esormais que $\phi$ n'est pas équivalent à un isomorphisme d\'efini sur $F$. 
En excluant les cas où $FC^{st}(\mathfrak{g}(F))=\{0\}$, il reste les 
 types suivants

$(A_{n-1},ram)$ avec $n$ pair, $B_{n}$, $C_{n}$, $(D_{n},nr)$ avec $n$ pair, $(D_{n},ram)$ avec $n$ impair, $E_{6}$ déployé, $E_{7}$. 

On peut pr\'eciser que, dans le cas o\`u $G$ est de type $(D_{4},nr)$, l'action du Frobenius sur ${\cal D}$ est d'ordre au plus $2$: un groupe trialitaire n'a pas de forme int\'erieure autre qu'elle-m\^eme.

Fixons un sommet $s\in S^{st}(G)$. Reprenons les constructions et notations de  \ref{unsousensemble}. On suppose que  l'\'epinglage $\mathfrak{E} $ et l'épinglage affine $\mathfrak{E}_{a}$ sont conservés par $\Gamma_{F}$. 
 En cons\'equence du (iii) du lemme \ref{uncorollaire}, quitte \`a remplacer $\phi$ par $Ad(x) \phi$ pour un \'el\'ement $x\in G(F^{nr})$ convenable, ce qui ne change pas $transfert_{\phi}$, on peut supposer que $\phi Fr(\phi)^{-1}$  est de la forme $Ad(\boldsymbol{\omega})$ pour un $\boldsymbol{\omega}\in \boldsymbol{\Omega}^{nr}$. Posons  alors $s_{H}= (\phi^{Imm})^{-1}(s)$. Puisque $(\phi Fr(\phi)^{-1})^{Imm}$ fixe $s$, $s_{H}$ est fix\'e par l'action galoisienne, donc $s_{H}\in S^{st}(H)$. Notons $\sigma\mapsto \boldsymbol{\omega}(\sigma)$ l'unique cocycle de $\Gamma_{F}^{nr}$ dans $\boldsymbol{\Omega}^{nr}$ tel que $\boldsymbol{\omega}(Fr)=\boldsymbol{\omega}$. On le rel\`eve par inflation en un cocycle de $\Gamma_{F}$ dans $\boldsymbol{\Omega}^{nr}$.  
On peut identifier $H$ \`a son image dans $G$ par $\phi$. L'action galoisienne sur $H$ est d\'efinie par  $\sigma\mapsto \sigma_{H} =Ad(\boldsymbol{\omega}(\sigma))\circ \sigma_{G}$. Par l'identification $\Gamma_{F}^{nr}\simeq \Gamma_{{\mathbb F}_{q}}$,  le cocycle $\sigma\mapsto \boldsymbol{\omega}(\sigma)$ s'identifie \`a un cocycle  d\'efini sur $\Gamma_{{\mathbb F}_{q}}$.   On peut identifier $H_{s_{H}}$ \`a $G_{s}$ muni de l'action galoisienne $\sigma\mapsto \sigma_{H}=Ad(\boldsymbol{\omega}(\sigma))_{s}\circ\sigma_{G}$, pour $\sigma\in \Gamma_{{\mathbb F}_{q}}$ (pour simplifier, on utilise les notations  $\sigma_{G}$ et $\sigma_{H}$ plut\^ot que $\sigma_{G_{s}}$ et $\sigma_{H_{s_{H}}}$).

Si $G$ est classique, l'espace $FC^{st}(\mathfrak{g}_{s}({\mathbb F}_{q}))$ est une droite et l'ensemble ${\bf FC}^{st}_{s,{\mathbb F}_{q}}$ est form\'e d'un unique faisceau-caract\`ere dont on note ${\cal O}$ le support. Si $G$ est de type $E_{6}$ déployé ou $E_{7}$, l'espace $FC^{st}(\mathfrak{g}_{s}({\mathbb F}_{q}))$, qui est suppos\'e non nul, est de dimension $2$. Rappelons en passant que la non-nullit\'e de l'espace pr\'ec\'edent impose les conditions $\delta_{3}(q-1)=1$ dans le cas $E_{6}$ et $\delta_{4}(q-1)=1$ dans le cas $E_{7}$, o\`u, pour $n\in {\mathbb N}_{>0}$, on note $\delta_{n}$ la fonction caract\'eristique du sous-groupe $n{\mathbb Z}$ de ${\mathbb Z}$ . L'ensemble ${\bf FC}^{st}_{s,{\mathbb F}_{q}}$ a deux \'el\'ements. Mais ces deux faisceaux sont support\'es par la m\^eme orbite nilpotente, que l'on note ${\cal O}$. On a d\'efini un sous-groupe  $R$ de $ G_{s}$,  Il est fix\'e par $Ad(k)_{s}$ pour tout \'el\'ement  $k\in \boldsymbol{\Omega}^{nr}$. En cons\'equence, les actions galoisiennes $\sigma\mapsto \sigma_{G}$ et $\sigma\mapsto \sigma_{H}$ co\"{\i}ncident sur $R$. On fixe $N\in {\cal O}\cap \mathfrak{r}({\mathbb F}_{q})$. Chaque ${\cal E}\in {\bf FC}^{st}_{s,{\mathbb F}_{q}}$ correspond \`a un caract\`ere $\xi_{{\cal E}}$ du groupe $Z_{G_{s}}(N)/Z_{G_{s}}(N)^0$. On pose $\Xi=\{\xi_{{\cal E}},{\cal E}\in {\bf FC}^{st}_{s,{\mathbb F}_{q}}\}$.  On munit un tel ${\cal E}$ de l'action galoisienne d\'etermin\'ee par le scalaire $r_{N} =1$. Puisque $N\in \mathfrak{r}$, on a $\boldsymbol{\omega}\in K_{s,N}^{\dag,nr}$. En reprenant les d\'efinitions, on voit que $\epsilon_{s,{\cal E},N}(\boldsymbol{\omega})=1$ sauf dans les deux cas ci-dessous.

(1)  Supposons que $G$ est de type $(B_{n})$.  On a $2n+1=k^2+h^2$ avec $k,h\in {\mathbb N}$, $k$ est pair et $\vert h-k\vert =1$. Le groupe $\boldsymbol{\Omega}^{nr} $ a deux \'el\'ements.  L'hypoth\`ese que $\phi$ n'est pas un isomorphisme sur $F$ modulo automorphismes int\'erieurs signifie que $\boldsymbol{\omega}$ est l'\'el\'ement non trivial de $\boldsymbol{\Omega}^{nr}$. Il r\'esulte des d\'efinitions que $\epsilon_{s,{\cal F},N}(\boldsymbol{\omega})=1$ si $h=k-1$ et $\epsilon_{s,{\cal E},N}(\boldsymbol{\omega})=-1$ si $h=k+1$.

(2) Supposons $G$ de 
  type $(A_{n-1},ram)$ avec $n$ pair. On a $n=h^2+k(k+1)$ avec $h,k\in {\mathbb N}$ et $h=k$ ou $h=k+1$.  Comme ci-dessus, $\boldsymbol{\omega}$ est l'\'el\'ement non trivial de $\boldsymbol{\Omega}^{nr}$. On a $\epsilon_{s,{\cal E},N}(\boldsymbol{\omega})=1$ si $h=k+1$ et $\epsilon_{s,{\cal E},N}(\boldsymbol{\omega})=-1$ si $h=k$.
  
  En tout cas, nous poserons simplement $\epsilon=\epsilon_{s,{\cal E},N}(\boldsymbol{\omega})$. Par l'identification $H_{s_{H}}=G_{s}$, ${\cal E}$ s'identifie \`a un faisceau ${\cal E}_{H}$ sur $\mathfrak{h}_{s_{H}}$. On le munit de l'action galoisienne tel que le Frobenius agisse par multiplication par $\epsilon$ sur la fibre au-dessus de $N$. En posant $\xi=\xi_{{\cal E}}$, on note $f_{\xi}^G$ la fonction caract\'eristique de ${\cal E}$ et $f_{\xi}^H$ celle de ${\cal E}_{H}$. Ces fonctions sont d\'efinies sur $\mathfrak{g}_{s}({\mathbb F}_{q})$, resp. $\mathfrak{h}_{s_{H}}({\mathbb F}_{q})$. On en d\'eduit par rel\`evement des fonctions sur $\mathfrak{g}(F)$, resp. $\mathfrak{h}(F)$, que l'on a souvent not\'ees de la m\^eme fa\c{c}on. Pour \'eviter les confusions, on les note ici ${\bf f}_{\xi}^G$ et ${\bf f}_{\xi}^H$. D'apr\`es la d\'efinition (3) de \ref{resolution}, on a l'\'egalit\'e
 $\iota_{\phi}({\bf f}^H_{\xi})=c(s){\bf f}^G_{\xi}$,
 o\`u $c(s)=\vert G_{s}({\mathbb F}_{q})\vert^{-1} \vert H_{s_{H}}({\mathbb F}_{q})\vert $.  Le théorème \ref{enonce} \'equivaut \`a la relation 
  
 (3) $transfert_{\phi}({\bf f}^{H}_{\xi})=c(s){\bf f}^G_{\xi}$ pour tout $\xi\in \Xi$.

\subsubsection{D\'ebut  de la preuve \label{debutdelapreuve}}
  Compl\'etons $N$ en un $\mathfrak{sl}(2)$-triplet $(f,h,N)$ contenu dans $\mathfrak{r}({\mathbb F}_{q})$. Il s'en d\'eduit un groupe à un paramètre $x_{*,h}$, cf. \ref{sl2triplets} puis  
   une graduation $(\mathfrak{g}_{s,i})_{i\in {\mathbb Z}}$ de $\mathfrak{g}_{s}$: pour $i\in {\mathbb Z}$, $\mathfrak{g}_{s,i}=\{X\in \mathfrak{g}_{s};  \forall t\in \bar{{\mathbb F}}_{q}^{\times}, \,\, Ad(x_{*,h}(t))X=t^{i}X\}$. 
  
  {\bf Remarque.} On v\'erifie cas par cas que notre hypoth\`ese $(Hyp)_{1}(p)$, qui implique $p\geq 3h(G)$, implique $p\geq 4h(G_{s})$. On voit alors que, pour tout $i$ tel que  $\mathfrak{g}_{s,i}\not=\{0\}$, $\mathfrak{g}_{s,i}$ est aussi l'ensemble des $X\in \mathfrak{g}_{s}$ tels que $ [h,X]=iX$.
  
  \bigskip
  En fait, on a $\mathfrak{g}_{s,i}=\{0\}$ pour $i$ impair. On note $P$ le sous-groupe parabolique de $G_{s}$ d'alg\`ebre de Lie $\oplus_{i\geq0} \mathfrak{g}_{s,i}$ et $M$ sa composante de Levi d'alg\`ebre de Lie $\mathfrak{g}_{s,0}$. L'action du groupe $M$  dans $\mathfrak{g}_{s,2}$ admet une unique orbite ouverte que l'on note $\tilde{\mathfrak{g}}_{s,2}$. 
    L'\'el\'ement $N$ appartient \`a cet ensemble. Parce que le $\mathfrak{sl}(2)$-triplet est contenu dans $\mathfrak{r}({\mathbb F}_{q})$, tous ces objets sont conserv\'es par les deux actions galoisiennes $\sigma\mapsto \sigma_{G}$ et $\sigma\mapsto \sigma_{H}$. Notons $\mathfrak{k}_{u}^{nr}$ l'image r\'eciproque dans $\mathfrak{k}_{s}^{nr}$ du radical nilpotent de l'alg\`ebre de Lie de $P$, c'est-\`a-dire de $\oplus_{i\geq 2}\mathfrak{g}_{s,i}$. Consid\'erons un groupe r\'eductif connexe $\tilde{G}$ d\'efini sur $F^{nr}$ et un homomorphisme surjectif $\tilde{G}\to G$ d\'efini sur $F^{nr}$ de noyau central dans $\tilde{G}$. On a alors  $G=\tilde{G}_{AD}$ et  il y a un homomorphisme surjectif $\tilde{G}_{s}\to G_{s}$. Notons $\tilde{P}$ l'image r\'eciproque de $P$ dans $\tilde{G}_{s}$ et notons $K^{\tilde{P},nr}_{s}$ le sous-groupe des $x\in K_{s}^{\tilde{G},\dag,nr}$, c'est-\`a-dire des $x\in \tilde{G}(F^{nr})$ qui fixent le sommet $s$, tels que l'action r\'eduite de $Ad(x)$ dans $\tilde{G}_{s}$ conserve $\tilde{P}$.    
 
\begin{lem}{Supposons construits des familles $(X^G_{\xi})_{\xi\in \Xi}$ et $(X^H_{\xi})_{\xi\in \Xi}$ v\'erifiant les conditions suivantes:

(i) pour tout $\xi\in \Xi$, $X^G_{\xi}$ est un \'el\'ement elliptique r\'egulier de $\mathfrak{g}(F)$,  $X^H_{\xi}$ est un \'el\'ement elliptique r\'egulier de $\mathfrak{h}(F)$ et les classes de conjugaison stable de ces \'el\'ements se correspondent;

(ii) pour tout $\xi\in \Xi$, $X^G_{\xi}$   appartient \`a $\mathfrak{k}_{u}^{nr}$; l'ensemble des \'el\'ements $x\in \tilde{G}(F^{nr})$ tels que $x^{-1}X^G_{\xi}x\in \mathfrak{k}_{u}^{nr}$ est \'egal \`a $K_{s}^{\tilde{P},nr}$;  $X^H_{\xi}$ appartient \`a $\mathfrak{k}_{u}^{nr}$ et l'ensemble des \'el\'ements $x\in \tilde{G}(F^{nr})$ tels que $x^{-1}X^H_{\xi}x\in \mathfrak{k}_{u}^{nr}$ est \'egal \`a $K_{s}^{\tilde{P},nr}$;

(iii)  les matrices carr\'ees  $({\bf f}^G_{\xi'}(X_{\xi}^G) )_{\xi,\xi'\in \Xi}$ et $({\bf f}_{\xi'}^H(X_{\xi}^H))_{\xi,\xi'\in \Xi}$ sont inversibles et \'egales.

Alors  on a l'\'egalit\'e $\iota_{\phi}=transfert_{\phi}$}\end{lem}

Preuve. Puisque $\tilde{G}_{AD}=G$, on a aussi $\tilde{G}_{SC}=G_{SC}$ et $\tilde{G}$ s'inscrit dans une suite d'homomorphismes $G_{SC}\to \tilde{G}\to G$. On note $K_{SC,s}^{P,nr}$ le groupe des $x\in K_{SC,s}^{0,nr}= K_{s}^{G_{SC},0,nr}$ tels que l'action r\'eduite de $Ad(x)$ dans $\mathfrak{g}_{s}$ conserve $P$. Remarquons que l'on peut remplacer $K_{SC,s}^{0,nr}$ par $K_{SC,s}^{\dag,nr}$: ces deux groupes sont \'egaux puisque $G_{SC}$ est simplement connexe. L'image r\'eciproque de $K_{s}^{\tilde{P},nr}$ dans $G_{SC}(F^{nr})$ est \'egale \`a ce groupe $K_{SC,s}^{P,nr}$. Alors l'hypoth\`ese (ii) reste v\'erifi\'ee si l'on remplace $\tilde{G}$ par $G_{SC}$. C'est-\`a-dire que  l'ensemble des \'el\'ements $x\in G_{SC}(F^{nr})$ tels que $x^{-1}X^G_{\xi}x\in \mathfrak{k}_{u}^{nr}$ est \'egal \`a $K_{s}^{G_{SC},0,nr}$
et de m\^eme pour $X^H_{\xi}$. Cela nous autorise \`a supposer que  $\tilde{G}=G_{SC}$.

 Fixons $\xi,\xi'\in \Xi$ et posons pour simplifier $f=f^G_{\xi'}$, ${\bf f}={\bf f}^G_{\xi'}$ et $X=X^G_{\xi}$.  Notons $T_{X}=G_{X}$. Avec nos d\'efinitions
de mesures, on v\'erifie que l'injection naturelle $T_{X,sc}(F)\backslash G_{SC}(F)\to T_{X}(F)\backslash G(F)$ pr\'eserve localement les mesures. Il en est de m\^eme si l'on remplace $X$ par un conjugu\'e stable de $X$. On en d\'eduit l'\'egalit\'e $S^G(X,{\bf f})=S^{G_{SC}}(X,{\bf f})$. 
  On note $\tilde{f}$ la fonction \'egale \`a $f$  sur $\mathfrak{u}_{P}({\mathbb F}_{q})=\oplus_{i\geq2}\mathfrak{g}_{s,i}({\mathbb F}_{q})$  et nulle hors de ce sous-ensemble. Notons $P_{sc}$ l'image r\'eciproque de $P$ dans $G_{SC,s}$. Il est muni  de la structure galoisienne issue de celle de $G_{SC,s}$. Pour $Y\in \mathfrak{g}_{s}({\mathbb F}_{q})$, on a l'\'egalit\'e 
    $$f(Y)=\sum_{g\in P_{sc}({\mathbb F}_{q})\backslash G_{SC,s}({\mathbb F}_{q})}\tilde{f}(gYg^{-1}).$$
     Notons  $\tilde{{\bf f}}$ la fonction sur $\mathfrak{g}(F)$ qui rel\`eve $\tilde{f}$. On obtient que l'image de ${\bf f}$ dans $I(\mathfrak{g}_{SC}(F))$ est \'egale \`a $c^G_{1}$ fois l'image de $\tilde{{\bf f}}$, o\`u 
     $$c^G_{1}=\vert G_{SC,s}({\mathbb F}_{q})\vert \vert P_{sc}({\mathbb F}_{q})\vert ^{-1} .$$ 
     D'o\`u 
    $$S^{G_{SC}}(X,{\bf f})=c_{1}^GS^{G_{SC}}(X,\tilde{{\bf f}}).$$
     L'homomorphisme $G_{SC,s}\to G_{s}$ est surjectif et son noyau est fini et central. On en d\'eduit  $\vert G_{SC,s}({\mathbb F}_{q})\vert=\vert G_{s}({\mathbb F}_{q})\vert$. D'o\`u aussi
      $$c^G_{1}=\vert G_{s}({\mathbb F}_{q})\vert \vert P_{sc}({\mathbb F}_{q})\vert ^{-1} .$$ 
          On calcule
  $$I^{G_{SC}}(X,\tilde{{\bf f}})= d^{G}(X)^{1/2}\int_{T_{X,sc}(F)\backslash G_{SC}(F)}\tilde{{\bf f}}(g^{-1}Xg)\, dg$$
  $$=d^{G}(X)^{1/2}mes(T_{X,sc}(F))^{-1}\int_{G_{SC}(F)}\tilde{{\bf f}}(g^{-1}Xg)\,dg.$$
   La fonction $\tilde{{\bf f}} $ est \`a support dans $\mathfrak{k}_{u}^{nr}\cap \mathfrak{g}(F)$.  En utilisant l'hypoth\`ese (ii) de  l'\'enonc\'e, on obtient
   $$I^{G_{SC}}(X,\tilde{{\bf f}})=d^{G}(X)^{1/2}mes(T_{X,sc}(F))^{-1}\int_{K^{P,nr}_{SC,s}\cap G_{SC}(F)}\tilde{{\bf f}}(g^{-1}Xg)\,dg.$$
   Mais la fonction $\tilde{{\bf f}}
$ est par construction invariante par conjugaison par $K^{P,nr}_{SC,s}\cap G_{SC}(F)$. D'o\`u
 $$I^{G_{SC}}(X,\tilde{{\bf f}})=c_{2}^G(X)  \tilde{{\bf f}}(X),$$
 o\`u
 $$c_{2}^G(X)=d^{G}(X)^{1/2}mes(T_{X,sc}(F))^{-1}mes(K^{P,nr}_{SC,s}\cap G_{SC}(F)).$$
  On a $mes(K_{SC,s}^{P,nr}\cap G_{SC}(F))=\vert P_{sc}({\mathbb F}_{q})\vert mes(\mathfrak{k}^{G, +}_{s})$. D'o\`u
  $$c_{2}^G(X)=d^{G}(X)^{1/2}mes(T_{X,sc}(F))^{-1} \vert P_{sc}({\mathbb F}_{q})\vert mes(\mathfrak{k}^{G, +}_{s}).$$

  Puisque $X\in \mathfrak{k}_{u}^{nr}$, on a $\tilde{{\bf f}}(X)={\bf f}(X)$. D'o\`u 
  $$I^{G_{SC}}(X,\tilde{{\bf f}})=c_{2}^G(X)  {\bf f}(X).$$
  La fonction ${\bf f}$ est cuspidale et  stable. Il en est de m\^eme de $\tilde{{\bf f}}$ puisque les images des deux fonctions dans $I(\mathfrak{g}_{SC}(F))$ sont proportionnelles. Les int\'egrales orbitales de $\tilde{{\bf f}}$ sont donc constantes sur toute classe de conjugaison stable elliptique.  
   Puisque $X$ est r\'egulier elliptique, on a donc
   $$S^{G_{SC}}(X,\tilde{{\bf f}})=c_{3}^G(X)I^{G_{SC}}(X,\tilde{{\bf f}}),$$
      o\`u $c_{3}(X)$ est le nombre de classes de conjugaison par $G_{SC}(F)$ dans la classe de conjugaison stable de $X$. 
    Finalement
    $$S^G(X,{\bf f})= c^G(X) {\bf f}(X) ,$$
    o\`u 
    $$c^G(X)=c_{1}^Gc_{2}^G(X)c_{3}^G(X)=\vert G_{s}({\mathbb F}_{q})\vert  d^{G}(X)^{1/2}mes(T_{X,sc}(F))^{-1}mes(\mathfrak{k}^{G, +}_{s})  c_{3}^G(X).$$
    
    En r\'etablissant les indices $\xi$ et $\xi'$ et les exposants $G$, on a plus pr\'ecis\'ement
     $$(1) \qquad  S^G(X^G_{\xi},{\bf f}^G_{\xi'})= c^G(X^G_{\xi}) {\bf f}^G_{\xi'}(X^G_{\xi}) .$$
     Un m\^eme calcul prouve que
     $$(2) \qquad S^H(X^H_{\xi},{\bf f}^H_{\xi'})= c^H(X^H_{\xi}) {\bf f}^H_{\xi'}(X^H_{\xi}) .$$
     
     La famille $({\bf f}^G_{\xi'})_{\xi'\in \Xi}$ est une base de $FC^{st}(\mathfrak{g}(F))$. La relation (1) et l'hypoth\`ese (iii) disant que la matrice $({\bf f}^G_{\xi'}(X_{\xi}^G) )_{\xi,\xi'\in \Xi}$ est inversible entra\^{\i}nent que les formes lin\'eaires $S^G(X^G_{\xi},.)$ s\'eparent les \'el\'ements de  $FC^{st}(\mathfrak{g}(F))$. Pour d\'emontrer la relation \ref{preparatifs}(3), 
     il suffit de prouver que, pour tout $\xi,\xi'\in \Xi$, on a l'\'egalit\'e
     $$(3) \qquad c(s)S^G(X_{\xi}^G,{\bf f}_{\xi'}^G))=S^G(X_{\xi}^G,transfert_{\phi}({\bf f}_{\xi'}^H)).$$
     Puisque les classes de conjugaison stable de $X^G_{\xi}$ et $X^H_{\xi}$ se correspondent, il r\'esulte de la d\'efinition du transfert endoscopique que
     $$S^G(X_{\xi}^G,transfert_{\phi}({\bf f}_{\xi'}^H))=S^H(X_{\xi}^H,{\bf f}_{\xi'}^H),$$
     d'o\`u, gr\^ace \`a (2):
     $$S^G(X_{\xi}^G,transfert_{\phi}({\bf f}_{\xi'}^H))=c^H(X^H_{\xi}){\bf f}^H_{\xi'}(X^H_{\xi}) .$$
       Gr\^ace \`a (1), on a aussi
     $$c(s)S^G(X_{\xi}^G,{\bf f}_{\xi'}^G))=c(s)c^G(X^G_{\xi}){\bf f}^G_{\xi'}(X^G_{\xi}) .$$
     En utilisant l'hypoth\`ese (iii), on voit que (3) r\'esulte de l'\'egalit\'e
$$(4) \qquad c(s)c^G(X^G_{\xi})=c^H(X^H_{\xi}),$$
qu'il nous reste \`a d\'emontrer. Parce que les classes stables de $X_{\xi}^G$ et $X_{\xi}^H$ se correspondent, on a $d^{G}(X_{\xi}^G)=d^{H}(X_{\xi}^H)$, $mes(T_{X^G_{\xi},sc}(F))=mes(T_{X^H_{\xi},sc}(F))$ et il  est connu que $c^G_{3}(X^G_{\xi})=c_{3}^H(X^H_{\xi})$ (en posant $T_{X_{\xi},sc}=T_{X^G_{\xi},sc}=T_{X^H_{\xi},sc}$, ces deux termes sont \'egaux au nombre d'\'el\'ements de $H^1(F,T_{X_{\xi},sc})$). Le terme $mes(\mathfrak{k}^{G, +}_{s})$ est une puissance de $q$ qui ne change pas quand on remplace $G$ par $H$.  Alors (4) r\'esulte imm\'ediatement des d\'efinitions. $\square$

Dans les paragraphes suivants, on construira cas par cas les familles $(X_{\xi}^G)_{\xi\in \Xi}$ etc... v\'erifiant les hypoth\`eses du lemme. Comme on l'a d\'ej\`a dit, dans le cas des groupes classiques, $\Xi$ est r\'eduit \`a un \'el\'ement et on supprime les $\xi$ de la notation.

\subsubsection{Preuve dans le cas $B_{n}$\label{Bn}}

 On va d\'ecrire certains groupes classiques en termes d'alg\`ebre lin\'eaire. Commen\c{c}ons par fixer les  notations pour le groupe $GL(n)$, o\`u $n$ est un entier strictement positif. On choisit pour tore maximal $T$ le tore diagonal et pour sous-groupe de Borel $B$ le groupe des matrices triangulaires sup\'erieures. Pour $i,j\in \{1,...,n\}$, on note $E_{i,j}$ la matrice dont tous les coefficients sont nuls sauf le $(i,j)$-i\`eme qui vaut $1$. Alors $(B,T,(E_{i,i+1})_{i=1,...,n-1})$ est l'\'epinglage "standard" de $GL(n)$. On note $\theta$ l'automorphisme ext\'erieur usuel de $GL(n)$:  il conserve l'\'epinglage, en envoyant $E_{i,i+1}$ sur $E_{n-i-1,n-i}$ et agit sur $Z(G)$ par $z\mapsto z^{-1}$. 
   Evidemment, si $V$ est un espace vectoriel sur $F$ de dimension $n$, muni d'une base $(e_{i})_{i=1,...,n}$ et si $F'$ est une extension de $F$,  on identifie \`a l'aide de cette base le groupe $GL_{F'}(V)$ des automorphismes $F'$-lin\'eaires de $V\otimes_{F}F'$ \`a $GL(n,F')$. Le groupe de Weyl $W$ s'identifie \`a celui des permutations des vecteurs de base, ou encore \`a celui des permutations de l'ensemble $\{1,...,n\}$.  
   
On suppose que $G$ est de type $B_{n}$ avec $n\geq2$. L'hypoth\`ese que $FC^{st}(\mathfrak{g}(F))\not=\{0\}$ entra\^{\i}ne que $2n+1=k^2+h^2$, o\`u $h,k\in {\mathbb N}$, $k$ est pair, $h$ est impair et  $h=k+\eta$, avec $\eta\in \{\pm 1\}$.  

Introduisons un espace $V$ sur $F$ de dimension $2n+1$ muni d'une base $(e_{i})_{i=1,...,2n+1}$ et d'une forme quadratique $q$ telle que $q(e_{i},e_{2n+2-i})=1$ et $q(e_{i},e_{j})=0$ pour $i+j\not=2n+2$. L'alg\`ebre $\mathfrak{g}$ s'identifie \`a l'ensemble des $X\in \mathfrak{gl}(V)$ tels que $q(Xv,v')+q(v,Xv')=0$ pour tous $v,v'\in V$.   On choisit l'\'epinglage $\mathfrak{ E}=(B,T,(E_{\alpha_{i}})_{i=1,...,n})$ de $G$ suivant: $B$ est le Borel triangulaire sup\'erieur, $T$ est le tore diagonal; pour $i=1,...,n$, $E_{\alpha_{i}}=E_{i,i+1}-E_{2n+1-i,2n+2-i}$. On le compl\`ete en un \'epinglage affine en fixant une uniformisante $\varpi_{F}$ de $F$ et en posant $E_{\alpha_{0}}=\varpi_{F}(-E_{2n,1}+E_{2n+1,2})$. Notons $\omega$ l'unique \'el\'ement non trivial de $\boldsymbol{\Omega}^{nr} $.  C'est l'\'el\'ement de $G(F)$ qui permute $e_{1}$ et $-\varpi_{F}e_{2n+1}$ et multiplie $e_{j}$ par $-1$ pour $j=2,...,2n$. Notons ${\cal S}$ le $\mathfrak{o}_{F}$-r\'eseau engendr\'e par les \'el\'ements $e_{i}$ pour $i=1,...,2n+1-k^2/2$ et $\varpi_{F}e_{i}$ pour $i=2n+2-k^2/2,...,2n+1$. Notons ${\cal S}^{\star}$ son dual, c'est-\`a-dire ${\cal S}^{\star}=\{v\in V; \forall v'\in V, \,q(v',v)\in \mathfrak{o}_{F}\}$. 
C'est le $\mathfrak{o}_{F}$-r\'eseau engendr\'e par les \'el\'ements $\varpi_{F}^{-1}e_{i}$ pour $i=1,...,k^2/2$ et $e_{i}$ pour $i=k^2/2+1,...,2n$.  L'ensemble 
$\mathfrak{k}_{s}$ est celui des éléments de $\mathfrak{g}(F)$  qui conservent ${\cal S}$ et ${\cal S}^{\star}$.   L'espace $\bar{V}'={\cal S}/\varpi_{F}{\cal S}^{\star}$ sur ${\mathbb F}_{q}$ est muni de la base $(\bar{e}_{1+k^2/2},...,\bar{e}_{2n+1-k^2/2})$, o\`u $\bar{e}_{i}$ est la r\'eduction de $e_{i}$, et de la forme quadratique $q'$ qui est la r\'eduction de $q$. L'espace $\bar{V}''={\cal S}^{\star}/{\cal S}$ est muni de la base $(\bar{e}_{2n+2-k^2/2},...,\bar{e}_{2n+1},\bar{e}_{1},...,\bar{e}_{k^2/2})$, o\`u, pour $i=2n+2-k^2/2,...,2n+1$, $\bar{e}_{i}$ est la r\'eduction de $e_{i}$ et, pour $i=1,...,k^2/2$, $\bar{e}_{i}$ est la r\'eduction de $\varpi_{F}^{-1}e_{i}$. Il est muni de la forme quadratique $q''$ qui est la r\'eduction de  $\varpi_{F}q$. 
L'alg\`ebre $\mathfrak{g}_{s}$ se d\'ecompose en $\mathfrak{g}'_{s}\oplus \mathfrak{g}''_{s}$, les facteurs \'etant les alg\`ebres de Lie des groupes sp\'eciaux orthogonaux $G'_{s}$ de $(\bar{V}',q')$ et  $G''_{s}$ de $(\bar{V}'',q'')$. On a ${\cal O}={\cal O}'\times {\cal O}''$, o\`u ${\cal O}'\subset \mathfrak{g}'_{s}$, resp. ${\cal O}''\subset \mathfrak{g}''_{s}$, est param\'etr\'ee par la partition $(2h-1,...,3,1)$, resp. $(2k-1,...,3,1)$.  

Posons $l=min(h,k)$. Pour $m\in \{1,...,l\}$, fixons un \'el\'ement $\alpha_{m}\in \bar{F}^{\times}$ de sorte que $\alpha_{m}^{4m}=2(-1)^{l+1}\varpi_{F}$. On pose $V_{m}=F(\alpha_{m})$. C'est un espace vectoriel  sur $F$ de dimension $4m$ et on le munit de la forme quadratique $q_{m}(v,v')=(-1)^{m+1}(4m)^{-1}trace_{F(\alpha_{m})/F}(\bar{v}v')$, o\`u $v\mapsto \bar{v}$ est la conjugaison galoisienne relative \`a l'extension $F(\alpha_{m})/F(\alpha_{m}^2)$. On note  ${\cal S}_{m}$ le $\mathfrak{o}_{F(\alpha_{m})}$-sous-module $\mathfrak{p}_{F(\alpha_{m})}^{m+1-l-h}$. Son dual ${\cal S}_{m}^{\star}$ est \'egal \`a $\mathfrak{p}_{F(\alpha_{m})}^{l+h-5m}$. On a $dim_{{\mathbb F}_{q}}({\cal S}_{m}/\mathfrak{p}_{F}{\cal S}_{m}^{\star})=2m-1+2(h-l)$ et $dim_{{\mathbb F}_{q}}({\cal S}_{m}^{\star}/{\cal S}_{m})=2m-1+2(k-l)$. 
On fixe une droite $V_{0}$ sur $F$, muni d'un \'el\'ement non nul $v_{0}$, et on le munit de la forme quadratique $q_{0}$ telle que $q_{0}(v_{0},v_{0})=1$ si $h>k$, $q_{0}(v_{0},v_{0})=-\frac{1}{2}\varpi_{F}$ si $h<k$. On pose ${\cal S}_{0}=\mathfrak{o}_{F}v_{0}$. L'espace quadratique $(V,q)$ est isomorphe \`a la somme directe des $(V_{m},q_{m})$ pour $m=0,...,l$. Il est un peu difficile d'\'ecrire explicitement l'isomorphisme car les $(V_{m},q_{m})$ ne sont pas d\'eploy\'es. Mais cela devient plus simple si l'on regroupe ces espaces deux par deux. Distinguons deux cas. 

Supposons d'abord $h>k$ donc $l=k$.   Pour $m=1,...,k/2$, on pose $W_{m}=V_{2m}\oplus V_{2m-1}$, ${\cal S}_{W_{m}}={\cal S}_{2m}\oplus {\cal S}_{2m-1}$,  $d'_{m}=8m$, $d''_{m}=8m-4$, $W_{0}=V_{0}$, ${\cal S}_{W_{0}}={\cal S}_{0}$.  On d\'ecompose $\{1,...,2n+1\}$ en intervalles croissants 
$$I''_{1}\cup... \cup I''_{k/2}\cup I'_{k/2}\cup...\cup I'_{1}\cup I_{0}\cup J'_{1}\cup...\cup J'_{k/2}\cup J''_{k/2}\cup...\cup J''_{1}$$
de sorte que $\vert I''_{m}\vert =\vert J''_{m}\vert =d''_{m}/2$ et $\vert I'_{m}\vert =\vert J'_{m}\vert =d'_{m}/2$ pour $m=1,...,k/2$ et $\vert I_{0}\vert =1$. On a $I_{0}=\{n+1\}$. On v\'erifie que, pour $m=1,...,k/2$, on peut identifier $W_{m}$ au sous-espace de $V$ engendr\'e par les $e_{i}$ pour $i\in I''_{m}\cup I'_{m}\cup J'_{m}\cup J''_{m}$ et que l'on peut identifier $W_{0}$  \`a la droite   engendr\'ee par  $e_{n+1}$, ces identifications \'etant telles que ${\cal S}$ soit la somme des ${\cal S}_{W_{m}}$. 

Supposons maintenant $h<k$ donc $l=h$.  Pour $m=1,...,k/2$, on pose $W_{m}=V_{2m-1}\oplus V_{2m-2}$, ${\cal S}_{W_{m}}={\cal S}_{2m-1}\oplus {\cal S}_{2m-2}$, $d'_{m}=8m-8$, $d''_{m}=8m-4$, sauf pour $m=1$ auquel cas $d'_{1}=1$. On d\'ecompose $\{1,...,2n+1\}$ en intervalles croissants $$ I''_{1}\cup... \cup I''_{k/2}\cup I'_{k/2}\cup...\cup I'_{2}\cup I_{1}\cup J'_{2}\cup...\cup J'_{k/2}\cup J''_{k/2}\cup...\cup J''_{1}$$
de sorte que $\vert I''_{m}\vert =\vert J''_{m}\vert =d''_{m}/2$ pour $m=1,...,k/2$, $\vert I'_{m}\vert =\vert J'_{m}\vert =d'_{m}/2$ pour $m=2,...,k/2$ et $\vert I_{1}\vert=d'_{1}=1$. On v\'erifie que, pour $m\geq2$, on peut identifier $W_{m}$ au sous-espace de $V$ engendr\'e par les $e_{i}$ pour $i\in I''_{m}\cup I'_{m}\cup J'_{m}\cup J''_{m}$ et que l'on peut identifier $W_{1}$ au sous-espace    engendr\'e par les $e_{i}$ pour $i\in I''_{1}\cup I_{1}\cup J''_{1}$, ces identifications \'etant telles que ${\cal S}$ soit la somme des ${\cal S}_{W_{m}}$. 

On note $X^G$ l'\'el\'ement de $\mathfrak{g}(F)$ qui conserve chaque $V_{m}$ pour $m=0,...,l$ et y agit par la multiplication par $\alpha_{m}$, avec l'exception $m=0$: $X^G$ agit par $0$ sur $V_{0}$. C'est un \'el\'ement elliptique r\'egulier, ses valeurs propres \'etant les conjugu\'es des $\alpha_{m}$, plus $0$. Il est clair qu'il conserve chaque ${\cal S}_{m}$ donc aussi ${\cal S}$, c'est-\`a-dire que $X^G\in \mathfrak{k}_{s}$. Notons $N$ l'image de $X^G$ dans $\mathfrak{g}_{s}({\mathbb F}_{q})$. Pour tout $m\in \{0,...,l\}$, l'action de $X_{G}$ sur $V_{m}$ se r\'eduit en des \'el\'ements nilpotents r\'eguliers de ${\cal S}_{m}/\mathfrak{p}_{F}{\cal S}_{m}^*$ et de ${\cal S}_{m}^*/{\cal S}_{m}$. Les partitions param\'etrant les composantes $N'\in \mathfrak{g}'_{s}$, resp. $N''\in \mathfrak{g}''_{s}$, sont donc form\'ees des dimensions des espaces ${\cal S}_{m}/\mathfrak{p}_{F}{\cal S}_{m}^*$, resp. ${\cal S}_{m}^*/{\cal S}_{m}$, c'est-\`a-dire $(2h-1,...,3,1)$ et $(2k-1,...,3,1)$. Donc 
 $N$ appartient \`a ${\cal O}$. On va prouver

(1) on peut choisir les isomorphismes ci-dessus de sorte que $N\in \mathfrak{r}({\mathbb F}_{q})$.

 On remarque qu'avec nos identifications, les vecteurs $e_{1}$ et $e_{2n+1}$ appartiennent \`a $W_{1}$. L'\'el\'ement  $\omega$  conserve chaque $W_{m}$ et, si $m\not=1$, il agit sur cet espace par multiplication par $-1$. Cette multiplication est centrale, donc la composante de $X^G$ dans $W_{m}$ pour $m\not=1$ est invariante par $Ad(\omega)$.  Donc seule importe  la composante de $X^G$ dans $W_{1}$. On peut pour simplifier supposer que $V=W_{1}$ si $h<k$, $V=W_{1}\oplus W_{0}$ si $h>k$, autrement dit que $n=2$, $k=2$ et $h=1$, ou que $n=6$, $k=2$ et $h=3$. On va pr\'eciser nos identifications dans ces cas.

Supposons $n=2$, $k=2$, $h=1$. On introduit les \'el\'ements suivants de $W_{1}=V_{0}\oplus V_{1}$:
$$e'_{1}=\frac{1}{2}\alpha_{1}^2-v_{0},\, e'_{2}=\alpha_{1}^1,\, e'_{3}=\alpha_{1}^0,\, e'_{4}=-\alpha_{1}^{-1},\, e'_{5}=\alpha_{1}^{-2}+\frac{1}{\varpi_{F}}v_{0}.$$
Ils forment une base de $W_{1}$. On identifie $W_{1}$ \`a $V$ par l'isomorphisme qui envoie $e'_{i}$ sur $e_{i}$ pour tout $i=1,...,5$. On voit qu'alors $\omega$  agit par  multiplication par $-1$ sur $V_{1}$ et par l'identit\'e de $V_{0}$. Cet \'el\'ement commute \`a $X^{G}$, c'est-\`a-dire que $Ad(\omega)(X^{G})=X^G$, a fortiori la r\'eduction $N$ de $X^G$ est fix\'ee par $Ad(\omega)_{s}$.

Supposons $n=6$, $k=2$, $h=3$. On introduit les \'el\'ements suivants de $W_{1}\oplus W_{0}=V_{0}\oplus V_{1}\oplus V_{2}$:
$$e'_{1}=\frac{1}{2}(\alpha_{2}^4-\alpha_{1}^2),\, e'_{2}=\alpha_{2}^3,\, e'_{3}=\alpha_{2}^2,\, e'_{4}=\alpha_{2}^1,\, e'_{5}=\alpha_{1}^1,\, e'_{6}=\frac{1}{2}(v_{0}-\alpha_{2}^0),\, e'_{7}=\alpha_{1}^0,$$
$$e'_{8}=v_{0}+\alpha_{2}^0,\, e'_{9}=-\alpha_{1}^{-1},\, e'_{10}=\alpha_{2}^{-1},\, e'_{11}=-\alpha_{2}^{-2},\, e'_{12}=\alpha_{2}^{-3},\, e'_{13}=-\alpha_{1}^{-2}-\alpha_{2}^{-4}.$$
Evidemment, $\alpha_{1}^0$ et $\alpha_{2}^0$ sont tous deux \'egaux \`a $1$ mais ils d\'esignent clairement l'un l'\'el\'ement  neutre de $F(\alpha_{1})^{\times}$, l'autre celui de $F(\alpha_{2})^{\times}$. 
La famille $(e'_{i})_{i=1,...,13}$ est une base de $W_{1}\oplus W_{0}$. On identifie $W_{1}\oplus W_{0}$ \`a $V$ par l'isomorphisme qui envoie $e'_{i}$ sur $e_{i}$ pour tout $i=1,...,13$. Notons $V_{1}^{\flat}$ le $F$ espace  engendr\'e par $\alpha_{1}^1,\alpha_{1}^0,\alpha_{1}^{-1}$. On voit que $\omega$ agit par multiplication par $-1$ sur $V_{0}\oplus V_{1}^{\flat}\oplus V_{2}$ et par l'identit\'e sur $F\alpha_{1}^2$. 
Cette fois, $X^G$ n'est plus invariant par $Ad(\omega)$. Mais on voit que la r\'eduction de $\omega$ agit par $-1$ sur  $\bar{V}'$ et agit sur $\bar{V}''$ par multiplication par $-1$ sur ${\cal S}_{2}^*/{\cal S}_{2}$ et par l'identit\'e sur ${\cal S}_{1}^*/{\cal S}_{1}$. Puisque $N''$ pr\'eserve ces deux sous-espaces de $\bar{V}''$, on a $Ad(\omega)_{s}(N)=N$. Cela d\'emontre (1).

On choisit l'\'el\'ement $N$ ainsi d\'efini dans la construction de \ref{preparatifs}.

Il nous reste \`a d\'efinir $X^H$. Si $h<k$, on a vu ci-dessus que $X^G$ \'etait invariant par $Ad(\omega)$. Alors $X^G$ appartient \`a $\mathfrak{h}(F)$ et on prend $X^H=X^G$.  Le terme $\epsilon$ d\'efini en \ref{preparatifs} vaut $1$ donc ${\bf f}^H(X^H)={\bf f}^G(X^G)$.   Les assertions (i) et (iii) du lemme \ref{debutdelapreuve} sont v\'erifi\'ees. L'assertion (ii) se d\'emontre comme en \cite{W7} 5.2(3) pour le groupe  $\tilde{G}=G$.   

Supposons $h>k$.   Rappelons que  $Fr_{H}=Ad(\omega)\circ Fr_{G}$.  On d\'efinit $X^H$ comme somme de composantes agissant dans $W_{m}$ pour $m=1,...,h$. Pour $m\geq2$, la composante de $X^G$ dans $W_{m}$  est invariante par $Ad(\omega)$ et on la conserve comme composante de $X^H$. Pour traiter le cas $m=1$, on simplifie les notations comme plus haut en supposant $n=6$, $k=2$ et $h=3$. Notons $\tau$ l'\'el\'ement de $G$ qui agit par multiplication par l'identit\'e sur $V_{0}\oplus   V_{2}$ et par multiplication par $-1$ sur $V_{1}$. L'automorphisme $Ad(\tau)$ fixe  $X^G$. D'apr\`es la description ci-dessus de $\omega$, on a $\omega \tau=\tau \omega$ et cet \'el\'ement agit par multiplication par $-1$ sur $V_{0}\oplus F\alpha_{1}^2\oplus V_{2}$ et par l'identit\'e sur $V_{1}^{\flat}$.  On voit que $\omega \tau$ fixe les vecteurs $e_{5},e_{7},e_{9}$ et multiplie les autres vecteurs de base par $-1$. Donc $\omega\tau$ appartient \`a $T(\mathfrak{o}_{F})$. Par le th\'eor\`eme de Lang, on peut fixer $t\in T(\mathfrak{o}_{F^{nr}})$ tel que $t^{-1}Fr_{H}(t)=\omega\tau$. Posons $X^{H}=Ad(t^{-1})(X^G)$. On a $Fr_{H}(X^H)=Ad(Fr_{H}(t))^{-1}\circ Fr_{H}(X^G)$. On a aussi $Fr_{H}(X^G)=Ad(\omega)\circ Fr_{G}(X^G)=Ad(\omega)(X^G)$. Parce que $Ad(\tau)$ fixe $X^G$, on a aussi $Fr_{H}(X^G)=Ad(\omega \tau)(X^G)=Ad(t^{-1}Fr_{H}(t))(X^G)$. Alors $Fr_{H}(X^H)=Ad(Fr_{H}(t))^{-1}\circ Ad(t^{-1}Fr_{H}(t))(X^G)=Ad(t^{-1})(X^G)=X^H$. Donc $X^H\in \mathfrak{h}(F)$. 
Puisque $X^G$ et $X^H$ sont conjugu\'es par d\'efinition, leurs classes de conjugaison stable se correspondent. L'assertion (ii) du lemme \ref{debutdelapreuve} pour $X^G$  et le groupe $\tilde{G}=G$ se d\'emontre  comme en \cite{W7} 5.2(3). Puisque $t$ appartient \`a $K_{s}^{P,nr}$, cette assertion pour $X^G$ entra\^{\i}ne la m\^eme assertion pour $X^H$. Il reste \`a d\'emontrer l'assertion (iii) du lemme \ref{debutdelapreuve}.   On a d\'efini $N$ comme la r\'eduction de $X^G$ dans $\mathfrak{g}_{s}({\mathbb F}_{q})$.  Notons $N^H$ la r\'eduction de $X^H$ dans $\mathfrak{h}_{s}({\mathbb F}_{q})$. Par d\'efinition ${\bf f}^G(X^G)=f^G(N)=1$, ${\bf f}^H(X^H)=f^H(N^H)$ et $f^H(N)=-1$ car le terme $\epsilon$ d\'efini en \ref{preparatifs} vaut $-1$.  Notons $\bar{t}$ la r\'eduction de $t$ dans $H_{s}$. On a $N^H=\bar{t}^{-1}N\bar{t}$ donc, par d\'efinition, $f^H(N^H)=\xi( \bar{z})^{-1}f^H(N)=-\xi( \bar{z})^{-1}$, o\`u $\xi$ est l'unique \'el\'ement de $\Xi$ et $\bar{z}=\bar{t}Fr_{H}(\bar{t})^{-1}$. D\'ecomposons $\bar{z}$ en $\bar{z}'\times \bar{z}''$, avec $\bar{z}'\in G'_{s}$ et $\bar{z}''\in G''_{s}$. On a d\'ecrit ci-dessus l'\'el\'ement $t^{-1}Fr_{H}(t)=\omega\tau$. On en d\'eduit que 
  $\bar{z}''=-1$ et que $\bar{z}'$ a  pour valeurs propres $+1$ avec multiplicit\'e $3$ et $-1$ avec multiplicit\'e $6$.  L'\'el\'ement $N'$ appartient \`a l'orbite nilpotente param\'etr\'ee par $(5,3,1)$.  On sait que son commutant dans le groupe orthogonal tout entier de $(\bar{V}',q')$ est produit d'un groupe unipotent et d'un groupe isomorphe \`a $({\mathbb Z}/2{\mathbb Z})^3$ engendr\'e par trois \'el\'ements $\bar{z}_{1}$, $\bar{z}_{3}$ et $\bar{z}_{5}$, chaque $\bar{z}_{i}$ poss\'edant pour valeurs propres $1$ avec multiplicit\'e $9-i$ et $-1$ avec pour multiplicit\'e $i$. Le groupe $Z_{G'_{s}}(N')$ est inclus dans le groupe pr\'ec\'edent et $\xi$ est la restriction du caract\`ere qui vaut $-1$ sur $\bar{z}_{1}$ et $\bar{z}_{5}$ et qui vaut $1$ sur $\bar{z}_{3}$. En comparant les valeurs propres, on voit que 
  $\bar{z}'$ est le produit de $\bar{z}_{1}\bar{z}_{5}$ par un \'el\'ement unipotent.   Donc $\xi(\bar{z}')=1$. Un calcul analogue montre que $\xi(\bar{z}'')=-1$. Donc $\xi(\bar{z})=-1$ et  $f^H(N^H)=f^G(N)$. Donc  ${\bf f}^H(X^H)={\bf f}^G(X^G)$, ce qui d\'emontre l'assertion (iii) du lemme
\ref{debutdelapreuve}. Cela ach\`eve la d\'emonstration. $\square$  

\subsubsection{Preuve dans le cas $C_{n}$\label{Cn}}
On suppose que $G$ est de type $C_{n}$ avec $n\geq2$. 
L'hypoth\`ese $FC^{st}(\mathfrak{g}(F))\not=\{0\}$ signifie que $n=k(k+1)$, avec $k\in {\mathbb N}$, $k\geq1$.

  Introduisons un espace $V$ sur $F$ de dimension $2n$ muni d'une base $(e_{i})_{i=1,...,2n}$ et d'une forme symplectique $q$ d\'efinie par
  $q(e_{i},e_{j})=0$ si $i+j\not=2n+1$, $q(e_{i},e_{2n+1-i})=(-1)^{i}$ si $i\in \{1,...,2n\}$. L'alg\`ebre $\mathfrak{g}$ s'identifie \`a l'ensemble des $X\in \mathfrak{gl}(V)$ tels que $q(Xv,v')+q(v,Xv')=0$ pour tous $v,v'\in V$. On choisit l'\'epinglage $\mathfrak{ E}=(B,T,(E_{\alpha_{i}})_{i=1,...,n})$ suivant: $B$ est le Borel triangulaire sup\'erieur, $T$ est le tore diagonal;  $E_{\alpha_{i}}=E_{i,i+1}+E_{2n-i,2n+1-i}$ pour $i=1,...,n-1$ et $E_{\alpha_{n}}=E_{n,n+1}$.   On le compl\`ete en un \'epinglage affine en fixant une uniformisante $\varpi_{F}$ de $F$ et en posant$E_{\alpha_{0}}=\varpi_{F}E_{2n,1}$. Notons $\omega$ l'unique \'el\'ement non trivial de $\boldsymbol{\Omega}^{nr}$. Il s'identifie \`a un \'el\'ement de $GL(V) (\bar{F})$ tel que $\omega(e_{i})=\lambda e_{n+i}$ pour $i=1,...,n$, $\omega(e_{i})=\lambda^{-1}e_{i-n}$ pour $i=n+1,...,2n$, o\`u $\lambda$ est une racine carr\'ee de $\varpi_{F}$ dans   $\bar{F}$ (rappelons que $G$ est le groupe adjoint du groupe symplectique, l'\'el\'ement ci-dessus se descend bien en un \'el\'ement de $G(F)$ qui est ind\'ependant du choix de $\lambda$). Notons ${\cal S}$ le $\mathfrak{o}_{F}$-r\'eseau engendr\'e par les \'el\'ements $e_{i}$ pour $i=1,...,3n/2$ et $\varpi_{F}e_{i}$ pour $i=3n/2+1,...,2n$. Notons ${\cal S}^{\star}$ son dual.  C'est le $\mathfrak{o}_{F}$-r\'eseau engendr\'e par les \'el\'ements $\varpi_{F}^{-1}e_{i}$ pour $i=1,...,n/2$ et $e_{i}$ pour $i=n/2+1,...,2n$.  Notons $\mathfrak{s}$ le $\mathfrak{o}_{F}$-module des \'el\'ements $X\in \mathfrak{gl}(V)(F)$ qui conservent ${\cal S}$ et ${\cal S}^{\star}$. Alors $\mathfrak{k}_{s}$ est l'intersection de $\mathfrak{g}(F)$ et de $\mathfrak{s}$. Les ${\mathbb F}_{q}$-espaces $\bar{V}'={\cal S}/\varpi_{F}{\cal S}^{\star}$ et $\bar{V}''={\cal S}^{\star}/{\cal S}$  sont de dimension $k(k+1)$ et sont naturellement munis de formes symplectiques. On en d\'eduit une d\'ecomposition $\mathfrak{g}_{s}=\mathfrak{sp}(k(k+1))\oplus \mathfrak{sp}(k(k+1))$ et ${\cal O}$ est la somme des deux orbites nilpotentes param\'etr\'ees par la partition $(2k,...,4,2)$. 
  
  L'\'egalit\'e $n=k(k+1)$ \'equivaut \`a $n=\sum_{h=1,...,k}2h$. Pour tout $h=1,...,k$, effectuons les m\^emes constructions que ci-dessus  en rempla\c{c}ant $n$ par $2h$. On affecte les objets d'un indice $h$. On identifie $V$ \`a $\oplus_{h=1,...,k}V_{h}$ de la fa\c{c}on suivante. On d\'ecompose l'ensemble $\{1,...,2n\}$ en r\'eunion croissante d'intervalles 
  $$I''_{1}\cup...\cup I''_{k}\cup I'_{k}\cup...\cup I'_{1}\cup J'_{1}\cup... \cup J'_{k}\cup J''_{k}\cup...\cup J''_{1}$$
  tels que, pour $h=1,...,k$, $\vert I''_{h}\vert =\vert I'_{h}\vert =\vert J'_{h}\vert =\vert J''_{h}\vert =h$. 
 Soit $h\in \{1,...,k\}$. Notons $a_{h}:\{1,...,4h\}\to I''_{h}\cup I'_{h}\cup J'_{h}\cup J''_{h}$ l'unique bijection croissante. Alors notre  identification envoie $e_{i,h}$ sur $(-1)^{h(h-1)/2}e_{ a_{h}(i)}$ pour $i=1,...,4h$. On en d\'eduit un plongement $\oplus_{h=1,...,k}\mathfrak{gl}(V_{h})\to \mathfrak{gl}(V)$.  
  Pour $h=1,...,k$, l'espace $\mathfrak{gl}(V_{h})$ est conserv\'e par l'action de $\omega$ et la restriction de cette action est celle de $\omega_{h}$. De m\^eme, $\mathfrak{s}\cap \mathfrak{gl}(V_{h})=\mathfrak{s}_{h}$. On peut alors construire $X^G$ et $N$ comme les sommes, en un sens \'evident, d'\'el\'ements $X_{h}$ et $N_{h}$ pour $h=1,...,k$. Fixons donc $h$. On pose $X^{h}=\sum_{j=0,...,2h}E_{\alpha_{j},h}$. Dans $\mathfrak{gl}(V_{h})$, c'est un \'el\'ement tel que $(X_{h})^{4h}=\varpi_{F}$, il est donc clairement semi-simple r\'egulier et elliptique. Il est fixe par l'action de $Ad(\omega_{h})$ et appartient \`a $\mathfrak{s}_{h}$. Sa r\'eduction est $N_{h}=\sum_{j=0,...,2h, j\not=h}\bar{E}_{\alpha_{j},h}$, avec une notation \'evidente.    Remarquons que $\mathfrak{g}_{s,h}$ s'identifie \`a $\mathfrak{sp}(2h)\oplus \mathfrak{sp}(2h)$ et $N_{h}$ est la somme de deux \'el\'ements nilpotents r\'eguliers. La somme $N=\oplus_{h=1,...,k}N_{h}$ est donc un \'el\'ement de $\mathfrak{g}_{s}({\mathbb F}_{q})\simeq \mathfrak{sp}(n)({\mathbb F}_{q})\oplus \mathfrak{sp}(n)({\mathbb F}_{q})$ dont chaque composante est param\'etr\'ee par la partition $(2k,...,4,2)$. C'est-\`a-dire que $N\in {\cal O}$. Posons $X^G=\oplus_{h=1,...,k}X_{h}$. Puisque $X^G$ est fixe par $Ad(\omega)$, $N$ est fixe par $Ad(\omega)_{s}$ donc appartient \`a $\mathfrak{r}({\mathbb F}_{q})$. On choisit cet \'el\'ement $N$ dans la construction de \ref{preparatifs}.   L'\'el\'ement $X^G$ est elliptique r\'egulier puisque les valeurs propres des $X_{h}$ sont toutes distinctes (la valuation d'une valeur propre de $X_{h}$ est $\frac{1}{4h}$). 
  L'assertion (ii) de \ref{debutdelapreuve} est v\'erifi\'ee pour $X^G$ et $\tilde{G}=G_{SC}$: la d\'emonstration est la m\^eme qu'en \cite{W7} 5.2(3), en interpr\'etant l'immeuble $Imm_{F^{nr}}(G)$ \`a l'aide de l'alg\`ebre lin\'eaire.  Puisque $X^G$ est fixe par $Ad(\omega)$, on peut choisir $X^H=X^G$, c'est bien un \'el\'ement de $\mathfrak{h}(F)$. L'assertion (i) du lemme \ref{debutdelapreuve} est claire, l'assertion (ii) pour $X^H$ est la m\^eme que pour $X^G$ et l'assertion (iii) est triviale, puisque le terme $\epsilon$ de \ref{preparatifs} vaut $1$. $\square$  
  
   \subsubsection{Preuve dans le cas $(D_{n},nr)$\label{Dnnr}}
  On suppose que $G$ est de type $D_{n}$ avec $n\geq4$, que $G$ est d\'eploy\'e sur $F^{nr}$ et que  l'action de $\Gamma_{F}^{nr}$ sur ${\cal D}$ se quotiente par un homomorphisme $\Gamma_{F}^{nr}\to \{1,\theta\}$ (pas forc\'ement surjectif).
L'hypoth\`ese $FC^{st}(\mathfrak{g}(F))\not=\{0\}$ signifie que $n=k^2$, avec $k\in {\mathbb N}$, $k$ pair.   

 Introduisons le groupe adjoint $\underline{G}$ de type $D_{n}$ et d\'eploy\'e sur $F$. Le groupe $G$ n'est pas forc\'ement isomorphe \`a $\underline{G}$, mais peut l'\^etre. Donc les objets d\'efinis pour $G$ ont aussi un sens pour ce groupe. D\'ecrivons-les. 
  Introduisons un espace $V$ sur $F$ de dimension $2n$ muni d'une base $(e_{i})_{i=1,...,2n}$ et d'une forme quadratique $q$ d\'efinie par
  $q(e_{i},e_{j})=0$ si $i+j\not=2n+1$ et $q(e_{i},e_{2n+1-i})=1$. L'alg\`ebre $\underline{\mathfrak{g}}$ s'identifie \`a l'ensemble des $X\in \mathfrak{gl}(V)$ tels que $q(Xv,v')+q(v,Xv')=0$ pour tous $v,v'\in V$. On choisit l'\'epinglage $\mathfrak{ E}=(B,T,(E_{\alpha_{i}})_{i=1,...,n})$ suivant: $B$ est le Borel triangulaire sup\'erieur, $T$ est le tore diagonal; pour $i=1,...,n-1$, $E_{\alpha_{i}}=E_{i,i+1}-E_{2n-i,2n+1-i}$  et $E_{\alpha_{n}}=E_{n-1,n+1}-E_{n,n+2}$.   On le compl\`ete en un \'epinglage affine en fixant une uniformisante $\varpi_{F}$ de $F$ et en posant$E_{\alpha_{0}}=\varpi_{F}(-E_{2n-1,1}+E_{2n,2})$. L'automorphisme $\theta$ de $\underline{G}$ est la conjugaison par  l'automorphisme de $V$, que nous notons encore $\theta$, qui permute $e_{n}$ et $e_{n+1}$ et fixe $e_{i}$ pour $i\not=n,n+1$. Notons $A$ l'automorphisme de $V$ qui permute $e_{1}$ et $\varpi_{F}e_{2n}$, $e_{n}$ et $e_{n+1}$ et qui fixe $e_{i}$ pour $i\not=1,n,n+1,2n$. Fixons un \'el\'ement $\lambda\in \bar{F}^{\times}$ tel que $\lambda^2=\varpi_{F}$. Notons $C$ l'automorphisme de $V$ d\'efini par $C(e_{i})=\lambda(-1)^{i}e_{n+i}$ pour $1\leq i\leq n$ et $C(e_{i})=\lambda^{-1}(-1)^{i+1}e_{i-n}$ pour $n+1\leq i\leq 2n$. Les automorphismes $A$ et $C$  appartiennent au groupe sp\'ecial orthogonal de $(V,q)$, leurs images dans le groupe adjoint $\underline{G}$ appartiennent \`a $\underline{G}(F)$ et ces images engendrent le groupe $\boldsymbol{\Omega}^{nr}$. 
 Notons ${\cal S}$ le $\mathfrak{o}_{F}$-r\'eseau engendr\'e par les \'el\'ements $e_{i}$ pour $i=1,...,3n/2$ et $\varpi_{F}e_{i}$ pour $i=3n/2+1,...,2n$. Son dual ${\cal S}^{\star}$ est le $\mathfrak{o}_{F}$-r\'eseau engendr\'e par les \'el\'ements $\varpi_{F}^{-1}e_{i}$ pour $i=1,...,n/2$ et $e_{i}$ pour $i=n/2+1,...,2n$.  L'ensemble $\underline{\mathfrak{k}}_{s}$ est celui des éléments de $\underline{\mathfrak{g}}(F)$  qui conservent ${\cal S}$ et ${\cal S}^{\star}$.    Les ${\mathbb F}_{q}$-espaces $\bar{V}'={\cal S}/\varpi_{F}{\cal S}^{\star}$ et $\bar{V}''={\cal S}^{\star}/{\cal S}$  sont de dimension $k^2$ et sont naturellement munis de formes quadratiques. On en d\'eduit une d\'ecomposition $\underline{\mathfrak{g}}_{s}=\mathfrak{so}(k^2)\oplus \mathfrak{so}(k^2)$ et ${\cal O}$ est la somme des deux orbites nilpotentes param\'etr\'ees par la partition $(2k-1,...,3,1)$. 
 
 Puisque $G$ est quasi-d\'eploy\'e, il y a pour $G$ deux cas possibles: ou bien $G=\underline{G}$, ou bien on peut identifier $G$ \`a $\underline{G}$ muni de l'action galoisienne telle que $\sigma_{G}=\sigma_{\underline{G}}$ pour $\sigma\in I_{F}$ et $Fr_{G}=\theta\circ Fr_{\underline{G}}$ pour tout rel\`evement $Fr$ du Frobenius dans $\Gamma_{F}$. 
 
Pour tout $m\in \{1,...,k/2\}$, introduisons un espace quadratique $(W_{m},Q_{m})$ similaire \`a $(V,q)$ mais de dimension $16m-8$. On note avec des indices $m$ les objets analogues pour cet espace \`a ceux introduits pour $(V,q)$.  D\'ecomposons l'ensemble $\{1,...,2n\}$ en r\'eunion d'intervalles croissants
  $$I''_{1}\cup...\cup I''_{k/2}\cup I'_{k/2}\cup...\cup I'_{1}\cup J'_{1}\cup...\cup J'_{k/2}\cup J''_{k/2}\cup...\cup J''_{1}$$
  de sorte que $\vert I''_{m}\vert =\vert I'_{m}\vert =\vert J'_{m}\vert =\vert J''_{m}\vert =4m-2$ pour tout $m=1,...,k/2$.  Notons $a_{m}:\{1,...,16m-8\}\to I''_{m}\cup I'_{m}\cup J'_{m}\cup J''_{m}$ l'unique bijection croissante. On identifie la somme directe des $(W_{m}, Q_{m})$ \`a $(V,q)$ en envoyant $e_{i,m}$ sur $e_{a_{m}(i)}$ pour tout $m=1,...,k/2$ et tout $i=1,...,16m-8$. 
 On en d\'eduit un plongement $\oplus_{m=1,...,k/2}\mathfrak{gl}(W_{m})\to \mathfrak{gl}(V)$. Pour tout $m=1,...,k/2$, $\mathfrak{gl}(W_{m})$ est conserv\'e par les actions de $\theta$, $A$ et $C$. Si $m=1$, ces actions co\"{\i}ncident avec celles de $\theta_{1}$, $A_{1}$ et $C_{1}$. Si $m\geq2$, les actions de $\theta$ et $A$ sur $\mathfrak{gl}(W_{m})$ sont triviales tandis que l'action de $C$ co\"{\i}ncide avec celle de $C_{m}$. On a aussi ${\cal S}\cap W_{m}={\cal S}_{m}$.

 Fixons un \'element $\mu\in \mathfrak{o}_{F}^{\times}-\mathfrak{o}_{F}^{\times,2}$ et une racine carr\'ee $\nu$ de $\mu$ dans $\mathfrak{o}_{F^{nr}}^{\times}$.  On pose $\mu_{G}=1$, $\nu_{G}=1$ si $G=\underline{G}$, $\mu_{G}=\mu$, $\nu_{G}=\nu$ si $G\not=\underline{G}$. 
 Pour tout $m\in \{1,...,k\}$, fixons un \'el\'ement $\alpha_{m}\in \bar{F}^{\times}$ de sorte que $\alpha_{m}^{4m-2}=\varpi_{F}$ pour $m\geq2$ tandis que $\alpha_{1}^2=\varpi_{F}\mu_{G}^{-1}$. 
  On pose $V_{m}=F(\alpha_{m})$. C'est un espace vectoriel sur $F$ de dimension $4m-2$ et on le munit de la forme quadratique $q_{m}(v,v')=\gamma_{m}\frac{1}{4m-2}trace_{F(\alpha_{m})/F}(\bar{v}v')$, o\`u $v\mapsto \bar{v}$ est la conjugaison galoisienne relative \`a l'extension quadratique $F(\alpha_{m})/F(\alpha_{m}^2)$ et o\`u $\gamma_{m}=\frac{(-1)^m}{2}$ pour $m\geq2$, tandis que $\gamma_{1}=- \frac{1}{2\mu_{G}}$. Notons par anticipation $(V^G,q^G)$ la somme directe des $(V_{m},q_{m})$.
    Si $G=\underline{G}$, l'espace quadratique $(V,q)$ est isomorphe sur $F$ \`a $(V^G,q^G)$.   Si $G\not=\underline{G}$, les deux espaces ne sont pas isomorphes sur $F$ mais le sont sur l'extension quadratique non ramifi\'ee $F(\nu)$ de $F$. Comme dans le cas d'un groupe de type $B_{n}$, regroupons nos espaces $(V_{m},q_{m})$ deux par deux en notant $(W'_{m},Q'_{m})$ la somme directe de $(V_{2m},q_{2m})$ et $(V_{2m-1},q_{2m-1})$. On va identifier $(W'_{m},Q'_{m})$ \`a $(W_{m},Q_{m})$. 
    
    Supposons d'abord $m\geq2$. Introduisons les \'el\'ements suivants de $W'_{m}$:
        $$(1) \qquad \left\lbrace \begin{array}{c}
   e'_{1}=\alpha_{2m}^{4m-1}+\alpha_{2m-1}^{4m-3},\, e'_{8m-4}=\alpha_{2m}^0+\alpha_{2m-1}^0, \\e'_{8m-3}=\alpha_{2m}^{0}-\alpha_{2m-1}^{0},\, e'_{16m-8}=-\alpha_{2m}^{1-4m}+\alpha_{2m-1}^{3-4m};\\
   
 \text{\,\, pour\,\,} i=2,...2m, \,e'_{i}=\alpha_{2m}^{4m-i}; \\ \text{\,\, pour\,\,} i=2m+1,...,6m-4,\, e'_{i}=\alpha_{2m-1}^{6m-3-i}; \\  \text{\,\, pour\,\,} i=6m-3,...,8m-5, \,e'_{i}=\alpha_{2m}^{8m-4-i}; \\ \text{\,\, pour\,\,} i=8m-2,...,10m-4, \,e'_{i}=(-\alpha_{2m})^{8m-3-i}; \\ \text{\,\, pour\,\,} i=10m-3,... 14m-8, \,e'_{i}=(-\alpha_{2m-1})^{10m-4-i}; \\ \text{\,\, pour\,\,} i=14m-7,...,16m-9, \,e'_{i}=(-\alpha_{2m})^{12m-7-i}.\\ \end{array}\right. $$

  Ces \'el\'ements forment une base de $W'_{m}$ et on  identifie $(W'_{m},Q'_{m})$ \`a $(W_{m},Q_{m})$ en envoyant $e'_{i}$ sur $e_{i,m}$ pour tout $i=1,...,16m-8$.  Modulo cette identification, on a ${\cal S}_{m}=\mathfrak{p}_{F(\alpha_{2m})}^{1-2m}\oplus \mathfrak{p}_{F(\alpha_{2m-1})}^{2-2m}$. On v\'erifie que l'automorphisme $C_{m}$ conserve $V_{2m}$, resp. $V_{2m-1}$, et que sa restriction \`a cet espace est la multiplication par $-\lambda \alpha_{2m}^{1-4m}$, resp. $\lambda \alpha_{2m-1}^{3-4m}$.  
  
 Supposons maintenant $m=1$. On introduit les \'el\'ements suivants de $W'_{1}$:
 $$e'_{1}=\alpha_{2}^{3}+\mu_{G}\alpha_{1}^1,\, e'_{2}=2\alpha_{2}^2,\, e'_{3}=\alpha_{2}^1,\, e'_{4}=\alpha_{2}^0+\nu_{G}\alpha_{1}^0,$$
 $$ e'_{5}=\alpha_{2}^0-\nu_{G}\alpha_{1}^0,\,e'_{6}=-2\alpha_{2}^{-1},\, e'_{7}=\alpha_{2}^{-2},\, e'_{8}=-\alpha_{2}^{-3}+\alpha_{1}^{-1}.$$ 
 On identifie comme pr\'ec\'edemment $(W'_{1},Q'_{1})$ \`a $(W_{1},Q_{1})$ \`a l'aide de cette base. Si $G=\underline{G}$, l'identification est la m\^eme que pr\'ec\'edemment et est d\'efinie sur $F$. 
 Si $G\not=\underline{G}$, cette identification est  d\'efinie sur  l'extension $F(\nu)$. Pr\'ecis\'ement, on voit que l'action sur $W'_{1}$ d'un \'el\'ement $\sigma\in \Gamma_{F}-\Gamma_{F(\nu)}$  s'identifie au compos\'e de $\theta$ et de l'action de $\sigma$ sur $W_{1}$.  On a encore  ${\cal S}_{1}=\mathfrak{p}_{F(\alpha_{2})}^{-1}\oplus \mathfrak{o}_{F(\alpha_{1})}$. L'automorphisme $A_{1}$, conserve $V_{2}$, resp. $V_{1}$, et sa restriction \`a cet espace est  la multiplication par  $-1$, resp. l'identit\'e. L'automorphisme $C_{1}$, conserve $V_{2}$, resp. $V_{1}$, et sa restriction \`a cet espace est la multiplication par $-\lambda\alpha_{2}^{-3}$, resp. $\nu_{G}^{-1}\lambda\alpha_{1}^{-1}$.

 Modulo ces isomorphismes, on a identifi\'e $(V,q)$ \`a  $(V^G,q^G)$. Si $G=\underline{G}$, cet isomorphisme est d\'efini sur $F$, autrement dit $G$ s'identifie au groupe adjoint du groupe sp\'ecial orthogonal de $(V^G,q^G)$. Si $G\not=\underline{G}$, cet isomorphisme est d\'efini sur  $F(\nu)$ et l'action sur $V^G$ d'un \'el\'ement $\sigma\in \Gamma_{F}-\Gamma_{F(\nu)}$  s'identifie au compos\'e de $\theta$ et de l'action de $\sigma$ sur $V$. De nouveau,   $G$ s'identifie au groupe adjoint du groupe sp\'ecial orthogonal de $(V^G,q^G)$.
  Notons $X^G$ l'\'el\'ement de $\underline{\mathfrak{g}}(F^{nr})$ qui pr\'eserve chaque $V_{m}$ et y agit par multiplication par $\alpha_{m}$. D'apr\`es ce que l'on vient de dire, on a $X^G\in \mathfrak{g}(F)$. Il pr\'eserve ${\cal S}$ et ${\cal S}^*$ donc appartient \`a $\mathfrak{k}_{s}$. D'apr\`es la description ci-dessus des automorphismes $A$ et $C$, $X^G$ est fix\'e par ces automorphismes. On d\'efinit $N$ comme la r\'eduction de $X^G$ dans $\mathfrak{g}_{s}({\mathbb F}_{q})$.  Alors $N$ est fix\'e par $A_{s}$ et $C_{s}$ donc appartient \`a $\mathfrak{r}({\mathbb F}_{q})$. On voit comme  en \ref{Bn} que $N\in {\cal O}$. On utilise cet \'el\'ement dans la construction de \ref{preparatifs}. On pose $X^H=X^G$. Puisque $X^G$ est fix\'e par $A$ et $C$, donc par $\boldsymbol{\Omega}^{nr}$, $X^H$ appartient \`a $\mathfrak{h}(F)$. 
   L'hypoth\`ese (i) du lemme \ref{debutdelapreuve} est claire et  l'hypoth\`ese (iii) aussi puisque  le terme $\epsilon$ de \ref{preparatifs} vaut $1$. On d\'emontre comme en \cite{W7} 5.2(3)   l'hypoth\`ese (ii)  pour le groupe $\tilde{G}$ \'egal au groupe sp\'ecial orthogonal de $(V,q)$.  $\square$

\subsubsection{Preuve dans le cas $(D_{n},ram)$}\label{Dnram}
Soit $E$ une extension quadratique ramifiée de $F$. On suppose que $G$ est de type $D_{n}$ avec $n\geq4$, que l'action de $\Gamma_{F}$ sur ${\cal D}$ a pour noyau $\Gamma_{E}$ et que cette action se quotiente en une  bijection $\Gamma_{E/F}\to\{1,\theta\}$. L'hypothèse $FC^{st}(\mathfrak{g}(F))\not=\{0\}$ signifie que $n=h^2$ avec $h\in {\mathbb N}$, $h$ impair. On fixe une uniformisante $\varpi_{E}$ de $E$ telle que $\varpi_{E}^2\in F^{\times}$ et on pose $\varpi_{F}=\varpi_{E}^2$. 

Introduisons comme en \ref{Dnnr} l'espace quadratique $(V,q)$, le groupe adjoint $\underline{G}$ du groupe spécial orthogonal de la forme $q$ et son épinglage $\mathfrak{E}$. Le groupe $G$ s'identifie à $\underline{G}$ sur $E$ mais, pour $\sigma\in \Gamma_{F}-\Gamma_{E}$, on a l'égalité $\sigma_{G}=\theta\circ \sigma_{\underline{G}}$. 
Pour cette structure, un épinglage affine  est formé des vecteurs suivants

$E_{\beta_{i}}=E_{\alpha_{i}}=E_{i,i+1}-E_{2n-i,2n+1-i}$ pour $i=1,...,n-2$;

$E_{\beta_{n-1,n}}=E_{\alpha_{n-1}}+E_{\alpha_{n}}=E_{n-1,n}+E_{n-1,n+1}-E_{n,n+2}-E_{n+1,n+2}$;

$E_{\beta_{0}}=\varpi_{E}(E_{n,1}-E_{n+1,1}+E_{2n,n}-E_{2n,n+1})$. 

Posons $v_{0}=\varpi_{E}(e_{n}-e_{n+1})$, $v_{1}=e_{n}+e_{n+1}$. Notons $V^G$ le sous-$F$-espace vectoriel de $V\otimes_{F}\bar{F}$ engendré par $v_{0}$, $v_{1}$ et les $e_{i}$ pour $i\in \{1,...,n-1\}\cup\{n+2,...,2n\}$. La forme $q$ se prolonge $\bar{F}$-linéairement à $V\otimes_{F}\bar{F}$ puis se restreint en une forme quadratique $q^G$ sur $V^G$ qui prend ses valeurs dans $F$. Le groupe $G$ s'identifie au groupe adjoint du groupe spécial orthogonal de $(V^G,q^G)$.

Fixons une racine primitive quatrième
de l'unité dans $\bar{F}^{\times}$, que l'on note $\zeta$. Posons $\lambda=\zeta\varpi_{E}$ et définissons l'automorphisme $C$ de $V$ par les formules

$C(e_{i})=\lambda(-1)^{i}e_{n+1+i}$ pour $i=1,...,n-1$, $C(e_{i})=\lambda^{-1}(-1)^{i+1}e_{i-n-1}$ pour $i=n+2,...,2n$;

$C(e_{n})=\zeta e_{n}$, $C(e_{n+1})=-\zeta e_{n+1}$. 

Il appartient au groupe spécial orthogonal de $V$  et son image  dans $G$ appartient à $G(F^{nr})$. On note encore $C$ cette image. On vérifie qu'il conserve l'épinglage affine. Alors $\boldsymbol{\Omega}^{nr}=\{1,C\}$.

Notons ${\cal S}$ le $\mathfrak{o}_{F}$-réseau de $V^G$ engendré par $v_{0}$, $v_{1}$, les $e_{i}$ pour $i\in \{1,...,n-1\}\cup\{n+2,...,(3n+1)/2\}$ et les $\varpi_{F}e_{i}$ pour $i\in\{(3n+3)/2,...,2n\}$. Son dual est le $\mathfrak{o}_{F}$-réseau engendré par $\varpi_{F}^{-1}v_{0}$, $v_{1}$, les $\varpi_{F}^{-1}e_{i}$ pour $i\in\{1,...,(n-1)/2\}$ et  les $e_{i}$ pour $i\in \{(n+1)/2,...,n-1\}\cup\{n+2,...,2n\}$. L'algèbre $\mathfrak{k}_{s}$ est l'ensemble des éléments de $\mathfrak{g}(F)$ qui conservent ${\cal S}$ et ${\cal S}^*$. Les ${\mathbb F}_{q}$-espaces $\bar{V}'={\cal S}/\varpi_{F}{\cal S}^*$ et $\bar{V}''={\cal S}^*/{\cal S}$ sont de dimension $h^2$ et sont naturellement munis de formes quadratiques. D'où une décomposition $\mathfrak{g}_{s}=\mathfrak{so}(h^2)\oplus \mathfrak{so}(h^2)$ et ${\cal O}$ est la somme des deux orbites nilpotentes paramétrées par la partition $(2h-1,...,3,1)$.

Pour tout $m\in \{1,...,h\}$, fixons un élément $\alpha_{m}\in \bar{F}^{\times}$ tel que $\alpha_{m}^{4m-2}=\varpi_{F}$.  On pose $V_{m}=F(\alpha_{m})$. C'est un espace vectoriel sur $F$ de dimension $4m-2$ et on le munit de la forme quadratique $q_{m}(v,v')=(-1)^{m+1}\frac{1}{8m-4}trace_{F(\alpha_{m})/F}(\bar{v}v')$, où $v\mapsto \bar{v}$ est la conjugaison galoisienne relative à l'extension quadratique $F(\alpha_{m})/F(\alpha_{m}^2)$. On va définir un isomorphisme d'espaces quadratiques de $\oplus_{m=1,...,h}(V_{m},q_{m})$ sur $(V^G,q^G)$.  Pour $m\in \{1,...,(h-1)/2\}$, notons $(W_{m},Q_{m})$ la somme directe de $(V_{2m},q_{2m})$ et de $(V_{2m+1},q_{2m+1})$. Introduisons les éléments $(e'_{i})_{i=1,...16m}$ de $W_{m}$  définis par les formules \ref{Dnnr} (1) où l'on remplace $m$ par $m+\frac{1}{2}$. Décomposons l'ensemble $\{1,...,n-1\}\cup \{n+2,...,2n\}$ en réunion croissante
$$I''_{1}\sqcup...\sqcup I''_{(h-1)/2}\sqcup I'_{(h-1)/2}\sqcup ...\sqcup I'_{1}\sqcup J'_{1}\sqcup...\sqcup J'_{(h-1)/2}\sqcup J''_{(h-1)/2}\sqcup ...\sqcup J''_{1}$$
de sorte que $\vert I''_{m}\vert=\vert I'_{m}\vert=\vert J'_{m}\vert=\vert J''_{m}\vert=4m$ pour tout $m=1,...,(h-1)/2$. Notons $a_{m}:\{1,...,16m\}\to I''_{m}\sqcup
I'_{m}\sqcup J'_{m}\sqcup J''_{m}$ l'unique bijection croissante. On identifie $W_{m}$ à un sous-espace de $V^G$ en identifiant $e'_{i}$ à $e_{a_{m}(i)}$ pour tout $i\in \{1,...,16m\}$. Il reste l'espace $V_{1}$ que l'on identifie au sous-espace de $V^G$ engendré par $v_{0}$ et $v_{1}$ en envoyant $\alpha_{1}^0$ sur $\frac{1}{2}v_{1}$ et $\alpha_{1}^1$ sur $-\frac{1}{2}v_{0}$. On a alors identifié $\oplus_{m=1,...,h}V_{m}$ à $V^G$ et on vérifie que cet isomorphisme transporte la forme $\oplus_{m=1,...,h}V_{m}$  en $q^G$. Modulo cet isomorphisme, ${\cal S}$ s'identifie à $\oplus_{m=1,...,h} \mathfrak{p}_{F(\alpha_{m})}^{1-m}$. On vérifie que l'automorphisme $C$ conserve chaque $V_{m}\otimes_{F}\bar{F}$ et que sa restriction à cet espace est la multiplication par $\lambda (-1)^m\alpha_{m}^{1-2m}$. 

Notons $X^G$ l'élément de $\mathfrak{g}(F)$ qui préserve chaque $V_{m}$ et y agit par multiplication par $\alpha_{m}$. Il préserve ${\cal S}$ et ${\cal S}^*$ donc appartient à $\mathfrak{k}_{s}$. D'après la description ci-dessus, il est fixé par $C$. On définit $N$ comme la réduction de $X^G$ dans $\mathfrak{g}_{s}({\mathbb F}_{q})$. Alors $N$ est fixé par $C_{s}$ donc appartient à $\mathfrak{r}({\mathbb F}_{q})$.  On voit comme  en \ref{Bn} que $N\in {\cal O}$. On utilise cet \'el\'ement dans la construction de \ref{preparatifs}. On pose $X^H=X^G$. Puisque $X^G$ est fix\'e par $C$, donc par $\boldsymbol{\Omega}^{nr}$, $X^H$ appartient \`a $\mathfrak{h}(F)$. 
   L'hypoth\`ese (i) du lemme \ref{debutdelapreuve} est claire et  l'hypoth\`ese (iii) aussi puisque  le terme $\epsilon$ de \ref{preparatifs} vaut $1$. On d\'emontre comme en \cite{W7} 5.2(3)   l'hypoth\`ese (ii)  pour le groupe $\tilde{G}$ \'egal au groupe sp\'ecial orthogonal de $(V^G,q^G)$.  $\square$

 \subsubsection{Preuve dans le  cas  o\`u $G$ est de type $E_{6}$  déployé ou $E_{7}$\label{exceptionnels}}
   Supposons que $G$ est déployé de type $E_{6}$.   L'hypoth\`ese $FC^{st}(\mathfrak{g}(F))\not=\{0\}$ \'equivaut \`a $\delta_{3}(q-1)=1$.  Pour simplifier la notation, on note $E_{i}=E_{\alpha_{i}}$ les \'el\'ements de l'\'epinglage affine. L'unique \'el\'ement non trivial $\theta$ de $\boldsymbol{Out}(G)$ fixe $E_{0}$, $E_{2}$ et $E_{4}$ et \'echange $E_{1}$ et $E_{6}$ ainsi que $E_{3}$ et $E_{5}$. Le groupe $\boldsymbol{\Omega}^{nr}$ a trois \'el\'ements et est engendr\'e par l'\'el\'ement qui fixe $E_{4}$ et agit  sur les autres \'el\'ements de l'\'epinglage affine par les    permutations $E_{1}\mapsto E_{6}\mapsto E_{0}\mapsto E_{1}$, $E_{3}\mapsto E_{5}\mapsto E_{2}\mapsto E_{3}$. Pour fixer la notation, supposons que l'\'el\'ement $\omega$  qui d\'efinit $H$  soit cet \'el\'ement.  
     
       Le groupe $G_{s,SC}$ s'identifie \`a $SL(3)^3$, les trois copies de $SL(3)$ ayant respectivement pour \'epinglages $\{\bar{E}_{1},\bar{E}_{3}\}$, $\{\bar{E}_{6},\bar{E}_{5}\}$, $\{\bar{E}_{0},\bar{E}_{2}\}$, o\`u $\bar{E}_{i}$ est la r\'eduction de $E_{i}$ dans $\mathfrak{g}_{s}$. Plus pr\'ecis\'ement, on identifie $G_{s,SC}$ \`a $SL(3)^3$ de sorte que l'action d\'eduite de $Ad(\omega)$ soit $(g_{1},g_{2},g_{3})\mapsto (g_{3},g_{1},g_{2})$.    L'orbite ${\cal O}$ est le produit des orbites r\'eguli\`eres. On choisit $N=\bar{E}_{1}+\bar{E}_{3}+\bar{E}_{6}+\bar{E}_{5}+\bar{E}_{0}+\bar{E}_{2}$. Il appartient \`a ${\cal O}$ et est fix\'e par $Ad(\omega')_{s}$ tout $\omega'\in \boldsymbol{\Omega}^{nr}$, donc appartient \`a ${\cal O}\cap \mathfrak{r}({\mathbb F}_{q})$.   On peut supposer que le   $\mathfrak{sl}(2)$-triplet  fix\'e en \ref{preparatifs} a pour \'el\'ement  semi-simple l'unique \'el\'ement possible appartenant \`a $\mathfrak{t}_{s}({\mathbb F}_{q})$ (on rappelle que  $T_{s}$ est  le sous-tore maximal de $G_{s}$ d\'eduit de $T$).   Notons $\bar{\boldsymbol{\zeta}}_{1/3}$ le groupe des racines cubiques de l'unit\'e dans $\bar{{\mathbb F}}_{q}^{\times}$.  Il est contenu dans ${\mathbb F}_{q}^{\times}$   d'apr\`es l'hypoth\`ese $\delta_{3}(q-1)=1$. Le groupe $Z_{G_{s,SC}}(N)/Z_{G_{s,SC}}(N)^0$ s'identifie \`a $(\bar{\boldsymbol{\zeta}}_{1/3})^3$. Le groupe $Z_{G_{s}}(N)/Z_{G_{s}}(N)^0$ est un quotient du groupe pr\'ec\'edent et l'ensemble $\Xi$ se rel\`eve en celui des caract\`eres de    $Z_{G_{s,SC}}(N)/Z_{G_{s,SC}}(N)^0$ de la forme $(z_{1},z_{2},z_{3})\mapsto \xi(z_{1}z_{2}z_{3})$ o\`u $\xi$ est un caract\`ere non trivial de $\bar{\boldsymbol{\zeta}}_{1/3}$. On identifie ainsi $\Xi$ \`a l'ensemble \`a deux \'el\'ements des caract\`eres non triviaux de $\bar{\boldsymbol{\zeta}}_{1/3}$. 
       
          Le groupe $P$ est simplement le groupe de Borel de $G_{s}$ déduit de l'alc\^ove $C^{nr}$, $K_{SC,s}^{P,nr}$ est un sous-groupe d'Iwahori  de $G_{SC}(F^{nr})$ et $ \mathfrak{k}_{u}^{nr}$ est le premier terme de la filtration usuelle d'une alg\`ebre d'Iwahori.

        Notons $E$ l'extension non ramifi\'ee de degr\'e $3$ de $F$, fixons un \'el\'ement $y\in \mathfrak{o}_{E}^{\times}$, posons $x=Norm_{E/F}(y)$. Posons
    $$X^G=E_{1}+xE_{3}+E_{4}+E_{5}+E_{6}+E_{2}+E_{0},\,\, X^H=E_{1}+yE_{3}+E_{4}+Fr(y)E_{5}+E_{6}+Fr^2(y)E_{2}+E_{0}.$$
    On a $X^G\in \mathfrak{g}(F)$ et $X^H\in \mathfrak{h}(F)$. 
   \begin{lem}{(i) L'\'el\'ement $X^G$, resp. $X^H$, appartient \`a $\mathfrak{k}_{u}^{nr}$ et est semi-simple r\'egulier et elliptique. Il reste elliptique dans $\mathfrak{g}(F^{nr})$, resp. $\mathfrak{h}(F^{nr})$. 
   
   (ii) L'ensemble des $g\in G_{SC}(F^{nr})$ tels que $g^{-1}X^Gg\in \mathfrak{k}_{u}^{nr}$, resp. $g^{-1}X^H g\in \mathfrak{k}_{u}^{nr}$, est \'egal \`a $K_{SC,s}^{P,nr}$.}\end{lem}
   
   Preuve. Que $X^G$ appartienne \`a $\mathfrak{k}_{u}^{nr}$ est clair.  Fixons une uniformisante $\varpi_{F}$ de $F$.  Rappelons que  $s_{0}$ est le sommet appartenant \`a $S(\bar{C}^{nr})$ associ\'e \`a la racine $\alpha_{0}$.   Posons $E^{\star}_{0}=\varpi_{F}^{-1}E_{0}$. Cet \'el\'ement appartient \`a $\mathfrak{k}_{s_{0}}$ et  sa r\'eduction $\bar{E}^{\star}_{0}$ dans $\mathfrak{g}_{s_{0}}({\mathbb F}_{q})$  est un \'el\'ement non nul du sous-espace radiciel  de cette alg\`ebre associ\'e \`a la racine $\alpha_{0}$. Fixons $\lambda\in \bar{F}^{\times}$ tel que $\lambda^{12}=\varpi_{F}$. Posons $t=\prod_{i=1,...,6}\check{\varpi}_{i}(\lambda)$. On calcule $tX^Gt^{-1}=\lambda Y$, o\`u $Y=E^{\star}_{0}+xE_{3}+\sum_{i=1,...,6,i\not=3}E_{i}$. L'\'el\'ement $Y$ appartient \`a $\mathfrak{k}_{s_{0}}$ et sa r\'eduction $\bar{Y}$ dans $\mathfrak{g}_{s_{0}}({\mathbb F}_{q})$ est un \'el\'ement cyclique selon la d\'efinition de \cite{Kos} 6.2. Un tel \'el\'ement est semi-simple r\'egulier, cf. \cite{Kos} lemme 6.3 et corollaire 6.4. On peut choisir un \'el\'ement $\bar{g}\in G_{s_{0}}$ de sorte que $\bar{g}\bar{Y}\bar{g}^{-1}\in \mathfrak{t}_{s_{0}}$. Posons $\bar{Z}= \bar{g}\bar{Y}\bar{g}^{-1}$. Un raisonnement standard permet de relever $\bar{g}$ en un \'el\'ement $g\in K^{0,nr}_{s_{0}}$ de sorte que $gYg^{-1}\in \mathfrak{t}(F^{nr})$. L'\'el\'ement $\bar{Y}$ \'etant semi-simple r\'egulier, $Y$ et  $X^G$ le sont aussi. Notons $T_{X}$ le commutant de $X^G$ dans $G$.  On a $Ad(gt)(T_{X})=T$. L'automorphisme $Ad(gt)$ identifie $T_{X}$ \`a $T$ muni de l'action galoisienne $\sigma\mapsto Ad( (gt)^{-1}\sigma_{G}(gt))\circ \sigma_{G}$. Puisque $T$ est d\'eploy\'e, $X_{*}(T_{X})$ s'identifie \`a $X_{*}(T)$ muni de l'action galoisienne $\sigma\mapsto w_{\sigma}$, o\`u $w_{\sigma}$ est l'image dans le groupe de Weyl  de $(gt)^{-1}\sigma_{G}(gt)$. Pour prouver que $X^G$ est elliptique dans $\mathfrak{g}(F^{nr})$, il suffit de prouver qu'il existe $\sigma\in I_{F}$ tel que $w_{\sigma}$ agisse sans point fixe non nul dans $X_{*}(T)$. Fixons une racine primitive $\zeta$ d'ordre $12$ de l'unit\'e dans $\bar{F}^{\times}$. On peut choisir $\sigma\in I_{F}$ de sorte que $\sigma(\lambda)=\zeta\lambda$. On a alors $(gt)^{-1}\sigma_{G}(gt)=t^{-1}\sigma_{G}(t)=\prod_{i=1,...,6}\check{\varpi}_{i}(\zeta)$. On a  $Ad(t^{-1}\sigma_{G}(t))(Y)=\zeta Y$. L'\'el\'ement $t^{-1}\sigma_{G}(t)$ appartient \`a $K^{0,nr}_{s_{0}}$. Notons $\bar{t}_{\sigma}$ sa r\'eduction dans $G_{s_{0}}$. L'action de $w_{\sigma}$ sur $X_{*}(T)$ s'identifie \`a l'action de $\bar{g}^{-1}\bar{t}_{\sigma}\bar{g}$ sur $X_{*}(T_{s_{0}})$. Or $Ad(\bar{g}\bar{t}_{\sigma}\bar{g}^{-1})(\bar{Z})=Ad(\bar{g}\bar{t}_{\sigma})(\bar{Y})=\zeta Ad(\bar{g})(\bar{Y})=\zeta \bar{Z}$. Le th\'eor\`eme 9.2 de \cite{Kos} implique que $w_{\sigma}$ est une transformation de Coxeter. Une telle transformation n'a pas de point fixe non nul dans $X_{*}(T)$, cf. \cite{Kos} lemme 8.1. Cela prouve le (i) de l'\'enonc\'e pour l'\'el\'ement $X^G$. 
   
   En \ref{racinesaffines}, on a identifi\'e l'appartement $App_{F^{nr}}(T)$ \`a $ X_{*}(T)\otimes_{{\mathbb Z}}{\mathbb R}$. En particulier,   $s_{0}$ s'identifie \`a $0$. Soit $z$ le point de $App_{F^{nr}}(T)$ tel que $\alpha_{i}(z)=\frac{1}{12}$  pour  tout $i=1,...,6$. Il r\'esulte de \cite{W7} lemme 8.2 que $\mathfrak{k}_{u}^{nr}=\mathfrak{k}^{nr}_{z,\frac{1}{12}}$. On v\'erifie que $K_{SC,s}^{P,nr}=K_{SC,z}^{0,nr}$. Alors l'assertion (ii) de l'\'enonc\'e pour l'\'el\'ement $X^G$ se d\'emontre comme en \cite{W7} 8.4. R\'ep\'etons l'argument. Posons $F'=F^{nr}(\lambda)$.  Alors $\lambda^{-1}\mathfrak{k}_{u}^{nr}=\lambda^{-1}\mathfrak{k}_{z,\frac{1}{12}}^{nr}\subset \lambda^{-1}\mathfrak{k}_{z,\frac{1}{12},F'}=\mathfrak{k}_{z,0,F'}$. 
   Soit $g\in G_{SC}(F^{nr})$ tel que $g^{-1}X^Gg\in \mathfrak{k}_{u}^{nr}$. Alors $\lambda^{-1}g^{-1}X^Gg\in \mathfrak{k}_{z,0,F'}$ d'o\`u  
     $\lambda^{-1}X^G\in \mathfrak{k}_{gz,0,F'}$. D'apr\`es la preuve de (i), $\lambda^{-1}X^G$ est un \'el\'ement entier et de reduction r\'eguli\`ere de $\mathfrak{t}_{X}(F')$. Il r\'esulte du lemme \cite{W7} 8.1 que $gz$ est un point de l'appartement associ\'e \`a $T_{X}$ dans l'immeuble $Imm_{F'}(G)$ (la preuve de (i) montre que $T_{X}$ est d\'eploy\'e sur $F'$). Mais il est aussi fixe par $I_{F}$. Puisque $T_{X}$ est elliptique sur $F^{nr}$, l'appartement associ\'e \`a $T_{X}$ poss\`ede un unique point fixe par $I_{F}$, qui n'est autre que $z$. Cela prouve $gz=z$, donc $g\in K_{SC,z}^{0,nr}=K_{SC,s}^{P,nr}$. Cela ach\`eve la preuve du lemme pour l'\'el\'ement $X^G$. Evidemment, la m\^eme preuve s'applique \`a $X^H$. $\square$
   
   Posons $u=\check{\varpi}_{3}(yx^{-1})\check{\varpi}_{5}(Fr(y))\check{\varpi}_{2}(Fr^2(y))$. On constate que $Ad(u)(X^G)=X^H$. Donc 
   
   (1) la classe de conjugaison stable de $X^H$ dans $H$ correspond \`a celle de $X^G$ dans $G$. 
   
   Notons $\boldsymbol{\zeta}_{1/3}$ le groupe des racines cubiques de l'unit\'e dans $\bar{F}^{\times}$. Il s'identife \`a $\bar{\boldsymbol{\zeta}}_{1/3}$ par l'application de  r\'eduction $z\mapsto \bar{z}$.   Notons ${\cal Z}= 
 \boldsymbol{\zeta}_{1/3}-\{1\}$. Les ensembles ${\cal Z}$ et $\Xi$ ont chacun  deux \'el\'ements.   Il est clair que
    
   (2) la matrice carr\'ee $(\xi(\bar{\zeta}))_{\xi\in \Xi,\zeta\in {\cal Z}}$ est inversible.
   
     Pour $\zeta\in {\cal Z}$, fixons un \'el\'ement $b_{\zeta}\in F^{nr,\times}$ tel que $b_{\zeta}^{-1}Fr^3(b_{\zeta})=\zeta^{-1}$ et posons $y_{\zeta}=b_{\zeta}^3$.   C'est un \'el\'ement de $E$. On pose $x_{\zeta}=Norme_{E/F}(y_{\zeta})$. On note 
   $X_{\zeta}^G$ et $X_{\zeta}^H$ les \'el\'ements $X^G$ et $X^H$ construits ci-dessus \`a l'aide de $x_{\zeta}$ et $y_{\zeta}$. L'hypoth\`ese (i) du lemme \ref{debutdelapreuve} r\'esulte de (1)  et du (i) du lemme ci-dessus tandis que l'hypoth\`ese (ii) du lemme \ref{debutdelapreuve} r\'esulte du (ii) du lemme ci-dessus. En vertu de (2) ci-dessus, l'hypoth\`ese (iii) du lemme \ref{debutdelapreuve} r\'esulte de
   
   (3) pour tous $\xi\in \Xi$ et $\zeta\in {\cal Z}$, on a les \'egalit\'es ${\bf f}^G_{\xi}(X^G_{\zeta})={\bf f}^H_{\xi}(X^H_{\zeta})=\xi(\bar{\zeta})$.
   
 Prouvons (3). L'\'el\'ement $\epsilon$ de \ref{preparatifs} vaut $1$. En cons\'equence $f_{\xi}^G(N)=f_{\xi}^H(N)=1$.  Posons simplement $x=x_{\zeta}$, $y=y_{\zeta}$, $b=b_{\zeta}$. Posons $a=bFr(b)Fr^2(b)$. Il r\'esulte des d\'efinitions que $x=a^3$ et que $a^{-1}Fr(a)=\zeta^{-1}$.  La r\'eduction $\bar{X}^G_{\zeta} $ dans $\mathfrak{g}_{s}({\mathbb F}_{q})$ est \'egale \`a $\bar{E}_{1}+\bar{a}^3\bar{E}_{3}+\bar{E}_{5}+\bar{E}_{6}+\bar{E}_{2}+\bar{E}_{0}$.  Notons $t_{1}$ l'\'el\'ement diagonal de $SL(3)(\bar{{\mathbb F}}_{q})$ de coefficients diagonaux $(\bar{a}^{-1},\bar{a}^{-1},\bar{a}^2)$. Posons $t=(t_{1},1,1)$. On v\'erifie que $t^{-1}Nt=\bar{X}_{\zeta}^G$. On a $tFr_{G}(t)^{-1}=(\bar{a}^{-1}Fr(\bar{a}),1,1)=(\bar{\zeta}^{-1},1,1)$. On a donc $\xi(Fr_{G}(t)t^{-1})=\xi(\bar{\zeta})$, d'o\`u 
 $${\bf f}^G_{\xi}(X^G_{\zeta})=f^G_{\xi}(\bar{X}^G_{\zeta})=\xi(\bar{\zeta})f^G_{\xi}(N)=\xi(\bar{\zeta}).$$

La r\'eduction $\bar{X}^H_{\zeta}$ dans $\mathfrak{h}_{s}({\mathbb F}_{q})$ est \'egale \`a  $\bar{E}_{1}+\bar{y}\bar{E}_{3}+Fr(\bar{y})\bar{E}_{5}+\bar{E}_{6}+Fr^2(\bar{y})\bar{E}_{2}+\bar{E}_{0}$. Notons $t'_{1}$, $t'_{2}$ et $t'_{3}$ les \'el\'ements diagonaux de $SL(3)(\bar{{\mathbb F}}_{q})$ de coefficients diagonaux $\bar{(b}^{-1},\bar{b}^{-1},\bar{b}^2)$, resp. $(Fr(\bar{b})^{-1},Fr(\bar{b})^{-1},Fr(\bar{b}^2))$, $(Fr^2(\bar{b})^{-1},Fr^2(\bar{b})^{-1},Fr^2(\bar{b}^2))$. Posons $t'=(t'_{1},t'_{2},t'_{3})$. On v\'erifie que $t^{'-1}Nt'=\bar{X}^H_{\zeta}$. On a $t'Fr_{H}(t')^{-1}= (\bar{b}^{-1}Fr^3(\bar{b}),1,1)=(\bar{\zeta}^{-1},1,1)$.  Donc, comme ci-dessus,
 $${\bf f}^H_{\xi}(X^H_{\zeta})=f^H_{\xi}(\bar{X}^H_{\zeta})=\xi(\bar{\zeta})f^H_{\xi}(N)=\xi(\bar{\zeta}).$$
 Cela d\'emontre (3) et ach\`eve la preuve dans le cas o\`u $G$ est de type $E_{6}$.

Une preuve similaire s'applique quand $G$ est de type $E_{7}$. $\square$

\subsubsection{Preuve dans le cas o\`u $G$ est   de type $(A_{n-1},ram)$ avec $n$ pair\label{An-1ram}}
 On suppose que $G$ est de type $A_{n-1}$ avec $n$ pair et $n\geq 4$, que $G$ n'est pas d\'eploy\'e sur $F$ mais l'est sur une extension quadratique  ramifi\'ee $E$ de $F$. L'hypoth\`ese que $FC^{st}(\mathfrak{g}(F))\not=\{0\}$ \'equivaut \`a ce  que $n=h^2+k(k+1)$ avec $h,k\in {\mathbb N}$ et $h=k$ ou $h=k+1$.  Remarquons que la parit\'e de $n$ entra\^{\i}ne celle de $h$. 
 Le groupe $G$ est le groupe adjoint d'un groupe unitaire que l'on note $\tilde{G}$. On affecte d'un $\tilde{}$ les objets relatifs \`a $\tilde{G}$. 

On fixe une uniformisante $\varpi_{E}$ de $E$ telle que $\varpi_{E}^2\in F^{\times}$. On introduit un espace $V$ sur $E$ de dimension $n$ muni d'une base $(e_{i})_{i=1,...,n}$ et d'une forme hermitienne $q$ (relativement \`a l'extension $E/F$) telle que $q(e_{i},e_{n+1-i})=(-1)^{i}\varpi_{E}^{-1}$ pour $i=1,...,n$ et $q(e_{i},e_{j})=0$ pour $i+j\not=n+1$ (la forme $q$ est bien hermitienne car $n$ est pair). Gr\^ace \`a cette base, on identifie $GL(V\otimes_{E}\bar{F})$ \`a $GL(n,\bar{F})$. On introduit l'\'epinglage   de $GL(n,\bar{F})$ que l'on note ici $\tilde{\mathfrak{E}}=(\tilde{B},\tilde{T},(E_{\alpha_{i}})_{i=1,...,n-1})$    et l'automorphisme $\theta$ d\'efinis au d\'ebut du paragraphe \ref{Bn}.   En notant $\sigma\mapsto \sigma_{GL(n)}$ l'action naturelle de $\Gamma_{F}$ sur $GL(n,\bar{F})$, on v\'erifie que le groupe unitaire $\tilde{G}$ de $V$ s'identifie \`a $GL(n,\bar{F})$ muni de l'action $\sigma\mapsto \sigma_{G}$ telle que $\sigma_{G}=\sigma_{GL(n)}$ si $\sigma\in \Gamma_{E}$ et $\sigma_{G}=\theta\circ \sigma_{GL(n)}$ si $\sigma\not\in \Gamma_{E}$.  On v\'erifie que $\theta(E_{n-1,1})=-E_{n,2}$. Rappelons que $\Delta^{nr}_{a}=\{\beta_{i,n-i}; i=1,...,n/2-1\}\cup \{\beta_{n/2},\beta_{0}\}$. Posons $E_{\beta_{i,n-i}}=E_{\alpha_{i}}+E_{\alpha_{n-i}}$ pour $i=1,...,n/2-1$, $E_{\beta_{n/2}}=E_{\alpha_{n/2}}$ et  $E_{\beta_{0}}=\varpi_{E}(E_{n-1,1}+E_{n,2})$.     Alors $\tilde{\mathfrak{ E}}_{a}=(\tilde{T}^{nr},(E_{\beta})_{\beta\in \Delta^{nr}_{a}})$ est un \'epinglage affine de $\tilde{G}$ conserv\'e par l'action galoisienne. On d\'eduit de cet \'epinglage affine un tel \'epinglage $\mathfrak{ E}_{a}=(T^{nr},(E_{\beta})_{\beta\in \Delta^{nr}_{a}})$ de $G$,  o\`u $T$ est l'image de $\tilde{T}$ dans $G=\tilde{G}_{AD}$. Remarquons que $\tilde{\mathfrak{g}}(F)\subset \mathfrak{gl}(n,E)$ et que $\mathfrak{g}(F)$ est le sous-espace des \'el\'ements de trace nulle dans $\tilde{\mathfrak{g}}(F)$. 

L'unique \'el\'ement non trivial $\omega\in \boldsymbol{\Omega}^{nr}$ se rel\`eve en  l'\'el\'ement  de $\tilde{G}$ qui permute $e_{1}$ et $\varpi_{E}e_{n}$ et fixe $e_{i}$ pour $i=2,...,n-1$.  Notons $F'$ l'extension quadratique non ramifi\'ee de $F$. Pour $\sigma\in \Gamma_{F}$, on pose $a(\sigma)=0$ si $\sigma\in \Gamma_{F'}$ et $a(\sigma)=1$ sinon. On a $\sigma_{H}=Ad(\omega)^{a(\sigma)}\circ \sigma_{G}$.  Posons $n'=h^2$ et $n''=k(k+1)$. Notons ${\cal S}$ le $\mathfrak{o}_{E}$-r\'eseau de $V$ engendr\'e par les $e_{i}$ pour $i=1,...,n'/2$ et les  $\varpi_{E}e_{i}$ pour $n'/2+1,...,n$ et notons ${\cal S}^{\star}$ son dual, c'est-\`a-dire ${\cal S}^{\star}=\{v\in V; \forall v'\in {\cal S},\,q(v',v)\in \mathfrak{o}_{E}\}$. C'est le $\mathfrak{o}_{E}$-r\'eseau engendr\'e par les $e_{i}$  pour $i=1,...,n-n'/2$ et les $\varpi_{E}e_{i}$ pour $i=n-n'/2+1,...,n$. Alors $\tilde{\mathfrak{k}}_{s}$ est le sous-ensemble des $X\in \tilde{\mathfrak{g}}(F)$ tels que $X({\cal S})\subset {\cal S}$ et on a $\mathfrak{k}_{s}=\mathfrak{g}(F)\cap \tilde{\mathfrak{k}}_{s}$. L'espace $\bar{V}'={\cal S}/\varpi_{E}{\cal S}^{\star}$ sur ${\mathbb F}_{q}$ est de dimension $h^2$ et est naturellement muni d'une forme orthogonale. On note $\tilde{G}'_{s}$ son groupe sp\'ecial orthogonal. L'espace $\bar{V}''={\cal S}^{\star}/{\cal S}$ sur ${\mathbb F}_{q}$ est de dimension $k(k+1)$ et est naturellement muni d'une forme symplectique, qui est la r\'eduction de $\varpi_{E}q$. On note $\tilde{G}''_{s}$ son groupe symplectique. L'extension $E/F$ \'etant ramifi\'ee, on a l'\'egalit\'e $\tilde{\mathfrak{g}}_{s}=\mathfrak{g}_{s}$ et 
on obtient la d\'ecomposition $\mathfrak{g}_{s}=\tilde{\mathfrak{g}}'_{s}\oplus \tilde{\mathfrak{g}}''_{s}$ et ${\cal O}={\cal O}'\times {\cal O}''$, o\`u ${\cal O}'$ est param\'etr\'ee par la partition  $(2h-1,...,3,1)$ et ${\cal O}''$ est param\'etr\'ee par la partition $(2k,...,4,2)$.

  On pose $\eta=h-k\in \{0,1\}$. Pour $m=1,...,h$, on pose $d_{m}=4m-1-2\eta$ et on introduit un \'el\'ement $\alpha_{m}\in \bar{F}^{\times}$ tel que $\alpha_{m}^{d_{m}}=\varpi_{E}/2$. On pose $E_{m}=E(\alpha_{m})$ et $F_{m}=F(\alpha_{m}^2)$. On v\'erifie que $E_{m}$ est la compos\'ee des deux extensions $E$ et $F_{m}$ de $F$. On munit $E_{m}$ de la forme hermitienne $q_{m}$ (relative \`a l'extension $E/F$) d\'efinie par $q_{m}(v,v')=(-1)^m\frac{1}{2d_{m}}trace_{E_{m}/E}(\bar{v}v')$, o\`u $v\mapsto \bar{v}$ est la conjugaison galoisienne associ\'ee \`a l'extension $E_{m}/F_{m}$. On peut identifier convenablement $(V,q)$ \`a la somme directe des $(E_{m},q_{m})$ de sorte que ${\cal S}$ s'identifie \`a la somme des $\mathfrak{p}_{E_{m}}^{1-m}$.  L'identification est un peu compliqu\'ee parce que les $E_{m}$ sont de dimension impaire sur $E$.  De nouveau, on regroupe les $E_{m}$ deux-\`a-deux. Pour $m=1,...,h/2$, posons $n''_{m}=8m-4$ et $n'_{m}=8m-2-4\eta$. On d\'ecompose l'ensemble $\{1,...,n\}$ en une suite croissante d'intervalles
$$I'_{1}\cup ...\cup I'_{h/2}\cup I''_{h/2}\cup...\cup I''_{1}\cup J''_{1}\cup...\cup J''_{h/2}\cup J'_{h/2}\cup... J'_{1}$$
de sorte que $\vert I''_{m}\vert =\vert J''_{m}\vert =n''_{m}$ et $\vert I'_{m}\vert =\vert J'_{m}\vert =n'_{m}$. On voit que l'on peut identifier $E_{2m-1}\oplus E_{2m}$ au sous-espace de $V$ engendr\'e par les $e_{i}$ pour $i\in I'_{m}\cup I''_{m}\cup J''_{m}\cup J'_{m}$.  On pr\'ecisera ci-dessous l'identification pour $m=1$.  On note $X^G$ l'\'el\'ement de $\mathfrak{g}(F)$ qui conserve chaque $E_{m}
$ et y agit par multiplication par $\alpha_{m}$, sauf dans le cas o\`u $m=1$ et $\eta=1$. Dans ce cas, $X^G$ annule $E_{1}=E$.

{\bf Remarque.} On fait cette exception pour que $X^G$ appartienne \`a $\mathfrak{g}(F)$, c'est-\`a-dire soit de trace nulle.

\bigskip

Pour tout $m$, la multiplication par $\alpha_{m}$ conserve $\mathfrak{p}_{E_{m}}^{1-m}$ donc $X^G\in \tilde{k}_{s}$. On note $N$ la r\'eduction de $X^G$ dans $\mathfrak{g}_{s}$. On voit comme toujours que $N\in {\cal O}$. On doit d\'efinir un \'el\'ement $X^H$ et prouver que $N$ et $X^H$ ont les propri\'et\'es requises. Notre \'el\'ement $X^H$ va \^etre somme de ses composantes sur chaque $E_{m-1}\oplus E_{m}$. 
Avec nos identifications,  $Ad(\omega)$ agit trivialement sur $E_{2m-1}\oplus E_{2m}$ d\`es que $m\geq 2$. Donc les composantes de $X^G$ sur ces espaces appartiennent \`a $\mathfrak{h}(F)$ et les r\'eductions de ces composantes dans $\mathfrak{g}_{s}$ sont fix\'ees par $Ad(\omega)_{s}$.  On impose $X^H=X^G$ sur $E_{2m-1}\oplus E_{2m}$.  Il est clair que tout se passe bien sur ces composantes et que, pour achever  la preuve, il reste \`a consid\'erer le  
  facteur $E_{1}\oplus E_{2}$. On peut simplifier les notations en supposant que $h=2$. On distingue les cas $k=2$ ($\eta=0$) et $k=1$ ($\eta=1$). 

Supposons $k=2$. On a $n=10$. Introduisons les \'el\'ements suivants de $E_{1}\oplus E_{2}$:
$$e'_{1}=\alpha_{1}^{0}+\alpha_{2}^0,\, e'_{2}=\alpha_{2}^{-1},\,e'_{3}=\alpha_{2}^{-2},\, e'_{4}=-\alpha_{1}^{-1},\,e'_{5}=\alpha_{2}^{-3},$$

$$ e'_{6}=-\alpha_{2}^{-4},\, e'_{7}=\alpha_{1}^{-2},\, e'_{8}=\alpha_{2}^{-5},\, e'_{9}=\alpha_{2}^{-6},\, e'_{10}=\varpi_{E}^{-1}(-\alpha_{1}^{0}+\alpha_{2}^0).$$
Comme en \ref{Bn}, $\alpha_{1}^{0}$ et $\alpha_{2}^{0}$ d\'esignent clairement les unit\'es, l'une de $E_{1}$, l'autre de $E_{2}$. 
On identifie $E_{1}\oplus E_{2}$ \`a $V$ par   l'application qui envoie $e'_{i}$ sur $e_{i}$ pour tout $i=1,...,10$, dont on v\'erifie qu'elle est une isom\'etrie. Notons $E_{1}^{\flat}$ le sous-$E$-espace de $E_{1}$ engendr\'e par $\alpha_{1}^{-1}$ et $\alpha_{1}^{-2}$. On voit que $\omega$ agit par l'identit\'e sur  $E_{1}^{\flat}\oplus E_{2}$ et par multiplication par $-1$ sur $E\alpha_{1}^0$. Cette action ne fixe par $X^G$. Notons $N$ la r\'eduction de $X^G$ dans $\mathfrak{g}_{s}$, que l'on d\'ecompose en $N=N'\oplus N''$, avec $N'\in \tilde{\mathfrak{g}}'_{s}$, $N''\in \tilde{\mathfrak{g}}''_{s}$. On voit que $ \omega$ se r\'eduit en  l'identit\'e de $\bar{V}''_{s}$ et en l'automorphisme de $\bar{V}'_{s}$ qui agit par l'identit\'e sur $\mathfrak{p}_{E_{2}}^{-1}/\mathfrak{p}_{E_{2}}^{2}$ et par multiplication par $-1$ sur $\mathfrak{o}_{E_{1}}/\mathfrak{p}_{E_{1}}$. Or $N'$ conserve ces deux sous-espaces donc $N$ est fix\'e par $Ad(\omega)_{s}$. 
 Cela d\'emontre que  $N$ appartient \`a $ \mathfrak{r}({\mathbb F}_{q})$.    
Notons $\tau$ l'\'el\'ement de $\tilde{G}$ qui agit par l'identit\'e sur $E_{2}$ et par multiplication par $-1$ sur $E_{1}$. Il commute \`a $X^G$. D'apr\`es la description ci-dessus de $\omega$, on a $\omega\tau=\tau\omega$ et cet \'el\'ement agit par l'identit\'e sur $F\alpha_{1}^0\oplus E_{2}$ et par multiplication par $-1$ sur $E_{1}^{\flat}$. Autrement dit, $\omega\tau$ multiplie les vecteurs $e_{4}$ et $e_{7}$ par $-1$ et fixe les autres vecteurs de base. Donc $\omega\tau\in \tilde{T}(\mathfrak{o}_{F})$. Par le th\'eor\`eme de Lang, on peut fixer $t\in \tilde{T}(\mathfrak{o}_{F^{nr}})$ tel que $t^{-1}Fr_{H}(t)=\omega\tau$. Posons $X^{H}=Ad(t^{-1})(X^G)$. Par le m\^eme calcul qu'en \ref{Bn}, on voit que $Fr_{H}(X^H)=X^H$. Donc $X^H\in \mathfrak{h}(F)$. Puisque $X^G$ et $X^H$ sont conjugu\'es, leurs classes de conjugaison stable se correspondent. L'assertion (ii) du lemme \ref{debutdelapreuve} pour $X^G$ et le groupe $\tilde{G}$ est l'assertion 5.2(3) de \cite{W7}. Puisque $t\in K_{s}^{\tilde{P},nr}$, cette assertion pour $X^G$ entra\^{\i}ne la m\^eme assertion pour $X^H$. Il reste \`a d\'emontrer l'assertion (iii) du lemme \ref{debutdelapreuve}. Notons $N^H$ la r\'eduction de $X^H$ dans $\mathfrak{h}_{s}$. Par d\'efinition, ${\bf f}^G(X^G)=f^G(N)=1$, ${\bf f}^H(X^H)=f^H(N^H)$ et $f^H(N)=-1$ car le terme $\epsilon$ d\'efini en \ref{preparatifs} vaut $-1$. Notons $\bar{t}$ la r\'eduction de $t$ dans $\tilde{G}_{s}$. On a $N^H=\bar{t}^{-1}N\bar{t}$, d'o\`u $f^H(N^H)=\xi(\bar{z})f^H(N)=-\xi(\bar{z})$, o\`u $\xi$ est l'unique \'el\'ement de $\Xi$ et $\bar{z}=Fr_{H}(\bar{t})\bar{t}^{-1}$. D\'ecomposons $\bar{z}$ en $\tilde{z}'\times \bar{z}''$, o\`u $\bar{z}'\in \tilde{G}'_{s}$ et $\bar{z}''\in\tilde{G}''_{s}$. D'apr\`es la description ci-dessus de  
$t^{-1}Fr_{H}(t)=\omega\tau$, $\bar{z}'$ est l'\'el\'ement neutre tandis que $\bar{z}''$ a quatre valeurs propres \'egale \`a $1$ et deux valeurs propres \'egales \`a $-1$. 
L'\'el\'ement $N''$ appartient \`a l'orbite nilpotente param\'etr\'ee par la partition $(4,2)$. Son commutant dans $\tilde{G}''_{s}$ est produit d'un groupe unipotent et d'un groupe isomorphe \`a $({\mathbb Z}/2{\mathbb Z})^2$ engendr\'e par deux \'el\'ements $\bar{z}_{2}$, resp. $\bar{z}_{4}$, poss\'edant pour valeurs propres $1$ avec multiplicit\'es $4$, resp. $2$,  et $-1$ avec multiplicit\'es $2$, resp. $4$. On en d\'eduit que $\bar{z}''$ est le produit de $\bar{z}_{2}$ par un \'el\'ement unipotent. Or $\xi$ vaut $-1$ sur $\bar{z}_{2}$, donc $\xi(\bar{z})=-1$. On en d\'eduit ${\bf f}^H(X^H)=1={\bf f}^G(X^G)$, ce qui d\'emontre le (iii) du lemme \ref{debutdelapreuve}.

Supposons maintenant $k=1$. On a $n=6$. Introduisons les \'el\'ements suivants de $E_{1}\oplus E_{2}$:
$$e'_{1}=\alpha_{1}^{0}+\alpha_{2}^0,\, e'_{2}=\alpha_{2}^{-1},\,e'_{3}=\alpha_{2}^{-2},\, e'_{4}=\alpha_{2}^{-3},\,e'_{5}=\alpha_{2}^{-4},\,\, e'_{6}=\varpi_{E}^{-1}(-\alpha_{1}^0+\alpha_{2}^{0}).$$
On identifie encore $E_{1}\oplus E_{2}$  à $V$ \`a l'aide de cette base. Alors $\omega$ agit par l'identit\'e sur $E_{2}$ et par multiplication par $-1$ sur $E_{1}$. Donc $X^G$ est fixe par $\omega$. On pose $X^H=X^G$. On a $X^H\in \mathfrak{h}(F)$. Cette fois, on a $\epsilon=1$.  La preuve ci-dessus se trivialise et conduit au r\'esultat cherch\'e.      $\square$

\subsection{Utilisation de la construction de Lusztig}

\subsubsection{Définition d'un caractère $\epsilon^{\flat}_{N}$}\label{caractereepsilonflat}

Soit $G$ un groupe réductif connexe défini sur $F^{nr}$. {\bf On suppose} $G$ {\bf simplement connexe}. On fixe un épinglage $\mathfrak{E}=(B,T,(E_{\alpha})_{\alpha\in \Delta})$  et un épinglage affine $\mathfrak{E}_{a}$ tous deux définis sur $F^{nr}$ et on en déduit les objets introduits en \ref{alcove} et \ref{racinesaffines}. 

Soient $M$ un $F^{nr}$-Levi de $G$ et $s$ un sommet de $Imm_{F^{nr}}(M_{AD})$.  Soit  ${\cal E}\in {\bf FC}^{M_{s}}(\mathfrak{m}_{SC,s})$. On note ${\cal O}$ l'orbite nilpotente de $\mathfrak{m}_{SC,s}$ supportant ${\cal E}$ et on fixe $N\in {\cal O}$. Le faisceau ${\cal E}$ est associé à une représentation irréductible  $\xi_{N}$ de $Z_{M_{s}}(N)/Z_{M_{s}}(N)^0$ dans un espace $V$. Il s'agit a priori d'un $\bar{{\mathbb Q}}_{l}$-espace  mais nous l'identifions à un ${\mathbb C}$-espace par l'isomorphisme $\bar{{\mathbb Q}}_{l}\simeq {\mathbb C}$ que l'on a fixé.

On a défini le groupe $Norm_{G(F^{nr})}(M,s)$ en \ref{normalisateur}. Pour tout élément $g$ de ce groupe, l'automorphisme $Ad(g)$ définit naturellement des automorphismes $Ad(g)_{s}$ de $M_{s}$ et $M_{AD,s}$. Prouvons que

(1) pour tout $g\in Norm_{G(F^{nr})}(M,s)$, $Ad(g)_{s}$ conserve ${\cal E}$.

Preuve.  Le procédé de \ref{reduction} nous ramène au cas au $G$ est absolument simple. 
 Le couple $(M,s)$ est excellent, cf. \ref{excellentscouples}(1), et on l'identifie à un couple $(M_{\Lambda},s_{\Lambda})$ où $\Lambda$ est un excellent sous-ensemble propre de $\Delta_{a}^{nr}$. Si $\Delta_{a}^{nr}-\Lambda$ n'a qu'un élément, on a $M=G$, $Norm_{G(F^{nr})}(M,s)=K_{s}^{M,0,nr}$ et l'assertion est claire. Supposons que $\Delta_{a}^{nr}-\Lambda$  a au moins deux éléments. D'après \ref{normalisateur} (3) et (4), on a la suite exacte
 $$(2) \qquad 1 \to K_{s}^{M,0,nr}\to Norm_{G(F^{nr})}(M,s)\to W^{aff}(\Lambda^{aff})\to 1.$$
 D'après la proposition \ref{excellentsensemblesLambda}, le groupe $W^{aff}(\Lambda^{aff})$ est engendré par les symétries $w_{\gamma}$ associées au éléments $\gamma\in \Delta_{a}^{nr}-\Lambda$. Il en résulte que $Norm_{G(F^{nr})}(M,s)$ est engendré par ses intersections $K_{{\cal F}_{\gamma}}^{0,nr}\cap  Norm_{G(F^{nr})}(M,s)$, cf. \ref{excellentsensemblesLambda} pour la définition de ${\cal F}_{\gamma}$. Mais, pour $g\in K_{{\cal F}_{\gamma}}^{0,nr}\cap  Norm_{G(F^{nr})}(M,s)$, $Ad(g)_{s}$ est l'action de la réduction $\bar{g}$ de $g$ dans $G_{{\cal F}_{\gamma}}$. L'assertion résulte alors de \ref{faisceauxcaracteres}(2). Cela prouve (1).

On note $Norm_{G(F^{nr})}(M,s,N)$ le sous-groupe des $g\in Norm_{G(F^{nr})}(M,s)$ tels que $Ad(g)_{s}$ fixe $N$. Montrons que l'on a l'isomorphisme

(3) $K_{s,N}^{M,0,nr}\backslash Norm_{G(F^{nr})}(M,s,N) \simeq K_{s}^{M,0,nr}\backslash Norm_{G(F^{nr})}(M,s).$

L'homomorphisme naturel du membre de gauche dans celui de droite est injectif. Soit $g\in Norm_{G(F^{nr})}(M,s)$. L'automorphisme  $Ad(g)_{s}$ conserve le  faisceau ${\cal E}$ donc aussi ${\cal O}$. Il existe donc $x\in M_{s}$ tel que $Ad(g)_{s}(N)=xNx^{-1}$. On relève $x$ en un élément $k\in K_{s}^{M,0,nr}$. Alors $k^{-1}g$ appartient à $Norm_{G(F^{nr})}(M,s,N) $. L'homomorphisme est donc surjectif, ce qui prouve (3).

On a défini l'ensemble $P(s)=p_{M}^{-1}(s)$ en \ref{reduction}. Soit $x\in P(s)$. Pour $g\in K_{x}^{0,nr}\cap Norm_{G(F^{nr})}(M,s)$, la réduction $\bar{g}$ de $g$ dans $G_{x}$ appartient au groupe $Norm_{G_{x}}(M_{s})$. On note 
$$c_{x}: K_{x}^{0,nr}\cap Norm_{G(F^{nr})}(M,s) \to Norm_{G_{x}}(M_{s})$$
cet homomorphisme de réduction. Il est surjectif. En effet, fixons un tore $T^{nr}\in {\cal T}_{max}^{M,nr}$ tel que $s$ appartienne à l'appartement associé à $T^{nr}$ dans $Imm_{F^{nr}}(M_{AD})$. Le point $x$ appartient à $App_{F^{nr}}(T^{nr})$ et il se déduit de $T^{nr}$ un sous-tore maximal $T_{x}$ de $G_{x}$. Le groupe $ Norm_{G_{x}}(M_{s})$ est engendré par $Norm_{G_{x}}(T_{x},M_{s})$ et $M_{s}$. Le premier   groupe est contenu dans l'image de $ c_{x}$ d'après le (iii) du lemme \ref{normalisateur}, le second est contenu dans l'image du sous-groupe $K_{s}^{M,0,nr}$  de l'espace de départ. Cela prouve la surjectivité de $ c_{x}$. On a défini le groupe $Norm_{G_{x}}(M_{s},N)$ en \ref{faisceauxcaracteres}. L'homomorphisme $c_{x}$ se restreint en un homomorphisme surjectif
$$c_{x,N}: K_{x}^{0,nr}\cap Norm_{G(F^{nr})}(M,s,N) \to Norm_{G_{x}}(M_{s},N).$$

A la suite de Lusztig, on a défini  en \ref{faisceauxcaracteres} une représentation $\epsilon_{N}^{\flat}$ de $Norm_{G_{x}}(M_{s},N)$ dans $V$. Nous noterons plus précisément $\epsilon_{x,N}^{\flat}$ cette représentation. Elle  prolonge $\xi_{N}$, ou plus exactement le relèvement de $\xi_{N}$ en une représentation  de $Z_{M_{s}}(N)$.

\begin{thm}{Il existe une unique représentation $\epsilon_{N}^{\flat}$ de $Norm_{G(F^{nr})}(M,s,N)$ dans $V$ telle que, pour tout $x\in P(s)$ et tout $g\in K_{x}^{0,nr}\cap Norm_{G(F^{nr})}(M,s,N)$ , on ait l'égalité
$\epsilon_{N}^{\flat}(g)=\epsilon_{x,N}^{\flat}\circ c_{x,N}(g)$.}\end{thm} 

Preuve. Le procédé de \ref{reduction} nous ramène au cas où $G$ est absolument simple. On identifie $(M,s)$ à un couple $(M_{\Lambda},s_{\Lambda})$ où $\Lambda$ est un excellent sous-ensemble propre de $\Delta_{a}^{nr}$. Si $\Delta_{a}^{nr}-\Lambda$ n'a qu'un élément, on a $M=G$, $Norm_{G(F^{nr})}(M,s)=K_{s}^{M,0,nr}$ et $P(s)=\{s\}$. L'application $\epsilon_{N}^{\flat}$ définie par $\epsilon_{N}^{\flat}(g)=\epsilon_{s,N}^{\flat}\circ c_{s,N}(g)$ pour $g\in K_{s,N}^{M,0,nr}$ vérifie les propriétés requises et c'est la seule possible. Supposons que $\Delta_{a}^{nr}-\Lambda$  a au moins deux éléments.  En vertu de (3), la suite exacte (2) se récrit
 $$1 \to K_{s,N}^{M,0,nr}\to Norm_{G(F^{nr})}(M,s,N)\to W^{aff}(\Lambda^{aff})\to 1.$$
 Notons $\boldsymbol{\xi}$ l'homomorphisme de $K_{s,N}^{M,0,nr}$ dans $Aut(V)$ qui relève $\xi_{N}$, autrement dit $\boldsymbol{\xi}(k)=\xi_{N}(\bar{k})$ pour tout $k\in K_{s,N}^{M,0,nr}$, où $\bar{k}$ est la réduction de $k$ dans $Z_{M_{s}}(N)$.   Posons $S=\{w_{\gamma}; \gamma\in \Delta_{a}^{nr}-\Lambda\}$. Le couple $(W^{aff}(\Lambda^{aff}),S)$ est un système  de Coxeter. On rappelle que $W^{aff}(\Lambda^{aff})$ est alors muni d'une longueur  que l'on note $l$ et qu'il y a une notion de décomposition réduite d'un élément de ce groupe. Soit  $\gamma\in \Delta_{a}^{nr}-\Lambda$. Posons $w=w_{\gamma}$. Par définition de cet élément, il appartient à l'image de $K_{{\cal F}_{\gamma}}^{0,nr}\cap Norm_{G(F^{nr})}(M,s,N)$. On fixe un relèvement $n_{w}$ de $w$ dans cet ensemble. Pour utiliser les notations déjà introduites, fixons $x\in {\cal F}$. Posons $\epsilon(w)=\epsilon_{x,N}^{\flat}\circ c_{x,N}(n_{w})$. C'est un élément de $Aut(W)$. Il résulte des définitions que
 
(4) $\epsilon(w)\boldsymbol{\xi}(k)\epsilon(w)^{-1}=\boldsymbol{\xi}(n_{w}kn_{w}^{-1})$ pour tout $k\in K_{s,N}^{M,0,nr}$. 

Soit $g\in Norm_{G(F^{nr})}(M,s,N)$, notons $w$ son image dans $W^{aff}(\Lambda^{aff})$ et notons $D(w)$ l'ensemble des décompositions réduites de $w$. On définit une application ${\bf F}:D(w)\to Aut(V)$ de la façon suivante. Soit ${\bf w}=(w_{1},...,w_{l})\in D(w)$, c'est-à-dire que $l=l(w)$, que les $w_{i}$ appartiennent à $S$ et que l'on a l'égalité 
$w=w_{1}...w_{l}$. Il existe un unique $k({\bf w})\in K_{s,N}^{M,0,nr}$ tel que $g=k({\bf w})n_{1}...n_{l}$, où, pour simplifier, on a posé $n_{i}=n_{w_{i}}$ pour tout $i=1,...,l$. On pose ${\bf F}({\bf w})=\boldsymbol{\xi}(k({\bf w}))\epsilon(w_{1})...\epsilon(w_{l})$. Remarquons que cette définition se symétrise: il existe un unique $h({\bf w})\in K_{s,N}^{M,0,nr}$ tel que $g= n_{1}...n_{l}h({\bf w})$ et on a ${\bf F}({\bf w})= \epsilon(w_{1})...\epsilon(w_{l})\boldsymbol{\xi}(h({\bf w}))$. Cela résulte de (4). On utilisera le résultat suivant. Notons $\gamma_{1}$,...,$\gamma_{l}$ les éléments de $\Delta_{a}^{nr}-\Lambda$ tels que $w_{i}=w_{\gamma_{i}}$ pour tout $i=1,...,l$. Soit $x\in \bar{C}^{nr}\cap P(s)$,  supposons que $x$ appartienne à l'adhérence  de ${\cal F}_{\gamma_{i}}$ pour tout $i$. On a  $K_{{\cal F}_{\gamma_{i}}}^{0,nr}\subset K_{x}^{0,nr}$ pour tout $i$ et il résulte des définitions que $g\in K_{x}^{0,nr}$. L'égalité ci-dessous  résulte alors des définitions et  de \ref{faisceauxcaracteres}(6):

(5) ${\bf F}({\bf w})=\epsilon_{x,N}^{\flat}\circ c_{x,N}(g)$.

Montrons que

(6) la fonction ${\bf F}$ est constante sur $D(w)$. 

On raisonne par récurrence sur la longueur de $w$. Evidemment, l'assertion est triviale si la longueur de $w$ est $0$ ou $1$ (il n'y a qu'une décomposition de $w$). On suppose cette longueur supérieure ou égale à $2$. D'après le lemme 4 de \cite{B} IV.1.5, il suffit de prouver l'égalité ${\bf F}({\bf w})={\bf F}({\bf w}')$ pour deux éléments ${\bf w}$ et ${\bf w}'$ de $D$ vérifiant l'une des hypothèses suivantes, où l'on note ${\bf w}=(w_{1},...,w_{l})$ et ${\bf w}'=(w'_{1},...,w'_{l})$:

(7)  $w_{l}=w'_{l}$;

(8) $w_{1}=w'_{1}$;

(9) il existe $v,v'\in S$ de sorte que ${\bf w}=(v,v',v,v',...)$ tandis que ${\bf w}'=(v',v,v',v,...)$. 

Supposons (7) vérifiée. On pose $\underline{g}=gn_{l}^{-1}$ et on note $\underline{w}$ son image dans $W^{aff}$. On a $\underline{w}w_{l}=w$ et la longueur de $\underline{w}$ est $l-1$. 
 Les suites ${\underline{\bf w}}=(w_{1},...,w_{l-1})$ et $\underline{{\bf w}}'=(w'_{1},...,w'_{l-1})$ sont deux décompositions réduites de $\underline{w}$. En utilisant l'élement $\underline{g}$ au lieu de $g$, on définit une fonction $\underline{{\bf F}}$ sur $D(\underline{w})$. Par définition, on a ${\bf F}({\bf w})=\underline{{\bf F}}(\underline{{\bf w}})\epsilon(w_{l})$ et ${\bf F}({\bf w}')=\underline{{\bf F}}(\underline{{\bf w}}')\epsilon(w_{l})$. En utilisant l'hypothèse de récurrence, on a $\underline{{\bf F}}(\underline{{\bf w}})=\underline{{\bf F}}(\underline{{\bf w}}')$, d'où l'égalité cherchée ${\bf F}({\bf w})={\bf F}({\bf w}')$.
 
 La preuve dans le cas où (8) vérifiée est similaire, en utilisant la définition "symétrique" de ${\bf F}$.
 
 Supposons (9) vérifiée. On a $v\not=v'$ puisque ${\bf w}$ est une décomposition réduite et que l'on a supposé $l\geq2$. L'existence de deux décompositions réduites distinctes de la forme indiquée entraîne que le groupe engendré par $v$ et $v'$ est fini. Ce groupe n'est pas $W_{aff}(\Lambda^{aff})$ tout entier car on a déjà vu plusieurs fois que $W^{aff}(\Lambda^{aff})$ était infini.   Cela entraîne que $\Delta_{a}^{nr}-\Lambda$ a au moins trois éléments. Notons $\gamma$ et $\gamma'$ les deux éléments de cet ensemble tels que $v=w_{\gamma}$ et $v'=w_{\gamma'}$, fixons un troisième élément $\delta$ et  notons $s_{\delta}$ le sommet de $\bar{C}^{nr}$ associé à  cet élément. Alors $s_{\delta}$ est adhérent aux facettes ${\cal F}_{\gamma}$ et ${\cal F}_{\gamma'}$.  D'après (5), on a
 $${\bf F}({\bf w})=\epsilon_{s_{\delta},N}^{\flat}\circ c_{s_{\delta},N}(g)={\bf F}({\bf w}').$$
 Cela achève la preuve de (6). 
 
 Gr\^ace à (6), on peut définir $\epsilon_{N}^{\flat}(g)$ comme la valeur constante de la fonction ${\bf F}$. Il résulte des définitions et de (4) que l'on a l'égalité
 
 (10) $\epsilon_{N}^{\flat}(g)\boldsymbol{\xi}(k)\epsilon_{N}^{\flat}(g)^{-1}=\boldsymbol{\xi}(gkg^{-1})$ pour tous $g\in Norm_{G(F^{nr})}(M,s,N)$ et $k\in K_{s,N}^{M,0,nr}$.
 
  Montrons que
 
 (11) $\epsilon_{N}^{\flat}$ est un homomorphisme de $Norm_{G(F^{nr})}(M,s,N)$ dans $Aut(V)$. 
 
 Soient $g,g'\in Norm_{G(F^{nr})}(M,s,N)$. Notons $w$ et $w'$ leurs images dans $W^{aff}(\Lambda^{aff})$. On veut prouver l'égalité $\epsilon_{N}^{\flat}(g)\epsilon_{N}^{\flat}(g')=\epsilon_{N}^{\flat}(gg')$. On raisonne par récurrence sur $l(w)+l(w')$.  Fixons des décompositions réduites $(w_{1},...,w_{l})$ et $(w'_{1},...,w'_{l'})$ de $w$ et $w'$ et écrivons $g=kn_{1}...,n_{l}$, $g'=k'n'_{1}...n'_{l'}$ avec $k,k'\in K_{s,N}^{M,0,nr}$. Supposons d'abord que $l(ww')=l(w)+l(w')$. La suite $(w_{1},...,w_{l},w'_{1},...,w'_{l'})$ est une décomposition réduite de $ww'$ et on a $gg'=k''n_{1}...n_{l}n'_{1}...,n'_{l'}$, avec $k''=gk'g^{-1}k$. L'égalité $\epsilon_{N}^{\flat}(g)\epsilon_{N}^{\flat}(g')=\epsilon_{N}^{\flat}(gg')$ résulte immédiatement des définitions et de (10). Supposons maintenant que $l(ww')<l(w)+l(w')$. Cela entraîne $l(w),l(w')\geq1$. Il existe  un unique $j\in \{1,...,l'\}$ tel que $l(ww'_{1}...,w'_{j-1})=l+j-1$ tandis que $l(ww'_{1},...,w'_{j})=l+j-2$. Supposons d'abord $j<l'$. 
 Posons $g'_{1}=k'n'_{1}...,n'_{j}$ et $g'_{2}=n'_{j+1}...n'_{l'}$. On a $g'=g'_{1}g'_{2}$. Il résulte des définitions  que
$\epsilon_{N}^{\flat}(g')=\epsilon_{N}^{\flat} (g'_{1}) \epsilon_{N}^{\flat}(g'_{2})$. Il résulte de l'hypothèse de récurrence que $\epsilon_{N}^{\flat}(g)\epsilon_{N}^{\flat}(g'_{1})=\epsilon_{N}^{\flat}(gg'_{1})$ puis que $\epsilon_{N}^{\flat}(gg'_{1})\epsilon_{N}^{\flat}(g'_{2})=\epsilon_{N}^{\flat}(gg'_{1}g'_{2})=\epsilon_{N}^{\flat}(gg')$. D'où l'égalité requise. Supposons maintenant $j=l'$. On pose $g'_{1}=k'n'_{1}...,n'_{l'-1}$ et $g_{1}=gg'_{1}$. On a $\epsilon_{N}^{\flat}(g')=\epsilon_{N}^{\flat}(g'_{1})\epsilon_{N}^{\flat}(n'_{l})$ et, par l'hypothèse de récurrence, $\epsilon_{N}^{\flat}(g)\epsilon_{N}^{\flat}(g'_{1})=\epsilon_{N}^{\flat}(gg'_{1})$. Il reste à prouver que  $\epsilon_{N}^{\flat}(gg'_{1})\epsilon_{N}^{\flat}(n'_{l})=\epsilon_{N}^{\flat}(gg'_{1}n'_{l})$. Autrement dit, on est ramené au cas où  $l(w')=1$, c'est-à-dire $w'\in S$,   et où $g'=n_{w'}$, avec toujours $l(ww')<l(w)+l(w')=l(w)+1$. Cela entraîne $l(ww')=l(w)-1$. Les images dans $W^{aff}(\Lambda^{aff})$ de $gg'$ et ${g'}^{-1}$ sont $ww'$ et $w'$. Gr\^ace à l'hypothèse de récurrence, on a $\epsilon_{N}^{\flat}(g)=\epsilon_{N}^{\flat}(gg')\epsilon_{N}^{\flat}({g'}^{-1})$.  Notons $\gamma'$ l'élément de $\Delta_{a}^{nr}-\Lambda$ tel que $w'=w_{\gamma'}$. En appliquant (5) à un élément $x\in {\cal F}_{\gamma'}$, on obtient l'égalité $\epsilon_{N}^{\flat}({g'}^{-1})\epsilon_{N}^{\flat}(g')=1$. D'où l'égalité cherchée $\epsilon_{N}^{\flat}(g)\epsilon_{N}^{\flat}(g')=\epsilon_{N}^{\flat}(gg')$. Cela prouve (11). 

 Montrons que  

(12) pour tout élément $x\in P(s)\cap \bar{C}^{nr}$ et tout $g\in K_{x}^{0,nr}\cap Norm_{G(F^{nr})}(M,s,N)$ , on a l'égalité
$\epsilon_{N}^{\flat}(g)=\epsilon_{x,N}^{\flat}\circ c_{x,N}(g)$.

 Le système de racines de $G_{x}$ est $\Lambda(x)$. D'après \cite{L2}, théorème 9.2, le groupe $Norm_{G_{x}}(M_{s})$ est un groupe de Coxeter engendré par des symétries associées aux éléments de $\Lambda(x)-\Lambda$.  Il en résulte que l'image de $g$ dans $W^{aff}(\Lambda^{aff})$ appartient au groupe engendré par les $w_{\gamma'}$ pour $\gamma'\in \Lambda(x)-\Lambda$. En appliquant (5) à une décomposition quelconque de cette image, on obtient l'égalité cherchée. 

Montrons enfin que

(13) pour tout élément $x\in P(s) $ et tout $g\in K_{x}^{0,nr}\cap Norm_{G(F^{nr})}(M,s,N)$ , on a l'égalité
$\epsilon_{N}^{\flat}(g)=\epsilon_{x,N}^{\flat}\circ c_{x,N}(g)$.

D'après la proposition \ref{excellentsensemblesLambda}, il existe $n\in Norm_{G(F^{nr})}(M,s)$ tel que $nx\in P(s)\cap \bar{C}^{nr}$. D'après (3), on peut supposer $n\in Norm_{G(F^{nr})}(M,s,N)$. Fixons un tel $n$, posons $y=nx$, $\epsilon_{N}^{\flat}(n)=A$ et  $h=ngn^{-1}$. On a  $h\in K_{y}^{0,nr}\cap Norm_{G(F^{nr})}(M,s,N)$. En appliquant (12), on a l'égalité $\epsilon_{N}^{\flat}(h)=\epsilon_{y,N}^{\flat}\circ c_{y,N}(h)$. D'après (11), on a $\epsilon_{N}^{\flat}(h)=A\epsilon_{N}^{\flat}(g)A^{-1}$. L'automorphisme $Ad(n)$ se descend en un isomorphisme $a:G_{x}\to G_{y}$. On a $c_{y,N}(h)=a\circ c_{x,N}(g)$. En posant $\bar{g}=c_{x,N}(g)$, on voit qu'il reste à prouver l'égalité $A^{-1}\epsilon_{y,N}^{\flat}(a(\bar{g}))A=\epsilon_{x,N}^{\flat}(\bar{g})$. D'après \ref{faisceauxcaracteres}(8) et (9), le membre de gauche est la valeur en $\bar{g}$ de la fonction analogue  à $\epsilon_{x,N}^{\flat}$ déduite du  faisceau-caractère cuspidal associé à la représentation $z\mapsto A^{-1}\xi_{N}(a(z))A^{-1}$ de $Z_{M_{s}}(N)/Z_{M_{s}}(N)^0$. D'après (10), cette représentation est identique à $\xi_{N}$. Cette fonction analogue est donc égale à $\epsilon^{\flat}_{x,N}$ et l'égalité cherchée s'en déduit. Cela prouve (13).

On a prouvé que $\epsilon_{N}^{\flat}$ vérifiait toutes les conditions requises. Son unicité est évidente, les groupes $K_{x}^{0,nr}\cap Norm_{G(F^{nr})}(M,s,N)$  engendrant $Norm_{G(F^{nr})}(M,s,N)$ quand $x$ décrit seulement les sommets de $P(s)\cap \bar{C}^{nr}$. $\square$ 

La représentation $\epsilon_{N}^{\flat}$ vérifie les propriétés (14) à (17) ci-dessous. Elles se démontrent toutes en  montrant que le membre de droite des égalités en question vérifie la propriété caractéristique du membre de gauche. Cela résulte essentiellement des propriétés similaires \ref{faisceauxcaracteres}(4) à (9) des représentations $\epsilon_{x,N}^{\flat}$. On laisse la mise au point au lecteur.

Soient $G'$ un groupe réductif connexe défini sur $F^{nr}$ et $\phi:G\to G'$ un isomorphisme défini sur $F^{nr}$. Posons $M'=\phi(M)$. De $\phi$ se déduit un isomorphisme de $Imm_{F^{nr}}(M_{AD})$ sur $Imm_{F^{nr}}(M'_{AD})$, on note $s'$ l'image de $s$ par cet isomorphisme. De $\phi$ se déduit un isomorphisme $\phi_{s}:M_{s}\to M'_{s'}$. Celui-ci transporte ${\cal E}$ en un faisceau ${\cal E}'\in {\bf FC}^{M'_{s'}}(\mathfrak{m}'_{SC,s'})$. Fixons $N\in {\cal O}$, posons $N'=\phi_{s}(N)$ et fixons un automorphisme $A\in Aut(V)$. Considérons que ${\cal E}'$ est associé à la représentation $\xi_{N'}$ de $Z_{M'_{s'}}(N')/Z_{M'_{s'}}(N')^0$ définie par $\xi_{N'}(x')=A\xi_{N}(\delta_{s}^{-1}(x'))A^{-1}$ pour tout $x'\in Z_{M'_{s'}}(N')/Z_{M'_{s'}}(N')^0$. On en déduit une représentation $\epsilon_{N'}^{\flat}$ de $ Norm_{G'(F^{nr})}(M',s',N')$. On a l'égalité

(14) $\epsilon_{N'}^{\flat}(k')=A\epsilon_{N}^{\flat}(\delta^{-1}(k'))A^{-1}$ pour tout $k'\in Norm_{G'(F^{nr})}(M',s',N')$.

Soit $L$ un $F^{nr}$-Levi de $G$ contenant $M$, notons $M_{L,sc}$ l'image réciproque de $M$ dans $L_{SC}$. Puisque $G$ est simplement connexe, on a $L_{SC}\subset L$ et $M_{L,sc}=M\cap L_{SC}$.   On définit une représentation $\epsilon_{N}^{L_{SC},\flat}$ de $Norm_{L_{SC}(F^{nr})}(M_{sc},s,N)$ en appliquant la construction à $L_{SC}$.  Le groupe $Norm_{L_{SC}(F^{nr})}(M_{sc},s,N)$ est inclus dans $ Norm_{G(F^{nr})}(M,s,N)$. On a l'égalité

(15) $\epsilon_{N}^{L_{SC},\flat}(g)=\epsilon_{N}^{\flat}(g)$ pour tout $g\in Norm_{L_{SC}(F^{nr})}(M_{sc},s,N)$.

Supposons que $G$ soit défini sur $F$, que $M$ soit un $F$-Levi, que $s$ appartienne à $Imm(M_{AD})$ et que   ${\cal E}$ appartienne à $ {\bf FC}^{M_{s}}_{{\mathbb F}_{q}}(\mathfrak{m}_{s,SC})$. L'orbite ${\cal O}$ est conservée par l'action galoisienne et on peut fixer $N\in {\cal O}\cap \mathfrak{m}_{s,SC}({\mathbb F}_{q})$. On peut aussi fixer un automorphisme $A$ de $V$ de sorte que $A\xi_{N}(Fr(z))A^{-1}=\xi_{N}(z)$ pour tout $z\in Z_{M_{s}}(N)/Z_{M_{s}}(N)^0$. On a l'égalité

(16) $\epsilon_{N}^{\flat}(g)=A\epsilon_{N}^{\flat}(Fr(g))A^{-1}$ pour tout $g\in Norm_{G(F^{nr})}(M,s,N)$.

\subsubsection{Cas stable}\label{epsilonflatcasstable}

On conserve les données du paragraphe précédent. On  suppose que ${\cal E}\in {\bf FC}^{st}(\mathfrak{m}_{SC,s})$. Cette hypothèse implique que la représentation $\xi_{N}$ est un caractère, autrement dit on a $V={\mathbb C}$. {\bf On suppose que} $M_{AD}$ {\bf est absolument simple}.   On a défini en \ref{constructionfamilleadmissible}  un caractère $\epsilon_{N}$ du groupe $K_{s,N}^{M_{AD},\dag,nr}$ (remarquons que la définition de ce caractère n'utilisait que la structure de $M_{AD}$ sur $F^{nr}$).   
 Le groupe $K_{s,N}^{M,\dag,nr}$ est contenu dans $Norm_{G(F^{nr})}(M,s,N)$ et s'envoie naturellement dans le groupe $K_{s,N}^{M_{AD},\dag,nr}$ et on note $\mu:K_{s,N}^{M,\dag,nr}\to K_{s,N}^{M_{AD},\dag,nr}$ cet homomorphisme naturel.

\begin{thm}{On a l'égalité $\epsilon_{N}^{\flat}(k)=\epsilon_{N}\circ \mu(k)$ pour tout $k\in K_{s,N}^{M,\dag,nr}$. }\end{thm}

La démonstration occupe la suite de la sous-section 4.3.  Faisons simplement ici quelques remarques préliminaires.  L'assertion de l'énoncé ne dépend pas du choix de l'élément $N\in {\cal O}$. En effet, remplaçons $N$ par $N'=\bar{g}N\bar{g}^{-1}$, où $g$ est un élément de $K_{s}^{M,0,nr}$. D'après \ref{caractereepsilonflat}(14)  appliqué à l'automorphisme $Ad(g)$ de $G$, on a $\epsilon_{N'}^{\flat}(k')=\epsilon_{N}^{\flat}(g^{-1}k'g)$ pour tout $k'\in K_{s,N'}^{M,\dag,nr}$. D'après 
l'hypothèse $Hyp(3,N,N')$ de \ref{familles}, on a aussi $\epsilon_{N'}\circ\mu(k')=\epsilon_{N}\circ\mu(g^{-1}k'g)$. L'égalité de l'énoncé est donc vérifiée pour $N$ si et seulement si elle l'est pour $N'$.

 En utilisant \ref{reduction}, on déduit  de $G$ des groupes $G'_{i}$, des Levi $M'_{i}$, des sommets $s_{i}$. L'hypothèse que $M_{AD}$ est absolument simple implique qu'il existe un indice $j$ tel que $M'_{j,AD}=M_{AD}$ tandis que $M'_{i}$ est un sous-tore maximal de $G'_{i}$ pour $i\not=j$. On a $\mathfrak{m}_{SC,s}= \mathfrak{m}'_{SC,j,s_{j}}$. 
 Le faisceau-caractère se décompose  en $\prod_{i=1,...,n}{\cal E}_{i}$ et on a ${\cal E}_{j}={\cal E}$ tandis que ${\cal E}_{i}$ est le faisceau constant sur $\mathfrak{m}'_{SC,i,s_{i}}=\{0\}$ pour $i\not=j$. On a $K_{s,N}^{M,\dag,nr}=\prod_{i=1,...,n}K_{s_{i},N_{i}}^{M'_{i},\dag,nr}$ et $\epsilon_{N}^{\flat}$ est égal au produit des $\epsilon_{N_{i}}^{\flat}$ sur ce groupe. L'assertion du théorème pour $G$ équivaut à la m\^eme assertion pour $G'_{j}$ et à l'assertion que, pour $i\not=j$, $\epsilon_{N_{i}}^{\flat}$ est trivial sur $K_{s_{i},N_{i}}^{M'_{i},\dag,nr}$. Cette dernière assertion résulte de \ref{faisceauxcaracteres}(10) et de  la construction de \ref{caractereepsilonflat}.  Cela nous ramène au cas où $G$ est absolument quasi-simple, ce que l'on suppose désormais. Remarquons que l'égalité de l'énoncé est vérifiée quand $k\in K_{s,N}^{M,0,nr}$:  les  deux membres de l'égalité sont égaux à $\xi_{N}(\bar{k})$, où $\bar{k}$ est la réduction de $k$ dans $Z_{M_{s}}(N)$. Le théorème est alors trivial quand $M=G$: on a $K_{s,N}^{M,\dag,nr}=K_{s,N}^{M,0,nr}$ puisque $M=G$ est simplement connexe. Supposons dorénavant que $M\not=G$. Il suffit de prouver l'égalité de l'énoncé pour $k$ appartenant à un sous-ensemble de $ K_{s,N}^{M,\dag,nr}$ dont l'image dans $ K_{s,N}^{M,0,nr}\backslash  K_{s,N}^{M,\dag,nr}$ engendre ce quotient. Un tel ensemble est fourni par le lemme \ref{descriptionKdag}. Ce lemme et l'assertion \ref{caractereepsilonflat}(15) nous permet de remplacer $G$ par $L_{SC}$, où $L$ est un $F^{nr}$-Levi de $G$ contenant $M$ comme $F^{nr}$-Levi propre maximal. Le groupe $L_{SC}$ n'est pas forcément absolument quasi-simple. S'il ne l'est pas, la maximalité de $M$ entraîne que $L_{SC}=L'\times M_{SC}$ où $L'$ est un groupe irréductible  sur $F^{nr}$ et $M_{sc}=T'\times M_{SC}$, où $T'$ est un $F^{nr}$-sous-tore maximal de $L'$. Comme on vient de le voir, l'égalité du théorème est vérifiée dans ce cas. Il reste à traiter le cas où $L_{SC}$ est absolument quasi-simple. On peut aussi bien remplacer $G$ par $L_{SC}$. En résumé, on peut supposer que $G$ est absolument quasi-simple et que $M$ est un $F^{nr}$-Levi propre maximal de $G$. Compte tenu des hypothèses
imposées à $M$  et puisqu'on a vu ci-dessus que l'on pouvait aussi exclure le cas où $M$ est un  tore, les possibilités sont les suivantes:

 $G$ est de type $(A_{n-1},ram)$, $M$ est de type $(A_{n-3},ram)$ avec   $n\geq5$;
 
 $G$ est de type $B_{n}$, $M$ est de type $B_{n-1}$ avec $n\geq 3$;
  
  $G$ est de type $C_{n}$, $M$ est de type $C_{n-1}$ avec $n \geq3$ et $n$ est impair;  
    
   $G$ est de type $(D_{n},nr)$,  $M$ est de type $(D_{n-1},nr)$ avec $n\geq 5$ et $n$ est impair;  
   
 $G$ est de type $(D_{n},ram)$, $M$ est de type $(D_{n-1},ram)$ avec $n\geq 10$ et $n$ est pair;  
 
  $G$ est de type $(E_{6},ram)$, $M$ est de type $(A_{5},ram)$;
 
  $G$ est de type $E_{7}$, $M$ est de type $E_{6}$ déployé;  

 $G$ est de type $E_{8}$, $M$ est de type $E_{7}$.
 
 On identifie $(M,s)$ à un couple $(M_{\Lambda},s_{\Lambda})$, où $\Lambda$ est un excellent sous-ensemble propre de $\Delta_{a}^{nr}$ tel que $\Delta_{a}^{nr}-\Lambda$ ait deux éléments. On note $\gamma_{1}$ et $\gamma_{2}$ les deux éléments de cet ensemble et on note $s_{1}$ et $s_{2}$ les sommets associés. On utilise l'assertion \ref{descriptionKdag}(9) qui décrit $ K_{s,N}^{M,0,nr}\backslash  K_{s,N}^{M,\dag,nr}$.  Dans  cette référence, on avait noté $w_{\gamma_{i}}$, pour $i\in \{1,2\}$,   une symétrie de $W^{aff}(\Lambda^{aff})$ provenant d'un élément de $K_{{\cal F}_{\gamma_{i}}}^{0,nr}$. Mais ${\cal F}_{\gamma_{i}}=s_{j}$ où $j$ est l'autre élément de $\{1,2\}$. Il est plus commode de noter $w_{j}=w_{\gamma_{i}}$.   Moyennant cette modification,  \ref{descriptionKdag}(9) 
 nous ramène à démontrer l'assertion suivante:
 
 (1) il existe  $N\in {\cal O}$ et,  pour $i=1,2$,  il existe un
 élément $n_{i}\in Norm_{G(F^{nr})}(M,s,N)\cap K_{s_{i}}^{0,nr}$ dont l'image dans $W^{aff}(\Lambda^{aff})$ soit la symétrie $w_{i}$, de sorte que l'égalité de l'énoncé soit vérifiée pour $k=n_{1}n_{2}$.
 
 C'est cette assertion qui sera vérifiée cas par cas dans les paragraphes suivants.  
 
  On a fix\'e un \'epinglage affine $\mathfrak{ E}_{a}$ de $G$. On doit en fixer un de $M_{AD}$. On a $T\subset M$ donc  l'appartement associ\'e \`a $T^{nr}$ dans $Imm_{F^{nr}}(G)$ est contenu dans $Imm_{F^{nr}}^G(M_{ad})$.  L'ensemble $p_{M}(C^{nr})$ est contenu dans une unique alc\^ove $C^{M,nr}$ de $Imm_{F^{nr}}(M_{AD})$. Notons $T^{M}=T/Z(M)$. Alors de $T^{M,nr}$ et de $C^{M,nr}$ se d\'eduit un ensemble de racines affines $\Delta_{a}^{M,nr}$.  
On conna\^{\i}t l'ensemble $S(\bar{C}^{M,nr})$ des sommets adh\'erents \`a $C^{M,nr}$: il  est \'egal \`a $\{s\}\cup\{p_{M}(s_{\beta}); \beta\in \Lambda\}$. Autrement dit, $\Delta_{a}^{M,nr}$ s'identifie \`a la r\'eunion de $\Lambda$ et de l'ensemble r\'eduit \`a la racine associ\'ee au sommet $s$, que l'on note $\beta_{s}$. Remarquons que, en num\'erotant les \'el\'ements de ces ensembles de racines affines comme dans \cite{B}, le num\'ero d'un \'el\'ement de $\Lambda$ n'est en g\'en\'eral pas le m\^eme selon qu'on le consid\`ere comme un \'el\'ement de $\Delta_{a}^{nr}$ ou de $\Delta_{a}^{M,nr}$. 
 On a d\'ej\`a fix\'e un \'el\'ement $E_{\beta}$ pour $\beta\in \Lambda$, qui appartient \`a $\mathfrak{m}$. On peut choisir  un \'el\'ement $E_{\beta_{s}}$ de sorte que $(T^{M,nr},(E_{\beta})_{\beta\in \Delta_{a}^{M,nr}})$ soit un \'epinglage affine pour $M_{AD}$.

  \subsubsection{Les cas faciles}\label{lescasfaciles}

    On suppose que $G$ est de type $C_{n}$, $(D_{n},nr)$ avec $n$ impair, $(D_{n},ram)$ avec $n$ pair où $E_{7}$. Alors l'ensemble $\Delta_{a}^{nr}-\Lambda$ est respectivement $\{\alpha_{(n-1)/2},\alpha_{(n+1)/2}\}$, $\{\alpha_{(n-1)/2},\alpha_{(n+1)/2}\}$, $\{\beta_{n/2-1},\beta_{n/2}\}$, $\{\alpha_{3},\alpha_{5}\}$. On introduit le groupe $\boldsymbol{\Omega}^{nr}$ relatif à $G$ et notre épinglage affine $\mathfrak{E}_{a}$. Rappelons que c'est un sous-groupe de $G_{AD}(F^{nr})$.  On voit qu'il existe un \'el\'ement $ g\in  \boldsymbol{\Omega}^{nr}$ qui \'echange les deux sommets de $\Delta_{a}^{nr}-\Lambda$.    Alors $Ad(g)$ est un automorphisme de $G$ défini sur $F^{nr}$ qui conserve $M=M_{\Lambda}$, $s=s_{\Lambda}$ et le faisceau ${\cal E}$   mais qui \'echange les deux \'el\'ements $\gamma_{1}$ et $\gamma_{2}$ de $\Delta_{a}^{nr}-\Lambda$.  Fixons $N\in {\cal O}$. En multipliant $g$ par un \'el\'ement convenable de $K_{s}^{M,0,nr}$, on  peut supposer que $Ad(g)_{s}$ conserve $N$. Fixons un élément $n_{1}\in Norm_{G(F^{nr})}(M,s,N)\cap K_{s_{1}}^{0,nr}$ dont l'image dans $W^{aff}(\Lambda^{aff})$ soit la symétrie $w_{1}$. Posons $n_{2}=Ad(g)(n_{1})^{-1}$. Il est clair que  $n_{2}\in Norm_{G(F^{nr})}(M,s,N)\cap K_{s_{2}}^{0,nr}$ et que son image dans $W^{aff}(\Lambda^{aff})$ est la symétrie $w_{2}$. En utilisant \ref{caractereepsilonflat}(14), on voit que $\epsilon_{N}^{\flat}(
 n_{1}n_{2})=1$.  Les automorphismes $Ad(n_{1})$, $Ad(n_{2})$ et  $Ad(g)$ conservent $M$ et se descendent  en des automorphismes $\nu_{1}$,$\nu_{2}$ et  $\delta$ de $M_{AD}$. Ils appartiennent tous trois à  $ Aut_{s,N}(M_{AD},F^{nr})$. On a $\mu(n_{1}n_{2})=\nu_{1}\nu_{2}$. On a aussi $\nu_{2}=\delta\nu_{1}^{-1}\delta^{-1}$. Alors le lemme \ref{prolongement} implique que $\epsilon_{N}\circ\mu(n_{1}n_{2})=1$. Cela démontre l'assertion (1) de \ref{epsilonflatcasstable}.

 \subsubsection{Preuve  du théorème \ref{epsilonflatcasstable} dans le cas $(A_{n-1},ram)$ avec $n$ impair}\label{preuveAn-1ram}
 On suppose que $G$ est de type $A_{n-1}$ avec $n$ impair et $n\geq 5$ et que $M$ est de type $A_{n-3}$. On suppose que $G$ n'est pas d\'eploy\'e sur $F^{nr}$ mais l'est sur l'extension quadratique  ramifi\'ee $E$ de $F^{nr}$.   
 L'hypoth\`ese $ {\bf FC}^{st}(\mathfrak{m}_{SC,s})\not=\emptyset$ entra\^{\i}ne   que $n-2=h^2+k(k+1)$ avec $h,k\in {\mathbb N}$ et $h=k$ ou $h=k+1$.  Remarquons que l'imparit\'e de $n$ entra\^{\i}ne celle de $h$. 
 Le groupe $G$ est un groupe spécial  unitaire.  

On fixe une uniformisante $\varpi_{E}$ de $E$ telle que $\varpi_{E}^2\in F^{nr,\times}$. On introduit un espace $V$ sur $E$ de dimension $n$ muni d'une base $(e_{i})_{i=1,...,n}$ et d'une forme hermitienne $q$ (relativement \`a l'extension $E/F^{nr}$) telle que $q(e_{i},e_{n+1-i})=(-1)^{i} $ pour $i=1,...,n$ et $q(e_{i},e_{j})=0$ pour $i+j\not=n+1$. Gr\^ace \`a cette base, on identifie $GL(V\otimes_{E}\bar{F})$ \`a $GL(n,\bar{F})$. On introduit l'\'epinglage   de $GL(n,\bar{F})$ que l'on note ici $\tilde{\mathfrak{E}}=(\tilde{B},\tilde{T},(E_{\alpha_{i}})_{i=1,...,n-1})$    et l'automorphisme $\theta$ d\'efinis au d\'ebut du paragraphe \ref{Bn}.   En notant $\sigma\mapsto \sigma_{GL(n)}$ l'action naturelle de $\Gamma_{F}$ sur $GL(n,\bar{F})$, on v\'erifie que le groupe spécial unitaire $G$ de $V$ s'identifie \`a $SL(n,\bar{F})$ muni de l'action $\sigma\mapsto \sigma_{G}$ telle que $\sigma_{G}=\sigma_{GL(n)}$ si $\sigma\in \Gamma_{E}$ et $\sigma_{G}=\theta\circ \sigma_{GL(n)}$ si $\sigma\in \Gamma_{F^{nr}}- \Gamma_{E}$.  On construit avec ces donn\'ees un \'epinglage 
affine $\mathfrak{ E}_{a}=(T^{nr},(E_{\beta})_{\beta\in \Delta^{nr}_{a}})$  de $G$   que nous n'avons pas besoin de pr\'eciser (on a $T=\tilde{T}\cap G$). Rappelons que $\Delta_{a}^{nr}=\{\beta_{0},\beta_{1,n-1},...,\beta_{(n-1)/2,(n+1)/2}\}$.   L'ensemble $\Delta_{a}^{nr}-\Lambda$ est \'egal \`a $\{\beta_{k(k+1)/2,n-k(k+1)/2},\beta_{k(k+1)/2+1,n-1-k(k+1)/2}\}$.  Le Levi   $M$   est le sous-groupe des \'el\'ements  de $G$ qui conservent les droites $Ee_{k(k+1)/2+1}$ et $Ee_{n-k(k+1)/2}$. Son groupe d\'eriv\'e 
 $M_{der}$ est le groupe sp\'ecial unitaire du sous-espace $V_{M}$ de $V$ engendr\'e par les vecteurs $e_{i}$ pour $i=1,...,n$, $i\not=k(k+1)/2+1, n-k(k+1)/2$.

 On pose $\gamma_{1}=\beta_{k(k+1)/2,n-k(k+1)/2}$ et $\gamma_{2}=\beta_{k(k+1)/2+1,n-1-k(k+1)/2}$. Consid\'erons l'un des sommets $s_{j}$ pour $j=1,2$. Si $j=1$, on pose $n''=k(k+1)$ et $n'=h^2+2$. Si $j=2$, on pose $n''=k(k+1)+2$, $n'=h^2$. 
 Notons ${\cal S}$ le $\mathfrak{o}_{E}$-r\'eseau de $V$ engendr\'e par les $e_{i}$ pour $i=1,...,n-n''/2$ et les  $\varpi_{E}e_{i}$ pour $n-n''/2+1,...,n$ et notons ${\cal S}^{\star}$ son dual, c'est-\`a-dire ${\cal S}^{\star}=\{v\in V; \forall v'\in {\cal S},\,q(v',v)\in \mathfrak{o}_{E}\}$. C'est le r\'eseau engendr\'e par les $\varpi_{E}^{-1}e_{i}$  pour $i=1,...,n''/2$ et les $e_{i}$ pour $i=n''/2+1,...,n$. Alors $\mathfrak{k}^{nr}_{s_{j}}$ est le sous-ensemble des $X\in \mathfrak{g}(F^{nr})$ tels que $X({\cal S})\subset {\cal S}$. L'espace $\bar{V}'={\cal S}/\varpi_{E}{\cal S}^{\star}$ sur $\bar{{\mathbb F}}_{q}$ est de dimension $n'$. Il est naturellement muni d'une base $\{\bar{e}_{i}; i=n''/2+1,..., n-n''/2\}$ o\`u $\bar{e}_{i}$ est la r\'eduction de $e_{i}$, et  d'une forme orthogonale. On note $G'_{s_{j}}$ son groupe sp\'ecial orthogonal. L'espace $\bar{V}''={\cal S}^{\star}/{\cal S}$ sur $\bar{{\mathbb F}}_{q}$ est de dimension $n''$. Il  est naturellement muni  d'une base $\{\bar{e}_{i}; i\in \{1,...,n''/2\}\cup\{n+1-n''/2,...,n\}\}$, o\`u $\bar{e}_{i}$ est la r\'eduction de $\varpi_{E}^{-1}e_{i}$ pour $i\in \{1,...,n''/2\}$ et $\bar{e}_{i}$ est la r\'eduction de $e_{i}$ si $i\in \{n+1-n''/2,...,n\}$. L'espace $\bar{V}''$ est muni d'une forme symplectique r\'eduction de $\varpi_{E}q$. On note $G''_{s_{j}}$ son groupe symplectique. On a l'\'egalit\'e $G_{s_{j}}=G'_{s_{j}}\times G''_{s_{j}}$. Le groupe $M_{s}$ se d\'ecompose conform\'ement en $M'_{s,j}\times M''_{s,j}$. 
 Si $j=1$, on a $M''_{s,1}=  G''_{s_{1}}$ et $M'_{s,1}$ est le Levi de $G'_{s_{1}}$ form\'e des \'el\'ements de ce groupe qui conservent les droites $\bar{{\mathbb F}}_{q}\bar{e}_{n''/2+1}$ et $\bar{{\mathbb F}}_{q}\bar{e}_{n-n''/2} $. Si $j=2$, on a $M'_{s,2}=G'_{s_{2}} $ et   $M''_{s,2}$ est  le Levi de $G''_{s_{2}}$ form\'e des \'el\'ements de ce groupe qui conservent les droites $\bar{{\mathbb F}}_{q}\bar{e}_{n''/2}$ et $\bar{{\mathbb F}}_{q}\bar{e}_{n+1-n''/2} $.  Evidemment, les groupes $M'_{s,j,SC}$ et $M''_{s,j,SC}$ ne dépendent pas de $j$. 
 On choisit un \'el\'ement $N\in {\cal O}$, qui se d\'ecompose en $N'\oplus N''$. 
 
 Supposons $j=2$. On note $n_{2}$ l'élément de  $G$   qui envoie $e_{k(k+1)/2+1}$ sur $\varpi_{E}e_{n-k(k+1)/2}$, $e_{n-k(k+1)/2}$ sur $-\varpi_{E}^{-1}e_{k(k+1)/2+1}$ et qui fixe $e_{i}$ pour $i\in \{1,...,n\}$, $i\not=k(k+1)/2+1, n-k(k+1)/2$. Cet \'el\'ement appartient \`a $K_{s_{2}}^{0,nr}\cap Norm_{G(F^{nr})}(M,s,N)$ et son image par $c_{s_{2}}$ est un élément de $Norm_{G_{s_{2}}}(M_{s})$ dont l'image  dans $W^{G_{s_{2}}}(M_{s})$ est l'\'el\'ement non trivial de ce groupe. Remarquons que l'action de $Ad(n_{2})$ sur $M$ est triviale. La r\'eduction $\bar{n}_{2}$ de $n_{2}$ dans $G_{s_{2}}$  se d\'ecompose en $(\bar{n}'_{2}, \bar{n}''_{2})$.  On a $\bar{n}'_{2}=1$ tandis que $\bar{n}''_{2}$ est l'\'el\'ement de $G''_{s_{2}}$ qui envoie $\bar{e}_{k(k+1)/2+1}$ sur $\bar{e}_{n-k(k+1)/2}$, $\bar{e}_{n-k(k+1)/2}$ sur $-\bar{e}_{k(k+1)/2+1}$ et qui fixe $\bar{e}_{i}$ pour les autres indices $i$. Dans \cite{W1} lemme VIII.9, on a calcul\'e $\epsilon_{s_{2},N}^{\flat}(\bar{n}_{2})$: ce terme vaut $(-1)^k$. On en d\'eduit
 
 (1)   $ \epsilon_{N}^{\flat}(n_{2})=(-1)^k$.
 
 Supposons $j=1$. Notons $\bar{V}'_{M}$ le sous-espace de $\bar{V}'$ engendr\'e par les vecteurs $\bar{e}_{i}$ pour $i\in \{k(k+1)/2+2,...,n-1-k(k+1)/2\}$ et notons $M'_{s,der}$ son groupe sp\'ecial orthogonal, qui est le groupe d\'eriv\'e de $M'_{s}$.
 L'\'el\'ement $N'$ appartient \`a $\mathfrak{m}'_{s,der}$ et 
 a pour partition associ\'ee $(2h-1,...,3,1)$. Donc le sous-espace $\bar{V}'_{M}$ se d\'ecompose en somme directe orthogonale de sous-espaces $\bar{V}'_{m}$, pour $m=1,...,h$, de sorte que $dim(\bar{V}'_{m})=2m-1$, que $\bar{V}'_{m}$ soit conserv\'e par $N'$ et que la restriction de $N'$ \`a $\bar{V}'_{m}$ soit un nilpotent r\'egulier. Notons $V'_{M}$ le sous-espace de $V$ engendr\'e par les vecteurs $e_{i}$ pour $i=k(k+1)/2+2,...,n-1-k(k+1)/2\}$. On peut relever la d\'ecomposition de $\bar{V}'_{M}$ en une d\'ecomposition orthogonale de $V'_{M}$ en sous-espaces $V'_{m}$ de sorte que ${\cal S}\cap V'_{M}=\oplus_{m=1,...,h}{\cal S}\cap V'_{m}$ et que $\bar{V}'_{m}$ soit la r\'eduction de ${\cal S}\cap V'_{m}$. On note $n_{1}$ l'élément de $G$ qui  permute $e_{k(k+1)/2+1}$ et $e_{n-k(k+1)/2}$, qui fixe $e_{i}$ pour $i\in \{1,...,k(k+1)/2\}\cup\{n+1-k(k+1)/2,...,n\}$, qui conserve chaque $V'_{m}$ et y agit par l'identit\'e si $m=1,...,h-1$ et par multiplication par $-1$ pour $m=h$. Cet \'el\'ement appartient \`a $K_{s_{1}}^{0,nr}\cap Norm_{G(F^{nr})}(M,s,N)$ et son image par $c_{s_{1}}$ est un élément de $Norm_{G_{s_{1}}}(M_{s})$ dont l'image  dans $W^{G_{s_{1}}}(M_{s})$ est l'\'el\'ement non trivial de ce groupe. La r\'eduction $\bar{n}_{1}$  de $n_{1}$ dans $G_{s_{1}}$ se d\'ecompose en $(\bar{n}_{1}', \bar{n}_{1}'')$.  On a $\bar{n}''_{1}=1$ tandis que $\bar{n}'_{1}$ est l'\'el\'ement de $G'_{s_{1}}$ qui permute $\bar{e}_{k(k+1)/2}$ et $\bar{e}_{n-k(k+1)/2}$, qui conserve chaque $\bar{V}'_{m}$ et y agit par l'identit\'e si $m\not=h$ et par multiplication par $-1$ pour $m=h$. Dans \cite{W1} lemme VIII.9, on a calcul\'e $\epsilon_{s_{1},N}^{\flat}(\bar{n}_{1})$: ce terme vaut $1$. D'où
 
 (2) $\epsilon_{N}^{\flat}(n_{1})=1$. 
 
 On sait que $n_{1}n_{2}$ appartient à $K_{s,N}^{M,\dag,nr}$. Pour calculer $\mu(n_{1}n_{2})$, il suffit de calculer l'action de $Ad(n_{1}n_{2})$ sur $M_{AD}$. On a déjà remarqué que celle de $Ad(n_{2})$ était triviale. 
  L'automorphisme $Ad(n_{1})_{M_{AD}}$  est la conjugaison par la restriction de $n_{1}$ \`a $V_{M}$.   Remarquons que $M_{AD,s}=M'_{s,1,der}\times M''_{s,1}$. La restriction de $n_{1}$ \`a $V_{M}$ n'appartient pas \`a $\tilde{M}_{der}$ car c'est un \'el\'ement de d\'eterminant $-1$. Il convient de la multiplier par la multiplication par $-1$ sur $V_{M}$ tout entier, ce qui ne change pas l'action par conjugaison de $n_{1}$ sur $M_{AD}$. L'\'el\'ement obtenu  s'envoie alors sur un élément  $y\in K_{s,N}^{M_{AD},0,nr}$. Notons $\bar{y}$ sa r\'eduction dans   $Z_{M_{AD,s}}(N)$.  Elle s'\'ecrit de fa\c{c}on \'evidente $\bar{y}=(\bar{y}',\bar{y}'')$.  On a $\bar{y}''=-1$  et  $\bar{y}'$ est l'\'el\'ement de $M'_{s,1,der}$ qui  agit trivialement sur $\bar{V}'_{h}$ et par multiplication par $-1$ sur $\bar{V}'_{m}$ pour $m=1,...,h-1$. On calcule facilement $\xi_{N}(\bar{y})=(-1)^{[(k+1)/2]+(h-1)/2}$. En se rappelant que $h=k$ ou $k+1$, on v\'erifie que $(-1)^{[(k+1)/2]+(h-1)/2}=(-1)^k$. Par d\'efinition, on a $\epsilon_{N}\circ \mu(n_{1}n_{2})=\xi_{N}(\bar{y})$. On en d\'eduit
 
  (3)   $\epsilon_{N}\circ \mu(n_{1}n_{2})=(-1)^k $.
  
  Les égalités (1), (2) et (3)  prouvent l'assertion (1) de \ref{epsilonflatcasstable}, donc le théorème.

 \subsubsection{Preuve  du théorème \ref{epsilonflatcasstable} dans le cas $(A_{n-1},ram)$ avec $n$ pair}\label{preuveAn-1rampair}
  On suppose que $G$ est de type $A_{n-1}$ avec $n$ pair et $n\geq 6$ et que $M$ est de type $A_{n-3}$. On suppose que $G$ n'est pas d\'eploy\'e sur $F^{nr}$ mais l'est sur l'extension quadratique  ramifi\'ee $E$ de $F^{nr}$.   
 L'hypoth\`ese  $ {\bf FC}^{st}(\mathfrak{m}_{SC,s})\not=\emptyset$ entra\^{\i}ne que $n-2=h^2+k(k+1)$ avec $h,k\in {\mathbb N}$ et $h=k$ ou $h=k+1$.  Remarquons que la parit\'e de $n$ entra\^{\i}ne celle de $h$ et que l'on a $n\geq8$. 
 Le groupe $G$ est un groupe spécial unitaire. En \ref{An-1ram}, on a r\'ealis\'e le groupe  unitaire tout entier comme celui d'un espace $(V,q)$ (le corps de base dans cette r\'ef\'erence \'etait $F$, on l'\'etend en $F^{nr}$). On reprend les d\'efinitions de ce paragraphe.  L'ensemble $\Delta_{a}^{nr}-\Lambda$ est \'egal \`a $\{\beta_{h^2/2,n-h^2/2},\beta_{h^2/2+1,n-1-h^2/2}\}$. Le Levi  $M$  est le sous-groupe des \'el\'ements  de $G$ qui conservent les droites $Ee_{h^2/2+1}$ et $Ee_{n-h^2/2}$.    Son groupe d\'eriv\'e $M_{der}$ est le groupe sp\'ecial unitaire du sous-espace $V_{M}$ de $V$ engendr\'e par les vecteurs $e_{i}$ pour $i=1,...,n$, $i\not=h^2/2+1, n-h^2/2$. On doit fixer un \'epinglage affine de $M$. Pour cela, on pose $e_{i}^M=e_{i}$ pour $i=1,...,h^2/2$, $e_{i}^M=-e_{i+1}$ pour $i=h^2/2+1,...,n-h^2/2-2$, $e_{i}^M=e_{i+2}$ pour $i=n-h^2/2-1,...,n-2$. Alors la base $(e_{i}^M)_{i=1,...,n-2}$ de $V_{M}$ v\'erifie les conditions de \ref{An-1ram} et on construit l'\'epinglage affine de $M$ comme dans ce paragraphe. Il v\'erifie les conditions de \ref{epsilonflatcasstable}, c'est-\`a-dire que, pour $\beta\in \Lambda$, les $E_{\beta}$ des \'epinglages affines de $G$ et de $M$ co\"{\i}ncident. 
  
 On pose $\gamma_{1}=\beta_{h^2/2,n-h^2/2}$ et $\gamma_{2}=\beta_{h^2/2+1,n-1-h^2/2}$. Consid\'erons  un sommet $s_{j}$ avec $j=1,2$. Si $j=1$, on pose $n'=h^2$ et $n''=k(k+1)+2$. Si $j=2$, on pose $n'=h^2+2$, $n''=k(k+1)$. 
 Notons ${\cal S}$ le $\mathfrak{o}_{E}$-r\'eseau de $V$ engendr\'e par les $e_{i}$ pour $i=1,...,n'/2$ et les  $\varpi_{E}e_{i}$ pour $n'/2+1,...,n$ et notons ${\cal S}^{\star}$ son dual, c'est-\`a-dire ${\cal S}^{\star}=\{v\in V; \forall v'\in {\cal S},\,q(v',v)\in \mathfrak{o}_{E}\}$. C'est le r\'eseau engendr\'e par les $e_{i}$  pour $i=1,...,n-n'/2$ et les $\varpi_{E}e_{i}$ pour $i=n-n'/2+1,...,n$. Alors $\mathfrak{k}^{nr}_{s}$ est le sous-ensemble des $X\in \mathfrak{g}(F^{nr})$ tels que $X({\cal S})\subset {\cal S}$. L'espace $\bar{V}'={\cal S}/\varpi_{E}{\cal S}^{\star}$ sur $\bar{{\mathbb F}}_{q}$ est de dimension $n'$. Il est naturellement muni d'une base $\{\bar{e}_{i}; i\in \{1,...,n'/2\}\cup\{n-n'/2+1,..., n\}\}$ o\`u $\bar{e}_{i}$ est la r\'eduction de $e_{i}$ si $i\in \{1,...,n'/2\}$, de $\varpi_{E}e_{i}$ si $i\in \{n-n'/2+1,...,n\}$, et  d'une forme orthogonale. On note $G'_{s_{j}}$ son groupe sp\'ecial orthogonal. L'espace $\bar{V}''={\cal S}^{\star}/{\cal S}$ sur $\bar{{\mathbb F}}_{q}$ est de dimension $n''$. Il  est naturellement muni  d'une base $\{\bar{e}_{i}; i\in \{n'/2+1,...,n-n'/2\}\} $, o\`u $\bar{e}_{i}$ est la r\'eduction de $ e_{i}$. L'espace $\bar{V}''$ est  muni d'une forme symplectique r\'eduction de $\varpi_{E}q$. On note $G''_{s_{j}}$ son groupe symplectique. On v\'erifie que $G_{s_{j}}=G'_{s_{j}}\times G''_{s_{j}}$. Le groupe $M_{s}$ se d\'ecompose conform\'ement en $M'_{s,j}\times M''_{s,j}$. 
 Si $j=1$, on a $M'_{s,1}=  G'_{s_{1}}$ et $M''_{s,1}$ est le Levi de $G''_{s_{1}}$ form\'e des \'el\'ements de ce groupe qui conservent les droites $\bar{{\mathbb F}}_{q}\bar{e}_{n'/2+1}$ et $\bar{{\mathbb F}}_{q}\bar{e}_{n-n'/2} $. Si $j=2$, on a $M''_{s,2}=G''_{s_{2}} $ et   $M'_{s,2}$ est  le Levi de $G'_{s_{2}}$ form\'e des \'el\'ements de ce groupe qui conservent les droites $\bar{{\mathbb F}}_{q}\bar{e}_{n'/2}$ et $\bar{{\mathbb F}}_{q}\bar{e}_{n+1-n'/2} $.  Les groupes d\'eriv\'es de $M'_{s,j}$ et $M''_{s,j}$ ne dépendent pas de $j$, on les note
 $M'_{s,der}$ et $M''_{s,der}$  . On v\'erifie que $M_{AD,s}=(M'_{s,der}\times M''_{s,der})/\{\pm 1\}$, où $\{\pm 1\}$ s'envoie diagonalement dans les centres des deux facteurs. 
   D'apr\`es la description de l'\'el\'ement $\omega\in \boldsymbol{\Omega}^{M,nr}$ que l'on a donn\'ee en \ref{An-1ram}, l'action $Ad(\omega)$ se r\'eduit en une action $ \boldsymbol{\omega}$ sur le groupe $M_{AD, s}$ qui  est triviale sur $M''_{s,der}$ et qui co\"{\i}ncide sur $M'_{s,der}$ avec l'action adjointe de l'\'el\'ement du groupe orthogonal de $\bar{V}'$ qui permute $\bar{e}_{1}$ et $\bar{e}_{n}$ et fixe les autres vecteurs de base.  On choisit un \'el\'ement $N\in {\cal O}\cap \mathfrak{r}$, qui se d\'ecompose en $N'\oplus N''$. 
 La condition que $N$ appartienne \`a $\mathfrak{r}$ \'equivaut \`a ce que $N'$ soit  fix\'e par  $\boldsymbol{\omega}$.  
   
  Supposons $j=1$. Notons  $n_{1}$   l'\'el\'ement de $G$ qui envoie $e_{h^2/2+1}$ sur $e_{n-h^2/2}$, $e_{n-h^2/2}$ sur $-e_{h^2/2+1}$ et qui fixe $e_{i}$ pour $i\in \{1,...,n\}$, $i\not=h^2/2+1, n-h^2/2$. L'action de $Ad(n_{1})$ sur $M$ est triviale donc $Ad(n_{1})_{s}$ fixe $N$.
Alors    $n_{1}$ appartient \`a $K_{s_{j}}^{0,nr}\cap Norm_{G(F^{nr})}(M,s,N)$ et son image par $c_{s_{1}}$ dans $Norm_{G_{s_{1}}}(M_{s})$ s'envoie sur l'élément  non trivial de $W^{G_{s_{1}}}(M_{s})$.  La r\'eduction $\bar{n}_{1}$ de $n_{1}$ dans $G_{s_{1}}$  se d\'ecompose en $(\bar{n}'_{1}, \bar{n}''_{1})$.  On a $\bar{n}'_{1}=1$ tandis que $\bar{n}''_{1}$ est l'\'el\'ement de $G''_{s_{1}}$ qui envoie $\bar{e}_{h^2/2+1}$ sur $\bar{e}_{n-h^2/2}$, $\bar{e}_{n-h^2/2}$ sur $-\bar{e}_{h^2/2+1}$ et qui fixe $\bar{e}_{i}$ pour les autres indices $i$. Dans \cite{W1} lemme VIII.9, on a calcul\'e $\epsilon^{\flat}_{s_{1},N}(\bar{n}_{1})$: ce terme vaut $(-1)^k$. On en d\'eduit
 
 (1) $ \epsilon_{N}^{\flat}(n_{1})=(-1)^k$.
 
  Supposons $j=2$. Notons $\bar{V}'_{M}$ le sous-espace de $\bar{V}'$ engendr\'e par les vecteurs $\bar{e}_{i}$ pour $i\in \{1,...,h^2/2\}\cup\{n+1-h^2/2,...,n\}$. Le groupe $M'_{s,der}$ est son groupe sp\'ecial orthogonal. L'\'el\'ement $N'$ appartient \`a $\mathfrak{m}'_{s,der}$ et 
 a pour partition associ\'ee $(2h-1,...,3,1)$. Donc le sous-espace $\bar{V}'_{M}$ se d\'ecompose en somme directe orthogonale de sous-espaces $\bar{V}'_{m}$, pour $m=1,...,h$, de sorte que $dim(\bar{V}'_{m})=2m-1$, que $\bar{V}'_{m}$ soit conserv\'e par $N'$ et que la restriction de $N'$ \`a $\bar{V}'_{m}$ soit un nilpotent r\'egulier. Notons $V'_{M}$ le sous-espace de $V$ engendr\'e par les vecteurs $e_{i}$ pour $i\in \{1,...,h^2/2\}\cup\{n+1-h^2/2,...,n\}$. On peut relever la d\'ecomposition de $\bar{V}'_{M}$ en une d\'ecomposition orthogonale de $V'_{M}$ en sous-espaces $V'_{m}$ de sorte que ${\cal S}\cap V'_{M}=\oplus_{m=1,...,h}{\cal S}\cap V'_{m}$ et que $\bar{V}'_{m}$ soit la r\'eduction de ${\cal S}\cap V'_{m}$. On note $n_{2}$ l'élément de  $G$   qui  permute $e_{h^2/2+1}$ et $\varpi_{E}e_{n-h^2/2}$, qui fixe $e_{i}$ pour $i\in \{ h^2/2+2,...,n-1-h^2/2\}$, qui conserve chaque $V'_{m}$ et y agit par l'identit\'e si $m=1,...,h-1$ et par multiplication par $-1$ pour $m=h$. Cet \'el\'ement appartient \`a $K_{s_{2}}^{0,nr}\cap Norm_{G(F^{nr})}(M,s,N)$ et son image par $c_{s_{2}}$ dans $Norm_{G_{s_{2}}}(M_{s})$ s'envoie sur l'élément  non trivial de $W^{G_{s_{2}}}(M_{s})$.  
  La r\'eduction $\bar{n}_{2}$  se d\'ecompose en $(\bar{n}'_{2}, \bar{n}''_{2})\in G'_{s_{2}}\times G''_{s_{2}}$.  On a $\bar{n}''_{2}=1$ tandis que $\bar{n}'_{2}$ est l'\'el\'ement de $G'_{s_{2}}$ qui permute $\bar{e}_{h^2/2+1}$ et $\bar{e}_{n-h^2/2}$, qui conserve chaque $\bar{V}'_{m}$ et y agit par l'identit\'e si $m\not=h$ et par multiplication par $-1$ pour $m=h$. Dans \cite{W1} lemme VIII.9, on a calcul\'e $\epsilon^{\flat}_{s_{2},N}(\bar{n})$: ce terme vaut $1$.  D'où
  
  (2) $ \epsilon_{N}^{\flat}(n_{2})=1$. 
  
  Pour calculer $\mu(n_{1}n_{2})$, il suffit de calculer l'action de $Ad(n_{1}n_{2})$ sur $M_{AD}$. On a déjà remarqué que celle de $Ad(n_{1})$ était triviale. 
 L'automorphisme $Ad(n_{2})_{M_{AD}}$  est la conjugaison par la restriction de $n_{2}$ \`a $V_{M}$.  On a d\'ecrit ci-dessus l'\'el\'ement $\boldsymbol{\omega}$. C'est une sym\'etrie  orthogonale \'el\'ementaire  de $\bar{V}'$ qui fixe $N'$. Or la seule  symétrie qui vérifie cette condition est celle qui fixe $\bar{V}'_{m}$ pour $m\geq2$ et agit sur $\bar{V}'_{1}$ par multiplication par $-1$.   On voit alors que $Ad(n_{2})_{M_{AD}}$ co\"{\i}ncide avec $Ad(y) \boldsymbol{\omega}$, o\`u $y$ est un \'el\'ement de $K_{s,N}^{M_{AD},0,nr}$. Notons $\bar{y}$ sa r\'eduction, qui est l'image dans $M_{AD,s}$ d'un couple $(\bar{y}',\bar{y}'')$.   On a $\bar{y}''=1$  tandis que  $\bar{y}'$ est l'\'el\'ement de $M'_{s,der}$ qui  conserve les espaces $\bar{V}'_{m}$ et y agit par l'identit\'e pour $m\not=1,h$ et par multiplication par $-1$ pour $m=1,h$. On calcule ais\'ement $\xi_{N}(\bar{y})=-1$. Par d\'efinition, on a $\epsilon_{N}\circ \mu(n_{1}n_{2})=\xi_{N}(\bar{y})\boldsymbol{\epsilon}(\boldsymbol{\omega})$ et $\boldsymbol{\epsilon}(\boldsymbol{\omega})=(-1)^{k+1}$, d'o\`u 
 
 (3) $\epsilon_{N}\circ \mu(n_{1}n_{2})=(-1)^k$.

   Les égalités (1), (2) et (3)  prouvent l'assertion (1) de \ref{epsilonflatcasstable}, donc le théorème.
   
  \subsubsection{Preuve  du théorème \ref{epsilonflatcasstable} dans le cas $B_{n}$}\label{preuveBn}
  
  On suppose que $G$ est de type $B_{n}$ avec $n\geq 3$ et que $M$ est de type $B_{n-1}$.  L'hypoth\`ese $ {\bf FC}^{st}(\mathfrak{m}_{SC,s})\not=\emptyset$ entra\^{\i}ne  que $2n-1=k^2+h^2$ avec $k,h\in {\mathbb N}$, $k$ pair et $\vert k-h\vert =1$.  Le groupe $G$ est un groupe spinoriel. On se rappelle que le faisceau ${\cal E}$ est $M_{AD}$-équivariant. Autrement dit, pour $N\in {\cal O}$, le caractère $\xi_{N}$ est trivial sur l'image de $Z(M)$ dans $Z_{M_{s}}(N)$, a fortiori sur l'image de $Z(G)$. Cela entra\^{\i}ne que les fonctions $\epsilon_{N}^{\flat}$ et $\epsilon_{N}\circ \mu$ sont invariantes par $Z(G)$. On peut donc remplacer le groupe $G$ par son groupe adjoint $G_{AD}$, à condition que les éléments $n_{1}$ et $n_{2}$ que nous devons construire  appartiennent à l'image dans ce groupe de $G(F^{nr})$. Pour ce paragraphe, nous supposons donc que $G$ est adjoint, autrement dit que c'est un groupe spécial orthogonal.

  En \ref{Bn}, on a r\'ealis\'e $G$ comme le groupe sp\'ecial orthogonal d'un espace quadratique $(V,q)$ (le corps de base dans cette r\'ef\'erence \'etait $F$, on l'\'etend en $F^{nr}$). On reprend les d\'efinitions de ce paragraphe. 
   L'ensemble $\Delta_{a}-\Lambda$ est \'egal \`a $\{\alpha_{k^2/2},\alpha_{k^2/2+1}\}$. Le Levi   $M$    est le sous-groupe des \'el\'ements  de $G$ qui conservent les droites $F^{nr}e_{k^2/2+1}$ et $F^{nr}e_{2n+1-k^2/2}$. Son groupe d\'eriv\'e $\tilde{M}_{der}$ est le groupe sp\'ecial orthogonal du sous-espace $V_{M}$ de $V$ engendr\'e par les vecteurs $e_{i}$ pour $i=1,...,2n+1$, $i\not=k^2/2+1, 2n+1-k^2/2$. On doit 
 fixer un \'epinglage affine de $M$.   Pour cela, on pose $e_{i}^M=e_{i}$ pour $i=1,...,k^2/2$, $e_{i}^M=e_{i+1}$ pour $i=k^2/2+1,...,2n-1-k^2/2$, $e_{i}^M=e_{i+2}$ pour $i=2n-k^2/2,...,2n-1$. Alors la base $(e_{i}^M)_{i=1,...,2n-1}$ de $V^M$ v\'erifie la condition de \ref{Bn} et on construit un \'epinglage affine de $M$ comme dans ce paragraphe.   
     
   On pose $\gamma_{1}=\alpha_{k^2/2}$ et $\gamma_{2}=\alpha_{k^2/2+1}$. Consid\'erons le sommet $s_{j}$ pour $j=1,2$. Si $j=1$, on pose $n'=(h^2+1)/2$ et $n''=k^2/2$. Si $j=2$, on pose $n'=(h^2-1)/2$, $n''=k^2/2+1$. 
 Notons ${\cal S}$ le $\mathfrak{o}_{F^{nr}}$-r\'eseau de $V$ engendr\'e par les $e_{i}$ pour $i=1,...,2n+1-n''$ et les  $\varpi_{F}e_{i}$ pour $2n+2-n'',...,n$ et notons ${\cal S}^{\star}$ son dual, c'est-\`a-dire ${\cal S}^{\star}=\{v\in V; \forall v'\in {\cal S},\,q(v',v)\in \mathfrak{o}_{F^{nr}}\}$. C'est le r\'eseau engendr\'e par les $\varpi_{F}^{-1}e_{i}$  pour $i=1,...,n''/2$ et les $e_{i}$ pour $i=n''+1,...,2n+1$. Alors $\mathfrak{k}^{nr}_{s_{j}}$ est le sous-ensemble des $X\in \mathfrak{g}(F^{nr})$ tels que $X({\cal S})\subset {\cal S}$. L'espace $\bar{V}'={\cal S}/\varpi_{F}{\cal S}^{\star}$ sur $\bar{{\mathbb F}}_{q}$ est de dimension $2n'+1$. Il est naturellement muni d'une base $\{\bar{e}_{i}; i =n''+1,...,2n+1-n''\}$ o\`u $\bar{e}_{i}$ est la r\'eduction de $e_{i}$. Il est muni  d'une forme quadratique. On note $G'_{s_{j}}$ son groupe spécial orthogonal. L'espace $\bar{V}''={\cal S}^{\star}/{\cal S}$ sur $\bar{{\mathbb F}}_{q}$ est de dimension $2n''$. Il  est naturellement muni  d'une base $\{\bar{e}_{i}; i\in \{1,...,n''\}\cup\{2n-n''+2,...,2n+1\} \}$, o\`u $\bar{e}_{i}$ est la r\'eduction de $ \varpi_{F}^{-1}e_{i}$ pour $i\in \{1,...,n''\}$ et de $e_{i}$ pour $i\in \{2n-n''+2,...,2n+1\}$. L'espace $\bar{V}''$ est muni d'une forme  quadratique r\'eduction de $\varpi_{F}q$. On note $G''_{s_{j}}$ son groupe spécial orthogonal. On v\'erifie que $G_{s_{j}}=G'_{s_{j}}\times G''_{s_{j}}$. Le groupe $M_{s}$ se d\'ecompose conform\'ement en $M'_{s,j}\times M''_{s,j}$. 
 Si $j=1$, on a $M''_{s,1}=  G''_{s_{1}}$ et $M'_{s,1}$ est le Levi de $G'_{s,1}$ form\'e des \'el\'ements de ce groupe qui conservent les droites $\bar{{\mathbb F}}_{q}\bar{e}_{n''+1}$ et $\bar{{\mathbb F}}_{q}\bar{e}_{2n+1-n''} $. Si $j=2$, on a $M'_{s,2}=G'_{s_{2}} $ et   $M''_{s,2}$ est  le Levi de $G''_{s_{2}}$ form\'e des \'el\'ements de ce groupe qui conservent les droites $\bar{{\mathbb F}}_{q}\bar{e}_{n''}$ et $\bar{{\mathbb F}}_{q}\bar{e}_{2n+2-n''} $.  Les groupes dérivés de $M'_{s,j}$ et $M''_{s,j}$ ne dépendent pas de $j$, on les note
   $M'_{s,der}$ et $M''_{s,der}$.  Plus précisément, les espaces $\bar{V}'$, resp. $\bar{V}''$, définis ci-dessus dépendent de $j$ mais leurs sous-espaces $\bar{V}'_{M}$  engendr\'e par les vecteurs $\bar{e}_{i}$ pour $i\in \{ k^2/2+2,...,2n-k^2/2\}$, resp. $\bar{V}''_{M}$   engendr\'e par les vecteurs $\bar{e}_{i}$ pour $i\in \{1,...,k^2/2\}\cup \{2n+2-k^2/2,...,2n+1\}$, peuvent s'identifier. Alors $M'_{s,der}$ et $M''_{s,der}$ sont les groupes spéciaux orthogonaux de ces espaces.    On v\'erifie que $M_{AD,s}=\tilde{M}'_{s,der}\times \tilde{M}''_{s,der}$. 
 
  On choisit un \'el\'ement $N\in {\cal O}$, qui se d\'ecompose en $N'\oplus N''$. On imposera ci-dessous une condition suppl\'ementaire \`a $N''$.   L'\'el\'ement $N'$ appartient \`a $\mathfrak{m}'_{s,der}$ et 
 a pour partition associ\'ee $(2h-1,...,3,1)$. Donc le sous-espace $\bar{V}'_{M}$ se d\'ecompose en somme directe orthogonale de sous-espaces $\bar{V}'_{m}$, pour $m=1,...,h$, de sorte que $dim(\bar{V}'_{m})=2m-1$, que $\bar{V}'_{m}$ soit conserv\'e par $N'$ et que la restriction de $N'$ \`a $\bar{V}'_{m}$ soit un nilpotent r\'egulier.   L'\'el\'ement $N''$ appartient \`a $\mathfrak{m}''_{s_{M},der}$ et 
 a pour partition associ\'ee $(2k-1,...,3,1)$. Donc le sous-espace $\bar{V}''_{M}$ se d\'ecompose en somme directe orthogonale de sous-espaces $\bar{V}''_{m}$, pour $m=1,...,k$, ayant des propri\'et\'es analogues \`a ci-dessus. On peut imposer et on impose  \`a $N''$ la condition: $\bar{V}''_{1}$ est la droite port\'ee par $\bar{e}_{1}-\bar{e}_{2n+1}$. Pour $j\in \{1,2\}$, on a défini le réseau ${\cal S}$. L'intersection ${\cal S}\cap V_{M}$ ne dépend pas de $j$. 
 On peut relever ces d\'ecompositions de $\bar{V}'_{M}$  et $\bar{V}''_{M}$ en une d\'ecomposition orthogonale de $V_{M}$ en sous-espaces $V'_{m}$ et $V''_{m}$ de sorte que ${\cal S}\cap V_{M}=(\oplus_{m=1,...,h}{\cal S}\cap V'_{m})\oplus(\oplus_{m=1,...,k}{\cal S}\cap V''_{m})$, que $\bar{V}'_{m}$ soit la r\'eduction de ${\cal S}\cap V'_{m}$ et que $\bar{V}''_{m}$ soit la r\'eduction de ${\cal S}^{\star}\cap V''_{m}$. On peut imposer et on impose la condition: $V''_{1}$ est la droite port\'ee par $e_{1}-\varpi_{F}e_{2n+1}$.    D'apr\`es la description de l'\'el\'ement $\omega\in \boldsymbol{\Omega}^{M}$ que l'on a donn\'ee en \ref{Bn}, cet \'el\'ement conserve les sous-espaces $V'_{m}$ et $V''_{m}$ et y agit par multiplication par $-1$, sauf sur $V''_{1}$ o\`u il agit par l'identit\'e. 
   L'action $Ad(\omega)$ se r\'eduit en une action $ \boldsymbol{\omega}$ sur le groupe $M_{AD, s}$ qui  est triviale sur $M'_{s,der}$ et qui co\"{\i}ncide sur $M''_{s,der}$ avec l'action adjointe de l'\'el\'ement du groupe orthogonal de $\bar{V}''$ qui conserve chaque $\bar{V}''_{m}$ et y agit par l'identit\'e pour $m\geq2$ et par multiplication par $-1$ pour $m=1$. Cet \'el\'ement fixe $N$ donc $N\in \mathfrak{r}$.   
 
   Supposons $j=1$.  Notons $ n_{1}$  l'\'el\'ement de $G$ qui permute $e_{k^2/2+1}$ et $e_{2n+1-k^2/2}$, conserve chaque espace $V'_{m}$ et $V''_{m}$ et y agit par l'identit\'e, sauf sur $V'_{h}$ o\`u il agit par multiplication par $-1$.  Il appartient  bien à l'image de $G_{SC}(F^{nr})$ dans $G$.     On a $n_{1}\in K_{s_{1}}^{0,nr}\cap Norm_{G(F^{nr})}(M,s,N)$ et l'image naturelle de $n_{1}$ dans $Norm_{G_{s_{1}}}(M_{s})$ s'envoie sur l'élément non trivial de $W^{G_{s_{1}}}(M_{s})$. 
    La r\'eduction $\bar{n_{1}}$  se d\'ecompose en $(\bar{n}'_{1}, \bar{n}''_{1})$.  On a $\bar{n}''_{1}=1$ tandis que $\bar{n}'_{1}$ est l'\'el\'ement de $G'_{s_{1}}$ qui permute  $\bar{e}_{k^2/2+1}$ et $\bar{e}_{2n+1-k^2/2}$, conserve les $\bar{V}'_{m}$ et y agit par l'identit\'e sauf sur $\bar{V}'_{h}$ o\`u il agit par multiplication par $-1$. Dans \cite{W1} lemme VIII.9, on a calcul\'e $\epsilon^{\flat}_{s_{1}, N}(\bar{n}_{1})$: ce terme vaut $1$. D'où
   
   (1) $\epsilon_{N}^{\flat}(n_{1})=1$.

 Supposons $j=2$. On note $n_{2}$  l'\'el\'ement de $G$ qui permute $e_{k^2/2+1}$ et $\varpi_{F}e_{2n+1-k^2/2}$, conserve chaque espace $V'_{m}$ et $V''_{m}$ et y agit par l'identit\'e, sauf sur $V''_{k}$ o\`u il agit par multiplication par $-1$.   Il appartient  bien à l'image de $G_{SC}(F^{nr})$ dans $G$.  On a  $n_{2}\in K_{s_{2}}^{0,nr}\cap Norm_{G(F^{nr})}(M,s,N)$  et  l'image naturelle de $n_{2}$ dans $Norm_{G_{s_{2}}}(M_{s})$ s'envoie sur l'élément non trivial de $W^{G_{s_{2}}}(M_{s})$.  La r\'eduction $\bar{n}_{2}$  se d\'ecompose en $(\bar{n}'_{2}, \bar{n}''_{2})$.  On a $\bar{n}'_{2}=1$ tandis que $\bar{n}''_{2}$ est l'\'el\'ement de $G''_{s_{2}}$ qui permute  $\bar{e}_{k^2/2+1}$ et $\bar{e}_{2n+1-k^2/2}$, conserve les $\bar{V}''_{m}$ et y agit par l'identit\'e sauf sur $\bar{V}''_{k}$ o\`u il agit par multiplication par $-1$. Dans \cite{W1} lemme VIII.9, on a calcul\'e $\epsilon^{\flat}_{s_{2},N}(\bar{n})$: ce terme vaut $1$. D'où
 
 (2) $\epsilon_{N}^{\flat}(n_{2})=1$. 
 
 Le produit $n_{1}n_{2}$ multiplie $e_{k^2/2+1}$ par $\varpi_{F}$ et $e_{2n+1-k^2/2}$ par $\varpi_{F}^{-1}$, il conserve les espaces $V'_{m}$ et $V''_{m}$ en y agissant par l'identité, sauf pour $V'_{h}$ et $V''_{k}$ où il agit par multiplication par $-1$. C'est bien un élément de $K_{s,N}^{\dag,nr}$ mais son image $\mu(n_{1}n_{2})$ dans $M_{AD}$ n'appartient pas à $K_{s,N}^{M_{AD},0,nr}$: ses réductions agissent dans $\bar{V}'_{M}$ et $\bar{V}''_{M}$ par des éléments de déterminant $-1$. Notons $y=\mu(n_{1}n_{2})\omega$. Alors $y$ appartient à $K_{s,N}^{M_{AD},0,nr}$. L'action de $\bar{y}$ conserve les $\bar{V}'_{m}$ et $\bar{V}''_{m}$ et il agit par multiplication sur $-1$ sur ces sous-espaces, à l'exception de $V'_{h}$, $V''_{k}$ et $V''_{1}$ où il agit par l'identité. On calcule $\xi_{N}(\bar{y})=(-1)^{(h-1)/2+k/2-1}$. On a aussi par définition $\epsilon_{N}(\omega)=(-1)^{(h+k+1)/2}$, d'où
 
 (3) $\epsilon_{N}\circ\mu(n_{1}n_{2})=1$. 
 
   Les égalités (1), (2) et (3)  prouvent l'assertion (1) de \ref{epsilonflatcasstable}, donc le théorème.
   
    \subsubsection{Début de la preuve  du théorème \ref{epsilonflatcasstable}  pour $G$ de type $E_{8}$\label{preuveE8}}
  On suppose $G$ de type $E_{8}$ et $M$ de type $E_{7}$. On note $(\alpha_{i})_{i=0,...,8}$ les \'el\'ements de $\Delta_{a}$.  L'ensemble $\Lambda$ est $\Delta_{a}-\{\alpha_{3},\alpha_{6}\}$, cf. \ref{couplesstables}. On doit d'abord identifier l'ensemble $\Delta_{a}^{M,nr}$ introduit en \ref{epsilonflatcasstable}. Posons
  $(\beta_{0},\beta_{1},\beta_{2},\beta_{3},\beta_{5},\beta_{6},\beta_{7})=(\alpha_{0},\alpha_{8},\alpha_{1},\alpha_{7},\alpha_{2},\alpha_{4},\alpha_{5})$ et $\beta_{4}=\alpha_{3}+\alpha_{4}+\alpha_{5}+\alpha_{6}$. On vérifie que
  $$(1) \qquad \beta_{0}=-(2\beta_{1}+2\beta_{2}+3\beta_{3}+4\beta_{4}+3\beta_{5}+2\beta_{6}+\beta_{7}).$$
  Les racines de $M$ sont celles qui sont combinaisons linéaires à coefficients rationnels d'elements de $\Lambda$. C'est clairement le cas des racines $\beta_{i}$ pour $i\not=4$ et la relation (1) entraîne que c'est aussi le cas de $\beta_{4}$.   On voit que les racines $\beta_{i}$ pour $i=1,...,7$ vérifient  les m\^emes relations de produits scalaires qu'une base du système $\Sigma^M$ de type $E_{7}$ de $M$. Il en résulte que cet ensemble de racines est une base de $\Sigma^M$. La relation (1) entraîne que $\beta_{0}$ est l'opposé de la plus grande racine pour cette base. Enfin, les racines affines $\beta_{i}$ pour $i=1,...,7$ et $\beta_{0}+1=\alpha_{0}+1$ sont positives sur $C^{nr}$ donc le sont aussi sur $C^{M,nr}$. Par définition de $\Delta_{a}^{M,nr}$, on peut donc supposer que $\Delta_{a}^{M,nr}=\{\beta_{i};i=0,...,7\}$.  La racine $\beta_{s}$ est $\beta_{4}$.

   Le groupe $\boldsymbol{Aut}({\cal D}_{a}^{M})=\boldsymbol{\Omega}^M$ a deux \'el\'ements. L'\'el\'ement non trivial $\boldsymbol{\omega}$ fixe   $E_{\beta_{2}}$ et $E_{\beta_{4}}$ et permute $E_{\beta_{0}}$ et $E_{\beta_{7}}$, $E_{\beta_{1}}$ et $E_{\beta_{6}}$, $E_{\beta_{3}}$ et $E_{\beta_{5}}$.
   On choisit $N=\sum_{\beta\in \Lambda}\bar{E}_{\beta}$, o\`u les $\bar{E}_{\beta}$ sont les r\'eductions des $E_{\beta}$ dans $\mathfrak{m}_{AD, s}$. Cet \'el\'ement est fixe par la r\'eduction de $\boldsymbol{\omega}$ donc appartient \`a $\mathfrak{r}$.  On a calcul\'e le groupe $M_{AD,s}$ en \cite{W7} 8.14. On a $M_{AD,s}=(SL(4)\times SL(4)\times SL(2))/ \bar{\zeta}_{1/4}$, o\`u $\bar{\zeta}_{1/4}$ est le groupe des racines $4$-i\`emes de l'unit\'e dans $\bar{{\mathbb F}}_{q}$. Fixons un g\'en\'erateur $\underline{i}$ de ce groupe. Posons  $z_{1}=\check{\beta}_{0}(\underline{i})\check{\beta}_{1}(\underline{i}^2)\check{\beta}_{3}(\underline{i}^3)$, $z_{2}=\check{\beta}_{7}(\underline{i})\check{\beta}_{6}(\underline{i}^2)\check{\beta}_{5}(\underline{i}^3)$. Le 
    plongement $ \bar{\zeta}_{1/4}\to SL(4)\times SL(4)\times SL(2) $ est
   $$\underline{i}\mapsto ( z_{1},z_{2}, \check{\beta}_{2}(-1)).$$
   Le groupe $Z_{M_{AD,s}}(N)/Z_{M_{AD,s}}(N)^0$ s'identifie au centre de $M_{AD,s}$, lequel est engendr\'e par $z_{1}$, $z_{2}$ et $\check{\beta}_{2}(-1)$ (ou seulement $z_{1}$ et $\check{\beta}_{2}(-1)$ puisque $z_{1}z_{2}\check{\beta}_{2}(-1)$ se projette sur $1$ dans $M_{AD,s}$). Il y a deux faisceaux-caract\`eres cuspidaux \`a support dans l'orbite ${\cal O}$ de $N$ et ${\cal E}$ est l'un d'eux. En notant ici $i$ l'\'el\'ement habituel de ${\mathbb C}^{\times}$, il existe $\eta=\pm 1$ tel que le caract\`ere $\xi_{N}$ v\'erifie $\xi_{N}(z_{1})=\xi_{N}(z_{2})=\eta i$, $\xi_{N}(\check{\beta}_{2}(-1))=-1$.  
   
 En  \ref{epsilonflatcasstable}, on a noté $\gamma_{1}$ et $\gamma_{2}$ les deux éléments de $\Delta_{a}^{nr}-\Lambda$. Il est plus simple ici de les noter simplement $\alpha_{3}$ et $\alpha_{6}$. Pour chacune de ces deux racines $\alpha$, on introduira dans les deux paragraphes suivants un élément  $n_{\alpha}\in K_{s_{\alpha}}^{0,nr}\cap Norm_{G(F^{nr})}(M,s,N)$ dont l'image dans $Norm_{G_{s_{\alpha}}}(M_{s})$ s'envoie sur l'élément non trivial de $W^{G_{s_{\alpha}}}(M_{s})$ et qui vérifie les propriétés suivantes
 
 (2) $\epsilon_{N}^{\flat}(n_{\alpha})=-1$;
 
 (3) l'action de $Ad(n_{\alpha})$ sur $M_{AD}$ conserve l'épinglage affine $\mathfrak{E}^M_{a}$.

  Il résulte de (2) que $\epsilon_{N}^{\flat}(n_{\alpha_{3}}n_{\alpha_{6}})=1$. Il résulte de (3) que $\mu(n_{\alpha_{3}}n_{\alpha_{6}})$ est un élément de $\boldsymbol{\Omega}^{nr}$, donc que $\epsilon_{N}\circ\mu(n_{\alpha_{3}}n_{\alpha_{6}})=1$.   Cela  démontre  l'assertion (1) de \ref{epsilonflatcasstable}, donc le théorème.
   
    Rappelons le résultat suivant, que nous utiliserons plusieurs fois: \`a toute racine $\alpha\in \Sigma$ est associ\'e le sous-groupe de $G$ dont l'alg\`ebre de Lie est engendr\'ee par $\mathfrak{u}_{\alpha}$, $\mathfrak{u}_{-\alpha}$ et par l'\'el\'ement $\check{\alpha}\in X_{*}(T)\subset \mathfrak{t}$. Ce sous-groupe est en g\'en\'eral isomorphe \`a un quotient de $SL(2)$. Ici, puisque $G$ est simplement connexe, il est isomorphe \`a $SL(2)$. On l'appelle le groupe $SL(2)$ associ\'e \`a $\alpha$.  Soit $g$ un \'el\'ement de ce groupe et soit $\beta\in \Sigma$. Supposons que $\alpha$ et $\beta$ sont orthogonales. Parce que toutes les racines de $G$ sont de m\^eme longueur, $Ad(g)$ fixe tout \'el\'ement de $\mathfrak{u}_{\beta}$. 
   
   \subsubsection{Sommet $s_{\alpha_{6}}$ pour $G$ de type $E_{8}$\label{alpha6}}

   Posons $x=s_{\alpha_{6}}$. On a $G_{x}=(Spin(10)\times SL(4))/ \bar{\zeta}_{1/4}$, o\`u le plongement de $ \bar{\zeta}_{1/4}$ dans $Spin(10)\times SL(4)$ est
   $$\underline{i}\mapsto (\check{\alpha}_{1}(-1)\check{\alpha}_{4}(-1)\check{\alpha}_{2}(-\underline{i})\check{\alpha}_{5}(\underline{i}), \check{\alpha}_{7}(-\underline{i})\check{\alpha}_{8}(-1)\check{\alpha}_{0}(\underline{i})).$$
   Notons $G_{*}=Spin(10)$ et $G_{**}=SL(4)$ les deux composantes  de $G_{x,SC}$. La paire de Borel $(B_{x},T_{x})$ de $G_{x}$ se rel\`eve en le produit de deux telles paires de $G_{*}$ et $G_{**}$ que l'on note $(B_{*},T_{*})$ et $(B_{**},T_{**})$. On note $W_{*}$ et $W_{**}$ les groupes de Weyl de $G_{*}$ et $G_{**}$. 
   L'image r\'eciproque de $M_{s}$ dans  $G_{x,SC}=G_{*}\times G_{**}$ est $M_{*}\times G_{**}$, o\`u $M_{*,SC}=SL(2)\times Spin(6)$; l'\'el\'ement $N$ s'\'ecrit $N_{*}\oplus N_{**}$ o\`u $N_{*}=\sum_{i=1,2,4,5} \bar{E}_{\alpha_{i}}$ et $N_{**}=\sum_{i=0,7,8}\bar{E}_{\alpha_{i}}$.  
   Il est commode de r\'eindexer les racines simples de  la composante $G_{*}$ selon la num\'erotation habituelle pour un groupe $Spin(10)$. On note $(\gamma_{i})_{i=1,...,5}$ ces racines simples, d\'efinies par $\gamma_{1}=\alpha_{1}$, $\gamma_{2}=\alpha_{3}$, $\gamma_{3}=\alpha_{4}$, $\gamma_{4}=\alpha_{5}$, $\gamma_{5}=\alpha_{2}$. On peut identifier le groupe de Weyl $W_{*}$ de $G_{*}$ au groupe des permutations $w$ de $\{1,...,10\}$ telles que $w(11-i)=11-w(i)$ pour tout $i$ et que  l'ensemble des $i\in \{1,...,5\}$ tels que $wi\geq6$ a un nombre pair d'\'el\'ements. 
   La sym\'etrie \'el\'ementaire $w_{i}$ associ\'ee \`a la  racine $\gamma_{i}$ pour $i=1,...,4$ permute $i$ et $i+1$ ainsi que $11-i$ et $10-i$. La sym\'etrie \'el\'ementaire $w_{5}$ associ\'ee \`a $\gamma_{5}$ permute $4$ et $6$ ainsi que $5$ et $7$. On note $l:W_{*}\to {\mathbb N}$ la longueur relative \`a ces g\'en\'erateurs.  Introduisons l'\'el\'ement $w_{*}$ de $W_{*}$ qui  permute $1$ et $9$ ainsi que $2$ et $10$. Compl\'etons-le par l'identit\'e de la deuxi\`eme composante $G_{**}$. On obtient un \'el\'ement  du groupe de Weyl $W^{G_{x}}$ de $G_{x}$ qui normalise $M_{s}$ et dont    l'image dans $W^{G_{x}}(M_{s})$ est non triviale.  Introduisons le sous-groupe $H$ de $G$ engendr\'e par  les $\check{\alpha}_{i}$ et les  $E_{\pm \alpha_{i}}$ pour $i=1,2,3,4,5$. Il est isomorphe \`a $Spin(10)$ (cette fois, il s'agit d'un groupe sur $F^{nr}$) et  est muni d'un \'epinglage $\mathfrak{E}^H=(B^H,T^H,(E_{\alpha_{i}})_{i=1,2,3,4,5})$, o\`u $B^H=B\cap H$ et $T^H=T\cap H$. Introduisons la section de Springer $S$ relative  \`a cet \'epinglage. 
    En fait, $H$ est le  groupe d\'eriv\'e d'un Levi standard de $G$ et $S$ est simplement la restriction de la section de Springer de $W$ relative \`a $\mathfrak{E}$. 
   Le groupe $W_{*}$ est aussi le groupe de Weyl de $H$.    Les g\'en\'erateurs $S(w_{i})$ de la section de Springer appartiennent par construction \`a $K_{x}^{0,nr}$ donc l'image de $S$ est contenue dans ce groupe.  Posons 
   
    $n_{\alpha_{6}}=\check{\varpi}_{6}(-1)S(w_{*})^{-1}$.
   
   On va prouver que cet élément vérifie les conditions de \ref{preuveE8}, en particulier les relations (2) et (3) de ce paragraphe. 
    Pour simplifier la notation, on pose  $n=n_{\alpha_{6}}$. 
   On a  $n\in K_{x}^{0,nr}$.  Puisque $w_{*}$ conserve $M_{s}$, $Ad(n)$ conserve $M$. Montrons que
   
   (1) $Ad(n)$ fixe l'\'epinglage affine $\mathfrak{ E}_{a}^{M}$.
   
   Pour  $i=0,7,8$, il est clair que $Ad(S(w_{*}))$ fixe $E_{\alpha_{i}}$. Pour  $i=1,2,4,5$, il r\'esulte de la d\'efinition de $w_{*}$ que $w_{*}(\alpha_{i})=\alpha_{i}$. La propri\'et\'e \ref{epinglages}(5) entra\^{\i}ne que   $Ad(S(w_{*}))$ fixe $E_{\alpha_{i}}$. D'autre part $\check{\varpi}_{6}(-1)$ agit trivialement sur $E_{\alpha_{i}}$ pour $i\in \{1,...,8\}-\{6\}$ et il agit aussi trivialement sur $E_{\alpha_{0}}$ car $<\alpha_{0},\check{\varpi}_{6}>=-4$. Donc $Ad(n)$ fixe $E_{\alpha}$ pour tout $\alpha \in  \Lambda$. Il  reste \`a prouver qu'il fixe aussi le dernier \'el\'ement $E_{\beta_{4}}$ de l'\'epinglage affine $\mathfrak{ E}_{a}^{M}$. Introduisons l'\'el\'ement $\delta\in W_{*}$ qui est la permutation 
   $(1,2,3,4,5,6,7,8,9,10)\mapsto (4,5,1,2,3,8,9,10,6,7)$. On constate que $w_{*}=\delta^{-1}w_{5}\delta$ et que $l(w_{*})=l(\delta^{-1})+l(w_{5})+l(\delta)$. Il en r\'esulte que $S(w_{*})=S(\delta^{-1})S(w_{5})S(\delta)$. A l'aide de la relation \ref{epinglages}(4), on calcule $S(\delta^{-1})=tS(\delta)^{-1}$, o\`u $t=\prod_{\gamma>0; \omega(\gamma)<0}\check{\gamma}(-1)$. L'ensemble des $\gamma>0$ tels que $\omega(\gamma)<0$ est $\{\gamma_{1}+\gamma_{2},\gamma_{1}+\gamma_{2}+\gamma_{3},\gamma_{1}+\gamma_{2}+\gamma_{3}+\gamma_{4}, \gamma_{2},\gamma_{2}+\gamma_{3},\gamma_{2}+\gamma_{3}+\gamma_{4}\}$. D'o\`u $t=\check{\gamma}_{1}(-1)=\check{\alpha}_{1}(-1)$. 
   Par construction, $S(w_{5})$ est contenu dans le groupe $SL(2)$ (de $H$ ou $G$) associ\'e \`a la racine $\gamma_{5}$. Il en r\'esulte que $S(\delta)^{-1}S(w_{5})S(\delta)$ est contenu dans le groupe $SL(2)$ associ\'e \`a la racine $\delta^{-1}(\gamma_{5})$. Cette racine est la plus grande racine positive de $H$, c'est-\`a-dire $\gamma_{1}+2\gamma_{2}+2\gamma_{3}+\gamma_{4}+\gamma_{5}$, autrement dit $\alpha_{1}+2\alpha_{3}+2\alpha_{4}+\alpha_{5}+\alpha_{2}$. On constate que cette racine est orthogonale \`a $\beta_{4}=\alpha_{3}+\alpha_{4}+\alpha_{5}+\alpha_{6}$. Il en r\'esulte que le groupe $SL(2)$ associ\'e \`a $\delta^{-1}(\gamma_{5})$ fixe $E_{\beta_{4}}$. Par contre,  $Ad(t)$ multiplie $E_{\beta_{4}}$ par $-1$, d'o\`u l'\'egalit\'e $Ad(S(w_{*}))(E_{\beta_{4}})=-E_{\beta_{4}}$. Mais on a aussi $Ad(\check{\varpi}_{6}(-1))(E_{\beta_{4}})=-E_{\beta_{4}}$, donc $Ad(n)$ fixe $E_{\beta_{4}}$. Cela prouve (1).
   
   D'apr\`es (1),  l'action de $Ad(n)$ sur $M_{AD}$ est triviale a fortiori sa réduction dans $M_{s}$ fixe $N$. Donc $n\in Norm_{G(F^{nr})}(M,s,N)$. 
        
   Pour calculer $\epsilon_{x,N}^{\flat}\circ c_{x}(n)$, remarquons que $\check{\varpi}_{6}(-1)$ appartient \`a $K_{s,N}^{0,nr}$. On a donc $\epsilon_{x,N}^{\flat}\circ c_{x}(n)=\xi_{N}(\overline{{\check{\varpi}}_{6}(-1)})\epsilon_{x,N}^{\flat}\circ c_{x}(S(w_{*})^{-1})$. On a l'\'egalit\'e $\check{\varpi}_{6}=-\check{\alpha}_{7}-2\check{\alpha}_{8}-3\check{\alpha}_{0}$ et on en d\'eduit l'\'egalit\'e $\overline{\check{\varpi}_{6}(-1)}=z_{2}^2$. D'o\`u
   
   (2) $\epsilon^{\flat}_{x,N}\circ c_{x}(n)=-\epsilon^{\flat}_{x,N}\circ c_{x}(S(w_{*})^{-1})$. 
   
   On va calculer  le dernier terme selon la m\'ethode expliqu\'ee en \ref{calculepsilonflat}.  Le calcul se concentre dans l'image dans $G_{x}$ de la composante $G_{*}$ de $G_{x,SC}$ et on peut aussi bien remplacer cette image par $G_{*}$ lui-m\^eme. On note $S_{*}$ la section de Springer de $G_{*}$ relative \`a l'\'epinglage fix\'e. L'\'el\'ement $c_{x}(S(w*))$ est l'image dans $G_{x}$ de $S_{*}(w_{*})$. 
   Le groupe $M_{*}$ est un Levi standard de $G_{*}$. Notons $P_{*}$ le sous-groupe parabolique standard de $G_{*}$ de composante de Levi $M_{*}$. On choisit l'\'el\'ement $\underline{N}_{*}=
    =\sum_{i=1,...,5}\bar{E}_{\alpha_{i}}$. Il appartient \`a $\mathfrak{u}_{B_{*}}$ et \`a l'orbite induite de la $M_{*}$-orbite de $N_{*}$.  Fixons un \'el\'ement $Z\in \mathfrak{z}(M_{*})_{reg}$ que l'on pr\'ecisera ci-dessous. Soit $\nu\in \bar{{\mathbb F}}_{q}^{\times}$. On doit calculer l'\'el\'ement $u_{\nu}\in U_{P_{*}}$ tel que
   
     $$(3) \qquad u_{\nu}^{-1}(\nu Z+\underline{N}_{*})u_{\nu}=\nu Z+N_{*}.$$
   Notons $B_{*}^{opp}$ le sous-groupe de Borel de $G_{*}$ contenant $T_{*}$ et oppos\'e \`a $B_{*}$. Quitte \`a exclure un nombre fini de valeurs de $\nu$, on a une \'egalit\'e 
   
   $$(4) \qquad u_{\nu}S_{*}(w_{*}) =\underline{v}_{\nu}t_{\nu}v_{\nu}$$
    o\`u 
    $\underline{v}_{\nu}\in U_{B_{*}^{opp}}$, $t_{\nu}\in T_{*}$,   $v_{\nu}\in U_{B_{*}}$.  On note ${\bf t}_{\nu}$ l'image de $t_{\nu}$ dans $T_{*}/Z(M_{*})^0$. Alors les applications alg\'ebriques $\nu\mapsto \underline{v}_{\nu}$ et $\nu\mapsto {\bf t}_{\nu}$ se prolongent en $\nu=0$. On a $\underline{v}_{0}=1$ et ${\bf t}_{0}={\bf z}$, o\`u ${\bf z}$ est l'image dans $T_{*}/Z(M_{*})^0$ d'un \'el\'ement $z\in Z(M_{*})$. D'apr\`es \ref{calculepsilonflat}(4) , on a alors 
    
    (5) $\epsilon_{x,N}^{\flat}(S_{*}(w_{*})^{-1})=\xi_{N}(z)^{-1}$. 
    
    Si $G_{*}$ \'etait un groupe $SO(10)$, tous les termes ci-dessus se calculeraient par de fastidieux mais simples calculs matriciels. Mais $G_{*}$ est le groupe $Spin(10)$. Sans entrer dans le d\'etail des calculs, indiquons-en la m\'ethode. 
 Introduisons un espace $V$ de dimension $10$ sur $\bar{{\mathbb F}}_{q}$, muni d'une base $(e_{i})_{i=1,...,10}$ et de la forme quadratique $q$ d\'efinie par $q(e_{i},e_{11-i})=(-1)^{i}$ pour $i\leq 5$ et $q(e_{i},e_{j})=0$ si $i+j\not=11$. Notons $\tilde{G}_{*}$ son groupe sp\'ecial orthogonal. Le groupe $G_{*}$ est le rev\^etement spinoriel de $\tilde{G}_{*}$.  
     L'espace $\mathfrak{g}_{*}=\tilde{\mathfrak{g}}_{*}$ s'identifie \`a celui des \'el\'ements $X\in End(V)$ tels que $(Xv,v')+(v,Xv')=0$ pour tous $v,v'\in V$.  Avec les notations de \ref{Bn}, on peut identifier $\bar{E}_{\gamma_{1}}=E_{1,2}+E_{9,10}$, $\bar{E}_{\gamma_{2}}=E_{2,3}+E_{8,9}$, $\bar{E}_{\gamma_{3}}=E_{3,4}+E_{7,8}$, $\bar{E}_{\gamma_{4}}=E_{4,5}+E_{6,7}$, $\bar{E}_{\gamma_{5}}=E_{4,6}+E_{5,7}$. On  choisit pour $Z$ l'\'el\'ement  $Z=E_{1,1}+E_{2,2}-E_{9,9}-E_{10,10}$. On note $\pi=G_{*}\to \tilde{G}_{*}$ la projection.   Posons $\tilde{S}_{*}=\pi\circ S_{*}$. C'est la section de Springer de $\tilde{G}_{*}$. 
    Dans \cite{W1} lemme X.4, on a calcul\'e explicitement la section de Springer  $\tilde{S}_{*}$. Pour $w\in W_{*}$ et $i\in \{1,...,10\}$, posons $r(i,w)=\vert \{k\in {\mathbb N}; i<k\leq10; w(i)> w(k), i+k\not=11\}$. Alors $\tilde{S}(w)(e_{i})=(-1)^{r(i,w)}e_{wi}$. Remarquons que la projection $\pi$ se restreint en une bijection entre les sous-ensembles des \'el\'ements unipotents des deux groupes $G_{*}$ et $\tilde{G}_{*}$. En particulier l'\'el\'ement $u_{\nu}$ v\'erifiant (3)  est d\'etermin\'e par un calcul matriciel dans $\tilde{G}_{*}$. Introduisons le sous-groupe de Levi standard $L_{*}$ de $G_{*}$ tel que $\Delta^{L_{*}}=\{\gamma_{1},\gamma_{2},\gamma_{3},\gamma_{4}\}$. Son groupe d\'eriv\'e $L_{*,der}$ est isomorphe \`a $SL(5)$. Notons $Q_{*}$ le sous-groupe parabolique standard de composante de Levi $L_{*}$. Pour simplifier la notation, fixons $\nu$ et abandonnons les indices $\nu$: on pose $u=u_{\nu}$. On \'ecrit $u=v_{1}v_{2}$ avec $v_{1}\in U_{B_{*}\cap L_{*}}$ et $v_{2}\in U_{Q_{*}}$. Ces \'el\'ements se d\'eterminent par un calcul matriciel dans $\tilde{G}_{*}$. On a $uS_{*}(\delta^{-1})=v_{1}S_{*}(\delta^{-1})v_{3}$, o\`u $v_{3}=S_{*}(\delta^{-1})^{-1}v_{2}S_{*}(\delta^{-1})$. Le terme $v_{3}$ se calcule aussi dans $\tilde{G}_{*}$: on peut remplacer $S_{*}(\delta^{-1})$ par $\tilde{S}_{*}(\delta^{-1})$ dans la d\'efinition de $v_{3}$. Les deux termes $v_{1}$ et $S_{*}(\delta^{-1})$ appartiennent \`a $L_{*,der}$. La restriction \`a $W_{*}^{L_{*}}$ de $S_{*}$ est la section de Springer du groupe $L_{*,der}\simeq SL(5)$. Un calcul matriciel dans $SL(5)$ calcule $v_{1}S_{*}(\delta^{-1})$ sous la forme $v_{1}S_{*}(\delta^{-1})=\underline{v}_{1}t_{1}v_{4}$, o\`u $\underline{v}_{1}\in U_{B_{*}^{opp}\cap L_{*}}$, $t_{1}\in T_{*}$, $v_{4}\in U_{B_{*}\cap L_{*}}$. On obtient $uS_{*}(\delta^{-1})=\underline{v}_{1}t_{1}v_{5}$, o\`u $v_{5}=v_{4}v_{3}\in U_{B_{*}}$. Notons $L'_{*}$ le groupe de Levi standard tel que $\Delta^{L'_{*}}=\{\gamma_{5}\}$. Son groupe d\'eriv\'e $L'_{*,der}$ est isomorphe \`a $SL(2)$. Notons $Q'_{*}$ le sous-groupe parabolique standard de composante de Levi $L'_{*}$. Par un calcul matriciel dans $\tilde{G}_{*}$, on \'ecrit $v_{5}=v_{6}v_{7}$ o\`u $v_{6}\in U_{B_{*}\cap L'_{*}}$ et $v_{7}\in U_{Q'_{*}}$. On a $v_{5}S_{*}(w_{5})=v_{6}S_{*}(w_{5})v_{8}$, o\`u $v_{8}=S_{*}(w_{5})^{-1}v_{7}S_{*}(w_{5})$. Comme ci-dessus, $v_{8}$ se calcule dans $\tilde{G}_{*}$. Les termes $v_{6}$ et $S_{*}(w_{5})$ appartiennent \`a $L'_{*,der}$. Dans ce groupe isomorphe \`a $SL(2)$, on calcule $v_{6}S_{*}(w_{5})=\underline{v}_{2}t_{2}v_{9}$, avec $\underline{v}_{2}\in U_{B_{*}^{opp}\cap L'_{*}}$, $t_{2}\in T_{*}$ et $v_{9}\in U_{B_{*}\cap L'_{*}}$. En posant $v_{10}=v_{9}v_{8}\in U_{B_{*}}$, on obtient $uS(\delta^{-1})S(w_{5})=\underline{v}_{1}t_{1}\underline{v}_{2}t_{2}v_{10}$. Ecrivons $v_{10}=v_{11}v_{12}$, avec $v_{11}\in U_{B_{*}\cap L_{*}}$ et $v_{12}\in U_{Q_{*}}$. On a $v_{10}S_{*}(\delta)=v_{11}S_{*}(\delta)v_{13}$, o\`u $v_{13}=S(\delta)^{-1}v_{12}S(\delta)$. Tous ces \'el\'ements se calculent dans $\tilde{G}_{*}$. Comme dans la premi\`ere \'etape du calcul, on calcule dans le groupe $L_{*,der}$ une \'egalit\'e $v_{11}S_{*}(\delta)=\underline{v}_{3}t_{3}v_{14}$, o\`u $\underline{v}_{3}\in U_{B_{*}^{opp}\cap L_{*}}$, $t_{3}\in T_{*}$ et $v_{14}\in U_{B_{*}\cap L_{*}}$. En posant $v_{15}=v_{14}v_{13}\in U_{B_{*}}$, on obtient $uS_{*}(\delta^{-1})S_{*}(w_{5})S_{*}(\delta)=\underline{v}_{1}t_{1}\underline{v}_{2}t_{2}\underline{v}_{3}t_{3}v_{15}$. Ou encore, en posant $\underline{v}_{4}=\underline{v}_{1}t_{1}\underline{v}_{2}t_{2}\underline{v}_{3}t_{2}^{-1}t_{1}^{-1}$ et $t=t_{1}t_{2}t_{3}$, 
    $$u_{\nu}S_{*}(w_{*})=\underline{v}_{4}tv_{15}.$$
    En comparant avec (4), on a $t_{\nu}=t$. En explicitant les calculs ci-dessus, on trouve
    $$t_{\nu}=\left(\check{\gamma}_{1}(2\nu^7)\check{\gamma}_{2}(4\nu^{14})\check{\gamma}_{3}(4\nu^{14})\check{\gamma}_{4}(2\nu^7)\check{\gamma}_{5}(2\nu^7)\right)^{-1}.$$
    Mais le groupe $X_{*}(Z(M_{*})^0$ est engendr\'e par l'\'el\'ement $\check{\gamma}_{1}+2\check{\gamma}_{2}+2\check{\gamma}_{3}+\check{\gamma}_{4}+\check{\gamma}_{5}$. Donc $t_{\nu}\in Z(M_{*})^0$ et ${\bf t}_{\nu}=1$. D'apr\`es (5), on a donc $\epsilon_{x,N}^{\flat}(S_{*}(w_{*})^{-1})=1$, puis, d'apr\`es (2), $\epsilon^{\flat}_{x,N}\circ c_{x}(n)=-1$.  
     
  On a bien démontré les assertions (2) et (3) de \ref{preuveE8}.

 \subsubsection{Sommet $s_{\alpha_{3}}$ pour $G$ de type $E_{8}$}\label{alpha3}
 
  Posons  $x=s_{\alpha_{3}}$. On a $G_{x}=(SL(8)\times SL(2))/\bar{\zeta}_{1/4}$, o\`u le plongement de $\bar{\zeta}_{1/4} $ dans $SL(8)\times SL(2)$ est
 $$\underline{i}\mapsto (\check{\alpha}_{2}(-\underline{i})\check{\alpha}_{4}(-1)\check{\alpha}_{5}(\underline{i})\check{\alpha}_{7}(-\underline{i})\check{\alpha}_{8}(-1)\check{\alpha}_{0}(\underline{i}),\check{\alpha}_{1}(-1)).$$
 On note $G_{*}=SL(8)$ et $G_{**}=SL(2)$ les deux composantes de $G_{x,SC}$.   L'image r\'eciproque de $M_{s}$ dans $G_{x,SC}$ est $M_{*}\times G_{**}$, o\`u $M_{*,SC}=SL(4)\times SL(4)$; l'\'el\'ement $N$ s'\'ecrit $N_{*}\oplus N_{**}$ o\`u $N_{*}=\sum_{i=2,4,5,7,8,0} \bar{E}_{\alpha_{i}}$ et $N_{**}= \bar{E}_{\alpha_{1}}$.
 
 Notons $H$ le sous-groupe de $G$ engendr\'e par  les $\check{\alpha}_{i}$ et les $E_{\pm\alpha_{i}}$ pour $i=2,4,5,6,7,8,0$. Introduisons un espace $V$ de dimension $8$ sur $ F$ muni d'une base $(e_{i})_{i=1,...,8}$. On identifie $H$ \`a $SL(V)$ de sorte que $T\cap H$ s'identifie au groupe diagonal de $SL(V)$ et   les \'el\'ements $E_{\alpha_{i}}$ pour $i=2,4,5,6,7,8,0$ s'identifient respectivement aux \'el\'ements de l'\'epinglage habituel $E_{1,2}$, $E_{2,3}$ etc...
   Notons $n_{\alpha_{3}}$ l'\'el\'ement de $H$ tel que $n_{\alpha_{3}}(e_{i})=-e_{i+4}$ pour $i=1,...,4$ et $n_{\alpha_{3}}(e_{i})=e_{i-4}$ pour $i=5,...,8$.  On va prouver que cet élément vérifie les conditions de \ref{preuveE8}, en particulier les relations (2) et (3) de \ref{preuveE8}. 
    Pour simplifier la notation, on pose simplement $n=n_{\alpha_{3}}$.

   L'action de $n$ sur l'appartement $App_{F^{nr}}(T)$ conserve le sous-groupe de racines affines engendr\'e par les $\alpha_{i}^{aff}$ pour $i=0,...,8$, $i\not=3$. Puisque le sommet $x$ est l'\'el\'ement de cet appartement d\'efini par les \'egalites $\alpha_{i}^{aff}(x)=0$ pour $i\not=3$, l'action de $n$  fixe $x$. Donc $n\in K_{x}^{0,nr}$. L'action $Ad(\bar{n})$ conserve $M_{s}$ donc $Ad(n)$ conserve $M$ par un raisonnement d\'ej\`a fait plusieurs fois. Montrons que
   
   (1) l'action de $Ad(n)$ sur $M_{AD}$ conserve l'\'epinglage affine $\mathfrak{ E}_{a}^{M}$ et y agit par $\boldsymbol{\omega}$.
   
   Par construction, l'automorphisme $Ad(n)$  permute les couples $E_{\alpha_{2}}$ et $E_{\alpha_{7}}$, $E_{\alpha_{4}}$ et $E_{\alpha_{8}}$, $E_{\alpha_{5}}$ et $E_{\alpha_{0}}$. Puisque $\alpha_{1}$ est orthogonale aux racines $\alpha_{i}$ pour $i=2,4,5,6,7,8,0$, tout \'el\'ement de $H$ fixe $E_{\alpha_{1}}$. Donc $Ad(n)$ fixe $E_{\alpha_{1}}$. 
  Par la traduction entre les racines  dans $G$ et dans $M$ que l'on a expliqu\'ee en \ref{preuveE8}, il reste \`a d\'emontrer que $Ad(n)$  fixe le dernier \'el\'ement $E_{\beta_{4}}$ de $\mathfrak{ E}_{a}^{M}$. Par construction, on a $n=n_{1}n_{2}n_{3}n_{4}$ o\`u, pour $i=1,...,4$, $n_{i}(e_{i})=-e_{i+4}$, $n_{i}(e_{i+4})=e_{i}$ et $n_{i}(e_{j})=e_{j}$ pour $j=1,...,8$, $j\not=i,i+4$. L'\'el\'ement $n_{i}$ appartient au groupe $SL(2)$ associ\'e \`a la racine $\gamma_{i}$, o\`u $\gamma_{1}=\alpha_{2}+\alpha_{4}+\alpha_{5}+\alpha_{6}$, $\gamma_{2}=\alpha_{4}+\alpha_{5}+\alpha_{6}+\alpha_{7}$, $\gamma_{3}=\alpha_{5}+\alpha_{6}+\alpha_{7}+\alpha_{8}$, $\gamma_{4}=\alpha_{6}+\alpha_{7}+\alpha_{8}+\alpha_{0}$. On constate que les $\gamma_{i}$ sont orthogonales \`a $\beta_{4}=\alpha_{3}+\alpha_{4}+\alpha_{5}+\alpha_{6}$. Il en r\'esulte que ces groupes $SL(2)$ fixent $E_{\beta_{4}}$. Donc $Ad(n)$ fixe aussi $E_{\beta_{4}}$. Cela d\'emontre (1).

  On calcule $\epsilon_{x,N}^{\flat}\circ c_{x}(n)$ par la m\^eme m\'ethode qu'en \ref{alpha6} et on utilise les m\^emes notations. De nouveau, le calcul se concentre dans la composante $G_{*}$. 
    On pose   $\underline{N}_{*}=\sum_{i=2,4,5,6,7,8,0}\bar{E}_{\alpha_{i}}$.   
 Notons $L$ le $\mathfrak{o}_{F^{nr}}$-r\'eseau de $V$ engendr\'e par la base $(e_{i})_{i=1,...,8}$ et posons $\bar{V}=L/\mathfrak{p}_{F^{nr}}L$. L'espace $\bar{V}$ est muni de la base r\'eduite $(\bar{e}_{i})_{i=1,...,8}$ et  le groupe $G_{*}$ s'identifie \`a $SL(\bar{V})$. Avec des notations \'evidentes, on  a $N_{*}=\sum_{i=1,...,7, i\not=4}\bar{E}_{i,i+1}$ et $\underline{N}_{*}=\sum_{i=1,...,7}\bar{E}_{i,i+1}$. On choisit $Z=\sum_{i=1,...,4}E_{i,i}-\sum_{i=5,...,8}E_{i,i}$. Pour $\nu\in \bar{{\mathbb F}}_{q}^{\times}$, des calculs matriciels faciles dans $G_{*}=SL(8)$ permettent de calculer l'\'el\'ement $u_{\nu}\in U_{P_{*}}$ tel que $u_{\nu}^{-1}(\nu Z+\underline{N}_{*})u_{\nu}=\nu Z+N_{*}$, puis d'\'ecrire $u_{\nu}\bar{n}^{-1}=\underline{v}_{\nu}t_{\nu}v_{\nu}$ avec $\underline{v}_{\nu}\in U_{B_{*}^{opp}}$, $t_{\nu}\in T_{*}$ et $v_{\nu}\in U_{B_{*}}$.  On obtient
  $$t_{\nu}=\left(\check{\alpha}_{2}(2\nu)^4\check{\alpha}_{4}(2\nu)^8\check{\alpha}_{5}(2\nu)^{12}\check{\alpha}_{6}(2\nu)^{16}\check{\alpha}_{7}(-(2\nu)^{12})\check{\alpha}_{8}(2\nu)^8\check{\alpha}_{0}(-(2\nu)^4)\right)^{-1}.$$
  L'\'el\'ement 
  $$\check{\alpha}_{2}+2\check{\alpha}_{4}+3\check{\alpha}_{5}+4\check{\alpha}_{6}+3\check{\alpha}_{7}+2\check{\alpha}_{8}+\check{\alpha}_{0} $$
  appartient \`a $X_{*}(Z(M_{*})^0)$. Donc $t_{\nu}\in zZ(M_{*})^0$, o\`u $z= \check{\alpha}_{7}(-1)\check{\alpha}_{0}(-1)$. On a $z\in Z(M_{*})$ et son image dans $Z(M_{s})$ est \'egale \`a   $z_{1}^2$, cf. \ref{preuveE8}. D'apr\`es \ref{calculepsilonflat}(4), on obtient  $\epsilon^{\flat}_{x,N}\circ c_{x}(n)=\xi_{N}(z_{1}^{-2}) $, d'où $\epsilon^{\flat}_{x,N}\circ c_{x}(n)=-1$.  
     
  On a bien démontré les assertions (2) et (3) de \ref{preuveE8}.

  \subsubsection{D\'ebut de la preuve du théorème \ref{epsilonflatcasstable} dans le cas $(E_{6},ram)$\label{preuveE6ram}}
  
  On suppose que $G$ est de type $E_{6}$, que  $M$ est de type $A_{5}$ et que $G$ n'est pas d\'eploy\'e sur $F^{nr}$ mais l'est sur l'extension quadratique ramifi\'ee $E$ de $F^{nr}$. On note $(\alpha_{i})_{i=1,...,6}$ les \'el\'ements de $\Delta$. Un \'el\'ement de $\Gamma_{F^{nr}}-\Gamma_{E}$ agit sur $\Delta$ par l'automorphisme non trivial $\theta$ de cet ensemble. On a $\Delta_{a}^{nr}=\{\beta_{0},\beta_{1,6},\beta_{3,5},\beta_{4},\beta_{2}\}$. La racine $\beta_{0}$ est \'egale \`a  la restriction commune de $\alpha'_{0}$ et $\alpha''_{0}$ o\`u $\alpha'_{0}=-(\alpha_{1}+\alpha_{2}+2\alpha_{3}+2\alpha_{4}+\alpha_{5}+\alpha_{6})$ et $\alpha''_{0}=-(\alpha_{1}+\alpha_{2}+\alpha_{3}+2\alpha_{4}+2\alpha_{5}+\alpha_{6})$.   L'ensemble $\Lambda$ est \'egal \`a $\{\beta_{0},\beta_{3,5},\beta_{2}\}$. On doit décrire l'ensemble $\Delta_{a}^{M,nr}$ introduit en \ref{epsilonflatcasstable}. Nous notons par des termes soulignés les racines de $M$. Posons $(\underline{\alpha}_{1},...,\underline{\alpha}_{5})=(\alpha_{3},\alpha_{4}+\alpha_{5}+\alpha_{6},\alpha_{2},\alpha_{1}+\alpha_{3}+\alpha_{4},\alpha_{5})$. Les restrictions à $T^{nr}$ de ces racines sont combinaisons linéaires à coefficients rationnels d'éléments de $\Lambda$, donc ces racines appartiennent à l'ensemble $\Sigma^M$ des racines de $T$ dans $\mathfrak{m}$. On voit que les racines $\underline{\alpha}_{i}$ pour $i=1,...,5$ vérifient  les m\^emes relations de produits scalaires qu'une base du système $\Sigma^M$ de type $A_{5}$. Il en résulte que cet ensemble de racines est une base de $\Sigma^M$. Il s'en déduit une base $\Delta^{M,res}$ de l'ensemble des restrictions $\Sigma^{M,res}$ que l'on note $(\underline{\beta}_{15},\underline{\beta}_{24},\underline{\beta}_{3})$. Pour obtenir un ensemble $\Delta_{a}^{M,nr}$, on ajoute la racine supplémentaire  $\underline{\beta}_{0}$ égale à la restriction commune de $\underline{\alpha}'_{0}$ et $\underline{\alpha}''_{0}$ o\`u $\underline{\alpha}'_{0}=-(\underline{\alpha}_{1}+\underline{\alpha}_{2}+\underline{\alpha}_{3}+\underline{\alpha}_{4})$ et $\underline{\alpha}''_{0}=-(\underline{\alpha}_{2}+\underline{\alpha}_{3}+\underline{\alpha}_{4}+\underline{\alpha}_{5})$. On constate que $\underline{\beta}_{0}=\beta_{0}$. Les racines affines $\underline{\beta}_{15}$, $\underline{\beta}_{24}$, $\underline{\beta}_{3}$ et $\underline{\beta}_{0}+\frac{1}{2}$  sont positives sur $C^{nr}$ donc le sont aussi sur $C^{M,nr}$. Par définition de l'ensemble $\Delta_{a}^{M,nr}$ de  \ref{epsilonflatcasstable}, on peut donc supposer que $\Delta_{a}^{M,nr}$ est celui que l'on vient de décrire.    La racine $\underline{\beta}_{s}$ associée au sommet $s$ est  est $\underline{\beta}_{24}$. 
 Le groupe $\boldsymbol{\Omega}^{M,nr}$ a deux \'el\'ements. L'\'el\'ement non trivial permute $\underline{\beta}_{0}$ et $\underline{\beta}_{1,5}$ et fixe $\underline{\beta}_{2,4}$ et $\underline{\beta}_{3}$. Les entiers $h$ et $k$ de \ref{constructionfamilleadmissible} valent $2$ et $1$. Le caract\`ere $\boldsymbol{\epsilon}$ de $\boldsymbol{ \Omega}^{M,nr}$ d\'efini dans ce paragraphe est donc trivial.

    On choisit $N=\sum_{\beta\in \Lambda}\bar{E}_{\beta}$. Notons $\{\pm 1\}^3_{1}$  le sous-groupe des $(z_{1},z_{2},z_{3})\in \{\pm 1\}^3$ tels que $z_{1}z_{2}z_{3}=1$. On a $M_{AD,s}=(SL(2)\times SL(2)\times SL(2))/\{\pm 1\}^3_{1}$, o\`u $\{\pm 1\}^3$ est identifi\'e au produit des centres des trois composantes $SL(2)$.    Le groupe $Z_{M_{AD,s}}(N)/Z_{M_{AD,s}}(N)^0$ a deux \'el\'ements et $\xi_{N}$ est son caract\`ere non trivial.

   Comme en \ref{preuveE8}, nous noterons simplement $\beta_{1,6}$ et $\beta_{4}$ les racines notées   $\gamma_{1}$ et $\gamma_{2}$  en  \ref{epsilonflatcasstable}.  Pour chacune de ces deux racines $\beta$, on introduira dans les deux paragraphes suivants un élément  $n_{\beta}\in K_{s_{\beta}}^{0,nr}\cap Norm_{G(F^{nr})}(M,s,N)$ dont l'image dans $Norm_{G_{s_{\beta}}}(M_{s})$ s'envoie sur l'élément non trivial de $W^{G_{s_{\beta}}}(M_{s})$ et qui vérifie les propriétés suivantes
 
 (1) $\epsilon_{N}^{\flat}(n_{\beta})=1$;
 
 (2) l'action de $Ad(n_{\beta})$ sur $M_{AD}$ conserve l'épinglage affine $\mathfrak{E}^M_{a}$.

  Il résulte de (1) que $\epsilon_{N}^{\flat}(n_{\beta_{16}}n_{\beta_{4}})=1$. Il résulte de (2) que $\mu(n_{\beta_{16}}n_{\beta_{4}})$ est un élément de $\boldsymbol{\Omega}^{M,nr}$, donc que $\epsilon_{N}\circ\mu(n_{\beta_{16}}n_{\beta_{4}})=1$.   Cela  démontre  l'assertion (1) de \ref{epsilonflatcasstable}, donc le théorème.

   \subsubsection{Sommet $s_{\beta_{1,6}}$ pour $G$ de type $(E_{6},ram)$\label{beta16} }
  
  Posons $x=s_{\beta_{1,6}}$. On a $G_{x}=(Spin(7)\times SL(2))/ \{\pm 1\}$, le groupe $\{\pm 1\}$ s'identifiant aux centres des deux composantes. Notons $G_{*}==Spin(7)$ et $G_{**}=SL(2)$ les deux composantes de $G_{s,SC}$. De la paire de Borel $(B_{x},T_{x})$ de $G_{x}$   se d\'eduisent de telles paires $(B_{*},T_{*})$ et $(B_{**},T_{**}
)$ de $G_{*}$ et $G_{**}$.  On note $W_{*}$ et $W_{**}$ les groupes de Weyl de $G_{*}$ et $G_{**}$. 
 L'image r\'eciproque de $M_{s}$ dans ce groupe est $M_{*}\times G_{**}$ o\`u $M_{*,SC}=SL(2)\times SL(2)$. L'\'el\'ement $N$ s'\'ecrit $N_{*}\oplus N_{**}$, o\`u $N_{*}=\bar{E}_{\beta_{3,5}}+\bar{E}_{\beta_{2}}$ et $N_{**}=\bar{E}_{\beta_{0}}$.  
  On r\'eindexe les racines de $G_{*}$ selon la num\'erotation habituelle pour un groupe $Spin(7)$ en posant $\gamma_{1}=\beta_{2}$, $\gamma_{2}=\beta_{4}$ et $\gamma_{3}=\beta_{3,5}$. Le groupe de Weyl $W_{*}$ de $G_{*}$ s'identifie au groupe de permutations $w$ de l'ensemble $ \{1,...,7\}$ telles que $w(8-i)=8-w(i)$ pour tout $i$. La sym\'etrie \'el\'ementaire $w_{*,i}$ associ\'ee \`a $\gamma_{i}$ pour $i=1,2 $ permute $i$ et $i+1$ ainsi que $8-i$ et $7-i$, tandis que la sym\'etrie \'el\'ementaire $w_{*,3}$ associ\'ee \`a $\gamma_{3}$ permute $3$ et $5$. Pour $i=1,...,6$, notons $w_{i}$ la sym\'etrie \'el\'ementaire de $W$ associ\'ee \`a $\alpha_{i}$. Introduisons la section de Springer de $G$ associ\'ee \`a notre \'epinglage $\mathfrak{E}$. On voit que les \'el\'ements $S(w_{2})$, $S(w_{4})$ et $S(w_{3})S(w_{5})$ appartiennent \`a $K_{x}^{0,nr}$. Ils se r\'eduisent en des \'el\'ements de $G_{x}$ qui normalisent $T_{x}$ et d\'efinissent donc des \'el\'ements de $W_{*}\times W_{**}$. Les images dans $W_{**}$ des trois \'el\'ements ci-dessus sont $1$ tandis que leurs images dans $W_{*}$ sont respectivement $w_{*,1}$, $w_{*,2}$, $w_{*,3}$.

  Introduisons l'\'el\'ement $w_{*}\in W_{*}$ qui permute $1$ et $6$, $2$ et $7$ et fixe $3,4,5$. On le compl\`ete par l'identit\'e de la seconde composante $G_{**}$. L'\'el\'ement obtenu conserve $M_{s}$ et son image dans   $W^{G_{x}}(M_{s})$ est l'\'el\'ement non trivial de ce groupe.  On a l'\'egalit\'e 
  $$(1) \qquad w_{*}=w_{*,2}w_{*,3}w_{*,2}w_{*,1}w_{*,2}w_{*,3}w_{*,2}.$$
  Introduisons l'\'el\'ement $w\in W$ d\'efini par
  $$w=w_{4}w_{3}w_{5}w_{4}w_{2}w_{4}w_{5}w_{3}w_{4}.$$
  On v\'erifie que c'est une d\'ecomposition de longueur minimale, c'est-\`a-dire qu'en notant $l$ la longueur de $W$ relative aux g\'en\'erateurs $(w_{i})_{i=1,...,6}$, on a $l(w)=9$. Posons
  $n_{\beta_{1,6}}=S(w)^{-1}$. Nous allons prouver  que cet élément  vérifie les conditions de \ref{preuveE6ram}, en particulier que les relations (1) et (2) de ce paragraphe sont vérifiées. Pour simplifier, on pose $n=n_{\beta_{1,6}}$.

   D'apr\`es ce qui pr\'ec\`ede,  $n$ appartient \`a $K_{x}^{0,nr}$ et se r\'eduit en un \'el\'ement  $\bar{n}$ de $Norm_{G_{x}}(T_{x})$ dont l'image dans $W_{*}\times W_{**}$ est  $(w_{*},1)$.  L'action $Ad(\bar{n})$ conserve $M_{s}$ donc    $Ad(n)$ conserve $M$. Montrons que
 
 (2) $Ad(n)$ fixe l'\'epinglage affine $\mathfrak{ E}_{a}^{M,nr}$.
  
  On v\'erifie que $w(\alpha_{i})=\alpha_{i}$ pour $i=2,3,5$. D'apr\`es \ref{epinglages} (5),  $S(w)$ fixe $E_{\alpha_{i}}$ pour ces $i$, donc fixe $E_{\underline{\beta}_{3}}=E_{\beta_{2}}=E_{\alpha_{2}}$ et $E_{\underline{\beta}_{1,5}}=E_{\beta_{3,5}}=E_{\alpha_{3}}+E_{\alpha_{5}}$. Posons $\omega=w_{4}w_{5}w_{3}w_{4}$. On a $\omega^2=1$ et $w=\omega w_{2}\omega$.  On a $l(w)=l(\omega)+l(w_{2})+l(\omega)$. On en d\'eduit $S(w)=S(\omega)S(w_{2})S(\omega)$. A l'aide de \ref{epinglages} (4), on calcule $S(\omega)=tS(\omega)^{-1}$, o\`u $t=\prod_{\alpha>0, \omega(\alpha)<0}\check{\alpha}(-1)$. L'ensemble des $\alpha>0$ tels que $\omega(\alpha)<0$ est $\{\alpha_{4},\alpha_{3}+\alpha_{4},\alpha_{4}+\alpha_{5},\alpha_{3}+\alpha_{4}+\alpha_{5}\}$. D'o\`u $t=1$. On obtient $S(w)=S(\omega)^{-1}S(w_{2})S(\omega)$. Par construction, $S(w_{2})$ est contenu dans le  sous-groupe $SL(2) $ de $G$  associ\'e \`a la racine $\alpha_{2}$. Donc $S(w)$ est contenu dans le groupe $SL(2)$ associ\'e \`a la racine $\omega(\alpha_{2})=\alpha_{2}+\alpha_{3}+2\alpha_{4}+\alpha_{5}$. On constate que cette racine est orthogonale \`a $\alpha'_{0}$, $\alpha''_{0}$, $\underline{\alpha}_{2}$ et $\underline{\alpha}_{4}$. Donc $S(w)$ fixe les espaces radiciels associ\'es \`a ces racines. Les \'el\'ements $E_{\underline{\beta}_{0}}=E_{\beta_{0}}$, resp. $E_{\underline{\beta}_{2,4}}$ appartiennent \`a la somme des espaces radiciels associ\'es aux racines $\alpha'_{0}$ et $\alpha''_{0}$, resp. $\underline{\alpha}_{2}$ et $\underline{\alpha}_{4}$. Donc $Ad(n)$ fixe $E_{\underline{\beta}_{0}}$ et $E_{\underline{\beta}_{2,4}}$. Cela prouve (1).
  
  D'apr\`es (2) et par construction de $N$, $n$ appartient à $Norm_{G(F^{nr})}(M,s,N)$.

  Calculons $\epsilon^{\flat}_{x,N}\circ c_{x}(n)$.   Introduisons un espace $V$ de dimension $7$ sur $\bar{{\mathbb F}}_{q}$, muni d'une base $(e_{i})_{i=1,...,7}$ et de la forme quadratique d\'efinie par $q(e_{i},e_{8-i})=(-1)^{i+1}$ pour tout $i\not=4$, $q(e_{4},e_{4})=-2$. On note $\tilde{G}_{*}$ son groupe sp\'ecial orthogonal. On peut identifier $G_{*}$ \`a son rev\^etement spinoriel  de sorte que $T_{*}$ s'identifie \`a l'image r\'eciproque du tore diagonal et  $\bar{E}_{\beta_{2}}$, $\bar{E}_{\beta_{4}}$ et $\bar{E}_{\beta_{3,5}}$ s'identifient respectivement \`a $E_{1,2}+E_{6,7}$, $E_{2,3}+E_{5,6}$, $E_{3,4}+E_{4,5}$, avec les notations de \ref{Bn}. On note $\pi:G_{*}\to \tilde{G}_{*}$ la projection.  Introduisons la section de Springer $S_{*}$ de $W_{*}$ \`a valeurs dans $G_{*}$ associ\'ee \`a l'\'epinglage que l'on vient de d\'ecrire. Il r\'esulte des constructions que les r\'eductions dans $G_{x}$ de $S(w_{4})$, resp. $S(w_{3})S(w_{5})$, $S(w_{2})$, sont les images dans ce groupe des \'el\'ements $S_{*}(w_{*,2})$, resp. $S_{*}(w_{*,3})$, $S(w_{*,1})$, de $G_{*}$. Donc $c_{x}(n)$ est l'image dans $G_{x}$ de l'\'el\'ement $S_{*}(w_{*})$ de $G_{*}$. On peut aussi bien remplacer $G_{x}$ par $G_{*}$ et calculer $\epsilon_{*,N}^{\flat}(S_{*}(w_{*}))$, où $\epsilon_{*,N}^{\flat}$ est la fonction de \ref{faisceauxcaracteres} pour le groupe $G_{*}$ . On effectue le calcul comme en \ref{alpha6}. 
  On choisit  $\underline{N}_{*}=\bar{E}_{\beta_{3,5}}+\bar{E}_{\beta_{4}}+\bar{E}_{\beta_{2}}$ et $Z=E_{1,1}+E_{2,2}-E_{6,6}-E_{7,7}$. Pour $\nu\in \bar{{\mathbb F}}_{q}^{\times}$, un calcul matriciel permet de calculer l'\'el\'ement $u_{\nu}\in U_{B_{*}}$ tel que
  $$u_{\nu}^{-1}(\nu Z+\underline{N}_{*})u_{\nu}=\nu Z+N_{*}.$$
  On calcule ensuite le produit $u_{\nu}S_{*}(w_{*})$ sous la forme $\underline{v}_{\nu}t_{\nu}v_{\nu}$, cf. \ref{alpha6}(4) en utilisant  l'égalité (1) et en calculant successivement
  $$u_{\nu}S_{*}(w_{*,2}),\, u_{\nu}S_{*}(w_{*,2})S(w_{*,3}),\, u_{\nu}S_{*}(w_{*,2})S(w_{*,3})S(w_{*,2})S(w_{*,1})S(w_{*,2}),$$
  $$ u_{\nu}S_{*}(w_{*,2})S(w_{*,3})S(w_{*,2})S(w_{*,1})S(w_{*,2})S(w_{*,3}),$$
  $$ u_{\nu}S_{*}(w_{*,2})S(w_{*,3})S(w_{*,2})S(w_{*,1})S(w_{*,2})S(w_{*,3})S(w_{*,2}).$$
  A chaque fois, le calcul se passe soit dans le groupe $\tilde{G}_{*}$, soit  dans un groupe $SL(2)$ ou $SL(3)$ pour le terme central. On obtient
  $$t_{\nu}=\left(\check{\gamma}_{1}(2\nu^5)\check{\gamma}_{2}(4\nu^{10})\check{\gamma}_{3}(2\nu^5)\right)^{-1}.$$
  Mais $\check{\gamma}_{1}+2\check{\gamma}_{2}+\check{\gamma}_{3}$ appartient \`a $X_{*}(Z(M_{*})^0)$. Donc $t_{\nu}\in Z(M_{*})^0$. D'apr\`es \ref{alpha6}(7), on en d\'eduit $\epsilon_{*,N}^{\flat}(S_{*}(w_{*}))=1$, d'où
  
  $\epsilon_{N}^{\flat}(n)=1$. 
  
    On a bien démontré les assertions (1) et (2) de \ref{preuveE6ram}. 
  
   \subsubsection{Sommet $s_{\beta_{4}}$ pour $G$ de type $(E_{6},ram)$ \label{beta4}}

 Posons $x=s_{\beta_{4}}$. On a $G_{x}=(SL(4)\times SL(2))/\bar{\zeta}_{1/4}$, le groupe $\bar{\zeta}_{1/4}$ s'envoyant dans $SL(4)\times SL(2)$ par
$$\underline{i}\mapsto( \check{\alpha}_{0}(\underline{i})\check{\alpha}_{1,6}(-1)\check{\alpha}_{3,5}(-\underline{i}), \check{\alpha}_{2}(-1)).$$
 Notons $G_{*}=SL(4)$ et $G_{**}=SL(2)$ les deux composantes de $G_{x,SC}$. On note $(B_{*},T_{*})$ et $(B_{**},T_{**})$ les paires de Borel de $G_{*}$ et $G_{**}$ qui se d\'eduisent de $(B_{x},T_{x})$ et on note $W_{*}$ et $W_{**}$ les groupes de Weyl. L'image r\'eciproque dans $G_{*}\times G_{**}$ de $M_{s}$ est $M_{*}\times G_{**}$, o\`u $M_{*,SC}=SL(2)\times SL(2)$. L'\'el\'ement $N$ s'\'ecrit $N_{*}\oplus N_{**}$, o\`u $N_{*}=\bar{E}_{\beta_{0}}+\bar{E}_{\beta_{3,5}}$ et $N_{**}=\bar{E}_{\beta_{2}}$.   On identifie le groupe  $W_{*}$   au groupe des permutations de l'ensemble $\{1,2,3,4\}$, les sym\'etries \'el\'ementaires associ\'ees aux racines $\beta_{0},\beta_{1,6},\beta_{3,5}$ s'identifiant respectivement aux permutations $w_{*,1}$ de $1$ et $2$, $w_{*,2}$ de $2$ et $3$, $w_{*,3}$ de $3$ et $4$. 
 Introduisons l'\'el\'ement $w_{*}\in W_{*}$ qui envoie $(1,2,3,4)$ sur $(3,4,1,2)$. Compl\'et\'e par l'\'el\'ement neutre de $W_{**}$, on obtient un \'el\'ement de $W_{*}\times W_{**}$ dont l'image dans $W^{G_{x}}$ normalise $M_{s}$ et  s'envoie elle-m\^eme sur l'\'el\'ement non trivial de  
  $W^{G_{x}}(M_{s})$.  On a l'\'egalit\'e $w_{*}=w_{*,2}w_{*,1}w_{*,3}w_{*,2}$. Notons $S_{*}$ la section de Springer du groupe $G_{*}$ relative à l'épinglage $\bar{E}_{\beta_{0}},\bar{E}_{\beta_{1,6}},\bar{E}_{\beta_{3,5}}$ de ce groupe. Posons $\bar{n}_{*}=S_{*}(w_{*})$. 
  
  Rappelons que l'élément $E_{\beta_{0}}$ a été construit de la façon suivante: on a fixé  deux générateurs $E_{\alpha'_{0}}$, resp. $E_{\alpha''_{0}}$, de $\mathfrak{u}_{\alpha'_{0},\mathfrak{o}_{F^{nr}}}$, resp. $\mathfrak{u}_{\alpha''_{0},\mathfrak{o}_{F^{nr}}}$, tels que $\theta(E_{\alpha'_{0}})=E_{\alpha''_{0}}$, on a fixé une uniformisante $\varpi_{E}$ de $E$ telle que $\varpi_{E}^2\in F^{nr}$ et on a posé $E_{\beta_{0}}=\varpi_{E}(E_{\alpha'_{0}}-E_{\alpha''_{0}})$.
  Pour $\alpha=\alpha'_{0}$ ou $\alpha''_{0}$, on note $E_{-\alpha}$ l'unique élément de $\mathfrak{u}_{-\alpha}$ tel que $[E_{\alpha},E_{-\alpha}]=\check{\alpha}$ et on pose $n_{\alpha}=exp(\nu\varpi_{E}E_{\alpha})exp(-\nu\varpi_{E}^{-1}E_{-\alpha})exp(\nu\varpi_{E}E_{\alpha})$, où $\nu=1$ si $\alpha=\alpha'_{0}$, $\nu=-1$ si $\alpha=\alpha''_{0}$. Introduisons la section de Springer $S$ de $G$. Rappelons que, pour tout $\alpha\in \Sigma$, on note $w_{\alpha}$ la symétrie associée à $\alpha$. Posons
  $$n_{1}=n_{\alpha'_{0}}n_{\alpha''_{0}},\, n_{2}=S(w_{\alpha_{1}})S(w_{\alpha_{6}}),\, n_{3}=S(w_{\alpha_{3}})S(w_{\alpha_{5}}).$$
   Ce sont des éléments de $K_{x}^{0,nr}$ dont les réductions dans $G_{x}$ sont les images dans ce groupe de $S_{*}(w_{*,1})$, $S_{*}(w_{*,2})$ et $S_{*}(w_{*,3})$. On pose $n_{\beta_{4}}=n_{2}n_{1}n_{3}n_{2}$. On va prouver que cet élément vérifie les conditions de \ref{preuveE6ram}, en particulier les relations (1) et (2). On pose simplement $n=n_{\beta_{4}}$. 
  
  L'élément $n$ appartient  par construction à $G(F^{nr})$ et plus précisément à $K_{x}^{0,nr}$. Son  image dans $G_{x}$ est la m\^eme que celle de $\bar{n}_{*}$. Cette image normalise $M_{s}$, donc $n$ appartient à $Norm_{G(F^{nr})}(M,s)$.   
  
    En \ref{epinglages}(3), on a introduit un couple $\{\pm \dot{E}_{\alpha}\}$  d'éléments de $\mathfrak{u}_{\alpha}$ pour tout $\alpha\in \Sigma$. Fixons $\dot{E}_{\alpha'_{0}}$ et $\dot{E}_{\alpha''_{0}}$. On suppose que $\dot{E}_{\alpha''_{0}}=\theta(\dot{E}_{\alpha'_{0}})$. Pour $\alpha=\alpha'_{0}$ ou $\alpha''_{0}$,  on note $\dot{E}_{-\alpha}$ l'unique élément de $\mathfrak{u}_{-\alpha}$ tel que $[\dot{E}_{\alpha},\dot{E}_{-\alpha}]=\check{\alpha}$ et on pose $\dot{n}_{\alpha}=exp(\dot{E}_{\alpha})exp(-\dot{E}_{-\alpha})exp(\dot{E}_{\alpha})$. Il existe un élément $\lambda\in \mathfrak{o}^{\times}_{F^{nr}}$ tel que $E_{\alpha'_{0}}=\lambda\dot{E}_{\alpha'_{0}}$. Puisque $E_{\alpha''_{0}}=\theta(E_{\alpha'_{0}})$, on a aussi $E_{\alpha''_{0}}=\lambda\dot{E}_{\alpha''_{0}}$. 
 On a donc l'égalité $n_{\alpha}=\check{\alpha}(\nu\lambda\varpi_{E})\dot{n}_{\alpha}$ pour $\alpha=\alpha'_{0}$ ou $\alpha=\alpha''_{0}$, où $\nu$ est comme ci-dessus. Posons $\dot{n}_{1}=\dot{n}_{\alpha'_{0}}\dot{n}_{\alpha''_{0}}$ et $\dot{n}=n_{2}\dot{n}_{1}n_{3}n_{2}$. Posons $\gamma'_{0}=n_{2}(\alpha'_{0})$, $\gamma''_{0}=n_{2}(\alpha''_{0})$. Alors 
    
    (1) $n=\check{\gamma}'_{0}(\lambda\varpi_{E})\check{\gamma}''_{0}(-\lambda\varpi_{E})\dot{n}$. 
    
    On calcule
    $\gamma'_{0}=-(\alpha_{1}+\alpha_{2}+2\alpha_{3}+2\alpha_{4}+\alpha_{5})$, $\gamma''_{0}=-(\alpha_{2}+\alpha_{3}+2\alpha_{4}+2\alpha_{5}+\alpha_{6})$.

 Montrons que
 
 (2) $Ad(n)$ conserve $M$ et  l'\'epinglage affine $\mathfrak{ E}_{a}^{M,nr}$.
 
 Puisque l'élément $w_{*}$ permute les racines $\beta_{0}$ et $\beta_{3,5}$, l'élément $n$ envoie le couple $(\alpha_{3},\alpha_{5})$ sur l'un des couples $(\alpha'_{0},\alpha''_{0})$ ou $(\alpha''_{0},\alpha'_{0})$. 
 L'\'el\'ement $\dot{n}$ appartient au groupe $T_{S}S(W)$ introduit en \ref{epinglages}. D'après \ref{epinglages}(3), $Ad(\dot{n})$ envoie $E_{\beta_{3,5}}=E_{\alpha_{3}}+E_{\alpha_{5}}$ sur un élément de la forme $\pm  \dot{E}_{\alpha'_{0}}\pm \dot{E}_{\alpha''_{0}}$. On a $<\alpha'_{0},\check{\gamma}'_{0}>=1$, $<\alpha'_{0},\check{\gamma}'_{0}>=1$,  $<\alpha'_{0},\check{\gamma}''_{0}>=0$, $<\alpha''_{0},\check{\gamma}'_{0}>=0$, $<\alpha''_{0},\check{\gamma}''_{0}>=1$. D'après (1),  $Ad(n)$ envoie $E_{\beta_{3,5}}$ sur $\pm \lambda\varpi_{E}\dot{E}_{\alpha'_{0}}\pm \lambda\varpi_{E}\dot{E}_{\alpha''_{0}}=\varpi_{E}(\nu_{1} E_{\alpha'_{0}}+\nu_{2} E_{\alpha''_{0}})$, où $\nu_{1},\nu_{2}\in \{\pm 1\}$. Cet élément doit appartenir à $\mathfrak{g}(F^{nr})$, ce qui implique $\nu_{2}=-\nu_{1}$, d'où $Ad(n)(E_{\beta_{3,5}})=\nu_{1}E_{\beta_{0}}$. 
 L'égalité ci-dessus se réduit dans $\mathfrak{g}_{x}$ en $\bar{n}_{*}(\bar{E}_{\beta_{3,5}})=\nu_{1}\bar{E}_{\beta_{0}}$. Mais $\bar{n}_{*}$ appartient à l'image de $S_{*}$. D'après \ref{epinglages}(2), on a $n_{*}(\bar{E}_{\beta_{3,5}})=\bar{E}_{\beta_{0}}$. Donc $\nu_{1}=1$ et $Ad(n)$ envoie $E_{\beta_{3,5}}$ sur $E_{\beta_{0}}$. Un calcul similaire que l'on laisse au lecteur prouve que $Ad(n)$ envoie $E_{\beta_{0}}$ sur $E_{\beta_{3,5}}$. En se rappelant que $E_{\underline{\beta}_{0}}=E_{\beta_{0}}$ et  $E_{\underline{\beta}_{1,5}}=E_{\beta_{3,5}}$, $Ad(n)$ permute $E_{\underline{\beta}_{0}}$ et $E_{\underline{\beta}_{1,5}}$.  
 
  Par construction, $n$  appartient au produit des sous-groupes $SL(2)$ de $G$ associ\'ees aux racines $\alpha'_{0}$ et $\alpha''_{0}$,  resp. $\alpha_{1}$ et $\alpha_{6}$, $\alpha_{3}$ et $\alpha_{5}$. Toutes ces racines sont orthogonales \`a $\alpha_{2}=\underline{\beta}_{3}$ donc tous ces sous-groupes $SL(2)$ fixent $E_{\underline{\beta}_{3}}$. Cela entraîne que $Ad(n)$ fixe $E_{\underline{\beta}_{3}}$. 
  
  Il nous reste à prouver que $Ad(n)$ fixe $E_{\underline{\beta}_{2,4}}=E_{\underline{\alpha}_{2}}+E_{\underline{\alpha}_{4}}$. Posons
  $\tau_{1}=n_{2}^{-1}\dot{n}_{1}n_{2}$, $\tau_{3}=n_{2}^{-1}n_{3}n_{2}$. D'après \ref{epinglages}(1), on a $n_{2}^{-1}=\check{\alpha}_{1}(-1)\check{\alpha}_{6}(-1)n_{2}$, d'où $n=t\tau_{1}\tau_{3}$, où $t=\check{\gamma}'_{0}(\lambda\varpi_{E})\check{\gamma}''_{0}(-\lambda\varpi_{E})\check{\alpha}_{1}(-1)\check{\alpha}_{6}(-1)$. On calcule $\underline{\alpha}_{2}(t)=\underline{\alpha}_{4}(t)=1$. Pour prouver que $Ad(n)$ fixe $E_{\underline{\beta}_{2,4}}$, il suffit donc de prouver
  
  (3) $Ad(\tau_{1})$ et $Ad(\tau_{3})$ multiplient $E_{\underline{\beta}_{2,4}}$ par $-1$.
  
  Posons $\tau'=n_{2}^{-1}\dot{n}_{\alpha'_{0}}n_{2}$, $\tau''=n_{2}^{-1}\dot{n}_{\alpha''_{0}}n_{2}$. On a $\tau_{1}=\tau'\tau''$.  L'élément $\tau'$, resp. $\tau''$, est une symétrie du groupe $SL(2) $ associé  à la racine $\gamma'_{0}$, resp. $\gamma''_{0}$. On a de plus $\tau''=\theta(\tau')$. Posons  $\delta'=-(\alpha_{1}+\alpha_{2}+2\alpha_{3}+3\alpha_{4}+2\alpha_{5}+\alpha_{6})$ et $\delta''=-(\alpha_{2}+\alpha_{3}+\alpha_{4}+\alpha_{5})$.   Le groupe $SL(2)$ associ\'e \`a $\gamma''_{0}$   conserve les plans $\mathfrak{u}_{\underline{\alpha}_{2}}\oplus \mathfrak{u}_{\delta''}$ et  $\mathfrak{u}_{\underline{\alpha}_{4}}\oplus \mathfrak{u}_{\delta'}$. L'\'el\'ement  $Ad(\tau'')$    permute les deux droites indiqu\'ees. Donc $Ad(\tau'')(E_{\underline{\alpha}_{2}})\in \mathfrak{u}_{\delta''}$ et $Ad(\tau'')(E_{\underline{\alpha}_{4}})\in \mathfrak{u}_{\delta'}$. Mais $\delta'$ et $\delta''$ sont fixées par $\theta$ donc $\theta$ agit par l'identité sur $\mathfrak{u}_{\delta'}$ et $\mathfrak{u}_{\delta''}$ et fixe $Ad(\tau'')(E_{\underline{\alpha}_{2}})$ et $Ad(\tau'')(E_{\underline{\alpha}_{4}})$. 
  On a alors
  $$Ad(\tau')Ad(\tau'')(E_{\underline{\alpha}_{2}})=\theta\circ Ad(\tau'')\circ \theta \circ Ad(\tau'')(E_{\underline{\alpha}_{2}})=\theta\circ Ad(\tau'')^2(E_{\underline{\alpha}_{2}})$$
  $$=\theta\circ Ad(\check{\gamma}''_{0}(-1))(E_{\underline{\alpha}_{2}})=-\theta(E_{\underline{\alpha}_{2}})=-E_{\underline{\alpha}_{4}}.$$
  De m\^eme $Ad(\tau')Ad(\tau'')(E_{\underline{\alpha}_{4}})=-E_{\underline{\alpha}_{2}}$. Donc $Ad(\tau_{1})$ multiplie $E_{\underline{\beta}_{2,4}}$ par $-1$. Un calcul similaire vaut pour $\tau_{3}$, les racines $\delta'$ et $\delta''$ étant remplacées par $\alpha_{1}+\alpha_{2}+\alpha_{3}+\alpha_{4}+\alpha_{5}+\alpha_{6}$ et $\alpha_{4}$. Cela prouve (3) et achève la preuve de (2).

 Le terme  $\epsilon_{x,N}^{\flat}\circ c_{x}(n)$  est \'egal à $\epsilon_{*, N}^{\flat}(S_{*}(w_{*}))$, où $\epsilon_{*,N}^{\flat}$ est la fonction relative au groupe   $G_{*}$. On utilise la m\^eme m\'ethode qu'en \ref{alpha6}, en  posant  $\underline{N}_{*}=\bar{E}_{\beta_{0}}+\bar{E}_{\beta_{1,6}}+\bar{E}_{\beta_{3,5}}$ et $Z=\check{\beta}_{0}+2\check{\beta}_{1,6}+\check{\beta}_{3,5}\in \mathfrak{t}_{x}$. Le calcul matriciel dans $G_{*}\simeq SL(4)$ est facile. Avec les notations de \ref{alpha6}, on trouve
 $$t_{\nu}=(\check{ \beta}_{0}+2\check{\beta}_{1,6}+\check{\beta}_{3,5})(2\nu)^{-2}.$$
  Mais $\check{\beta}_{0}+2\check{\beta}_{1,6}+\check{\beta}_{3,5} $ appartient \`a $X_{*}(Z(M_{*})^0)$.  On en d\'eduit
$\epsilon_{x,N}^{\flat}\circ c_{x}(n)=1$, ou encore

$\epsilon_{N}^{\flat}(n)=1$. 

On a bien démontré les assertions (1) et (2) de \ref{preuveE6ram}. Cela achève la preuve du théorème \ref{epsilonflatcasstable}.

\subsection{Isomorphismes d'espaces $FC^{st}(\mathfrak{h}_{SC}(F))$} 

\subsubsection{Définition des isomorphismes}\label{isomorphismesFC}

 {\bf On suppose que} $G$ {\bf est  défini sur} $F$, {\bf qu'il est adjoint et absolument simple.}
  
 Soit $H,H'\in {\cal L}_{F}^{nr,st}$. On suppose que ces deux groupes sont conjugués par un élément de $G(F^{nr})$. Puisqu'il s'agit de $F^{nr}$-Levi, ils sont m\^eme conjugués par un élément de $G_{SC}(F^{nr})$. 
  Nous allons définir un isomorphisme
 $$(1) \qquad \iota_{H',H}:FC^{st}(\mathfrak{h}_{SC}(F))\to FC^{st}(\mathfrak{h}'_{SC}(F)).$$
 Avec les m\^emes notations qu'en \ref{resolution}, les espaces ci-dessus sont respectivement isomorphes à 
  $\oplus_{s\in \underline{S}^{st}(H_{AD})}FC^{st}(\mathfrak{h}_{SC,s}({\mathbb F}_{q}))$ et $\oplus_{s'\in \underline{S}^{st}(H'_{AD})}FC^{st}(\mathfrak{h}'_{SC,s'}({\mathbb F}_{q}))$. Fixons un élément $x\in G_{SC}(F^{nr})$ tel que $xHx^{-1}=H'$. 
   De $x$ se déduit un isomorphisme $Ad(x)^{Imm}:Imm_{F^{nr}}(H_{AD})\to Imm_{F^{nr}}(H'_{AD})$. Comme en \ref{resolution}, il s'en déduit une bijection $\underline{Ad(x)}:\underline{S}^{st}(H)\to \underline{S}^{st}(H')$ de sorte que, pour $s\in \underline{S}^{st}(H)$, $\underline{Ad(x)}(s)$ soit l'unique élément de $\underline{S}^{st}(H')$ qui soit conjugué à $Ad(x)^{Imm}(s)$ par un élément de $H'_{AD}(F^{nr})$. Soit $s\in \underline{S}^{st}(H)$, posons $s'=\underline{Ad(x)}(s)$.  Quitte à multiplier $x$ à gauche par un élément de $H'_{SC}(F^{nr})$, on peut supposer que $Ad(x)^{Imm}(s)=s'$.  De $Ad(x)$ se déduit alors un isomorphisme
$Ad(x)_{s}:H_{AD,s}\to H'_{AD,s'}$. Soit ${\cal E}\in {\bf FC}^{st}_{{\mathbb F}_{q}}(\mathfrak{h}_{SC,s})$. Par $Ad(x)_{s}$, ${\cal E}$ se transforme en un faisceau ${\cal E}'$ qui appartient à ${\bf FC}^{st}_{{\mathbb F}_{q}}(\mathfrak{h}'_{SC,s'})$. Fixons des éléments $N\in {\cal O}_{{\cal E}}^{\Gamma_{{\mathbb F}_{q}}}$ et $N'\in {\cal O}_{{\cal E}'}^{{\mathbb F}_{q}}$. Quitte à multiplier encore $x$ à gauche par un élément de $K_{s'}^{H'_{SC},0,nr}$, on peut supposer que $Ad(x)_{s}(N)=N'$. Parce que tous nos objets sont définis sur $F$ ou ${\mathbb F}_{q}$, on voit que $Fr(x)x^{-1}$ appartient à $Norm_{G_{SC}(F^{nr})}(H',s',N')$. Fixons une action de Frobenius sur ${\cal E}'$. Elle est déterminée par un scalaire $r_{N'}\in {\mathbb C}^{\times}$. On définit une action de Frobenius sur ${\cal E}$ en utilisant le scalaire 
$$(2) \qquad r_{N}=r_{N'}\epsilon_{N'}^{\flat}(Fr(x)x^{-1}).$$
Notons $f_{{\cal E}}$ et $f_{{\cal E}'}$ les fonctions caractéristiques de ${\cal E}$ et ${\cal E}'$.  Faisons maintenant varier le faisceau ${\cal E}$. La famille formée des fonctions $f_{{\cal E}} $ est une base de $FC^{st}(\mathfrak{h}_{SC, s}({\mathbb F}_{q}))$ tandis que la famille formée des fonctions $f_{{\cal E}'} $ est une base de $FC^{st}(\mathfrak{h}'_{SC,s'}({\mathbb F}_{q}))$. On d\'efinit un isomorphisme  
$$\iota_{s}:FC^{st}(\mathfrak{h}_{SC,s}({\mathbb F}_{q}))\to FC^{st}(\mathfrak{h}'_{SC, s'}({\mathbb F}_{q}))$$ 
par l'\'egalit\'e
$$(3) \qquad \iota_{s}(f_{{\cal E}})=\vert H'_{AD,s'}({\mathbb F}_{q})\vert^{-1}\vert H_{AD,s_{H}}({\mathbb F}_{q})\vert f_{{\cal E}'}.$$
  En faisant varier $s$, on obtient l'isomorphisme (1). 

\begin{lem}{La définition de $\iota_{H',H}$ ne dépend pas des choix effectués.}\end{lem} 

Preuve. 
Les choix sont ceux des ensembles $\underline{S}(H_{AD})$ et $\underline{S}(H'_{AD})$ et, pour $s$ et ${\cal E}$ fixés,  des éléments $N$ et $N'$, de l'élément $x\in G_{SC}(F^{nr})$ et du scalaire $r_{N'}$. Il est clair d'après (2) que le choix de $r_{N'}$ ne modifie pas la définition de $\iota_{H',H}$.   Quant à $x$, il y a un double choix. D'un premier choix d'un élément qui conjugue $H$ en $H'$, on déduit une bijection $\underline{Ad(x)}$. Ensuite, après avoir choisi un sommet $s$, dont on note $s'$ l'image par $\underline{Ad(x)}$, un faisceau ${\cal E}$ et des éléments $N$ et $N'$, on modifie $x$ de sorte qu'il appartienne à $Norm_{G_{SC}(F^{nr})}(H',s',N')$. Montrons 
d'abord que la bijection $\underline{Ad(x)}$ ne dépend pas du choix de $x$.  A ce point, on ne peut changer $x$ qu'en le multipliant à gauche par un élément de $Norm_{G_{SC}(F^{nr})}(H')$. Pour $n\in Norm_{G_{SC}(F^{nr})}(H')$ et $s'\in \underline{S}^{st}(H'_{AD})$, l'unique élément de $\underline{S}^{st}(H'_{AD})$ qui est conjugué à $ns'$ par un élément de $H'_{AD}(F^{nr})$ est $s'$ lui-m\^eme d'après le lemme \ref{lemmeautomorphismes}. Donc $\underline{Ad(nx)}=\underline{Ad(x)}$. Cela étant, 
  les termes $s$, $s'$, $N$ et $N'$ étant fixés, le deuxième choix de $x$  est unique à multiplication à gauche près par un élément $y\in Norm_{G_{SC}(F^{nr})}(H',s',N')$. Remplacer $x$ par $yx$ pour un tel $y$ remplace 
$\epsilon_{N'}^{\flat}(Fr(x)x^{-1})$ par $\epsilon_{N'}^{\flat}(Fr(y)Fr(x)x^{-1}y^{-1})$. Parce que $\epsilon_{N'}^{\flat}$ est un caractère, ceci vaut $\epsilon_{N'}^{\flat}(Fr(y))\epsilon_{N'}^{\flat}(Fr(x)x^{-1})\epsilon_{N'}^{\flat}(y)^{-1}$. Parce que $\epsilon_{N'}^{\flat}$ est invariant par Frobenius, c'est encore égal à $\epsilon_{N'}^{\flat}(Fr(x)x^{-1})$  et la définition n'a pas changé. L'indépendance des choix de $N$ et $N'$ résulte d'un calcul formel similaire à celui de la preuve de la proposition \ref{resolution}. On le laisse au lecteur.    Si on modifie l'ensemble $\underline{S}(H_{AD})$, un sommet $s$  est remplacé  par $h(s)$   pour un élément $h\in H_{SC}(F)$ et le sommet $s'$ associé ne change pas. On peut remplacer $x$ par $xh^{-1}$, ce qui ne modifie pas $\epsilon_{N'}^{\flat}(Fr(x)x^{-1})$ puisque $h\in H_{SC}(F)$. On voit que le nouvel isomorphisme $\iota_{h(s)}$ est le composé de  $\iota_{s}$ et de l'isomorphisme de transport de structure $FC^{st}(\mathfrak{h}_{SC,h(s)}({\mathbb F}_{q}))\to  FC^{st}(\mathfrak{h}_{SC,s}({\mathbb F}_{q}))$ défini par $Ad(h)^{-1}$. Ce dernier est inessentiel puisqu'il devient l'identité quand on envoie les deux espaces dans $FC^{st}(\mathfrak{h}_{SC}(F))\subset I(\mathfrak{h}_{SC}(F))$. Si on modifie $\underline{S}(H'_{AD})$, pour un sommet $s$ fixé, $s'$ est remplacé par un sommet $h'(s')$ pour un $h'\in H'_{SC}(F)$. L'élément $h'$ détermine un isomorphisme $Ad(h')_{s'}:H'_{AD,s'}\to H'_{AD,h'(s')}$ qui définit un isomorphisme de transport de structure $FC^{st}(\mathfrak{h}'_{SC,s'}({\mathbb F}_{q}))\to  FC^{st}(\mathfrak{h}'_{SC,h'(s')}({\mathbb F}_{q}))$. Pour réaliser cet isomorphisme, on doit remplacer le faisceau ${\cal E}'$ par $Ad(h')_{s'}({\cal E}')$,  $N'$ par $Ad(h')_{s'}(N')$ et poser $r_{Ad(h')_{s'}(N')}=r_{N'}$. Pour construire le nouveau $\iota_{s}$, on peut remplacer $x$ par $h'x$. 
Le terme $r_{N'}\epsilon_{N'}^{\flat}(Fr(x)x^{-1})$ est alors remplacé par $r_{Ad(h')_{s'}(N')}\epsilon_{ Ad(h')_{s'}(N')}^{\flat}(Fr(h')Fr(x)x^{-1}{h'}^{-1})=r_{N'}\epsilon_{ Ad(h')_{s'}(N')}^{\flat}(h'Fr(x)x^{-1}{h'}^{-1})$ puisque $h'$ est fixé par $Fr$.  D'après \ref{caractereepsilonflat}(14),  on a l'égalité $\epsilon_{ Ad(h')_{s'}(N')}^{\flat}(h'Fr(x)x^{-1}{h'}^{-1})=\epsilon_{N'}^{\flat}(Fr(x)x^{-1})$. On voit que le nouveau $\iota_{s}$ est le composé  de l'isomorphisme de transport de structure $FC^{st}(\mathfrak{h}'_{SC,s'}({\mathbb F}_{q}))\to  FC^{st}(\mathfrak{h}'_{SC,h'(s')}({\mathbb F}_{q}))$  et de l'ancien $\iota_{s}$. Encore une fois, le premier isomorphisme  est inessentiel. Cela achève la preuve. $\square$

\subsubsection{Composition d'isomorphismes $\iota_{H',H}$}\label{compositionderechef}
On considère trois éléments $H_{i}\in {\cal L}_{F}^{nr,st}$, pour $i\in \{1,2,3\}$, que l'on suppose conjugués deux à deux par un élément de $G(F^{nr})$. 
\begin{lem}{On a l'égalité $\iota_{H_{3},H_{2}}\circ\iota_{H_{2},H_{1}}=\iota_{H_{3},H_{1}}$.}\end{lem}  

Preuve. Elle est analogue à celle du lemme \ref{composition}. Brièvement, considérons  pour $i=1,2,3$ un sommet $s_{i}\in Imm(H_{i,AD})$, un faisceau ${\cal E}_{i}\in {\bf FC}^{st}_{{\mathbb F}_{q}}(\mathfrak{h}_{i,SC,s_{i}})$ et un élément $N_{i}\in {\cal O}_{{\cal E}_{i}}^{\Gamma_{{\mathbb F}_{q}}}$. Considérons aussi deux éléments $x,y\in G_{SC}(F^{nr})$ et supposons que les septuplets $(s_{1},{\cal E}_{1},N_{1},s_{2},{\cal E}_{2},N_{2},x)$, resp.  $(s_{2},{\cal E}_{2},N_{2},s_{3},{\cal E}_{3},N_{3},y)$, vérifient les conditions imposées dans la construction de l'isomorphisme $\iota_{H_{2},H_{1}}$, resp. $\iota_{H_{3},H_{2}}$. Alors le septuplet $(s_{1},{\cal E}_{1},N_{1},s_{3},{\cal E}_{3},N_{3},yx)$ vérifie les conditions imposées dans la construction de l'isomorphisme $\iota_{H_{3},H_{1}}$. On fixe une action de Frobenius sur ${\cal E}_{3}$ déterminée par un scalaire $r_{3}$ et on applique les constructions. Le faisceau ${\cal E}_{1}$ se retrouve muni de deux actions de Frobenius: l'une obtenue par la construction de $\iota_{H_{3},H_{1}}$, déterminée par un scalaire $r_{1}$, l'autre obtenue par la composition de $\iota_{H_{3},H_{2}}$ et de $\iota_{H_{2},H_{1}}$, déterminée par un scalaire que l'on note $\tilde{r}_{1}$. Comme en \ref{composition}, le lemme équivaut à l'égalité $r_{1}=\tilde{r}_{1}$. Avec des notations compréhensibles, on  a
$$\tilde{r}_{1}=r_{2}\epsilon_{N_{2}}^{\flat}(Fr(x)x^{-1})=r_{3}\epsilon^{\flat}_{N_{3}}(Fr(y)y^{-1})\epsilon_{N_{2}}^{\flat}(Fr(x)x^{-1}).$$
D'après \ref{caractereepsilonflat}(14) appliqué à l'isomorphisme $Ad(y)$, on a $\epsilon_{N_{2}}^{\flat}(Fr(x)x^{-1})=\epsilon_{N_{3}}(yFr(x)x^{-1}y^{-1})$, d'où
$$\tilde{r}_{1}=r_{3}\epsilon_{N_{3}}(Fr(y)y^{-1})\epsilon_{N_{3}}(yFr(x)x^{-1}y^{-1})=r_{3}\epsilon_{N_{3}}(Fr(y)Fr(x)x^{-1}y^{-1}).$$
Mais ce dernier terme est égal à $r_{1}$, d'où l'égalité $\tilde{r}_{1}=r_{1}$ cherchée. Cela achève la démonstration. $\square$

\subsubsection{Comparaison avec le transfert endoscopique}\label{comparaison}

 \begin{prop}{Soient  $H,H'\in {\cal L}_{F}^{nr,st}$. Supposons que $H$ et $H'$ soient stablement conjugués. Alors $\iota_{H,H'}$ coïncide avec l'isomorphisme de transfert défini en \ref{automorphismes}.}\end{prop}

 Preuve. L'assertion est triviale si $H=G$. On suppose que $H\not=G$. Fixons $g\in G_{SC}(F^{nr})$ tel que $gHg^{-1}=H'$ et $Fr(g)g^{-1}\in H'(F^{nr})$. L'homomorphisme $Ad(g)$ se quotiente en un torseur intérieur $\phi:H_{AD}\to H'_{AD}$. 
   A l'aide de $\phi$, on a défini en \ref{resolution} un autre isomorphisme  $\iota_{\phi}$,  et on a prouvé qu'il coïncidait avec l'isomorphisme de transfert, cf.  théorème \ref{enonce}. Il nous suffit donc de prouver l'égalité
 
 (1)  $\iota_{H',H}=\iota_{\phi}$.
 
 Soit $s\in S^{st}(H_{AD})$. Quitte à multiplier $g$ à gauche par un élément de $H'_{SC}(F^{nr})$, on peut supposer que $g$ transporte $s$ en un sommet $s'\in S^{st}(H'_{AD})$. Alors $Ad(g)$ définit un isomorphisme $Ad(g)_{s}:H_{AD,s}\to H'_{AD,s'}$. Soit ${\cal E}\in {\bf FC}^{st}_{{\mathbb F}_{q}}(\mathfrak{h}_{SC,s})$, posons ${\cal E}'=Ad(g)_{s}({\cal E})$. Fixons $N\in {\cal O}_{{\cal E}}^{\Gamma_{{\mathbb F}_{q}}}$ et $N'\in 
{\cal O}_{{\cal E}'}^{\Gamma_{{\mathbb F}_{q}}}$. Quitte à multiplier encore $g$ à gauche par un élément de $H'_{SC}(F^{nr})$, on peut supposer que $Ad(g)_{s}(N)=N'$. Munissons ${\cal E}'$ d'une action de Frobenius déterminée par un scalaire $r_{N'}\in {\mathbb C}^{\times}$. La construction de $\iota_{H',H}$ munit ${\cal E}$ d'une action de Frobenius déterminée par le scalaire $r_{N}=r_{N'}\epsilon_{N'}^{\flat}(Fr(g)g^{-1})$. La construction de $\iota_{\phi}$ munit ${\cal E}$ d'une action de Frobenius déterminée par le scalaire $\tilde{r}_{N}=r_{N'}\epsilon_{s',{\cal E}',N'}(k_{\phi})$. On note $f_{{\cal E}'}$, resp. $f_{{\cal E}}$, resp.  $\tilde{f}_{{\cal E}}$, les fonctions caractéristiques de ${\cal E}'$ muni de son action de Frobenius, resp. de ${\cal E}$ muni de l'action définie par $r_{N}$, resp. de ${\cal E}$ muni de l'action définie par $\tilde{r}_{N}$. On a par définition
$$ \iota_{H',H}(f_{{\cal E}})=\vert H'_{AD,s'}({\mathbb F}_{q})\vert^{-1}\vert H_{AD,s}({\mathbb F}_{q})\vert f_{{\cal E}'},$$
$$\iota_{\phi}(\tilde{f}_{{\cal E}})=\vert H'_{AD,s'}({\mathbb F}_{q})\vert^{-1}\vert H_{AD,s}({\mathbb F}_{q})\vert f_{{\cal E}'}.$$
L'égalité (1) résulte donc de l'égalité $\tilde{f}_{{\cal E}}=f_{{\cal E}}$, ou encore de l'égalité $r_{N} =\tilde{r}_{N}$, ou encore de l'égalité

(2) $\epsilon_{N'}^{\flat}(Fr(g)g^{-1})=\epsilon_{s',{\cal E}',N'}(k_{\phi})$. 

Puisque  $\phi$ se déduit de $Ad(g)$, on a $Fr(\phi)\phi^{-1}=Ad(k_{\phi})$, où $k_{\phi}$ est l'image naturelle dans $H'_{AD}$ de $Fr(g)g^{-1}$. L'élément $Fr(g)g^{-1}$ appartient à $Norm_{G_{SC}(F^{nr})}(H',s',N')\cap H'_{sc}(F^{nr})=K_{s',N'}^{H'_{sc},\dag,nr}$. Avec les notations de \ref{epsilonflatcasstable}, on a donc $k_{\phi}=\mu(Fr(g)g^{-1})$. Alors le théorème \ref{epsilonflatcasstable} implique (2). Cela achève la démonstration. $\square$

\subsubsection{Comparaison de deux éléments de $FC^{st}(\mathfrak{h}_{SC}(F))$}\label{HsHderechef}

On suppose que $G$ est quasi-déployé sur $F$. Soit $\Lambda$ un sous-ensemble propre de $\Delta_{a}^{nr}$ conservé par $\Gamma_{F}^{nr}$. On suppose que le sommet $s_{\Lambda}$ de $Imm(M_{\Lambda,AD})$ appartient à $S^{st}(M_{\Lambda,AD})$. Pour simplifier, on pose $(M,s)=(M_{\Lambda},s_{\Lambda})$. Soient $H\in {\cal L}_{F}^{nr,st}$ et $s_{H}$ un sommet de $Imm(H_{AD})$. On suppose que les couples $(M,s)$ et $(H,s_{H})$ sont conjugués par un élément de $G(F^{nr})$.
 Soit $x\in \bar{C}$. Supposons  que $x\in Imm^G(H)$ et $p_{H}(x)=s_{H}$. Le lemme \ref{HsH} dit que $\Lambda\subset \Lambda(x)$ et que les Levi $H_{s_{H}}$ et $M_{s}$ de $G_{x}$ sont conjugués. Soit ${\cal E}\in {\bf FC}^{st}_{{\mathbb F}_{q}}(\mathfrak{m}_{SC,s})$ muni d'une action de Frobenius. 
On note $f$ sa fonction caractéristique et ${\bf f}$ la fonction sur $\mathfrak{m}_{SC}(F)$ qui s'en déduit.  Puisque  les Levi $H_{s_{H}}$ et $M_{s}$ de $G_{x}$ sont conjugués, la construction de Lusztig rappelée en \ref{fonctionscaracteristiques} déduit de ${\cal E}$ un faisceau-caractère ${\cal E}_{H}$ sur $\mathfrak{h}_{SC,s_{H}}$ muni d'une action de Frobenius. On note 
$f_{H}$ sa fonction caractéristique et ${\bf f}_{H}$ la fonction sur $\mathfrak{h}_{SC}(F)$ qui s'en déduit. Les fonctions ${\bf f}$ et ${\bf f}_{H}$ appartiennent respectivement à $FC^{st}(\mathfrak{m}_{SC}(F))$ et $FC^{st}(\mathfrak{h}_{SC}(F))$. 

\begin{lem}{On a l'égalité $\iota_{H,M}({\bf f})=\vert H_{AD,s_{H}}({\mathbb F}_{q})\vert^{-1}\vert M_{AD,s}({\mathbb F}_{q})\vert {\bf f}_{H}$.}\end{lem}
  
Preuve. Puisque les couples $(H,s_{H})$ et $(M,s)$ sont conjugués par un élément de $K_{x}^{nr,0}$ d'après le lemme \ref{HsH}, le lemme est une conséquence immédiate des définitions et du théorème \ref{epsilonflatcasstable} appliqué à un sommet de $Imm_{F^{nr}}(G)$ adhérent à la facette à laquelle appartient $x$. $\square$

 \section{Relation avec les faisceaux-caract\`eres unipotents sur un groupe fini}
  
 \subsection{Quelques constructions}
 
 \subsubsection{Le groupe $G$\label{legroupeG}}
 
 Dans cette section, notre objet d'\'etude est le groupe ${\bf G}$ introduit en \ref{discretisation}. Il est plus commode d'oublier la provenance $p$-adique de ce groupe et de le consid\'erer directement comme un groupe d\'efini sur ${\mathbb F}_{q}$. Toutefois, pour certains points,  il est utile de r\'eintroduire  de fa\c{c}on occulte le groupe $p$-adique initial. C'est la raison des constructions de ce premier paragraphe.

 Soit $G$ un groupe r\'eductif connexe d\'efini sur ${\mathbb F}_{q}$.  On impose l'hypoth\`ese $(Hyp)_{1}(p)$.

   Pour cette premi\`ere sous-section , {\bf on suppose que} $G$ {\bf est adjoint et absolument simple.}   On fixe une paire de Borel \'epingl\'ee $\mathfrak{E}=(B,T,(E_{\alpha})_{\alpha\in \Delta})$ d\'efinie sur ${\mathbb F}_{q}$. On compl\`ete l'\'epinglage en un \'epinglage affine $\mathfrak{E}_{a}=(T,(E_{\alpha})_{\alpha\in \Delta_{a}})$ lui-aussi d\'efini sur ${\mathbb F}_{q}$, c'est-\`a-dire que l'on fixe un \'el\'ement non nul $E_{\alpha_{0}}\in \mathfrak{u}_{\alpha_{0}}$ fix\'e par $\Gamma_{{\mathbb F}_{q}}$.  
 
 Introduisons le groupe r\'eductif connexe $G_{F}$ d\'efini et quasi-d\'eploy\'e sur $F$  qui a m\^eme donn\'ees de racines que $G$. Cela signifie que l'on peut fixer une paire de Borel \'epingl\'ee $\mathfrak{E}_{F}=(B_{F},T_{F},(E_{F,\alpha})_{\alpha\in \Delta})$ de $G_{F}$ d\'efinie sur $F$ et des isomorphismes duaux $X_{*}(T_{F})\to X_{*}(T)$ et $X^{*}(T_{F})\to X^*(T)$ qui transportent les racines et racines simples, resp. coracines et coracines simples, de $T_{F}$ sur celles de $T$. De plus, $G_{F}$ est d\'eploy\'e sur $F^{nr}$ et les isomorphismes ci-dessus sont \'equivariants pour les actions de $\Gamma_{F}^{nr}\simeq \Gamma_{{\mathbb F}_{q}}$. En \ref{discretisation}, on a associ\'e \`a ces donn\'ees un groupe ${\bf G}$ muni d'un \'epinglage et il est clair que ${\bf G}$ s'identifie \`a $G$, l'\'epinglage en question s'identifiant \`a $\mathfrak{E}$. L'espace $\mathfrak{g}$ s'identifie  \`a $\mathfrak{g}_{F,s_{0}}$. On peut relever $E_{\alpha_{0}}$ en un \'el\'ement  $E'_{F,\alpha_{0}}\in \mathfrak{k}_{F,s_{0}}^{nr}\cap \mathfrak{u}_{F,\alpha_{0}}(F)$. On fixe une uniformisante $\varpi_{F}$ de $F$ et on pose $E_{F,\alpha_{0}}=\varpi_{F}E'_{F,\alpha_{0}}$. Alors $\mathfrak{E}_{F,a}=(T_{F},(E_{F,\alpha})_{\alpha\in \Delta_{a}})$ est un \'epinglage affine de $G_{F}$ d\'efini sur $F$. 
 
 Cette construction nous permet d'utiliser les d\'efinitions de \ref{racinesaffines} et \ref{discretisation}. En particulier, on dispose de l'alc\^ove $C^{nr}\subset App_{F^{nr}}(T_{F})\subset Imm_{F^{nr}}(G_{F})$, du groupe ${\mathbb T}$, de son plongement $a_{T}:{\mathbb T}\to App_{F^{nr}}(T_{F})$ et de l'homomorphisme surjectif $j_{T}:{\mathbb T}\to T$. Puisque $G_{F}$ est adjoint, le plongement $a_{T}$ est injectif et, comme en \ref{dansboldsymbolg}, on note ${\bf j}_{T}:a_{T}({\mathbb T})\to T$ l'application $j_{T}\circ a_{T}^{-1}$.  
 Les groupes $G$ et $G_{F}$ ont m\^eme groupe de Weyl $W$. Notons $\tilde{W}$ le groupe de transformations de $App_{F^{nr}}(T_{F})$ engendr\'e par $W$ et le groupe des translations par $X_{*}(T)$.  Il respecte le sous-ensemble $a_{T}({\mathbb T})$ et, pour tout $x\in a_{T}({\mathbb T})$, il existe $\tilde{w}\in \tilde{W}$ tel que $\tilde{w}(x)\in a_{T}({\mathbb T})\cap \bar{C}^{nr}$. On fait agir $\tilde{W}$ sur $T$: l'action de $W$ est l'action naturelle et les translations par $X_{*}(T)$ agissent trivialement. L'homomorphisme ${\bf j}_{T}$ est \'equivariant pour les actions de $\tilde{W}$. Donc
 
 (1) pour tout $t\in T$, il existe $x\in \bar{C}^{nr}\cap a_{T}({\mathbb T})$ et $w\in W$ de sorte que ${\bf j}_{T}(x)=w(t)$. 
 
 Soit $x\in \bar{C}^{nr}\cap a_{T}({\mathbb T})$. On note $\Lambda(x)$ l'ensemble des $\alpha\in \Delta_{a}$ telles que $\alpha^{aff}(x)=0$. Posons $t={\bf j}_{T}(x)$. Montrons que 
 
 (2) le syst\`eme de racines de $G_{t}$ relatif au tore $T$ a pour base $\Lambda(x)$.
 
 En effet, notons ${\cal F}$ la facette de $Imm_{F^{nr}}(G_{F})$ \`a laquelle appartient $x$. Le syst\`eme de racines du groupe $G_{F,{\cal F}}$ a pour base $\Lambda(x)$. On a introduit en \ref{lesgroupesGF} un plongement $\iota_{x,{\cal F}}:G_{F,{\cal F}}\to G$. Son image est justement le groupe $G_{t}$. Le plongement \'etant d\'eduit d'une conjugaison par un \'el\'ement de $T_{F}$, il  identifie les ensembles de racines de ces groupes. L'assertion (2) s'ensuit.
 
 Pour $\alpha\in \Delta_{a}$, le sommet $s_{\alpha}$ associ\'e \`a $\alpha$ appartient \`a $a_{T}({\mathbb T})$. On pose $t_{\alpha}={\bf j}_{T}(s_{\alpha})$. Introduisons la base $(\check{\varpi}_{\alpha})_{\alpha\in \Delta}$ de $X_{*}(T)$ duale de $\Delta$. On se rappelle la relation
 $$(3) \qquad \sum_{\alpha\in \Delta_{a}}d(\alpha)\alpha=0$$
 o\`u $d(\alpha_{0})=1$. On calcule ais\'ement $t_{\alpha_{0}}=1$ et $t_{\alpha}=\check{\varpi}_{\alpha}(\bar{\zeta}_{1/d(\alpha)})$ pour $\alpha\in \Delta$. 
 
 \subsubsection{Conjugaison des \'el\'ements $t_{\alpha}$\label{talpha}}
 
   Soient $\alpha,\alpha'\in \Delta_{a}$. Notons $Aut(G,t_{\alpha}\mapsto t_{\alpha'}) $ l'ensemble des automorphismes de $G$ qui envoient $t_{\alpha}$ sur $t_{\alpha'}$ (cet ensemble peut \^etre vide). Posons $Int(G,t_{\alpha}\mapsto t_{\alpha'})=\{g\in G; Ad(g)\in 
     Aut(G,t_{\alpha}\mapsto t_{\alpha'})\}$.  On a introduit un sous-groupe $\boldsymbol{\Omega}$ de $G$ en \ref{epinglages}. Notons $\boldsymbol{\Omega}(\alpha\mapsto \alpha')$ l'ensemble des \'el\'ements $\omega\in \boldsymbol{\Omega}$ tels que $\omega(\alpha)=\alpha'$. Si $G$ n'est pas de type $A_{n-1}$ avec $n$ impair, on a introduit un groupe ${\bf Aut}({\cal D}_{a})$ en \ref{epinglages}. Notons   ${\bf Aut}({\cal D}_{a},\alpha\mapsto \alpha')$ l'ensemble des $\delta\in {\bf Aut}({\cal D}_{a})$ tels que $\delta(\alpha)=\alpha'$.  Dans le cas o\`u $\alpha=\alpha'$, on a simplement $Int(G,t_{\alpha}\mapsto t_{\alpha})=Z_{G}(t_{\alpha})$ et on 
    simplifie les autres notations en $Aut(G,t_{\alpha})$, $\boldsymbol{\Omega}(\alpha)$ et ${\bf Aut}({\cal D}_{a},\alpha)$.

   \begin{lem}{Soient $\alpha,\alpha'\in \Delta_{a}$.

   (i) Soit $\delta\in Aut(G,t_{\alpha}\mapsto t_{\alpha'})$. Supposons que $\delta$ conserve $T$ et  envoie l'ensemble de racines $\Delta_{a}-\{\alpha\}$ sur $\Delta_{a}-\{\alpha'\}$. Alors $\delta$ envoie  $\alpha$ sur $\alpha'$. 
   
   (ii)  L'application  produit $Z_{G}(t_{\alpha'})^0\times \boldsymbol{\Omega}(\alpha\mapsto \alpha')\to G$ est injective et son image  est $Int(G,t_{\alpha}\mapsto t_{\alpha'})$.

    (iii) Supposons que  $G$ n'est pas de type $A_{n-1}$ avec $n$ impair.  L'application produit $Ad(Z_{G}(t_{\alpha'})^0)\times {\bf Aut}({\cal D}_{a},\alpha\mapsto \alpha')\to Aut(G)$ est injective et son image est $Aut(G,t_{\alpha}\mapsto t_{\alpha'})$. 
     }\end{lem}
    
  Preuve.         
     Soit $\delta$ v\'erifiant les conditions de (i). Puisque $t_{\alpha}$, resp. $t_{\alpha'}$,  est un \'el\'ement de $G $ d'ordre $d(\alpha)$, resp. $d(\alpha')$, l'\'egalit\'e $\delta(t_{\alpha})=t_{\alpha'}$ implique $d(\alpha)=d(\alpha')$. Supposons d'abord $d(\alpha')=d(\alpha)=1$.   Il existe alors $\omega,\omega'\in \boldsymbol{\Omega}$ tels que $\omega(\alpha_{0})=\alpha$ et $\omega'(\alpha')=\alpha_{0}$, cf. \ref{epinglages}(2). L'automorphisme $\omega'\delta\omega$ conserve $T$ et l'ensemble $\Delta_{a}-\{\alpha_{0}\}=\Delta$. Donc il  conserve le sous-ensemble $\Sigma^+$ des racines positives, donc conserve la plus grande racine $-\alpha_{0}$. Cela entra\^{\i}ne que $\delta$ envoie $\alpha$ sur $\alpha'$. Supposons maintenant que $d(\alpha')=d(\alpha)\not=1$, a fortiori $\alpha,\alpha'\not=\alpha_{0}$. Appliquons $\delta$ \`a la relation \ref{legroupeG}(3). On obtient
     $$d(\alpha)\delta(\alpha)+\sum_{\gamma\in \Delta_{a}-\{\alpha\}}d(\gamma)\delta(\gamma)=0,$$
     ou encore
     $$\delta(\alpha)=-\sum_{\gamma\in \Delta_{a}-\{\alpha'\}}d(\alpha)^{-1}d(\delta^{-1}(\gamma))\gamma.$$
     Dans le deuxi\`eme membre, on exprime la racine $\alpha_{0}\in \Delta_{a}-\{\alpha'\}$ gr\^ace \`a \ref{legroupeG}(3). En tenant compte de l'\'egalit\'e $d(\alpha)=d(\alpha')$, on obtient
     $$(1) \qquad \delta(\alpha)=d(\delta^{-1}(\alpha_{0}))\alpha'+\sum_{\gamma\in \Delta-\{\alpha'\}}d(\alpha)^{-1}(d(\delta^{-1}(\alpha_{0}))d(\gamma)-d(\delta^{-1}(\gamma)))\gamma.$$
     C'est l'\'ecriture de la racine $\delta(\alpha)$ dans la base $\Delta$. Donc les coefficients du membre de droite sont entiers, tous de m\^eme signe, et la valeur absolue du  coefficient d'une racine $\gamma\in \Delta$ est major\'ee par $d(\gamma)$. En appliquant cela au coefficient de $\alpha'$, on obtient  $d(\delta^{-1}(\alpha_{0}))\leq d(\alpha')$. 
     Puisque $\delta$ transporte $(\alpha,t_{\alpha})$ sur $(\delta(\alpha),t_{\alpha'})$, on a $\delta(\alpha)(t_{\alpha'})=\alpha(t_{\alpha})=\bar{\zeta}_{1/d(\alpha)}=\bar{\zeta}_{1/d(\alpha')}$.  D'apr\`es (1), on a $<\delta(\alpha),\check{\varpi}_{\alpha'}>=d(\delta^{-1}(\alpha_{0}))$.  De la d\'efinition de $t_{\alpha'}$ r\'esulte l'\'egalit\'e $\delta(\alpha)(t_{\alpha'})=\bar{\zeta}_{d(\delta^{-1}(\alpha_{0}))/d(\alpha')}$. En comparant les deux \'egalit\'es pr\'ec\'edentes, on obtient $d(\delta^{-1}(\alpha_{0}))\equiv 1\,\,mod\,\,d(\alpha'){\mathbb Z}$. Avec l'in\'egalit\'e $d(\delta^{-1}(\alpha_{0}))\leq d(\alpha')$, cela entra\^{\i}ne $d(\delta^{-1}(\alpha_{0}))=1$. Pour $\gamma\in \Delta-\{\alpha'\}$, le coefficient de $\gamma$ dans l'expression (1) est donc $d(\alpha)^{-1}(d(\gamma)-d(\delta^{-1}(\gamma)))$. Montrons que
     
     (2) la somme de ces coefficients est nulle. 
     
     Parce que $d(\delta^{-1}(\alpha_{0}))=1=d(\alpha_{0})$, on a 
     $$\sum_{\gamma\in \Delta-\{\alpha'\}}(d(\gamma)-d(\delta^{-1}(\gamma)))=\sum_{\gamma\in \Delta_{a}-\{\alpha'\}}(d(\gamma)-d(\delta^{-1}(\gamma)))$$
     $$=(\sum_{\gamma\in \Delta_{a}-\{\alpha'\}}d(\gamma))-(\sum_{\gamma\in \Delta_{a}-\{\alpha'\}}d(\delta^{-1}(\gamma)))$$
     $$(\sum_{\gamma\in \Delta_{a}-\{\alpha'\}}d(\gamma))-(\sum_{\gamma\in \Delta_{a}-\{\alpha\}}d(\gamma))=
      (\sum_{\gamma\in \Delta_{a}}d(\gamma))-d(\alpha')- (\sum_{\gamma\in \Delta_{a}}d(\gamma))+d(\alpha).$$
      Ceci est nul puisque $d(\alpha)=d(\alpha')$. Cela prouve (2). 
      
      Puisque tous les coefficients sont de m\^eme signe, (2) entra\^{\i}ne qu'ils sont nuls. Alors (1) se simplifie en $\delta(\alpha)=\alpha'$, ce qui prouve l'assertion (i) de l'\'enonc\'e. 
          
            Soit $\omega\in \boldsymbol{\Omega}(\alpha\mapsto \alpha')$. Alors $ \omega$ conserve $T$ et  $\Delta_{a}$ et envoie $\alpha$ sur $\alpha'$ donc aussi $\Delta_{a}-\{\alpha\}$ sur $\Delta_{a}-\{\alpha'\}$. D'apr\`es \ref{epinglages} (3), on a $d(\alpha)=d(\alpha')$. On voit alors que  $Ad(\omega)(t_{\alpha})$  v\'erifie les relations: $\alpha'(Ad(\omega)(t_{\alpha}))=\bar{\zeta}_{1/d(\alpha')}$ et $\gamma(Ad(\omega)(t_{\alpha}))=1$ pour $\gamma\in \Delta_{a}-\{\alpha'\}$. Ces relations caract\'erisent $t_{\alpha'}$ donc 
            $Ad(\omega)(t_{\alpha})=t_{\alpha'}$. 
         Cela prouve  que $\omega\in Int(G,t_{\alpha}\mapsto t_{\alpha'})$. Il en r\'esulte que l'image de l'application de (ii) est contenue dans cet ensemble. Inversement, soit $g \in Int(G,t_{\alpha}\mapsto t_{\alpha'})$. Alors $Ad(g)$ envoie $Z_{G}(t_{\alpha})^0$ sur $Z_{G}(t_{\alpha'})^0$. Cet automorphisme transporte $T$ en un sous-tore maximal $T'$ de $Z_{G}(t_{\alpha'})^0$ et transporte la base $\Delta_{a}-\{\alpha\}$ du syst\`eme de racines de $Z_{G}(t_{\alpha})^0$ relatif \`a $T$ en une base $\Delta(\alpha')$ du syst\`eme de racines de $Z_{G}(t_{\alpha'})^0$ relatif \`a $T'$. D'apr\`es \ref{legroupeG}(2), il existe $h\in Z_{G}(t_{\alpha'})^0$ tel que $Ad(h)$ transporte $T$  en $T'$ et  $\Delta_{a}-\{\alpha'\}$ en $\Delta(\alpha')$. Alors $Ad(h^{-1}g)$ conserve $T$ et envoie $\Delta_{a}-\{\alpha\}$ sur $\Delta_{a}-\{\alpha'\}$. Cet automorphisme  envoie $t_{\alpha}$ sur $t_{\alpha'}$ d'apr\`es l'hypoth\`ese sur $g$. D'apr\`es (i), il envoie $\alpha$ sur $\alpha'$ . Son image dans le groupe de Weyl appartient donc \`a $\Omega$. Introduisons l'\'el\'ement $\omega\in \boldsymbol{\Omega}$ qui a la m\^eme image dans le groupe de Weyl. On a $\omega\in \boldsymbol{\Omega}(\alpha\mapsto \alpha')$.  Posons $t=h^{-1}g\omega^{-1}$. Alors $Ad(t)$ conserve $T$ et agit trivialement sur $\Delta_{a}$. A fortiori, il agit trivialement sur $\Delta$ donc conserve $B$. Alors $t\in T\subset Z_{G}(t_{\alpha'})^0$. On a $g=ht\omega \in Z_{G}(t_{\alpha'})^0\boldsymbol{\Omega}(\alpha\mapsto \alpha')$, ce qui d\'emontre la premi\`ere assertion de (ii). Soient $z_{1},z_{2}\in Z_{G}(t_{\alpha'})^0$ et $\omega_{1},\omega_{2}\in \boldsymbol{\Omega}(\alpha\mapsto \alpha')$. Supposons $z_{1}\omega_{1}=z_{2}\omega_{2}$. Posons $g=z_{2}^{-1}z_{1}=\omega_{2}\omega_{1}^{-1}$. Alors $g\in Z_{G}(t_{\alpha'})^0\cap \boldsymbol{\Omega}(\alpha')$.  L'automorphisme $Ad(g)$ conserve $T$, $\Delta_{a}$ et $\alpha'$. Donc cet automorphisme conserve $\Delta_{a}-\{\alpha'\}$. Puisque $g\in Z_{G}(t_{\alpha'})^0$, $Ad(g)$ se restreint \`a $Z_{G}(t_{\alpha'})^0$ en un automorphisme int\'erieur. Puisque  $\Delta_{a}-\{\alpha'\}$ est une base du syst\`eme de racines de ce groupe, $Ad(g)$ agit par l'identit\'e sur $\Delta_{a}-\{\alpha'\}$, donc aussi sur $\Delta_{a}$. Puisque $g\in \boldsymbol{\Omega}$, cela entra\^{\i}ne $g=1$. Cela d\'emontre (ii). 
     
     La preuve de (iii) est similaire. $\square$

     Dans le cas o\`u $\alpha=\alpha'$, les assertions (ii) et (iii) deviennent

        (3) on a les \'egalit\'es $Z_{G}(t_{\alpha})=Z_{G}(t_{\alpha})^0\boldsymbol{\Omega}(\alpha)$ et $Z_{G}(t_{\alpha})^0\cap \boldsymbol{\Omega}(\alpha)=\{1\}$;
   
   (4) si $G$ n'est pas de type $A_{n-1}$ avec $n$ impair, on a les \'egalit\'es $Aut(G,t_{\alpha})=Ad(Z_{G}(t_{\alpha})^0) {\bf Aut}({\cal D}_{a},\alpha)$ et $Ad(Z_{G}(t_{\alpha})^0) \cap {\bf Aut}({\cal D}_{a},\alpha)=\{1\}$.
   
     \subsubsection{Construction de Levi\label{constructiondeLeviG}}
    Soit $\Lambda$ un sous-ensemble de $\Delta_{a}$ tel que $\Lambda\not=\Delta_{a}$. Notons $A_{\Lambda}$ la composante neutre du sous-groupe des $t\in T$ tels que $\alpha(t)=1$ pour tout $\alpha\in \Lambda$. C'est un sous-tore de $T$. Le groupe $X_{*}(A_{\Lambda})$ est celui des $x_{*}\in X_{*}(T)$ tels que $<\alpha,x_{*}>=0$ pour tout $\alpha\in \Lambda$. Notons $M_{\Lambda}$ le commutant de $A_{\Lambda}$. C'est un Levi de $G$ dont le centre est \'egal \`a $A_{\Lambda}$. Notons ${\mathbb Z}[\Lambda]$ le sous-${\mathbb Z}$-module de $X^*(T)$ engendr\'e par $\Lambda$ et notons ${\mathbb Q}[\Lambda]$ le sous-${\mathbb Q}$-espace vectoriel de  $X^*_{{\mathbb Q}}(T)$ engendr\'e par $\Lambda$.  L'ensemble $\Sigma^{M_{\Lambda}}$ des racines de $T$ dans $M_{\Lambda}$ est $\Sigma\cap {\mathbb Q}[\Lambda]$.  Notons aussi $L_{\Lambda}$ le sous-groupe connexe de $G$ contenant $T$ tel que l'ensemble de racines $\Sigma^{L_{\Lambda}}$ de $T$ dans $L_{\Lambda}$ soit $\Sigma\cap {\mathbb Z}[\Lambda]$. L'ensemble  $\Lambda$ est une base de $\Sigma^{L_{\Lambda}}$. Le groupe $L_{\Lambda}$ est un sous-groupe de $M_{\Lambda}$ de m\^eme rang semi-simple $\vert \Lambda\vert $. La composante neutre du centre de $L_{\Lambda}$ est $A_{\Lambda}$.  
    
      \begin{lem}{ Pour tout $x\in \bar{C}^{nr}$ tel que $\Lambda\subset \Lambda(x)$, 
   on a l'\'egalit\'e $L_{\Lambda}=Z_{G}({\bf j}_{T}(x))^0\cap M_{\Lambda}$. 
  Il existe un unique \'el\'ement $t_{\Lambda}\in M_{\Lambda,AD}$ tel que, pour tout $x$ comme ci-dessus,  la projection de ${\bf j}_{T}(x)$ dans $M_{\Lambda,AD}$ soit $t_{\Lambda}$.   On a $L_{\Lambda}=Z_{M_{\Lambda}}(t_{\Lambda})^0$.}\end{lem}
  
  La preuve est similaire \`a celle de \ref{constructiondeLevi}(2). $\square$

 \subsection{Faisceaux-caract\`eres unipotents}

  \subsubsection{Rappels sur la repr\'esentation de Springer g\'en\'eralis\'ee\label{Springer}}
  Dans cette section, {\bf on l\`eve l'hypoth\`ese que }$G$ {\bf est adjoint et absolument simple. }
   
 On note ${\cal I}(G)$ l'ensemble des couples $({\cal O},{\cal L})$, o\`u ${\cal O}\in \mathfrak{g}_{nil}/conj$ et ${\cal L}$ est un syst\`eme local irr\'eductible sur ${\cal O}$.    Il s'identifie \`a celui des classes de conjugaison par $G$ de couples $(N,\xi)$, o\`u $N\in \mathfrak{g}_{nil}$ et $\xi$ est une repr\'esentation irr\'eductible de $Z_{G}(N)/Z_{G}(N)^0$ dans un $\bar{{\mathbb Q}}_{l}$-espace.
 
 {\bf Remarque.} Comme en  \ref{fonctionscaracteristiques}, on suppose fixé un isomorphisme $\bar{{\mathbb Q}}_{l}\simeq {\mathbb C}$. Toutes les fonctions caractéristiques que l'on construira seront identifiées à des fonctions à valeurs complexes.
 \bigskip
 
 Notons ${\cal I}_{0}^{FC}(G)$ l'ensemble des classes de conjugaison par $G$ de couples $(M,{\cal E})$, o\`u 
 $M$ est un  Levi de $G$ et ${\cal E}\in {\bf FC}^M(\mathfrak{m}_{SC})$.   Le groupe $G$ agit par conjugaison sur l'ensemble des triplets $(M,{\cal E},\rho)$ o\`u $M$ et ${\cal E}$ sont comme ci-dessus et $\rho\in Irr(W(M))$. On note ${\cal I}^{FC}(G)$ l'ensemble des classes de conjugaison. La repr\'esentation de Springer g\'en\'eralis\'ee \'etablit une bijection
 $$Spr:{\cal I}^{FC}(G)\to {\cal I}(G).$$

  Le groupe $\Gamma_{{\mathbb F}_{q}}$ agit naturellement sur les ensembles  ${\cal I}(G)$, ${\cal I}_{0}^{FC}(G)$ et ${\cal I}^{FC}(G)$. On note ${\cal I}_{{\mathbb F}_{q}}(G)$, ${\cal I}_{0,{\mathbb F}_{q}}^{FC}(G)$ et ${\cal I}^{FC}_{{\mathbb F}_{q}}(G)$ ces ensembles de points fixes.  Fixons une paire de Borel $(B,T)$ d\'efinie sur ${\mathbb F}_{q}$.  Un \'el\'ement de ${\cal I}_{0,{\mathbb F}_{q}}^{FC}(G)$ est repr\'esent\'e par un couple $(M,{\cal E})$ tel que $M$ soit un ${\mathbb F}_{q}$-Levi  standard de $G$ et que ${\cal E}$ soit conserv\'e par l'action de Frobenius.   D'apr\`es \ref{faisceauxcaracteres}(1), un ${\mathbb F}_{q}$-Levi standard $M$ tel que ${\bf FC}^M(\mathfrak{m}_{SC})$ est non vide n'est conjugu\'e \`a aucun autre ${\mathbb F}_{q}$-Levi standard que lui-m\^eme et  tout \'el\'ement de $Norm_{G}(M)$  agit par l'identité sur l'ensemble ${\bf FC}^M(\mathfrak{m}_{SC})$. Il en r\'esulte  que deux couples  $(M,{\cal E})$ v\'erifiant les conditions ci-dessus ne sont conjugu\'es que s'ils sont \'egaux. Donc ${\cal I}_{0,{\mathbb F}_{q}}^{FC}(G)$ s'identifie \`a l'ensemble de ces couples $(M,{\cal E})$.  Pour un tel couple, on fixe une action de Frobenius sur ${\cal E}$. L'ensemble ${\cal I}_{{\mathbb F}_{q}}^{FC}(G)$ est celui des triplets $(M,{\cal E},\rho)$ tels que $(M,{\cal E})\in {\cal I}_{0,{\mathbb F}_{q}}^{FC}(G)$ et que la classe de $\rho$ soit conserv\'ee par $\Gamma_{{\mathbb F}_{q}}$. Comme en \ref{fonctionscaracteristiques}, on note $Irr_{{\mathbb F}_{q}}(W(M))$ cet ensemble de repr\'esentations et, pour $\rho 
   \in Irr_{{\mathbb F}_{q}}(W(M))$, on fixe un prolongement continu $\rho^{\flat}$ de $\rho$ au  produit semi-direct $W(M)\rtimes \Gamma_{{\mathbb F}_{q}}$.   On a d\'efini en loc. cit. la fonction caract\'eristique $\boldsymbol{\chi}_{M,{\cal E},\rho}$. 
  La bijection $Spr$ est \'equivariante pour les actions galoisiennes. En particulier, elle se restreint en une bijection  de ${\cal I}^{FC}_{{\mathbb F}_{q}}(G)$ sur ${\cal I}_{{\mathbb F}_{q}}(G)$.

 Soit $(M,{\cal E},\rho)\in {\cal I}_{{\mathbb F}_{q}}^{FC}(G)$.  Notons $({\cal O},{\cal L})$  son image par la repr\'esentation de Springer g\'en\'eralis\'ee.  Elle est conserv\'ee par l'action galoisienne.  Munissons ${\cal L}$ d'une  action de Frobenius  dont l'action sur chaque fibre d\'efinie sur ${\mathbb F}_{q}$ est d'ordre fini. Introduisons le prolongement d'intersection   du syst\`eme local ${\cal L}[dim({\cal O})]$ sur ${\cal O}$. Il est encore muni d'une action de Frobenius,  notons $\chi_{M,{\cal E},\rho}$ sa fonction caract\'eristique. 
Posons 
 $$b(M,{\cal E},\rho)=\frac{1}{2}(dim({\cal O}_{{\cal E}})-dim({\cal O})+dim(\mathfrak{g})-dim(\mathfrak{m})).$$
  D'apr\`es \cite{L4} 6.3, il  existe alors un unique choix de l'action de Frobenius sur ${\cal L}$ de sorte que l'on ait l'\'egalit\'e
 
  (1) $\chi_{M,{\cal E},\rho}(X)=(-1)^{dim(Z(M)^0)}q^{-b(M,{\cal E},\rho)} \boldsymbol{\chi}_{M,{\cal E},\rho}(X)$ pour tout $X\in \mathfrak{g}_{nil}({\mathbb F}_{q})$.

On fixe cette normalisation. De la formule ci-dessus et de \ref{fonctionscaracteristiques}(1) et (2), on d\'eduit

(2) $\chi_{M,{\cal E},\rho}=q^{-b(M,{\cal E},\rho)}\vert W(M)\vert ^{-1}\sum_{w\in W(M)}trace(\rho^{\flat}(w^{-1}Fr))Q_{M,{\cal E},w}$.

On note ${\cal Y}_{{\cal O},{\cal L}}$ la restriction \`a   ${\cal O}$ de $\chi_{M,{\cal E},\rho}$. C'est la fonction caract\'eristique de ${\cal L}$. 
On note aussi $\chi^{\natural}_{ M,{\cal E},\rho}$ la fonction caract\'eristique de l'inverse de l'action de Frobenius du prolongement d'intersection   du syst\`eme local ${\cal L}[dim({\cal O})]$ sur ${\cal O}$. Sa restriction \`a ${\cal O}$ est la conjugu\'ee $\bar{{\cal Y}}_{{\cal O},{\cal L}}$. Pour $w\in W(M)$, il existe une unique fonction $Q^{\natural}_{M,{\cal E},w}$ sur $\mathfrak{g}_{nil}({\mathbb F}_{q})$ de sorte que, pour tout $\rho\in Irr_{{\mathbb F}_{q}}(W(M))$, on ait l'\'egalit\'e

(3) $\chi^{\natural}_{M,{\cal E},\rho}=q^{b(M,{\cal E},\rho)}\vert W(M)\vert ^{-1}\sum_{w\in W(M)}trace(\rho^{\flat}(w Fr^{-1}))
Q^{\natural}_{M,{\cal E},w}$. 

  \subsubsection{Faisceaux-caract\`eres unipotents cuspidaux sur $G$\label{faisceauxcaracteresunipotentscuspidaux}}
         
   Lusztig a d\'efini la notion de faisceau-caract\`ere unipotent cuspidal sur $G$. Signalons que "unipotent" se r\'ef\`ere \`a la classification de ces faisceaux en termes de groupe dual de $G$, cela ne signifie nullement qu'ils sont \`a support unipotent.    Un tel faisceau ${\cal E}$ est  port\'e par une classe de conjugaison ${\cal X}_{{\cal E}}$ dans $G$. Fixons un \'el\'ement de cette classe, que l'on \'ecrit $tu$, $t$ \'etant sa composante semi-simple et $u$ sa composante unipotente. Il y a alors un caract\`ere $\xi_{{\cal E},tu}$, ou simplement $\xi_{tu}$, de $Z_{G}(tu)/Z_{G}(tu)^0$ tel que ${\cal E}$ se restreigne \`a ${\cal X}$ en le syst\`eme local $G/Z_{G}(tu)^0\times_{\xi_{tu}}\bar{{\mathbb Q}}_{l}$ d\'efini comme en \ref{faisceauxcaracterescuspidaux}. On a l'\'egalit\'e $Z_{G_{t}}(u)^0=Z_{G}(tu)^0$ d'o\`u  l'inclusion $Z_{G_{t}}(u)/Z_{G_{t}}(u)^0\subset Z_{G}(tu)/Z_{G}(tu)^0$. La classe de conjugaison ${\cal X}$ peut porter plusieurs faisceaux-caract\`eres unipotents et cuspidaux, correspondant \`a plusieurs caract\`eres $\xi_{tu}$. Deux tels caract\`eres distincts ont des restrictions distinctes \`a $Z_{G_{t}}(u)/Z_{G_{t}}(u)^0$. On a la propri\'et\'e suivante, cf. \cite{L2} proposition 2.7 et \cite{L3} proposition 7.11:
   
   (1)   soit $t$ la partie semi-simple d'un \'el\'ement de ${\cal X}_{{\cal E}}$; alors l'ensemble des \'el\'ements unipotents $u\in G_{t}$ tels que $tu\in {\cal X}_{{\cal E}}$ est une unique classe de conjugaison par $G_{t}$;  $t$ est \'el\'ement isol\'e de $G$, c'est-\`a-dire que la composante neutre de $Z(G_{t})$ est $Z(G)^0$.

    On note ${\bf FC}_{u}(G)$ l'ensemble des faisceaux-caract\`eres unipotents cuspidaux sur $G$. Notons $\underline{\pi}:G\to G_{AD}$ l'homomorphisme naturel. L'application image r\'eciproque par $\underline{\pi}$ identifie ${\bf FC}_{u}(G_{AD})$ et ${\bf FC}_{u}(G)$. En particulier les caract\`eres $\xi_{tu}$ ci-dessus sont triviaux sur $Z(G)$. 
   
{\bf Supposons} $G$ {\bf adjoint et absolument simple} et introduisons les objets de \ref{legroupeG}.  L'ensemble ${\bf FC}_{u}(G)$ est partiellement d\'ecrit dans \cite{DLM} p. 495 \`a 509. On voit que,  pour tout ${\cal E}\in {\bf FC}_{u}(G)$, il existe une  racine $\alpha\in \Delta_{a}$ telle que ${\cal X}_{{\cal E}}$ contienne un \'el\'ement de partie semi-simple $t_{\alpha}$. De plus,  cette racine $\alpha$ est fix\'ee par tout \'el\'ement de $Aut({\cal D}_{a})$. Il r\'esulte du lemme \ref{talpha} que cette racine $\alpha$ est unique: le lemme \ref{talpha} implique que deux \'el\'ements $t_{\alpha}$ et $t_{\alpha'}$ sont conjugu\'es dans $G$ si et seulement si $\alpha$ et $\alpha'$ sont conjugu\'es par l'action du groupe $\boldsymbol{\Omega}$. On dira que $\alpha$ est associ\'ee \`a ${\cal E}$. 

Soit ${\cal E}\in {\bf FC}_{u}(G)$ et soit $\alpha\in \Delta_{a}$ sa racine associ\'ee. Notons ${\cal O}_{{\cal E}}$ l'ensemble des $N\in \mathfrak{g}_{t_{\alpha},nil}$ tel que $t_{\alpha}exp(N)\in {\cal X}_{{\cal E}}$. D'apr\`es (1), c'est une unique orbite pour l'action de $G_{t_{\alpha}}$. Pour $N\in {\cal O}_{{\cal E}}$, le groupe $Z_{G_{t_{\alpha}}}(N)/Z_{G_{t_{\alpha}}}(N)^0$ s'identifie \`a $Z_{G_{t_{\alpha}}}(exp(N))/Z_{G_{t_{\alpha}}}(exp(N))^0$. Ainsi, le caract\`ere $\xi_{t_{\alpha}exp(N)}$ se restreint en un caract\`ere de $Z_{G_{t_{\alpha}}}(N)/Z_{G_{t_{\alpha}}}(N)^0$ qui d\'efinit un syst\`eme local ${\cal E}_{t_{\alpha}}$ sur ${\cal O}_{{\cal E}}$. Comme on l'a dit dans la preuve de \ref{legroupeG}(2), le groupe $G_{t_{\alpha}}$ s'identifie \`a $G_{F, s_{\alpha}}$. On peut identifier ${\cal E}_{t_{\alpha}}$ \`a un syst\`eme local sur $\mathfrak{g}_{F,s_{\alpha},nil}$. En comparant la classification de  \cite{W7} paragraphe 9 avec les r\'esultats de \cite{DLM}, on obtient les r\'esultats suivants:

(2) pour tout ${\cal E}\in {\bf FC}_{u}(G)$, le syst\`eme local ${\cal E}_{t_{\alpha}}$ appartient \`a ${\bf FC}^{st}(\mathfrak{g}_{F, s_{\alpha}})$, o\`u $\alpha\in \Delta_{a}$ est la racine associ\'ee \`a ${\cal E}$;

(3) l'application ${\cal E}\mapsto (s_{\alpha},{\cal E}_{t_{\alpha}})$ est une bijection de ${\bf FC}_{u}(G)$ sur l'ensemble des couples $(s,\tilde{{\cal E}})$ o\`u $s\in S^{nr,st}(\bar{C}^{nr}) $ et $\tilde{{\cal E}}\in {\bf FC}^{st}(\mathfrak{g}_{F, s})$. 

Tout \'el\'ement de ${\bf FC}_{u}(G)$ est conserv\'e par l'action galoisienne.  La correspondance (3) n'est pas \'equivariante pour les actions galoisiennes (parce qu'en g\'en\'eral,  $t_{\alpha}$ n'est pas fix\'e par l'action galoisienne).

   Pour des raisons de r\'ecurrence, on doit consid\'erer le groupe $G=\{1\}$. Dans ce cas ${\bf FC}_{u}(G)$ est r\'eduit au syst\`eme local trivial de rang $1$ sur $G$. 
  
  \subsubsection{Faisceaux-caract\`eres unipotents sur $G$\label{faisceauxcaracteresunipotents}}
   Soit $M$ un Levi de $G$ et ${\cal E}\in {\bf FC}_{u}(M_{AD})$. On simplifie les notations en supprimant les termes ${\cal E}$ des notations du paragraphe pr\'ec\'edent: ${\cal X}={\cal X}_{{\cal E}}$ etc... Fixons $x=tu\in {\cal X}\subset M_{AD}$. Notons $Z_{M}(x)$ l'image r\'eciproque dans $M$ de $Z_{M_{AD}}(x)$. On a la suite exacte
  $$1\to Z(M)\to Z_{M}(x)\to Z_{M_{AD}}(x)\to 1.$$
  On rel\`eve le caract\`ere $\xi_{x}$ en un caract\`ere de $Z_{M}(x)$ trivial sur $Z(M)$. Notons $S_{x}$ l'image r\'eciproque de $x$ dans $M$ et $S_{x,reg}$ le sous-ensemble des \'el\'ements de $S_{x}$ dont le commutant dans $G$ est contenu dans $M$.  
     
   Posons ${\cal V}_{reg}=\{(y,gZ_{M}(x))\in  G\times G/Z_{M}(x); g^{-1}yg\in S_{x,reg}\}$. On note $p_{reg}:{\cal V}_{reg}\to G$ la projection $(y,gZ_{M}(x))\mapsto y$, on note $\nabla_{reg} $  son image et $\nabla$ l'adh\'erence de $\nabla_{reg}$. Posons $\tilde{{\cal V}}_{reg}= \{(y,gZ_{M}(x)^0)\in  G\times G/Z_{M}(x)^0; g^{-1}yg\in S_{x,reg}\}$. Le groupe $Z_{M}(x)/Z_{M}(x)^0$ agit par multiplication \`a droite sur $\tilde{{\cal V}}_{reg}$ et le quotient de $\tilde{{\cal V}}_{reg}$ par cette action est ${\cal V}_{reg}$. On d\'efinit le syst\`eme local ${\cal C}_{reg}=\tilde{{\cal V}}_{reg}\otimes_{\xi_{x}}\bar{{\mathbb Q}}_{l}$ sur ${\cal V}_{reg}$ comme en \ref{faisceauxcaracterescuspidaux}. On introduit le syst\`eme local $p_{reg,*}({\cal C}_{reg})$ sur $\nabla_{reg}$. On note $K$ le complexe prolongement d'intersection de $p_{reg,*}({\cal C}_{reg})[dim(\nabla)]$.
    C'est un  faisceau pervers sur $G$ \`a support dans $\nabla$.  
  
     Lusztig et Shoji ont prouv\'e que l'alg\`ebre d'endomorphismes du complexe $K$ \'etait isomorphe \`a l'alg\`ebre du groupe $W(M)$. Concr\`etement, cela signifie ce qui suit. Pour $n\in Norm_{G}(M)$, l'automorphisme $Ad(n)$ conserve $M$ et se descend en un automorphisme de $M_{AD}$ que l'on note encore $Ad(n)$. Cet automorphisme  conserve  le faisceau ${\cal E}$ donc conserve ${\cal X}$.
  Notons $Norm_{G}(M,x)$ le sous-groupe des $n\in Norm_{G}(M)$ tels que $Ad(n)(x)=x$.  Pour  $n\in Norm_{G}(M)$, puisque $Ad(n)$ conserve ${\cal X}$,   il existe  $m\in M$ tel que $nm\in Norm_{G}(M,x)$. Donc l'inclusion $Norm_{G}(M,x)\to Norm_{G}(M)$ se quotiente en une bijection $Norm_{G}(M,x)/Z_{M}(x)\simeq W(M)$. Alors il existe un caract\`ere $\epsilon_{x}$ de $Norm_{G}(M,x)$ dont la restriction \`a $Z_{M}(x)$ est \'egale \`a  $\xi_{x}$.  Fixons un tel caract\`ere. On d\'efinit alors  l'action de $W(M)$ sur $K$ comme en \ref{faisceauxcaracteres}.

  Pour tout $\rho\in Irr(W(M))$, fixons une r\'ealisation de $\rho$ dans un espace $V_{\rho}$. Il existe un unique faisceau pervers irr\'eductible  $K_{M,{\cal E},\rho}$ de sorte que 
  $$K\simeq \oplus_{\rho\in Irr(W(M))}V_{\rho}\otimes K_{M,{\cal E},\rho},$$
  l'isomorphisme entrela\c{c}ant l'action de $W(M)$ sur $K$ en l'action \'evidente sur le membre de droite.   Les faisceaux-caract\`eres unipotents sur $G$ sont les faisceaux pervers $K_{M,{\cal E},\rho}$  quand on fait varier $M$, ${\cal E}$ et $\rho$. 
  
  Le caract\`ere $\epsilon_{x}$ n'est pas unique, on peut le multiplier par un caract\`ere de $W(M)$. Cela modifie l'indexation par $\rho$ des faisceaux-caract\`eres. 
  Soit $x'$ un autre point de ${\cal X}$. Fixons $m$ tel que $x'=Ad(m^{-1})(x)$. Alors $Norm_{G}(M,x')=
  m^{-1}Norm_{G}(M,x)m$. Le caract\`ere $\epsilon_{x}$ \'etant fix\'e, d\'efinissons le caract\`ere $\epsilon_{x'}$ de $Norm_{G}(M,x')$ par $\epsilon_{x'}(m^{-1}nm)=\epsilon_{x}(n)$ pour tout $n\in Norm_{G}(M)$.   Cette d\'efinition ne d\'epend pas du choix de $m$. En effectuant les m\^emes constructions o\`u l'on remplace $x$ et $\epsilon_{x}$ par $x'$ et $\epsilon_{x'}$, on obtient un complexe $K'$ muni d'une action de $W(M)$ et il existe un isomorphisme $W(M)$-\'equivariant de $K'$ sur $K$. 
  
  On appellera famille admissible  (pour $M$ et ${\cal E}$) une famille $(\epsilon_{x})_{x\in {\cal X}}$ telle que: 
  
  pour tout $x\in {\cal X}$, $\epsilon_{x}$ est un caract\`ere de $Norm_{G}(M,x)$ dont la restriction \`a $Z_{M}(x)$ est \'egale \`a  $\xi_{x}$;
  
  pour tous $x,x'\in {\cal X}$ et tout $m\in M$ tels que $x'=Ad(m^{-1})(x)$, on a l'\'egalit\'e $\epsilon_{x'}=\epsilon_{x}\circ Ad(m)$.
  
  Supposons que $M$ soit un ${\mathbb F}_{q}$-Levi et que ${\cal E}$ soit conserv\'e par l'action de $\Gamma_{{\mathbb F}_{q}}$ (comme on l'a dit en \ref{faisceauxcaracteresunipotentscuspidaux}, cette deuxi\`eme condition est automatique si $M_{AD}$ est absolument simple). Soit  $(\epsilon_{x})_{x\in {\cal X}}$  une famille admissible. Pour $x\in {\cal X}$, posons $^{Fr}\epsilon_{x}=\epsilon_{Fr(x)}\circ Fr$. Alors $(^{Fr}\epsilon_{x})_{x\in {\cal X}}$ est une famille admissible. Les conditions suivantes sont \'equivalentes:
  
  $^{Fr}\epsilon_{x}=\epsilon_{x}$ pour tout $x\in {\cal X}$;
  
  il existe $x\in {\cal X}^{\Gamma_{{\mathbb F}_{q}}}$ tel que  $\epsilon_{x}\circ Fr=\epsilon_{x}$;
  
  pour tout $x\in {\cal X}^{\Gamma_{{\mathbb F}_{q}}}$, on a  $\epsilon_{x}\circ Fr=\epsilon_{x}$.
  
  Si elles sont v\'erifi\'ees, on dira que la famille admissible 
  $(\epsilon_{x})_{x\in {\cal X}}$ est d\'efinie sur ${\mathbb F}_{q}$.
 Supposons donn\'ee une telle famille.  Soit $x\in {\cal X}^{\Gamma_{{\mathbb F}_{q}}}$.  Munissons le faisceau ${\cal E}$ d'une action de Frobenius, c'est-\`a-dire d'un isomorphisme $Fr^*({\cal E})\to {\cal E}$. Comme en \ref{fonctionscaracteristiques},   cet isomorphisme est de la forme
  $$\begin{array}{ccc}Fr^*(M/Z_{M}(x)^0\times_{\xi_{x}}\bar{{\mathbb Q}}_{l})&\to&M/Z_{M}(x)^0\times_{\xi_{x}}\bar{{\mathbb Q}}_{l}\\ (mZ_{M}(x)^0,v)&\mapsto&(Fr^{-1}(m)Z_{M}(x)^0,r_{x}v),\\ \end{array}$$
  o\`u $r_{x}$ est un  élément de $\bar{{\mathbb Q}}_{l}^{\times}$ dont on suppose  que c'est une racine de l'unité. On en d\'eduit une action de Frobenius sur $K$.   
    Pour tout $\rho\in Irr_{{\mathbb F}_{q}}(W(M))$, on choisit un prolongement continu $\rho^{\flat}$ de $\rho$ \`a $W(M)\rtimes \Gamma_{{\mathbb F}_{q}}$. Evidemment, si $\Gamma_{{\mathbb F}_{q}}$ agit trivialement sur $W(M)$, on suppose que $\rho^{\flat}$ est triviale sur $\Gamma_{{\mathbb F}_{q}}$. Comme en \ref{fonctionscaracteristiques}, on r\'ecup\`ere une action de Frobenius  sur $K_{M,{\cal E},\rho}$. On note $\boldsymbol{\chi}_{M,{\cal E},\rho}$ la fonction caract\'eristique de $K_{M,{\cal E},\rho}$. 
  
  Pour $w\in W(M)$, fixons un repr\'esentant $n$ de $w$ dans $Norm_{G}(M,x)$ et un \'el\'ement $g\in G$ tel que $gFr(g)^{-1}=n^{-1}$. On pose $M_{w}=g^{-1}Mg$, ${\cal X}_{w}=Ad(g^{-1})({\cal X})$, $x_{w}=Ad(g^{-1})(x)$. Comme en \ref{fonctionscaracteristiques}, on d\'efinit le faisceau-caract\`ere unipotent et cuspidal ${\cal E}_{w}$ sur $M_{w,AD}$. Il est muni de l'action de Frobenius associ\'ee comme ci-dessus au scalaire $r_{x_{w}}=r_{x}\epsilon_{x}(n)$. On reprend les constructions ci-dessus en rempla\c{c}ant $M$ et ${\cal E}$ par $M_{w}$ et ${\cal E}_{w}$. On obtient un complexe $K_{w}$ muni d'un Frobenius. On note $\boldsymbol{\chi}_{M,{\cal E},w}$ sa fonction caract\'eristique. Pour $\rho\in Irr_{{\mathbb F}_{q}}(W(M))$, on a alors l'\'egalit\'e
  
  (1) $\boldsymbol{\chi}_{M,{\cal E},\rho}=\vert W(M)\vert ^{-1}\sum_{w\in W(M)}trace(\rho^{\flat}(w^{-1}Fr))\boldsymbol{\chi}_{M,{\cal E},w}$. 
  
  {\bf Remarques.}(2)     Les  d\'efinitions ci-dessus d\'ependent du choix de la famille admissible $(\epsilon_{x})_{x\in {\cal X}}$ d\'efinie sur ${\mathbb F}_{q}$. Nous verrons en \ref{calculdorbites} qu'il existe une  telle famille  "canonique".

  \subsubsection{ Le caract\`ere $\mu_{t}$\label{mut}}
Soient $ M$ un Levi de $G$ et ${\cal E}\in {\bf FC}_{u}(M_{AD})$. On simplifie les notations en posant  ${\cal X}={\cal X}_{{\cal E}}$ etc.... Notons $S_{ss}$ l'ensemble des \'el\'ements de $M$ dont l'image dans $M_{AD}$ est la partie semi-simple d'un \'el\'ement de ${\cal X}$.   Soit $t\in S_{ss}$.
  Le groupe $G_{t}\cap M$ est un Levi de $G_{t}$: c'est le commutant de $Z(M)^0$ dans $G_{t}$. On a \'evidemment $M_{t}=Z_{M}(t)^0\subset G_{t}\cap M\subset Z_{M}(t)$. Puisque $G_{t}\cap M$ est un Levi, il est connexe donc $G_{t}\cap M=M_{t}$. Notons $t_{ad}$ l'image de $t$ dans $M_{AD}$ et notons $Z_{M}(t_{ad})$ le sous-groupe des $m\in M$ tels que $Ad(m)(t_{ad})=t_{ad}$. Alors $ M_{t}=Z_{M}(t_{ad})^0$. En particulier, on peut fixer un \'el\'ement unipotent $u\in M_{t}$ tel que $t_{ad}u\in {\cal X}$ (on identifie les \'el\'ements unipotents de $M$ et $M_{AD}$). On note $N$ l'\'el\'ement de $\mathfrak{m}_{t,SC}$ tel que $u=exp(N)$. 
Le caract\`ere $\xi_{t_{ad}u}$ de $Z_{M}(t_{ad}u)$ se restreint en un caract\`ere $\xi_{N}$ de $Z_{M_{t}}(N)=Z_{M_{t}}(u)$. On en d\'eduit un syst\`eme local ${\cal E}_{t}$ port\'e par l'orbite ${\cal O}_{N}$ de $N$ dans $\mathfrak{m}_{t,SC}$. Ce syst\`eme local est cuspidal. On peut appliquer les constructions de \ref{faisceauxcaracteres} au groupe $G_{t}$, \`a son Levi $M_{t}$ et au syst\`eme local ${\cal E}_{t}$. En particulier, on  d\'efinit le groupe $Norm_{G_{t}}(M_{t},N)$. En posant $W_{t}(M_{t})=Norm_{G_{t}}(M_{t})/M_{t}$, on a l'isomorphisme 
$Norm_{G_{t}}(M_{t},N)/Z_{M_{t}}(N)\simeq W_{t}(M_{t})$.  Montrons que
 
 (1) $Norm_{G_{t}}(M_{t})\subset Norm_{G}(M)$ et $Norm_{G_{t}}(M_{t},N)\subset Norm_{G}(M,t_{ad}u)$; ces  inclusions se quotientent en un homomorphisme injectif
 $\iota_{t}:W_{t}(M_{t})\to W(M)$. 
 
 Soit $n\in Norm_{G_{t}}(M_{t})$. Alors $Ad(n)$ conserve $M_{t}$ donc aussi $Z(M_{t})^0$. Mais $t_{ad}$ est un \'el\'ement isol\'e de $M_{AD}$, ce qui entra\^{\i}ne $Z(M_{t})^0=Z(M)^0$. Donc $Ad(n)$ conserve $Z(M)^0$ et aussi $M$, c'est-\`a-dire $n\in Norm_{G}(M)$. Supposons $n\in Norm_{G_{t}}(M_{t},N)$.  L'action de $Ad(n)$ dans $M_{AD}$ conserve $t_{ad}$ et $u=exp(N)$, donc conserve $t_{ad}u$. Cela prouve la premi\`ere assertion de (1). Ces inclusions se quotientent  \'evidemment en des homomorphismes injectifs 
 $$Norm_{G(t)}(M_{t})/(Norm_{G_{t}}(M_{t})\cap M)\to Norm_{G}(M)/(Norm_{G}(M)\cap M),$$
 $$Norm_{G(t)}(M_{t},N)/(Norm_{G_{t}}(M_{t},N)\cap M)\to Norm_{G}(M,t_{ad}u)/(Norm_{G}(M,t_{ad}u)\cap M).$$
  Mais les ensembles de départ sont $W_{t}(M_{t})$ et ceux d'arrivée sont  $W(M)$. Cela prouve (1).

 Fixons une famille admissible $(\epsilon_{x})_{x\in {\cal X}}$ pour le couple $(M,{\cal E})$. 
 Le caract\`ere $\epsilon_{t_{ad}u}$ de $Norm_{G}(M,t_{ad}u)$ se  restreint en un caract\`ere de $Norm_{G_{t}} (M_{t},N)$. On dispose aussi du caract\`ere canonique $\epsilon_{N}^{\flat}$ de $Norm_{G_{t}}(M_{t},N)$, cf. \ref{faisceauxcaracteres}. Ces deux caract\`eres  co\"{\i}ncident avec $\xi_{N}$ sur $Z_{M_{t}}(N)$. Il y a donc un unique caract\`ere $\mu_{t}$ de $W_{t}(M_{t})$ tel que $\mu_{t}(n)\epsilon_{N}^{\flat}(n)=\epsilon_{t_{ad}u}(n)$ pour tout $n\in  Norm_{G_{t}}(M_{t},N)$, o\`u, par abus de notation, on a identifi\'e $\mu_{t}$ \`a un caract\`ere de $ Norm_{G_{t}}(M_{t},N)$. On v\'erifie que ce caract\`ere $\mu_{t}$ est ind\'ependant du choix de l'\'el\'ement $u$.

 Montrons que
 
 (2) $W(M)$ est engendré par les images des homomorphismes $\iota_{t}$ quand $t$ décrit $S_{ss}$. 
 
   L'assertion est vide si $M=G$. Supposons $M\not=G$. D'après \cite{L2}, théorème 9.2, $W(M)$ est engendré par les sous-groupes $W^L(M)$ quand $L$ décrit les Levi de $G$ contenant $M$ comme sous-groupe de Levi propre maximal. Il suffit de démontrer la m\^eme assertion pour un tel sous-groupe. Fixons un sous-tore maximal $T$ contenu dans $M$ et notons $\Sigma$, resp. $\Sigma^L$,  $\Sigma^M$, l'ensemble des racines de $T$ dans $\mathfrak{g}$, resp. $\mathfrak{l}$, $\mathfrak{m}$. Fixons $\beta\in \Sigma^L-\Sigma^M$. 
    Puisque $S_{ss}$ est stable par multiplication par $Z(M)^0$, il est clair que l'on peut trouver $t\in S_{ss}$ tel que $\alpha(t)\not=1$ pour tout $\alpha\in \Sigma-\Sigma^L$ tandis que $\beta(t)=1$.  Pour un tel $t$, on a $G_{t}= L_{t}$ mais $G_{t}\not=M_{t}$. Alors $W_{t}(M_{t})$ n'est pas trivial et son image par $\iota_{t}$ est contenue dans $W^L(M)$. Puisque $M$ est maximal propre dans $L$, $W^L(M)$ n'a que deux éléments. Puisque $\iota_{t}$ est injective, elle est forcément surjective.   Cela prouve (2).

     \subsubsection{L'ensemble $\Lambda$ associ\'e \`a un couple $(M,{\cal E})$\label{Lambda}}
   
   Pour la fin de la sous-section 5.2, {\bf on suppose que} $G$ {\bf est adjoint et absolument simple}. A tout sous-ensemble propre $\Lambda\subset \Delta_{a}$, on a associ\'e en \ref{constructiondeLeviG} un Levi $M_{\Lambda}$ et un \'el\'ement semi-simple $t_{\Lambda}\in M_{\Lambda,AD}$.
   
    \begin{prop}{Soient $M$ un Levi de $G$ et ${\cal E}\in {\bf FC}_{u}(M_{AD})$. 
    
    (i) Il existe un unique sous-ensemble propre $\Lambda\subset \Delta_{a}$  tel  qu'il existe $g\in G$ de sorte que $gMg^{-1}=M_{\Lambda}$  et que $Ad(g)$ transporte ${\cal X}_{{\cal E}}$ en un sous-ensemble de $M_{\Lambda,AD}$ contenant un \'el\'ement de partie semi-simple  $t_{\Lambda}$. L'ensemble $\Lambda$ est conservé par l'action du groupe $\Omega$. 
   
   (ii) Soit $t$ un \'el\'ement de $M$ dont l'image dans $M_{AD}$ est la partie semi-simple d'un \'el\'ement de ${\cal X}_{{\cal E}}$. Alors il existe $g\in G$ et $x\in \bar{C}^{nr}\cap {\mathbb T}$ tels que $gtg^{-1}={\bf j}_{T}(x)$, $\Lambda\subset \Lambda(x)$, $gMg^{-1}=M_{\Lambda}$ et que  l'image de $gtg^{-1}$ dans $M_{\Lambda,AD}$ soit $t_{\Lambda}$. }\end{prop} 
   
   Preuve. Les classes de conjugaison de Levi de $G$ correspondent bijectivement
aux classes de conjugaison de Levi de $G_{F}$. Il r\'esulte de    \ref{faisceauxcaracteresunipotentscuspidaux}(3) que les classes de conjugaison de couples $(M,{\cal E})$ form\'es d'un Levi de $G$ et d'un faisceau ${\cal E}\in {\bf FC}_{u}(M_{AD})$ correspondent bijectivement aux les classes de conjugaison de triplets $(M_{F},s_{M_{F}},{\cal E}_{F})$, o\`u $M_{F}$ est un Levi de $G_{F}$, $s_{M_{F}}\in S^{nr,st}(M_{F})$ et ${\cal E}_{F}\in {\bf FC}^{st}(M_{F,s_{M_{F}},SC})$. Les assertions d'existence de la  proposition  sont les analogues de la proposition \ref{couples} et l'assertion d'unicité est l'analogue de l'assertion (1) de \ref{couplesstables} . La proposition se d\'emontre par les m\^emes arguments. On laisse la mise au point au lecteur. $\square$

{\bf Remarques.} (1) De m\^eme qu'en \ref{couplesstables}, l'ensemble $\Lambda$ ne d\'ependait que de $M_{F}$ et $s_{M_{F}}$ et pas de ${\cal E}_{F}$, de m\^eme l'ensemble $\Lambda$ de la proposition ne d\'epend que de $M$ et de  l'orbite ${\cal X}_{{\cal E}}$.

(2) L'ensemble $\Lambda$ se  d\'ecrit cas par cas par les m\^emes formules qu'en \ref{couplesstables}.

 \subsubsection{Restriction au commutant d'un \'el\'ement semi-simple\label{restriction}}
 Comme dans le paragraphe précédent, on suppose que $G$ est adjoint, absolument simple. {\bf On suppose de plus que }$G$  {\bf est d\'eploy\'e sur }${\mathbb F}_{q}$. Soit $M$ un ${\mathbb F}_{q}$-Levi de $G$ et soit ${\cal E}\in {\bf FC}_{u}(M_{AD})$. On simplifie les notations en abandonnant les indices ${\cal E}$ des objets associ\'es \`a ce faisceau. Soit $t$ un \'el\'ement semi-simple  de $M({\mathbb F}_{q})$, notons $t_{ad}$ son image dans $M_{AD}$, soit $u$ un \'el\'ement unipotent de $M_{AD,t_{ad}}({\mathbb F}_{q})$, supposons $t_{ad}u\in {\cal X}$. On suppose que le groupe  $G_{t}$ est d\'eploy\'e sur ${\mathbb F}_{q}$.  
  
 On utilise les constructions de \ref{mut}. 
 Pour $\rho\in Irr(W(M))$, posons $\rho_{t}=\rho\circ\iota_{t}$. Cette repr\'esentation  se d\'ecompose en somme de repr\'esentations irr\'eductibles de $W_{t}(M_{t})$. Pour $\rho'\in Irr(W_{t}(M_{t}))$, on note $<\rho_{t},\rho'>$ la multiplicit\'e de $\rho'$ dans $\rho_{t}$. 
 
 Fixons une famille admissible $(\epsilon_{x})_{x\in {\cal X}}$ pour $(M,{\cal E})$, que l'on suppose d\'efinie sur ${\mathbb F}_{q}$. On fixe une action de Frobenius sur ${\cal E}$, dont on d\'eduit une telle action sur ${\cal E}_{t}$, puis diverses fonctions caract\'eristiques, cf. \ref{faisceauxcaracteres} et \ref{faisceauxcaracteresunipotents}. 
  
  \begin{prop}{Pour tous $\rho\in Irr(W(M))$ et $N'\in  \mathfrak{g}_{t,nil}({\mathbb F}_{q})$, on a l'\'egalit\'e
  $$\boldsymbol{\chi}_{M,{\cal E},\rho}(t\,exp(N'))=\sum_{\rho'\in Irr(W_{t}(M_{t}))}<\rho_{t},\mu_{t}\rho'>\boldsymbol{\chi}_{M_{t},{\cal E}_{t},\rho'}(N').$$}\end{prop}
  
  La preuve est due \`a Lusztig (voir aussi \cite{AA} th\'eor\`eme 6.8). On la reprend en \ref{preuverestriction} pour faire appara\^{\i}tre clairement le caract\`ere $\mu_{t}$.  
  
  \subsubsection{Un lemme de restriction}\label{lemmerestriction}
  On conserve les données du paragraphe précédent. Pour $w\in W(M)$, soit $W_{t}(M_{t})[w]$ l'ensemble des $w'\in W_{t}(M_{t})$ tels que $\iota_{t}(w') $ soit conjugué à $w$ par un élément de $W(M)$.

  \begin{lem}{Pour tous $w\in W(M)$ et $N'\in  \mathfrak{g}_{t,nil}({\mathbb F}_{q})$, on a l'\'egalit\'e
  $$  \boldsymbol{\chi}_{M,{\cal E},w}(t\,exp(N'))=\vert Z_{W(M)}(w)\vert \vert W_{t}(M_{t})\vert ^{-1}\sum_{w'\in W_{t}(M_{t})[w]}\mu_{t}(w')\boldsymbol{\chi}_{M_{t},{\cal E}_{t},w'}(N').$$}\end{lem}
  
 Preuve.     Pour tout $w'\in W_{t}(M_{t})$, on choisit un repr\'esentant $\nu_{w'}\in Norm_{G_{t}}(M_{t},N)$ de $w'$ et un \'el\'ement $\gamma_{w'}\in G_{t}$ tel que $\gamma_{w'}Fr(\gamma_{w'})^{-1}=\nu_{w'}^{-1}$. On  construit le Levi $M_{t,w'}=\gamma_{w'}^{-1}M_{t}\gamma_{w'}$ de $G_{t}$ et le faisceau ${\cal E}_{t,w'}$ sur $\mathfrak{m}_{t,w',SC}$. D'o\`u la fonction $\boldsymbol{\chi}_{M_{t},{\cal E}_{t},w'}$ sur $\mathfrak{g}_{t}({\mathbb F}_{q})$.     De m\^eme, pour $w\in W(M)$, on choisit un repr\'esentant $n_{w}\in Norm_{G}(M,t_{ad}u)$ de $w$ et un \'el\'ement $g_{w}\in G$ tel que $g_{w}Fr(g_{w})^{-1}=n_{w}^{-1}$. On construit le Levi $M_{w}=g_{w}^{-1}Mg_{w}$ de $G$ et le faisceau ${\cal E}_{w}$ sur $M_{w,AD}$. D'o\`u la  fonction $\boldsymbol{\chi}_{M,{\cal E},w}$.    On note $S$, resp. $S_{w}$, l'image r\'eciproque dans $M$, resp.  $M_{w}$, du support de ${\cal E}$, resp. ${\cal E}_{w}$, et $S_{ss}$, resp. $S_{w,ss}$, l'ensemble des parties semi-simples d'\'el\'ements de $S$, resp. $S_{w}$. On note $\pi^M:M\to M_{AD}$, resp. $\pi^{M_{w}}=M_{w}\to M_{w,AD}$ les projections naturelles. Pour $y\in M$, resp. $y\in M_{w}$, on pose simplement $y_{ad}=\pi^M(y)$, resp. $y_{ad}=\pi^{M_{w}}(y)$. Fixons $w\in W(M)$. Notons ${\cal K}_{w}$ l'ensemble des $k\in G({\mathbb F}_{q})$ tels que $ktk^{-1}\in S_{w,ss}$. Pour $k\in {\cal K}_{w}$, posons $M[k]=k^{-1}M_{w}k\cap G_{t}$. C'est un Levi de $G_{t}$ d\'efini sur ${\mathbb F}_{q}$: c'est le commutant dans $G_{t}$ du sous-tore $k^{-1}Z(M_{w})^0k$ de $G_{t}$. 
Notons ${\cal O}[k]$ l'ensemble des $Y\in \mathfrak{m}[k]_{nil}$ tels que $kt\,exp(Y)k^{-1}\in S_{w}$.  Posons $t[k]=ktk^{-1}$. L'application $Ad(k)$ identifie $M[k]$ à $M_{w}\cap G_{t[k]}=M_{w,t_[k]}$ et ${\cal O}[k]$ à l'ensemble des $Y\in \mathfrak{m}_{w,t[k]}$ tels que $t[k]exp(Y)\in S_{w}$. Ce dernier ensemble est une unique orbite nilpotente sous $M_{w,t[k]}$, donc ${\cal O}[k]$ est une unique orbite nilpotente sous $M[k]$.  On munit ${\cal O}[k]$ d'un syst\`eme local  ${\cal E}[k]$ muni d'une action de Frobenius: c'est l'image r\'eciproque de ${\cal E}_{w}$ par l'application $Y\mapsto \pi^{M_{w}}(kt\,exp(Y)k^{-1})$. 
  On en d\'eduit comme en \ref{fonctionscaracteristiques} un complexe $K[k]$ sur $G_{t}$ et sa fonction caract\'eristique que l'on note $\boldsymbol{\chi}[k]$. D'apr\`es \cite{L3} th\'eor\`eme 8.5, pour $N'\in  \mathfrak{g}_{t,nil}({\mathbb F}_{q})$,  on a l'\'egalit\'e
 $$(1) \qquad \boldsymbol{\chi}_{M,{\cal E},w}(t \,exp(N'))=\vert G_{t}({\mathbb F}_{q})\vert ^{-1}\vert M_{w}({\mathbb F}_{q})\vert ^{-1}\sum_{k\in {\cal K}_{w}}\vert M[k]({\mathbb F}_{q})\vert  \boldsymbol{\chi}[k](N').$$

 Montrons que 
 
 (2) pour tout $k\in {\cal K}_{w}$, il existe  $\gamma\in G_{t}$ tels que $M[k]=\gamma M_{t}\gamma^{-1}$ et  que  $g_{w}k\gamma \in Norm_{G}(M)$. 
 
  Soit $k\in {\cal K}_{w}$. Alors $ktk^{-1}\in S_{w,ss}$, d'o\`u $g_{w}ktk^{-1}g_{w}^{-1}\in S_{ss}$.  L'ensemble $S_{ss}$ est celui des $t'\in M$ dont l'image $t'_{ad}$ dans $M_{AD}$ est conjugu\'e \`a $t_{ad}$. On peut donc fixer $m\in M$ tel qu'en posant $h=mg_{w}k$ et $t'=hth^{-1}$, on ait $t'\in M$ et $t'_{ad}=t_{ad}$. Posons $M_{t'}=G_{t'}\cap M$.  Comme on l'a dit en \ref{mut}, on a $M_{t'}=Z_{M}(t'_{ad})^0=Z_{M}(t_{ad})^0$ et ce groupe   est un Levi de $G_{t'}$.  On a $Ad(h)^{-1}(M_{t'})=M[k]$ d'apr\`es les d\'efinitions de $h$ et $t'$. Introduisons le sous-ensemble  $\Lambda\subset \Delta_{a}$ associ\'e \`a $M$ et ${\cal E}$, cf. \ref{Lambda}. D'apr\`es la proposition \ref{Lambda}, on peut conjuguer la situation par un \'el\'ement de $G$ de sorte que
   $M=M_{\Lambda}$, $t_{ad}=t_{\Lambda}$ et qu'il existe $x\in \bar{C}^{nr}\cap a_{T}({\mathbb T})$ de sorte que $t={\bf j}_{T}(x)\in T$ et $\Lambda\subset \Lambda(x)$. 
  Le syst\`eme de racines de $G_{t}$ a pour base $\Lambda(x)$. Le Levi $ M[k]$ est conjugu\'e par un \'el\'ement de $G_{t}$ \`a un Levi standard de $G_{t}$. C'est-\`a-dire qu'il existe $\Lambda'\subset \Lambda(x)$ et $\gamma\in G_{t}$ tels que $Ad(\gamma^{-1}h^{-1})(M_{t'})=L_{\Lambda'}$. Autrement dit, puisque $M_{t'}=Z_{M}(t_{ad})^0=L_{\Lambda}$, on a $Ad(\gamma^{-1}h^{-1})(L_{\Lambda})=L_{\Lambda'}$.  Donc $Ad(\gamma^{-1}h^{-1})$ envoie la composante neutre du centre de $ L_{\Lambda}$ sur la composante neutre du centre de $L_{\Lambda'}$.  Ces composantes neutres $A_{\Lambda}$ et $A_{\Lambda'}$ sont aussi celles des centres de $M_{\Lambda}$ et $M_{\Lambda'}$. Donc $Ad(\gamma^{-1}h^{-1})(M_{\Lambda})=M_{\Lambda'}$. D'autre part, par construction, $Ad(\gamma^{-1}h^{-1})$ envoie $t'$ sur $t$, donc l'image $t'_{ad}=t_{ad}=t_{\Lambda}$ de $t'$ dans $M_{\Lambda}$ sur l'image de $t$ dans $M_{\Lambda',AD}$, c'est-\`a-dire $t_{\Lambda'}$. L'unicit\'e de la proposition \ref{Lambda} implique l'\'egalit\'e $\Lambda'=\Lambda$.  En cons\'equence, on a $L_{\Lambda'}=L_{\Lambda}=M_{t}$ et $M[k]=\gamma M_{t}\gamma^{-1}$. L'\'el\'ement $h\gamma=mg_{w}k\gamma$ normalise $M$ donc $g_{w}k\gamma$ aussi. Cela prouve (2).

  Fixons un ensemble ${\cal W}_{t}$ de repr\'esentants des classes de conjugaison dans $W_{t}(M_{t})$. Montrons que

 (3) pour tout $k\in {\cal K}_{w}$, il existe un unique $w'\in {\cal W}_{t}$ et il existe $\kappa\in G_{t}({\mathbb F}_{q})$ tels que $k^{-1}M_{w}k=\kappa^{-1}\gamma_{w'}^{-1}M\gamma_{w'}\kappa$.
 
  Soit $k\in {\cal K}_{w}$.  Fixons $\gamma$ v\'erifiant (2). Puisque $Ad(\gamma)$ conjugue $M_{t}$ en $M[k]$ et que ces Levi sont d\'efinis sur ${\mathbb F}_{q}$, on voit ais\'ement qu'il existe  $w'\in {\cal W}_{t}$,  $\nu\in Norm_{G_{t}}(M_{t})$ et $\kappa\in G_{t}({\mathbb F}_{q})$ de sorte que $\gamma=\kappa^{-1}\gamma_{w'}^{-1}\nu^{-1}$. On  a $Norm_{G_{t}}(M_{t})\subset Norm_{G}(M)$.  Puisque  $g_{w}k\gamma$ appartient aussi à $ Norm_{G}(M)$, on a $g_{w}k\kappa^{-1}\gamma_{w'}^{-1}\in Norm_{G}(M)$. Donc $Ad(\kappa^{-1}\gamma_{w'}^{-1})(M)=Ad(k^{-1}g_{w}^{-1})(M)$, ou encore $k^{-1}M_{w}k=\kappa^{-1}\gamma_{w'}^{-1}M\gamma_{w'}\kappa$. D'o\`u l'existence de $w'$ et $\kappa$ v\'erifiant les conditions de (3). Soient $\underline{w}'\in {\cal W}_{t}$ et $\underline{\kappa}\in G_{t}({\mathbb F}_{q})$ tels que $k^{-1}M_{w}k=\underline{\kappa}^{-1}\gamma_{\underline{w}'}^{-1}M\gamma_{\underline{w}'}\underline{\kappa}$. On a alors $\kappa^{-1}\gamma_{w'}^{-1}M\gamma_{w'}\kappa\cap G_{t}=\underline{\kappa}^{-1}\gamma_{\underline{w}'}^{-1}M\gamma_{\underline{w}'}\underline{\kappa}\cap G_{t}$, c'est-\`a-dire $\kappa^{-1}M_{t,w'}\kappa=\underline{\kappa}^{-1}M_{t,\underline{w}'}\underline{\kappa}$. Par d\'efinition de ${\cal W}_{t}$, les Levi $M_{t,w'}$ et $M_{t,\underline{w}'}$ ne peuvent \^etre conjugu\'es par un \'el\'ement de $G_{t}({\mathbb F}_{q})$ que si $w'=\underline{w}'$. Cela d\'emontre l'unicit\'e de $w'$, d'o\`u (3).
  
  Pour tout $w'\in {\cal W}_{t}$, notons ${\cal K}_{w,w'}$ le sous-ensemble des $k\in {\cal K}_{w}$ tels que $k^{-1}M_{w}k$ et $\gamma_{w'}^{-1}M\gamma_{w'}$ soient conjugu\'es par un \'el\'ement de $G_{t}({\mathbb F}_{q})$. D'apr\`es (3), ${\cal K}_{w}$ est union disjointe des ${\cal K}_{w,w'}$ quand $w'$ parcourt ${\cal W}_{t}$. Montrons que
  
  (4) soient $w'\in {\cal W}_{t}$ et $k\in {\cal K}_{w,w'}$; alors on a l'\'egalit\'e $ \boldsymbol{\chi}[k]=\mu_{t}(w')\boldsymbol{\chi}_{M_{t},{\cal E}_{t},w'}$. 
  
 Fixons  $\kappa\in G_{t}({\mathbb F}_{q})$ tel que $k^{-1}M_{w}k=\kappa^{-1}\gamma_{w'}^{-1}M\gamma_{w'}\kappa$. Posons $r=\gamma_{w'}\kappa k^{-1}g_{w}^{-1}$. L'\'egalit\'e  précédente \'equivaut \`a $r\in Norm_{G}(M)$.  On peut fixer $m\in M$ de sorte que, en posant $n=mr$, on ait $n\in Norm_{G}(M,t_{ad}u)$.  Posons $m_{w}=g_{w}^{-1}mg_{w}$. C'est un élément de $M_{w}$. On a le diagramme commutatif
 $$(5)\quad \begin{array}{ccccc} M[k]&\stackrel{Ad(k)}{\to}&M_{w}&\stackrel{Ad(m_{w})}{\to}&M_{w}\\ Ad(\kappa)\downarrow\qquad&&Ad(g_{w})\downarrow\qquad\quad&& \\M_{t,w'}&&M&& Ad(g_{w})\downarrow\qquad\quad\\ Ad(\gamma_{w'})\downarrow\qquad\quad&Ad(r)\nearrow\quad&&Ad(m)\searrow\quad&\\ M&&\stackrel{Ad(n)}{\to}&&M\\ \end{array}$$
 
  L'orbite ${\cal O}[k]$ est l'ensemble des $Y\in \mathfrak{m}[k]_{nil}$ tels que $kt\,exp(Y)k^{-1}\in S_{w}$ et ${\cal E}[k]$ est l'image r\'eciproque de ${\cal E}_{w}$ par l'application $Y\mapsto \pi^{M_{w}}(kt\,exp(Y)k^{-1})$. La commutativité du diagramme ci-dessus entraîne que ${\cal O}[k]$ est aussi l'ensemble des $Y\in \mathfrak{m}[k]_{nil}$ tels que $\gamma_{w'}\kappa t\,exp(Y)\kappa^{-1}\gamma_{w'}^{-1}\in S$ et que ${\cal E}[k]$ est l'image réciproque de ${\cal E}$ par l'application $Y\mapsto \pi^M(\gamma_{w'}\kappa t\,exp(Y)\kappa^{-1}\gamma_{w'}^{-1})$.    Posons $u=exp(N)$, $N_{t,w'}=Ad(\gamma_{w'})^{-1}(N)$, notons ${\cal O}$ l'orbite de $N$ dans $\mathfrak{m}_{t}$ et ${\cal O}_{t,w'}$ celle de $N_{t,w'}$ dans $\mathfrak{m}_{t,w'}$. Puisque $\kappa $ appartient à $G_{t}$, on obtient que ${\cal O}[k]= {\cal O}_{t,w'}$ et que ${\cal E}[k] ={\cal E}_{t,w'}$.    Il en r\'esulte que les fonctions $\boldsymbol{\chi}[k]$ et $\boldsymbol{\chi}_{M_{t},{\cal E}_{t},w'}$ sont proportionnelles. Pour calculer la constante de proportionnalit\'e, on doit comparer les actions de Frobenius  sur nos différents faisceaux. Posons   $Y=Ad(\kappa)^{-1}(N_{t,w'})$. L'action de Frobenius initiale sur ${\cal E}$ est d\'etermin\'ee par un scalaire $r_{x_{ad}}\in {\mathbb C}^{\times}$. L'action de Frobenius sur ${\cal E}_{t,w'}$ est déterminée par le scalaire $r_{N_{t,w'}}=r_{x_{ad}}\epsilon_{N}^{\flat}(Fr(\gamma_{w'})\gamma_{w'}^{-1})$. Puisque $\kappa\in G_{t}({\mathbb F}_{q})$,  on a $N\in \mathfrak{m}[k]_{nil}({\mathbb F}_{q})$ et cette action est aussi déterminée par le  m\^eme scalaire au point $N$, que l'on note $r_{{\cal E}_{t,w'},N}=r_{x_{ad}}\epsilon_{N}^{\flat}(Fr(\gamma_{w'})\gamma_{w'}^{-1})$. 
  Posons $x_{w}=g_{w}^{-1}xg_{w}$. L'action de Frobenius sur ${\cal E}_{w}$ est déterminée par le scalaire $r_{x_{w,ad}}=r_{x_{ad}}\epsilon_{x_{ad}}(Fr(g_{w})g_{w}^{-1})$. Posons $x'_{w}= m_{w}^{-1}x_{w}m_{w}$ et $x''_{w}=kt\,exp(N)k^{-1}$. La commutativité du diagramme ci-dessus et le fait que $Ad(n)_{M_{AD}}$ fixe $x_{ad}$  impliquent que $x'_{w,ad}=x''_{w,ad}$.  Puisque $x''_{w}\in G({\mathbb F}_{q})$, cela entra\^{\i}ne que  $Fr(m_{w})m_{w}^{-1}$ appartient à $Z_{M_{w}}(x_{w,ad})$  et que 
   l'action de Frobenius sur ${\cal E}_{w}$ est aussi déterminée par le scalaire $r_{x''_{w,ad}}=r_{w,ad}\xi_{x_{w,ad}}(Fr(m_{w})m_{w}^{-1})=r_{x_{ad}}\epsilon_{x_{ad}}(Fr(g_{w})g_{w}^{-1})\xi_{x_{w,ad}}(Fr(m_{w})m_{w}^{-1})$. Par définition de ${\cal E}[k]$, son action de Frobenius est déterminée par le m\^eme scalaire en $N$, que l'on note $r_{{\cal E}[k],N}=r_{x_{ad}}\epsilon_{x_{ad}}(Fr(g_{w})g_{w}^{-1})\xi_{x_{w,ad}}(Fr(m_{w})m_{w}^{-1})$. Pour démontrer (4), nous devons prouver l'égalité $r_{{\cal E}[k],N}=
 \mu_{t}(w')r_{{\cal E}_{t,w'},N}$, autrement dit
 $$r_{x_{ad}}\epsilon_{x_{ad}}(Fr(g_{w})g_{w}^{-1})\xi_{x_{w,ad}}(Fr(m_{w})m_{w}^{-1})=\mu_{t}(w')r_{x_{ad}}\epsilon_{N}^{\flat}(Fr(\gamma_{w'})\gamma_{w'}^{-1}),$$
 ou encore 
 
 (6) $\epsilon_{x_{ad}}(Fr(g_{w})g_{w}^{-1})\xi_{x_{w,ad}}(Fr(m_{w})m_{w}^{-1})=\mu_{t}(w')\epsilon_{N}^{\flat}(Fr(\gamma_{w'})\gamma_{w'}^{-1})$. 
 
 Par définition, on a l'égalité $\xi_{x_{w,ad}}(Fr(m_{w})m_{w}^{-1})=\xi_{x_{ad}}(g_{w}Fr(m_{w})m_{w}^{-1}g_{w}^{-1})$, puis, par définition de la famille admissible $(\epsilon_{x'})_{x'\in {\cal X}}$, $\xi_{x_{w,ad}}(Fr(m_{w})m_{w}^{-1})=\epsilon_{x_{ad}}(g_{w}Fr(m_{w})m_{w}^{-1}g_{w}^{-1})$. Le membre de gauche de (6) est donc égal à  $\epsilon_{x_{ad}}(Fr(g_{w}m_{w})(g_{w}m_{w})^{-1})$. Puisque $k\in G({\mathbb F}_{q})$, on a $Fr(g_{w}m_{w})(g_{w}m_{w})^{-1}=Fr(g_{w}m_{w}k)(g_{w}m_{w}k)^{-1}$. La commutativité du diagramme (5) montre que $g_{w}m_{w}k=n\gamma_{w'}\kappa$. Puisque $\kappa$ appartient à $G_{t}({\mathbb F}_{q})$, on obtient $Fr(g_{w}m_{w})(g_{w}m_{w})^{-1}=Fr(n\gamma_{w'})(n\gamma_{w'})^{-1}$ et le membre de gauche de (6)  est égal à $\epsilon_{x_{ad}}(Fr(n\gamma_{w'})(n\gamma_{w'})^{-1})$. Mais $n$ appartient à $Norm_{G}(M,x_{ad})$. L'expression précédente est donc égale à $\epsilon_{x_{ad}}(Fr(n))\epsilon_{x_{ad}}(Fr(\gamma_{w'})\gamma_{w'}^{-1})\epsilon_{x_{ad}}(n)^{-1}$. Les termes extr\^emes disparaissent car notre famille admissible est définie sur ${\mathbb F}_{q}$. Le membre de gauche de (6) vaut donc $\epsilon_{x_{ad}}(Fr(\gamma_{w'})\gamma_{w'}^{-1})$. Alors l'égalité (6) résulte de la définition du caractère $\mu_{t}$. Cela prouve (6) et (4).

Prouvons que

(7) soit $w'\in {\cal W}_{t}$; alors ${\cal K}_{w,w'}$ est l'ensemble des \'el\'ements $k\in G({\mathbb F}_{q})$ tels que $k^{-1}M_{w}k$ soit conjugu\'e \`a $\gamma_{w'}^{-1}M\gamma_{w'}$ par un \'el\'ement de $G_{t}({\mathbb F}_{q})$. 

Notons ${\cal K}'_{w,w'}$ l'ensemble des \'el\'ements $k\in G({\mathbb F}_{q})$ tels que $k^{-1}M_{w}k$ soit conjugu\'e \`a $\gamma_{w'}^{-1}M\gamma_{w'}$ par un \'el\'ement de $G_{t}({\mathbb F}_{q})$. Par d\'efinition de ${\cal K}_{w,w'}$,  on a l'inclusion ${\cal K}_{w,w'}\subset {\cal K}'_{w,w'}$. Inversement, toujours par d\'efinition de ${\cal K}_{w,w'}$, pour d\'emontrer que ${\cal K}'_{w,w'}$ est contenu dans ${\cal K}_{w,w'}$, il suffit de d\'emontrer qu'il est contenu dans ${\cal K}_{w}$. Soit $k\in {\cal K}'_{w,w'}$, fixons $\kappa\in G_{t}({\mathbb F}_{q})$ tel que $k^{-1}M_{w}k=\kappa^{-1}\gamma_{w'}^{-1}M\gamma_{w'}\kappa$. Alors $g_{w}k\kappa^{-1}\gamma_{w'}^{-1}$ normalise $M$, donc $Ad(g_{w}k\kappa^{-1}\gamma_{w'}^{-1})$ conserve $S$. Donc $Ad(g_{w}k\kappa^{-1}\gamma_{w'}^{-1})(t)\in S_{ss}$ et $Ad(k\kappa^{-1}\gamma_{w'}^{-1})(t)\in S_{w,ss}$. Mais $Ad(\kappa)$ et $Ad(\gamma_{w'})$ fixent $t$ donc $Ad(k)(t)\in S_{w,ss}$. Par d\'efinition, cela signifie que $x\in {\cal K}_{w}$. D'o\`u (7).

L'ensemble $W_{t}(M_{t})[w]$ est celui  des $w'\in W_{t}(M_{t})$ tels que $\gamma_{w'}^{-1}M\gamma_{w'}$ soit conjugué à $M_{w}$ par un élément de $G({\mathbb F}_{q})$.     
Posons ${\cal W}_{t}[w]={\cal W}_{t}\cap W_{t}(M_{t})[w]$. Si $w'\in {\cal W}_{t}-{\cal W}_{t}[w]$, l'ensemble ${\cal K}_{w,w'}$ est vide. Soit $w'\in {\cal W}_{t}[w]$. Fixons $k_{0}\in G({\mathbb F}_{q})$ tel que $ k_{0}^{-1}M_{w}k_{0}=\gamma_{w'}^{-1}M\gamma_{w'}$. Alors ${\cal K}_{w,w'}=Norm_{G({\mathbb F}_{q})}(M_{w})k_{0}G_{t}({\mathbb F}_{q})$. L'intersection $k_{0}^{-1}Norm_{G({\mathbb F}_{q})}(M_{w})k_{0}\cap G_{t}({\mathbb F}_{q})$ est \'egale \`a $Norm_{G_{t}({\mathbb F}_{q})}(\gamma_{w'}^{-1}M\gamma_{w'})$. Par le m\^eme raisonnement qu'en \ref{mut}(1), ce groupe est \'egal \`a $Norm_{G_{t}({\mathbb F}_{q})}(M_{t,w'})$. On obtient
$$\vert {\cal K}_{w,w'}\vert =\vert Norm_{G({\mathbb F}_{q})}(M_{w})\vert \vert G_{t}({\mathbb F}_{q})\vert \vert Norm_{G_{t}({\mathbb F}_{q})}(M_{t,w'})\vert ^{-1}.$$
Pour tout $k\in {\cal K}_{w,w'}$, on a $\vert M[k]({\mathbb F}_{q})\vert =\vert M_{t,w'}({\mathbb F}_{q})\vert $ et $\boldsymbol{\chi}[k]=\mu_{t}(w')\boldsymbol{\chi}_{M_{t},{\cal E}_{t},w'}$. L'\'egalit\'e (1) devient
$$(8) \qquad \boldsymbol{\chi}_{M,{\cal E},w}(t\,exp(N'))=\vert Norm_{G({\mathbb F}_{q})}(M_{w})\vert \vert M_{w}({\mathbb F}_{q})\vert ^{-1}\sum_{w'\in {\cal W}_{t}[w]}\vert M_{t,w'}({\mathbb F}_{q})\vert$$
$$ \vert Norm_{G_{t}({\mathbb F}_{q})}(M_{t,w'})\vert ^{-1}\mu_{t}(w')\boldsymbol{\chi}_{M_{t},{\cal E}_{t},w'}(N').$$
Par le th\'eor\`eme de Lang, on a l'\'egalit\'e
$$Norm_{G({\mathbb F}_{q})}(M_{w})/M_{w}({\mathbb F}_{q})=(Norm_{G}(M_{w})/M_{w})^{\Gamma_{{\mathbb F}_{q}}}.$$L'application $Ad(g_{w})$ identifie $Norm_{G}(M_{w})/M_{w}$ muni de son action de Frobenius \`a $W(M)$ muni de l'action $v\mapsto w^{-1}vw$. Il en r\'esulte que $\vert Norm_{G({\mathbb F}_{q})}(M_{w})/M_{w}({\mathbb F}_{q})\vert=\vert Z_{W(M)}(w)\vert$. De m\^eme $\vert Norm_{G_{t}({\mathbb F}_{q})}(M_{t,w'})/M_{t,w'}({\mathbb F}_{q})\vert =\vert Z_{W_{t}(M_{t})}(w')\vert $. Dans la formule (8), on peut remplacer la sommation sur $w'\in {\cal W}_{t}[w]$ par une sommation sur $W_{t}(M_{t})[w]$, \`a condition de multiplier le terme index\'e par $w'$ par l'inverse du nombre d'\'el\'ements de l'orbite de $w'$, c'est-\`a-dire $\vert W_{t}(M_{t})\vert^{-1} \vert Z_{W_{t}(M_{t})}(w')\vert $. Cette formule devient celle de l'énoncé. $\square$

\subsubsection{Preuve de la proposition \ref{restriction}}\label{preuverestriction}
Soient $\rho\in Irr(W(M))$ et $N'\in \mathfrak{g}_{t,nil}({\mathbb F}_{q})$. 
 D'apr\`es \ref{faisceauxcaracteresunipotents}(1), on a l'\'egalit\'e
$$\boldsymbol{\chi}_{M,{\cal E},\rho}(t \,exp(N'))=\vert W(M)\vert ^{-1}\sum_{w\in W(M)}trace(\rho(w)^{-1})\boldsymbol{\chi}_{M,{\cal E},w}(N').$$
En utilisant le lemme \ref{lemmerestriction}, on obtient
 $$\boldsymbol{\chi}_{M,{\cal E},\rho}(t \,exp(N'))=\vert W(M)\vert ^{-1}\sum_{w\in W(M)}trace(\rho(w)^{-1})\vert Z_{W(M)}(w)\vert \vert W_{t}(M_{t})\vert ^{-1}$$
 $$\sum_{w'\in W_{t}(M_{t})[w]}\mu_{t}(w')\boldsymbol{\chi}_{M_{t},{\cal E}_{t},w'}(N').$$
 Pour $w'\in W_{t}(M_{t})$,  notons $Cl_{W(M)}(w')$ la classe de conjugaison par $W(M)$ de $\iota_{t}(w')$.   On obtient
  $$\boldsymbol{\chi}_{M,{\cal E},\rho}(t \,exp(N'))=\vert W_{t}(M_{t})\vert ^{-1}\sum_{w'\in W_{t}(M_{t})}\mu_{t}(w')\boldsymbol{\chi}_{M_{t},{\cal E}_{t},w'}(N')$$
  $$\sum_{w\in Cl_{W(M)}(w')}\vert W(M)\vert ^{-1}\vert Z_{W(M)}(w)\vert trace(\rho(w)^{-1}).$$
 Pour  $w'\in W_{t}$ et $w\in Cl_{W}(w')$, on a  les égalités $trace(\rho(w)^{-1})=trace(\rho_{t}(w')^{-1})$ et $\vert W(M)\vert ^{-1}\vert Z_{W(M)}(w)\vert=\vert Cl_{W(M)}(w')\vert ^{-1}$. La derni\`ere somme de la formule ci-dessus vaut donc $trace(\rho_{t}(w')^{-1})$. D'o\`u
  $$\boldsymbol{\chi}_{M,{\cal E},\rho}(t \,exp(N'))=\vert W_{t}(M_{t})\vert ^{-1}\sum_{w'\in W_{t}(M_{t})}\mu_{t}(w')trace(\rho_{t}(w')^{-1})\boldsymbol{\chi}_{M_{t},{\cal E}_{t},w'}(N').$$
  En d\'ecomposant $\rho_{t}$ en composantes simples, on obtient
  $$\boldsymbol{\chi}_{M,{\cal E},\rho}(t \,exp(N'))= \sum_{\rho'\in Irr(W_{t}(M_{t}))}<\rho_{t},\rho'>\vert W_{t}(M_{t})\vert ^{-1}$$
  $$\sum_{w'\in W_{t}}\mu_{t}(w')trace(\rho'(w')^{-1})\boldsymbol{\chi}_{M_{t},{\cal E}_{t},w'}(N').$$
  Puis, en utilisant \ref{fonctionscaracteristiques}(1),
   $$\boldsymbol{\chi}_{M,{\cal E},\rho}(t \,exp(N'))= \sum_{\rho'\in Irr(W_{t}(M_{t}))}<\rho_{t},\rho'>\boldsymbol{\chi}_{M_{t},{\cal E},\rho'\mu_{t}^{-1}}.$$
     Par le  changement de variables $\rho'\mapsto \rho'\mu_{t}$, la formule ci-dessus devient celle de la proposition \ref{restriction}. $\square$

  \subsection{Param\'etrages}
   
     \subsubsection{L'orbite sp\'eciale associ\'ee \`a un triplet $(M,{\cal E},\rho)$\label{orbitespeciale}}
  Pour tout couple $(M,{\cal E})$ form\'e d'un Levi $M$ de $G$ et d'un faisceau ${\cal E}\in {\bf FC}_{u}(M_{AD})$, fixons une famille admissible $(\epsilon_{x})_{x\in {\cal X}_{{\cal E}}}$.  
  
    On note ${\cal A}(G)$ l'ensemble des classes de conjugaison par $G$ de triplets $(M,{\cal E},\rho)$ o\`u $M$ est un Levi de $G$, ${\cal E}\in {\bf FC}_{u}(M_{AD})$ et $\rho\in Irr(W(M))$. Ainsi qu'on l'a d\'ej\`a fait, nous noterons une classe de conjugaison par un triplet la repr\'esentant.  On a d\'efini  le faisceau-caract\`ere $K_{M,{\cal E},\rho}$. 
  
  Lusztig a d\'efini la notion d'orbite unipotente sp\'eciale. Notons ${\cal U}_{sp}$ l'ensemble de ces orbites. 
  A tout  faisceau-caract\`ere unipotent sur $G$,  il associe  une telle orbite (c'est une cons\'equence de \cite{L5} corollaire 16.7). Pour  ${\cal O}\in {\cal U}_{sp}$, notons ${\cal A}(G,{\cal O})$ le sous-ensemble des $(M,{\cal E},\rho)$ tels que l'orbite associ\'ee \`a $K_{M,{\cal E},\rho}$ soit ${\cal O}$.

  \begin{thm}{ Soient ${\cal O}\in {\cal U}_{sp}$, $(M,{\cal E},\rho)\in {\cal A}(G,{\cal O})$ et $x=tu\in G$, $t$ \'etant la partie semi-simple de $x$ et $u$ sa partie unipotente. Notons ${\cal O}_{u}$ l'orbite de $u$ dans $G$. Supposons que la fibre de $K_{M,{\cal E},\rho}$ en $x$ soit non nulle; alors $dim({\cal O}_{u})< dim({\cal O})$ ou ${\cal O}_{u}={\cal O}$.}\end{thm}
  
  Cf. \cite{L4} th\'eor\`eme 10.7.

  \subsubsection{Calcul de l'orbite sp\'eciale associ\'ee \`a un triplet $(M,{\cal E},\rho)$\label{calculdorbites}}
  
  On suppose ici que $G$ est d\'eploy\'e sur ${\mathbb F}_{q}$. Notons $Unip(G)$ l'ensemble des classes d'isomorphismes de repr\'esentations unipotentes irr\'eductibles de $G({\mathbb F}_{q})$. Notons $Unip_{cusp}(G)$ le sous-ensemble des repr\'esentations  dans $Unip(G)$ qui sont cuspidales.  Remarquons que l'application image r\'eciproque par l'application $\pi:G({\mathbb F}_{q})\to G_{AD}({\mathbb F}_{q})$ est une bijection de $Unip(G_{AD})$ sur $Unip(G)$, resp. $Unip_{cusp}(G_{AD})$ sur $Unip_{cusp}(G)$. 
   Notons ${\cal A}Unip(G)$ l'ensemble des classes de conjugaison par $G$ de triplets $(M,\pi_{cusp},\rho)$, o\`u $M$ est un ${\mathbb F}_{q}$-Levi de $G$, $\pi_{cusp}\in Unip_{cusp}(M_{AD})$ et $\rho\in Irr(W(M))$.    Soit $(M,\pi_{cusp},\rho)\in {\cal A}Unip(G)$. La repr\'esentation $\pi_{cusp}$ se rel\`eve en une repr\'esentation de $M({\mathbb F}_{q})$ encore not\'ee $\pi_{cusp}$.  Soit $P$ un sous-groupe parabolique de $G$ de composante de Levi $M$. On introduit la repr\'esentation induite $Ind_{P}^G(\pi_{cusp})$ de $G({\mathbb F}_{q})$. Son  alg\`ebre d'entrelacements est canoniquement isomorphe \`a l'alg\`ebre du groupe $W(M)$. La repr\'esentation $Ind_{P}^G(\pi_{cusp})$ se d\'ecompose selon les repr\'esentations irr\'eductibles de ce groupe. Ainsi, \`a la repr\'esentation $\rho$, on associe une composante irr\'eductible $\pi_{M,\pi_{cusp},\rho}$ de la repr\'esentation induite. C'est une repr\'esentation unipotente de $G({\mathbb F}_{q})$. On sait que l'application $(M,\pi_{cusp},\rho)\mapsto \pi_{M,\pi_{cusp},\rho}$ est une bijection de ${\cal A}Unip(G)$ sur $Unip(G)$. 
  
  A toute repr\'esentation $\pi\in Unip(G)$, Lusztig a associ\'e une orbite unipotente sp\'eciale. D'o\`u par composition une application ${\cal A}Unip(G)\to {\cal U}_{sp}$. Pour toute ${\cal O}\in {\cal U}_{sp}$, on note ${\cal A}Unip(G,{\cal O})$ la fibre de cette application au-dessus de ${\cal O}$. On sait d\'ecrire ces ensembles ${\cal A}Unip(G,{\cal O})$. Pour cela, on se ram\`ene imm\'ediatement au cas o\`u $G$ est adjoint et absolument simple. Dans le cas ces groupes classiques, ces ensembles  se d\'ecrivent \`a l'aide de "symboles", cf. \cite{C} 13.8. Dans le cas des groupes exceptionnels, ils sont d\'ecrits dans les tables de \cite{C}, p. 478-488.
  
  {\bf Remarque.} (1) Dans ces tables, Curtis regroupe les repr\'esentations selon leur orbite associ\'ee mais ne d\'esigne pas cette orbite. On r\'ecup\`ere cette orbite de la fa\c{c}on suivante. Dans la description par Curtis d'un ensemble ${\cal A}Unip(G,{\cal O})$, le premier terme est de la forme $(T,{\bf 1},\rho)$, o\`u $T$ est un tore d\'eploy\'e maximal, ${\bf 1}$ est l'unique repr\'esentation de $T_{AD}({\mathbb F}_{q})=\{1\}$ et $\rho\in Irr(W)$. La repr\'esentation $\rho$ est sp\'eciale. Notons ${\cal E}$ l'unique \'el\'ement de ${\bf FC}^T(\mathfrak{t}_{SC})$. Le triplet $(T,{\cal E},\rho)$ appartient \`a l'ensemble ${\cal I}^{FC}(G)$. Son image par la repr\'esentation de Springer est le couple $({\cal O},{\cal L}_{1})$, o\`u ${\cal L}_{1}$ est le syst\`eme local trivial sur ${\cal O}$. La repr\'esentation de Springer \'etant d\'ecrite dans les tables de \cite{C}, p. 427-433, cela permet le calcul de ${\cal O}$. 
  \bigskip

 En \cite{L10} th\'eor\`emes 3.3 et 3.7, Lusztig d\'efinit une bijection $\tau_{cusp}:Unip_{cusp}(G)\to {\bf FC}_{u}(G)$. Remarquons que, puisque $G$ est d\'eploy\'e, il n'y a pas de diff\'erence entre les classes de conjugaison de Levi et de ${\mathbb F}_{q}$-Levi. On d\'efinit une bijection $\tau:{\cal A}Unip(G)\to {\cal A}(G)$ en envoyant un triplet $(M,\pi_{cusp},\rho)$ sur $(M,\tau_{cusp}^{M_{AD}}(\pi_{cusp}),\rho)$. Cela \'etant, le corollaire 3.9 et les assertions 3.10 de \cite{L10} affirment la propri\'et\'e suivante. Pour tout couple $(M,{\cal E})$ form\'e d'un Levi $M$ de $G$ et d'un faisceau ${\cal E}\in {\bf FC}_{u}(M_{AD})$, il existe une unique famille admissible $(\epsilon_{x})_{x\in {\cal X}_{{\cal E}}}$ de sorte qu'en utilisant celle-ci pour construire les faisceaux-caract\`eres $K_{M,{\cal E},\rho}$, le diagramme suivant soit commutatif:
 $$\begin{array}{ccccc}{\cal A}Unip(G)&&\stackrel{\tau}{\to}&&{\cal A}(G)\\ &\searrow&&\swarrow&\\ &&{\cal U}_{sp}&&\\ \end{array}$$
 D\'esormais, nous utiliserons exclusivement ces familles admissibles.  Puisqu'on a dit ci-dessus que l'application ${\cal A}Unip(G)\to {\cal U}_{sp}$ \'etait connue, il en est de m\^eme de l'application ${\cal A}(G)\to {\cal U}_{sp}$.
 
 {\bf Remarques.} (2) Dans \cite{L10}, Lusztig ne dit pas comment il caract\'erise l'unique famille admissible. Mais, comme il a bien  voulu me l'indiquer, la d\'efinition est donn\'ee dans   \cite{L9} 1.4. On peut la traduire ainsi. Soit $(M,{\cal E})$ comme ci-dessus. Notons $\dot{{\cal E}}$ le faisceau sur $M$ image r\'eciproque de ${\cal E}$ par la projection $\pi^M:M\to M_{AD}$. Autrement dit $\dot{{\cal E}}=K^M_{M,{\cal E},1}$, o\`u $1$ est l'unique repr\'esentation de $W^M(M)=\{1\}$. Soit ${\cal O}^M$ l'image de $\dot{{\cal E}}$ dans ${\cal U}_{sp}^M$. Notons $\underline{{\cal O}}^G $ l'orbite unipotente de $G$ induite  de ${\cal O}^M$. Quelle que soit la famille admissible utilis\'ee, il existe une unique repr\'esentation $\rho\in W(M)$ telle que l'orbite associ\'ee \`a $K_{M,{\cal E},\rho}$ soit $\underline{{\cal O}}^G$ et cette repr\'esentation est de dimension $1$. Alors la famille canonique est celle pour laquelle cette repr\'esentation $\rho$ est triviale. 
 
 (3) Les assertions ci-dessus r\'esultent aussi de la preuve par Shoji de la conjecture de Lusztig, cf. \cite{Shoji} et \cite{Shoji2}. En effet, modulo un bon  choix de familles admissibles (il est vrai un peu confus), Shoji montre que, pour $(M,\pi_{cusp},\rho)\in {\cal A}Unip(G)$, dont on note $(M,{\cal E},\rho)$ l'image par $\tau$, la fonction caract\'eristique de $K_{M,{\cal E},\rho}$ est \'egale \`a une constante pr\`es au caract\`ere "fant\^ome" asoci\'e \`a $(M,\pi_{cusp},\rho)$. On peut en d\'eduire que les orbites sp\'eciales associ\'ees aux deux triplets sont \'egales. 
 
 (4) On a suppos\'e que $G$ \'etait d\'eploy\'e sur ${\mathbb F}_{q}$. Mais, pour calculer l'application ${\cal A}(G)\to {\cal U}_{sp}$, qui est ind\'ependante du corps de base, on peut \'etendre ce corps de base en un corps ${\mathbb F}_{q^n}$ sur lequel $G$ est d\'eploy\'e. 
 \bigskip
 
 Consid\'erons le cas d'un faisceau-caract\`ere unipotent cuspidal ${\cal E}$ de $G_{AD}$. Comme ci-dessus, notons $\dot{{\cal E}}$ son image r\'eciproque sur $G$. Soit $x=tu$ un \'el\'ement de $G$ au-dessus duquel la fibre de $\dot{{\cal E}}$ est non nulle.
 Alors
 
 (5) l'orbite unipotente sp\'eciale associ\'ee \`a $\dot{{\cal E}}$ est l'orbite de $u$. 
 
 Preuve. Cela r\'esulte imm\'ediatement du th\'eor\`eme \ref{parametrage3} ci-dessous mais on peut donner une preuve d\`es maintenant. On peut supposer $G$ adjoint et absolument simple. Notons ${\cal O}'_{{\cal E}}$ l'orbite de $u$ et ${\cal O}_{{\cal E}}$ l'orbite sp\'eciale associ\'ee \`a $\dot{{\cal E}}$. Enum\'erons tous les faisceaux-caract\`eres unipotents cuspidaux de $G_{AD}$: ${\cal E}_{1}$,...,${\cal E}_{n}$. Ils sont d\'ecrits explicitement, on conna\^{\i}t donc les orbites ${\cal O}'_{{\cal E}_{1}}$,...,${\cal O}'_{{\cal E}_{n}}$. On est rest\'e \'evasifs sur la bijection $\tau_{cusp}$ mais son existence implique que la famille ${\cal O}_{{\cal E}_{1}}$,...,${\cal O}_{{\cal E}_{n}}$ est \'egale, \`a permutation pr\`es, \`a la famille des orbites sp\'eciales associ\'ees aux \'el\'ements de $Unip_{cusp}(G)$. Cette famille est elle-aussi connue. On constate alors cas par cas que les familles ${\cal O}'_{{\cal E}_{1}}$,...,${\cal O}'_{{\cal E}_{n}}$ et ${\cal O}_{{\cal E}_{1}}$,...,${\cal O}_{{\cal E}_{n}}$ sont \'egales \`a permutation pr\`es. Le th\'eor\`eme \ref{orbitespeciale} implique que $dim({\cal O}'_{{\cal E}_{i}})<dim({\cal O}_{{\cal E}_{i}})$ ou ${\cal O}'_{{\cal E}_{i}}={\cal O}_{{\cal E}_{i}}$ pour tout $i=1,...,n$. De ces deux propri\'et\'es r\'esulte imm\'ediatement l'\'egalit\'e ${\cal O}'_{{\cal E}_{i}}={\cal O}_{{\cal E}_{i}}$ pour tout $i$.

 \subsubsection{Une propri\'et\'e des familles admissibles canoniques\label{normalisation}}
 Soient $ M$ un Levi de $G$ et ${\cal E}\in {\bf FC}_{u}(M_{AD})$. On simplifie les notations en posant  ${\cal X}={\cal X}_{{\cal E}}$ etc....   Notons $S_{ss}$ l'ensemble des \'el\'ements de $M$ dont l'image dans $M_{AD}$ est la partie semi-simple d'un \'el\'ement de ${\cal X}$. Comme on l'a dit, on utilise la famille admissible canonique dans toutes les constructions pr\'ec\'edentes. 
 
 \begin{prop}{ 
 (i) Pour tout $t\in S_{ss}$, le caract\`ere $\mu_{t}$ de \ref{mut} est \'egal \`a $1$.
 
 (ii) Supposons  que  $M$ soit  un ${\mathbb F}_{q}$-Levi et que ${\cal E}$ soit conserv\'e par l'action de $\Gamma_{{\mathbb F}_{q}}$. Alors la famille admissible canonique  est d\'efinie sur ${\mathbb F}_{q}$. }\end{prop}
 
 Preuve. Pour la premi\`ere assertion, on se ram\`ene ais\'ement au cas o\`u $G$ est adjoint et absolument simple. On peut supposer que $M$ est un Levi propre et que $W_{t}(M_{t})\not=\{1\}$ sinon l'assertion est vide. Notons $\dot{{\cal E}}$ l'image r\'eciproque de ${\cal E}$ dans $M$, soit ${\cal O}^M\in {\cal U}_{sp}^M$ l'orbite unipotente dans $M$ associ\'ee \`a $\dot{{\cal E}}$. Notons ${\cal O}^G$ l'orbite de $G$ engendr\'ee par ${\cal O}^M$ (et non pas l'orbite induite). En g\'en\'eral, l'orbite engendr\'ee par une orbite sp\'eciale de $M$ ne reste pas sp\'eciale dans $G$. Mais, dans notre cas o\`u ${\cal O}^{M}$ est associ\'ee \`a un faisceau-caract\`ere cuspidal, on v\'erifie cas par cas que ${\cal O}^G$ est sp\'eciale.  
 Identifions par l'exponentielle les orbites  unipotentes dans un groupe \`a des orbites nilpotentes dans son alg\`ebre de Lie.  Par d\'efinition de $S_{ss}$, $\mathfrak{m}_{t}$ contient un \'el\'ement $N$ tel que le faisceau $\dot{{\cal E}}$ soit non nul en $t\,exp(N)$. Notons ${\cal O}_{t}^M$ l'orbite de $N$ dans $\mathfrak{m}_{t}$. D'apr\`es \ref{calculdorbites}(5), l'orbite ${\cal O}^M$ (descendue en une orbite nilpotente) contient $N$, autrement dit, ${\cal O}^M$ est l'image de ${\cal O}_{t}^M$ par l'application naturelle $\mathfrak{m}_{t,nil}/conj\to \mathfrak{m}_{nil}/conj$. Notons ${\cal O}_{t}^G$ l'orbite de $\mathfrak{g}_{t}$ engendr\'ee par ${\cal O}_{t}^M$. Alors ${\cal O}^G$ est l'image de ${\cal O}_{t}^G$ par l'application naturelle $\mathfrak{g}_{t,nil}/conj\to \mathfrak{g}_{nil}/conj$. Pour $\rho'\in W_{t}(M_{t})$, notons $({\cal O}_{t,\rho'},{\cal L}_{t,\rho'})\in {\cal I}(G_{t})$ l'image de $(M_{t},{\cal E}_{t},\rho')$ par la correspondance de Springer g\'en\'eralis\'ee. Notons ${\cal O}_{\rho'}$ l'orbite nilpotente de $\mathfrak{g}$ engendr\'ee par ${\cal O}_{t,\rho'}$.  Une propri\'et\'e g\'en\'erale de la correspondance de Springer g\'en\'eralis\'ee est que ${\cal O}_{t}^G\leq{\cal O}_{t,\rho'}$ pour l'ordre usuel des orbites, avec \'egalit\'e si et seulement si $\rho'$ est le caract\`ere signature $sgn_{t}$ de $W_{t}(M_{t})$. Puisque l'application $\mathfrak{g}_{t,nil}/conj\to \mathfrak{g}_{nil}/conj$ est strictement croissante, on en d\'eduit ${\cal O}^G\leq{\cal O}_{\rho'}$ avec \'egalit\'e si et seulement si $\rho'=sgn_{t}$. 
 
  Notons $sgn$ le caract\`ere signature de $W(M)$.
 D'apr\`es \ref{calculdorbites}, on sait calculer l'orbite sp\'eciale associ\'ee au triplet $(M,{\cal E},sgn)$. On constate cas par cas que c'est ${\cal O}^G$. Appliquons la proposition \ref{restriction} \`a $(M,{\cal E},sgn)$. Quitte \`a \'etendre le corps de base ${\mathbb F}_{q}$ en un corps ${\mathbb F}_{q^n}$, on peut supposer v\'erifi\'ees les hypoth\`eses de cette proposition: ces hypoth\`eses se r\'esumaient \`a l'invariance par l'action de Frobenius de divers objets, ce que l'on peut assurer par changement de base.  La restriction \`a $W_{t}(M_{t})$ du caract\`ere $sgn$ est le caract\`ere $sgn_{t}$.  Il y a une seule repr\'esentation $\rho'\in Irr(W_{t}(M_{t}))$ telle que $<sgn_{t},\mu_{t}\rho'>\not=0$, \`a savoir $\rho'=\mu_{t}sgn_{t}$ (remarquons que $\mu_{t}=\mu_{t}^{-1}$ ainsi que c'est le cas pour tout caractère de $W(M)$). Puisque $\boldsymbol{\chi}_{M_{t},{\cal E}_{t},\mu_{t}sgn_{t}}$ est non nulle sur ${\cal O}_{t,\mu_{t}sgn_{t}}$, la proposition \ref{restriction}  implique que la fonction $N'\mapsto \boldsymbol{\chi}_{M,{\cal E},\rho}(t\,exp(N'))$ est non nulle sur ${\cal O}_{\mu_{t}sgn_{t}}$. Puisque l'orbite sp\'eciale associ\'ee \`a $(M,{\cal E},sgn)$ est ${\cal O}^G$, le th\'eor\`eme \ref{orbitespeciale} implique que $dim({\cal O}_{\mu_{t}sgn_{t}})\leq dim({\cal O}^G)$. D'apr\`es les rappels ci-dessus, cela entra\^{\i}ne $\mu_{t}=1$. 
   Cela prouve le (i) de l'\'enonc\'e. 
 
  Remarquons que la propri\'et\'e (i) caract\'erise la famille admissible canonique. En effet, deux familles admissibles diff\`erent par un caract\`ere de $W(M)$. Si ces deux familles v\'erifient (i), ce caract\`ere est trivial sur l'image de l'homomorphisme $\iota_{t}:W_{t}(M_{t})\to W(M)$, cela pour tout $t\in S_{ss}$.  Puisque ces images     engendrent le groupe $W(M)$, cf. \ref{mut}(2),  le caract\`ere est trivial et les deux familles sont \'egales.

 Prouvons l'assertion (ii). On suppose donc que $M$ est un ${\mathbb F}_{q}$-Levi  et que ${\cal E}$ est conserv\'e par l'action de $\Gamma_{{\mathbb F}_{q}}$.  On consid\`ere la famille admissible canonique $(\epsilon_{x})_{x\in {\cal X}}$, qui v\'erifie  (i). On a d\'efini en \ref{faisceauxcaracteres} la famille admissible $(^{Fr}\epsilon_{x})_{x\in {\cal X}}$. Il s'agit de d\'emontrer que ces deux familles sont \'egales. D'apr\`es l'unicit\'e affirm\'ee  ci-dessus, il suffit de v\'erifier que la famille $(^{Fr}\epsilon_{x})_{x\in {\cal X}}$ v\'erifie encore la propri\'et\'e (i). Soit $t\in S_{ss}$. Choisissons un \'el\'ement unipotent $u\in M_{t}$ tel que $t_{ad}u\in {\cal X}$ et d\'efinissons $N\in \mathfrak{m}_{t,nil}$ par $exp(N)=u$. Avec les notations de \ref{mut}, on  doit prouver que $^{Fr}\epsilon_{t_{ad}u}(n)=\epsilon^{\flat}_{N}(n)$ pour tout $n\in Norm_{G_{t}}(M_{t},N)$. On a $^{Fr}\epsilon(t_{ad}u)(n)=\epsilon_{Fr(t_{ad}u)}(Fr(n))$. L'\'el\'ement $Fr(t)$  appartient \`a $S_{ss}$. En utilisant (i) pour cet \'el\'ement, on a $\epsilon_{Fr(t_{ad}u)}(Fr(n))=\epsilon^{\flat}_{Fr(N)}(Fr(n))$, o\`u $\epsilon^{\flat}_{Fr(N)}$ est le caract\`ere canonique de 
 $Norm_{G_{Fr(t)}}(M_{Fr(t)},Fr(N))$. D'après \ref{faisceauxcaracteres}(9), on a  l'\'egalit\'e $\epsilon^{\flat}_{Fr(N)}(Fr(n))=\epsilon^{\flat}_{N}(n)$. D'o\`u l'\'egalit\'e cherch\'ee $^{Fr}\epsilon_{t_{ad}u}(n)=\epsilon^{\flat}_{N}(n)$, ce qui d\'emontre la seconde assertion de la proposition. $\square$
 
   \subsubsection{Majoration d'orbites\label{premiertheoreme}}
 Pour la fin de la section 5, on suppose $G$ adjoint et absolument simple.  Soient ${\cal O}\in {\cal U}_{sp}$ et $(M,{\cal E},\rho)\in {\cal A}(G,{\cal O})$   On note ${\cal X}={\cal X}_{{\cal E}}$ et $S_{ss}$ l'ensemble des \'el\'ements de $M$ dont l'image dans $M_{AD}$ est la partie semi-simple d'un \'el\'ement de ${\cal X}$.      Soit $t\in S_{ss}$. On d\'efinit   la repr\'esentation $\rho_{t}=\rho\circ\iota_{t}$ de $W_{t}(M_{t})$ comme en   \ref{restriction}.  Soit $\rho'\in Irr(W_{t}(M_{t}))$. Par la correspondance de Springer g\'en\'eralis\'ee, le triplet $(M_{t},{\cal E}_{t},\rho')$ correspond \`a un couple $({\cal O}_{t,\rho'},{\cal L}_{t,\rho'})\in {\cal I}(G_{t})$. On note ${\cal O}_{M_{t},{\cal E}_{t},\rho'}$ ou simplement ${\cal O}_{\rho'}$ l'orbite nilpotente de $\mathfrak{g}$ engendr\'ee par ${\cal O}_{t,\rho'}$. On identifie par l'exponentielle les orbites nilpotentes de $\mathfrak{g}$ et les orbites unipotentes de $G$. 

\begin{prop}{Supposons $<\rho_{t},\rho'>>0$. Alors $dim({\cal O}_{\rho'})< dim({\cal O})$ ou ${\cal O}_{\rho'}={\cal O}$.}\end{prop}

Preuve. Comme on l'a d\'ej\`a dit, quitte \`a remplacer le corps de base ${\mathbb F}_{q}$ par ${\mathbb F}_{q^n}$ pour un entier $n$ assez grand, on peut supposer v\'erifi\'ees les hypoth\`eses de \ref{restriction}.  Notons ${\cal R}$ l'ensemble des $\rho''\in Irr(W_{t}(M_{t}))$ telles que $<\rho_{t},\rho''>>0$ et munissons-le de l'ordre d\'efini par $\rho''_{1}\leq \rho''_{2}$ si et seulement si ${\cal O}_{t,\rho''_{1}}$ est incluse dans l'adh\'erence de Zariski de ${\cal O}_{t,\rho''_{2}}$. Remarquons que cette condition implique que ${\cal O}_{\rho''_{1}}$  est incluse dans l'adh\'erence de Zariski de ${\cal O}_{\rho''_{2}}$. Soit $\rho''$ un \'el\'ement maximal  de ${\cal R}$ tel que $\rho'\leq \rho''$. On voit que la conclusion du th\'eor\`eme pour $\rho''$ implique la m\^eme conclusion pour $\rho'$. On peut donc supposer que $\rho'=\rho''$, autrement dit que $\rho'$ est maximale. Notons ${\cal R}'$ le sous-ensemble des $\rho''\in {\cal R}$ telles que ${\cal O}_{t,\rho''}={\cal O}_{t,\rho'}$. Pour $\rho''\in {\cal R}-{\cal R}'$, la restriction \`a ${\cal O}_{t,\rho'}$ de la fonction $\boldsymbol{\chi}_{M_{t},{\cal E}_{t},\rho''}$ est nulle. Pour $\rho''\in {\cal R}'$, cette restriction est proportionnelle \`a ${\cal Y}_{{\cal O}_{t,\rho'},{\cal L}_{t,\rho''}}$. Quand $\rho''$ parcourt ${\cal R}'$, ces fonctions sont lin\'eairement ind\'ependantes. Il en r\'esulte que la fonction
$$\sum_{\rho''\in Irr(W_{t}(M_{t}))}<\rho_{t},\rho''>\boldsymbol{\chi}_{M_{t},{\cal E}_{t},\rho''}$$
est non nulle sur ${\cal O}_{t,\rho'}$. La proposition \ref{restriction} (o\`u $\mu_{t}=1$ d'apr\`es la proposition \ref{normalisation}) entra\^{\i}ne que la fonction $N'\mapsto \boldsymbol{\chi}_{M,{\cal E},\rho}(t\,exp(N'))$ n'est pas nulle sur ${\cal O}_{t,\rho'}$. Le th\'eor\`eme \ref{orbitespeciale} entra\^{\i}ne alors la conclusion de la proposition. $\square$

\subsubsection{Param\'etrage de ${\cal A}(G,{\cal O})$\label{parametrage3}}

      Consid\'erons les triplets $(N,d, \mu)$ o\`u $N\in \mathfrak{g}_{nil}$,   $d\in \bar{A}(N)$ et $\mu$ est une repr\'esentation irr\'eductible de $Z_{\bar{A}(N)}(d)$. Le groupe $G$ agit par conjugaison sur ces triplets. On note  $\bar{{\cal C}}^{Irr}$ l'ensemble des classes de conjugaison.  Comme toujours, nous noterons une classe de conjugaison par un triplet la représentant.
 Pour une orbite ${\cal O}\in \mathfrak{g}_{nil}/conj$, on note  $\bar{{\cal C}}^{Irr}({\cal O})$ le sous-ensemble des classes de conjugaison   de  triplets $(N,d,\mu)$  où $N\in {\cal O}$.
    
    Soit $(N,d,\mu)\in \bar{{\cal C}}^{Irr}$.  Soit $t$ un élément semi-simple de $Z_{G}(N)$ dont l'image dans $\bar{A}(N)$ est \'egale \`a $d$. On a la suite d'homomorphismes
  $$Z_{G_{t}}(N)/Z_{G_{t}}(N)^0\to Z_{G}(N)/Z_{G}(N)^0=A(N)\to \bar{A}(N),$$
  et l'image de $Z_{G_{t}}(N)/Z_{G_{t}}(N)^0$ dans $\bar{A}(N)$ est contenue dans $Z_{\bar{A}(N)}(d)$. 
  Ainsi $\mu$ se rel\`eve en une repr\'esentation de $Z_{G_{t}}(N)/Z_{G_{t}}(N)^0$. Cette repr\'esentation n'est pas forc\'ement irr\'eductible puisque l'homomorphisme $Z_{G_{t}}(N)/Z_{G_{t}}(N)^0\to Z_{\bar{A}(N)}(d)$ n'est pas forc\'ement surjectif. Posons $u=exp(N)$ et notons
  ${\cal X}_{t,u}$ la classe de conjugaison par $G_{t}$ de $tu$. De la repr\'esentation pr\'ec\'edente  se d\'eduit un syst\`eme local  pas forc\'ement irr\'eductible sur ${\cal X}_{t,u}$ que l'on note ${\cal E}_{t,\mu}$. 
     
    Disons que ${\cal O}$ est exceptionnelle dans les cas suivants:
  
  $G$ est de type $E_{7}$ et ${\cal O}$ est de type $A_{4}+A_{1}$;
  
  $G$ est de type $E_{8}$ et ${\cal O}$ est de type $A_{4}+A_{1}$ ou $E_{6}(a_{1})+A_{1}$. 
  
 \begin{thm}{ Soit ${\cal O}\in {\cal U}_{sp}$. Supposons que    ${\cal O}$ ne soit pas exceptionnelle. Il existe une unique application $\nabla_{{\cal O}}:{\cal A}(G,{\cal O})\to \bar{{\cal C}}^{Irr}({\cal O})$ v\'erifiant les conditions suivantes. Soit $(M,{\cal E},\rho)\in {\cal A}(G,{\cal O})$, posons $\nabla_{{\cal O}}(M,{\cal E},\rho)=(N,d,\mu)$. Soit $t$ un élément semi-simple appartenant à $Z_{G}(N)$, notons $d_{t}$ son image dans $\bar{A}(N)$. On a
 
 (i) si $d_{t}$ n'est pas conjugu\'e \`a $d$ par un \'el\'ement de $\bar{A}(N)$, la restriction \`a ${\cal X}_{t,u}$ de $K_{M,{\cal E},\rho}$ est nulle;
 
 (ii) supposons $d_{t}=d$; alors la restriction \`a ${\cal X}_{t,u}$ de $K_{M,{\cal E},\rho}$ est \'egale au syst\`eme local ${\cal E}_{t,\mu}$ \`a un d\'ecalage pr\`es.
 
 L'application $\nabla_{{\cal O}}$ est bijective. }\end{thm}
 
 Cf. \cite{L8} th\'eor\`eme 2.4.
 
 {\bf Remarque.}  On doit signaler que Lusztig ne donne pas de vraie d\'emonstration de ce théorème, il se contente de donner quelques indications. De plus, 
  l'\'enonc\'e de Lusztig est erron\'e. Il oublie l'hypoth\`ese que ${\cal O}$ n'est pas exceptionnelle (l'\'enonc\'e est faux pour les orbites exceptionnelles). M\^eme dans le cas d'une orbite qui n'est pas exceptionnelle, Lusztig  \'ecrit un \'enonc\'e plus fort o\`u les groupes $G_{t}$ sont remplac\'es par $Z_{G}(t)$.   Il ne nous semble pas que la d\'emonstration esquiss\'ee par Lusztig permette d'atteindre cet \'enonc\'e plus fort. Dans le cas des groupes classiques, l'\'enonc\'e ci-dessus est une cons\'equence de la proposition XI.29 de \cite{W1}. Nous n'expliquerons pas ici comment il se d\'eduit de ce r\'esultat car cela deviendra clair en \ref{preuveDnram} et \ref{preuveAn-1ram} lorsque nous reprendrons les m\'ethodes de \cite{W1}  pour d\'emontrer un \'enonc\'e analogue \`a celui ci-dessus dans le cas des groupes classiques ramifi\'es. Pour les groupes exceptionnels, la m\'ethode indiqu\'ee par Lusztig ram\`ene \`a des calculs de restrictions de repr\'esentations de groupes de Weyl et de correspondances de Springer g\'en\'eralis\'ees.  On dispose de tables pour ces restrictions et  les calculs  peuvent certainement \^etre traités  complètement.  
  
  \bigskip
  
   Supposons que ${\cal O}$ soit exceptionnelle. On d\'efinit encore une bijection $\nabla_{{\cal O}}:{\cal A}(G,{\cal O})\to \bar{{\cal C}}^{Irr}({\cal O})$ par les formules imitant la classification des repr\'esentations unipotentes de $G({\mathbb F}_{q})$, que l'on va rappeler.  On fixe $N\in {\cal O}$. Alors $\bar{A}(N)\simeq {\mathbb Z}/2{\mathbb Z}$. On note $1,d_{2}$ les deux \'el\'ements de ce groupe et $1,sgn$ ses deux caract\`eres. 
   
   Supposons que $G$ soit de type $E_{7}$ et que ${\cal O}$ soit de type $A_{4}+A_{1}$.  L'ensemble ${\cal A}(G,{\cal O})$ a quatre \'el\'ements. Deux sont  associ\'es \`a des triplets $(M,{\cal E},\rho)$ pour lesquels $M$ est un tore, ${\cal E}$ est le syst\`eme trivial (qui est l'unique syst\`eme unipotent cuspidal pour un tore) et $\rho$ est l'une des deux repr\'esentations $\rho_{512,11}$ et $\rho_{512,12}$ de $W=W(M)$ (avec les notations de \cite{C} p. 483 et, ci-dessous, p. 485 et 486). On note simplement  ces triplets $\rho_{512,11}$ et $\rho_{512,12}$. Deux autres \'el\'ements sont associ\'es \`a des triplets $(M,{\cal E},\rho)$ o\`u $M=G$, ${\cal E}$ est l'un des deux \'el\'ements de ${\bf FC}_{u}(G)$ et $\rho$ est la repr\'esentation triviale de $W(M)=\{1\}$. Nous noterons ces deux \'el\'ements $E_{7}[i]$ et $E_{7}[-i]$, sans nous soucier de savoir  exactement \`a quel faisceau-caract\`ere cuspidal ils correspondent. Alors $\nabla_{{\cal O}}$ est d\'efinie par les formules:
   $$\nabla_{{\cal O}}(\rho_{512,11})=(N,1,1 ),\,\,\nabla_{{\cal O}}(\rho_{512,12})=(N,1,sgn),$$
   $$\nabla_{{\cal O}}(E_{7}[i])=(N,d_{2},1),\,\, \nabla_{{\cal O}}(E_{7}[-i])=(N,d_{2},sgn).$$

 Supposons que $G$ soit de type $E_{8}$ et que ${\cal O}$ soit de type $A_{4}+A_{1}$. L'ensemble ${\cal A}(G,{\cal O})$ a quatre \'el\'ements. Deux sont  associ\'es \`a des triplets $(M,{\cal E},\rho)$ pour lesquels $M$ est un tore, ${\cal E}$ est le syst\`eme trivial  et $\rho$ est l'une des deux repr\'esentations $\rho_{4096,26}$ et $\rho_{4096,27}$ de $W=W(M)$. On note simplement  ces triplets $\rho_{4096,26}$ et $\rho_{4096,27}$. Deux autres \'el\'ements sont associ\'es \`a des triplets $(M,{\cal E},\rho)$ o\`u $M_{AD}$ est de type $E_{7}$, ${\cal E}$ est l'un des deux \'el\'ements de ${\bf FC}_{u}(M_{AD})$ et $\rho$ est la repr\'esentation non triviale de $W(M)={\mathbb Z}/2{\mathbb Z}$. Nous noterons ces deux \'el\'ements $(E_{7}[i],sgn)$ et $(E_{7}[-i],sgn)$. L'application  $\nabla_{{\cal O}}$ est d\'efinie par les formules:
   $$\nabla_{{\cal O}}(\rho_{4096,26})=(N,1,1 ),\,\,\nabla_{{\cal O}}(\rho_{4096,27})=(N,1,sgn),$$
   $$\nabla_{{\cal O}}(E_{7}[i],sgn)=(N,d_{2},1),\,\, \nabla_{{\cal O}}(E_{7}[-i],sgn)=(N,d_{2},sgn).$$

  Supposons que $G$ soit de type $E_{8}$ et que ${\cal O}$ soit de type $E_{6}(a_{1})+A_{1}$. L'ensemble ${\cal A}(G,{\cal O})$ a quatre \'el\'ements. Deux sont associ\'es \`a des triplets $(M,{\cal E},\rho)$ pour lesquels $M$ est un tore, ${\cal E}$ est le syst\`eme trivial   et $\rho$ est l'une des deux repr\'esentations $\rho_{4096,11}$ et $\rho_{4096,12}$ de $W=W(M)$. On note simplement  ces triplets $\rho_{4096,11}$ et $\rho_{4096,12}$. Deux autres \'el\'ements sont associ\'es \`a des triplets $(M,{\cal E},\rho)$ o\`u $M_{AD}$ est de type $E_{7}$, ${\cal E}$ est l'un des deux \'el\'ements de ${\bf FC}_{u}(M_{AD})$ et $\rho$ est la repr\'esentation  triviale de $W(M)={\mathbb Z}/2{\mathbb Z}$. Nous noterons ces deux \'el\'ements $(E_{7}[i],1)$ et $(E_{7}[-i],1)$. L'application 
 $\nabla_{{\cal O}}$ est d\'efinie par les formules:
   $$\nabla_{{\cal O}}(\rho_{4096,11})=(N,1,1 ),\,\,\nabla_{{\cal O}}(\rho_{4096,12})=(N,1,sgn),$$
   $$\nabla_{{\cal O}}(E_{7}[i],1)=(N,d_{2},1),\,\, \nabla_{{\cal O}}(E_{7}[-i],1)=(N,d_{2},sgn).$$
   
   On note $\bar{{\cal C}}^{Irr}_{sp}$ la réunion des $\bar{{\cal C}}^{Irr}({\cal O})$ sur les orbites ${\cal O}\in {\cal U}_{sp}$. 
   La réunion des bijections $\nabla_{{\cal O}}$ sur les ${\cal O}\in {\cal U}_{sp}$ définit une bijection notée $\nabla:{\cal A}(G)\to \bar{{\cal C}}^{Irr}_{sp}$.

   \subsubsection{Rel\`evement d'\'el\'ements de $\bar{A}(N)$\label{relevement}}
   Considérons l'ensemble des couples $(N,d)$ o\`u $N\in \mathfrak{g}_{nil}$ et  $d\in \bar{A}(N)$. Le groupe $G$ agit par conjugaison sur ces triplets. On note  $\bar{{\cal C}}$ l'ensemble des classes de conjugaison. 
   Le groupe $\Gamma_{{\mathbb F}_{q}}$ agit sur $\mathfrak{g}_{nil}/conj$ en conservant le sous-ensemble des orbites sp\'eciales. Cette action est triviale, sauf dans le cas o\`u  $G$ est de type $D_{n}$  et n'est pas d\'eploy\'e, cf. \ref{actiongaloisienne}. On note $\bar{{\cal C}}_{{\mathbb F}_{q}}$ le sous-ensemble des $(N,d)\in \bar{{\cal C}}$ tels que l'orbite de $N$ soit conservée par $\Gamma_{{\mathbb F}_{q}}$. 
   
   Notons ${\cal J}$ l'ensemble des couples $(\alpha,{\cal O}')$ où $\alpha\in \Delta_{a}$ et ${\cal O}'\in \mathfrak{g}_{t_{\alpha},nil}/conj$. 
   Pour un tel couple, soit $N'\in {\cal O}'$. On a $t_{\alpha}\in Z_{G}(N')$ et on note $d_{\alpha}$ l'image de $t_{\alpha}$ dans $\bar{A}(N')$. Alors la classe de conjugaison par $G$ du couple $(N',d_{\alpha})$ appartient à $\bar{{\cal C}}$.  Cela définit une application $J:{\cal J}\to \bar{{\cal C}}$. Avec les m\^emes notations, le groupe $Z_{G_{t_{\alpha}}}(N')/Z_{G_{t_{\alpha}}}(N')^0$ s'envoie naturellement dans $Z_{\bar{A}(N')}(d_{\alpha})$. On note ${\cal J}_{max}$ le sous-ensemble des $(\alpha,{\cal O}')\in {\cal J}$ tels que cette application soit surjective. 
   
   Soit $\alpha\in \Delta_{a}$, supposons que $\alpha$ est fix\'ee par l'action galoisienne. L'\'el\'ement $t_{\alpha}$ d\'efini en \ref{talpha} n'est pas forc\'ement fix\'e par cette action: il l'est si $\alpha=\alpha_{0}$ mais, si $\alpha\in \Delta$, on a $t_{\alpha}=\check{\varpi}_{\alpha}(\bar{\zeta}_{1/d(\alpha)})$ d'o\`u $Fr(t_{\alpha})=\check{\varpi}_{\alpha}(Fr(\bar{\zeta}_{1/d(\alpha)}))=\check{\varpi}_{\alpha}(\bar{\zeta}_{q/d(\alpha)})$. Toutefois,  on a l'\'egalit\'e $G_{Fr(t_{\alpha})}=G_{t_{\alpha}}$ donc $G_{t_{\alpha}}$ est conserv\'e par l'action galoisienne.  On note ${\cal J}_{{\mathbb F}_{q}}$, resp. ${\cal J}_{max,{\mathbb F}_{q}}$, le sous-ensemble des $(\alpha,{\cal O}')\in {\cal J}$, resp. ${\cal J}_{max}$,  tels que $\alpha$ et ${\cal O}'$ soient fixés par l'action galoisienne.

   \begin{lem}{Il existe une section $S_{J}:\bar{{\cal C}}\to {\cal J}_{max}$ de l'application $J$ qui envoie $\bar{{\cal C}}_{{\mathbb F}_{q}}$ dans ${\cal J}_{max,{\mathbb F}_{q}}$. } 
   \end{lem}
   
   Preuve. On utilise le groupe $G_{F}$, cf. \ref{legroupeG}. Soit $(N,d)\in \bar{{\cal C}}$.  Posons $\bar{R}= Z_{\bar{A}(N)}(d)$. Le triplet $(N,d,\bar{R})$ appartient \`a l'ensemble $\bar{{\cal R}}_{max}$ d\'efini en \ref{auxiliaire}. D'apr\`es le (ii) du lemme de \ref{auxiliaire}, on peut fixer $\alpha\in \Delta_{a}$ et une orbite ${\cal O}''\in \mathfrak{g}_{F,s_{\alpha},nil}/conj$ de sorte que $\bar{r}_{s_{\alpha}}({\cal O}'')=(N,d,\bar{R})$. On se rappelle que $\iota_{s_{\alpha}}(G_{F,s_{\alpha}})=G_{t_{\alpha}}$. On pose ${\cal O}'=\iota_{s_{\alpha}}({\cal O}'')$.   Il r\'esulte de la d\'efinition de l'application $\bar{r}_{s_{\alpha}}$  que 
   $(\alpha,{\cal O}')$ appartient à ${\cal J}_{max}$ et que $J(\alpha,{\cal O}')=(N,d)$ dans $\bar{{\cal C}}$. 
   
   Supposons $(N,d)\in \bar{{\cal C}}_{{\mathbb F}_{q}}$. 
      On peut conjuguer le couple $(\alpha,{\cal O}'')$ de la construction ci-dessus par un \'el\'ement de $\boldsymbol{\Omega}$. Il suffit de v\'erifier  qu'une telle conjugaison bien choisie  transforme le couple $(\alpha,{\cal O}'')$ en un  élément de ${\cal J}_{{\mathbb F}_{q}}$. C'est exactement la d\'emonstration que l'on a faite  en  \ref{surjectivite}. $\square$

 \subsubsection{S\'eparation des \'el\'ements de ${\cal A}(G)$\label{deuxiemetheoreme}}
 On a d\'efini en \ref{orbitespeciale} l'ensemble ${\cal A}(G)$ form\'e de classes de conjugaison de triplets $(M,{\cal E},\rho)$. A tout couple $(M,{\cal E})$ intervenant dans un tel triplet, on a associ\'e en \ref{Lambda} un sous-ensemble $\Lambda\subset \Delta_{a}$. Notons $\boldsymbol{\Lambda}$ l'ensemble de ces sous-ensembles $\Lambda$ images de couples $(M,{\cal E})$.  Notons ${\cal A}'(G)$ l'ensemble des triplets $(\Lambda,{\cal E},\rho)$ o\`u $\Lambda\in \boldsymbol{\Lambda}$, ${\cal E}\in {\bf FC}_{u}(M_{\Lambda,AD})$ et $\rho\in Irr(W(M_{\Lambda}))$. Cet ensemble s'identifie \`a un sous-ensemble de l'ensemble des triplets pr\'ec\'edents par l'application $(\Lambda,{\cal E},\rho)\mapsto (M_{\Lambda},{\cal E},\rho)$.  Ce sous-ensemble est un ensemble de repr\'esentants des classes de conjugaison dans l'ensemble des triplets. Donc ${\cal A}'(G)$ s'identifie \`a ${\cal A}(G)$.  Modulo cette bijection, on peut remplacer ${\cal A}(G)$ par ${\cal A}'(G)$ dans diverses d\'efinitions.  
 
 Pour une orbite ${\cal O}\in \mathfrak{g}_{nil}/conj$ on note $\bar{{\cal C}}({\cal O})$ l'ensemble des $(N,d)\in \bar{{\cal C}}$ tels que $N\in {\cal O}$. Si on fixe $N\in {\cal O}$, $\bar{{\cal C}}({\cal O})$ s'identifie à l'ensemble des classes de conjugaison dans $\bar{A}(N)$. 
   On a une application d'oubli $(N,d,\mu)\to (N,d)$ de $\bar{{\cal C}}^{Irr}$ dans $\bar{{\cal C}}$, qui est surjective. Soit ${\cal O}\in {\cal U}_{sp}$. On a une bijection $\nabla_{{\cal O}}:{\cal A}'(G,{\cal O})\to \bar{{\cal C}}^{Irr}({\cal O})$ d'où, par composition, une surjection ${\cal A}'(G,{\cal O})\to \bar{{\cal C}}({\cal O})$. Pour ${\bf d}\in \bar{{\cal C}}({\cal O})$, on note ${\cal A}'(G,{\cal O},{\bf d})$ la fibre de cette surjection au-dessus de ${\bf d}$.

 \begin{thm}{Soit ${\cal O}\in {\cal U}_{sp}$. Il existe une application injective $b:\bar{{\cal C}}({\cal O})\to {\cal J}$ telle que les propri\'et\'es suivantes soient v\'erifi\'ees.  

 (i)  $J\circ b$ est l'identité de $\bar{{\cal C}}({\cal O})$. 
 
 (ii) Si ${\cal O}$ est conserv\'ee par l'action galoisienne,  $b$ prend ses valeurs dans ${\cal J}_{{\mathbb F}_{q}}$. 
 
 Pour les deux propriétés suivantes, on fixe ${\bf d}\in \bar{{\cal C}}({\cal O})$ et on pose $b({\bf d})=(\alpha,{\cal O}')$.  
 
 (iii) Soit $(\Lambda,{\cal E},\rho)\in {\cal A}'(G,{\cal O},{\bf d})$. Alors $\alpha\not\in \Lambda$. 
 Il existe une unique repr\'esentation $\rho'\in Irr(W_{t_{\alpha}}(M_{\Lambda,t_{\alpha}}))$ telle que $<\rho_{t_{\alpha}},\rho'>>0$ et que son image  $Spr(M_{\Lambda,t_{\alpha}},{\cal E}_{t_{\alpha}},\rho')$ par la correspondance de Springer g\'en\'eralis\'ee soit port\'ee par    ${\cal O}'$. Pour cette repr\'esentation $\rho'$, on a $<\rho_{t_{\alpha}},\rho'>=1$. Posons
   $Spr(M_{\Lambda,t_{\alpha}},{\cal E}_{t_{\alpha}},\rho')=({\cal O}',{\cal L}_{\alpha;\Lambda,{\cal E},\rho})$. L'application $(\Lambda,{\cal E},\rho)\mapsto {\cal L}_{\alpha;\Lambda,{\cal E},\rho}$ d\'efinie sur  ${\cal A}'(G,{\cal O},d)$ est injective.
 
 (iv)  Soit $(\Lambda,{\cal E},\rho)\in {\cal A}'(G,{\cal O})-{\cal A}'(G,{\cal O},{\bf d})$. Supposons $\alpha\not\in \Lambda$. 
 Alors il n'existe pas de repr\'esentation $\rho'\in Irr(W_{t_{\alpha}}(M_{\Lambda,t_{\alpha}}))$ telle que $<\rho_{t_{\alpha}},\rho'>>0$ et que son image  $Spr(M_{\Lambda,t_{\alpha}},{\cal E}_{t_{\alpha}},\rho')$ par la correspondance de Springer g\'en\'eralis\'ee soit port\'ee par   ${\cal O}'$. }\end{thm}
 
   Preuve. Supposons d'abord que ${\cal O}$ n'est pas exceptionnelle. On fixe une section $S_{J}:\bar{{\cal C}}\to {\cal J}_{max}$ comme dans le lemme \ref{relevement}. On définit $b$ comme la restriction de $S_{J}$ à $\bar{{\cal C}}({\cal O})$. Les propriétés (i) et (ii) de l'énoncé résultent de celles de la section $S_{J}$. 
   
   Soit ${\bf d}\in \bar{{\cal C}}({\cal O})$, posons $b({\bf d})=(\alpha,{\cal O}')$. Fixons $N\in {\cal O}'$ et notons $d_{\alpha}$ l'image de $t_{\alpha}$ dans $\bar{A}(N)$. D'après (i), ${\bf d}$ est la classe de conjugaison de $(N,d_{\alpha})$ et $N$ appartient à ${\cal O}$. 
   
   Soit 
  $(\Lambda,{\cal E},\rho)\in {\cal A}(G,{\cal O},{\bf d})$.   On a $\nabla(\Lambda,{\cal E},\rho)=(N,d_{\alpha},\mu)$, où $\mu$ est une représentation irréductible de $Z_{\bar{A}(N)}(d_{\alpha})$. Puisque $(N,d_{\alpha})\in {\cal J}_{max}$, l'homomorphisme 
  $$Z_{G_{t_{\alpha}}}(N)/Z_{G_{t_{\alpha}}}(N)^0\to Z_{\bar{A}(N)}(d_{\alpha}) $$
   est surjectif et $\mu$ se relève en une représentation  irréductible  $\mu_{\alpha}$ de $Z_{G_{t_{\alpha}}}(N')/Z_{G_{t_{\alpha}}}(N')^0$. Cette repr\'esentation d\'etermine un syst\`eme local irr\'eductible  sur ${\cal O}'$ que l'on note ${\cal L}_{\alpha;\Lambda,{\cal E},\rho}$. Puisque $\nabla$ est injective, l'application $(\Lambda,{\cal E},\rho)\mapsto \mu$ est injective sur $ {\cal A}(G,{\cal O},{\bf d})$ et il en est de m\^eme de l'application $(\Lambda,{\cal E},\rho)\mapsto {\cal L}_{\alpha;\Lambda,{\cal O},\rho}$. Posons   ${\cal L}={\cal L}_{\alpha;\Lambda,{\cal E},\rho}$, $u=exp(N)$ et notons ${\cal X}_{t_{\alpha},u}$ l'orbite de $t_{\alpha}u$ pour l'action de $G_{t_{\alpha}}$. L'assertion (ii) du th\'eor\`eme \ref{parametrage3} et la construction disent que la restriction du complexe $K_{M_{\Lambda},{\cal E},\rho}$ \`a  ${\cal X}_{t_{\alpha},u}$ est, \`a un d\'ecalage pr\`es, le syst\`eme local d\'eduit de ${\cal L}$ par l'exponentielle. Il en r\'esulte que $t_{\alpha}u$ appartient au support de  $K_{M_{\Lambda},{\cal E},\rho}$. Par construction de ce complexe, la partie semi-simple $t_{\alpha}$ de $t_{\alpha}u$ est donc conjugu\'ee \`a un \'el\'ement de $M_{\Lambda}$ dont l'image dans $M_{\Lambda,AD}$ est la composante semi-simple d'un élément  du support de ${\cal E}$. D'apr\`es la proposition \ref{Lambda}, il existe $x\in \bar{C}^{nr}\cap a_{T}({\mathbb T})$ tel que $\Lambda\subset \Lambda(x)$ et que $t_{\alpha}$  soit conjugu\'e \`a ${\bf j}_{T}(x)$. Les groupes $G_{t_{\alpha}}$ et $G_{{\bf j}_{T}(x)}$ sont conjugu\'es. Leurs ensembles de racines ont 
  pour ensembles de racines simples $\Delta_{a}-\{\alpha\}$ et $\Lambda(x)$. Ces ensembles ont donc m\^eme nombre d'\'el\'ements et cela entra\^{\i}ne que $x=s_{\beta}$ pour une racine $\beta\in \Delta_{a}$. La condition $\Lambda\subset \Lambda(x)$ entra\^{\i}ne que $\beta\not\in \Lambda$. Puisque $t_{\alpha}$ et ${\bf j}_{T}(x)=t_{\beta}$ sont conjugu\'es, le lemme \ref{talpha} entra\^{\i}ne que $\alpha$ et $\beta$ sont conjugu\'es par l'action du groupe $\boldsymbol{\Omega}$. Or ce groupe conserve $\Lambda$, cf. lemme \ref{Lambda}(i). Donc $\alpha\not\in \Lambda$.  

Quitte \`a \'etendre le corps de base ${\mathbb F}_{q}$ on peut supposer que les hypoth\`eses de \ref{restriction} sont v\'erifi\'ees. En particulier, on dispose de la fonction caract\'eristique $\boldsymbol{\chi}_{M_{\Lambda},{\cal E},\rho}$. Puisque 
 la restriction du complexe $K_{M_{\Lambda},{\cal E},\rho}$ \`a  ${\cal X}_{t_{\alpha},u}$ est, \`a un d\'ecalage pr\`es, le syst\`eme local d\'eduit de ${\cal L}$, il existe un  nombre complexe $c$ de module $1$  tel que l'on ait l'\'egalit\'e
 
 (1) $\boldsymbol{\chi}_{M_{\Lambda},{\cal E},\rho}(t_{\alpha}exp(N'))=c{\cal Y}_{{\cal O}',{\cal L}}(N')$ pour tout $N'\in ({\cal O}')^{\Gamma_{{\mathbb F}_{q}}}$.
  
 D'apr\`es la proposition \ref{restriction}, on a aussi
 $$(2) \qquad \boldsymbol{\chi}_{M_{\Lambda},{\cal E},\rho}(t_{\alpha}exp(N'))=\sum_{\rho'\in Irr(W_{t_{\alpha}}(M_{\Lambda,t_{\alpha}}))}<\rho_{t},\rho'> \boldsymbol{\chi}_{M_{\Lambda,t_{\alpha}},{\cal E}_{t_{\alpha}},\rho'}(N').$$
 Pour $\rho'\in Irr(W_{t_{\alpha}}(M_{\Lambda,t_{\alpha}}))$, notons $({\cal O}'_{\rho'},{\cal L}_{\rho'})$ l'image de $(M_{\Lambda,t_{\alpha}},{\cal E}_{t_{\alpha}},\rho')$ par la repr\'esentation de Springer g\'en\'eralis\'ee et notons ${\cal O}_{\rho'}$ l'orbite de $\mathfrak{g}$ engendr\'ee par ${\cal O}'_{\rho'}$. Le support de $\boldsymbol{\chi}_{M_{\Lambda,t_{\alpha}},{\cal E}_{t_{\alpha}},\rho'}$ est contenu dans l'adh\'erence de ${\cal O}'_{\rho'}$. Pour $N'\in {\cal O}'$, on peut donc restreindre l'ensemble de sommation de (2) \`a   l'ensemble ${\cal R}$ des  $\rho'\in Irr(W_{t_{\alpha}}(M_{\Lambda,t_{\alpha}}))$ tels que $<\rho_{t},\rho'>>0$ et que l'adh\'erence de ${\cal O}'_{\rho'}$ contienne ${\cal O}'$. Prouvons que
 
 (3) pour $\rho'\in {\cal R}$, on a ${\cal O}'_{\rho'}={\cal O}'$. 

 Supposons que ${\cal O}'_{\rho'}\not={\cal O}'$. Puisque ${\cal O}'$ est contenue dans l'adh\'erence de ${\cal O}'_{\rho'}$ et que ${\cal O}$ est l'orbite de $\mathfrak{g}$ engendr\'ee par ${\cal O}'$, il r\'esulte de \ref{orbitesengendrees} que ${\cal O}$ est contenue dans l'adh\'erence de ${\cal O}_{\rho'}$ et que ${\cal O}\not={\cal O}_{\rho'}$. Ces deux propri\'et\'es contredisent la proposition  \ref{premiertheoreme}. Cela prouve (3). 
 
 En cons\'equence, pour $\rho'\in {\cal R}$, la restriction de $\boldsymbol{\chi}_{M_{\Lambda,t_{\alpha}},{\cal E}_{t_{\alpha}},\rho'}$ \`a ${\cal O}'$ est \'egale \`a ${\cal Y}_{{\cal O}',{\cal L}_{\rho'}}$. L'\'egalit\'e (2) devient
 $$(4) \qquad \boldsymbol{\chi}_{M_{\Lambda},{\cal E},\rho}(t_{\alpha}exp(N'))=\sum_{\rho'\in  {\cal R}}<\rho_{t},\rho'> {\cal Y}_{{\cal O}',{\cal L}_{\rho'}}(N').$$
 En comparant avec (1), on obtient que ${\cal R}$ est r\'eduit \`a un unique \'el\'ement. Pour cet \'el\'ement $\rho'$, on a ${\cal L}_{\rho'}={\cal L}$ et $<\rho_{t},\rho'>=1$. Cela d\'emontre le (iii) de l'\'enonc\'e. 
 
  Soit  $(\Lambda,{\cal E},\rho)\in {\cal A}(G,{\cal O})  $.  Supposons $\alpha\not\in \Lambda$ et supposons qu'il existe  une repr\'esentation $\rho'\in Irr(W_{t_{\alpha}}(M_{\Lambda,t_{\alpha}}))$ telle que $<\rho_{t_{\alpha}},\rho'>>0$ et que son image    par la correspondance de Springer g\'en\'eralis\'ee soit port\'ee par   ${\cal O}'$. On d\'efinit l'ensemble ${\cal R}$ comme ci-dessus. L'hypoth\`ese pr\'ec\'edente  entra\^{\i}ne que ${\cal R}$ n'est pas vide.  Les assertions (3) et (4) restent vraies. Il en r\'esulte que la fonction $N'\mapsto \boldsymbol{\chi}_{M_{\Lambda},{\cal E},\rho}(t_{\alpha}exp(N'))$ est non nulle sur ${\cal O}'$. A fortiori, la restriction \`a  ${\cal X}_{t_{\alpha},u}$ du complexe $K_{M_{\Lambda},{\cal E},\rho}$ est non nulle.   Le th\'eor\`eme \ref{parametrage3} dit que $\nabla(\Lambda,{\cal E},\rho)$ est de la forme $(N,d_{\alpha},\mu)$. Par d\'efinition, on a alors $(\Lambda,{\cal E},\rho)\in {\cal A}(G,{\cal O},d)$. Cela prouve (iv).

    Il reste \`a traiter le cas des orbites exceptionnelles. Pour une telle orbite ${\cal O}$, on fixe encore un \'el\'ement $N\in {\cal O}$. On a $\bar{A}(N)=A(N)={\mathbb Z}/2{\mathbb Z}$ et on utilise les m\^emes notations qu'en \ref{parametrage3}. Les éléments de $\bar{{\cal C}}({\cal O})$ sont $(N,1)$ et $(N,d_{2})$ et on les note simplement ${\bf 1}$ et ${\bf d}_{2}$.  On note ${\cal L}_{1}$ et ${\cal L}_{sgn}$  les syst\`emes locaux sur ${\cal O}$ correspondant   aux deux caract\`eres $1$ et $sgn$ de $A(N)$ . Remarquons que le groupe $G$ est déployé donc l'assertion (ii) de l'énoncé est triviale. 
    
    Supposons que $G$ soit  de type $E_{7}$ et que ${\cal O}$ soit de type $A_{4}+A_{1}$. D'apr\`es 
     \ref{parametrage3}, on a ${\cal A}(G,{\cal O},{\bf 1})=\{\rho_{512,11},\rho_{512,12}\}$ et ${\cal A}(G,{\cal O},{\bf d}_{2})=
      \{E_{7}[i], E_{7}[-i]\}$.   Remarquons que  $t_{\alpha_{0}}=1$ donc $G_{t_{\alpha_{0}}}=G$ et ${\cal O}$ est aussi une orbite de $\mathfrak{g}_{t_{\alpha_{0}},nil}$. On pose $b({\bf 1})=(\alpha_{0},{\cal O})$. On a bien $J\circ b({\bf 1})={\bf 1}$.  Le  groupe $G_{t_{\alpha_{4}}}$ est de type $A_{3}\times A_{3}\times A_{1}$. Notons ${\cal O}'$ l'orbite nilpotente r\'eguli\`ere de $\mathfrak{g}_{t_{\alpha_{4}}}$. On pose $b({\bf d}_{2})=(\alpha_{4},{\cal O}')$.  On doit prouver que $J\circ b({\bf d}_{2})={\bf d}_{2}$. L'orbite engendrée par ${\cal O}'$ est ${\cal O}$: avec la normalisation usuelle des produits scalaires, les normes $\vert\vert x_{*,{\cal O}'}\vert\vert$ et $\vert\vert x_{*,{\cal O}}\vert\vert$ sont toutes deux égales à $42$ et on vérifie  que ${\cal O}$ est la seule orbite nilpotente de $G$ pour laquelle cette norme a cette valeur. On peut supposer $N\in {\cal O}'$. Par définition, $J\circ b({\bf d}_{2})$ est alors égal au couple $(N,d_{\alpha_{4}})$, où $d_{\alpha_{4}}$ est l'image de $t_{\alpha_{4}}$ dans $ \bar{A}(N)$.   Si $d_{\alpha_{4}}=1$, alors $t_{\alpha_{4}}$ appartient à $Z_{G}(N)^0$. Ce groupe est un tore de dimension $2$, cf. \cite{C} p. 404. Puisque $t_{\alpha_{4}}$ est semi-simple, on a $Z_{G}(N)^0\subset G_{t_{\alpha_{4}}}$, donc $Z_{G}(N)^0\subset Z_{G_{t_{\alpha_{4}}}}(N)$. Or ce dernier groupe est fini, contradiction. Donc $d_{\alpha_{4}}=d_{2}$ et $J\circ b({\bf d}_{2})={\bf d}_{2}$. 
 
       Il r\'esulte de \cite{C} p. 430 que les images de $\rho_{512,11}$ et $\rho_{512,12}$ par la correspondance de Springer sont $({\cal O},{\cal L}_{1})$ et $({\cal O},{\cal L}_{sgn})$. L'assertion (iii) pour ${\bf d}={\bf 1}$ en r\'esulte. 
  Supposons $N\in {\cal O}'$. Les \'el\'ements 
   $E_{7}[i]$ et $E_{7}[-i]$ sont des faisceaux cuspidaux port\'es par l'orbite de $t_{\alpha_{4}} exp(N)$. On en d\'eduit par restriction et via l'exponentielle des syst\`emes locaux sur ${\cal O}'$.  Ces syst\`emes locaux sont distincts comme on l'a dit en 
   \ref{faisceauxcaracteresunipotentscuspidaux}. L'assertion (iii) pour ${\bf d}={\bf d}_{2}$ en r\'esulte. 
   Pour  ${\bf d}={\bf 1}$, l'assertion (iv) est vide: les \'el\'ements de $A(G,{\cal O},{\bf d}_{2}) $ ont pour ensemble $\Lambda$ associ\'e l'ensemble $\Delta_{a}-\{\alpha_{4}\}$, qui contient $\alpha_{0}$.    Il reste \`a prouver cette assertion (iv) pour ${\bf d}={\bf d}_{2}$. Il s'agit de prouver ce qui suit. Soit $\rho=\rho_{512,11}$ ou $\rho=\rho_{512,12}$. Notons $\rho_{t_{\alpha_{4}}}$ la restriction de $\rho$ au groupe de Weyl $W_{t_{\alpha_{4}}}$ de $G_{t_{\alpha_{4}}}$. Soit $\rho'$ une composante irr\'eductible de $\rho_{t_{\alpha_{4}}}$. Alors l'image de $\rho'$ par la correspondance de Springer est port\'ee par une orbite  distincte de ${\cal O}'$.  En fait, pour notre groupe $G_{t_{\alpha_{4}}}$ de type $A_{3}\times A_{3}\times A_{1}$, l'unique repr\'esentation irr\'eductible de $W_{t_{\alpha_{4}}}$ dont l'image par la correspondance de Springer est port\'ee par ${\cal O}'$ 
   est la repr\'esentation triviale. On doit donc prouver que la repr\'esentation triviale n'intervient pas dans $\rho_{t_{\alpha_{4}}}$. Mais $\rho_{t_{\alpha_{4}}}$ est calcul\'e en \cite{Alvis} table 32 et on voit que cette propri\'et\'e est v\'erifi\'ee.

 Supposons que $G$ soit de type $E_{8}$ et que ${\cal O}$ soit de type $A_{4}+A_{1}$.  On a ${\cal A}(G,{\cal O},{\bf 1})=\{\rho_{4096,26},\rho_{4096,27}\}$ et ${\cal A}(G,{\cal O},{\bf d}_{2})=\{(E_{7}[i],sgn),(E_{7}[-i],sgn)\}$. On pose $b({\bf 1})=(\alpha_{0},{\cal O})$. On a bien $J\circ b({\bf 1})={\bf 1}$. Le groupe $G_{t_{\alpha_{3}}}$ est de type $A_{1}\times A_{7}$. On note ${\cal O}'$ l'orbite nilpotente de $\mathfrak{g}_{t_{\alpha_{3}}}$ qui est le produit de l'orbite r\'eguli\`ere du facteur $A_{1}$ et de l'orbite du facteur $A_{7}$ param\'etr\'ee par la partition $(4^2)$. On pose $b({\bf d}_{2})=(\alpha_{3},{\cal O}')$. On doit prouver que $J\circ b({\bf d}_{2})=({\bf d}_{2})$. L'orbite engendrée par ${\cal O}'$ est ${\cal O}$: avec la normalisation usuelle des produits scalaires, les normes $\vert\vert x_{*,{\cal O}'}\vert\vert$ et $\vert\vert x_{*,{\cal O}}\vert\vert$ sont toutes deux égales à $42$ et on vérifie  que ${\cal O}$ est la seule orbite nilpotente de $G$ pour laquelle cette norme a cette valeur. On peut supposer $N\in {\cal O}'$, on a alors  $J\circ b({\bf d}_{2})=(N,d_{\alpha_{3}})$.   Si $d_{\alpha_{3}}=1$, alors $t_{\alpha_{3}}\in Z_{G}(N)^0$. Ce groupe est  de rang $3$, cf. \cite{C} p. 406. Puisque $t_{\alpha_{3}}$ est semi-simple, il est contenu dans un sous-tore maximal $T'$ de $Z_{G}(N)^0$. On a $T'\subset G_{t_{\alpha_{3}}}$, donc $T'\subset Z_{G_{t_{\alpha_{3}}}}(N)$. Or un sous-tore maximal de ce dernier groupe est de dimension $1$, contradiction. Donc $d_{\alpha_{3}}=d_{2}$ et $J\circ b({\bf d}_{2})=({\bf d}_{2})$.

  Les images  par la correspondance de Springer de $\rho_{4096,26}$ et $\rho_{4096,27}$ sont $({\cal O},{\cal L}_{1})$ et $({\cal O},{\cal L}_{sgn})$. L'assertion (iii) pour ${\bf d}={\bf 1}$ en r\'esulte.  
  Les termes $(E_{7}[\pm i],sgn)$ sont de la forme $(\Lambda,{\cal E},sgn)$ o\`u $\Lambda=\Delta_{a}-\{\alpha_{3},\alpha_{6}\}$, $M_{\Lambda}$ est de type $E_{7}$, ${\cal E}$ est l'un des deux \'el\'ements de ${\bf FC}_{u}(M_{\Lambda,AD})$ et $sgn$ est le caract\`ere non trivial de $W(M_{\Lambda})={\mathbb Z}/2{\mathbb Z}$.  La repr\'esentation $sgn_{t_{\alpha_{3}}}$ est le caract\`ere non trivial  de $W_{t_{\alpha_{3}}}(M_{\Lambda,t_{\alpha_{3}}})={\mathbb Z}/2{\mathbb Z}$.   Les images par la correspondance de Springer g\'en\'eralis\'ee des deux triplets  $(M_{\Lambda,t_{\alpha_{3}}},{\cal E}_{t_{\alpha_{3}}},sgn_{t_{\alpha_{3}}})$ sont port\'ees par ${\cal O}'$ et associ\'ees \`a deux syst\`emes locaux distincts sur cette orbite. L'assertion (iii) pour ${\bf d}={\bf d}_{2}$ en r\'esulte.   De nouveau, l'assertion (iv) est  vide pour ${\bf d}={\bf 1}$. Pour ${\bf d}={\bf d}_{2}$, elle est \'equivalente \`a l'assertion suivante.  Soit $\rho=\rho_{4096,26}$ ou $\rho=\rho_{4096,27}$. Notons $\rho_{t_{\alpha_{3}}}$ la restriction de $\rho$ au groupe de Weyl $W_{t_{\alpha_{3}}}$ de $G_{t_{\alpha_{3}}}$. Soit $\rho'$ une composante irr\'eductible de $\rho_{t_{\alpha_{3}}}$. Alors l'image de $\rho'$ par la correspondance de Springer est port\'ee par une orbite distincte de ${\cal O}'$. Cela r\'esulte de \cite{Alvis} table 42, o\`u les repr\'esentations $\rho_{4096,26}$ et $\rho_{4096,27}$ sont not\'ees respectivement $4096^*_{x}$ et $4096^*_{z}$. 
 
 Supposons que $G$ soit de type $E_{8}$ et que ${\cal O}$ soit de type $E_{6}(a_{1})+A_{1}$.   On a ${\cal A}(G,{\cal O},{\bf 1})=\{\rho_{4096,11},\rho_{4096,12}\}$, ${\cal A}(G,{\cal O},{\bf d}_{2})=\{(E_{7}[i],1),(E_{7}[-i],1)\}$. On pose $b({\bf 1})=(\alpha_{0},{\cal O})$. On a bien $J\circ b({\bf 1})={\bf 1}$. Notons ${\cal O}''$ le produit des deux orbites nilpotentes r\'eguli\`eres des composantes de $G_{t_{\alpha_{3}}}$ de type $A_{1}$ et $A_{7}$. On pose $b({\bf d}_{2})=(\alpha_{3},{\cal O}'')$.  On doit prouver que $J\circ b({\bf d}_{2})=({\bf d}_{2})$. L'orbite engendrée par ${\cal O}''$ est ${\cal O}$: avec la normalisation usuelle des produits scalaires, les normes $\vert\vert x_{*,{\cal O}''}\vert\vert$ et $\vert\vert x_{*,{\cal O}}\vert\vert$ sont toutes deux égales à $170$ et on vérifie  que ${\cal O}$ est la seule orbite nilpotente de $G$ pour laquelle cette norme a cette valeur. On peut supposer $N\in {\cal O}''$, on a alors  $J\circ b({\bf d}_{2})=(N,d_{\alpha_{3}})$.   Si $d_{\alpha_{3}}=1$, alors $t_{\alpha_{3}}\in Z_{G}(N)^0$. Ce groupe est  un tore de dimension $1$, cf. \cite{C} p. 407. On obtient une contradiction comme dans le premier cas ci-dessus.   Donc $d_{\alpha_{3}}=d_{2}$ et $J\circ b({\bf d}_{2})=({\bf d}_{2})$.

Les images  par la correspondance de Springer de $\rho_{4096,11}$ et $\rho_{4096,12}$ sont $({\cal O},{\cal L}_{1})$ et $({\cal O},{\cal L}_{sgn})$. L'assertion (iii) pour ${\bf d}={\bf 1}$ en r\'esulte.  
Les termes $(E_{7}[\pm i],1)$ sont de la forme $(\Lambda,{\cal E},1)$ o\`u $\Lambda=\Delta_{a}-\{\alpha_{3},\alpha_{6}\}$, $M_{\Lambda}$ est de type $E_{7}$, ${\cal E}$ est l'un des deux \'el\'ements de ${\bf FC}_{u}(M_{\Lambda,AD})$ et $1$ est le caract\`ere  trivial de $W(M_{\Lambda})={\mathbb Z}/2{\mathbb Z}$.  La repr\'esentation $1_{t_{\alpha_{3}}}$ est le caract\`ere  trivial  de $W_{t_{\alpha_{3}}}(M_{\Lambda,t_{\alpha_{3}}})={\mathbb Z}/2{\mathbb Z}$.   Les images par la correspondance de Springer g\'en\'eralis\'ee des deux triplets  $(M_{\Lambda,t_{\alpha_{3}}},{\cal E}_{t_{\alpha_{3}}},1)$ sont port\'ees par ${\cal O}''$ et associ\'ees \`a deux syst\`emes locaux distincts sur cette orbite. L'assertion (iii) pour ${\bf d}={\bf d}_{2}$ en r\'esulte.  L'assertion (iv) est vide pour ${\bf d}={\bf 1}$. Pour ${\bf d}={\bf d}_{2}$, elle se  d\'eduit encore de \cite{Alvis} table 42, o\`u $\rho_{4096,11}$ et $\rho_{4096,12}$ sont not\'ees respectivement $4096_{z}$ et $4096_{x}$. Cela ach\`eve la preuve du th\'eor\`eme. $\square$
 
 {\bf Remarque.}  Nous n'utiliserons le th\'eor\`eme \ref{parametrage3} que via le th\'eor\`eme ci-dessus dont il est une cons\'equence. Mais ce dernier th\'eor\`eme est beaucoup plus faible que le th\'eor\`eme \ref{parametrage3} et beaucoup plus simple \`a prouver. Par exemple, supposons que  $A(G,{\cal O})$ est r\'eduit \`a un seul \'el\'ement (ce qui est souvent les cas pour les groupes exceptionnels). Cet \'el\'ement est alors de la forme $(M,{\cal E},\rho)$ o\`u $M=T$, ${\cal E}$ est le syst\`eme local trivial sur $M_{AD}$ et $\rho$ est une repr\'esentation irr\'eductible sp\'eciale de $W$.  L'image de $\rho$ par la repr\'esentation de Springer est $({\cal O},{\cal L}_{1})$, o\`u ${\cal L}_{1}$ est le syst\`eme local trivial sur ${\cal O}$. On pose alors $b(1)=(\alpha_{0},{\cal O})$ et toutes
  les assertions du th\'eor\`eme sont v\'erifi\'ees. Pour les autres orbites des groupes exceptionnels, la d\'emonstration du th\'eor\`eme devrait r\'esulter comme dans le cas ci-dessus des orbites exceptionnelles des tables de \cite{Alvis} mais nous n'avons pas eu la patience de traiter toutes les orbites.

 \section{  L'espace des combinaisons lin\'eaires stables d'int\'egrales orbitales nilpotentes }

\subsection{Int\'egrales orbitales nilpotentes et espace ${\cal D}(\mathfrak{g}(F))$}

\subsubsection{Normalisation des int\'egrales orbitales nilpotentes\label{integralesnilpotentes}}
Dans cette section finale, on suppose que $G$ est un groupe r\'eductif connexe d\'efini sur $F$,  {\bf  adjoint, absolument simple et quasi-d\'eploy\'e}. On impose   l'hypoth\`ese $(Hyp)_{2}(p)$  de \ref{leshypothesessurp}.

On fixe un \'epinglage $\mathfrak{E}$ et un \'epinglage affine $\mathfrak{E}_{a}$ comme en \ref{alcove} et \ref{racinesaffines}, tous deux conserv\'es par l'action galoisienne. On utilise les constructions de ces paragraphes. 

 Rappelons qu'en \ref{lespaceFC}, on a fix\'e un caract\`ere $\psi$ de $F$, une forme bilin\'eaire sym\'etrique non d\'eg\'en\'er\'ee $<.,.>$ sur $\mathfrak{g}(F)$ et la transformation de Fourier $f\mapsto \hat{f}$ qui s'en d\'eduit. 
 
Soit ${\mathbb O}\in \mathfrak{g}_{nil}(F)/conj$. Pour $N\in {\mathbb O}$, l'espace tangent \`a ${\mathbb O}$ au point $N$ s'identifie \`a $\mathfrak{g}/Z_{\mathfrak{g}}(N)$. Cet espace est muni de la forme symplectique $(X,Y)\mapsto < N,[X,Y]>$. De $\psi$ et de cette forme se d\'eduit une mesure autoduale sur $\mathfrak{g}/Z_{\mathfrak{g}}(N)$. La collection de ces mesures sur les espaces tangents en tout point de ${\mathbb O}$ d\'efinit une mesure $d_{{\mathbb O}}$ sur cette orbite ${\mathbb O}$. Elle est invariante par conjugaison par $G(F)$. On d\'efinit l'int\'egrale orbitale $I_{{\mathbb O}}$ sur $C_{c}^{\infty}(\mathfrak{g}(F))$ par $$I_{{\mathbb O}}(f)=\int_{{\mathbb O}}f(N)\,d_{{\mathbb O}}N$$
pour toute $f\in C_{c}^{\infty}(\mathfrak{g}(F))$. La famille $(I_{{\mathbb O}})_{{\mathbb O}\in \mathfrak{g}_{nil}(F)/conj}$ est une base de l'espace $I(\mathfrak{g}(F))_{nil}^*$ des distributions invariantes \`a support nilpotent. 

Pour $f\in C_{c}^{\infty}(\mathfrak{g}(F))$ et $\lambda\in F^{\times}$, notons $f^{\lambda}$ la fonction d\'efinie par $f^{\lambda}(X)=f(\lambda X)$. Soit ${\mathbb O}\in \mathfrak{g}_{nil}(F)/conj$. L'orbite ${\mathbb O}$ est invariante par l'homoth\'etie $N\mapsto \lambda^2N$. On a l'\'egalit\'e

(1) $I_{{\mathbb O}}(f^{\lambda^2})=\vert \lambda\vert _{F}^{-dim({\mathbb O})}I_{{\mathbb O}}(f)$.

Pour $d\in {\mathbb N}$, notons $I(\mathfrak{g}(F))_{nil}^*[d]$ le sous-espace de $I(\mathfrak{g}(F))_{nil}^*$ engendr\'e par les $I_{{\mathbb O}}$ pour les orbites ${\mathbb O}$ telles que $dim({\mathbb O})=d$. L'espace $I(\mathfrak{g}(F))_{nil}^*$ est somme directe des sous-espaces $I(\mathfrak{g}(F))_{nil}^*[d]$. Pour $J\in I(\mathfrak{g}(F))_{nil}^*$, notons $J[d]$ sa composante sur $I(\mathfrak{g}(F))_{nil}^*[d]$. Montrons que

(2) $J$ est stable si et seulement si $J[d]$ l'est pour tout $d\in {\mathbb Z}$. 

L'assertion (1) montre que, pour $f\in I(\mathfrak{g}(F))$, $J[d](f)$ se calcule par interpolation des valeurs $J(f^{\lambda^2})$ pour $\lambda\in F^{\times}$. L'espace $I^{inst}(\mathfrak{g}(F))$ d\'efini en \ref{lespacecalD} est conserv\'e par l'application $f\mapsto f^{\lambda}$. Si $J$ est stable et que  $f\in I^{inst}(\mathfrak{g}(F))$, on a $J(f^{\lambda^2})=0$ pour tout $\lambda$, donc $J[d](f)=0$. Cela prouve (2).

\subsubsection{La fonction-test associ\'ee \`a un couple $(s,{\cal O})$\label{fonctiontest}}

Soit $s$ un  sommet de l'immeuble $Imm(G)$ et soit ${\cal O}\in \mathfrak{g}_{s,nil}({\mathbb F}_{q})/conj$.  Fixons un $\mathfrak{sl}(2)$-triplet $(f,h,e)$ de $\mathfrak{g}_{s}({\mathbb F}_{q})$ tel que $e\in {\cal O}$. Remarquons que $-f$ appartient aussi \`a ${\cal O}$ (le $\mathfrak{sl}(2)$-triplet d\'etermine un homomorphisme $SL(2)\to G_{s}$ d\'efini sur ${\mathbb F}_{q}$ et $-f$ est conjugu\'e \`a $e$ par un \'el\'ement de l'image de $SL(2,{\mathbb F}_{q})$). Au triplet est associ\'e une graduation $(\mathfrak{g}_{s,i})_{i\in {\mathbb Z}}$ de $\mathfrak{g}_{s}$: $\mathfrak{g}_{s,i}=\{X\in \mathfrak{g}_{s}; \forall \,\, \lambda\in \bar{{\mathbb F}}_{q}^{\times},\,\,Ad(x_{*,h}(\lambda))(X)=\lambda^{i}X\}$. On pose $\mathfrak{g}_{s,\geq i}=\oplus_{j\geq i} \mathfrak{g}_{s,j}$. Notons $\underline{h}'_{s,{\cal O}}$ la fonction caract\'eristique de $-f+\mathfrak{g}_{s,\geq-1}$ et posons
$$\underline{h}_{s,{\cal O}}= q^{-dim(\mathfrak{g}_{s,1})/2)}\vert G_{s}({\mathbb F}_{q})\vert ^{-1}\sum_{x\in G_{s}({\mathbb F}_{q})}{^x\underline{h}'}_{s,{\cal O}}.$$
 Cette fonction ne d\'epend pas du choix du $\mathfrak{sl}(2)$-triplet. 
 On rel\`eve $\underline{h}_{s,{\cal O}}$ en une fonction sur $\mathfrak{k}_{s}$ invariante par translations par $\mathfrak{k}_{s}^+$, que l'on \'etend en une fonction sur $\mathfrak{g}(F)$ nulle hors de $\mathfrak{k}_{s}$. On note $h_{s,{\cal O}}$ cette derni\`ere fonction. Dans le cas où ${\cal O}=\{0\}$, cette fonction est la fonction caractéristique de $\mathfrak{k}_{s}$. 
 
On a introduit en \ref{lespacecalD} le sous-espace ${\cal H}\subset I(\mathfrak{g}(F))$. 

\begin{prop}{(i) La fonction $h_{s,{\cal O}}$ appartient \`a ${\cal H}$.

(ii) Soit ${\mathbb O}\in \mathfrak{g}_{nil}(F)/conj$. Si $rel_{s,F}({\cal O})$ n'est pas incluse dans l'adh\'erence de ${\mathbb  O}$ pour la topologie $p$-adique, alors $I_{{\mathbb  O}}(h_{s,{\cal O}})=0$.

(iii) Soit ${\mathbb O}=rel_{s,F}({\cal O})$. Alors $I_{{\mathbb O}}(h_{s,{\cal O}})= 1$.}\end{prop}

 Preuve.  La fonction $\underline{h}'_{s,{\cal O}}$ est invariante par translations par $\mathfrak{g}_{s,\geq0}$ et  cet ensemble contient une sous-alg\`ebre d'Iwahori. Il en résulte que la fonction sur $\mathfrak{g}(F)$ déduite de $\underline{h}'_{s,{\cal O}}$ appartient à ${\cal H}$. La fonction $h_{s,{\cal O}}$  est combinaison linéaire de conjuguées de la fonction précédente et appartient donc  elle aussi à ${\cal H}$.  On renvoie à  \cite{DK} lemme 7.3.1  pour les assertions (ii) et (iii). $\square$
 
 Quand $s$ d\'ecrit les sommets de $Imm(G)$, ou m\^ eme seulement les éléments de $S(\bar{C})$,  et ${\cal O}$ d\'ecrit $\mathfrak{g}_{s,nil}({\mathbb F}_{q})/conj$, les orbites $rel_{s,F}({\cal O})$ d\'ecrivent $\mathfrak{g}_{nil}(F)/conj$ tout entier, cf. \ref{debacker}(8). Il r\'esulte de la proposition que les fonctions $h_{s,{\cal O}}$ s\'eparent les int\'egrales orbitales nilpotentes. 
 
 \subsubsection{ Calcul de la fonction $\underline{h}_{s,{\cal O}}$\label{calcul}}
On conserve les donn\'ees du paragraphe pr\'ec\'edent. On a d\'efini l'ensemble ${\cal I}^{FC}_{0,{\mathbb F}_{q}}(G_{s})$ en \ref{Springer}. Pour $(M,{\cal E})\in {\cal I}^{FC}_{0,{\mathbb F}_{q}}(G_{s})$, on note $W_{s}(M)=Norm_{G_{s}}(M)/M$. 

 \begin{prop}{Pour tout $X\in \mathfrak{g}_{s,nil}({\mathbb F}_{q})$, on a l'\'egalit\'e
 $$\underline{h}_{s,{\cal O}}(X)=\vert G_{s}({\mathbb F}_{q})\vert ^{-1}q^{dim(\mathfrak{g}_{s})-dim({\cal O})/2}\sum_{(M,{\cal E})\in {\cal I}^{FC}_{0,{\mathbb F}_{q}}(G_{s})}\vert W_{s}(M)\vert ^{-1}q^{-dim(Z(M)^0)}$$
 $$\sum_{w\in W_{s}(M)} \vert Z(M_{w})^0({\mathbb F}_{q})\vert Q_{M,{\cal E},w}^{\natural}(e)Q_{M,{\cal E},w}(-X).$$}\end{prop}

Preuve. La proposition est cons\'equence de la proposition 6.12 de \cite{L4}. On applique les constructions de cette r\'ef\'erence au $\mathfrak{sl}(2)$-triplet $(-f,h,-e)$. Lusztig introduit en \cite{L4} 2.4 une fonction $\boldsymbol{\Gamma}_{\phi}$ sur $\mathfrak{g}_{s,nil}({\mathbb F}_{q})$ d\'efinie par
$$\boldsymbol{\Gamma}_{\phi}(X)=q^{dim(\mathfrak{g}_{s,1})/2-dim(\mathfrak{g}_{s,\geq1})}\sum_{g\in G_{s}({\mathbb F}_{q})}\psi(<-f,gXg^{-1}>){\bf 1}_{\mathfrak{g}_{s,\geq 1}}(gXg^{-1})$$
pour $X\in \mathfrak{g}_{s,nil}({\mathbb F}_{q})$, 
o\`u ${\bf 1}_{\mathfrak{g}_{s,\geq1}}$ est la fonction caract\'eristique de $\mathfrak{g}_{s,\geq1}({\mathbb F}_{q})$.  En appliquant nos d\'efinitions, on a
$$\hat{\underline{h}}'_{s,{\cal O}}(X)=q^{dim(\mathfrak{g}_{s,\geq -1})-dim(\mathfrak{g}_{s})/2}\psi(<-f,X>){\bf 1}_{\mathfrak{g}_{s,\geq1}}(X),$$
puis
$$\hat{\underline{h}}_{s,{\cal O}}(X)=q^{-dim(\mathfrak{g}_{s,1})/2+dim(\mathfrak{g}_{s,\geq-1})-dim(\mathfrak{g}_{s})/2}\vert G_{s}({\mathbb F}_{q})\vert ^{-1}$$
$$\sum_{g\in G_{s}({\mathbb F}_{q})}\psi(<-f,gXg^{-1}>){\bf 1}_{\mathfrak{g}_{s,\geq1}}(gXg^{-1}).$$
D'o\`u 
$$(1) \qquad \hat{\underline{h}}_{s,{\cal O}}(X)=c \boldsymbol{\Gamma}_{\phi}(X)$$
pour $X\in \mathfrak{g}_{s,nil}({\mathbb F}_{q})$, o\`u 
$$c=q^{-dim(\mathfrak{g}_{s,1})+dim(\mathfrak{g}_{s,\geq-1})+dim(\mathfrak{g}_{s,\geq1})-dim(\mathfrak{g}_{s})/2}\vert G_{s}({\mathbb F}_{q})\vert ^{-1},$$
ou encore
$$c=q^{dim(\mathfrak{g}_{s})/2}\vert G_{s}({\mathbb F}_{q})\vert ^{-1}.$$
Les fonctions $\hat{\underline{h}}_{s,{\cal O}}$ et $\boldsymbol{\Gamma}_{\phi}$ sont toutes deux \`a support nilpotent, l'\'egalit\'e (1) est donc v\'erifi\'ee pour tout $X\in \mathfrak{g}_{s}({\mathbb F}_{q})$. Par inversion de Fourier, on en d\'eduit
$$(2) \qquad \underline{h}_{s,{\cal O}}(X)=c\hat{\boldsymbol{\Gamma}}_{\phi}(-X)$$
pour tout $X\in \mathfrak{g}_{s}({\mathbb F}_{q})$. 

La proposition 6.12 de \cite{L4} calcule la restriction de $\hat{\boldsymbol{\Gamma}}_{\phi}$ \`a $\mathfrak{g}_{s,nil}({\mathbb F}_{q})$ sous la forme
$$(3) \qquad \hat{\boldsymbol{\Gamma}}_{\phi}(X)=\sum_{(M,{\cal E})\in  {\cal I}^{FC}_{0,{\mathbb F}_{q}}(G_{s})}\boldsymbol{\Psi}_{M,{\cal E}}(X)$$
pour $X\in \mathfrak{g}_{s,nil}({\mathbb F}_{q})$. On doit prendre garde \`a  quelques points techniques. Lusztig ne normalise pas la transformation de Fourier comme nous, la sienne est la n\^otre multipli\'ee par $q^{dim(\mathfrak{g}_{s})/2}$. Comme on l'a dit, son triplet est $(-f,h,-e)$, le terme $-N$ apparaissant dans sa formule est $e$ pour nous. Il suppose que le groupe $G_{s}$ est d\'eploy\'e. Ses arguments valent dans le cas g\'en\'eral, mais il appara\^{\i}t des Frobenius dans les formules. Enfin, on peut simplifier la somme en $\iota'$ apparaissant dans sa formule en vertu de l'\'egalit\'e 6.6(b) de \cite{L4}. En tenant compte de ces remarques et en traduisant les notations de Lusztig en les n\^otres, on obtient la formule suivante. Soit $(M,{\cal E})\in {\cal I}^{FC}_{0,{\mathbb F}_{q}}(G_{s})$. Pour $X\in \mathfrak{g}_{s,nil}({\mathbb F}_{q})$, on a l'\'egalit\'e
$$(4) \qquad \boldsymbol{\Psi}_{M,{\cal E}}(X)=\sum_{\rho,\rho_{1}\in Irr_{{\mathbb F}_{q}}(W_{s}(M))}q^{b(M,{\cal E},\rho_{1})-b(M,{\cal E},\rho)-dim({\cal O})/2+dim(\mathfrak{g}_{s})/2-dim(Z(M)^0)} $$
$$\vert W_{s}(M)\vert ^{-1}\sum_{w\in W_{s}(M)}trace(\rho^{\flat}(w^{-1}Fr))trace(\rho_{1}^{\flat}(Fr^{-1}w))\vert Z(M_{w})^0({\mathbb F}_{q})\vert \chi^{\natural}_{M,{\cal E},\rho}(e)\chi_{M,{\cal E},\rho_{1}}(X)$$
$$=q^{-dim({\cal O})/2+dim(\mathfrak{g}_{s})/2-dim(Z(M)^0)} \vert W_{s}(M)\vert ^{-1}\sum_{w\in W_{s}(M)} \vert Z(M_{w})^0({\mathbb F}_{q})\vert S^{\natural}(e,w)S(X,w),$$
où
 $$S^{\natural}(e,w)= \sum_{\rho\in Irr_{{\mathbb F}_{q}}(W_{s}(M))}q^{-b(M,{\cal E},\rho)}trace(\rho^{\flat}(w^{-1}Fr)) \chi^{\natural}_{M,{\cal E},\rho}(e)$$
et
$$S(X,w)=\sum_{\rho_{1}\in Irr_{{\mathbb F}_{q}}(W_{s}(M))}q^{b(M,{\cal E},\rho_{1})}trace(\rho_{1}^{\flat}(Fr^{-1}w))\chi_{M,{\cal E},\rho_{1}}(X).$$
Dans la définition de $S^{\natural}(e,w)$, utilisons la formule  \ref{Springer}(3).  Alors
$$S^{\natural}(e,w)=\vert W_{s}(M)\vert ^{-1}\sum_{v\in W_{s}(M)}Q^{\natural}_{M,{\cal E},v}(e)\sum_{\rho\in Irr_{{\mathbb F}_{q}}(W_{s}(M))}trace(\rho^{\flat}(v Fr^{-1}))trace(\rho^{\flat}(w^{-1}Fr)) .$$
Notons $Cl_{Fr}(w)$ la classe de $Fr$-conjugaison de $w$. Fixons $v\in W_{s}(M)$. Selon les formules d'orthogonalit\'e pour le groupe $W_{s}(M)\rtimes \Gamma_{{\mathbb F}_{q}}$, la somme en $\rho$ vaut $0$ si $v\not\in Cl_{Fr}(w)$ et vaut $\vert W_{s}(M)\vert \vert Cl_{Fr}(w)\vert ^{-1}$ si $v\in Cl_{Fr}(w)$. Dans ce dernier cas, on a $Q^{\natural}_{M,{\cal E},v}(e)=Q^{\natural}_{M,{\cal E},w}(e)$ et la somme de ces termes sur $v\in Cl_{Fr}(w)$ vaut $\vert Cl_{Fr}(w)\vert  Q^{\natural}_{M,{\cal E},w}(e)$.  Donc $S^{\natural}(e,w)=Q^{\natural}_{M,{\cal E},w}(e)$. Un calcul analogue utilisant la formule \ref{Springer}(2) prouve que   $S(X,w)=Q_{M,{\cal E},w}(X)$. Alors la formule (4) devient
$$\boldsymbol{\Psi}_{M,{\cal E}}(X)=q^{-dim({\cal O})/2+dim(\mathfrak{g}_{s})/2-dim(Z(M)^0)}\vert W_{s}(M)\vert ^{-1}$$
$$\sum_{w\in W_{s}(M)}
\vert Z(M_{w})^0({\mathbb F}_{q})\vert Q^{\natural}_{M,{\cal E},w}(e)Q_{M,{\cal E},w}(X).$$
L'\'enonc\'e r\'esulte de cette \'egalit\'e et de (2) et (3). $\square$

\subsubsection{Description de $SI(\mathfrak{g}(F))^*_{nil}$\label{descriptionnil}}
On a introduit en \ref{lespacecalD} un espace  ${\cal D}(\mathfrak{g}(F))$ et un homomorphisme $D^G:{\cal D}(\mathfrak{g}(F))\to I(\mathfrak{g}(F))^*$ qui est antilinéaire. On a note $\hat{D}^G$ la compos\'ee de cette application et de la transformation de Fourier dans $I(\mathfrak{g}(F))^*$. D'apr\`es le lemme \ref{lespacecalD},  l'application $d\mapsto res_{{\cal H}}(\hat{D}^G(d))$ est injective. Il en est de m\^eme de l'application $J\mapsto res_{{\cal H}}(J)$ de $I(\mathfrak{g}(F))^*_{nil}$ dans ${\cal H}^*$. De plus, ces deux applications ont m\^eme image dans ${\cal H}^*$: cette assertion est \'equivalente par transformation de Fourier \`a  la proposition 5.5 de \cite{W5}. Il existe donc un unique isomorphisme antilin\'eaire $\delta:I(\mathfrak{g}(F))^*_{nil}\to {\cal D}(\mathfrak{g}(F))$ de sorte que $res_{{\cal H}}(J)=res_{{\cal H}}(\hat{D}^G(\delta(J))$ pour tout $J\in I(\mathfrak{g}(F))_{nil}^*$.  On a d\'efini le sous-espace ${\cal D}^{st}(\mathfrak{g}(F))$ de ${\cal D}(\mathfrak{g}(F))$. 
On note $SI(\mathfrak{g}(F))^*_{nil}=I(\mathfrak{g}(F))^*_{nil}\cap SI(\mathfrak{g}(F))^*$ l'espace des distributions stables \`a support nilpotent.

\begin{lem}{ L'isomorphisme $\delta$ se restreint en un isomorphisme (antilin\'eaire) de $SI(\mathfrak{g}(F))^*_{nil}$ sur ${\cal D}^{st}(\mathfrak{g}(F))$.  }\end{lem}

Preuve. Notons $S{\cal H}^*$ l'image de $SI(\mathfrak{g}(F))^*$ par $res_{{\cal H}}$. D'apr\`es le lemme \ref{lespacecalD}, pour $d\in {\cal D}(\mathfrak{g}(F))$, $res_{{\cal H}}(\hat{D}^G(d))$ appartient \`a $S{\cal H}^*$ si et seulement si $d$ appartient \`a ${\cal D}^{st}(\mathfrak{g}(F))$. Il suffit donc de prouver que, pour $J\in I(\mathfrak{g}(F))^*_{nil}$, on a $res_{{\cal H}}(J)\in S{\cal H}^*$ si et seulement si $J$ est stable. Le sens "si" est \'evident. Supposons $res_{{\cal H}}(J)\in S{\cal H}^*$. Fixons ${\bf J}\in SI(\mathfrak{g}(F))^*$ tel que $res_{{\cal H}}(J)=res_{{\cal H}}({\bf J})$. Soit $f\in I^{inst}(\mathfrak{g}(F))$.  Il existe un entier $n\in {\mathbb Z}$ tel que $f^{\lambda}\in {\cal H}$ pour tout $\lambda$ tel que $val_{F}(\lambda)\geq n$. Pour un tel $\lambda$, on a $J(f^{\lambda})={\bf J}(f^{\lambda})$ et ceci est nul puisque ${\bf J}$ est stable. La distribution $J$ est somme d'int\'egrales orbitales nilpotentes, lesquelles v\'erifient la propri\'et\'e  \ref{integralesnilpotentes}(2): la fonction $\lambda\mapsto J(f^{\lambda^2})$ est combinaison lin\'eaire de fonctions $\lambda\mapsto \vert \lambda\vert _{F}^{d}$ o\`u $d\in {\mathbb Z}$. Une telle fonction qui est nulle pour tout $\lambda$ tel que $val_{F}(\lambda^2)\geq n$ est nulle pour tout $\lambda$. Donc $J(f^{\lambda^2})=0$ pour tout $\lambda\in F^{\times}$. En particulier, pour $\lambda=1$, $J(f)=0$. Donc $J$ annule $I^{inst}(\mathfrak{g}(F))$ et $J$ est stable. $\square$.

  \subsubsection{L'ensemble ${\cal B}(G)$\label{calBG}}
  Soient $M$ un $F$-Levi de $G$ et $s_{M}$ un sommet de $Imm(M_{AD})$ tel que $FC^{st}(\mathfrak{m}_{SC,s_{M}}({\mathbb F}_{q}))$ ne soit pas nul.  Cet espace $FC^{st}(\mathfrak{m}_{SC,s_{M}}({\mathbb F}_{q}))$ a pour base les fonctions caract\'eristiques des \'el\'ements de  ${\bf FC}^{st}_{{\mathbb F}_{q}}(\mathfrak{m}_{SC,s_{M}})$, qui sont des faisceaux-caract\`eres cuspidaux sur $\mathfrak{m}_{SC,s_{M}}$, \`a support nilpotent et conserv\'es par l'action de Frobenius. On fixe une telle base ${\cal B}(M,s_{M})$. Pour $\varphi\in {\cal B}(M,s_{M})$, on note ${\cal E}_{\varphi}$ le faisceau-caract\`ere associ\'e. Il est muni d'une action de Frobenius. On note $\bar{{\cal O}}_{\varphi}$ l'orbite nilpotente de $\mathfrak{m}_{SC,s_{M}}$ supportant ${\cal E}_{\varphi}$ et on en fixe un \'el\'ement $N_{\varphi}$. Au faisceau ${\cal E}_{\varphi}$ est associ\'e un caract\`ere $\xi_{\varphi}$ de $Z_{M_{AD,s_{M}}}(N)/Z_{M_{AD,s_{M}}}(N)^0$. De $\varphi$ est issue une fonction sur $\mathfrak{m}_{SC}(F)$. On la  note  encore $\varphi$ si cela ne  cr\'ee pas de confusion, ou $\varphi^M$ s'il para\^{\i}t utile de distinguer les deux fonctions. Par cette application de rel\`evement, ${\cal B}(M,s_{M})$ s'identifie \`a un sous-ensemble de $FC^{st}(\mathfrak{m}_{SC}(F))$. Ce dernier ensemble  a pour base la r\'eunion des ${\cal B}(M,s_{M})$ quand $s_{M}$ parcourt un ensemble de repr\'esentants des orbites de l'action de $M(F)$ dans l'ensemble des sommets $s_{M}$ tels que $FC^{st}(\mathfrak{m}_{SC,s_{M}}({\mathbb F}_{q}))\not=\{0\}$.
  
Remarquons que, pour un \'el\'ement $g\in G(F)$ tel que l'action $Ad(g)$ conserve le couple $(M,s_{M})$, cette action est triviale sur $FC^{st}(\mathfrak{m}_{SC,s_{M}}({\mathbb F}_{q}))$, cf. \ref{automorphismes}. On peut donc supposer que, si deux couples $(M,s_{M})$ et $(M',s_{M'})$ v\'erifient tous deux les hypoth\`eses pr\'ec\'edentes et s'ils sont conjugu\'es par un \'el\'ement de $G(F)$, cette conjugaison envoie ${\cal B}(M,s_{M})$ sur ${\cal B}(M',s_{M'})$. 
  On note ${\cal B}(G)$ l'ensemble des classes de conjugaison par $G(F)$ de triplets $(M,s_{M},\varphi)$, o\`u $M$ et $s_{M}$ v\'erifient les hypoth\`eses ci-dessus et $\varphi\in {\cal B}(M,s_{M})$.  
  
 Une variante de cet ensemble nous sera utile.  On a vu en \ref{couplesstables} que l'ensemble des classes de conjugaison de couples $(M,s_{M})$ v\'erifiant les hypoth\`eses ci-dessus \'etait en bijection avec un sous-ensemble de l'ensemble des sous-ensembles propres $\Lambda\subset \Delta_{a}^{nr}$. Notons $\boldsymbol{\Lambda}$ ce sous-ensemble et ${\cal B}'(G)$ l'ensemble des couples $(\Lambda,\varphi)$, o\`u $\Lambda\in \boldsymbol{\Lambda}$ et $\varphi\in {\cal B}(M_{\Lambda},s_{\Lambda})$. Alors l'application $(M,s_{M})\mapsto \Lambda$ d\'efinit une bijection de ${\cal B}(G)$ sur ${\cal B}'(G)$. 
 
 Pour $(M,s_{M})$ comme ci-dessus, on a d\'efini en \ref{L*stFsd} l'ensemble $W^{I_{F}}(M)/Fr-conj$. On note $Irr_{F}(W^{I_{F}}(M))$ l'ensemble des repr\'esentations irr\'eductibles de $W^{I_{F}}(M)$ qui sont conserv\'ees par l'action  de $\Gamma_{F}^{nr}$. Pour tout $\rho\in Irr_{F}(W^{I_{F}}(M))$, on fixe un prolongement continu $\rho^{\flat}$ de $\rho$ au produit semi-direct $W^{I_{F}}(M)\rtimes \Gamma_{F}^{nr}$. 
 On note ${\cal B}_{W}(G)$, resp. ${\cal B}_{Irr}(G)$,  l'ensemble des classes de conjugaison par $G(F)$ de quadruplets $(M,s_{M},\varphi,w)$, o\`u $M$ et $s_{M}$ v\'erifient les hypoth\`eses ci-dessus, $\varphi\in {\cal B}(M,s_{M})$ et $w\in W^{I_{F}}(M)/Fr-conj$, resp. $\rho\in Irr_{F}(W^{I_{F}}(M))$.   On a de nouveau les variantes \'evidentes ${\cal B}'_{W}(G)$, resp. ${\cal B}'_{Irr}(G)$, form\'ees de classes de triplets $(\Lambda,\varphi,w)$, resp. $(\Lambda,\varphi,\rho)$. 
  
  \subsubsection{Deux bases de ${\cal D}^{st}(\mathfrak{g}(F))$\label{deuxbases}}
Soit $(M,s_{M},\varphi,w)\in {\cal B}_{W}(G)$. D'apr\`es la proposition \ref{Levistandard}, $M$ est conjugu\'e par un \'el\'ement de $G(F)$ \`a un unique Levi standard et on ne perd rien \`a supposer que $M$ est ce Levi.  Au couple $(M,w)$, on a associ\'e en \ref{parametrage} une classe de conjugaison stable dans ${\cal L}_{F}^{nr,st}$ (le passage au groupe $G^*$ de cette r\'ef\'erence dispara\^{\i}t puisque $G$ est quasi-d\'eploy\'e). Soit $H$ un \'el\'ement de cette classe. On a d\'efini un isomorphisme canonique $\iota_{H,M}:FC^{st}(\mathfrak{m}_{SC}(F))\to FC^{st}(\mathfrak{h}_{SC}(F))$ en \ref{isomorphismesFC}.  On note $\varphi_{H}=\iota_{H,M}(\varphi)$. Par d\'efinition de l'isomorphisme $\iota_{H,M}$, on peut fixer un sommet $s_{H}\in Imm(H_{AD})$ de sorte que les couples $(H,s_{H})$ et $(M,s_{M})$ soient conjugu\'es par un \'el\'ement de $G(F^{nr})$ et que $\varphi_{H}$ provienne d'un \'el\'ement de $ FC^{st}(\mathfrak{h}_{SC,s_{H}}({\mathbb F}_{q}))$. On a d\'efini la distribution $D^G_{X_{H},\varphi_{H}}$ en \ref{descriptioncasgeneral}. 

 On peut consid\'erer $\varphi_{H}$ comme un \'el\'ement de ${\cal K}^{st}(\mathfrak{g}(F))$, cf. \ref{lespaceKstablegeneral}. Gr\^ace qux propositions \ref{compositionderechef} et \ref{comparaison}, la proposition \ref{descriptionfinale} \'equivaut \`a dire que, quand $(M,s_{M},\varphi,w)$ d\'ecrit $ {\cal B}_{W}(G)$, les distributions $k^{G,st}(\varphi_{H})$ d\'ecrivent une base de $D^G({\cal D}^{st}(\mathfrak{g}(F)))$. 

On a noté  $W^{nr}(H)=H(F^{nr})\backslash Norm_{G(F^{nr})}(H)$ et 
$W_{F}(H)=H(F)\backslash Norm_{G(F)}(H)$. Fixons un ensemble de repr\'esentants ${\cal L}_{F}^{nr,st}(M,w)$ des classes de conjugaison par $G(F)$ dans la classe de conjugaison stable dans  ${\cal L}_{F}^{nr,st}$ associ\'ee \`a $(M,w)$. Pour $H\in {\cal L}_{F}^{nr,st}$, on a d\'efini en \ref{uncalculdemesures} une mesure $m(H)$ dont on a prouv\'e qu'elle ne d\'ependait que de la classe de conjugaison stable de $H$. Elle est donc constante sur ${\cal L}_{F}^{nr,st}(M,w)$, on note $m(w)$ cette valeur constante.  Posons
$$c(w)=m(w)^{-1} \vert M_{AD,s_{M}}({\mathbb F}_{q})\vert q^{dim(M_{SC,s_{M}})-dim(\bar{{\cal O}}_{\varphi})/2}.$$
Introduisons la distribution
 $$D(M,s_{M}, \varphi,w)=c(w) \sum_{H\in {\cal L}_{F}^{nr,st}(M,w)}\vert W^{nr}(H)^{\Gamma_{F}^{nr}}\vert \vert W_{F}(H)\vert ^{-1}D^G_{X_{H},\varphi_{H}}.$$
  Notons $k(M,s_{M}, \varphi,w)$ la restriction \`a $\mathfrak{g}_{tn}(F)$ de la transform\'ee de Fourier de la distribution $D(M,s_{M}, \varphi,w)$.
 Il r\'esulte des d\'efinitions de \ref{lespaceKstablegeneral} que, pour $H\in {\cal L}_{F}^{nr,st}(M,w)$, on a l'\'egalit\'e  $k^{G,st}(\varphi_{H})=c(w)^{-1}k(M,s_{M},\varphi,w)$. Donc $(k(M,s_{M},\varphi,w))_{(M,s_{M},\varphi,w)\in {\cal B}_{W}(G)}$ est une base de  $D^G({\cal D}^{st}(\mathfrak{g}(F)))$.
 
   Soit $(M,s_{M},\varphi,\rho)\in {\cal B}_{Irr}(G)$. On d\'efinit la distribution
$$k(M,s_{M},\varphi,\rho)=\vert W^{I_{F}}(M)\vert ^{-1}\sum_{w\in W^{I_{F}}(M)}trace(\rho^{\flat}(wFr^{-1})) k(M,s_{M},\varphi,w).$$
La famille $(k(M,s_{M},\varphi,\rho))_{(M,s_{M},\varphi,\rho)\in {\cal B}_{Irr}(G)}$ est encore une base de  $D^G({\cal D}^{st}(\mathfrak{g}(F)))$.

 De nouveau, on peut remplacer dans ces d\'efinitions le couple $(M,s_{M})$ par $(M_{\Lambda},s_{\Lambda})$ o\`u $\Lambda$ est le sous-ensemble de $\Delta_{a}^{nr}$ associ\'e \`a $(M,s_{M})$. On notera aussi $D(\Lambda,\varphi,w)$, $k(\Lambda,\varphi,w)$ et $k(\Lambda,\varphi,\rho)$ les distributions $D(M,s_{M},\varphi,w)$, $k(M,s_{M},\varphi,w)$ et $k(M,s_{M},\varphi,\rho)$ ci-dessus, modulo cette correspondance $(M,s_{M})\leftrightarrow \Lambda$. 
      
\subsubsection{Evaluation de la distribution $\hat{k}(\Lambda,\varphi,w)$\label{evaluationw}}
Soit $(\Lambda,\varphi,w)\in {\cal B}'_{W}(G)$. Pour simplifier la notation, on pose $M=M_{\Lambda}$. On note $Cl_{Fr}(w)$ la classe de $Fr$-conjugaison de $w$ dans $W^{I_{F}}(M)$. On pose
$$c(\Lambda,\varphi)=q^{dim(M_{s_{\Lambda}})/2-dim(\bar{{\cal O}}_{\varphi})/2}.$$
Soit $s\in S(\bar{C})$. On affecte d'indices $s$ les objets relatifs au groupe $G_{s}$. 
 Le sommet $s$ est associ\'e \`a une orbite   de l'action de $\Gamma_{F}^{nr}$ sur $\Delta_{a}^{nr}$. C'est le compl\'ementaire de l'ensemble $\Lambda(s)$, cf. \ref{constructiondeLevi}. Rappelons que $\Lambda$ est conserv\'e par l'action de $\Gamma_{F}^{nr}$, on a donc $\Lambda\subset \Lambda(s)$ ou $\Lambda\cup \Lambda(s)=\Delta_{a}^{nr}$.  
 Si $ \Lambda\subset \Lambda(s)$, on a $s\in Imm^G(M_{ad})$ et $p_{M}(s)=s_{\Lambda}$. Le groupe $M_{s_{\Lambda}}$ est alors un ${\mathbb F}_{q}$-Levi de $G_{s}$.  D'après le (iii) du lemme \ref{normalisateur}, de l'inclusion naturelle $Norm_{G(F^{nr})}(T^{nr},M,s_{\Lambda})\cap K_{s}^{nr,0}\to Norm_{G(F^{nr})}(M)$ se déduit un homomorphisme injectif $j_{s}:W_{s}(M_{s_{\Lambda}})\to W^{I_{F}}(M)$,  qui est compatible aux actions de Frobenius. Pour $w'\in W_{s}(M_{s_{\Lambda}})$, on note $w'\sim w$ si et seulement si $j_{s}(w')\in Cl_{Fr}(w)$.
 
 On rappelle que pour toute fonction $f$ sur $\mathfrak{g}(F)$ ou $\mathfrak{g}_{s}({\mathbb F}_{q})$, on note $f^-$ la fonction d\'efinie par $f^-(X)=f(-X)$. 
  
Soit $f\in {\bf C}(\mathfrak{g}_{s}({\mathbb F}_{q}))$.   On d\'eduit de $f$ une fonction sur $\mathfrak{g}(F)$ que l'on a souvent encore not\'ee $f$.  Pour plus de pr\'ecision, notons-la ici $f^G$. 

\begin{prop}{On suppose que $\hat{f}$ est \`a support nilpotent. 

(i) Si $ \Lambda\cup \Lambda(s)=\Delta_{a}^{nr}$, $\hat{k}(\Lambda,\varphi,w)(f^G)=0$.

(ii) Si $ \Lambda\subset \Lambda(s)$, 
$$\hat{k}(\Lambda,\varphi,w)(f^G)=c(\Lambda,\varphi)mes(K_{s}^0)\vert W_{s}(M_{s_{\Lambda}})\vert ^{-1}\vert W^{I_{F}}(M)\vert \vert Cl_{Fr}(w)\vert^{-1}$$
$$\sum_{w'\in W_{s}(M_{s_{\Lambda}}); w'\sim w}(Q_{M_{s_{\Lambda}},\varphi,w'}^-,f) .$$}\end{prop}

Preuve.  
Posons $f^{inv}=\vert G_{s}({\mathbb F}_{q})\vert ^{-1}\sum_{g\in G_{s}({\mathbb F}_{q})}{^gf}$. Remplacer $f$ par $f^{inv}$ ne change pas les  assertions de l'\'enonc\'e. En cons\'equence, on peut supposer $f$ invariante par conjugaison par $G_{s}({\mathbb F}_{q})$. Posons  $k=k(\Lambda,\varphi,w)$ et $D=D(\Lambda,\varphi,w)$. Par d\'efinition,
 $\hat{k}(f^G)=k(\hat{f}^G)$. Puisque $\hat{f}$ est \`a support nilpotent,  $\hat{f}^G$ est \`a support topologiquement nilpotent. Par d\'efinition de $k$,  on a $k(\hat{f}^G)=\hat{D}(\hat{f}^G)=D(f^{-,G})$. D'o\`u
 $$(1) \qquad \hat{k}(f^G)= c(w)\sum_{H\in {\cal L}_{F}^{nr,st}(M,w)}\vert W^{nr}(H)^{\Gamma_{F}^{nr}}\vert \vert W_{F}(H)\vert ^{-1}D^G_{X_{H},\varphi_{H}}(f^{-,G}).$$
 Fixons $H$ intervenant dans cette somme. Par d\'efinition,
 $$D^G_{X_{H},\varphi_{H}}(f^{-,G})=\int_{H(F)\backslash G(F)}X(g)\,dg,$$
 o\`u
 $$X(g)=\int_{A_{H}(F)\backslash H(F)} \int_{\mathfrak{h}_{SC}(F)}{^g(f^{-,G})}(X_{H}+h^{-1}Zh)\bar{\varphi}^{H_{SC}}_{H}(Z)\,dZ\,dh.$$
 Notons ${\cal X}$ l'ensemble des $g\in G(F)$ tels que $gs\in Imm^G(H)$ et $p_{H}(gs)=s_{H}$.  On va montrer
 
 (2) pour $g\in G(F)-H(F){\cal X}$, on a $X(g)=0$. 
 
 Fixons $g\in G(F)$  et supposons $X(g)\not=0$. L'application $Ad(g)$ transporte la fonction $f$ en une fonction $^gf\in C(\mathfrak{g}_{gs}({\mathbb F}_{q}))$ et on a $^g(f^{-G})=(^gf)^{-G}$. On note simplement $^gf^{-G}$ cette fonction. 
  Puisque $\varphi^{H_{SC}}_{H}$ est \`a support topologiquement nilpotent et que $X(g)\not=0$, le support de la fonction $^gf^{-,G}$  coupe $X_{H}+\mathfrak{h}_{SC,tn}(F)$. A fortiori, $\mathfrak{k}_{gs}$ coupe $X_{H}+\mathfrak{h}_{SC,tn}(F)$. D'apr\`es le lemme \ref{ImmGHad},  $gs$ appartient \`a $ Imm^G(H)$. Le sommet $gs$ se projette en un point $(gs)^H\in Imm(H_{AD})$ qui appartient \`a une facette ${\cal F}_{g}^H$ de cet immeuble. On a les inclusions $\mathfrak{h}_{SC,{\cal F}_{g}^H}\subset \mathfrak{h}_{{\cal F}_{g}^H}\subset \mathfrak{g}_{gs}$. 
  Notons $\bar{X}_{H}$ la projection de $X_{H}$ dans $\mathfrak{g}_{gs}$ et d\'efinissons une fonction $f_{1}\in C(\mathfrak{h}_{SC,{\cal F}_{g}^H})$ par $f_{1}(Z)={^gf^-}(\bar{X}_{H}+Z)$ pour tout $Z\in \mathfrak{h}_{SC,{\cal F}_{g}^H}({\mathbb F}_{q})$. On a alors l'\'egalit\'e  $^gf^{-,G}(X_{H}+Z)=f_{1}^{H_{SC}}(Z)$ pour tout $Z\in \mathfrak{h}_{SC,tn}(F)$.  On en d\'eduit
 $$X(g)=\int_{A_{H}(F)\backslash H(F)}\int_{\mathfrak{h}_{SC}(F)}f_{1}^{H_{SC}}(h^{-1}Zh)\bar{\varphi}^{H_{SC}}_{H}(Z)\, dZ\, dh.$$
   Pour $h\in H(F)$, d\'efinissons une fonction $f_{2,h}$ sur $\mathfrak{h}_{SC,s_{H}}({\mathbb F}_{q})$ par
 $$f_{2,h}(Z)=\int_{\mathfrak{k}^{H_{SC},+}_{s_{H}}}f_{1}^{H_{SC}}(h^{-1}(\dot{Z}+Z')h)\,dZ',$$
 o\`u $\dot{Z}$ est un rel\`evement de $Z$ dans $\mathfrak{k}^{H_{SC}}_{s_{H}}$. On a l'\'egalit\'e
 $$\int_{\mathfrak{h}_{SC}(F)}f_{1}^{H_{SC}}(h^{-1}Zh)\bar{\varphi}^{H_{SC}}_{H}(Z)\, dZ=
 \sum_{Z\in \mathfrak{h}_{SC,s_{H}}({\mathbb F}_{q})}f_{2,h}(Z)\bar{\varphi}_{H}(Z).$$
 Posons ${\cal F}'=h({\cal F}_{g}^H)$. Par construction, la fonction $f_{2,h}$ est support\'ee par l'image dans $\mathfrak{h}_{SC,s_{H}}({\mathbb F}_{q})$  de $\mathfrak{k}_{{\cal F}'}^{H_{SC}}\cap \mathfrak{k}_{s_{H}}^{H_{SC}}$ et est invariante par l'image de $\mathfrak{k}_{{\cal F}'}^{H_{SC},+}\cap \mathfrak{k}_{s_{H}}^{H_{SC}}$. Si ${\cal F}'$ n'est pas \'egale \`a $\{s_{H}\}$,  il existe un sous-groupe parabolique propre $P$  d\'efini sur ${\mathbb F}_{q}$ de $H_{SC,s_{H}}$ tel que ces deux ensembles soient respectivement $P({\mathbb F}_{q})$ et $U_{P}({\mathbb F}_{q})$.  La cuspidalit\'e de $\varphi_{H}$ entra\^{\i}ne la nullit\'e de la somme ci-dessus. Puisque $X(g)\not=0$ par hypoth\`ese, il existe $h\in H(F)$ tel que  cette somme  soit non nulle, donc tel que $h({\cal F}_{g}^H)=s_{H}$. Mais alors $hg\in {\cal X}$.  Cela prouve (2). 
 
  Puisque l'application $g\mapsto X(g)$ est invariante \`a gauche par $H(F)$,  il nous suffit de calculer $X(g)$ pour $g\in {\cal X}$. Soit donc $g\in {\cal X}$. 
 Posons $H'= g^{-1}Hg$ et notons $s_{H'}$ l'image de $s_{H}$ dans $Imm(H'_{AD})$ par $Ad(g)^{-1}$. Alors $s\in Imm^G(H')$ et la projection de $s$ dans $Imm(H'_{AD})$ est $s_{H'}$. Le (i) du lemme \ref{HsH} dit que  $\Lambda\subset \Lambda(s)$.  Si cette condition n'est pas v\'erifi\'ee, on a donc ${\cal X}=\emptyset$, d'o\`u $\hat{k}(\Lambda,\varphi,w)(f^G)=0$ gr\^ace \`a (2). 
 Cela d\'emontre le (i) de l'\'enonc\'e. 

 Supposons $\Lambda\subset \Lambda(s)$. Le (ii) du lemme \ref{HsH} affirme l'existence de $k\in K_{s}^{nr,0}$ tel que $Ad(k)$ transporte $(H',s_{H'})$ sur $(M_{\Lambda},s_{\Lambda})=(M,s_{\Lambda})$. En notant $\bar{k}$ la r\'eduction de $k$ dans $G_{s}$, on a aussi $Ad(\bar{k})(H'_{s_{H'}})=M_{s_{\Lambda}}$. Fixons un tel $k$. L'\'el\'ement $Fr(\bar{k})\bar{k}^{-1}$ de $G_{s}$ normalise $M_{s_{\Lambda}}$ et d\'efinit un \'el\'ement $w'\in W_{s}(M_{s_{\Lambda}})$. Cet \'el\'ement d\'epend de l'\'el\'ement $k$ choisi mais sa classe de $Fr$-conjugaison n'en d\'epend pas. Posons $x=kg^{-1}$. On a $xHx^{-1}=M$. Puisque $H$ est param\'etr\'e par l'image de $w$ dans  $W^{I_{F}}(M)/Fr-conj$, l'image de $Fr(x)x^{-1}$ dans $W^{I_{F}}(M)$ est $Fr$-conjugu\'e \`a $w$. Cela \'equivaut \`a la relation $w'\sim w$. 
  L'isomorphisme $Ad(g)^{-1}$ transporte $\varphi_{H}$ en $\varphi_{H'}\in C(\mathfrak{h}'_{SC,s_{H'}}({\mathbb F}_{q}))$. Posons $X_{H'}=g^{-1}X_{H}g$. C'est un \'el\'ement de $A_{H'}^{nr}(F)$ entier et de r\'eduction r\'eguli\`ere.  Un changement de variables \'evident conduit \`a l'\'egalit\'e
 $$X(g)=\int_{A_{H'}(F)\backslash H'(F)}\int_{\mathfrak{h}'_{SC}(F)}f^{-,G}(X_{H'}+h^{-1}Zh)\bar{\varphi}_{H'}^{H'_{SC}}(Z)\,dZ\,dh.$$
 En reprenant la preuve de (2), on voit que l'int\'egrale int\'erieure est support\'ee par le groupe des $h'\in H'(F)$ tels que $h's_{H'}=s_{H'}$, c'est-\`a-dire par $K_{s_{H'}}^{H',\dag}$. Cette int\'egrale est constante sur cet ensemble car la fonction $\varphi_{H'}$ est invariante par ce groupe. On obtient
 $$X(g)=mes(A_{H'}(F)\backslash K_{s_{H'}}^{H',\dag})\int_{\mathfrak{h}'_{SC}(F)}f^{-,G}(X_{H'}+Z)\bar{\varphi}_{H'}^{H'_{SC}}(Z)\,dZ$$
 $$=mes(A_{H'}(F)\backslash K_{s_{H'}}^{H',\dag})mes(k_{s_{H'}}^{H'_{SC},+})\sum_{\bar{Z}\in \mathfrak{h}'_{SC,s_{H'}}({\mathbb F}_{q})}f^-(\bar{X}_{H'}+\bar{Z})\bar{\varphi}_{H'}(\bar{Z}),$$
 o\`u $\bar{X}_{H'}$ est la r\'eduction de $X_{H'}$.

 La fonction $\varphi$ est la fonction caract\'eristique du faisceau caract\`ere ${\cal E}_{\varphi}$ sur $\mathfrak{m}_{SC,s_{\Lambda}}$. Il correspond \`a ce faisceau un tel faisceau ${\cal E}_{\varphi,H'}$ sur $\mathfrak{h}'_{SC,s_{H'}}$ muni d'une action de Frobenius d\'efinie comme en \ref{fonctionscaracteristiques}. Notons $\varphi_{w'}$ sa fonction caractéristique. D'après le lemme \ref{HsHderechef}, on a l'égalité  $\varphi_{H'}= \vert M_{AD,s_{\Lambda}}({\mathbb F}_{q})\vert \vert H'_{AD, s_{H'}}({\mathbb F}_{q})\vert ^{-1}\varphi_{w'}$.  En introduisant la fonction $\varphi_{w',\bar{X}_{H'}}^{\flat}$ de \ref{fonctionscaracteristiques} associ\'ee \`a l'\'el\'ement $\bar{X}_{H'}$,  on obtient
 $$X(g)=c_{1}(H')\sum_{\bar{X}\in \mathfrak{g}_{s}({\mathbb F}_{q})}f^{-}(\bar{X})\bar{\varphi}_{w',\bar{X}_{H'}}^{\flat}(\bar{X}),$$
 o\`u 
 $$ c_{1}(H')= mes(A_{H'}(F)\backslash K_{s_{H'}}^{H',\dag})mes(k_{s_{H'}}^{H'_{SC},+})\vert M_{AD,s_{\Lambda}}({\mathbb F}_{q})\vert \vert H'_{AD, s_{H'}}({\mathbb F}_{q})\vert ^{-1}.$$
  Puisqu'on a suppos\'e $f$ invariante par $G_{s}({\mathbb F}_{q})$, on peut aussi bien moyenner $\varphi_{w',\bar{X}_{H'}}^{\flat}$ par ce groupe. Cela remplace cette fonction par $\vert H'_{s_{H'}}({\mathbb F}_{q})\vert \vert G_{s}({\mathbb F}_{q})\vert ^{-1}\varphi_{w',\bar{X}_{H'}}$, o\`u $\varphi_{w',\bar{X}_{H'}}$ est d\'efinie en \ref{fonctionscaracteristiques}. On obtient alors
 $$X(g)=  c_{2}(H')(\varphi_{w',\bar{X}_{H'}},f^-).$$
 o\`u 
 $$(3) \qquad c_{2}(H')=\vert H'_{s_{H'}}({\mathbb F}_{q})\vert c_{1}(H')= mes(A_{H'}(F)\backslash K_{s_{H'}}^{H',\dag})mes(k_{s_{H'}}^{H'_{SC},+})$$
 $$\vert M_{AD, s_{\Lambda}}({\mathbb F}_{q})\vert \vert H'_{AD, s_{H'}}({\mathbb F}_{q})\vert ^{-1} \vert H'_{s_{H'}}({\mathbb F}_{q})\vert.$$
Notons simplement $Q_{w'}$ la fonction $Q_{M_{s_{\Lambda}},{\cal E}_{\varphi},w'}$.  Puisque $\hat{f}$ est \`a support nilpotent par hypoth\`ese, la relation \ref{fonctionscaracteristiques}(5) implique  la deuxi\`eme \'egalit\'e ci-apr\`es 
 $$(\varphi_{w',\bar{X}_{H'}},f^-)=(\hat{\varphi}_{w',\bar{X}_{H'}},\hat{f}^-)=(\hat{Q}_{w'},\hat{f}^-)=(Q_{w'},f^-)=(Q_{w'}^-,f).$$  
   D'o\`u
 $$(4) \qquad X(g)=c_{2}(H)(Q^-_{w'},f).$$
  
 Notons plus pr\'ecis\'ement $w'_{g}$ l'\'el\'ement $w'$ de la construction ci-dessus. Pour $w'\in W_{s}(M_{s_{\Lambda}})$, notons ${\cal X}(w')$ l'ensemble des $g\in {\cal X}$  tels que  $w'_{g}$   soit $Fr$-conjugu\'e  dans   $W_{s}(M_{s_{\Lambda}})$ \`a  $w'$. On a vu que, si ${\cal X}(w')$ n'est pas vide, alors $w'\sim w$. 
  L'ensemble ${\cal X}$ est r\'eunion disjointe des ${\cal X}(w')$ quand $w'$ d\'ecrit $W_{s}(M_{s_{\Lambda}})/Fr-conj$. On remarque que, pour  $g\in {\cal X}$ et $h\in H(F)$ tel que $hg\in {\cal X}$, les \'el\'ements $w'_{g}$  et $w'_{hg}$ sont $Fr$-conjugu\'es: ils ne d\'ependent que du groupe  $H'$ qui ne change pas quand on remplace $g$ par $hg$ (on a $H'=g^{-1}Hg$). Donc $H(F){\cal X}$ est aussi l'union disjointe des ensembles $H(F){\cal X}(w')$. D'apr\`es sa d\'efinition, l'ensemble ${\cal X}$ est invariant par multiplication \`a droite par $K_{s}^0$. Une telle multiplication remplace $H'_{s_{H'}}$ par un ensemble conjugu\'e par un   \'el\'ement de $G_{s}({\mathbb F}_{q})$ et ne modifie pas  l'\'el\'ement $w'$. En cons\'equence, $ {\cal X}(w')$ est invariant par multiplication \`a droite par   $K_{s}^0$. Il en r\'esulte que $H(F){\cal X}(w')$ est  ouvert et ferm\'e.  Pour tout \'el\'ement $g\in H(F){\cal X}(w')$, $X(g)$ est calcul\'e par la formule (4). On obtient
   $$(5) \qquad D_{X_{H},\varphi_{H}^G}(f^{-,G})=c_{2}(H)\sum_{w'\in W_{s}(M_{s_{\Lambda}})/Fr-conj}mes(H(F)\backslash H(F){\cal X}(w'))(Q_{w'}^-,f).$$
   
   Fixons un \'el\'ement 
    $w'\in W_{s}(M_{s_{\Lambda}})$. Fixons un $\bar{{\mathbb F}}_{q}$-Levi $M_{s_{\Lambda},w'}$ de $G_{s}$, d\'efini sur ${\mathbb F}_{q}$,  conjugu\'e \`a $M_{s_{\Lambda}}$ et param\'etr\'e par la classe de $Fr$-conjugaison de $w'$.  On  rel\`eve $M_{s_{\Lambda},w'}$  en un $F^{nr}$-Levi de $G$ que l'on note $H_{w'}$.  Notons $Cl_{s,Fr}(w')$ la classe de $Fr$-conjugaison de $w'$. 
        On va prouver la relation 
    
    (6) si $H_{w'}$ n'est pas conjugu\'e \`a $H$ par un \'el\'ement de $G(F)$, ${\cal X}(w')=\emptyset$;
    si $H_{w'}$ est conjugu\'e \`a $H$ par un \'el\'ement de $G(F)$,
   $$mes(H(F)\backslash H(F){\cal X}(w'))=\vert W_{F}(H)\vert \vert  W_{s}(M_{s_{\Lambda}})\vert ^{-1} \vert Cl_{s,Fr}(w')\vert $$
   $$  \vert H_{s_{H}}({\mathbb F}_{q})\vert^{-1} q^{dim(M_{s_{M}})/2}mes (K_{s}^0).$$
     
     Le sommet $s$ d\'etermine un sommet de $Imm(H_{w',AD})$ que l'on note $s_{H_{w'}}$. On a $H_{w',s_{H_{w'}}}=M_{s_{\Lambda},w'}$.  Notons ${\cal Y}(w')$  l'ensemble des $g\in G(F)$ tels qu'il existe $y\in K_{s}^{0}$ de sorte que $Ad(gy)(H_{w'})=H$. 
Montrons qu'on a l'\'egalit\'e
   
   (7) $H(F){\cal X}(w')={\cal Y}(w')$.  
      
   Supposons $g\in  {\cal X}(w')$ et reprenons les constructions pr\'ec\'edentes, en particulier celle du  groupe $H'$. Alors $H'_{s_{H'}}$ est conjugu\'e \`a $ H_{w',s_{H_{w'}}}$ par un \'el\'ement de $G_{s}({\mathbb F}_{q})$. Par un proc\'ed\'e de rel\`evement  d\'ej\`a plusieurs fois utilis\'e, il existe $y\in K_{s}^0$ tel que $Ad(y)(H_{w'})=H'$. Puisque $Ad(g)(H')=H$, cela entra\^{\i}ne que $g\in {\cal Y}(w')$. D'autre part, ${\cal Y}(w')$ est invariant \`a gauche par $H(F)$. Donc $H(F){\cal X}(w')\subset {\cal Y}(w')$.   Inversement, soit $g\in {\cal Y}(w')$. Fixons 
    $y\in K_{s}^0$ de sorte que $Ad(gy)(H_{w'})=H$. Posons $H'_{w'}=Ad(y)(H_{w'})$, qui v\'erifie les m\^emes conditions que $H_{w'}$. On a $s\in Imm^G(H'_{w'})$. Notons $s_{H'_{w'}}$ l'image de $s$ dans $Imm(H'_{w',AD})$ et $s'_{H}$ l'image de $s_{H'_{w'}}$ dans $Imm(H_{AD})$ par $Ad(g)$. Les couples suivants sont conjugu\'es par un \'el\'ement de $G(F^{nr})$:
   
   $(H,s_{H})$,  $(M,s_{\Lambda})$, $(H'_{w'},s_{H'_{w'}})$, $(H,s'_{H})$.
   
   Il existe donc un \'el\'ement $n\in Norm_{G(F^{nr})}(H)$ tel que l'action de $n$ sur $Imm_{F^{nr}}(H_{AD})$ envoie $s_{H}$ sur $s'_{H}$. D'apr\`es le lemme \ref{lemmeautomorphismes} appliqu\'e au groupe $H_{SC}$, il existe $h\in H(F)$ (en fait provenant d'un \'el\'ement de $H_{SC}(F)$) tel que $hs'_{H}=s_{H}$. Posons $g'=hgy$.  Alors  $Ad(g')$ transporte le couple $(H_{w'},s_{H_{w'}})$ sur $(H,s_{H})$. Puisque $s\in Imm^G(H_{w'})$ et $p_{H_{w'}}(s)= s_{H_{w'}}$, on a $g's\in Imm^G(H)$ et $p_{H}(s)=s_{H}$. Donc $g'\in {\cal X}$. Le groupe $H'$ de notre construction  associ\'e \`a cet \'el\'ement $g'$ est $H_{w'}$. Donc $g'\in {\cal X}(w')$. On a d\'ej\`a remarqu\'e que cet ensemble \'etait invariant par multiplication \`a droite par $K_{s}^0$. Donc $g=h^{-1}g'y^{-1}\in H(F){\cal X}(w')$. Cela prouve (7).

 Si $H_{w'}$ n'est pas conjugu\'e \`a $H$ par un \'el\'ement de $G(F)$, ${\cal Y}(w')$ est vide donc ${\cal X}(w')$ aussi, ce qui d\'emontre la premi\`ere assertion de (6). Supposons que $H_{w'}$ soit conjugu\'e \`a $H$ par un \'el\'ement de $G(F)$, fixons $x\in G(F)$ tel que $Ad(x)(H_{w'})=H$.  Par d\'efinition, on a $ {\cal Y}(w')=Norm_{G(F)}(H)xK_{s}^0$. Alors $H(F)\backslash  {\cal Y}(w')=H(F)\backslash Norm_{G(F)}(H)xK_{s}^0\simeq H_{w'}(F)\backslash Norm_{G(F)}(H_{w'})K_{s}^0$. L'ensemble  $H_{w'}(F)\backslash Norm_{G(F)}(H_{w'})K_{s}^0/K_{s}^0$ est en bijection avec $W_{F}(H_{w'})/W_{s,{\mathbb F}_{q}}(H_{w',s_{H_{w'}}})$, o\`u $W_{F}(H_{w'})=H_{w'}(F)\backslash Norm_{G(F)}(H_{w'})$ et 
 $$W_{s,{\mathbb F}_{q}}(H_{w',s_{H_{w'}}})=  H_{w',s_{H_{w'}}}({\mathbb F}_{q})\backslash Norm_{G_{s}({\mathbb F}_{q})}(H_{w',s_{H_{w'}}})$$
 $$\simeq (H_{w'}(F)\cap K_{s}^0)) \backslash (Norm_{G(F)}(H_{w'})\cap K_{s}^0).$$
  Chaque double classe est isomorphe \`a $(K_{s}^0\cap H_{w'}(F))\backslash K_{s}^0=K_{s_{H_{w'}}}^{H_{w'},0}\backslash K_{s}^0$. On a $W_{F}(H_{w'})\simeq W_{F}(H)$ puisque $H$ et $H_{w'}$ sont conjugu\'es par un \'el\'ement de $G(F)$. On obtient
 $$mes(H(F)\backslash {\cal Y}(w'))=\vert W_{F}(H)\vert \vert W_{s,{\mathbb F}_{q}}(H_{w',s_{H_{w'}}})\vert ^{-1}mes(K_{s_{H_{w'}}}^{H_{w'},0})^{-1}mes(K_{s}^0).$$
 Le th\'eor\`eme de Lang appliqu\'e \`a $H_{w',s_{H_{w'}}}$ entra\^{\i}ne que $W_{s,{\mathbb F}_{q}}(H_{w',s_{H_{w'}}})$ s'identifie \`a l'ensemble de points fixes de l'action de $Fr$ dans $H_{w',s_{H_{w'}}}\backslash Norm_{G_{s}}(H_{w',s_{H_{w'}}})$. Ce quotient s'identifie \`a $W_{s}(M_{s_{\Lambda}})$ muni d'une action de Frobenius tordue par $w'$. On voit alors que  $W_{s,{\mathbb F}_{q}}(H_{w',s_{H_{w'}}})$  s'identifie au groupe  des $u\in W_{s}(M_{s_{\Lambda}})$ tels que $Fr(u) w' u^{-1}=w'$. Notons $Z_{W_{s} (M_{s_{\Lambda}}),Fr}(w')$ ce groupe.  L'application $u\mapsto Fr(u)w' u^{-1}$ est une bijection de $W_{s}(M_{s_{\Lambda}})/Z_{W_{s} (M_{s_{\Lambda}}),Fr}(w')$ sur $Cl_{s,Fr}(w')$. On en d\'eduit
 $$(8) \qquad \vert W_{s,{\mathbb F}_{q}}(H_{w',s_{H_{w'}}})\vert=\vert W_{s}(M_{s_{\Lambda}})\vert \vert Cl_{s,Fr}(w')\vert ^{-1}.$$ D'autre part, on a $K_{s_{H_{w'}}}^{H_{w'},0}\simeq K_{s_{H}}^{H,0}$ et
 $$mes(K_{s_{H}}^{H,0})= \vert  H_{s_{H}}({\mathbb F}_{q})\vert  q^{-dim(H_{s_{H}})/2}=\vert  H_{s_{H}}({\mathbb F}_{q})\vert q^{-dim(M_{s_{\Lambda}})/2}.$$
 D'o\`u
  $$mes(H(F)\backslash {\cal Y}(w'))=\vert W_{F}(H)\vert \vert W_{s}(M_{s_{\Lambda}})\vert ^{-1}\vert \vert Cl_{s,Fr}(w')\vert \vert H_{s_{H}}({\mathbb F}_{q})\vert^{-1} q^{dim(M_{s_{M}})/2}mes (K_{s}^0).$$
  Cela d\'emontre la deuxi\`eme assertion de (6). 
  
  Notons $W_{s}(M_{s_{\Lambda}})[H]$ l'ensemble des $w'\in W_{s}(M_{s_{\Lambda}})$ tels que les groupes $H_{w'}$  et $H$ soient conjugués par un \'el\'ement de $G(F)$. Posons 
  $$(9) \qquad c_{3}(H)= c_{2}(H)  \vert W_{s}(M_{s_{\Lambda}})\vert ^{-1}\vert   H_{s_{H}}({\mathbb F}_{q})\vert^{-1} q^{dim(M_{s_{M}})/2}mes (K_{s}^0).$$
   
  En utilisant (6), l'\'egalit\'e (5) se transforme en
  $$D_{X_{H},\varphi_{H}^G}(f^{-,G})=c_{3}(H)\vert W_{F}(H)\vert \sum_{w'\in W_{s}(M_{s_{\Lambda}})[H]/Fr-conj }\vert  Cl_{s,Fr}(w')\vert  (Q_{w'}^-,f).$$
  Sommer sur les classes de $Fr$-conjugaison de $w'$ un terme affect\'e du coefficient $\vert  Cl_{s,Fr}(w')\vert $ revient \`a sommer sur $w'$ le m\^eme terme o\`u l'on supprime ce coefficient. D'o\`u
  $$(10) \qquad D_{X_{H},\varphi_{H}^G}(f^{-,G})=c_{3}(H)\vert W_{F}(H)\vert \sum_{w'\in W_{s}(M_{s_{\Lambda}})[H]}   (Q_{w'}^-,f).$$
  
  Posons
  $$c_{3}=  m(w) \vert M_{AD,s_{\Lambda}}({\mathbb F}_{q})\vert  \vert W_{s}(M_{s_{\Lambda}})\vert ^{-1}  q^{dim(M_{s_{\Lambda}})/2-dim(M_{SC,s_{\Lambda}})}mes (K_{s}^0).$$
  Montrons que
  
  (11) $c_{3}(H)= c_{3}$.
  
  Notons $\pi:H_{SC}\to H$ l'homomorphisme naturel. Montrons d'abord que
  
  (12) $\pi(H_{SC}(F))K_{s_{H}}^{H,\dag}=H(F)$.
  
  Soit $h\in H(F)$. D'apr\`es le lemme \ref{lemmeautomorphismes} appliqu\'e \`a $H_{SC}$, il existe $h'_{sc}\in H_{SC}(F)$ tel que $\pi(h'_{sc})s_{H}=hs_{H}$. Donc $\pi(h'_{sc})^{-1}h\in K_{s_{H}}^{H,\dag}$, ce qui d\'emontre (12). 
  
  On a rappel\'e plus haut la mesure $m(H)$ introduite en \ref{uncalculdemesures}. Par d\'efinition, ce nombre est la mesure totale de $A_{H}(F) \pi(H_{SC}(F))\backslash H(F)$, cet ensemble \'etant muni de la mesure invariante \`a droite telle que, pour toute fonction $\phi$ sur $A_{H}(F)\backslash H(F)$, on ait l'\'egalit\'e
  $$\int_{A_{H}(F) \pi(H_{SC}(F))\backslash H(F)}\int_{H_{SC}(F)}\phi(\pi(h_{sc})x)\,dh_{sc}\,dx=\int_{A_{H}(F)\backslash H(F)}\phi(h)\,dh.$$
  Appliquons cela \`a la fonction caract\'eristique $\phi$ du groupe $A_{H}(F)\backslash K_{s_{H}}^{H,\dag}$. A droite, on  obtient $mes(A_{H}(F)\backslash K_{s_{H}}^{H,\dag})$. A gauche, l'assertion (11) montre que l'int\'egrale int\'erieure ne d\'epend pas de $x$. Puisque l'image r\'eciproque par $\pi$ de $K_{s_{H}}^{H,\dag}$ est $K_{s_{H}}^{H_{SC},0}$, cette int\'egrale vaut $mes(K_{s_{H}}^{H_{SC},0})$. Le membre de gauche vaut donc $m(H)mes(K_{s_{H}}^{H_{SC},0})$, d'o\`u l'\'egalit\'e
  $$mes(A_{H}(F)\backslash K_{s_{H}}^{H,\dag})=m(H)mes(K_{s_{H}}^{H_{SC},0}).$$
  On a aussi
  $$mes(K_{s_{H}}^{H_{SC},0})=\vert H_{SC,s_{H}}({\mathbb F}_{q})\vert mes(\mathfrak{k}_{s_{H}}^{H_{SC},+}).$$
  Enfin, l'homomorphisme $H_{SC,s_{H}}\to H_{AD,s_{H}}$ est surjectif et de noyau fini central. On en d\'eduit l'\'egalit\'e
  $$\vert H_{SC,s_{H}}({\mathbb F}_{q})\vert =\vert H_{AD,s_{H}}({\mathbb F}_{q})\vert .$$
  En appliquant ces \'egalit\'es et les d\'efinitions (9) et (3), on obtient
  $$c_{3}(H)=  m(H) mes(k_{s_{H}}^{H_{SC},+})^2\vert M_{AD,s_{\Lambda}}({\mathbb F}_{q})\vert  \vert  W_{s}(M_{s_{\Lambda}})\vert ^{-1}  q^{dim(M_{s_{\Lambda}})/2}mes (K_{s}^0).$$
  On a aussi $mes(k_{s_{H}}^{H_{SC},+})^2=q^{-dim(H_{SC,s_{H}})}=q^{-dim(M_{SC,s_{\Lambda}})}$ et  $m(w)=m(H)$ par d\'efinition de $m(w)$. 
   Alors l'\'egalit\'e ci-dessus devient (11). 
  
  Rempla\c{c}ons $c_{3}(H)$ par $c_{3}$ dans l'\'egalit\'e (10) et reportons l'\'egalit\'e obtenue dans l'expression (1). Remarquons que les termes $\vert W^{nr}(H)^{\Gamma_{F}^{nr}}\vert $ intervenant dans cette expression sont en fait ind\'ependants de $H$ puisque tous les $H$ intervenant sont stablement conjugu\'es. En fait, pour $H\in {\cal L}_{F}^{nr,st}(M,w)$, le groupe $W^{nr}(H)$ muni de son action de Frobenius s'identifie \`a $W^{I_{F}}(M)$ muni de l'action tordue $u\mapsto w^{-1} Fr(u)w$. Comme plus haut, cf. \'egalit\'e (8), on en d\'eduit
  $$\vert W^{nr}(H)^{\Gamma_{F}^{nr}}\vert =\vert W^{I_{F}}(M)\vert \vert Cl_{Fr}(w)\vert ^{-1}.$$
  L'expression (1) devient
  $$\hat{k}(f^G)=c(w)c_{3}\vert W^{I_{F}}(M)\vert\vert Cl_{Fr}(w)\vert ^{-1}\sum_{H\in {\cal L}^{nr,st}_{F}(M,w)} \sum_{w'\in W_{s}(M_{s_{\Lambda}})[H]}   (Q_{w'}^-,f).$$
  Pour $w'\in W_{s}(M_{s_{\Lambda}})$, l'existence de $H\in {\cal L}^{nr,st}_{F}(M,w)$ tel que $w'\in W_{s}(M_{s_{\Lambda}})[H]$ implique $w'\sim w$. Inversement, si $w'\sim w$, il existe un unique $H\in {\cal L}^{nr,st}_{F}(M,w)$ tel que $H_{w'}$ soit conjugu\'e \`a $H$ par un \'el\'ement de $G(F)$, c'est-\`a-dire tel que $w'\in W_{s}(M_{s_{\Lambda}})[H]$. On obtient alors
  $$\hat{k}(f^G)=c(w)c_{3}\vert W^{I_{F}}(M)\vert\vert Cl_{Fr}(w)\vert ^{-1} \sum_{w'\in W_{s}(M_{s_{\Lambda}}); w'\sim w}(Q^-_{w'},f).$$
 En explicitant les constantes, on voit que $c(w)c_{3}=c(\Lambda,\varphi)mes(K_{s}^0)\vert W_{s}(M_{s_{\Lambda}})\vert ^{-1}$. L'\'egalit\'e ci-dessus devient celle du 
 (ii) de l'\'enonc\'e. $\square$
 
 \subsubsection{Evaluation de la distribution $\hat{k}(\Lambda,\varphi,\rho)$\label{evaluationrho}}
 Soit $(\Lambda,\varphi,\rho)\in {\cal B}'_{Irr}(G)$ et soit $s\in S(\bar{C})$. On utilise les m\^emes notations que dans le paragraphe pr\'ec\'edent.

  \begin{prop}{Soit $f\in {\bf C}(\mathfrak{g}_{s}({\mathbb F}_{q}))$. On suppose que $\hat{f}$ est \`a support nilpotent.    
  
  (i) Si $ \Lambda\cup \Lambda(s)=\Delta_{a}^{nr}$, $\hat{k}(\Lambda,\varphi,\rho)(f^G)=0$.
  
  (ii) Si $ \Lambda\subset \Lambda(s)$,
  $$\hat{k}(\Lambda,\varphi,\rho)(f^G)=c(\Lambda,\varphi)mes(K_{s}^0)\vert W_{s}(M_{s_{\Lambda}})\vert ^{-1}\sum_{w'\in W_{s}(M_{s_{\Lambda}})}trace(\rho^{\flat}(j_{s}(w')Fr^{-1}))(Q^-_{M_{s_{\Lambda}},\varphi,w'},f).$$}\end{prop}
  
  Preuve. Posons $k=k(\Lambda,\varphi,\rho)$. Par  d\'efinition, on a
  $$\hat{k}(f^G)=\vert W^{I_{F}}(M)\vert^{-1}\sum_{w\in W^{I_{F}}(M)}trace(\rho^{\flat}(wFr^{-1}))\hat{k}(\Lambda,\varphi,w)(f^G).$$
  On applique la proposition pr\'ec\'edente. Le (i) de l'\'enonc\'e est imm\'ediat. Supposons $ \Lambda\subset \Lambda(s)$. Alors
  $$\hat{k}(f^G)=c\sum_{w\in W^{I_{F}}(M)} \vert Cl_{Fr}(w)\vert ^{-1}trace(\rho^{\flat}(wFr^{-1}))\sum_{w'\in W_{s}(M_{s_{\Lambda}});w'\sim w }(Q^-_{M_{s_{\Lambda}},\varphi,w'},f),$$
  o\`u
  $$c=c(\Lambda,\varphi)mes(K_{s}^0)\vert W_{s}(M_{s_{\Lambda}})\vert ^{-1}.$$
  Pour tout $w'\in W_{s}(M_{s_{\Lambda}})$, l'ensemble des $w\in W^{I_{F}}(M)$ tels que $w'\sim w$ est $Cl_{Fr}(j_{s}(w'))$. Pour tout \'el\'ement $w$ de cet ensemble, on a $trace(\rho^{\flat}(wFr^{-1}))=trace(\rho^{\flat}(j_{s}(w')Fr^{-1}))$. Alors la formule ci-dessus se transforme en celle de l'\'enonc\'e. $\square$

  \subsubsection{Calcul de $\hat{k}(\Lambda,\varphi,\rho)( h_{s,{\cal O}})$\label{calculcrucial}}
  Soient $(\Lambda,\varphi,\rho)\in {\cal B}'_{Irr}(G)$, $s\in S(\bar{C})$ et ${\cal O}\in \mathfrak{g}_{s,nil}({\mathbb F}_{q})/conj$.  On pose $M=M_{\Lambda}$ et on  fixe $N\in {\cal O}$. {\bf On suppose que} $ \Lambda\subset \Lambda(s)$.
  
  Posons $W_{s}(M_{s_{\Lambda}})^+=W_{s}(M_{s_{\Lambda}})\rtimes \Gamma_{{\mathbb F}_{q}}$ et notons $\tilde{W}_{s}(M_{s_{\Lambda}})$ son sous-ensemble $W_{s}(M_{s_{\Lambda}})\rtimes\{Fr^{-1}\}$. 
   Consid\'erons une repr\'esentation  continue $\pi^+$ de $W_{s}(M_{s_{\Lambda}})^+$, qui n'est pas suppos\'ee irr\'eductible (rappelons que la continuit\'e signifie que $\pi^+$ est triviale sur un sous-groupe ouvert de $\Gamma_{{\mathbb F}_{q}}$). Notons $\pi$ sa restriction \`a $W_{s}(M_{s_{\Lambda}})$. On d\'ecompose ces repr\'esentations en sommes de repr\'esentations irr\'eductibles
   $$\pi^+=\sum_{\sigma^+\in Irr(W_{s}(M_{s_{\Lambda}})^+)}m(\pi^+,\sigma^+)\sigma^+,$$
    $$\pi=\sum_{\rho'\in Irr(W_{s}(M_{s_{\Lambda}}))}m(\pi,\rho')\rho',$$
    avec des multiplicit\'es $m(\pi^+,\sigma^+),m(\pi,\rho')\in {\mathbb N}$. Consid\'erons une repr\'esentation $\sigma^+\in Irr(W_{s}(M_{s_{\Lambda}})^+)$. Si sa restriction $\sigma^+_{\vert  W_{s}(M_{s_{\Lambda}})}$ n'est pas irr\'eductible, la trace de $\sigma^+$ est nulle sur $\tilde{W}_{s}(M_{s_{\Lambda}})$. Supposons que $\sigma^+_{\vert  W_{s}(M_{s_{\Lambda}})}$ soit irr\'eductible. Notons $\rho'$ cette repr\'esentation. On a $\rho'\in Irr_{{\mathbb F}_{q}}(W_{s}(M_{s_{\Lambda}}))$. On a fix\'e un prolongement $\rho^{'\flat}$ de $\rho'$ \`a $W_{s}(M_{s_{\Lambda}})^+$. La repr\'esentation $\sigma^+$ n'est pas forc\'ement \'egale \`a ce prolongement mais il existe une racine de l'unit\'e $\zeta(\sigma^+)\in {\mathbb C}^{\times}$ de sorte que $trace(\sigma^+)$ co\"{\i}ncide sur $\tilde{W}_{s}(M_{s_{\Lambda}})$ avec $\zeta(\sigma^+)trace(\rho^{'\flat})$. Pour tout $\rho'\in  Irr_{{\mathbb F}_{q}}(W_{s}(M_{s_{\Lambda}}))$, posons
    $$m(\pi^+,\rho^{'\flat})=\sum_{\sigma^+}m(\pi^+,\sigma^+)\zeta(\sigma^+),$$
    o\`u on somme sur les $\sigma^+\in Irr(W_{s}(M_{s_{\Lambda}})^+)$ tels que $\sigma^+_{\vert  W_{s}(M_{s_{\Lambda}})}=\rho'$. C'est un nombre complexe. Remarquons que, si l'on supprime les racines de l'unit\'e $\zeta(\sigma^+)$ dans la somme ci-dessus, on obtient $m(\pi,\rho')$. On en d\'eduit:

    (1) $\vert m(\pi^+,\rho^{'\flat})\vert \leq m(\pi,\rho')$ et $\vert m(\pi^+,\rho^{'\flat}\vert =1$ si $m(\pi,\rho')=1$.
    
    Notre raisonnement conduit \`a la relation
    
    (2) $trace(\pi^+)$ co\"{\i}ncide sur $\tilde{W}_{s}(M_{s_{\Lambda}})$  avec 
    $$\sum_{\rho'\in Irr_{{\mathbb F}_{q}}(W_{s}(M_{s_{\Lambda}}))}m(\pi^+,\rho^{'\flat})trace(\rho^{'\flat}).$$
    
    Notons $\rho_{s}$, resp. $\rho_{s}^{\flat}$, les repr\'esentations de $W_{s}(M_{s_{\Lambda}})$, resp. $W_{s}(M_{s_{\Lambda}})^+$, d\'efinies par $\rho_{s}(w')=\rho(j_{s}(w'))$, resp. $\rho^{\flat}_{s}(w' \sigma)=\rho^{\flat}(j_{s}(w')\sigma)$ pour tout $\sigma\in \Gamma_{{\mathbb F}_{q}}$. Pour $\rho'\in Irr(W_{s}(M_{s_{\Lambda}}))$, notons $(\bar{{\cal O}}_{\rho'},{\cal L}_{\rho'})$ l'image de $(M_{s_{\Lambda}},{\cal E}_{\varphi},\rho')$ par la repr\'esentation de Springer dans $G_{s}$.

  \begin{prop}{On a l'\'egalit\'e
  $$\hat{k}(\Lambda,\varphi,\rho)( h_{s,{\cal O}})=\sum_{\rho'\in  Irr_{{\mathbb F}_{q}}(W_{s}(M_{s_{\Lambda}}))}q^{dim(\bar{{\cal O}}_{\rho'})/2-dim({\cal O})/2}m(\rho_{s}^{\flat},\rho^{'\flat})\chi^{\natural}_{M_{s_{\Lambda}},{\cal E}_{\varphi},\rho'}(N).
  $$}\end{prop}
  
 Preuve. La fonction $h_{s,{\cal O}}$ est d\'eduite de la fonction $\underline{h}_{s,{\cal O}}$ sur $\mathfrak{g}_{s}({\mathbb F}_{q})$ et la transform\'ee de Fourier de celle-ci est \`a support nilpotent. On peut donc utiliser la proposition \ref{evaluationrho}, dans laquelle on remplace la fonction $\underline{h}_{s,{\cal O}}$ par sa valeur calcul\'ee par la proposition \ref{calcul}. On obtient
$$\hat{k}(\Lambda,\varphi,\rho)( h_{s,{\cal O}})=c_{1}\vert W_{s}(M_{s_{\Lambda}})\vert ^{-1}\sum_{w''\in W_{s}(M_{s_{\Lambda}})}trace(\rho^{\flat}(j_{s}(w'')Fr^{-1}))\sum_{(L',{\cal E}')\in {\cal I}^{FC}_{0,{\mathbb F}_{q}}(G_{s})} $$
$$\vert W_{s}(L')\vert ^{-1}q^{-dim(Z(L')^0)} \sum_{w'\in W_{s}(L')}\vert Z(L'_{w'})^0({\mathbb F}_{q})\vert Q^{\natural}_{L',{\cal E}',w'}(N)(Q^-_{M_{s_{\Lambda}},{\cal E}_{\varphi},w''},Q^-_{L',{\cal E}',w'}),$$
o\`u
$$c_{1}=c(\Lambda,\varphi)mes(K_{s}^0) \vert G_{s}({\mathbb F}_{q})\vert ^{-1}q^{dim(\mathfrak{g}_{s})-dim({\cal O})/2}.$$
On calcule le produit $(Q^-_{M_{s_{\Lambda}},{\cal E}_{\varphi},w''},Q^-_{L',{\cal E}',w'})=(Q_{M_{s_{\Lambda}},{\cal E}_{\varphi},w''},Q_{L',{\cal E}',w'})$ gr\^ace \`a \ref{fonctionscaracteristiques}(4). Il n'est non nul que si $(L',{\cal E}')=(M_{s_{\Lambda}},{\cal E}_{\varphi})$ et $w'$ et $w''$ sont $Fr$-conjugu\'es, c'est-\`a-dire $w''\in Cl_{s,Fr}(w')$. Remarquons que la premi\`ere condition fait dispara\^{\i}tre la somme en $(L',{\cal E}')$ de l'expression ci-dessus. Si les conditions pr\'ec\'edentes sont satisfaites, le produit ci-dessus vaut 
$$\vert W_{{\mathbb F}_{q}}(M_{s_{\Lambda},w'})\vert \vert Z(M_{s_{\Lambda},w'})^0({\mathbb F}_{q})\vert ^{-1}q^{dim({\bar{\cal O}}_{\varphi})-dim(M_{SC,s_{\Lambda}})}.$$
On obtient
$$\hat{k}(\Lambda,\varphi,\rho)( h_{s,{\cal O}})=c_{2}\vert W_{s}(M_{s_{\Lambda}})\vert ^{-2}\sum_{w'\in W_{s}(M_{s_{\Lambda}})}\vert Cl_{s,Fr}(w')\vert \vert W_{{\mathbb F}_{q}}(M_{s_{\Lambda},w'})\vert $$
$$ trace(\rho^{\flat}(j_{s}(w')Fr^{-1}))Q^{\natural}_{M_{s_{\Lambda}},{\cal E}_{\varphi},w'}(N),$$
o\`u
$$c_{2}=c_{1}q^{-dim(Z(M_{s_{\Lambda}})^0)+dim(\bar{{\cal O}}_{\varphi})-dim(M_{SC,s_{\Lambda}})}.$$
  Comme on l'a vu en \ref{evaluationw}(8), on a l'\'egalit\'e
$$\vert Cl_{s,Fr}(w')\vert \vert W_{{\mathbb F}_{q}}(M_{s_{\Lambda},w'})\vert =\vert W_{s}(M_{s_{\Lambda}})\vert$$
pour tout $w'\in W_{s}(M_{s_{\Lambda}})$, ce qui simplifie l'expression ci-dessus en
$$\hat{k}(\Lambda,\varphi,\rho)( h_{s,{\cal O}})=c_{2}\vert W_{s}(M_{s_{\Lambda}})\vert ^{-1}\sum_{w'\in W_{s}(M_{s_{\Lambda}})}  trace(\rho^{\flat}(j_{s}(w')Fr^{-1}))Q^{\natural}_{M_{s_{\Lambda}},{\cal E}_{\varphi},w'}(N).$$
On utilise (2) ci-dessus et on obtient
$$\hat{k}(\Lambda,\varphi,\rho)( h_{s,{\cal O}})=c_{2}\sum_{\rho'\in Irr_{{\mathbb F}_{q}}(W_{s}(M_{s_{\Lambda}}))}m(\rho_{s}^{\flat},\rho^{'\flat}) $$
$$\vert W_{s}(M_{s_{\Lambda}})\vert ^{-1}\sum_{w'\in W_{s}(M_{s_{\Lambda}})}  trace(\rho^{'\flat}(w'Fr^{-1}))Q^{\natural}_{M_{s_{\Lambda}},{\cal E}_{\varphi},w'}(N).$$
D'apr\`es \ref{Springer}(3), la derni\`ere ligne ci-dessus vaut $q^{-b(M_{s_{\Lambda}},{\cal E}_{\varphi},\rho')}\chi^{\natural}_{M_{s_{\Lambda}},{\cal E}_{\varphi},\rho'}(N)$. D'o\`u
$$\hat{k}(\Lambda,\varphi,\rho)( h_{s,{\cal O}})=\sum_{\rho'\in Irr_{{\mathbb F}_{q}}(W_{s}(M_{s_{\Lambda}}))}c_{3}(\rho')
 \chi^{\natural}_{M_{s_{\Lambda}},{\cal E}_{\varphi},\rho'}(N),$$
 o\`u
 $$c_{3}(\rho')=c_{2}q^{-b(M_{s_{\Lambda}},{\cal E}_{\varphi},\rho')}m(\rho_{s}^{\flat},\rho^{'\flat}).$$
 On a l'\'egalit\'e $mes(K_{s}^0)=q^{-dim(\mathfrak{g}_{s})/2}\vert G_{s}({\mathbb F}_{q})\vert $. En se rappelant la d\'efinition
 $$c(\Lambda,\varphi)=q^{dim(M_{s_{\Lambda}})/2-dim(\bar{{\cal O}}_{\varphi})/2},$$
 on calcule la constante 
 $$c_{3}(\rho')=q^{dim({\bar{\cal O}}_{\rho'})/2-dim({\cal O})/2}m(\rho_{s}^{\flat},\rho^{'\flat}).$$
La formule ci-dessus devient celle de l'\'enonc\'e. $\square$

\subsection{ Deux th\'eor\`emes de maximalit\'e}

\subsubsection{Param\'etrage de l'ensemble ${\cal B}_{Irr}(G)$\label{parametragecalB}}
Considérons l'ensemble des triplets $(N,d,\mu)$ où 

$N\in \boldsymbol{\mathfrak{g}}_{nil}/conj$ et l'orbite de $N$ est conservée par $I_{F}$;

$d\in \tilde{\bar{A}}(N)$;

$\mu$ est une représentation irréductible de $Z_{\bar{A}(N)}(d)$. 

Le groupe ${\bf G}$ agit sur cet ensemble, on note $\bar{{\cal C}}^{Irr}$ l'ensemble des classes de conjugaison. Pour une orbite $\boldsymbol{{\cal O}}\in (\boldsymbol{\mathfrak{g}}_{nil}/conj)^{I_{F}}$, on note $\bar{{\cal C}}^{Irr}(\boldsymbol{{\cal O}})$ le sous-ensemble des classes de triplets $(N,d,\mu)$ tels que $N\in \boldsymbol{{\cal O}}$. On note $\bar{{\cal C}}^{Irr}_{sp}$, resp. $\bar{{\cal C}}^{Irr}_{sp,F}$, la réunion des $\bar{{\cal C}}^{Irr}(\boldsymbol{{\cal O}})$  pour $\boldsymbol{{\cal O}}\in {\cal U}_{sp}$, resp. et, de plus, $\boldsymbol{{\cal O}}$ est conservée par $\Gamma_{F}^{nr}\simeq \Gamma_{{\mathbb F}_{q}}$.

Nous allons d\'efinir une application injective
$$(1) \qquad \boldsymbol{\nabla}_{F}:{\cal B}_{Irr}(G)\to \bar{{\cal C}}^{Irr}_{sp,F}.$$

Dans ce paragraphe, nous donnons la d\'efinition dans le cas o\`u $G$ est d\'eploy\'e sur $F^{nr}$. Quand $G$ n'est pas déployé sur $F^{nr}$, nous donnerons la définition dans les paragraphes \ref{An-1ramSI} à \ref{E6ramSI}. 

Supposons d'abord que $G$ est d\'eploy\'e sur $F$ et que $q-1$ soit divisible par $3\times 4\times 5$. On a dit en \ref{faisceauxcaracteresunipotentscuspidaux} que l'ensemble ${\bf FC}_{u}({\bf G})$ \'etait en bijection avec celui des couples $(s,{\cal E})$ o\`u $s\in S^{nr,st}(\bar{C}^{nr})$ et ${\cal E}\in {\bf FC}^{st}(\mathfrak{g}_{s})$. Sous les hypoth\`eses que l'on vient de poser, on a $S^{nr,st}(\bar{C}^{nr})=S^{st}(\bar{C})$ et ${\bf FC}^{st}(\mathfrak{g}_{s})={\bf FC}_{{\mathbb F}_{q}}^{st}(\mathfrak{g}_{s})$ pour tout $s\in S^{st}(\bar{C})$ (cette derni\`ere \'egalit\'e r\'esulte de  \cite{W7} paragraphe 9).  En fixant une action de Frobenius sur tout \'el\'ement de ${\bf FC}_{{\mathbb F}_{q}}^{st}(\mathfrak{g}_{s})$, cet ensemble s'identifie \`a ${\cal B}(G,s)$. On obtient que ${\bf FC}_{u}({\bf G})$ s'identifie au sous-ensemble des quadruplets $(M,s_{M},\varphi,\rho)\in {\cal B}_{Irr}(G)$ tels que $M=G$ (et donc $\rho=1$). Il y a une bijection naturelle entre les classes de conjugaison de $F$-Levi de $G$ et les classes de conjugaison de ${\mathbb F}_{q}$-Levi de ${\bf G}$.  Soient $(M,{\bf M})$ un couple de Levi se correspondant. Le r\'esultat pr\'ec\'edent vaut en rempla\c{c}ant $(G,{\bf G})$ par $(M,{\bf M})$. On a $W^{I_{F}}(M)\simeq W({\bf M})$ et $Irr_{{\mathbb F}_{q}}(W^{I_{F}}(M))\simeq Irr(W({\bf M}))$ puisque $G$ est d\'eploy\'e.     On d\'efinit alors une bijection de  ${\cal B}_{Irr}(G)$ sur ${\cal A}({\bf G})$: \`a un quadruplet $(M,s_{M},\varphi,\rho)\in {\cal B}_{Irr}(G)$, on associe le triplet 
$({\bf M},\boldsymbol{{\cal E}},\rho)\in {\cal A}({\bf G})$ o\`u ${\bf M}$ correspond \`a $M$ comme ci-dessus et $\boldsymbol{{\cal E}}\in {\bf FC}_{u}({\bf M})$ correspond \`a $(s_{M},\varphi)$.

Affaiblissons les hypoth\`eses sur $G$ en supposant seulement que $G$ est d\'eploy\'e sur $F^{nr}$. On peut fixer une extension finie non ramifi\'ee $F'$ de $F$ telle que les hypoth\`eses pr\'ec\'edentes soient v\'erifi\'ees sur le corps de base $F'$. L'ensemble ${\cal A}({\bf G})$ est insensible \`a un tel changement de corps de base. L'ensemble ${\cal B}_{Irr}(G)$ l'est. Mais notons plus pr\'ecis\'ement ${\cal B}_{Irr,F}(G)$ et ${\cal B}_{Irr,F'}(G)$ ces ensembles relatifs aux deux corps de base, on a une injection naturelle ${\cal B}_{Irr,F}(G)\to {\cal B}_{Irr,F'}(G)$. 
En effet, ${\cal B}_{Irr,F}(G)$ appara\^{\i}t comme le sous-ensemble des \'el\'ements de ${\cal B}_{Irr,F'}(G)$ dont tous les termes v\'erifient une certaine propri\'et\'e de conservation par l'action de $\Gamma_{F'/F}$ (dans le cas des faisceaux-caract\`eres cuspidaux, cette conservation se traduit par une condition de divisibilit\'e de $q-1$, cf. \cite{W7} paragraphe 9). La compos\'ee de cette injection et de la bijection  $ {\cal B}_{Irr,F'}(G)\to {\cal A}({\bf G})$ fournit une injection naturelle ${\cal B}_{Irr,F}(G)\to {\cal A}({\bf G})$. 

Puisque $G$ est déployé sur $F^{nr}$, il s'identifie au groupe $G_{F}$ introduit en \ref{legroupeG}. Alors $\bar{{\cal C}}^{Irr}_{sp}$ est le m\^eme ensemble qu'en \ref{parametrage3}. On a défini dans ce paragraphe une bijection $\nabla:{\cal A}({\bf G})\to \bar{{\cal C}}_{sp}^{Irr}$. On définit$\boldsymbol{\nabla}_{F}$ comme le composé de l'injection ${\cal B}_{Irr}(G)\to {\cal A}({\bf G})$ et de $\nabla$. 

  Il nous reste \`a prouver que

(1) $\boldsymbol{\nabla}_{F}$ prend ses valeurs dans $ \bar{{\cal C}}^{Irr}_{sp,F}$. 

 Autrement dit, soit $(N,d,\mu)$ un élément de  l'image de $\boldsymbol{\nabla}_{F}$. On doit prouver que l'orbite $\boldsymbol{{\cal O}}$ de $N$ est conserv\'ee par l'action de $\Gamma_{{\mathbb F}_{q}}$. 
D'apr\`es \ref{actiongaloisienne}, c'est toujours le cas si $G$ n'est pas de type $D_{n}$ ou si $G$ est d\'eploy\'e. Supposons que $G$ soit de type $D_{n}$ avec $n\geq4$ et que le Frobenius $Fr$ agisse sur ${\cal D}$ par l'automorphisme habituel   d'ordre $2$. Supposons que $\boldsymbol{{\cal O}}$ ne soit pas conservée par $\Gamma_{{\mathbb F}_{q}}$. Alors tous les termes de la partition associée à $\boldsymbol{{\cal O}}$  sont pairs. On sait que ${\cal A}({\bf G},\boldsymbol{{\cal O}})$ est r\'eduit \`a un seul \'el\'ement qui est de la forme $({\bf T}, \boldsymbol{{\cal E}},\rho)$ o\`u $\boldsymbol{{\cal E}}$ est le syst\`eme local trivial sur ${\bf T}_{AD}=\{1\}$ et $\rho$ est une repr\'esentation irr\'eductible de $W$ qui n'est pas conserv\'ee par l'action galoisienne. Cette derni\`ere propri\'et\'e interdit \`a un tel triplet d'appartenir \`a l'image de l'injection ${\cal B}_{Irr,F}(G)\to {\cal A}({\bf G})$. D'o\`u  (1). Il reste le cas o\`u $G$ est de type $D_{4}$ et o\`u le Frobenius agit sur ${\cal D}$ par un automorphisme d'ordre $3$. Notons $\lambda$ la partition associée à $\boldsymbol{{\cal O}}$.  Puisque $\boldsymbol{{\cal O}}$ est spéciale, on voit à l'aide de  \ref{actiongaloisienne} que $\lambda$ est l'une des partitions $(51^3)$, $(4^2)$, $(31^5)$, $(2^4)$. Les partitions   $(4^2)$ et $(2^4)$ s'\'eliminent par le m\^eme argument que ci-dessus. Supposons $\lambda=(51^3)$, resp. $(31^5)$. Le calcul de ${\cal A}({\bf G},\boldsymbol{{\cal O}})$ résulte de  \cite{C} p. 409, on trouve que cet ensemble n'a qu'un \'el\'ement qui est de la forme  $({\bf T}, \boldsymbol{{\cal E}},\rho)$, o\`u $\boldsymbol{{\cal E}}$ est comme ci-dessus et $\rho$ est une repr\'esentation irr\'eductible de $W$.  D\'ecrivons cette repr\'esentation. Les repr\'esentations irr\'eductibles de $W$ sont param\'etr\'ees par des couples $(\nu,\nu')\in {\cal P}_{2}(4)$, les couples $(\nu,\nu')$ et $(\nu',\nu)$ \'etant identifi\'es (en fait, si $\nu=\nu'$, le couple param\`etre deux repr\'esentations irr\'eductibles mais peu nous importe). Ici le calcul donne que $\rho$ est param\'etr\'ee par $((31),\emptyset)$, resp. $((21^2),\emptyset)$.  Ces repr\'esentations ne sont pas invariantes par l'action galoisienne: en notant $(w_{i})_{i=1,...,4}$ les g\'en\'erateurs de $W$ associ\'es aux racines simples, on calcule $trace(\rho(w_{3}w_{4}))=3$ et $trace(\rho(w_{1}w_{3}))=-1$, alors que $w_{3}w_{4}$ et $w_{1}w_{3}$ sont conjugu\'es par l'action galoisienne. Comme ci-dessus, cela interdit au triplet  $({\bf T}, \boldsymbol{{\cal E}},\rho)$ d'appartenir \`a l'image de l'injection ${\cal B}_{Irr,F}(G)\to {\cal A}({\bf G})$. D'o\`u (1).

Cela ach\`eve la d\'efinition de $\boldsymbol{\nabla}_{F}$ quand $G$ est d\'eploy\'e sur $F^{nr}$. 
Signalons une propri\'et\'e de l'application ${\cal B}_{irr}(G)\to {\cal A}({\bf G})$ que  l'on a d\'efinie ci-dessus. Soit $(M,s_{M},\varphi,\rho)\in {\cal B}_{Irr}(G)$, notons 
$({\bf M},\boldsymbol{{\cal E}},\rho)\in {\cal A}({\bf G})$ son image. En \ref{couples}, on a associ\'e \`a $(M,s_{M})$ un sous-ensemble $\Lambda\subset \Delta_{a}$. En \ref{Lambda}, on a associ\'e \`a $({\bf M},\boldsymbol{{\cal E}})$ un sous-ensemble analogue. Montrons que

(2) ces deux ensembles $\Lambda$ sont \'egaux.

On peut supposer que $M$ est standard. On choisit pour alc\^ove $C^{nr,M}$ de $Imm_{F^{nr}}(M_{AD})$ celle qui contient $p_{M}(C^{nr})$. Il s'en d\'eduit un ensemble de racines affines $\Delta_{a}^M$ et on peut supposer que $s_{M}$ est associ\'e \`a une racine $\alpha\in \Delta_{a}^M$. Ainsi, l'ensemble   $\Sigma^M$ des racines de $T$ dans $M$ et l'ensemble $\Sigma^{M_{s_{M}}}$ des racines de $T_{s_{M}}$ dans $M_{s_{M}}$ sont des sous-ensembles de $\Sigma$. Par d\'efinition de l'application de ${\cal B}_{irr}(G)$ dans $ {\cal A}({\bf G})$, on peut supposer que ${\bf M}$ est le Levi standard de ${\bf G}$ dont l'ensemble de racines est $\Sigma^M$ et que la racine associ\'ee \`a $\boldsymbol{{\cal E}}$ est $\alpha$. L'ensemble de racines de ${\bf T}$ dans ${\bf M}_{t_{\alpha}}$ est alors $\Sigma^{M_{s_{M}}}$. Alors les deux ensembles $\Lambda$ sont caract\'eris\'es par la m\^eme propri\'et\'e. A savoir qu'il existe un \'el\'ement $w\in W$ qui conjugue $\Sigma^M$, resp. $\Sigma^{M_{s_{M}}}$, en l'ensemble des \'el\'ements de $\Sigma$ qui sont combinaisons lin\'eaires \`a coefficients rationnels, resp. entiers, d'\'el\'ements de $\Lambda$. Cette caract\'erisation commune entra\^{\i}ne leur \'egalit\'e. D'o\`u (2).

   \subsubsection{Le cas  de type $(A_{n-1},ram)$\label{An-1ramSI}}
   
   Dans ce paragraphe et les suivants, on donne la d\'efinition de $\boldsymbol{\nabla}_{F}$ dans le cas o\`u $G$ n'est pas d\'eploy\'e sur $F^{nr}$. La d\'efinition copie la classification des représentations unipotentes des groupes finis non déployés. Elle est faite pour que les th\'eor\`emes \ref{premiertheoremeSI} et \ref{deuxiemetheoremeSI} soient v\'erifi\'es. 
   
     On suppose ici que $G$ est de type $(A_{n-1},ram)$.  Puisque $A(N)=\{1\}$ pour tout \'el\'ement nilpotent $N\in \boldsymbol{\mathfrak{g}}_{nil}$ et que toute orbite nilpotente est sp\'eciale et conserv\'ee par $\Gamma_{F}$, l'ensemble $\bar{{\cal C}}^{Irr}_{sp,F}$ s'identifie simplement \`a $\boldsymbol{\mathfrak{g}}_{nil}/conj
$ ou encore \`a ${\cal P}(n)$.

Pour $(M,s_{M},\varphi)\in {\cal B}(G)$, $M_{SC}$ est de type $A_{u-1}$ o\`u $u$ est un nombre triangulaire,  c'est-\`a-dire $u=v(v+1)/2$ pour un $v\in {\mathbb N}$ et $u$ est de m\^eme parit\'e que $n$. Le groupe $M$, ou encore l'entier $u$, d\'etermine un unique couple $(s_{M},\varphi)$.  Le groupe de Weyl $W(M)$ est celui d'un syst\`eme de racines de type $B_{(n-u)/2}$ (ou $C_{(n-u)/2}$). Ses repr\'esentations sont param\'etr\'ees par l'ensemble ${\cal P}_{2}((n-u)/2)$. On obtient une bijection 
  $(M,s_{M},\varphi,\rho)\mapsto (v,\alpha,\beta)$   de ${\cal B}_{Irr}(G)$ sur l'ensemble ${\cal M}(n)$ des triplets $(v,\alpha,\beta)$ tels que

$v\in {\mathbb N}$ et $v(v+1)/2$ est de m\^eme parit\'e que $n$;

$\alpha$, $\beta$ sont des partitions;

$v(v+1)/2+2S(\alpha)+2S(\beta)=n$. 

L'application $\boldsymbol{\nabla}_{F}$ s'identifie \`a une application de ${\cal M}(n)$ sur ${\cal P}(n)$. Sa d\'efinition est issue de \cite{AMR} paragraphe 5.B. 
Soit $(v,\alpha,\beta)\in {\cal M}(n)$. Fixons un entier $t\geq n$ de m\^eme parit\'e que $v+1$. On a $l(\alpha),l(\beta)\leq (t-v-1)/2$. Posons $X(\alpha)
 =\{2\alpha_{j}+t+v+1-2j; j=1,...,(t+v+1)/2\}$, $Y(\beta)=\{2\beta_{j}+t-v-2j; j=1,...,(t-v-1)/2\}$.  Pour une partition $\lambda\in {\cal P}(n)$, posons ${\bf Z}(\lambda)=\{\lambda_{j}+t-j; j=1,...,t\}$. Alors $\boldsymbol{\nabla}_{F}(v,\alpha,\beta)$ est l'unique partition $\lambda\in {\cal P}(n)$ telle que ${\bf Z}(\lambda)=X(\alpha)\cup Y(\beta)$. On v\'erifie que $\boldsymbol{\nabla}_{F}$ est injective.

\subsubsection{Le cas  de type $(D_{n},ram)$\label{DnramSI}}
On suppose que $G$ est de type $(D_{n},ram)$.   Pour $\lambda\in {\cal P}^{orth}(2n)$, on dit que $\lambda$ sp\'eciale  si et seulement si $\lambda_{2j-1}\equiv \lambda_{2j}\,mod\,2{\mathbb Z}$ pour tout entier $j\geq1$.  Consid\'erons une orbite  dans $\boldsymbol{\mathfrak{g}}_{nil}$ param\'etr\'ee par une partition $\lambda\in {\cal P}^{orth}(2n)$. L'orbite est sp\'eciale si et seulement si $\lambda$ l'est. Supposons qu'il en soit ainsi. On associe \`a $\lambda$ un symbole $({\bf X}(\lambda),{\bf Y}(\lambda))$ de la fa\c{c}on suivante.  On introduit la partition $\lambda+(2n-1,2n-2,...,0)$ et on consid\`ere qu'elle a $2n$ termes (\'eventuellement le dernier peut \^etre nul). On constate que $n$ termes sont pairs, on les note $2x_{1}>...>2x_{n}$, et que $n$ termes sont impairs, on les note $2y_{1}+1>...>2y_{n}+1$. Alors ${\bf X}(\lambda)=(x_{1},...,x_{n})$ et ${\bf Y}(\lambda)=(y_{1},...,y_{n})$. On a d\'efini en \ref{parametrageorbitesnilpotentes} l'ensemble des intervalles $Int(\lambda)\subset Jord_{bp}(\lambda)$. Pour $\Delta\in Int(\lambda)$, notons $j_{min}(\Delta)$, resp. $j_{max}(\Delta)$, le plus petit, resp. grand, indice $j\geq1$ tel que $\lambda_{j}\in \Delta$.  L'indice $j_{min}(\Delta)$ est impair tandis que $j_{max}(\Delta)$ est pair. On v\'erifie que
$${\bf X}(\lambda)-({\bf X}(\lambda)\cap {\bf Y}(\lambda))=\{x_{(j_{min}(\Delta)+1)/2}; \Delta\in Int(\lambda)\},$$
$${\bf Y}(\lambda)-({\bf X}(\lambda)\cap {\bf Y}(\lambda))=\{y_{j_{max}(\Delta)/2}; \Delta\in Int(\lambda)\}.$$
 Consid\'erons l'ensemble des  paires $(X,Y)$ v\'erifiant les conditions suivantes:
 
   $X,Y\subset {\mathbb N}$, $X\cap Y={\bf X}(\lambda)\cap {\bf Y}(\lambda)$, $X\cup Y={\bf X}(\lambda)\cup {\bf Y}(\lambda)$.
   
   \noindent Dans la derni\`ere \'egalit\'e, on consid\`ere les deux membres comme des partitions, ou encore des ensembles "avec multiplicit\'es". Disons que 
  deux paires $(X,Y)$ sont \'equivalentes si et seulement si $(X',Y')=(X,Y)$ ou $(X',Y')=(Y,X)$. L'ensemble des classes d'\'equivalence est appel\'ee la famille de $\lambda$ et est not\'ee 
 $Fam(\lambda)$.  Un \'el\'ement de cette famille sera not\'ee comme   une paire $(X,Y)$ repr\'esentant sa classe d'\'equivalence. On note $Fam(\lambda)_{pair}$, resp. $Fam(\lambda)_{imp}$, l'ensemble des $(X,Y)\in Fam(\lambda)$ tels que $\vert X\vert \equiv\vert Y\vert \equiv n\,mod\,2{\mathbb Z}$, resp. $\vert X\vert \equiv\vert Y\vert \equiv n+1\,mod\,2{\mathbb Z}$. Notons $\bar{A}(\lambda)^{\vee}$ le groupe dual de $\bar{A}(\lambda)$. On va d\'efinir une bijection
 $$(1) \qquad \phi_{\lambda}:Fam(\lambda)_{imp}\to \tilde{\bar{A}}(\lambda)\times \bar{A}(\lambda)^{\vee}.$$
 
 Notons ${\bf -1}$ l'\'el\'ement de $\{\pm 1\}^{Int(\lambda)}$ dont toutes les composantes valent $-1$. 
   On d\'efinit d'abord une bijection 
 $$(2) \qquad Fam(\lambda)\to (\{\pm 1\}^{Int(\lambda)}\times \{\pm 1\}^{Int(\lambda)})/\{1,( {\bf -1},{\bf -1})\}.$$
  Pour $(X,Y) \in Fam(\lambda)$, il existe deux uniques sous-ensembles $Int^X,Int^Y\subset Int(\lambda)$ tels que $X=({\bf X}(\lambda)-\{x_{(j_{min}(\Delta)+1)/2};\Delta\in Int^X\})\cup \{y_{j_{max}(\Delta)/2};\Delta\in Int^Y\}$, $Y=({\bf Y}(\lambda)-\{y_{j_{max}(\Delta)/2};\Delta\in Int^Y\})\cup \{x_{(j_{min}(\Delta)+1)/2};\Delta\in Int^X\}$. On associe \`a $Int^X$ l'\'el\'ement $\nu^X=(\nu^X_{\Delta})_{\Delta\in Int(\lambda)}\in \{\pm 1\}^{Int(\lambda)}$ tel que 
 $$\nu^X_{\Delta}=\left\lbrace\begin{array}{cc}1,& \text{ si } \Delta\not\in Int^X,\\ -1,& \text{ si }\Delta\in Int^X.\\ \end{array}\right.$$
 On d\'efinit de m\^eme $\nu^Y$. Si l'on remplace $(X,Y)$ par $(Y,X)$, le couple $(\nu^X,\nu^Y)$ est remplac\'e par $(\nu^X ({\bf -1}),\nu^Y{(\bf -1}))$. On voit alors que l'application $(X,Y)\mapsto (\nu^X,\nu^Y)$ se quotiente en la bijection (2).    Rappelons que l'ensemble $Int(\lambda)$ est naturellement ordonn\'e. 
 D\'efinissons une bijection
 $$(3) \begin{array}{ccc} (\{\pm 1\}^{Int(\lambda)}\times \{\pm 1\}^{Int(\lambda)})/\{1,( {\bf -1},{\bf -1})\}&\to& \{\pm 1\}^{Int(\lambda)}\times(\{\pm 1\}^{Int(\lambda)}/\{1,{\bf -1}\})\\ (\nu^X,\nu^Y)&\mapsto&(d,\mu)\\ \end{array}$$
 par les formules  $d_{\Delta}=\nu^X_{\Delta}\nu^Y_{\Delta}$, $\mu_{\Delta}=(\prod_{\Delta'\geq \Delta}\nu^X_{\Delta'})\prod_{\Delta'> \Delta}\nu^Y_{\Delta'}$ pour tout $\Delta\in Int(\lambda)$ (comme plus haut, on voit que l'application d\'efinie par ces formules se quotiente en la bijection (3)). La compos\'ee des bijections (2) et (3) envoie $Fam(\lambda)_{imp}$ sur  $\{\pm 1\}^{Int(\lambda)}_{-1}\times(\{\pm 1\}^{Int(\lambda)}/\{1,{\bf -1}\})$.
 Mais  $\{\pm 1\}^{Int(\lambda)}_{-1}=\tilde{\bar{A}}(\lambda)$ et  $\{\pm 1\}^{Int(\lambda)}/\{1,{\bf -1}\} $ s'identifie naturellement au dual de $\{\pm 1\}^{Int(\lambda)}_{1}$, autrement dit au dual $\bar{A}(\lambda)^{\vee}$ de $\bar{A}(\lambda)$. On a ainsi \'etabli la bijection $\phi_{\lambda}$ de la relation (1).

Pour $(M,s_{M},\varphi)\in {\cal B}(G)$, ou bien $M=T$, ou bien 
$M_{SC}$ est  
de type $ D_{h^2}$ o\`u $h$ est un entier impair tel que $3\leq h$ et $h^2\leq n$. . Par convention, nous posons $h=1$ dans le cas o\`u $M=T$. Le groupe $M$, ou encore l'entier $h$, d\'etermine un unique couple $(s_{M},\varphi)$.  Le groupe de Weyl $W^{I_{F}}(M)$ est celui d'un syst\`eme de racines de type $B_{n-h^2}$ ou $C_{n-h^2}$. Ses repr\'esentations sont param\'etr\'ees par l'ensemble ${\cal P}_{2}(n-h^2)$ (si $n-h^2=1$, le groupe de Weyl a deux éléments, son caractère trivial, resp. non trivial, correspond à l'élement $(1,\emptyset)$, resp. $(\emptyset,1)$, de ${\cal P}_{2}(1)$). 
L'ensemble ${\cal B}_{Irr}(G)$ s'identifie \`a celui des triplets $(h,\alpha,\beta)$ o\`u $h$ est comme ci-dessus et $(\alpha,\beta)\in {\cal P}_{2}(n-h^2)$. Pour un tel triplet, on peut consid\'erer que $\alpha$ a $n+h$ termes et que $\beta$ en a $n-h$. Posons 

 $X(\alpha)=\alpha+(n+h-1,...,0)$, $Y(\beta)=\beta+(n-h-1,...,0)$. 

\noindent Il existe une unique partition  $\lambda\in {\cal P}^{orth}(2n)$ qui est sp\'eciale et  telle que $(X(\alpha),Y(\beta))\in Fam(\lambda)$. La non-nullit\'e de $h$ entra\^{\i}ne que $\lambda$ a au moins un terme impair, donc l'orbite $\boldsymbol{{\cal O}}$ associée à  $\lambda$  est conservée par $\Gamma_{{\mathbb F}_{q}}$. L'imparit\'e de $h$ entra\^{\i}ne que $(X(\alpha),Y(\beta))\in Fam_{imp}(\lambda)$.  Posons $(d,\mu)=\phi_{\lambda}(X(\alpha),Y(\beta))$. On d\'efinit alors $\boldsymbol{\nabla}_{F}(h,\alpha,\beta)=(N,d,\mu)$, où $N$ est un élément de $\boldsymbol{{\cal O}}$ et $d$, resp.  $\mu$, est identifié à un élément de $\tilde{\bar{A}}(N)$, resp. à une représentation de $\bar{A}(N)$ .  L'application $\boldsymbol{\nabla}_{F}$ ainsi d\'efinie est injective (en fait bijective).

  \subsubsection{Le cas $(D_{4},3-ram)$\label{D3ramSI}}
  On suppose que $G$ est de type $(D_{4},3-ram)$. 
Les couples $(M,s_{M})$ tel que ${\cal B}(M,s_{M})\not=\emptyset$ sont

$M=G$, $s_{M}=s_{0}$; $M=G$, $s_{M}=s_{134}$; $M=T$ et $s_{M}$ est l'unique point de $Imm(M_{AD})$. 

Dans chaque cas,  ${\cal B}(M,s_{M})$ est r\'eduit \`a un \'el\'ement. Dans le cas o\`u $M=G$, on a $W^{I_{F}}(G)=\{1\}$ et au couple $(G,s_{G})$
 est associ\'e un unique \'el\'ement de ${\cal B}_{Irr}(G)$ que l'on note simplement $(G,s_{G})$.   Si $M=T$, $W^{I_{F}}(M)$ s'identifie au groupe de Weyl de $ G_{s_{0}}$, qui est de type $G_{2}$. On note ses repr\'esentations comme dans \cite{C} page 412 et on note $(T,\rho)$ l'\'el\'ement de ${\cal B}_{Irr}(G)$ associ\'e \`a une repr\'esentation $\rho$ de ce groupe.  
 
 L'ensemble des orbites sp\'eciales de ${\bf G}$ conserv\'ees par $I_{F}$ est param\'etr\'e par l'ensemble de partitions   $\{(71),(53),(3^21^2) ,(2^21^4),(1^8)\}$, cf. \ref{actiongaloisienne}. On l'identifie \`a cet ensemble de partitions. Sauf dans le cas de la partition $\lambda=(3^21^2)$, on a $A(\lambda)=\{1\}$ et $\lambda$ se compl\`ete en un unique \'el\'ement de $\bar{{\cal C}}^{Irr}_{sp,F}$ que l'on note encore   $\lambda$. Si $\lambda=(3^21^2)$, on choisit un représentant $N$ de l'orbite associée à $\lambda$  qui appartient à l'ensemble ${\cal N}(\Gamma_{F})$ de \ref{calCF}. On  a $\tilde{\bar{A}}(N)=A(N)\gamma= \{\pm 1\}\gamma$. On note $1$, resp. $sgn$, le caract\`ere trivial, resp.  non trivial, de $ \{\pm 1\}$. Pour $d\in \tilde{\bar{A}}(N)$ et pour un caractère $\mu$ de $\bar{A}(N)$, on note $(3^21^2,d,\mu)$ le triplet $(N,d,\mu)$. Cette définition ne dépend pas du choix de $N$: d'après la preuve de \ref{preuvedecalCF}, on ne peut modifier $N$ que par conjugaison par un élément de ${\bf G}^{I_{F}}$ et une telle conjugaison fixe $\gamma$. 

  On d\'efinit $\boldsymbol{\nabla}_{F}$ par le tableau suivant:
  
  $$\begin{array}{cc}(M,s_{M},\varphi,\rho)&\boldsymbol{\nabla}_{F}(M,s_{M},\varphi,\rho)\\ &\\ (G,s_{0})& (3^21^2,\gamma,sgn)\\ 
  (G,s_{134})&(3^21^2,-\gamma,sgn) \\ (T,\phi_{1,0})&(71)\\ (T,\phi'_{1,3})&(53) \\ (T,\phi_{21})&(3^21^2,\gamma,1)\\ (T,\phi_{2,2})&(3^21^2,-\gamma,1)\\ (T,\phi_{1,3}'')&(2^21^4)\\ (T,\phi_{1,6})&1^8\\ \end{array}.$$ 
  
  \subsubsection{Le cas $(E_{6},ram)$\label{E6ramSI}}
  
  On suppose que $G$ est de type $(E_{6},ram)$. Notons $M(A_{5})$ l'unique $F$-Levi standard $M$ tel que $M_{SC}$ soit de type $A_{5}$.   
  Les couples $(M,s_{M})$ tel que ${\cal B}(M,s_{M})\not=\emptyset$ sont
  
  $M=G$, $s_{M}=s_{0}$; 
  
  si $\delta_{3}(q-1)=1$, $M=G$  et $s_{M}=s_{35}$; 
  
  $M=M(A_{5})$, $s_{M}$ est le sommet tel que $M_{s_{M},SC}$ soit de type $A_{1}\times A_{1}\times A_{1}$;  
    
  $M=T$ et $s_{M}$ est l'unique point de $Imm(M_{AD})$. 

Sauf dans le cas o\`u $M=G$ et $s_{M}=s_{35}$, ${\cal B}(M,s_{M})$ a un seul \'el\'ement et on l'omet de la notation. Dans le cas o\`u $M=G$ et $s_{M}=s_{35}$, ${\cal B}(M,s_{M})$ a deux \'el\'ements que l'on note $\varphi_{1}$ et $\varphi_{2}$. Si $M=G$, $W^{I_{F}}(M)=\{1\}$ et on omet son unique repr\'esentation de la notation. Si $M=M(A_{5})$,  $W^{I_{F}}(M)$ a deux \'el\'ements, on note $1$ et $sgn$ ses deux caract\`eres. Si $M=T$, $W^{I_{F}}(M)$ s'identifie au groupe de Weyl de $G_{s_{0}}$, qui est de type $F_{4}$. On note ses repr\'esentations comme dans \cite{C} p. 412 et on note $(T,\rho)$ l'\'el\'ement de ${\cal B}_{Irr}(G)$ associ\'e \`a une repr\'esentation $\rho$ de ce groupe. 

 Toute orbite nilpotente de $\boldsymbol{\mathfrak{g}}$ est conservée par $\Gamma_{F}$, cf. \ref{iotanilram}. Soit $N$ un élément de l'ensemble ${\cal N}(\Gamma_{F})$ de \ref{calCF} et dont l'orbite $\boldsymbol{{\cal O}}$ est spéciale.  
  On a en g\'en\'eral $A(N)=\{1\}$. Alors $N$ se compl\`ete en un unique \'el\'ement de $ \bar{{\cal C}}^{Irr}_{sp,F}$ que l'on note simplement  comme  $\boldsymbol{{\cal O}}$ est d\'esign\'ee dans les tables de \cite{C}. Si $A(N)\not=\{1\}$, $N$ est fixe par $\gamma$ d'apr\`es le choix de notre ensemble ${\cal N}$ et on a $\tilde\bar{{A}}(N)=\bar{A}(N)\gamma$, o\`u $\gamma$ agit trivialement sur $\bar{A}(N)$. 
Si $N$ est de type $A_{2}$ ou $ E_{6}(a_{3})$, on a $\tilde{\bar{A}}(N)=\{\pm 1\}\gamma$ et on note $1$ et $sgn$ les deux caract\`eres de $ \{\pm 1\}$. Si $N$ est de type  $ D_{4}(a_{1})$, on a $\tilde{\bar{A}}(N)\simeq \mathfrak{S}_{3}\gamma$. On note $1,g_{2},g_{3}$ des repr\'esentants des trois classes de conjugaison dans $\mathfrak{S}_{3}$, $g_{i}$ \'etant d'ordre $i$ pour $i=2,3$.   Le groupe $Z_{\bar{A}(N)}(\gamma)$ est $\mathfrak{S}_{3}$ dont on note les repr\'esentations irr\'eductibles $1$, $sgn$, $r$ ($r$ est de dimension $2$). Le groupe $Z_{\bar{A}(N)}(g_{2}\gamma)$ est  d'ordre $2$ et  on note ses caract\`eres $1$ et $sgn$. Le groupe $Z_{\bar{A}(N)}(g_{3}\gamma)$ est  d'ordre $3$ et on note ses caract\`eres $1$, $\theta$, $\theta^2$. Pour $d\in \tilde{\bar{A}}(N)$ et pour une représentation irréductible  $\mu$ de $Z_{\bar{A}(N)}(d)$, on note $(\star,d,\mu)$ le triplet $(N,d,\mu)$, où $\star$ est la désignation de $\boldsymbol{{\cal O}}$ dans les tables de \cite{C} p.412. Comme en \ref{D3ramSI}, cette définition ne dépend pas du choix de $N$, modulo la remarque qui suit.  

{\bf Remarque.} Il n'est pas clair pour nous de savoir qui est $\theta$ et qui est $\theta^2$, ni de savoir  qui est $\varphi_{1}$ et qui est $\varphi_{2}$. Dans les formules qui suivent, il y a donc une ambigu\"{\i}t\'e, qui ne sera pas importante pour la suite.
\bigskip

On d\'efinit $\boldsymbol{\nabla}_{F}$ par le tableau suivant:

$$\begin{array}{ccc}\text{conditions}&(M,s_{M},\varphi,\rho)&\boldsymbol{\nabla}_{F}(M,s_{M},\varphi,\rho)\\ &&\\ 
&(G,s_{0})&(D_{4}(a_{1}),\gamma,sgn)\\ \delta_{3}(q-1)=1 &(G,s_{35},\varphi_{1}),(G,s_{35},\varphi_{2}) &(D_{4}(a_{1}),g_{3}\gamma,\theta) ,(D_{4}(a_{1}),g_{3}\gamma,\theta^2)\\ &
 (M(A_{5}),s_{M(A_{5})},1)&D_{5}(a_{1})\\ &(M(A_{5}),s_{M(A_{5})},sgn)&A_{2}+A_{1}\\ &(T,\phi_{1,0})&E_{6}\\ &(T,\phi'_{2,4})&E_{6}(a_{1})\\ &(T,\phi_{4,1})&D_{5}\\ &(T,\phi_{9,2})&(E_{6}(a_{3}),\gamma,1)\\ &
 (T,\phi''_{2,4})&(E_{6}(a_{3}),\gamma,sgn)\\ &(T,\phi'_{8,3})&(E_{6}(a_{3}),-\gamma,1)\\ &(T,\phi'_{1,12})&(E_{6}(a_{3}),-\gamma,sgn)\\ &(T,\phi'_{4,7})&A_{4}+A_{1}\\ &(T,\phi''_{8,3})&D_{4}\\ &(T,\phi'_{9,6})&A_{4}\\ &(T,\phi_{12,4})&(D_{4}(a_{1}),\gamma,1)\\ &(T,\phi''_{6,6})&(D_{4}(a_{1}),\gamma,r)\\ &(T,\phi_{16,5})&(D_{4}(a_{1}),g_{2}\gamma,1)\\ &(T,\phi_{4,8})&(D_{4}(a_{1}),g_{2}\gamma,sgn)\\ &(T,\phi'_{6,6})&(D_{4}(a_{1}),g_{3}\gamma,1)\\ &(T,\phi''_{9,6})&A_{3}\\ &(T,\phi''_{4,7})&A_{2}+2A_{1}\\ &(T,\phi'_{8,9})&2A_{2}\\ &(T,\phi''_{8,9})&(A_{2},\gamma,1)\\ &(T,\phi''_{1,12})&(A_{2},\gamma,sgn)\\ &(T,\phi_{9,10})&(A_{2},-\gamma,1)\\ &(T,\phi'_{2,16})&(A_{2},-\gamma,sgn)\\ &(T,\phi_{4,13})&2A_{1}\\ &(T,\phi''_{2,16})&A_{1}\\ &(T,\phi_{1,24})&1\\ \end{array}$$.

\subsubsection{Un premier th\'eor\`eme\label{premiertheoremeSI}}
Pour $\boldsymbol{{\cal O}}\in {\cal U}_{sp}^{\Gamma_{F}}$, on note ${\cal B}_{Irr}(G,\boldsymbol{{\cal O}})$ l'image réciproque de $\bar{{\cal C}}^{Irr}(\boldsymbol{{\cal O}})$ par l'application $\boldsymbol{\nabla}_{F}$.  Ce sous-ensemble de ${\cal B}_{Irr}(G)$ s'identifie par la bijection ${\cal B}_{Irr}(G)\simeq {\cal B}'_{Irr}(G)$ \`a un sous-ensemble ${\cal B}'_{Irr}(G,\boldsymbol{{\cal O}})$ de ${\cal B}'_{Irr}(G)$.

 Soient $\boldsymbol{{\cal O}}\in {\cal U}_{sp}^{\Gamma_{F}}$,  $(\Lambda,\varphi,\rho)\in {\cal B}'_{Irr}(G,\boldsymbol{{\cal O}})$ et  $s\in S(\bar{C})$. On suppose que $ \Lambda\subset \Lambda(s)$.  On pose $M=M_{\Lambda}$.  
 Soit $\rho'\in Irr(W_{s}(M_{s_{\Lambda}}))$. Le triplet $(M_{s_{\Lambda}},{\cal E}_{\varphi},\rho')$ appartient \`a ${\cal I}^{FC}(G_{s})$. On note $(\bar{{\cal O}}_{\rho'},{\cal L}_{\rho'})$ son image par la repr\'esentation de Springer g\'en\'eralis\'ee de $G_{s}$. On a d\'efini l'application $\iota_{s,nil}:\mathfrak{g}_{s,nil}/conj\to (\boldsymbol{\mathfrak{g}}_{nil}/conj)^{I_{F}} $ en \ref{surFnr}. Posons  $\iota_{s,nil}(\bar{{\cal O}}_{\rho'})=\boldsymbol{{\cal O}}_{\rho'}$. On a aussi d\'efini la multiplicit\'e $m(\rho_{s},\rho')$ de $\rho'$ dans $\rho_{s}=\rho\circ j_{s}$.

 \begin{thm}{ Supposons $m(\rho_{s},\rho')>0$. Alors ou bien $dim(\boldsymbol{{\cal O}}_{\rho'})< dim(\boldsymbol{{\cal O}})$, ou bien $\boldsymbol{{\cal O}}=\boldsymbol{{\cal O}}_{\rho'}$.}\end{thm}
 
 On d\'emontrera ce th\'eor\`eme en m\^eme temps que  le suivant dans les paragraphes \ref{casnonramSI} \`a \ref{preuveE6rambis}.
 
 \subsubsection{Un deuxi\`eme th\'eor\`eme\label{deuxiemetheoremeSI}}
 
 Fixons $\boldsymbol{{\cal O}}\in {\cal U}_{sp}^{\Gamma_{F}}$. Notons $\bar{{\cal C}}(\boldsymbol{{\cal O}})$ l'ensemble des classes de conjugaison par ${\bf G}$ dans l'ensemble des couples $(N,d)$ où $N\in \boldsymbol{{\cal O}}$ et $d\in \tilde{\bar{A}}(N)$.  On d\'efinit une application
 $$\begin{array}{ccccc}{\cal B}_{Irr}(G,\boldsymbol{{\cal O}})&&\to&&\bar{{\cal C}}(\boldsymbol{{\cal O}})\\ (M,s_{M},\varphi,\rho)&\mapsto&\boldsymbol{\nabla}_{F}(M,s_{M},\varphi,\rho)=(N,d,\mu)&\mapsto&(N,d)\\ \end{array}$$
 Pour ${\bf d}\in  \bar{{\cal C}}(\boldsymbol{{\cal O}})$, on note 
 ${\cal B}_{Irr}(G,\boldsymbol{{\cal O}},{\bf d})$ la fibre de cette application au-dessus de ${\bf d}$. On a aussi la variante ${\cal B}'_{Irr}(G,\boldsymbol{{\cal O}},{\bf d})$. 
 
 Notons ${\cal J}_{F}$ l'ensemble des couples $(s,{\bar{\cal O}})$ o\`u $s\in S(\bar{C})$ et $\bar{{\cal O}}\in \mathfrak{g}_{s,nil}/conj$ est une orbite conserv\'ee par l'action de $\Gamma_{{\mathbb F}_{q}}$. 
 
 \begin{thm}{Il existe une application injective $b_{F}: \bar{{\cal C}}(\boldsymbol{{\cal O}})\to {\cal J}_{F}$ telle que les propri\'et\'es suivantes soient v\'erifi\'ees. Soit $ {\bf d}\in \bar{{\cal C}}(\boldsymbol{{\cal O}})$, posons $b_{F}({\bf d})=(s,\bar{{\cal O}})$. 
 
 (i)  On a $\iota_{s,nil}(\bar{{\cal O}})=\boldsymbol{{\cal O}}$.
 
 (ii) Soit  $(\Lambda,\varphi,\rho)\in {\cal B}'_{Irr}(G,\boldsymbol{{\cal O}},{\bf d})$.   Alors $\Lambda\subset \Lambda(s)$. 
  Il existe une unique repr\'esentation $\rho'\in Irr(W_{s}(M_{s_{\Lambda}}))$ telle que $<\rho_{s},\rho')>0$ et que l'image de $(M_{s_{\Lambda}},{\cal E}_{\varphi},\rho')$ par la correspondance de Springer g\'en\'eralis\'ee pour $G_{s}$ soit port\'ee par   $\bar{{\cal O}}$. Cette repr\'esentation $\rho'$ est conserv\'ee par l'action galoisienne de $\Gamma_{{\mathbb F}_{q}}$ et  on a  $<\rho_{s},\rho'>=1$. Posons  $Spr(M_{s_{\Lambda}},{\cal E}_{\varphi},\rho')=({\bar{\cal O}},{\cal L}_{\Lambda,\varphi,\rho})$. L'application $(\Lambda,\varphi,\rho)\mapsto {\cal L}_{\Lambda,\varphi,\rho}$ d\'efinie sur ${\cal B}'_{Irr}(G,\boldsymbol{{\cal O}},{\bf d})$ est injective.
  
  (iii) Soit $(\Lambda,\varphi,\rho)\in {\cal B}'_{Irr}(G,\boldsymbol{{\cal O}})- {\cal B}'_{Irr}(G,\boldsymbol{{\cal O}},{\bf d})$. Supposons $\Lambda\subset \Lambda(s)$. Alors il n'existe pas de repr\'esentation $\rho'\in Irr(W_{s}(M_{s_{\Lambda}}))$ telle que $<\rho_{s},\rho')>0$ et que l'image de $(M_{s_{\Lambda}},{\cal E}_{\varphi},\rho')$ par la correspondance de Springer g\'en\'eralis\'ee pour $G_{s}$ soit port\'ee par  $\bar{{\cal O}}$.}\end{thm}

   \subsubsection{Preuve des th\'eor\`emes \ref{premiertheoremeSI}  et \ref{deuxiemetheoremeSI} dans le cas o\`u $G$ est d\'eploy\'e sur $F^{nr}$\label{casnonramSI}}
On suppose que $G$  est d\'eploy\'e sur $F^{nr}$. La construction de l'application $\boldsymbol{\nabla}_{F}$ en \ref{parametragecalB}  associe \`a $(\Lambda,\varphi,\rho)\in {\cal B}'_{Irr}(G,\boldsymbol{{\cal O}})$ un \'el\'ement $({\bf M},\boldsymbol{{\cal E}},\rho)$ de ${\cal A}({\bf G},\boldsymbol{{\cal O}})$. D'apr\`es \ref{parametragecalB}(2), l'ensemble $\Lambda$ est celui associ\'e \`a ${\bf M}$ comme en \ref{Lambda}. Le triplet  
$({\bf M},\boldsymbol{{\cal E}},\rho)$ de ${\cal A}({\bf G},\boldsymbol{{\cal O}})$ s'identifie au triplet  $(\Lambda,\boldsymbol{{\cal E}},\rho)$ de ${\cal A}'({\bf G},\boldsymbol{{\cal O}})$. 

Soit $s\in S(\bar{C})$ tel que $\Lambda\subset \Lambda(s)$. Alors l'\'el\'ement $t={\bf j}_{T}(s)\in {\bf T}$ v\'erifie les hypoth\`eses de \ref{premiertheoreme}. On a $W_{t}({\bf M}_{\Lambda,t})=W_{s}(M_{\Lambda,s_{\Lambda}})$ et le th\'eor\`eme \ref{premiertheoremeSI}  est une traduction de la proposition \ref{premiertheoreme}.

Puisque $G$ est d\'eploy\'e sur $F^{nr}$, l'ensemble $\bar{{\cal C}}(\boldsymbol{{\cal O}})$ est le m\^eme qu'en \ref{deuxiemetheoreme}.   On a d\'efini  dans ce paragraphe une application $b: \bar{{\cal C}}(\boldsymbol{{\cal O}})\to {\cal J}$. Soit $ {\bf d}\in \bar{{\cal C}}(\boldsymbol{{\cal O}})$, posons $b({\bf d})=(\alpha,\boldsymbol{{\cal O}}_{\alpha})$.  Puisque $\boldsymbol{{\cal O}}$ est par hypothèse  conserv\'ee par $\Gamma_{F}^{nr}$, l'assertion (ii) du th\'eor\`eme \ref{deuxiemetheoreme} dit  que $\alpha\in \Delta_{a}$ est fixe par $\Gamma_{F}^{nr}$  et  que $\boldsymbol{{\cal O}}_{\alpha}$ est conserv\'ee par l'action galoisienne.   Il  est associ\'e à $\alpha$ un sommet $s_{\alpha}\in S(\bar{C})$. L'homomorphisme $\iota_{s_{\alpha}}:\mathfrak{g}_{s_{\alpha}}\to \boldsymbol{\mathfrak{g}}$ identifie $\mathfrak{g}_{s_{\alpha}}$ \`a $\boldsymbol{\mathfrak{g}}_{t_{\alpha}}$. Ainsi $\boldsymbol{{\cal O}}_{\alpha}$ s'identifie \`a une orbite $\bar{{\cal O}}\in \mathfrak{g}_{s_{\alpha},nil}/conj$. 
Alors $(s_{\alpha},\bar{{\cal O}})$ appartient \`a ${\cal J}_{F}$. On d\'efinit $b_{F}$ par $b_{F}({\bf d})=(s_{\alpha},\bar{{\cal O}})$. Le groupe $W_{s_{\alpha}}(M_{\Lambda,s_{\Lambda}})$ est \'egal \`a $W_{t_{\alpha}}({\bf M}_{\Lambda,t_{\alpha}})$ et ${\cal E}_{\varphi}$ est \'egal \`a $\boldsymbol{{\cal E}}_{t_{\alpha}}$. Alors le th\'eor\`eme \ref{deuxiemetheoremeSI} est une simple traduction du th\'eor\`eme \ref{deuxiemetheoreme}, \`a ceci pr\`es que celui-ci n'affirme pas que l'unique repr\'esentation $\rho'$ de l'assertion (iii) est conserv\'ee par l'action galoisienne. Mais cela r\'esulte de son unicit\'e. En effet, consid\'erons cette repr\'esentation $\rho'$ et soit $\sigma\in \Gamma_{{\mathbb F}_{q}}$. Il fait partie des hypoth\`eses que $\rho$ est conserv\'ee par l'action de $\Gamma_{{\mathbb F}_{q}}$. Donc $\rho_{s_{\alpha}}$ est conserv\'ee par $\sigma$. Puisque $\rho'$ intervient dans $\rho_{s_{\alpha}}$, $\sigma(\rho')$ intervient aussi. Il fait aussi partie des hypoth\`eses que ${\cal E}_{\varphi}$ est conserv\'e par l'action de $\Gamma_{{\mathbb F}_{q}}$. La correspondance de Springer g\'en\'eralis\'ee est \'equivariante pour les actions galoisiennes, donc l'image par cette correspondance de $(M_{s_{\Lambda}},{\cal E}_{\varphi},\sigma(\rho'))$ est port\'ee par $\sigma(\bar{{\cal O}})=\bar{{\cal O}}$. Alors $\sigma(\rho')$ v\'erifie les conditions qui caract\'erisent $\rho'$, donc $\sigma(\rho')=\rho'$. Cela ach\`eve la preuve du th\'eor\`eme  \ref{deuxiemetheoremeSI} dans le cas o\`u $G$ est d\'eploy\'e sur $F^{nr}$.

\subsubsection{Rappels sur la correspondance de Springer g\'en\'eralis\'ee pour les groupes classiques\label{groupesclassiques}}
Pour ce paragraphe, le corps de base est ${\mathbb F}_{q}$. On suppose que $G$ est un groupe symplectique ou sp\'ecial orthogonal d\'efini  et d\'eploy\'e sur ${\mathbb F}_{q}$, c'est-\`a-dire $G=Sp(N)$ ou $G=SO(N)$.  Dans le cas symplectique, on suppose que $N\geq2$ est pair. Dans le cas sp\'ecial orthogonal, on suppose $N\geq 3$. Puisque $G$ est d\'eploy\'e, on a ${\cal I}(G)={\cal I}_{{\mathbb F}_{q}}(G)$, ${\cal I}^{FC}(G)={\cal I}^{FC}_{{\mathbb F}_{q}}(G)$. 
On utilise les notations suivantes. Pour deux entiers $a\geq b\geq0$, on pose $[a,b]=(a,a-1,...,b+1,b)$. Si de plus, $a$ et $b$ sont de m\^eme parit\'e, on pose $[a,b]_{2}=(a,a-2,...,b+2,b)$. Pour une  partition $\mu=(\mu_{1},\mu_{2},...)$, on note $2\mu=(2\mu_{1},2\mu_{2},...)$. 

Supposons $G=Sp(N)$. On pose $N=2n$. L'ensemble ${\cal I}(G)$ s'identifie \`a celui des couples $(\lambda,\epsilon)$ o\`u $\lambda\in {\cal P}^{symp}(2n)$ et $\epsilon\in \{\pm1\}^{Jord_{bp}(\lambda)}$. L'ensemble ${\cal I}^{FC}(G)$ s'identifie \`a celui des triplets $(k,\alpha,\beta)$, o\`u $k\in {\mathbb N}$, $k(k+1)\leq 2n$ et $(\alpha,\beta)\in {\cal P}_{2}(n-k(k+1)/2)$. L'entier $k$ correspond \`a un Levi $M_{k}$ de $G$ tel que $M_{k,SC}$ est de type $C_{k(k+1)/2}$ et $(\alpha,\beta)$ d\'etermine une repr\'esentation irr\'eductible $\rho_{\alpha,\beta}$ de $W(M_{k})$ qui est un groupe de Weyl de type $C_{n-k(k+1)/2}$. 

 Consid\'erons une partition $\lambda\in {\cal P}^{symp}(2n)$.    On consid\`ere qu'elle \`a $2n$ termes, on pose $mult_{\lambda}(0)=2n-l(\lambda)$ et on note $Jord_{*}(\lambda)$ l'ensemble des $i\in {\mathbb N}$ tels que  $mult_{\lambda}(i)\geq1$, c'est-\`a-dire que $Jord_{*}(\lambda)$ peut contenir $0$.
On v\'erifie que $\lambda+[2n-1,0]$ poss\`ede $n$ termes pairs et $n$ termes impairs. On  introduit les deux partitions $z=(z_{1},...,z_{n})$ et $z'=(z'_{1},...,z'_{n})$ telles que $\lambda'+[2n-1,0]=\{2z_{i};i=1,...,n\}\cup \{2z'_{i}+1; i=1,...,n\}$. On pose ${\bf A}(\lambda)=(z'+[n+1,2])\cup \{0\}$, ${\bf B}(\lambda)=z+[n,1]$. Notons ${\bf A}(\lambda)=(a_{1},...,a_{n+1})$, ${\bf B}(\lambda)=(b_{1},...,b_{n})$. On v\'erifie que $a_{j}\geq a_{j+1}+2$ pour $j\leq n$ et $b_{j}\geq b_{j+1}+2$ pour $j\leq n-1$ et que $a_{j}\geq b_{j}\geq a_{j+1}$. Posons ${\bf A}(\lambda)\cup {\bf B}(\lambda)={\bf C}(\lambda)=(c_{1},...,c_{2n+1})$. 
  Soit $i\in Jord_{*}(\lambda)$. On note $j_{min}(i)$, resp. $j_{max}(i)$, le plus petit, resp. grand, $j\in \{1,...,2n+1\}$ tel que $\lambda_{j}=i$. Le nombre $i$  
  contribue \`a ${\bf A}(\lambda)\cup {\bf B}(\lambda)$ par  une sous-partition ${\bf C}(\lambda,i)=(c_{j_{min}(i)},...,c_{j_{max}(i)})$. On pose ${\bf A}(\lambda,i)={\bf A}(\lambda)\cap {\bf C}(\lambda,i)$, ${\bf B}(\lambda,i)={\bf B}(\lambda)\cap {\bf C}(\lambda,i)$. Supposons que $i$ est pair. Alors $ {\bf C}(\lambda,i)$ est l'intervalle $[i/2+2n+1-j_{min}(i),i/2+2n+1-j_{max}(i)]$. Ces \'el\'ements interviennent alternativement dans ${\bf A}(\lambda)$ et ${\bf B}(\lambda)$: si $c_{j}\in {\bf A}(\lambda)$, resp. $c_{j}\in {\bf B}(\lambda)$, alors $c_{j+1}\in {\bf B}(\lambda)$, resp $c_{j+1}\in {\bf A}(\lambda)$.  Le plus grand terme $c_{j_{min}(i)}$ est $a_{(j_{min}(i)+1)/2}$ si $j_{min}(i)$ est impair, $b_{j_{min}(i)/2}$ si $j_{min}(i)$ est pair. Le plus petit terme $c_{j_{max}(i)}$ est $a_{(j_{max}(i)+1)/2}$ si $j_{max}(i)$ est impair, $b_{j_{max}(i)/2}$ si $j_{max}(i)$ est pair. Supposons maintenant $i$ impair. La multiplicit\'e $mult_{\lambda}(i)$ est paire, donc $j_{min}(i)$ et $j_{max}(i)$ ne sont pas de m\^eme parit\'e. Alors ${\bf C}(\lambda,i)$ est la r\'eunion de deux ensembles \'egaux \`a $[(i+1)/2+2n-j_{min}(i),(i+1)/2+2n+1-j_{max}(i)]_{2}$. Si $j_{min}(i)$ est impair, on a $c_{j_{min}(i)}=c_{j_{min}(i)+1}=a_{(j_{min}(i)+1)/2}=b_{(j_{min}(i)+1)/2}$,  $c_{j_{max}(i)-1}=c_{j_{max}(i)}=a_{j_{max}(i)/2}=b_{j_{max}(i)/2}$.  Si $j_{min}(i)$ est pair, on a $c_{j_{min}(i)}=c_{j_{min}(i)+1}=a_{j_{min}(i)/2+1} =b_{j_{min}(i)/2}$, $c_{j_{max}(i)-1}=c_{j_{max}(i)}=a_{(j_{max}(i)+1)/2} =b_{(j_{max}(i)-1)/2}$. 
  
  Soit $m\in \{0,...,2n\}$.  Posons $\nu_{m}(\lambda)=1$ si   $S_{m}(\lambda) $ est impair  et $\nu_{m}(\lambda)=0$ sinon. Remarquons que la condition $S_{m}(\lambda)$ impair implique $\lambda_{m}$ impair et $\lambda_{m}=\lambda_{m+1}$. 
  A l'aide des formules ci-dessus,   on v\'erifie que
 $$(1) \qquad S_{m}(\lambda)=2S_{m}({\bf C}(\lambda))+m^2-m-4nm+\nu_{m}(\lambda).$$
  
  Soit $\epsilon\in \{\pm 1\}^{Jord_{bp}(\lambda)} $.   On pose 
  $${\bf A}(\lambda,\epsilon)=(\cup_{i\in Jord_{*}(\lambda)-Jord_{bp}(\lambda)}{\bf A}(\lambda,i))\cup(\cup_{i\in Jord_{bp}(\lambda); \epsilon_{i}=1}{\bf A}(\lambda,i))\cup(\cup_{i\in Jord_{bp}(\lambda); \epsilon_{i}=-1}{\bf B}(\lambda,i));$$
   $${\bf B}(\lambda,\epsilon)=(\cup_{i\in Jord_{*}(\lambda)-Jord_{bp}(\lambda)}{\bf B}(\lambda,i))\cup(\cup_{i\in Jord_{bp}(\lambda); \epsilon_{i}=1}{\bf B}(\lambda,i))\cup(\cup_{i\in Jord_{bp}(\lambda); \epsilon_{i}=-1}{\bf A}(\lambda,i)).$$

Consid\'erons maintenant un triplet $(k,\alpha,\beta)\in {\cal I}^{FC}(G)$. Parce que $(\alpha,\beta)\in {\cal P}_{2}(n-k(k+1)/2)$, on a $l(\alpha),l(\beta)\leq n-[\frac{k+1}{2}]$. On consid\`ere que $\alpha$ a $n+1+[\frac{k}{2}]$ termes et que $\beta$ en a $n-[\frac{k}{2}]$. Si $k$ est pair, on pose $A(\alpha)=\alpha+[2n+k,0]_{2}$, $B(\beta)=\beta+[2n-k-1,1]_{2}$. Si $k$ est impair, on pose $A(\alpha)=\alpha+[2n+k, 1]_{2}$, $B(\beta)=\beta+[2n-k-1,0]_{2}$. 
 Alors $(\lambda,\epsilon)=Spr(k,\alpha,\beta)$ est caract\'eris\'e par les propri\'et\'es suivantes:   $\lambda$   est l'unique \'el\'ement de ${\cal P}^{symp}(2n)$ telle que ${\bf A}(\lambda)\cup {\bf B}(\lambda)=A(\alpha)\cup B(\beta)$; $\epsilon$ est l'unique \'el\'ement de  $\{\pm 1\}^{Jord_{bp}(\lambda)}$ tel que 
 $$({\bf A}(\lambda,\epsilon),{\bf B}(\lambda,\epsilon))= \left\lbrace\begin{array}{cc}(A(\alpha),B(\beta)),&\text{ si }k\text{ est pair,}\\ (B(\beta),A(\alpha)),&\text{ si }k\text{ est  impair.}\\ \end{array}\right.$$
 
 Inversement, partons de $\epsilon\in \{\pm 1\}^{Jord_{bp}(\lambda)} $. Il existe un unique triplet $(k,\alpha,\beta)\in {\cal I}^{FC}(G)$ tel que $(\lambda,\epsilon)=Spr(k,\alpha,\beta)$ et il est utile de calculer l'entier $k$. Posons $J=\{i\in Jord_{bp}(\lambda); mult_{\lambda}(i)\text{ est impair et }\epsilon_{i}=-1\}$. Posons
$$D=\sum_{i\in J}(-1)^{mult_{\lambda}(\geq i)}.$$
Il r\'esulte de \cite{W1} lemme XI.4 que 

(2) $k=\vert 2D+\frac{1}{2}\vert -\frac{1}{2}$.

\noindent Cette \'egalit\'e implique que  $D\geq 0$ si $k$ est pair et $D<0$ si $k$ est impair.

Supposons maintenant $G=SO(N)$ avec $N$ impair. On pose $N=2n+1$. L'ensemble  ${\cal I}(G)$ s'identifie \`a celui des couples $(\lambda,\epsilon)$, o\`u $\lambda\in {\cal P}^{orth}(2n+1)$ et $\epsilon\in \{\pm 1\}^{Jord_{bp}(\lambda)}/\{1,-{\bf 1}\}$, en notant  $-{\bf 1}$ l'\'el\'ement de $\{\pm 1\}^{Jord_{bp}(\lambda)}$ dont toutes les composantes valent $-1$. L'ensemble ${\cal I}^{FC}(G)$ s'identifie \`a celui des triplets $(h,\alpha,\beta)$, o\`u $h\in {\mathbb N}$ est un entier impair tel que $h^2\leq 2n+1$ et $(\alpha,\beta)\in {\cal P}_{2}(n+(1-h^2)/2)$. L'entier $h$ correspond \`a un Levi $M_{h}$ de $G$ tel que $M_{h,SC}$ est de type $B_{(h^2-1)/2}$ et $(\alpha,\beta)$ d\'etermine une repr\'esentation irr\'eductible $\rho_{\alpha,\beta}$ de $W(M_{h})$ qui est un groupe de Weyl de type $B_{n+(1-h^2)/2}$.  

Consid\'erons une partition $\lambda\in {\cal P}^{orth}(2n+1)$. On consid\`ere qu'elle \`a $2n+1$ termes, on pose $mult_{\lambda}(0)=2n+1-l(\lambda)$ et on note $Jord_{*}(\lambda)$  l'ensemble des $i\in {\mathbb N}$ tels que  $mult_{\lambda}(i)\geq1$.
On v\'erifie que $\lambda+[2n,0]$ poss\`ede $n$ termes pairs et $n+1$ termes impairs. On  introduit les deux partitions $z=(z_{1},...,z_{n})$ et $z'=(z'_{1},...,z'_{n+1})$ telles que $\lambda+[2n,0]=\{2z_{i};i=1,...,n\}\cup \{2z'_{i}+1; i=1,...,n+1\}$.
 On pose ${\bf A}(\lambda)=z'+[n,0]$, ${\bf B}(\lambda)=z+[n-1,0]$. Notons ${\bf A}(\lambda)=(a_{1},...,a_{n+1})$, ${\bf B}(\lambda)=(b_{1},...,b_{n})$. On v\'erifie que $a_{j}\geq a_{j+1}+2$ pour $j\leq n$ et $b_{j}\geq b_{j+1}+2$ pour $j\leq n-1$ et que $a_{j}\geq b_{j}\geq a_{j+1}$. Posons ${\bf A}(\lambda)\cup {\bf B}(\lambda)={\bf C}(\lambda)=(c_{1},...,c_{2n+1})$. 
  Soit $i\in Jord_{*}(\lambda)$. On d\'efinit les termes 
   $j_{min}(i)$,  $j_{max}(i)$,   ${\bf C}(\lambda,i)$, ${\bf A}(\lambda,i) $ et ${\bf B}(\lambda,i) $ comme dans le cas symplectique. Supposons que $i$ est impair. Alors $ {\bf C}(\lambda,i)$ est l'intervalle $[(i+1)/2+2n-j_{min}(i),(i+1)/2+2n-j_{max}(i)]$. Ces \'el\'ements interviennent alternativement dans ${\bf A}(\lambda)$ et ${\bf B}(\lambda)$.  Le plus grand terme $c_{j_{min}(i)}$ est $a_{(j_{min}(i)+1)/2}$ si $j_{min}(i)$ est impair, $b_{j_{min}(i)/2}$ si $j_{min}(i)$ est pair. Le plus petit terme $c_{j_{max}(i)}$ est $a_{(j_{max}(i)+1)/2}$ si $j_{max}(i)$ est impair, $b_{j_{max}(i)/2}$ si $j_{max}(i)$ est pair. Supposons maintenant $i$ pair. La multiplicit\'e $mult_{\lambda}(i)$ est paire, donc $j_{min}(i)$ et $j_{max}(i)$ ne sont pas de m\^eme parit\'e. Alors ${\bf C}(\lambda,i)$ est la r\'eunion de deux ensembles \'egaux \`a $[i/2+2n-j_{min}(i),i/2+2n+1-j_{max}(i)]_{2}$. Si $j_{min}(i)$ est impair, on a $c_{j_{min}(i)}=c_{j_{min}(i)+1}=a_{(j_{min}(i)+1)/2}=b_{(j_{min}(i)+1)/2}$, $c_{j_{max}(i)-1}=c_{j_{max}(i)}=a_{j_{max}(i)/2}=b_{j_{max}(i)/2}$. Si $j_{min}(i)$ est pair, on a $c_{j_{min}(i)}=c_{j_{min}(i)+1}= a_{j_{min}(i)/2+1}= b_{j_{min}(i)/2}$, $c_{j_{max}(i)-1}=c_{j_{max}(i)}= a_{(j_{max}(i)+1)/2}=b_{(j_{max}(i)-1)/2}$. 
  
  Soit $m\in \{1,...,2n+1\}$.  Posons $\nu_{m}(\lambda)=1$ si $\lambda_{m}$ est pair et $m\equiv S_{m}(\lambda)+1\,\,mod\,\,2{\mathbb Z}$ et $\nu_{m}(\lambda)=0$ sinon. Notons que, si $\nu_{m}(\lambda)=1$, on a $\lambda_{m}=\lambda_{m+1}$. On pose aussi $\nu_{0}(\lambda)=0$.  Pour $m\in \{0,...,2n+1\}$,    on v\'erifie  \`a l'aide des formules ci-dessus que
 $$(3) \qquad S_{m}(\lambda)=2S_{m}({\bf C}(\lambda))+m^2-4nm+\nu_{m}(\lambda).$$
  
  Soit $\epsilon\in \{\pm 1\}^{Jord_{bp}(\lambda)}/\{1,-{\bf 1}\}$. Relevons $\epsilon$ en un \'el\'ement de $\{\pm 1\}^{Jord_{bp}(\lambda)}$. On d\'efinit alors ${\bf A}(\lambda,\epsilon)$ et ${\bf B}(\lambda,\epsilon)$ par les m\^emes formules que dans le cas symplectique.       Changer de rel\`evement $\epsilon$ revient \`a permuter ${\bf A}(\lambda,\epsilon)$ et ${\bf B}(\lambda,\epsilon)$.

Consid\'erons maintenant un triplet $(h,\alpha,\beta)\in {\cal I}^{FC}(G)$. Parce que $(\alpha,\beta)\in {\cal P}_{2}(n+(1-h^2)/2)$, on a $l(\alpha),l(\beta)\leq n+(1-h)/2$. On consid\`ere que $\alpha$ a $n+(1+h)/2$ termes et que $\beta$ en a $n+(1-h)/2$. On pose $A(\alpha)=\alpha+[2n-1+h,0]_{2}$, $B(\beta)=\beta+[2n-1-h,0]_{2}$.  Alors $(\lambda,\epsilon)=Spr(h,\alpha,\beta)$ est caract\'eris\'e par les propri\'et\'es suivantes:   $\lambda$   est l'unique \'el\'ement de ${\cal P}^{orth}(2n+1)$ telle que ${\bf A}(\lambda)\cup {\bf B}(\lambda)=A(\alpha)\cup B(\beta)$; $\epsilon$ est l'unique \'el\'ement de  $\{\pm 1\}^{Jord_{bp}(\lambda)}/\{1,-{\bf 1}\}$ tel que $(A(\alpha),B(\beta))$ soit \'egal \`a  $({\bf A}(\lambda,\epsilon),{\bf B}(\lambda,\epsilon))$ \`a permutation pr\`es des deux facteurs.   

 Inversement, partons de $\epsilon\in \{\pm 1\}^{Jord_{bp}(\lambda)}/\{1,-{\bf 1}\} $. Il existe un unique triplet $(h,\alpha,\beta)\in {\cal I}^{FC}(G)$ tel que $(\lambda,\epsilon)=Spr(h,\alpha,\beta)$. L'entier $h$ se calcule de la fa\c{c}on suivante. On rel\`eve $\epsilon$ en un \'el\'ement de $\{\pm 1\}^{Jord_{bp}(\lambda)}$. On d\'efinit l'ensemble $J$ et le nombre $D$ comme dans le cas symplectique. D'apr\`es \cite{W1} lemme XI.4, on a
 
 (4) $h=\vert 2D+1\vert $.
 
 \noindent Si $D\geq 0$, on a $({\bf A}(\lambda,\epsilon) ,{\bf B}(\lambda,\epsilon))=(A(\alpha),B(\beta))$. Si  $D<0$,on a $({\bf A}(\lambda,\epsilon) ,{\bf B}(\lambda,\epsilon))=(B(\beta),A(\alpha))$.

 Supposons enfin $G=SO(N)$ avec $N$ pair. On pose $N=2n$. L'ensemble ${\cal I}(G)$ s'envoie surjectivement sur celui des couples $(\lambda,\epsilon)$ o\`u $\lambda\in {\cal P}^{orth}(2n)$ et $\epsilon\in \{\pm 1\}^{Jord_{bp}(\lambda)}/\{1,-{\bf 1}\}$.  Les fibres de cette application ont un seul \'el\'ement sauf au-dessus d'une partition $\lambda$ dont tous les termes sont pairs, auquel cas la fibre a deux \'el\'ements (pour une telle partition, $Jord_{bp}(\lambda)=\emptyset$ et le terme $\epsilon$ dispara\^{\i}t). 
 L'ensemble ${\cal I}^{FC}(G)$ s'envoie surjectivement sur  celui des classes d'\'equivalence de triplets $(h,\alpha,\beta)$, o\`u $h\in {\mathbb N}$, $h$ est  pair, $h^2\leq 2n$ et $(\alpha,\beta)\in {\cal P}_{2}(n-h^2/2)$. Deux triplets $(h_{1},\alpha_{1},\beta_{1})$ et $(h_{2},\alpha_{2},\beta_{2})$ sont \'equivalents si et seulement s'ils sont \'egaux ou si $h_{1}=h_{2}=0$ et $(\alpha_{1},\beta_{1})=(\beta_{2},\alpha_{2})$. 
  Les fibres de cette application   ont un seul \'el\'ement sauf au-dessus d'un triplet de la forme $(0,\alpha,\alpha)$, auquel cas la fibre a deux \'el\'ements. 
Si $h\geq2$,  $h$ correspond \`a un Levi $M_{h}$ de $G$ tel que $M_{h,SC}$ est de type  $D_{h^2/2}$  et $(\alpha,\beta)$ d\'etermine une repr\'esentation irr\'eductible $\rho_{\alpha,\beta}$ de $W(M_{h})$ qui est un groupe de Weyl de type $B_{n-h^2/2}$. Si $h=0$, le Levi $M_{0}$ est un tore maximal d\'eploy\'e  $T$ de $G$. Le groupe de Weyl $W(M_{0})=W$ est de type $D_{n}$. Si $\alpha\not=\beta$, le couple $(\alpha,\beta)$ d\'etermine une repr\'esentation irr\'eductible $\rho_{\alpha,\beta}$ de $ W$ qui est la restriction de la repr\'esentation similaire d'un groupe de Weyl de type $B_{n}$. On a $\rho_{\alpha,\beta}=\rho_{\beta,\alpha}$.  Si $\alpha=\beta$, cette restriction se d\'ecompose en deux repr\'esentations irr\'eductibles de  $ W$. Chacune d'elles correspond \`a un \'el\'ement de ${\cal I}^{FC}(G)$.  Comme on le voit, il y a dans ces descriptions des \'el\'ements exceptionnels dont on peut dire qu'ils se d\'edoublent. Pour ne pas alourdir la r\'edaction, on identifiera ${\cal I}(G)$, resp. ${\cal I}^{FC}(G)$, \`a l'ensemble des paires $(\lambda,\epsilon)$, resp. des triplets $(h,\alpha,\beta)$ comme ci-dessus, en consid\'erant que certaines paires ou certains triplets sont d\'edoubl\'es. On esp\`ere que cela ne cr\'eera pas de confusion.  

  Consid\'erons une partition $\lambda\in {\cal P}^{orth}(2n)$.    On consid\`ere qu'elle \`a $2n$ termes, on pose $mult_{\lambda}(0)=2n-l(\lambda)$ et on note $Jord_{*}(\lambda)$ l'ensemble des $i\in {\mathbb N}$ tels que  $mult_{\lambda}(i)\geq1$. On v\'erifie que $\lambda+[2n-1,0]$ poss\`ede $n$ termes pairs et $n$ termes impairs. On  introduit les deux partitions $z=(z_{1},...,z_{n})$ et $z'=(z'_{1},...,z'_{n})$ telles que $\lambda+[2n-1,0]=\{2z_{i};i=1,...,n\}\cup \{2z'_{i}+1; i=1,...,n\}$. 
   On pose ${\bf A}(\lambda)=(z+[n-1,0])$, ${\bf B}(\lambda)=z'+[n-1,0]$. Notons ${\bf A}(\lambda)=(a_{1},...,a_{n})$, ${\bf B}(\lambda)=(b_{1},...,b_{n})$. On v\'erifie que $a_{j}\geq a_{j+1}+2$ pour $j\leq n-1$ et $b_{j}\geq b_{j+1}+2$ pour $j\leq n-1$ et que $a_{j}\geq b_{j}\geq a_{j+1}$. Posons ${\bf A}(\lambda)\cup {\bf B}(\lambda)={\bf C}(\lambda)=(c_{1},...,c_{2n})$. 
  Soit $i\in Jord_{*}(\lambda)$. On d\'efinit $j_{min}(i)$,   $j_{max}(i)$,   ${\bf C}(\lambda,i) $, ${\bf A}(\lambda,i) $ et ${\bf B}(\lambda,i) $ comme dans le cas symplectique. Supposons que $i$ est impair. Alors $ {\bf C}(\lambda,i)$ est l'intervalle $[(i-1)/2+2n-j_{min}(i),(i-1)/2+2n-j_{max}(i)]$. Ces \'el\'ements interviennent alternativement dans ${\bf A}(\lambda)$ et ${\bf B}(\lambda)$.   Le plus grand terme $c_{j_{min}(i)}$ est $a_{(j_{min}(i)+1)/2}$ si $j_{min}(i)$ est impair, $b_{j_{min}(i)/2}$ si $j_{min}(i)$ est pair. Le plus petit terme $c_{j_{max}(i)}$ est $a_{(j_{max}(i)+1)/2}$ si $j_{max}(i)$ est impair, $b_{j_{max}(i)/2}$ si $j_{max}(i)$ est pair. Supposons maintenant $i$ pair. La multiplicit\'e $mult_{\lambda}(i)$ est paire, donc $j_{min}(i)$ et $j_{max}(i)$ ne sont pas de m\^eme parit\'e. Alors ${\bf C}(\lambda,i)$ est la r\'eunion de deux ensembles \'egaux \`a $[i/2+2n-1-j_{min}(i),i/2+2n-j_{max}(i)]_{2}$. Si $j_{min}(i)$ est impair, on a $c_{j_{min}(i)}=c_{j_{min}(i)+1}=a_{(j_{min}(i)+1)/2}=b_{(j_{min}(i)+1)/2}$,  $c_{j_{max}(i)-1}=c_{j_{max}(i)}=a_{j_{max}(i)/2}=b_{j_{max}(i)/2}$.  Si $j_{min}(i)$ est pair, on a $c_{j_{min}(i)}=c_{j_{min}(i)+1}=a_{j_{min}(i)/2+1} =b_{j_{min}(i)/2}$, $c_{j_{max}(i)-1}=c_{j_{max}(i)}=a_{(j_{max}(i)+1)/2} =b_{(j_{max}(i)-1)/2}$.

 Soit $m\in \{1,...,2n\}$.  Posons $\nu_{m}(\lambda)=1$ si   $\lambda_{m}$ est pair et $S_{m}(\lambda)\equiv m+1\, mod\,2{\mathbb Z} $   et $\nu_{m}(\lambda)=0$ sinon. Si $\nu_{m}(\lambda)=1$, on a $\lambda_{m+1}=\lambda_{m}$.  On pose aussi $\nu_{0}(\lambda)=0$. Pour $m\in \{0,...,2n\}$, on v\'erifie \`a l'aide des formules ci-dessus que
 $$(5) \qquad S_{m}(\lambda)=2S_{m}({\bf C}(\lambda))+m^2+2m-4nm+\nu_{m}(\lambda).$$
 
  Soit $\epsilon\in \{\pm 1\}^{Jord_{bp}(\lambda)}/\{1,-{\bf 1}\} $. On rel\`eve $\epsilon$ en un \'el\'ement de  $\{\pm 1\}^{Jord_{bp}(\lambda)}$.   On d\'efinit ${\bf A}(\lambda,\epsilon)$ et ${\bf B}(\lambda,\epsilon)$ comme dans le cas symplectique.  
    Changer de rel\`evement $\epsilon$ \'echange ${\bf A}(\lambda,\epsilon)$ et ${\bf B}(\lambda,\epsilon)$.

 Consid\'erons maintenant un triplet $(h,\alpha,\beta)\in {\cal I}^{FC}(G)$. Parce que $(\alpha,\beta)\in {\cal P}_{2}(n-h^2/2)$, on v\'erifie que $l(\alpha),l(\beta)\leq n- h/2$. On consid\`ere que $\alpha$ a $n+h/2$ termes et que $\beta$ en a $n-h/2$.  On pose $A(\alpha)=\alpha+[2n+ h-2,0]_{2}$, $B(\beta)=\beta+[2n- h-2,0]_{2}$.   Alors $(\lambda,\epsilon)=Spr(h,\alpha,\beta)$ est caract\'eris\'e par les propri\'et\'es suivantes:  $\lambda$ est l'unique \'el\'ement de ${\cal P}^{orth}(2n)$ tel que $A(\alpha)\cup B(\beta)={\bf C}(\lambda)$; $\epsilon$ est l'unique \'el\'ement de $\{\pm 1\}^{Jord_{bp}(\lambda)}/\{1,-{\bf 1}\} $ tel que $(A(\alpha),B(\beta))$ soit \'egal \`a $(({\bf A}(\lambda,\epsilon),{\bf B}(\lambda,\epsilon))$ ou à $({\bf B}(\lambda,\epsilon), {\bf A}(\lambda,\epsilon))$. Remarquons que les termes qui se d\'edoublent se correspondent. Pr\'ecis\'ement, si $h=0$ et $\alpha=\beta$, on a $Spr(0,\alpha,\alpha)=2\alpha$ (on a d\'ej\`a dit que, pour une partition $\lambda=2\alpha$, le terme $\epsilon$ dispara\^{\i}t). 
 
  Inversement, partons de $\epsilon\in \{\pm 1\}^{Jord_{bp}(\lambda)}/\{1,-{\bf 1}\} $. Il existe un unique triplet $(h,\alpha,\beta)\in {\cal I}^{FC}(G)$ tel que $(\lambda,\epsilon)=Spr(h,\alpha,\beta)$. L'entier $h$ se calcule de la fa\c{c}on suivante. On rel\`eve $\epsilon$ en un \'el\'ement de $\{\pm 1\}^{Jord_{bp}(\lambda)}$.  On d\'efinit l'ensemble $J$ et le nombre $D$ comme dans le cas symplectique.  D'apr\`es \cite{W1} lemme XI.4, on a

 (6) $h=\vert 2D\vert $.
 
 \noindent Si $D>0$, on a $({\bf A}(\lambda,\epsilon) ,{\bf B}(\lambda,\epsilon))=(A(\alpha),B(\beta))$. Si  $D<0$,on a $({\bf A}(\lambda,\epsilon) ,{\bf B}(\lambda,\epsilon))=(B(\beta),A(\alpha))$. Si $D=0$, seule compte la paire $(\alpha,\beta)$ \`a permutation pr\`es et on peut donc supposer que $({\bf A}(\lambda,\epsilon) ,{\bf B}(\lambda,\epsilon))=(A(\alpha),B(\beta))$.

 \subsubsection{Preuve des th\'eor\`emes \ref{premiertheoremeSI} et \ref{deuxiemetheoremeSI} dans le cas $(D_{n},ram)$ \label{preuveDnram}}

On revient \`a notre groupe $G$ d\'efini sur $F$. On suppose que $G$ est de type $(D_{n},ram)$  avec $n\geq4$. En \ref{DnramSI}, on  a identifi\'e ${\cal B}_{Irr}(G)$ \`a l'ensemble des triplets $(h,\alpha,\beta)$,o\`u $h\in {\mathbb N}$ est un entier impair tel que $h^2\geq n$ et $(\alpha,\beta)\in {\cal P}_{2}(n-h^{2})$. L'entier $h$ d\'etermine un Levi de $G$ auquel est associ\'e un  ensemble $\Lambda_{h}\subset \Delta_{a}^{nr}$ et $(\alpha,\beta)$ d\'etermine une repr\'esentation irr\'eductible $\rho_{\alpha,\beta}$ de $W^{I_{F}}(M_{\Lambda_{h}})$ qui est un groupe de Weyl de type $B_{n-h^2}$. 
Pour un tel triplet, on a explicit\'e $\boldsymbol{\nabla}_{F}(h,\alpha,\beta)=(\lambda,d,\mu)$ o\`u $\lambda\in {\cal P}^{orth}(2n)$ est une partition orthogonale sp\'eciale ayant au moins un terme impair, $d\in \{\pm 1\}^{Int(\lambda)}_{-1}$ et $\mu\in \{\pm 1\}^{Int(\lambda)}/\{1,-{\bf 1}\}$. 

L'ensemble $S(\bar{C})$ s'identifie \`a celui des couples $(n',n'')\in {\mathbb N}^2$ tels que $n'+n''=n-1$. Pour le sommet $s$ param\'etr\'e par $(n',n'')$, on a $G_{s}=SO(2n'+1)\times SO(2n''+1)$. Notons $G'_{s}$ et $G''_{s}$ ces deux composantes. On a explicit\'e en \ref{groupesclassiques} les ensembles ${\cal I}(G'_{s})$, ${\cal I}^{FC}(G'_{s})$ et les ensembles analogues pour $G''_{s}$. On ajoute des indices $s$ et des $'$ ou $''$ aux notations de cette r\'ef\'erence. 
  L'ensemble $\mathfrak{g}_{s,nil}/conj$ est param\'etr\'e par ${\cal P}^{orth}(2n'+1)\times {\cal P}^{orth}(2n''+1)$. On a d\'efini l'application $\bar{c}_{s}:\mathfrak{g}_{s,nil}/conj\to \bar{{\cal C}}$ en \ref{surFnr}. D'apr\`es \ref{csram}, l'image de 
$(\lambda', \lambda'')\in {\cal P}^{orth}(2n'+1)\times {\cal P}^{orth}(2n''+1)$ par $\bar{c}_{s}$ est $(\lambda,d)$, o\`u $\lambda=\lambda'\cup \lambda''$ et $ d\in \{\pm 1\}^{Int(\lambda)}_{-1}$ est l'\'el\'ement d\'efini  par la formule suivante: pour $\Delta\in Int(\lambda)$,   $d_{\Delta}=(-1)^{mult_{\lambda''}(\Delta)}$. Remarquons que l'on a aussi $ d_{\Delta}=(-1)^{mult_{\lambda'}(\Delta)}$ puisque $mult_{\lambda'}(\Delta)+mult_{\lambda''}(\Delta)=mult_{\lambda}(\Delta)$ qui est pair.

Pour un entier $h$ impair tel que $h^2\leq n$, on a $\Lambda_{h}\subset \Lambda(s)$ si et seulement si $h^2\leq 2n'+1$ et $h^2\leq 2n''+1$. Si ces conditions sont v\'erifi\'ees, $M_{\Lambda_{h},s_{\Lambda_{h}}}$ s'identifie \`a $M'_{s,h}\times M''_{s,h}$ et le groupe $W'_{s}(M'_{s,h})\times W''_{s}(M''_{s,h})$ appara\^{\i}t naturellement comme un sous-groupe de $W^{I_{F}}(M_{\Lambda_{h}})$. 

\begin{prop}{Soient $(h,\alpha,\beta)\in {\cal B}_{Irr}(G)$ et $(n',n'')$ un couple param\'etrant un sommet $s\in S(\bar{C})$. On pose $\boldsymbol{\nabla}_{F}(h,\alpha,\beta)=(\lambda,d,\mu)$.

(i) Supposons $\Lambda_{h}\subset \Lambda(s)$, soient $(h,\alpha',\beta')\in {\cal I}^{FC}(G'_{s})$ et $(h,\alpha'',\beta'')\in {\cal I}^{FC}(G''_{s})$, notons $(\lambda',\epsilon')$ et $(\lambda'',\epsilon'')$ leurs images par la correspondance de Springer g\'en\'eralis\'ee. Supposons que $\rho_{\alpha',\beta'}\otimes \rho_{\alpha'',\beta''}$ soit une composante irr\'eductible de la restriction de $\rho_{\alpha,\beta}$ \`a $W'_{s}(M'_{s,h})\times W''_{s}(M''_{s,h})$. Alors $\lambda'\cup \lambda''\leq \lambda$ et, si $\lambda'\cup \lambda''=\lambda$, on a l'\'egalit\'e $\bar{c}_{s}(\lambda',\lambda'')=(\lambda,d)$. 

(ii) Inversement, soit $(\lambda',\lambda'')\in {\cal P}^{orth}(2n'+1)\times {\cal P}^{orth}(2n''+1)$, supposons $\bar{c}_{s}(\lambda',\lambda'')=(\lambda,d)$. Alors $\Lambda_{h}\subset \Lambda(s)$ et il  existe d'uniques triplets $(h,\alpha',\beta')\in {\cal I}^{FC}(G'_{s})$ et $(h,\alpha'',\beta'')\in {\cal I}^{FC}(G''_{s})$ v\'erifiant les conditions suivantes:

(a) leurs images par la correspondance de Springer g\'en\'eralis\'ee sont de la forme $(\lambda',\epsilon')$ et $(\lambda'',\epsilon'')$;

(b) $\rho_{\alpha',\beta'}\otimes \rho_{\alpha'',\beta''}$ est une composante irr\'eductible de la restriction de $\rho_{\alpha,\beta}$ \`a $W'_{s}(M'_{s,h})\times W''_{s}(M''_{s,h})$.

Pour ces triplets, $\rho_{\alpha',\beta'}\otimes \rho_{\alpha'',\beta''}$  intervient avec  multiplicit\'e $1$ dans la restriction de $\rho_{\alpha,\beta}$.  Pour $i\in Jord_{bp}(\lambda')$, resp $i\in Jord_{bp}(\lambda'')$, on a l'\'egalit\'e $\epsilon'_{i}=\mu_{\Delta}$, resp.  $\epsilon''_{i}=\mu_{\Delta}$, o\`u $\Delta$ est l'intervalle de $\lambda$ tel que $i\in \Delta$.}\end{prop}

{\bf Remarque.} Comme on l'a dit,  les termes $\epsilon'$, $\epsilon''$ et $\mu$ appartiennent \`a des quotients par des groupes $\{1,-{\bf 1}\}$. La derni\`ere   assertion du th\'eor\`eme signifie que l'on peut choisir des rel\`evements de ces termes qui v\'erifient les formules indiqu\'ees. 

  \bigskip

On a besoin d'un pr\'eliminaire.  Explicitons les ensembles ${\bf X}(\lambda)$ et ${\bf Y}(\lambda)$ associ\'es en \ref{DnramSI} \`a la partition sp\'eciale $\lambda\in {\cal P}^{orth}(2n)$.    On consid\`ere que $\lambda$ a $2n$ termes et que $mult_{\lambda}(0)=2n-l(\lambda)$. On note  $Jord_{*}(\lambda)$ l'ensemble des $i\in {\mathbb N}$ tels que $mult_{\lambda}(i)>0$. 
L'ensemble $Jord_{bp}(\lambda)$ se d\'ecompose en union disjointe d'intervalles et on a not\'e $Int(\lambda)$ l'ensemble de ces intervalles. Posons $Int_{*}(\lambda)=Int(\lambda)\cup (Jord_{*}(\lambda)-Jord_{bp}(\lambda))$. 
Soit $\Delta\in Int_{*}(\lambda)$. Notons $j_{min}(\Delta)$, resp. $j_{max}(\Delta)$, le plus petit, resp. grand, entier $j\geq1$ tel que $\lambda_{j}\in \Delta$. Les entiers $j_{min}(\Delta)$ et $j_{max}(\Delta)$ sont respectivement  impair et pair. Si $\Delta\not\in Int(\lambda)$, c'est-\`a-dire $\Delta=\{i\}$ avec $i$ pair, on a $\lambda_{j}=i$ pour tout $j\in \{j_{min}(\Delta),...,j_{max}(\Delta)\}$. Supposons $\Delta\in Int(\lambda)$. Alors, 
pour tout entier $j$ pair tel que $j_{min}(\Delta)<j<j_{max}(\Delta)$, on a $\lambda_{j}=\lambda_{j+1}$  et ce terme appartient \`a $\Delta$.    Tout \'el\'ement $\Delta\in Int_{*}(\lambda)$ contribue \`a ${\bf X}(\lambda)$, resp. ${\bf Y}(\lambda)$, par la sous-partition ${\bf X}(\lambda,\Delta)=(x_{(j_{min}(\Delta)+1)/2},...,x_{j_{max}(\Delta)/2})$, resp. ${\bf Y}(\lambda,\Delta)=(y_{(j_{min}(\Delta)+1)/2},...,y_{j_{max}(\Delta)/2})$. 
Si $\Delta=\{i\}$, avec $i$ pair, on a $x_{j}=y_{j}=i/2+n-j$ pour tout $j\in \{(j_{min}(\Delta)+1)/2,...,j_{max}(\Delta)/2\}$. Supposons $\Delta\in Int(\lambda)$. Alors $x_{j}=\frac{\lambda_{2j-1}+1}{2}+n-j$ et $y_{j}=\frac{\lambda_{2j}}{2}+n-1-j$  pour tout $j\in \{(j_{min}(\Delta)+1)/2,...,j_{max}(\Delta)/2\}$. On constate que
$$x_{(j_{min}(\Delta)+1)/2}> y_{(j_{min}(\Delta)+1)/2}=x_{(j_{min}(\Delta)+1)/2+1}>... $$
$$> y_{j_{max}(\Delta)/2-1}=x_{j_{max}(\Delta)/2}>y_{j_{max}(\Delta)/2}.$$
Soit $m\in \{0,...,2n\}$. Posons $\boldsymbol{\nu}_{m}(\lambda)=1$ si $\lambda_{m}$ est pair et $m$ est impair, $\boldsymbol{\nu}_{m}(\lambda)=0$ sinon. A l'aide des formules ci-dessus, on calcule
 $$(1) \qquad S_{m}(\lambda)=2S_{m}({\bf X}(\lambda)\cup {\bf Y}(\lambda))+[\frac{m}{2}]-m(4n-1-m)/2+\boldsymbol{\nu}_{m}(\lambda).$$

{\bf Preuve de la proposition.} On consid\`ere les donn\'ees de l'assertion  (i). 
On veut prouver que $\lambda'\cup \lambda''\leq \lambda$. C'est \'equivalent \`a l'assertion suivante. Soient $m'\in \{0,...,2n'+1\}$ et  $m''\in \{0,...,2n''+1\}$. Posons $m=m'+m''$. Alors
 
 (2) $S_{m'}(\lambda')+S_{m''}(\lambda'')\leq S_{m}(\lambda)$. 
 
 Le nombre $S_{m'}(\lambda')$ est calcul\'e par \ref{groupesclassiques}(3). 
 Puisque ${\bf C}(\lambda')=A(\alpha')\cup B(\beta')$, il existe deux entiers $x'\in \{0,...,n'+(1+h)/2\}$ et $y'\in \{0,...,n'+(1-h)/2\}$ tels que $x'+y'=m'$ et $S_{m'}({\bf C}(\lambda'))=S_{x'}(A(\alpha'))+S_{y'}(B(\beta'))$. On calcule
 $$S_{x'}(A(\alpha'))=S_{x'}(\alpha')+x'(2n'+h-x'),\, S_{y'}(B(\beta'))=S_{y'}(\beta')+y'(2n'-h-y').$$
 D'o\`u
 $$S_{m'}(\lambda')=2S_{x'}(\alpha')+2S_{y'}(\beta')+e',$$
 o\`u
 $$e'=2x'(2n'+h-x')+2y'(2n'-h-y')+
 {m'}^2-4n'm'+\nu_{m'}(\lambda')$$
 $$=2x'(h-x')+2y'(-h-y')+
 {m'}^2+\nu_{m'}(\lambda').$$
 En introduisant des termes analogues pour la partition $\lambda''$, on obtient
 $$(3) \qquad S_{m'}(\lambda')+S_{m''}(\lambda'')=2S_{x'}(\alpha')+2S_{y'}(\beta')+2S_{x''}(\alpha'')+2S_{y''}(\beta'')+
  e'+e''.$$
 L'hypoth\`ese que $\rho_{\alpha',\beta'}\otimes \rho_{\alpha'',\beta''}$ intervient dans la restriction de $\rho_{\alpha,\beta}$ implique que $\alpha\geq \alpha'\cup \alpha''$ et $\beta\geq \beta'\cup \beta''$, cf. par exemple \cite{W3} 1.13(38). En posant $x=x'+x''$ et $y=y'+y''$, on a alors 
 $$(4) \qquad S_{x'}(\alpha')+S_{x''}(\alpha'')\leq S_{x}(\alpha),\, S_{y'}(\beta')+S_{y''}(\beta'')\leq S_{y}(\beta),$$
 d'o\`u
 $$(5) \qquad S_{m'}(\lambda')+S_{m''}(\lambda'')\leq 2S_{x}(\alpha)+2S_{y}(\beta)+ e'+e''.$$
 
  Le nombre $S_{m}(\lambda)$ est calcul\'e par (1). 
 Puisque ${\bf X}(\lambda)\cup {\bf Y}(\lambda)=X(\alpha)\cup Y(\beta)$ et que $x+y=m$ par construction de ces termes, on a
 $$(6) \qquad S_{m}({\bf X}(\lambda)\cup {\bf Y}(\lambda))\geq S_{x}(X(\alpha))+S_{y}(Y(\beta)).$$
 On calcule
 $$S_{x}(X(\alpha))=S_{x}(\alpha)+ x(n+h)-x(x+1)/2,\, S_{y}(Y(\beta))=S_{y}(\beta)+y(n-h)-y(y+1)/2.$$
 D'o\`u
 $$S_{m}(\lambda)\geq 2S_{x}(\alpha)+2S_{y}(\beta)+e,$$
 o\`u
 $$e=2x(n+h)-x(x+1)+2y(n-h)-y(y+1)+[\frac{m}{2}]-m(4n-1-m)/2+\boldsymbol{\nu}_{m}(\lambda)$$
 $$=2xh-x^2-2yh-y^2+[\frac{m}{2}]-m/2+m^2/2+\boldsymbol{\nu}_{m}(\lambda).$$
 Gr\^ace \`a (5), on en d\'eduit 
 $$(7) \qquad S_{m'}(\lambda')+S_{m''}(\lambda'')\leq S_{m}(\lambda)+E,$$
 o\`u $E=e'+e''-e$. 
 On calcule 
 $$(8) \qquad E=-\frac{(x'-y'-x''+y'')^2}{2}+\frac{m}{2}-[\frac{m}{2}]+\nu_{m'}(\lambda')+\nu_{m''}(\lambda'')-\boldsymbol{\nu}_{m}(\lambda).$$
 On a $m=x'+y'+x''+y''$ donc $x'-y'-x''+y''$ est de m\^eme parit\'e que $m$ et la somme des trois premiers termes de $E$ est n\'egative ou nulle. Si $\nu_{m'}(\lambda')=\nu_{m''}(\lambda'')=0$, alors $E\leq 0$ et (2) r\'esulte de (7). Supposons $\nu_{m'}(\lambda')=1$. Par d\'efinition de ce terme, on a alors $1\leq m'\leq 2n'$, $\lambda'_{m'}$ est pair, $\lambda'_{m'+1}=\lambda'_{m'}$ et $\nu_{m'-1}(\lambda')=\nu_{m'+1}(\lambda')=0$. Montrons que
 
 (9) si les analogues de (2) pour les couples $(m'-1,m'')$ et $(m'+1,m'')$ sont v\'erifi\'ees, alors (2) est v\'erifi\'ee pour le couple $(m',m'')$. 
 
 Supposons que (2) ne soit pas v\'erifi\'ee pour $(m',m'')$, c'est-\`a-dire $S_{m'}(\lambda')+S_{m''}(\lambda'')> S_{m}(\lambda)$. On a $S_{m'}(\lambda')=S_{m'-1}(\lambda')+\lambda'_{m'}$ et $S_{m}(\lambda)=S_{m-1}(\lambda)+\lambda_{m}$. Puisque  $S_{m'-1}(\lambda')+S_{m''}(\lambda'')\leq S_{m-1}(\lambda)$ par hypoth\`ese, on a $\lambda'_{m'}> \lambda_{m}$. A fortiori $\lambda'_{m'+1}=\lambda'_{m'}> \lambda_{m+1}$. On a $S_{m'+1}(\lambda')=S_{m'}(\lambda')+\lambda'_{m'+1}$ et $S_{m+1}(\lambda)=S_{m}(\lambda)+\lambda_{m+1}$. L'in\'egalit\'e $S_{m'}(\lambda')+S_{m''}(\lambda'')> S_{m}(\lambda)$ entra\^{\i}ne alors $S_{m'+1}(\lambda')+S_{m''}(\lambda'')> S_{m+1}(\lambda)$. C'est contraire \`a l'hypoth\`ese. Cette contradiction d\'emontre (9).
 
 Supposons $\nu_{m'}(\lambda')=1$ et $\nu_{m''}(\lambda'')=0$ (une preuve similaire vaut si $\nu_{m'}(\lambda')=0$ et $\nu_{m''}(\lambda'')=1$). Alors les deux couples $(m'-1,m'')$ et $(m'+1,m'')$ v\'erifient $\nu_{m'-1}(\lambda')=\nu_{m''}(\lambda'')=0$ et $\nu_{m'+1}(\lambda')=\nu_{m''}(\lambda'')=0$. Comme on l'a dit ci-dessus, la relation (2) est v\'erifi\'ee pour ces deux couples. L'assertion (9) entra\^{\i}ne qu'elle l'est pour $(m',m'')$. Supposons maintenant que $\nu_{m'}(\lambda')=\nu_{m''}(\lambda'')=1$. Alors les deux couples $(m',m''-1)$ et $(m',m''+1)$ v\'erifient les hypoth\`eses pr\'ec\'edentes:  $\nu_{m'}(\lambda')=1$ et $\nu_{m''-1}(\lambda'')=\nu_{m''+1}(\lambda'')=0$. On vient de prouver que (2) \'etait v\'erifi\'ee pour ces deux couples. L'assertion sym\'etrique \`a (9) o\`u l'on \'echange les composantes $G'_{s}$ et $G''_{s}$ entra\^{\i}ne que (2) est aussi v\'erifi\'ee pour le couple $(m',m'')$. Cela ach\`eve la preuve de (2).

Avant de prouver la deuxi\`eme assertion du (i) de l'\'enonc\'e, examinons ce que devient le calcul ci-dessus dans le cas o\`u $\lambda'\cup \lambda''=\lambda$. Pour tout entier $m\in \{1,...,2n\}$, on peut choisir un couple $(m',m'')\in \{0,...,2n'+1\}\times \{0,...,2n''+1\}$ tel que $m=m'+m''$ et 
$$(10) \qquad S_{m}(\lambda)=S_{m'}(\lambda')+S_{m''}(\lambda'').$$
Ce couple n'est pas unique: par exemple, si $(m',m'')$ convient avec $m'\geq1$  et $m''\leq 2n''$  et que $\lambda'_{m'}=\lambda''_{m''+1}$, le couple $(m'-1,m''+1)$ convient aussi. Parmi les couples possibles, on choisit celui tel que $m'$ soit maximal. Les entiers $m'$ et $m''$  se calculent alors de la fa\c{c}on suivante. Posons $i=\lambda_{m}$. L'ensemble des $m'$ plus grands termes de $\lambda'$ contient tous ceux qui sont strictement sup\'erieurs \`a $i$, ils sont en nombre $mult_{\lambda'}(\geq i+1)$. Il y a en plus $m'-mult_{\lambda'}(\geq i+1)$ termes \'egaux \`a $i$. De m\^eme pour $\lambda''$. De plus, si $m-mult_{\lambda}(\geq i+1)\leq mult_{\lambda'}(i)$, on a $m'-mult_{\lambda'}(\geq i+1)=m-mult_{\lambda}(\geq i+1)$ et $m''=mult_{\lambda''}(\geq i+1)$, tandis que si $m-mult_{\lambda}(\geq i+1)> mult_{\lambda'}(i)$, on a $m'=mult_{\lambda'}(\geq i)$ et $m''-mult_{\lambda''}(\geq i+1)= m-mult_{\lambda}(\geq i+1)-mult_{\lambda'}(i)$.
Posons $E_{0}=\nu_{m'}(\lambda')+\nu_{m''}(\lambda'')-\boldsymbol{\nu}_{m}(\lambda)$. Montrons que

(11)  $E_{0}=0$. 

Posons $i=\lambda_{m}$. 
Supposons d'abord $i$ impair donc  $\boldsymbol{\nu}_{m}(\lambda)=0$. Il r\'esulte de la description ci-dessus que ou bien $\lambda'_{m'}=i$, ou bien $m'\leq 2n'$ et $\lambda'_{m'}>i\geq \lambda'_{m'+1}$, ou bien $m'=2n'+1$. Dans les trois cas, $\nu_{m'}(\lambda')=0$. De m\^eme $\nu_{m''}(\lambda'')=0$, d'o\`u $E_{0}=0$. Supposons maintenant $i$ pair. Alors $\boldsymbol{\nu}_{m}(\lambda)=1$ si $m$ est impair et $\boldsymbol{\nu}_{m}(\lambda)=0$ si $m$ est pair. Supposons $m-mult_{\lambda}(\geq i+1)\leq mult_{\lambda'}(i)$. Alors ou bien $m''\leq 2n''$ et $\lambda''_{m''}>i\geq \lambda''_{m''+1}$, ou  bien $m''=2n''+1$. Dans les deux cas, $\nu_{m''}(\lambda'')=0$. On a $\lambda'_{m'}=i$ qui est pair. Parce que $\lambda'$ est orthogonale, on a 
$$S_{m'}(\lambda') \equiv S_{mult_{\lambda'}(\geq i+1)}(\lambda')\equiv mult_{\lambda'}(\geq i+1)\,\, mod\,\ 2{\mathbb Z}.$$
D'o\`u $m'-S_{m'}(\lambda')\equiv m'-mult_{\lambda'}(\geq i+1)\,\, mod\,\, 2{\mathbb Z}$. On a $m'-mult_{\lambda'}(\geq i+1)= m-mult_{\lambda}(\geq i+1)$ mais $mult_{\lambda}(\geq i+1)$ est pair car $\lambda$ est sp\'eciale.  Puisque $\nu_{m'}(\lambda')$ vaut $1$ si 
$m'-mult_{\lambda'}(\geq i+1)$ est impair et $0$ si $m'-mult_{\lambda'}(\geq i+1)$ est pair, on voit que $\nu_{m'}(\lambda')=\boldsymbol{\nu}_{m}(\lambda)$, d'o\`u $E_{0}=0$. On laisse au lecteur le calcul similaire dans le cas o\`u $m-mult_{\lambda}(\geq i+1)> mult_{\lambda'}(i)$. Cela prouve (11). 
  
L'entier $E$ de la formule (7) devient $E=-\frac{(x'-y'-x''+y'')^2}{2}+\frac{m}{2}-[\frac{m}{2}]$, qui est toujours n\'egatif ou nul. En comparant (7) et (10), on obtient d'abord que $E=0$, autrement dit

(12) $x'-y'-x''+y''=0$ si $m$ est pair; $\vert x'-y'-x''+y''\vert =1$ si $m$ est impair.

Une seconde cons\'equence est que les in\'egalit\'es (4) et (6) que l'on a utilis\'ees sont forc\'ement des \'egalit\'es.

D\'emontrons la seconde assertion du (i) de l'\'enonc\'e. On suppose $\lambda'\cup \lambda''=\lambda$. Soit $\Delta\in Int(\lambda)$. Effectuons les constructions ci-dessus pour chacun des entiers $m_{1}=j_{min}(\Delta)-1$ et $m_{2}=j_{max}(\Delta)$. On en d\'eduit divers nombres que l'on indexe par $1$ et $2$: $m'_{1}$ et $m'_{2}$, $x'_{1}$ et $x'_{2}$ etc... (dans le cas particulier o\`u $m_{1}=0$, tous les nombres index\'es par $1$ sont nuls). On note par un indice $0$ leurs diff\'erences: $m'_{0}=m'_{2}-m'_{1}$ etc... Puisque $m_{1}$ et $m_{2}$ sont pairs, l'assertion (12) appliqu\'ee aux deux entiers $m_{1}$ et $m_{2}$ entra\^{\i}ne que
$$(13) \qquad x_{0}'-y_{0}'-x_{0}''+y_{0}''=0.$$
Par construction, on a 

(14) $x'_{0}+y'_{0}=mult_{\lambda'}(\Delta)$ et $x''_{0}+y''_{0}=mult_{\lambda''}(\Delta)$. 

\noindent D'autre part, la relation (6), qui est maintenant une \'egalit\'e, appliqu\'ee à chacun des  entiers $m_{1}$ et $m_{2}$, entra\^{\i}ne que $x_{0}=x'_{0}+x''_{0}$ est le nombre d'\'el\'ements de $X(\alpha)\cap ({\bf X}(\lambda,\Delta)\cup {\bf Y}(\lambda,\Delta))$ et que $y_{0}=y'_{0}+y''_{0}$ est le nombre d'\'el\'ements de $Y(\beta)\cap ({\bf X}(\lambda,\Delta)\cup {\bf Y}(\lambda,\Delta))$. On peut supposer que $X(\alpha)$ et $Y(\beta)$ se d\'eduisent comme en \ref{DnramSI} d'\'el\'ements $\nu^X$ et $\nu^Y$. On voit que
$$(15)\qquad  x_{0}=\vert X(\alpha)\cap ({\bf X}(\lambda,\Delta)\cup {\bf Y}(\lambda,\Delta))\vert =\left\lbrace\begin{array}{cc}m_{0}/2,&\text{ si }\nu^X_{\Delta}=\nu^Y_{\Delta}\\ m_{0}/2+1,&\text{ si }\nu^X_{\Delta}=1\text{ et }\nu^Y_{\Delta}=-1\\ m_{0}/2-1,&\text{ si }\nu^X_{\Delta}=-1\text{ et }\nu^Y_{\Delta}=1\\ \end{array}\right.$$
$$(16)\qquad  y_{0} =\vert Y(\beta)\cap ({\bf X}(\lambda,\Delta)\cup {\bf Y}(\lambda,\Delta))\vert =\left\lbrace\begin{array}{cc}m_{0}/2,&\text{ si }\nu^X_{\Delta}=\nu^Y_{\Delta}\\ m_{0}/2+1,&\text{ si }\nu^X_{\Delta}=-1\text{ et }\nu^Y_{\Delta}=1\\ m_{0}/2-1,&\text{ si }\nu^X_{\Delta}=1\text{ et }\nu^Y_{\Delta}=-1.\\ \end{array}\right.$$
En se rappelant que $d_{\Delta}=\nu^X_{\Delta}\nu^Y_{\Delta}$, on obtient

(17) $x_{0}=y_{0}$ si $ d_{\Delta}=1$, $\vert x_{0}-y_{0}\vert =2$ si $d_{\Delta}=-1$. 

On a
$$x''_{0}+y''_{0}\equiv x''_{0}-y''_{0}\equiv (x_{0}-y_{0})/2-(x'_{0}-y'_{0}-x''_{0}+y''_{0})/2\,\,mod\,\,2{\mathbb Z}.$$
Il r\'esulte alors  de (13), (14) et (17) que $mult_{\lambda''}(\Delta)$ est pair si $d_{\Delta}=1$ et est impair si $d_{\Delta}=-1$. Ces \'egalités pour tout $\Delta\in Int(\lambda)$ \'equivalent \`a l'\'egalit\'e $\bar{c}_{s}(\lambda',\lambda'')=(\lambda,d)$. Cela ach\`eve la preuve du (i) de l'\'enonc\'e.

Tournons-nous vers le (ii). Soit $(\lambda',\lambda'')\in {\cal P}^{orth}(2n'+1)\times {\cal P}^{orth}(2n''+1)$, supposons $\bar{c}_{s}(\lambda',\lambda'')=(\lambda,d)$.  Pour $\eta^X,\eta^Y\in \{\pm 1\}$, notons $Int(\lambda;\eta^X,\eta^Y)$ l'ensemble des $\Delta\in Int(\lambda)$ tels que $\nu^X_{\Delta}=\eta^X$ et $\nu^Y_{\Delta}=\eta^Y$. L'\'egalite (15)  entra\^{\i}ne que
$\vert  X(\alpha)\vert =n+\vert Int(\lambda;1,-1)\vert -\vert Int(\lambda;-1,1)\vert $. Mais $\vert X(\alpha)\vert =n+h$. Donc 

(18) $h= \vert Int(\lambda;1,-1)\vert - \vert Int(\lambda;-1,1)\vert $.

\noindent A fortiori
$h\leq  \vert Int(\lambda;1,-1)\vert $. Pour $\Delta\in Int(\lambda;1,-1)$, on a $d_{\Delta}=-1$ et l'\'egalit\'e $\bar{c}_{s}(\lambda',\lambda'')=(\lambda,d)$ implique que $mult_{\lambda''}(\Delta)$ et $mult_{\lambda'}(\Delta)$ sont impairs.  Donc $\lambda'$ et $\lambda''$ contiennent toutes deux un \'el\'ement de $\Delta$ pour tout $\Delta\in Int(\lambda;1,-1)$. Ces \'el\'ements sont distincts et impairs et il y en a au moins $h$. Leurs somme est donc minor\'ee par $1+3+...+(2h-1)=h^2$. Cette somme est major\'ee par $S(\lambda')=2n'+1$ et par $S(\lambda'')=2n''+1$. Donc $h^2\leq 2n'+1$ et $h^2\leq 2n''+1$, ce qui signifie que $\Lambda_{h}\subset \Lambda(s)$. C'est la premi\`ere assertion de (ii). 

Consid\'erons des triplets $(h,\alpha',\beta')$ et $(h,\alpha'',\beta'')$ v\'erifiant les conditions de (ii). On pose $Spr(h,\alpha',\beta')=(\lambda',\epsilon')$ et $Spr(h,\alpha'',\beta'')=(\lambda'',\epsilon'')$. Plus pr\'ecis\'ement,  on rel\`eve $\epsilon'$ en  l'\'el\'ement de
$\{\pm 1\}^{Jord_{bp}(\lambda')}$ tel que $A(\alpha')={\bf A}(\lambda',\epsilon')$ et $B(\beta')={\bf B}(\lambda',\epsilon')$.  De m\^eme pour $\epsilon''$.   

On reprend notre calcul en choisissant pour tout $m\in \{0,...,2n\}$ le couple $(m',m'')$ d\'efini plus haut, que l'on note $(m'_{m},m''_{m})$. On en d\'eduit des entiers $x'$ etc... que l'on note $x'_{m}$ etc... Ces entiers sont nuls dans le cas $m=0$. On peut supposer et on suppose que les suites $m\mapsto x'_{m}$ etc... sont croissantes. 
Pour $m\geq1$, on pose ${\bf m}'_{m}=m'_{m}-m'_{m-1}$, ${\bf x}'_{m}=x'_{m}-x'_{m-1}$ etc... 
Pour tout $m\geq1$, le couple $({\bf m}'_{m},{\bf m}''_{m})$ est \'egal \`a $( 1,0)$ ou \`a $( 0,1)$.  Disons que $m$ provient de $\lambda'$ dans le premier cas, de $\lambda''$ dans le second. 
Fixons un intervalle $\Delta\in Int(\lambda)$.   Posons $u'=mult_{\lambda'}(\Delta)$ et $u''=mult_{\lambda''}(\Delta)$.  Notons $k'_{1}< k'_{2}<...<k'_{u'}$ les \'el\'ements de $\{ j_{min}(\Delta),...,j_{max}(\Delta)\}$ qui proviennent de $\lambda'$ et $k''_{1}<k''_{2}<...<k''_{u''}$ ceux qui proviennent de $\lambda''$.  Ajoutons \`a ces suites les termes $k'_{0}=k''_{0}=j_{min}(\Delta)-1$. Pour $r=1,...,u'$, posons $\underline{x}'[r]=x'_{k'_{r}}-x'_{k'_{r-1}}={\bf x}'_{k'_{r}}$, $\underline{y}'[r]=y'_{k'_{r}}-y'_{k'_{r-1}}={\bf y}'_{k'_{r}}$. On a $\underline{x}'[r]+\underline{y}'[r] =1$, donc 
le couple  $(\underline{x}'[r],\underline{y}'[r])$ est \'egal \`a $(1,0)$ ou $(0,1)$. On pose des d\'efinitions analogues pour $\lambda''$, c'est-\`a-dire en rempla\c{c}ant les $'$ par des $''$. On va prouver

(19) quand $r$ d\'ecrit $\{1,...,u'\}$, les couples $(\underline{x}'[r],\underline{y}'[r])$ sont alternativement $(1,0)$ et $(0,1)$; le premier couple $(\underline{x}'[1],\underline{y}'[1])$ vaut $(1,0)$ si $\nu^X_{\Delta}=1$ et $(0,1)$ si $\nu^X_{\Delta}=-1$; une assertion analogue vaut pour les termes relatifs \`a $\lambda''$. 

Supposons pour fixer les id\'ees que $\nu^X_{\Delta}=1$. L'assertion devient $(\underline{x}'[r],\underline{y}'[r])=$

\noindent $(\underline{x}''[r],\underline{y}''[r])=(1,0)$ si $r$ est impair et $(\underline{x}'[r],\underline{y}'[r])=(\underline{x}''[r],\underline{y}''[r])=(0,1)$ si $r$ est pair. 
On raisonne par r\'ecurrence sur $m\in \{ j_{min}(\Delta),...,j_{max}(\Delta)\}$ en supposant que les assertions relatives \`a $\lambda'$, resp. $\lambda''$,  sont v\'erifi\'ees pour tout $r$ tel que $k'_{r}<m$, resp. $k''_{r}<m$. Cette hypoth\`ese est vide pour $m=j_{min}(\Delta)$.   Posons ${\bf X}(\lambda)\cup {\bf Y}(\lambda)=(z_{1},...,z_{2n})$.  Supposons pour fixer les id\'ees que $m=k'_{r}$ et notons $s$ le plus grand \'el\'ement de $\{0,...,u''\}$ tel que $k''_{s}<m$. On a $r+s=m-j_{min}(\Delta)+1$. Supposons d'abord $m$ impair. On a $z_{m}> z_{m+1}$ et $z_{j_{min}(\Delta)-1}> z_{j_{min}(\Delta)}$.   La relation (6), qui est une \'egalit\'e,  d\'etermine enti\`erement  les couples $(x_{m},y_{m})$ et $(x_{j_{min}(\Delta)-1},y_{j_{min}(\Delta)-1})$. L'entier $x_{m}-x_{j_{min}(\Delta)-1}$, resp. $y_{m}-y_{j_{min}(\Delta)-1}$, est le nombre d'\'el\'ements de $X(\alpha)\cap \{z_{j_{min}(\Delta)},...,z_{m}\}$, resp. de $Y(\alpha)\cap  \{z_{j_{min}(\Delta)},...,z_{m}\}$. D'apr\`es la construction et l'hypoth\`ese $\nu^X_{\Delta}=1$, ces nombres sont $(m-j_{min}(\Delta))/2+1$ et $(m-j_{min}(\Delta))/2$. Mais on a aussi
$$x_{m}-x_{j_{min}(\Delta)-1}=(\sum_{t'=1,...,r}\underline{x}'[t'])+(\sum_{t''=1,...,s}\underline{x}''[t'']),$$
$$ y_{m}-y_{j_{min}(\Delta)-1}=(\sum_{t'=1,...,r}\underline{y}'[t'])+(\sum_{t''=1,...,s}\underline{y}''[t'']).$$
L'hypoth\`ese de r\'ecurrence permet de calculer ces sommes: par exemple, 
$\sum_{t'=1,...,r-1}\underline{x}'[t']$ est le nombre de termes impairs dans l'ensemble $\{1,...,r-1\}$. En tenant compte de l'\'egalit\'e $r+s=m-j_{min}(\Delta)+1$ et de l'imparit\'e de $m$ et $j_{min}(\Delta)$, on obtient
$$x_{m}-x_{j_{min}(\Delta)-1}=\underline{x}'[r]+(\sum_{t'=1,...,r-1}\underline{x}'[t'])+(\sum_{t''=1,...,s}\underline{x}''[t''])$$
$$=\left\lbrace\begin{array}{cc}\underline{x}'[r]+(m-j_{min}(\Delta))/2,&\text{ si }r\text{ est impair},\\\underline{x}'[r]+ (m-j_{min}(\Delta))/2+1,&\text{ si }r\text{ est pair}.\\ \end{array}\right.$$
En comparant avec l'\'egalit\'e $x_{m}-x_{j_{min}(\Delta)-1}=(m-j_{min}(\Delta))/2+1$, on obtient $\underline{x}'[r]=1$ si $r$ est impair et $\underline{x}'[r]
=0$ si $r$ est pair.  D'o\`u forc\'ement $\underline{y}'[r]=0$ si $r$ est impair et $\underline{y}'[r]
=1$ si $r$ est pair. Cela prouve notre assertion dans le cas o\`u $m$ est impair.  Supposons $m$ pair. L'assertion (12) appliqu\'ee aux entiers $m$ et $m-2$ entra\^{\i}ne
$${\bf x}'_{m}-{\bf y}'_{m} -{\bf x}''_{m} +{\bf y}''_{m}=
-({\bf x}'_{m-1}-{\bf y}'_{m-1}-{\bf x}''_{m-1}+{\bf y}''_{m-1}).$$
Le membre de gauche est $\underline{x}'[r]-\underline{y}'[r]$. Le membre de droite est $-(\underline{x}'[r-1]-\underline{y}'[r-1])$ si $m-1=k'_{r-1}$, $\underline{x}''[s]-\underline{y}''[s]$ si $m-1=k''_{s}$. Puisque $m$ est pair, $r+s$ l'est aussi et l'hypoth\`ese de r\'ecurrence entra\^{\i}ne que ce membre de droite vaut en tout cas $(-1)^{r-1}$. On a alors forc\'ement $(\underline{x}'[r],\underline{y}'[r])=(1,0)$ si $r$ est impair, $(\underline{x}'[r],\underline{y}'[r])=(0,1)$ si $r$ est pair, ce qui ach\`eve la preuve de (19). 

Posons ${\bf C}(\lambda',\Delta)=\cup_{i\in \Delta}{\bf C}(\lambda',i)$  et notons ${\bf A}(\lambda',\Delta)$, ${\bf B}(\lambda',\Delta)$, $A(\alpha',\Delta)$, $B(\beta',\Delta)$ les intersections de ${\bf A}(\lambda')$ etc... avec ${\bf C}(\lambda',\Delta)$. Notons $mult_{\lambda'}(>\Delta)$ la somme des $mult_{\lambda'}(i)$ pour les entiers $i$ strictement sup\'erieurs \`a tout \'el\'ement de $\Delta$. Autrement dit, le plus grand terme de ${\bf C}(\lambda',\Delta)$ est $c_{mult_{\lambda'}(>\Delta)+1}$ avec les notations utilis\'ees en \ref{groupesclassiques}.  L'ensemble ${\bf C}(\lambda',\Delta)$ est form\'e d'entiers distincts  intervenant alternativement dans ${\bf A}(\lambda',\Delta)$ et ${\bf B}(\lambda',\Delta)$. Le plus grand terme de ${\bf C}(\lambda',\Delta)$ appartient \`a ${\bf A}(\lambda',\Delta)$ si $mult_{\lambda'}(>\Delta)$ est pair, \`a ${\bf B}(\lambda',\Delta)$ si $mult_{\lambda'}(>\Delta)$ est impair. On a $A(\alpha',\Delta)\cup B(\beta',\Delta)={\bf C}(\lambda',\Delta)$ et l'assertion (19) dit que ${\bf C}(\lambda',\Delta)$ est form\'e d'entiers   intervenant alternativement dans $A(\alpha',\Delta)$ et $B(\beta',\Delta)$, le plus grand terme intervenant dans $A(\alpha',\Delta)$ si $\nu^X_{\Delta}=1$, dans $B(\beta',\Delta)$ si $\nu^X_{\Delta}=-1$. En cons\'equence, on a $(A(\alpha',\Delta),B(\beta',\Delta))=({\bf A}(\lambda',\Delta), {\bf B}(\lambda',\Delta))$ si $\nu^X_{\Delta}=(-1)^{mult_{\lambda'}(>\Delta)}$ et $(A(\alpha',\Delta),B(\beta',\Delta))=({\bf B}(\lambda',\Delta), {\bf A}(\lambda',\Delta))$ si $\nu^X_{\Delta}=(-1)^{1+mult_{\lambda'}(>\Delta)}$.  En se rappelant la construction de $\epsilon'$, cela signifie que $\epsilon'_{i}=\nu^X_{\Delta}(-1)^{mult_{\lambda'}(>\Delta)}$ pour tout $i\in \Delta\cap Jord_{bp}(\lambda')$. 
Par hypoth\`ese, on a $\bar{c}_{s}(\lambda',\lambda'')=(\lambda,d)$. Donc, pour tout $\Delta'\in Int(\lambda)$, on a $(-1)^{mult_{\lambda'}(\Delta')}=(-1)^{mult_{\lambda''}(\Delta')}=d_{\Delta'}=\nu^X_{\Delta'}\nu^Y_{\Delta'}$. On a alors
$$\nu^X_{\Delta}(-1)^{mult_{\lambda'}(>\Delta)}=(\prod_{\Delta'\geq \Delta}\nu^X_{\Delta'})(\prod_{\Delta'> \Delta}\nu^Y_{\Delta'})=\mu_{\Delta}.$$
 Cela d\'emontre que $\epsilon'_{i}=\mu_{\Delta}$ pour tout $i\in \Delta\cap Jord_{bp}(\lambda')$. On prouve de m\^eme que $\epsilon''_{i}=\mu_{\Delta}$ pour tout $i\in \Delta\cap Jord_{bp}(\lambda'')$. C'est la derni\`ere assertion du (ii) de l'\'enonc\'e. Cela entra\^{\i}ne qu'il y a au plus un couple de triplets $(h,\alpha',\beta')$ et $(h,\alpha'',\beta'')$ v\'erifiant les conditions impos\'ees. 

 On doit maintenant prouver l'existence de ce couple. On d\'efinit des \'el\'ements $\epsilon'\in \{\pm 1\}^{Jord_{bp}(\lambda')}$ et $\epsilon''\in \{\pm 1\}^{Jord_{bp}(\lambda'')}$ par les formules de l'\'enonc\'e. Il existe d'uniques \'el\'ements $(h',\alpha',\beta')\in {\cal I}^{FC}(G'_{s})$ et $(h'',\alpha'',\beta'')\in {\cal I}^{FC}(G''_{s})$ de sorte que $Spr(h',\alpha',\beta')=(\lambda',\epsilon')$ et $Spr(h'',\alpha'',\beta'')=(\lambda'',\epsilon'')$. Prouvons que
 
 (20) $h'=h''=h$. 
 
 En adaptant les notations de  \ref{groupesclassiques}, on note   $J'$ l'ensemble des $i\in Jord_{bp}(\lambda')$ tels que  $mult_{\lambda'}(i)$  est impair et $\epsilon'_{i}=-1$. On pose  
 $$D'=\sum_{i\in J'}(-1)^{mult_{\lambda'}(\geq i)}.$$
 On peut supprimer de la d\'efinition de $J'$ la condition $\epsilon'_{i}=-1$ en introduisant dans la somme d\'efinissant $D'$ le terme $\frac{1-\epsilon'_{i}}{2}$. On peut supprimer la condition $mult_{\lambda'}(i)$ impair en rempla\c{c}ant le terme $(-1)^{mult_{\lambda'}(\geq i)}$ par $\frac{(-1)^{mult_{\lambda'}(\geq i)}-(-1)^{mult_{\lambda'}(> i)}}{2}$. Alors
 $$D'=\frac{1}{4}\sum_{i\in Jord_{bp}(\lambda')}(1-\epsilon'_{i})((-1)^{mult_{\lambda'}(\geq i)}-(-1)^{mult_{\lambda'}(> i)}).$$
 Puisque $\epsilon'_{i}$ est constant sur les intervalles de $\lambda'$, on peut regrouper les $i$ intervenant dans chaque intervalle. En utilisant la d\'efinition de $\epsilon'$, on obtient
 $$D'=\frac{1}{4}\sum_{\Delta\in Int(\lambda')}(1-\mu_{\Delta})((-1)^{mult_{\lambda'}(\geq \Delta)}-(-1)^{mult_{\lambda'}(> \Delta)}),$$
 o\`u $mult_{\lambda'}(\geq \Delta)$ est le somme des $mult_{\lambda'}(i)$ sur les entiers $i$ sup\'erieurs ou \'egaux au plus petit terme de $\Delta$. 
 Puisque $\bar{c}_{s}(\lambda',\lambda'')=(\lambda,d)$, on a
 $$(-1)^{mult_{\lambda'}(\geq \Delta)}=\prod_{\Delta'\geq \Delta}\nu^X_{\Delta'}\nu^Y_{\Delta'},\,\, (-1)^{mult_{\lambda'}(> \Delta)}=\prod_{\Delta'> \Delta}\nu^X_{\Delta'}\nu^Y_{\Delta'}.$$
 En se rappelant la d\'efinition de $\mu_{\Delta}$, on calcule
 $$\mu_{\Delta}((-1)^{mult_{\lambda'}(\geq \Delta)}-(-1)^{mult_{\lambda'}(> \Delta)})=\nu^Y_{\Delta}-\nu^X_{\Delta},$$
 puis
 $$D'=\frac{1}{4}(\sum_{\Delta\in Int(\lambda')}(-1)^{mult_{\lambda'}(\geq \Delta)}-(-1)^{mult_{\lambda'}(> \Delta)})+\frac{1}{4}(\sum_{\Delta\in Int(\lambda')}\nu^X_{\Delta}-\nu^Y_{\Delta}).$$
 La premi\`ere somme entre parenth\`eses vaut $(-1)^{mult_{\lambda'}(\geq \Delta_{min})}-1$, o\`u $\Delta_{min}$ est le plus petit intervalle. Puisque $\lambda'\in {\cal P}^{orth}(2n'+1)$, le nombre $mult_{\lambda'}(\geq\Delta_{min})$ est impair et cette premi\`ere somme  vaut $-2$. Avec les notations de (18), le sommand de la deuxi\`eme somme entre parenth\`eses vaut $2$ si $\Delta\in Int(\lambda;1,-1)$, $-2$ si $\Delta\in Int(\lambda;-1,1)$ et $0$ dans les autres cas. Cette deuxi\`eme somme vaut donc $2(\vert Int(\lambda;1,-1)\vert -\vert  Int(\lambda;-1,1)\vert )$. D'apr\`es (18), c'est $2h$. D'o\`u
 $D'=\frac{1}{2}(h-1)$. D'apr\`es \ref{groupesclassiques}(4), on a $h'=\vert 2D'+1\vert $, donc  $h'=h$. On d\'emontre de m\^eme que $h''=h$, ce qui prouve (20). 
 
 On applique, du moins partiellement, nos constructions  aux triplets $(h,\alpha',\beta')$ et $(h,\alpha'',\beta'')$. Pour $m\in \{0,...,2n\}$, on choisit selon la recette  expliqu\'ee ci-dessus des entiers $m'_{m}\in \{0,...,2n'+1\}$ et $m''_{m}\in \{0,...,2n''+1\}$ de sorte que $m=m'_{m}+m''_{m}$ et $S_{m}(\lambda)=S_{m'_{m}}(\lambda')+S_{m''_{m}}(\lambda'')$. A l'aide de ces entiers $m'_{m}$ et $m''_{m}$, on choisit des entiers $x'_{m}$, $y'_{m}$, $x''_{m}$ et $y''_{m}$ comme au d\'ebut de la preuve et on obtient l'\'egalit\'e (4) (o\`u les termes sont maintenant affect\'es d'indices $m$). On pose $x_{m}=x'_{m}+x''_{m}$, $y_{m}=y'_{m}+y''_{m}$. Pour $m\geq1$, on pose comme plus haut ${\bf m}'_{m}=m'_{m}-m'_{m-1}$  etc... 
 Prouvons que
 
 (21) on a l'\'egalit\'e $S_{m}({\bf X}(\lambda)\cup {\bf Y}(\lambda))=S_{x_{m}}(X(\alpha))+S_{y_{m}}(Y(\beta))$;
 
 (22) on a l'\'egalit\'e $x'_{m}-y'_{m}-x''_{m}+y''_{m}=0$ si $m$ est pair et $\vert x'_{m}-y'_{m}-x''_{m}+y''_{m}\vert =1$ si $m$ est impair.
 
 On prouve ces relations par r\'ecurrence sur $m$. Remarquons d'abord que, pour $m$ impair, la relation (22) se d\'eduit de la m\^eme relation pour $m-1$. En effet, en tenant compte de cette relation pour $m-1$,  on a $x'_{m}-y'_{m}-x''_{m}+y''_{m}={\bf x}'_{m}-{\bf y}'_{m}-{\bf x}''_{m}+{\bf y}''_{m}$. Par construction, un et un seul des entiers ${\bf x}'_{m}$ etc.... vaut $1$, les autres valant $0$. D'o\`u l'\'egalit\'e (22) pour $m$. 
 
 Supposons d'abord que $\lambda_{m}$ est pair et que $m$ est impair.    Par r\'ecurrence, on a  
 
 \noindent $S_{m-1}({\bf X}(\lambda)\cup {\bf Y}(\lambda))=S_{x_{m-1}}(X(\alpha))+S_{y_{m-1}}(Y(\beta))$. On  note encore ${\bf X}(\lambda)\cup {\bf Y}(\lambda)=\{z_{1},...,z_{2n}\}$. Les hypoth\`eses sur $m$ entra\^{\i}nent $z_{m-1}> z_{m}$ et l'\'egalit\'e pr\'ec\'edente d\'etermine enti\`erement $x_{m-1}$ et $y_{m-1}$. De plus, le $x_{m-1}+1$-i\`eme terme de $X(\alpha)$ est \'egal au $y_{m-1}+1$-i\`eme terme de $Y(\beta)$, ces deux termes \'etant \'egaux \`a $z_{m}=z_{m+1}$. En cons\'equence, on a les \'egalit\'es
 $$S_{m}({\bf X}(\lambda)\cup {\bf Y}(\lambda))=S_{x_{m-1}+1}(X(\alpha))+S_{y_{m-1}}(Y(\beta))=S_{x_{m-1}}(X(\alpha))+S_{y_{m-1}+1}(Y(\beta)).$$
 Puisque $(x_{m},y_{m})$ est \'egal par construction soit \`a $(x_{m-1}+1,y_{m-1})$, soit \`a $(x_{m-1},y_{m-1}+1)$, on obtient dans les deux cas l'\'egalit\'e (21) pour $m$.  Supposons maintenant que $\lambda_{m}$ est pair et que $m$ est pair. Supposons pour fixer les id\'ees que $m$ provient de $\lambda'$, c'est-\`a-dire ${\bf m}'_{m}=1$ et ${\bf m}''_{m}=0$. Posons $i=\lambda_{m}$. On a note ${\bf C}(\lambda',i)=(c_{j_{min}(i)},...,c_{j_{max}(i)})$ et $m'_{m}$ est un \'el\'ement de $\{j_{min}(i),...,j_{max}(i)\}$. 
 Il  r\'esulte de la  construction des suites $(m'_{m})_{m\in \{0,...,2n\}}$ et $(m''_{m})_{m\in \{0,...,2n\}}$ que $m'_{m}-j_{min}(i)+1=m-mult_{\lambda}(\geq i+1)$. L'entier $m$ est pair par hypoth\`ese et $mult_{\lambda}(\geq i+1)$ l'est parce que $\lambda$ est sp\'eciale. Donc 
 l'entier $m'_{m}-j_{min}(i)+1$ est pair. A  fortiori $m'_{m}> j_{min}(i)$.  On a alors $m'_{m-1}=m'_{m}-1\geq j_{min}(i)$. D'autre part, puisque $m'_{m}-j_{min}(i)+1$ et $i$ sont pairs, on a  $c_{m'_{m}-1}=c_{m'_{m}}$, c'est-\`a-dire   $c_{m'_{m-1}}=c_{m'_{m}}$.   Ces deux termes appartiennent l'un \`a $A(\alpha')$, l'autre \`a $B(\beta')$. Cela entra\^{\i}ne que les deux couples $({\bf x}'_{m-1},{\bf y}'_{m-1})$ et $({\bf x}'_{m},{\bf y}'_{m})$ sont distincts, l'un valant $(1,0)$, l'autre $(0,1)$. Donc $(x_{m},y_{m})=(x_{m-2}+1,y_{m-2}+1)$. Alors l'\'egalit\'e (21) pour $m$ se d\'eduit comme ci-dessus de la m\^eme \'egalit\'e pour $m-2$. En utilisant (22) pour $m-2$, on a aussi
 $$x'_{m}-y'_{m}-x''_{m}+y''_{m}={\bf x}'_{m-1}+{\bf x}'_{m}-{\bf y}'_{m-1}-{\bf y}'_{m}$$
 et ceci vaut $0$ puisque les couples $({\bf x}'_{m-1},{\bf y}'_{m-1})$ et $({\bf x}'_{m},{\bf y}'_{m})$ sont distincts. D'o\`u (22) pour $m$. 
 
 Supposons maintenant que $\lambda_{m}$ est impair. Notons $\Delta$  l'intervalle de $\lambda$ tel que $\lambda_{m}\in \Delta$. Introduisons les m\^emes notations qu'avant l'assertion (19): les suites $k'_{1}<...<k'_{u'}$ et $k''_{1}<...<k''_{u''}$ et, pour $r\in \{1,...,u'\}$, resp. $r\in \{1,...,u''\}$ les \'el\'ements $\underline{x}'[r]=x'_{k'_{r}}-x'_{k'_{r-1}}={\bf x}'_{k'_{r}}$ etc...  
 Montrons que
 
 (23) l'assertion (19) reste v\'erifi\'ee. 
 
 En effet, parce que l'application $i\mapsto \epsilon'_{i}$ est
constante sur $Jord_{bp}(\lambda')\cap \Delta$,  le couple $(A(\alpha',\Delta), B(\beta',\Delta))$ est \'egal \`a $({\bf A}(\lambda',\Delta),{\bf B}(\lambda',\Delta))$ ou a $({\bf B}(\lambda',\Delta),{\bf A}(\lambda',\Delta))$. Les entiers de ${\bf C}(\lambda',\Delta)$ interviennent alternativement dans ${\bf A}(\lambda',\Delta)$ et ${\bf B}(\lambda',\Delta)$, donc interviennent aussi alternativement dans $A(\alpha',\Delta)$ et $B(\beta',\Delta)$. Cela \'equivaut \`a la premi\`ere assertion de (19). Dans le paragraphe qui suit (19), on a calcul\'e la valeur constante de $\epsilon'_{i}$ pour $i\in Jord_{bp}(\lambda')\cap \Delta$ en utilisant la valeur du premier couple $(\underline{x}'[1],\underline{y}'[1])$. Dans la situation pr\'esente o\`u l'on conna\^{\i}t cette valeur constante de $\epsilon'_{i}$, le m\^eme calcul d\'etermine ce premier couple. C'est la deuxi\`eme assertion de (19). Cela prouve (23). 

Pour fixer les id\'ees, on suppose que $\nu^X_{\Delta}=1$ et que $m$ provient de $\lambda'$. On note $r$ l'\'el\'ement de $\{1,...,u'\}$ tel que $m=k'_{r}$ et $s$ le plus grand \'el\'ement de $\{0,...,u''\}$ tel que $k''_{s}<m$. 
Supposons d'abord  $m=j_{min}(\Delta)$. Donc $m$ est impair et  $r=1$.  Comme on l'a dit, la relation (22) pour $m$ se d\'eduit   de la m\^eme relation pour $m-1$. On a $z_{m-1}> z_{m}> z_{m+1}$. La relation (21) pour $m-1$ d\'etermine enti\`erement le couple $(x_{m-1},y_{m-1})$ et, compte tenu de la construction de $X(\alpha)$ et $Y(\beta)$ et de notre hypoth\`ese $\nu^X_{\Delta}=1$, la relation (21) pour $m$ (que l'on doit d\'emontrer) \'equivaut \`a l'assertion $x_{m}=x_{m-1}+1$ et $y_{m}=y_{m-1}$. Compte tenu de notre hypoth\`ese $m=k'_{1}$, cela \'equivaut \`a $(\underline{x}'[1],\underline{y}'[1])=(1,0)$. Mais cela r\'esulte de (19) et cela d\'emontre (21) pour $m$. Supposons maintenant que $m$ est impair et $m> j_{min}(\Delta)$. On a d\'ej\`a v\'erifi\'e (22) pour $m$. On a $z_{m-1}=z_{m}$. Compte tenu de la relation (21) pour $m-2$, cette relation pour $m$ \'equivaut \`a l'\'egalit\'e $(x_{m},y_{m})=(x_{m-2}+1,y_{m-2}+1)$, ou encore \`a ${\bf x}_{m}+{\bf x}_{m-1}={\bf y}_{m}+{\bf y}_{m-1}$, ou encore \`a ${\bf x}_{m}-{\bf y}_{m}={\bf y}_{m-1}-{\bf x}_{m-1}$. Le membre de gauche vaut $\underline{x}'[r]-\underline{y}'[r]=(-1)^{r-1}$ d'apr\`es (19) et notre hypoth\`ese $\nu^X_{\Delta}=1$. Si $m-1=k'_{r-1}$, le membre de droite vaut de m\^eme $-(-1)^{r-2}$, d'o\`u l'\'egalit\'e cherch\'ee. Si $m-1=k''_{s}$, ce membre de gauche vaut $-(-1)^{s-1}$. On a $r+s=m-j_{min}(\Delta)+1$. C'est un nombre  impair et on obtient encore l'\'egalit\'e cherch\'ee.
Supposons maintenant que $m$ est pair. Compte tenu des assertions (22) pour $m-2$ et $m-1$, cette assertion pour $m$ \'equivaut \`a
$${\bf x}'_{m}-{\bf y}'_{m}-{\bf x}''_{m}+{\bf y}''_{m}=-({\bf x}'_{m-1}-{\bf y}'_{m-1}-{\bf x}''_{m-1}+{\bf y}''_{m-1}).$$
Le membre de gauche est $\underline{x}'[r]-\underline{y}'[r]$. D'apr\`es (19) et notre hypoth\`ese $\nu^X_{\Delta}=1$, ceci vaut $(-1)^{r-1}$. Si $m-1=k'_{r-1}$, le membre de droite vaut de m\^eme $-(-1)^{r-2}$, d'o\`u l'\'egalit\'e requise. Si $m=k''_{s}$, ce membre de droite vaut $(-1)^{s-1}$. Mais $r+s=m-j_{min}(\Delta)+1$  est pair donc $(-1)^{s-1}=(-1)^{r-1}$, d'o\`u de nouveau l'\'egalit\'e requise.
  Pour prouver (21), supposons d'abord que  $m\leq j_{max}(\Delta)-1$, donc $m\leq j_{max}(\Delta)-2$ puisque ces deux entiers sont pairs.  On a alors $z_{m-1}<z_{m}=z_{m+1}$. L'assertion (21) pour $m$ se d\'eduit alors de la m\^eme assertion pour $m-1$ comme dans le cas ci-dessus o\`u $\lambda_{m}$ est pair et $m$ est impair: l'assertion (21) est v\'erifi\'ee que le couple $(x_{m},y_{m})$ soit $(x_{m-1}+1,y_{m-1})$ ou $(x_{m-1},y_{m-1}+1)$.    Supposons enfin $m=j_{max}(\Delta)$, donc $r=u'$. On a $z_{m-1}< z_{m}< z_{m+1}$. La relation (21) pour $m-1$ d\'etermine enti\`erement le couple $(x_{m-1},y_{m-1})$.   Compte tenu de la construction de $X(\alpha)$ et $Y(\beta)$, cette relation pour $m$ \'equivaut aux \'egalit\'es $({\bf x}_{m},{\bf y}_{m})=(1,0)$ si $\nu^Y_{\Delta}=-1$, $({\bf x}_{m},{\bf y}_{m})=(0,1)$ si $\nu^Y_{\Delta}=1$. Ou encore \`a l'\'egalit\'e ${\bf x}_{m}-{\bf y}_{m}=-\nu^Y_{\Delta}$. On a ${\bf x}_{m}-{\bf y}_{m}=\underline{x}'[u']-\underline{y}'[u']=(-1)^{u'-1}$. Compte tenu de notre hypoth\`ese $\nu^X_{\Delta}=1$, on a $-\nu^Y_{\Delta}=-d_{\Delta}$. Mais $u'=mult_{\lambda'}(\Delta)$ et l'hypoth\`ese $\bar{c}'_{s}(\lambda',\lambda'')=(\lambda,d)$ entra\^{\i}ne que $(-1)^{mult_{\lambda'}(\Delta)}=d_{\Delta}$. D'o\`u l'\'egalit\'e cherch\'ee, qui d\'emontre (21) pour $m$. Cela ach\`eve la preuve de (21) et (22). 
  
 En utilisant (21), les m\^emes calculs qu'au d\'ebut de la preuve conduisent \`a l'\'egalit\'e
 $$S_{m}(\lambda)=2S_{x_{m}}(\alpha)+2S_{y_{m}}(\beta)+e$$
 avec la m\^eme valeur de $e$. Posons $T_{m}=S_{x_{m}}(\alpha)-S_{x'_{m}}(\alpha')-S_{x''_{m}}(\alpha'')$, $U_{m}=S_{y_{m}}(\beta)-S_{y'_{m}}(\beta')-S_{y''_{m}}(\beta'')$. Jointe \`a (3), l'\'egalit\'e pr\'ec\'edente conduit \`a
 $$(24) \qquad S_{m}(\lambda)=2T_{m}+2U_{m}+S_{m'_{m}}(\lambda')+S_{m''_{m}}(\lambda'')-E$$
 o\`u $E$ est comme en (8). L'assertion (11) reste v\'erifi\'ee: elle r\'esulte seulement de notre choix des suites $(m'_{m})_{m=1,...,2n}$ et $(m''_{m})_{m=1,...,2n}$. L'assertion (22) entra\^{\i}ne que la somme des trois premiers termes de (8) est nulle. Donc $E=0$. Mais on a par hypoth\`ese l'\'egalit\'e $S_{m}(\lambda)=S_{m'_{m}}(\lambda')+S_{m''_{m}}(\lambda'')$. Alors (24) entra\^{\i}ne
 $T_{m}+U_{m}=0$. On a $(x_{m},y_{m})=(x_{m-1}+1,y_{m-1})$ ou $(x_{m-1},y_{m-1}+1)$. D'o\`u $U_{m}=U_{m-1}$ dans le premier cas et $T_{m}=T_{m-1}$ dans le second. L'\'egalit\'e $T_{m}+U_{m}=0$ entra\^{\i}ne alors par r\'ecurrence les deux \'egalites $T_{m}=U_{m}=0$. Cela prouve les \'egalit\'es
 $$S_{x_{m}}(\alpha)=S_{x'_{m}}(\alpha')+S_{x''_{m}}(\alpha'')\,\, S_{y_{m}}(\beta)=S_{y'_{m}}(\beta')+S_{y''_{m}}(\beta'')$$
 pour tout $m\in \{1,...,2n\}$.  Par construction, quand $m$ d\'ecrit $\{0,...,2n\}$, les entiers $x'_{m}$, resp. $x''_{m}$ d\'ecrivent enti\`erement les ensembles $\{0,...,n' +(1+h)/2\}$, resp. $\{0,...,n''+(1+h)/2\}$. Donc $x_{m}$ d\'ecrit enti\`erement l'ensemble $\{0,...,n+h\}$. 
 Notons $\ell_{1}<... <\ell_{n+h}$ les entiers $m$ tels que ${\bf x}_{m}=1$.  Pour $r=1,...,n+h$, on a $x_{\ell_{r}}=r$. Posons $x'[r]=x'_{\ell_{r}}$, $x''[r]=x''_{\ell_{r}}$. Alors $x'[r]+x''[r]=r$ et $S_{r}(\alpha)=S_{x'[r]}(\alpha')+S_{x''[r]}(\alpha'')$. Ces relations entra\^{\i}nent par r\'ecurrence sur $r$ que $\alpha'$ est la partition extraite de $\alpha$ form\'ee des $\alpha_{r}$ pour les $r$ tels que $x'[r]=x'[r-1]+1$ tandis que $\alpha''$ est la partition extraite de $\alpha$ form\'ee des $\alpha_{r}$ pour les $r$ tels que $x''[r]=x''[r-1]+1$. En cons\'equence $\alpha'\cup \alpha''=\alpha$. On prouve de m\^eme que $\beta'\cup \beta''=\beta$. Mais alors on sait que la multiplicit\'e de $\rho_{\alpha',\beta'}\otimes \rho_{\alpha'',\beta''}$ dans la restriction de $\rho_{\alpha,\beta}$ vaut $1$, cf. \cite{W3} 1.13(38). Cela prouve que nos triplets $(h,\alpha',\beta')$ et $(h,\alpha'',\beta'')$ v\'erifient les conditions du (ii) de l'\'enonc\'e et cela prouve en m\^eme temps l'assertion concernant la multiplicit\'e de  $\rho_{\alpha',\beta'}\otimes \rho_{\alpha'',\beta''}$ dans la restriction de $\rho_{\alpha,\beta}$. Cela ach\`eve la preuve de la proposition. $\square$

D\'eduisons les th\'eor\`emes \ref{premiertheoremeSI} et \ref{deuxiemetheoremeSI}  de la proposition. Le th\'eor\`eme \ref{premiertheoremeSI} est une cons\'equence imm\'ediate du (i) de l'\'enonc\'e ci-dessus. Fixons une partition  sp\'eciale $\lambda\in {\cal P}^{orth}(2n)$ dont au moins un terme est impair et fixons $d\in \tilde{\bar{A}}(\lambda)/conj=\{\pm 1\}^{Int(\lambda)}_{-1}$. Fixons un sommet $s\in S(\bar{C})$ param\'etr\'e par un couple $(n',n'')$ et fixons $(\lambda',\lambda'')\in {\cal P}^{orth}(2n'+1)\times {\cal P}^{orth}(2n''+1)$ de sorte que  $\bar{c}_{s}(\lambda',\lambda'')=(\lambda,d)$. L'existence d'un tel triplet $(s,\lambda',\lambda'')$ r\'esulte de \ref{surjectivite} et se v\'erifie d'ailleurs par un calcul simple. En identifiant $(\lambda',\lambda'')$ \`a un \'el\'ement de $\mathfrak{g}_{s,nil}/conj$, on pose $b_{F}(\lambda,d)=(s,\lambda',\lambda'')$. Cela d\'efinit l'application $b_{F}$ du th\'eor\`eme \ref{deuxiemetheoremeSI}. L'assertion (i) de ce th\'eor\`eme est imm\'ediate et (iii) r\'esulte de la derni\`ere assertion du (i) de la proposition. Les assertions du (ii) du th\'eor\`eme \ref{deuxiemetheoremeSI} se d\'eduisent de celles du (ii) de la proposition. En particulier, l'injectivit\'e de l'application $(\Lambda,\varphi,\rho)\mapsto {\cal L}_{\Lambda,\varphi,\rho}$ est une traduction du fait que les caract\`eres $\epsilon'$ et $\epsilon''$ du (ii) du th\'eor\`eme ci-dessus se d\'eduisent injectivement du terme $\mu$.

\subsubsection{Preuve des th\'eor\`emes \ref{premiertheoremeSI} et \ref{deuxiemetheoremeSI} dans le cas $(A_{n-1},ram)$ \label{preuveAn-1ram}}

On suppose que $G$ est de type $(A_{n-1},ram)$ (ce qui suppose $n\geq3$). En \ref{An-1ramSI}, on a identifi\'e ${\cal B}_{Irr}(G)$ \`a l'ensemble des triplets $(v,\alpha,\beta)$ o\`u $v\in {\mathbb N}$, $v(v+1)/2\leq n$, $v(v+1)/2$ est de m\^eme parit\'e que $n$ et $(\alpha,\beta)\in {\cal P}_{2}(n/2-v(v+1)/4)$. L'entier $v$ d\'etermine un Levi de $G$ auquel est associ\'e un ensemble $\Lambda_{v}\subset \Delta_{a}^{nr}$ et $(\alpha,\beta)$ d\'etermine une repr\'esentation irr\'eductible $\rho_{\alpha,\beta}$ de $W^{I_{F}}(M_{\Lambda_{v}})$, qui est un groupe de Weyl de type $B_{n/2-v(v+1)/4}$. Pour un tel triplet, on a explicit\'e $\boldsymbol{\nabla}_{F}(v,\alpha,\beta)=\lambda\in {\cal P}(n)$. Introduisons les notations suivantes. Pour $\lambda\in {\cal P}(n)$ et $i\in Jord(\lambda)$, on pose
$$mult_{\lambda}(>i;imp)=\sum_{j>i; j \,\,impair}mult_{\lambda}(j),\, mult_{\lambda}(>i;pair)=\sum_{j>i; j \,\,pair}mult_{\lambda}(j).$$

L'ensemble $S(\bar{C})$ s'envoie surjectivement sur l'ensemble des couples $(N',N'')\in {\mathbb N}^2$ tels que $N'+N''=n$, $N'$ est pair et $N''\not=2$. Cette application est bijective si $n$ est impair. Si $n$ est pair, la fibre au-dessus d'un \'el\'ement $(N',N'')$ n'a qu'un \'el\'ement sauf dans le cas o\`u  $n$ est pair et $(N',N'')=(n,0)$ o\`u elle en a deux. Ces deux sommets sont conjugu\'es par l'action de l'\'el\'ement non trivial de $\boldsymbol{\Omega}^{nr}$. Pour un sommet $s$ de couple associ\'e $(N',N'')$, posons $G'_{s}=Sp(N')$ et $G''_{s}=SO(N'')$, plus exactement $G''_{s}$ est le groupe $SO(N'')$ d\'eploy\'e. On a $G_{s}=G'_{s}\times G''_{s}$ si $n$ est impair. Si $n$ est pair et $N'N''\not=0$, on a $G_{s}=(G'_{s}\times G''_{s})/\{\pm 1\}$, o\`u $\{\pm 1\}$ est identifi\'e au centre des deux composantes. Si $n$ est pair et $N'=0$, resp. $N''=0$, on a $G_{s}=G''_{s}/\{\pm 1\}$, resp. $G_{s}=G'_{s}/\{\pm 1\}$. On a d\'ecrit en \ref{groupesclassiques}  les ensembles ${\cal I}(G'_{s})$, ${\cal I}^{FC}(G'_{s})$ et leurs analogues pour $G''_{s}$. On  ajoute des indices $s$ et des $'$ et $''$ aux notations de \ref{groupesclassiques} et on adopte la m\^eme convention que dans ce paragraphe   quant au probl\`eme de d\'edoublement de certains termes dans le cas o\`u $N''$ est pair. 
  Si $n$ est impair, on a ${\cal I}^{FC}(G_{s})={\cal I}^{FC}(G'_{s})\times {\cal I}^{FC}(G''_{s})$. Si $n$ est pair, un couple $((k,\alpha',\beta'),(h,\alpha'',\beta''))\in {\cal I}^{FC}(G'_{s})\times {\cal I}^{FC}(G''_{s})$ se quotiente en un \'el\'ement de ${\cal I}^{FC}(G_{s})$ si et seulement si $h+k(k+1)$ est divisible par $4$. 

L'ensemble $\mathfrak{g}_{s,nil}/conj
$ s'envoie surjectivement sur ${\cal P}^{symp}(N')\times {\cal P}^{orth}(N'')$, les nombres d'\'el\'ements des fibres de cette application se d\'eduisant de ce que l'on a dit  en \ref{groupesclassiques}. L'application $\bar{c}_{s}:\mathfrak{g}_{s,nil}/conj\to \bar{{\cal C}}={\cal P}(n)$ envoie $(\lambda',\lambda'')$ sur $\lambda'\cup \lambda''$.

Pour un entier $v\in {\mathbb N}$ tel que $v(v+1)/2\leq n$ et $v(v+1)/2$ est de m\^eme parit\'e que $n$, on pose $k=[\frac{v}{2}]$, $h=[\frac{v+1}{2}]$. On a l'\'egalit\'e $h^2+k(k+1)=v(v+1)/2$. On a  $\Lambda_{v}\subset \Lambda(s)$ si et seulement si $k(k+1)\leq N'$ et $h^2\leq N''$. Si ces conditions sont v\'erifi\'ees, $M_{\Lambda_{v},s_{\Lambda_{v}}}$ s'identifie \`a l'image dans $G_{s}$ de $M'_{s,k}\times M''_{s,h}$ et le groupe $W'(M'_{s,k})\times W''(M''_{s,h})$ s'identifie naturellement \`a un sous-groupe de $W^{I_{F}}(M_{\Lambda_{v}})$. Remarquons que, si $n$ est pair, $v(v+1)/2$ l'est aussi et un calcul simple montre que  le couple $(k,h)$ v\'erifie la condition \'evoqu\'ee plus haut: $h+k(k+1)$ est divisible par $4$.

\begin{prop}{Soient $(v,\alpha,\beta)\in {\cal B}_{Irr}(G)$ et $s\in S(\bar{C})$. On note $(N',N'')$ le couple param\'etrant $s$,  on pose $\boldsymbol{\nabla}_{F}(v,\alpha,\beta)=\lambda$, $k=[\frac{v}{2}]$ et  $h=[\frac{v+1}{2}]$.

(i) Supposons $\Lambda_{v}\subset \Lambda(s)$, soient $(k,\alpha',\beta')\in {\cal I}^{FC}(G'_{s})$ et $(h,\alpha'',\beta'')\in {\cal I}^{FC}(G''_{s})$, notons $(\lambda',\epsilon')$ et $(\lambda'',\epsilon'')$ leurs images par la correspondance de Springer g\'en\'eralis\'ee. Supposons que $\rho_{\alpha',\beta'}\otimes \rho_{\alpha'',\beta''}$ soit une composante irr\'eductible de la restriction de $\rho_{\alpha,\beta}$ \`a $W'(M'_{s,k})\times W''(M''_{s,h})$. Alors $\lambda'\cup \lambda''\leq \lambda$.

(ii) Soient $(\lambda',\lambda'')\in {\cal P}^{symp}(N')\times {\cal P}^{orth}(N'')$, supposons $\lambda'\cup \lambda''=\lambda$. Alors $\Lambda_{v}\subset \Lambda(s)$. Il existe d'uniques triplets $(k,\alpha',\beta')\in {\cal I}^{FC}(G'_{s})$ et $(h,\alpha'',\beta'')\in {\cal I}^{FC}(G''_{s})$ v\'erifiant les conditions suivantes:

(a) leurs images par la  correspondance de Springer g\'en\'eralis\'ee sont de la forme $(\lambda',\epsilon')$ et $(\lambda'',\epsilon'')$;

(b) $\rho_{\alpha',\beta'}\otimes \rho_{\alpha'',\beta''}$ est une composante irr\'eductible de la restriction de $\rho_{\alpha,\beta}$ \`a $W'(M'_{s,k})\times W''(M''_{s,h})$.

Pour ces triplets, $\rho_{\alpha',\beta'}\otimes \rho_{\alpha'',\beta''}$ intervient avec multiplicit\'e $1$ dans la restriction de $\rho_{\alpha,\beta}$. Pour $i\in Jord_{bp}(\lambda')$, on a $\epsilon'_{i}=(-1)^{n+mult_{\lambda}(>i;imp)}$. Pour $i\in Jord_{bp}(\lambda'')$, on a $\epsilon''_{i}=(-1)^{mult_{\lambda}(>i;pair)}$.
}\end{prop}

{\bf Remarque.} Comme en \ref{preuveDnram}, la derni\`ere assertion signifie que l'on peut choisir un rel\`evement $\epsilon''$ v\'erifiant les \'egalit\'es indiqu\'ees.
\bigskip

Preuve. On consid\`ere les donn\'ees de l'assertion  (i). Comme en \ref{preuveDnram}, on doit prouver 
  l'assertion suivante. Soient $m'\in \{0,...,N'\}$ et  $m''\in \{0,...,N''\}$. Posons $m=m'+m''$. Alors
 
 (1) $S_{m'}(\lambda')+S_{m''}(\lambda'')\leq S_{m}(\lambda)$. 
 
 Le nombre $S_{m'}(\lambda')$ est calcul\'e par \ref{groupesclassiques}(1). 
 Puisque ${\bf C}(\lambda')=A(\alpha')\cup B(\beta')$, il existe deux entiers $x'\in \{0,...,n'+1+[\frac{k}{2}]\}$ et $y'\in \{0,...,n'-[\frac{k}{2}]\}$ tels que $x'+y'=m'$ et $S_{m'}({\bf C}(\lambda'))=S_{x'}(A(\alpha'))+S_{y'}(B(\beta'))$. On calcule
 $$S_{x'}(A(\alpha'))=S_{x'}(\alpha')+x'(2n'+k+1-x'),\, S_{y'}(B(\beta'))=S_{y'}(\beta')+y'(2n'-k-y').$$
 D'o\`u
 $$S_{m'}(\lambda')=2S_{x'}(\alpha')+2S_{y'}(\beta')+e',$$
 o\`u
 $$e'=2x'(2n'+k+1-x')+2y'(2n'-k-y')+
 {m'}^2-m'-4n'm'+\nu_{m'}(\lambda')$$
 $$=x'(2k-2x'+1)+y'(-2k-2y'-1)+
 {m'}^2+\nu_{m'}(\lambda').$$
 
 Il existe de m\^eme des entiers $x''\in \{0,...,n''+[\frac{h+1}{2}]\}$ et $y''\in \{0,...n''-[\frac{h}{2}]\}$ tels que $x''+y''=m''$ et $S_{m''}({\bf C}(\lambda''))=S_{x''}(A(\alpha''))+S_{y''}(B(\beta''))$. Un calcul analogue
\`a celui ci-dessus (que l'on a d\'ej\`a fait en \ref{preuveDnram} dans le cas o\`u $N''$ est impair) conduit \`a l'\'egalit\'e
$$ S_{m''}(\lambda'')=2S_{x''}(\alpha'')+2S_{y''}(\beta'')+e'',$$
o\`u
$$e''=2x''(h-x'')+2y''(-h-y'')+{m''}^2+\nu_{m''}(\lambda'').$$
 
 D'o\`u
 $$ S_{m'}(\lambda')+S_{m''}(\lambda'')=2S_{x'}(\alpha')+2S_{y'}(\beta')+2S_{x''}(\alpha'')+2S_{y''}(\beta'')+
  e'+e''.$$
 L'hypoth\`ese que $\rho_{\alpha',\beta'}\otimes \rho_{\alpha'',\beta''}$ intervient dans la restriction de $\rho_{\alpha,\beta}$ implique que $\alpha\geq \alpha'\cup \alpha''$ et $\beta\geq \beta'\cup \beta''$, cf. \cite{W3} 1.13(38). En posant $x=x'+x''$ et $y=y'+y''$, on a alors 
 $$ S_{x'}(\alpha')+S_{x''}(\alpha'')\leq S_{x}(\alpha),\, S_{y'}(\beta')+S_{y''}(\beta'')\leq S_{y}(\beta),$$
 d'o\`u
 $$(2) \qquad S_{m'}(\lambda')+S_{m''}(\lambda'')\leq 2S_{x}(\alpha)+2S_{y}(\beta)+ e'+e''.$$
 
 Utilisons les constructions de \ref{An-1ramSI}. On a choisi un entier $t\geq n$ de m\^eme parit\'e que $v+1$ et d\'efini ${\bf Z}(\lambda)=\lambda+[t-1,0]$. On a $S_{m}(\lambda)=S_{m}({\bf Z}(\lambda))-m(2t-1-m)/2$. Puisque ${\bf Z}(\lambda)=X(\alpha)\cup Y(\beta)$ et $x+y=m$, on a
 $S_{m}({\bf Z}(\lambda))\geq S_{x}(X(\alpha))+S_{y}(Y(\beta))$. On calcule
 $$S_{x}(X(\alpha))=2S_{x}(\alpha)+x(t+v-x),\, S_{y}(Y(\beta))=2S_{y}(\beta)+y(t-v-1-y).$$
  D'o\`u
 $$S_{m}(\lambda)\geq 2S_{x}(\alpha)+2S_{y}(\beta)+e,$$
 o\`u
 $$e=-m(2t-1-m)/2+x(t+v-x)+y(t-v-1-y)=xv-x^2-y(v+1)-y^2+m(m+1)/2.$$
  Gr\^ace \`a (2), on en d\'eduit 
 $$(3) \qquad S_{m'}(\lambda')+S_{m''}(\lambda'')\leq S_{m}(\lambda)+E,$$
 o\`u $E=e'+e''-e$. 
 En se rappelant la d\'efinition $k=[\frac{v}{2}]$ et $h=[\frac{v+1}{2}]$, on calcule 
 $$(4) \qquad E=-\frac{(x'-y'-x''+y'')(x'-y'-x''+y''+(-1)^{v+1})}{2}
 +\nu_{m'}(\lambda')+\nu_{m''}(\lambda'').$$
 Le premier terme de $E$ est n\'egatif ou nul. 
  Si $\nu_{m'}(\lambda')=\nu_{m''}(\lambda'')=0$, alors $E\leq 0$ et (1) r\'esulte de (3). La m\^eme preuve que l'on a utilis\'ee dans \ref{preuveDnram} (cf. l'assertion (9) de cette r\'ef\'erence et le paragraphe qui suit sa d\'emonstration) montre que (1) est aussi v\'erifi\'ee en tout cas. Cela prouve cette assertion (1), c'est-\`a-dire le (i) de l'\'enonc\'e. 
  
 Passons \`a la preuve de (ii). Soient $(\lambda',\lambda'')\in {\cal P}^{symp}(N')\times {\cal P}^{orth}(N'')$, supposons $\lambda'\cup \lambda''=\lambda$. On consid\`ere que $\lambda$ a $t$ termes, on pose $mult_{\lambda}(0)=t-l(\lambda)$ et on note $Jord_{*}(\lambda)$
l'ensemble des $i\in {\mathbb N}$ tels que $mult_{\lambda}(i)>0$. Pour tout $i\in Jord_{*}(\lambda)$, on note $\{j_{min}(i),...,j_{max}(i)\}$ l'intervalle des $j\in \{1,...,t\}$ tels que $\lambda_{j}=i$. 
Notons $i_{1}>....> i_{s}$ l'ensemble des $i\in Jord_{*}(\lambda)$ tels que $mult_{\lambda}(i)$ soit impair. On a $s\equiv t\,\, mod\,\, 2{\mathbb Z}$. On voit par r\'ecurrence que $j_{min}(i_{r})\equiv j_{max}(i_{r})\equiv r\,\,mod\,\,2{\mathbb Z}$ pour tout $r=1,...,s$. Par construction, $X(\alpha)$ est le sous-ensemble des \'el\'ements pairs de ${\bf Z}(\lambda)$ et $Y(\beta)$ est le sous-ensemble des \'el\'ements impairs. On a aussi $v+1=\vert X(\alpha)\vert -\vert Y(\beta)\vert $. Un \'el\'ement $i\in Jord_{*}(\lambda)$ contribue \`a ${\bf Z}(\lambda)$ par les \'el\'ements $i+t-j
$ o\`u $j$ parcourt $\{j_{min}(i),...,j_{max}(i)\}$. Si $mult_{\lambda}(i)$ est pair, il y a dans cet ensemble autant d'\'el\'ements pairs que d'impairs. Si $mult_{\lambda}(i)$ est impair et $i+t-j_{min}(i)$ est impair, il y a $(mult_{\lambda}(i)+1)/2$ \'el\'ements impairs et $(mult_{\lambda}(i)-1)/2$ \'el\'ements pairs. Si $mult_{\lambda}(i)$ est impair et $i+t-j_{min}(i)$ est pair, il y a $(mult_{\lambda}(i)+1)/2$ \'el\'ements pairs et $(mult_{\lambda}(i)-1)/2$ \'el\'ements impairs. Ces consid\'erations conduisent \`a l'\'egalit\'e
$$ v+1=\vert \{r\in \{1,...,s\}; i_{r}+t-r\,\, pair\}\vert -\vert \{r\in \{1,...,s\}; i_{r}+t-r\,\, impair\}\vert .$$
Pour $\eta,\tau\in \{\pm 1\}$, posons
$$e(\eta,\tau)=\vert \{r\in \{1,...,s\}; (-1)^{i_{r}+t-r}=\eta,(-1)^{r}=\tau\}\vert . $$
 L'\'egalit\'e pr\'ec\'edente se r\'ecrit
 $$v+1=e(1,1)+e(1,-1)-e(-1,1)-e(-1,-1).$$
 La somme $e(1,1)+e(-1,1)$, resp. $e(1,-1)+e(-1,-1)$, est le nombre de termes pairs, resp. impairs, dans $\{1,...,s\}$, c'est-\`a-dire $[\frac{s}{2}]$, resp. $[\frac{s+1}{2}]$. L'\'egalit\'e ci-dessus se r\'ecrit sous les deux formes
 $$v+1-[\frac{s+1}{2}]+[\frac{s}{2}]=2e(1,1)-2e(-1,-1),\,\, v+1+[\frac{s+1}{2}]-[\frac{s}{2}]=2e(1,-1)-2e(-1,1).$$
 Parce que $s\equiv t\equiv v+1\,\,mod\,\,2{\mathbb Z}$, on a $[\frac{s+1}{2}]-[\frac{s}{2}]=[\frac{v+2}{2}]-[\frac{v+1}{2}]$. On a aussi $v+1=[\frac{v+2}{2}]+[\frac{v+1}{2}]$ et les \'egalit\'es pr\'ec\'edentes se r\'ecrivent
  $$(5) \qquad [\frac{v+1}{2}]=e(1,1)-e(-1,-1),\,\, [\frac{v+2}{2}]=e(1,-1)-e(-1,1).$$
 A fortiori
 $$(6) \qquad e(1,1)\geq [\frac{v+1}{2}], \,\, e(1,-1)\geq [\frac{v+2}{2}].$$
 Supposons $t$ pair. Par d\'efinition, $e(1,-1)$ est inf\'erieur ou \'egal au nombre des $r\in \{1,...,s\}$ tels que $i_{r}$ est  impair. Pour un tel entier, la multiplicit\'e  $mult_{\lambda'}(i_{r})$ est paire puisque $\lambda'$ est symplectique. Puisque $mult_{\lambda}(i_{r})$ est impaire par d\'efinition de la suite $(i_{1},...,i_{s})$, on a $mult_{\lambda''}(i_{r})>0$. D'o\`u $e(1,-1)\leq \vert Jord_{bp}(\lambda'')\vert $. Puisque $t$ est pair, $v$ est impair et $ [\frac{v+2}{2}]= [\frac{v+1}{2}]=h$. Alors (6) implique $h\leq \vert Jord_{bp}(\lambda'')\vert $. La somme des \'el\'ements de $Jord_{bp}(\lambda'')$ est sup\'erieure ou \'egale \`a $1+3+.... +2h-1=h^2$. D'o\`u $h^2\leq N''$. Par d\'efinition, $e(1,1)$ est inf\'erieur ou \'egal au nombre des $r\in \{1,...,s\}$ tels que $i_{r}$ est  pair. Il faut faire  attention d'exclure le cas o\`u $i_{s}=0$. En tout cas, $e(1,1)-1$ est inf\'erieur ou \'egal au nombre des $r\in \{1,...,s\}$ tels que $i_{r}$ est  pair et non nul. On prouve comme ci-dessus que ce dernier nombre est inf\'erieur ou \'egal \`a $\vert Jord_{bp}(\lambda')\vert $. D'o\`u $e(1,1)-1\leq \vert  Jord_{bp}(\lambda')\vert $. Puisque $v$ est impair, on a $[\frac{v+1}{2}]=[\frac{v}{2}]+1=k+1$. Alors (6) implique $k\leq \vert  Jord_{bp}(\lambda')\vert $.  La somme des \'el\'ements de $Jord_{bp}(\lambda')$ est sup\'erieure ou \'egale \`a $2+4+... +2k=k(k+1)$. D'o\`u $k(k+1)\leq N'$. Cela d\'emontre que $\Lambda_{v}\subset \Lambda(s)$. On laisse au lecteur le cas similaire o\`u $t$ est impair. Cela prouve la premi\`ere assertion de (ii). 

  Consid\'erons des triplets $(k,\alpha',\beta')$ et $(h,\alpha'',\beta'')$ v\'erifiant les conditions (a) et (b) du (ii) de l'\'enonc\'e.  Pour un entier $m\in \{0,...,n\}$, on choisit un couple d'entiers $(m'_{m},m''_{m})\in \{0,...,N'\}\times\{0,...,N''\}$ tels que $m=m'_{m}+m''_{m}$ et $S_{m}(\lambda)=S_{m'_{m}}(\lambda')+S_{m''_{m}}(\lambda'')$. Cette condition ne d\'etermine pas enti\`erement le couple. On impose de plus la condition suivante: si $\lambda_{m}$ est pair,  on choisit parmi les couples possibles celui pour lequel $m'_{m}$ est maximal; si $\lambda_{m}$ est impair,  on choisit parmi les couples possibles celui pour lequel $m''_{m}$ est maximal. On reprend nos constructions en indexant par $m$ les entiers que l'on a introduits: $x'_{m}$, $y'_{m}$ etc... Pour $m\geq1$, on pose ${\bf m}'_{m}=m'_{m}-m'_{m-1}$, ${\bf x}'_{m}=x'_{m}-x'_{m-1}$ etc... Soit $i\in Jord_{bp}(\lambda')$. Posons $M=mult_{\lambda}(> i)$, $M'=mult_{\lambda'}(>i)$, $M''=mult_{\lambda''}(>i)$ et ${\bf M}'=mult_{\lambda'}(i)$.  Il r\'esulte de nos constructions que $m'_{M}=M'$, $m''_{M}=M''$ et que, pour tout $\ell\in \{1,...,{\bf M}'\}$, on a
$m'_{M+\ell}=M'+\ell$, $m''_{M+\ell}=M''$. Soit $m\in \{M,...,M+{\bf M}'\}$. On a   $\lambda''_{m''_{m}}=\lambda''_{M}> \lambda''_{m''_{m}+1}$ donc $\nu_{m''_{m}}(\lambda'')=0$.  On a soit $\lambda'_{m'_{m}}> \lambda'_{m'_{m}+1}$ si $m=M$, soit  $\lambda'_{m'_{m}}=i\in Jord_{bp}(\lambda')$ si $m\geq M+1$, en tout cas $\nu_{m'_{m}}(\lambda')=0$. Alors $E\leq 0$.  Nos hypoth\`eses et la relation (3) entra\^{\i}nent
que $E=0$ et que cette relation est une \'egalit\'e. Toutes les in\'egalit\'es utilis\'ees dans le calcul sont donc elles-aussi des \'egalit\'es. La nullit\'e de $E$ entra\^{\i}ne que $x'_{m}-y'_{m}-x''_{m}+y''_{m}$ vaut $0$ ou $(-1)^v$. Mais $x'_{m}-y'_{m}-x''_{m}+y''_{m}$ est de m\^eme parit\'e que $x'_{m}+y'_{m}+x''_{m}+y''_{m}=m$. D'o\`u
$x'_{m}-y'_{m}-x''_{m}+y''_{m}=0$ si $m$ est pair, $x'_{m}-y'_{m}-x''_{m}+y''_{m}=(-1)^v$ si $m$ est impair. 
Par diff\'erence, cela entra\^{\i}ne que   ${\bf x}'_{M+1}-{\bf y}'_{M+1}-{\bf x}''_{M+1}+{\bf y}''_{M+1}=(-1)^{v+M}$. Mais ${\bf m}''_{M+1}=0$ donc ${\bf x}''_{M+1}={\bf y}''_{M+1}=0$ et $({\bf x}'_{M+1},{\bf y}'_{M+1})$ est \'egal \`a $(1,0)$ ou $(0,1)$. D'o\`u
$$ ({\bf x}'_{M+1},{\bf y}'_{M+1})=\left\lbrace\begin{array}{cc}(1,0),&\text{ si }v+M\text{ est pair},\\ (0,1),&\text{ si }v+M\text{ est impair}\\ \end{array}\right.$$
Par d\'efinition de $x'_{M+1}$ et $y'_{M+1}$, cela signifie que le plus grand terme de l'intervalle ${\bf C}(\lambda',i)$ intervient dans $A(\alpha')$ si $v+M$ est pair et dans $B(\beta')$ si $v+M$ est impair. En se rappelant comment sont reli\'es $A(\alpha')$, $B(\beta')$, ${\bf A}(\lambda',\epsilon')$ et ${\bf B}(\lambda',\epsilon')$,  on voit que le plus grand terme de l'intervalle ${\bf C}(\lambda',i)$ intervient dans ${\bf A}(\lambda',\epsilon')$ si $v+M+k$ est pair, dans ${\bf B}(\lambda',\epsilon')$ si $v+M+k$ est impair. On se rappelle, cf. \ref{groupesclassiques}, que le couple $({\bf A}(\lambda',\epsilon'),{\bf B}(\lambda',\epsilon'))$ se déduit d'un couple $({\bf A}(\lambda'),{\bf B}(\lambda'))$ et le plus grand terme de ${\bf C}(\lambda',i)$
intervient dans ${\bf A}(\lambda')$ si $M'$ est pair et dans ${\bf B}(\lambda')$ si $M'$ est impair.  Les définitions entra\^{\i}nent alors que  $\epsilon'_{i}=(-1)^{v+M+k+M'}$.  On a
$$M+M'=mult_{\lambda}(>i)+mult_{\lambda'}(>i)=\sum_{j>i}(mult_{\lambda}(j)+mult_{\lambda'}(j)).$$
Modulo $2{\mathbb Z}$, on a $mult_{\lambda'}(j)\equiv 0$ si $j$ est impair, puisque $\lambda'$ est symplectique. Si $j$ est pair, on a $mult_{\lambda}(j)+mult_{\lambda'}(j)\equiv mult_{\lambda''}(j)\equiv 0$ puisque $\lambda''$ est orthogonale. Donc $M+M'\equiv mult_{\lambda}(>i;imp)\,\,mod\,\,2{\mathbb Z}$. R\'ecrivons cette relation, qui nous servira plus loin:

(7) $ mult_{\lambda}(>i)+mult_{\lambda'}(>i)\equiv mult_{\lambda}(>i;imp)\,\,mod\,\,2{\mathbb Z}$.

\noindent On a aussi   $v=h+k$ d'o\`u $v+k\equiv h\equiv N''\equiv n\,\, mod\,\, 2{\mathbb Z}$. D'o\`u 
 $\epsilon'_{i}=(-1)^{n+mult_{\lambda}(>i;imp)}$. C'est une des \'egalit\'es affirm\'ees dans le (ii) de l'\'enonc\'e. 
 
 Pour le caract\`ere $\epsilon''$, on commence par le relever en un \'el\'ement de $\{\pm 1\}^{Jord_{bp}(\lambda'')}$ tel que ${\bf A}(\lambda'',\epsilon'')=A(\alpha'')$, ${\bf B}(\lambda'',\epsilon'')=B(\beta'')$. Un calcul analogue, qu'on laisse au lecteur, conduit \`a l'\'egalit\'e 
  $\epsilon''_{i}=(-1)^{v+1+mult_{\lambda}(>i;pair)}$ pour $i\in Jord_{bp}(\lambda'')$. A multiplication \'eventuelle pr\`es par $-{\bf 1}$, cela \'equivaut \`a la formule de l'\'enonc\'e. 
 
\bigskip 
Puisque les caract\`eres $\epsilon'$ et $\epsilon''$ sont uniquement d\'etermin\'es, il y a au plus un couple de triplets $(k,\alpha',\beta')$ et $(h,\alpha'',\beta'')$ v\'erifiant les conditions (a) et (b) du (ii) de l'\'enonc\'e. Il faut maintenant prouver son existence.  On introduit les \'el\'ements $\epsilon'\in \{\pm 1\}^{Jord_{bp}(\lambda')}$ et $\epsilon''\in \{\pm 1\}^{Jord_{bp}(\lambda'')}$ d\'efinis par $\epsilon'_{i}=(-1)^{n+mult_{\lambda}(>i;imp)}$ pour tout $i\in Jord_{bp}(\lambda')$ et $\epsilon''_{i}=(-1)^{v+1+mult_{\lambda}(>i;pair)}$ pour tout $i\in Jord_{bp}(\lambda'')$. Introduisons les triplets $(\underline{k},\alpha',\beta')\in {\cal I}^{FC}(G'_{s})$ et $(\underline{h},\alpha'',\beta'')\in {\cal I}^{FC}(G''_{s})$ dont les images par la correspondance de Springer g\'en\'eralis\'ee sont $(\lambda',\epsilon')$, $(\lambda'',\epsilon'')$. On a introduit en \ref{groupesclassiques} pour chacune de ces paires un entier not\'e $D$, que l'on note ici $D'$  pour la paire $(\lambda',\epsilon')$ et $D''$ pour la paire $(\lambda'',\epsilon'')$.  Montrons que

(8) $\underline{k}=k$, $\underline{h}=h$ et $D''\geq0$.

 Comme en \ref{groupesclassiques},  posons $J'=\{i\in Jord_{bp}(\lambda'); mult_{\lambda'}(i)\text{ est impair et }\epsilon'_{i}=-1\}$. On a
 $$D'=\sum_{i\in J'}(-1)^{mult_{\lambda'}(\geq i)}.$$
 Posons ${\bf J}'=\{i\in Jord_{bp}(\lambda'); mult_{\lambda'}(i)\text{ est impair}\}$. 
  On peut remplacer la somme en $i\in J'$ par une somme en $i\in {\bf J}'$, en glissant un terme $\frac{1}{2}(1-\epsilon'_{i})$ dans la somme. On obtient
$$D'=\frac{1}{2}\sum_{i\in {\bf J}'}(1-\epsilon'_{i})(-1)^{mult_{\lambda'}(\geq i)}=D'_{1}+D'_{2},$$
o\`u
$$D'_{1}=\frac{1}{2}\sum_{i\in {\bf J}'}(-1)^{mult_{\lambda'}(\geq i)},$$
$$D'_{2}=-\frac{1}{2}\sum_{i\in {\bf J}'}\epsilon'_{i}(-1)^{mult_{\lambda'}(\geq i)}.$$
Pour $i\in {\bf J}'$, on a $(-1)^{mult_{\lambda'}(\geq i)}=\frac{1}{2}((-1)^{mult_{\lambda'}(\geq i)}-(-1)^{mult_{\lambda'}(> i)})$. Pour $i\in Jord(\lambda')-{\bf J}'$, cette derni\`ere expression est nulle. Donc
$$D'_{1}=\frac{1}{4}\sum_{i\in Jord(\lambda')}((-1)^{mult_{\lambda'}(\geq i)}-(-1)^{mult_{\lambda'}(> i)}).$$
Cette  somme se simplifie et on obtient
$$D'_{1}=\frac{1}{4}((-1)^{mult_{\lambda'}(\geq i_{min})}-1),$$
o\`u $i_{min}$ est le plus petit terme de $\lambda'$. On a $mult_{\lambda'}(\geq i_{min})=l(\lambda')=l(\lambda)-l(\lambda'')$ (rappelons par exemple que $l(\lambda)$ est le nombre de termes non nuls de $\lambda$). Puisque $\lambda''$ est orthogonale, on a $l(\lambda'')\equiv N''\equiv n\,\,mod\,\,2{\mathbb Z}$, d'o\`u
$$(9)\qquad D'_{1}=\frac{1}{4}((-1)^{n+l(\lambda)}-1).$$
 D'apr\`es  (7), pour $i\in {\bf J}'$, on a $\epsilon'_{i}=(-1)^{n+mult_{\lambda}(>i)+mult_{\lambda'}(>i)}$. On a aussi $(-1)^{mult_{\lambda'}(\geq i)}=(-1)^{mult_{\lambda'}(> i)+1}$ et $D'_{2}$ se r\'ecrit
 $$D'_{2}=\frac{1}{2}(-1)^n\sum_{i\in {\bf J}'}(-1)^{mult_{\lambda}(>i)}.$$
 On a introduit plus haut la suite $i_{1}> ... > i_{s}$. On voit que ${\bf J}'$ est l'ensemble des $i_{r}$ pour $r=1,...,s$, tels que $i_{r}$ est pair et $i_{r}\not=0$: cela r\'esulte du fait que, pour un entier $i$ pair strictement positif, la multiplicit\'e $mult_{\lambda''}(i)$ est paire donc $mult_{\lambda'}(i)$ est impaire si et seulement si $mult_{\lambda}(i)$ l'est. Pour $r=1,...,s$, posons $\eta_{r}=(-1)^{i_{r}+t-r}$, $\tau_{r}=(-1)^r$. Alors $i_{r}$ est pair si et seulement si $\eta_{r}\tau_{r}=(-1)^t$. On a remarqu\'e plus haut que $j_{min}(i_{r})$ \'etait de la m\^eme parit\'e que $r$. On a  $mult_{\lambda}(>i_{r})=j_{min}(i_{r})-1$, d'o\`u $(-1)^{mult_{\lambda}(>i_{r})}=(-1)^{r-1}=-\tau_{r}$. On obtient
 $$D'_{2}=\frac{1}{2}(-1)^{n+1}\sum_{r=1,...,s; \eta_{r}\tau_{r}=(-1)^t,\, i_{r}\not=0}\tau_{r}.$$ 
  Ceci est \'egal \`a la m\^eme somme o\`u l'on supprime la restriction $i_{r}\not=0$, dont on retire la contribution de l'éventuel  terme $i_{s}=0$.  Ce terme existe si $mult_{\lambda}(0)$ est impaire, c'est-\`a-dire $t-l(\lambda)$ est impair. Dans ce cas, sa contribution est   $\tau_{s}=(-1)^s=(-1)^t$. En tout cas, cette contribution vaut $\frac{1}{2}(-1)^t(1-(-1)^{l(\lambda)+t})=\frac{1}{2}((-1)^{t}-(-1)^{l(\lambda)})$.  
  D'o\`u
  $D'_{2}=D'_{3}+D'_{4}$, o\`u
  $$D'_{3}=  \frac{1}{2}(-1)^{n+1}\sum_{r=1,...,s; \eta_{r}\tau_{r}=(-1)^t}\tau_{r},$$
  $$D'_{4}= \frac{1}{4}(-1)^n((-1)^t-(-1)^{l(\lambda)}).$$
  Gr\^ace \`a (9), on a 
  $$D'_{1}+D'_{4}=\frac{1}{4}((-1)^{n+t}-1).$$
  Avec les notations introduites plus haut, on a
  $$D'_{3}=\frac{1}{2}(-1)^{n+1}(e((-1)^t,1)-e((-1)^{t+1},-1)).$$
  En se rappelant que $t$ est de la m\^eme parit\'e que $v+1$ et que $k=[\frac{v}{2}]$, il r\'esulte de (5) que
  $$D'_{3}=\frac{1}{2}(-1)^{n+t+1}(k+1),$$
  d'o\`u
  $$D'=D'_{1}+D'_{3}+D'_{4}=\frac{1}{2}(-1)^{n+t+1}(k+1)+\frac{1}{4}((-1)^{n+t}-1).$$
  Alors $2D'+\frac{1}{2}=(-1)^{n+t+1}(k+\frac{1}{2})$, puis $\underline{k}=\vert 2D'+\frac{1}{2}\vert -\frac{1}{2}=k$.
  
  On laisse au lecteur la preuve similaire des assertions concernant $\underline{h}$ et $D''$.   Cela prouve  (8). 
  
  Si $h>0$, la relation $D''\geq0$ entra\^{\i}ne que ${\bf A}(\lambda'',\epsilon'')=A(\alpha'')$ et ${\bf B}(\lambda'',\epsilon'')=B(\beta'')$. Si $h=0$, on a le droit de permuter $\alpha''$ et $\beta''$ et on suppose les \'egalit\'es pr\'ec\'edentes v\'erifi\'ees. 
  
  Reprenons nos calculs pour ces triplets $(k,\alpha',\beta')$ et $(h,\alpha'',\beta'')$. Pour $m\in \{0,...,n\}$, on choisit selon la recette expliqu\'ee plus haut des entiers $m'_{m}\in \{0,...,N'\}$ et $m''_{m}\in \{0,...,N''\}$ tels que $m=m'_{m}+m''_{m}$ et $S_{m}(\lambda)=S_{m'_{m}}(\lambda')+S_{m''_{m}}(\lambda'')$. On choisit ensuite des entiers $x'_{m}$, $y'_{m}$, $x''_{m}$, $y''_{m}$ de sorte que
  $$S_{m'_{m}}({\bf C}(\lambda'))=S_{x'_{m}}(A(\alpha'))+S_{y'_{m}}(B(\beta')),\,\, S_{m''_{m}}({\bf C}(\lambda''))=S_{x''_{m}}(A(\alpha'))+S_{y''_{m}}(B(\beta')).$$
  On impose \`a ces termes des conditions suppl\'ementaires.  On suppose que les suites $m\mapsto x'_{m}$ etc... sont croissantes. Supposons $m'_{m}=m'_{m-1}+1$. Si $\lambda'_{m'_{m}}$ est pair, il n'y a pas le choix pour $x'_{m}$ et $y'_{m}$. Supposons $\lambda'_{m'_{m}}$ impair. Si $m'_{m}-mult_{\lambda'}(>\lambda'_{m'_{m}})$ est pair, il n'y a de nouveau pas le choix et on a $x'_{m}=x'_{m-2}+1$, $y'_{m}=y'_{m-2}+1$.  Mais, si
     $m'_{m}-mult_{\lambda'}(> \lambda'_{m'_{m}})$ est impair, les deux couples $(x'_{m-1}+1,y'_{m-1})$ et $(x'_{m-1},y'_{m-1}+1)$ conviennent. On choisit le couple $(x'_{m},y'_{m}) $ de sorte que $x'_{m}-y'_{m}=x'_{m-1}-y'_{m-1}+(-1)^{v+m}$. En revenant au cas o\`u  $m'_{m}-mult_{\lambda'}(>\lambda'_{m'_{m}})$ est pair, les \'egalit\'es $x'_{m}=x'_{m-2}+1$ et $y'_{m}=y'_{m-2}+1$ et l'\'egalit\'e pr\'ec\'edente pour $m-1$, c'est-\`a-dire $x'_{m-1}-y'_{m-1}=x'_{m-2}-y'_{m-2}+(-1)^{v+m-1}$ impliquent $x'_{m}-y'_{m}=x'_{m-1}-y'_{m-1}+(-1)^{v+m}$. Cette \'egalit\'e est donc v\'erifi\'ee pour tout $m$ tel que $m'_{m}=m'_{m-1}+1$ et $\lambda'_{m'_{m}}$ est impair. 
     On fait un choix analogue pour les entiers $x''_{m}$ et $y''_{m}$, au signe pr\`es. C'est-\`a-dire que si $m''_{m}=m''_{m-1}+1$ et si $\lambda''_{m''_{m}}$ est pair,    $(x''_{m},y''_{m}) $ est tel que $x''_{m}-y''_{m}=x''_{m-1}-y''_{m-1}+(-1)^{v+m+1}$. Pour $m\geq1$, on d\'efinit ${\bf m}'_{m}=m'_{m}-m'_{m-1}$, ${\bf x}'_{m}=x'_{m}-x'_{m-1}$ etc... 
  
  On d\'efinit $E_{m}$ par la formule (4), o\`u l'on ajoute des indices $m$ \`a tous les termes. Prouvons que, pour tout $m\in \{0,...,n\}$ on a les \'egalit\'es
  
  (10) $S_{m}({\bf Z}(\lambda))=S_{x_{m}}(X(\alpha))+S_{y_{m}}(Y(\beta))$;
  
  (11) $E_{m}=0$. 
  
  La d\'emonstration se fait par r\'ecurrence sur $m$.  Posons $i=\lambda_{m}$, supposons $i$ pair et ${\bf m}'_{m}=1$, ${\bf m}''_{m}=0$. D'apr\`es notre choix des suites $(m'_{m})_{m=0,...,n}$ et $(m''_{m})_{m=0,...,n}$, cela signifie que $m\in \{M+1,...,M+{\bf M}'\}$, o\`u $M=mult_{\lambda}(>i)$ et ${\bf M}'=mult_{\lambda'}(i)$. Pour $\ell=1,...,{\bf M}'$, les couples $({\bf x}'_{M+\ell},{\bf y}'_{M+\ell})$ sont alternativement $(1,0)$ et $(0,1)$. 
   Dans la preuve ci-dessus de l'unicit\'e des triplets, on a d\'eduit la valeur de $\epsilon'_{i}$ \`a partir de la connaissance des couples $({\bf x}'_{M+\ell},{\bf y}'_{M+\ell})$. La m\^eme preuve s'inverse: ici on conna\^{\i}t la valeur de $\epsilon'_{i}$ et on en d\'eduit que
   $$(12) \qquad ({\bf x}'_{m},{\bf y}'_{m})=\left\lbrace\begin{array}{cc}(1,0),& \text{ si }v+m\text{ est impair},\\ (0,1),&\text{ si }v+m\text{ est pair.}\\ \end{array}\right.$$
   Par d\'efinition de $X(\alpha)$ et $Y(\beta)$ et d'apr\`es la relation (10) pour $m-1$, cette relation (10) pour $m$ (que l'on doit prouver)  \'equivaut aux \'egalit\'es $({\bf x}_{m},{\bf y}_{m})=(1,0)$ si $\lambda_{m}+t-m$ est pair, $({\bf x}_{m},{\bf y}_{m})=(0,1)$ si $\lambda_{m}+t-m$ est impair. Puisque $({\bf x}_{m},{\bf y}_{m})=({\bf x}'_{m},{\bf y}'_{m})$, que $\lambda_{m}=i$ est pair et que $t$ et $v$ sont de parit\'e oppos\'ee, ces \'egalit\'es r\'esultent de (12). D'o\`u (10) pour $m$. Par un raisonnement d\'ej\`a fait, on a $\nu_{m'_{m}}(\lambda')=\nu_{m'_{m-1}}(\lambda')=\nu_{m''_{m}}(\lambda'')=\nu_{m''_{m-1}}(\lambda'')=0$. On se rappelle que $x'_{m}-y'_{m}-x''_{m}+y''_{m}$ est de la m\^eme parit\'e que $m$. La relation (11) pour $m-1$ signifie donc que $x'_{m-1}-y'_{m-1}-x''_{m-1}+y''_{m-1}$ vaut $0$ si $m-1$ est pair et $(-1)^v$ si $m-1$ est impair. Compte tenu de cette relation, la relation (11) pour $m$ \'equivaut aux \'egalit\'es ${\bf x}'_{m}-{\bf y}'_{m}=(-1)^{v+1}$ si $m$ est pair et ${\bf x}'_{m}-{\bf y}'_{m}=(-1)^v$ si $m$ est impair. Cela r\'esulte encore de (12), d'o\`u la relation (11) pour $m$.  
   
   Supposons encore   ${\bf m}'_{m}=1$ et ${\bf m}''_{m}=0$ mais supposons $i$ impair. De nouveau, la relation (10) pour $m$ \'equivaut aux \'egalit\'es $({\bf x}'_{m},{\bf y}'_{m})=(1,0)$ si $\lambda_{m}+t-m$ est pair, $({\bf x}'_{m},{\bf y}'_{m})=(0,1)$ si $\lambda_{m}+t-m$ est impair. Mais ici $\lambda_{m}+t$ est de m\^eme parit\'e que $v$ et on a choisi les termes $x'_{m}$ et $y'_{m}$ de sorte que les conditions pr\'ec\'edentes soient v\'erifi\'ees. On a $\nu_{m''_{m}}(\lambda'')=0$. Supposons que $m'_{m}-mult_{\lambda'}(>i)$ soit pair. Alors $\nu_{m'_{m}}(\lambda')=0$. Comme ci-dessus, la  relation (11) pour $m$ \'equivaut aux relations: $x'_{m}-y'_{m}-x''_{m}+y''_{m}$ vaut $0$ si $m$ est pair et $(-1)^v$ si $m$ est impair. Mais on a aussi $\nu_{m'_{m-2}}(\lambda')=0$ car, ou bien $m'_{m-2}-mult_{\lambda'}(>i)>0$ et alors ce nombre est pair, ou bien $m'_{m-2}=mult_{\lambda'}(>i)$ et alors $\lambda'_{m'_{m-2}}> \lambda'_{m'_{m-2}+1}=i$. Les relations pr\'ec\'edentes sont donc vraies quand on remplace $m$ par $m-2$. Puisque $(x'_{m},y'_{m})=(x'_{m-2}+1,y'_{m-2}+1)$, elles restent vraies pour $m$, d'o\`u (11) pour $m$. Supposons maintenant que $m'_{m}-mult_{\lambda'}(>i)$ soit impair. Alors $\nu_{m'_{m}}(\lambda')=1$. L'\'egalit\'e $E_{m}=0$ \'equivaut alors \`a
   $$\frac{1}{2}(x'_{m}-y'_{m}-x''_{m}+y''_{m})(x'_{m}-y'_{m}-x''_{m}+y''_{m} +(-1)^{v+1})=1.$$
   Mais l'\'egalit\'e $E_{m-1}=0$ \'equivaut comme ci-dessus aux relations: $x'_{m-1}-y'_{m-1}-x''_{m-1}+y''_{m-1}$ vaut $0$ si $m-1$ est pair et $(-1)^v$ si $m-1$ est impair. On en d\'eduit   que l'\'egalit\'e $E_{m}=0$ \'equivaut aux \'egalit\'es ${\bf x}'_{m}-{\bf y}'_{m}=(-1)^{v}$ si $m$ est pair, ${\bf x}'_{m}-{\bf y}'_{m}=(-1)^{v+1}$ si $m$ est impair. Mais on a choisi les termes $x'_{m}$ et $y'_{m}$ de sorte qu'il en soit ainsi. Cela prouve (11) pour $m$. 
   
   On laisse au lecteur la preuve analogue dans le cas o\`u ${\bf m}'_{m}=0$ et ${\bf m}''_{m}=1$. Cela prouve (10) et (11). 
   
A l'aide de ces relations, on ach\`eve la preuve de la proposition comme celle de la proposition \ref{preuveDnram}. On pose $T_{m}=S_{x_{m}}(\alpha)-S_{x'_{m}}(\alpha')-S_{x''_{m}}(\alpha'')$, $U_{m}=S_{y_{m}}(\beta)-S_{y'_{m}}(\beta')-S_{y''_{m}}(\beta'')$. Les m\^emes calculs qu'au d\'ebut de la preuve, joints \`a l'\'egalit\'e (10), conduisent \`a l'\'egalit\'e
$$S_{m}(\lambda)=2T_{m}+2U_{m}+S_{m'_{m}}(\lambda')+S_{m''_{m}}(\lambda'')-E_{m}.$$
Gr\^ace \`a (11), cela entra\^{\i}ne $T_{m}+U_{m}=0$. On en d\'eduit par r\'ecurrence $T_{m}=U_{m}=0$, puis, toujours par r\'ecurrence, $\alpha=\alpha'\cup \alpha''$, $\beta=\beta'\cup \beta''$. On en d\'eduit que nos triplets $(k,\alpha',\beta')$ et $(h,\alpha'',\beta'')$ v\'erifient les conditions requises et que la multiplicit\'e de $\rho_{\alpha',\beta'}\otimes \rho_{\alpha'',\beta''}$ dans la restriction de $\rho_{\alpha,\beta}$ vaut $1$. Cela ach\`eve la preuve de la proposition. $\square$

{\bf Remarque.} Supposons que le triplet $(h,\alpha'',\beta'')$ de l'assertion (ii) se d\'edouble, c'est-\`a-dire $h=0$ et $\alpha''=\beta''$.  Il y a deux repr\'esentations $\rho_{\alpha'',\alpha''}$. Le groupe $W''_{s}$ est de type $D_{n''}$ et est un sous-groupe d'indice $2$ dans un groupe de Weyl de type $B_{n''}$.  Les deux composantes sont \'echang\'ees par l'automorphisme de ce groupe d\'eduit de l'action d'un \'el\'ement de ce groupe de Weyl de type $B_{n''}$. La restriction de $\rho_{\alpha,\beta}$ \`a $W'_{s}(M'_{s,k})\times W''_{s}$ est conserv\'ee par cet automorphisme. Donc le 
 produit tensoriel  de $\rho_{\alpha',\beta'}$ et  de chacune des deux repr\'esentations $\rho_{\alpha'',\beta''}$   intervient dans la restriction de $\rho_{\alpha,\beta}$.   Cela ne perturbe pas notre r\'esultat d'unicit\'e car la partition $\lambda''$ se d\'edouble elle aussi et chacune des repr\'esentations $\rho_{\alpha'',\beta''}$ correspond \`a une et une seule des partitions d\'edoubl\'ees.
 
 \bigskip

 Les th\'eor\`emes  \ref{premiertheoremeSI} et \ref{deuxiemetheoremeSI} se d\'eduisent de la proposition ci-dessus comme en \ref{preuveDnram}.
 
 \subsubsection{Preuve des th\'eor\`emes \ref{premiertheoremeSI} et \ref{deuxiemetheoremeSI} dans le cas $(D_{4},3-ram)$ \label{preuveD43ram}}
 
 On suppose que $G$ est de type $(D_{4},3-ram)$. On a d\'ecrit en \ref{D3ramSI} l'ensemble ${\cal B}_{Irr}(G)$ et on a d\'efini l'application $\boldsymbol{\nabla}_{F}$. On va traiter successivement tous les \'el\'ements $(M,s_{M},\varphi,\rho)\in {\cal B}_{Irr}(G)$. Pour chacun d'eux, on a d\'ecrit  l'ensemble $\Lambda$ tel que $(M,s_{M})$ soit conjugu\'e \`a $(M_{\Lambda},s_{\Lambda})$, cf. \ref{couplesstables}. On note $s_{0},s_{134},s_{2}$ les sommets de $S(\bar{C})$, correspondant aux racines $\beta_{0},\beta_{134},\beta_{2}$. Dans le tableau ci-dessous, seuls figurent les sommets $s$   tels que $\Lambda\subset \Lambda(s)$. Pour chacun d'eux,  on calcule la restriction de $\rho$ \`a $W_{s}(M_{\Lambda,s_{\Lambda}})$ puis, pour chaque composante irr\'eductible $\rho'$ de cette restriction, on calcule l'image $(\bar{{\cal O}}_{\rho'},{\cal L}_{\rho'})$ de $(M_{\Lambda,s_{\Lambda}},\rho')$ par la correspondance de Springer g\'en\'eralis\'ee.  Expliquons la notation de cette image. On peut identifier $(\bar{{\cal O}}_{\rho'},{\cal L}_{\rho'})$ \`a un couple $(N',\xi')$, o\`u $N'\in \bar{{\cal O}}_{\rho'}$ et $\xi'$ est une repr\'esentation irr\'eductible de $Z_{G_{s}}(N')/Z_{G_{s}}(N')^0$. Quand ce groupe est trivial ou quand $\xi'=1$, on remplace le couple $(N',\xi')$ par le symbole ou la partition qui param\`etre ${\bar{\cal O}}_{\rho'}$. Dans le cas o\`u $G_{s}$ est classique et $\xi'\not=1$, on ajoute la mention: $\xi'\not=1$. Dans le cas o\`u $s=s_{0}$, auquel cas $G_{s}$ est de type $G_{2}$, et o\`u $\bar{{\cal O}}_{\rho'}$ est l'orbite $G_{2}(a_{1})$, auquel cas $Z_{G_{s}}(N')/Z_{G_{s}}(N')^0\simeq \mathfrak{S}_{3}$, on \'ecrit plus pr\'ecis\'ement $\xi'$ qui est l'une des trois repr\'esentations irr\'eductibles $1,sgn,r$ de ce groupe, $r$ \'etant de dimension $2$. 
 
 Expliquons comment s'effectuent les calculs. Dans le cas o\`u $M=G$, il y a par construction un seul sommet $s$ tel que $\Lambda\subset \Lambda(s)$. C'est l'unique sommet qui figure dans la table ci-dessous.  Le syst\`eme local  ${\cal L}_{\rho'}$ est simplement le faisceau-caract\`ere cuspidal dont est issu la fonction $\varphi$. Il est d\'ecrit en \cite{W7} 8.3. Dans le cas o\`u $M=T$, on a $\Lambda=\emptyset$ et les trois sommets de $S(\bar{C})$ v\'erifient la condition $ \Lambda\subset \Lambda(s)$. Le groupe $W^{I_{F}}(T)=W^{I_{F}}$ est simplement le groupe de Weyl de $G_{s_{0}}$. Pour le sommet $s=s_{0}$, la restriction identifie $W^{I_{F}}$ \`a $W_{s_{0}}$  et $(\bar{{\cal O}}_{\rho},{\cal L}_{\rho})$ est simplement l'image de $\rho$ par la repr\'esentation de Springer de ce groupe de type $G_{2}$. Elle est d\'ecrite en \cite{C} p. 427. Pour les deux autres sommets, la restriction de $\rho$ \`a $W_{s}$ est d\'ecrite dans les tables 21 et 22 de \cite{Alvis}. Le groupe $G_{s}$ est produit de groupes de type $A_{n-1}$ pour lesquels la description de la correspondance de Springer est imm\'ediate: les repr\'esentations des groupes de Weyl comme les orbites nilpotentes sont param\'etr\'ees par des partitions et, dans ce param\'etrage, la correspondance de Springer est l'identit\'e. On doit signaler une subtilit\'e dans l'utilisation des tables de \cite{Alvis}. Le groupe $W(A_{2})$ de cette r\'ef\'erence est engendr\'e par deux r\'eflexions associ\'ees \`a des racines longues. Ici, notre groupe $W_{s_{2}}$ est bien un groupe $W(A_{2})$ mais il est engendr\'e par deux r\'eflexions associ\'ees \`a des racines courtes. La raison en est que notre diagramme ${\cal D}_{a}^{nr}$ n'est pas le diagramme de Dynkin compl\'et\'e associ\'e \`a un groupe de type $G_{2}$. Pour calculer les restrictions \`a $W_{s_{2}}$, on doit composer les tables de \cite{Alvis} avec l'automorphisme de $W^{I_{F}}$ qui permute les deux g\'en\'erateurs standards, \'echangeant ainsi les r\'eflexions associ\'ees aux racines courtes et celles associ\'ees aux racines longues. Pour les repr\'esentations irr\'eductibles, cela \'echange $\phi_{1,3}'$ et $\phi_{1,3}''$. 
 Par ailleurs, Alvis n'utilise pas les m\^emes notations que Curtis. La traduction est: $\chi_{1,1}=\phi_{1,0}$, $\chi_{1,2}=\phi_{1,6}$, $\chi_{1,3}=\phi_{1,3}'$, $\chi_{1,4}=\phi_{1,3}''$, $\chi_{2,1}=\phi_{2,1}$, $\chi_{2,2}=\phi_{2,2}$.
 
 Les r\'esultats sont les suivants. Quand les restrictions de $\rho$ ont plusieurs composantes, on note chacune d'elles en les s\'eparant par des $+$. Il s'av\`ere que les multiplicit\'es de chaque composante sont toujours \'egales \`a $1$. On expliquera plus loin pourquoi certains termes sont soulign\'es. La deuxi\`eme ligne du tableau rappelle le type de $G_{s,SC}$.
 
 $$\begin{array}{ccccc}(M,s_{M},\varphi,\rho)&s_{0}&s_{134}&s_{2}&\boldsymbol{\nabla}_{F}(M,s_{M},\varphi,\rho)\\ &G_{2}&\tilde{A}_{1}\times A_{1}&\tilde{A}_{2}&\\&&&&\\
 
 (G,s_{0})&\underline{(G_{2}(a_{1}),sgn)}&&&(3^21^2,\gamma,sgn)\\(G,s_{134})&\underline{((2,2),\xi'\not=1)}&&&(3^21^2,-\gamma,sgn)\\ (T,\phi_{1,0})&\underline{G}_{2}&(2,2)&(3)&(71)\\ (T,\phi'_{1,3})&(G_{2}(a_{1}),r)&(2,1^2)&\underline{(3)}&(53)\\ (T,\phi_{2,1})&\underline{(G_{2}(a_{1}),1)}&(2,1^2)+(1^2,2) &(21)&(3^21^2,\gamma,1)\\ (T,\phi_{2,2})&\tilde{A}_{1}&\underline{(2,2)}+(1^2,1^2)&(21)&(3^21^2,-\gamma,1)\\ (T,\phi''_{1,3})&\underline{A}_{1}&\underline{(1^2,2)}&1^3&(2^21^4)\\ (T,\phi_{1,6})&\underline{1}&\underline{(1^2,1^2)}&\underline{(1^3)}&(1^8)\\ \end{array}$$

   Pour d\'emontrer le \ref{premiertheoremeSI},  soit $(M,s_{M},\varphi,\rho)\in {\cal B}_{Irr}(G)$. Notons $\boldsymbol{{\cal O}}\in \boldsymbol{\mathfrak{g}}_{nil}/conj$ l'orbite figurant dans  $\boldsymbol{\nabla}_{F}(M,s_{M},\varphi,\rho)$. Soit $(s,\bar{{\cal O}}_{\rho'},{\cal L}_{\rho'})$ un triplet intervenant dans la ligne associ\'ee \`a $(M,s_{M},\varphi,\rho)$ du tableau ci-dessus. Posons $\boldsymbol{{\cal O}}_{\rho'}=\iota_{s}(\bar{{\cal O}}_{\rho'})\in \boldsymbol{\mathfrak{g}}_{nil}/conj$. Le th\'eor\`eme \ref{premiertheoremeSI} r\'esulte de l'assertion plus forte: $\boldsymbol{{\cal O}}_{\rho'}\leq \boldsymbol{{\cal O}}$ pour l'ordre usuel des partitions. Mais $\boldsymbol{{\cal O}}_{\rho'}$ est calcul\'e en \ref{iotanilram} et on constate que cette assertion est v\'erifi\'ee. Les termes soulign\'es dans le tableau ci-dessus sont ceux pour lesquels ${\boldsymbol{\cal O}}_{\rho'}=\boldsymbol{{\cal O}}$.

Fixons une orbite $\boldsymbol{{\cal O}}\in {\cal U}_{sp}$ conservée par $\Gamma_{F}$ et un élément $N$ de cette orbite appartenant à notre ensemble de représentants  ${\cal N}(\Gamma_{F})$. Dans  le cas o\`u $A(N)=\{1\}$, l'ensemble ${\cal B}_{Irr}(G,\boldsymbol{{\cal O}})$ est r\'eduit \`a un \'el\'ement $(M,s_{M},\varphi,\rho)$. Le th\'eor\`eme  \ref{deuxiemetheoremeSI}  se r\'esume \`a l'assertion suivante: il existe un triplet
 $(s,\bar{{\cal O}}_{\rho'},{\cal L}_{\rho'})$ intervenant dans la ligne associ\'ee \`a $(M,s_{M},\varphi,\rho)$ tel que $\boldsymbol{{\cal O}}_{\rho'}=\boldsymbol{{\cal O}}$. Autrement dit, il y a dans cette ligne un terme soulign\'e. C'est clair. Il reste le cas o\`u $\boldsymbol{{\cal O}}$ est param\'etr\'ee par $(3^21^2)$. On a alors $\tilde{\bar{A}}(N)/conj\simeq \{\pm \gamma\}$. On d\'efinit l'application $b_{F}$ du th\'eor\`eme \ref{deuxiemetheoremeSI}  par $b_{F}(N,\gamma)=(s_{0},G_{2}(a_{1}))$ et $b_{F}(N,-\gamma)=(s_{134},(2,2))$. L'assertion (i) de ce th\'eor\`eme \'equivaut  aux \'egalit\'es $\iota_{s_{0}}(G_{2}(a_{1}))=\iota_{s_{134}}(2,2)=(3^21^2)$, cf. \ref{iotanilram}. L'assertion (ii) pour $\gamma$ signifie que pour les deux \'el\'ements $(M,s_{M},\varphi,\rho)$ tels que $\boldsymbol{\nabla}_{F}(M,s_{M},\varphi,\rho)=(3^21^2,\gamma,1)$ ou $(3^21^2,\gamma,sgn)$ le couple $(s_{0},G_{2}(a_{1}))$ appara\^{\i}t dans la ligne associ\'ee \`a $(M,s_{M},\varphi,\rho)$ et que le caract\`ere $\xi'$ qui compl\`ete ce couple est diff\'erent pour les deux \'el\'ements $(M,s_{M},\varphi,\rho)$. C'est clair. L'assertion pour $-\gamma$ est similaire. L'assertion (iii) du th\'eor\`eme pour $\gamma$ signifie que, pour les  deux \'el\'ements $(M,s_{M},\varphi,\rho)$ tels que $\boldsymbol{\nabla}_{F}(M,s_{M},\varphi,\rho)=(3^21^2,-\gamma,1)$ ou $(3^21^2,-\gamma,sgn)$ le couple $(s_{0},G_{2}(a_{1}))$ n'appara\^{\i}t  pas dans la ligne associ\'ee \`a $(M,s_{M},\varphi,\rho)$. C'est clair.  L'assertion pour $-\gamma$ est similaire.  Cela d\'emontre le th\'eor\`eme \ref{deuxiemetheoremeSI}.
 
  \subsubsection{Preuve des th\'eor\`emes \ref{premiertheoremeSI} et \ref{deuxiemetheoremeSI} dans le cas $(E_{6},ram)$ \label{preuveE6rambis}}
 On suppose que $G$ est de type $(E_{6},ram)$. On d\'emontre les th\'eor\`emes \ref{premiertheoremeSI} et \ref{deuxiemetheoremeSI} par la m\^eme m\'ethode que dans le cas $(D_{4},3-ram)$. On utilise les m\^emes notations. On \'ecrira ci-dessous le tableau dont les th\'eor\`emes r\'esultent. 
 
 Concernant les repr\'esentations irr\'eductibles d'un groupe de Weyl de type $F_{4}$, le dictionnaire entre les notations de Alvis et de Curtis est le suivant:
 $$\chi_{1,1}=\phi_{1,0},\, \chi_{1,2}=\phi_{1,12}'',\, \chi_{1,3}=\phi_{1,12}',\, \chi_{1,4}=\phi_{1,24},\, \chi_{2,1}=\phi_{2,4}'',\, \chi_{2,2}=\phi_{2,16}',\, \chi_{2,3}=\phi_{2,4}',$$
 $$\chi_{2,4}=\phi_{2,16}'',\, \chi_{4,1}=\phi_{4,8},\, \chi_{9,1}=\phi_{9,2},\, \chi_{9,2}=\phi_{9,6}'',\, \chi_{9,3}=\phi_{9,6}',\, \chi_{9,4}=\phi_{9,10},\, \chi_{6,1}=\phi_{6,6}',$$
  $$ \chi_{6,2}=\phi_{6,6}'',\, \chi_{12,1}=\phi_{12,4},\, \chi_{4,2}=\phi_{4,1},\, \chi_{4,3}=\phi_{4,7}'',\, \chi_{4,4}=\phi_{4,7}',\, \chi_{4,5}=\phi_{4,13},\, \chi_{8,1}=\phi_{8,3}'',$$
 $$ \chi_{8,2}=\phi_{8,9}',\, \chi_{8,3}=\phi_{8,3}',\, \chi_{8,4}=\phi_{8,9}'',\, \chi_{16,1}=\phi_{16,5}.$$
 
 Ici encore, notre diagramme ${\cal D}_{a}^{nr}$ n'est pas le diagramme compl\'et\'e d'un groupe de type $F_{4}$. Pour calculer les restrictions \`a des groupes $W_{s}$ des repr\'esentations irr\'eductibles du groupe $W^{I_{F}}$ de type $F_{4}$,  on doit composer ces repr\'esentations avec l'automorphisme de ce groupe qui \'echange les g\'en\'erateurs correspondant aux racines courtes avec ceux  correspondant aux racines longues. Cet automorphisme agit sur les repr\'esentations irr\'eductibles en conservant celles qui, dans la notation de Curtis, ne sont pas affect\'ees de $'$ ou $''$. Il \'echange celles affect\'ees d'un $'$ avec celles affect\'ees d'un $''$, par exemple il permute $\phi_{1,12}'$ et $\phi_{1,12}''$, sauf pour le couple  form\'e de $\phi_{6,6}'$ et $\phi_{6,6}''$: ces deux repr\'esentations sont fix\'ees par l'automorphisme.

  Consid\'erons un \'el\'ement $(M,s_{M},\varphi,\rho)\in {\cal B}_{Irr}(G)$. Si $M=G$, le calcul que nous devons faire est imm\'ediat comme en \ref{preuveD43ram}. Il y a deux \'el\'ements $(G,s_{35},\varphi_{i})$ pour $i=1,2$ qui n'existent que dans le cas o\`u $\delta_{3}(q-1)=1$. Sous cette condition, les faisceaux-caract\`eres cuspidaux \`a support nilpotent dans $\mathfrak{g}_{s_{35}}$ dont sont issus ces deux triplets sont port\'es par la m\^eme orbite $(3,3)$ et correspondent aux deux caract\`eres d'ordre $3$ du groupe $Z_{G_{s_{35}}}(N')/Z_{G_{s_{35}}}(N')^0$, o\`u $N'$ est un \'el\'ement de cette orbite. On a not\'e ces caract\`eres $\theta$ et $\theta^2$. 
  
   Si $M=T$, on
a besoin de conna\^{\i}tre   la correspondance de Springer  pour les groupes $G_{s}$. Puisque $M=T$, cette correspondance est insensible \`a une isog\'enie, on peut remplacer $G_{s}$ par son groupe adjoint ou le rev\^etement simplement connexe de celui-ci.  Le groupe $G_{s_{0}}$ est de type $F_{4}$ et la correspondance de Springer est calcul\'ee en \cite{C} p. 428. Il y a dans $\mathfrak{g}_{s_{0},nil}$ une orbite de type $F_{4}(a_{3})$ telle que, pour $N'$ dans cette orbite, le groupe  $Z_{G_{s_{0}}}(N')/Z_{G_{s_{0}}}(N')^0$ soit isomorphe \`a $\mathfrak{S}_{4}$.  On adopte les notations de \cite{C} pour  les repr\'esentations irr\'eductibles de ce groupe (qui d\'eterminent les syst\`emes locaux sur l'orbite). 
Pour $s=s_{35}$ et $s=s_{4}$, $G_{s,AD}$ est produit de groupes de type $A_{n-1}$ et le calcul de cette correspondance est imm\'ediat. Si $s=s_{16}$, $G_{s,AD}$ contient un facteur de type $B_{3}$, c'est-\`a-dire un facteur $SO(7)$. Si $s=s_{2}$, $G_{s,AD}$ est de type $C_{4}$ et on peut remplacer $G_{s}$ par $G_{s,SC}=Sp(8)$. On a donc besoin du calcul de la correspondance de Springer pour les groupes $SO(7)$ et $Sp(8)$. On a expliqu\'e ce calcul en \ref{groupesclassiques}. Le r\'esultat est le suivant, avec les notations de cette r\'ef\'erence. Dans la colonne de droite, le premier terme est une partition $\lambda$, orthogonale ou symplectique. Le deuxi\`eme terme est un caract\`ere $\epsilon\in \{\pm 1\}^{Jord_{bp}(\lambda)}$ qui d\'etermine le faisceau port\'e par l'orbite param\'etr\'ee par $\lambda$. Ce caract\`ere est  omis si c'est le caract\`ere unit\'e. Dans les autres cas, les \'el\'ements de $Jord_{bp}(\lambda)$ sont rang\'es en ordre d\'ecroissant. Par exemple, si $\lambda=(42^2)\in {\cal P}^{symp}(8)$, le caract\`ere $\epsilon=(1,-1)$ a pour composantes $\epsilon_{4}=1$, $\epsilon_{2}=-1$. 

$$\begin{array}{cccccc}Sp(8)&&&&SO(7)&\\(\alpha,\beta)&Spr(0,\alpha,\beta)&&&(\alpha,\beta)&Spr(1,\alpha,\beta)\\ (4,\emptyset)&(8)&&&(3,\emptyset)&(7)
\\(31,\emptyset) &(61^2)&&&(21,\emptyset)&(51^2,(1,-1))\\(2^2,\emptyset)&(42^2,(1,-1))&&&(1^3,\emptyset)&(31^4,(1,-1))\\(21^2,\emptyset)&(41^4)&&&(2,1)&(51^2)\\(1^4,\emptyset)&(21^6)&&&(1^2,1)&(32^2)\\(3,1)&(62)&&&(1,2)&(3^21)\\(21,1)&(42^2)&&&(1,1^2)&(31^4)
\\(1^3,1)&(2^31^2)&&&(\emptyset,3)&(3^21,(-1,1))\\(2,2)&(4^2)&&&(\emptyset,21)&(2^21^3)\\(1^2,2)&(3^22)&&&(\emptyset,1^3)&(1^7)\\(2,1^2)&(421^2)&&&&\\(1^2,1^2)&(2^4)&&&&\\(1,3)&(4^2,-1)&&&&\\(1,21)&(3^21^2)&&&&\\(1,1^3)&(2^21^4)&&&&\\(\emptyset,4)&(62,(-1,-1))&&&&\\(\emptyset,31)&(421^2,(-1,-1))&&&&\\(\emptyset,2^2)&(2^4,-1)&&&&\\(\emptyset,21^2)&(2^21^4,-1)&&&&\\(\emptyset,1^4)&(1^8)&&&&\\ \end{array}$$
 
 Dans le cas o\`u $M=M(A_{5})$ avec la notation de \ref{E6ramSI}, c'est-\`a-dire $M_{SC}$ est de type $A_{5}$, l'ensemble $\Lambda$ associ\'e \`a $M$ est $\{s_{0},s_{35},s_{2}\}$ et l'ensemble des $s\in S(\bar{C})$ tels que $\Lambda\subset \Lambda(s)$ est $\{s_{16},s_{4}\}$. Pour ces sommets $s=s_{16}$ et $s_{4}$, on a besoin de conna\^{\i}tre partiellement la correspondance de Springer g\'en\'eralis\'ee pour le groupe $G_{s}$, plus exactement l'image par cette correspondance d'\'el\'ements de ${\cal I}^{FC}(G_{s})$ de la forme $(M_{ s_{M}},\rho')$. Ces \'el\'ements ne se descendent pas au groupe adjoint $G_{s,AD}$. Lusztig et Spaltenstein ont calculé cette  correspondance de Springer généralisée. Toutefois, nous n'aurons pas besoin du calcul complet de $Spr(M_{s_{M}},\rho')=({\bar{\cal O}}_{\rho'},{\cal L}_{\rho'})$, mais seulement de l'orbite $\bar{{\cal O}}_{\rho'}$. Or celle-ci se d\'eduit de la normalisation de Lusztig, cf. \cite{L2} th\'eor\`eme 9.2 et proposition 9.5. 
 Le groupe $W_{s}(M_{s_{M}})$ a deux \'el\'ements, donc $\rho'$ est l'un des deux caract\`eres $1$ ou $sgn$.  L'espace $\mathfrak{m}_{s_{M},SC}$ porte un faisceau-caract\`ere cuspidal dont le support est ici l'orbite nilpotente r\'eguli\`ere de cet espace, notons-la $\bar{{\cal O}}_{s_{M}}$. Alors $Spr(M_{s_{M}},1)$ est port\'e par l'orbite induite de $\bar{{\cal O}}_{s_{M}}$ \`a $\mathfrak{g}_{s}$, c'est-\`a-dire l'orbite r\'eguli\`ere de cet espace, et $Spr(M_{s_{M}},sgn)$ est port\'e par l'orbite engendr\'ee, qui se calcule ais\'ement. Dans le tableau ci-dessous, on ne fait pas figurer le  faisceau ${\cal L}_{\rho'}$, que nous n'avons pas calcul\'e. 
 
  Les restrictions \`a des groupes $W_{s}$ de repr\'esentations $\rho$ de $W^{I_{F}}$ ont beaucoup de composantes irr\'eductibles. Pour conserver \`a notre tableau une taille raisonnable, on ne fait figurer que les couples $({\cal O}_{\rho'},{\cal L}_{\rho'})$ dont l'orbite ${\cal O}_{\rho'}$ est maximale. Comme l'application $\iota_{s,nil}$ est strictement croissante d'apr\`es le tableau de \ref{iotanilram}, on voit en effet qu'il suffit
 de consid\'erer ces couples pour d\'emontrer nos th\'eor\`emes. Signalons qu'il 
  s'av\`ere que toutes les composantes irr\'eductibles interviennent  avec multiplicit\'e $1$. 

On obtient le tableau suivant, que l'on a divis\'e en quatre pour des raisons typographiques. 
 
 $$\begin{array}{ccccccc}(M,s_{M},\varphi,\rho)&s_{0}&s_{16}&s_{35}&s_{4}&s_{2}&\boldsymbol{\nabla}_{F}(M,s_{M},\varphi,\rho)\\&F_{4}&\tilde{A}_{1}\times B_{3}&\tilde{A}_{2}\times A_{2}&\tilde{A}_{3}\times A_{1}&C_{4}&\\&&&&&&\\(G,s_{0})&\underline{(F_{4}(a_{3}),sgn)}&&&&&(D_{4}(a_{1}),\gamma,sgn)\\ (G,s_{35},\varphi_{1})&&&\underline{((3,3),(\theta,\theta))}&&&(D_{4}(a_{1}),g_{3}\gamma,\theta)\\ (G,s_{35},\varphi_{2})&&&\underline{((3,3),(\theta^2,\theta^2))}&&&(D_{4}(a_{1}),g_{3}\gamma,\theta^2)\\ \end{array}$$
 
 $$\begin{array}{ccccccc}(M,s_{M},\varphi,\rho)&s_{0}&s_{16}&s_{35}&s_{4}&s_{2}&\boldsymbol{\nabla}_{F}(M,s_{M},\varphi,\rho)\\&F_{4}&\tilde{A}_{1}\times B_{3}&\tilde{A}_{2}\times A_{2}&\tilde{A}_{3}\times A_{1}&C_{4}&\\&&&&&&\\
 (M(A_{5}),s_{M(A_{5})},1)&&\underline{(2,7)}&&(4,2)&&D_{5}(a_{1})\\(M(A_{5}),s_{M(A_{5})},sgn)&&\underline{(2,32^2)}&&\underline{(2^2,2)}&&A_{2}+A_{1}\\ \end{array}$$
 
 $$\begin{array}{ccccc}(M,s_{M},\varphi,\rho) &s_{0}&s_{16}&s_{35}&\boldsymbol{\nabla}_{F}(M,s_{M},\varphi,\rho)\\&F_{4}&\tilde{A}_{1}\times B_{3}&\tilde{A}_{2}\times A_{2}&\\&&&&\\(T,\phi_{1,0})&\underline{F_{4}}&(2,7)&(3,3)&E_{6}
 
 \\(T,\phi'_{2,4})&(F_{4}(a_{1}),sgn)&(2,(51^2,(1,-1)))&(3,21)&E_{6}(a_{1})
 
 \\ (T,\phi_{4,1})&\underline{(F_{4}(a_{1}),1)}&(2,51^2) &(21,3) &D_{5}
 
 \\&&+(1^2,7)&+(3,21)&
 
 \\ (T,\phi_{9,2})&\underline{(F_{4}(a_{2}),1)}&(2,7)&(3,3)&(E_{6}(a_{3}),\gamma,1)
 
 \\&&&&
 
 \\(T,\phi''_{2,4})&\underline{(F_{4}(a_{2}),sgn)}&(2,7)&(21,3)&(E_{6}(a_{3}),\gamma,sgn)

 \\(T,\phi'_{8,3})&C_{3}&(2,51^2) &(3,3)&(E_{6}(a_{3}),-\gamma,1)

 \\(T,\phi'_{1,12})&(F_{4}(a_{3}),sgn)&(2,(31^4,(1,-1)))&(3,1^3)&(E_{6}(a_{3}),-\gamma,sgn)
 
 \\(T,\phi'_{4,7})&(C_{3}(a_{1}),sgn)&(2,32^2) &(3,21)&A_{4}+A_{1}

 \\(T,\phi''_{8,3})&\underline{B_{3}}&(2,51^2) &(3,3)&D_{4}
 
 \\&&+\underline{(1^2,7)}&&
 
 \\(T,\phi'_{9,6})&(F_{4}(a_{3}),\psi_{31})&(2,(51^2,(1,-1))) &(3,21)&A_{4}
 
 \\&&&&
 
 \\ (T,\phi_{12,4})&\underline{(F_{4}(a_{3}),1)}&(2,3^21) &(21,3) &(D_{4}(a_{1}),\gamma,1)
 
 \\&&+(1^2,51^2)&+(3,21)&

 \\(T,\phi''_{6,6})&\underline{(F_{4}(a_{3}),\psi_{22})}&(2,31^4) &(1^3,3) &(D_{4}(a_{1}),\gamma,r)
 
 \\&&+(1^2,51^2)&+(21,21)+(3,1^3)&
 
 \\ (T,\phi_{16,5})&(C_{3}(a_{1}),1)&\underline{(2,51^2)}&(21,3) &(D_{4}(a_{1}),g_{2}\gamma,1)
 
 \\&&&+(3,21)&
 
 \\ (T,\phi_{4,8})&(B_{2},sgn)&\underline{(2,(51^2,(1,-1)))}&(21,21)&(D_{4}(a_{1}),g_{2}\gamma,sgn)
 
 \\ (T,\phi'_{6,6})&\tilde{A}_{2}+A_{1}&(2,3^21) &\underline{(3,3)}& (D_{4}(a_{1}),g_{3}\gamma,1)

 \\ (T,\phi''_{9,6})&(B_{2},1)&(2,3^21) &(21,3)&A_{3}
 
 \\&&+(1^2,51^2)&&
 
 \\ (T,\phi''_{4,7})&\underline{A_{2}+\tilde{A}_{1}}&\underline{(2,(3^21,(-1,1)))}&(21,3)&A_{2}+2A_{1}
 
 \\ (T,\phi'_{8,9})&\underline{\tilde{A}_{2}}&(2,32^2)&(21,21) &2A_{2}
 
 \\&&&+\underline{(3,1^3)}&
 
 \\ (T,\phi''_{8,9})&\underline{(A_{2},1)}&(2,2^21^3) &\underline{(1^3,3)} &(A_{2},\gamma,1)
 
 \\&&+\underline{(1^2,3^21)}&+(21,21)&
 
 \\ (T,\phi''_{1,12})&\underline{(A_{2},sgn)}&(1^2,(3^21,(-1,1)))&\underline{(1^3,3)}&(A_{2},\gamma,sgn)
 
 \\ (T,\phi_{9,10})&A_{1}+\tilde{A}_{1}&\underline{(2,31^4) }&(21,21)&(A_{2},-\gamma,1)
 
 \\&&+(1^2,32^2)&&
 
 \\ (T,\phi'_{2,16})&(\tilde{A}_{1},sgn)&\underline{(2,(31^4,(1,-1)))}&(21,1^3)&(A_{2},-\gamma,sgn)
 
 \\ (T,\phi_{4,13})&\underline{(\tilde{A}_{1},1)}&\underline{(2,1^7) }&(1^3,21) &2A_{1}
 
 \\&&+\underline{(1^2,31^4)}&+\underline{(21,1^3)}&
 
 \\ (T,\phi''_{2,16})&\underline{A_{1}}&\underline{(1^2,2^21^3)}&\underline{(1^3,21)}&A_{1}
 
 \\ (T,\phi_{1,24})&\underline{1}&\underline{(1^2,1^7)}&\underline{(1^3,1^3)}&1\\ \end{array}$$

 $$\begin{array}{cccc}(M,s_{M},\varphi,\rho) &s_{4}&s_{2}& \boldsymbol{\nabla}_{F}(M,s_{M},\varphi,\rho)\\&\tilde{A}_{3}\times A_{1}&C_{4}&\\&&&\\(T,\phi_{1,0})&(4,2)&(8)&E_{6}
 
 \\(T,\phi'_{2,4})&(4,2)&\underline{(8)}&E_{6}(a_{1})
 
 \\ (T,\phi_{4,1}) &(31,2) +(4,1^2)&(62)&D_{5}

 \\ (T,\phi_{9,2})&(4,2)&(61^2) +(4^2)&(E_{6}(a_{3}),\gamma,1)

 \\(T,\phi''_{2,4})&(2^2,2)&(42^2,(1,-1))&(E_{6}(a_{3}),\gamma,sgn)

 \\(T,\phi'_{8,3})&(4,2)&\underline{(62)}&(E_{6}(a_{3}),-\gamma,1)

 \\(T,\phi'_{1,12})&(4,1^2)&\underline{(62,(-1,-1))}&(E_{6}(a_{3}),-\gamma,sgn)
 
 \\(T,\phi'_{4,7})&\underline{(4,2)}&(4^2,-1)&A_{4}+A_{1}

 \\(T,\phi''_{8,3})&(31,2)&(42^2)&D_{4}

 \\(T,\phi'_{9,6})&(31,2)+\underline{(4,1^2)} &\underline{(4^2)}&A_{4}

 \\ (T,\phi_{12,4})&(31,2)&(421^2)+(3^22) &(D_{4}(a_{1}),\gamma,1)

 \\(T,\phi''_{6,6})&(21^2,2) +(31,1^2)&(421^2)&(D_{4}(a_{1}),\gamma,r)

 \\ (T,\phi_{16,5})&(31,2)&\underline{(42^2)}&(D_{4}(a_{1}),g_{2}\gamma,1)

 \\ (T,\phi_{4,8})&(2^2,2)&\underline{(42^2,(1,-1))}&(D_{4}(a_{1}),g_{2}\gamma,sgn)
 
 \\ (T,\phi'_{6,6})& (31,2)&(3^22)&(D_{4}(a_{1}),g_{3}\gamma,1)

 \\ (T,\phi''_{9,6})&(2^2,2)&\underline{(41^4)} +(2^4)&A_{3}

 \\ (T,\phi''_{4,7})&(21^2,2)&(2^31^2)&A_{2}+2A_{1}
 
 \\ (T,\phi'_{8,9}) &(2^2,2) +\underline{(31,1^2)}&\underline{(3^21^2)}&2A_{2}

 \\ (T,\phi''_{8,9})&(21^2,2)&(2^31^2)&(A_{2},\gamma,1)

 \\ (T,\phi''_{1,12})&(1^4,2)&(21^6)&(A_{2},\gamma,sgn)
 
 \\ (T,\phi_{9,10})&(21^2,2) +\underline{(2^2,1^2)}&\underline{(2^4)}&(A_{2},-\gamma,1)

 \\ (T,\phi'_{2,16})&\underline{(2^2,1^2)}&\underline{(2^4,-1)}&(A_{2},-\gamma,sgn)
 
 \\ (T,\phi_{4,13})&(1^4,2) +\underline{(21^2,1^2)}&\underline{(2^21^4)}&2A_{1}

 \\ (T,\phi''_{2,16})&\underline{(1^4,2)}&\underline{(21^6)}&A_{1}
 
 \\ (T,\phi_{1,24})&\underline{(1^4,1^2)}&\underline{(1^8)}&1\\ \end{array}$$

 Le th\'eor\`eme  \ref{premiertheoremeSI} se d\'emontre comme en \ref{preuveD43ram} en utilisant les tableaux ci-dessus et le calcul de $\iota_{s,nil}$ fait en \ref{iotanilram} pour tous les sommets $s\in S(\bar{C})$. Les termes soulign\'es ont la m\^eme signification qu'en \ref{preuveD43ram}. 
 
  Fixons une orbite $\boldsymbol{{\cal O}}\in {\cal U}_{sp}$  et un élément $N$ de cette orbite appartenant à notre ensemble de représentants  ${\cal N}\Gamma_{F})$.  Dans  le cas o\`u $A(N)=\{1\}$, l'ensemble ${\cal B}_{Irr}(G,\boldsymbol{{\cal O}})$ est r\'eduit \`a un \'el\'ement $(M,s_{M},\varphi,\rho)$. Le th\'eor\`eme \ref{deuxiemetheoremeSI} se r\'eduit \`a l'assertion que, dans la ligne des tableaux ci-dessus correspondant \`a cet \'el\'ement, il y a au moins un terme soulign\'e (on doit \'evidemment r\'eunir les deux derniers tableaux...). Supposons $A(N)\not=\{1\}$. Il y a trois \'el\'ements $N$ qui v\'erifient cette condition. On d\'efinit dans chaque cas l'application $b_{F}$ du th\'eor\`eme \ref{deuxiemetheoremeSI} par les formules suivantes:
  
  supposons $\boldsymbol{{\cal O}}$ de type $E_{6}(a_{3})$; alors $\tilde{\bar{A}}(N)/conj\simeq \{\pm \gamma\}$; on pose 
  
  \noindent $b_{F}(N,\gamma)=(s_{0},F_{4}(a_{2}))$, $b_{F}(N,-\gamma)=(s_{2},(62))$;
  
  supposons $\boldsymbol{{\cal O}}$ de type $A_{2}$; alors $\tilde{\bar{A}}(N)/conj\simeq \{\pm \gamma\}$; on pose $b_{F}(N,\gamma)=(s_{0},A_{2})$,  $b_{F}(N,-\gamma)=(s_{2},(2^4))$;
 
supposons $\boldsymbol{{\cal O}}$ de type $D_{4}(a_{1})$; alors $\tilde{\bar{A}}(N)/conj\simeq \{ \gamma,g_{2}\gamma,g_{3}\gamma\}$; on pose $b_{F}(N,\gamma)=(s_{0},F_{4}(a_{3}))$, $b_{F}(N,g_{2}\gamma)=( s_{16},(2,51^2))$,   $b_{F}(n,g_{3}\gamma)=(s_{35},(3,3))$.

\noindent Comme en \ref{preuveD43ram}, les propri\'et\'es requises de cette application se lisent sur les tableaux ci-dessus.

\subsection{Une nouvelle base de $SI(\mathfrak{g}(F))^*_{nil}$}

\subsubsection{Un calcul plus pr\'ecis de $\hat{k}(\Lambda,\varphi,\rho)(h_{s,{\cal O}_{s}})$\label{uncalculplusprecis}}
 
Pour $s\in S(\bar{C})$ et ${\cal O}_{s}\in \mathfrak{g}_{s,nil}({\mathbb F}_{q})/conj$, on note $\bar{{\cal O}}_{s}$ l'orbite dans $\mathfrak{g}_{s,nil}$ contenant ${\cal O}_{s}$. 

Soient $\boldsymbol{{\cal O}}\in {\cal U}_{sp}^{\Gamma_{F}}$,  ${\bf d}\in \bar{C}(\boldsymbol{{\cal O}})$ 
 et $(\Lambda,\varphi,\rho)\in {\cal B}'_{Irr}(G,\boldsymbol{{\cal O}},{\bf d})$. Posons $b_{F}({\bf d})=(s,\bar{{\cal O}}_{{\bf d}})$.  Le (ii) du th\'eor\`eme \ref{deuxiemetheoremeSI} affirme l'existence et l'unicit\'e d'une certaine repr\'esentation $\rho'\in Irr_{{\mathbb F}_{q}}(W_{s}(M_{\Lambda,s_{\Lambda}}))$. On la note ici $\rho'_{\Lambda,\varphi,\rho}$ et on pose $Spr(M_{s_{\Lambda}},{\cal E}_{\varphi},\rho'_{\Lambda,\varphi,\rho})=(\bar{{\cal O}}_{{\bf d}},{\cal L}_{\Lambda,\varphi,\rho})$. 

\begin{prop}{Soient $s\in S(\bar{C})$, ${\cal O}_{s}\in \mathfrak{g}_{s,nil}({\mathbb F}_{q})/conj$, $\boldsymbol{{\cal O}}\in {\cal U}_{sp}^{\Gamma_{F}}$, ${\bf d}\in \bar{C}(\boldsymbol{{\cal O}})$ et   $(\Lambda,\varphi,\rho)\in {\cal B}'_{Irr}(G,\boldsymbol{{\cal O}})$.  Fixons un \'el\'ement $N_{s}\in {\cal O}_{s}$. 

(i) Si $\hat{k}(\Lambda,\varphi,\rho)(h_{s,{\cal O}_{s}})\not=0$, alors $dim(\iota_{s,nil}(\bar{{\cal O}}_{s}))<dim(\boldsymbol{{\cal O}})$ ou $ \iota_{s,nil}(\bar{{\cal O}}_{s})=\boldsymbol{{\cal O}}$.

(ii) Supposons $(s,\bar{{\cal O}}_{s})=b_{F}({\bf d})$. Si $(\Lambda,\varphi,\rho)\not\in {\cal B}'_{Irr}(G,\boldsymbol{{\cal O}},{\bf d})$, $\hat{k}(\Lambda,\varphi,\rho)(h_{s,{\cal O}_{s}})=0$. Si $(\Lambda,\varphi,\rho)\in {\cal B}'_{Irr}(G,\boldsymbol{{\cal O}},{\bf d})$,  on a $\Lambda\subset \Lambda(s)$ et 
$$\hat{k}(\Lambda,\varphi,\rho)(h_{s,{\cal O}_{s}})=m(\rho_{s}^{\flat},\rho^{'\flat}_{\Lambda,\varphi,\rho}) \bar{{\cal Y}}_{\bar{{\cal O}}_{{\bf d}},{\cal L}_{\Lambda,\varphi,\rho}}(N_{s}).$$}\end{prop}

Preuve. Le terme $\hat{k}(\Lambda,\varphi,\rho)(h_{s,{\cal O}_{s}})$  est calcul\'e par les propositions \ref{evaluationrho} et \ref{calculcrucial}. D'o\`u

 (1) si $\Lambda\not\subset \Lambda(s)$,  $\hat{k}(\Lambda,\varphi,\rho)(h_{s,{\cal O}_{s}})=0$:
 
 (2) si $\Lambda\subset \Lambda(s)$,
$$\hat{k}(\Lambda,\varphi,\rho)(h_{s,{\cal O}_{s}}) =\sum_{\rho'\in  Irr_{{\mathbb F}_{q}}(W_{s}(M_{s_{\Lambda}}))}q^{dim(\bar{{\cal O}}_{\rho'})/2-dim(\bar{{\cal O}}_{s})/2}m(\rho_{s}^{\flat},\rho^{'\flat})\chi^{\natural}_{M_{s_{\Lambda}},{\cal E}_{\varphi},\rho'}(N_{s}).$$

Supposons $\hat{k}(\Lambda,\varphi,\rho)(h_{s,{\cal O}_{s}})\not=0$. Alors $\Lambda\subset \Lambda(s)$ et 
 il existe $\rho'\in Irr_{{\mathbb F}_{q}}(W_{s}(M_{\Lambda,s_{\Lambda}}))$ tel que $m(\rho_{s}^{\flat},\rho^{'\flat})\not=0$ et que $\chi^{\natural}_{M_{s_{\Lambda}},{\cal E}_{\varphi},\rho'}$ ne soit pas nul sur ${\cal O}_{s}$. La premi\`ere propri\'et\'e implique $m(\rho_{s},\rho')\not=0$. Le th\'eor\`eme \ref{premiertheoremeSI} implique que 

(3) $dim(\boldsymbol{{\cal O}}_{\rho'})< dim(\boldsymbol{{\cal O}})$ ou $\boldsymbol{{\cal O}}_{\rho'}=\boldsymbol{{\cal O}}$. 
 
 \noindent Puisque $\chi^{\natural}_{M_{s_{\Lambda}},{\cal E}_{\varphi},\rho'}$ n'est pas nulle sur ${\cal O}_{s}$,    $\bar{{\cal O}}_{s}$ est contenue dans l'adh\'erence de Zariski de  $\bar{{\cal O}}_{\rho'}$. On utilise la propri\'et\'e de croissance de l'application $\iota_{s,nil}$: cela entra\^{\i}ne que $\iota_{s,nil}(\bar{{\cal O}}_{s}) $ est contenue dans l'adh\'erence de $\iota_{s,nil}(\bar{{\cal O}}_{\rho'})=\boldsymbol{{\cal O}}_{\rho'}$.  D'où
 
(4) $dim( \iota_{s,nil}(\bar{{\cal O}}_{s}))< dim(\boldsymbol{{\cal O}}_{\rho'})$ ou $\iota_{s,nil}(\bar{{\cal O}}_{s})=\boldsymbol{{\cal O}}_{\rho'}$.  

\noindent L'assertion (i) de l'\'enonc\'e r\'esulte de (3) et (4).

Supposons $(s,\bar{{\cal O}}_{s})=b_{F}({\bf d})$. On a alors $ \iota_{s,nil}(\bar{{\cal O}}_{s})=\boldsymbol{{\cal O}}$. Sous l'hypoth\`ese $\Lambda\subset \Lambda(s)$, le m\^eme calcul montre que, dans la formule (2), on peut se limiter aux $\rho'$ tels que $\bar{{\cal O}}_{s}$ est contenue dans l'adh\'erence de Zariski de  $\bar{{\cal O}}_{\rho'}$ et $ \iota_{s,nil}(\bar{{\cal O}}_{s})=\boldsymbol{{\cal O}}_{\rho'}$. La croissance stricte de l'application $\iota_{s,nil}$, cf. \ref{orbitesengendrees},  implique que ces conditions \'equivalent \`a $\bar{{\cal O}}_{\rho'}=\bar{{\cal O}}_{s}$. Le th\'eor\`eme \ref{deuxiemetheoremeSI} dit qu'une telle repr\'esentation n'existe pas si $(\Lambda,\varphi,\rho)\not\in {\cal B}'_{Irr}(G,\boldsymbol{{\cal O}},{\bf d})$. Dans ce cas $\hat{k}(\Lambda,\varphi,\rho)(h_{s,{\cal O}_{s}})=0$. Puisque cette nullit\'e est aussi v\'erifi\'ee si $\Lambda\not\subset \Lambda(s)$, on obtient la premi\`ere assertion du (ii) de l'\'enonc\'e. Supposons $(\Lambda,\varphi,\rho)\in {\cal B}'_{Irr}(G,\boldsymbol{{\cal O}},{\bf d})$. Alors le th\'eor\`eme \ref{deuxiemetheoremeSI} dit que $\Lambda\subset \Lambda(s)$ et qu'il n'y a qu'une repr\'esentation $\rho'$ v\'erifiant la condition ci-dessus,
 \`a savoir $\rho'_{\Lambda,\varphi,\rho}$.  Le terme $q^{dim({\bar{\cal O}}_{\rho'})/2-dim(\bar{{\cal O}}_{s})/2}$ dispara\^{\i}t de la formule (2) pour cette repr\'esentation puisque $\bar{{\cal O}}_{\rho'}=\bar{{\cal O}}_{s}$. Puisque $Spr(M_{s_{\Lambda}},{\cal E}_{\varphi},\rho'_{\Lambda,\varphi,\rho})=(\bar{{\cal O}}_{{\bf d}},{\cal L}_{\Lambda,\varphi,\rho})$ et $N_{s}\in \bar{{\cal O}}_{s}=\bar{{\cal O}}_{{\bf d}}$, on a l'\'egalit\'e $\chi^{\natural}_{M_{s_{\Lambda}},{\cal E}_{\varphi},\rho_{\Lambda,\varphi,\rho}'}(N_{s})=\bar{{\cal Y}}_{\bar{{\cal O}}_{{\bf d}},{\cal L}_{\Lambda,\varphi,\rho}}(N_{s})$, cf. \ref{Springer}. 
 Donc (2) se transforme en la formule du (ii) de l'\'enonc\'e. $\square$

 \subsubsection{Construction d'une famille de distributions\label{constructiondelabase}}
 La famille $(k(\Lambda,\varphi,\rho))_{(\Lambda,\varphi,\rho)\in {\cal B}'_{Irr}(G)}$  est une base de $D^G({\cal D}^{st}(\mathfrak{g}(F)))$. Pour $(\Lambda,\varphi,\rho)\in {\cal B}'_{Irr}(G)$, notons $I(\Lambda,\varphi,\rho)$ l'\'el\'ement de $SI(\mathfrak{g}(F))^*_{nil}$ tel que $D^G\circ\delta(I(\Lambda,\varphi,\rho)))=k(\Lambda,\varphi,\rho)$. Alors $(I(\Lambda,\varphi,\rho))_{(\Lambda,\varphi,\rho)\in {\cal B}'_{Irr}(G)}$ est une base de $SI(\mathfrak{g}(F))^*_{nil}$. Pour $(\Lambda,\varphi,\rho)\in {\cal B}'_{Irr}(G)$, la distribution $I(\Lambda,\varphi,\rho)$ appartient a fortiori \`a $I(\mathfrak{g}(F))^*_{nil}$   et on peut \'ecrire  cet \'el\'ement dans la base $(I_{{\mathbb  O}})_{{\mathbb O}\in \mathfrak{g}_{nil}(F)/conj}$ de cet espace:
  $$(1) \qquad I(\Lambda,\varphi,\rho)=\sum_{{\mathbb  O}\in \mathfrak{g}_{nil}(F)/conj}c(\Lambda,\varphi,\rho,{\mathbb  O})I_{{\mathbb O}}.$$
  Pour ${\mathbb O}\in \mathfrak{g}_{nil}(F)/conj$, on note $\bar{{\mathbb O}}$ son image naturelle dans  $\mathfrak{g}_{nil}/conj$.

\begin{lem}{Soient  $ \boldsymbol{{\cal O}}\in {\cal U}_{sp}^{\Gamma_{F}}$, $(\Lambda,\varphi,\rho)\in {\cal B}'_{Irr}(G,\boldsymbol{{\cal O}})$ et ${\mathbb  O}\in \mathfrak{g}_{nil}(F)/conj$. Supposons $c(\Lambda,\varphi,\rho,{\mathbb O})\not=0$. Alors $dim({\mathbb  O})< dim(\boldsymbol{{\cal O}})$ ou ${\bf Rel}( \boldsymbol{{\cal O}})=\bar{{\mathbb O}}$. }\end{lem}

Preuve.  Notons ${\cal V}$ l'ensemble des ${\mathbb O}'\in \mathfrak{g}_{nil}(F)/conj$ telles que $c(\Lambda,\varphi,\rho,{\mathbb O}')\not=0$. On ordonne de fa\c{c}on usuelle les \'el\'ements de $\mathfrak{g}_{nil}(F)/conj$: ${\mathbb O}'\leq {\mathbb O}''$ si et seulement si ${\mathbb O}'\subset Cl_{F}({\mathbb O}'')$, o\`u $Cl_{F}({\mathbb O}'')$ est l'adh\'erence de ${\mathbb O}''$ pour la topologie $p$-adique. 

Supposons  d'abord que  ${\mathbb  O}$ est un \'el\'ement maximal de ${\cal V}$. 
 Fixons $s\in S(\bar{C})$ et ${\cal O}_{s}\in \mathfrak{g}_{s,nil}({\mathbb F}_{q})/conj$ tels que $rel_{s,F}({\cal O}_{s})={\mathbb O}$, cf. \ref{debacker}(9). Calculons $I(\Lambda,\varphi,\rho)(h_{s,{\cal O}_{s}})$. Pour ${\mathbb O}'\in \mathfrak{g}_{nil}(F)/conj$, la proposition \ref{fonctiontest} dit que  si $I_{{\mathbb O}'}(h_{s,{\cal O}_{s}})\not=0$, alors ${\mathbb  O}\leq {\mathbb  O}'$. Puisque ${\mathbb O}$ est suppos\'e maximal dans ${\cal V}$, on a $c(N,\varphi,\rho,{\mathbb  O}')=0$ si ${\mathbb  O}< {\mathbb  O}'$.  
  Il r\'esulte alors de (1) que 
$$ I(\Lambda,\varphi,\rho)(h_{s,{\cal O}_{s}})= c(\Lambda,\varphi,\rho,{\mathbb  O})I_{{\mathbb O}}(h_{s,{\cal O}_{s}}).$$
Ce dernier terme  est non nul d'apr\`es la proposition \ref{fonctiontest}, donc $ I(\Lambda,\varphi,\rho)(h_{s,{\cal O}_{s}})\not=0$. 

Puisque $h_{s,{\cal O}_{s}}\in {\cal H}$, la d\'efinition de l'application $\delta$ implique que
$$I(\Lambda,\varphi,\rho)(h_{s,{\cal O}_{s}})=\hat{k}(\Lambda,\varphi,\rho)(h_{s,{\cal O}_{s}}).$$
Puisque cette expression n'est pas nulle, le (i) de la proposition \ref{uncalculplusprecis} dit que  la conclusion du lemme est v\'erifi\'ee (notons que $dim({\mathbb O})=dim(\iota_{s,nil}(\bar{{\cal O}}_{s})$).
Cela d\'emontre le lemme sous l'hypoth\`ese que ${\mathbb  O}$ est maximale dans ${\cal V}$. Si ${\mathbb O}$ n'est pas maximale, on peut fixer un \'el\'ement maximal ${\mathbb  O}'$ de ${\cal V}$ tel que ${\mathbb O}\subsetneq
Cl_{F}({\mathbb  O}')$. On a alors $dim({\mathbb O})< dim({\mathbb O}')$ et l'assertion du lemme d\'ej\`a prouv\'ee pour ${\mathbb O}'$ entra\^{\i}ne $dim({\mathbb O})< dim(\boldsymbol{{\cal O}})$. $\square$

Pour $ \boldsymbol{{\cal O}}\in {\cal U}_{sp}^{\Gamma_{F}}$ et $(\Lambda,\varphi,\rho)\in {\cal B}'_{Irr}(G,\boldsymbol{{\cal O}})$, on note  $I^{max}(\Lambda,\varphi,\rho)$ la distribution d\'efinie par 
$$I^{max}(\Lambda,\varphi,\rho)=\sum_{{\mathbb O}\in \mathfrak{g}_{nil}(F)/conj, {\bf Rel}(\boldsymbol{{\cal O}})=\bar{{\mathbb O}}}c(\Lambda,\varphi,\rho,{\cal O})I_{{\mathbb  O}}.$$
D'apr\`es le lemme pr\'ec\'edent, on peut aussi bien remplacer la condition ${\bf Rel}(\boldsymbol{{\cal O}})=\bar{{\mathbb O}}$ par $dim({\mathbb O})=dim(\boldsymbol{{\cal O}})$. Donc, d'apr\`es \ref{integralesnilpotentes}(3), $I^{max}(\Lambda,\varphi,\rho)$ appartient \`a $SI(\mathfrak{g}(F))^*_{nil}$. 
D'autre part, $I^{max}(\Lambda,\varphi,\rho)$ est support\'ee par l'orbite  ${\bf Rel}(\boldsymbol{{\cal O}})\in\mathfrak{g}_{nil}/conj$.

\subsubsection{Ind\'ependance lin\'eaire des \'el\'ements de la famille \label{independance}}
\begin{thm}{La famille de distributions $(I^{max}(\Lambda,\varphi,\rho))_{(\Lambda,\varphi,\rho)\in {\cal B}'_{Irr}(G)}$ est une base de l'espace $SI(\mathfrak{g}(F))^*_{nil}$.}\end{thm}

Preuve. Cette famille est index\'ee par l'ensemble ${\cal B}'_{Irr}(G)$ qui indexe aussi la base $(I(\Lambda,\varphi,\rho))_{(\Lambda,\varphi,\rho)\in {\cal B}'_{Irr}(G)}$ de $SI(\mathfrak{g}(F))^*_{nil}$. C'est donc une base si et seulement si cette famille est lin\'eairement ind\'ependante.   On d\'ecompose la famille en sous-familles index\'ees par ${\cal U}_{sp}^{\Gamma_{F}}$. La sous-famille associ\'ee \`a $ \boldsymbol{{\cal O}}\in {\cal U}_{sp}^{\Gamma_{F}}$ est $(I^{max}(\Lambda,\varphi,\rho))_{(\Lambda,\varphi,\rho)\in {\cal B}'_{Irr}(G,\boldsymbol{{\cal O}})}$. Ses \'el\'ements sont des distributions port\'ees par l'orbite ${\bf Rel}(\boldsymbol{{\cal O}})$. Il en r\'esulte que la famille $(I^{max}(\Lambda,\varphi,\rho))_{(\Lambda,\varphi,\rho)\in {\cal B}'_{Irr}(G)}$ est lin\'eairement ind\'ependante si et seulement si chaque sous-famille l'est.

Fixons  $ \boldsymbol{{\cal O}}\in {\cal U}_{sp}^{\Gamma_{F}}$. Consid\'erons une combinaison lin\'eaire
$$D= \sum_{(\Lambda,\varphi,\rho)\in {\cal B}'_{Irr}(G,\boldsymbol{{\cal O}})}x(\Lambda,\varphi,\rho)I^{max}(\Lambda,\varphi,\rho),$$
avec des coefficients $x(\Lambda,\varphi,\rho)\in {\mathbb C}$. Supposons $D=0$. On veut prouver que $x(\Lambda,\varphi,\rho)=0$ pour tout $(\Lambda,\varphi,\rho)\in {\cal B}_{Irr}'(G,\boldsymbol{{\cal O}})$. Puisque ${\cal B}'_{Irr}(G,\boldsymbol{{\cal O}})=\sqcup_{{\bf d}\in \bar{{\cal C}}(\boldsymbol{{\cal O}}}{\cal B}'_{Irr}(G,\boldsymbol{{\cal O}},{\bf d})$, on est ramen\'e \`a prouver

(1) pour tout ${\bf d}\in \bar{{\cal C}}(\boldsymbol{{\cal O}})$, on a $x(\Lambda,\varphi,\rho)=0$ pour tout $(\Lambda,\varphi,\rho)\in {\cal B}_{Irr}'(G,\boldsymbol{{\cal O}},{\bf d})$.

Fixons  ${\bf d}\in \bar{{\cal C}}(\boldsymbol{{\cal O}})$.  Posons $b_{F}({\bf d})=(s,\bar{{\cal O}}_{{\bf d}})$. Rappelons que $\bar{{\cal O}}_{{\bf d}}$ est une orbite dans $\mathfrak{g}_{s,nil}$ qui est conserv\'ee par l'action galoisienne. Soit ${\cal O}_{s}\in \mathfrak{g}_{s,nil}({\mathbb F}_{q})/conj$ une orbite telle que ${\cal O}_{s}\subset \bar{{\cal O}}_{{\bf d}}$. Posons ${\mathbb  O}=rel_{s,F}({\cal O}_{s})$. Puisque $\iota_{s,nil}(\bar{{\cal O}}_{{\bf d}})=N$, on a $ {\bf Rel}(\boldsymbol{{\cal O}})=\bar{{\mathbb O}}$.  Soit $(\Lambda,\varphi,\rho)\in {\cal B}_{Irr}'(G,\boldsymbol{{\cal O}})$. Calculons $I^{max}(\Lambda,\varphi,\rho)(h_{s,{\cal O}_{s}})$.  La diff\'erence $I^{max}(\Lambda,\varphi,\rho)-I(\Lambda,\varphi,\rho)$ est combinaison lin\'eaire d'int\'egrales orbitales $I_{{\mathbb  O}'}$ pour des ${\mathbb  O}'\in \mathfrak{g}_{nil}(F)/conj$ telles  que $dim({\mathbb O}')< dim(\boldsymbol{{\cal O}})=dim({\mathbb O})$. La proposition \ref{fonctiontest} implique que $I^{max}(\Lambda,\varphi,\rho)-I(\Lambda,\varphi,\rho)$ annule $h_{s,{\cal O}_{s}}$. On a aussi $I(\Lambda,\varphi,\rho)(h_{s,{\cal O}_{s}})=\hat{k}(\Lambda,\varphi,\rho)(h_{s,{\cal O}_{s}})$ par d\'efinition de l'application $\delta$.  D'o\`u $I^{max}(\Lambda,\varphi,\rho)(h_{s,{\cal O}_{s}})=\hat{k}(\Lambda,\varphi,\rho)(h_{s,{\cal O}_{s}})$. On applique la proposition \ref{uncalculplusprecis}. L'expression pr\'ec\'edente est nulle si $(\Lambda,\varphi,\rho)\not\in {\cal B}'_{Irr}(G,\boldsymbol{{\cal O}},{\bf d})$. Si $(\Lambda,\varphi,\rho)\in {\cal B}'_{Irr}(G,\boldsymbol{{\cal O}},{\bf d})$, elle vaut $m(\rho_{s}^{\flat},\rho^{'\flat}_{\Lambda,\varphi,\rho}) \bar{{\cal Y}}_{\bar{{\cal O}}_{{\bf d}},{\cal L}_{\Lambda,\varphi,\rho}}(N_{s})$. On en d\'eduit
$$(2) \qquad D(h_{s,{\cal O}_{s}})=\sum_{(\Lambda,\varphi,\rho)\in {\cal B}'_{Irr}(G,\boldsymbol{{\cal O}},{\bf d})}x(\Lambda,\varphi,\rho)m(\rho_{s}^{\flat},\rho^{'\flat}_{\Lambda,\varphi,\rho}) \bar{{\cal Y}}_{\bar{{\cal O}}_{{\bf d}},{\cal L}_{\Lambda,\varphi,\rho}}(N_{s}).$$
Ceci est nul puisque $D=0$ par hypoth\`ese. 
Pour $(\Lambda,\varphi,\rho)\in {\cal B}'_{Irr}(G,\boldsymbol{{\cal O}},{\bf d})$, le th\'eor\`eme \ref{deuxiemetheoremeSI} dit que  $m(\rho_{s},\rho'_{\Lambda,\varphi,\rho})=1$. Il en r\'esulte que $m(\rho_{s}^{\flat},\rho^{'\flat}_{\Lambda,\varphi,\rho})\not=0$ (cette multiplicit\'e est  \'egale \`a une racine de l'unit\'e qui n'est pas forc\'ement $1$). D'apr\`es le m\^eme th\'eor\`eme, l'application $(\Lambda,\varphi,\rho)\mapsto {\cal L}_{\Lambda,\varphi,\rho}$ est injective sur ${\cal B}'_{Irr}(G,\boldsymbol{{\cal O}},{\bf d})$. Donc les fonctions $\bar{{\cal Y}}_{\bar{{\cal O}}_{{\bf d}},{\cal L}_{\Lambda,\varphi,\rho}}$ qui interviennent ci-dessus, qui sont port\'ees par $\bar{{\cal O}}_{{\bf d}}^{\Gamma_{{\mathbb F}_{q}}}$,  sont lin\'eairement ind\'ependantes. En faisant varier ${\cal O}_{s}$ parmi les \'el\'ements de $\mathfrak{g}_{s,nil}({\mathbb F}_{q})/conj$  qui sont contenus dans  $\bar{{\cal O}}_{{\bf d}}$, l'\'el\'ement $N_{s}$ varie parmi toutes les classes de conjugaison par $G_{s}({\mathbb F}_{q})$ dans $\bar{{\cal O}}_{{\bf d}}^{\Gamma_{{\mathbb F}_{q}}}$. La nullit\'e de l'expression (2) pour toutes ces orbites ${\cal O}_{s}$ entra\^{\i}ne $x(\Lambda,\varphi,\rho)=0$ pour tout $(\Lambda,\varphi,\rho)\in {\cal B}_{Irr}'(G,\boldsymbol{{\cal O}},{\bf d})$. Cela d\'emontre (1) et le th\'eor\`eme. $\square$

\subsubsection{Nombres d'\'el\'ements\label{nombresdelements}}
Fixons $\boldsymbol{{\cal O}}\in {\cal U}_{sp}^{\Gamma_{F}}$. Posons  ${\bf Rel}(\boldsymbol{{\cal O}})=\bar{{\mathbb  O}}\in \mathfrak{g}_{nil}/conj$. L'orbite $\bar{{\mathbb  O}}$ est conserv\'ee par l'action galoisienne. Soit $e\in \bar{{\mathbb  O}}\cap \mathfrak{g}(F)$. D'apr\`es les constructions du paragraphe  \ref{cohomologie}, les ensembles $H^1(\Gamma_{F},\bar{A}(e))$ et $\bar{{\cal C}}^{\sharp}_{F}(\boldsymbol{{\cal O}})$ sont en bijection. Notons 
 $h^1(\boldsymbol{{\cal O}})$ leur nombre d'éléments.   Rappelons que l'on a d\'efini la notion d'orbite exceptionnelle en \ref{parametrage3}.

\begin{lem}{Supposons que $\boldsymbol{{\cal O}}$ ne soit pas exceptionnelle ou que $\delta_{4}(q-1)=1$.  Alors le nombre d'\'el\'ements de ${\cal B}'_{Irr}(G,\boldsymbol{{\cal O}})$ est $h^1(\boldsymbol{{\cal O}})$. Supposons que  $\boldsymbol{{\cal O}}$ soit exceptionnelle et que $\delta_{4}(q-1)=0$. Alors ce nombre d'\'el\'ements est $2$ tandis que $h^1(\boldsymbol{{\cal O}})=4$.}\end{lem}

Preuve. Notons $N$ l'élément de notre ensemble de représentants ${\cal N}(\Gamma_{F})$ qui appartient à $\boldsymbol{{\cal O}}$. On a d\'efini l'ensemble $\bar{{\cal C}}^{Irr}(\boldsymbol{{\cal O}})$ en \ref{parametragecalB}. C'est celui des classes de conjugaison par $\bar{A}(N)$ dans l'ensemble des couples $(d,\mu)$ o\`u $d\in \tilde{\bar{A}}(N)$ et $\mu$ est une repr\'esentation irr\'eductible de $Z_{\bar{A}(N)}(d)$.   Fixons un ensemble de repr\'esentants $\tilde{D}$ des classes de conjugaison par $\bar{A}(N)$ dans $\tilde{\bar{A}}(N)$. Alors le nombre d'\'el\'ements de $\bar{{\cal C}}^{Irr}(\boldsymbol{{\cal O}})$ est la somme sur $d\in \tilde{D}$ de $\vert Irr(Z_{\bar{A}(N)}(d))\vert $. Notons $\bar{{\cal C}}'_{F}(\boldsymbol{{\cal O}})$ l'ensemble des classes de conjugaison par $\bar{A}(N)$ dans l'ensemble des couples $(d,v)\in \tilde{\bar{A}}(N)\times \bar{A}(N)$ tels que $vdv^{-1}=d$. Pour tout $d\in \tilde{\bar{A}}(N)$, notons $Z_{\bar{A}(N)}(d)/conj$ l'ensemble des classes de conjugaison dans $Z_{\bar{A}(N)}(d)$. Le nombre d'\'el\'ements de $\bar{{\cal C}}'_{F}(\boldsymbol{{\cal O}})$ est la somme sur $d\in \tilde{D}$ de $\vert Z_{\bar{A}(N)}(d)/conj\vert $. Mais $\vert Z_{\bar{A}(N)}(d)/conj\vert =\vert Irr(Z_{\bar{A}(N)}(d))\vert $. Donc

(1) les ensembles $\bar{{\cal C}}^{Irr}(\boldsymbol{{\cal O}})$ et $\bar{{\cal C}}'_{F}(\boldsymbol{{\cal O}})$ ont m\^eme nombre d'\'el\'ements. 

L'application $\boldsymbol{\nabla}_{F}$ d\'efinie en \ref{parametragecalB} est une injection de ${\cal B}'_{Irr}(G,\boldsymbol{{\cal O}})$ dans $\bar{{\cal C}}^{Irr}(\boldsymbol{{\cal O}})$. Elle n'est pas toujours surjective. La construction de $\boldsymbol{\nabla}_{F}$ s'appuie sur les descriptions de \cite{W7} des espaces $FC^{st}(\mathfrak{m}_{SC}(F))$ pour les $F$-Levi $M$ de $G$. Dans ces descriptions apparaissent des conditions sur la divisiblit\'e de $q-1$. On peut dire que certaines fonctions disparaissent quand $G$ est exceptionnel et que $q-1$ n'est pas divisible par $3$, $4$ ou $5$. Supposons que $G$ ne soit pas exceptionnel ou que $q-1$ soit divisible par $3\times 4\times 5$. On s'aper\c{c}oit alors que $\boldsymbol{\nabla}_{F}$ se restreint en une bijection de ${\cal B}'_{Irr}(G,\boldsymbol{{\cal O}})$ sur $\bar{{\cal C}}^{Irr}(\boldsymbol{{\cal O}})$.

Si $G$ est classique ou de type $(D_{4},3-ram)$, les ensembles  $\bar{{\cal C}}_{F}^{\sharp}(\boldsymbol{{\cal O}})$ et $\bar{{\cal C}}'_{F}(\boldsymbol{{\cal O}})$ sont tous deux égaux à $\tilde{\bar{A}}(N)\times \bar{A}(N)$, cf. \ref{calCF} remarque (4). Supposons $G$ exceptionnel. Alors 
  $\bar{{\cal C}}_{F}^{\sharp}(\boldsymbol{{\cal O}})$ est l'ensemble  des classes de conjugaison par $\bar{A}(N)$ dans l'ensemble des couples $(d,v)\in \tilde{\bar{A}}(N)\times \bar{A}(N)$ qui v\'erifient la condition $vdv^{-1}=d^q$. On a $\tilde{\bar{A}}(N)=\bar{A}(N)$ ou $\tilde{\bar{A}}(N)=\bar{A}(N)\times {\mathbb Z}/2{\mathbb Z}$. Alors $\bar{{\cal C}}_{F}^{\sharp}(\boldsymbol{{\cal O}})$ s'identifie à l'ensemble  des classes de conjugaison par $\bar{A}(N)$ dans l'ensemble des couples $(d,v)\in \bar{A}(N)\times \bar{A}(N)$ qui v\'erifient la condition $vdv^{-1}=d^q$. De m\^eme, $\bar{{\cal C}}'_{F}(\boldsymbol{{\cal O}})$ s'identifie à l'ensemble  des classes de conjugaison par $\bar{A}(N)$ dans l'ensemble des couples $(d,v)\in \bar{A}(N)\times \bar{A}(N)$ qui v\'erifient la condition $vdv^{-1}=d$. Supposons que $q-1$ soit divisible par $3\times 4\times 5$. Puisque $\bar{A}(N)=\mathfrak{S}_{n}$ pour un entier $n\in \{1,...,5\}$, les deux conditions précédentes sont identiques.   Ces consid\'erations d\'emontrent l'\'egalit\'e $\vert {\cal B}'_{Irr}(G,\boldsymbol{{\cal O}})\vert =\vert\bar{{\cal C}}^{\sharp}_{F}(\boldsymbol{{\cal O}})\vert $, c'est-\`a-dire  $\vert {\cal B}'_{Irr}(G,\boldsymbol{{\cal O}})\vert =h^1(\boldsymbol{{\cal O}})$ sous l'hypoth\`ese que $G$ n'est pas exceptionnel ou que $q-1$ soit divisible par $3\times 4\times 5$. 

Supposons maintenant que $G$ soit exceptionnel et levons l'hypothèse de. divisibilité de $q-1$.  En vertu de (1), l'\'egalit\'e affirm\'ee par l'\'enonc\'e du lemme  \'equivaut \`a

(2) $\vert \bar{{\cal C}}^{Irr}(\boldsymbol{{\cal O}})\vert -\vert {\cal B}'_{Irr}(G,\boldsymbol{{\cal O}})\vert =\vert \bar{{\cal C}}'_{F}(\boldsymbol{{\cal O}})\vert-\vert \bar{{\cal C}}^{\sharp}_{F}(\boldsymbol{{\cal O}})\vert$.  

 On fixe un ensemble de repr\'esentants $D$ des classes de conjugaison dans $\bar{A}(N)$. Comme on l'a dit ci-dessus, on peut remplacer les $\tilde{\bar{A}}(N)$ par $\bar{A}(N)$ dans les définitions et on obtient
$$ \vert \bar{{\cal C}}'_{F}(\boldsymbol{{\cal O}})\vert =\sum_{d\in D}\vert Z_{\bar{A}(N)}(d)/conj\vert .$$
Pour $d\in D$, notons $Z_{\bar{A}(N)}(d;q)$ l'ensemble des $v\in \bar{A}(N)$  tels que $vdv^{-1}=d^q$. Il n'est jamais vide (cf. \ref{calCF} remarque (6))  et, si l'on en fixe un \'el\'ement $v$, on a $Z_{\bar{A}(N)}(d;q)=Z_{\bar{A}(N)}(d)v=vZ_{\bar{A}(N)}(d)$. Notons $Z_{\bar{A}(N)}(d;q)/conj$ l'ensemble des classes de conjugaison par $Z_{\bar{A}(N)}(d)$ dans $Z_{\bar{A}(N)}(d;q)$. Alors 
$$ \vert \bar{{\cal C}}^{\sharp}_{F}(\boldsymbol{{\cal O}})\vert =\sum_{d\in D}\vert Z_{\bar{A}(N)}(d;q)/conj\vert .$$
On calcule ais\'ement le membre de droite de (2) \`a l'aide des deux \'egalit\'es ci-dessus. Traitons par exemple le cas o\`u $G$ est de type $E_{8}$ et $\boldsymbol{{\cal O}}$ est de type $E_{8}(a_{7})$. On a alors $\bar{A}(N)\simeq \mathfrak{S}_{5}$. 
L'ensemble $D$ est param\'etr\'e par ${\cal P}(5)$, on note $\lambda(d)$ la partition param\'etrant $d\in D$.    Si $ \lambda(d)$ est \'egale \`a $(1^5)$, $(21^3)$, $(2^21)$, on a $d^2=1$ donc $d^q=d$   et  $Z_{\bar{A}(N)}(d;q)=Z_{\bar{A}(N)}(d)$. Supposons que $\lambda(d)=(31^2)$ ou $(32)$.   Alors $Z_{\bar{A}(N)}(d)\simeq {\mathbb Z}/6{\mathbb Z}$. Pour $v\in Z_{\bar{A}(N)}(d;q)$, la conjugaison par $v$ agit dans ${\mathbb Z}/6{\mathbb Z}$ par multiplication par $q$. On en d\'eduit que $Z_{\bar{A}(N)}(d;q)/conj$ a m\^eme nombre d'\'el\'ements que ${\mathbb Z}/(6{\mathbb Z}+(q-1){\mathbb Z})$. D'o\`u
$$\vert Z_{\bar{A}(N)}(d)/conj\vert -\vert Z_{\bar{A}(N)}(d;q)/conj\vert =4 (1-\delta_{3}(q-1)).$$
Si $\lambda(d)=(41)$ ou $\lambda(d)=(5)$, on trouve de m\^eme
 $$\vert Z_{\bar{A}(N)}(d)/conj\vert -\vert Z_{\bar{A}(N)}(d;q)/conj\vert =2(1-\delta_{4}(q-1)),\,\, \text{ resp. }4(1-\delta_{5}(q-1)).$$
 D'o\`u, dans notre exemple:
 $$(3) \qquad \vert \bar{{\cal C}}'_{F}(\boldsymbol{{\cal O}})\vert-\vert \bar{{\cal C}}^{\sharp}_{F}(\boldsymbol{{\cal O}})\vert=8(1-\delta_{3}(q-1))+2(1-\delta_{4}(q-1))+4(1-\delta_{5}(q-1)).$$
 
 Le membre de gauche de (2) calcule le d\'efaut de surjectivit\'e de l'injection de ${\cal B}'_{Irr}(G,\boldsymbol{{\cal O}})$ dans $\bar{{\cal C}}^{Irr}(\boldsymbol{{\cal O}})$ d\'eduite de $\boldsymbol{\nabla}_{F}$. Comme on l'a dit ci-dessus, il s'agit de voir quels sont les \'el\'ements $(M,s_{M},\varphi,\rho)$ qui disparaissent quand $q-1$ n'est pas divisible par $3\times 4\times 5$. Dans le cas o\`u $G$ est de type $(E_{6},ram)$, cela se voit directement sur la d\'efinition de $\boldsymbol{\nabla}_{F}$, cf. \ref{E6ramSI}: il y a deux termes qui disparaissent si $\delta_{3}(q-1)=0$, qui concernent l'orbite $\boldsymbol{{\cal O}}$ de type $D_{4}(a_{1})$. Le membre de gauche de (2) vaut donc $2(1-\delta_{3}(q-1))$ pour cette orbite et $0$ pour les autres orbites. 
  Supposons maintenant $G$ d\'eploy\'e sur $F^{nr}$. Les \'el\'ements $(M,s_{M},\varphi,\rho)$ qui disparaissent se  d\'eduisent des r\'esultats de  \cite{W7}  
paragraphe 9. Leurs images $\boldsymbol{{\cal O}}\in {\cal U}_{sp}^{\Gamma_{F}}$ se calculent en associant d'abord \`a $(M,s_{M},\varphi,\rho)$ un triplet $({\bf M}, \boldsymbol{{\cal E}},\rho)$ comme en \ref{parametragecalB}, puis en utilisant les tables de \cite{C}. Donnons le r\'esultat dans le m\^eme cas que plus haut: $G$ est de type $E_{8}$ et $\boldsymbol{{\cal O}}$ est de type $E_{8}(a_{7})$. Si $\delta_{3}(q-1)=0$, il y a quatre termes $(M,s_{M},\varphi,\rho)$ qui disparaissent avec $M$ de type $E_{6}$. Pr\'ecisement, le couple $(M,s_{M})$ est unique et son ensemble $\Lambda$ associ\'e est $\Delta_{a}-\{\alpha_{2},\alpha_{4},\alpha_{7}\}$. Ce couple se compl\`ete par deux fonctions $\varphi$ possibles et deux repr\'esentations $\rho$. Les termes suivants disparaissent aussi  (notons que le terme $\rho$ est inexistant si $M=G$):
  
  si $\delta_{3}(q-1)=0$,  deux termes $(G,s_{\alpha_{4}},\varphi)$ et deux termes $(G,s_{\alpha_{7}},\varphi)$;
   
  si  $\delta_{4}(q-1)=0$,  deux termes $(G,s_{\alpha_{6}},\varphi)$;
  
  si $\delta_{5}(q-1)=0$, quatre termes $(G,s_{\alpha_{5}},\varphi)$. 
  
  On constate alors que $\vert \bar{{\cal C}}'_{F}(\boldsymbol{{\cal O}})\vert-\vert \bar{{\cal C}}^{\sharp}_{F}(\boldsymbol{{\cal O}})\vert$ est \'egal au membre de droite de (3), ce qui d\'emontre (2) pour notre exemple. 
  
  Des calculs analogues valent dans tous les cas o\`u $\boldsymbol{{\cal O}}$ n'est pas exceptionnelle.
  
  Dans les cas o\`u $\boldsymbol{{\cal O}}$ est exceptionnelle, on a $\bar{A}(N)={\mathbb Z}/2{\mathbb Z}$ et $\vert \bar{{\cal C}}^{\sharp}_{F}(N)\vert=4$, d'o\`u $h^1(\boldsymbol{{\cal O}})=4$. Mais, si $\delta_{4}(q-1)=0$, il y a deux termes qui disparaissent de ${\cal B}'_{Irr}(G,\boldsymbol{{\cal O}})$: ceux de la forme $(M,s_{M},\varphi,\rho)$ o\`u $M$ est de type $E_{7}$ et $s_{M}=s_{\alpha_{4}}$ (pour la num\'erotation des racines de $M$).  Dans ce cas, ${\cal B}'_{Irr}(G,\boldsymbol{{\cal O}})$ n'a que deux \'el\'ements. $\square$
  
\subsubsection{L'\'enonc\'e final\label{final}}
Pour toute orbite $\bar{{\mathbb O}}\in (\mathfrak{g}_{nil}/conj)^{\Gamma_{F}}$, notons $SI(\mathfrak{g}(F))^*_{\bar{{\mathbb O}}}$ le sous-espace des \'el\'ements de $SI(\mathfrak{g}(F))^*_{nil}$ qui sont combinaisons lin\'eaires d'int\'egrales orbitales $I_{{\mathbb O}}$  avec ${\mathbb O}\subset \bar{{\mathbb O}}$. 

\begin{thm}{(i) L'espace $SI(\mathfrak{g}(F))^*_{nil}$ est la somme directe des $SI(\mathfrak{g}(F))^*_{\bar{{\mathbb  O}}}$ quand $\bar{{\mathbb \cal O}}$ d\'ecrit  $(\mathfrak{g}_{nil}/conj)^{\Gamma_{F}}$.

(ii) Pour $\bar{{\mathbb  O}}\in (\mathfrak{g}_{nil}/conj)^{\Gamma_{F}}$, l'espace $SI(\mathfrak{g}(F))^*_{\bar{{\mathbb  O}}}$ est nul si $\bar{{\mathbb  O}}$ n'est pas sp\'eciale.

(iii) Soit $\bar{{\mathbb O}}\in (\mathfrak{g}_{nil}/conj)^{\Gamma_{F}}$, supposons $\bar{{\mathbb  O}}$ sp\'eciale et fixons $e\in \bar{{\mathbb  O}}^{\Gamma_{F}}$. Si $\bar{{\mathbb O}}$ n'est pas exceptionnelle ou si $\delta_{4}(q-1)=1$, la dimension de $SI(\mathfrak{g}(F))^*_{\bar{{\mathbb O}}}$ est \'egale au nombre d'\'el\'ements de $H^1(\Gamma_{F},\bar{A}(e))$. Si $\bar{{\mathbb  O}}$ est exceptionnelle et $\delta_{4}(q-1)=0$, la dimension de $SI(\mathfrak{g}(F))^*_{\bar{{\mathbb  O}}}$ est $2$ tandis que le nombre d'\'el\'ements de $H^1(\Gamma_{F},\bar{A}(e))$ est $4$. }\end{thm}

C'est une cons\'equence directe du th\'eor\`eme \ref{independance} et du lemme \ref{nombresdelements}.

jean-loup.waldspurger@imj-prg.fr

CNRS- Institut de Math\'ematiques de Jussieu- Paris rive gauche

4 place Jussieu

Bo\^{\i}te courrier 247

75252 Paris Cedex 05

\end{document}